\documentclass[a4paper,12pt]{book}
\usepackage{amsthm,amsfonts,amssymb,euscript}
\usepackage{latexsym, multicol, fancybox}
\usepackage{graphicx}
\usepackage{color}
\usepackage{amsmath, amsthm, amssymb, bm}
\usepackage{epstopdf}
\usepackage{caption}
\usepackage{psfrag}

\numberwithin{equation}{section}

\newtheorem{theorem}{Theorem}[section]
\newtheorem{lemma}[theorem]{Lemma}
\newtheorem{proposition}[theorem]{Proposition}
\newtheorem{corollary}[theorem]{Corollary}
\newtheorem{definition}[theorem]{Definition}
\newtheorem{remark}[theorem]{Remark}

\newtheorem*{thmmain}{Main Theorem}
\newtheorem*{thmM0}{Theorem M0}
\newtheorem*{thmM1}{Theorem M1}
\newtheorem*{thmM2}{Theorem M2}
\newtheorem*{thmM3}{Theorem M3}
\newtheorem*{thmM4}{Theorem M4}
\newtheorem*{thmM5}{Theorem M5}
\newtheorem*{thmM6}{Theorem M6}
\newtheorem*{thmM7}{Theorem M7}
\newtheorem*{thmM8}{Theorem M8}

\newcommand{\Red}{\textcolor{red}}

\newcommand{\bea}{\begin{eqnarray}}
\newcommand{\eea}{\end{eqnarray}}
\def\beaa{\begin{eqnarray*}}
\def\eeaa{\end{eqnarray*}}
\def\ba{\begin{array}}
\def\ea{\end{array}}
\def\be#1{\begin{equation} \label{#1}}
\def \eeq{\end{equation}}

\newcommand{\nn}{\nonumber}

\newcommand{\ds}{\displaystyle}

\def\a{{\alpha}}

\def\b{{\beta}}
\def\be{{\beta}}
\def\ga{\gamma}
\def\Ga{\Gamma}
\def\de{\delta}
\def\De{\Delta}
\def\ep{\epsilon}
\def\ka{\kappa}

\def\la{\lambda}
\def\La{\Lambda}

\def\Si{\Sigma}
\def\om{\omega}
\def\Om{\Omega}

\def\vphi{\varphi}

\def\th{\theta}

\def\ze{\zeta}
\def\ka{\kappa}
\def\nab{\nabla}
\def\vsi{{\varsigma}}
\def\Up{\Upsilon}

\def\pr{{\partial}}
\def\les{\lesssim}
\def\c{\cdot}

\def\AA{{math\cal A}}
\def\BB{{\mathcal B}}
\def\CC{{\mathcal C}}
\def\MM{{\mathcal M}}
\def\NN{{\mathcal N}}

\def\LL{{\mathcal L}}
\def\II{{\mathcal I}}

\def\FF{{\mathcal F}}

\def\HH{{\mathcal H}}

\def\TT{{\mathcal T}}

\def\SS{{\mathcal S}}
\def\UU{{\mathcal U}}
\def\JJ{{\mathcal J}}
\def\KK{{\mathcal K}}
\def\Lie{{\mathcal L}}

\def\DD{{\mathcal D}}

\def\RR{{\mathcal R}}

\def\AA{{\mathcal A}}

\def\HH{{\mathcal H}}

\def\Lie{{\mathcal L}}

\def\kab{\underline{\ka}}

\def\lap{{\De}}

\def\D{{\bf D}}

\def\O{{\bf O}}

\def\R{{\bf R}}

\def\S{{\bf S}}

\def\T{{\bf T}}
\def\Z{{\bf Z}}
\def\g{{\bf g}}

\def\SSS{{\Bbb S}}
\def\RRR{{\Bbb R}}

\def\f12{{\frac 1 2}}

\def\dual{{\,\,^*}}
\def\div{{\mbox div\,}}
\def\curl{{\mbox curl\,}}
\def\lot{\mbox{ l.o.t.}}

\def\Lb{{\,\underline{l}}}
\def\Hb{\,\underline{H}}
\def\Lb{{\,\underline{L}}}

\def\Omb{\underline{\Om}}

\def\Cb{\underline{C}}

\def\Xh{\,^{(h)}X}

\def\trch{{\mbox tr}\, \chi}
\def\chih{{\widehat \chi}}
\def\chib{{\underline \chi}}
\def\chibh{{\underline{\chih}}}

\def\etab{{\underline \eta}}
\def\omb{{\underline{\om}}}
\def\bb{{\underline{\b}}}
\def\aa{\protect\underline{\a}}
\def\xib{{\underline \xi}}

\def\Xib{\underline{\Xi}}

\def\Ab{\underline{A}}
\def\Bb{\underline{B}}
\def\Xb{\underline{X}}

\newcommand{\BBb}{\underline{\mathcal{\BB}}}

\def\Xh{\widehat{X}}
\def\Xbh{\widehat{\Xb}}

\def\ub{\underline{u}}

\def\dds{ \, d  \hspace{-2.4pt}\dual    \mkern-19mu /}
\def\ddd{ \,  d \hspace{-2.4pt}    \mkern-6mu /}

\def\tr{\mbox{tr}}
\def\atr{\,^{(a)}\mbox{tr}}

\def\trchb{{\tr \,\chib}}

\def\Div{\mbox{\textbf{Div}}}

\def\atrch{\atr\chi}
\def\atrchb{\atr\chib}

\def\hot{\widehat{\otimes}}
\def\rhod{\,\dual\hspace{-2pt}\rho}

\def\piX{\, ^{(X)}\pi}

\def\fb{\protect\underline{f}}
\def\err{\mbox{Err}}
\def\ov{\overline}

\newcommand{\hch}{\widehat{\chi}}

\def\f12{\frac 1 2}
\def\lab{\label}

\def\bsplit{\begin{split}}

\newcommand{\Mext}{\,{}^{(ext)}\mathcal{M}}
\def\MMext{\,{}^{(ext)}\mathcal{M}}
\newcommand{\Mint}{\,{}^{(int)}\mathcal{M}}
\def\Mtop{ \, ^{(top)} \MM}

\newcommand{\rint}{\,{}^{(int)}r}
\newcommand{\rext}{\,{}^{(ext)}r}

\newcommand{\Lext}{\,{}^{(ext)}\mathcal{L}_0}
\newcommand{\Lint}{\,{}^{(int)}\mathcal{L}_0}

\def\qf{\frak{q}}
\def\qfb{\underline{\frak{q}}}

\def\Bk{\mathfrak{B}}

\def\Dk{\mathfrak{D}}
\def\Nk{\mathfrak{N}}

\def\Rk{\mathfrak{R}}

\def\Ik{\mathfrak{I}}
\def\Jk{\mathfrak{J}}
\def\hk{\mathfrak{h}}

\def\dkb{ \, \mathfrak{d}     \mkern-9mu /}
\def\dk{\mathfrak{d}}
\def\Dkext{\,{}^{(ext)}\mathfrak{D}}

\DeclareFontFamily{U}{mathx}{\hyphenchar\font45}
\DeclareFontShape{U}{mathx}{m}{n}{
      <5> <6> <7> <8> <9> <10>
      <10.95> <12> <14.4> <17.28> <20.74> <24.88>
      mathx10
      }{}
\DeclareSymbolFont{mathx}{U}{mathx}{m}{n}
\DeclareFontSubstitution{U}{mathx}{m}{n}
\DeclareMathAccent{\widecheck}{0}{mathx}{"71}

\def\Zc{\widecheck{Z}}
\def\Hc{\widecheck{H}}
\def\Hbc{\widecheck{\Hb}}
\def\trXc{\widecheck{\tr X}}
\def\trXbc{\widecheck{\tr\Xb}}
\def\Pc{\widecheck{P}}
\def\ombc{\widecheck{\omb}}
\def\Gac{\widecheck{\Ga}}
\def\Rc{\widecheck R}
\def\Mc{\widecheck M}
\def\rhoc{\widecheck{\rho}}

\def\kac{\widecheck \ka}
\def\kabc{\widecheck{\underline{\ka}}}

\def\omc{\widecheck \omega}
\def\ombc{\underline{\widecheck{\omega}}}

\def\trchc{\widecheck{\tr\chi}}
\def\trchbc{\widecheck{\tr\chib}}
\def\yc{\widecheck{y}}
\def\zc{\widecheck{z}}
\def\zec{\widecheck{\ze}}
\def\etac{\widecheck{\eta}}

\def\muc{\widecheck{\mu}}
\def\rhodc{\widecheck{\rhod}}

\newcommand{\deh}{\delta_{\mathcal{H}}}
\newcommand{\dec}{\delta_{dec}}
\newcommand{\dee}{\delta_{extra}}
\newcommand{\dt}{\delta_B}

\newcommand{\idl}{\LL_0}

\newcommand{\rintl}{\,{}^{(int)}r_{\idl}}
\newcommand{\rextl}{\,{}^{(ext)}r_{\idl}}

\def\epg{\overset{\circ}{\ep}}
\def\dg{\overset{\circ}{\de}}
\def\Jp{J^{(p)} }

\def\piT{{^{(\T)}\pi}}

\def\MMextend{\MM^{(extend)} }

\def\ovS{\overset{\circ}{ S}}
\def\ug{\overset{\circ}{u}}   
\def\sg{\overset{\circ}{s}}   
\def\rg{\overset{\circ}{r}}   
\def\ovPhi{\protect\overset{\circ}{\Phi}\,}
\def\ovphi{\protect\overset{\circ}{\phi}\,}

\newcommand{\Lieb}{\Lie \mkern-10mu /\,}

\newcommand{\mr}{\textrm{Match}}

\def\jS{j^{Sw}}

\newcommand{\Go}{\textrm{Good}}

\def\DDov{\ov{\DD}}

\newcommand{\Gabb}{{\Ga \mkern-11mu /\,}}

\def\GabbX{{^{(X)} \Gabb}}

\def\ovu{\overset{\circ}{ u}}

\def\ovs{\overset{\circ}{ s}}

\def\ovr{\overset{\circ}{ r}}
\def\ovla{\protect\overset{\circ}{ \la}\,}

\def\epg{\protect\overset{\circ}{\ep}}  

\def\ug{\overset{\circ}{u}}   
\def\sg{\overset{\circ}{s}}               
\def\rg{\overset{\circ}{r}}   
 
\def\mg{\overset{\circ}{m}}   

\def\ovu{\overset{\circ}{ u}}

\def\ovs{\overset{\circ}{ s}}

\def\ovr{{\overset{\circ}{ r}\,}}
\def\ovm{{\overset{\circ}{ m}\,}}
\def\ovS{\overset{\circ}{ S}}

\def\ovg{\overset{\circ}{ g}}

\def\ovGa{\overset{\circ}{ \Ga}}
\def\ovPhi{\protect\overset{\circ}{\Phi}\,}
\def\ovphi{\protect\overset{\circ}{\phi}\,}
\def\JpS{{J^{(p, \S)}}}

\def\lapzero{{\overset{\circ}{ \lap}}}
\def\nabzero{{\overset{\circ}{ \nab}}}

\def\Cb{\underline{C}}
\def\Cbp{\Cb^{(p)}} 
\def\Mp{M^{(p)}}

\def\CbpS{\Cb^{(\S, p)}} 
\def\MpS{M^{(\S, p)}}

\def\JpSS{J^{(p, \S)}}

\def\CbpS{\Cb^{(\S, p)}} 
\def\MpS{M^{(\S, p)}}
\def\undB{\underline{B}}
\def\Psih{\widehat{\Psi}}

\def\kadot{\dot{\ka}}
\def\kabdot{\dot{\kab}}
\def\mudot{\dot{\mu}}
\def\Lab{\underline{\La}}

\def\ut{\widetilde{u}}
\def\rt{\widetilde{r}}
\def\et{\widetilde{e}}
\def\Et{\widetilde{E}}

\def\nabt{\widetilde{\nab}}

\def\at{\widetilde{\a}} 
\def\bt{\widetilde{\b}} 
\def\bbt{\widetilde{\bb}} 
\def\aat{\widetilde{\aa}} 
\def\rhot{\widetilde{\rho}}
\def\rhodt{\dual\widetilde{\rho}}
\def\trcht{\widetilde{\trch}}
\def\trchbt{\widetilde{\trchb}}
\def\chiht{\widetilde{\chih}}

\def\chibht{\widetilde{\chibh}}

\def\zet{\widetilde{\ze}}
\def\xit{\widetilde{\xi}}
\def\xibt{\widetilde{\xib}}
\def\omt{\widetilde{\om}}
\def\ombt{\widetilde{\omb}} 
\def\etat{\widetilde{\eta}}
\def\etabt{\widetilde{\etab}}

\def\divt{\widetilde{\div}}
\def\curlt{\widetilde{\curl}}

\def\Jt{\widetilde{J}}

\def\Lextt{\,^{(ext) }\widetilde{\LL_0} }
\def\Ikt{\,^{(ext)}\widetilde{\Ik}}

\def\st{\widetilde{s}}
\def\St{\widetilde{S}}
\def\kat{\widetilde{\ka}}

\def\Fb{\underline{F}}

\def\xt{\widetilde{x}}
\def\Xt{\widetilde{X}}

\def\Rhat{{\widehat R}} 
\def\Nhat{{\widehat{N}}}

\def\Rkstar{\,^\star\Rk}
\def\Rkstarr{\,^{(\Si_*)}\Rk}

\def\Sint{{\,^{(int)}\Si}}

\def\Sk{\frak G}
\def\Skstar{\, ^\star\frak{G}}

\def\Skstarr{\, ^{(\Si_*)}\frak{G}}

\def\Rkext{\,^{(ext)} \Rk}
\def\Skext{\,^{(ext)} \Sk}
\def\Skint{\,^{(int)} \Sk}
\def\Rkint{\,^{(int)} \Rk}
\def\Gacint{\, ^{(int)} \Gac}

\def\Rktop{\,^{(top)} \Rk}
\def\Sktop{\,^{(top)} \Sk}

\def\ehat{{\widehat{e}} }

\def\IkPT{\,^{(PT) }\Ik}

\def\IkPText{\,^{(PT-ext) }\Ik}
\def\IkPTint{\,^{(PT- int) }\Ik}

\def\Kh{\,^{(h)}K}

\begin{document}

\title{Kerr stability  for  small angular momentum}
\author{Sergiu Klainerman, J\'er\'emie Szeftel}

\maketitle

{\bf Abstract.} 
\textit{This is our  main paper in a series   in which we   prove the  full, unconditional,  nonlinear  stability   of  the   Kerr   family  $Kerr(a, m)$ for small angular momentum, i.e. $|a|/m\ll 1$, in the context of
 asymptotically flat solutions of the Einstein vacuum equations (EVE). Three papers in the series, \cite{KS-GCM1} and \cite{KS-GCM2}   and \cite{GKS1} have already been released. We expect that  the  remaining ones
 \cite{GKS2}, \cite{KS:Kerr-B} and \cite{Shen} will appear  shortly. Our   work  extends  the strategy developed in \cite{KS}, in which only axial polarized perturbations of Schwarzschild were treated,  by   developing new  geometric and analytic  ideas on how to deal with   with  general  perturbations  of Kerr.   We note   that the   restriction to  small angular momentum appears  only in connection to Morawetz type estimates in \cite{GKS2} and \cite{KS:Kerr-B}.}

\tableofcontents

%%%%%%%%%%%%%%%%%%%%%%%%%%%%%%%%%%%%%%%%%%%%%

\chapter{Introduction}

%%%%%%%%%%%%%%%%%%%%%%%%%%%%%%%%%%%%%%%%%%%%%

This is our  main paper in a series   in which we   prove the  full  nonlinear  stability   of  the   Kerr   family  $Kerr(a, m)$ for small angular momentum, i.e. $|a|/m\ll 1$, in the context of asymptotically flat solutions of the Einstein vacuum equations (EVE),
 \bea
 \lab{eq:EVE-intro}
 \R_{\a\b}=0.
 \eea
  We recall that the Kerr  family,  discovered  by  R. Kerr \cite{Kerr}  in 1963,  consists of explicit, stationary, asymptotically flat,  solutions of  EVE.  It is considered by physicists  and astrophysicists to be the main mathematical model of a  black hole.

%%%%%%%%%%%%%%%%%%%%%%%%%%%%%%% 
  
\section{Kerr stability conjecture}

%%%%%%%%%%%%%%%%%%%%%%%%%%%%%%%

The discovery of black holes, first  as  explicit solutions of  EVE and later as possible explanations of astrophysical phenomena, 
 has not only revolutionized our understanding of the universe,   it  also gave mathematicians a monumental task: to test the physical reality of these solutions. This may seem nonsensical since physics tests the reality of its objects by experiments and observations and, as such, needs mathematics to formulate the theory and make quantitative predictions, not to test it. The problem, in this case, is that black holes are by definition non-observable and thus no direct experiments are possible.  Astrophysicists ascertain the presence of such objects through indirect observations\footnote{The most recent Nobel  prize in Physics  was  awarded to R.  Penrose  for his theoretical    foundations and to   R. Genzel and A. Ghez    for providing  observational   evidence for the presence of super massive black holes in the center of our galaxy.}  and numerical experiments, but  both are limited in scope to the range of possible observations or 
 the  specific initial conditions  in which numerical  simulations are conducted.  One can rigorously check that the Kerr solutions have vanishing Ricci curvature, that is, their mathematical reality is undeniable.   But to be real in a physical sense, they have to satisfy certain properties that can be neatly formulated in unambiguous mathematical language.
 Chief among them\footnote{Other  such properties concern the  rigidity of the Kerr family  or  the dynamical formations  of black holes from regular configurations. } is the problem of stability, that is, to show that if the precise initial data corresponding to Kerr are perturbed a bit, the basic features of the corresponding solutions do not change much\footnote{If  the Kerr family would be  unstable  under  perturbations, black holes  would be nothing more than mathematical artifacts.}.  This leads naturally to 
 the following  conjecture.

   {\bf Conjecture} (Stability of Kerr conjecture).\,\,{\it  Vacuum, asymptotically flat,  initial data sets, 
   sufficiently close to  $Kerr(a,m)$, $|a|/m <1$,  initial data,  have  maximal developments with complete
   future null infinity  and with   domain of outer communication\footnote{This presupposes the existence of an event horizon. Note that the  existence of such an event horizon   can only be established  upon the completion of the proof of the conjecture.}   which
   approaches  (globally)  a nearby Kerr solution.}

In this section, we provide a brief introduction to the current state of the art concerning the Kerr stability conjecture. For a more in depth introduction to the problem,  we refer the reader to \cite{Daf} or  the introduction    of  \cite{KS}.

%%%%%%%%%%%%%%%%%%%%%%%%%%%%%%%% 
  
\subsection{The Kerr stability problem in the physics literature}

%%%%%%%%%%%%%%%%%%%%%%%%%%%%%%%%

The  nonlinear stability of the Kerr    family  has become, ever since its   discovery  by  R. Kerr \cite{Kerr}  in 1963, a central topic in  general relativity. The first stability  results obtained by  physicists   in the context  of   the linearized EVE  near   a fixed member of the  Kerr family were    mode stability results. The metric perturbation point of view  was initiated by  Regge-Wheeler \cite{Re-W} who discovered the master Regge-Wheeler equation for   odd-parity perturbations. An alternative approach via the Newman-Penrose  (NP) formalism      was  first  undertaken  by Bardeen-Press \cite{Bar-Press}. This latter type of analysis was later extended to the Kerr family by Teukolsky  \cite{Teuk} who made the  important  discovery that  the extreme  curvature components, relative   to  a principal null frame, satisfy  decoupled, separable,  wave equations.   These extreme curvature components    also turn out to be \textit{gauge invariant}  in the sense that small perturbations of the frame lead  to quadratic errors in  their  expression.     The  full  extent   of what could be done  by  mode analysis, in both approaches,    can be found in    Chandrasekhar's   book  \cite{Chand}. Chandrasekhar  also introduced (see \cite{Chand2}) a  transformation theory relating  the two approaches.  More precisely, he found  a transformation which connects  the Teukolsky  equations  to the       Regge-Wheeler  one. The full  mode stability, i.e. lack of exponentially growing modes,  for the Teukolsky equation in Kerr is  due  to  Whiting \cite{Whit} (see also \cite{Yacov} for a stronger quantitive version).

%%%%%%%%%%%%%%%%%%%%%%%%%%%%%%%% 
  
\subsection{The scalar linear wave equation in Kerr}

%%%%%%%%%%%%%%%%%%%%%%%%%%%%%%%%

 Mode stability     is far from establishing  even boundedness of solutions to the linearized equations and  falls thus far short of what is needed to understand nonlinear stability. To achieve that and, in addition, to   derive realistic decay estimates, one needs  an entirely different approach based on  a far reaching  extension of the classical  vectorfield  method\footnote{Method based on the symmetries of Minkowski space   to derive uniform, robust,   decay  for nonlinear wave equations, see  \cite{Kl-vectorfield},  \cite{Kl-null}, \cite{Kl-vect2},  \cite{Ch-Kl0}.} used in the proof of  the nonlinear stability of Minkowski \cite{Ch-Kl}.

 The new   method, which    has emerged  in the last 18 years   in connection to the study of boundedness and decay  for the    scalar wave equation in the Kerr space $\KK(a,m)$,  compensates  for   the lack of  enough Killing and conformal Killing vectorfields   on a Schwarzschild  or  Kerr  background  by   introducing  new vectorfields  whose deformation tensors have    coercive properties    in various, not necessarily causal, regions of spacetime.  
  The starting and most  demanding   part of the new method is  the derivation of  a global,  simultaneous,  \textit{Energy-Morawetz}   estimate  which degenerates   in the trapping region. This task is  somewhat  easier in Schwarzschild,  or    for axially symmetric solutions in Kerr,  where the trapping region  is  restricted to a    smooth hypersurface.  The first  such estimates, in Schwarzschild,    were proved  by Blue and Soffer in \cite{B-S1}, \cite{B-S2} followed by  a long sequence of further  improvements in \cite{B-St}, \cite{DaRo1}, \cite{MaMeTaTo} etc. 
  
     In the absence of axial symmetry,   the derivation of    an    Energy-Morawetz estimate   in    $Kerr(a,m)$, $|a/m|\ll 1$, requires  a more refined  analysis involving either      Fourier decompositions,   see   \cite{DaRo2}, \cite{Ta-Toh},  or  a systematic use of the second order Carter operator, see \cite{A-B}.    The   derivation of    such an estimate  in the full sub-extremal case $|a|<m$  is even more subtle  and   was  achieved  by Dafermos, Rodnianski and Shlapentokh-Rothman \cite{mDiRySR2014}  by  combining mode decomposition with   the vectorfield method.
     
 Once the energy-Morawetz  estimate is derived,  one can combine it  with  local estimates near the horizon, based on  its red shift properties, as introduced in  \cite{DaRo1}, and $r^p$ weighted  estimates, first  introduced  in \cite{Da-Ro3},   to derive  realistic  uniform decay properties of the solutions.

%%%%%%%%%%%%%%%%%%%%%%%%%%%%%%%% 
  
\subsection{Stability of Schwarzschild}

%%%%%%%%%%%%%%%%%%%%%%%%%%%%%%%%

      The first application of the new vectorfield method to   the linearized Einstein equation  near Schwarzschild space, due to Dafermos, Holzegel   and Rodnianski, appeared in  \cite{D-H-R}.  
      The paper is the first to  introduce  and  make use  of a physical space version of  Chandrasekhar's  transformation  to provide realistic boundedness and decay  of  solutions of the Teukolsky equations using the new vectorfield method. This method, of  estimating   the extreme curvature  components    by passing from Teukolsky to a Regge-Wheeler  type equation, to which the vectorfield method can be applied,  is 
       important in all future developments in the subject.  
      
      The first  nonlinear stability result of the Schwarzschild space  appears  in  \cite{KS}. In its simplest version, the  result can be stated as follows.

\begin{theorem}[Klainerman-Szeftel \cite{KS}] 
\lab{MainThm-firstversion-KS}
The future globally hyperbolic   development  of  an   \emph{axially symmetric, polarized},  asymptotically  flat   initial data set, sufficiently close  (in a specified  topology)  to a Schwarzschild  initial data set   of  mass $m_0>0$,  has a complete    future null infinity  $\II^+$ and converges 
in  its causal past  $\JJ^{-1}(\II^{+})$  to another  nearby  Schwarzschild solution of mass $m_{\infty}$ close to $m_0$. 
 \end{theorem}

      The restriction to  axial polarized perturbations  is the  simplest  assumption   which insures  that the final state is itself Schwarzschild  and thus avoids the additional complications  of  the Kerr stability problem which we discuss below. We note that in a  just released  preprint,  the authors in \cite{DHRT}  dispense of any  symmetry assumptions   by properly  preparing a co-dimension 3 subset of the initial data  such that the final state is  still Schwarzschild.

%%%%%%%%%%%%%%%%%%%%%%%%%%%%%%%% 
  
\subsection{The case of Kerr with small angular momentum}

%%%%%%%%%%%%%%%%%%%%%%%%%%%%%%%%      

The first  breakthrough       result on the linear stability of Kerr,      for $|a|/m\ll 1$,  is due independently to \cite{Ma} and \cite{D-H-R-Kerr}. Both results extend the method  of \cite{D-H-R}, mentioned above,  by providing estimates to the extreme linearized curvature components via    a similar Chandrasekhar transformation which   takes the Teukolsky equations   to a generalized Regge-Wheeler (gRW)  equation.   The   passage to a  tensorial  version of  gRW equation, in the fully nonlinear setting,  plays an essential role in our work, see \cite{GKS1}, \cite{GKS2} and the discussion below. The  result  of  \cite{D-H-R-Kerr}  was  recently extended to the full subextremal range in   the  outstanding    paper of   Y.  Shlapentokh-Rothman and  R. Teixeira da Costa  \cite{Y-R}.

The first    linear stability results for the full linearized   Einstein vacuum equations   near   $Kerr(a, m)$,   for $|a|/m\ll 1$,  appear in    \cite{ABBMa2019}  and  \cite{HHV}.  The first paper,  based on the NP formalism,     builds  on the results of  \cite{Ma} and \cite{D-H-R-Kerr}  while the second paper is based on a version of  the   metric formalism.    Though the ultimate  relevance  of these  papers to  nonlinear stability  remains open  they are  both  remarkable results  in so far as they deal with difficulties that looked insurmountable even ten years ago.

  Though  it does not quite   fit in the framework  of our discussion,  we  would like to end this quick  survey 
    of results   by mentioning the  striking achievement  of Hintz and Vasy \cite{HVas} on  the nonlinear   stability   of the stationary part  of  Kerr-de Sitter with small angular momentum, see \cite{HVas}. The result  does not concern EVE  but rather   the Einstein vacuum equation  with  a strictly positive cosmological 
constant
\bea
\lab{eq:Einstein-cosmological}
\R_{\a\b}  +\Lambda \g_{\a\b}=0, \qquad \La>0.
\eea
 It is important to note that,  despite  the   fact that,  formally,  \eqref{eq:EVE-intro} is the limit\footnote{To pass  to the limit  requires one  to  understand  all global in time solutions of \eqref{eq:Einstein-cosmological} with $\La=1$, not only  those which are small perturbations of   Kerr-de Sitter,  treated by \cite{HVas}.} of \eqref{eq:Einstein-cosmological} as $\La\to 0$, the  global  behavior  of  the corresponding solutions   is   radically different\footnote{Major   differences between   formally close equations occur in many other   contexts. For example, the  incompressible Euler equations   are formally the limit of the Navier-Stokes equations as the viscosity tends to zero. Yet, at fixed  viscosity, the global properties of the Navier-Stokes equations are radically different   from that  of the  Euler equations.}.

  The main  simplification in the case of stationary solutions of  \eqref{eq:Einstein-cosmological} is that  the  expected decay   rates  of perturbations  near   Kerr-de Sitter     is exponential,   while in the case $\La=0$ the decay is  lower degree polynomial\footnote{While there is exponential decay in the stationary part treated in \cite{HVas}, note that lower degree polynomial decay is expected in connection to the  stability of the complementary causal region  (called   cosmological   or expanding)  of the full  Kerr-de Sitter  space, see e.g. \cite{Vo}.},  with various components of tensorial quantities  decaying at different  rates,  and the slowest decaying rate\footnote{Responsible  for carrying  gravitational waves at large distances so that they are  detectable.} being no better than $t^{-1}$. Despite  this major simplification,  the work of Hintz and Vasy   is the first  general nonlinear  stability  result  in GR  where    one has to prove  asymptotic stability towards a family of solutions, i.e. full quantitative   convergence  to a final state  close,  but  different from the  initial one\footnote{ The nonlinear stability  of Schwarzschild result in \cite{KS}  is, one the other hand,  the first such result in the  more demanding  case  of asymptotically flat solutions  of EVE.}.  It is  also fair to say that  the work of Hintz-Vasy deals  with  some of the geometric features  of the black hole stability problem   without  having to worry about the  considerable  analytic difficulties of the    physically relevant  Kerr stability problem. On the other hand, as  it  is  apparent in our work here,  the geometric and analytic difficulties   of the Kerr stability problem  are highly entangled and cannot be neatly separated as in the $\La>0$ case.  Thus  the geometric framework   of our work is  very  different from that of \cite{HVas}.

%%%%%%%%%%%%%%%%%%%%%%%%%%%%%%%

\section{Kerr stability  for  small angular momentum}   

%%%%%%%%%%%%%%%%%%%%%%%%%%%%%%%

%%%%%%%%%%%%%%%%%%%%%%%%%%%%%%%
   
\subsection{Simplest version of our main  theorem}

%%%%%%%%%%%%%%%%%%%%%%%%%%%%%%%

The simplest  version of  our main theorem can be stated as  follows.
\begin{theorem}[Main Theorem, first version] 
\lab{MainThm-firstversion}
The future globally hyperbolic   development  of  a general,   asymptotically  flat,    initial data set, sufficiently close  (in a suitable  topology)  to a   $Kerr(a_0, m_0) $   initial data set,  
  for sufficiently  small $a_0/m_0$,  has a complete    future null infinity  $\II^+$ and converges 
in  its causal past  $\JJ^{-1}(\II^{+})$  to another  nearby Kerr spacetime $Kerr(a_\infty, m_\infty )$ with parameters    $(a_{\infty}, m_\infty)$ close to the initial ones $(a_0, m_0)$.
 \end{theorem}
\begin{figure}[h!]
\centering
\includegraphics[scale=0.7]{kerr_5.pdf}
\caption{The Penrose diagram of the  final   space-time in  the Main Theorem with complete future null infinity $\II^+$ and future event horizon $\HH^+$. }
\lab{fig0-introd}
\end{figure}

Our proof  rests on     the following   major    ingredients.

\begin{enumerate}
\item   A  formalism  to derive   tensorial versions of the  Teukolsky  and Regge-Wheeler type  equations  in the   full nonlinear setting.
\item An analytic mechanism to derive  estimates  for solutions of  these.

\item A dynamical mechanism  to  identify the  final  values of   $(a_\infty, m_\infty)$.

\item  A   dynamical mechanism for    finding the right gauge conditions   in which convergence to the final state  takes place.

\item  A precisely formulated continuity argument,    based on a grand bootstrap   scheme,   which assigns to all geometric quantities  involved in the process   specific  decay rates,  which can   be   dynamically recovered  from the initial  conditions   by a long series  of estimates,   and thus ensure  convergence to a   final Kerr  state. 
\item  The continuity  argument is based  on the crucial concept of   finite,  GCM admissible spacetimes $\MM=\Mext\cup \Mint\cup \Mtop$, see Figure \ref{fig1-introd}, whose
 defining characteristic is  its spacelike, GCM  boundary $\Si_*$.   Note that the boundaries  $\Mext\cap\Mtop $ and $\Mint\cap \Mtop $ are timelike\footnote{Asymptotically null   as we pass to the limit.} and that $\Mtop$ is  needed   to have the entire space $\MM$  causal.   The regions $\Mext$ and $\Mint $ are separated by the timelike hypersurface $\TT$ and the spacelike boundary $\AA$ is beyond the future horizon $\HH_+$ of the limiting space.       Finally the region $\LL_0$,  is the initial data  layer   in which $\MM$ is  prescribed  as a solution of the Einstein vacuum equations.
\end{enumerate}

\begin{figure}[h!]
\centering
\includegraphics[scale=1.1]{Kerr_1.pdf}
\caption{The GCM admissible space-time $\mathcal{M}$}
\lab{fig1-introd}
\end{figure} 
\begin{remark}
As in \cite{KS}  we  construct spacetimes  starting    from  the initial  layer $\LL_0$, see Figure \ref{fig1-introd}.   The initial layers we consider are those which arise  from the evolution of  asymptotically flat initial data sets\footnote{  As constructed  in  the works
  \cite{KlNi},  \cite{KlNi2} and  \cite{Ca-Ni}.},  supported on a spacelike hypersurface $\Si_0$.  Thus the   future development of an initial  layer $\LL_0$  should be interpreted as a   future development of the corresponding    initial data set,  see Definition  \ref{definition:developmentofLL_0}. 
\end{remark}

\begin{remark}
\lab{remark-Mtop-intro}
As  mentioned above the region $\Mtop$ is only needed  as causal completion to $\Mext\cup\Mint$ and can be easily determined by a standard 
local existence result once  the  geometry of $\Mext\cup\Mint$  is controlled. For that reason we will mostly ignore it  in this introduction. We also note, as in \cite{KS}, 
that $\Mext$ is by far the harder region to  control, even though $\Mint $ contains the  degenerate   region of trapped null geodesics.
\end{remark}

Here is a short summary of how we deal with these issues.

  \begin{itemize}
 \item In \cite{KS-GCM1} and \cite{KS-GCM2} we  have provided a   framework for  
 dealing  with  the issue  (4), by  constructing generalized notions of  generally covariant modulated (GCM)  spheres\footnote{Generalizing those   used in the  nonlinear stability of  Schwarzschild  in the polarized case, see \cite{KS}.} in the asymptotic region of a general perturbation of Kerr.  The paper     \cite{KS-GCM2}   also  contains   a definition of  the   angular momentum  for GCM spheres. These results are needed here in connection to the construction of the essential boundary $\Si_*$, see also\footnote{The result in \cite{Shen},  where $\Si_*$ is actually  constructed from these GCM pieces,  generalizes the construction of GCMH   from \cite{KS}  to the non-polarized   case  needed here.}  \cite{Shen}.

   \item In  \cite{GKS1}  we   deal with  issue (1) by   developing a geometric formalism  of non-integrable  horizontal structures,  well adapted to perturbations of Kerr, and  use    it  to derive  the   generalized Regge-Wheeler (gRW)  equation   in the context of  general   perturbations of Kerr.  In the linear case,    complex scalar versions  of such equations were first  derived independently  in  \cite{Ma} and  \cite{D-H-R-Kerr}, based on an extension of the  physical space Chandrasekhar type transformation   introduced  in \cite{Chand2} and first  exploited   in \cite{D-H-R},      in the context of the linearized Einstein vacuum  equations near Schwarzschild space.

    \item In  the forthcoming  paper  \cite{GKS2}  we deal  with issue (2)  by deriving estimates for  gRW using an    extension of the classical  vectorfield  method,   based on commutation with second order operators.  In the context of the standard scalar wave equation in Kerr,  such an approach was developed by Andersson and Blue  in their important paper  \cite{A-B}. We  note  that the results  on decay in  \cite{Ma} and  \cite{D-H-R-Kerr}, on the other hand, depend heavily on  mode decompositions for the linearized  gRW equations  in Kerr, an approach whose
    generalization to the full nonlinear setting  seems to present substantial difficulties. Such decompositions were also  essential in     the recent  remarkable  result  \cite{Y-R}  which   derives decay estimates for solutions of the gRW equation  in Kerr(a,m) for the full subextremal case $|a|<m$.
    
    \item The nonlinear terms  present  in the  full version of  the  gRW  equation derived in \cite{GKS1}, as well as those generated 
     by commutation   with vectorfields and  second order Carter operator, are treated  in a similar spirit   as the treatment of the nonlinear terms in  \cite{KS},  by   showing that they verify  a favorable null type structure.
    
    \item In the present paper     we state  a precise version of our main Theorem  \ref{MainThm-firstversion},  define the main objects  and  provide a roadmap for the entire proof.
    We also   deal, in  detail,  with the issues (3) and (5) as follows:
    \begin{itemize}
    \item  We introduce the concept of PG structures (Chapter 2), which allows us to  extend, in perturbations of Kerr, the main features embodied  by   the principal null frames in Kerr.
  
    \item   We define  (Chapter 3) the  notion of  finite,  GCM admissible,  spacetimes  $\MM$,    whose  defining feature, as mentioned above,    is given by   their future, spacelike     boundary    $\Si_* $, see  Figure  \ref{fig1-introd}. 
    This hypersurface  is  foliated by  GCM spheres,  as  defined in  \cite{KS-GCM1}, \cite{KS-GCM2}, and    is used to initialize  the    basic  PG structure and  sphere foliations\footnote{We note that the  null frames  of the PG structure are not   adapted to the sphere foliation,  in the same way that the principal null frame  in   Kerr is not adapted to the  $S(t, r) $ spheres in the    Boyer-Lindquist coordinates. They do however  verify specific compatibility assumptions   described in this  paper in connection to  what we call  principal geodesic structures, see  section \ref{section:PGstructuresKerr}.}  of $\MM$.

      \item We provide   a full set of bootstrap assumptions (Chapter 3)  on these admissible spacetimes. These are of two types: assumptions on decay, involving  derivatives up to order $k_{small} $ for all components of Ricci and  curvature coefficients, relative to the adapted frame,  and assumptions on boundedness, involving derivatives up to $k_{large}= 2k_{small}+1$.

    \item  Relying  on the   estimates for the  extreme components of the curvature, derived in the forthcoming paper \cite{GKS2}, and the GCM conditions on $\Si_*$ we  derive here  complete decay estimates for  all other 
    Ricci and curvature components,  thus improving\footnote{By showing that they     depend only on the smallness of the initial perturbation.}  the bootstrap assumptions on decay.  
   
   \item  To  improve the bootstrap  assumptions on boundedness, we  cannot rely  on the  PG frame, which loses derivatives, but need instead to  use a different frame, which we  call  principal temporal (PT).  In Chapter  9 of this paper,  we  show how  to control the PT frame  at the highest level of derivatives, conditional  on  boundedness estimates   for the curvature.  The estimates for the latter, hyperbolic in nature\footnote{The estimates for the PT frame, assuming  the curvature as given,  are based  on the GCM assumptions on $\Si_*$   and  transport equations. The curvature  estimates, derived in \cite{KS:Kerr-B}, are based instead on Energy-Morawetz and $r^p$-weighted estimates as  well as  by treating the null  Bianchi   equations as a Maxwell type system.},   are delayed to the forthcoming  paper \cite{KS:Kerr-B}. 
   
   \item Finally, we  show how  $\MM$  can be extended to a strictly larger  GCM admissible spacetime $\widetilde{\MM}$ and  thus complete the continuation argument  mentioned in  item (5) above. 
    \end{itemize}
  \end{itemize}

%%%%%%%%%%%%%%%%%%%%%%%%%%%%%%%%%%%%%%  
  
\subsection{Short comparison with Theorem \ref{MainThm-firstversion-KS}}

%%%%%%%%%%%%%%%%%%%%%%%%%%%%%%%%%%%%%%

Our  proof follows the main outline of \cite{KS} in which we   have settled the conjecture   in the restricted class of polarized  perturbations,  see Theorem  \ref{MainThm-firstversion-KS}.

Besides fixing the angular momentum  to be zero,  the   polarization assumption made in \cite{KS}   led to
 important  conceptual and technical simplifications. The most important  challenges to extend    the result in \cite{KS}  to  unconditional perturbations  of Kerr are as follows.

   \begin{itemize}
   \item[a.] The lack of  integrability\footnote{We note that the PT frames, used in Chapter 9, are also non-integrable.}  of the PG   structures of  $\MM$,   which inherits the lack of integrability of the principal null frames in Kerr.
   \item[b.] The  structure and derivation of the gRW equations are considerable more complex.
   \item[c.]  The vectorfield approach  used in \cite{KS} is no longer appropriate  to   Morawetz estimates in perturbations of Kerr.
   \item[d.] The construction  of GCM surfaces    in the  general  setting    is both conceptually and technically  more involved than in the polarized case. For a comprehensive 
    discussion of these  we refer the reader to the introduction of \cite{KS-GCM1} and   \cite{KS-GCM2}.
   \item[e.] The derivation of decay estimates  in the general setting   setting   is both conceptually and technically  more involved  than in the polarized case.
   This is  ultimately due  to  the  lack  of integrability of the PG structures  which  are incompatible with    nonlocal estimates, such as integration of Hodge  type elliptic systems on $S$-foliations, see Remark \ref{remark:Sfolaitions-intro}  below.   To avoid this difficulty in our work  we  need  to construct a secondary  integrable structure and a mechanism  to go back and forth  from the integrable to the non-integrable one.  
   \item[f.] Unlike  in \cite{KS},  where both the decay and  boundedness estimates are based on the same integrable  frame,   we use here two different\footnote{In fact,  we use yet another  frame, namely the integrable frame  associated to $\Si_*$. } types of non-integrable frames: PG frames  for decay and PT  frames for boundedness. 
   \end{itemize}
   We refer the reader  to the introduction of \cite{GKS1}   for a thorough discussion of the items b) and c) and  \cite{KS-GCM1},  \cite{KS-GCM2} for the item d).
   \begin{remark}
   \lab{remark:Sfolaitions-intro}
   In connection to point e) above it is important to remark that   various types of $S$-foliations and their adapted null frames  play a  a fundamental role in many of  the major  mathematical  results in GR, 
    starting with  \cite{Ch-Kl} but also    \cite{Chr-BH},     \cite{KlNi}, \cite{KRS}, \cite{D-H-R}, \cite{Daf-Luk},  \cite{KS} and others. $S$-foliations also play an important role in   applications to fluids such as 
    pioneered by Christodoulou in   \cite{Ch-SW}. Our work here is the first where $S$-foliations  are replaced  by the more complex geometric structures  as mentioned in point e).
   \end{remark}

In what follows we  describe the main conceptual innovations to deal with  a) and e) in this paper. We start by describing the geometric  properties of our  admissible spacetime $\MM$ in Figure  \ref{fig1-introd}.

%%%%%%%%%%%%%%%%%%%%%%%%%%%%%%%   
   
\section{Main geometric structures}
 
%%%%%%%%%%%%%%%%%%%%%%%%%%%%%%%  
 
 As mentioned above, both the results of \cite{Ch-Kl} on the nonlinear stability of the Minkowski space and the result of \cite{KS}  on the nonlinear stability of Schwarzschild under polarized perturbations
  rely on a geometric  formalism   based  on \textit{S-foliations}, i.e. foliations by topological  $2$-spheres,   and   \textit{adapted  null fames}   $(e_3 , e_4, \HH) $,  with $e_3, e_4 $ forming a null pair  and $\HH$,   the horizontal 
  space of vectors orthogonal to both,  tangent  to the $S$-foliation. In both works,   this geometric structure    was   constructed such that it most resembles   the situation
   in the unperturbed case. Thus,  for example,  in  the proof of  stability of the Minkowski case  \cite{Ch-Kl},  all  components of the curvature tensor,  decomposed relative to the frame, converge to zero -- albeit at different rates.
  The same holds true  in \cite{KS}, after   the $\rho$ components of  the curvature  is properly  normalized by subtracting  its  Schwarzschild value.
  
   By contrast, the principal null  vectors  $(e_3, e_4)$ in Kerr, relative to which the  curvature tensor  takes a simple form,   do not lead to integrable horizontal structures, i.e.  the \textit{ horizontal }  space of vectors $\HH$  perpendicular to $(e_3, e_4)$ is not integrable  in the sense of Frobenius.  Thus a  geometric formalism based on $S$-foliations and adapted frames, as developed in   \cite{Ch-Kl}  and used in many other important works   in mathematical GR (see Remark \ref{remark:Sfolaitions-intro}),     is no longer appropriate in perturbations of Kerr.  The   Newman-Penrose (NP),   see \cite{NP},  circumvents this difficulty by working  with principal null pairs $(e_3, e_4)$
       and a specified  basis\footnote{Or rather  the complexified vectors $m=e_1+ ie_2$ and $\ov{m} = e_1-i e_2$.}   $(e_1, e_2)$  for $\HH$.  It thus reduces all calculations to  equations  involving    the  Christoffel  symbols  of the frame.  This  un-geometric feature of the formalism makes it difficult to use it in  the nonlinear setting of  the Kerr  stability  problem.  Indeed  complex calculations, such those needed  to derive the nonlinear analogue of gRW, mentioned above,   depend    on higher  derivatives  of  all   connection coefficients of the NP  frame rather than only  those which are geometrically significant.  This  seriously affects  and complicates  the structure of non-linear corrections and makes  it difficult to avoid  artificial gauge type singularities.

%%%%%%%%%%%%%%%%%%%%%%%%%%%%%%%       
       
\subsection{General horizontal   formalism}       
 
%%%%%%%%%%%%%%%%%%%%%%%%%%%%%%% 
 
  In our work we rely instead on a  tensorial approach, based on horizontal structures which closely mimics the calculations  done  in  integrable  settings while maintaining the important diagonalizable  properties of the  principal  directions.    This allows us to maintain, with minimal  changes, the  geometric formalism  of \cite{Ch-Kl}  widely used today  in mathematical GR.  The formalism,  developed  in detail in \cite{GKS1}, is 
  succinctly reviewed in section \ref{subsection:review-horiz.structures}.  It  was used in \cite{GKS1}  to derive a tensorial, nonlinear version of the  gRW equation of \cite{Ma} and \cite{D-H-R-Kerr}. 
  The idea is very simple: we   define  Ricci coefficients   exactly  as in \cite{Ch-Kl},  relative to an arbitrary   basis of  vectors $(e_1, e_2) $ of $\HH$,  
    \beaa
 \begin{split}
\chib_{ab}&=\g(\D_ae_3, e_b),\qquad\quad\,\,\chi_{ab}=\g(\D_ae_4, e_b),\qquad\quad \ze_a=\frac 1 2 \g(\D_{a}e_4,  e_3),\\
\etab_a&=\frac 1 2 \g(\D_4 e_3, e_a),\qquad \quad \eta_a=\frac 1 2 \g(\D_3 e_4, e_a),\\
\xib_a&=\frac 1 2 \g(\D_3 e_3 , e_a),\qquad \quad\xi_a=\frac 1 2 \g(\D_4 e_4, e_a),\\
\omb&=\frac 1 4 \g(\D_3e_3 , e_4),\qquad\quad\,\, \om=\frac 1 4 \g(\D_4 e_4, e_3),
 \end{split}
\eeaa
and remark that, due to the  lack of integrability of  $\HH$,  the null fundamental  forms  $\chi$ and $\chib$ are no longer symmetric.  They can be both decomposed as follows
\beaa
\chi_{ab}=\frac 1 2 \trch\de_{ab}+\frac1 2 \in_{ab}\atrch +\chih_{ab}, \qquad  \chib_{ab}=\frac 1 2 \trchb\de_{ab}+\frac1 2 \in_{ab} \atrchb +\chibh_{ab},
\eeaa
 where  the new scalars $\atrch$, $\atrchb$ measure the lack of integrability  of  the horizontal structure. The null curvature components are also defined 
 as in \cite{Ch-Kl}, 
\beaa
\a_{ab}=\R_{a4b4},\,\,\,\, \b_a=\frac 12 \R_{a434}, \,\,\,\,    \bb_a=\frac 1 2 \R_{a334},  \,\,\,\,  \aa_{ab}=\R_{a3b3}, \,\,\,\,  \rho=\frac 1 4 \R_{3434}, \,\,\,\,  \rhod=\frac 1 4\dual \R_{3434}.
\eeaa
   The  null structure and null Bianchi equations can then be derived  as  in the integrable case. The only new features are the presence of the scalars $\atrch, \atrchb $ in the  equations. Finally we note that the equations acquire additional simplicity if we pass to complex notations\footnote{The dual here is taken with respect to the fully antisymmetric  horizontal $1$ tensor $\in_{ab}$.}, 
   \beaa
   \bsplit
A &:= \a+i\dual\a, \quad B:=\b+i\dual\b, \quad\, P:=\rho+i\dual\rho,\quad \Bb:=\bb+i\dual\bb, \quad \Ab:=\aa+i\dual\aa,
\\
 X&:= \chi+i\dual\chi, \quad \Xb:=\chib+i\dual\chib, \quad H:=\eta+i\dual \eta, \quad \Hb:=\etab+i\dual \etab, \quad 
 Z :=\ze+i\dual\ze, \\
\Xi&:= \xi + i \dual \xi, \quad\,\,\, \Xib:= \xib+ i \dual \xib.
\end{split}
\eeaa    
  Note that, in particular, $\tr X=\trch-i\atrch,  \, \tr \Xb=\trchb-i\atrchb$.

%%%%%%%%%%%%%%%%%%%%%%%%%%%%%%%  
  
\subsection{Principal geodesic structures} 
\lab{section:PGS-introd}

%%%%%%%%%%%%%%%%%%%%%%%%%%%%%%%

 The geometric formalism  based on   these  non-integrable frames, though perfectly adapted to calculations,  is insufficient to derive   estimates,  which   often  involves  the integration  of  Hodge type elliptic systems on $S$-foliations. It is for this reason that we develop here  a more complex formalism  which  combines $S$-foliations with  non-integrable   frames.   This  approach   requires  in fact  two pairs of frames, the non-integrable   one  which most resemble the principal frame  of Kerr, and  a secondary  one  which is adapted to the $S$-foliation.   To estimate various quantities  we need to constantly   pass from one frame to the other. This is done according to the  general change of frames formula
  \bea
 \lab{eq:framechange-integr-PG-intro}
\bsplit
 \la^{-1}   e_4'&= e_4 + f^b  e_b +\frac 1 4 |f|^2  e_3,\\
  e_a'&= \left(\de_a^b +\frac{1}{2}\fb_af^b\right) e_b +\frac 1 2  \fb_a  e_4 +\left(\frac 1 2 f_a +\frac{1}{8}|f|^2\fb_a\right)   e_3,\\
\la  e_3'&=  \left(1+\frac{1}{2}f\c\fb  +\frac{1}{16} |f|^2  |\fb|^2\right) e_3 + \left(\fb^b+\frac 1 4 |\fb|^2f^b\right) e_b  + \frac 1 4 |\fb|^2 e_4,
\end{split}
\eea
 where  $f, \fb$ are arbitrary $1$ forms and  $\la$ is an arbitrary real scalar, see Lemma \ref{Lemma:Generalframetransf}.

  The  transformation formulas \eqref{eq:framechange-integr-PG-intro}  provide the most general  
 way of passing between two different null frames. They play an essential role all through our work, most prominently in the construction of GCM surfaces in \cite{KS-GCM1}, \cite{KS-GCM2}.

 At the heart of this   dual geometric  formalism  lies the following crucial  definition, see Definition \ref{def:PGstructure}.

\begin{definition}[PG structure]
An outgoing  principal geodesic  (PG)  structure consists  of a null pair $(e_3, e_4)$  and  the induced horizontal structure $\HH$, together with a scalar function $r$ such that
 \begin{enumerate}
\item $e_4$ is  a    null outgoing geodesic  vectorfield, i.e.  $\D_4 e_4=0$,  

\item $r$ is  an affine parameter,  i.e. $e_4(r)=1$,

\item the   gradient  of $r$,  given by  $N=\g^{\a\b} \pr_\b r \pr_\a$,  is perpendicular  to  $\HH$.
\end{enumerate}
A similar concept of incoming  PG structure is defined by interchanging the roles  of $e_3, e_4$
\end{definition}
 
 %%%%%%%%%%%%%%%%%%%%%%%%%%%%%%%
 
 \subsection{Initialization of PG structures} 
 
 %%%%%%%%%%%%%%%%%%%%%%%%%%%%%%%
 
 Such structures are     initialized in our work on  the  boundary  $\Si_*$, see Figure  \ref{fig1-introd}. This leads to the following  definitions, see details in section \ref{subsection:InitializationSi*}.

\begin{definition}\lab{def:framedhypersurface-intro}
 A framed hypersurface consists of   a set  $\big(\Si, r, (\HH,  e_3, e_4)\big)$ where
 \begin{enumerate}
 \item $\Si$ is  smooth a hypersurface in $\MM$,
 
 \item  $( e_3, e_4)$ is a null  pair on $\Si$ such that $e_4$ is transversal to $\Si$, and $\HH$, the  horizontal space   perpendicular on $e_3, e_4$,    is  tangent to $\Si$,
 
 \item  the function $r:\Si\to \RRR$ is a regular  function on $\Si$   such that
$   \HH(r) =0$. 
\end{enumerate}
 \end{definition}
 
 To define an appropriate  initial data set we need also to prescribe an additional horizontal $1$-form $f$ as follows.
\begin{definition}[PG-data set]
\lab{definition:PGdataset-intro}
The boundary data of a PG structure   (PG-data set)    consists of  
\begin{enumerate}
\item  a framed hypersurface   $\big(\Si, r, (\HH, e_3, e_4)\big)$   as in Defintion \ref{def:framedhypersurface-intro},
 
 \item  a fixed  1-form $f$ on the spheres $S$ of the $r$-foliation of $\Si$ verifying  the condition
 \beaa
 b_\Si |f|^2<4\quad \text{on}\quad \Si,
 \eeaa
 where $b_\Si$ is such that $\nu= e_3 + b_\Si e_4$ is tangent to  $\Si$.
 \end{enumerate}
\end{definition}

The following is precisely Proposition \ref{Prop:Initialization of PG structures} in the main text.
 \begin{proposition}
 \lab{Prop:Initialization of PG structures:intro}
 Given a     PG data  set $\big(\Si,  r,  (\HH, e_3, e_4),  f\big) $ as in Definition \ref{definition:PGdataset-intro},  there exists a unique  PG   structure  $\big(r',( \HH', e'_3, e'_4) \big)$   defined in a neighborhood of  $\Si$   such that the following hold true
\begin{enumerate}
\item   The function $r'$ is prescribed on $\Si$ by  $r'=r$.
  
\item    Along $\Si$,  the  restriction  of  the spacetime PG null frame  $(\HH', e'_3, e' _4)$  and  the   given null frame
 $( \HH,  e_3, e_4 )$ on $\Si$ are related by the   transformation formulas \eqref{eq:framechange-integr-PG-intro} with transition coefficients $(f, \fb, \la)$,  where $(e_1, e_2)$ is a fixed, arbitrary,  orthonormal  basis of  $\HH$,  where $f$ is part of the PG-data set, and where $\fb$ and $\la$ are given by 
\beaa
\la=1, \qquad \fb  = -\frac{(\nu(r)-b_\Si)}{1-\frac{1}{4}b_\Si |f|^2}f.  
\eeaa
 \end{enumerate} 
 \end{proposition}

 %%%%%%%%%%%%%%%%%%%%%%%%%%%%%%%
 
 \section{GCM initial data sets}
 
 %%%%%%%%%%%%%%%%%%%%%%%%%%%%%%%
 
 The hypersurface $\Si_*$  in Figure \ref{fig1-introd}  is not only a framed hypersurface.  It also verifies  crucial   general covariant modulated (GCM) conditions. Given the importance  of these conditions 
  we  describe below the main ingredients needed in their definitions.  We concentrate  first on the boundary $S_*$  of $\Si_*$, see Figure \ref{fig1-introd},    on which  various quantities   are initialized and transported  along  $\Si_*$.

%%%%%%%%%%%%%%%%%%%%%%%%%%%%%%%  
   
  \subsection{Last sphere   $S_*$   of $\Si_*$}
 
%%%%%%%%%%%%%%%%%%%%%%%%%%%%%%% 

  To  define the geometry of $S_*$   we need the effective  uniformization  results derived in \cite{KS-GCM2}, which  we review  in section \ref{subsection:effective-unifrmization}.  Based on these   results, we endow  $S_*$ with coordinates  $(\th, \vphi)$ such that the following  conditions are verified.
  
   \begin{enumerate}
 \item[i.] The    induced metric $g$ on $S_*$  takes the form
 \bea\lab{eq:formofthemetriconSstarusinguniformization:0-intro}
 g= e^{2\phi} r^2\Big( (d\th)^2+ \sin^2 \th (d\vphi)^2\Big).
 \eea
 
 \item[ii.] The  functions 
 \bea
 J^{(0)} :=\cos\th, \qquad J^{(-)} :=\sin\th\sin\vphi, \qquad  J^{(+)} :=\sin\th\cos\vphi,
 \eea
 verify  the balanced  conditions 
 \bea
 \int_{S_*}  J^{(p)} =0, \quad p=0,+,-.
 \eea
 \end{enumerate}

   Recall that $\Si_*$ is  assumed to be a framed hypersurface in the  sense of  Definition \ref{def:framedhypersurface-intro} and thus endowed with a frame $(e_3, e_4, \HH)$ 
   and function $r$ on it such that $\HH(r)=0$.

\begin{definition}
\lab{define:am-onSi:intro}
We define the parameters $(m, a) $  of $S_*$  by the formulas
\bea
\lab{eq:Hawkingmass-intro}
\frac{2m}{r}=1+\frac{1}{16\pi}\int_{S_{*}}\trch \trchb,
\eea
and 
\bea
\lab{eq:angular momentum:intro}
a &:=& \frac{r^3}{8\pi m}\int_{S_*} J^{(0)}\curl\b.
\eea
\eqref{eq:Hawkingmass-intro}  is the usual Hawking  mass  of $S_*$ while  \eqref{eq:angular momentum:intro}  was   introduced in \cite{KS-GCM2}.
\end{definition}

%%%%%%%%%%%%%%%%%%%%%%%%%%%%%%%
 
 \subsection{GCM conditions  for  $\Si_*$}

%%%%%%%%%%%%%%%%%%%%%%%%%%%%%%% 

  The coordinates  $(\th, \vphi)$ on $S_*$   and  the $\ell=1$ basis $\Jp$    are extended to  $\Si_*$ by setting
  \bea\lab{eq:canonical-ell=1modesonSi:0-intro}
  \nu(\th)=\nu(\vphi)=0, \qquad    \nu(\Jp)=0, \quad p=0,+,-,
  \eea
  where $\nu=e_3 + b_{*} e_4$ is tangent to $\Si_*$ and  normal   to  the  $r$-foliation on $\Si_*$.  
 We also extend the parameters $(a, m)$ to be constant along $\Si_*$.

  We are now ready to define  the crucial concept of a  GCM hypersurface
   \begin{definition}[GCM hypersurface]
 \lab{def:CanonicalGCM-hypersurface-intro}
 Consider a   framed hypersurface $\Si_*$   with    end sphere $S_*$,   coordinates $(\th, \vphi)$,   
   and functions $J^{(0)}$, $J^{(+)}$ and $J^{(-)}$ defined as in \eqref{eq:formofthemetriconSstarusinguniformization:0-intro}--\eqref{eq:canonical-ell=1modesonSi:0-intro}.  $\Si_*$ is called  a GCM  hypersurface if in addition the following  conditions\footnote{The scalar $\mu:=-\div \ze -\rho +\frac 1 4 \chih\c \chibh$   is the familiar  mass aspect function, as in \cite{Ch-Kl} and \cite{KS}.}  are verified.
 \begin{enumerate}
 \item On any sphere $S$  of the $r$-foliation of $\Si_*$, the following holds
 \bea
\begin{split}
\lab{eq:Si^*-GCM1-intro}
\trch &=\frac{2}{r},\quad\\
 \trchb &=-\frac{2(1-\frac{2m}{r} )}{r}+\underline{C}_0+\sum_{p=0, +,-}\underline{C}_pJ^{(p)},\quad\\
  \mu &=\frac{2m}{r^3}+M_0+\sum_{p=0, +,-}M_p J^{(p)},\\
\int_S J^{(p)}\div\eta &=0, \qquad \int_S J^{(p)}\div\xib=0,\qquad p=0,+,-, \\\
b_{*}\big|_{SP}&=-1-\frac{2m}{r}, 
\end{split}
\eea
where  $\underline{C}_0$, $\underline{C}_p$, $M_0$, $M_p$ are scalar functions on $\Si_*$ constant  along the leaves of the foliation, and  $SP$ denotes the south poles of the spheres on $\Si_*$, i.e. $\th=\pi$.  

\item  In addition, we have on  the  last sphere $S_*$ of $\Si_*$
\bea
\lab{eq:Si^*-GCM2-intro}
\trchb=-\frac{2(1-\frac{2m}{r})}{r},\qquad \int_{S_*} \Jp \, \div\b =0,\qquad  p=0,+,-,
\eea
  as well as 
\bea\lab{eq:Si^*-GCM3-intro}
\int_{S_*}J^{(+)}\curl\b=0, \qquad \int_{S_*}J^{(-)}\curl\b=0.
\eea
 \end{enumerate}
 \end{definition}

\begin{remark}
Given the  five degrees of freedom of the transition parameters  $(f, \fb, \la)$   in  the general  change of frame formula \eqref{eq:framechange-integr-PG-intro} we expect to  be able  to impose   five  GCM conditions on  a sphere $S\subset\Si_*$.   Since  the   frame of $\Si_*$  is tangent  to  its $S$-foliation we implicitly have $\atrch=\atrchb=0$. It would be natural to impose  Schwarzschildian values  for $\trch, \trchb$ and $\mu$, to account for the remaining three degrees of freedom. This would lead however to  a   differential  system in $(f, \fb, \la)$  which  is not solvable, due to the presence of a kernel and a co-kernel at the level of $\ell=1$ modes.  We are thus obliged to     relax  these conditions  by  imposing, in the case of $\trchb$ and $\mu$, Schwarzschildian values only for  the
$\ell\ge 2$ modes, see \eqref{eq:Si^*-GCM1-intro}.   The remaining degrees of  freedom allow us  to  prescribe  also the $\ell=1$ modes of $\div\xib$ and $\div\eta$, as in  \eqref{eq:Si^*-GCM1-intro}.
These conditions on the $\ell=1$ modes correspond in fact at the level of $(f, \fb, \la)$ to ODEs for the $\ell=1$ modes of $\div f$ and $\div\fb$ along\footnote{Our first GCM   result, in \cite{KS-GCM1}, 
 is based  in fact  on prescribing the $\ell=1$ modes of $\div f$ and $\div \fb$.}   $\Si_*$. As a consequence, we can freely prescribe these $\ell=1$ modes on $S_*$, which allows us to obtain \eqref{eq:Si^*-GCM2-intro} on $S_*$. Using the additional freedom   of rigid rotations  for  frames on $S_*$  we  can also  insure  that \eqref{eq:Si^*-GCM3-intro} holds.   The remaining condition  on $b_*$   is related to the  freedom to  choose  the hypersurface $\Si_*$.   
\end{remark}

%%%%%%%%%%%%%%%%%%%%%%%%%%%%%%%%%%%%%%%%%%%%%
    
\section{GCM admissible spacetimes} 

%%%%%%%%%%%%%%%%%%%%%%%%%%%%%%%%%%%%%%%%%%%%%

We are now ready to define  our  GCM admissible spacetime,  concept  of fundamental importance in  our proof. As can be seen in Figure  \ref{fig1-introd},  $\MM=\Mext\cup \Mint\cup\Mtop$.
Each of the   domains $\Mext, \Mint $ and $\Mtop$ are endowed with  a PG structure, all ultimately induced by  $\Si_*$.  The crucial  structure is that of $\Mext$.  Once it is fixed, those of $\Mint$ and $\Mtop$  can be easily  derived.

%%%%%%%%%%%%%%%%%%%%%%%%%%%%

\subsection{The GCM-PG data set on $\Si_*$}
\lab{sec:admissibleGMCPGdatasetonSigmastar:intro}

%%%%%%%%%%%%%%%%%%%%%%%%%%%%

To initialize the PG structure of $\Mext$,  according to Proposition \ref{Prop:Initialization of PG structures:intro}, 
   we  assume not only that $\Si_*$,  in Figure \ref{fig1-introd},    is a GCM hypersurface, 
as in  Definition \ref{def:CanonicalGCM-hypersurface-intro}, but also  that  it is endowed  with a $1$-form $f$ which makes it into a GCM-PG
 data set  $\big(\Si_*, r, (e_3,e_4, \HH), f\big)$. In addition $\Si_*$ is specified by a function $u$    such that $u= c_*- r$, for some constant $c_*$ to be specified.

 Here are  therefore  the  main features of the boundary $\Si_*$.
 \begin{itemize}
 \item[-]  $\big(\Si_*, r, (e_3,e_4, \HH), f\big)$ is a   GCM-PG data set,  in the sense of Definitions \ref{definition:PGdataset-intro} and  \ref{def:CanonicalGCM-hypersurface-intro}, with $r$ 
  decreasing from its value $r_*$ on $S_*$.
 \item[-] the parameters   $a, m$ are  defined by \eqref{eq:Hawkingmass-intro} and  \eqref{eq:angular momentum:intro}, 
\item[-]  the  transition parameter $f$ is given by $f=\frac{a}{r}  d \vphi$     on $S_*$ and  transported to $\Si_*$ by 
 $\nab_\nu(r f)=0$. 
 
\item[-]  Along  $\Si_*$ we have  $u= c_*- r $  with  $c_*=1+r(S_1)$ where $S_1=\Si_*\cap \BB_1$, see  Figure \ref{fig1-introd}.

 \item[-] The function $r$ verifies a  dominance condition on $S_*$,   see  \eqref{eq:behaviorofronSstar-rough}, 
\bea
\lab{eq:behaviorofronSstar-rough:intro}
r_* \sim  u_*^{1+\dec},
\eea
where $u_*$ and $r_*$ denote respectively the value of  $u$ and $r$ on $S_*$.
\end{itemize}

%%%%%%%%%%%%%%%%%%%%%%%%%%%%%%%%%%%%

\subsection{The  PG structures  and $S$ foliations of $\Mext,\Mint$}  
\lab{section:PGstructures:intro}

%%%%%%%%%%%%%%%%%%%%%%%%%%%%%%%%%%%%

\begin{itemize}
\item[-] The  outgoing  PG structure on  $\Mext $ is fixed   from  the GCM-PG data  set  of     $\Si_*$,   with the help of  Proposition \ref{Prop:Initialization of PG structures:intro}. 
$\Mext$ is also endowed with the  $S(u, r) $ foliation where $u$ is extended from $\Si_*$ by setting $e_4(u)=0$.   The hypersurfaces of constant $u$  are timelike\footnote{They become  null at infinity.}. Note also
 that $u=u_*$ is the hypersurface separating  $\Mext$ from $\Mtop$ while  $u=u_1$ is the boundary $\BB_1$.

\item[-] $\Mext$  terminates at the inner boundary  $\TT=\Mint\cap\Mext$.    $\Mint$  is endowed  with  an ingoing  PG structure  initialized at $\TT$, defined starting  by renormalizing 
$e_3$  on $\TT$ and extending it  geodesically in $\Mint$.  We can also  extend $r$   from $\TT$ in $\Mint$ by setting $e_3(r)=-1$.  We define $\ub$ in  $\Mint $  such that it coincides with 
$u$ on $\TT$ and   $e_3(\ub)=0$.  The corresponding hypersurfaces are timelike.

\item[-]     Note that        $\Mint\cup\Mext $ in Figure \ref{fig1-introd} is not a causal region.    This is ultimately due to the fact       that the functions $u, \ub$  are not null but time-like.
   Thus,  see Remark \ref{remark-Mtop-intro},  the region   $ \Mtop$   is needed  as a completion of  $  \Mint\cup\Mext $ to a causal region.

\item[-] The  black hole parameters $(a, m)$ are  extended everywhere in $\MM$ to be constant. We also define an ingoing  PG structure on $\Mtop$ suitably initialized from the outgoing PG structure of $\Mext$ on $\{u=u_*\}$.

\end{itemize}
     
\begin{remark}
 it is important to note that $\Mext$ comes equipped  not only with  the PG frame $(e_3, e_4, \HH)$ but also with the secondary, integrable,  frame $(e_3',  e_4', \HH')$ adapted to    the spheres $S(u, r)$, i.e. $\HH'$ is tangent to 
  the  $S$ spheres.  We also have precise formulas\footnote{ The passage from the  PG  frame $(e_3, e_4, \HH)$  to the integrable  one  $(e_3', e_4', \HH')$   is obtained by the transformation formulas
   \eqref{General-frametransformation}  with  parameters $ (f, \fb, \la)$ given by \eqref{def:transition-functs:ffbla}.}   to pass form  one frame to the other whenever needed. 
\end{remark}

%%%%%%%%%%%%%%%%%%%%%%%%%%%%%%%

\subsection{GCM admissible spacetimes}

%%%%%%%%%%%%%%%%%%%%%%%%%%%%%%%

We are now ready to define our central  concept which, in addition to the geometric specifications made above  for $\Si_*$, $\Mext$, $\Mint$ and $\Mtop$, contains information 
 about decay and boundedness of  the linearized\footnote{Linearization consists for scalar quantities in subtracting Kerr values, but is slightly more subtle for 1-forms. See Definition \ref{def:renormalizationofallnonsmallquantitiesinPGstructurebyKerrvalue}.}   Ricci and curvature coefficients. As in \cite{KS}, we divide these into  the sets  we denote by   $\Ga_g, \Ga_b$.  For example, $\Ga_g$ includes in particular $\trXc, \trXbc, \Xh,  Z $, as well as  the curvature  components\footnote{In fact $A, B$ behave  even better, see \eqref{equation:defboudednessnormsMext:chap3} \eqref{equation:defdecaynormsMext:chap3}.}   $r A, rB,  r P  $. The set  $\Ga_b$ contains in particular   the Ricci coefficient $  \Xbh, H, \widecheck{\omb}  $ and the slow decaying curvature components $\Ab$ and $r\Bb$.   We refer the reader to  Definition \ref{definition.Ga_gGa_b}  for the precise definition of $\Ga_b$ and $\Ga_g$.

\begin{definition}
\lab{definition:GCMadmissible-intro}
A finite space $\MM=\Mext\cup\Mint\cup\Mtop$ as in Figure \ref{fig1-introd}   is called a GCM  admissible  spacetime  with parameters $(a, m)$ if  the following hold true.

\begin{enumerate}
\item The boundary $\Si_*$  is endowed with the  PG-GCM  data set described in section \ref{sec:admissibleGMCPGdatasetonSigmastar:intro}.

\item The domains $\Mext, \Mint, \Mtop$ are  endowed with the PG data sets and $S$ foliations  described in section \ref{section:PGstructures:intro}.

\item The linearized\footnote{Obtained for scalars by subtracting their   Kerr   values,  expressed in term of the scalar functions $(r, \th)$. The case of 1-forms is slightly more subtle.}   Ricci and curvature coefficients     verify  bootstrap  assumptions $ {(\bf BA)}_\ep$, in  $\Mext$, $\Mint$ and $\Mtop$,  measured in terms of  a small parameter $\ep$ with $\ep\gg \ep_0$, with $\ep_0$ the size of the original perturbation. The bootstrap assumptions are expressed in terms of: 
  \begin{itemize}
  \item[-] uniform decay norms  denoted  here by  $\NN^{(Dec)}_{k_{small}}$, for a maximum of $k_{small}  $ derivatives,
  
  \item[-]       $r^p$-weighted  supremum norms  denoted by
  $\NN^{(Sup)}_{k_{large}}$ for a maximum of $k_{large }$ derivatives, 
  
  \item[-]   the number $k_{small}$ is  sufficiently large  and  $k_{large} =2 k_{small}+1$.
  \end{itemize}
  Thus    $ {(\bf BA)}_\ep$ can be expressed in the form
  \bea
\NN^{(Sup)}_{k_{large}}+\NN^{(Dec)}_{k_{small}}\leq \ep.
\eea 
\end{enumerate}
\end{definition}

\begin{remark} 
The bootstrap assumptions  for decay  $\NN^{(Dec)}_{k_{small}}\leq \ep$ imply in particular  the following decay rates\footnote{Here $\dec$ is a small positive  constant.} in $\Mext$.
\beaa
|\Ga_g | &\leq &\ep r^{-2}  u^{-\frac 12 -\dec}, \qquad 
| \nab_3\Ga_g| \leq  \ep r^{-2} u^{-1-\dec}, \qquad 
|\Ga_b | \leq\ep r^{-1}  u^{-1 -\dec}.
\eeaa
In addition  each  derivatives  $\nab, \nab_4$ improve  the decay in $r$ while  each additional $\nab_3$  derivative keeps the decay unchanged. 
We express  this schematically in the  form
\beaa
|\dk^{\le k} \Ga_g | &\leq &\ep r^{-2}  u^{-\frac 12 -\dec}, \qquad 
|\dk^{\le k-1} \nab_3\Ga_g| \leq \ep r^{-2} u^{-1-\dec}, \qquad 
|\dk^{\le k} \Ga_b | \leq \ep r^{-1}  u^{-1 -\dec},
\eeaa
where $\dk=(\nab_3, r\nab_4, r\nab)$ and $\dk^{\le k} $ refers to derivatives up to order  $k\le k_{small}$.  
\end{remark}

%%%%%%%%%%%%%%%%%%%%%%

\section{Principal  temporal  frames}
  
%%%%%%%%%%%%%%%%%%%%%%  

 As mentioned earlier,  the PG structures  are  adequate  for  deriving  decay estimates  but  deficient in terms of loss of derivatives and thus inadequate   for deriving  boundedness estimates for the top derivatives of the Ricci coefficients.  Indeed  the  
 $\nab_4$ equations  for $\tr \Xb$, $\Xbh$  and $\Xib$ in Proposition    \ref{prop-nullstrandBianchi:complex:outgoing}   contain angular derivatives\footnote{This loss can be overcome for integrable foliations such as geodesic foliations and double null foliations relying on elliptic Hodge systems on 2-spheres of the foliation, but not for non integrable structures such as PG structures.} of other Ricci coefficients. Similarly,  the same situation occurs for  ingoing PG  structures where  the $ \nab_3$ equations for $\tr X$, $\Xh$, and $\Xi$   are manifestly losing derivatives. Thus, in order to derive  boundedness estimates for the top derivatives of the Ricci coefficients, we are forced to  introduce  new frames   which we call principal temporal  (PT).  These frames are  used only  in  Chapter 9 where they play an essential  role.

%%%%%%%%%%%%%%%%%%%%%%%%%%%%%%%%%%%%%%%%%%%%%

\subsection{Outgoing PT structures}

%%%%%%%%%%%%%%%%%%%%%%%%%%%%%%%%%%%%%%%%%%%%%

\begin{definition}
\lab{definition:outgoingPT-intro}
An outgoing  PT structure    $\{ (e_3, e_4, \HH), r, \th, \Jk\}$     on $\MM$ consists of  a null pair $(e_3, e_4)$, the induced horizontal structure   $\HH$,   functions $(r, \th)$, and a horizontal  1-form $\Jk$ such that the following hold true:
\begin{enumerate}
\item   $e_4$ is geodesic.

\item We have
\bea
e_4(r)=1,\qquad    e_4(\th)=0, \qquad \nab_4 (q \Jk)=0, \qquad  q= r+a i \cos \th.
\eea

\item We have 
\bea
\Hb=-\frac{a\ov{q}}{|q|^2}\Jk.
\eea
\end{enumerate}
An extended outgoing PT structure possesses,  in addition,  a scalar function $u$ verifying $e_4(u)=0$. 
\end{definition}

\begin{definition}\lab{def:initialoutgoingPTdataset-intro}
An outgoing PT  initial data set  consists of a hypersurface $\Si$ transversal to $e_4$ together with a null pair $(e_3, e_4)$, the induced horizontal structure $\HH$,  scalar   functions $(r, \th)$,  and a horizontal  1-form $\Jk$,  all defined on $\Si$.
\end{definition}

The following is precisely Lemma \ref{lemma:constructionoutgoingPTframes} in the main text.
\begin{lemma}
\lab{lemma:constructionoutgoingPTframes-intor}
Any outgoing PT  initial data set, as in Definition \ref{def:initialoutgoingPTdataset-intro}, can be  locally extended to   an outgoing PT structure.
\end{lemma}

%%%%%%%%%%%%%%%%%%%%
  
\subsection{Ingoing PT structures}

%%%%%%%%%%%%%%%%%%%%

\begin{definition}
\lab{definition:ingoingPT-intro}
An ingoing  PT structure    $\{ (e_3, e_4, \HH), r, \th, \Jk\}$     on $\MM$ consists of  a null pair $(e_3, e_4)$, the induced horizontal structure   $\HH$,  functions $(r, \th)$, and a horizontal  1-form $\Jk$ such that the following hold true:
\begin{enumerate}
\item   $e_3$ is geodesic.

\item We have
\bea
e_3(r)=-1,\qquad    e_3(\th)=0, \qquad \nab_3 (\ov{q}\Jk)=0, \qquad  q= r+a i \cos \th.
\eea

\item We have 
\bea
H=\frac{aq}{|q|^2}\Jk.
\eea
\end{enumerate}
An extended  ingoing PT structure possesses,  in addition,  a function $\ub$ verifying $e_3(\ub)=0$. 
\end{definition}

\begin{definition}\lab{def:initialingoingPTdataset-intro}
An ingoing PT  initial data set  consists of a hypersurface $\Si$ transversal to $e_3$ together with a null pair $(e_3, e_4)$, the induced horizontal structure $\HH$,  scalar   functions $(r, \th)$,  and a horizontal  1-form $\Jk$,  all defined on $\Si$.
\end{definition}

\begin{lemma}
\lab{lemma:constructioningoingPTframes-intro}
Any ingoing PT  initial data set, as in Definition \ref{def:initialingoingPTdataset-intro}, can be  locally extended to   an ingoing PT structure.
\end{lemma}

%%%%%%%%%%%%%%%%%%%%%%%%%%%%%%%  

\section{Outline of  the proof of the main theorem}
\lab{section:outlineproof-intro}

%%%%%%%%%%%%%%%%%%%%%%%%%%%%%%%

The detailed version of the main Theorem  is found in  section \ref{section:MainTheorem}.
We sketch below  the main steps in our proof. We refer the reader to  sections  \ref{sec:mainintermediateresults:chap3} and \ref{sec:endoftheproofofthemaintheorem:chap3} for more  details.
We also give an outline of the main conclusions of the Theorem.

%%%%%%%%%%%%%%%%%%%%%%%%%%%%%%%%%%%%%

\subsection{Control of the initial data} 

%%%%%%%%%%%%%%%%%%%%%%%%%%%%%%%%%%%%%
  
The  main results on  the initial data  is  stated   in  Theorem   M0   and     proved in   section \ref{sec:proofofTheoremM0},    based on the initial data  and bootstrap assumptions  in  the initial  layer   $\LL_0$.
The result  provides estimates for the main linearized quantities  restricted to the past  boundary $\BB_1\cup \BBb_1$ of our GCM admissible spacetime,  see Figure \ref{fig1-introd}.   It is important to note 
that  $\BB_1, \BBb_1$ are   not  causal,  but rather timelike, with   $\BB_1$ asymptotically null.  They are thus not to be  regarded as   fixed  hypersurfaces where  the initial data is prescribed. In fact  they change  throughout  the  continuation argument  at the heart of the  proof, while  remaining constrained to the boundary layer $\LL_0$. As in  \cite{KS} the proof of Theorem M0 is  quite subtle  due  to the fact that the spheres of the  foliation induced by $\Mext$ differ substantially   from spheres of the initial data layer $\Lext$ along  the outgoing  direction. This  anomalous behavior  reflects  the difference between  the center of mass  frames of the final and initial  Kerr states  and is as such an important feature of our result.

%%%%%%%%%%%%%%%%%%%%%%%%%%%%%%%  
  
 \subsection{Theorems M1--M5}

%%%%%%%%%%%%%%%%%%%%%%%%%%%%%%%
 
 Given a   GCM admissible spacetime, Theorems M1--M5, stated in section \ref{sec:mainintermediateresults:chap3},     improve the    decay   estimates  for $k\le k_{small} $  of  the  bootstrap assumptions $ {(\bf BA)}_\ep$ (see Definition \ref{definition:GCMadmissible-intro}),     i.e.  derive  estimates in which $\ep$ is replaced\footnote{Thus establishing  that  the bounds depend only on the initial conditions and universal constants.} by $\ep_0$.

{\bf  Theorem  M1.} \textit{Improved decay estimates for $\qf$ and $A$.}    This is our main result  concerning the  improved  decay estimates for 
$A$.  This is achieved as follows:
\begin{itemize}
\item[-] In \cite{GKS1}  we derive a  tensorial nonlinear version of  the gRW     equation.  This is  a  tensorial wave  equation for   a $2$-tensor $\qf$,   derived  from $A$ by a Chandrasekhar  transformation  of the form $\qf= \nab_3^2 A + C_1\nab_3 A  +C_2 A$, for specific  scalar  functions $C_1, C_2$.   The  wave equation for $\qf$  still contains linear terms in  $A$. Thus, in reality, we have to deal with a coupled wave-transport system  for the variables  $(\qf, A)$. The linear  theory  for  such systems,  in a  fixed Kerr background,  was derived in \cite{Ma} and \cite{D-H-R-Kerr}.  The first physical space version of the Chandrasekhar  transformation  has appeared in \cite{D-H-R}, in linear perturbations of Schwarzschild.   An adapted  nonlinear version  of the transformation plays an important role in \cite{KS}.

\item[-]  It is important to note that, as in \cite{KS},  the construction of $\qf$ and the estimates   for  $(\qf, A)$ mentioned below,   need to be done in a    global  frame for $\MM$  in which the component $H $ has   better 
decay in $\Mext$ then  the same component in the   PG frame of $\Mext$.  Simple transformation formulas allow us to  transfer  results obtained in  the global frame to results in  the   original PG  frames and vice-versa.

\item[-]  In the forthcoming paper \cite{GKS2}, we  derive  boundedness and decay estimates  for the coupled system mentioned above. The most demanding  part  is  the  derivation of  a  Morawetz type  estimate  for the coupled system $(\qf, A)$, a step  which requires  a nonlinear adaptation of  the   Anderson-Blue    \cite{A-B}       extension of the vectorfield method, mentioned earlier.   The papers   \cite{Ma}, \cite{D-H-R-Kerr}   derive the corresponding  estimate, in a fixed Kerr,  by  appealing  to  a mode decomposition.
\end{itemize}

{\bf Theorem M2.} \textit{ Improved  estimates for $\Ab$ on $\Si_*$ and $\Mint$.}  At a linear level,  $\Ab$  can be treated  in a similar manner  as $A$, i.e. we  can pass   from the Teukolsky equation for $\Ab$  to a gRW equations for a $2$ tensor $\qfb$ derived from $\Ab$   by a similar  second order transformation  formula as for $A$, with $e_3$ replaced by $e_4$. The difficulty is that the  nonlinear terms  in the gRW  equation are  not so easy to control in view of their low decay in powers of $r$. 
In  \cite{KS}  we relied  instead on a nonlinear version of the  well known   Teukolsky-Starobinsky identity  which  relates  $\nab^2_3$  derivatives of $\qf$ to four angular derivatives  of    $\aa$, see Proposition 2.3.15 in \cite{KS}, from which we can, in principle, recover $\aa$. The   non-integrable situation treated here requires in fact that we use both the gRW  equations for $\qfb$  and   an appropriate version of the  Teukolsky-Starobinsky identities.  The details   will appear in  \cite{GKS2}.

{\bf Theorem M3.}  \textit{Improved estimates for $(\Ga_g, \Ga_b$) on $\Si_*$.} Theorem M3, proved in Chapter  5 of this work, makes use of the improved  estimates for $\a$, $\aa$, and $\qf$ of Theorems M1 and M2, to   derive  improved  estimates for all other  Ricci and curvature components  restricted to $\Si_*$.   Together with Theorem M4, this  is the most  subtle part of the entire proof in that  it  depends crucially on the  properties of $\Si_*$, mentioned above,  
 and the  difficult estimates of  $\a, \aa$,   and  uses  in fact   almost all   other elements of our  overall scheme. 
Here are some of the key ideas in the proof. 
\begin{itemize}

\item[-] To derive  decay estimates for all other quantities  along $\Si_*$ it is natural    to make use of  the   transport equations  along  $\nu= e_3+b_* e_4$   induced on $\Si_*$ 
by the null structure and null Bianchi equations.

\item[-]  Integrating these  transport  equations  starting from $\BB_1\cap\Si_*$,  where we have smallness  information  in terms of $\ep_0$,     is   prohibitive  since  such an integration loses   all  decay with respect to the  $u$ factor. 
To integrate in the opposite direction, starting  from  $S_*$, we need initial conditions on $S_*$.  This is, in a nut-shell, the very  reason   our   GCM conditions were introduced. 

\item[-] Using  the propagation equations  along $\Si_*$, the GCM conditions, in particular those on the final sphere  $S_*$,  the Hodge type equations on the   $S$ spheres and   
the information already derived for $\a, \aa, \qf$,  one can  derive improved estimates  for    all  linearized Ricci and curvature coefficients  $(\Ga_g, \Ga_b)$   on $\Si_*$.
\end{itemize}

{\bf Theorem M4.} \textit{Improved estimates for $\Ga_g, \Ga_b $  in $\Mext$.}   Theorem M4,  proved in Chapter  6
of this work, extends the  estimates proved of Theorem M3  on $\Si_*$   to the entire region $\Mext$.    
There are two type of difficulties. The first, type  already encountered   in \cite{KS},  is  to derive  sufficient decay  for $\Ga_g$ quantities  in the regions near the black hole where $r$ is just bounded.  The second  type of difficulties,   are  due to the lack of integrability of the PG  structure of $\Mext$.  Here are some of the key ideas in the proof.  For a more comprehensive discussion of this step  we refer to  section \ref{section:MainResultsandStrategy}.

\begin{itemize}
\item[-]  Ideally one would use the null structure and Bianchi  equations in the $e_4$ direction to  transport information from $\Si_*$ to $\Mext$. Unfortunately, as it turns out,  many   of these equations   are strongly  overshooting in $r$.  As in  \cite{KS} we   devise new  renormalized quantities which  verify useful transport equations  which can be integrated  from $\Si_*$ in the  $e_4$ direction.  
\item[-] In \cite{KS} we  were able to  combine   these transport equations  with  elliptic Hodge systems  on the leaves of the $S$-foliation to derive estimates for the remaining quantities.  This becomes a problem  in our case due to the lack of integrability of  the  PG structure.  What we do instead is  to go back and forth between the  PG frame and  the integrable frame associated to the $S(u, r)$ spheres, and  perform  our elliptic estimates on these $S$-spheres.

 \item[-] The process  generates additional derivatives in the  direction of  the    vectorfield $\T$, analogous to the time translation of Kerr,
  which  turns out to be almost Killing.   Fortunately  the equations  obtained by  commutations with $\T$ are no longer overshooting and thus can be integrated directly from $\Si_*$. 
  \item[-] We combine  all these ingredients,  making  use of the fact that in $\Mext $ the  defining function  $ r$ is also sufficiently large,  to derive estimates for all
  elements  of $(\Ga_g, \Ga_b)$ in $\Mext$.
\end{itemize}

{\bf Theorem M5.}   \textit{Improved estimates for $\Ga_g, \Ga_b $  in $\Mint \cup \Mtop.$} This step, proved in Chapter 7, is significantly easier than Theorem M4 due to the fact that $\Mint$ is bounded in $r$ and $\Mtop$ is a local existence region. We first control the foliation of $\Mint$ and $\Mtop$ from the one of $\Mext$ respectively on $\TT$ and $\{u=u_*\}$, and then propagate this control, using transport equations along $e_3$, respectively to $\Mint$ and $\Mtop$ thanks to the equations of the corresponding ingoing PG structures.

%%%%%%%%%%%%%%%%%%%%%%%%%%%%%%%

\subsection{Extension  of GCM admissible spacetimes}    

%%%%%%%%%%%%%%%%%%%%%%%%%%%%%%%

We end the proof by invoking a continuity argument as  in \cite{KS}, see section  \ref{sec:endoftheproofofthemaintheorem:chap3}. The  argument requires  a definition  of  a set   $\UU(  u_*) $ of  GCM admissible spacetimes verifying the   bootstrap assumptions ${ \bf  BA}_\ep$ such that  $\ep$ and   the values $(r_*, u_*)$  of $(r, u) $  on $S_*$ verify
\bea
\lab{constraints:ep,r_*,u_*}
  \ep =\ep_0^{\frac{2}{3}},\qquad r_* =\de_*\ep_0^{-1}u_*^{1+\dec},
\eea
where $\de_*>0$ is a small constant satisfying $\de_*\gg \ep$. 

{{\bf Theorem M6.} \textit{The set $\UU(u_*) $ is not empty}.  More precisely,  we show that there exists $\de_0>0$ small enough such that, for sufficiently small constants $\ep_0>0$ and $\ep>0$ satisfying   the constraint in \eqref{constraints:ep,r_*,u_*},  
  \beaa
  [1, 1+\de_0]\subset \UU(u_*).
  \eeaa

   Once the  estimates assumed by ${\bf (BA)}_\ep$  have  been   improved we extend $\MM$  and its  foliation  to a larger GCM admissible  spacetime 
$\widetilde{\MM} $. This is achieved as follows.

{\bf Theorem M7.}  \textit{Extension  argument}.  We show  that any  GCM admissible spacetime in $\UU(u_*)$ for some $0<u_*<+\infty$ 
has a GCM admissible extension in in $\UU(u'_*)$    for some  $u_*'>u_*$, initialized by Theorem M0, which verifies  
the improved decay bootstrap assumptions.

The  main steps in the extension are, as in \cite{KS}:
\begin{itemize} 
\item[-]  First extend  $\MM$ and its foliation  to a strictly larger space $\MM'$.
\item[-]  To  make sure that the  extended spacetime is  GCM admissible, one has to construct  a new  GCM    hypersurface $\widetilde{\Si}_* $ in $\MM'\setminus\MM$  and  use it to  define a  new extended GCM admissible spacetime $\widetilde{\MM}$.  It is at this stage  that we have to prove the existence  of GCM spheres  in  $\MM'\setminus \MM$.  More precisely, using the bounds on  the Ricci and curvature coefficients  on $\MM'$,  defined by the extended foliation, we have to construct  GCM spheres  in  $\MM'\setminus \MM$.

\item[-]  The GCM  spheres   mentioned above  are used  as building blocks for  the new spacelike hypersurface  $\widetilde{\Si}_* $.   
The construction of   $\widetilde{\Si}_* $, similar to  that in \cite{KS},   is explicitly done in our context in \cite{Shen}.  
Once this is done we can also a construct  a new    GCM-PG data set  on $\widetilde{\Si}_* $ and 
   use it   construct  thus  the desired GCM admissible extension   $\widetilde{\MM}$.

\item[-]  One needs to  check that relative to the new  structure we     improve the original  bootstrap assumption for decay, i.e.  $\NN^{(Dec)}_{k_{small}} \les \ep_0.$
\end{itemize}

 {\bf Theorem M8.}  \textit{ Estimates for the top order derivatives.}

 The new  admissible spacetime  $\widetilde{\MM}$ is strictly larger that $\MM$ and verifies $ \NN^{(Dec)}_{k_{small}}\les \ep_0$. It still  remains to improve  the second half of the bootstrap assumptions concerning $\NN^{(Sup)}_{k_{large}}$ and show that $\NN^{(Sup)}_{k_{large}}\les \ep_0$. Part  of the argument for this is provided in Chapter 9 where,  using  the PT frame,   we derive boundedness estimates for the top derivatives of the  Ricci coefficients in terms of  bounds for the  top derivatives of the   curvature coefficients.  These latter bounds will be  derived  in a forthcoming paper \cite{KS:Kerr-B} by taking advantage of  energy-Morawetz type estimates and the   Maxwell like character    of the Bianchi identities.  Both type of estimates   are proved by an induction argument   starting with the improved  estimates for       $k\le  k_{small}$, i.e. $ 
 \NN^{(Dec)}_{k_{small}}\les \ep_0$, derived in  Theorems M0--M7.

%%%%%%%%%%%%%%

\subsection{Conclusions}

%%%%%%%%%%%%%%

The  precise version of our  main  theorem, see section \ref{section:MainTheorem}, states  a few  important  conclusions. Here are some of them.
\begin{enumerate}
\item[-]  The future null infinity  $\II_{+}$ of  the limiting space   $\MM_\infty$ is  complete.  The other future boundary of $\MM_\infty$ is given by the spacelike hypersurface $\AA$, which  can be shown  to  belong  to the complement of $\JJ^{-}(\II^{+})$.   In particular this  establishes the existence of  the event horizon $\HH^+$.
\item[-] The  spheres  $S(u, r)$  converge to   round spheres,  i.e.  $\lim_{r\to \infty} r^2 K(u, r) =1$, where $K$ is the  Gauss curvature of $S$.
\item[-] The quantities $a_\infty, m_\infty$ can  be determined by taking limits of well defined quasi-local quantities  which we define below.
\end{enumerate}

%%%%%%%%%%%%%%%%%%%%%%%%%%%%%%%%%

\subsubsection{Limits   of quasi-local  quantities on $\II^+$}    

%%%%%%%%%%%%%%%%%%%%%%%%%%%%%%%%%

The quasi-local  quantities appearing   below  are defined relative to the  integrable frame of $\Mext_\infty$, i.e. the frame $ (e_3', e_4', \HH')$ with $\HH'$ tangent to the spheres $S(u, r)$.  The first  quantity is the well known Hawking mass. The second quantity was  first defined in  \cite{KS-GCM2}.

\begin{definition}
We define the following quasi-local quantities on a given sphere  $S=S(u, r)\subset\Mext_\infty$ and its  integrable  frame $(e_3', e_4', \HH')$.
\begin{enumerate}
\item 
We define the Hawking mass of  $S$  to be 
\bea
m_H(u, r) &=&   \frac{|S(u, r)|^{1/2}}{ 4\pi^{1/2} }     \left(1+\frac{1}{16\pi}\int_{S(u, r)} \trch'\trchb'\right).
\eea
where  $ \trch', \trchb'$ are  calculated with respect to the integrable frame of $\Mext_\infty$. 
\item We define  the quasi-local  angular momentum of $S$ to be the triplet 
\bea
{\mathfrak j} _{\ell=1, p}(u, r) :=    \frac{1}{r^2}   \int _{S(u, r)} ( \curl' \b')\Jp, \qquad p=-, 0, +.
\eea
where  $\b'$,  $\curl' \b'$  are defined with respect to the integrable frame of $\Mext$ and  the triplet $\Jp$ is defined by\footnote{ Here   $\th$ and $\vphi$   are such that $ e_4(\th)= e_4(\vphi)=0$ initialized   at $\II^+$   to be  standard spherical  coordinates.}
\beaa
J^0=\cos \th, \qquad J^+=\sin\th \cos\vphi, \qquad J^{-}=\sin\th \sin \vphi
\eeaa
\end{enumerate}
\end{definition}

\begin{proposition}
The following statements hold true
\begin{enumerate}
\item  The Hawking mass has a limit as $r\to \infty$, called the Bondi mass 
\beaa
M_B(u)=\lim_{r\to \infty}  m_H(u,r).
\eeaa
 \item The Bondi mass has a limit as $u \to \infty$
 \beaa
\lim_{u\to \infty} M_B(u)=m_\infty. 
 \eeaa
 \item The  quasi-local  angular momentum   ${\mathfrak j} _{\ell=1, p}(u, r)$ of $S$ has a limit as $ r\to \infty$
  \beaa
 \JJ_{\ell=1, p}(u) &=&\lim_{r\to \infty} {\mathfrak j} _{\ell=1, p}(u, r).
 \eeaa
 \item The  triplet   $ \JJ_{\ell=1, p}(u)  $ has a limit as $u\to \infty$ and 
\beaa
\lim_{u\to \infty} \JJ_{\ell=1,0}(u)= 2 a_\infty m_\infty, \qquad   \lim_{u\to \infty} \JJ_{\ell=1,\pm}(u) =0,
\eeaa
which defines $a_\infty$.
\end{enumerate}
\end{proposition}
We also note that other definitions of angular momentum have been  proposed in the literature, see 
\cite{Sz} for a comprehensive  review, and    \cite{Rizzi} and \cite{Chen} for other  interesting proposals.

%%%%%%%%%%%%%%%%%%%%%%%%%%%%%%%

\section{Comments on the  full subextremal case}  

%%%%%%%%%%%%%%%%%%%%%%%%%%%%%%%

Though our result is  restricted to  small angular  momentum, there are reasons to hope that a full stability result, for the full subextremal  case,   is conceivable in the near future.   To start with,  the only   limitation in our work  to small values of $|a|/m$  comes from  the proof of the Morawetz type estimates for the gRW  wave equations in  \cite{GKS2} and \cite{KS:Kerr-B}. On the other hand, Morawetz estimates  for  gRW in Kerr,  in the full subextremal range,  have  been recently  derived  by  R.  Shlapentokh-Rothman and  R. Teixeira da Costa in  \cite{Y-R}. Their work rests however on  mode decompositions,   which rely strongly on the specific structure of Kerr. Thus the only   remaining obstacle, while important,  seems to be more  of a technical nature rather than conceptual.

%%%%%%%%%%%%%%%%%%
  
\section{Organization of the paper}

%%%%%%%%%%%%%%%%%%

In Chapter   2  we provide  a full descriptions of the main geometric structures needed in our work.  Chapter 3   contains   the precise version of our main theorem, its main conclusions,  as well as a full strategy of its proof,  divided in the  nine supporting intermediate results, Theorems M0--M8.   In Chapters 4 to 8,  we then give complete proofs of  Theorems, M0 and   M3--M7. Finally, we provide in Chapter 9 a  proof of Theorem M8 by assuming the curvature estimates in \cite{KS:Kerr-B}.

%%%%%%%%%%%%%%%%%%%
 
\section{Acknowledgements} 

%%%%%%%%%%%%%%%%%%

As we have pointed out in our introduction to \cite{KS}, our  results would  be inconceivable without the  remarkable achievements  obtained, first, during the  so called golden age  of black hole physics,  and then the   equally golden age   of mathematical GR in the last 30--40  years.  In addition to the references made in \cite{KS}, which have influenced our work on the nonlinear  stability of Schwarzschild,    we  have  to single  out the   works of  Andersson-Blue in  \cite{A-B}  and the works Ma in  \cite{Ma}  and Dafermos-Holzegel- Rodnianski in  \cite{D-H-R-Kerr} which play an important role  in our approach.  We    thank    E. Giorgi, our  co-author  of \cite{GKS1}-\cite{GKS2},   for her many useful comments  about this work.  We also thank  D. Shen, the author  of \cite{Shen}, for  reading part of the  manuscript.  And, again,   we thank our wives Anca and  Emilie   for their  continuing, priceless,   patience, understanding and support. We also thank Anca and Emilie  for their   help with  the  main drawings  our  papers.

The first author   has been supported  by the  NSF grant  DMS 1800841.  He would like to thank the Laboratoire Jacques-Louis Lions and    IHES   for their  hospitality during  his many visits in  Paris.  The second author is supported by the ERC grant  ERC-2016 CoG 725589 EPGR.

%%%%%%%%%%%%%%%%%%%%%%%%%%%%%%%%%%%%%%%%%%%%%

%%%%%%%%%%%%%%%%%%%%%%%%%%%%%%%%%%%%%%%%%%%%%

%%%%%%%%%%%%%%%%%%%%%%%%%%%%%%%%%%%%%%%%%%%%%

\chapter{Preliminaries}
\lab{chapter:preliminaries}

%%%%%%%%%%%%%%%%%%%%%%%%%%%%%%%%%%%%%%%%%%%%%

%%%%%%%%%%%%%%%%%%%%%%%%%%%%%%%%%%%%%%%%%%%%%

\section{A general formalism}

%%%%%%%%%%%%%%%%%%%%%%%%%%%%%%%%%%%%%%%%%%%%%

We  review  the general formalism we have introduced in \cite{GKS1}.

%%%%%%%%%%%%%%%%%%%%%%%%%

\subsection{Null pairs and horizontal structures}
\lab{subsection:review-horiz.structures}

%%%%%%%%%%%%%%%%%%%%%%%%%

Considered a fixed    null pair $e_3, e_4$, i.e. $\g(e_3, e_3)=\g(e_4, e_4)=0$,   $\g(e_3, e_4)=-2,$ and
 denote  by  $\O(\MM)$ the vectorspace  of horizontal vectorfields $X$  on $\MM$, i.e.  $\g(e_3, X)= \g(e_4, X)=0$.
  Given a fixed   orientation  on $\MM$,  with corresponding  volume form  $\in$,  we define  the induced 
 volume form on   $\O(\MM)$ by,  
 \beaa
  \in(X, Y):=\frac 1 2\in(X, Y, e_3, e_4). 
  \eeaa
 A null  frame on $\MM$ consists of a choice of horizontal vectorfields  $e_1, e_2$, such that\footnote{We use greek 
 indices for $1,2,3,4$ and latin indices $a,b$ for $1,2$.}
 \beaa
 \g(e_a, e_b)=\de_{ab}\qquad  a, b=1,2.
 \eeaa  
 The commutator $[X,Y]$ of two horizontal vectorfields
may fail however to be horizontal. We say that the pair $(e_3, e_4 )$ is integrable if   $\O(\MM)$  forms an integrable distribution,
i.e. $X, Y\in\O(\MM) $ implies that $[X,Y]\in\O(\MM)$. As  known  the  principal null pair in Kerr fails to be integrable.
Given an arbitrary vectorfield $X$ we denote by $^{(h)}X$
its  horizontal projection, $^{(h)}X=X+ \frac 1 2 \g(X,e_3)e_4+ \frac 1 2   \g(X,e_4) e_3$. A  $k$-covariant tensor-field $U$ is said to be horizontal,  $U\in \O_k(\MM)$,
if  for any $X_1,\ldots X_k$ we have 
\beaa
U(X_1,\ldots X_k)=U( ^{(h)} X_1,\ldots  ^{(h)} X_k).
\eeaa

 We define the left and right duals of horizontal 1-forms $\om$ and covariant  2-tensors $U$,
\beaa
\dual \om_{a}&=&\in_{ab}\om_b,\quad  \om^*\, _{a}=\om_b\in_{ba},\quad 
(\dual U)_{ab}=\in_{ac} U_{cb},\qquad (U^*)_{ab}=
U_{ac}\in_{cb}.
\eeaa
Note that  $\dual(\dual \om)=-\om, \,  \dual\om=-\om^*, \,  \dual (\dual U)=-U$.
Moreover  if $U$ is symmetric, then $\dual U=- U^*$ and  if  $U=\widehat U$ is symmetric, traceless,
then, $\dual \widehat U=-\widehat U^*$ is also symmetric 
traceless.

\begin{definition}\label{definition-SS-real}
We denote by $\SS_1=\SS_1(\MM)$ the  set of horizontal 1-forms  on $\MM$, and by $\SS_2=\SS_2(\MM)$
  the set of symmetric traceless   horizontal $2$-forms on $\MM$.
  \end{definition}
  
Given  $\xi, \eta\in\SS_1 $,  $U, V\in \SS_2$  we denote\footnote{Note that the definition  of $ \xi\hot \eta$ is the same as in \cite{Ch-Kl} and  differs from the one  in \cite{GKS1} by a factor $1/2$.}
\beaa
\xi\c \eta&:=&\de^{ab} \xi_a\eta_b,\qquad 
\xi\wedge\eta:=\in^{ab} \xi_a\eta_b=\xi\c\dual \eta,\qquad 
(\xi\hot \eta)_{ab}:=\xi_a \eta_b +\xi_b \eta_a-\de_{ab} \xi\c \eta,\\
(\xi\c U)_a&:=&\de^{bc} \xi_b U_{ac}, \qquad  (U\wedge V)_{ab} := \ep^{ab}U_{ac}V_{cb}.
\eeaa

For any $ X, Y\in \O(\MM)$ we define  the induced metric $g(X, Y)=\g(X, Y)$ and the null second fundamental forms
\bea
\chi(X,Y)=\g(\D_Xe_4 ,Y), \qquad \chib(X,Y)=\g(\D_X\Lb,Y).
\eea
Observe that  $\chi$ and $\chib$  are  symmetric if and only if   the horizontal structure is integrable. Indeed this  follows easily from the formulas,
 \beaa
 \chi(X,Y)-\chi(Y,X)&=&\g(\D_X e_4, Y)-\g(\D_Ye_4,X)=-\g(e_4, [X,Y]),\\
 \chib(X,Y)-\chib(Y,X)&=&\g(\D_X e_3, Y)-\g(\D_Ye_3,X)=-\g(e_3, [X,Y]).
\eeaa
  Note  that  we  can view  $\chi$ and $\chib$ as horizontal 2-covariant tensor-fields
 by extending their definition to arbitrary vectorfields  $X, Y$  by setting  $\chi(X, Y)= \chi( ^{(h)}X, ^{(h)}Y)$,  $\chib(X, Y)= \chib( ^{(h)}X, ^{(h)}Y)$.
 Given an horizontal 2-tensor $U$  we define its trace $\tr U$  and anti-trace $\atr U$
\beaa
\tr (U):=\de^{ab}U_{ab}, \qquad \atr U=\in^{ab} U_{ab}.
\eeaa
Accordingly we  decompose $\chi, \chib$ as follows,
\beaa
\chi_{ab}&=&\chih_{ab} +\frac 1 2 \de_{ab} \trch+\frac 1 2 \in_{ab}\atrch,\\
\chib_{ab}&=&\chih_{ab} +\frac 1 2 \de_{ab} \trch+\frac 1 2 \in_{ab}\atrchb.
\eeaa

We define the horizontal covariant operator $\nab$ as follows. Given $X, Y\in \O(\MM)$
 \bea
 \nab_X Y&:=&^{(h)}(\D_XY)=\D_XY- \frac 1 2 \chib(X,Y)e_4 -  \frac 1 2 \chi(X,Y) e_3.
 \eea
Note that,
\beaa
 \nab_X Y-\nab_Y X   &=&[X, Y]-\frac 1 2 (\atrchb\,  e_4+\atrch \, e_3)\in(X, Y).
\eeaa
In particular,
 \bea
 [X, Y]^\perp&=&\frac 1 2(\atrchb\,  e_4+\atrch \, e_3)\in(X, Y).
 \eea
 Also, for  all  $X,Y, Z\in \O(\MM)$,
 \beaa
 Z g (X,Y)=g(\nab_Z X, Y)+ g(X, \nab_ZY).
 \eeaa
 {\bf Remark.}\quad In the integrable case, $\nab$ coincides with the Levi-Civita connection
 of the metric induced on the integral surfaces of   $\O(\MM)$.  
 Given $X$ horizontal, $\D_4X$ and $\D_3 X$ are in general not horizontal.
 We define $\nab_4 X$ and $\nab_3 X$  to be the horizontal projections
 of the former.  More precisely,
 \beaa
 \nab_4 X&:=&^{(h)}(\D_4 X)=\D_4 X- \frac 1 2 \g(X, \D_4 e_3 ) e_4- \frac 1 2  \g(X, \D_4 e_4)  e_3 ,\\
 \nab_3 X&:=&^{(h)}(\D_3 X)=\D_3 X-   \frac 1 2 \g(X, \D_3e_3) e_3 - \frac 1 2   \g(X, \D_3 e_4 ) e_3. 
 \eeaa
The definition can be easily extended to arbitrary  $  \O_k(\MM) $ tensor-fields  $U$ 
\beaa
 \nab_4U(X_1,\ldots, X_k)&=&e_4 (U(X_1,\ldots, X_k))- \sum_i U( X_1,\ldots, \nab_4 X_i, \ldots X_k),\\
  \nab_3 U(X_1,\ldots, X_k)&=&e_3 (U(X_1,\ldots, X_k)) -\sum_i U( X_1,\ldots, \nab_3 X_i, \ldots X_k).\\
 \eeaa

%%%%%%%%%%%%%%%%%%%%%%%%%

\subsection{Ricci and curvature  coefficients}

%%%%%%%%%%%%%%%%%%%%%%%%%

Given a null frame $e_1, e_2, e_3, e_4$ we define  the general  connection coefficients,
\bea
\lab{eq:Lacoefficients}
(\La_\mu)_{\a\b} :&=& \g\big(\D_{e_\mu} e_\b,  e_\a\big)
\eea
   and  the  special  ones
 \bea
 \begin{split}
\chib_{ab}&=\g(\D_ae_3, e_b),\qquad \,\,\,\,\,\,\,\chi_{ab}=\g(\D_ae_4, e_b),\\
\xib_a&=\frac 1 2 \g(\D_3 e_3 , e_a),\qquad\,\,\,\,\, \xi_a=\frac 1 2 \g(\D_4 e_4, e_a),\\
\omb&=\frac 1 4 \g(\D_3e_3 , e_4),\qquad\,\,\,\,\,\,\, \om=\frac 1 4 \g(\D_4 e_4, e_3),\qquad \\
\etab_a&=\frac 1 2\g(\D_4 e_3, e_a),\qquad \quad \eta_a=\frac 1 2 \g(\D_3 e_4, e_a),\qquad\\
 \ze_a&=\frac 1 2 \g(\D_{e_a}e_4,  e_3).
 \end{split}
\eea
Note that these   account for all the  connection coefficients except  
\beaa
(\La_\mu)_{ab}:=\g(\D_{e_\mu} e_b, e_a), \qquad \mu=1,2,3,4,\quad a, b=1,2.
\eeaa

We have the Ricci formulas
\bea
\lab{eq:Ricciformula}
\D_a e_b&=&\nab_a e_b+\frac 1 2 \chi_{ab} e_3+\frac 1 2  \chib_{ab}e_4,\nn\\
\D_a e_4&=&\chi_{ab}e_b -\ze_a e_4,\nn\\
\D_a e_3&=&\chib_{ab} e_b +\ze_ae_3,\nn\\
\D_3 e_a&=&\nab_3 e_a +\eta_a e_3+\xib_a e_4,\nn\\
\D_3 e_3&=& -2\omb e_3+ 2 \xib_b e_b,\label{ricci}\\
\D_3 e_4&=&2\omb e_4+2\eta_b e_b,\nn\\
\D_4 e_a&=&\nab_4 e_a +\etab_a e_4 +\xi_a e_3,\nn\\
\D_4 e_4&=&-2 \om e_4 +2\xi_b e_b,\nn\\
\D_4 e_3&=&2 \om e_3+2\etab_b e_b.\nn
\eea 

For a given horizontal   1-form $\xi$, we  define the frame independent   operators\footnote{Note that the definition of $\nab\hot $ differs from the given in \cite{GKS1} by a factor $1/2$.}  
\beaa
\div\xi&=&\de^{ab}\nab_b\xi_a,\qquad 
\curl\xi=\in^{ab}\nab_a\xi_b,\qquad 
(\nab\hot \xi)_{ba}=\nab_b\xi_a+\nab_a  \xi_b-\de_{ab}( \div \xi).
\eeaa
We also define the usual  curvature components, see \cite{Ch-Kl}, 
\bea
\bsplit
\a_{ab}&=\R_{a4b4},\quad \b_a=\frac 12 \R_{a434}, \quad \rho=\frac 1 4 \R_{3434}, \quad\rhod=\frac 1 4 \dual\R_{3434},\\
\bb_a& =\frac 1 2 \R_{a334}, \quad \aa_{ab}=\R_{a3b3},
\end{split}
\eea
where $\dual\R$ denotes the Hodge dual of $\R$.

%%%%%%%%%%%%%%%%%%%%      
     
\subsection{Commutation formulas}

%%%%%%%%%%%%%%%%%%%%%

  \begin{lemma}
   \lab{lemma:comm-gen}
Let $U_{A}= U_{a_1\ldots a_k} $ be a general $k$-horizontal  tensorfield.
\begin{enumerate}
\item  We have
\bea
\bsplit
[\nab_3, \nab_b] U_A&= -\chib_{bc}\nab_c U_{A} +(\eta_b-\ze_b)  \nab_3 U_A  +\sum_{i=1}^k \big( \chib_{a_i  b}  \eta_c-\chib_{bc} \eta_{a_i} \big) U_{a_1\ldots }\,^ c \,_{\ldots a_k} \\
&+\err_{3bA}[U],\\
\err_{3bA}[U]&=\sum_{i=1}^k \Big( \chi_{a_i  c}  \xib_c-\chi_{bc}\, \xib_{a_i}-\in_{a_i c} \dual\bb_b  \Big) U_{a_1\ldots}\,^ c\,_{\ldots a_k} +\xib_b \nab_4 U_{A }.
\end{split} 
\eea

\item We have
\bea
\bsplit
[\nab_4, \nab_b] U_A&= -\chi_{bc}\nab_c U_{A}  +(\etab_b+\ze_b)  \nab_4 U_A+\sum_{i=1}^k\big(\chi_{a_i  b}  \etab_c-\chi_{bc} \etab_{a_i} \big) U_{a_1\ldots c\ldots a_k} \\
&+\err_{4bA}[U],\\
\err_{4bA}[U]&=\sum_{i=1}^k \Big( \chib_{a_i  c}  \xi_c-\chib_{bc}\, \xi_{a_i}+\in_{a_i c} \dual\b_b  \Big) U_{a_1\ldots}\,^ c\,_{\ldots a_k} +\xi_b \nab_3  U_{A }.
\end{split} 
\eea

\item We have
\bea
\bsplit
\, [\nab_4, \nab_3] U_A&= 2(\etab_b-\eta_b ) \nab_b U_A + 2\sum_{i=1}^k\big( \eta_{a_i} \etab_b-\etab_{a_i} \eta_b- \in_{a_i b}\dual \rho) U_{a_1\ldots}\,^b\,_{\ldots a_k} \\
&+ 2 \om \nab_3 U_A -2\omb \nab_4 U_A +\err_{43A},\\
\err_{43A}&= 2\sum_{i=1}^k \big( \xib_{a_i}  \xi_b- \xi_{a_i}  \xib_b )U_{a_1\ldots} \,^b\,_{\ldots a_k}.
\end{split}
\eea
\end{enumerate}
\end{lemma}
 
\begin{proof}
See Lemma 3.8 in \cite{GKS1}.
\end{proof}

%%%%%%%%%%%%%%%%%%%%%%%%%%

\subsection{Null structure  and Bianchi equations}

%%%%%%%%%%%%%%%%%%%%%%%%%%

The null structure  equations  are given in the following proposition, see \cite{GKS1},  and   \cite{Ch-Kl}  for the  integrable case\footnote{Note that the term  $\xib\wedge(-\eta+\etab+2\ze)$ in $\nab_3\atrchb$  differs from that in  \cite{Ch-Kl}, see (7.4.1c) on page 165.}.

\begin{proposition}
\label{prop-nullstr}
We have
\beaa
\nab_3\trchb&=&-|\chibh|^2-\frac 1 2 \big( \trchb^2-\atrchb^2\big)+2\div\xib  - 2\omb \trchb +  2 \xib\c(\eta+\etab-2\ze),\\
\nab_3\atrchb&=&-\trchb\atrchb +2\curl \xib -2\omb\atrchb+ 2 \xib\wedge(-\eta+\etab+2\ze),\\
\nab_3\chibh&=&-\trchb\,  \chibh+ \nab\hot \xib- 2\omb\, \chibh+   \xib\hot(\eta+\etab-2\ze)-\aa,
\eeaa
\beaa
\nab_3\trch
&=& -\chibh\c\chih -\frac 1 2 \trchb\trch+\frac 1 2 \atrchb\atrch    +   2   \div \eta+ 2 \omb \trch + 2 \big(\xi\c \xib +|\eta|^2\big)+ 2\rho,\\
\nab_3\atrch
&=&-\chibh\wedge\chih-\frac 1 2(\atrchb \trch+\trchb\atrch)+ 2 \curl \eta + 2 \omb \atrch + 2 \xib\wedge\xi  -  2 \dual \rho,\\
\nab_3\chih
&=&-\frac 1 2 \big( \trch \chibh+\trchb \chih\big)-\frac 1 2 \big(-\dual \chibh \, \atrch+\dual \chih\,\atrchb\big)
+\nab\hot \eta +2 \omb \chih\\
&+&\xib\hot\xi +\eta\hot\eta,
\eeaa
\beaa
\nab_4\trchb
&=& -\chih\c\chibh -\frac 1 2 \trch\trchb+\frac 1 2 \atrch\atrchb    +  2   \div \etab+ 2 \om \trchb + 2\big( \xi\c \xib +|\etab|^2\big)+2\rho,\\
\nab_4\atrchb
&=&-\chih\wedge\chibh-\frac 1 2(\atrch \trchb+\trch\atrchb)+ 2 \curl \etab + 2 \om \atrchb + 2 \xi\wedge\xib+2 \dual \rho,\\
\nab_4\chibh
&=&-\frac 1 2 \big( \trchb \chih+\trch \chibh\big)-\frac 1 2 \big(-\dual \chih \, \atrchb+\dual \chibh\,\atrch\big)
+\nab\hot \etab +2 \om \chibh\\
&+&\xi\hot\xib +\etab\hot\etab,
\eeaa
\beaa
\nab_4\trch&=&-|\chih|^2-\frac 1 2 \big( \trch^2-\atrch^2\big)+ 2 \div\xi  - 2 \om \trch + 2   \xi\c(\etab+\eta+2\ze),\\
\nab_4\atrch&=&-\trch\atrch + 2 \curl \xi - 2 \om\atrch+ 2 \xi\wedge(-\etab+\eta-2\ze),\\
\nab_4\chih&=&-\trch\,  \chih+\nab\hot \xi- 2 \om \chih+    \xi\hot(\etab+\eta+2\ze)-\a.
\eeaa
\beaa
\nab_3 \ze+2\nab\omb&=& -\chibh\c(\ze+\eta)-\frac{1}{2}\trchb(\ze+\eta)-\frac{1}{2}\atrchb(\dual\ze+\dual\eta)+ 2 \omb(\ze-\eta)\\
&&+\hch\c\xib+\frac{1}{2}\trch\,\xib+\frac{1}{2}\atrch\dual\xib +2 \om \xib -\bb,
\\
\nab_4 \ze -2\nab\om&=& \chih\c(-\ze+\etab)+\frac{1}{2}\trch(-\ze+\etab)+\frac{1}{2}\atrch(-\dual\ze+\dual\etab)+2 \om(\ze+\etab)\\
&& -\chibh\c\xi -\frac{1}{2}\trchb\,\xi-\frac{1}{2}\atrchb\dual\xi -2 \omb \xi -\b,
\\
\nab_3 \etab -\nab_4\xib &=& -\chibh\c(\etab-\eta) -\frac{1}{2}\trchb(\etab-\eta)+\frac{1}{2}\atrchb(\dual\etab-\dual\eta) -4 \om \xib  +\bb, \\
\nab_4 \eta    -    \nab_3\xi &=& -\chih\c(\eta-\etab) -\frac{1}{2}\trch(\eta-\etab)+\frac{1}{2}\atrch(\dual\eta-\dual\etab)-4\omb \xi -\b.
\eeaa
\beaa
\nab_3\om+\nab_4\omb&=&   \rho  +4\om\omb +\xi\c \xib +(\eta-\etab)\c\ze -\eta\c\etab.
\eeaa

We also have the Codazzi equations
\beaa
\div\chih +\ze\c\chih &=& \frac{1}{2}\nab\trch+\frac{1}{2}\trch\ze -\frac{1}{2}\dual\nab\atrch-\frac{1}{2}\atrch\dual\ze -\atrch\dual\eta-\atrch\dual\xi -\b,\\
\div\chibh -\ze\c\chibh &=& \frac{1}{2}\nab\trchb-\frac{1}{2}\trchb\ze -\frac{1}{2}\dual\nab\atrchb+\frac{1}{2}\atrchb\dual\ze -\atrchb\dual\etab-\atrchb\dual\xib +\bb,
\eeaa
\beaa
\curl\ze&=&-\frac 1 2 \chih\wedge\chibh   +\frac 1 4 \big(  \trch\atrchb-\trchb\atrch   \big)+\om \atrchb -\omb\atrch+\dual \rho.
\eeaa
\end{proposition}

The null Bianchi equations  are given below, see  \cite{GKS1},  and     \cite{Ch-Kl}  for the  integrable case.
    \begin{proposition}\label{prop:bianchi} 
    We have
    \beaa
    \nab_3\a- \nab\hot \b&=&-\frac 1 2 \big(\trchb\a+\atrchb\dual \a)+4\omb \a+
 (\ze+4\eta)\hot \b - 3 (\rho\chih +\rhod\dual\chih),\\
\nab_4\beta - \div\a &=&-2(\trch\beta-\atrch \dual \b) - 2  \om\b +\a\c  (2 \ze +\etab) + 3  (\xi\rho+\dual \xi\rhod),\\
     \nab_3\b - (\nab\rho+\dual\nab\rhod) &=&-(\trchb \b+\atrchb \dual \b)+2 \omb\,\b+2\bb\c \chih+3 (\rho\eta+\rhod\dual \eta)+    \a\c\xib,\\
 \nab_4 \rho-\div \b&=&-\frac 3 2 (\trch \rho+\atrch \rhod)+(2\etab+\ze)\c\b-2\xi\c\bb-\frac 1 2 \chibh \c\a,\\
   \nab_4 \rhod+\curl\b&=&-\frac 3 2 (\trch \rhod-\atrch \rho)-(2\etab+\ze)\c\dual \b-2\xi\c\dual \bb+\frac 1 2 \chibh \c\dual \a, \\
     \nab_3 \rho+\div\bb&=&-\frac 3 2 (\trchb \rho -\atrchb \rhod) -(2\eta-\ze) \c\bb+2\xib\c\b- \frac{1}{2}\chih\c\aa,
 \\
   \nab_3 \rhod+\curl\bb&=&-\frac 3 2 (\trchb \rhod+\atrchb \rho)- (2\eta-\ze) \c\dual \bb-2\xib\c\dual\b-\frac 1 2 \chih\c\dual \aa,\\
     \nab_4\bb + \nab\rho-\dual\nab\rhod &=&-(\trch \bb+ \atrch \dual \bb)+ 2\om\,\bb+2\b\c \chibh
    - 3 (\rho\etab-\rhod\dual \etab)-    \aa\c\xi,\\
     \nab_3\bb +\div\aa &=&-2(\trchb\,\bb-\atrchb \dual \bb)- 2  \omb\bb-\aa\c(-2\ze+\eta) - 3  (\xib\rho-\dual \xib \rhod),\\
     \nab_4\aa+ \nab\hot \bb&=&-\frac 1 2 \big(\trch\aa +\atrch\dual \aa)+4\om \aa+
 (\ze-4\etab)\hot \bb - 3  (\rho\chibh -\rhod\dual\chibh).
\eeaa
    \end{proposition}

%%%%%%%%%%%%%%%%%%%%%%%%

\subsection{Main equations in complex form}

%%%%%%%%%%%%%%%%%%%%%%%%

In this section, we recall the complex notations introduced in \cite{GKS1} that will allow to simplify the main equations.
\begin{definition}
\beaa
A:=\a+i\dual\a, \quad B:=\b+i\dual\b, \quad P:=\rho+i\dual\rho,\quad \Bb:=\bb+i\dual\bb, \quad \Ab:=\aa+i\dual\aa,
\eeaa      
 and   
\beaa
&& X=\chi+i\dual\chi, \quad \Xb=\chib+i\dual\chib, \quad H=\eta+i\dual \eta, \quad \Hb=\etab+i\dual \etab, \quad Z=\ze+i\dual\ze, \\ 
&& \Xi=\xi+i\dual\xi, \quad \Xib=\xib+i\dual\xib.
\eeaa    
In particular, note that 
\beaa
\tr X = \trch-i\atrch, \quad \Xh=\chih+i\dual\chih, \quad \tr\Xb = \trchb -i\atrchb, \quad \Xbh=\chibh+i\dual\chibh.
\eeaa
\end{definition}

\begin{remark}
$A$, $B$, $\Ab$, $\Bb$, $X$, $\Xb$, $H$, $\Hb$, $\Xi$, $\Xib$ and $Z$ are anti-self dual tensors, i.e. they verify $\dual U=-iU$. 
\end{remark}

\begin{definition}
We define derivatives of complex quantities as follows
\begin{itemize}
\item For two scalar functions $a$ and $b$, we define
\beaa
\DD(a+ib) &:=& (\nabla+i\dual\nabla)(a+ib).
\eeaa

\item For a 1-form $f$, we define
\beaa
\DD\c(f+i\dual f) &:=& (\nabla+i\dual\nabla)\c(f+i\dual f)
\eeaa
and  
\beaa
\DD\hot(f+i\dual f) &:=& (\nabla+i\dual\nabla)\hot(f+i\dual f).
\eeaa

\item For a symmetric traceless 2-form $u$, we define
\beaa
\DD\c(u+i\dual u) &:=& (\nabla+i\dual\nabla)\c(u+i\dual u).
\eeaa
\end{itemize}
\end{definition}

These complex notations allow us to rewrite the null structure equations as follows, see \cite{GKS1}.  
\begin{proposition}
\label{prop-nullstr:complex} 
The following hold true\footnote {The convention  for $\hot$ is the one  from \cite{Ch-Kl} and differs by a factor of $2$  from the one in \cite{GKS1}.}
\beaa
\nab_3\tr\Xb +\frac{1}{2}(\tr\Xb)^2+2\omb\,\tr\Xb &=& \DD\c\ov{\Xib}+\Xib\c\ov{\Hb}+\ov{\Xib}\c(H-2Z)-\frac{1}{2}\Xbh\c\ov{\Xbh},\\
\nab_3\Xbh+\Re(\tr\Xb) \Xbh+ 2\omb\,\Xbh&=& \frac{1}{2}\DD\hot \Xib+  \frac{1}{2} \Xib\hot(H+\Hb-2Z)-\Ab,
\eeaa
\beaa
\nab_3\tr X +\frac{1}{2}\tr\Xb\tr X-2\omb\tr X &=& \DD\c\ov{H}+H\c\ov{H}+2P+\Xib\c\ov{\Xi}-\frac{1}{2}\Xbh\c\ov{\Xh},\\
\nab_3\widehat{X} +\frac{1}{2}\tr\Xb\, \widehat{X} -2\omb\widehat{X} &=& \frac{1}{2}\DD\hot H  +\frac{1}{2}H\hot H -\frac{1}{2}\ov{\tr X} \widehat{\Xb}+\frac{1}{2}\Xib\hot\Xi,
\eeaa
\beaa
\nab_4\tr\Xb +\frac{1}{2}\tr X\tr\Xb -2\om\tr\Xb &=& \DD\c\ov{\Hb}+\Hb\c\ov{\Hb}+2\ov{P}+\Xi\c\ov{\Xib}-\frac{1}{2}\Xh\c\ov{\Xbh},\\
\nab_4\widehat{\Xb} +\frac{1}{2}\tr X\, \widehat{\Xb} -2\om\widehat{\Xb} &=& \frac{1}{2}\DD\hot\Hb  +\frac{1}{2}\Hb\hot\Hb -\frac{1}{2}\ov{\tr\Xb} \widehat{X}+\frac{1}{2}\Xi\hot\Xib,
\eeaa
\beaa
\nab_4\tr X +\frac{1}{2}(\tr X)^2+2\om\tr X &=& \DD\c\ov{\Xi}+\Xi\c\ov{H}+\ov{\Xi}\c(H+2Z)-\frac{1}{2}\Xh\c\ov{\Xh},\\
\nab_4\Xh+\Re(\tr X)\Xh+ 2\om\Xh&=& \frac{1}{2}\DD\hot \Xi+  \frac{1}{2} \Xi\hot(\Hb+H+2Z)-A.
\eeaa
Also,
\beaa
\nab_3Z +\frac{1}{2}\tr\Xb(Z+H)-2\omb(Z-H) &=& -2\DD\omb -\frac{1}{2}\widehat{\Xb}\c(\ov{Z}+\ov{H})\\
&&+\frac{1}{2}\tr X\Xib+2\om\Xib -\Bb+\frac{1}{2}\ov{\Xib}\c\Xh,\\
\nab_4Z +\frac{1}{2}\tr X(Z-\Hb)-2\om(Z+\Hb) &=& 2\DD\om +\frac{1}{2}\widehat{X}\c(-\ov{Z}+\ov{\Hb})\\
&&-\frac{1}{2}\tr\Xb\Xi-2\omb\Xi -B-\frac{1}{2}\ov{\Xi}\c\Xbh,\\
\nab_3\Hb -\nab_4\Xib &=&  -\frac{1}{2}\ov{\tr\Xb}(\Hb-H) -\frac{1}{2}\Xbh\c(\ov{\Hb}-\ov{H}) -4\om\Xib+\Bb,\\
\nab_4H -\nab_3\Xi &=&  -\frac{1}{2}\ov{\tr X}(H-\Hb) -\frac{1}{2}\Xh\c(\ov{H}-\ov{\Hb}) -4\omb\Xi-B,
\eeaa
and
\beaa
\nab_3\om+\nab_4\omb -4\om\omb -\xi\c \xib -(\eta-\etab)\c\ze +\eta\c\etab&=&   \rho.
\eeaa
Also,
\beaa
\frac{1}{2}\ov{\DD}\c\Xh +\frac{1}{2}\Xh\c\ov{Z} &=& \frac{1}{2}\DD\ov{\tr X}+\frac{1}{2}\ov{\tr X}Z-i\Im(\tr X)(H+\Xi)-B,\\
\frac{1}{2}\ov{\DD}\c\Xbh -\frac{1}{2}\Xbh\c\ov{Z} &=& \frac{1}{2}\DD\ov{\tr\Xb}-\frac{1}{2}\ov{\tr\Xb}Z-i\Im(\tr\Xb)(\Hb+\Xib)+\Bb,
\eeaa
and,
\beaa
\curl\ze&=&-\frac 1 2 \chih\wedge\chibh   +\frac 1 4 \big(  \trch\atrchb-\trchb\atrch   \big)+\om \atrchb -\omb\atrch+\dual \rho.
\eeaa
\end{proposition}
    
The complex notations allow us to rewrite the Bianchi identities as follows, see \cite{GKS1}.  
  \begin{proposition}\label{prop:bianchi:complex} 
    We have\footnote{ Here $\hot$  is the usual   standard one from  \cite{Ch-Kl}.}, 
 \beaa
 \nab_3A -\frac{1}{2}\DD\hot B &=& -\frac{1}{2}\tr\Xb A+4\omb A +\frac{1}{2}(Z+4H)\hot B -3\ov{P}\Xh,\\
\nab_4B -\frac{1}{2}\ov{\DD}\c A &=& -2\ov{\tr X} B -2\om B +\frac{1}{2}A\c  (\ov{2Z +\Hb})+3\ov{P} \,\Xi,\\
\nab_3B-\DD\ov{P} &=& -\tr\Xb B+2\omb B+\ov{\Bb}\c \Xh+3\ov{P}H +\frac{1}{2}A\c\ov{\Xib},\\
\nab_4P -\frac{1}{2}\DD\c \ov{B} &=& -\frac{3}{2}\tr X P +\frac{1}{2}(2\Hb+Z)\c\ov{B} -\ov{\Xi}\c\Bb -\frac{1}{4}\Xbh\c \ov{A}, \\
\nab_3P +\frac{1}{2}\ov{\DD}\c\Bb &=& -\frac{3}{2}\ov{\tr\Xb} P -\frac{1}{2}(\ov{2H-Z})\c\Bb +\Xib\c \ov{B} -\frac{1}{4}\ov{\Xh}\c\Ab, \\
\nab_4\Bb+\DD P &=& -\tr X\Bb+2\om\Bb+\ov{B}\c \Xbh-3P\Hb -\frac{1}{2}\Ab\c\ov{\Xi},\\
\nab_3\Bb +\frac{1}{2}\ov{\DD}\c\Ab &=& -2\ov{\tr\Xb}\,\Bb -2\omb\,\Bb -\frac{1}{2}\Ab\c (\ov{-2Z +H})-3P \,\Xib,\\
\nab_4\Ab +\frac{1}{2}\DD\hot\Bb &=& -\frac{1}{2}\tr X \Ab+4\om\Ab +\frac{1}{2}(Z-4\Hb)\hot \Bb -3P\Xbh.
\eeaa
    \end{proposition}

%%%%%%%%%%%%%%%%%%%%%%%%%%%%%%%%%%%%%%

\section{Null frame transformations}

%%%%%%%%%%%%%%%%%%%%%%%%%%%%%%%%%%%%%%

Consider two null frames  $(e_4, e_3, e_1, e_2)$ and $(e_4', e_3', e_1', e_2')$ on $\MM$ with $H=\{e_3, e_4\}^\perp$ and $H'=\{e_3', e_4'\}^\perp$  the corresponding horizontal structures. We denote by $\Ga', \Ga$ the connection coefficients relative to the two frames.  We  denote by $\nab, \nab \hot, \div, \curl , \nab_3, \nab_4$    the   standard operators   corresponding to   to $H$ and by  $\nab',  \nab' \hot,  \div', \curl', \nab'_3, \nab_4 ' $ those corresponding to $H'$.  The goal  is to establish transition formulas between the Ricci  and curvature  coefficients  of the  two frames.

%%%%%%%%%%%%%%%%%%%%%%%%%%%%%%%%%

\subsection{Transformation between two null frames}

%%%%%%%%%%%%%%%%%%%%%%%%%%%%%%%%%

\begin{lemma}
\lab{Lemma:Generalframetransf}
The following transformation formulas hold true.
\begin{enumerate}
\item A general null transformation\footnote{ In full generality, one could also rotate $e_1', e_2'$, but this would not change the horizontal structure and, as it  turns out, is in fact  not needed. The dot product and magnitude  $|\c |$ are taken with respect to the standard euclidian norm of $\RRR^2$.}  between two  null frames $(e_4, e_3, e_1, e_2)$ and $(e_4', e_3', e_1', e_2')$ on $\MM$
  can be written in   the form,
 \bea
 \lab{General-frametransformation}
 \bsplit
  e_4'&=\la\left(e_4 + f^b  e_b +\frac 1 4 |f|^2  e_3\right),\\
  e_a'&= \left(\de_a^b +\frac{1}{2}\fb_af^b\right) e_b +\frac 1 2  \fb_a  e_4 +\left(\frac 1 2 f_a +\frac{1}{8}|f|^2\fb_a\right)   e_3,\\
 e_3'&=\la^{-1}\left( \left(1+\frac{1}{2}f\c\fb  +\frac{1}{16} |f|^2  |\fb|^2\right) e_3 + \left(\fb^b+\frac 1 4 |\fb|^2f^b\right) e_b  + \frac 1 4 |\fb|^2 e_4 \right),
 \end{split}
 \eea
 where $\la$, $f_a= f^a , \fb_a= \fb^a $  are    scalar functions, called the transition coefficients of the change of frame. Note, in particular,
  \beaa
   e_a'&=& e_a +\frac 1 2  \fb_a \la^{-1} e_4'  +\frac 1 2 f_a e_3,\\
    e_3'&=&\la^{-1}\left(  e_3 +  \fb^ae_a' -\frac 1 4 |\fb|^2\la^{-1} e_4'\right). 
 \eeaa

  \item The inverse transformation is given by the formulas
  \bea
 \lab{General-frametransformation'}
 \bsplit
  e_4&=\la'\left(e'_4 + f_b'  e'_b +\frac 1 4 |f'|^2  e'_3\right),\\
  e_a&= \left(\de_a^b  +\frac{1}{2}\fb'_a f'^b \right) e'_b +\frac 1 2  \fb'_a  e'_4 +\left(\frac 1 2 f'_a +\frac{1}{8}|f'|^2\fb'_a\right)   e'_3,\\
 e_3&=(\la')^{-1}\left( \left(1+\frac{1}{2}f'\c\fb'  +\frac{1}{16} |f'|^2  |\fb'|^2\right) e'_3 + \left(\fb'^b+\frac 1 4 |\fb'|^2f'^b\right) e'_b  + \frac 1 4 |\fb'|^2 e'_4 \right),
 \end{split}
 \eea
 where
\bea
\lab{relations:laffb-to-primes}
\bsplit
\la' &= \la^{-1} \left(1+\frac{1}{2}f\c\fb  +\frac{1}{16} |f|^2  |\fb|^2\right),\\
f_a'  &= -\frac{\la}{1+\frac{1}{2}f\c\fb  +\frac{1}{16} |f|^2  |\fb|^2}\left(f_a +\frac{1}{4}|f|^2\fb_a\right),\\
\fb_a' &= -\la^{-1}\left(\fb_a+\frac 1 4 |\fb|^2f_a\right).
\end{split}
\eea
Moreover
\bea
\lab{relations:ffbf'fb'}
\bsplit
\fb'_a f_b'=\fb_b f_a, \qquad  \la' |f'|^2 = \la |f|^2, \qquad 
(\la')^{-1}|\fb'|^2= \la^{-1} |\fb|^2
\end{split}
\eea
 Denoting    $F=(f, \fb,  \la-1)$,   we also write, for small $|F|$,      
\bea
\lab{eq:f-f'transf-simplified}
\bsplit
\la' &= \la^{-1} \left(1+\frac{1}{2}f\c\fb \right)+ O(|F|^3),\\
f'_a&=-\la f_a + O(|F|^3),\\
\fb'_a&=  -\la^{-1}\fb_a + O(|F|^3).
\end{split}
\eea
  \end{enumerate}
  \end{lemma}
 
\begin{proof}
 For the first part of the lemma, see  section 3.1   in \cite{KS-GCM1}. To check the second part  we make use of   \eqref{General-frametransformation} to deduce
\beaa
\g(e_a', e_4) &=& - f_a -\frac{1}{4}|f|^2\fb_a,\\ 
\g(e_3', e_4)&=& -2\la^{-1} \left(1+\frac{1}{2}f\c\fb  +\frac{1}{16} |f|^2  |\fb|^2\right),\\
\g(e_3', e_a)&=&\la^{-1}\left(\fb_a+\frac 1 4 |\fb|^2f_a\right).
\eeaa
Similarly,   using   \eqref{General-frametransformation'}, we deduce
\beaa
\g(e_4, e_a') &=& \la' f_a',\\ 
\g(e_4, e_3')&=& -2\la',\\
\g(e_a, e_3')&=&-\fb_a'.
\eeaa
We infer
\beaa
- f_a -\frac{1}{4}|f|^2\fb_a &=& \g(e_a', e_4)=\la' f_a',\\
-2\la^{-1} \left(1+\frac{1}{2}f\c\fb  +\frac{1}{16} |f|^2  |\fb|^2\right) &=& \g(e_3', e_4)=-2\la',\\
\la^{-1}\left(\fb_a+\frac 1 4 |\fb|^2f_a\right) &=& \g(e_3', e_a)=-\fb_a',
\eeaa
and hence
\beaa
\la' f_a' &=& - f_a -\frac{1}{4}|f|^2\fb_a,\\
\la' &=& \la^{-1} \left(1+\frac{1}{2}f\c\fb  +\frac{1}{16} |f|^2  |\fb|^2\right),\\
\fb_a' &=& -\la^{-1}\left(\fb_a+\frac 1 4 |\fb|^2f_a\right).
\eeaa
In particular
\beaa
f_a'  &=& -{\la'}^{-1}\left(f_a +\frac{1}{4}|f|^2\fb_a\right)\\
&=& -\frac{\la}{1+\frac{1}{2}f\c\fb  +\frac{1}{16} |f|^2  |\fb|^2}\left(f_a +\frac{1}{4}|f|^2\fb_a\right).
\eeaa
The identities \eqref{relations:ffbf'fb'}  follow from  the relations  induced by $\g(e_a, e_b')$, $ \g(e_3. e_3') $, $\g(e_4, e_4')$.
This concludes the proof of the lemma.
\end{proof}

\begin{remark}
As a consequence of the above lemma  both $f$ and $f'$ can be regarded
 as horizontal  vectors on both $H=\{e_3, e_4\}^\perp$ and $H'=\{e_3', e_4'\}^\perp$.
\end{remark}

%%%%%%%%%%%%%%%%%%%%%%%%%%%%%%%%%%%%%%%%%%%%%

\subsection{Transformation formulas for Ricci and curvature  coefficients}

%%%%%%%%%%%%%%%%%%%%%%%%%%%%%%%%%%%%%%%%%%%%%

\begin{proposition}
\lab{Proposition:transformationRicci}
The following transformation formulas between $\Ga$ and $\Ga'$ hold true.
\begin{itemize}
\item The transformation formula for $\xi$ is given by 
\bea
\bsplit
\la^{-2}\xi' &= \xi +\frac{1}{2}\nab_{\la^{-1}e_4'}f+\frac{1}{4}(\trch f -\atrch\dual f)+\om f +\err(\xi,\xi'),\\
\err(\xi,\xi') &= \frac{1}{2}f\c\chih+\frac{1}{4}|f|^2\eta+\frac{1}{2}(f\c \ze)\,f -\frac{1}{4}|f|^2\etab \\
&+ \la^{-2}\left( \frac{1}{2}(f\c\xi')\,\fb+ \frac{1}{2}(f\c\fb)\,\xi'   \right)  +\lot
       \end{split}
\eea

\item The transformation formula for $\xib$ is given by 
\bea
\bsplit
\la^2\xib' &= \xib + \frac{1}{2}\la\nab_3'\fb' +    \omb\,\fb + \frac{1}{4}\trchb\,\fb - \frac{1}{4}\atrchb\dual\fb +\err(\xib, \xib'),\\
\err(\xib, \xib') &=   \frac{1}{2}\fb\c\chibh - \frac{1}{2}(\fb\c\ze)\fb +  \frac 1 4 |\fb|^2\etab  -\frac 1 4 |\fb|^2\eta'+\lot
       \end{split}
\eea

\item The transformation formulas for $\chi $ are  given by 
\bea
\bsplit
\la^{-1}\trch' &= \trch  +  \div'f + f\c\eta + f\c\ze+\err(\trch,\trch')\\
\err(\trch,\trch') &= \fb\c\xi+\frac{1}{4}\fb\c\left(f\trch -\dual f\atrch\right) +\om (f\c\fb)  -\omb |f|^2 \\
& -\frac{1}{4}|f|^2\trchb -  \frac 1 4 ( f\c\fb) \la^{-1}\trch' +\frac 1 4  (\fb\wedge f) \la^{-1}\atrch'+\lot,
\end{split}
\eea
\bea
\bsplit
\la^{-1}\atrch' &= \atrch  +  \curl'f + f\wedge\eta + f\wedge\ze +\err(\atrch,\atrch'),\\
\err(\atrch,\atrch') &= \fb\wedge\xi+\frac{1}{4}\left(\fb\wedge f\trch +(f\c\fb)\atrch\right) +\om f\wedge\fb   \\
& -\frac{1}{4}|f|^2\atrchb -  \frac 1 4 ( f\c\fb) \la^{-1}\atrch' +\frac 1 4   \la^{-1}(f\wedge \fb)\trch'+\lot,
\end{split}
\eea
\bea
\bsplit
\la^{-1}\chih' &= \chih  +  \nab'\hot f + f\hot\eta + f\hot\ze+\err(\chih,\chih')\\
\err(\chih,\chih') &=\fb\hot\xi+\frac{1}{4}\fb\hot\left(f\trch -\dual f\atrch\right) +\om f\hot\fb  -\omb f\hot f -\frac{1}{4}|f|^2\atrchb\\
&  +\frac 1 4  (f\hot\fb) \la^{-1}\trch' +\frac 1 4  (\dual f\hot\fb) \la^{-1}\atrch' +\frac 1 2  \fb\hot (f\c\la^{-1}\chih')+\lot
\end{split}
\eea

\item The transformation formulas for $\chib $ are  given by 
\bea
\bsplit
\la\trchb' &= \trchb +\div'\fb +\fb\c\etab  -  \fb\c\ze +\err(\trchb, \trchb'),\\
\err(\trchb, \trchb') &= \frac{1}{2}(f\c\fb)\trchb+f\c\xib -|\fb|^2\om + (f\c\fb)\omb   -\frac 1 4 |\fb|^2\la^{-1}\trch'+\lot,
\end{split}
\eea
\bea
\bsplit
\la\atrchb' &= \atrchb +\curl'\fb +\fb\wedge\etab  -  \ze\wedge\fb+\err(\atrchb, \atrchb'),\\
\err(\atrchb, \atrchb') &= \frac{1}{2}(f\c\fb)\atrchb+f\wedge\xib  + (f\wedge\fb)\omb   -\frac 1 4 |\fb|^2\la^{-1}\atrch'\\
&+\lot,
\end{split}
\eea
\bea
\bsplit
\la\chibh' &= \chibh +\nab'\hot\fb +\fb\hot\etab  -  \fb\hot\ze +\err(\chibh, \chibh'),\\
\err(\chibh, \chibh') &= \frac{1}{2}(f\hot\fb)\trchb  +f\hot\xib -(\fb\hot\fb)\om + (f\hot\fb)\omb   -\frac 1 4 |\fb|^2\la^{-1}\chih'+\lot
\end{split}
\eea

\item  The transformation formula for $\ze$ is given by 
\bea
\bsplit
\ze' &= \ze -\nab'(\log\la)  -\frac{1}{4}\trchb f +\frac{1}{4}\atrchb \dual f +\om\fb -\omb f +\frac{1}{4}\fb\trch\\
&+\frac{1}{4}\dual\fb\atrch+\err(\ze, \ze'),\\
\err(\ze, \ze') &= -\frac{1}{2}\chibh\c f + \frac{1}{2}(f\c\ze)\fb -  \frac{1}{2}(f\c\etab)\fb +\frac{1}{4}\fb(f\c\eta) + \frac{1}{4}\fb(f\c\ze)  \\
& + \frac{1}{4}\dual\fb(f\wedge\eta) + \frac{1}{4}\dual\fb(f\wedge\ze) +  \frac{1}{4}\fb\div'f  + \frac{1}{4}\dual\fb \curl'f +\frac{1}{2}\la^{-1}\fb\c\chih' \\
&  -\frac{1}{16}(f\c\fb)\fb\la^{-1}\trch' +\frac{1}{16}  (\fb\wedge f) \fb\la^{-1}\atrch'  -  \frac{1}{16}\dual\fb ( f\c\fb) \la^{-1}\atrch'\\
& +\frac{1}{16}\dual\fb \la^{-1}(f\wedge \fb)\trch' +\lot
\end{split}
\eea

\item   The transformation formula for $\eta$ is given by 
\bea
\bsplit
\eta' &= \eta +\frac{1}{2}\la \nab_3'f  +\frac{1}{4}\fb\trch -\frac{1}{4}\dual\fb\atrch -\omb\, f +\err(\eta, \eta'),\\
\err(\eta, \eta') &= \frac{1}{2}(f\c\fb)\eta +\frac{1}{2}\fb\c\chih
+\frac{1}{2}f(\fb\c\ze)  -  (\fb\c f)\eta'+ \frac{1}{2}\fb (f\c\eta') +\lot
\end{split}
\eea

\item   The transformation formula for $\etab$ is given by 
\bea
\bsplit
\etab' &= \etab +\frac{1}{2}\nab_{\la^{-1}e_4'}\fb +\frac{1}{4}\trchb f - \frac{1}{4}\atrchb\dual f -\om\fb +\err(\etab, \etab'),\\
\err(\etab, \etab') &=  \frac{1}{2}f\c\chibh + \frac{1}{2}(f\c\etab)\fb-\frac 1 4  (f\c\ze)\fb  -\frac 1 4 |\fb|^2\la^{-2}\xi'+\lot
\end{split}
\eea

\item   The transformation formula for $\om$ is given by
\bea
\bsplit
\la^{-1}\om' &=  \om -\frac{1}{2}\la^{-1}e_4'(\log\la)+\frac{1}{2}f\c(\ze-\etab) +\err(\om, \om'),\\
\err(\om, \om') &=   -\frac{1}{4}|f|^2\omb - \frac{1}{8}\trchb |f|^2+\frac{1}{2}\la^{-2}\fb\c\xi' +\lot
\end{split}
\eea

\item   The transformation formula for $\omb$ is given by
\bea
\bsplit
\la\omb' &= \omb+\frac{1}{2}\la e_3'(\log\la)  -\frac{1}{2}\fb\c\ze -\frac{1}{2}\fb\c\eta +\err(\omb,\omb')\\
\err(\omb,\omb') &= f\c\fb\,\omb-\frac{1}{4} |\fb|^2\om  +\frac{1}{2}f\c\xib + \frac{1}{8}(f\c\fb)\trchb + \frac{1}{8}(\fb\wedge f)\atrchb \\
& -\frac{1}{8}|\fb|^2\trch  -\frac{1}{4}\la \fb\c\nab_3'f    +\frac{1}{2}  (\fb\c f)(\fb\c\eta')- \frac{1}{4}|\fb|^2 (f\c\eta')+\lot
\end{split}
\eea

\end{itemize}
where, for the transformation formulas of the Ricci coefficients above, $\lot$ denote expressions of the type
\beaa
         \lot&=O((f,\fb)^3)\Ga +O((f,\fb)^2) \Gac
\eeaa
containing no derivatives of $f$, $\fb$, $\Ga$ and $\Gac$. 

\begin{itemize}
\item The transformation formula for $\a, \aa $  are  given by
\bea
\bsplit
\la^{-2} \a'&=\a +\err(\a, \a')\\
\err(\a, \a')&=  \big(  f\hot \b  -\dual f \hot \dual  \b )+ \left( f\hot f-\frac 1 2  \dual f \hot   \dual f \right) \rho
+  \frac 3  2 \big(  f \hot  \dual  f\big) \rhod +\lot,
\end{split}
\eea
\bea
\bsplit
\la^2\aa'&=\aa +\err(\a, \a')\\
\err(\aa, \aa')&=  \big(  \fb \hot \bb  -\dual \fb \hot \dual  \bb )+ \big( \fb \hot \fb-\frac 1 2  \dual \fb \hot   \dual \fb \big) \rho
+  \frac 3  2 \big(  \fb \hot  \dual  \fb\big) \rhod +\lot
\end{split}
\eea

\item   The transformation formula for $\b , \bb $  are  given by
   \bea
  \bsplit
\la^{-1}   \b'&=\b +\frac 3 2\big(  f \rho+\dual  f  \rhod\big)+\err(\b, \b') \\
  \err(\b, \b')&= \frac 1 2 \a\c\fb+\lot,
  \end{split}
  \eea
  \bea
  \bsplit
  \la\bb'&=\bb +\frac 3 2\big(  \fb \rho+\dual  \fb  \rhod\big)+\err(\bb, \bb') \\
  \err(\bb, \bb')&= \frac 1 2  \aa\c f +\lot 
  \end{split}
  \eea
  \item The transformation formula for $\rho$ and $\rhod $  are  given by
  \bea
  \bsplit
 \rho' &= \rho +\err(\rho, \rho'),\\
\err(\rho, \rho') &= \fb\c\b - f\c\bb +\frac{3}{2}\rho(f\c\fb) -\frac{3}{2}\rhod (f\wedge\fb) +\lot
   \end{split}
  \eea
  \bea
  \bsplit
  \rhod' &= \rhod +\err(\rhod, \rhod'),\\
\err(\rhod, \rhod') &= -\fb\c\dual\b - f\c\dual\bb +\frac{3}{2}\rhod(f\c\fb) +\frac{3}{2}\rho (f\wedge\fb) +\lot
   \end{split}
  \eea  

\end{itemize}
where, for the transformation formulas of the Ricci coefficients above, $\lot$ denote expressions of the type
\beaa
         \lot&=O((f,\fb)^3)(\rho, \rhod) +O((f,\fb)^2)(\a,\b,\aa, \b)
\eeaa
containing no derivatives of $f$, $\fb$, $\a$, $\b$, $(\rho, \rhod)$, $\bb$, and $\aa$. 
\end{proposition} 

\begin{proof}
See  section 3.2 in \cite{KS-GCM1}.
\end{proof}

%%%%%%%%%%%%%%%%%%%%%%%%%%

\subsection{Transport equations for $(f, \fb, \la)$}
\lab{section:transport(f,fb,la)}

%%%%%%%%%%%%%%%%%%%%%%%%%%

\begin{corollary}\lab{cor:transportequationine4forchangeofframecoeff:simplecasefirst}
Under the assumption 
\beaa
\xi'=0, \qquad \om'=0, \qquad \etab'+\ze'=0, 
\eeaa
we  have the following transport equations for  $(f, \fb, \la)$
\beaa
\nab_{\la^{-1}e_4'}f+\frac{1}{2}(\trch f -\atrch\dual f)+2\om f &=& -2\xi - f\c\chih+E_1(f, \Ga),\\
\la^{-1}e_4'(\log\la) &=& 2\om+f\c(\ze-\etab)+E_2(f, \Ga),\\
\nab_{\la^{-1}e_4'}\fb+\frac{1}{2}(\trch\fb  +\atrch\dual\fb) &=& -2(\etab+\ze)   +2\nab'(\log\la)  +2\omb f \\ && +E_3({\nab'}^{\leq 1}f, \fb, \Ga, \la^{-1}\chi'),
\eeaa
where $E_1(f, \Ga)$, $E_2(f, \Ga)$ and $E_3({\nab'}^{\leq 1}f, \fb, \Ga, \la^{-1}\chi')$ are given by
\beaa
E_1(f, \Ga) &=& -(f\c\ze)f-\frac{1}{2}|f|^2\eta +\frac{1}{2}|f|^2\etab +O(f^3\Ga),\\
E_2(f, \Ga) &=& - \frac{1}{2}|f|^2\omb  -\frac{1}{4}\trchb|f|^2+O(f^3\Ga+f^2\chibh),
\eeaa
and
\beaa
E_3({\nab'}^{\leq 1}f, \fb, \Ga, \la^{-1}\chi') &=&  \frac{1}{4}(f\c\eta)\fb + \frac{1}{2}(f\c\ze)\fb    + \frac{1}{4}\dual\fb(f\wedge\eta) + \frac{1}{4}\dual\fb(f\wedge\ze)\\
&& +  \frac{1}{4}\fb\div'f + \frac{1}{4}\dual\fb \curl'f +\frac{1}{2}\la^{-1}\fb\c\chih' \\
&&+O\Big((\la^{-1}\trch', \la^{-1}\atrch')(f, \fb)^3+(f, \fb)^3\Ga\Big).
\eeaa
\end{corollary}

\begin{proof}
See section \ref{proofof:cor:transportequationine4forchangeofframecoeff:simplecasefirst}.
\end{proof}

\begin{corollary}\lab{cor:transportequationine4forchangeofframecoeffinformFFbandlamba}
Assume that we have
\beaa
\Xi'=0, \qquad \om'=0, \qquad \Hb'+Z'=0.
\eeaa
We introduce
\beaa
F:=f+i\dual f, \qquad \underline{F}:=\fb+i\dual \fb.
\eeaa
Then, we have
\beaa
\nab_{\la^{-1}e_4'}F+\frac{1}{2}\ov{\tr X} F+2\om F &=& -2\Xi -\chih\c F+E_1(f, \Ga),\\
\la^{-1}\nab_4'(\log\la) &=& 2\om+f\c(\ze-\etab)+E_2(f, \Ga),\\
\nab_{\la^{-1}e_4'}\underline{F}+\frac{1}{2}\tr X\underline{F} &=& -2(\Hb+Z)   +2\DD'(\log\la)  +2\omb F  +E_3({\nab'}^{\leq 1}f, \fb, \Ga, \la^{-1}\chi').
\eeaa

Moreover, introducing a complex valued scalar function $q$ satisfying $e_4(q)=1$, we have
\beaa
\nab_{\la^{-1}e_4'}(\ov{q}F) &=& -2\ov{q}\om F  -2\ov{q}\Xi +E_4(f, \Ga),\\
\nab_{\la^{-1}e_4'}\left[q\Big(\underline{F}-2q\DD'(\log\la)\Big)+\ov{q}e_3(r)F\right] &=&  -2q(\Hb+Z)  -2q^2\DD'\Big(2\om+f\c(\ze-\etab)\Big)\\
&& +e_3(r)\left(-2\ov{q}\om F  -2\ov{q}\Xi\right) +2\omb(q-\ov{q})F \\
&&+E_5({\nab'}^{\leq 1}f,\fb, {\nab'}^{\leq 1}\la, \D^{\leq 1}\Ga),
\eeaa
where 
\beaa
E_4(f, \Ga) &:=& -\frac{1}{2}\ov{q}\left(\ov{\tr X} -\frac{2}{\ov{q}}\right)F -\ov{q}\chih\c F+\ov{q}E_1(f, \Ga)+f\c\nab(\ov{q})F+\frac{1}{4}|f|^2e_3(\ov{q})F
\eeaa
and
\beaa
E_5({\nab'}^{\leq 1}f,\fb, {\nab'}^{\leq 1}\la, \D^{\leq 1}\Ga) &=&-\frac{q}{2}\left(\tr X -\frac{2}{q}\right)\underline{F}       +qE_3({\nab'}^{\leq 1}f, \fb, \Ga, \la^{-1}\chi')\\
&& -2q^2\DD'\left(E_2(f, \Ga)\right)  -2q[\nab_{\la^{-1}e_4'},q\DD']\log\la\\
&&+\left(f\c\nab(q)+\frac{1}{4}|f|^2e_3(q)\right)\Big(\underline{F}-2q\DD'(\log\la)\Big)\\
&&+e_3(r)E_4(f, \Ga).
\eeaa
\end{corollary}

\begin{proof}
See section \ref{proofof:cor:transportequationine4forchangeofframecoeffinformFFbandlamba}.
\end{proof}

\begin{remark}
In practice, we will integrate first the transport equations for $F$, then for $\la$. Finally, we will integrate the renormalized transport equation for $\underline{F}$ and we recover one less derivative than for $F$ and $\log(\la)$. Note that the renormalization of the transport equation for $\underline{F}$ is needed in order to avoid a potential log-loss due to the terms $\DD'(\log\la)$ and $\omb F$ on the RHS. 
\end{remark}

%%%%%%%%%%%%%%%%%%%%%%%%%

\section{Principal  geodesic  structures}
\lab{section:outgoingprinc.structure}

%%%%%%%%%%%%%%%%%%%%%%%%%

%%%%%%%%%%%%%%%%%%%%%%%%%%%%

\subsection{Principal outgoing geodesic structures}
\lab{subsection:Principaloutgoinggeodesicstructures}

%%%%%%%%%%%%%%%%%%%%%%%%%%%%

\begin{definition}[PG structure]\lab{def:PGstructure}
An outgoing PG structure consists  of a null pair $(e_3, e_4)$  and  the induced horizontal structure $\HH=\O(\MM)$, together with a scalar function $r$ such that

 \begin{enumerate}
\item $e_4$ is  a    null outgoing geodesic  vectorfield, i.e.  $\D_4 e_4=0$,  

\item $r$ is  an affine parameter,  i.e. $e_4(r)=1$,

\item the   gradient  of $r$,  given by  $N=\g^{\a\b} \pr_\b r \pr_\a$,  is perpendicular  to  $\HH$.

\end{enumerate}
\end{definition}

\begin{lemma}
\lab{Lemma:outgoing-geodesic-horizontal-structure}
Given a PG  structure  as above,    we have
\beaa
 \om=0, \qquad \xi=0, \qquad \etab+\ze=0.
\eeaa
\end{lemma}

\begin{proof}
Since $e_4$ is geodesic, we have $\xi=\om=0$. Also, in view of the Ricci formulas \eqref{eq:Ricciformula}, 
\beaa
\,[e_a, e_4] &=& \chi_{ab}e_b-\nab_4e_a-(\ze+\etab)_a e_4,
\eeaa
with $(e_a)_{a=1,2}$ an orthonormal  basis of $\HH$. Applying  the commutator  formula  to $r$,  and using  $e_4(r)=1$, $e_1(r)=e_2(r)=0$, and  $(\nab_4 e_a)(r)=0 $,  we infer that
\beaa
\etab+\ze &=& 0
\eeaa
as stated.
\end{proof}

In view of the above, the following relations hold for PG structure
\bea
\xi =0, \qquad\quad\,\, \om=0, \qquad \etab+\ze=0,\qquad e_4(r)=1,\qquad \nab(r)=0.
\eea

The following lemma shows how to initialize a PG structure on a hypersurface of $\MM$  transversal to $e_4$.

\begin{lemma}\lab{lemma:initializationofaPGstructureonacausalhypersurface} 
Let a  hypersurface  $\Si$, a scalar function $r$ and   null pair   $(e_4,  e_3)$  both  defined on $\Si$ and let $\HH$ the  corresponding horizontal space. Assume that $e_4$ is transversal to $\Si$, and impose the transversality condition
 \beaa
 e_4(r)=1\quad\textrm{ on }\,\,\Si.
 \eeaa
Under this transversality condition\footnote{Note that we need to take the transversality condition into account    since $\HH$ is not necessarily included in the tangent space of $\Si$.}, we assume that\footnote{That is $X(r)=0$ for any $X\in  \HH$.}   $\HH(r)=0$. 
Then, we can extend $r$  and the null frame $(e_4,  e_3)$  uniquely   to a  PG structure.
\end{lemma}

\begin{proof}
Since $e_4$ is transversal to $\Si$, we may extend it geodesically, i.e. $\D_{e_4} e_4=0$, to a neighborhood of $\Si$.  We then  extend    $r$ in a neighborhood of $\Si$   such that $ e_4(r)=1$.  Since $\HH(r)=0$  on $\Si$  we have $\HH$ orthogonal to $  N=\g^{\a\b} \pr_\b r \pr_\a$  on $\Si$.  Since $e_4$ and $N$ are well  defined  in a neighborhood of $\Si$, we  extend $\HH$ by  choosing it to be orthogonal to both.  We then choose $e_3$  the unique null vector  perpendicular to $\HH$  and such that $\g(e_3, e_4)=-2$.  Thus $r$ and the frame $(e_3, e_4)$  define a PG structure in a neighborhood of  $\Si$ coinciding with the given one on $\Si$.
\end{proof}

%%%%%%%%%%%%%%%%%%%%%%%%%%%%%%%%%%%%%%%%%%%%%

\subsection{Null structure and Bianchi identities for an outgoing PG  structure}

%%%%%%%%%%%%%%%%%%%%%%%%%%%%%%%%%%%%%%%%%%%%%

In an outgoing PG structure, Propositions \ref{prop-nullstr:complex}  and \ref{prop:bianchi:complex}  take the following form.
\begin{proposition}
\lab{prop-nullstrandBianchi:complex:outgoing}
\beaa
\nab_4\tr X +\frac{1}{2}(\tr X)^2 &=& -\frac{1}{2}\Xh\c\ov{\Xh},\\
\nab_4\Xh+\Re(\tr X)\Xh &=& -A,
\eeaa
\beaa
\nab_4\tr\Xb +\frac{1}{2}\tr X\tr\Xb &=& -\DD\c\ov{Z}+Z\c\ov{Z}+2\ov{P}-\frac{1}{2}\Xh\c\ov{\Xbh},\\
\nab_4\widehat{\Xb} +\frac{1}{2}\tr X\, \widehat{\Xb}  &=& -\frac{1}{2}\DD\hot Z  +\frac{1}{2}Z\hot Z -\frac{1}{2}\ov{\tr\Xb} \widehat{X},
\eeaa
\beaa
\nab_3\tr\Xb +\frac{1}{2}(\tr\Xb)^2+2\omb\,\tr\Xb &=& \DD\c\ov{\Xib}-\Xib\c\ov{Z}+\ov{\Xib}\c(H-2Z)-\frac{1}{2}\Xbh\c\ov{\Xbh},\\
\nab_3\Xbh+\Re(\tr\Xb) \Xbh+ 2\omb\,\Xbh&=& \frac{1}{2}\DD\hot \Xib+  \frac{1}{2} \Xib\hot(H-3Z)-\Ab,
\eeaa
\beaa
\nab_3\tr X +\frac{1}{2}\tr\Xb\tr X-2\omb\tr X &=& \DD\c\ov{H}+H\c\ov{H}+2P-\frac{1}{2}\Xbh\c\ov{\Xh},\\
\nab_3\widehat{X} +\frac{1}{2}\tr\Xb\, \widehat{X} -2\omb\widehat{X} &=& \frac{1}{2}\DD\hot H  +\frac{1}{2}H\hot H -\frac{1}{2}\ov{\tr X} \widehat{\Xb}.
\eeaa
Also,
\beaa
\nab_4Z +\tr X Z &=&  -\widehat{X}\c\ov{Z} -B,\\
\nab_4H +\frac{1}{2}\ov{\tr X}(H+Z) &=&   -\frac{1}{2}\Xh\c(\ov{H}+\ov{Z}) -B,\\
\eeaa
\beaa
\nab_3Z +\frac{1}{2}\tr\Xb(Z+H)-2\omb(Z-H) &=& -2\DD\omb -\frac{1}{2}\widehat{\Xb}\c(\ov{Z}+\ov{H})+\frac{1}{2}\tr X\Xib -\Bb+\frac{1}{2}\ov{\Xib}\c\Xh,\\
\nab_3Z +\nab_4\Xib &=&  -\frac{1}{2}\ov{\tr\Xb}(Z+H) -\frac{1}{2}\Xbh\c(\ov{Z}+\ov{H}) -\Bb,\\
\eeaa
and
\beaa
\nab_4\omb  -(2\eta+\ze)\c\ze &=&   \rho.
\eeaa
Also,
\beaa
\frac{1}{2}\ov{\DD}\c\Xh +\frac{1}{2}\Xh\c\ov{Z} &=& \frac{1}{2}\DD\ov{\tr X}+\frac{1}{2}\ov{\tr X}Z-i\Im(\tr X)H -B,\\
\frac{1}{2}\ov{\DD}\c\Xbh -\frac{1}{2}\Xbh\c\ov{Z} &=& \frac{1}{2}\DD\ov{\tr\Xb}-\frac{1}{2}\ov{\tr\Xb}Z-i\Im(\tr\Xb)(-Z+\Xib)+\Bb.
\eeaa
Also,
 \beaa
 \nab_3A -\frac{1}{2}\DD\hot B &=& -\frac{1}{2}\tr\Xb A+4\omb A +\frac{1}{2}(Z+4H)\hot B -3\ov{P}\Xh,\\
\nab_4B -\frac{1}{2}\ov{\DD}\c A &=& -2\ov{\tr X} B  +\frac{1}{2}A\c  \ov{Z},\\
\nab_3B-\DD\ov{P} &=& -\tr\Xb B+2\omb B+\ov{\Bb}\c \Xh+3\ov{P}H +\frac{1}{2}A\c\ov{\Xib},\\
\nab_4P -\frac{1}{2}\DD\c \ov{B} &=& -\frac{3}{2}\tr X P -\frac{1}{2}Z\c\ov{B}  -\frac{1}{4}\Xbh\c \ov{A}, \\
\nab_3P +\frac{1}{2}\ov{\DD}\c\Bb &=& -\frac{3}{2}\ov{\tr\Xb} P -\frac{1}{2}(\ov{2H-Z})\c\Bb +\Xib\c \ov{B} -\frac{1}{4}\ov{\Xh}\c\Ab, \\
\nab_4\Bb+\DD P &=& -\tr X\Bb+\ov{B}\c \Xbh+3P Z,\\
\nab_3\Bb +\frac{1}{2}\ov{\DD}\c\Ab &=& -2\ov{\tr\Xb}\,\Bb -2\omb\,\Bb -\frac{1}{2}\Ab\c (\ov{-2Z +H})-3P \,\Xib,\\
\nab_4\Ab +\frac{1}{2}\DD\hot\Bb &=& -\frac{1}{2}\ov{\tr X} \Ab +\frac{5}{2}Z\hot \Bb -3P\Xbh.
\eeaa
    \end{proposition} 
    
\begin{proof}
We are simply making use of the fact that, in an outgoing PG structure, we have $\om=\xi=0$, $\etab=-\ze$, and hence also $\Xi=0$, $\Hb=-Z$.
\end{proof}

%%%%%%%%%%%%%%%%%%%%%%%%%%%%%%%%%%%%%%%%%%%%%

\subsection{Coordinates  associated to an outgoing PG structure}
\lab{section:coordinates.outgoingPGstr}

%%%%%%%%%%%%%%%%%%%%%%%%%%%%%%%%%%%%%%%%%%%%%

\begin{definition}
Assume given an outgoing principal   geodesic structure $\{r, ( e_3, e_4), \HH\}$. 
In addition to $r$, we define scalar  functions $(u, \th, \vphi)$ such that
\bea
e_4(u)=e_4(\th)=e_4(\vphi)=0.
\eea
\end{definition}

\begin{proposition}
\lab{prop:e_4(xyz)}
The following equations hold true for the coordinates $(u,r,\th,\vphi)$ associated to an outgoing  PG structure
\beaa
e_4(e_3(r)) &=&  -2\omb,
\eeaa
\beaa
\nab_4\DD u +\frac{1}{2}\tr X\DD u &=& -\frac{1}{2}\Xh\c\ov{\DD}u,\\
e_4(e_3(u)) &=& -\Re\Big((Z+H)\c\ov{\DD} u\Big),
\eeaa
\beaa
\nab_4(\DD \cos \th) +\frac{1}{2}\tr X (\DD \cos \th) &=& -\frac{1}{2}\Xh\c\ov{\DD}(\cos \th),\\
e_4(e_3(\cos \th)) &=& -\Re\Big((Z+H)\c\ov{\DD} ( \cos \th) \Big).
\eeaa
\end{proposition}

\begin{proof}
For a scalar function $f$, we have the commutator formulas 
\beaa
\,[\nab_4, \DD]f &=&  -\frac{1}{2}\tr X\DD f -\frac{1}{2}\Xh\c\ov{\DD}f,\\
\,[\nab_4, \nab_3]f &=&  -\Re\Big((Z+H)\c\ov{\DD} f\Big) -2\omb \nab_4f.
\eeaa
Hence, for a scalar function $f$ such that $\nab_4f=1$ of $\nab_4f=0$, we have
\beaa
\nab_4\DD f +\frac{1}{2}\tr X\DD f &=& -\frac{1}{2}\Xh\c\ov{\DD}f,\\
e_4(e_3(f)) &=& -\Re\Big((Z+H)\c\ov{\DD} f\Big) -2\omb \nab_4f.
\eeaa
We then apply these transport equations respectively with the choice $f=r$, $f=u$ and $ f = cos \th$.
\end{proof}

%%%%%%%%%%%%%%%%%%%%%%%%%%%%%%%%%%%%%%%%%%%%%

\subsection{Integrable frame adapted to a PG structure}
\lab{section:adapted-integr-frame}

%%%%%%%%%%%%%%%%%%%%%%%%%%%%%%%%%%%%%%%%%%%%%

Given a PG structure\footnote{Here, $(e_1, e_2)$  denotes an arbitrary orthonormal basis of the horizontal space  $\HH$.}   $\{r,(e_3, e_4, e_1, e_2)\}$, see Definition  \ref{def:PGstructure},  and a coordinate $u$ transported by $e_4(u)=0$, we construct in the following lemma a  new frame $(e_3', e_4' , e_1', e_2')$ adapted to the topological spheres $S(u, r)$, i.e. with  $(e_1', e_2')$ 
  tangent to $S(u, r)$.   

\begin{lemma}
\lab{Lemma:Transformation-principal-to-integrable frames}
Consider  a PG structure   $\{r,(e_3, e_4, e_1, e_2)\}$, and a coordinate $u$ transported by $e_4(u)=0$, and 
assume that
\beaa
e_3(u)>0, \qquad (e_3(u))^2 +4|\nab u|^2e_3(r)>0.
\eeaa
Then, there exists a    frame transformation of type \eqref{General-frametransformation},  taking $(e_1, e_2, e_3, e_4)$ into a frame  $(e_3', e_4',  e_1', e_2') $, and  verifying the following conditions
\begin{enumerate}
\item  The horizontal vectors $(e_1', e_2')$ are tangent to $S(u, r) $.
\item We have $\g(e_3', e_4)=-2$.
\end{enumerate}
Moreover, the coefficients $(\la, f, \fb)$ of the  frame transformation are given by\footnote{Note that the formula for $\la$ implies $\la>0$. Indeed, using $e_3(u)>0$ and $(e_3(u))^2 +4|\nab u|^2e_3(r)>0$, and the above formulas for $\la$, $f$ and $\fb$, we have
\beaa
\la &=& 1-\frac{4|\nab u|^2e_3(r)}{(e_3(u)+ \sqrt{(e_3(u))^2 +4|\nab u|^2e_3(r)})\sqrt{(e_3(u))^2 +4|\nab u|^2e_3(r)}}+\frac{1}{16}|f|^2|\fb|^2\\
&\geq& \frac{(e_3(u))^2}{(e_3(u))^2 +4|\nab u|^2e_3(r)}>0.
\eeaa}  
\bea\lab{def:transition-functs:ffbla}
\bsplit
f  &= -\frac{4}{e_3(u)+ \sqrt{(e_3(u))^2 +4|\nab u|^2e_3(r)}}\nab u,\\
\fb &= \frac{2e_3(r)}{\sqrt{(e_3(u))^2 +4|\nab u|^2e_3(r)}}\nab u,\\
\la &=1+\frac{1}{2}f\c\fb  +\frac{1}{16}|f|^2|\fb|^2.
\end{split}
\eea
\end{lemma}

\begin{remark}
In Kerr, we have, see section \ref{section:PGstructuresKerr}, 
\beaa
|\nab u|^2=\frac{a^2(\sin\th)^2}{|q|^2}, \qquad e_3(r)=-\frac{\De}{|q|^2}, \qquad e_3(u)=\frac{2(r^2+a^2)}{|q|^2},
\eeaa
and hence
\beaa
e_3(u)>0, \qquad (e_3(u))^2 +4|\nab u|^2e_3(r)= \frac{4\Si^2}{|q|^4}>0
\eeaa
as well as
\beaa
f =  -\frac{2|q|^2}{r^2+a^2+ \Si} \nab u, \qquad \fb=-\frac{\De}{\Si}\nab u.
\eeaa
\end{remark}
\begin{proof}
Recall the frame transformation \eqref{General-frametransformation}     
\beaa
 \bsplit
  e_4'&=\la\left(e_4 + f^b  e_b +\frac 1 4 |f|^2  e_3\right),\\
  e_a'&= \left(\de_a^b +\frac{1}{2}\fb_af^b\right) e_b +\frac 1 2  \fb_a  e_4 +\left(\frac 1 2 f_a +\frac{1}{8}|f|^2\fb_a\right)   e_3,\\
 e_3'&=\la^{-1}\left( \left(1+\frac{1}{2}f\c\fb  +\frac{1}{16} |f|^2  |\fb|^2\right) e_3 + \left(\fb^b+\frac 1 4 |\fb|^2f^b\right) e_b  + \frac 1 4 |\fb|^2 e_4 \right).
 \end{split}
 \eeaa
In particular, the choice for $\la$ is equivalent to the normalization $\g(e_3', e_4)=-2$.  It thus suffices to look for $(f, \fb)$ such that $(e_1', e_2')$ is tangent to $S(u,r)$, which is equivalent to $e_a'(r)=e_a'(u)=0$. Now, recall that
\beaa
e_4(r)=1, \qquad e_4(u)=0, \qquad e_1(r)=e_2(r)=0.
\eeaa
Thus, we have $e_a'(r)=e_a'(u)=0$ if and only if
\beaa
0 &=& \fb_a +\left(f_a+\frac{1}{4}|f|^2\fb_a\right)e_3(r),\\
0 &=& \left(\de_a^b+\frac{1}{2}\fb_af^b\right) e_b(u)+\left(\frac{1}{2}f_a+\frac{1}{8}|f|^2\fb_a\right)e_3(u),
\eeaa
which is equivalent to 
\beaa
\fb &=&  -\frac{e_3(r)}{1+\frac{1}{4}|f|^2e_3(r)}f,\\
0 &=& \nab u+\frac{1}{2}(f\c\nab u)\fb+\left(\frac{1}{2}f+\frac{1}{8}|f|^2\fb\right)e_3(u).
\eeaa
Plugging the first equation in the second, we infer
\beaa
-(f\c\nab u)e_3(r)f+e_3(u)f+ 2\left(1+\frac{1}{4}|f|^2e_3(r)\right)\nab u &=& 0.
\eeaa
We look for $f$ under the form
\beaa
f=h\nab u
\eeaa
where $h$ is a scalar function to be determined. Plugging in the equation above, we obtain  
\beaa
\left(-h^2|\nab u|^2e_3(r)+he_3(u)+ 2\left(1+\frac{1}{4}h^2|\nab u|^2e_3(r)\right)\right)\nab u &=& 0.
\eeaa
Thus, $(e_1', e_2')$ is tangent to $S(u,r)$ if and only if $h$ satisfies 
\beaa
-\frac{1}{2}|\nab u|^2e_3(r)h^2+e_3(u)h+2=0.
\eeaa
The discriminant is given by
\beaa
Disc &=& (e_3(u))^2 +4|\nab u|^2e_3(r)>0,
\eeaa
where the strict positivity comes from the assumptions, and the roots are given by
\beaa
h_{\pm} &=& \frac{-e_3(u)\pm \sqrt{(e_3(u))^2 +4|\nab u|^2e_3(r)}}{-|\nab u|^2e_3(r)}.
\eeaa

We choose the root $h_+$, i.e.
\beaa
h &=& \frac{-e_3(u)+ \sqrt{(e_3(u))^2 +4|\nab u|^2e_3(r)}}{-|\nab u|^2e_3(r)}\\
&=& \frac{4|\nab u|^2e_3(r)}{-|\nab u|^2e_3(r)\big(e_3(u)+ \sqrt{(e_3(u))^2 +4|\nab u|^2e_3(r)}\big)}\\
&=& -\frac{4}{e_3(u)+ \sqrt{(e_3(u))^2 +4|\nab u|^2e_3(r)}},
\eeaa
where we recall that assumption $e_3(u)>0$.  We infer
\beaa
f &=& h\nab u= -\frac{4}{e_3(u)+ \sqrt{(e_3(u))^2 +4|\nab u|^2e_3(r)}}\nab u
\eeaa
and 
\beaa
\fb &=&  -\frac{e_3(r)}{1+\frac{1}{4}|f|^2e_3(r)}f= -\frac{he_3(r)}{1+\frac{1}{4}h^2|\nab u|^2e_3(r)}\nab u= -\frac{2he_3(r)}{4+e_3(u)h}\nab u\\
&=& \frac{2e_3(r)}{\sqrt{(e_3(u))^2 +4|\nab u|^2e_3(r)}}\nab u
\eeaa
as desired.
\end{proof}

\begin{definition}\lab{def:associatedintegrableframetoPGstructure}
 We refer  to the   frame introduced in Lemma \ref{Lemma:Transformation-principal-to-integrable frames} 
 as the associated integrable frame  to the PG structure.
\end{definition}

\begin{corollary}
Let $(e_1, e_2, e_3, e_4)$ a  PG frame and let $(e_3', e_4',  e_1', e_2') $ the associated integrable frame. Consider the frame transformation of type \eqref{General-frametransformation},  taking  $(e_3', e_4',  e_1', e_2') $ into  $(e_1, e_2, e_3, e_4)$, and let $(\la', f', \fb')$ the corresponding coefficients. Then, we have $\la'=1$.  
\end{corollary}

\begin{proof}
The frame transformation of type \eqref{General-frametransformation},  taking  $(e_3', e_4',  e_1', e_2') $ into  $(e_1, e_2, e_3, e_4)$, is the inverse transformation of the one of Lemma \ref{Lemma:Transformation-principal-to-integrable frames}. Thus, the proof follows immediately from the properties of $(\la, f, \fb)$ in Lemma \ref{Lemma:Transformation-principal-to-integrable frames} and the identities \eqref{relations:laffb-to-primes} relating $(\la, f, \fb)$ and $(\la', f', \fb')$.
\end{proof}

%%%%%%%%%%%%%%%%%%%%%%%%%%

\section{Canonical  outgoing  PG structure in Kerr}
\lab{section:PGstructuresKerr}

%%%%%%%%%%%%%%%%%%%%%%%%%%

%%%%%%%%%%%%%%%%%%%%%%%

\subsection{Boyer-Lindquist coordinates}

%%%%%%%%%%%%%%%%%%%%%%%

The Kerr metric in standard Boyer-Lindquist coordinates  $(t, r, \th, \phi) $ is given by 
$$\g=-\frac{|q|^2\Delta}{\Sigma^2}(dt)^2+\frac{\Sigma^2(\sin\theta)^2}{|q|^2}\left(d\phi-\frac{2amr}{\Sigma^2}dt\right)^2+\frac{|q|^2}{\Delta}(dr)^2+|q|^2(d\theta)^2,$$
where 
$$
\left\{\ba{lll}
q &=& r+ i a \cos\th,\\
\Delta &=& r^2-2mr+a^2,\\
\Sigma^2 &=& (r^2+a^2)|q|^2+2mra^2(\sin\theta)^2=(r^2+a^2)^2-a^2(\sin\theta)^2\Delta.
\ea\right.
$$
We also   note   that  
\bea
\T=\pr_t , \qquad \Z=\pr_\phi,
\eea
 are both Killing and $\T$  is only time-like in the     complement of the ergoregion,
 i.e. in $|q|^2> 2 mr$. 
 The    domain of outer communication    is given by, 
\begin{equation*}
\RR=\{(\theta,r,t,\vphi)\in(0,\pi)\times(r_+,\infty)\times\mathbb{R}\times\mathbb{S}^1\},
\end{equation*}
where     $r_+:=m+\sqrt{m^2-a^2}$,  the larger root of $\De$, corresponds to the      event horizon.

\begin{lemma}
The following  principal null  directions    are canonical in Kerr.
\begin{enumerate}
\item The          null pair  $(e_3, e_4 )$        for which $e_4$ is geodesic   is given by
\bea
\lab{eq:Out.PGdirections-Kerr}
 e_4=\frac{r^2+a^2}{\Delta}\pr_t+\pr_r+\frac{a}{\Delta}\pr_\phi,\quad e_3=\frac{r^2+a^2}{|q|^2}\pr_t-\frac{\Delta}{|q|^2}\pr_r+\frac{a}{|q|^2}\pr_\phi.
\eea

\item  The  null pair  $(e_3, e_4 )$ for which $e_3$ is geodesic is given by
\bea
\lab{eq:Inc.PGdirections-Kerr}
e_3=\frac{r^2+a^2}{\Delta}\pr_t-\pr_r+\frac{a}{\Delta}\pr_\phi,\qquad  e_4 =\frac{r^2+a^2}{|q|^2}\pr_t+\frac{\Delta}{|q|^2}\pr_r+\frac{a}{|q|^2}\pr_\phi.
\eea

\item In both cases,  a principal null frame can be obtained by adding the vectorfields 
\beaa
 e_1=\frac{1}{|q|}\pr_\th,\quad e_2=\frac{a\sin\th}{|q|}\pr_t+\frac{1}{|q|\sin\th}\pr_\phi.
\eeaa
\end{enumerate}
\end{lemma}

\begin{proof}
Straightforward verification.
\end{proof}

\begin{definition}
We refer to  the null frame obtained by adding $(e_1, e_2)$ to  $(e_3, e_4)$ in  \eqref{eq:Out.PGdirections-Kerr} as the canonical outgoing principal  null frame.  We refer to  the null frame obtained by adding $(e_1, e_2)$ to  $(e_3, e_4)$ in  \eqref{eq:Inc.PGdirections-Kerr} as the 
canonical ingoing  principal  null frame.
\end{definition}

\begin{remark}
The canonical ingoing  principal  null frame, i.e. the null frame obtained by adding $(e_1, e_2)$ to  $(e_3, e_4)$ in  \eqref{eq:Inc.PGdirections-Kerr}, is regular towards the future for all $r>0$.
\end{remark}

\begin{remark}
In the remaining of section \ref{section:PGstructuresKerr}, we    will only considers the outgoing PG structure of Kerr.  
\end{remark}

%%%%%%%%%%%%%%%%%%%%%%%%%%%%%%%%%%%%%%%%%%%%%

\subsection{Canonical complex 1-form  $\Jk$} 
\lab{section:formsJk}

%%%%%%%%%%%%%%%%%%%%%%%%%%%%%%%%%%%%%%%%%%%%%

\begin{definition}
\lab{def:JkandJ}
We define the  following  complex  1-form  $\Jk$ in Kerr as follows
\bea
\Jk:=j+i\dual j,
\eea
where the real 1-form $j$ is defined by
\bea
j_1=0, \qquad  j_2=\frac{\sin\th}{|q|}.
\eea
Hence
\beaa
\Jk_1=\frac{i\sin\th}{|q|}, \qquad \Jk_2=\frac{\sin\th}{|q|}.
\eeaa
\end{definition}

\begin{remark}
The relevance of the complex 1-form $\Jk$ is due to its link to the complex 1-forms $H$, $\Hb$ and $Z$, see  \eqref{eq:RicciinKerr-out}.
\end{remark}

\begin{lemma}\lab{lemma:propertiesofhorizontal1formJkinKerr}
 The  following identities hold  true.
 \begin{itemize}
 \item We have
 \bea
 \dual \Jk =- i \Jk, \qquad  \Jk\c\ov{\Jk}=\frac{2(\sin\th)^2}{|q|^2}.
 \eea
 
 \item We   have
\bea
\lab{transport-Jk-Kerr}
\nab_4\Jk+\frac 1 2 \tr X \Jk=0, \qquad \nab_3\Jk+\frac 1 2 \tr \Xb  \Jk=0.
\eea
 
 \item  The complex 1-form $\Jk$ verifies
\bea
\ov{\DD}\c\Jk = \frac{4i(r^2+a^2)\cos\th}{|q|^4},\qquad \DD\hot\Jk = 0.
\eea

\item  We have
\bea 
\DD (q)= -a\Jk, \qquad \DD(\ov{q})= a \Jk.
\eea
\end{itemize}
\end{lemma}

\begin{proof}
Straightforward computations.   Note that \eqref{transport-Jk-Kerr}   holds true in both the ingoing and outgoing frame. In the  outgoing frame it becomes
\beaa
\nab_4\Jk &=& -\frac{1}{q}\Jk,     \qquad
\nab_3\Jk = \frac{\De q}{|q|^4}\Jk.
\eeaa
To check the third identity it helps to first check  the following identities for   $j=\Re(\Jk)$
\beaa
\div j=0, \quad \curl j=\frac{2(r^2+a^2)\cos\th}{|q|^4}, \quad \nab\hot j=0,\quad \nab_a j_b =\frac{(r^2+a^2)\cos\th}{|q|^4}\in_{ab}.
\eeaa
The last identity is a straightforward verification.
\end{proof}

\begin{lemma}
The Killing vectorfield given on BL coordinates by   $\T=\pr_t$ can be expressed  in terms of the  outgoing PG  frame  by 
\bea
\T &=& \frac 1 2\left( e_3+\frac{\De}{|q|^2} e_4 -2a\Re(\Jk)^be_b\right).
\eea
\end{lemma}

\begin{proof}
Using
\beaa
\frac{\sin\th}{|q|} e_2 &=& j^b e_b=\Re(\Jk)^be_b
\eeaa
yields the formula.
\end{proof}

%%%%%%%%%%%%%%%%%%%%%%%%%%%%%%%%%%%%%%%%%%%%%

\subsection{Canonical outgoing  PG structure   in Kerr}
\lab{subs:Kerrvalues1}

%%%%%%%%%%%%%%%%%%%%%%%%%%%%%%%%%%%%%%%%%%%%%

 Consider the canonical  principal outgoing  null  frame  in BL coordinates
\bea
\lab{eq:Out.PGframe-Kerr}
\bsplit
e_4&=\frac{r^2+a^2}{\Delta}\pr_t+\pr_r+\frac{a}{\Delta}\pr_\phi,\quad e_3=\frac{r^2+a^2}{|q|^2}\pr_t-\frac{\Delta}{|q|^2}\pr_r+\frac{a}{|q|^2}\pr_\phi, \\
 e_1&=\frac{1}{|q|}\pr_\th,\quad e_2=\frac{a\sin\th}{|q|}\pr_t+\frac{1}{|q|\sin\th}\pr_\phi.
 \end{split}
\eea

\begin{lemma}
\lab{Lemma:OutPG.Kerr}
The functions     $r$  together with the principal frame  \eqref{eq:Out.PGframe-Kerr}   defines an outgoing PG  structure   in Kerr. Moreover  
\begin{enumerate}
 \item
 The complex curvature   null  components w.r.t. the frame  are given by
\beaa
A=B=\Bb=\Ab=0,\quad P=-\frac{2m}{q^3}.
\eeaa
\item 
The  vanishing  complex Ricci coefficients in Kerr are
\beaa
 \Xh=\Xbh=\Xi=\Xib=\om=0.
\eeaa
\item The non-vanishing   complex Ricci  coefficients are\footnote{Note that all complex 1-forms are  expressed with respect to $\Jk$.}
\bea
\lab{eq:RicciinKerr-out}
\bsplit
\tr X&=\frac{2}{q}, \qquad\quad\,\,\,\, \tr\Xb=-\frac{2\Delta q}{|q|^4},\\
\Hb&= -\frac{a\ov{q}}{|q|^2}\Jk, \qquad\, H=\frac{aq}{|q|^2}\Jk,\qquad Z=\frac{a\ov{q}}{|q|^2}\Jk,\\
\omb&=\frac{1}{2}\pr_r\left(\frac{\De}{|q|^2}\right),
\end{split}
\eea
with   the 1-form $\Jk$   of Definition \ref{def:JkandJ}.

 \item We have 
 \bea
 \lab{eq:Kerr.La_abc}
 \bsplit
 (\La_1)_{12}&=0, \qquad \qquad \qquad  (\La_2)_{12} =-\frac{r^2+a^2}{|q|^3} \cot \th, \,\, \\
  (\La_3)_{12}&= - \frac{a\De\cos \th}{|q|^4}, \qquad  (\La_4)_{12} =- \frac{a \cos \th}{|q|^2}.
  \end{split}
 \eea
  \end{enumerate}
\end{lemma}

\begin{proof}
Straightforward computations.  See also the section on Kerr in \cite{GKS1}.
 \end{proof}

%%%%%%%%%%%%%%%%%%%%%%%%%%%%%%%%%

\subsection{Outgoing  PG  coordinates in Kerr}
\lab{sec:outgoingPGcoordinatesinKerr:chap2} 

%%%%%%%%%%%%%%%%%%%%%%%%%%%%%%%%%

Note that relative to the BL coordinates we have 
\beaa
e_4(r)=1,\qquad  e_4(\th)=0, \qquad e_4(\phi)=\frac{a}{\De}.
\eeaa
To  derive the coordinate system  associated to the  outgoing PG structure, we need to introduce in addition to $(r,\th)$
\beaa
u:= t- f(r), \quad f'(r)=\frac{r^2+ a^2}{ \De}, \qquad \vphi:= \phi -h(r), \quad h'(r)=\frac{a}{\Delta},
\eeaa
which satisfies
\beaa
e_4(u)=0, \qquad e_4(\vphi) =0.
\eeaa

\begin{lemma}\lab{lemma:urthetavphicoordinatesinKerrchap2}
The  coordinate system  given by $(u, r, \th, \vphi)$, with
\beaa
u=t-f(r), \quad f'(r)=\frac{r^2+a^2}{\Delta}, \qquad \vphi = \phi -h(r), \quad h'(r)=\frac{a}{\Delta},
\eeaa
 is   the canonical  coordinate system associated to the  canonical outgoing PG structure  of Kerr, called the outgoing Eddington-Finkelstein  (EF) coordinate system of Kerr.
 Moreover,
\begin{enumerate}
\item  The action of  the  outgoing PG  frame on the coordinates $(u, r, \th, \vphi) $ is given by
 \bea
\lab{eq:OPG-valuesKerr}
\begin{array}{llll}
e_4(r)=1, & e_4(u)=0, & e_4(\th)=0, & e_4(\vphi)=0,\\[2mm]
\displaystyle e_3(r)=-\frac{\Delta}{|q|^2}, & \displaystyle e_3(u)=\frac{2(r^2+a^2)}{|q|^2},  & e_3(\th)=0,  & \displaystyle e_3(\vphi)=\frac{2a}{|q|^2}, \\[2mm] 
e_1(r)=0, & e_1(u)=0,  & \displaystyle e_1(\th)=\frac{1}{|q|},  & e_1(\vphi)=0, \\[2mm] 
e_2(r)=0, & \displaystyle e_2(u)=\frac{a\sin\th}{|q|},  & e_2(\th)=0,  & \displaystyle e_2(\vphi)=\frac{1}{|q|\sin\th}.
\end{array}
\eea

\item In the outgoing EF coordinates,  the metric takes the form 
 \beaa
\g &=& -\left(1-\frac{2mr}{|q|^2}\right)(du)^2 -2dr du +2a(\sin\th)^2dr d\vphi\\
&&-\frac{4mra(\sin\th)^2}{|q|^2}du d\vphi+|q|^2(d\th)^2+\frac{\Si^2(\sin\th)^2}{|q|^2}(d\vphi)^2.
\eeaa

\item  In the outgoing EF coordinates, the determinant on the metric is given by 
\bea\lab{eq:determinantofKerrmetricinoutgoingEFcoord}
\det(\g) &=& -|q|^4(\sin\th)^2.
\eea

\item  In the outgoing EF coordinates, the inverse metric coefficients $\g^{\a\b}$ are given by
\beaa
&&\g^{rr} =\frac{\Delta}{|q|^2},\qquad \g^{ru}=-\frac{r^2+a^2}{|q|^2},\qquad \g^{r\th} = 0,\qquad \g^{r\vphi}= -\frac{a}{|q|^2},\\
&& \g^{uu} = \frac{a^2(\sin\th)^2}{|q|^2},\qquad \g^{u\th} = 0,\qquad \g^{u\vphi} = \frac{a}{|q|^2},\\
&& \g^{\th\th} = \frac{1}{|q|^2},\qquad \g^{\th\vphi} = 0,\qquad \g^{\vphi\vphi} = \frac{1}{|q|^2(\sin\th)^2}.
\eeaa
\end{enumerate}
\end{lemma} 

\begin{proof}
Straightforward verification.
\end{proof}

%%%%%%%%%%%%%%%%%%%%%%%%%%%%

\subsection{Canonical basis of $\ell=1$ modes in Kerr}
\lab{section:cabell=1basisKerr}

%%%%%%%%%%%%%%%%%%%%%%%%%%%%

\begin{definition}
\lab{definition:canonical.ell=1basisKerr}
We define the canonical basis of $\ell=1$ modes of Kerr to be the triplet of scalar functions\footnote{Note that  they are regular  everywhere including the axis of symmetry.} $\Jp$, $p\in\{0, +, -\}$, 
 \bea
 J^{(0)} :=\cos\th, \qquad J^{(-)} :=\sin\th\sin\vphi, \qquad  J^{(+)} :=\sin\th\cos\vphi.
 \eea
\end{definition}

\begin{remark}
Note that we have
\bea
e_4(\Jp)=0.
\eea
\end{remark}

\begin{lemma}
\lab{Lemma:lapJ^p-Kerr}
The following identities hold true.
\bea
\lab{eq:lapJ^p-Kerr1}
\bsplit
\lap J^{(0)}&=-\frac{2(r^2+  a^2)}{|q|^4}J^{(0)},\\
\lap  J^{(\pm)}&=- \frac{2 r^2 }{|q|^4} J^{(\pm)}.
\end{split}
\eea
Also
\bea
\lab{eq:lapJ^p-Kerr2}
\DD\hot \DD \Jp=0.
\eea
\end{lemma}

\begin{remark}
To avoid any confusion concerning the notation $\De$  in sections \ref{section:cabell=1basisKerr} and  \ref{section:complex1formJkplusminusinKerr}, notice that:
\begin{itemize}
\item In the statement of Lemma \ref{Lemma:lapJ^p-Kerr} and in its proof, as well as in the proof of Lemma \ref{Lemma:propertiesJk_pm}, $\De$ denotes the horizontal Laplacian. 

\item In the statement of Lemma \ref{Le:nab_4nab_3Psi} and in its proof, as well as in the statement of Lemma \ref{Lemma:propertiesJk_pm}, $\De$ denotes the scalar function given by $r^2-2mr+a^2$. 
\end{itemize}
\end{remark}

\begin{proof}
Given a scalar function $h$ we  have, in view of $(\La_1)_{12}=0$ and the formula for $(\La_2)_{12}$ in \eqref{eq:Kerr.La_abc}, 
\beaa
\lap h &=& e_1(e_1 h)+ e_2(e_2 h) +\frac{r^2+a^2}{|q|^3} \cot \th  e_1(h)
\eeaa
and \eqref{eq:lapJ^p-Kerr1} follows by choosing $h=\Jp$ and recalling    the relations, see \eqref{eq:OPG-valuesKerr}, 
\beaa
e_1(\th)=\frac{1}{|q|},  \quad e_1(\vphi)=0, \quad   e_2(\th)=0,  \quad   e_2(\vphi)=\frac{1}{|q|\sin\th}.
\eeaa

Next, we focus on \eqref{eq:lapJ^p-Kerr2}. First, note that for a scalar function $h$, we have
\beaa
\DD\hot\DD h=2\nab\hot\nab h +2i\dual(\nab\hot\nab h)
\eeaa
so that it suffices to prove $\nab\hot\nab\Jp=0$. 
Also, given a scalar function $h$ we  have, in view of $(\La_1)_{12}=0$ and the formula for $(\La_2)_{12}$ in \eqref{eq:Kerr.La_abc}, 
\beaa
(\nab\hot\nab h)_{11} &=& e_1(e_1 h) - e_2(e_2 h) - \frac{r^2+a^2}{|q|^3} \cot \th  e_1(h),\\
(\nab\hot\nab h)_{12} &=& e_1(e_2 h) + e_2(e_1 h) - \frac{r^2+a^2}{|q|^3} \cot \th  e_2(h)
\eeaa
and we obtain $\nab\hot\nab\Jp=0$ by choosing $h=\Jp$ so that \eqref{eq:lapJ^p-Kerr2} follows.  
\end{proof}

\begin{lemma}
The following  identities hold true
\bea
 e_3(J^{(0)}) =0, \qquad e_3 (J^{(+)} )+ \frac{2a}{|q|^2}J^{(-)}=0, \qquad e_3 (J^{(- )}) - \frac{2a}{|q|^2}J^{(+)}=0.
\eea
\end{lemma}

\begin{proof}
Straightforward calculation.
\end{proof}

\begin{lemma}
\lab{Lemma:DDJ^{(+)}}
We have
\bea
\DD J^{(0)} =i \Jk.
\eea
\end{lemma}

\begin{proof}
We  check that
 \beaa
  \nab J^{(0)}=- \Im(\Jk), \qquad  \dual \nab J^{(0)} =\Re(\Jk).
  \eeaa
 Indeed, recalling the definition of $\Jk$, 
\beaa
e_1 (J^{(0)}) &=& -\frac{\sin\th }{|q|}=-\Im(\Jk_1), \qquad  e_2 (J^{(0)}) =0=- \Im(\Jk_2),\\
\dual \nab J^{(0)} &=&-\Im( \dual \Jk) =- \Im(- i  \Re \Jk) =\Re(\Jk),
\eeaa
Therefore
\beaa
\DD J^{(0)} &=&  \nab J^{(0)}+i \dual \nab J^{(0)}=-\Im(\Jk) +i \Re(\Jk)=i \Jk
\eeaa
as stated.
\end{proof}
 
It remains to  also calculate  $\DD J^{(\pm)}$. To do that, we introduce in the next section the  complex 1-forms $\Jk_\pm$, as counterparts to $\Jk$.

%%%%%%%%%%%%%%%%%%%%%%%%%%%

\subsection{Canonical complex 1-forms $\Jk_\pm$}
\lab{section:complex1formJkplusminusinKerr}

%%%%%%%%%%%%%%%%%%%%%%%%%%%

\begin{definition}
\lab{definition:Jk_pm} 
We define the complex 1-forms $\Jk_{\pm}$ as follows  
\bea
\Jk_{\pm}=j_{\pm}+ i \dual j_{\pm},
\eea
 where the real 1-forms $j_{\pm}$ are given by
 \bea
 \bsplit
 (j_+)_1&=\frac{1}{|q|} \cos\th\cos\vphi, \qquad (j_+)_2=-\frac{1}{|q|} \sin\vphi,\\
 (j_-)_1& =\frac{1}{|q|} \cos\th\sin\vphi, \qquad (j_-)_2=\frac{1}{|q|}  \cos\vphi.
 \end{split}
 \eea
\end{definition}

We can now obtain the following analog of Lemma \ref{Lemma:DDJ^{(+)}}.
\begin{lemma}
\lab{Lemma:DD J^{(pm)} =Jk_pm}
The complex 1-forms $\Jk_\pm$  are anti-selfadjoint, i.e.  $\dual\Jk_{\pm}=- i\Jk_{\pm}$, and verify
\bea
 \DD J^{(+)} =\Jk_+, \qquad   \DD J^{(-)} =\Jk_-.
\eea
\end{lemma}

\begin{proof}
 To   prove  $\DD J^{(\pm)} =\Jk_\pm$ we      check the following
\beaa
\nab J^{(\pm)}=\Re( \Jk_\pm),  \qquad \dual  \nab J^{(\pm)}=\Im( \Jk_\pm).
\eeaa
Indeed, in view of Definition \ref{definition:Jk_pm}  and Definition \ref{definition:canonical.ell=1basisKerr},
\beaa
e_1 (J^{(+)})& =& e_1(\sin\th \cos \vphi)=\frac{1}{|q|} \cos \th \cos\vphi=(j_+)_1=\big(\Re(\Jk_{+})\big)_1,\\
e_2 (J^{(+)})& =& e_2(\sin\th \cos \vphi)=-\frac{1}{|q|}\sin \vphi=(j_+)_2=\big(\Re(\Jk_{+})\big)_2,
\eeaa
\beaa
e_1 (J^{(-)})& =& e_1(\sin\th \sin \vphi)=\frac{1}{|q|} \cos \th \sin\vphi= (j_-)_1=\big(\Re(\Jk_{-})\big)_1,\\
e_2 (J^{(-)})& =& e_2(\sin\th \sin \vphi)=\frac{1}{|q|}\cos \vphi= (j_-)_2=\big(\Re(\Jk_{-})\big)_2.
\eeaa
Also
\beaa
\dual \nab J^{(+)}&=&\Re( \dual \Jk_+ )=\Re (-i \Jk_+)=\Im (\Jk_+),\\
\dual \nab J^{(-)}&=&\Re( \dual \Jk_- )=\Re (-i \Jk_-)=\Im (\Jk_-).
\eeaa
Thus
\beaa
\DD J^{(+)} &=&  \nab J^{(+)}+i \dual \nab J^{(+)}=  \Re(\Jk_+)+i \Im(\Jk_+)=\Jk_+,\\
\DD J^{(-)} &=&  \nab J^{(-)}+i \dual \nab J^{(-)}=  \Re(\Jk_-)+i \Im(\Jk_-)=\Jk_-,
\eeaa
as stated.
\end{proof}  

\begin{lemma}
\lab{Le:nab_4nab_3Psi}
Let  $\Psi$ be a  complex 1-form in  Kerr of the form
\beaa
\Psi=\frac{1}{|q|} \big(\psi+ i \dual \psi\big)
\eeaa
where the real 1-form $\psi$ is such that $e_4(\psi_1)=e_4(\psi_2) = 0$. Then
\beaa
\nab_4 \Psi+\frac{1}{q} \Psi=0.
\eeaa
Also, we have
\beaa
(\nab_3 \Psi)_a  =\frac{\De q}{|q|^4}\Psi _a+\frac{1}{|q|}  e_3 \big(\psi_a+i   \dual \psi_a\big), \quad a=1,2. 
\eeaa
\end{lemma}

\begin{proof}
Since $(\La_4)_{12}= -\frac{a\cos\th}{|q|^2}$ and $(\La_3)_{12}= - \frac{a\De\cos \th}{|q|^4}$ in view of  \eqref{eq:Kerr.La_abc}, we easily check, using also $e_4(\psi_1)=e_4(\psi_2) = 0$, 
\beaa
\nab_4 \psi=-\frac{a\cos\th}{|q|^2} \dual \psi, \qquad \nab_3\psi_a=-\frac{a\De\cos \th}{|q|^4}\dual \psi_a+e_3(\psi_a),
\eeaa
and thus
\beaa
\nab_4 (\psi+i\dual \psi) &=& i \frac{a\cos\th}{|q|^2} (\psi+i \dual\psi), \\
\nab_3(\psi+i\dual \psi)_a &=& i \frac{a\De\cos\th}{|q|^4} (\psi+i \dual\psi)_a+e_3(\psi_a+i\dual \psi_a).
\eeaa
The conclusion then easily follows using $e_4(r)=1$, $e_4(\th)=e_3(\th)=0$ and $e_3(r)=-\frac{\De}{|q|^2}$ to compute $e_4(|q|^2)$ and $e_3(|q|^2)$.
\end{proof}

\begin{lemma}
\lab{Lemma:propertiesJk_pm} 
The complex 1-forms $\Jk_\pm$  verify 
\bea
\lab{eq:propertiesJk_pm1}
\nab_4 \Jk_\pm =-   \frac{1}{q} \Jk_{\pm}, \qquad 
\nab_3 \Jk_\pm =\frac{\De q}{|q|^4} \, \Jk_{\pm} \,\mp  \,\frac{2a}{|q|^2}  \Jk_{\mp},
\eea
\bea
\lab{eq:propertiesJk_pm2}
\DD\hot \Jk_{\pm}=0, \qquad \ov{\DD}\c \Jk_\pm= - \frac{4r^2 }{|q|^4} J^{(\pm)}\mp \frac{4ia^2\cos\th}{|q|^4}J^{(\mp)},
\eea
\bea
\lab{eq:propertiesJk_pm3}
 \bsplit
 \Jk_+\c\ov{\Jk_+}&=\frac{2(\cos\th)^2(\cos\vphi)^2+2(\sin\vphi)^2}{|q|^2},\\
  \Jk_-\c\ov{\Jk_-}&=\frac{2(\cos\th)^2(\sin\vphi)^2+2(\cos\vphi)^2}{|q|^2},
  \end{split}
 \eea
and
\bea
\lab{eq:propertiesJk_pm4}
\Re(\Jk_+)\c\Re(\Jk)=-\frac{1}{|q|^2}J^{(-)}, \qquad \Re(\Jk_-)\c\Re(\Jk)=\frac{1}{|q|^2}J^{(+)}.
\eea
\end{lemma}

\begin{proof}
\eqref{eq:propertiesJk_pm1} is an immediate consequence of Lemma \ref{Le:nab_4nab_3Psi} using that $e_4(\th)=e_4(\vphi)=e_3(\th)=0$ and $e_3(\vphi)=\frac{2a}{|q|^2}$.

To check \eqref{eq:propertiesJk_pm2}, recall from  Lemma \ref{Lemma:DD J^{(pm)} =Jk_pm} that $\Jk_\pm=\DD J^{(\pm)}$ which yields
\beaa
\DD\hot\Jk_{\pm} &=& \DD\hot\DD J^{(\pm)}=0,\\
\ov{\DD}\c\Jk_{\pm} &=& \ov{\DD}\c\DD J^{(\pm)}=2\De J^{(\pm)}+2i(\nab_1\nab_2-\nab_2\nab_1)J^{(\pm)}\\
&=& - \frac{4r^2 }{|q|^4} J^{(\pm)}+2i(\nab_1\nab_2-\nab_2\nab_1)J^{(\pm)}
\eeaa
where we used \eqref{eq:lapJ^p-Kerr1} and \eqref{eq:lapJ^p-Kerr2}. This yields the first identity in \eqref{eq:propertiesJk_pm2}, while the second one follows from the above together with the following identity for a scalar function $h$ only depending on $(\th, \vphi)$
\beaa
(\nab_1\nab_2-\nab_2\nab_1)h &=& e_1(e_2h)-e_2(e_1h)+\frac{r^2+a^2}{|q|^3}\cot\th e_2(h)= \frac{2a^2\cos\th}{|q|^4}\pr_\vphi(h)
\eeaa
which is then applied to $h=J^{(\pm)}$ noticing that $\pr_\vphi(J^{(\pm)})=\mp J^{(\mp)}$.

The remaining estimates \eqref{eq:propertiesJk_pm3} and \eqref{eq:propertiesJk_pm4} are straightforward verifications from the definitions.
\end{proof}

%%%%%%%%%%%%%%%%%%%%%%%%%%%%%%%%%%%%%%%%%%%%%%%%%

\subsection{Additional coordinates system}
\lab{sec:outgoingPGcoordinatesinKerr:additionalcoordsystems:chap2} 

%%%%%%%%%%%%%%%%%%%%%%%%%%%%%%%%%%%%%%%%%%%%%%%%%%

Recall from \eqref{eq:determinantofKerrmetricinoutgoingEFcoord} that the determinant of the Kerr metric in the $(u,r,\th, \vphi)$ coordinates is given by $\det(\g) = -|q|^4(\sin\th)^2$ which illustrates that the $(r,u, \th, \vphi)$ coordinates system is singular at the axis, i.e. at $\th=0$ and $\th=\pi$. In this section, we introduce a coordinates system which is regular at the axis
\bea
(u,r,x^1, x^2), \qquad x^1:=J^{(+)}, \qquad x^2:=J^{(-)}.
\eea

We start by showing that the coordinates $(x^1, x^2)$ defined above form a coordinates system of the spheres $S(u,r)$ which is regular at the axis. 

\begin{lemma}\lab{lemma:regularx1x2coordinatessystemonSurx10x20}
The metric $g$ induced by $\g$ on $S(u,r)$ has the following form in Kerr in the $(x^1, x^2)$ coordinates system
\beaa
g &=& |q|^2\Bigg[\left(\frac{1-(x^2)^2}{1-|x|^2}+ \frac{a^2(x^2)^2}{|q|^2}\left(1+\frac{2mr}{|q|^2}\right)\right)(dx^1)^2\\
&&+\left(\frac{2x^1x^2}{1-|x|^2}- \frac{2a^2x^1x^2}{|q|^2}\left(1+\frac{2mr}{|q|^2}\right)\right)dx^1 dx^2\\
&&+\left(\frac{1-(x^1)^2}{1-|x|^2}+ \frac{a^2(x^1)^2}{|q|^2}\left(1+\frac{2mr}{|q|^2}\right)\right)(dx^2)^2\Bigg].
\eeaa
\end{lemma}

\begin{proof}
Straightforward verification.
\end{proof}
 
\begin{remark}
The coordinates $(x^1, x^2)$ verify the following properties
\begin{enumerate}
\item We have $|(x^1, x^2)|\leq 1$, with $|(x^1, x^2)|=1$ corresponding to the equator, i.e. $\th=\frac{\pi}{2}$, and $|(x^1, x^2)|=0$ corresponding to the poles.

\item $(x^1, x^2)$ are regular coordinates away from the equator, i.e. for $|(x^1, x^2)|<1$, or equivalently for $\th\neq\frac{\pi}{2}$

\item Since the $(\th, \vphi)$ coordinates system is regular for $\th\neq 0, \pi$, and the $(x^1, x^2)$ coordinates system is regular for $\th\neq\frac{\pi}{2}$, the spheres $S(u,r)$ are covered by these two coordinates systems (or rather 3 since one needs one $(x^1, x^2)$ coordinates system per hemisphere). 
\end{enumerate}
\end{remark}

\begin{remark}
Since we have $e_4(\th)=e_4(\vphi)=0$, $(x^1, x^2)$ satisfy 
\beaa
e_4(x^1)=e_4(x^2)=0.
\eeaa
\end{remark}

Using the above coordinates on $S(u,r)$, we consider the coordinates $(u,r,x^1,x^2)$ on Kerr in the following lemma. 
\begin{lemma}\lab{lemma:urthetaJplusJminuscoordinatesinKerrchap2}
Let $(x^1, x^2)$ coordinates on $S(u,r)$ given by $(x^1, x^2)=(J^{(+)}, J^{(-)})$. Then, 
\begin{enumerate}
\item The Kerr metric is given in the $(u,r,x^1,x^2)$ coordinates system by   
\beaa
\g &=& -\left(1-\frac{2mr}{|q|^2}\right)(du)^2 -2dr du +2adr\big(-x^2 dx^1+x^1dx^2\big)\\
&&  -\frac{4mra}{|q|^2}du \big(-x^2 dx^1+x^1dx^2\big)+g,
\eeaa
where  $g$ denotes the induced metric by $\g$ on $S(u,r)$, which is given by Lemma \ref{lemma:regularx1x2coordinatessystemonSurx10x20} in the $(x^1, x^2)$ coordinates system.

\item  In the $(u,r,x^1,x^2)$ coordinates, the determinant on the metric is given by 
\bea
\det(\g) &=& \frac{|q|^4}{(\cos\th)^2}
\eea

\item  In the $(u,r,x^1,x^2)$ coordinates, the inverse metric coefficients $\g^{\a\b}$ are given by
\beaa
&&\g^{rr} =\frac{\Delta}{|q|^2},\qquad \g^{ru}=-\frac{r^2+a^2}{|q|^2},\qquad \g^{rx^1}=\frac{ax^2}{|q|^2},\qquad \g^{rx^2}=-\frac{ax^1}{|q|^2},\\
&& \g^{uu} = \frac{a^2(\sin\th)^2}{|q|^2},\qquad \g^{ux^1} = -\frac{ax^2}{|q|^2},\qquad \g^{ux^2} = \frac{ax^1}{|q|^2},\\
&& \g^{x^1x^1} = \frac{(\cos\th\cos\vphi)^2+(\sin\vphi)^2}{|q|^2},\qquad \g^{x^1x^2} = \frac{((\cos\th)^2-1)\sin\vphi\cos\vphi}{|q|^2},\\ 
&& \g^{x^2x^2} = \frac{(\cos\th\sin\vphi)^2+(\cos\vphi)^2}{|q|^2}.
\eeaa
\end{enumerate}
\end{lemma}

\begin{proof}
Straightforward verification.
\end{proof}

\begin{remark}
In view of the above lemma, away from the horizon, the outer domain of communication of Kerr can be covered by the following three coordinates systems
\begin{itemize}
\item $(u,r,\th, \vphi)$ away from the poles $\th=0$ and $\th=\pi$, e.g. in $\frac{\pi}{4}<\th<\frac{3\pi}{4}$,

\item two copies of $(u,r,J^{(+)}, J^{(-)})$ away from the equator $\th=\frac{\pi}{2}$, e.g. one copy in $0\leq \th<\frac{\pi}{3}$ and another  copy  in $\frac{2\pi}{3}<\th\leq\pi$. 
\end{itemize}
\end{remark}

%%%%%%%%%%%%%%%%%%%%%%%%%%%%%%%%%

\subsection{Asymptotic for the outgoing PG structure in Kerr}

%%%%%%%%%%%%%%%%%%%%%%%%%%%%%%%%%

We denote by $\jS$ the following real 1-form 
\beaa
\jS_1=0, \qquad \jS_2=\frac{\sin\th}{r}. 
\eeaa

 \begin{corollary} 
 \lab{corollary:Kerr.PG.asymptotics}
In Kerr, relative to its  canonical    outgoing    PG frame (see section \ref{subs:Kerrvalues1}),  the following identities hold asymptotically for  large $r$:
\begin{enumerate}
\item We have,
\beaa
j &=&\big(1+O(r^{-2}) \big) \jS.
\eeaa

\item We have
\beaa
\nab u= a j =a\big(1+O(r^{-2}) \big) \jS, 
\eeaa
and
\beaa
 e_3(r)=-\Up  +O\big(r^{-2} \big), \qquad e_3(u)=2+O\big(r^{-2}\big),
\eeaa
where we have used the notation 
\beaa
\Up &:=& 1-\frac{2m}{r}.
\eeaa

 \item The real components of  the non-vanishing  Ricci   scalars are given by        
\beaa
\trch &=&  \frac 2 r +O(r^{-3}), \\
  \trchb&=& -\frac{2\Up}{r}+O (r^{-3} ),\\
\atrch&=& \frac{2a \cos \th }{r^2} +O\big(r^{-4}\big),\\
  \atrchb&=&  \frac{2a\Up  \cos \th}{r^2}+O\big(r^{-4}\big).
\eeaa
Also
\beaa
 \omb&=& \frac{m}{r^2} +O(r^{-3}).
 \eeaa
 
 \item  The Ricci coefficient 1-forms are given by
 \beaa
 \etab=-\ze&=& -\frac{a}{r^2}  \Big(\jS+\frac{a\cos\th}{r}\dual\jS\Big)\big(1+O(r^{-2} ) \big) ,
 \\
 \eta&=& \frac{a}{r^2} \Big(\jS-\frac{a\cos\th}{r}\dual\jS\Big)\big(1+O(r^{-2} ) \big), 
 \eeaa
and, relative to the frame $e_1, e_2$,
\beaa
\sin \th (\La_2)_{12}&=&-\frac{\cos \th}{r}+ O(r^{-3}), \quad (\La_3)_{21} =\frac{a\Up  \cos\th}{r^2}+O(r^{-4}), \quad \\
(\La_4)_{21} &=& \frac{a \cos \th}{r^2}+O(r^{-4} ).
\eeaa 
 \item The real non-vanishing curvature components are given by
\beaa
\rho&=& -\frac{2m}{r^3} +O(r^{-5} ),\\
\rhod&=& \frac{6am \cos \th}{r^4}+ O(r^{-5} ).
\eeaa
\end{enumerate}
\end{corollary}

\begin{proof}
Straightforward  verification.  
\end{proof}

%%%%%%%%%%%%%%%%%%%%%%%%%%%%%%%%%%%

\subsection{Asymptotic of the  associated integrable  frame  in Kerr} 

%%%%%%%%%%%%%%%%%%%%%%%%%%%%%%%%%%%

We consider now the associated integrable frame  $(e_3', e_4', e_1', e_2')$ as defined by Lemma \ref{Lemma:Transformation-principal-to-integrable frames}.  
\begin{lemma}
\lab{Lemma:Transformation-principal-to-integrable-Kerr}
Consider the change of frame transformation in Kerr from the  PG frame  $(e_3, e_4, e_1, e_2)$ to
 the associated integrable frame $(e_3', e_4', e_1', e_2')$ as defined by Lemma \ref{Lemma:Transformation-principal-to-integrable frames}. The following   statements hold true for large $r$.
 \begin{enumerate}
 \item  We have $\la=1+O(r^{-2})$, $f_1=\fb_1=0$, and
 \bea
f_2 =-\frac{a\sin\th}{r}+ O\big( \sin\th r^{-3} \big), \qquad \fb_2=- \frac{a\sin\th \Up}{r}+ O\big( \sin\th r^{-3} \big).
\eea

\item  The  connection coefficients  $( \La'_\mu)_{12}=\g\big( \D_{e_\mu'} e_2',  e_1'\big)  $, $\mu=1,2,3,4$,  verify
\bea
(\La_1')_{12}=O( \sin \th r^{-4}), \qquad  \sin \th  (\La_2')_{21}= \frac{\cos \th}{r} +O(r^{-3} ),
\eea
and 
\bea
(\La_3')_{21}=  O\big(r^{-4}\big), \qquad   (\La_4')_{21}=     O\big(r^{-4}\big).
\eea

\item  The curvature components are given by
\bea\lab{eq:curvaturecomponentsintegrableframeinKerrasymptoticchap2}
\bsplit
\a'&= O(\sin^2\th r^{-5} ),\\
\b'_1&= O(\sin \th r^{-5}), \qquad \quad  \,\,\b_2'= \frac{3 a m \sin\th}{r^4} +O(\sin \th r^{-5}),\\
\rho'&=-\frac{2m}{r^3} + O(r^{-5}), \qquad 
\rhod'=\frac{6 am\cos \th }{r^4}  + O(\sin^2 \th r^{-5}),\\
\bb'&= O\big(\sin \th  r^{-4} \big),\\
\aa'&= O\big(\sin^2 \th  r^{-5} \big).
\end{split}
\eea 

\item We have,
\bea
\bsplit
\div' \b'&=O(r^{-6}),\\
\curl' \b' &=  \frac{6 am \cos \th}{r^5}+O(r^{-6}).
\end{split}
\eea

\item We have,
\bea
\bsplit
 \div' f&=O(\sin^2\th r^{-5} ), \qquad   \nab'\hot f=O(\sin \th r^{-5} ),\\
  \div' \fb&=O(\sin^2 \th r^{-5} ), \qquad   \nab'\hot \fb=O(\sin \th r^{-5} ).
  \end{split}
 \eea 
 Also,
 \bea
 \bsplit
 \nab'_3(  f_1, \fb_1)  +\frac 1 2 \trchb( f_1, \fb_1)&=  O \big(\sin\th  r^{-5} \big), \\
 \nab'_3 ( f_2, \fb_2) +\frac 1 2 \trchb (f_2, \fb_2) &= O \big(\sin\th  r^{-3} \big),
\\
  \nab'_4(  f_1, \fb_1)  +\frac 1 2 \trch( f_1, \fb_1)&=  O \big(\sin\th  r^{-5} \big), \\
 \nab'_4 ( f_2, \fb_2) +\frac 1 2 \trch (f_2, \fb_2) &= O \big(\sin\th  r^{-3} \big).
\end{split}
\eea

\item The connection coefficients  behave as follows
\bea
\bsplit
\trch'& = \frac 2 r  +O(r^{-3} ), \qquad \qquad 
 \chih'=   O(r^{-3} ), \\
\trchb'&= -\frac{2 \Up }{r} + O(r^{-3} ), \qquad \,\,\chibh'=     O(r^{-3}), \\
 \ze'&=O(\sin\th  r^{-3} ), \qquad  \quad \quad\eta'=O(\sin\th  r^{-3} ).
\end{split}
\eea
Also,
\bea
\bsplit
\xi'&=O\big(\sin \th r^{-3} \big), \qquad \,\, \xib'=O\big(\sin \th r^{-3} \big),\\
\om'&=O\big(r^{-3}\big), \qquad\qquad  \, \omb'=\frac{m}{r^2} +O\big(r^{-3}\big).
\end{split}
\eea
 \end{enumerate}
\end{lemma} 

\begin{proof}
 Straightforward verification  in view of  Lemma \ref{Lemma:Transformation-principal-to-integrable frames}, Corollary \ref{corollary:Kerr.PG.asymptotics} and the transformation formulas of Proposition   \ref{Proposition:transformationRicci}.  
\end{proof}

\begin{remark}
Recall from Lemma \ref{Lemma:OutPG.Kerr} that $\a=\aa=\b=\bb=0$ in the PG structure of Kerr. Notice from  \eqref{eq:curvaturecomponentsintegrableframeinKerrasymptoticchap2} that this does not hold for the associated integrable frame, denoted by prime.
\end{remark}

We end this section  with the following proposition.
\begin{proposition}
\lab{prop:modesforwidecheck{trch, trchb}-Kerr}
Consider the linearized  quantities
\beaa
\widecheck{\trch}':=\trch' -\frac{2}{r'}, \qquad \widecheck{\trchb}':=\trchb' +\frac{2\left(1-\frac{2m}{r'}\right)}{r'}
\eeaa
where $r'$ denotes  the area radius of the spheres $S(u, r)$ in Kerr.
We have, for large $r$,
\beaa
\int_{S(u,r)}\widecheck{\trch}'J^{(p)} &=& O\left(\frac{1}{r^2}\right),\quad p=0,+,-,\\
\int_{S(u,r)}\widecheck{\trchb}'J^{(p)} &=& O\left(\frac{1}{r^2}\right),\quad p=0,+,-.
\eeaa
where $J^{(p)}$, $p=0,+,-$,  denote the standard $\ell=1$ spherical harmonics, i.e.
\beaa
J^{(-)}=\sin\th \sin \vphi, \qquad J^{(0)}=\cos \th, \qquad  J^{(+)}=\sin\th \cos \vphi.
\eeaa
\end{proposition}

\begin{proof}
See section \ref{appendix-proofofprop:modesforwidecheck{trch, trchb}-Kerr}.
\end{proof}

%%%%%%%%%%%%%%%%%%%%%%%%%%%%%%%%%%%%%%%%%%%%%

\section{Initialization of PG  structures on  a hypersurface}
\lab{subsection:InitializationSi*}

%%%%%%%%%%%%%%%%%%%%%%%%%%%%%%%%%%%%%%%%%%%%%

%%%%%%%%%%%%%%%%%%%%

\subsection{Framed hypersurfaces}

%%%%%%%%%%%%%%%%%%%%

\begin{definition}\lab{def:framedhypersurface}
 A framed hypersurface consists of   a set  $\big(\Si, r, (\HH,  e_3, e_4)\big)$ where
 \begin{enumerate}
 \item $\Si$ is  smooth a hypersurface in $\MM$,
 
 \item  $( e_3, e_4)$ is a null  pair on $\Si$ such that $e_4$ is transversal to $\Si$, and $\HH$, the  horizontal space   perpendicular on $e_3, e_4$,    is  tangent to $\Si$,
 
 \item  the function $r:\Si\to \RRR$ is a regular function on $\Si$   such that
$   \HH(r) =0$.

\end{enumerate}
 \end{definition}

For a given framed hypersurface $\big(\Si, r, (\HH, e_3, e_4)\big)$, we denote by  $\nu$ the  vectorfield tangent to $\Si$, normal to the $r$-foliation,  and normalized by the condition\footnote{This is always possible since $e_4$ is transversal to $\Si$.} $\g(\nu, e_4)=-2$. Thus,
  \bea
  \lab{eq:definition-nuSi}
\nu&=&e_3+ b_\Si e_4
\eea
with $b_\Si$ a  given  scalar function  on $\Si$.

%%%%%%%%%%%%%%%%%%%%%%%%%%%%%%%%%

\subsection{Initialization of PG structures}  

%%%%%%%%%%%%%%%%%%%%%%%%%%%%%%%%%

   In what follows,  we consider  natural  initial  data  structures  on hypersurfaces  $\Si$  of $\MM$ 
   which  generate  PG structures.

\begin{definition}[PG-data set]
\lab{definition:PGdataset}
The boundary data of a PG structure   (PG-data set)    consists of  
\begin{enumerate}
\item  a framed hypersurface   $\big(\Si, r, (\HH, e_3, e_4)\big)$   as in Defintion \ref{def:framedhypersurface},
 
 \item  a fixed  1-form $f$ on the spheres $S$ of the $r$-foliation of $\Si$ verifying  the condition
 \bea\lab{eq:|f|-sqrt|y| <2}
 b_\Si |f|^2<4\quad \text{on}\quad \Si,
 \eea
 where $b_\Si$ appears in \eqref{eq:definition-nuSi}.
 \end{enumerate}
\end{definition}

 \begin{proposition}
 \lab{Prop:Initialization of PG structures}
 Given a     PG data  set $\big(\Si,  r,  (\HH, e_3, e_4),  f\big) $ as in Definition \ref{definition:PGdataset},  there exists a unique  PG   structure  $\big(r',( \HH', e'_3, e'_4) \big)$   defined in a neighborhood of  $\Si$   such that the following hold true
 \begin{enumerate}
  \item   The function $r'$ is prescribed on $\Si$ by  $r'=r$.
\item    Along $\Si$,  the  restriction  of  the spacetime null frame  $(\HH', e'_3, e' _4)$  and  the   given null frame
 $( \HH,  e_3, e_4 )$ on $\Si$ are related by the transformation formulas,  where $(e_1, e_2)$ is a fixed orthonormal  basis of  $\HH$,
  \bea
 \lab{eq:framechange-integr-PG}
\bsplit
   e_4'&=e_4 + f^b  e_b +\frac 1 4 |f|^2  e_3,\\
  e_a'&= \left(\de_a^b +\frac{1}{2}\fb_af^b\right) e_b +\frac 1 2  \fb_a  e_4 +\left(\frac 1 2 f_a +\frac{1}{8}|f|^2\fb_a\right)   e_3,\\
 e_3'&= \left(1+\frac{1}{2}f\c\fb  +\frac{1}{16} |f|^2  |\fb|^2\right) e_3 + \left(\fb^b+\frac 1 4 |\fb|^2f^b\right) e_b  + \frac 1 4 |\fb|^2 e_4,
\end{split}
\eea
where $\fb$ is chosen such that\footnote{Note that $\fb$ is well defined thanks to the condition \eqref{eq:|f|-sqrt|y| <2}.}
 \bea
\lab{relations:laffb-to-primes-Si^*}
\bsplit
\fb &= -\frac{(\nu(r)-b_\Si)}{1-\frac{1}{4}b_\Si |f|^2}f.
\end{split}
\eea

 \end{enumerate} 
 \end{proposition} 

\begin{proof} 
First, note that in view of \eqref{eq:definition-nuSi}, we have on $\Si$
 \beaa
 e_4+ f^b  e_b+\frac 1 4 |f|^2  e_3 &=& \left(1-\frac{1}{4}b_\Si |f|^2\right)e_4+ f^b  e_b+\frac{1}{4}|f|^2\nu.
 \eeaa
Since $\nu$ and $e_2$ are tangent to $\Si$, and $e_4$ is transversal to $\Si$,  under the condition \eqref{eq:|f|-sqrt|y| <2}, the vectorfield $e_4+ f^b  e_b+\frac 1 4 |f|^2  e_3$ is transversal to $\Si$. Hence, in view of \eqref{eq:framechange-integr-PG}, $e_4'$ is  transversal to $\Si$. 

$e_4'$ being transversal to $\Si$, we introduce the following transversality condition on $\Si$ for $r'$
 \bea\lab{eq:transversality-rprimeanduprime}
 e_4'(r')=1\quad \text{on}\quad \Si.
 \eea
The transversality condition \eqref{eq:transversality-rprimeanduprime} allows us to specify all derivatives of $r'$ on $\Si$.  In particular, we would like to compute $e_1'(r)$  and $e_2'(r)$ on $\Si$. In view of \eqref{eq:framechange-integr-PG}, we have on $\Si$
\beaa
 e_a'(r') &=& \left(\de_a^b +\frac{1}{2}\fb_af^b\right) e_b(r') +\frac 1 2  \fb_a  e_4(r') +\left(\frac 1 2 f_a +\frac{1}{8}|f|^2\fb_a\right)   e_3(r').
\eeaa
 Since $(e_1, e_2)$ are tangent to $\Si$, since $r'=r$ on $\Si$, and since $e_1(r)=e_2(r)=0$ on $\Si$, we have on $\Si$
 \beaa
 e_a(r') &=& e_a(r)=0
 \eeaa
 and hence
\beaa
 e_a'(r') &=& \frac 1 2  \fb_a  e_4(r') +\left(\frac 1 2 f_a +\frac{1}{8}|f|^2\fb_a\right)   e_3(r')\\
 &=& \frac 1 2\left(e_4(r') +\frac{1}{4}|f|^2   e_3(r')\right)\fb_a+\frac 1 2e_3(r') f_a.
\eeaa  
Together with the definition \eqref{relations:laffb-to-primes-Si^*} of $\fb$, we obtain on $\Si$
\bea\lab{eq:intermediraycomputatonofe2primeofrprime}
 e_a'(r')  &=& \frac 1 2\left(   e_3(r')-\frac{\nu(r)-b_\Si}{1-\frac{1}{4}b_\Si |f|^2}  \left(e_4(r') +\frac{1}{4}|f|^2   e_3(r')\right)\right) f_a.
\eea

In view of \eqref{eq:intermediraycomputatonofe2primeofrprime}, we need  to compute $e_4(r')$ and $e_3(r')$.  In view of \eqref{eq:framechange-integr-PG} and the transversality condition \eqref{eq:transversality-rprimeanduprime}, we have on $\Si$
\beaa
1 &=& e_4'(r')=e_4(r') + f^b  e_b(r') +\frac 1 4 |f|^2  e_3(r')\\
&=& e_4(r') +\frac 1 4 |f|^2  e_3(r').
\eeaa 
 Also, since $\nu$ is tangent to $\Si$ and $r'=r$ on $\Si$, we have
 \beaa
 \nu(r) &=& \nu(r')=e_3(r')+b_\Si e_4(r').
 \eeaa
We infer from both identities, using also \eqref{eq:|f|-sqrt|y| <2},  
\beaa
e_4(r')=\frac{1-\frac{1}{4}\nu(r) |f|^2}{1-\frac{1}{4}b_\Si |f|^2}, \qquad e_3(r') = \frac{\nu(r)-  b_\Si}{1-\frac{1}{4}b_\Si |f|^2}, \quad \text{on}\quad \Si.
\eeaa
Together with \eqref{eq:intermediraycomputatonofe2primeofrprime}, we infer on $\Si$
\beaa
 e_a'(r')  &=& \frac 1 2\left(   e_3(r')-\frac{\nu(r)-b_\Si}{1-\frac{1}{4}b_\Si |f|^2}  \left(e_4(r') +\frac{1}{4}|f|^2   e_3(r')\right)\right) f_a\\
 &=&  \frac 1 2\left(   e_3(r')-\frac{\nu(r)-b_\Si}{1-\frac{1}{4}b_\Si |f|^2}\right)f_a\\
&=& 0.
\eeaa
Hence, we have finally obtained 
\bea\lab{eq:computatonofe2primeofrprimee1primeofrprimee1primeofuprime}
e_1'(r')=0, \qquad e_2'(r')=0,\quad \text{on}\quad \Si.
\eea 
 
Finally, on the hypersurface $\Si$, we are given a scalar function $r'$ and a null frame $( e'_1, e'_2, e'_3, e'_4)$, such that $e_4'$ is transversal to $\Si$ and $r'$ satisfies the  transversality condition \eqref{eq:transversality-rprimeanduprime} on $\Si$. Also, under this transversality condition, \eqref{eq:computatonofe2primeofrprimee1primeofrprimee1primeofuprime} holds. We are thus in position to apply Lemma \ref{lemma:initializationofaPGstructureonacausalhypersurface}, according to which there exists a unique  PG   structure  $\big(r',( e'_1, e'_2, e'_3, e'_4) \big)$   defined in a neighborhood of  $\Si$. This concludes the proof of the proposition.
\end{proof}

%%%%%%%%%%%%%%%%%%%%%%%%%%%%%%%%

\subsection{GCM hypersurfaces} 
  \lab{section:canonical-coord-BPG}

%%%%%%%%%%%%%%%%%%%%%%%%%%%%%%%%%

In this section, we introduce the concept of GCM hypersurfaces and GCM-PG data sets which are at the core of the construction of our bootstrap regions in Chapter 3, see section \ref{sec:admissibleGMCPGdatasetonSigmastar}.  The construction starts with framed hypersurfaces   $\big(\Si_*, r,  (\HH, e_3, e_4) \big)$ which terminate in a  boundary $S_*$    on which  the given function $r$ is constant, i.e. $S_*$ is a leaf of the $r$-foliation of $\Si_*$. 
     
The definition of GCM hypersurfaces below includes in particular conditions on the Ricci coefficients $\eta$ and $\xib$, see \eqref{eq:Si^*-GCM1}. Now, a priori, the only well defined Ricci coefficients on $\Si_*$ are $\trch$, $\trchb$,  $\chih$, $\chibh$ and $\ze$. To make sense of all Ricci coefficients on $\Si_*$, including in particular $\eta$ and $\xib$, we need to choose transversality conditions on $\Si_*$. We choose them to be compatible with an outgoing geodesic foliation initialized on $\Si_*$, i.e., we choose the following transversality conditions on $\Si_*$
\bea\lab{eq:transversalityconditionsonSigma:firsttime}
\xi=0, \qquad \om=0, \qquad \etab=-\ze, \qquad \textrm{on}\quad\Si_*.
\eea   
   
The GCM hypersurfaces below also involve a basis of $\ell=1$ modes $J^{(p)}$, $p=0,+,-$. This basis is defined as follows
 \begin{enumerate}
 \item There exist  coordinates $(\th, \vphi)$ on $S_*$   such that the   induced metric $g$ on $S_*$  takes the form (see \cite{KS-GCM2} for results on effective uniformization and canonical $\ell=1$ modes. A short review of the main results  on this topic  can also be found in section \ref{subsection:effective-unifrmization}. ).
 \bea\lab{eq:formofthemetriconSstarusinguniformization:0}
 g= e^{2\phi} r^2\Big( (d\th)^2+ \sin^2 \th (d\vphi)^2\Big).
 \eea
 
 \item The  functions 
 \bea
 J^{(0)} :=\cos\th, \qquad J^{(-)} :=\sin\th\sin\vphi, \qquad  J^{(+)} :=\sin\th\cos\vphi,
 \eea
 verify  the balanced  conditions 
 \bea
 \int_{S_*}  J^{(p)} =0, \qquad p=0,+,-.
 \eea
 
  \item We choose $(e_1, e_2)$ on  $S_*$   as follows
  \bea \lab{eq:canonical-e1ande2onSstar:0}
  e_1:=\frac{1}{r e^\phi} \pr_\th, \qquad e_2:=\frac{1}{r\sin \th e^\phi} \pr_\vphi.
  \eea
  
  \item The above $(\th, \vphi)$ coordinates on $S_*$ are extended to  $\Si_*$ by setting
  \bea
  \lab{eq:canonical-corrdsSi:0}
  \nu(\th)=\nu(\vphi)=0,
  \eea
  where $\nu$          is the   normal to  the  $r$-foliation on $\Si_*$   defined in \eqref{eq:definition-nuSi}. 
  
  \item We also extend   the $\Jp$  functions to $\Si_*$ by setting
  \bea\lab{eq:canonical-ell=1modesonSi:0}
  \nu(\Jp)=0, \quad p=0,+,-.
  \eea
\end{enumerate}

   \begin{definition}[GCM hypersurface]
 \lab{def:CanonicalGCM-hypersurface}
 Consider a   framed hypersurface  with    end sphere $S_*$,  transversality conditions \eqref{eq:transversalityconditionsonSigma:firsttime},   coordinates $(\th, \vphi)$  
   and functions $J^{(0)}$, $J^{(+)}$ and $J^{(-)}$ defined as in \eqref{eq:formofthemetriconSstarusinguniformization:0}--\eqref{eq:canonical-ell=1modesonSi:0}.  We call  such a framed hypersurface   a GCM  hypersurface if in addition the following  GCM conditions are verified.
 \begin{enumerate}
 \item We have on any sphere of the  foliation induced by $r$,
 \bea
\begin{split}
\lab{eq:Si^*-GCM1}
& \trch=\frac{2}{r},\quad\\
 \trchb&=-\frac{2\Up}{r}+\underline{C}_0+\sum_{p=0, +,-}\underline{C}_pJ^{(p)},\quad\\
  \mu&=\frac{2m}{r^3}+M_0+\sum_{p=0, +,-}M_p J^{(p)},\\
&\int_S J^{(p)}\div\eta=0, \qquad \int_S J^{(p)}\div\xib=0,\qquad p=0,+,-, \\\
& b_{\Si_*}\big|_{SP}=-1-\frac{2m}{r}, 
\end{split}
\eea
where  $\underline{C}_0$, $\underline{C}_p$, $M_0$, $M_p$ are scalar functions on $\Si_*$ constant  along the leaves of the foliation, and  $SP$ denotes the south poles of the spheres on $\Si_*$, i.e. $\th=\pi$.  

\item  In addition, we have on  the  end  sphere $S_*$
\bea
\lab{eq:Si^*-GCM2}
\widecheck{\trchb}=0,\qquad \int_{S_*} \Jp \, \div\b =0,\qquad  p=0,+,-,
\eea
  as well as 
\bea\lab{eq:Si^*-GCM3}
\int_{S_*}J^{(+)}\curl\b=0, \qquad \int_{S_*}J^{(-)}\curl\b=0.
\eea
 \end{enumerate}
 \end{definition}

\begin{definition}
\lab{define:am-onSi}
We define the parameters $(m, a) $  of a GCM  hypersurface to be constant on $\Si_*$ with  $m$ being  the Hawking mass of $S_*$, i.e.
\bea
\frac{2m}{r}=1+\frac{1}{16\pi}\int_{S_{*}}\trch \trchb,
\eea
and with $a$ given by
\bea
a &:=& \frac{r^3}{8\pi m}\int_{S_*} J^{(0)}\curl\b.
\eea
\end{definition}

\begin{definition}
\lab{definition:GCM-datasets}
A GCM-PG data set  is a PG data set  $\big(\Si_*, r, (\HH, e_3, e_4), f \big)$, as in Definition \ref{definition:PGdataset}  such that the framed hypersurface   $\big(\Si_*, r, (\HH, e_3, e_4)\big)$ is a  GCM hypersurface, see Definition \ref{def:CanonicalGCM-hypersurface}.
\end{definition}

%%%%%%%%%%%%%%%%%%%%%%%%%%%%%%%%%%%%%%%%%%%%%

\section{Linearization of outgoing PG structures}

%%%%%%%%%%%%%%%%%%%%%%%%%%%%%%%%%%%%%%%%%%%%%

In this section, we assume given  an outgoing  principal geodesic structure with  associated PG coordinates  $(u, r, \th, \vphi)$.  To compare various quantities with  the corresponding ones in Kerr  we also assume given   two constants $m$ and  $|a|<m$.   With these  fixed  values we define, as  in Kerr,  the following functions of $(r, \th)$
\bea\lab{eq:quantities-qDeSi}
\bsplit
q &:= r+ia\cos\th,\\
\Delta &:= r^2+a^2 -2mr,\\
\Si^2 &:= (r^2+a^2)|q|^2+2mra^2(\sin\th)^2.
\end{split}
\eea

%%%%%%%%%%%%%%%%%%%%%%%%%%

\subsection{Adapted basis of $\ell=1$ modes}
\lab{section:can.ell=1basis.Mext}

%%%%%%%%%%%%%%%%%%%%%%%%%%

Recall that  the coordinates  $(\th, \vphi)$ verify $e_4(\th)=e_4(\vphi)=0$.  We define\footnote{See definition \ref{definition:canonical.ell=1basisKerr}  in Kerr.}   the basis of $\ell=1$ modes $\Jp$, $p\in\{0, +, -\}$, adapted to the PG structure according to
 \bea
 J^{(0)} :=\cos\th, \qquad  J^{(+)} :=\sin\th\cos\vphi, \qquad J^{(-)} :=\sin\th\sin\vphi.
 \eea
Clearly,
\bea
e_4(\Jp)=0.
\eea

The $\ell=1$ modes of a scalar function on $S(u,r)$ are defined as follows.
\begin{definition}
\lab{Definition:ell=1modesofascalarfunction}
Given a scalar function $f$ on a sphere $S=S(u, r)$,  we define   the $\ell=1$ modes of $f$  to be the triplet of numbers
\beaa
(f)_{\ell=1} &:=& \left(\frac{1}{r^2}\int_{S}fJ^{(0)}, \frac{1}{r^2}\int_{S}fJ^{(+)}, \frac{1}{r^2}\int_{S}fJ^{(-)}\right).
\eeaa
\end{definition}

%%%%%%%%%%%%%%%%%%%%%%%%%%%%%%%%%%%%%%%%%%%%%

\subsection{The auxiliary complex 1-forms $\Jk$, $\Jk_\pm$}
\lab{section:auxilliaryformsMext}

%%%%%%%%%%%%%%%%%%%%%%%%%%%%%%%%%%%%%%%%%%%%%

To linearize 1-forms, we    rely on a complex horizontal 1-form $\Jk$  verifying the following  properties
\bea\lab{eq:basicpropertiesofJkinMextusedtolinearize}
\nab_4\Jk = -\frac{1}{q}\Jk, \qquad \dual\Jk = -i\Jk,\qquad \Jk\c\ov{\Jk}=\frac{2(\sin\th)^2}{|q|^2}.
\eea

\begin{remark}
In Kerr, we have, see Definition \ref{def:JkandJ},
\beaa
\Jk_1=\frac{i\sin\th}{|q|}, \qquad \Jk_2=\frac{\sin\th}{|q|}.
\eeaa
\end{remark}

In addition to $\Jk$, we also introduce  the complex 1-forms $\Jk_\pm$ verifying
\bea
\lab{eq:basicpropertiesofJk_pminMextusedtolinearize}
\nab_4\Jk_\pm = -\frac{1}{q}\Jk_\pm , \qquad \dual\Jk_\pm = -i\Jk_\pm,
\eea
and
\bea
\lab{eq:propertiesJk_pm3-again}
 \bsplit
 \Jk_+\c\ov{\Jk_+}&=\frac{2(\cos\th)^2(\cos\vphi)^2+2(\sin\vphi)^2}{|q|^2},\\
  \Jk_-\c\ov{\Jk_-}&=\frac{2(\cos\th)^2(\sin\vphi)^2+2(\cos\vphi)^2}{|q|^2}.
  \end{split}
 \eea
 
\begin{remark}
In Kerr, see Definition \ref{definition:Jk_pm}, we have $\Jk_{\pm}=j_{\pm}+ i \dual j_{\pm}$ where 
 \beaa
 \bsplit
 (j_+)_1&=\frac{1}{|q|} \cos\th\cos\vphi, \quad (j_+)_2=-\frac{1}{|q|} \sin\vphi,\quad 
 (j_-)_1 =\frac{1}{|q|} \cos\th\sin\vphi, \quad (j_-)_2=\frac{1}{|q|}  \cos\vphi.
 \end{split}
 \eeaa
\end{remark}

The following lemma shows that,  upon a suitable choice for $\Jk, \Jk_\pm$ on a hypersurface $\Si$ transversal to $e_4$, it suffices in fact to transport  $\Jk, \Jk_\pm$  using the equations
\bea\lab{eq:transportequationforthe1formJk}
\nab_4\Jk = -\frac{1}{q}\Jk, \qquad  \nab_4\Jk _\pm= -\frac{1}{q}\Jk_\pm.
\eea

\begin{lemma}
\lab{Lemma:Propagation-Jk} 
Assume $\Jk, \Jk_\pm$ are  initialized  on a hypersurface $\Si$ transversal to $e_4$ such that we have, along $\Si$,
\bea
\lab{eq:basicpropertiesofJkinMextusedtolinearize-half}
\dual\Jk = -i\Jk,\qquad \Jk\c\ov{\Jk}=\frac{2(\sin\th)^2}{|q|^2},
\eea
\bea
\lab{eq:basicpropertiesofJk_pminMextusedtolinearize-half}
\bsplit
\dual\Jk_\pm &= -i\Jk_\pm,\qquad\\
 \Jk_+\c\ov{\Jk_+}&=\frac{2(\cos\th)^2(\cos\vphi)^2+2(\sin\vphi)^2}{|q|^2},\\
  \Jk_-\c\ov{\Jk_-}&=\frac{2(\cos\th)^2(\sin\vphi)^2+2(\cos\vphi)^2}{|q|^2}.
\end{split}
\eea
and
\bea
\lab{eq:propertiesJk_pm4_Si}
 \Re(\Jk_+)\c\Re(\Jk)=-\frac{1}{|q|^2}J^{(-)}, \qquad \Re(\Jk_-)\c\Re(\Jk)=\frac{1}{|q|^2}J^{(+)},
 \eea
and  extend $\Jk, \Jk_\pm$ to a neighborhood of $\Si$ according to   \eqref{eq:transportequationforthe1formJk}.
Then,  the identities  \eqref{eq:basicpropertiesofJkinMextusedtolinearize-half}, \eqref{eq:basicpropertiesofJk_pminMextusedtolinearize-half} and
\eqref{eq:propertiesJk_pm4_Si} are propagated to a neighborhood of $\Si$.
\end{lemma}

\begin{proof}
In view of the transport equation along $e_4$ satisfied by $\Jk$, we have
\beaa
\nab_4(\dual\Jk+i\Jk) +\frac{1}{q}(\dual\Jk+i\Jk) &=& 0
\eeaa
and, using also $e_4(q)=e_4(\ov{q})=1$ and $e_4(\th)=0$, 
\beaa
\nab_4\left(\Jk\c\ov{\Jk}-\frac{2(\sin\th)^2}{|q|^2}\right) +\left(\frac{1}{q}+\frac{1}{\ov{q}}\right)\left(\Jk\c\ov{\Jk}-\frac{2(\sin\th)^2}{|q|^2}\right) &=& 0
\eeaa
from which \eqref{eq:basicpropertiesofJkinMextusedtolinearize-half} is propagated to a neighborhood of $\Si$ as stated.  The identities  \eqref{eq:basicpropertiesofJk_pminMextusedtolinearize-half} can be propagated to a neighborhood of $\Si$ in the same manner. 

Next, we focus on \eqref{eq:propertiesJk_pm4_Si}. We write, setting $\Re(\Jk)  =j, \Re(\Jk_\pm) =j_\pm$, 
\beaa
\nab_4 j_+ &=&\nab_4\Re(\Jk_+)=\Re \left(\nab_4 \Jk_+ \right)=- \Re\left( \frac{1}{q} \Jk\right)= -\frac{1}{|q|^2}  \Re\big( \ov{q}  \Jk_+\big)=-\frac{1}{|q|^2} ( rj_{+} + a \cos \th  \dual j_+ ),\\
\nab_4 j &=& \nab_4\Re(\Jk)= -\frac{1}{|q|^2}  \Re\big( \ov{q}  \Jk\big)=-\frac{1}{|q|^2} ( rj  + a \cos \th  \dual j ).
\eeaa
We deduce
\beaa
\nab_4(j_+ \c j)&=& -\frac{1}{|q|^2} ( rj_{+} + a \cos \th  \dual j_+ )\c j  -\frac{1}{|q|^2} ( rj  + a \cos \th  \dual j )\c j_+=-\frac{2 r }{|q|^2}  j_+\c j.
\eeaa
Hence, since $e_4(J^{(-)})=0$,
\beaa
\nab_4\Big(j_+\c j + \frac{1}{|q|^2} J^{(-)}\Big)&=& -\frac{2 r }{|q|^2} j _{+} \c j-\frac{2 r}{|q|^4} J^{(-)}=  -\frac{2 r }{|q|^2} \Big( j_+ \c  j+ \frac{1}{|q|^2} J^{(-)}\Big).
\eeaa
We deduce
\beaa
\nab_4\Big(\Re(\Jk_+)\c\Re(\Jk)+ \frac{1}{|q|^2} J^{(-)}\Big)+\frac{2 r}{|q|^4} \Big(\Re(\Jk_+)\c\Re(\Jk)+ \frac{1}{|q|^2} J^{(-)}\Big)=0.
\eeaa
Thus, the first identity of \eqref{eq:propertiesJk_pm4_Si} is propagated to a neighborhood of $\Si$ as stated. The second identity of \eqref{eq:propertiesJk_pm4_Si} can be propagated to a neighborhood of $\Si$ in the same manner. 
\end{proof}

\begin{remark}
Note that the relations of Lemma \ref{Lemma:Propagation-Jk} are verified in Kerr, see Lemma \ref{Lemma:propertiesJk_pm}.
\end{remark}

%%%%%%%%%%%%%%%%%%%%%%%%%%%%%%%%%%%%%%%%%%%%%

\subsection{Definition of linearized quantities for an  outgoing  PG structure}

%%%%%%%%%%%%%%%%%%%%%%%%%%%%%%%%%%%%%%%%%%%%%

Recall that, given an outgoing PG  structure  with associated coordinates $(u, r, \th, \vphi)$  the following
hold true.
\begin{enumerate}
\item The following identities hold
 \beaa
 \xi=0, \quad \om=0, \quad e_4(r)=1, \quad e_4(u)= e_4(\th)=e_4(\vphi)=0, \quad  e_1(r)=e_2(r)=0.
 \eeaa
 In addition, we have
 \beaa
 \Hb=-Z.
 \eeaa
 \item The  quantities
 \beaa
\Xh, \quad \Xbh, \quad \Xib, \quad A, \quad B, \quad \Bb, \quad \Ab,
\eeaa
  vanish  in      Kerr  and therefore are   small in perturbations.
\end{enumerate}
We renormalize below all other quantities, not vanishing in Kerr\footnote{Since  $\Hb=-Z$,  $\Hb$ does not need to be included in Definition \ref{def:renormalizationofallnonsmallquantitiesinPGstructurebyKerrvalue}.}, by subtracting their  $\mbox{Kerr}(a, m) $  values  for suitably chosen constants\footnote{The precise values of these constants  will be defined  in section \ref{sec:defintioncanonicalspacetime}, see also Definition \ref{define:am-onSi}.} $(a, m)$.

\begin{definition}\lab{def:renormalizationofallnonsmallquantitiesinPGstructurebyKerrvalue} 
Let a complex horizontal 1-form $\Jk$ satisfying  \eqref{eq:basicpropertiesofJkinMextusedtolinearize}.
We  define  the following renormalizations, for given constants $(a, m)$, 
\begin{enumerate}
\item Linearization of Ricci and curvature coefficients.
\bea
\bsplit
\trXc &:= \tr X-\frac{2}{q}, \qquad\,\qquad \,\,\,\,    \trXbc := \tr\Xb+\frac{2q\Delta}{|q|^4},\\
\\
\Zc &:= Z-\frac{a\ov{q}}{|q|^2}\Jk,\qquad \qquad \quad 
\Hc := H-\frac{aq}{|q|^2}\Jk,\\
\\
  \ombc& := \omb  - \frac 1 2 \pr_r\left(\frac{\De}{|q|^2} \right),\qquad \,\, \Pc := P+\frac{2m}{q^3}.\\
\end{split}
\eea

\item Linearization of derivatives of $r, q, u$.
\bea
\bsplit
\widecheck{\DD q} &:=\DD q+a\Jk, \qquad\quad \widecheck{\DD \ov{q}} :=\DD \ov{q}-a\Jk,\\
\widecheck{e_3(r)} &:= e_3(r)+\frac{\Delta}{|q|^2},\\
\widecheck{\DD u} &:= \DD u -a\Jk,\qquad \widecheck{e_3(u)} := e_3(u) -\frac{2(r^2+a^2)}{|q|^2}.
\end{split}
\eea

\item Linearization for $\Jk$ and $\Jk_\pm$.
\bea
\bsplit
\widecheck{\ov{\DD}\c\Jk}& := \ov{\DD}\c\Jk-\frac{4i(r^2+a^2)\cos\th}{|q|^4}, \qquad \widecheck{\nab_3\Jk}:=\nab_3\Jk -\frac{\De q}{|q|^4}\Jk,\\
 \widecheck{\ov{\DD}\c\Jk_\pm}&:=\ov{\DD}\c \Jk_\pm+\frac{4}{r^2} J^{(\pm)} \pm \frac{4ia^2\cos\th}{|q|^4}J^{(\mp)}, \\  
 \widecheck{\nab_3 \Jk_\pm } &:= \nab_3 \Jk_\pm - \frac{\De q}{|q|^4} \, \Jk_{\pm} \pm   \,\frac{2a}{|q|^2}  \Jk_{\mp}. 
  \end{split}
\eea

\item Linearization for $\Jp$.
\bea
\bsplit
\widecheck{\DD J^{(0)}}  &:= \DD J^{(0)} -i\Jk, \qquad \qquad  \qquad \widecheck{\DD(J^{(\pm)})}:= \DD(J^{(\pm)})-\Jk_{\pm},
\\
\widecheck{e_3(J^{(+)})}&:=e_3(J^{(+)}) + \frac{2a}{|q|^2}J^{(-)}, \qquad
 \widecheck{e_3(J^{(-)})}:=e_3(J^{(-)}) - \frac{2a}{|q|^2}J^{(+)}.
 \end{split}
\eea
 \end{enumerate}
\end{definition}

%%%%%%%%%%%%%%%%%%%%%%%%%%%%%%%%%%%%%%%%%%%%%

\subsection{Definition of the notations $\Ga_b$ and $\Ga_g$ for error terms}

%%%%%%%%%%%%%%%%%%%%%%%%%%%%%%%%%%%%%%%%%%%%%

\begin{definition}
\lab{definition.Ga_gGa_b}
The set of all linearized quantities is of the form $\Ga_g\cup \Ga_b$ with  $\Ga_g,  \Ga_b$
 defined as follows.
 \begin{enumerate}
\item 
 The set   $\Ga_g$   with
 \bea
 \bsplit
 \Ga_g &= \Big\{\trXc,\quad  \Xh,\quad \Zc,\quad \trXbc , \quad r\Pc, \quad  rB, \quad  rA\Big\}.
 \end{split}
 \eea
 
 \item  The set  $\Ga_b=\Ga_{b,1}\cup \Ga_{b, 2}\cup \Ga_{b,3} \cup \Ga_{b, 4}$   with
 \bea
 \bsplit
 \Ga_{b,1}&= \Big\{\Hc, \quad \Xbh, \quad \ombc, \quad \Xib,\quad  r\Bb, \quad \Ab\Big\},\\
  \Ga_{b, 2}&= \Big\{r^{-1}\widecheck{e_3(r)}, \quad \widecheck{\DD q}, \quad \widecheck{\DD \ov{q}}, \quad \widecheck{\DD u},  \quad r^{-1}\widecheck{e_3(u)}\Big\},\\
  \Ga_{b, 3}&= 
   \Big\{ \widecheck{\DD(J^{(0)})}, \quad \widecheck{\DD(J^{(\pm)})}, \quad e_3(J^{(0)}), \quad    \widecheck{e_3(J^{(\pm)})}\Big\}, \\
   \Ga_{b,4}&=\Bigg\{ r\,\widecheck{\ov{\DD}\c\Jk}, \quad r\,\DD\hot\Jk, \quad r\,\widecheck{\nab_3\Jk}, \quad
   r\,\widecheck{\ov{\DD}\c\Jk_\pm}, \quad r\,\DD\hot\Jk_\pm, \quad r\,\widecheck{\nab_3\Jk_\pm} \Bigg\}. 
   \end{split}
 \eea
\end{enumerate}
We also define, with the help of  the weighted derivatives $\dk=\{\nab_3, r\nab_4, \dkb=r\nab \}$,
\beaa
\Ga_g^{(s)} = \dk^{\leq s}\Ga_g, \qquad  \Ga_b^{(s)} = \dk^{\leq s}\Ga_b.
\eeaa
\end{definition}

\begin{remark}
The justification for the above decompositions has to do with the expected  decay properties of the linearized  components of the outgoing PG structure. More precisely, we expect that,    see sections \ref{section:main-norms}
 and \ref{section:Bootstrap assumptions} for details,
\bea\lab{eq:expectedbehaviorGabGag:chap2}
\bsplit
\big|\Ga^{(s)}_g|&\les \ep \min\Big\{ r^{-2 } u^{-1/2-\dec},  \, r^{-1} u^{-1-\dec} \Big\}, \\
\big|\nab_3 \Ga^{(s-1)}_g| &\les \ep  r^{-2 } u^{-1-\dec},\\
\big|\Ga^{(s)}_b\big| &\les \ep  r^{-1 } u^{-1-\dec},
\end{split}
\eea
for a small constant $\dec>0$.
\end{remark}

%%%%%%%%%%%%%%%%%%%%%%%%%%%%

\subsection{Approximate Killing vectorfield  $\T$} 
\lab{section:Definition-boldT}

%%%%%%%%%%%%%%%%%%%%%%%%%%%%

Given an outgoing PG structure  on $\MM$  with adapted coordinates $(u, r, \th, \vphi)$, i.e. $e_4(r)=1$ and $e_4(u)=e_4(\th)=e_4(\vphi)=0$,  we define a vectorfield $\T$ as follows.
\begin{definition}
The vectorfield $\T$ is defined  by
\beaa
\T &:=& \frac{1}{2}\left(e_3+\frac{\Delta}{|q|^2}e_4 -2a\Re(\Jk)^be_b\right).
\eeaa
\end{definition} 

Note that we have
\bea
\bsplit
&\T(u) = 1+\frac{1}{2}\left(\widecheck{e_3(u)} -2a\Re(\Jk)\c\widecheck{ \nab u }\right),\qquad \T(r) = \frac{1}{2}\widecheck{e_3(r)},\\
&\T(\cos\th) = \frac{1}{2}\left(e_3(\cos\th) - 2a\Re(\Jk)\c\widecheck{ \nab \cos\th }\right),
\end{split}
\eea
and\footnote{We use in particular the identity 
\beaa
|\Re(\Jk)|^2= \frac{(\sin\th)^2}{|q|^2}.
\eeaa} 
\bea
\g(\T, \T) &=& -1+\frac{2mr}{|q|^2}.
\eea

The following proposition shows that $\T$ is an approximate Killing vectorfield. 
\begin{proposition}
\lab{Proposition:deftensorT}
We have $\piT_{44} =0$,  $ \piT_{4a}\in \Ga_g $ and all  other components of  $\piT$ are in $\Ga_b$.  Moreover
\beaa
g^{ab}\,\piT_{ab} &=& \Ga_g.
\eeaa
In addition
\beaa
\g(\D_a \T, e_4),\,  \g(\D_4 \T, e_a) \in \Ga_g, \qquad 
\g(\D_a \T, e_3),\,  \g(\D_3 \T, e_a) \in \Ga_b,
\eeaa
and
\beaa
 \g(\D_a\T,  e_b)= -\frac{2amr\cos\th}{|q|^4}\in_{ab}+\Ga_b.
 \eeaa
\end{proposition}

\begin{proof}
See section \ref{appendix:Proof-Proposition:deftensorT}. 
\end{proof}

%%%%%%%%%%%%%%%%%%%%%%%%%%%%%%%%%%%%%%%%%%%%

\section{Ingoing PG structures}\lab{sec:Principalingoinggeodesicstructures}

%%%%%%%%%%%%%%%%%%%%%%%%%%%%%%%%%%%%%%%%%%%%%

%%%%%%%%%%%%%%%%%%%%%%%%%%%%%%%%%%%%%%%%%%%%

\subsection{Definition of ingoing PG structures}
\lab{sec:basicdefinitionforingoingPGstruct}

%%%%%%%%%%%%%%%%%%%%%%%%%%%%%%%%%%%%%%%%%%%%%

Ingoing PG structures are PG structures where  the roles of $e_3$ and $e_4$ are reversed 
compared to outgoing ones. In particular, we have for ingoing PG structures   
 \bea
 \D_3e_3=0, \qquad e_3(r)=-1, \qquad \nab(r)=0,
 \eea
 and hence
 \bea
 \xib=0, \qquad \omb=0, \qquad \eta=\ze.
 \eea
 
Also, in addition to $r$, we define the ingoing PG coordinates  $(\ub, \th, \vphi)$ such that 
\bea
e_3(\ub)=e_3(\th)= e_3(\vphi)=0.
\eea

Finally, we introduce horizontal complex 1-forms $\Jk$, $\Jk_\pm$ satisfying 
\bea\lab{eq:transportequationforthe1formJk:ingoingcase}
\nab_3\Jk = \frac{1}{\ov{q}}\Jk, \qquad  \nab_3\Jk _\pm= \frac{1}{\ov{q}}\Jk_\pm,
\eea
and we define the adapted basis of $\ell=1$ modes $\Jp$, $p=0,+,-$, by 
\bea
J^{(0)}=\cos\th, \qquad J^{(+)}=\sin\th\cos\vphi, \qquad J^{(-)}=\sin\th\sin\vphi,
\eea
so that we have in particular $e_3(\Jp)=0$ for $p=0,+,-$.

All the equations of Proposition \ref{prop-nullstrandBianchi:complex:outgoing} and Proposition \ref{prop:e_4(xyz)} 
for outgoing PG structures have a counterpart for ingoing PG structures. The equations can be easily deduced from the ones of Proposition \ref{prop-nullstrandBianchi:complex:outgoing} and Proposition \ref{prop:e_4(xyz)} by performing the following substitutions
\beaa
&& u\to \ub, \quad r\to r, \quad \th\to \th, \quad \vphi\to \vphi, \quad e_4\to e_3, \quad e_3\to e_4, \quad e_a\to e_a, \\
&& \a\to \aa, \quad \b\to -\bb, \quad \rho\to \rho, \quad \rhod\to -\rhod, \quad \bb\to -\b,\quad \aa\to \a,\\
&& \xi\to \xib, \quad \om\to \omb, \quad \chi\to \chib, \quad \eta\to\etab, \quad \etab\to \eta, \quad\ze\to -\ze, \quad \chib\to \chi,\quad \omb\to\om, \quad \xib\to\xi.
\eeaa

%%%%%%%%%%%%%%%%%%%%%%%%%%%%%%%%%%%%%%%%%%%%%%%%

\subsection{Ingoing PG structures in Kerr}\lab{sec:KerrvaluesiningoingPGstructure:chap2}

%%%%%%%%%%%%%%%%%%%%%%%%%%%%%%%%%%%%%%%%%%%%%%%%

 In Kerr, relative to  BL coordinates $(t, r, \th, \phi) $, the ingoing principal null pair is given by, see \eqref{eq:Inc.PGdirections-Kerr}, 
 we have
 \beaa
 e_4 = \frac{r^2+a^2}{|q|^2} \pr_t +\frac{\De}{|q|^2} \pr_r +\frac{a}{|q|^2} \pr_{\phi}, \qquad 
 e_3=\frac{r^2+a^2}{\De} \pr_t -\pr_r +\frac{a}{\De} \pr_{\phi}.
\eeaa
Also,  the functions $(\ub, \vphi)$ are  given by
\beaa
\ub :=t+f(r), \quad f'(r)=\frac{r^2+a^2}{\Delta}, \qquad \vphi:= \phi+h(r), \quad h'(r)=\frac{a}{\Delta}. 
\eeaa

\begin{remark}
Note that we have indeed
\beaa
e_3(r)=-1, \qquad e_3(\ub)= e_3(\th)= e_3 (\vphi)=0.
\eeaa
\end{remark}

$\Jk$, $\Jk_\pm$ and $\Jp$ are still defined according respectively to Definition \ref{def:JkandJ}, Definition \ref{definition:Jk_pm} and Definition \ref{definition:canonical.ell=1basisKerr}.

%%%%%%%%%%%%%%%%%%%%%%%%%%%%%%%%%%%%%%%%%%%%

\subsection{Linearization of ingoing PG structures}

%%%%%%%%%%%%%%%%%%%%%%%%%%%%%%%%%%%%%%%%%%%%%

To linearize 1-forms, we    rely on a complex horizontal 1-form $\Jk$  verifying the following  properties
\bea\lab{eq:basicpropertiesofJkinMintusedtolinearize}
\nab_3\Jk = \frac{1}{\ov{q}}\Jk, \qquad \dual\Jk = -i\Jk,\qquad \Jk\c\ov{\Jk}=\frac{2(\sin\th)^2}{|q|^2}.
\eea

\begin{definition}\lab{def:renormalizationofallnonsmallquantitiesinPGstructurebyKerrvalue:ingoingcase}
Let a complex horizontal 1-form $\Jk$ satisfying  \eqref{eq:basicpropertiesofJkinMintusedtolinearize}.
We  define  the following renormalizations, for given constants $(a, m)$,
\begin{enumerate}
\item Linearization of Ricci and curvature coefficients.
\bea
\bsplit
\trXc &:= \tr X-\frac{2\ov{q}\De}{|q|^4}, \qquad\,\qquad     \trXbc := \tr\Xb+\frac{2}{\ov{q}},\\
\\
\Zc &:= Z-\frac{aq}{|q|^2}\Jk,\qquad \qquad \quad\,\,\,
\Hbc := \Hb+\frac{a\ov{q}}{|q|^2}\Jk,\\
\\
\omc & := \om  + \frac{1}{2}\pr_r\left(\frac{\De}{|q|^2} \right),\qquad \quad\, \Pc := P+\frac{2m}{q^3}.\\
\end{split}
\eea

\item Linearization of derivatives of $r, q, \ub$.
\bea
\bsplit
\widecheck{\DD q} &:=\DD q+a\Jk, \qquad \widecheck{\DD \ov{q}} :=\DD \ov{q}-a\Jk,\\
\widecheck{e_4(r)} &:= e_4(r)-\frac{\Delta}{|q|^2},\\
\widecheck{\DD\ub} &:= \DD\ub -a\Jk,\qquad \widecheck{e_4(\ub)} := e_4(\ub) -\frac{2(r^2+a^2)}{|q|^2}.
\end{split}
\eea

\item Linearization for $\Jk$ and $\Jk_\pm$.
\bea
\bsplit
\widecheck{\ov{\DD}\c\Jk}& := \ov{\DD}\c\Jk-\frac{4i(r^2+a^2)\cos\th}{|q|^4}, \qquad \widecheck{\nab_4\Jk}:=\nab_4\Jk +\frac{\De \ov{q}}{|q|^4}\Jk,\\
\widecheck{\ov{\DD}\c\Jk_\pm}&:=\ov{\DD}\c \Jk_\pm+\frac{4}{r^2} J^{(\pm)}  \pm \frac{4ia^2\cos\th}{|q|^4}J^{(\mp)}, \\
\widecheck{   \nab_4 \Jk_\pm } &:= \nab_4 \Jk_\pm + \frac{\De \ov{q}}{|q|^4} \, \Jk_{\pm} \pm   \,\frac{2a}{|q|^2}  \Jk_{\mp}.
\end{split}
\eea

\item Linearization for $\Jp$.
\bea
\bsplit
\widecheck{\DD J^{(0)}}  &:= \DD J^{(0)} -i\Jk, \qquad \qquad  \qquad \widecheck{\DD(J^{(\pm)})}:= \DD(J^{(\pm)})-\Jk_{\pm},
\\
\widecheck{e_4(J^{(+)})}&:=e_4(J^{(+)}) + \frac{2a}{|q|^2}J^{(-)}, \qquad
\widecheck{e_4(J^{(-)})}:=e_4(J^{(-)}) - \frac{2a}{|q|^2}J^{(+)}.
\end{split}
\eea
\end{enumerate}
\end{definition}

\begin{definition}
\lab{definition.Ga_gGa_b:ingoingcase}
The set of all linearized quantities is of the form $\Ga_g\cup \Ga_b$ with  $\Ga_g,  \Ga_b$
 defined as follows.
 
\begin{enumerate}
\item 
 The set   $\Ga_g=\Ga_{g,1}\cup \Ga_{g, 2}$   with
 \bea
 \bsplit
 \Ga_{g,1} &= \Big\{\Xi, \quad \omc, \quad \trXc,\quad  \Xh,\quad \Zc, \quad \Hbc,\quad \trXbc , \quad r\Pc, \quad  rB, \quad  rA\Big\},\\
   \Ga_{g, 2}&= \left\{r\widecheck{e_4(r)}, \quad r\widecheck{e_4(\ub)}, \quad re_4(J^{(0)}), \quad r\widecheck{e_4(J^{(\pm)})}, \quad r^2\widecheck{\nab_4\Jk}, \quad r^2\widecheck{\nab_4\Jk_\pm} \right\}.
 \end{split}
 \eea
 
 \item  The set  $\Ga_b=\Ga_{b,1}\cup \Ga_{b, 2}$   with
 \bea
 \bsplit
 \Ga_{b,1}&= \Big\{\Xbh,\quad  r\Bb, \quad \Ab, \quad \widecheck{\DD q}, \quad \widecheck{\DD \ov{q}}, \quad \widecheck{\DD\ub}\Big\},\\
  \Ga_{b, 2}&=\left\{ \widecheck{\DD(J^{(0)})}, \quad \widecheck{\DD(J^{(\pm)})}, \quad r\,\widecheck{\ov{\DD}\c\Jk}, \quad r\,\DD\hot\Jk, \quad
   r\,\widecheck{\ov{\DD}\c\Jk_\pm}, \quad r\,\DD\hot\Jk_\pm\right\}.
 \end{split}
 \eea
\end{enumerate}
\end{definition}

%%%%%%%%%%%%%%%%%%%%%%%%%%%%%%%%%%%%%%%%%%%%%

\section{Principal temporal structures}
\lab{sec:introductionofPTstrctures}

%%%%%%%%%%%%%%%%%%%%%%%%%%%%%%%%%%%%%%%%%%%%%

The PG structures we have studied so far are perfectly adequate  for  deriving  decay estimates  but are  deficient in terms of loss of derivatives and thus inadequate   for deriving  boundedness estimates for the top derivatives of the Ricci coefficients.  Indeed  the  
 $\nab_4$ equations  for $\tr \Xb$, $\Xbh$  and $\Xib$ in Proposition    \ref{prop-nullstrandBianchi:complex:outgoing}   contain angular derivatives\footnote{This loss can be overcome for integrable foliations such as geodesic foliations and double null foliations relying on elliptic Hodge systems on 2-spheres of the foliation, but not for non integrable structures such as PG structures.} of other Ricci coefficients. Similarly,  the same situation occurs for  ingoing PG  structures where  the $ \nab_3$ equations for $\tr X$, $\Xh$, and $\Xi$   are manifestly losing derivatives. Thus, in order to derive  boundedness estimates for the top derivatives of the Ricci coefficients, we are forced to  introduce  new frames   which we call principal temporal  (PT). We first introduce outgoing PT structures, and then ingoing PT structures.

%%%%%%%%%%%%%%%%%%%%%%%%%%%%%%%%%%%%%%%%%%%%%

\subsection{Outgoing PT structures}
\lab{sec:outgoingPTstructures:chap2} 

%%%%%%%%%%%%%%%%%%%%%%%%%%%%%%%%%%%%%%%%%%%%%

\begin{definition}
\lab{definition:outgoingPT}
An outgoing  PT structure    $\{ (e_3, e_4, \HH), r, \th, \Jk\}$     on $\MM$ consists of  a null pair $(e_3, e_4)$, the induced horizontal structure   $\HH$,   functions $(r, \th)$, and a horizontal  1-form $\Jk$ such that the following hold true:
\begin{enumerate}
\item   $e_4$ is geodesic.

\item We have
\bea
e_4(r)=1,\qquad    e_4(\th)=0, \qquad \nab_4 (q \Jk)=0, \qquad  q= r+a i \cos \th.
\eea

\item We have 
\bea
\Hb=-\frac{a\ov{q}}{|q|^2}\Jk.
\eea
\end{enumerate}
An extended outgoing PT structure possesses,  in addition,  a scalar function $u$ verifying $e_4(u)=0$. 
\end{definition}

\begin{definition}\lab{def:initialoutgoingPTdataset}
An outgoing PT  initial data set  consists of a hypersurface $\Si$ transversal to $e_4$ together with a null pair $(e_3, e_4)$, the induced horizontal structure $\HH$,  scalar   functions $(r, \th)$,  and a horizontal  1-form $\Jk$,  all defined on $\Si$.
\end{definition}

\begin{lemma}
\lab{lemma:constructionoutgoingPTframes}
Any outgoing PT  initial data set, as in Definition \ref{def:initialoutgoingPTdataset}, can be  locally extended to   an outgoing PT structure.
\end{lemma} 

\begin{proof}
Extend first $e_4$  in a neighborhood of $\Si$ such that $\D_{e_4} e_4 =0$. Also,  extend  $e_a$, $a=1,2$, such that $\D_{e_4}e_a=0$. Finally, extend $e_3$  to  be the unique null companion of $e_4$ orthogonal to $(e_1, e_2)$. Thus,
\beaa
\xi=0, \qquad \om=0, \qquad \etab=0.
\eeaa
We also extend $(r, \th)$ and $\Jk$ so that,
\beaa
e_4(r)=1, \qquad  e_4(\th)=0, \qquad \nab_4(q\Jk)=0, \qquad q=r+i a \cos \th.
\eeaa

Since $\etab=0$, $(e_1, e_2, e_3, e_4)$ is not a PT frame. We  look for a  new frame  $(e_1', e_2', e_3', e_4')$     related  to  $(e_1, e_2, e_3, e_4)$ by the formula  \eqref{General-frametransformation} 
with transition parameters $(f=0, \fb, \la=1)$, i.e. 
\bea
\bsplit
  e_4'&=e_4,\\
  e_a'&= e_a +\frac 1 2  \fb_a  e_4,\\
 e_3'&=e_3 + \fb^b e_b  + \frac 1 4 |\fb|^2 e_4,
 \end{split}
 \eea
 where $\fb$ will be chosen later. Since $e_4'=e_4$, and since $e_4$ is geodesic, note first that we have in the frame $(e_1', e_2', e_3', e_4')$
\bea
\xi'=0, \qquad \om'=0.
\eea

Next, we define the functions $(r', \th')$ and  the horizontal 1-form $\Jk'$, for the horizontal structure $\HH'$ induced by $(e_3', e_4')$, by    
\bea
 r'=r, \qquad \th'=\th, \qquad q'=r'+i\cos(\th')=q, \qquad  \Jk'_{e_a'} =\Jk_{e_a}, \quad a=1,2.
 \eea 
Since $e_4'=e_4$, we have in particular 
\bea
e_4'(r')=1, \qquad e_4'(\th')=0.
\eea
Also, we compute $\nab_4'(q'\Jk')$. Since $ \Jk'_{e_a'} =\Jk_{e_a}$, $e_4'=e_4$, $q'=q$, and $\nab_4(q\Jk)=0$, we have
\beaa
\nab_4'(q'\Jk')_a &=& e_4'(q'\Jk'_a)-\g(\D_{e_4'}e_a', e_b')\Jk_b' = e_4(q\Jk_a)-\g(\D_{e_4'}e_a', e_b')\Jk_b\\
&=& \nab_4(q\Jk)_a -\Big(\g(\D_{e_4'}e_a', e_b')-\g(\D_{e_4}e_a, e_b)\Big)\Jk_b\\
&=& -\Big(\g(\D_{e_4'}e_a', e_b')-\g(\D_{e_4}e_a, e_b)\Big)\Jk_b.
\eeaa
Since 
\beaa
\g(\D_{e_4'}e_a', e_b') &=& \g\left(\D_{e_4}\left(e_a +\frac 1 2  \fb_a  e_4\right), e_b +\frac 1 2  \fb_b  e_4\right)\\
&=& \g(\D_{e_4}e_a, e_b) -\fb_b\xi_a+\fb_a\xi_b=\g(\D_{e_4}e_a, e_b)
\eeaa
where we used the fact that $\xi=0$, we infer
\bea
\nab_4'(q'\Jk') &=& 0.
\eea

In view of the above, in order for $\{ (e_3', e_4', \HH'), r', \th', \Jk'\}$ to satisfy all the properties of an outgoing PT structure, it remains to obtain the desired identity for $\etab'$ which is related to the choice of $\fb$. To this end, we compute\footnote{We use here a more precise transformation formula for $\etab$ than the one derived in Proposition \ref{Proposition:transformationRicci}.} 
\beaa
2\etab_a' &=& \g(\D_{e_4'}e_3', e_a') = \g\left(\D_{e_4}\left(e_3 + \fb^b e_b  + \frac 1 4 |\fb|^2 e_4\right), e_a +\frac 1 2  \fb_a  e_4\right)\\
&=& \g\left(\D_{e_4}e_3, e_a +\frac 1 2  \fb_a  e_4\right)+e_4(\fb_a)+\fb^b\g\left(\D_{e_4}e_b, e_a +\frac 1 2  \fb_a  e_4\right)\\
&&+\frac 1 4 |\fb|^2\g\left(\D_{e_4}e_4, e_a +\frac 1 2  \fb_a  e_4\right)\\
&=& 2\etab_a -2\om\fb_a+\nab_4\fb_a -\fb_b\fb_a\xi_b+\frac{1}{2}|\fb|^2\xi_a
\eeaa
and since $\xi=0$, $\om=0$ and $\etab=0$, we infer
\bea
2\etab' &=& \nab_4\fb.
\eea

We now fix $\fb$, and hence the frame $(e_1', e_2', e_3', e_4')$, as the solution of the following transport equation
\bea
\nab_4\fb = -2\Re\left(\frac{a\ov{q'}}{|q'|^2}\Jk'\right), \qquad \fb_{|_{\Si}}=0.
\eea
In view of the above transformation formula for $\etab'$, we infer $\etab' =- \Re(\frac{a\ov{q'}}{|q'|^2}\Jk')$, and hence
\bea
\Hb' =- \frac{a\ov{q'}}{|q'|^2}\Jk'.
\eea
Thus, $\{ (e_3', e_4', \HH'), r, \th, \Jk'\}$ satisfies all the properties of an outgoing PT structure, and, since $\fb=0$ on $\Si$,  coincides with the initial outgoing PT data set on $\Si$. This ends the proof of Lemma \ref{lemma:constructionoutgoingPTframes}.
\end{proof}

%%%%%%%%%%%%%%%%%%%%%%%%%%%%%%%%%%%%%%%%%%%%%

\subsection{Null structure equations in an outgoing  PT frame}

%%%%%%%%%%%%%%%%%%%%%%%%%%%%%%%%%%%%%%%%%%%%%

\begin{proposition}
\lab{Prop:NullStr-outgoingPTframe}
Consider an outgoing PT structure. Then,  the   equations in the $e_4$ direction for the Ricci coefficients  of the  outgoing PT frame take the form
\beaa
\nab_4\tr X +\frac{1}{2}(\tr X)^2 &=& -\frac{1}{2}\Xh\c\ov{\Xh},\\
\nab_4\Xh+\Re(\tr X)\Xh &=& -A,
\\
\nab_4\tr\Xb +\frac{1}{2}\tr X\tr\Xb  &=& -\DD\c\left(\frac{aq}{|q|^2}\ov{\Jk}\right) +  \frac{a^2}{|q|^2}|\Jk|^2 +2\ov{P} -\frac{1}{2}\Xh\c\ov{\Xbh},\\
\nab_4\widehat{\Xb} +\frac{1}{2}\tr X\, \widehat{\Xb}  &=& -\DD\hot\left(\frac{a\ov{q}}{|q|^2}\Jk\right)  +   \frac{a^2(\ov{q})^2}{|q|^4}\Jk\hot\Jk -\frac{1}{2}\ov{\tr\Xb} \widehat{X},
\\
\nab_4Z +\frac{1}{2}\tr X Z &=&  -\frac{1}{2}\tr X\frac{a\ov{q}}{|q|^2}\Jk-\frac{1}{2}\widehat{X}\c\left(\ov{Z}+\frac{aq}{|q|^2}\ov{\Jk}\right) -B,\\
\ \nab_4\Xib &=& -\nab_3\left(\frac{a\ov{q}}{|q|^2}\Jk\right) -\frac{1}{2}\ov{\tr\Xb}\left(\frac{a\ov{q}}{|q|^2}\Jk+H\right) -\frac{1}{2}\Xbh\c\left(\frac{aq}{|q|^2}\ov{\Jk}+\ov{H}\right) -\Bb,\\
\nab_4H &=&  -\frac{1}{2}\ov{\tr X}\left(H+\frac{a\ov{q}}{|q|^2}\Jk\right) -\frac{1}{2}\Xh\c\left(\ov{H} + \frac{aq}{|q|^2}\ov{\Jk}\right) -B,
\\
\nab_4\omb &=&   \left(\eta+\Re\left(\frac{a\ov{q}}{|q|^2}\Jk\right)\right)\c\ze +\eta\c\Re\left(\frac{a\ov{q}}{|q|^2}\Jk\right) +\rho.
\eeaa
\end{proposition}

\begin{remark}\lab{rmk:whythePTframedoesnotloosederivatives}
The main feature of the PT gauge choice is that treating all equations in Proposition \ref{Prop:NullStr-outgoingPTframe}  as transport equations in $e_4$ does not loose derivatives. Indeed,  the RHS of all the equations depend only on   Ricci and curvature coefficients, as well as first order derivatives of $q$ and $\Jk$ (in the RHS of the equations for $\tr \Xb$, $\Xbh$  and $\Xib$). The point is that first order derivatives of $q$ and $\Jk$ can be controlled at the same level of regularity\footnote{In the case of PG frames, these terms are replaced by first order derivatives of $Z$ leading to a loss of derivative.} than the Ricci coefficients of the PT frame. 
\end{remark}

\begin{proof}
These equations follow immediately from plugging the identities 
\beaa
\Xi=0, \qquad \om=0,\qquad  \Hb = -   \frac{a\ov{q}}{|q|^2}\Jk,
\eeaa
satisfied by an outgoing PT frame in  the general formulas of Proposition \ref{prop-nullstr:complex}.  
\end{proof}

%%%%%%%%%%%%%%%%%%%%%%%%%%%%%%%%%

\subsection{Linearized quantities for outgoing PT structures}

%%%%%%%%%%%%%%%%%%%%%%%%%%%%%%%%%

Given an extended  outgoing PT  structure   $\{ (e_3, e_4, \HH), r, \th,  u, \Jk\}$, the following   holds true:
\begin{enumerate}
\item We have
 \beaa
 \xi=\om=0, \quad e_4(r)=1, \quad e_4(u)= e_4(\th)=0, \qquad \nab_4(r\Jk)=0.
 \eeaa
 In addition, we have
 \beaa
 \Hb=-\frac{a\ov{q}}{|q|^2}\Jk.
 \eeaa
 
 \item The  quantities
 \beaa
\Xh, \quad \Xbh, \quad \Xib, \quad A, \quad B, \quad \Bb, \quad \Ab, \quad \DD r, \quad e_3(\cos\th), \quad \DD\hot\Jk,
\eeaa
  vanish  in      Kerr  and therefore are   small in perturbations.
\end{enumerate}
We renormalize below all other quantities, not vanishing in Kerr\footnote{Since  $ \Hb=-\frac{a\ov{q}}{|q|^2}\Jk$,  $\Hb$ does not need to be included in Definition \ref{def:renormalizationofallnonsmallquantitiesinPTstructurebyKerrvalue}.}, by subtracting their  $\mbox{Kerr}(a, m) $  values.

\begin{definition}\lab{def:renormalizationofallnonsmallquantitiesinPTstructurebyKerrvalue} 
We consider  the following renormalizations, for given constants $(a, m)$, 
\bea
\begin{split}
\trXc &:= \tr X-\frac{2}{q}, \qquad\qquad\,\,\,  \trXbc := \tr\Xb+\frac{2q\Delta}{|q|^4},\\
\Zc &:= Z -\frac{a\ov{q}}{|q|^2}\Jk,\qquad\qquad\,\,\,\, \Hc := H-\frac{aq}{|q|^2}\Jk,\\
\ombc &:= \omb-\frac 1 2 \pr_r\left(\frac{\De}{|q|^2} \right), \qquad \Pc := P+\frac{2m}{q^3},
\end{split}
\eea
as well as 
\bea
\begin{split}
\widecheck{e_3(r)} &:= e_3(r)+\frac{\Delta}{|q|^2},\qquad \widecheck{\DD(\cos\th)} := \DD(\cos(\th)) -i\Jk,\\
\widecheck{\DD u} &:=\DD u- a\Jk, \qquad\qquad\,\,\,\,  \widecheck{e_3(u)}:=e_3(u)-\frac{2(r^2+a^2)}{|q|^2},
\end{split}
\eea
and 
\bea
\widecheck{\ov{\DD}\c\Jk} := \ov{\DD}\c\Jk-\frac{4i(r^2+a^2)\cos\th}{|q|^4}, \qquad \widecheck{\nab_3\Jk}:=\nab_3\Jk -\frac{\De q}{|q|^4}\Jk.
\eea
\end{definition}

%%%%%%%%%%%%%%%%%%%%%%%%%%

\subsection{Transport equations for $(f, \fb, \la)$}
\lab{section:transport(f,fb,la):PTcase}

%%%%%%%%%%%%%%%%%%%%%%%%%%

We will need to compare PG structures with PT structures. To this end, we will control coefficients $(f, \fb, \la)$ corresponding to the change of frame between PG frames and PT frames. We derive in this section transport equations for $(f, \fb, \la)$ in the case where the second frame is an outgoing PT frame\footnote{This corresponds to the analog of section \ref{section:transport(f,fb,la)} where transport equations for $(f, \fb, \la)$ are derived in the case where the second frame is an outgoing PG frame.}.

The following is the analog of Corollary \ref{cor:transportequationine4forchangeofframecoeff:simplecasefirst}.
\begin{corollary}\lab{cor:transportequationine4forchangeofframecoeff:simplecasefirst:PTcase}
Under the assumption 
\beaa
\xi'=0, \qquad \om'=0, 
\eeaa
we  have the following transport equations for  $(f, \fb, \la)$
\beaa
\nab_{\la^{-1}e_4'}f+\frac{1}{2}(\trch f -\atrch\dual f)+2\om f &=& -2\xi - f\c\chih+E_1(f, \Ga),\\
\la^{-1}e_4'(\log\la) &=& 2\om+f\c(\ze-\etab)+E_2(f, \Ga),\\
\nab_{\la^{-1}e_4'}\fb &=& 2(\etab'-\etab)  -\frac{1}{2}(\trchb f -\atrchb\dual f) +2\om\fb\\
&&-f\c\chibh  +E_6(f, \fb, \Ga),
\eeaa
where $E_1(f, \Ga)$, $E_2(f, \Ga)$ and $E_6(f, \fb, \Ga)$ are given by
\beaa
E_1(f, \Ga) &=& -(f\c\ze)f-\frac{1}{2}|f|^2\eta +\frac{1}{2}|f|^2\etab +O(f^3\Ga),\\
E_2(f, \Ga) &=& - \frac{1}{2}|f|^2\omb  -\frac{1}{4}\trchb|f|^2+O(f^3\Ga+f^2\chibh),
\eeaa
and
\beaa
E_6(f, \fb, \Ga) &=& -(f\c\etab)\fb +\frac 1 2(f\c\ze)\fb +O((f, \fb)^3\Ga).
\eeaa
\end{corollary}

\begin{proof}
The transport equations for $f$ and $\la$ have already been derived in Corollary \ref{cor:transportequationine4forchangeofframecoeff:simplecasefirst}. We thus focus on the one for $\fb$. Assuming that $\xi'=0$, we have in view of the frame transformation \eqref{General-frametransformation}
\beaa
2\etab_a' &=&  \g(\D_{e_4'}e_3', e_a')= \g(\D_{\la^{-1}e_4'}(\la e_3'), e_a') = \g\left(\D_{\la^{-1}e_4'}\left(e_3 +  \fb^be_b' -\frac 1 4 |\fb|^2\la^{-1} e_4'\right), e_a'\right)\\
&=& \g\left(\D_{\la^{-1}e_4'}e_3 , e_a'\right) +\g\left(\D_{\la^{-1}e_4'}\left( \fb^be_b' \right), e_a'\right) -\frac 1 2 |\fb|^2\la^{-2}\xi'_a\\
&=& \g\left(\D_{\la^{-1}e_4'}e_3 , e_a'\right) +\g\left(\D_{\la^{-1}e_4'}\left( \fb^be_b' \right), e_a'\right).
\eeaa
We compute the terms on the right-hand side
\beaa
\g\left(\D_{\la^{-1}e_4'}e_3 , e_a'\right) &=& \g\left(\D_{\la^{-1}e_4'}e_3 , \left(\de_a^b +\frac{1}{2}\fb_af^b\right) e_b +\frac 1 2  \fb_a  e_4 \right)\\
&=& \left(\de_a^b +\frac{1}{2}\fb_af^b\right)\g\left(\D_{\la^{-1}e_4'}e_3 ,  e_b  \right)+\frac 1 2  \fb_a\g\left(\D_{\la^{-1}e_4'}e_3 ,   e_4 \right)\\
&=& \left(\de_a^b +\frac{1}{2}\fb_af^b\right)\g\left(\D_{e_4+f^ce_c+\frac{1}{4}|f|^2e_3}e_3 ,  e_b  \right)+\frac 1 2  \fb_a\g\left(\D_{e_4+f^be_b+\frac{1}{4}|f|^2e_3}e_3 ,   e_4 \right)\\
&=& \left(\de_a^b +\frac{1}{2}\fb_af^b\right)\left(2\etab_b+f^c\chib_{cb}+\frac{1}{2}|f|^2\xib_b\right)+\frac 1 2  \fb_a\big(-4\om-f\c\ze \big)\\
&&+O((f, \fb)^3\Ga)\\
&=& 2\etab_a+(f\c\etab)\fb_a +f^c\chib_{ca}+\frac{1}{2}|f|^2\xib_a+\frac 1 2  \fb_a\big(-4\om-f\c\ze \big)+O((f, \fb)^3\Ga)\\
\eeaa
and
\beaa
\g\left(\D_{\la^{-1}e_4'}\left( \fb^be_b' \right), e_a'\right) &=& \la^{-1}e_4'(\fb_a)+ \fb^b\g\left(\D_{\la^{-1}e_4'}e_b', e_a'\right) = \la^{-1}e_4'(\fb_a) - \fb^b\g\left(\D_{\la^{-1}e_4'}e_a', e_b'\right)\\
&=& \nab_{\la^{-1}e_4'}\fb_a.
\eeaa
We infer
\beaa
2\etab_a' &=& \g\left(\D_{\la^{-1}e_4'}e_3 , e_a'\right) +\g\left(\D_{\la^{-1}e_4'}\left( \fb^be_b' \right), e_a'\right)\\
&=& 2\etab_a+\nab_{\la^{-1}e_4'}\fb_a +\frac{1}{2}\trchb f_a- \frac{1}{2}\atrchb\dual f_a -2\om\fb_a+f^c\chibh_{ca} + (f\c\etab)\fb_a -\frac 1 2  \fb_a (f\c\ze) \\
&&+O((f, \fb)^3\Ga)
\eeaa
and hence
\beaa
\nab_{\la^{-1}e_4'}\fb &=& 2(\etab'-\etab)  -\frac{1}{2}(\trchb f -\atrchb\dual f) +2\om\fb-f\c\chibh  +E_6(f, \fb, \Ga)
\eeaa
where 
\beaa
E_6(f, \fb, \Ga) &=& -(f\c\etab)\fb +\frac 1 2(f\c\ze)\fb +O((f, \fb)^3\Ga)
\eeaa
as desired. This concludes the proof of Corollary \ref{cor:transportequationine4forchangeofframecoeff:simplecasefirst:PTcase}.
\end{proof}

The following is the analog of Corollary \ref{cor:transportequationine4forchangeofframecoeffinformFFbandlamba}. 
\begin{corollary}\lab{cor:transportequationine4forchangeofframecoeffinformFFbandlamba:PTcase}
Assume that we have
\beaa
\Xi'=0, \qquad \om'=0, \qquad \Hb'=-\frac{a\ov{q'}}{|q'|^2}\Jk'.
\eeaa
We introduce
\beaa
F:=f+i\dual f, \qquad \underline{F}:=\fb+i\dual \fb.
\eeaa
Then, we have
\beaa
\nab_{\la^{-1}e_4'}F+\frac{1}{2}\ov{\tr X} F+2\om F &=& -2\Xi -\chih\c F+E_1(f, \Ga),\\
\la^{-1}\nab_4'(\log\la) &=& 2\om+f\c(\ze-\etab)+E_2(f, \Ga),\\
\nab_{\la^{-1}e_4'}\underline{F} &=& -2\left(\frac{a\ov{q'}}{|q'|^2}\Jk'-\frac{a\ov{q}}{|q|^2}\Jk\right) -2\left(\Hb+\frac{a\ov{q}}{|q|^2}\Jk\right) -\frac{1}{2}\ov{\tr X}F +2\om\underline{F}\\
&& -F\c\chibh  +E_6(f, \fb, \Ga).
\eeaa

Moreover, introducing a complex valued scalar function $q$ satisfying $e_4(q)=1$, we have
\beaa
\nab_{\la^{-1}e_4'}(\ov{q}F) &=& -2\ov{q}\om F  -2\ov{q}\Xi +E_4(f, \Ga),
\eeaa
where 
\beaa
E_4(f, \Ga) &=& -\frac{1}{2}\ov{q}\left(\ov{\tr X} -\frac{2}{\ov{q}}\right)F -\ov{q}\chih\c F+\ov{q}E_1(f, \Ga)+f\c\nab(\ov{q})F+\frac{1}{4}|f|^2e_3(\ov{q})F.
\eeaa
\end{corollary}

\begin{proof}
The transport equations for $F$ and $\la$ have already been derived in Corollary \ref{cor:transportequationine4forchangeofframecoeffinformFFbandlamba}. We thus focus on the one for $\Fb$. Since $\xi'=0$, recall from Corollary \ref{cor:transportequationine4forchangeofframecoeff:simplecasefirst:PTcase}  that we have 
\beaa
\nab_{\la^{-1}e_4'}\fb &=& 2(\etab'-\etab)  -\frac{1}{2}(\trchb f -\atrchb\dual f) +2\om\fb-f\c\chibh  +E_6(f, \fb, \Ga).
\eeaa
In view of the definition of $F$ and $\underline{F}$, this yields
\beaa
\nab_{\la^{-1}e_4'}\underline{F} &=& 2(\Hb'-\Hb)  -\frac{1}{2}\ov{\tr X}F +2\om\underline{F} -F\c\chibh  +E_6(f, \fb, \Ga).
\eeaa
Plugging $\Hb'=-\frac{a\ov{q'}}{|q'|^2}\Jk'$, we obtain 
\beaa
\nab_{\la^{-1}e_4'}\underline{F} &=& -2\left(\frac{a\ov{q'}}{|q'|^2}\Jk'-\frac{a\ov{q}}{|q|^2}\Jk\right) -2\left(\Hb+\frac{a\ov{q}}{|q|^2}\Jk\right) -\frac{1}{2}\ov{\tr X}F +2\om\underline{F} -F\c\chibh \\
&& +E_6(f, \fb, \Ga)
\eeaa
as desired.
\end{proof}

\begin{remark}
In practice, we will integrate first the transport equations for $F$, then the one for $\la$, and finally the one for $\underline{F}$. Note that the transport equation for $\Fb$ in Corollary \ref{cor:transportequationine4forchangeofframecoeffinformFFbandlamba:PTcase} is at the same regularity level that the one for $F$ and $\la$, while the one for $\Fb$ in Corollary \ref{cor:transportequationine4forchangeofframecoeffinformFFbandlamba} looses one derivative. This is another manifestation of the fact that, unlike the PT frame, the PG frame exhibits a loss of derivative.   
\end{remark}

%%%%%%%%%%%%%%%%%%%%%%%%%%%%%%%%%%%%%%%%%%%%%

\subsection{Ingoing PT structures}
\lab{sec:ingoingPTstructures:chap2}  

%%%%%%%%%%%%%%%%%%%%%%%%%%%%%%%%%%%%%%%%%%%%%

\begin{definition}
\lab{definition:ingoingPT}
An ingoing  PT structure    $\{ (e_3, e_4, \HH), r, \th, \Jk\}$     on $\MM$ consists of  a null pair $(e_3, e_4)$, the induced horizontal structure   $\HH$,  functions $(r, \th)$, and a horizontal  1-form $\Jk$ such that the following hold true:
\begin{enumerate}
\item   $e_3$ is geodesic.

\item We have
\bea
e_3(r)=-1,\qquad    e_3(\th)=0, \qquad \nab_3 (\ov{q}\Jk)=0, \qquad  q= r+a i \cos \th.
\eea

\item We have 
\bea
H=\frac{aq}{|q|^2}\Jk.
\eea
\end{enumerate}
An extended  ingoing PT structure possesses,  in addition,  a function $\ub$ verifying $e_3(\ub)=0$. 
\end{definition}

\begin{definition}\lab{def:initialingoingPTdataset}
An ingoing PT  initial data set  consists of a hypersurface $\Si$ transversal to $e_3$ together with a null pair $(e_3, e_4)$, the induced horizontal structure $\HH$,  scalar   functions $(r, \th)$,  and a horizontal  1-form $\Jk$,  all defined on $\Si$.
\end{definition}

\begin{lemma}
\lab{lemma:constructioningoingPTframes}
Any ingoing PT  initial data set, as in Definition \ref{def:initialingoingPTdataset}, can be  locally extended to   an ingoing PT structure.
\end{lemma} 

\begin{proof}
Straightforward adaptation of the proof of Lemma \ref{lemma:constructionoutgoingPTframes}.
\end{proof}

%%%%%%%%%%%%%%%%%%%%%%%%%%%%%%%

\subsection{Null structure equations in an ingoing  PT frame}

%%%%%%%%%%%%%%%%%%%%%%%%%%%%%%%

\begin{proposition}\lab{Prop:NullStr-ingoingPTframe}
Consider an ingoing PT structure. Then,  the   equations in the $e_3$ direction for the Ricci coefficients  of the  ingoing PT frame take the form
\beaa
\nab_3\tr\Xb +\frac{1}{2}(\tr\Xb)^2 &=& -\frac{1}{2}\Xbh\c\ov{\Xbh},\\
\nab_3\Xbh+\Re(\tr\Xb) \Xbh&=& -\Ab,
\\
\nab_3\tr X +\frac{1}{2}\tr\Xb\tr X &=& \DD\c\left(\frac{a\ov{q}}{|q|^2}\ov{\Jk}\right)+\frac{a^2}{|q|^2}|\Jk|^2+2P-\frac{1}{2}\Xbh\c\ov{\Xh},\\
\nab_3\widehat{X} +\frac{1}{2}\tr\Xb\, \widehat{X}  &=& \frac{1}{2}\DD\hot\left(\frac{a q}{|q|^2}\Jk\right)  +\frac{1}{2}
\frac{a^2q^2}{|q|^4}\Jk\hot\Jk -\frac{1}{2}\ov{\tr X} \widehat{\Xb},
\\
\nab_3Z +\frac{1}{2}\tr\Xb Z &=& -\frac{1}{2}\tr\Xb \frac{a q}{|q|^2}\Jk -\frac{1}{2}\widehat{\Xb}\c\left(\ov{Z}+\frac{a \ov{q}}{|q|^2}\ov{\Jk}\right) -\Bb,\\
\nab_3\Hb &=&  -\frac{1}{2}\ov{\tr\Xb}\left(\Hb- \frac{a q}{|q|^2}\Jk\right) -\frac{1}{2}\Xbh\c\left(\ov{\Hb}-\frac{a\ov{q}}{|q|^2}\ov{\Jk}\right) +\Bb,\\
\nab_3\Xi&=&\nab_4\left(\frac{a q}{|q|^2}\Jk\right)  +\frac{1}{2}\ov{\tr X}\left(\frac{a q}{|q|^2}\Jk-\Hb\right) +\frac{1}{2}\Xh\c\left(\frac{a\ov{q}}{|q|^2}\ov{\Jk} -\ov{\Hb}\right) -B,
\\
\nab_3\om&=& \left(\Re\left(\frac{a q}{|q|^2}\Jk\right)-\etab\right)\c\ze -\Re\left(\frac{a q}{|q|^2}\Jk\right)\c\etab+   \rho.
\eeaa
\end{proposition}

\begin{proof}
These equations follow immediately from plugging the identities 
\beaa
\Xib=0, \qquad \omb=0,\qquad  H=\frac{a q}{|q|^2}\Jk. 
\eeaa
satisfied by an ingoing PT frame in  the general formulas of Proposition \ref{prop-nullstr:complex}.  
\end{proof}

%%%%%%%%%%%%%%%%%%%%%%%%%%%%%%%%

\subsection{Linearized quantities in   an ingoing PT frame}

%%%%%%%%%%%%%%%%%%%%%%%%%%%%%%%%%

Given an extended  ingoing PT  structure   $\{ (e_3, e_4, \HH), r, \th,  \ub, \Jk\}$   hold true:
\begin{enumerate}
\item We have
 \beaa
 \xib=\omb=0, \quad e_3(r)=-1, \quad e_3(\ub)= e_3(\th)=0, \qquad \nab_3(\ov{q}\Jk)=0.
 \eeaa
 In addition, we have
 \beaa
 H=\frac{a q}{|q|^2}\Jk.
 \eeaa
 \item The  quantities
 \beaa
\Xh, \quad \Xbh, \quad \Xi, \quad A, \quad B, \quad \Bb, \quad \Ab, \quad \DD r, \quad e_4(\cos\th), \quad \DD\hot\Jk,
\eeaa
  vanish  in      Kerr  and therefore are   small in perturbations.
\end{enumerate}
We renormalize below all other quantities, not vanishing in Kerr\footnote{Since  $H=\frac{aq}{|q|^2}\Jk$,  $H$ does not need to be included in Definition \ref{def:linearizedPT-ingoingcase:chap2}.}, by subtracting their  $\mbox{Kerr}(a, m) $  values.

\begin{definition}\lab{def:linearizedPT-ingoingcase:chap2}
We consider the following renormalizations, for given constants $(a, m)$,
\bea
\bsplit
\trXc &:= \tr X-\frac{2\ov{q}\De}{|q|^4}, \qquad\,\qquad     \trXbc := \tr\Xb+\frac{2}{\ov{q}},\\
\Zc &:= Z-\frac{aq}{|q|^2}\Jk,\qquad \qquad \quad\,\,\,
\Hbc := \Hb+\frac{a\ov{q}}{|q|^2}\Jk,\\
\omc & := \om  + \frac{1}{2}\pr_r\left(\frac{\De}{|q|^2} \right),\qquad \quad\, \Pc := P+\frac{2m}{q^3},
\end{split}
\eea
as well as 
\bea
\bsplit
\widecheck{e_4(r)} &:= e_4(r)-\frac{\Delta}{|q|^2},\qquad  \widecheck{\DD(\cos\th)} := \DD(\cos(\th)) -i\Jk,\\
\widecheck{\DD\ub} &:= \DD\ub -a\Jk,\qquad\qquad\,\,\,\, \widecheck{e_4(\ub)} := e_4(\ub) -\frac{2(r^2+a^2)}{|q|^2},
\end{split}
\eea
and
\bea
\bsplit
\widecheck{\ov{\DD}\c\Jk}& := \ov{\DD}\c\Jk-\frac{4i(r^2+a^2)\cos\th}{|q|^4}, \qquad \widecheck{\nab_4\Jk}:=\nab_4\Jk +\frac{\De \ov{q}}{|q|^4}\Jk.
\end{split}
\eea
\end{definition}

%%%%%%%%%%%%%%%%%%%%%%%%%%%%%%%%%%%%%%%%%%%%%
%%%%%%%%%%%%%%%%%%%%%%%%%%%%%%%%%%%%%%%%%%%%%

%%%%%%%%%%%%%%%%%%%%%%

\chapter{GCM admissible spacetimes}
\lab{chapter-GCM admissible.spacetimes}

%%%%%%%%%%%%%%%%%%%%%%

In this chapter we introduce the crucial notion of  general covariant modulated (GCM) admissible spacetimes, define  our main norms and state the main results.

%%%%%%%%%%%%%%%%%%%%%%

\section{Initial  data layer}
\lab{sec:defintionoftheinitialdatalayer}

%%%%%%%%%%%%%%%%%%%%%%

We consider a spacetime region $(\LL_0, \g)$, sketched  below in figure  \ref{fig0}, where
\begin{itemize}
\item The  Lorentzian spacetime metric $\g$ is close to the metric\footnote{We will write $\LL_0(a_0, m_0)$ whenever we need to emphasize the dependence on $(a_0, m_0)$.} of  $Kerr(a_0, m_0)$, $ |a_0|<m_0$,  in a suitable topology\footnote{This topology will be specified in our initial data layer assumptions, see \eqref{def:initialdatalayerassumptions} as well as section \ref{sec:initialdatalayernorm}.}. 

\item $\LL_0=\Lext\cup\Lint$. 

\item The intersection   $\Lext\cap\Lint$ is non trivial.
\end{itemize}

$(\LL_0, g)$ is called the initial data layer if it satisfies the properties in sections \ref{sec:definitioninitialdatalayer:boundaries}--\ref{sec:definitioninitialdatalayer:coordLextp} below.

\begin{figure}[h!]
\centering
\includegraphics[scale=0.9]{Kerr_2_april1.pdf}
\caption{The initial data layer $\LL_0$}
\label{fig0}
\end{figure}

%%%%%%%%%%%%%%%%%%%%%%%

\subsection{Boundaries}
\lab{sec:definitioninitialdatalayer:boundaries}

%%%%%%%%%%%%%%%%%%%%%%%

The future and past  boundaries  of $\LL_0$ are  given by
\beaa
\pr^+\LL_0&=& \AA_0\cup \BB_{(3,0)} \cup \underline{\BB}_{(3,0)},\\
\pr^- \LL_0 &=&\BB_{(0,0)}\cup \underline{\BB}_{(0,0)},
\eeaa
where
\begin{enumerate}
  \item  The  past non spacelike outgoing boundary of the far region $\Lext$ is denoted by $\BB_{(0,0)}$. 
 
    \item  The  past non   spacelike incoming   boundary of the near region  $\Lint$ is denoted by $\underline{\BB}_{(0,0)}$.  
 
 \item  $\Lext$ is unbounded in the future outgoing  directions.
   
  \item  The  future non   spacelike outgoing  boundary of the far region $\Lext$ is denoted by $\BB_{(3,0)}$. 
   
   \item  The    future non spacelike outgoing   boundary of the near region  $\Lint$ is denoted by $\underline{\BB}_{(3,0)}$.  
 
  \item  The  future spacelike boundary of the near region $\Lint$ is denoted by $\AA_0$.
 \end{enumerate}

\begin{remark}\lab{remark:preliminaryremarkonconstants:chap3}
Below, we make use of  the following constants:
\begin{itemize}
\item  the constants  $m_0>0$ and $|a_0|<m_0$  are the given  mass and the angular momentum of the  perturbed  Kerr solution,
\item $\ep_0>0$ is a small constant measuring the size of the perturbation of the initial data,
\item  $\deh>0$ and $\de_*>0$ are sufficiently small constants,
\item $r_0$ is a sufficiently large constant.    
\end{itemize}
See section \ref{sec:discussionofsmallnessconstantforthemaintheorem} for the specification of these constants.
\end{remark}

%%%%%%%%%%%%%%%%%%%%%%%%%%%%%%

\subsection{Foliations of $\LL_0$ and adapted null frames}

%%%%%%%%%%%%%%%%%%%%%%%%%%%%%%%

The spacetime $\LL_0=\Lext\cup \Lint$ is foliated as follows.

%%%%%%%%%%%%%%%%%%%%%%%%%%%%%%

\subsubsection{Foliations of $\Lext$}

%%%%%%%%%%%%%%%%%%%%%%%%%%%%%%%

The  far region $\Lext$ is covered by an outgoing PG structure, i.e. it is foliated by two functions $(u_{\idl}, \rextl)$ such that
\begin{itemize}
\item We have 
\beaa
L_0( u_{\idl})=0, \qquad L_0( \rextl)=1,
\eeaa
where the vectorfield $L_0$ is null  and satisfies $\D_{L_0}L_0=0$.

\item We denote by $(\,{}^{(ext)}(e_0)_3, \,{}^{(ext)}(e_0)_4, \,{}^{(ext)}(e_0)_1, \,{}^{(ext)}(e_0)_2)$ the null frame satisfying 
on $\Lext$
\beaa
&&{}^{(ext)}(e_0)_4=L_0,\qquad {}^{(ext)}(e_0)_1(\rextl)={}^{(ext)}(e_0)_2(\rextl)=0.
\eeaa

\item The    outgoing  future non spacelike boundary $\BB_{(3,0)}$ and the   past outgoing  non spacelike boundary  $\BB_{(0,0)}$ are given by
\beaa
\BB_{(3,0)}=\left\{u_{\idl}=3\right\}, \qquad \BB_{(0,0)}=\left\{u_{\idl}=0\right\}.
\eeaa

\item The foliation  by $u_{\idl}$ of $\Lext$ terminates at the timelike boundary 
$$\Big\{ \rextl= r_0-1\Big\},$$
where $r_0$ has been introduced in Remark \ref{remark:preliminaryremarkonconstants:chap3}.
\end{itemize}

%%%%%%%%%%%%%%%%%%%%%%%%%%%%%%

\subsubsection{Foliations of $\Lint$}

%%%%%%%%%%%%%%%%%%%%%%%%%%%%%%%

The near region $\Lint$ is foliated by two functions $(\ub_{\idl}, \rintl)$ such that  
\begin{itemize}
\item We have 
\beaa
\Lb_0( \ub_{\idl})=0, \qquad \Lb_0( \rintl)=-1,
\eeaa
where the vectorfield $\Lb_0$ is null  and satisfies $\D_{\Lb_0}\Lb_0=0$.

\item We denote by $(\,{}^{(int)}(e_0)_3, \,{}^{(int)}(e_0)_4, \,{}^{(int)}(e_0)_1, \,{}^{(int)}(e_0)_2)$ the null frame satisfying 
on $\Lint$
\beaa
{}^{(int)}(e_0)_3=\Lb_0, \qquad {}^{(int)}(e_0)_1(\rintl)={}^{(int)}(e_0)_2(\rintl)=0.
\eeaa

\item The $(\ub_{\idl},\rint)$ foliation is initialized on $\rextl=r_0$ as it will be made precise below. 

\item The foliation  by $\ub_{\idl}$, of $\Lint$ terminates at the space like boundary 
$$\AA_0=\left\{\rintl= \Big(m_0+\sqrt{m_0^2-a_0^2}\Big)(1-2\deh)\right\}$$
  where  $m_0$, $a_0$ and $\deh$ have been introduced in Remark \ref{remark:preliminaryremarkonconstants:chap3}.

\item The    future non spacelike ingoing   boundary $\underline{\BB}_{(3,0)}$ and the  past incoming  non  spacelike boundary $\underline{\BB}_{(0,0)}$ are given by
\beaa
\underline{\BB}_{(3,0)}=\left\{\ub_{\idl}=3\right\}, \qquad \underline{\BB}_{(0,0)}=\left\{\ub_{\idl}=0\right\}.
\eeaa

\item The foliation  by $\ub_{\idl}$ of $\Lint$ terminates at the time like boundary 
$$\Big\{\rintl= r_0+1\Big\}$$
where $r_0$ has been introduced in Remark \ref{remark:preliminaryremarkonconstants:chap3}.
\end{itemize}

%%%%%%%%%%%%%%%%%%%%%%%%%%%%%%%%%%%%%

\subsection{Definition of additional scalars and 1-forms in $\LL_0$}
\lab{sec:definitioninitialdatalayer:coordLextp}

%%%%%%%%%%%%%%%%%%%%%%%%%%%%%%%%%%%%%

We introduce the following scalars and 1-forms adapted to the  initial data layer $\LL_0$ defined above.
\begin{enumerate}
\item  In $\Lext $, we consider coordinates   $({}^{(ext)}\th_{\idl},  {}^{(ext)}\vphi_{\idl})$ satisfying 
 \bea
 \lab{eq:definethandvphi-ext:IDLayer}
 \,{}^{(ext)}(e_0)_4({}^{(ext)}\th_{\idl})=\,{}^{(ext)}(e_0)_4({}^{(ext)}\vphi_{\idl})=0.
 \eea
We also consider scalar functions ${}^{(ext)}J^{(p)}$, $p=0,+,-$ defined by
 \bea
 \bsplit
 \,{}^{(ext)}J^{(0)}&=\cos\left({}^{(ext)}\th_{\idl}\right), \quad {}^{(ext)}J^{(+)}=\sin\left({}^{(ext)}\th_{\idl}\right)\cos\left({}^{(ext)}\vphi_{\idl}\right), \\ 
 {}^{(ext)}J^{(-)}&=\sin\left({}^{(ext)}\th_{\idl}\right)\cos\left({}^{(ext)}\vphi_{\idl}\right),
 \end{split}
 \eea
and a complex 1-form $\,{}^{(ext)}\Jk$ satisfying 
\bea
\nab_{\,{}^{(ext)}(e_0)_4}{}^{(ext)}\Jk=0.
\eea

 \item  In $\Lint $, we consider coordinates   $({}^{(int)}\th_{\idl},  {}^{(int)}\vphi_{\idl})$ satisfying 
 \bea
 \lab{eq:definethandvphi-int:IDLayer}
 \,{}^{(int)}(e_0)_3({}^{(int)}\th_{\idl})=\,{}^{(int)}(e_0)_3({}^{(int)}\vphi_{\idl})=0.
 \eea
 We also consider scalar functions ${}^{(int)}J^{(p)}$, $p=0,+,-$ defined by
 \bea
 \bsplit
 \,{}^{(int)}J^{(0)}&=\cos\left({}^{(ext)}\th_{\idl}\right), \quad {}^{(int)}J^{(+)}=\sin\left({}^{(int)}\th_{\idl}\right)\cos\left({}^{(int)}\vphi_{\idl}\right), \\ 
 {}^{(int)}J^{(-)}&=\sin\left({}^{(int)}\th_{\idl}\right)\cos\left({}^{(int)}\vphi_{\idl}\right),
 \end{split}
 \eea
and a complex 1-form $\,{}^{(int)}\Jk$ satisfying 
\bea
\nab_{\,{}^{(int)}(e_0)_3}{}^{(int)}\Jk=0.
\eea
 \end{enumerate}

%%%%%%%%%%%%%%%%%%%%%%%%%%%%%%%

\subsection{Initializations  of the foliation on $\Lint$}

%%%%%%%%%%%%%%%%%%%%%%%%%%%%%%%

The $(\ub_{\idl},\rintl)$ foliation is initialized on $\rextl=r_0$ 
by setting,
\beaa
\ub_{\idl}&=&u_{\idl}, \qquad \rintl=\rextl
\eeaa
and, 
\beaa
\,{}^{(int)}(e_0)_4=\la_0\,{}^{(ext)}(e_0)_4,\quad  {}^{(int)}(e_0)_3=\la_0^{-1}\,{}^{(ext)}(e_0)_3,\quad 
\,{}^{(int)}(e_0)_b=\,{}^{(ext)}(e_0)_b, \,\, b=1, 2,
\eeaa
where
\beaa
{}^{(ext)}\la_0=\frac{r_0^2-2m_0r_0+a_0^2}{r_0^2+a_0^2(\cos({}^{(ext)}\th_{\idl}))^2}.
\eeaa
We also initialize the coordinates $({}^{(int)}\th_{\idl},  {}^{(int)}\vphi_{\idl})$ on $\rextl=r_0$ as follows
\beaa
{}^{(int)}\th_{\idl}={}^{(ext)}\th_{\idl}, \qquad {}^{(int)}\vphi_{\idl}={}^{(ext)}\vphi_{\idl}.
\eeaa

%%%%%%%%%%%%%%%%%%%%%%%%%%%%%%%%%%%%%%%%%%

\section{GCM admissible spacetimes}\lab{sec:defintioncanonicalspacetime}
\lab{section:GCMadmissible-spacetimes}

%%%%%%%%%%%%%%%%%%%%%%%%%%%%%%%%%%%%%%%%%

We consider a spacetime $(\MM, \g)$, sketched  below in figure  \ref{fig1}, where
\begin{itemize}
\item The Lorentzian spacetime metric $\g$ is close to Kerr in a suitable topology\footnote{This topology will be specified in our bootstrap assumptions, see section \ref{section:Bootstrap assumptions} as well as section \ref{section:main-norms}.}. 

\item $\MM=\Mext\cup\Mint\cup {}^{(top)}\MM$. 

\item $\TT=\Mext\cap\Mint$ is a time-like hyper-surface.

\item $\Mext\cap {}^{(top)}\MM$ and $\Mint\cap {}^{(top)}\MM$ are essentially  time-like hyper-surfaces.
\end{itemize}
$(\MM, \g)$ is called a general covariant modulated admissible (or shortly GCM-admissible) spacetime  if it is defined as in sections \ref{sec:boundariesadmissiblespacetime}--\ref{sec:admissibledefamandthetavphionMM}   below.

\begin{figure}[h!]
\centering
\includegraphics[scale=1.1]{Kerr_1.pdf}
\caption{The GCM admissible space-time $\mathcal{M}$}
\lab{fig1}
\end{figure}

%%%%%%%%%%%%%%%%%%%%%%%

\subsection{Boundaries}
\lab{sec:boundariesadmissiblespacetime}

%%%%%%%%%%%%%%%%%%%%%%%

The future and past  boundaries  of $\MM$ are  given by
\beaa
\pr^+\MM&=& \AA \cup  \,^{(top)} \Si \cup\Sigma_*,\\
\pr^- \MM&=& \BB_1\cup \underline{\BB}_1,
\eeaa
where,  see figure  \ref{fig1},
\begin{enumerate}
\item The past boundary  $\BB_1\cup \underline{\BB}_1$ is included in  the  initial data   layer\footnote{Recall that $\LL_0$, defined in section \ref{sec:defintionoftheinitialdatalayer}, is a spacetime region in which the metric on $\MM$ is specified to be  a small perturbation of the Kerr data.} $\mathcal{L}_0$. 
  
 \item  The  future spacelike boundary of the far region $\Mext$ is denoted by $\Sigma_*$.
  
  \item  The  future  spacelike     boundary of the  top region ${}^{(top)}\MM$ is denoted by $\,^{(top)} \Si$.

  \item  The  future spacelike boundary of the near region $\Mint$ is denoted by $\AA$.
  
  \item The time-like hyper-surface  $\TT$, separating  $\Mext$ from $\Mint$,  starts at $\underline{\BB}_1\cap\BB_1$ and 
   terminates at $\Mext\cap\Mint\cap {}^{(top)}\MM$.
 \end{enumerate}

%%%%%%%%%%%%%%%%%%%%%%%%%%%%

\subsection{Principal geodesic structures on  $\MM$}
\lab{sec:admissiblePGstructuresonMM}

%%%%%%%%%%%%%%%%%%%%%%%%%%%%

The spacetime region $\MM=\Mext\cup \Mint\cup {}^{(top)}\MM$ admits the following  PG structures.

%%%%%%%%%%%%%%%%%%%%%%%%%%%%%

\subsubsection{Principal geodesic structure on  $\Mext$}

%%%%%%%%%%%%%%%%%%%%%%%%%%%%%

The  far region $\Mext$ is endowed with an outgoing PG structure   defined by a scalar  functions $\rext$ and null frame
\beaa
(\,{}^{(ext)}e_3, \,{}^{(ext)}e_4, \,{}^{(ext)}e_1, \,{}^{(ext)}e_2)
\eeaa
with ${}^{(ext)}e_4$ null geodesic outgoing and such that we have on $\Mext$, see section \ref{subsection:Principaloutgoinggeodesicstructures},
\beaa
{}^{(ext)}e_4(\rext)=1,\qquad  {}^{(ext)}e_1(\rext)={}^{(ext)}e_2(\rext)=0.
\eeaa
Moreover
\begin{enumerate}
\item We introduce in addition a scalar function $u$ satisfying on $\Mext$
\beaa
{}^{(ext)}e_4(u)=0.
\eeaa

\item The $(u,\rext)$ foliation is initialized on $ \Sigma_*$  as it will be made precise below. 

\item  The outgoing non spacelike past boundary  $\BB_1$  corresponds precisely to  $u=1$.

\item The foliation  by $u$ of $\Mext$ terminates at the timelike boundary 
$$\TT=\Big\{\rext= r_0\Big\},$$
 where $r_0$ has been introduced in Remark \ref{remark:preliminaryremarkonconstants:chap3}.
\end{enumerate}

%%%%%%%%%%%%%%%%%%%%%%%%%%%%%

\subsubsection{Principal geodesic structure on  $\Mint$}

%%%%%%%%%%%%%%%%%%%%%%%%%%%%%

The  near  region $\Mint$ is endowed with an ingoing PG  structure   defined by a scalar  function $\rint$ and null frame
\beaa
(\,{}^{(int)}e_3, \,{}^{(int)}e_4, \,{}^{(int)}e_1, \,{}^{(int)}e_2)
\eeaa
with ${}^{(int)}e_3$ null geodesic ingoing and such that we have on $\Mint$, see section \ref{sec:Principalingoinggeodesicstructures},
\beaa
{}^{(int)}e_3(\rint)=-1, \qquad {}^{(int)}e_1(\rint)={}^{(int)}e_2(\rint)=0.
\eeaa
Moreover
\begin{enumerate}
\item We introduce in addition a scalar function $\ub$ satisfying on $\Mint$
\beaa
{}^{(int)}e_3(\ub)=0.
\eeaa

\item The $(\ub,\rint)$ foliation is initialized on $\TT$ as it will be made precise below. 

\item  The ingoing non spacelike past boundary  $\underline{\BB}_1$  corresponds precisely to  $\ub=1$.

\item The foliation  by $\ub$ of $\Mint$ terminates at the space like boundary 
$$\AA=\left\{\rint= \left(m_0+\sqrt{m_0^2-a_0^2}\right)(1-\deh)\right\}$$
where $m_0 $ and $\deh$ have been defined above.

\item We have $\ub=u_*$ on $\AA\cap\,^{(top)} \Si$.
\end{enumerate}

%%%%%%%%%%%%%%%%%%%%%%%%%%%%%%%%

\subsubsection{Principal geodesic structure on ${}^{(top)}\MM$}

%%%%%%%%%%%%%%%%%%%%%%%%%%%%%%%%

The region ${}^{(top)}\MM$  is endowed with an ingoing  PG structure   defined by a scalar function
 ${}^{(top)}r$  and a null frame
 \beaa(\,{}^{(top)}e_3, \,{}^{(top)}e_4, \,{}^{(top)}e_1, \,{}^{(top)}e_2)
 \eeaa
 where  ${}^{(top)}e_3$ is null ingoing geodesic  and such that we have on ${}^{(top)}\MM$
\beaa
{}^{(top)}e_3({}^{(top)}r)=-1, \qquad {}^{(top)}e_1({}^{(top)}r)={}^{(top)}e_2({}^{(top)}r)=0.
\eeaa
Moreover
\begin{enumerate}
\item We introduce in addition a scalar function $\ub$ satisfying on ${}^{(top)}\MM$
\beaa
{}^{(top)}e_3(\ub)=0.
\eeaa

\item $\Mint\cap {}^{(top)}\MM=\big\{\ub=u_*\big\}$.

\item The function $\ub$  is continuous across $\Mint\cap {}^{(top)}\MM$. 

\item   The function  ${}^{(top)}r$    extends the function $\rint$  continuously  across  $\Mint\cap {}^{(top)}\MM$. 

\item The null frame $(\,{}^{(top)}e_3, \,{}^{(top)}e_4, \,{}^{(top)}e_1, \,{}^{(top)}e_2)$ is a continuous extension of the frame \newline $(\,{}^{(int)}e_3, \,{}^{(int)}e_4, \,{}^{(int)}e_1, \,{}^{(int)}e_2)$ across $\Mint\cap {}^{(top)}\MM$.

\item The $(\ub, {}^{(top)}r)$ foliation is initialized on $\{u=u_*\}$ as it will be made precise below. 

\item The foliation  by $\ub$ of ${}^{(top)}\MM$ terminates at the  boundary 
\beaa
{}^{(top)}\Si &=& \Big\{\ub+\sigma_{top}\left({}^{(top)}r\right) = u_*\Big\},
\eeaa
where the function $\sigma_{top}$ may be chosen such that\footnote{The particular choice of the function $\sigma_{top}(r)$ satisfying the desired constraints is irrelevant for the proof. One could for example make the following suitable choice 
\beaa
\sigma_{top}(r)=\frac{m^2}{r}+c_1\quad\textrm{for}\quad r\leq r_0, \quad\,\, \sigma_{top}(r)=-2(r-r_0)-4m\log\left(\frac{r}{r_0}\right)-\frac{m^2}{r}+c_2\quad\textrm{for}\quad r\geq r_0+m,
\eeaa
pick the constants $c_1$ and $c_2$ such that ${}^{(top)}\Si$ starts at $\AA\cap\{\ub=u_*\}$ and terminates at $S_*$, and smoothly extend $\sigma_{top}(r)$ to $(r_0, r_0+m)$ so that ${}^{(top)}\Si $ is everywhere spacelike.  See section \ref{sec:proofofprop:propertiesoftauusefulfortheoremM8:chap9}, and in particular Proposition \ref{prop:propertiesoftauusefulfortheoremM8:appendix} and Lemma \ref{lemma:defintionoff1andf2andtheirbasicpropertiesusefullater}, for explicit computations in Kerr.} 
\begin{enumerate}
\item ${}^{(top)}\Si$ is spacelike,

\item ${}^{(top)}\Si$ starts at $\AA\cap\{\ub=u_*\}$ and terminates at $S_*$,

\item denoting by $r_{+,top}(\ub)$ and $r_{-,top}(\ub)$ respectively the maximum and the minimum of ${}^{(top)}r$ along a level hypersurface of $\ub$ in ${}^{(top)}\MM(r\geq r_0)$, there holds
\bea\lab{eq:upperboundrpubminusrmubonMtop:chap3}
0\leq r_{+,top}(\ub) - r_{-,top}(\ub)\leq 4m_0
\eea
uniformly in $\ub$. 
\end{enumerate}
\end{enumerate}

%%%%%%%%%%%%%%%%%%%%%%%%%%%%

\subsection{The GCM-PG data set on $\Si_*$}
\lab{sec:admissibleGMCPGdatasetonSigmastar}

%%%%%%%%%%%%%%%%%%%%%%%%%%%%

To initialize the PG structure of $\Mext$ on its future spacelike boundary $\Si_*$,  we assume  given a  GCM-PG data set  $\big(\Si_*, r, (e_1,e_2,e_3,e_4), f\big)$ as defined in Definition \ref{definition:GCM-datasets}, i.e.  
\begin{enumerate}
\item $r$ is a scalar function on $\Si_*$ whose level sets are 2-spheres foliating $\Si_*$, $(e_1, e_2, e_3, e_4)$ is a null frame defined on $\Si_*$, and $f$ is a 1-form tangent to  the spheres of the $r$-foliation of $\Si_*$.

\item $e_4$ is transversal\footnote{This is in fact automatic since $\Si_*$ is spacelike.} to $\Si_*$, $(e_1, e_2)$ are tangent to $\Si_*$, and $e_1(r)=e_2(r)=0$, so that $(e_1, e_2)$ are tangent to the leaves of the $r$-foliation.  

\item The sphere $S_*$ is the final  sphere on $\Si_*$. 

\item The frame $(e_1, e_2, e_3, e_4)$ satisfies the transversality conditions \eqref{eq:transversalityconditionsonSigma:firsttime}. 

\item We are given on $\Si_*$ coordinates $(\th, \vphi)$ and a basis of $\ell=1$ modes $J^{(p)}$, $p=0,+,-$, that are defined as follows
\begin{enumerate}
\item $(\th, \vphi)$ and $J^{(p)}$, $p=0,+,-$, are initialized  on $S_*$ as in  section \ref{section:canonical-coord-BPG}, 

\item $(\th, \vphi)$ and $J^{(p)}$, $p=0,+,-$, are propagated  to $\Si_*$ by $\nu(\th)=\nu(\vphi)=0$, and $\nu(J^{(p)})=0$, $p=0,+,-$, where $\nu=e_3+b_* e_4$  denotes the unique  vectorfield   on $\Si_*$  orthogonal to the $r$- foliation such that $\g(\nu, e_4)=-2$. 
\end{enumerate}

\item The GCM conditions \eqref{eq:Si^*-GCM1}-\eqref{eq:Si^*-GCM3} of Definition \ref{def:CanonicalGCM-hypersurface} are verified.

\item The constants  $(a, m)$ are specified  according to Definition \ref{define:am-onSi}.
\end{enumerate} 

In addition, we assume the following.
\begin{enumerate}
\item[(a)] Let $r_*$  such that $S_*=S(r_*)$. Then, $r$  is   monotonically increasing from   $r_*$. 

\item[(b)] The function $r$ verifies the dominance condition\footnote{A precise condition will be given later  in
\eqref{eq:behaviorofronS-star}.} on $S_*$
\bea
\lab{eq:behaviorofronSstar-rough}
r_* \gg u_*^{1+\dec},
\eea
where $u_*$ denotes the value of the function $u$ on $S_*$, with $u$ is specified below, see \eqref{eq:choiceofuonSigmastar}.  

\item[(c)] The 1-form $f$ on $\Si_*$ is chosen by
\bea\lab{eq:GCMPGdatasetfprescribedforadmissblespacetime}
f_1 =0, \qquad f_2 =\frac{a\sin\th}{r}, \quad\textrm{on}\quad S_*, \qquad \nab_\nu(rf)=0\quad\textrm{on}\quad\Si_*,
\eea
where $(e_1, e_2)$ are specified on $S_*$ by \eqref{eq:canonical-e1ande2onSstar:0}.
\end{enumerate}

%%%%%%%%%%%%%%%%%%%%%%%%%%%%%

\subsection{Definition of $(m,a)$ in $\MM$}
\lab{sec:definitionofamthetandvphiadmissible}

%%%%%%%%%%%%%%%%%%%%%%%%%%%%%

Given a GCM admissible spacetime $\MM$, we define the values $(a, m)$ associated to $\MM$   to be the constants    $(a, m)$  associated to $S_*$ according to Definition \ref{define:am-onSi}. Thus   each  GCM admissible spacetime $\MM$ is naturally  equipped  with  constants $(a, m)$. Note that  these constants   depend on $\MM$, i.e. two different GCM admissible spacetimes have in general different  constants $(a, m)$ associated to them.

%%%%%%%%%%%%%%%%%%%%%%%%%%%%%%%

\subsection{Initialization of the  PG structures  of $\MM$}
\lab{sec:initalizationadmissiblePGstructure}

%%%%%%%%%%%%%%%%%%%%%%%%%%%%%%%

%%%%%%%%%%%%%%%%%%%%%%%%%%%%%%%

\subsubsection{Initialization of the  PG structure  of $\Mext$}
\lab{section:InitializationofPG structure-Mext}

%%%%%%%%%%%%%%%%%%%%%%%%%%%%%%%

 The PG structure of $\Mext$   is  initialized on $\Si_*$  by the above  GCM-PG data set according to  Proposition \ref{Prop:Initialization of PG structures}, i.e. 
 \begin{enumerate}
\item    Along $\Si_*$,  the  restriction  of  the PG frame  $(\,{}^{(ext)}e_3, \,{}^{(ext)}e_4, \,{}^{(ext)}e_1, \,{}^{(ext)}e_2)$ of $\Mext$ is prescribed by the transformation formulas
\bea
 \bsplit
   \,{}^{(ext)}e_4 &=e_4 + f^b  e_b +\frac 1 4 |f|^2  e_3,\\
  \,{}^{(ext)}e_a&= \left(\de_a^b +\frac{1}{2}\fb_af^b\right) e_b +\frac 1 2  \fb_a  e_4 +\left(\frac 1 2 f_a +\frac{1}{8}|f|^2\fb_a\right)   e_3,\\
 \,{}^{(ext)}e_3&= \left(1+\frac{1}{2}f\c\fb  +\frac{1}{16} |f|^2  |\fb|^2\right) e_3 + \left(\fb^b+\frac 1 4 |\fb|^2f^b\right) e_b  + \frac 1 4 |\fb|^2 e_4,
 \end{split}
 \eea
where $(e_1,e_2,e_3,e_4)$ is the null frame of the GCM-PG data set of $\Si_*$, and the 1-forms $f$ and $\fb$ are given respectively by \eqref{eq:GCMPGdatasetfprescribedforadmissblespacetime} and \eqref{relations:laffb-to-primes-Si^*}, i.e.
\bea\lab{eq:1formfusedintheinitializationofthePGframeonSigmastar}
f_1 =0, \qquad f_2 =\frac{a\sin\th}{r}, \quad\textrm{on}\quad S_*, \qquad \nab_\nu(rf)=0\quad\textrm{on}\quad\Si_*,
\eea
and\footnote{The fact that $\fb$ is well defined will follow from our bootstrap assumptions.}
\bea
\fb = -\frac{(\nu(r)-b_*)}{1-\frac{1}{4}b_* |f|^2}f\quad\textrm{on}\quad\Si_*,
\eea
where the scalar functions $(r, \th, b_*)$ on $\Si_*$, the constant $a$, and the vectorfield $\nu$ tangent to $\Si_*$ are part of the above GCM-PG data set. 

\item   The functions $\rext$ is prescribed on $\Si_*$ by  $\rext=r$ where $r$ belongs to the GCM-PG data set of $\Si_*$. 

\item The function $u$ is prescribed on $\Si_*$ by
\bea\lab{eq:choiceofuonSigmastar}
u= c_*-\rext,
\eea 
where $c_*$ is a constant that will be fixed in Remark \ref{rmk:calibrationofu}. We then set $u_*$ to be the value of $u$ on $S_*$.
 \end{enumerate}

%%%%%%%%%%%%%%%%%%%%%%%%%%%%%%%

\subsubsection{Initialization of the  PG structure  of $\Mint$}

%%%%%%%%%%%%%%%%%%%%%%%%%%%%%%%

 \begin{enumerate}
\item The $(\ub,\rint)$ foliation is initialized on $\TT$  such that,
\bea
\ub=u, \qquad   \rint=\rext.
\eea
In particular $S(\ub, \rint)$ coincide on $\TT$ with $S(u, \rext)$.

\item The null frame $(\,{}^{(int)}e_3, \,{}^{(int)}e_4, \,{}^{(int)}e_1, \,{}^{(int)}e_2)$ is defined on $\TT$ by the following  renormalization,
\bea
\,{}^{(int)}e_4=  \la{}^{(ext)}e_4,\,\,\, {}^{(int)}e_3= \la^{-1}{}^{(ext)}e_3, \,\,\, {}^{(int)}e_a={}^{(ext)}e_a, \, a=1,2,\,\,\textrm{ on }\TT
\eea
where
\bea
\la=\,{}^{(ext)}\la=\frac{\,{}^{(ext)}\Delta}{|\,{}^{(ext)}q|^2}.
\eea
\end{enumerate}

%%%%%%%%%%%%%%%%%%%%%%%%%%%%%%%%%%

\subsubsection{Initialization of the  PG structure  of $^{(top)}\MM$}

%%%%%%%%%%%%%%%%%%%%%%%%%%%%%%%%%%

 \begin{enumerate}
 \item 
The $(\ub, {}^{(top)}r)$-foliation of ${}^{(top)}\MM$  is initialized on $\{u=u_*\}$  such that,
\bea\lab{eq:continuityofubandronuequalustarforinitializationMtop}
\ub=u+2\int_{r_0}^{\rext}\frac{{\tilde{r}}^2+a^2}{{\tilde{r}}^2-2m\tilde{r}+a^2}d\tilde{r}, \qquad  {}^{(top)}r=\rext.
\eea
In particular, the 2-spheres $S(\ub,  {}^{(top)}r)$ coincide on $\{u=u_*\}$ with $S(u, \rext)$.

\begin{remark}
The initialization of $\ub$ in \eqref{eq:continuityofubandronuequalustarforinitializationMtop} agrees with the relation between $\ub$ and $u$ in Kerr, see section \ref{sec:outgoingPGcoordinatesinKerr:chap2} and section \ref{sec:KerrvaluesiningoingPGstructure:chap2}. Note also that $\ub=u$ at $r=r_0$. 
\end{remark}

\item Moreover, the null frame $(\,{}^{(top)}e_3, \,{}^{(top)}e_4, \,{}^{(top)}e_1, \,{}^{(top)}e_2)$ is prescribed on $\{u=u_*\}$  by the transformation formulas
\bea
 \bsplit
   \,{}^{(top)}e_4 &=\la\, {}^{(ext)}e_4,\\
  \,{}^{(top)}e_a&= {}^{(ext)}e_b +\frac 1 2  \fb_a \, {}^{(ext)}e_4,\\
 \,{}^{(top)}e_3&= \la^{-1}\left({}^{(ext)}e_3 + \fb^b\, {}^{(ext)}e_b  + \frac 1 4 |\fb|^2\, {}^{(ext)}e_4\right),
 \end{split}
 \eea
where
\bea
\la=\,{}^{(ext)}\la=\frac{\,{}^{(ext)}\Delta}{|\,{}^{(ext)}q|^2}, \qquad \fb=\,{}^{(ext)}h \widecheck{\,{}^{(ext)}e_3(r)}\,{}^{(ext)}\nab(u),
\eea
with the scalar function $\,{}^{(ext)}h$ given by\footnote{The fact that $\,{}^{(ext)}h$ is well defined will follow from our bootstrap assumptions.} 
\bea\lab{eq:choiceofhininitializationingoingPGstructureMtop:chap3}
h &=& \frac{4}{e_3(u)+\sqrt{(e_3(u))^2+4|\nab u|^2\widecheck{e_3(r)}}},
\eea
where we dropped $\,{}^{(ext)}$ for simplicity in the formula for $h$. 
\end{enumerate}

\begin{remark}
The choice \eqref{eq:choiceofhininitializationingoingPGstructureMtop:chap3} above ensures $\,{}^{(top)}\nab(\,{}^{(top)}r)=0$ on $\{u=u_*\}$, see Lemma \ref{lemma:necessaryidentitiesonuequalustarforthecontrolofMtopnuequlustar}, which is a necessary condition for the ingoing foliation of $\Mtop$ to be a PG structure.
\end{remark}

%%%%%%%%%%%%%%%%%%%%%%%%%%%%%%%

\subsection{Definition of coordinates $(\th, \vphi)$ in $\MM$}
\lab{sec:admissibledefamandthetavphionMM}

%%%%%%%%%%%%%%%%%%%%%%%%%%%%%%%

We introduce the following $(\th, \vphi)$ coordinates adapted to the  PG structures defined above.
\begin{enumerate}
\item  In $\Mext $, we initialize  $({}^{(ext)}\th,  {}^{(ext)}\vphi)$ on $\Si_*$ by
\bea
{}^{(ext)}\th=\th, \qquad {}^{(ext)}\vphi=\vphi,
\eea
where $(\th, \vphi)$ are associated to the GCM-PG data set as above, and propagate it to $\Mext$ by
 \bea
 \lab{eq:definethandvphi-ext}
 \,{}^{(ext)}e_4({}^{(ext)}\th)=\,{}^{(ext)}e_4({}^{(ext)}\vphi)=0.
 \eea
 
 \item  In $\Mint $, we initialize  $({}^{(int)}\th,  {}^{(int)}\vphi)$ on $\TT$ by
 \bea
{}^{(int)}\th={}^{(ext)}\th, \qquad {}^{(int)}\vphi={}^{(ext)}\vphi,
\eea
 and propagate it to $\Mint$ by
 \bea
 \lab{eq:definethandvphi-int}
 \,{}^{(int)}e_3({}^{(int)}\th)=\,{}^{(int)}e_3({}^{(int)}\vphi)=0.
 \eea
 
 \item In ${}^{(top)}\MM$, we initialize  ${}^{(top)}\th$ on $\{u=u_*\}$ by
 \bea
{}^{(top)}\th={}^{(ext)}\th.
\eea
and propagate it to ${}^{(top)}\MM$ by 
 \bea
 \lab{eq:definethandvphi-top}
 \,{}^{(top)}e_3({}^{(top)}\th)= 0.
 \eea
\end{enumerate}

\begin{remark}
Note that we do not need to define a coordinate ${}^{(top)}\vphi$. Indeed, the coordinate $\vphi$, in the various regions of $\MM$, plays an auxiliary role in the proof, and  is in fact only needed to control coordinates system on $\MM$ in the limit $u_*\to +\infty$ where the region $\Mtop$ actually disappears.
\end{remark}

%%%%%%%%%%%%%%%%%%%%%%%%%

\section{Main norms}\label{section:main-norms}

%%%%%%%%%%%%%%%%%%%%%%%%%

We define our main norms  on a  given GCM admissible spacetime $\MM$. The norms involve in particular the constants  $(a, m)$ specified in section \ref{sec:definitionofamthetandvphiadmissible}.

%%%%%%%%%%%%%%%%%%%%%%%%%%%%%%%%%%%%%

\subsection{Main norms on $\Si_*$}
\label{section:main-normsSigmastar}

%%%%%%%%%%%%%%%%%%%%%%%%%%%%%%%%%%%%%

All quantities appearing in this section are defined relative to the  GCM frame of $\Si_*$, the scalar function $r$ of  $\Si_*$, and the constant $m$  introduced in section \ref{sec:admissibleGMCPGdatasetonSigmastar}. 

 We introduce the function $u$ as in \eqref{eq:choiceofuonSigmastar}, i.e. we have
\bea
u=c_*-r
\eea
where $c_*$ is a fixed constant. Also, recall that we  have the following properties for angular derivatives $r$ on $\Si_*$
\bea
\lab{eq:relations-ru-onSi*}
e_1 (r)= e_2(r)=0.
\eea
Finally, recall that the frame $(e_1, e_2, e_3, e_4)$ of $\Si_*$ satisfies the transversality conditions \eqref{eq:transversalityconditionsonSigma:firsttime}, to which we add transversality conditions for $e_4(r)$ and $e_4(u)$
\bea\lab{eq:tranversalityconditionforthefoliationonSi*}
\xi=0, \qquad \om=0, \qquad \etab=-\ze, \qquad e_4(r)=1, \qquad e_4(u)=0, \qquad\textrm{ on }\quad\Si_*.
\eea

\begin{remark}
The transversality conditions \eqref{eq:tranversalityconditionforthefoliationonSi*} allow us to make sense of all the Ricci coefficients for the frame of $\Si_*$ and all first order derivatives of $r$ and $u$.
\end{remark}

\begin{definition}\lab{def:renormalizedquantitiesGCMfoliationofSigma*}
For the GCM foliation of     the boundary   $\Si_*$ of $ \MM$,              the linearized  quantities are defined as follows
 \bea
 \ba{cccccc}
\widecheck{\trch} &:=&\displaystyle \trch-\frac{2}{r}, \qquad &
\widecheck{\trchb} &:=&\displaystyle \trchb+\frac{2\Up}{r},\\
\\
\ombc &:=&\displaystyle \omb-\frac{m}{r^2},\qquad &\rhoc &:=&\displaystyle \rho +\frac{2m}{r^3},\\
\\
\widecheck{e_3(r)} &:=& e_3(r)+\Up,\qquad &\widecheck{e_3(u)} &:=&\displaystyle e_3(u) -2,\\[3mm]
\widecheck{b_*} &:=& \displaystyle b_*+1+\frac{2m}{r},
\ea
\eea
where $\Up = 1-\frac{2m}{r}$.
\end{definition}

With these normalizations we define   the sets $\Ga_g^*, \Ga^*_b$ as follows.
\begin{definition}
\lab{Definition:linearizedquantitiesSi*}
The set of all linearized quantities on $\Si_*$ is of the form $\Ga^*_g\cup \Ga^*_b$  defined as follows.
\begin{enumerate}
\item The set $\Ga^*_g$ contains
\bea
\bsplit
\Ga^*_g&:=\Big\{ \widecheck{\trch}, \, \chih, \, \ze,\,\widecheck{\trchb}, \,  r\a,\,  r\b, \,  r\big(\rhoc, \, \rhod\big)\Big\}.
\end{split}
\eea
\item The set $\Ga^*_b$ contains 
\bea
\Ga^*_b&:=&\Big\{ \eta,\,  \chibh, \, \widecheck{\omb}, \, \xib,\, r\bb,\, \aa,  \,   r^{-1}\widecheck{e_3(r)}\,   r^{-1}\widecheck{e_3(u)}, \,   r^{-1}\widecheck{b_*}\Big\}.
\eea
\end{enumerate}
\end{definition}

To define higher order derivatives norms on $\Si_*$, we use the weighted derivative operators $\dk_*$ tangential to $\Si_*$ defined as follows
\bea
\dk_* := \{\nab_\nu, \dkb\}. 
\eea

%%%%%%%%%%%%%%%%%%%%%%%%%%
 
\subsubsection{Boundedness  norms  on  $\Si_*$}

%%%%%%%%%%%%%%%%%%%%%%%%%%

For any $k\geq 0$, we introduce the following  norms
\bea
\,^{*}\mathfrak{B}_k := \sup_{\Si_*} \Big\{ r^2| \dk_*^{\leq k}\Ga^*_g|  + r| \dk_*^{\leq k}\Ga^*_b|\Big\}.
\eea

%%%%%%%%%%%%%%%%%%%%%%%%%%%%%%%%

\subsubsection{Decay norms on $\Si_*$}

%%%%%%%%%%%%%%%%%%%%%%%%%%%%%%%

Let $\dec>0$ a small constant to be specified later.  For any $k\geq 1$, we introduce the following  norms
\bea
\,^{*}\mathfrak{D}_k := \sup_{\Si_*}\Big\{    r^2u^{\frac 1 2 +\dec} | \dk_*^{\leq k}\Ga^*_g|+r^2u^{1 +\dec} | \dk_*^{\leq k-1}\nab_\nu\Ga^*_g|  + r u^{1+\dec} | \dk_*^{\leq k}\Ga^*_b|   \Big\}.
\eea

%%%%%%%%%%%%%%%%%%%%%%%%%%%%%%%%%%%%%

\subsection{Main norms in $\Mext $}\label{section:main-normsextregion}

%%%%%%%%%%%%%%%%%%%%%%%%%%%%%%%%%%%%%

All quantities appearing in this section are defined relative to the  outgoing PG structure of  $\Mext$.  As there is  no danger of confusion we will drop  the       prefixes $^{(ext)}$ in what follows.

%%%%%%%%%%%%%%%%%%%%%%%%%%%%%%%%%%%%%%
 
\subsubsection{Definition of complex horizontal 1-form $\Jk$ on $\Mext$}

%%%%%%%%%%%%%%%%%%%%%%%%%%%%%%%%%%%%%%

Recall that the PG frame of    $\Mext$  satisfies 
\beaa
\Xi=\om=0,\qquad \Hb=-Z
\eeaa
as well as,
\beaa
e_4(r)=1,\quad   e_4(u)=e_4(\th)= e_4(\vphi)=0, \quad e_1(r)=e_2(r)=0.
\eeaa

To define linearized quantities in $\Mext$ as in Definition \ref{def:renormalizationofallnonsmallquantitiesinPGstructurebyKerrvalue}, we need to use the complex horizontal 1-form $\Jk$ introduced in section \ref{section:auxilliaryformsMext}. Recall that $\Jk$ is defined in that section up to an initialization on a hypersurface transversal to $e_4$. We provide such initialization on the hypersurface $\Si_*$ below\footnote{Note that $q= r+ i a\cos\th$ is defined with respect to the constant $a$ of $\Mext$.}.

\begin{definition}[Definition of $\Jk$ in $\Mext$]\lab{def:definitionofJkonMext} 
We define the complex horizontal 1-form $\Jk$  on $\Mext$ by\footnote{Note that $\Jk$ is related on $\Si_*$ to the 1-form $f$ introduced in \eqref{eq:1formfusedintheinitializationofthePGframeonSigmastar} as follows
\beaa
\Jk=\frac{r}{|q|}(f+i\dual f)\quad\textrm{on}\quad\Si_*.
\eeaa}
\bea\lab{eq:basicdefinitionofJkonSstarSigmastarMextstar}
\bsplit
&\Jk_1=\frac{i\sin\th}{|q|}, \quad \Jk_2=\frac{\sin\th}{|q|},\quad\textrm{on}\quad S_*, \\
&\nab_\nu(|q|\Jk)=0\quad\textrm{on}\quad \Si_*,\\
& \nab_4\Jk = -\frac{1}{q}\Jk\quad\textrm{on}\quad \Mext.
\end{split}
\eea 
\end{definition}

The complex horizontal 1-form $\Jk$ satisfies the following lemma.
\begin{lemma}
The complex horizontal 1-form $\Jk$ of Definition \ref{def:definitionofJkonMext} satisfies on $\Mext$ the identities \eqref{eq:basicpropertiesofJkinMextusedtolinearize}, i.e.
\bea\lab{eq:basicpropertiesofJkinMextusedtolinearize:bis}
\nab_4\Jk = -\frac{1}{q}\Jk, \qquad \dual\Jk = -i\Jk,\qquad \Jk\c\ov{\Jk}=\frac{2(\sin\th)^2}{|q|^2}.
\eea
\end{lemma}

\begin{remark}
Since $\Jk$ satisfies  \eqref{eq:basicpropertiesofJkinMextusedtolinearize}, it can be used to define the linearized quantities in $\Mext$ as in Definition \ref{def:renormalizationofallnonsmallquantitiesinPGstructurebyKerrvalue}.
\end{remark}

\begin{proof}
The first identity holds true by Definition \ref{def:definitionofJkonMext}. The two other identities hold true on $S_*$ by Definition \ref{def:definitionofJkonMext}., and are then immediately transported to $\Si_*$ using the fact that $\nu(|q|\Jk)=0$ and $\nu(\th)=0$ on $\Si_*$. Finally, they are transported to $\Mext$ using $\nab_4\Jk = -q^{-1}\Jk$ and Lemma \ref{Lemma:Propagation-Jk}. 
\end{proof}

Recall\footnote{We stress  the fact that the linearized quantities appearing in the definitions of  $\Ga_g$ and $\Ga_b$  are done with respect to the constants $(a, m)$ of $\Mext$.} the set of quantities  $\Ga_g$, $\Ga_b$, see Definition \ref{definition.Ga_gGa_b}. Finally, we use the        
   weighted  derivative operator  
\bea
\dk:=\{\nab_3, r\nab_4, \dkb\}.
\eea

%%%%%%%%%%%%%%%%%%%%%%%%%%
 
\subsubsection{Boundedness  norms  in  $\Mext$}

%%%%%%%%%%%%%%%%%%%%%%%%%%

Let $\dt>0$ a small constant to be specified later.  For any $k\geq 1$, we introduce the following  norms
\bea\lab{equation:defboudednessnormsMext:chap3}
\bsplit
\,{}^{(ext)}\mathfrak{B}_k &:= \sup_{\Mext} \Big\{ r^2| \dk^{\leq k}\Ga_g|  + r| \dk^{\leq k}\Ga_b|+r^{\frac{7}{2}+\frac{\dt}{2}}\big(|\dk^{\leq k}A|+|\dk^{\leq k}B|\big)\\
&+r^{\frac{9}{2}+\dec}|\dk^{\le k-1}\nab_3 A|+r^4|\dk^{\le k-1} \nab_3B|\Big\}+\left(\int_{\Si_*}|\dk^{\le k}\Hc|^2\right)^{\frac{1}{2}}.
\end{split}
\eea

%%%%%%%%%%%%%%%%%%%%%%%%%%%%%%%%

\subsubsection{Decay norms in $\Mext$}

%%%%%%%%%%%%%%%%%%%%%%%%%%%%%%%

Let $\dec>0$ a small constant to be specified later. We define for $k\ge 1$,
\bea\lab{equation:defdecaynormsMext:chap3}
\bsplit
\,{}^{(ext)}\mathfrak{D}_k &:= \sup_{\Mext}\Big(ru^{1+\dec}+r^2u^{\frac{1}{2}+\dec}\Big)|\dk^{\le k} \Ga_g|+\sup_{\Mext}ru^{1+\dec}|\dk^{\le k}\Ga_b|\\
&+\sup_{\Mext}r^4u^{\frac{1}{2}+\dec}\Big(|\dk^{\le k-1}\nab_3 A|+|\dk^{\le k-1} \nab_3B|\Big)\\
&+\sup_{\Mext}r^2 u^{1+\dec} |\dk^{\le k-1} \nab_3\Ga_g|+\left(\int_{\Sigma_*}u^{2+2\dec}|\dk^{\le k}\Hc|^2\right)^{\frac{1}{2}}.
\end{split}
\eea

\begin{remark}
The integral bootstrap assumption on $\Sigma_*$ for $\Hc$ will only be needed in the proof of Proposition \ref{prop:constructionsecondframeinMext} and recovered in Proposition  \ref{Prop.Flux-bb-vthb-eta-xib}. In fact, other components satisfy an analog integral estimate  on $\Sigma_*$: this is the case of $\Xbh$, $\Xib$ and $r\Bb$, see  Proposition  \ref{Prop.Flux-bb-vthb-eta-xib}. But $\Hc$ is the only component for which we need to make this type of  bootstrap assumption. 
\end{remark}

%%%%%%%%%%%%%%%%%%%%%%%%%%%%%%%%%%%%%%%%%%%%%%%%

\subsection{Main norms in $\Mint$}\label{section:main-normsintregion}

%%%%%%%%%%%%%%%%%%%%%%%%%%%%%%%%%%%%%%%%%%%%%%%%

All quantities appearing in this section are defined relative to the  ingoing PG structure of  $\Mint$.  As there is  no danger of confusion we will drop  the       prefixes $^{(int)}$ in what follows.

%%%%%%%%%%%%%%%%%%%%%%%%%%%%%%%%%%%%%%
 
\subsubsection{Definition of complex horizontal 1-form $\Jk$ on $\Mint$}

%%%%%%%%%%%%%%%%%%%%%%%%%%%%%%%%%%%%%%

To define linearized quantities in $\Mint$ as in Definition \ref{def:renormalizationofallnonsmallquantitiesinPGstructurebyKerrvalue:ingoingcase}, we need to use the complex horizontal 1-form $\Jk$ introduced in section \ref{sec:Principalingoinggeodesicstructures}. Recall that $\Jk$ is defined in that section up to an initialization on a hypersurface transversal to $e_3$. We provide such initialization on the hypersurface $\TT$ below.

\begin{definition}[Definition of $\Jk$ in $\Mint$]\lab{def:definitionofJkonMint} 
We define the complex horizontal 1-form $\Jk$  on $\Mint$ by
\bea\lab{eq:basicdefinitionofJkonSstarSigmastarMintstar}
\bsplit
&\Jk={}^{(ext)}\Jk\quad\textrm{on}\quad \TT, \\
& \nab_3\Jk = \frac{1}{\ov{q}}\Jk\quad\textrm{on}\quad \Mint.
\end{split}
\eea 
\end{definition}

Recall also, see Definition \ref{definition.Ga_gGa_b:ingoingcase},   the   set of quantities $\Ga_g$, $\Ga_b$ for ingoing PG structures.

%%%%%%%%%%%%%%%%%%%%%%%%%%
 
\subsubsection{Boundedness  norms  in  $\Mint$}

%%%%%%%%%%%%%%%%%%%%%%%%%%

For any $k\geq 1$, we introduce the following  norms
\bea
\bsplit
\,{}^{(int)}\mathfrak{B}_k &:= \sup_{\Mint} \Big\{| \dk^{\leq k}\Ga_g|  + | \dk^{\leq k}\Ga_b|\Big\}.
\end{split}
\eea

%%%%%%%%%%%%%%%%%%%%%%%%%%%%%%%%

\subsubsection{Decay norms in $\Mint$}

%%%%%%%%%%%%%%%%%%%%%%%%%%%%%%%

We define for $k\ge 1$,
\bea
\bsplit
\,{}^{(int)}\mathfrak{D}_k &:= \sup_{\Mint}\ub^{1+\dec}\Big(|\dk^{\le k} \Ga_g|+|\dk^{\le k}\Ga_b|\Big).
\end{split}
\eea

%%%%%%%%%%%%%%%%%%%%%%%%%%%%%%%%%%%%%%%%%%%%%%%%

\subsection{Main norms in $\Mtop$ }\label{section:main-normstopregion}

%%%%%%%%%%%%%%%%%%%%%%%%%%%%%%%%%%%%%%%%%%%%%%%%

All quantities appearing in this section are defined relative to the  ingoing PG structure of  $\Mtop$.  As there is  no danger of confusion we will drop  the       prefixes $^{(top)}$ in what follows.

%%%%%%%%%%%%%%%%%%%%%%%%%%%%%%%%%%%%%%
 
\subsubsection{Definition of complex horizontal 1-form $\Jk$ on $\Mtop$}

%%%%%%%%%%%%%%%%%%%%%%%%%%%%%%%%%%%%%%

To define linearized quantities in $\Mtop$ as in Definition \ref{def:renormalizationofallnonsmallquantitiesinPGstructurebyKerrvalue:ingoingcase}, we need to use the complex horizontal 1-form $\Jk$ introduced in section \ref{sec:Principalingoinggeodesicstructures}. Recall that $\Jk$ is defined in that section up to an initialization on a hypersurface transversal to $e_3$. We provide such initialization on the hypersurface $\{u=u_*\}$ below.

\begin{definition}[Definition of $\Jk$ in $\Mtop$]\lab{def:definitionofJkonMtop} 
We define the complex horizontal 1-form $\Jk$  on $\Mtop$ by
\bea\lab{eq:basicdefinitionofJkonSstarSigmastarMtopstar}
\bsplit
&\Jk={}^{(ext)}\Jk\quad\textrm{on}\quad \{u=u_*\}, \\
& \nab_3\Jk = \frac{1}{\ov{q}}\Jk\quad\textrm{on}\quad \Mtop.
\end{split}
\eea 
\end{definition}

\begin{remark}\lab{rmk:noJpmandnoJkpminMtop:chap3}
We do not introduce the scalar function $J^{(\pm)}$ and the complex 1-forms $\Jk_\pm$ in $\Mtop$. Thus, the quantities $\Ga_g$, $\Ga_b$ in $\Mtop$ correspond to the ones in Definition \ref{definition.Ga_gGa_b:ingoingcase} where all linearized quantities based on $J^{(\pm)}$ and  $\Jk_\pm$ have been removed. 
\end{remark}

%%%%%%%%%%%%%%%%%%%%%%%%%%
 
\subsubsection{Boundedness  norms  in  $\Mtop$}

%%%%%%%%%%%%%%%%%%%%%%%%%%

For any $k\geq 1$, we introduce the following  norms
\bea
\bsplit
\,{}^{(top)}\mathfrak{B}_k &:= \sup_{\Mtop} \Big\{ r^2| \dk^{\leq k}\Ga_g|  + r| \dk^{\leq k}\Ga_b|+r^{\frac{7}{2}+\frac{\dt}{2}}\big(|\dk^{\leq k}A|+|\dk^{\leq k}B|\big)\\
&\qquad\qquad\qquad\qquad\,\,+r^{\frac{9}{2}+\dec}|\dk^{\le k-1}\nab_3 A|+r^4|\dk^{\le k-1} \nab_3B|\Big\}.
\end{split}
\eea

%%%%%%%%%%%%%%%%%%%%%%%%%%%%%%%%

\subsubsection{Decay norms in $\Mtop$}

%%%%%%%%%%%%%%%%%%%%%%%%%%%%%%%

To describe the decay norms in $\Mtop(r\geq r_0)$, we introduce the scalar function ${}^{(top)}u$ as follows
\bea
{}^{(top)}u &:=& \ub -2\int_{r_0}^r\frac{{\tilde{r}}^2+a^2}{{\tilde{r}}^2-2m\tilde{r}+a^2}d\tilde{r}.
\eea

\begin{remark}
Note in view of \eqref{eq:continuityofubandronuequalustarforinitializationMtop} that ${}^{(top)}u=u$ on $u=u_*$.
\end{remark}

We define for $k\ge 1$,
\bea
\bsplit
\,{}^{(top)}\mathfrak{D}_k &:=\,{}^{(top)}\mathfrak{D}_k^{\leq r_0}  + \,{}^{(top)}\mathfrak{D}_k^{\geq r_0},\\
\,{}^{(top)}\mathfrak{D}^{\leq r_0}_k &:= \sup_{\Mtop(r\leq r_0)}\ub^{1+\dec}\Big(|\dk^{\le k} \Ga_g|+|\dk^{\le k}\Ga_b|\Big),\\
\,{}^{(top)}\mathfrak{D}^{\geq r_0}_k &:= \sup_{\Mtop(r\geq r_0)}\Big(r({}^{(top)}u)^{1+\dec}+r^2({}^{(top)}u)^{\frac{1}{2}+\dec}\Big)|\dk^{\le k} \Ga_g|\\
&+\sup_{\Mtop(r\geq r_0)}r({}^{(top)}u)^{1+\dec}|\dk^{\le k}\Ga_b|+\sup_{\Mtop(r\geq r_0)}r^2 ({}^{(top)}u)^{1+\dec} |\dk^{\le k-1} \nab_3\Ga_g|\\
&+\sup_{\Mtop(r\geq r_0)}r^4({}^{(top)}u)^{\frac{1}{2}+\dec}\Big(|\dk^{\le k-1}\nab_3 A|+|\dk^{\le k-1} \nab_3B|\Big).
\end{split}
\eea

%%%%%%%%%%%%%%%%%%%%%%%%%%%%%%%%%%%%%

\subsection{Combined norms}\lab{sec:definitionofconcatenatednorm}

%%%%%%%%%%%%%%%%%%%%%%%%%%%%%%%%%%%%%

We define the following norms $\MM$ by combining our above norms on $\Si_*$, $\Mext$, $\Mint$ and ${}^{(top)}\MM$
\beaa
\Nk^{(Sup)}_k & :=& \,^{*}\mathfrak{B}_k+\,{}^{(ext)}\Bk_k+\,{}^{(int)}\Bk_k+\,{}^{(top)}\Bk_k,\\
\Nk^{(Dec)}_k & :=& \,^{*}\Dk_k+\,{}^{(ext)}\Dk_k+\,{}^{(int)}\Dk_k +\,{}^{(top)}\Dk_k. 
\eeaa

%%%%%%%%%%%%%%%%%%%%%%%%%%%%%%%%%%%%%%%

\subsection{Initial layer norm}\lab{sec:initialdatalayernorm}

%%%%%%%%%%%%%%%%%%%%%%%%%%%%%%%%%%%%%%%

Recall the notations of section \ref{sec:defintionoftheinitialdatalayer} concerning the initial data layer $\LL_0=\LL_0(a_0, m_0)$. Recall that the constants $m_0>0$ and $|a_0|< m_0$ are the mass and angular momentum of the initial Kerr  spacetime relative to which our initial perturbation is measured. We define   the initial layer norm to be\footnote{Recall that the initial data layer foliations satisfy $\etab+\ze=0$, as well as $\xi=\om=0$ on $\Lext$, and $\eta=\ze$ as well as $\xib=\omb=0$ on $\Lint$.}
\beaa
\Ik_k &:=&\, ^{(ext)} \Ik_k+\, ^{(int)} \Ik_k+\Ik_k',
\eeaa 
where
\beaa
\bsplit
\, ^{(ext)} \Ik_0 &:= \sup_{\Lext} \Big\{ r^2|\Ga_g|  + r|\Ga_b|+r^{\frac{7}{2}+\frac{\dt}{2}}\big(|A|+|B|\big)\Big\},
\end{split}
\eeaa
\beaa
\, ^{(int)} \Ik_0&:=&  \sup_{\Lint}\Big\{| \Ga_g|  + | \Ga_b|\Big\},
\eeaa
\beaa
\Ik_0' &:=&  \sup_{\Lint\cap\Lext} \left(  |f|+|\fb| +|\log(\la_0^{-1}\la)|   \right),\qquad \la_0= {}^{(ext)}\la_0=1-\frac{2m_0}{\rextl},
\eeaa
and 
\begin{itemize}
\item for ${}^{(ext)}\Ik_0$, $(\Ga_g, \Ga_b)$ is given by Definition \ref{definition.Ga_gGa_b} for the outgoing PG structure of $\Lext$,

\item  for ${}^{(int)}\Ik_0$, $(\Ga_g, \Ga_b)$ is given by Definition \ref{definition.Ga_gGa_b:ingoingcase} for the ingoing PG structure of $\Lint$,

\item for $\Ik_0'$, $(f, \fb, \la)$ denote the transition coefficients of Lemma \ref{Lemma:Generalframetransf}  from the frame of the outgoing PG structure of  $\Lext$ to the frame of the ingoing PG structure of $\Lint$  in the region $\Lint\cap\Lext$,
\end{itemize}
with $\Ik_k$ the corresponding  higher derivative norms obtained by replacing   each   component by $\dk^{\le k}$ of it. 

\begin{remark}  
 Note that  in the definition of   $\, ^{(ext)} \Ik_k$ we allow a higher power  of $r$   in front $\a$, $\b$ and their derivatives than what    it is consistent with the results  of  \cite{Ch-Kl} and \cite{KlNi}. The additional $r^{\dt}$ power, for $\dt$ small,   is consistent instead with the result  of   \cite{KlNi2}. See also Remark \ref{remark:resultsKl-Ni}.
\end{remark}

%%%%%%%%%%%%%%%%%%%%%%%%%%%%%%%%%%%%%%%%

\section{Main theorem}

%%%%%%%%%%%%%%%%%%%%%%%%%%%%%%%%%%%%%%%%

%%%%%%%%%%%%%%%%%%%%%%%%%%%%%%%%%%%%%%%%%%%%%%%%%%%%

\subsection{Smallness constants}\lab{sec:discussionofsmallnessconstantforthemaintheorem}

%%%%%%%%%%%%%%%%%%%%%%%%%%%%%%%%%%%%%%%%%%%%%%%%%%%%

Before stating our main theorems, we first introduce the following constants that will be involved in its statement.
\begin{itemize}
\item The constants $m_0>0$ and $|a_0|<m_0$ are the mass and the angular momentum of the Kerr solution relative to which our initial perturbation is measured. 

\item The integer $k_{large}$ which corresponds to the maximum number of derivatives of the solution.

\item The size of the initial data layer norm is measured by $\ep_0>0$. 

\item The size of the bootstrap assumption norms are measured by $\ep>0$.

\item $r_0>0$ is tied to the definition of $\TT$, i.e. $\TT=\{r=r_0\}$. 

\item $\deh>0$ measures the width of the region $|r-m_0-\sqrt{m_0^2-a_0^2}|\leq 2m_0\deh$ where the redshift estimate holds.

\item $\dec$ is tied to decay estimates  in $u$, $\ub$  for $\check{\Gamma}$ and $\Rc$.

\item $\dt$ is involved in the $r$-power of the supremum estimates for high derivatives of $\a$ and $\b$.

\item $\de_*$ is involved in the behavior of $r$ on $S_*$, see \eqref{eq:behaviorofronS-star} below.
\end{itemize}

  In what follows, $m_0$ and $a_0$ are fixed constants with $0\leq |a_0|<m_0$, $\deh$, $\dt$, and $\dec$ are fixed, sufficiently small, universal constants, and $r_0$ and $k_{large}$ are  fixed, sufficiently large, universal constant, chosen such that
\bea\lab{eq:constraintsonthemainsmallconstantsepanddelta}
\bsplit
& 0<\deh,\,\, \dec, \,\,\dt, \,\,\de_* \ll \min\{m_0, 1\},\qquad \dt> 2\dec, \\ 
& r_0\gg \max\{m_0,1\},\qquad k_{large}\gg \frac{1}{\dec}.
\end{split}
\eea
Then, $\ep$ and $\ep_0$ are chosen such that
\bea\lab{eq:constraintsonthemainsmallconstantsepanddelta:bis}
0<\ep_0, \ep\ll \min\left\{\deh, \dec, \dt, \de_*, \frac{1}{r_0}, \frac{1}{k_{large}}, m_0-|a_0|, 1\right\},
\eea
\bea\lab{eq:constraintofep0wrta0}
\ep_0, \ep\ll |a_0|\quad \textrm{ in the case }a_0\neq 0,
\eea
and
\bea\lab{eq:constraintbetweenepep*ep0}
\ep=\ep_0^{\frac{2}{3}}.
\eea

\begin{remark}
Note that we may always assume \eqref{eq:constraintofep0wrta0}, even if $0<|a_0|\les\ep_0$. Indeed, in that case, an initial data layer assumption of the type
\beaa
\Ik_k\les\ep_0,
\eeaa
see section \ref{sec:initialdatalayernorm} for the definition of the initial data layer norm $\Ik_k$, remains true by setting $a_0=0$.
\end{remark}

Using the definition of $\ep_0$, we  can  now make  precise the  condition  \eqref{eq:behaviorofronSstar-rough} of $r$  on $S_*$
\bea
\lab{eq:behaviorofronS-star}
r_*= \de_*  \ep_0^{-1} u_*^{1+\dec},
\eea
where we recall that $r_*$ and $u_*$ denote respectively the value of $r$ and $u$ on $S_*$. 

Also, we introduce the integer $k_{small}$ which corresponds to the  number of derivatives  for which the solution satisfies decay estimates. It is related to $k_{large}$ by
\bea\lab{eq:choiceksmallmaintheorem}
k_{small}=\left \lfloor\frac 1 2 k_{large}\right \rfloor +1.
\eea

From now on, in the rest of the paper, $\lesssim$ means bounded by a constant depending only on geometric universal constants (such as Sobolev embeddings, elliptic estimates,...) as well as the constants 
$$m_0,\, a_0, \, \deh,\, \dec, \,\dt, \, \de_*, \, r_0, \, k_{large}$$
\textit{but not on} $\ep$ and $\ep_0$.

%%%%%%%%%%%%%%%%%%%%%%%%%%%%%%%%

\subsection{Admissible  future null  complete spacetimes}

%%%%%%%%%%%%%%%%%%%%%%%%%%%%%%%%

We introduce in this section the notion of admissible  future null  complete spacetimes which corresponds 
formally to the limit as $u_*\to +\infty$ of the GCM admissible spacetimes of section \ref{section:GCMadmissible-spacetimes}.

\begin{definition}\lab{def:ep0kadmissibleLL0region}
We say that an initial data layer $\LL_0=\LL_0(a_0, m_0)$, defined  as in  section \ref{sec:defintionoftheinitialdatalayer}, is admissible   if it lies in the future of an asymptotically flat initial data set,  supported  on a spacelike hypersurface  $\Si_0$,   of ADM 
 mass $ m_0  $ and angular momentum  $a_0$. For a representation  see the lower part of the Figure  \ref{fig:penrosediagramfuturecomplete}.

 In addition,    we say that   $\LL_0$ is   $(\ep_0, k)$-admissible    if  it verifies the bounds\footnote{One expects in principal the weaker bound $\Ik_{k}\leq \ep_0$. See Remark \ref{rmk:whyonearthdoweneedep0squareandnotep0forinitaldata} for an explanation of the need of the stronger bound $\Ik_{k}\leq \ep_0^2$.}
\beaa
\Ik_{k}\leq \ep_0^2,
\eeaa
with $\Ik_k $ defined  as in  section \ref{sec:initialdatalayernorm}.
\end{definition}

\begin{remark}
\lab{remark:resultsKl-Ni}
The results in \cite{KlNi},  \cite{KlNi2} and  \cite{Ca-Ni}   identify large classes of initial data sets  on $\Si_0$  which generate  $(\ep_0, k)$-admissible  initial data layers.
\end{remark}

We define  a development of an  admissible initial data layer   $\LL_0$ as follows.

\begin{definition}
\lab{definition:developmentofLL_0}
We say that a spacetime is   a future development of  an  admissible initial data layer  $\LL_0$  if  it  is in fact a future development of   the   initial data set supported on    $\Si_0$.
\end{definition}

\begin{figure}[h!]
\centering
\includegraphics[scale=1]{kerr_4.pdf}
\caption{ Penrose diagram of  an admissible future complete spacetime}
\label{fig:penrosediagramfuturecomplete}
\end{figure}

\begin{definition}
\lab{Definition:admissible.futurenullcomplete} 
We say that an asymptotically flat Einstein vacuum spacetime  $\MM$, as in Figure \ref{fig:penrosediagramfuturecomplete},  is  admissible  future null  complete,   if it verifies the following properties.
\begin{itemize}
\item  It is   a  future  development   of  an admissible   initial data layer  set $\LL_0(a_0, m_0)$, in the sense  of  Definition \ref{definition:developmentofLL_0}.

 \item The future null infinity  $\II^{+}$ of   $\MM$  is  complete.  The other future boundary of $\MM$ is given by the spacelike hypersurface $\AA$ and belongs to the complement of $\JJ^{-}(\II^{+})$.

 \item $\MM=\Mext\cup\Mint$.
 
 \item $\TT=\Mext\cap\Mint$ is  time-like.
 
 \item  $\Mext$ comes equipped with  an outgoing PG structure, as in Definition \ref{def:PGstructure}, given by the function $r$  and  null frame\footnote{ Recall that  $\HH(r)=0$, $ e_4(r)=1$.}   $(e_3, e_4, \HH)$,  and a function $u$   verifying $e_4(u)=0$.
 
  \item Each   $S=S(u, r)$ sphere  in $\Mext$   comes equipped  with an adapted, integrable,  frame $(e_3', e_4', \HH')$. The passage from the  PG  frame to the integrable  one  is obtained by the transformation formulas
   \eqref{General-frametransformation}  with  parameters $ (f, \fb, \la)$ given by \eqref{def:transition-functs:ffbla}.
   
  \item The timelike surface $\TT$ is given  by  $\{r=r_0\}$, and $\Mint$ comes  equipped  with an ingoing PG structure  and  function $\ub$   initialized on $\TT$ as in  section \ref{sec:initalizationadmissiblePGstructure}.
\end{itemize}
\end{definition}

\begin{definition}
Given   an  admissible  future null  complete,  spacetime $\MM$  as in Definition \ref{Definition:admissible.futurenullcomplete}, and constants $(a_\infty, m_\infty)$, $|a_\infty|<m_\infty$, we define  the norms\footnote{Note that $\MM$ does not contain a region $\Mtop$, so that there are no norms corresponding to $\Mtop$ in the definition of $\Nk^{(Sup)}_{k_{large}}(a_\infty, m_\infty )$ and $\Nk^{(Dec)}_{k_{small}}(a_\infty, m_\infty )$.} $\Nk^{(Sup)}_{k_{large}}(a_\infty, m_\infty )$, $\Nk^{(Dec)}_{k_{small}}(a_\infty, m_\infty )$ as in  section \ref{sec:definitionofconcatenatednorm},  with  the constants $(a, m)$, which appear  in  the definition of the  linearized quantities replaced by 
 $(a_\infty, m_\infty)$.
\end{definition}

%%%%%%%%%%%%%%%%%%%%%%%%%%%%%%%%%%%%%%%%%%%%%%%%%%%%

\subsection{Statement of the main theorem} 
\lab{section:MainTheorem}

%%%%%%%%%%%%%%%%%%%%%%%%%%%%%%%%%%%%%%%%%%%%%%%%%%%%

We are now ready to give the following precise version of our main theorem. 

\begin{thmmain}[Main theorem, version 2]
Let $\LL=\LL(a_0, m_0) $  be an     $(\ep_0, k_{large}+10)$-admissible initial data layer as in Definition \ref{def:ep0kadmissibleLL0region},  with      $|a_0|/ m_0$ sufficiently small,  $ k_{large}$ sufficiently 
large,    and   $\ep_0>0$  sufficiently small. In particular, we assume\footnote{One expects in principal the weaker bound $\Ik_{k_{large}+10}\leq \ep_0$. See Remark \ref{rmk:whyonearthdoweneedep0squareandnotep0forinitaldata} for an explanation of the need of the stronger bound $\Ik_{k_{large}+10}\leq \ep_0^2$.}
\bea\lab{def:initialdatalayerassumptions}
\Ik_{k_{large}+10}\leq \ep_0^2.
\eea
 Then   $\LL=\LL(a_0, m_0) $   possesses an    admissible future complete       development        $\MM_\infty$  as in Definition  \ref{Definition:admissible.futurenullcomplete}.
Moreover:
\begin{enumerate}
\item   There exist constants $(a_\infty, m_\infty)$,  $|a_\infty| \ll m_\infty$,         such that the following  estimates hold true,      relative to the norms     $\Nk^{(Sup)}_{k_{large}}=\Nk^{(Sup)}_{k_{large}}(a_\infty, m_\infty ), \, \Nk^{(Dec)}_{k_{small}}=\Nk^{(Dec)}_{k_{small}}(a_\infty, m_\infty ) $ defined above, 
\bea\lab{def:bootstrapasumptionsglobalnorms}
\Nk^{(Sup)}_{k_{large}}+\Nk^{(Dec)}_{k_{small}} +| a_\infty- a_0| +|m_\infty- m_0| \le C\ep_0
\eea
where $C$ is a universal constant sufficiently large  and $k_{small}=\left \lfloor\frac 1 2 k_{large}\right \rfloor +1$.

\item  The space $\MM_\infty$ is  a  limit of finite  GCM admissible spacetimes\footnote{More precisely they are the limits of the   GCM family  $\UU(u_*)$, as  in Definition  \ref{definition:UU(u_*)}.  Note also that $\Mtop$ disappears in the limit.}.
\end{enumerate}

 The estimates \eqref{def:bootstrapasumptionsglobalnorms} imply, in particular:
\begin{itemize}
\item[-] On $\Mext_\infty$, we have
\beaa
 |\a|, |\b|&\les& \min\left\{ \frac{\ep_0}{r^3(u+2r)^{\frac{1}{2}+\dec}},\,   \frac{\ep_0}{r^2(u+2r)^{1+\dec}}\right\},\\
 |\rhoc|, |\rhodc| &\les & \min\left\{ \frac{\ep_0}{r^3 u ^{\frac{1}{2}+\dec}},\,   \frac{\ep_0}{r^2 u^{1+\dec}}\right\},\\
  |\bb| &\les& \frac{\ep_0}{r^2u^{1+\dec}},\\
   |\aa| &\les&\frac{\ep_0}{ru^{1+\dec}},
\eeaa
and 
\beaa
|\trchc|, |\widecheck{\atrch}| &\les& \frac{\ep_0}{r^2u^{1+\dec}},\\
|\chih|, |\widecheck{\ze}|, |\trchbc|, |\widecheck{\atrchb}| &\les& \min\left\{ \frac{\ep_0}{r^2u^{\frac{1}{2}+\dec}}\frac{\ep_0}{ru^{1+\dec}}\right\},\\
|\widecheck{\eta}|, |\chibh|, |\ombc|, |\xib| & \les& \frac{\ep_0}{ru^{1+\dec}}.
\eeaa

\item[-] On $\Mint_\infty$   we have, for any linearized quantity  $\widecheck{\psi}$,
\beaa
| \widecheck{\psi}|  \les \frac{\ep_0}{\ub^{1+\dec}}.
\eeaa
\end{itemize}
Note that analog statements of the above estimates also hold for $\dk^k$ derivatives with $k\leq k_{small}$.

 Moreover the following other statements  hold true. 
\begin{enumerate}
\item  Let  $m_H(u,r) $ denote the Hawking mass     adapted to the spheres $S=S(u, r)$ of $\Mext_\infty$, i.e.
\beaa
m_H(u, r) &=&   \sqrt{\frac{|S(u, r)|}{4\pi}}    \left(1+\frac{1}{16\pi}\int_{S(u, r)} \trch'\trchb'\right).
\eeaa
where  $ \trch', \trchb'$ are  calculated with respect to the integrable\footnote{ See Definition \ref{Definition:admissible.futurenullcomplete}. The integrable frame is such that $(e_1', e_2')$ is tangent to $S(u,r)$.  Recall, on the other hand, that the outgoing PG frame  of $\Mext$ is non integrable.} frame $(e_1', e_2', e_3', e_4')$ of $\Mext$. 
Then:
\begin{itemize}
\item[-]
The Bondi mass exists and is given by
\beaa
M_B(u):=\lim_{r\to \infty}  m_H(u,r).
\eeaa

\item[-] $M_B(u) $ has a limit as $u\to \infty$  and
\bea
\lim_{u\to \infty} M_B(u) = m_\infty,
\eea
i.e. $m_\infty$ coincides with the final Bondi mass.
\item[-] The Bondi mass law formula \eqref{eq:Bondi-masslaw}  holds true. In particular $m_\infty< m_0$.
\end{itemize}
 
 \item We define the  quasi-local angular momentum  for   a  sphere $S(u, r)$  to be the triplet 
 \beaa
 {\mathfrak j} _{\ell=1, p}(u, r) :=    \frac{1}{|S(u,r)|}   \int _{S(u, r)} ( \curl' \b')\Jp, \qquad p=0, +,-.
 \eeaa
 with $\curl' \b'$ defined relative the the integrable frame of $\Mext$.  Then
 \begin{itemize}
 \item[-]  The   triplet    ${\mathfrak j} _{\ell=1, p}(u, r)$ has a limit  as $r \to \infty$ at fixed $u$ given by 
 \bea
 \JJ_{\ell=1, p}(u) &=&\lim_{r\to \infty} {\mathfrak j} _{\ell=1, p}(u, r).
 \eea
 \item[-]  The  triplet   $ \JJ_{\ell=1, p}(u)  $ has a limit as $u\to \infty$ and 
\beaa
\lim_{u\to \infty} \JJ_{\ell=1,0}(u)= 2 a_\infty m_\infty, \qquad   \lim_{u\to \infty} \JJ_{\ell=1,\pm}(u) =0.
\eeaa
\end{itemize}

\item $\Mext$ is covered by three regular coordinates patches:
\begin{itemize}
\item[-] in the $(u,r,\th,\vphi)$ coordinates system, we have, for $\frac{\pi}{4}<\th<\frac{3\pi}{4}$, 
\beaa
\g &=& \g_{a_\infty,m_\infty}+\Big(du, dr, rd\th, r\sin\th d\vphi\Big)^2O\left(\frac{\ep_0}{u^{1+\dec}}\right),
\eeaa

\item[-] in the $(u,r,x^1,x^2)$ coordinates system,  with $x^1=J^{(+)}$ and $x^2=J^{(-)}$, we have, for $0\leq\th<\frac{\pi}{3}$ and for   $\frac{2\pi}{3}<\th\leq\pi$, 
\beaa
\g &=& \g_{a_\infty,m_\infty}+\Big(du, dr, rdx^1, rdx^2\Big)^2O\left(\frac{\ep_0}{u^{1+\dec}}\right),
\eeaa
\end{itemize}
where in each case, $\g_{a_\infty,m_\infty}$ denotes the Kerr metric expressed in the corresponding coordinates system of Kerr, see Lemma \ref{lemma:urthetavphicoordinatesinKerrchap2} and Lemma \ref{lemma:urthetaJplusJminuscoordinatesinKerrchap2}.

\item $\Mint$ is covered by three regular coordinates patches:
\begin{itemize}
\item[-] in the $(\ub,r,\th,\vphi)$ coordinates system, we have, for $\frac{\pi}{4}<\th<\frac{3\pi}{4}$, 
\beaa
\g &=& \g_{a_\infty,m_\infty}+\Big(d\ub, dr, rd\th, r\sin\th d\vphi\Big)^2O\left(\frac{\ep_0}{\ub^{1+\dec}}\right),
\eeaa

\item[-] in the $(\ub,r,x^1,x^2)$ coordinates system, with $x^1=J^{(+)}$ and $x^2=J^{(-)}$, we have, for $0\leq\th<\frac{\pi}{3}$ and for   $\frac{2\pi}{3}<\th\leq\pi$, 
\beaa
\g &=& \g_{a_\infty,m_\infty}+\Big(d\ub, dr, rdx^1, rdx^2\Big)^2O\left(\frac{\ep_0}{\ub^{1+\dec}}\right),
\eeaa
\end{itemize}
where in each case, $\g_{a_\infty,m_\infty}$ denotes the Kerr metric expressed in the corresponding coordinates system of Kerr, i.e. the analog for ingoing PG structures of Lemma \ref{lemma:urthetavphicoordinatesinKerrchap2} and Lemma \ref{lemma:urthetaJplusJminuscoordinatesinKerrchap2}.
\end{enumerate}
\end{thmmain}

\begin{remark}
By far the  most  demanding   part  in the proof of  the theorem is to establish  the  crucial estimates \eqref{def:bootstrapasumptionsglobalnorms}.   The  items  1--5  are important corollaries   of those   estimates, to be treated as conclusions in  section \ref{section:conclusion-mainThm}.  In  section  \ref{section:ProofMainThm}, we discuss the main  intermediary  results, Theorems  M0--M8,   in the proof of the estimates \eqref{def:bootstrapasumptionsglobalnorms}.
\end{remark}

\begin{remark}\lab{rmk:whyonearthdoweneedep0squareandnotep0forinitaldata}
The assumption \eqref{def:initialdatalayerassumptions} for the initial data layer norm $\Ik_{k_{large}+10}$ can in fact be replaced by its following weaker analog
\bea\lab{def:initialdatalayerassumptions:waekeranalog}
\Ik_{k_{large}+10}\leq \ep_0, \qquad {}^{(ext)}\Ik_{3}\leq\ep_0^2,
\eea
see also Remark \ref{rmk:whyonearthdoweneedep0squareandnotep0forinitaldata:thereturn}. The estimate ${}^{(ext)}\Ik_{3}\leq\ep_0^2$ in \eqref{def:initialdatalayerassumptions:waekeranalog} is then used at two instances in the proof of the main theorem:
\begin{itemize}
\item the proof of Theorem M0, see section \ref{sec:proofofTheoremM0}, where is it used to infer the $O(\ep_0)$ control of the curvature components of the PG structure of $\Mext$ and $\Mint$ respectively on the past boundaries $\BB_1$ and $\underline{\BB}_1$ of $\MM$,

\item the proof of Theorem M6, see section \ref{sec:proofofTheoremM6}, where it is used to exhibit in the initial data layer a spacetime $\MM$ satisfying our bootstrap assumptions, hence initiating the continuity argument in the proof of the main theorem. 
\end{itemize}
The precise places where we need the estimate ${}^{(ext)}\Ik_{3}\leq\ep_0^2$ in \eqref{def:initialdatalayerassumptions:waekeranalog}  are outlined in Remark \ref{remark:needforstrongerinitialdatalayernorminproofofThmM0:1} for the proof of Theorem M0, and Remark \ref{rem:thisistheonlyplacewhereweneedstrongercontrolofinitialdata} for the proof of Theorem M6. 
\end{remark}

%%%%%%%%%%%%%%%%%%%%%%%%%%%%%%%%%%%%%%%%

\section{Main bootstrap assumptions}
\lab{section:Bootstrap assumptions}
\label{section:bootstrap}

%%%%%%%%%%%%%%%%%%%%%%%%%%%%%%%%%%%%%%%%

 Given a GCM admissible  spacetime  $\MM$,  as defined in section \ref{section:GCMadmissible-spacetimes},  we  assume that the combined  norms\footnote{Recall that the  norms   are defined 
  with respect  to linearized quantities   which involve the  constants  $(a, m)$ specified in section \ref{sec:definitionofamthetandvphiadmissible}.} $\Nk^{(Sup)}_k$ and $\Nk^{(Dec)}_k$,  defined    in   section \ref{section:main-norms},  verifies   the following bounds

{\it{\bf BA-B} (Bootstrap Assumptions on  $r$-weighted  sup norms)
\bea\lab{def:bootstrapasumptionsglobalnormonenergie}
\Nk^{(Sup)}_{k_{large}} +|m-m_0|+|a-a_0|\le \ep,
\eea

{\bf BA-D} (Bootstrap Assumptions on decay)
\bea\lab{def:bootstrapasumptionsglobalnormondecay}
\Nk^{(Dec)}_{k_{small}}\le \ep.
\eea}
We shall often refer to them in the text  as ${\bf BA_\ep}$.

%%%%%%%%%%%%%%%%%%%%%%%%%%%%%

\section{The global null frame}\lab{section-globaleframe}

%%%%%%%%%%%%%%%%%%%%%%%%%%%%%

%%%%%%%%%%%%%%%%%%%%%%%%%%%%%

\subsection{The quantity $\qf$}
\lab{sec:defintionofquantityqf}

%%%%%%%%%%%%%%%%%%%%%%%%%%%%%

We will need to control a complex horizontal symmetric traceless 2-tensor $\qf$ of the form 
\beaa
\qf &=& q \ov{q}^{3}\Big((\nab_3-2\omb)(\nab_3-4\omb)A + C_1(\nab_3-4\omb)A  +C_2 A\Big), 
\eeaa
for specific complex scalar  functions\footnote{See \eqref{finalchoicefordefinition-C} for the choice of $C_1$ and $C_2$.} $C_1$, $C_2$. $\qf$ has been introduced in \cite{GKS1}, where a Regge Wheeler type equation is derived for it. Based on this Regge Wheeler type equation, we derive estimates for $\qf$ and $A$ in Theorem M1, see section  \ref{sec:mainintermediateresults:chap3}.

As in \cite{KS},  the quantity $\qf$ has to be defined   relative to a  global frame, smooth in the entire  region  $\MM$,  in which the quantity\footnote{This was $\eta$ in \cite{KS}.}  $\Hc$ behaves like $\Ga_g$.  To achieve this  in a way which does not affect the other $\Ga_g$ quantities  we proceed  in two steps:

\begin{enumerate}
\item We construct in section \ref{subsection:constructionsecondframeinMext} a second frame on $\Mext$ for which $\Hc$ belongs to $\Ga_g$.

\item We use this second  frame of $\Mext$  to construct in section \ref{section-globaleframe} 
a global frame  in $\MM$ for which $\Hc$ belongs to $\Ga_g$.
\end{enumerate} 
It is this global frame  which is used to analyze  the decay properties of $\qf$ in Theorem M1.

%%%%%%%%%%%%%%%%%%%%%%%%%%%%%

\subsection{Construction of a second frame in $\Mext$}
\label{subsection:constructionsecondframeinMext}

%%%%%%%%%%%%%%%%%%%%%%%%%%%%%

We denote in this section by $(e_1, e_2, e_3, e_4)$  the outgoing PG frame of the region $\Mext$, and by $(e_1', e_2', e_3', e_4')$ the second frame of $\Mext$ constructed below. The primed frame of $\Mext$ is obtained by performing  a transformation of the form
\bea
\lab{secondframe-Mext1}
\begin{split}
e_4'&= e_4 + f^a e_a +\frac 1 4 |f|^2  e_3,\\
e_a' &= e_a +\frac 1 2 f_a e_3,\quad a=1,2,\\
e_3'&=  e_3,
\end{split}
\eea
such that 
\bea
\lab{secondframe-Mext2}
\widecheck{H'}\in \Ga_g'
\eea
where $\Ga_g', \Ga_b'$ are defined below.

\begin{definition}\lab{definition.Ga_gGa_b:withprimseforsecondframeMext}
Let $(e_1', e_2', e_3', e_4')$ the second frame of $\Mext$ constructed in Proposition \ref{prop:constructionsecondframeinMext} below. With respect to that frame, we introduce the notations $\Ga_g', \Ga_b'$ as follows:
\begin{itemize}
\item the linearized quantities for the frame $(e_1', e_2', e_3', e_4')$ are defined  in the same way  as Definition \ref{def:renormalizationofallnonsmallquantitiesinPGstructurebyKerrvalue} for the outgoing PG frame of $\Mext$, with respect to the   coordinates $(r, \th)$ and the complex 1-form $\Jk$ of the PG structure\footnote{Thus, for example, $\widecheck{\tr X}' =\tr X'-\frac{2}{q}$, $\Hc'=H'-\frac{aq}{|q|^2}\Jk$, $\widecheck{e_3'(r)} = e_3'(r)-  \frac{\De}{|q|^2}$, $\widecheck{\DD' J^{(0)}}=\DD' J^{(0)}-i\Jk$, $\widecheck{\nab_3'\Jk}=\nab_3'\Jk-\frac{\De q}{|q|^4}\Jk$,...},  

\item in addition, we introduce the following linearized quantities which are trivial for an outgoing PG structure\footnote{Except $\Hbc$ which satisfies instead $\Hb=-Z$.}
\beaa
\Hbc':=\Hb'+\frac{a\ov{q}}{|q|^2}\Jk, \qquad \widecheck{e_4'(r)}:=e_4'(r)-1, \qquad \widecheck{\nab_4'\Jk}:=\nab_4'\Jk+\frac{1}{q}\Jk,
\eeaa

\item the notation $\Ga_b'$ is the one of Definition \ref{definition.Ga_gGa_b}, except that $\Hc'$ does not belong to $\Ga_b'$,

\item the notation $\Ga_g'$ is given by 
\beaa
\Ga_g' &=& \Ga_{g,1}'\cup\Ga_{g,2}',
\eeaa
where $\Ga_{g,1}'$ is the one of Definition \ref{definition.Ga_gGa_b}, and where $\Ga_{g,2}'$ is given by\footnote{Note that all quantities in $\Ga_{g,2}'$ vanish identically in the case of an outgoing PG structure except $\Hc'$ and $\Hbc'$.}  
\beaa
\Ga_{g,2}' := \Big\{\om',\,\,\,  \Xi',\,\,\, \Hbc',\,\,\, \Hc', \,\,\, \widecheck{e_4'(r)}, \,\,\, e_4'(u), \,\,\, e_4'(J^{(0)}), \,\,\, r^{-1}\nab'(r), \,\,\, \widecheck{\nab_4'\Jk}\Big\}.
\eeaa
\end{itemize}
\end{definition}

Using these notations we can  define the decay norms  $\Dk'_k$ exactly as the norms $\Dk_k$  in section \ref{section:main-normsextregion}.  In the proposition below, however,  we derive estimates for  $ \Dk'_k$ norms  for  values of $k\le k_{small} +129$, with the particular choice 129 being  sufficiently large to absorb possible losses of derivatives later on in the process of proving our main theorem.   This requires an interpolation between the decay norms, for $k\le k_{small}$ and the    boundedness norms for $k\leq k_{large}$. The interpolation  leads to a slight loss of decay  which affects the small constant $\dec$ in the definition of  the decay norms.
Thus, in the   statement of the proposition below, the norms  $ \Dk'_k$  are defined exactly as  $\Dk_k$ with       $\dec$ replaced by $\dec'=\dec-2\de_0$. The precise definitions of   $\de_0$ is as follows\footnote{Note that we have in view of \eqref{eq:constraintsonthemainsmallconstantsepanddelta} and \eqref{eq:choiceksmallmaintheorem}
\beaa
\dec(k_{large}-k_{small})\geq \frac{1}{2}\dec k_{large}-\dec \gg 1,
\eeaa
and we may thus assume $\dec(k_{large}-k_{small})\geq 390$ so that we have indeed $\de_0\leq\frac{\dec}{3}$.}
\bea
\lab{eq:constraintsonklossandde0forsecondframeofMext}
\de_0:=\frac{130}{k_{large}-k_{small}}, \qquad 0<\de_0\leq \frac{\dec}{3}.
\eea

\begin{proposition}\lab{prop:constructionsecondframeinMext}
Let $(e_4, e_3, e_1, e_2)$ the outgoing PG  frame of $\Mext$. There exists another frame $(e_4', e_3', e_1', e_2')$ of $\Mext$   given  by \eqref{secondframe-Mext1}
such that:
\begin{enumerate}
\item  The 1-form $f$ appearing in \eqref{secondframe-Mext1} vanishes in $\Mext(r\leq u^{\frac{1}{2}})$. In particular, the frame $(e_4', e_3', e_1', e_2')$ coincides with the outgoing PG frame $(e_4, e_3, e_1, e_2)$ of $\Mext$ in $\Mext(r\leq u^{\frac{1}{2}})$. 

\item The 1-form $f$ appearing in \eqref{secondframe-Mext1} vanishes also on $\{u=u_*\}$. In particular, the frame $(e_4', e_3', e_1', e_2')$ coincides with the outgoing PG frame $(e_4, e_3, e_1, e_2)$ of $\Mext$ on $\{u=u_*\}$. 

\item  The norms $\Dkext'_k$  defined as  $\Dkext_k$, with $\Ga_g, \Ga_b $ replaced by $\Ga_g',  \Ga_b'$, $\dk$ replaced by $\dk'$, $\dec$ replaced by $\dec'=\dec-2\de_0$, and $\Ga_g',  \Ga_b'$ given by Definition \ref{definition.Ga_gGa_b:withprimseforsecondframeMext}, verify the estimates
\bea
\max_{0\leq k\leq k_{small}+129}\Dkext'_k \les \ep. 
\eea

\item  The  horizontal 1-form $f$ verifies, for $k\leq k_{small}+130$  on $\Mext$,
\bea\lab{eq:estimateforfinconstructionsecondframeinMext}
\begin{split}
|(\dk')^kf| &\les \frac{\ep}{ru^{\frac{1}{2}+\dec'}+u^{1+\dec'}}, \\
 |(\dk')^{k-1}\nab_3'f| &\les \frac{\ep}{ru^{1+\dec'}}.
\end{split}
\eea

\item In addition to the control induced by $\xi'\in \Ga_g'$, we have, for $k\leq k_{small}+129$  on $\Mext$, 
\bea
|(\dk')^k\xi'| \les \frac{\ep}{r^{3+\dec'}}, \qquad   |(\dk')^{k-1}\nab_3'\xi'| \les \frac{\ep}{r^3u^{\frac{1}{2}+\dec'}}.
\eea
\end{enumerate}
\end{proposition}

\begin{remark}
The crucial point of Proposition \ref{prop:constructionsecondframeinMext} is that in the new frame $(e_4', e_3', e_1', e_2')$ of $\Mext$, $\Hc'$ belongs to $\Ga_g'$, see Definition \ref{definition.Ga_gGa_b:withprimseforsecondframeMext}, and thus displays a better decay in $r^{-1}$ than $\Hc$ corresponding to the outgoing PG frame $(e_4, e_3, e_1, e_2)$ of $\Mext$. Also, since $\Ga_b'$ and $\Ga_g'$ satisfy the same estimates in $\Mext(r\les u^{\frac{1}{2}})$, $\Hc$ displays the correct behavior in such regions. This is why we may choose  the frame $(e_4', e_3', e_1', e_2')$ to coincide with the outgoing PG frame $(e_4, e_3, e_1, e_2)$ of $\Mext$ in $\Mext(r\leq u^{\frac{1}{2}})$. 
\end{remark}

%%%%%%%%%%%%%%%%%%%%%%%%%%%%%%%%%%%%%%%

\subsection{Construction of the global null frame}

%%%%%%%%%%%%%%%%%%%%%%%%%%%%%%%%%%%%%%%

We start with the definition of the region where the frames of $\Mint$, ${}^{(top)}\MM$, and a conformal renormalization of the second frame of $\Mext$ (i.e. the one of Proposition \ref{prop:constructionsecondframeinMext}) will be matched.
\begin{definition}\lab{def:cutofffunctionforthematchingregion:bis}
We define the matching region as the spacetime region
\bea
\nn\mr &:=& \Big(\Mext\cap\{\rext\geq r_0+1\}\cap\{u_*-5\leq u\leq u_*\}\Big)\\
&&\cup\Big(\Mint\cap\{\rint\leq r_0-1\}\cap\{u_*-5\leq\ub\leq u_*\}\Big).
\eea
\end{definition}

Also, we introduce the notations ${}^{(glo)}\Ga_g, {}^{(glo)}\Ga_b$.

\begin{definition}\lab{definition.Ga_gGa_b:globalframeofMM}
Let $({}^{(glo)}e_4, {}^{(glo)}e_3, {}^{(glo)}e_1, ^{(glo)}e_2)$ the global frame of $\MM$ constructed in Proposition \ref{prop:existenceandestimatesfortheglobalframe:bis} below, together with the corresponding pair of scalars $({}^{(glo)}r, {}^{(glo)}J^{(0)})$ and the  complex 1-form ${}^{(glo)}\Jk$. With respect to that frame, we introduce the notations ${}^{(glo)}\Ga_g, {}^{(glo)}\Ga_b$ as follows, where we drop ${}^{(glo)}$ after the first item   to ease the notations:
\begin{itemize}
\item the linearized quantities for the frame $({}^{(glo)}e_4, {}^{(glo)}e_3, {}^{(glo)}e_1, ^{(glo)}e_2)$ are defined  in the same way  as Definition \ref{def:renormalizationofallnonsmallquantitiesinPGstructurebyKerrvalue:ingoingcase}  for ingoing PG structures, with respect to the scalars $({}^{(glo)}r, {}^{(glo)}J^{(0)})$ and the  complex 1-form ${}^{(glo)}\Jk$,

\item in addition, we introduce the following linearized quantities which are trivial for an ingoing PG structure\footnote{Except $\Hc$ which satisfies instead $H=-Z$.}
\beaa
\Hc:=H-\frac{aq}{|q|^2}\Jk, \qquad \widecheck{e_3(r)}:=e_3(r)+1, \qquad \widecheck{\nab_3\Jk}:=\nab_3\Jk-\frac{1}{\ov{q}}\Jk,
\eeaa

\item the notation $\Ga_b$ is given by
\beaa
\Ga_b &=& \Ga_{b,1}\cup\Ga_{b,2},
\eeaa
$\Ga_{b,1}$ is the one of Definition \ref{definition.Ga_gGa_b:ingoingcase}, and where $\Ga_{b,2}$ is given by\footnote{Note that all quantities in $\Ga_{b,2}$ vanish identically in the case of an ingoing PG structure.}  
\beaa
\Ga_{b,2} := \Big\{\omb,\,\,\,  \Xib, \,\,\, r^{-1}\widecheck{e_3(r)}, \,\,\, e_3(J^{(0)}), \,\,\, r\widecheck{\nab_3\Jk}\Big\}.
\eeaa

\item the notation $\Ga_g$ is given by 
\beaa
\Ga_g &=& \Ga_{g,1}\cup\{\Hc,\,\,\, r^{-1}\nab(r)\},
\eeaa
where $\Ga_{g,1}$ is the one of Definition \ref{definition.Ga_gGa_b:ingoingcase}.
\end{itemize}
\end{definition}

Here is our main proposition concerning our global frame.
\begin{proposition}\lab{prop:existenceandestimatesfortheglobalframe:bis}
Let the small constant $\de_0>0$ satisfying \eqref{eq:constraintsonklossandde0forsecondframeofMext}. There exists 
\begin{itemize}
\item a global null frame $({}^{(glo)}e_4, {}^{(glo)}e_3, {}^{(glo)}e_1, ^{(glo)}e_2)$,

\item a pair of scalars $({}^{(glo)}r, {}^{(glo)}J^{(0)})$, and  a complex 1-form ${}^{(glo)}\Jk$,  
\end{itemize}
all defined on $\MM$ such that:
\begin{itemize}
\item[(a)] In $\Mext\setminus\mr$, we have
\beaa
({}^{(glo)}e_4, {}^{(glo)}e_3, {}^{(glo)}e_1, ^{(glo)}e_2)= \left({}^{(ext)}\la\,{}^{(ext)}e_4', {}^{(ext)}\la^{-1}{}^{(ext)}e_3', {}^{(ext)}e_1', ^{(ext)}e_2'\right),
\eeaa
where $({}^{(ext)}e_4', {}^{(ext)}e_3', {}^{(ext)}e_1', ^{(ext)}e_2')$ denotes the second frame of $\Mext$, i.e. the fame of Proposition \ref{prop:constructionsecondframeinMext}, and where ${}^{(ext)}\la:=\frac{{}^{(ext)}\De}{|{}^{(ext)}q|^2}$. We also have ${}^{(glo)}r=\rext$, ${}^{(glo)}J^{(0)}=\cos({}^{(ext)}\th)$, and ${}^{(glo)}\Jk={}^{(ext)}\Jk$.

\item[(b)] In $\Mint\setminus\mr$, we have
\beaa
({}^{(glo)}e_4, {}^{(glo)}e_3, {}^{(glo)}e_1, ^{(glo)}e_2) = ({}^{(int)}e_4, {}^{(int)}e_3, {}^{(int)}e_1, {}^{(int)}e_2),
\eeaa
as well as ${}^{(glo)}r=\rint$, ${}^{(glo)}J^{(0)}=\cos({}^{(int)}\th)$, and ${}^{(glo)}\Jk={}^{(int)}\Jk$.

\item[(c)] In ${}^{(top)}\MM$, we have
\beaa
({}^{(glo)}e_4, {}^{(glo)}e_3, {}^{(glo)}e_1, ^{(glo)}e_2) = ({}^{(top)}e_4, {}^{(top)}e_3, {}^{(top)}e_1, {}^{(top)}e_2),
\eeaa
as well as ${}^{(glo)}r={}^{(top)}r$, ${}^{(glo)}J^{(0)}=\cos({}^{(top)}\th)$, and ${}^{(glo)}\Jk={}^{(top)}\Jk$.
 
\item[(d)] In the matching region, we have
\beaa
\max_{0\leq k\leq k_{small}+125}\sup_{\mr\cap\Mint}u^{1+\dec-2\de_0}\left|\dk^k({}^{(glo)}\Ga_g, {}^{(glo)}\Ga_b)\right| &\les& \ep,
\eeaa
and
\beaa
\max_{0\leq k\leq k_{small}+125}\sup_{\mr\cap\Mext}\Bigg\{\Big(r^2u^{\frac{1}{2}+\dec-2\de_0}+ru^{1+\dec-2\de_0}\Big)\left|\dk^k{}^{(glo)}\Ga_g\right|\qquad\quad\\
+ru^{1+\dec-2\de_0}\left|\dk^k{}^{(glo)}\Ga_b\right| +r^2u^{1+\dec-2\de_0}\left|\dk^{k-1}\nab_{{}^{(glo)}e_3}{}^{(glo)}\Ga_g\right|\\
+r^{\frac{7}{2}+\frac{\dt}{2}}\left|\dk^k({}^{(glo)}A, {}^{(glo)}B)\right| +r^4u^{\frac{1}{2}+\dec-2\de_0}\left|\dk^{k-1}\nab_{{}^{(glo)}e_3}{}^{(glo)}B\right|\\
+\Big(r^{\frac{9}{2}+\frac{\dt}{2}}+r^4u^{\frac{1}{2}+\dec-2\de_0}\Big)\left|\dk^{k-1}\nab_{{}^{(glo)}e_3}{}^{(glo)}A\right|\Bigg\} &\les& \ep,
\eeaa
where ${}^{(glo)}\Ga_g$ and ${}^{(glo)}\Ga_b$ are given by Definition \ref{definition.Ga_gGa_b:globalframeofMM}. 

\item[(e)] In the matching region, we have
\beaa
\max_{0\leq k\leq k_{small}+125}\sup_{\mr\cap\Mint}u^{1+\dec-2\de_0}\Big|\dk^k({}^{(glo)}r-\rint, \, {}^{(glo)}J^{(0)}-\cos({}^{(int)}\th),\\
 {}^{(glo)}\Jk-{}^{(int)}\Jk)\Big| &\les& \ep,
\eeaa
and 
\beaa
\max_{0\leq k\leq k_{small}+125}\sup_{\mr\cap\Mext}u^{1+\dec-2\de_0}\Big(\left|\dk^k({}^{(glo)}r-\rext)\right|\\
+r\left|\dk^k({}^{(glo)}J^{(0)}-\cos({}^{(ext)}\th))\right|+r^2\left|\dk^k( {}^{(glo)}\Jk-{}^{(ext)}\Jk)\right|\Big) &\les& \ep.
\eeaa

\item[(f)] Let $(f, \fb, \la)$ denote the change of frame coefficients from the frame of $\Mint$ to the global frame, and let  $(f', \fb', \la')$ denote the change of frame coefficients from the second frame of $\Mext$, i.e. the fame of Proposition \ref{prop:constructionsecondframeinMext}, to the global frame. Then,  in the matching region, we have
\beaa
\max_{0\leq k\leq k_{small}+125}\sup_{\mr\cap\Mint}u^{1+\dec-2\de_0}\left|\dk^k(f, \fb, \la)\right| &\les& \ep,
\eeaa
and
\beaa
\max_{0\leq k\leq k_{small}+125}\sup_{\mr\cap\Mext}ru^{1+\dec-2\de_0}\left|\dk^k\left(f', \fb', \log\left(\frac{|{}^{(ext)}q|^2}{{}^{(ext)}\De}\la'\right)\right)\right| &\les& \ep.
\eeaa

\item[(g)] In addition to the control induced by ${}^{(glo)}\xi\in {}^{(glo)}\Ga_g$, we have, for $k\leq k_{small}+125$  on $\Mext$, 
\beaa
|\dk^k{}^{(glo)}\xi| \les \frac{\ep}{r^{3+\dec-2\de_0}}, \qquad   |\dk^{k-1}\nab_{{}^{(glo)}e_3}{}^{(glo)}\xi| \les \frac{\ep}{r^3u^{\frac{1}{2}+\dec-2\de_0}},
\eeaa
and on $\Mtop$
\beaa
|\dk^k{}^{(glo)}\xi| \les \frac{\ep}{r^{3+\dec-2\de_0}}, \qquad   |\dk^{k-1}\nab_{{}^{(glo)}e_3}{}^{(glo)}\xi| \les \frac{\ep}{r^3({}^{(top)}u)^{\frac{1}{2}+\dec-2\de_0}}.
\eeaa
\end{itemize}
\end{proposition}

\begin{remark}
\label{remark-globalframeforM1:bis}
Recall that the global frame on $\MM$ of Proposition \ref{prop:existenceandestimatesfortheglobalframe:bis} will be needed to derive decay estimates for the quantity $\qf$ in Theorem M1 (stated in section  \ref{sec:mainintermediateresults:chap3}),  see also the discussion in  section \ref{section:outlineproof-intro}.
\end{remark}

%%%%%%%%%%%%%%%%%%%%%%%%%%%%%%%%%%%%%%%%%%%%%%%%

\section{Proof of the main theorem}
\lab{section:ProofMainThm}

%%%%%%%%%%%%%%%%%%%%%%%%%%%%%%%%%%%%%%%%%%%%%%%%

%%%%%%%%%%%%%%%%%%%%%%%%%%%%%%%%%%%%%%%%%%%%%%%%

\subsection{Main intermediate results}
\lab{sec:mainintermediateresults:chap3} 

%%%%%%%%%%%%%%%%%%%%%%%%%%%%%%%%%%%%%%%%%%%%%%%%

We are ready to   state our main intermediary results.  

\begin{thmM0}
Assume that the initial data layer  $\LL_0$, as   defined  in section \ref{sec:defintionoftheinitialdatalayer},  satisfies
\beaa
\Ik_{k_{large}+10}\leq \ep_0.
\eeaa
Then under the bootstrap assumptions {\bf BA-B} and {\bf BA-D}, the following holds true on the initial data hypersurface $\BB_1\cup\BBb_1$,
\beaa
\max_{0\leq k\leq k_{large-2}}&&\Bigg\{ \sup_{\BB_1}  \left[r^{\frac{7}{2}  +\de_B}\left( |\dk^k\,{}^{(ext)}A| + |\dk^k\,{}^{(ext)}B|\right)+r^{\frac{9}{2}  +\de_B}|\dk^{k-1}\nab_3(\,{}^{(ext)}A)|   \right]\\
\nn&&+ \sup_{\BB_1} \left[r^3\ \left|\dk^k\left(\,{}^{(ext)}P+\frac{2m}{q^3}\right)\right|+r^2|\dk^k\,{}^{(ext)}\Bb|+r|\dk^k\,{}^{(ext)}\Ab|\right]\Bigg\} \les \ep_0,
\eeaa
\beaa
\nn\max_{0\leq k\leq k_{large-2}}\sup_{\BBb_1}\Bigg[ |\dk^k\,{}^{(int)}A| + |\dk^k\,{}^{(int)}B|+ \left|\dk^k\left(\,{}^{(int)}P +\frac{2m}{q^3}\right)\right| &&\\
 +|\dk^k\,{}^{(int)}\Bb|+|\dk^k\,{}^{(int)}\Ab|\Bigg] &\les& \ep_0,
\eeaa
and
\beaa
\sup_{\BB_1\cup\BBb_1}\Big(\left|m-m_0\right|+\left|a-a_0\right|\Big)\les \ep_0.
\eeaa
\end{thmM0}

\begin{thmM1}
Assume given a GCM admissible spacetime $\MM$  as defined in section \ref{sec:defintioncanonicalspacetime} 
verifying the bootstrap   assumptions {\bf BA-B} and {\bf BA-D}  for  some sufficiently small   $\ep>0$. Then, if $\ep_0>0$ is sufficiently small, there exists $\dee>\dec$ such that we have the following   estimates in  $\MM$,
\begin{enumerate}
\item  The quantity $\qf$ verifies the estimate
\beaa
 \max_{0\leq k\leq k_{small}+100}\sup_{\Mext}\left\{\Big(ru^{\frac{1}{2}+\dee}+u^{1+\dee}\Big)|\dk^k\qf|+ru^{1+\dee}|\dk^{k-1}\nab_3\qf|\right\}\\
+\max_{0\leq k\leq k_{small}+100}\sup_{\Mint\cup\Mtop}\ub^{1+\dee}|\dk^k\qf|    &\les&  \ep_0.
\eeaa
Moreover, $\qf$ also satisfies the following estimate
\beaa
\max_{0\leq k\leq k_{small}+100}u^{2+2\dee}\int_{\Sigma_*(\geq u)}|\dk^{k-1}\nab_3\qf|^2 \les \ep_0^2.
\eeaa

\item The quantity $A$ verifies  the estimate, for all  $ k\le  k_{small}+100$, 
\beaa
\sup_{\Mext}\Big(\frac{r^2(2r+u)^{1+\dee}}{\log(1+u)}+r^3(2r+u)^{\frac{1}{2}+\dee}\Big)\Big(|\dk^k\,{}^{(ext)}A|+r|\dk^{k-1}\nab_3\,{}^{(ext)}A|\Big) \les \ep_0.
\eeaa
\end{enumerate}
\end{thmM1}

\begin{thmM2}
Under  the  same assumptions as above 
 we  have the following decay estimates for $\aa$
\beaa
\,{}^{(int)}\Dk_{k_{small}+80}[\aa] +\,{}^{(top)}\Dk_{k_{small}+80}[\aa]     \les  \ep_0, \qquad \max_{0\leq k\leq k_{small}+80}\int_{\Sigma_*}u^{2+2\dec}|\dk^k\aa|^2 \les \ep_0^2.
\eeaa
\end{thmM2}

\begin{thmM3}
Under  the  same assumptions as above 
 we  have the following decay estimates  on $\Si_*$
 \beaa
\,^{*}\mathfrak{D}_{k_{small}+60} \les \ep_0.
 \eeaa
\end{thmM3}

\begin{thmM4}
Under  the  same assumptions as above 
 we  have the following decay estimates  on $\Mext$
 \beaa
\,{}^{(ext)}\mathfrak{D}_{k_{small}+40} \les \ep_0.
 \eeaa
\end{thmM4}

\begin{thmM5}
Under  the  same assumptions as above we also have the following decay estimates in    $\Mint$ and $\Mtop$
\beaa
 \,{}^{(int)}\Dk_{k_{small}+20}+\,{}^{(top)}\Dk_{k_{small}+20} \les  \ep_0.
\eeaa
\end{thmM5}

Note that, as an immediate consequence of Theorem M1 to Theorem M5, we have obtained, under  the  same assumptions as above, the following improvement of our bootstrap assumptions on decay
\bea\lab{eq:finalimprovmentofabsolutlyallbootstrapasumptionseverywhereonMMforallderivatives}
\Nk^{(Dec)}_{k_{small}+20}\les \ep_0.
\eea

%%%%%%%%%%%%%%%%%%%%%%%%%%%%%%%%%%%%%%%%%%%%%%%%

\subsection{End of the proof of the main theorem}
\lab{sec:endoftheproofofthemaintheorem:chap3}

%%%%%%%%%%%%%%%%%%%%%%%%%%%%%%%%%%%%%%%%%%%%%%%%

We end the proof by invoking a continuity argument as  in \cite{KS}. The  argument  is based on Definition \ref{definition:UU(u_*)}  below   of  the set   $\UU(  u_*) $
of  GCM admissible spacetimes verifying  bootstrap assumptions ${ \bf  BA}_\ep$  with  $\ep=\ep_0^{\frac{2}{3}}$.

\begin{definition}[Definition of $\UU(u_*)$]
\lab{definition:UU(u_*)}
Let $\ep_0>0$ and $\ep=\ep_0^{\frac{2}{3}}$ be given small constants,  and let $\UU(u_*)$ be the set of all GCM admissible spacetimes $\MM$ defined in section \ref{sec:defintioncanonicalspacetime} such that
\begin{itemize}
\item $u_*$ is the value of $u$ on the last  sphere $S_*$ of $\Si_*$, 

\item $u_*$ satisfies, see \eqref{eq:behaviorofronS-star}, 
\beaa
r_*= \de_*  \ep_0^{-1} u_*^{1+\dec},
\eeaa

\item relative to the combined norms defined in   section \ref{sec:definitionofconcatenatednorm}, we have\footnote{i.e. the bootstrap assumptions \eqref{def:bootstrapasumptionsglobalnormonenergie} \eqref{def:bootstrapasumptionsglobalnormondecay} hold true.}
\beaa
\Nk^{(Sup)}_{k_{large}}\leq \ep,\qquad \Nk^{(Dec)}_{k_{small}}\le \ep.
\eeaa
\end{itemize}
\end{definition}

\begin{definition}
We define  $\UU$ to  be the set of all values of $u_*\geq 0$  for which the spacetime $\UU(u_*)$ exists. 
\end{definition}

The following theorem shows that $\UU$ is not empty. 
\begin{thmM6}\lab{th:M6}
There exists $\de_0>0$ small enough such that for a sufficiently small constant $\ep_0>0$  we have $[1,1+\de_0]\subset\UU$. 
\end{thmM6}

In view of Theorem M6, we may define $U_*$ as the supremum over all value of $u_*$ that belongs to $\UU$
\beaa
U_* :=\sup_{u_*\in\UU}u_*.
\eeaa

Assume by contradiction that 
$$U_*<+\infty.$$
Then, by the continuity of the flow, $U_*\in \UU$. According to  \eqref{eq:finalimprovmentofabsolutlyallbootstrapasumptionseverywhereonMMforallderivatives},   the bootstrap assumptions on decay \eqref{def:bootstrapasumptionsglobalnormondecay} on any spacetime of $\aleph(U_*)$ are improved by
\beaa
\Nk^{(Dec)}_{k_{small}+20}\les \ep_0.
\eeaa
To reach a contradiction, we still need an extension procedure for spacetimes in $\aleph(u_*)$ to larger values of $u$, as well as to improve our bootstrap assumptions on boundedness   \eqref{def:bootstrapasumptionsglobalnormonenergie}. This is done in two steps.

\begin{thmM7}
Any GCM admissible spacetime in $\aleph(u_*)$ for some $0<u_*<+\infty$ such that 
\beaa
\Nk^{(Dec)}_{k_{small}+20}\les \ep_0
\eeaa
has a GCM admissible extension  verifying \eqref{eq:behaviorofronS-star},    with  $u_*'>u_*$, initialized by Theorem M0, which verifies  
\beaa
\Nk^{(Dec)}_{k_{small}}\les \ep_0.
\eeaa
\end{thmM7}

\begin{thmM8}
The GCM admissible spacetime exhibited in Theorem M7 satisfies in addition 
\beaa
\Nk^{(Sup)}_{k_{large}}\les \ep_0
\eeaa
and therefore belongs to $\aleph(u_*')$. In particular $u_*'$ belongs to $\UU$. 
\end{thmM8}

In view of Theorem M8, we have reached a contradiction, and hence
$$U_*=+\infty$$
so that the spacetime may be continued forever. This concludes the proof of the main theorem.\\

%%%%%%%%%%%%%%%%%%%%%%%%%%%%%%%%%%%%%%%%%%%%%%%%

\section{Conclusions}
\lab{section:conclusion-mainThm}

%%%%%%%%%%%%%%%%%%%%%%%%%%%%%%%%%%%%%%%%%%%%%%%%

We denote by $\MM$ our global spacetime obtained in the limit $u_*\to +\infty$. Note that in the limit $u_*\to +\infty$, the region ${}^{(top)}\MM$ disappears\footnote{Indeed, recall  that $u\geq u_*$ and $\ub\geq u_*$ on ${}^{(top)}\MM$ by construction of our GCM admissible spacetime.} so that 
\bea
\MM=\Mint\cup\Mext,
\eea
where $\Mint$ is covered by a ingoing PG structure, and $\Mext$ by an outgoing PG structure. In particular, note that 
\bea
\bsplit
\Mint &= \{1\leq\ub<+\infty\}\cup\{r_+-\deh\leq\rint\leq r_0\},\\
\Mext &= \{1\leq u<+\infty\}\cup\{r_0\leq \rext<+\infty\},
\end{split}
\eea
where we recall that 
\bea
\lab{eq:definitionr_+}
r_+=m_0+\sqrt{m_0^2-a_0^2.}
\eea
Also, as a consequence of Theorem M0, the parameters $(m_\infty, a_\infty)$ obtained in the limit $u_*\to +\infty$ satisfy 
\bea
|m_\infty -m_0|+|a_\infty -a_0| \les \ep_0. 
\eea
This implies in particular $m_\infty>0$ and $|a_\infty|\ll m_\infty$.

%%%%%%%%%%%%%%%%%%%%%%%%%%%%

\subsection{The Penrose diagram of $\MM$}

%%%%%%%%%%%%%%%%%%%%%%%%%%%

{\bf Complete future null infinity.} We first deduce from our estimate that our spacetime $\MM$ has a complete future null infinity $\II_+$. Let us denote by $(e_4, e_3, e_1, e_2)$ and $(u,r)$ the null frame and scalar functions associated to our outgoing PG structure of $\Mext$. Recall also that the outgoing PG structure of $\Mext$ comes together with a scalar function $\th$ and a complex 1-form $\Jk$.  The scalar function $u$ of $\MM$ satisfies, using $e_4(u)=0$,  
\beaa
\g(\D u, \D u) &=& |\nab u|^2.
\eeaa
Recalling that (see Definition \ref{def:renormalizationofallnonsmallquantitiesinPGstructurebyKerrvalue})  $\nab u = a\Re(\Jk)+\widecheck{\nab u}=a\Re(\Jk)+\Ga_b$, as well as the identity $\Re(\Jk)\c\Re(\Jk)=\frac{(\sin\th)^2}{|q|^2}$, we infer on $\Mext$, using also our control for $\Ga_b$ induced by  the estimates $\Nk^{(Dec)}_{k_{small}}\les \ep_0$ of our main theorem, 
\beaa
\g(\D u, \D u) &=& \frac{a^2(\sin\th)^2}{|q|^2}+O\left(\frac{\ep_0}{r^2u^{2+2\dec}}\right).
\eeaa
In particular, the leaves of the $u$-foliation of $\Mext$ are asymptotically null as $r\to +\infty$, so that the portion of null infinity of $\MM$ corresponds to the limit $r\to+\infty$ along the leaves of the $u$-foliation of the outgoing PG structure of $\Mext$. As these leaves exist for all $u\geq 1$ with suitable estimates, it suffices to prove that $u$ is an affine parameter of $\II_+$. To this end, recall from our main theorem that the estimates $\Nk^{(Dec)}_{k_{small}}\les \ep_0$ hold which implies in particular
\bea\lab{eq:estimate1intheconclusionsforxibombandvsi}
\sup_{\Mext}u^{1+\dec}\left(r|\xib|+r\left|\omb-\frac{1}{2}\pr_r\left(\frac{\De}{|q|^2}\right)\right|+\left|\widecheck{e_3(u)}\right|\right) &\les&\ep_0.
\eea
We infer that 
\beaa
\lim_{r\to +\infty}\xib, \omb =0\textrm{ for all }1\leq u<\infty.
\eeaa
In view of the identity 
\beaa
\D_3e_3 &=& -2\omb e_3 + 2\xib_b e_b,
\eeaa
we infer that $e_3$ is a null geodesic generator of $\II_+$. Since we have 
\beaa
e_3(u)=\frac{2(r^2+a^2)}{|q|^2}+\widecheck{e_3(u)}=\frac{2(r^2+a^2)}{|q|^2}+O\left(\frac{\ep_0}{u^{1+\dec}}\right)=2+O\left(\frac{1}{r^2}+\ep_0\right)
\eeaa 
in view of \eqref{eq:estimate1intheconclusionsforxibombandvsi}, $u$ is an  affine parameter of $\II_+$ so that $\II_+$ is indeed complete.

{\bf Existence of a future event horizon.} Let us denote by $(e_4, e_3, e_1, e_2)$ and $(\ub,r)$ the null frame and scalar functions associated to our ingoing PG structure of $\Mint$. Note that the estimates $\Nk^{(Dec)}_{k_{small}}\les \ep_0$  imply
\beaa
\sup_{\Mint}\ub^{1+\dec}\left(\left|e_3(r)+1\right|+\left|e_4(r)-\frac{\De_\infty}{|q|^2}\right|\right) &\les&\ep_0.
\eeaa
where $\De_\infty=\De(a_\infty, m_\infty)$.

In particular, considering the spacetime region $r\leq r_+(1-\deh/2)$ of $\Mint$, and in view of the estimate $|m_\infty-m_0|\les\ep_0$, we infer, for all $r\le r_+(1-\deh/2)$, that 
\beaa
\De_\infty &=&\De(a_\infty, m_\infty) = r^2-2m_\infty r+a_\infty^2=(r-r_+(a_\infty, m_\infty))(r-r_-(a_\infty, m_\infty))\\
&=& (r-r_+)(r-r_-)+O(m_\infty -m_0)+O(a_\infty -a_0)\\
&\leq& -\frac{\de_\HH}{2}r_+\left(r_+-r_- -\frac{\de_\HH r_+}{2}\right)+O(\ep_0)\\
&\les& -\deh<0
\eeaa
and hence
\beaa
e_3(r)\leq -\frac{1}{2}<0, \qquad e_4(r)\les -\deh<0 \quad\textrm{on}\quad \Mint\left(r\leq r_+(1-\deh/2)\right).
\eeaa
Consider now $\gamma(s)$ any future directed null geodesic emanating from a point of the region $\Mint(r\leq r_+(1-\deh/2))$.  Since $\dot{\ga}$ is  a null vector, there exists at any point of $\ga(s)$ in $\Mint$ a scalar $\la$ and a 1-form $f$ such that 
\beaa
\dot{\ga} &=& \la\left(e_4+f^be_b+\frac{1}{4}|f|^2e_3\right),
\eeaa
where $\la>0$  (since $\dot{\ga}$ is future directed). Since $\nab(r)=0$, we infer
\beaa
\frac{dr}{ds}=\D_{\dot{\ga}}r &=& \la\left(e_4+f^be_b+\frac{1}{4}|f|^2e_3\right)r=\la\left(e_4(r)+\frac{1}{4}|f|^2e_3(r)\right).
\eeaa
 Since $e_3(r)<0$ and $e_4(r)<0$ in $\Mint(r\leq r_+(1-\deh/2))$ in view of the above, and since $\la>0$ and $|f|^2\geq 0$, we deduce that $r$ decreases along $\ga(s)$ so that, in particular,  $\ga(s)$ 
 cannot  reach $\II^+$.  Thus the past of $\II^+$  does not contain this region  and hence $\MM$ contains the event horizon $\HH_+$ of a black hole in its interior.
  Moreover, since any point on the timelike hyper surface $\TT=\Mint\cap\Mext$ is on an outgoing null geodesic in $\Mext$ of geodesic generator $e_4$ with $e_4(r)=1$ and defined for all $r\geq r_0$, $\TT$  is in the past of $\II_+$. Hence, since $\TT$ is one of the boundaries of $\Mint$, $\HH_+$ is actually located in the interior of the region $\Mint$.  

{\bf Asymptotic stationarity of $\MM$.} Recall that we have introduced a vectorfield $\T$ in $\Mext$ as well as one in $\Mint$ by
\beaa
\T &=& \frac{1}{2}\left(e_3+\frac{\Delta}{|q|^2}e_4 -2a\Re(\Jk)^be_b\right)\quad\textrm{in}\quad\Mext,\\ 
\T &=& \frac{1}{2}\left(e_4+\frac{\Delta}{|q|^2}e_3 -2a\Re(\Jk)^be_b\right)\quad\textrm{in}\quad\Mint.
\eeaa
Also,  see Proposition \ref{Proposition:deftensorT} , all components of $\piT$ belong, at least, to  $\Ga_b$. Thus, making  use of  the estimate $\Nk^{(Dec)}_{k_{small}}\les \ep_0$ 
of our main theorem, we deduce,
\beaa
|{}^{(\T)}\pi|\les \frac{\ep_0}{ru^{1+\dec}}\quad\textrm{in}\quad\Mext\quad\textrm{and}\quad|{}^{(\T)}\pi|\les \frac{\ep_0}{\ub^{1+\dec}}\quad\textrm{in}\quad\Mint.
\eeaa
In particular, $\T$ is an asymptotically Killing vectorfield and hence our spacetime $\MM$ is asymptotically stationary.\\

The above conclusions regarding $\II_+$ and $\HH_+$ allow us to draw the Penrose diagram of $\MM$, see figure \ref{fig:penrosediagramconclusionmaintheorem} below.
\begin{figure}[h!]
\centering
\includegraphics[scale=0.3]{kerr_3.pdf}
\caption{The Penrose diagram of the space-time $\mathcal{M}$ with past boundary $\BBb_1\cup\BB_1$}
\label{fig:penrosediagramconclusionmaintheorem}
\end{figure}

\subsection{Limits at $\II^+$}

%%%%%%%%%%%%%%%%%%%%%%%%%%%%

\subsubsection{Integrable frame}

%%%%%%%%%%%%%%%%%%%%%%%%%%%%

We  denote by $(e_4, e_3, e_1, e_2)$ and $(u,r)$ the null frame and  defining  functions associated to our outgoing PG structure of $\Mext$. Recall that $(e_1, e_2)$ is not tangent  to the spheres $S(u,r)$ since $\nab(u)\neq 0$. Recall that we have exhibited in  Lemma \ref{Lemma:Transformation-principal-to-integrable frames} a    frame transformation  taking $(e_3, e_4, e_1, e_2)$ into an integrable null frame  $(e_3', e_4',  e_1', e_2')$, i.e. such that  the horizontal vectors $(e_1', e_2')$ are tangent to $S(u, r)$. The corresponding frame coefficients  $(f, \fb, \la)$ satisfy, see \eqref{def:transition-functs:ffbla}  \eqref{eq:formulaforchangeofframecoeffandfbarforframetangenttosphereSur}, 
\bea\lab{eq:changeofframecoefficientsfromPGtointegralframeconclusions}
\bsplit
f &= -\Big(1+O(r^{-2})+r\Ga_b\Big)\nab u,\\ 
\fb &= -\left(1-\frac{2m}{r}+O(r^{-2})+r\Ga_b\right)\nab u,\\
\la &= 1+O(r^{-2})+r^{-1}\Ga_b.
\end{split}
\eea
Then, relying on the frame transformations formulas of Proposition \ref{Proposition:transformationRicci} and the above control of $(f, \fb, \la)$, we obtain\footnote{See Lemma \ref{Lemma:Transformation-principal-to-integrable-Kerr} where the corresponding passage from the PG frame to the integrable frame is done in details in the particular case of Kerr.} for the frame $(e_3', e_4',  e_1', e_2')$ 
\bea\lab{eq:behavioroftheintegrableframeforlarger1}
\bsplit
\trch'& = \frac 2 r  +O(r^{-3})+\dk^{\leq 1}\Ga_g, \qquad \qquad 
 \chih'=   O(r^{-3})+\dk^{\leq 1}\Ga_g, \\
\trchb'&= -\frac{2\left(1-\frac{2m_\infty}{r}\right)}{r} + O(r^{-3})+\dk^{\leq 1}\Ga_g, \qquad \,\,\chibh'=     O(r^{-3})+\dk^{\leq 1}\Ga_b, \\
 \ze'&=O(r^{-3})+\dk^{\leq 1}\Ga_g,
\end{split}
\eea
\bea\lab{eq:behavioroftheintegrableframeforlarger2}
\bsplit
\xi'&=O\big(r^{-3}\big)+r^{-1}\dk^{\leq 1}\Ga_b, \qquad\qquad\,\,\,  \xib'=O\big(r^{-3}\big)+\dk^{\leq 1}\Ga_b,\\
\om'&=O\big(r^{-3}\big)+r^{-1}\dk^{\leq 1}\Ga_b, \qquad\qquad  \,\, \omb'=\frac{m_\infty}{r^2} +O\big(r^{-3}\big)+\dk^{\leq 1}\Ga_b,\\
 \etab'&=O(r^{-3})+\dk^{\leq 1}\Ga_g, \qquad \qquad \qquad  \eta'=O(r^{-3})+\dk^{\leq 1}\Ga_b,\\
\end{split}
\eea
and
\bea\lab{eq:behavioroftheintegrableframeforlarger3}
\bsplit
\a'&= \a+O(r^{-5})+r^{-2}\Ga_g,\qquad \quad  \,\,\,\b'=\b+ O(r^{-4})+r^{-2}\Ga_g,\\
\rho'&=-\frac{2m_\infty}{r^3} + O(r^{-5})+r^{-1}\Ga_g, \qquad 
\rhod'=O(r^{-4})  + r^{-1}\Ga_g,\\
\bb'&= O\big(r^{-4}\big)+r^{-1}\Ga_b,\qquad\qquad\quad\,\,\aa'= O\big(r^{-5} \big)+\Ga_b.
\end{split}
\eea

%%%%%%%%%%%%%%%%%%%%%%%%%%%%%%%%%%%%%%

\subsubsection{The spheres at null infinity are round}

%%%%%%%%%%%%%%%%%%%%%%%%%%%%%%%%%%%%%%

Recalling that the integrable frame $(e_3', e_4',  e_1', e_2')$ is such that  the horizontal vectors $(e_1', e_2')$ are tangent to $S(u, r)$, the Gauss curvature $K$ of the spheres $S(u,r)$  is given by the Gauss equation
\beaa
K &=& -\rho'-\frac{1}{4}\trch'\trchb' +\frac{1}{2}\chih'\c\chibh'.
\eeaa
In view of \eqref{eq:behavioroftheintegrableframeforlarger1} and \eqref{eq:behavioroftheintegrableframeforlarger3}, we infer
\beaa
K-\frac{1}{r^2} &=& O(r^{-4})+r^{-1}\Ga_g.
\eeaa
Thus, in view of our estimates in $\Mext$ for $\Ga_g$, we deduce
\beaa
\left|K -\frac{1}{r^2}\right| &\les& \frac{1}{r^4}+\frac{\ep_0}{r^3u^{\frac{1}{2}+\dec}}
\eeaa
so that 
\beaa
\lim_{r\to +\infty}r^2K =1.
\eeaa
In particular the spheres at null infinity are round. 

Also, using Gauss-Bonnet, we have 
\beaa
4\pi &=& \int_SK=\int_S\left(\frac{1}{r^2}+O\left(\frac{1}{r^4}+\frac{\ep_0}{r^3u^{\frac{1}{2}+\dec}}\right)\right)\\
&=& \frac{|S|}{r^2}+O\left(\frac{|S|}{r^4}+\frac{|S|\ep_0}{r^3u^{\frac{1}{2}+\dec}}\right)
\eeaa
and hence
\bea\lab{eq:risalmostthearearadiusofthespehresSpfuandr}
\sqrt{\frac{|S|}{4\pi}} &=& r\left(1+O(r^{-2})+O(\ep_0r^{-1}u^{-\frac{1}{2}-\dec})\right)
\eea
which shows that $r$ is a good approximation of the area radius of the spheres $S(u,r)$.

%%%%%%%%%%%%%%%%%%%%%%%%%%%%

\subsubsection{Limits  at null infinity and Bondi mass}

%%%%%%%%%%%%%%%%%%%%%%%%%%%%

Recall the definition of the Hawking mass, associated here to the spheres $S(u,r)$
\beaa
m_H &=& \frac{\sqrt{\frac{|S|}{4\pi}}}{2}\left(1+\frac{1}{16\pi}\int_S\trch'\trchb'\right).
\eeaa
In view of \eqref{eq:risalmostthearearadiusofthespehresSpfuandr} and \eqref{eq:behavioroftheintegrableframeforlarger1}, we infer
\bea\lab{eq:estimatebetweenmHandmforconclusions}
m_H &=& m_\infty\Big(1+O(r^{-1})+O(\ep_0u^{-\frac{1}{2}-\dec})\Big).
\eea
 Also, we differentiate the identity for $m_H$ w.r.t. the vectorfield $e_4$ of the PG structure of $\Mext$ and find
\beaa
e_4(m_H) &=& \frac{m_H}{2|S|}e_4\left(|S|\right) +\frac{1}{32\pi}\sqrt{\frac{|S|}{4\pi}}e_4\left(\int_S\trch'\trchb'\right). 
\eeaa
Next, we make use of the following corollary of  Lemma \ref{lemma:e4derivativeofintegralonS}.
\begin{lemma}
\lab{lemma:e_4(int_S h-chap3)}
Given a scalar function $h$ we have
\beaa
e_4\left(\int_{S(r,u)}h\right) &=& \int_{S(r,u)}\Big(e_4(h)+\trch' h+f\c\nab' h\Big)+O(r^{-1})h.
\eeaa
\end{lemma}

\begin{proof}
 According to  Lemma \ref{lemma:e4derivativeofintegralonS} 
 \beaa 
e_4\left(\int_{S(r,u)}h\right) &=& \int_{S(r,u)}\Big(e_4(h)+\de^{ab}\g(\D_{e_a'}e_4, e_b')h\Big).
\eeaa 
In view of the precise formula  for $\de^{ab}\g(\D_{e_a'}e_4, e_b')$ in that lemma,  the estimates  for $f, \fb$ in \eqref{eq:changeofframecoefficientsfromPGtointegralframeconclusions}
 \beaa
\de^{ab}\g(\D_{e_a'}e_4, e_b') &=&\big( 1+ O(r^{-2}) \big) \trch + O(r^{-3} ) +O(r^{-2})  \dk^{\leq 1}\Ga_b.
\eeaa
Next we  need to  replace  $\trch$ with $ \trch'$. To this end, we make use\footnote{Note that  the formula $ \trch  =\big( 1+ O(r^{-2}) \big) \trch'  + O(r^{-3} ) +O(r^{-2})\dk^{\leq 1} \Ga_b  $ is not good enough.} of the transformation formula  for $\trch'$ in Proposition \ref{Proposition:transformationRicci}. Together with \eqref{eq:changeofframecoefficientsfromPGtointegralframeconclusions}   we deduce
\beaa
\trch' &= \trch  +  \div'f    +O(r^{-3}).
\eeaa
Consequently
\beaa
\de^{ab}\g(\D_{e_a'}e_4, e_b') &=&\trch'-\div' f +O(r^{-3}).
\eeaa
Hence
\beaa
e_4\left(\int_{S(r,u)}h\right) &=& \int_{S(r,u)}\Big(e_4(h) +\trch' h  -( \div'  f) h\Big)+ O(r^{-1} )  h \\
&=& \int_{S(r,u)}\Big(e_4(h) +\trch' h  + f \nab'  h\Big)+ O(r^{-1} )   h
\eeaa 
as stated.
\end{proof}

Using the above identity with the choice $h=1$ and $h=\trch'\trchb'$, we infer, using also  \eqref{eq:behavioroftheintegrableframeforlarger1},  $f=O(r^{-1})$ and  the gain in powers of $r$   for  $\nab'( \trch'\trchb')$,
\beaa 
e_4\left(|S|\right) &=& \int_{S}\trch'+O(r^{-1})=\frac{2|S|}{r}+O(1),\\
e_4\left(\int_S\trch'\trchb'\right) &=& \int_{S}\Big(e_4(\trch'\trchb')+(\trch')^2\trchb'\Big)+O(r^{-3}).
\eeaa
Plugging in the above identity for $e_4(m_H)$, we infer
\beaa
e_4(m_H) &=& \frac{m_H}{r} +\frac{1}{32\pi}\sqrt{\frac{|S|}{4\pi}}\int_{S}\Big(e_4(\trch'\trchb')+(\trch')^2\trchb'\Big)+O(r^{-2}). 
\eeaa
Now, we use the following computation, see the proof of Lemma \ref{lemma:transportequationforHawkingmass}, 
\beaa
e_4'(\trch'\trchb') &=& -{\trch'}^2\trchb' +2\trch'\rho' +2\trchb'\div'\xi' +2\trch'\div'\etab'\\
&&  +\trchb'\left(  2 \xi'\c(\etab'+3\ze')-|\chih'|^2\right)+\trch'\left(  2|\etab'|^2 -\chih'\c\chibh'\right).
\eeaa
In view of \eqref{eq:behavioroftheintegrableframeforlarger1}, \eqref{eq:behavioroftheintegrableframeforlarger2} and \eqref{eq:behavioroftheintegrableframeforlarger3}, this yields
\beaa
e_4'(\trch'\trchb') &=& -{\trch'}^2\trchb' +\frac{4}{r}\rho' -\frac{4}{r}\div'\xi' +\frac{4}{r}\div'\etab' -\frac{2}{r}\chih'\c\chibh'+O(r^{-5}).
\eeaa
Also, since $\la=1+O(r^{-2})$ and $f=O(r^{-1})$, we have
\beaa
e_4'(\trch'\trchb') &=& \la\left(e_4+f\c\nab+\frac{1}{4}|f|^2e_3\right)(\trch'\trchb')\\
&=& e_4(\trch'\trchb')+O(r^{-2})e_4(\trch'\trchb')+O(r^{-1})\nab(\trch'\trchb')\\
&&+O(r^{-2})e_3(\trch'\trchb')\\
&=& e_4(\trch'\trchb')+O(r^{-5})
\eeaa
and hence
\beaa
e_4(\trch'\trchb') &=& -{\trch'}^2\trchb' +\frac{4}{r}\rho' -\frac{4}{r}\div'\xi' +\frac{4}{r}\div'\etab' -\frac{2}{r}\chih'\c\chibh'+O(r^{-5}).
\eeaa
Integrating on $S$, and integrating the divergences by parts, we obtain 
\beaa
\int_{S}\Big(e_4(\trch'\trchb')+(\trch')^2\trchb'\Big) &=& \int_S\left(\frac{4}{r}\rho' -\frac{2}{r}\chih'\c\chibh'\right)+O(r^{-3}).
\eeaa
We deduce from the above
\beaa
e_4(m_H) &=& \frac{m_H}{r} +\frac{1}{8\pi r}\sqrt{\frac{|S|}{4\pi}}\int_S\left(\rho' -\frac{1}{2}\chih'\c\chibh'\right)+O(r^{-2}). 
\eeaa
Now, in view of the Gauss equation, and using Gauss-Bonnet and the definition of the Hawking mass, we have
\beaa
\int_S\left(\rho' -\frac{1}{2}\chih'\c\chibh'\right) &=& \int_S\left(-K-\frac{1}{4}\trch'\trchb'\right)= -4\pi-\frac{1}{4}\int_S\trch'\trchb' = -\frac{8\pi m_H}{\sqrt{\frac{|S|}{4\pi}}}
\eeaa
and hence 
\bea\lab{eq:estimatefore4derivativeofHawkingmassinMext}
|e_4(m_H)| &\les& r^{-2}.
\eea
Since $r^{-2}$ is integrable, we infer the existence of a limit to $m_H$ as $r\to +\infty$ along the leaves  of the $u$-foliation of $\Mext$
\beaa
M_B(u) &=& \lim_{r\to +\infty} m_H(u, r)\textrm{ for all }1\leq u<+\infty,
\eeaa
where $M_B(u)$ is the so-called Bondi mass.
 
Next, we have the following null structure equation in $\Mext$, see Proposition \ref{prop-nullstr} and recall that the frame $(e_4', e_3', e_1', e_2')$ is integrable so that $\atrch'=\atrchb'=0$, 
\beaa
\nab_4'\chibh'
&=&-\frac 1 2 \big( \trchb' \chih'+\trch' \chibh'\big)+\nab'\hot \etab' +2 \om' \chibh' +\xi'\hot\xib' +\etab'\hot\etab'.
\eeaa
In view of $\Nk^{(Dec)}_{k_{small}}\les \ep_0$ and \eqref{eq:behavioroftheintegrableframeforlarger1} \eqref{eq:behavioroftheintegrableframeforlarger2}, we  deduce
\beaa
|\nab_4(r\chibh')| &\les& \frac{1}{r^2}.
\eeaa
Since $r^{-2}$ is integrable, we infer the existence of a limit to $r\chibh'$ as $r\to +\infty$ along the leaves  of the $u$-foliation of $\Mext$
\beaa
\underline{\Theta}(u,\cdot) &=& \lim_{r\to +\infty}r\chibh'(r, u, \cdot)\textrm{ for all }1\leq u<+\infty.
\eeaa
  On the other hand, in view of     $\Nk^{(Dec)}_{k_{small}}\les \ep_0$ and \eqref{eq:behavioroftheintegrableframeforlarger1}  again,   
 \beaa
 r|\chibh'| &\les& \frac{\ep_0}{u^{1+\dec}}, \qquad \mbox{on} \quad \Mext.
 \eeaa
 We infer that
 \beaa
 |\underline{\Theta}(u,\cdot)| &\les& \frac{\ep_0}{u^{1+\dec}}\textrm{ for all }1\leq u<+\infty.
 \eeaa

%%%%%%%%%%%%%%%%%%%%%% 
 
\subsubsection{A Bondi mass formula} 

%%%%%%%%%%%%%%%%%%%%%%

We use the following computation, see the proof of Lemma \ref{lemma:transportequationforHawkingmass}, 
\beaa
e_3'(\trch'\trchb') &=& -\trch'{\trchb'}^2 +2\trchb'\rho'+2\trch'\div'\xib' +2\trchb'\div'\eta'\\
&& +\trch'\left(  2 \xib'\c(\eta'-3\ze')-|\chibh'|^2\right)+\trchb'\left(  2|\eta'|^2 -\chibh'\c\chih'\right).
\eeaa
Together with \eqref{eq:behavioroftheintegrableframeforlarger1} \eqref{eq:behavioroftheintegrableframeforlarger2} \eqref{eq:behavioroftheintegrableframeforlarger3} and the estimates $\Nk^{(Dec)}_{k_{small}}\les \ep_0$, we deduce
 \beaa
\left|e_3'(m_H) + \frac{r}{32\pi}\int_S\trch'|\chibh'|^2\right| &\les& \frac{1}{r}
 \eeaa
 and hence
 \beaa
\left|e_3'(m_H) + \frac{1}{4|S|}\int_S|r\chibh'|^2\right| &\les& \frac{1}{r}.
 \eeaa 
Letting $r\to +\infty$ along  the leaves  of the $u$-foliation of $\Mext$, and using that the spheres at null infinity are round, we infer in view of the definition of $M_B$ and $\underline{\Theta}$
\beaa
e_3'(M_B)(u) &=& -\frac{1}{4}\int_{\mathbb{S}^2}|\underline{\Theta}|^2(u,\cdot)\textrm{ for all }1\leq u<+\infty.
\eeaa
Since 
\beaa
e_3'(u)=e_3(u)+O(r^{-1})=2+\widecheck{e_3(u)}+O(r^{-1})
\eeaa
and $e_3'$ is orthogonal to the spheres foliating $\II_+$, we infer $e_3'=(2+\widecheck{e_3(u)})\pr_u$. Thus, we obtain the following Bondi mass type formula 
\bea
\lab{eq:Bondi-masslaw}
\pr_uM_B(u) &=& -\frac{1}{8(1+\frac{1}{2}\widecheck{e_3(u)})}\int_{\mathbb{S}^2}\underline{\Theta}^2(u,\cdot)\textrm{ for all }1\leq u<+\infty,
\eea
with $\widecheck{e_3(u)}$ satisfying \eqref{eq:estimate1intheconclusionsforxibombandvsi}.

%%%%%%%%%%%%%%%%%

\subsubsection{Final Bondi mass} 

%%%%%%%%%%%%%%%%%

In view of the estimate
 \beaa
 |\underline{\Theta}(u,\cdot)| &\les& \frac{\ep_0}{u^{1+\dec}}\textrm{ for all }1\leq u<+\infty,
 \eeaa
and the control for $\widecheck{e_3(u) }$ in \eqref{eq:estimate1intheconclusionsforxibombandvsi}, we infer that
\beaa
|\pr_uM_B(u)| &\les& \frac{\ep_0^2}{u^{2+2\dec}}\textrm{ for all }1\leq u<+\infty.
\eeaa 
In particular, since $u^{-2-2\dec}$ is integrable, the limit along $\II_+$ exists
\beaa
M_B(+\infty)=\lim_{u\to +\infty}M_B(u)
\eeaa
and is the so-called final Bondi mass. 

Also, recall \eqref{eq:estimatebetweenmHandmforconclusions}
\beaa
m_H &=& m_\infty\Big(1+O(r^{-1})+O(\ep_0u^{-\frac{1}{2}-\dec})\Big).
\eeaa
Fixing $u$ and letting $r\to +\infty$, we infer on $\II_+$, in view of the definition of the Bondi mass,
\beaa
M_B(u) &=& m_\infty\Big(1+O(\ep_0u^{-\frac{1}{2}-\dec})\Big).
\eeaa
Then, letting $u\to +\infty$ along $\II_+$, and in view  of the definition of the final Bondi mass, we infer 
\beaa
m_\infty=M_B(+\infty),
\eeaa
i.e. the final mass $m_\infty$ coincides with the final Bondi mass.

%%%%%%%%%%%%%%%%%%%%%%%%%%%%%%%%%%%

\subsection{The final angular momentum $a_\infty$}

%%%%%%%%%%%%%%%%%%%%%%%%%%%%%%%%%%%

We exhibit below a geometric quantity converging to the final angular momentum $a_\infty$ along $\II_+$ as $u\to +\infty$. Relying on the frame transformations formulas of Proposition \ref{Proposition:transformationRicci} and the above control of $(f, \fb, \la)$ in \eqref{eq:changeofframecoefficientsfromPGtointegralframeconclusions}, have
\beaa
\curl'\b'=\curl\b -\frac{3m_\infty}{r^3}\curl(f)+r^{-3}\Ga_g+O(r^{-6}).
\eeaa
Using again \eqref{eq:changeofframecoefficientsfromPGtointegralframeconclusions} for $f$, this yields
\beaa
\curl'\b'=\curl\b -\frac{3m_\infty}{r^3}\curl(\nab u)+r^{-3}\Ga_g+O(r^{-6}).
\eeaa
We have
\beaa
\curl(\nab u) &=& \in_{ab}\nab_a\nab_b u=\in_{ab}\left(\D_a\D_bu -\frac{1}{2}\chi_{ab}e_3(u)-\frac{1}{2}\chib_{ab}e_4(u)\right)\\
&=& -\frac{1}{2}\atrch e_3(u) -\frac{1}{2}\atrchb e_4(u)
\eeaa
which together with the fact that $e_4(u)=0$ implies 
\beaa
\curl'\b' &=& \curl\b +\frac{3m_\infty}{2r^3}\atrch e_3(u)+r^{-3}\Ga_g+O(r^{-6}).
\eeaa
Recall, see Definition \ref{def:renormalizationofallnonsmallquantitiesinPGstructurebyKerrvalue} and the definitions of $ \Ga_g, \Ga_b$,
\beaa
\tr X=\trch-i \atrch=\frac{2}{q}+\Ga_g =\frac{2(r-i \cos\th)}{r^2+a^2}+\Ga_g, \qquad  e_3(u)=\frac{2(r^2+a^2)}{|q|^2}+ r\Ga_b.
\eeaa
Hence
\beaa
\curl'\b' &=& \curl\b +\frac{6a_\infty m_\infty\cos\th}{r^5}+r^{-3}\Ga_g+O(r^{-6}).
\eeaa
For a scalar function $h$ on $\Mext$, we introduce, see section  \ref{section:can.ell=1basis.Mext},
\beaa
(h)_{\ell=1} &:=& \frac{1}{|S|}\left(\int_SJ^{(0)}, \int_SJ^{(+)}, \int_SJ^{(-)}\right),
\eeaa
where 
\beaa
J^{(0)}=\cos\th, \qquad J^{(+)}=\sin\th\cos\vphi, \qquad J^{(-)}=\sin\th\sin\vphi.
\eeaa
Now, in view of Lemma \ref{lemma:propertiesintegralofJpandJpJqonMext}, we have on $\Mext$
\beaa
\int_{S}\Jp J^{(q)} &=& \frac{4\pi}{3}r^2\de_{pq}+O\left(1+\ep_0 ru^{-\frac{1}{2}-\dec}\right),
\eeaa
from which we infer
\bea
\lab{eq:conclusions-AngMom}
\bsplit
r^5(\curl'\b')_{\ell=1,0} &= r^5(\curl\b)_{\ell=1,0} +2a_\infty m_\infty+r^2\Ga_g+O(r^{-1}),\\
r^5(\curl'\b')_{\ell=1,\pm} &= r^5(\curl\b)_{\ell=1,\pm} +r^2\Ga_g+O(r^{-1}).
\end{split}
\eea
We  next  appeal to the first  identity \eqref{eq:Proposition.Identity.ell=1Div B1} of Proposition \ref{Proposition:Identity.ell=1Div B}
according to which, in view of   the estimate $\Nk^{(Dec)}_{k_{small}}\les \ep_0$ of the main theorem,
\beaa
\nab_4 \left(\int_{S(r,u)} \frac{rJ^{(0)}}{\Si}   [\ov{D}]_{ren}    \left( r^4[B]_{ren} \right)    \right)=  O\left(\frac{\ep_0}{r^{\frac{3}{2}}u^{\frac{1}{2}+\dec}}\right)
\eeaa
and
\beaa
\nab_4 \left(\int_{S(r,u)} \frac{rJ^{(\pm)}}{\Si}   [\ov{D}]_{ren}    \left( r^4[B]_{ren} \right)    \right)\mp  \frac{a}{r^2}  \int_{S(r,u)} \frac{rJ^{(\mp)}}{\Si}   [\ov{D}]_{ren}    \left( r^4[B]_{ren}    \right)  =  O\left(\frac{\ep_0}{r^{\frac{3}{2}}u^{\frac{1}{2}+\dec}}\right),
\eeaa
where
\beaa
\bsplit
[B]_{ren}&:=B - \frac{3a}{2}\ov{\Pc}\Jk -\frac{a}{4}\ov{\Jk}\c A, \qquad 
[\ov{D}]_{ren}:= \ov{\DD}\c -\frac{a}{2}\ov{\Jk}\c\nab_4 -\frac{a}{2}\ov{\Jk},
\end{split}
\eeaa
are renormalized quantities. Recalling that $\Si^2=(r^2+a^2)|q|^2+2mra^2(\sin\th)^2$, see  \eqref{eq:quantities-qDeSi}, and  the estimates for $A$, $B$, $P$ and   $\Ga_g$  provided by the main theorem, we deduce
\beaa
\nab_4\left(r^3\int_S\Big(\ov{\DD}\c B+r^{-3}\dk^{\leq 1}\Ga_g\Big)\Jp\right) &=& O\left(\frac{\ep_0}{r^{\frac{3}{2}}u^{\frac{1}{2}+\dec}}\right).
\eeaa
Integrating forward from $\TT$,  using again $\Nk^{(Dec)}_{k_{small}}\les \ep_0$, we obtain
\beaa
\sup_{\Mext}u^{\frac{1}{2}+\dec}r^3\left|\int_S(\ov{\DD}\c B)\Jp\right| &\les& \ep_0. 
\eeaa
Since $\ov{\DD}\c B=(\nab-i\dual \nab)\c(\b+i\dual\b)=2\div\b+2i\curl\b$, this yields
\beaa
\sup_{\Mext}u^{\frac{1}{2}+\dec}r^5|(\curl\b)_{\ell=1}| &\les& \ep_0
\eeaa
and hence, back to \eqref{eq:conclusions-AngMom},
\bea\lab{eq:firstestimateforcurlbetaprimeellequal1modeforconclusions}
\bsplit
r^5(\curl'\b')_{\ell=1,0} &= 2a_\infty m_\infty+O\left(\frac{1}{r}+\frac{\ep_0}{u^{\frac{1}{2}+\dec}}\right),\\
r^5(\curl'\b')_{\ell=1,\pm} &= O\left(\frac{1}{r}+\frac{\ep_0}{u^{\frac{1}{2}+\dec}}\right).
\end{split}
\eea

To derive a limit  for  $r^5(\curl'\b')_{\ell=1,0} $ as $r\to \infty$, at constant $u$, we  need 
to  estimate the quantity $e_4(r^5(\curl'\b')_{\ell=1})$. We start with the following Bianchi identity for the integrable frame, in view of  Proposition \ref{prop-nullstr} with $\atrch=0$,
\beaa
\nab_4'\b'+(2\trch'+2\om')\b' &=& \div'\a'+\a'\c(2\ze'+\etab')+3\rho'\xi'+3\rhod' \dual \xi'.
\eeaa
Differentiating with $\curl'$, this yields
\beaa
e_4'(\curl'\b') -[\nab_4', \curl']\b' + 2\trch'\curl'\b' &=& \curl'\div'\a' -2\nab'\trch'\c\dual\b'-\curl'(\om'\b')\\
&&+\curl'\Big(\a'\c(2\ze'+\etab')+3\rho'\xi'\Big).
\eeaa
By standard commutation formulas (see \cite{GKS1} or \cite{Ch-Kl})
\beaa
[\nab_4', \curl']\b' &=& -\frac{1}{2}\trch'\curl'\b' -\chih'\c\nab'\dual\b'+(\etab'+\ze')\c\nab_4'\dual\b'+\xi'\c\nab_3'\dual\b'\\
&& +\xi'\c\chib'\c\b'+\etab'\c\chi'\c\b',
\eeaa
we infer, together with  \eqref{eq:behavioroftheintegrableframeforlarger1} \eqref{eq:behavioroftheintegrableframeforlarger2} \eqref{eq:behavioroftheintegrableframeforlarger3},
\bea\lab{eq:intermediarytransportequationcurlbetatintegrableconclusionsforainfty}
\nn e_4'(\curl'\b')  + \frac{5}{2}\trch'\curl'\b' &=& \curl'\div'\a' +\xi'\c\nab_3'\dual\b'+3\curl'\big(\rho'\xi'\big)\\
&&+O\Big(r^{-7}+\ep_0 r^{-\frac{13}{2}}\Big).
\eea
The  estimate for $\xi' =O\big(r^{-3}\big)+r^{-1}\dk^{\leq 1}\Ga_b$ provided by \eqref{eq:behavioroftheintegrableframeforlarger2} is not enough to derive a suitable estimate for the term involving $3\curl'\big(\rho'\xi'\big)$ on the RHS of \eqref{eq:intermediarytransportequationcurlbetatintegrableconclusionsforainfty}. To derive a sharper estimate for $\xi'$, we use again   the frame transformations formulas of Proposition \ref{Proposition:transformationRicci} and the control \eqref{eq:changeofframecoefficientsfromPGtointegralframeconclusions} for $(f, \fb, \la)$. We obtain 
\beaa
\xi' &=& \frac{1}{2}\nab_4f+\frac{1}{4}\trch f+O(r^{-3})+r^{-2}\dk^{\leq 1}\Ga_b\\
&=& \frac{1}{2r}\nab_4(rf)+O(r^{-3})+r^{-2}\dk^{\leq 1}\Ga_b
\eeaa
where we have kept the explicit form of the a priori problematic term. We now show that $\nab_4(rf)$ behaves better than expected from  \eqref{eq:changeofframecoefficientsfromPGtointegralframeconclusions}. In view of the explicit formula for $f$
\beaa
  f= -\frac{4}{e_3(u)+ \sqrt{(e_3(u))^2 +4|\nab u|^2e_3(r)}}\nab u
  \eeaa
 in \eqref{def:transition-functs:ffbla}, we have\footnote{Note that  $e_3(u) = 2 + r\Ga_b$,  $e_3(r)=-1+ r\Ga_b +O(r^{-2} )$          and $\nab u =\Ga_b +O(r^{-1} )$.}
\beaa
\nab_4(rf) &=& -4\big(1+r\Ga_b+O(r^{-2})\big)\Big(\nab_4(r\nab u)-2e_4(e_3(u))r\nab u\Big)+O(r^{-4})+r^{-3}\Ga_b.
\eeaa
The desired gain comes from the  following transport equations in $e_4$ of Proposition \ref{prop:e_4(xyz)}
\beaa
\nab_4\DD u +\frac{1}{2}\tr X\DD u &=& -\frac{1}{2}\Xh\c\ov{\DD}u,\qquad 
e_4(e_3(u)) = -\Re\Big((Z+H)\c\ov{\DD} u\Big),
\eeaa
from which
\beaa
\nab_4(r\nab u) &=& O(r^{-2})+\Ga_g,\qquad 
e_4(e_3(u)) = O(r^{-2})+r^{-1}\Ga_b.
\eeaa
Therefore
\beaa
\nab_4(rf) &=& O(r^{-2})+\Ga_g,
\eeaa
i.e. $\nab_4(rf)$ behaves indeed better than expected from \eqref{eq:changeofframecoefficientsfromPGtointegralframeconclusions}. Plugging in the above formula for $\xi'$, we obtain the following improvement of \eqref{eq:behavioroftheintegrableframeforlarger2} for $\xi'$
\bea
\lab{eq:improvedxi'}
\xi' &=& O(r^{-3})+r^{-1}\Ga_g+r^{-2}\dk^{\leq 1}\Ga_b.
\eea
Plugging in \eqref{eq:intermediarytransportequationcurlbetatintegrableconclusionsforainfty}, and using \eqref{eq:behavioroftheintegrableframeforlarger3} for $\rho'$, and the Bianchi identity for $\nab_3'\b'$  (see   Proposition
\ref{prop-nullstr} with $\atrch=0$), we deduce
\bea
e_4'(\curl'\b')  + \frac{5}{2}\trch'\curl'\b' &=& \curl'\div'\a' +O\Big(r^{-7}+\ep_0 r^{-\frac{13}{2}}\Big).
\eea
On the other hand, using the fact that $ f=\Ga_b +O(r^{-1})$ and $ \la = 1+O(r^{-2})+r^{-1}\Ga_b$ in view of  \eqref{eq:changeofframecoefficientsfromPGtointegralframeconclusions}, we have
\beaa
e_4'(r^3) -\frac{3}{2}r^3\trch' &=& 3r^2\la\left(e_4+f\c\nab+\frac{1}{4}|f|^2e_3\right)r -\frac{3}{2}r^3\trch'\\
&=& 3r^2\la -\frac{3}{2}r^3\trch' +O(1)+r\Ga_b= -\frac{3}{2}r^3\left(\trch' -\frac{2}{r}\right)+O(1)+r\Ga_b\\
& =& r^3\dk^{\leq 1}\Ga_g +O(1)+r\Ga_b,
\eeaa
where we used \eqref{eq:behavioroftheintegrableframeforlarger1} for $\trch'$. Using also the control for $\b'$ in \eqref{eq:behavioroftheintegrableframeforlarger3}, 
i.e. $\b'=\b+ O(r^{-4})+r^{-2}\Ga_g,$
we deduce
\bea
\lab{eq:transportcurlb'}
e_4'(r^3\curl'\b')  + \trch' r^3\curl'\b' &=& r^3\curl'\div'\a' +O\Big(r^{-4}+\ep_0 r^{-\frac{7}{2}}\Big).
\eea
Next we appeal to    Lemma \ref{lemma:e_4(int_S h-chap3)}
\beaa
e_4\left(\int_{S(r,u)}h\right) &=& \int_{S(r,u)}\Big(e_4(h)+\trch' h+f\c\nab' h\Big)+O(r^{-1})h.
\eeaa
Since $e_4(\Jp)=0$ for $p=0,+,-$, we infer
\beaa 
e_4\left(\int_{S}r^3\curl'\b'\Jp\right) &=& \int_{S}\Big(e_4(r^3\curl'\b')+\trch'r^3\curl'\b'\Big)\Jp+O(r^3){\dkb'}^{\leq 1}\curl'\b'.
\eeaa
Since
\beaa
e_4'(r^3\curl'\b') &=& \la\left(e_4+f\c\nab+\frac{1}{4}|f|^2\right)(r^3\curl'\b')= e_4(r^3\curl'\b')+O(r)\dk^{\leq 1}\curl'\b',
\eeaa
we deduce
\beaa 
e_4\left(\int_{S}r^3\curl'\b'\Jp\right) &=&\int_Se_4'(r^3\curl'\b') +\trch'  (r^3\curl'\b')  +O(r^3)\dk^{\leq 1}\curl'\b'. 
\eeaa
Thus, in view of \eqref{eq:transportcurlb'},
\beaa
 e_4\left(\int_{S}r^3\curl'\b'\Jp\right)&=& r^3\int_{S}\curl'\div'\a'\Jp+O(r^3)\dk^{\leq 1}\curl'\b'+O\Big(r^{-2}+\ep_0 r^{-\frac{3}{2}}\Big).
\eeaa
Together with the control for  $\b'=\b+ O(r^{-4})+r^{-2}\Ga_g$ in \eqref{eq:behavioroftheintegrableframeforlarger3}  and recalling the definition of $\ddd_2, \ddd_1$
\beaa
e_4\left(\int_{S}r^3\curl'\b'\Jp\right)&=& r^3\int_{S}\ddd_1'\, \ddd_2' \a' \c  (0,       \Jp) +O\Big(r^{-2}+\ep_0 r^{-\frac{3}{2}}\Big).
\eeaa
Integrating by parts we deduce
\beaa 
e_4\left(\int_{S}r^3\curl'\b'\Jp\right) &=& r^3\int_{S}\a'\c\dds_2'\dds_1'(0,\Jp)+O\Big(r^{-2}+\ep_0 r^{-\frac{3}{2}}\Big).
\eeaa
Now, according to Proposition \ref{prop:controlofDDprimehotDDprimeJonMext}, we have
\beaa
|\dds_2'\dds_1'(0,\Jp)|\les \frac{\ep_0}{r^3u^{\frac{1}{2}+\dec}}+\frac{1}{r^4}.
\eeaa
Together with  the control for $\a'$ in \eqref{eq:behavioroftheintegrableframeforlarger3}      we deduce
\beaa 
e_4\left(\int_{S}r^3\curl'\b'\Jp\right) &=& O\Big(r^{-2}+\ep_0 r^{-\frac{3}{2}}\Big).
\eeaa
Together with the definition of $(\curl'\b')_{\ell=1}$, we finally obtain the following control for $e_4((\curl'\b')_{\ell=1})$
\beaa 
e_4\left((\curl'\b')_{\ell=1}\right) &=& O\Big(r^{-2}+\ep_0 r^{-\frac{3}{2}}\Big).
\eeaa
Since $r^{-3/2}$ in integrable, we infer from the above control of $e_4((\curl'\b')_{\ell=1})$ the existence of a limit to $(\curl'\b')_{\ell=1}$ as $r\to +\infty$ along the leaves of the $u$-foliation of $\Mext$,  for all $ 1\leq u<+\infty,$
\bea
\JJ_{\ell=1, p} (u):=\lim_{r\to +\infty}(\curl'\b')_{\ell=1, p}(u,r),  \qquad p=  0, +,-.
\eea
Fixing $u$ and letting $r\to +\infty$ in \eqref{eq:firstestimateforcurlbetaprimeellequal1modeforconclusions}, we infer on $\II_+$, in view of the definition of $\JJ_{\ell=1, p} (u)$,
\beaa
\bsplit
\JJ_{\ell=1, 0} (u)  &= 2a_\infty m_\infty+O\left(\frac{\ep_0}{u^{\frac{1}{2}+\dec}}\right),\\
\JJ_{\ell=1,\pm}(u)  &= O\left(\frac{\ep_0}{u^{\frac{1}{2}+\dec}}\right).
\end{split}
\eeaa
Then, letting $u\to +\infty$ along $\II_+$, we deduce that $\JJ_{\ell=1,p}(u)$ admit limits, which we denote by $\JJ_{\ell=1,p}(+\infty)$, i.e.
\beaa
\JJ_{\ell=1,p}(+\infty):=\lim_{u\to +\infty}\JJ_{1,p}(u), \qquad p=0,+,-,
\eeaa
and that these limits satisfy 
\beaa
a_\infty=\frac{1}{2m_\infty}\JJ_{\ell=1,0}(+\infty), \qquad \JJ_{\ell=1,\pm}(+\infty)  &= 0.
\eeaa

%%%%%%%%%%%%%%%%%%

\subsection{Other conclusions} 

%%%%%%%%%%%%%%%%%%%

%%%%%%%%%%%%%%%%%%%%%%%%%%%%%%%

\subsubsection{Coordinates systems on $\Mext$ and $\Mint$}

%%%%%%%%%%%%%%%%%%%%%%%%%%%%%%%

In view of Proposition \ref{prop:3coordinatessysteminMextchap4}, $\Mext$ is covered by three regular coordinates patches:
\begin{itemize}
\item in the $(u,r,\th,\vphi)$ coordinates system, we have, for $\frac{\pi}{4}<\th<\frac{3\pi}{4}$, 
\beaa
\g &=& \g_{a_\infty,m_\infty}+\Big(du, dr, rd\th, r\sin\th d\vphi\Big)^2O\left(\frac{\ep_0}{u^{1+\dec}}\right),
\eeaa

\item in the $(u,r,x^1,x^2)$ coordinates system,  with $x^1=J^{(+)}$ and $x^2=J^{(-)}$, we have, for $0\leq\th<\frac{\pi}{3}$ and for   $\frac{2\pi}{3}<\th\leq\pi$, 
\beaa
\g &=& \g_{a_\infty,m_\infty}+\Big(du, dr, rdx^1, rdx^2\Big)^2O\left(\frac{\ep_0}{u^{1+\dec}}\right),
\eeaa
\end{itemize}
where in each case, $\g_{a_\infty,m_\infty}$ denotes the Kerr metric expressed in the corresponding coordinates system of Kerr, see Lemma \ref{lemma:urthetavphicoordinatesinKerrchap2} and Lemma \ref{lemma:urthetaJplusJminuscoordinatesinKerrchap2}.

Also, in view of Proposition \ref{prop:3coordinatessysteminMintchap4}, $\Mint$ is covered by three regular coordinates patches:
\begin{itemize}
\item in the $(\ub,r,\th,\vphi)$ coordinates system, we have, for $\frac{\pi}{4}<\th<\frac{3\pi}{4}$, 
\beaa
\g &=& \g_{a_\infty,m_\infty}+\Big(d\ub, dr, rd\th, r\sin\th d\vphi\Big)^2O\left(\frac{\ep_0}{\ub^{1+\dec}}\right),
\eeaa

\item in the $(\ub,r,x^1,x^2)$ coordinates system, with $x^1=J^{(+)}$ and $x^2=J^{(-)}$, we have, for $0\leq\th<\frac{\pi}{3}$ and for   $\frac{2\pi}{3}<\th\leq\pi$, 
\beaa
\g &=& \g_{a_\infty,m_\infty}+\Big(d\ub, dr, rdx^1, rdx^2\Big)^2O\left(\frac{\ep_0}{\ub^{1+\dec}}\right),
\eeaa
\end{itemize}
where in each case, $\g_{a_\infty,m_\infty}$ denotes the Kerr metric expressed in the corresponding coordinates system of Kerr, i.e. the analog for ingoing PG structures of Lemma \ref{lemma:urthetavphicoordinatesinKerrchap2} and Lemma \ref{lemma:urthetaJplusJminuscoordinatesinKerrchap2}.

%%%%%%%%%%%%%%%%%%%%%%%%%%%%%%%

\subsubsection{Asymptotic of the future event horizon}

%%%%%%%%%%%%%%%%%%%%%%%%%%%%%%%

Let
\beaa
r_{\pm,\infty} &:=& m_\infty\pm\sqrt{m_\infty^2-a_\infty^2}.
\eeaa
We show below that $\HH_+$ is located in the following region of $\Mint$
\bea\lab{eq:asymptoticlocationofthehorizonintheconlusionofthemaintheorem}
r_{+,\infty}\left(1 -\frac{\sqrt{\ep_0}}{\ub^{1+\dec}}\right)\leq r\leq r_{+,\infty}\left(1 +\frac{\sqrt{\ep_0}}{\ub^{1+\frac{\dec}{2}}}\right)\textrm{ on }\HH_+\textrm{ for any }1\leq\ub<+\infty.
\eea

We consider first the lower bound. Let us denote by $(e_4, e_3, e_1, e_2)$ and $(\ub,r)$ the null frame and scalar functions associated to our ingoing PG structure of $\Mint$. The estimates $\Nk^{(Dec)}_{k_{small}}\les \ep_0$  imply
\beaa
\sup_{\Mint}\ub^{1+\dec}\left(\left|e_3(r)+1\right|+\left|e_4(r)-\frac{\De}{|q|^2}\right|\right) &\les&\ep_0.
\eeaa
In particular, we have for all $r\leq r_{+,\infty}\left(1 -\frac{\sqrt{\ep_0}}{\ub^{1+\dec}}\right)$
\beaa
\De &=& r^2-2m_\infty r+a_\infty^2=(r-r_{+,\infty})(r-r_{-,\infty})\\
&\leq& -\frac{\sqrt{\ep_0}}{2\ub^{1+\dec}}r_{+,\infty}\left(r_{+,\infty}-r_{-,\infty} -\frac{\sqrt{\ep_0}}{2\ub^{1+\dec}}\right)\\
&\les& -\frac{\sqrt{\ep_0}}{\ub^{1+\dec}}<0
\eeaa
and hence
\beaa
e_3(r)\leq -\frac{1}{2}<0, \qquad e_4(r)\les -\frac{\sqrt{\ep_0}}{\ub^{1+\dec}}<0 \quad\textrm{on}\quad \Mint\left(r\leq r_{+,\infty}\left(1 -\frac{\sqrt{\ep_0}}{\ub^{1+\dec}}\right)\right).
\eeaa
Consider now $\gamma(s)$ any future directed null geodesic emanating from a point of the region $\Mint\left(r\leq r_{+,\infty}\left(1 -\frac{\sqrt{\ep_0}}{\ub^{1+\dec}}\right)\right)$. $\dot{\ga}$ being a null vector, there exists at any point of $\ga(s)$ in $\Mint$ a scalar $\la$ and a 1-form $f$ such that 
\beaa
\dot{\ga} &=& \la\left(e_4+f^be_b+\frac{1}{4}|f|^2e_3\right),
\eeaa
where $\la>0$ since $\dot{\ga}$ is future directed. Since $\nab(r)=0$, we infer
\beaa
\frac{dr}{ds}=\D_{\dot{\ga}}r &=& \la\left(e_4+f^be_b+\frac{1}{4}|f|^2e_3\right)r=\la\left(e_4(r)+\frac{1}{4}|f|^2e_3(r)\right).
\eeaa
 Since $e_3(r)<0$ and $e_4(r)<0$ in $\Mint\left(r\leq r_{+,\infty}\left(1 -\frac{\sqrt{\ep_0}}{\ub^{1+\dec}}\right)\right)$ in view of the above, and since $\la>0$ and $|f|^2\geq 0$, we deduce that $r$ decreases along $\ga(s)$ so that $\ga(s)$ either stays in $\Mint\left(r\leq r_{+,\infty}\left(1 -\frac{\sqrt{\ep_0}}{\ub^{1+\dec}}\right)\right)$ for all $s$, or exists $\MM$ through $r=r_+(1-\deh)$. Thus,  $\Mint\left(r\leq r_{+,\infty}\left(1 -\frac{\sqrt{\ep_0}}{\ub^{1+\dec}}\right)\right)$ lies strictly inside the black hole and hence, $\HH_+$ must lie in $r\geq r_{+,\infty}\left(1 -\frac{\sqrt{\ep_0}}{\ub^{1+\dec}}\right)$ which concludes the lower bound. 
 
 Next, we focus on proving the upper bound. We need to show that any 2-sphere 
\bea\lab{eq:asymptoticspheretothehorizon}
S(\ub_1):=S\left(\ub_1, r=r_{+,\infty}\left(1 +\frac{\sqrt{\ep_0}}{\ub_1^{1+\frac{\dec}{2}}}\right)\right),\quad 1\leq \ub_1<+\infty,
\eea
is in the past of $\II_+$. Since $\Mext$ lies in the past of $\II_+$, is suffices to show that from any point of $S(\ub_1)$ there exists a future directed null geodesic reaching $\Mext$ in finite time. We will in fact show that the future directed null geodesics from $S(\ub_1)$ with initial speed $e_4$ reach $\Mext$ in finite time. 
Assume, by contradiction, that there exists a null geodesic  from $S(\ub_1)$ with initial speed $e_4$, denoted by $\ga$, that does not reach $\Mext$ in finite time. Let $e_4'$ be the geodesic generator of $\ga$. In view of Lemma \ref{Lemma:Generalframetransf} on general null frame transformation, and denoting by $(e_4, e_3, e_1, e_2)$ the null frame\footnote{Recall that we assume by contradiction that $\ga$ does not reach $\Mext$ and hence stays in $\Mint$.} of $\Mint$, we look for $e_4'$ under the form 
\beaa
e_4' &=& \la\left(e_4 + f^ae_a +\frac 1 4 |f|^2  e_3\right).
\eeaa
Also, let 
\beaa
F:=f+i\dual f. 
\eeaa
Then, the fact that $e_4'$ is geodesic implies the following transport equations along $\ga$ for $F$ and $\la$ in view of  Corollary \ref{cor:transportequationine4forchangeofframecoeffinformFFbandlamba}
\beaa
\nab_{\la^{-1}e_4'}F+\frac{1}{2}\tr X F+2\om F &=& -2\Xi -\chih\c F+E_1(f, \Ga),\\
\la^{-1}\nab_4'(\log\la) &=& 2\om+f\c(\ze-\etab)+E_2(f, \Ga),
\eeaa
where $E_1(f, \Ga)$ and $E_2(f, \Ga)$ contain expressions of the type $O(\Ga f^2)$ with no derivatives and $\Ga$ denotes the Ricci coefficients w.r.t. the original null frame $(e_3, e_4, e_1, e_2)$ of $\Mint$.  

We then proceed as follows
\begin{enumerate}
\item First, since $e_4$ is the initial speed of $\ga$ on $S(\ub_1)$, $f$ and $\la$ satisfy 
\beaa
f=0,\quad \la=1\textrm{ on }\ga\cap S(\ub_1).
\eeaa

\item Then, we initiate a continuity argument by assuming for some 
$$\ub_1<\ub_2<\ub_1+\left(\frac{\ub_1}{\ep_0}\right)^{\frac{\dec}{2}}$$ 
that we have
\bea\lab{eq:bootstrapasumptiononfandUpfortheasymptoticlocationofthehorizon}
|f|\leq \frac{\sqrt{\ep_0}}{\ub_1^{\frac{1}{2}+\dec}},\quad \frac{\De}{|q|^2}\geq \frac{c_{\infty}\sqrt{\ep_0}}{2\ub_1^{1+\frac{\dec}{2}}},\quad 0<\la<+\infty\textrm{ on }\ga(\ub_1, \ub_2)\cap\Mint
\eea
where $\ga(\ub_1, \ub_2)$ denotes the portion of $\ga$ in $\ub_1\leq \ub\leq\ub_2$ and where the strictly positive constant $c_{\infty}$ is given by 
\beaa
c_{\infty}:=\frac{r_{+,\infty}(r_{+,\infty}-r_{-,\infty})}{(r_{+,\infty})^2+a_\infty^2}.
\eeaa

\item We have
\beaa
\la^{-1}e_4'(\ub) &=& e_4(\ub)+f\c\nab(\ub)+\frac{1}{4}|f|^2e_3(\ub)=\frac{2(r^2+a^2)}{|q|^2}+\widecheck{e_4(\ub)}+f\c\nab(\ub).
\eeaa
Relying on our control of the ingoing geodesic foliation of $\Mint$, the above assumption for $f$ and the transport equation for $F$, we obtain on $\ga(\ub_1, \ub_2)\cap\Mint$
\beaa
\sup_{\ga(\ub_1, \ub_2)\cap \Mint}|f| &\les& \frac{\ep_0}{\ub_1^{1+\dec}}(\ub_2-\ub_1)\\
&\les& \frac{\ep_0^{1-\frac{\dec}{2}}}{\ub_1^{1+\frac{\dec}{2}}}
\eeaa
which improves our assumption in \eqref{eq:bootstrapasumptiononfandUpfortheasymptoticlocationofthehorizon} on $f$.

\item We have in view of the control of $f$
\beaa
\la^{-1}e_4'(r) &=& e_4(r)+\frac{1}{4}|f|^2e_3(r)=\frac{\De}{|q|^2}+O\left(\frac{\ep_0}{\ub_1^{1+\dec}}\right),\\
\la^{-1}e_4'(\cos\th) &=& e_4(\cos\th)+f\c\nab(\cos\th)=O\left(\frac{\sqrt{\ep_0}}{\ub_1^{\frac{1}{2}+\dec}}\right).
\eeaa
This yields
\beaa
\la^{-1}e_4'\left(\log\left(\frac{\De}{|q|^2}\right)\right) &=& \la^{-1}e_4'(r)\pr_r\log\left(\frac{\De}{|q|^2}\right) -\la^{-1}e_4'(\cos\th)\pr_{\cos\th}\log(|q|^2)\\
&=& \pr_r\left(\frac{\De}{|q|^2}\right)+\frac{|q|^2}{\De}O\left(\frac{\ep_0}{\ub_1^{1+\dec}}\right)+O\left(\frac{\sqrt{\ep_0}}{\ub_1^{\frac{1}{2}+\dec}}\right).
\eeaa
Thanks to our assumption on the lower bound of $\frac{\De}{|q|^2}$, we infer
\beaa
\la^{-1}e_4'\left(\log\left(\frac{\De}{|q|^2}\right)\right) &=& \pr_r\left(\frac{\De}{|q|^2}\right)(1+O(\sqrt{\ep_0}))
\eeaa
and since we are in $\Mint$, $r\geq r_0$ and hence
\beaa
\la^{-1}e_4'\left(\log\left(\frac{\De}{|q|^2}\right)\right) &\geq & \frac{m_0}{2r_0^2}.
\eeaa
Integrating from $\ub=\ub_1$, we deduce
\beaa
\frac{\De}{|q|^2} &\geq& \big(1+O(\sqrt{\ep_0})\big)\frac{c_{\infty}\sqrt{\ep_0}}{\ub_1^{1+\frac{\dec}{2}}}\exp\left(\frac{m_0}{2r_0^2}(\ub-\ub_1)\right)
\eeaa
which is an improvement of our assumption in \eqref{eq:bootstrapasumptiononfandUpfortheasymptoticlocationofthehorizon} on $\frac{\De}{|q|^2}$. 

\item In view of the control of $f$ and of the ingoing geodesic foliation of $\Mint$, we rewrite the transport equation for $\la$ as
\beaa
\la^{-1}\nab_4'(\log\la) &=& 2\om+f\c(\ze-\etab)+E_2(f, \Ga)\\
&=& -\pr_r\left(\frac{\De}{|q|^2}\right)+O\left(\frac{\ep_0^{1-\frac{\dec}{2}}}{\ub_1^{1+\frac{\dec}{2}}}\right).
\eeaa
Since we have obtained above the other hand 
\beaa
\la^{-1}e_4'\left(\log\left(\frac{\De}{|q|^2}\right)\right) &=& \pr_r\left(\frac{\De}{|q|^2}\right)(1+O(\sqrt{\ep_0}))
\eeaa
we immediately infer
\beaa
\la^{-1}e_4'\left(\log\left(\la\left(\frac{\De}{|q|^2}\right)^2\right)\right)>0,\quad \la^{-1}e_4'\left(\log\left(\la\sqrt{\frac{\De}{|q|^2}}\right)\right)<0.
\eeaa
Integrating from $\ub=\ub_1$, this yields
\beaa
 \left(\big(1+O(\sqrt{\ep_0})\big)\frac{c_{\infty}\sqrt{\ep_0}}{\ub_1^{1+\frac{\dec}{2}}}\right)^2\left(\frac{\De}{|q|^2}\right)^{-2}\leq\la\leq \left(\big(1+O(\sqrt{\ep_0})\big)\frac{c_{\infty}\sqrt{\ep_0}}{\ub_1^{1+\frac{\dec}{2}}}\right)^\frac{1}{2}\left(\frac{\De}{|q|^2}\right)^{-\frac{1}{2}}.
\eeaa
Since $\frac{\De}{|q|^2}$ has an explicit lower bounded in view of our previous estimate, as well as an explicit upper bound since we are in $\Mint$, this yields an improvement of our assumptions in \eqref{eq:bootstrapasumptiononfandUpfortheasymptoticlocationofthehorizon} for $\la$.

\item Since we have improved all our bootstrap assumptions \eqref{eq:bootstrapasumptiononfandUpfortheasymptoticlocationofthehorizon}, we infer by a continuity argument the following bound
\beaa
\frac{\De}{|q|^2} &\geq& \big(1+O(\sqrt{\ep_0})\big)\frac{c_{\infty}\sqrt{\ep_0}}{\ub_1^{1+\frac{\dec}{2}}}\exp\left(\frac{m_0}{2r_0^2}(\ub-\ub_1)\right)\textrm{ on }\ga\left(\ub_1, \ub_1+\left(\frac{\ub_1}{\ep_0}\right)^{\frac{\dec}{2}}\right)\cap\Mint.
\eeaa
Now, in this $\ub$ interval, we may choose
\beaa
\ub_3:=\ub_1+\frac{2r_0^2}{m_0}\left(1+\frac{\dec}{2}\right)\log\left(\frac{\ub_1}{\ep_0}\right)
\eeaa
for which we have $\frac{\De}{|q|^2}\geq 1$. This is a contradiction since $\frac{\De}{|q|^2}<1$ in $\Mint$ (and even in $\MM$). Thus, we deduce that $\ga$ reaches $\Mext$ before $\ub=\ub_3$, a contradiction to our assumption on $\ga$. This concludes the proof of  \eqref{eq:asymptoticlocationofthehorizonintheconlusionofthemaintheorem}. 
\end{enumerate}

%%%%%%%%%%%%%%%%%%%%%%%%%%%%%%%%%%%%%%%%%%%%%%%%%%%%%

\section{Structure of the rest of the paper}

%%%%%%%%%%%%%%%%%%%%%%%%%%%%%%%%%%%%%%%%%%%%%%%%%%%%%

The rest of this paper is devoted to the proof of Theorem M0-M8. Note that the following results will be proved in separate papers:
\begin{itemize}
\item Theorem M1 concerning decay estimates for $\qf$ and $A$, and Theorem M2 concerning decay estimates for $\Ab$, will be proved in \cite{GKS2}. 

\item The control of the curvature components in Theorem M8, concerning top order boundedness estimates, will be proved in \cite{KS:Kerr-B}.  
\end{itemize}

In this paper, we prove the remaining results, i.e. Theorem M0 and Theorems M3-M7, as well as\footnote{Assuming the control of top order derivatives of  the curvature components in \cite{KS:Kerr-B}.} Theorem M8. More precisely:
\begin{enumerate}
\item We prove first consequences of the bootstrap assumptions, i.e. the control of coordinates systems and of the global frame of section \ref{section-globaleframe}, in Chapter 4.

\item Theorem M3 is proved in Chapter 5.

\item Theorem M4 is proved in Chapter 6.

\item Theorem M5 is proved in Chapter 7.

\item Theorems M0, M6 and M7 are proved in Chapter 8.

\item Theorem M8 is proved in Chapter 9 assuming the control of the curvature components in \cite{KS:Kerr-B}.
\end{enumerate}

\begin{remark}\lab{rmk:whyisTheoremM0onlyinchapter8}
Note that Theorem M0 should be proved first but has been postponed to chapter 8 for convenience as its proofs has a similar flavor to the ones of  Theorems M6 and M7. In particular, while its proof is in chapter 8, it relies only on the bootstrap assumptions and on the assumptions on the initial data layer.
\end{remark}

%%%%%%%%%%%%%%%%%%%%%%%%%%%%%%%%%%%

%%%%%%%%%%%%%%%%%%%%%%%%%%%%%%%%%%%%%%%%

\chapter{First consequences of the bootstrap assumptions}
\lab{Chapter:Firstconsequences-BA}

%%%%%%%%%%%%%%%%%%%%%%%%%%%%%%%%%%%%%%%%

In this chapter, we derive first consequences of our bootstrap assumptions on decay and boundedness. We start with the control of coordinates systems on $\Mext$ and $\Mint$ in section \ref{sec:controlofcoordinatessystems:chap4}. Then, we prove Proposition \ref{prop:constructionsecondframeinMext} on the construction of a second frame of $\Mext$ in section \ref{sec:proofofprop:constructionsecondframeinMext}. Finally, we prove Proposition \ref{prop:existenceandestimatesfortheglobalframe:bis} on the construction of a global frame on $\MM$ in section 
\ref{sec:proofofprop:existenceandestimatesfortheglobalframe:bis}.

%%%%%%%%%%%%%%%%%%%%%%%%%%%%%%%%%%%%%%%%

\section{Control of coordinates systems}
\lab{sec:controlofcoordinatessystems:chap4}

%%%%%%%%%%%%%%%%%%%%%%%%%%%%%%%%%%%%%%%%

In this section, we show that the metric in $\Mext$ and $\Mint$ is close to Kerr in 
suitable coordinates systems by relying on the bootstrap assumptions.  

\begin{proposition}\lab{prop:3coordinatessysteminMextchap4}
$\Mext$ is covered by three regular coordinates patches:
\begin{itemize}
\item in the $(u,r,\th,\vphi)$ coordinates system, we have, for $\frac{\pi}{4}<\th<\frac{3\pi}{4}$, 
\beaa
\g &=& \g_{a,m}+\Big(du, dr, rd\th, r\sin\th d\vphi\Big)^2O\left(\frac{\ep}{u^{1+\dec}}\right),
\eeaa

\item in the $(u,r,x^1,x^2)$ coordinates system,  with $x^1=J^{(+)}$ and $x^2=J^{(-)}$, we have, for $0\leq\th<\frac{\pi}{3}$ and for   $\frac{2\pi}{3}<\th\leq\pi$, 
\beaa
\g &=& \g_{a,m}+\Big(du, dr, rdx^1, rdx^2\Big)^2O\left(\frac{\ep}{u^{1+\dec}}\right),
\eeaa
\end{itemize}
where in each case, $\g_{a,m}$ denotes the Kerr metric expressed in the corresponding coordinates system of Kerr, see Lemma \ref{lemma:urthetavphicoordinatesinKerrchap2} and Lemma \ref{lemma:urthetaJplusJminuscoordinatesinKerrchap2}.
\end{proposition}

\begin{proposition}\lab{prop:3coordinatessysteminMintchap4}
$\Mint$ is covered by three regular coordinates patches:
\begin{itemize}
\item in the $(\ub,r,\th,\vphi)$ coordinates system, we have, for $\frac{\pi}{4}<\th<\frac{3\pi}{4}$, 
\beaa
\g &=& \g_{a,m}+\Big(d\ub, dr, rd\th, r\sin\th d\vphi\Big)^2O\left(\frac{\ep}{\ub^{1+\dec}}\right),
\eeaa

\item in the $(\ub,r,x^1,x^2)$ coordinates system, with $x^1=J^{(+)}$ and $x^2=J^{(-)}$, we have, for $0\leq\th<\frac{\pi}{3}$ and for   $\frac{2\pi}{3}<\th\leq\pi$, 
\beaa
\g &=& \g_{a,m}+\Big(d\ub, dr, rdx^1, rdx^2\Big)^2O\left(\frac{\ep}{\ub^{1+\dec}}\right),
\eeaa
\end{itemize}
where in each case, $\g_{a,m}$ denotes the Kerr metric expressed in the corresponding coordinates system of Kerr, i.e. the analog for ingoing PG structures of Lemma \ref{lemma:urthetavphicoordinatesinKerrchap2} and Lemma \ref{lemma:urthetaJplusJminuscoordinatesinKerrchap2}.
\end{proposition}

The proof of Proposition \ref{prop:3coordinatessysteminMextchap4} and Proposition \ref{prop:3coordinatessysteminMintchap4} being similar, we focus in the rest of this section on the proof of Proposition \ref{prop:3coordinatessysteminMextchap4}.

It will be in fact easier to control first the coefficients of the inverse metric. To this end, we rely on the following simple  lemma.
\begin{lemma}\lab{lemma:basicandgeneralformulaforinversemetricincoordinatessystem}
In a coordinates system $(x^\a)$, we have
\beaa
\g^{\a\b} &=& -\frac{1}{2}e_4(x^\a)e_3(x^\b)-\frac{1}{2}e_3(x^\a)e_4(x^\b)+\nab(x^\a)\c\nab(x^\b).
\eeaa
\end{lemma}

\begin{proof}
For a scalar function $h$, we have
\beaa
\D h &=& \g^{\a\b}e_\a(h) e_\b=-\frac{1}{2}e_4(h)e_3 -\frac{1}{2}e_3(h)e_4+\nab(h).
\eeaa
We infer
\beaa
\g(\D x^\a, \D x^\b) &=& -\frac{1}{2}e_4(x^\a)e_3(x^\b)-\frac{1}{2}e_3(x^\a)e_4(x^\b)+\nab(x^\a)\c\nab(x^\b). 
\eeaa
Since
\beaa
\g(\D x^\a, \D x^\b) &=& \g\left(\g^{\a\mu}\pr_{x^\mu}, \g^{\b\nu}\pr_{x^\nu}\right)= \g^{\a\mu}\g^{\b\nu}\g_{\mu\nu}=\g^{\a\b},
\eeaa
we deduce
\beaa
\g^{\a\b} &=& -\frac{1}{2}e_4(x^\a)e_3(x^\b)-\frac{1}{2}e_3(x^\a)e_4(x^\b)+\nab(x^\a)\c\nab(x^\b)
\eeaa
as stated.
\end{proof}

\begin{lemma}\lab{eq:inversemetricinperturbatonofKerrMext}
We have in $\Mext$
\beaa
\g^{rr} = -e_3(r),\qquad \g^{r\a}=-\frac{1}{2}e_3(x^\a)\quad\textrm{for}\quad x^\a=u, \,\th, \,\vphi, \, x^1, \, x^2,
\eeaa
and
\beaa
\g^{\a\b} &=& \nab(x^\a)\c\nab(x^\b)\quad\textrm{for}\quad x^\a,\, x^\b=u, \,\th, \,\vphi, \, x^1, \, x^2.
\eeaa
\end{lemma}

\begin{proof}
Recall that $(u,r,\th,\vphi)$ and $J^{(\pm)}$ verify on $\Mext$
\beaa
e_4(r)=1, \qquad \nab(r)=0, \qquad e_4(u)=e_4(\th)=e_4(\vphi)=e_4(J^{(+)})=e_4(J^{(-)})=0.
\eeaa
In view of Lemma \ref{lemma:basicandgeneralformulaforinversemetricincoordinatessystem}, we infer
\beaa
\g^{rr} = -e_3(r),\qquad \g^{r\a}=-\frac{1}{2}e_3(x^\a)\quad\textrm{for}\quad x^\a=u, \,\th, \,\vphi, \, x^1, \, x^2,
\eeaa
and
\beaa
\g^{\a\b} &=& \nab(x^\a)\c\nab(x^\b)\quad\textrm{for}\quad x^\a,\, x^\b=u, \,\th, \,\vphi, \, x^1, \, x^2,
\eeaa
as stated.
\end{proof}

\begin{corollary}\lab{cor:controlofinversemetricinMext}
In $\Mext$, we have
\beaa
&&\g^{rr} =\frac{\Delta}{|q|^2}+r\Ga_b,\quad \g^{ru}=-\frac{r^2+a^2}{|q|^2}+r\Ga_b,\quad \g^{rx^1}=\frac{ax^2}{|q|^2}+\Ga_b,\quad \g^{rx^2}=-\frac{ax^1}{|q|^2}+\Ga_b,\\
&& \g^{r\th} = (\sin\th)^{-1}\Ga_b,\quad \g^{r\vphi}= -\frac{a}{|q|^2}+(\sin\th)^{-2}\Ga_b,\\
&& \g^{uu} = \frac{a^2(\sin\th)^2}{|q|^2}+r^{-1}\Ga_b,\quad \g^{ux^1} = -\frac{ax^2}{|q|^2}+r^{-1}\Ga_b,\quad \g^{ux^2} = \frac{ax^1}{|q|^2}+r^{-1}\Ga_b,\\
&& \g^{u\th} = r^{-1}(\sin\th)^{-1}\Ga_b,\quad \g^{u\vphi} = \frac{a}{|q|^2}+r^{-1}(\sin\th)^{-2}\Ga_b,\\
&& \g^{x^1x^1} = \frac{(\cos\th\cos\vphi)^2+(\sin\vphi)^2}{|q|^2}+r^{-1}\Ga_b,\quad \g^{x^1x^2} = \frac{((\cos\th)^2-1)\sin\vphi\cos\vphi}{|q|^2}+r^{-1}\Ga_b,\\ 
&& \g^{x^2x^2} = \frac{(\cos\th\sin\vphi)^2+(\cos\vphi)^2}{|q|^2}+r^{-1}\Ga_b,\\
&& \g^{\th\th} = \frac{1}{|q|^2}+r^{-1}(\sin\th)^{-2}\Ga_b,\quad \g^{\th\vphi} = r^{-1}(\sin\th)^{-3}\Ga_b,\quad \g^{\vphi\vphi} = \frac{1}{|q|^2(\sin\th)^2}+r^{-1}(\sin\th)^{-4}\Ga_b.
\eeaa
\end{corollary}

\begin{proof}
Recall the notations
\beaa
\widecheck{e_3(r)} = e_3(r)+\frac{\Delta}{|q|^2},\qquad \widecheck{\nab u} = \nab u -a\Re(\Jk),\qquad \widecheck{e_3(u)} = e_3(u) -\frac{2(r^2+a^2)}{|q|^2},
\eeaa
and
\beaa
\widecheck{\nab J^{(0)}}=\nab J^{(0)}+\Im(\Jk), \qquad \widecheck{\nab J^{(\pm)}}=\nab J^{(\pm)}-\Re(\Jk_\pm), \qquad \widecheck{e_3(J^{(+)})}=e_3(J^{(\pm)})\pm\frac{2a}{|q|^2}J^{(\mp)}.
\eeaa
Also, recall that 
\beaa
&& \widecheck{e_3(r)}, \,\, \widecheck{e_3(u)}\in r\Ga_b, \qquad \qquad \widecheck{\nab u}\in \Ga_b,\\
&& e_3(J^{(0)}),\,\, \widecheck{e_3(J^{(\pm)})}\in\Ga_b, \qquad \widecheck{\nab J^{(p)}}\in\Ga_b,\,\, p=0,+,-.
\eeaa
We infer
\beaa
&& e_3(r)=-\frac{\Delta}{|q|^2}+r\Ga_b, \qquad e_3(u) =\frac{2(r^2+a^2)}{|q|^2}+r\Ga_b,\qquad \nab u =a\Re(\Jk)+\Ga_b, \\
&& \nab J^{(0)}=-\Im(\Jk)+\Ga_b, \qquad \nab J^{(\pm)}=\Re(\Jk_\pm)+\Ga_b,\\
&& e_3(J^{(0)})=\Ga_b, \qquad e_3(J^{(\pm)})=\mp\frac{2a}{|q|^2}J^{(\mp)}+\Ga_b. 
\eeaa

Also, we observe 
\beaa
&&\sin\th\nab(\th)=-\nab(J^{(0)}), \qquad \sin\th e_3(\th)=- e_3(J^{(0)}), \\
&&(\sin\th)^2\nab(\vphi) = -J^{(-)}\nab(J^{(+)})+J^{(+)}\nab(J^{(-)}), \\
&& (\sin\th)^2e_3(\vphi) = -J^{(-)}e_3(J^{(+)})+J^{(+)}\nab(J^{(-)}),
\eeaa  
so that we have in view of the above 
\beaa
&&\sin\th\nab(\th) = \Im(\Jk)+\Ga_b, \qquad (\sin\th)^2\nab(\vphi) = -J^{(-)}\Re(\Jk_+)+J^{(+)}\Re(\Jk_-)+\Ga_b,\\
&&\sin\th e_3(\th)=\Ga_b, \qquad (\sin\th)^2e_3(\vphi) =\frac{2a(\sin\th)^2}{|q|^2}+\Ga_b.
\eeaa
In view of Lemma \ref{eq:inversemetricinperturbatonofKerrMext}, the proof of the corollary follows then easily.
\end{proof}

We are now ready to prove Proposition \ref{prop:3coordinatessysteminMextchap4}.

\begin{proof}[Proof of Proposition \ref{prop:3coordinatessysteminMextchap4}]
In view of Corollary \ref{cor:controlofinversemetricinMext} for the inverse metric coefficients in $\Mext$, and  Lemma \ref{lemma:urthetavphicoordinatesinKerrchap2}  for the inverse metric coefficients of the Kerr metric, we have in the $(u,r,\th, \vphi)$ coordinates system of $\Mext$
\beaa
\g &=& \g_{a,m}+\Big(du, dr, rd\th, r\sin\th d\vphi\Big)^2(\sin\th)^{-2}r\Ga_b
\eeaa
and the conclusion follows from the control of $\Ga_b$ and the fact that $\sin\th> \frac{\sqrt{2}}{2}$ in the range  $\frac{\pi}{4}<\th<\frac{3\pi}{4}$. 

Similarly,  in view of Corollary \ref{cor:controlofinversemetricinMext} for the inverse metric coefficients in $\Mext$, and  Lemma \ref{lemma:urthetaJplusJminuscoordinatesinKerrchap2} for the inverse metric coefficients of the Kerr metric, we have in the $(u,r,x^1, x^2)$ coordinates system of $\Mext$
\beaa
\g &=& \g_{a,m}+\Big(du, dr, rdx^1, rdx^2\Big)^2(\cos\th)^{-2}r\Ga_b
\eeaa
and the conclusion follows from the control of $\Ga_b$ and the fact that $\cos\th> \frac{1}{2}$ in the range  $0\leq\th<\frac{\pi}{3}$ and $\cos\th<-\frac{1}{2}$ in the range  $\frac{2\pi}{3}<\th<\pi$. This concludes the proof of the Proposition \ref{prop:3coordinatessysteminMextchap4}.
\end{proof}

%%%%%%%%%%%%%%%%%%%%%%%%%%%%%%%%%%%%%%%%

\section{Proof of Proposition \ref{prop:constructionsecondframeinMext}}
\lab{sec:proofofprop:constructionsecondframeinMext}

%%%%%%%%%%%%%%%%%%%%%%%%%%%%%%%%%%%%%%%%

Let $(e_4, e_3, e_1, e_2)$ the outgoing PG frame of $\Mext$. We will exhibit another frame $(e_4', e_3', e_1', e_2')$ of $\Mext$  provided by
\bea\lab{eq:formulachangeofframeforsecondframeMext}
\begin{split}
e_4'&= e_4 + f^b e_b +\frac 1 4 |f|^2  e_3,\\
e_a'&= e_a +\frac 1 2 f_a e_3,\quad a=1,2,\\
e_3'&=  e_3,
\end{split}
\eea
where $f$ is such that 
\bea\lab{eq:propertiesoffforsecondframeofMext}
f=0\textrm{ on }S_*, \qquad \widecheck{\eta}'=0\textrm{ on }\Sigma_*, \qquad \xi'=0\textrm{ on }\Mext.
\eea
The desired estimates for the  Ricci coefficients and curvature components with respect to the new frame $(e_4', e_3', e_1', e_2')$ of $\Mext$ will be obtained in the region $\Mext(r\geq u^{\frac{1}{2}})$ using:
\begin{itemize}
\item the change of frame formulas of Proposition \ref{Proposition:transformationRicci}, applied to the change of frame from $(e_4, e_3, e_1, e_2)$ to $(e_4', e_3', e_1', e_2')$,

\item the estimates for $f$ on $\Mext$, and the fact that $\fb=0$ and $\la=1$ in the null frame transformation \eqref{eq:formulachangeofframeforsecondframeMext}, 

\item the estimates for the Ricci coefficients and curvature components with respect to the outgoing PG frame $(e_4, e_3, e_1, e_2)$ of $\Mext$ provided by the bootstrap assumptions on decay and boundedness.
\end{itemize}

Now, as it turns out, the frame $(e_4', e_3', e_1', e_2')$ does not satisfy the desired estimates in the region $\Mext(r\leq u^{\frac{1}{2}})$. We will thus introduce a third frame $(e_4'', e_3'', e_1'', e_2'')$ on $\Mext$, agreeing with $(e_4, e_3, e_1, e_2)$ on $\Mext(r\leq u^\frac{1}{2})$,  and with $(e_4', e_3', e_1', e_2')$ on $\Mext(r\geq u^{\frac{1}{2}})$, and satisfying all desired properties of Proposition \ref{prop:constructionsecondframeinMext}. 
In Steps 1--3 below, we study the properties of $(e_4', e_3', e_1', e_2')$. We then introduce the frame $(e_4'', e_3'', e_1'', e_2'')$ in Step 4,  and conclude,  in Steps 4--6,  the proof of Proposition \ref{prop:constructionsecondframeinMext}.

{\bf Step 1}. We start by deriving an equation for $f$ on $\Mext$. In view of the condition $\xi'=0$ on $\Mext$,  
see \eqref{eq:propertiesoffforsecondframeofMext}, in view of $\xi=\om=0$ satisfied by the outgoing PG structure of $\Mext$, and in view of Corollary \ref{cor:transportequationine4forchangeofframecoeffinformFFbandlamba}, we have
\bea\lab{eq:equationfore4primefforsecondframeofMext}
\nab_4'F+\frac{1}{2}\tr X F &=& -\chih\c F+E_1(f, \Ga)\,\,\textrm{ on }\Mext,
\eea
where $E_1(f, \Ga)$  contains expressions of the type $O(\Ga f^2)$ with no derivatives, and where
\beaa
F:=f+i\dual f.
\eeaa

We also derive an equation for $f$ on $\Sigma_*$. In view of  the change of frame formulas of Proposition \ref{Proposition:transformationRicci} in the particular case where $\la=1$ and $\fb=0$, we have on $\Mext$
\beaa 
\nab_3'F &=& 2H' -2H  +2\omb F +\err[\nab_3F],
\eeaa
where the lower order term $\err[\nab_3F]$ contains expressions of the type $O(f^2\Ga_b)$ with no derivatives. Now, we have\footnote{Recall that for the second frame, the Ricci coefficients and curvature components are also linearized using the scalar function $r$ and $\th$ and the complex 1-form $\Jk$ attached to the principal frame of $\Mext$.}
\beaa
H = \frac{aq}{|q|^2}\Jk+\Hc, \qquad H' =  \frac{aq}{|q|^2}\Jk+\Hc',
\eeaa
and hence
\beaa
H'-H &=&  \Hc'-\Hc,
\eeaa
so that, together with the condition $\Hc'=0$ on $\Sigma_*$,  
see \eqref{eq:propertiesoffforsecondframeofMext}, we infer
\bea\lab{eq:equationfore3primefforsecondframeofMext}
\nab_3'F &=&  -2\Hc  +2\omb F +\err[\nab_3F]\quad\textrm{on}\quad\Sigma_*.
\eea
Now, since $u+r$ is constant on $\Sigma_*$, the following vectorfield 
\beaa
\nu_{\Si_*}' := e_3'+b'e_4', \qquad b' := -\frac{e_3'(u+r)}{e_4'(u+r)},
\eeaa
is tangent to $\Sigma_*$. We compute in view of the above
\beaa
\nab_{\nu_{\Si_*}'}'F &=& \nab_3'F+b'\nab_4'F\\
&=&  -2\Hc  +2\omb F +\err[\nab_3F] +b'\left(-\frac{1}{2}\tr X F -\chih\c F + E_1(f, \Ga)\right).
\eeaa
Using \eqref{eq:formulachangeofframeforsecondframeMext}, as well as 
\beaa
e_4(r)=1, \qquad e_4(u)=0, \qquad \nab(r)=0,
\eeaa
we have
\beaa
b' &=& -\frac{e_3'(u+r)}{e_4'(u+r)}=  -\frac{e_3(u+r)}{\left(e_4 + f\c\nab +\frac 1 4 |f|^2  e_3\right)(u+r)}\\
&=& -\frac{e_3(u)+e_3(r)}{1+f\c\nab(u) +\frac 1 4 |f|^2(e_3(u)+e_3(r))}.
\eeaa
Recalling the following linearizations
\beaa
\widecheck{e_3(r)} = e_3(r)+\frac{\Delta}{|q|^2},\qquad \widecheck{\DD u} = \DD u -a\Jk,\qquad \widecheck{e_3(u)} = e_3(u) -\frac{2(r^2+a^2)}{|q|^2},
\eeaa
we deduce
\beaa
b' &=&  -\frac{1+\frac{2mr+a^2(\sin\th)^2}{|q|^2}+\widecheck{e_3(r)}+\widecheck{e_3(u)}}{1+f\c\nab(u) +\frac 1 4 |f|^2(e_3(u)+e_3(r))}
\eeaa
and hence
\bea\lab{eq:equationsatsfiedbyfonSigma*forsecondframeMext}
\nab_{\nu_{\Si_*}'}'F &=&   -2\Hc  +2\omb F +\err[\nab_3F]  \\
\nn&& -\frac{1+\frac{2mr+a^2(\sin\th)^2}{|q|^2}+\widecheck{e_3(r)}+\widecheck{e_3(u)}}{1+f\c\nab(u) +\frac 1 4 |f|^2(e_3(u)+e_3(r))}\left(-\frac{1}{2}\tr X F -\chih\c F + E_1(f, \Ga)\right)\,\,\textrm{on}\,\,\Si_*. 
\eea

{\bf Step 2}. Next, we estimate $f$ on $\Sigma_*$. In view of \eqref{eq:constraintsonthemainsmallconstantsepanddelta} and \eqref{eq:choiceksmallmaintheorem}, we have
\beaa
\dec(k_{large}-k_{small})\geq \frac{1}{2}\dec k_{large}-\dec \gg 1,
\eeaa
and we may thus assume from now on
\beaa
\frac{\dec}{3}(k_{large}-k_{small})\geq 130. 
\eeaa
Also, we introduce  a small constant $\de_0>0$ satisfying 
\beaa
\de_0=\frac{130}{k_{large}-k_{small}}\leq \frac{\dec}{3}.
\eeaa
Note in particular, from the bootstrap assumptions on decay and boundedness for the outgoing PG frame $(e_4, e_3, e_1, e_2)$ of $\Mext$ that 
\beaa
|\dk^{\leq k_{small}+130}\Ga_g| &\les&  |\dk^{\leq k_{small}}\Ga_g|^{1-\frac{130}{k_{large}-k_{small}}}|\dk^{\leq k_{large}}\Ga_g|^{\frac{130}{k_{large}-k_{small}}} \\
&\les& \min\left[\frac{\ep}{r^2}\left(\frac{1}{u^{\frac{1}{2}+\dec}}\right)^{1-\frac{130}{k_{large}-k_{small}}}, \,\frac{\ep}{r}\left(\frac{1}{u^{1+\dec}}\right)^{1-\frac{130}{k_{large}-k_{small}}}\right]\\
|\dk^{\leq k_{small}+130}\Ga_b| &\les&  |\dk^{\leq k_{small}}\Ga_b|^{1-\frac{130}{k_{large}-k_{small}}}|\dk^{\leq k_{large}}\Ga_b|^{\frac{130}{k_{large}-k_{small}}} \\
&\les& \frac{\ep}{r}\left(\frac{1}{u^{1+\dec}}\right)^{1-\frac{130}{k_{large}-k_{small}}}.
\eeaa
Since we have, in view of the definition of $\de_0$, 
\beaa
\frac{130}{k_{large}-k_{small}}(1+\dec)  &=&  \left(1+\dec\right)\de_0\leq 2\de_0, \\ 
\frac{130}{k_{large}-k_{small}}\left(\frac{1}{2}+\dec\right) &=&  \left(\frac{1}{2}+\dec\right)\de_0\leq \de_0,
\eeaa
we infer
\bea\lab{eq:estimateforframeMexttobeusedlater:chap4secondframeMext}
\nn\sup_{\Mext}\Big(r^2u^{\frac{1}{2}+\dec-2\de_0}+ru^{1+\dec-2\de_0}\Big)|\dk^{\leq k_{small}+130}\Ga_g| \\
+\sup_{\Mext}ru^{1+\dec-2\de_0}|\dk^{\leq k_{small}+130}\Ga_b| &\les& \ep.
\eea

Next, we assume the following local bootstrap assumption for $f$ on $\Si_*$
\bea\lab{eq:localbootstraponSigma*forsecondframeMext}
|\dk^{\leq k_{small}+130}f| \leq \frac{\sqrt{\ep}}{ru^{\frac{1}{2}+\dec-2\de_0}}\,\,\,\,\textrm{ on }u_1\leq u\leq u_*
\eea
where 
\beaa
1\leq u_1<u_*.
\eeaa
Since $f=0$ on $S_*$ in view of  \eqref{eq:propertiesoffforsecondframeofMext}, \eqref{eq:localbootstraponSigma*forsecondframeMext} holds for $u_1$ close enough to $u_*$, and our goal is to prove that we may in fact choose $u_1=1$ and replace $\sqrt{\ep}$ with $\ep$ in \eqref{eq:localbootstraponSigma*forsecondframeMext}. 

In view of  the bootstrap assumptions on boundedness  for the Ricci coefficients and curvature components with respect to the outgoing principal frame $(e_4, e_3, e_1, e_2)$ of $\Mext$, \eqref{eq:equationsatsfiedbyfonSigma*forsecondframeMext} yields
\beaa
\nab_{\nu_{\Si_*}'}'F &=&   -2\Hc +h, \qquad |\dk^kh| \les r^{-1}(|\dk^{\leq k}f|+|\dk^{\leq k}f|^4)\textrm{ for }k\leq k_{large}.
\eeaa
To differentiate this transport equation, we introduce the derivation 
\beaa
\widetilde{\nab} &:=& \nab' -\frac{\nab'(u+r)}{e_4'(u+r)}\nab_4'
\eeaa
so that $(\nab_{\nu_{\Si_*}'}', \widetilde{\nab})$ span all tangential derivatives to $\Si_*$. Note that 
\beaa
\widetilde{\nab} &=& \nab' -\frac{\nab(u+r)+\frac 1 2 f e_3(u+r)}{e_4(u+r)+f\c\nab(u+r)+\frac{1}{4}|f|^2e_3(u+r)}\nab_4'\\
&=& \nab' -\frac{\nab(u)+\frac 1 2 f e_3(u+r)}{1+f\c\nab(u)+\frac{1}{4}|f|^2e_3(u+r)}\nab_4'
\eeaa
 We introduce the following weighted derivatives on $\Si_*$:
 \beaa
 \widetilde{\dkb}:=r\widetilde{\nab}, \qquad \widetilde{\dk}=(\nab_{\nu_{\Si_*}'}', \widetilde{\dkb}).
 \eeaa
Using commutator identities, using also \eqref{eq:equationfore4primefforsecondframeofMext} and \eqref{eq:equationfore3primefforsecondframeofMext}, and in view of \eqref{eq:localbootstraponSigma*forsecondframeMext}, we infer
\beaa
|\nab_{\nu_{\Si_*}'}'\widetilde{\dkb}^kF| &\les& |\dk^{\leq k}\Hc| + \frac{\sqrt{\ep}}{r^2u^{\frac{1}{2}+\dec-2\de_0}}\textrm{ for }k\leq k_{small}+130, \,\, u_1\leq u\leq u_*.
\eeaa

Since $f=0$ on $S_*$ in view of  \eqref{eq:propertiesoffforsecondframeofMext}, and since $\nu_{\Si_*}'$ is tangent to $\Si_*$, we deduce on $\Si_*$, integrating along the integral curve of $\nu_{\Si_*}'$
\beaa
|\widetilde{\dkb}^kF| &\les& \int_u^{u_*}|\dk^{\leq k}\Hc|+\frac{\sqrt{\ep}}{u^{\frac{1}{2}+\dec-2\de_0}}\int_u^{u_*}\frac{1}{\nu_{\Si_*}'(u')r^2}\textrm{ for }k\leq k_{small}+130, \,\, u_1\leq u\leq u_*.
\eeaa
Since 
\beaa
\nu_{\Si_*}'(u) &=& e_3'(u)+b'e_4'(u)\\
&=& e_3(u)   -\frac{1+\frac{2mr+a^2(\sin\th)^2}{|q|^2}+\widecheck{e_3(r)}+\widecheck{e_3(u)}}{1+f\c\nab(u) +\frac 1 4 |f|^2(e_3(u)+e_3(r))}\left(e_4 + f^b e_b +\frac 1 4 |f|^2  e_3\right)u\\
&=& \frac{2(r^2+a^2)}{|q|^2}+\widecheck{e_3(u)}   -\frac{1+\frac{2mr+a^2(\sin\th)^2}{|q|^2}+\widecheck{e_3(r)}+\widecheck{e_3(u)}}{1+f\c\nab(u) +\frac 1 4 |f|^2(e_3(u)+e_3(r))}\left(f\c\nab(u) +\frac 1 4 |f|^2e_3(u)\right)
\eeaa
we have
\beaa
\nu_{\Si_*}'(u) &=& 2+O\left(\frac{1}{r}+\ep\right)
\eeaa
and hence, we have on $\Si_*$
\beaa
|\widetilde{\dkb}^kF| &\les& \int_u^{u_*}|\dk^{\leq k}\Hc|+\frac{\sqrt{\ep}}{u^{\frac{1}{2}+\dec-2\de_0}}\int_u^{u_*}\frac{1}{r^2}\textrm{ for }k\leq k_{small}+130, \,\,u_1\leq u\leq u_*.
\eeaa
Together with the behavior \eqref{eq:behaviorofronS-star} of $r$ on $\Sigma_*$, we infer
\beaa
|\widetilde{\dkb}^kF| &\les& \int_u^{u_*}|\dk^{\leq k}\Hc|+\frac{\ep}{ru^{\frac{1}{2}+\dec-2\de_0}}\textrm{ for }k\leq k_{small}+130,\,\, u_1\leq u\leq u_*.
\eeaa

Next, we estimate $\Hc$. We have by interpolation
\beaa
\|\dk^{\leq k_{small}+132}\Hc\|_{L^2(S)} &\les& \|\dk^{\leq k_{small}}\Hc\|_{L^2(S)}^{1-\frac{132}{k_{large}-k_{small}}}\|\dk^{\leq k_{large}}\Hc\|_{L^2(S)}^{\frac{132}{k_{large}-k_{small}}},
\eeaa
and hence, using $\de_0>0$, we have
\beaa
&& \int_{\Si_*(\geq u)}|\dk^{\leq k_{small}+132}\Hc| \\
&\les& \left(\int_{\Si_*(\geq u)}{u'}^{1+\de_0}|\dk^{\leq  k_{small}+132}\Hc|^2\right)^{\frac{1}{2}}\\
&\les& \frac{1}{u^{\frac{1}{2}+\dec - 2\de_0}}\left(\int_{\Si_*}{u'}^{2+2\dec}|\dk^{\leq k_{small}}\Hc|^2\right)^{\frac{1}{2}-\frac{132}{2(k_{large}-k_{small})}}\left(\int_{\Si_*}|\dk^{\leq k_{large}}\Hc|^2\right)^{\frac{132}{2(k_{large}-k_{small})}},
\eeaa
where we have used the fact, in view of the definition of $\de_0$, that 
\beaa
\frac{132}{k_{large}-k_{small}}(1+\dec)+\frac{\de_0}{2} & = & \left(\left(1+\frac{2}{130}\right)(1+\dec)+\frac{1}{2}\right)\de_0\leq 2\de_0
\eeaa
and
\beaa
\frac{1}{2}+\dec - 2\de_0 \geq  \frac{1}{2}+\dec -\frac{2}{3}\dec=\frac{1}{2}+ \frac{\dec}{3}>0.
\eeaa
Now, recall from the bootstrap assumptions on decay and boundedness for $\Hc$ on $\Mext$ that we have 
\beaa
\int_{\Si_*}u^{2+2\dec}|\dk^{\leq k_{small}}\Hc|^2 + \int_{\Si_*}|\dk^{\leq k_{large}}\Hc|^2 &\leq& \ep^2. 
\eeaa
We deduce
\beaa
\int_{\Si_*(\geq u)}|\dk^{\leq k_{small}+132}\Hc| &\les& \frac{\ep}{u^{\frac{1}{2}+\dec - 2\de_0}}.
\eeaa
Together with the Sobolev embedding on the 2-spheres $S$ foliating $\Si_*$, as well as  the behavior \eqref{eq:behaviorofronS-star} of $r$ on $\Sigma_*$, we infer 
\beaa
\int_u^{u_*}|\dk^{\leq k_{small}+130}\Hc| &\les& \frac{\ep}{ru^{\frac{1}{2}+\dec-2\de_0}}.
\eeaa
Plugging in the above estimate  for $F$, we infer on $\Si_*$
\beaa
|\widetilde{\dkb}^kF| &\les& \frac{\ep}{ru^{\frac{1}{2}+\dec-2\de_0}}\textrm{ for }k\leq k_{small}+130,\,\, u_1\leq u\leq u_*.
\eeaa
Together with \eqref{eq:equationfore4primefforsecondframeofMext} and \eqref{eq:equationfore3primefforsecondframeofMext}, we recover $e_4$ and $e_3$ derivatives to deduce
\beaa
|\dk^kF| &\les& \frac{\ep}{ru^{\frac{1}{2}+\dec-2\de_0}}\textrm{ for }k\leq k_{small}+130,\,\, u_1\leq u\leq u_*
\eeaa
and hence, since $F=f+i\dual f$, 
\beaa
|\dk^kf| &\les& \frac{\ep}{ru^{\frac{1}{2}+\dec-2\de_0}}\textrm{ for }k\leq k_{small}+130,\,\, u_1\leq u\leq u_*.
\eeaa
This is an improvement of the bootstrap assumption \eqref{eq:localbootstraponSigma*forsecondframeMext}. Thus, we may choose $u_1=1$, and $f$ satisfies the following estimate
\beaa
|\dk^kf| &\les& \frac{\ep}{ru^{\frac{1}{2}+\dec-2\de_0}}\textrm{ for }k\leq k_{small}+130\textrm{ on }\Sigma_*.
\eeaa
Together with \eqref{eq:equationfore3primefforsecondframeofMext}, as well as  the behavior \eqref{eq:behaviorofronS-star} of $r$ on $\Sigma_*$, and the control of $\Hc$ provided by \eqref{eq:estimateforframeMexttobeusedlater:chap4secondframeMext}, we infer 
\beaa
|\dk^{k-1}\nab_3'f| &\les& |\dk^{k-1}\Hc|+\frac{\ep}{r^2}\\
&\les& \frac{\ep}{ru^{1+\dec-2\de_0}}\textrm{ for }k\leq k_{small}+130\textrm{ on }\Sigma_*.
\eeaa
Collecting the two above estimates, we obtain
\begin{equation}\lab{eq:controlonSigma*forsecondframeMext}
|\dk^kf| \les \frac{\ep}{ru^{\frac{1}{2}+\dec-2\de_0}}, \quad |\dk^{k-1}\nab_3'f| \les \frac{\ep}{ru^{1+\dec-2\de_0}}\textrm{ for }k\leq k_{small}+130\textrm{ on }\Sigma_*.
\end{equation}

{\bf Step 3}. Next, we estimate $f$ on $\Mext$. We assume the following local bootstrap assumption
\bea\lab{eq:localbootstraponMextforsecondframeMext}
|\dk^{\leq k_{small}+130}f| \leq \frac{\sqrt{\ep}}{ru^{\frac{1}{2}+\dec-2\de_0}}\,\,\,\,\textrm{ on }r\geq r_1,
\eea
where $r_1\geq r_0$. In view of the control of $f$ on $\Si_*$ provided by \eqref{eq:controlonSigma*forsecondframeMext},  \eqref{eq:localbootstraponMextforsecondframeMext} holds for $r_1$ sufficiently large, and our goal is to prove that we may in fact choose $r_1=r_0$ and replace $\sqrt{\ep}$ with $\ep$ in \eqref{eq:localbootstraponMextforsecondframeMext}. 

From Corollary \ref{cor:transportequationine4forchangeofframecoeffinformFFbandlamba}, we may rewrite \eqref{eq:equationfore4primefforsecondframeofMext} as
\beaa
\nab_4'(qF) &=& E_4(f, \Ga)\,\,\textrm{ on }\Mext,
\eeaa
where 
\beaa
E_4(f, \Ga) &=& -\frac{1}{2}q\trXc F -q\chih\c F+ q E_1(f, \Ga) +f\c\nab(q)F+\frac{1}{4}|f|^2e_3(q)F.
\eeaa
In view of  the estimate \eqref{eq:estimateforframeMexttobeusedlater:chap4secondframeMext} for the Ricci coefficients and curvature components with respect to the outgoing principal  frame $(e_4, e_3, e_1, e_2)$ of $\Mext$, and in view of the form of $E_1(f, \Ga)$, we have 
\beaa
\left|\dk^kE_4(f, \Ga)\right| &\les& \ep r^{-1}u^{-\frac{1}{2}}|\dk^{\leq k}f|+|\dk^{\leq k}f|^2+|\dk^{\leq k}f|^4\textrm{ for }k\leq k_{small}+130.
\eeaa
Using commutator identities, and in view of \eqref{eq:localbootstraponMextforsecondframeMext}, we infer\footnote{Note that we have 
\beaa
\dec-2\de_0 \geq \dec - \frac{2}{3}\dec \geq \frac{\dec}{3}>0.
\eeaa}
\beaa
\nab_4'(\dkb,  \Lieb_\T)^k(qF)  &\leq& \frac{\ep}{r^2u^{1+\dec-2\de_0}}\textrm{ for }k\leq k_{small}+130,\,\, r\geq r_1.  
\eeaa
Integrating backwards from $\Sigma_*$ where we have \eqref{eq:controlonSigma*forsecondframeMext}, we deduce
\beaa
|(\dkb, \Lieb_\T)^kf|  &\leq& \frac{\ep}{ru^{\frac{1}{2}+\dec-2\de_0}}\,\,\,\textrm{ for }k\leq k_{small}+130,\,\, r\geq r_1.  
\eeaa
Together with \eqref{eq:equationfore4primefforsecondframeofMext}, we recover the $e_4$ derivatives and obtain
\beaa
|\dk^kf|  &\leq& \frac{\ep}{ru^{\frac{1}{2}+\dec-2\de_0}}\,\,\,\textrm{ for }k\leq k_{small}+130,\,\, r\geq r_1.  
\eeaa
This is an improvement of the bootstrap assumption \eqref{eq:localbootstraponMextforsecondframeMext}. Thus, we may choose $r_1=r_0$, and we have
\beaa
|\dk^kf| &\les& \frac{\ep}{ru^{\frac{1}{2}+\dec-2\de_0}}\,\,\,\textrm{ for }k\leq k_{small}+130\textrm{ on }\Mext.
\eeaa
Also, commuting once \eqref{eq:equationfore4primefforsecondframeofMext} with $e_3'$, using the schematic commutator identity 
\beaa
[\nab_3',\nab_4'] &=& 2\omb'\nab_4'-2\om' \nab_3'+2(\eta'-\etab')\c\nab'+\Big(\xib'\xi', \,\eta'\etab', \,\dual\rho'\Big),
\eeaa 
and proceeding as above to integrate backward from $\Sigma_*$ where $\nab_3'f$ is under control from \eqref{eq:controlonSigma*forsecondframeMext}, we also obtain 
\beaa
|\dk^{k-1}\nab_3'f| &\les& \frac{\ep}{ru^{1+\dec-2\de_0}}+\frac{\ep}{r^2u^{\frac{1}{2}+\dec-2\de_0}}\,\,\,\textrm{ for }k\leq k_{small}+130\textrm{ on }\Mext.
\eeaa
Collecting the two above estimates, we obtain in the region $r\geq u^{\frac{1}{2}}$
\bea\lab{eq:controlonMextforsecondframeMext}
\begin{split}
|\dk^kf| &\les \frac{\ep}{ru^{\frac{1}{2}+\dec-2\de_0}}, \,\,\,\textrm{ for }k\leq k_{small}+130\textrm{ on }\Mext(r\geq u^{\frac{1}{2}}),\\
 |\dk^{k-1}\nab_3'f| &\les \frac{\ep}{ru^{1+\dec-2\de_0}}\,\,\,\textrm{ for }k\leq k_{small}+130\textrm{ on }\Mext(r\geq u^{\frac{1}{2}}).
\end{split}
\eea
Also, since $f=0$ on $S_*$ in view of \eqref{eq:propertiesoffforsecondframeofMext}, and in view of the transport equation \eqref{eq:equationfore4primefforsecondframeofMext} for $F$, and the fact that $F=f+i\dual f$, $f$ satisfies the following addition property 
\bea\lab{eq:extrapropertyoffwhichvanishesonuequlustar}
f=0\quad\textrm{on}\quad\{u=u_*\}.
\eea

{\bf Step 4.} We now introduce a third frame $(e_4'', e_3'', e_1'', e_2'')$ on $\Mext$ given by 
\bea\lab{eq:formulachangeofframeforsecondframeMext:bis}
\begin{split}
e_4''&= e_4 + {f'}^b e_b +\frac 1 4 |f'|^2  e_3,\\
e_a''&= e_a +\frac 1 2 {f'}_a e_3,\quad a=1,2,\\
e_3''&=  e_3,
\end{split}
\eea
where the horizontal 1-form $f'$ is defined by
\bea\lab{eq:propertiesoffforsecondframeofMext:bis}
f' := \psi\left(\frac{u^{\frac{1}{2}}}{r}\right)f,
\eea
with $\psi$ a smooth function on $\mathbb{R}$ such that $\psi=1$ for $s\leq \frac{1}{2}$ and $\psi=0$ for $s\geq 1$. Note in particular that $(e_4'', e_3'', e_1'', e_2'')$ coincides with $(e_4, e_3, e_1, e_2)$ in $\Mext(r\leq u^\frac{1}{2})$.

We estimate $f'$. Note that 
\beaa
re_4\left(\log\left(\frac{u^{\frac{1}{2}}}{r}\right)\right) &=& -1,\\
r\nab\left(\log\left(\frac{u^{\frac{1}{2}}}{r}\right)\right) &=& \frac{r}{2u}\nab(u)=\frac{ar}{2u}\Re(\Jk)+\frac{r}{2u}\Ga_b,\\
e_3\left(\log\left(\frac{u^{\frac{1}{2}}}{r}\right)\right) &=& \frac{1}{2}\frac{e_3(u)}{u}+\frac{e_3(r)}{r}=\frac{r^2+a^2}{u|q|^2}-\frac{\De}{r|q|^2} + \frac{r}{2u}\Ga_b+\Ga_b.
\eeaa
Thus, in view of the bootstrap assumptions on boundedness, and since $\frac{1}{2}\leq \frac{u^{\frac{1}{2}}}{r}\leq 1$ on the support of $\psi'$, we infer, for $0\leq k\leq k_{large}$, 
\beaa
\left|\dk^k\left(\psi\left(\frac{u^{\frac{1}{2}}}{r}\right)\right)\right| \les 1, \qquad \left|\dk^{k-1}e_3\left(\psi\left(\frac{u^{\frac{1}{2}}}{r}\right)\right)\right| \les \frac{1}{r}.
\eeaa
Together with the estimates \eqref{eq:controlonMextforsecondframeMext} for $f$ on $\Mext(r\geq u^{\frac{1}{2}})$, 
the property \eqref{eq:extrapropertyoffwhichvanishesonuequlustar} for $f$, and the definition \eqref{eq:propertiesoffforsecondframeofMext:bis} for $f'$, we infer
\bea\lab{eq:controlonMextforsecondframeMext:bis}
\begin{split}
f' &= 0,\,\,\,\, \textrm{ on }\Mext\left(r\leq u^{\frac{1}{2}}\right)\quad\textrm{and on}\quad\{u=u_*\},\\
|\dk^kf'| &\les \frac{\ep}{ru^{\frac{1}{2}+\dec-2\de_0}}, \,\,\,\textrm{ for }k\leq k_{small}+130\textrm{ on }\Mext\left(r\geq u^{\frac{1}{2}}\right),\\
 |\dk^{k-1}\nab_3''f'| &\les \frac{\ep}{ru^{1+\dec-2\de_0}}\,\,\,\textrm{ for }k\leq k_{small}+130\textrm{ on }\Mext\left(r\geq u^{\frac{1}{2}}\right),
\end{split}
\eea
which are the desired estimates for $f'$.

{\bf Step 5.} In view of Proposition \ref{Proposition:transformationRicci} applied to our particular case, i.e. a triplet $(f', \fb', \la')$ with $\fb'=0$ and $\la'=1$, and the fact that $(e_4, e_3, e_1, e_2)$ is an outgoing PG frame, we have 
\beaa
\xi'' &=& \frac{1}{2}\nab_4''f'+\frac{1}{4}(\trch f' -\atrch\dual f') +\frac{1}{2}f'\c\chih+\frac{1}{4}|f'|^2\eta+\frac{1}{2}(f'\c \ze)\,f' -\frac{1}{4}|f'|^2\etab   +\lot,\\
\xib'' &=& \xib,\\
\etab'' &=& -\ze+\frac{1}{4}(\trchb f'-\atrchb\dual f')+\frac{1}{2}f'\c\chibh+\frac{|f'|^2}{4}\xib,\\
\eta'' &=& \eta +\frac{1}{2}\nab_3f' -\omb f' -\frac{1}{4}|f'|^2\xib,\\
\ze'' &=& \ze -\frac{1}{4}(\trchb f' +\atrchb\dual f') -\omb f' -\frac{1}{2}f'\c\chibh  - \frac{1}{2}(f'\c\xib) f',
\eeaa
\beaa
\trchb'' &=& \trchb+f'\c\xib,\\
\atrchb'' &=& \atrchb+f'\wedge\xib,\\
\chibh'' &=& \chibh+\frac{1}{2}f'\hot\xib,
\eeaa
\beaa
\trch'' &=& \trch+\div''f' +f'\c(\ze+\eta) -\omb |f'|^2 -\frac{|f'|^2}{4}\left(\trchb+f'\c\xib\right),\\
\atrch'' &=& \atrch+\curl''f' +f'\wedge(\ze-\eta)  -\frac{|f'|^2}{4}\left(\atrchb+f'\wedge\xib\right),\\
\chih'' &=& \chih +\frac{1}{2}\nab''\hot f'+\frac{1}{2}f'\hot(\ze+\eta) -\frac{1}{2}\omb f'\hot f'  -\frac{|f'|^2}{8}\left(\chibh+f'\hot\xib\right),
\eeaa
\beaa
\omb'' &=& \omb+\frac{1}{2}f'\c\xib,\\
{\om''} &=& \ze\c f' -\frac{1}{2}|f'|^2\omb  -\frac{1}{2}f'\c f'\c\chib -\frac{|f'|^2}{4}\xib\c f',
\eeaa
and
\beaa
\a''&=&\a + \big(  f'\hot \b  -\dual f' \hot \dual  \b )+ \big( f'\hot f'-\frac 1 2  \dual f' \hot   \dual f' \big) \rho
+  \frac 3  2 \big(  f' \hot  \dual  f'\big) \rhod +\lot,\\
\aa''&=&\aa,\\
\b''&=& \b +\frac 3 2\big(  f' \rho+\dual  f'  \rhod\big)+\lot,  \\
  \bb''&=&\bb + \frac 1 2  \aa\c f' +\lot, \\
  \rho'' &=& \rho -f'\c\bb+\lot,\\
  \rhod'' &=& \rhod - f'\wedge\bb+\lot,
 \eeaa 
 where   the  lower order terms denoted by $\lot$  are linear  with respect   to    the Ricci coefficients and curvature components of the outgoing PG structure of $\Mext$,  and quadratic or  higher order in $f'$, and do not contain derivatives of the latter. Also, we have
 \beaa
 e_4'(r) &=& 1+\frac{1}{4}|f'|^2e_3(r),\\
  e_4'(u) &=& f'\c\nab(u)+\frac{1}{4}|f'|^2e_3(u),\\
 e_4'(J^{(0)}) &=& f'\c\nab(J^{(0)})+\frac{1}{4}|f'|^2e_3(J^{(0)}),
 \eeaa
 \beaa
 \nab'(r) &=& \frac{1}{2}e_3(r)f',\\
 \nab'(u) &=& \nab(u)+\frac{1}{2}e_3(u)f',\\
 \nab'(J^{(0)}) &=& \nab(J^{(0)})+\frac{1}{2}e_3(J^{(0)})f',
 \eeaa
 \beaa
 \nab_4'\Jk &=& -\frac{1}{q}\Jk+f'\c\nab\Jk+\frac{1}{4}|f'|^2\nab_3\Jk,\\
 \nab_a'\Jk &=& \nab_a\Jk+\frac{1}{2}f_a\nab_3\Jk,\,\, a=1,2,
 \eeaa
 while the derivatives of $r$, $u$, $J^{(0)}$ and $\Jk$ in the $e_3$ direction are the same as for the PG structure of $\Mext$ since $e_3''=e_3$. Together with  the estimates \eqref{eq:controlonMextforsecondframeMext:bis} for $f'$ on $\Mext$, and 
 the estimates for the PG structure of $\Mext$ provided by\footnote{We need additional estimates for $\nab_3\Ga_g$, $\a$, $\b$, $\nab_3\a$ and $\nab_3\b$ compared to \eqref{eq:estimateforframeMexttobeusedlater:chap4secondframeMext}. They are easily obtained in the same way, i.e by interpolation between the bootstrap assumptions on decay and boundedness for the outgoing PG frame $(e_4, e_3, e_1, e_2)$ of $\Mext$.}   \eqref{eq:estimateforframeMexttobeusedlater:chap4secondframeMext}, we immediately infer
  \bea\lab{eq:recoveryofalldesiredestimatesexceptforeta'forsecondframeMext}
\nn\max_{0\leq k\leq k_{small}+129}\sup_{\Mext}\!\!\!\!&&\Bigg\{\Big(r^2u^{\frac{1}{2}+\dec-2\de_0}+ru^{1+\dec-2\de_0}\Big)|\dk^k(\Ga_g''\setminus\{\Hc''\})|+ru^{1+\dec-2\de_0}|\dk^k\Ga_b''|\\
\nn&&+r^2u^{1+\dec-2\de_0}\left|\dk^{k-1}\nab_3''\left(\Ga_g''\setminus\{\Hc''\}\right)\right|+r^{\frac{7}{2}+\frac{\dt}{2}}\Big(|\dk^kA''|+|\dk^kB''|\Big)\\
\nn&& +\left(r^{\frac{9}{2}+\frac{\dt}{2}}+r^4u^{\frac{1}{2}+\dec-2\de_0}\right)|\dk^{k-1}\nab_3''A''|\\
&&+r^4u^{\frac{1}{2}+\dec-2\de_0}|\dk^{k-1}\nab_3''B''|\Bigg\} \les \ep,
\eea
where, according to Definition \ref{definition.Ga_gGa_b:withprimseforsecondframeMext}, $\Ga_g'', \Ga_b''$ are defined as follows:
\begin{itemize}
\item the linearized quantities for the frame $(e_1'', e_2'', e_3'', e_4'')$ are defined  in the same way  as Definition \ref{def:renormalizationofallnonsmallquantitiesinPGstructurebyKerrvalue} for the outgoing PG frame of $\Mext$, with respect to the   coordinates $(r, \th)$ and the complex 1-form $\Jk$ of the PG structure\footnote{Thus, for example, $\widecheck{\tr X}'' =\tr X''-\frac{2}{q}$, $\Hc''=H''-\frac{aq}{|q|^2}\Jk$, $\widecheck{e_3''(r)} = e_3''(r)-  \frac{\De}{|q|^2}$, $\widecheck{\DD''J^{(0)}}=\DD''J^{(0)}-i\Jk$, $\widecheck{\nab_3''\Jk}=\nab_3''\Jk-\frac{\De q}{|q|^4}\Jk$,...},  

\item in addition, we introduce the following linearized quantities which are trivial for an outgoing PG structure\footnote{Except $\Hbc$ which satisfies instead $\Hb=-Z$.}
\beaa
\Hbc''=\Hb''+\frac{a\ov{q}}{|q|^2}\Jk, \qquad \widecheck{e_4''(r)}=e_4''(r)-1, \qquad \widecheck{\nab_4''\Jk}=\nab_4''\Jk+\frac{1}{q}\Jk,
\eeaa

\item the notation $\Ga_b''$ is the one of Definition \ref{definition.Ga_gGa_b}, except that $\Hc''$ does not belong to $\Ga_b''$,

\item the notation $\Ga_g'$ is given by 
\beaa
\Ga_g'' &=& \Ga_{g,1}''\cup\Ga_{g,2}'',
\eeaa
where $\Ga_{g,1}''$ is the one of Definition \ref{definition.Ga_gGa_b}, and where $\Ga_{g,2}''$ is given by\footnote{Note that all quantities in $\Ga_{g,2}''$ vanish identically in the case of an outgoing PG structure except  $\Hc''$ and $\Hbc''$.}  
\beaa
\Ga_{g,2}'' = \Big\{\om'',\,\,\,  \Xi'',\,\,\, \Hbc'',\,\,\, \Hc'', \,\,\, \widecheck{e_4''(r)}, \,\,\, e_4''(u), \,\,\, e_4''(J^{(0)}), \,\,\, r^{-1}\nab''(r), \,\,\, \widecheck{\nab_4''\Jk}\Big\}.
\eeaa
\end{itemize}

Furthermore, recall that $\xi=0$, and that $\xi'=0$ by the construction of $f$, see \eqref{eq:propertiesoffforsecondframeofMext}. Since we have, by the choice of $f'$, $\xi''=\xi$ on $\Mext(r\leq u^\frac{1}{2})$ and  $\xi''=\xi'$ on $\Mext(r\geq 2u^{\frac{1}{2}})$, we infer
\beaa
\xi''=0\quad\textrm{on}\quad\Mext\setminus\left\{\frac{r}{2}\leq u^{\frac{1}{2}}\leq r\right\}.
\eeaa 
Together with the fact that $\xi''\in\Ga_g''\setminus\{\Hc''\}$ and \eqref{eq:recoveryofalldesiredestimatesexceptforeta'forsecondframeMext}, this implies in particular, for $k\leq k_{small}+129$,  
\beaa
|\dk^k\xi''| &\les& \frac{\ep}{r^2u^{\frac{1}{2}+\dec-2\de_0}}{\bf 1}_{\frac{r}{2}\leq u^{\frac{1}{2}}\leq r}\les  \frac{\ep}{r^{3+2(\dec-2\de_0)}},\\
|\dk^{k-1}\nab_3''\xi''| &\les&  \frac{\ep}{r^2u^{1+\dec-2\de_0}}{\bf 1}_{\frac{r}{2}\leq u^{\frac{1}{2}}\leq r}\les  \frac{\ep}{r^3u^{\frac{1}{2}+\dec-2\de_0}},
\eeaa
and thus, in addition to the estimates for $\xi''$ and $\nab_3\xi''$ provided by \eqref{eq:recoveryofalldesiredestimatesexceptforeta'forsecondframeMext}, we have
\bea\lab{eq:xidoublepimeissupportedinMextforrlikeuonehalf}
\max_{0\leq k\leq k_{small}+129}\sup_{\Mext}\Big[r^{3+2(\dec-2\de_0)}|\dk^k\xi''|+r^3u^{\frac{1}{2}+\dec-2\de_0}|\dk^{k-1}\nab_3''\xi''|\Big] \les \ep.
\eea 

Finally, since we have, by the choice of $f'$, $\eta''=\eta$ on $\Mext(r\leq u^\frac{1}{2})$, and since  $\Hc\in\Ga_b$ satisfies \eqref{eq:estimateforframeMexttobeusedlater:chap4secondframeMext}, we infer, on $\Mext(r\leq u^\frac{1}{2})$, 
\bea\lab{eq:secondfromofMextcruialfirstestimateforetadoubleprime}
\max_{0\leq k\leq k_{small}+130}\sup_{\Mext(r\leq u^\frac{1}{2})}ru^{1+\dec-2\de_0}|\dk^k\Hc''| \les \ep.
\eea

{\bf Step 6.} Notice that \eqref{eq:secondfromofMextcruialfirstestimateforetadoubleprime} yields the desired estimate for $\Hc''$ in  $\Mext(r\leq u^\frac{1}{2})$. We now focus on estimating $\Hc''$ in the region $\Mext(r\geq u^{\frac{1}{2}})$. Proceeding as for the other Ricci coefficients would yield for $\Hc''$ the same behavior than $\Hc$ and hence a loss of $r^{-1}$ compared to the desired estimate. Instead, we rely on the following null structure equations in Proposition \ref{prop-nullstr:complex}  
\beaa
\nab_4''Z'' +\frac{1}{2}\tr X''(Z''-\Hb'')-2\om''(Z''+\Hb'') &=& 2\DD''\om'' +\frac{1}{2}\widehat{X}''\c(-\ov{Z}''+\ov{\Hb}'') -B''\\
&& -\frac{1}{2}\tr\Xb''\Xi''-2\omb''\Xi'' -\frac{1}{2}\ov{\Xi}''\c\Xbh'',\\
\nab_4''H'' +\frac{1}{2}\ov{\tr X}''(H''-\Hb'') &=&   -\frac{1}{2}\Xh''\c(\ov{H}''-\ov{\Hb}'') -B'' +\nab_3''\Xi''   -4\omb''\Xi''.
\eeaa
This yields
\beaa
&& \nab_4''(H''-Z'') + \frac{1}{2}\ov{\tr X}''(H''-Z'') +\frac{1}{2}\Xh''\c(\ov{H''-Z''}) \\
&=& - 2\DD''\om''   +\frac{1}{2}(\tr X'' - \ov{\tr X}'')(Z''-\Hb'')  -2\om''(Z''+\Hb'')\\
&&+\nab_3''\Xi''   -4\omb''\Xi''  +\frac{1}{2}\tr\Xb''\Xi''+2\omb''\Xi'' +\frac{1}{2}\ov{\Xi}''\c\Xbh''.
\eeaa
Since we have
\beaa
H''=\frac{aq}{|q|^2}\Jk+\Hc'', \qquad Z''=\frac{a\ov{q}}{|q|^2}\Jk+\Zc'', \qquad \Hb''=-\frac{a\ov{q}}{|q|^2}\Jk+\Hbc'', 
\eeaa
we deduce
\beaa
&& \nab_4''\left(\Hc''-\Zc''\right) + \frac{1}{2}\ov{\tr X}''\left(\Hc''-\Zc''\right)  \\
&=& - \nab_4''\left(\frac{a(q-\ov{q})}{|q|^2}\Jk\right) - \frac{1}{\ov{q}}\left(\frac{a(q-\ov{q})}{|q|^2}\Jk\right) +\frac{1}{2}\left(\frac{2}{q} - \frac{2}{\ov{q}}\right)\left(\frac{2a\ov{q}}{|q|^2}\Jk\right)\\
&& - 2\DD''\om''  -\frac{1}{2}\Xh''\c(\ov{H''-Z''}) +\frac{1}{2}(\tr X'' - \ov{\tr X}'')\left(\Zc''-\Hbc''\right)  -2\om''(Z''+\Hb'')\\
&& - \frac{1}{2}\ov{\trXc}''\left(\frac{a(q-\ov{q})}{|q|^2}\Jk\right)+\frac{1}{2}(\trXc'' - \ov{\trXc}'')\left(\frac{2a\ov{q}}{|q|^2}\Jk\right)\\
&&+\nab_3''\Xi''   -4\omb''\Xi''  +\frac{1}{2}\tr\Xb''\Xi''+2\omb''\Xi'' +\frac{1}{2}\ov{\Xi}''\c\Xbh''.
\eeaa
Also, using $e_4(q)=1$ and $\nab_4\Jk=-q^{-1}\Jk$, we have
\beaa
&&\nab_4''\left(\frac{a(q-\ov{q})}{|q|^2}\Jk\right)\\ 
&=& \left(\nab_4+f\c\nab+\frac{|f|^2}{4}\nab_3\right)\left(\frac{a(q-\ov{q})}{|q|^2}\Jk\right)\\
&=& \frac{a(q-\ov{q})}{|q|^2}\nab_4\Jk+e_4\left(\frac{a(q-\ov{q})}{|q|^2}\right)\Jk+\left(f\c\nab+\frac{|f|^2}{4}\nab_3\right)\left(\frac{a(q-\ov{q})}{|q|^2}\Jk\right)\\
&=& -\frac{1}{q}\frac{a(q-\ov{q})}{|q|^2}\Jk+a\left(-\frac{1}{\ov{q}^2}+\frac{1}{q^2}\right)\Jk+\left(f\c\nab+\frac{|f|^2}{4}\nab_3\right)\left(\frac{a(q-\ov{q})}{|q|^2}\Jk\right)
\eeaa
and hence
\beaa
&& \nab_4''\left(\Hc''-\Zc''\right) + \frac{1}{2}\ov{\tr X}''\left(\Hc''-\Zc''\right)  \\
&=&  - 2\DD''\om''  -\frac{1}{2}\Xh''\c(\ov{H''-Z''}) +\frac{1}{2}(\tr X'' - \ov{\tr X}'')\left(\Zc''-\Hbc''\right)  -2\om''(Z''+\Hb'')\\
&& - \frac{1}{2}\ov{\trXc}''\left(\frac{a(q-\ov{q})}{|q|^2}\Jk\right)+\frac{1}{2}(\trXc'' - \ov{\trXc}'')\left(\frac{2a\ov{q}}{|q|^2}\Jk\right) -\left(f\c\nab+\frac{|f|^2}{4}\nab_3\right)\left(\frac{a(q-\ov{q})}{|q|^2}\Jk\right)\\
&&+\nab_3''\Xi''   -4\omb''\Xi''  +\frac{1}{2}\tr\Xb''\Xi''+2\omb''\Xi'' +\frac{1}{2}\ov{\Xi}''\c\Xbh''.
\eeaa

Next, 
\begin{itemize}
\item we commute with $\dkb''$ and $\Lieb_\T''$, and we rely on the corresponding commutator identities,

\item we use the above equation for $\nab_4''(\Hc''-\Zc'')$ to recover the $e_4''$ derivatives, 

\item we rely on the estimates \eqref{eq:recoveryofalldesiredestimatesexceptforeta'forsecondframeMext}, as well as  the  estimate \eqref{eq:xidoublepimeissupportedinMextforrlikeuonehalf} for $\xi''$,
\end{itemize}
which allows us to derive, for $k\leq k_{small}+129$, 
\beaa
  \left|\nab_4''(\dk^k(\Hc''-\Zc''))+\frac 12 \ov{\tr X}''\dk^k(\Hc'' -\Zc'')\right| &\les& \frac{\ep}{r^3u^{\frac{1}{2}+\dec-2\de_0}}+\frac{\ep}{r^2}|\dk^{\leq k}(\Hc'' -\Zc'')|.
\eeaa
Since we have, by the choice of $f'$, $\eta''=\eta'$ on $\Mext(r\geq 2u^{\frac{1}{2}})$, and since $\Hc'=0$ on $\Si_*$, see \eqref{eq:propertiesoffforsecondframeofMext}, we have $\Hc''=0$ on $\Si_*$. Thus, integrating backwards from $\Sigma_*$, and using the control $\Zc''$ provided by \eqref{eq:recoveryofalldesiredestimatesexceptforeta'forsecondframeMext}, we infer
\beaa
&&\max_{0\leq k\leq k_{small}+129}\sup_{\Mext}r^2u^{\frac{1}{2}+\dec-2\de_0}|\dk^k\Hc''| \\
&\les& \ep+\max_{0\leq k\leq k_{small}+129}\sup_{\Mext}r^2u^{\frac{1}{2}+\dec-2\de_0}|\dk^k\Zc''|\\
&\les& \ep.
\eeaa
Also, commuting first the equation for $\nab_4''(\Hc''-\Zc'')$ with $\nab_3''$, using the schematic commutator identity 
\beaa
[\nab_3'',\nab_4''] &=& 2\omb''\nab_4''-2\om'' \nab_3''+2(\eta''-\etab'')\c\nab''+\Big(\xib''\xi'', \,\eta''\etab'', \,\dual\rho''\Big),
\eeaa 
and proceeding as above to integrate backward from $\Sigma_*$, we also obtain 
\beaa
&&\max_{0\leq k\leq k_{small}+129}\sup_{\Mext}r^2u^{1+\dec-2\de_0}|\dk^{k-1}\nab_3''\Hc''| \\
&\les& \ep+\max_{0\leq k\leq k_{small}+129}\sup_{\Mext}r^2u^{1+\dec-2\de_0}|\dk^{k-1}\nab_3''\Zc''|\\
&\les& \ep.
\eeaa
In view of the control of $\Hc''$ on $\Mext(r\leq u^\frac{1}{2})$ provided by \eqref{eq:secondfromofMextcruialfirstestimateforetadoubleprime}, this yields 
 \beaa
\nn\max_{0\leq k\leq k_{small}+129}\sup_{\Mext}\Big\{\Big(r^2u^{\frac{1}{2}+\dec-2\de_0}+ru^{1+\dec-2\de_0}\Big)|\dk^k\Hc''|\\
+r^2u^{1+\dec-2\de_0}\left|\dk^{k-1}\nab_3''\Hc''\right|\Big\} &\les& \ep.
\eeaa
Thus, together with \eqref{eq:recoveryofalldesiredestimatesexceptforeta'forsecondframeMext}, we infer
  \beaa
\nn\max_{0\leq k\leq k_{small}+129}\sup_{\Mext}&&\Bigg\{\Big(r^2u^{\frac{1}{2}+\dec-2\de_0}+ru^{1+\dec-2\de_0}\Big)|\dk^k\Ga_g''|+ru^{1+\dec-2\de_0}|\dk^k\Ga_b''|\\
\nn&&+r^2u^{1+\dec-2\de_0}\left|\dk^{k-1}\nab_3''\Ga_g''\right|+r^{\frac{7}{2}+\frac{\dt}{2}}\Big(|\dk^kA''|+|\dk^kB''|\Big)\\
\nn&& +\left(r^{\frac{9}{2}+\frac{\dt}{2}}+r^4u^{\frac{1}{2}+\dec-2\de_0}\right)|\dk^{k-1}\nab_3''A''|\\
&&+r^4u^{\frac{1}{2}+\dec-2\de_0}|\dk^{k-1}\nab_3''B''|\Bigg\} \les \ep.
\eeaa
Together with the control of $f'$ provided by \eqref{eq:controlonMextforsecondframeMext:bis} and the control of $\xi''$ provided by \eqref{eq:xidoublepimeissupportedinMextforrlikeuonehalf}, this concludes the proof of Proposition \ref{prop:constructionsecondframeinMext}.

%%%%%%%%%%%%%%%%%%%%%%%%%%%%%%%%%%%%%%%%%%%

\section{Proof of Proposition \ref{prop:existenceandestimatesfortheglobalframe:bis}}
\lab{sec:proofofprop:existenceandestimatesfortheglobalframe:bis}

%%%%%%%%%%%%%%%%%%%%%%%%%%%%%%%%%%%%%%%%%%%

Recall the small constant $\de_0>0$ introduced in the proof of Proposition \ref{prop:constructionsecondframeinMext}
\bea\lab{eq:introductionofconstantdelta0fortheproofofglobalframe:chap4}
\de_0=\frac{130}{k_{large}-k_{small}}\leq \frac{\dec}{3}.
\eea
In order to produce a global frame on $\MM$, we will proceed in several steps. First, we extend  the second frame of $\Mext$ constructed in Proposition \ref{prop:constructionsecondframeinMext} slightly inside $\Mint$. 

\begin{lemma}\lab{lemma:extensionofframesofMextintoMint}
Let  $({}^{(ext)}e_4', {}^{(ext)}e_3', {}^{(ext)}e_1', ^{(ext)}e_2')$ the second frame of $\Mext$ constructed in  Proposition \ref{prop:constructionsecondframeinMext}. We may extend this frame into $\Mint\cap\{\rint\geq r_0-2\}\cap\{u\leq u_*\}$. Furthermore:
\begin{itemize}
\item the linearized quantities $({}^{(ext)}\Ga_g', {}^{(ext)}\Ga_b')$ are defined 
using the extension of $(\rext, {}^{(ext)}\th)$ and ${}^{(ext)}\Jk$, associated to the outgoing PG structure or $\Mext$,  to the region $\Mint\cap\{\rint\geq r_0-2\}\cap\{u\leq u_*\}$ where they satisfy 
\beaa
\max_{0\leq k\leq k_{small}+129}\sup_{\Mint\cap\{\rint\geq r_0-2\}\cap\{u\leq u_*\}}u^{1+\dec-2\de_0}\Big|\dk^k(\rext-\rint, \\
\cos({}^{(ext)}\th)-\cos({}^{(int)}\th), {}^{(ext)}\Jk - {}^{(int)}\Jk)\Big| &\les& \ep,
\eeaa

\item we have
\beaa
\max_{0\leq k\leq k_{small}+127}\sup_{\Mint\cap\{\rint\geq r_0-2\}\cap\{u\leq u_*\}}u^{1+\dec-2\de_0}\left|\dk^k({}^{(ext)}\Ga_g', {}^{(ext)}\Ga_b')\right| &\les& \ep,
\eeaa
where  the Ricci coefficients and curvature components are the ones associated to  the frame $({}^{(ext)}e_4', {}^{(ext)}e_3', {}^{(ext)}e_1', ^{(ext)}e_2')$ of $\Mext$ extended into $\Mint$, 

\item we also have
\beaa
\max_{0\leq k\leq k_{small}+128}\sup_{\Mint\cap\{\rint\geq r_0-2\}\cap\{u\leq u_*\}}u^{1+\dec-2\de_0}\left|\dk^k\left(f, \fb, \log\left(\frac{{}^{(int)}\De}{|{}^{(int)}q|^2}\la\right)\right)\right| \les \ep,
\eeaa
where $(f, \fb, \la)$ denotes the change of frame coefficients from the ingoing PG frame  $({}^{(int)}e_4, {}^{(int)}e_3, {}^{(int)}e_1, ^{(int)}e_2)$ of $\Mint$ to $({}^{(ext)}e_4', {}^{(ext)}e_3', {}^{(ext)}e_1', ^{(ext)}e_2')$.
\end{itemize}
\end{lemma}

\begin{proof}
See section \ref{sec:proofoflemma:extensionofframesofMextintoMint}.
\end{proof}

To continue, we introduce the following definition.
\begin{definition}
For $j=1,2$,  we denote by $u_*^{(j)}$ the value of $\ub$ on the south pole of the sphere $\{u=u_*\}\cap\{\rext=r_0-j\}$ of the extension of the PG frame of $\Mext$ in $\Mint$. 
\end{definition}

\begin{remark}\lab{rmk:controlofvalueustarjgluingregion} 
Along $\ub=u^{(j)}$, $j=1,2$, in the region $r\sim r_0$, we have 
\beaa
\frac{du}{dr}=\frac{e_3(u)}{e_3(r)}=-2+O(r_0^{-1})
\eeaa
and hence, since $u_*^{(j)}$ is the value of $\ub$ on the south pole of $\{u=u_*\}\cap\{\rint=r_0-j\}$, we infer
\beaa
u_*^{(j)}=u_*-2j+O(r_0^{-1}).
\eeaa
\end{remark}

We now glue a renormalization of the  the second frame of $\Mext$ constructed in Proposition \ref{prop:constructionsecondframeinMext} to the ingoing PG frame of $\Mint$ in the matching region 
\bea
\mr_1 &:=& \Mint\cap\{r_0-1\leq\rint\leq r_0\}\cap\{\ub\leq u_*^{(1)}\}.
\eea  

\begin{lemma}\lab{lemma:matchingof3rdframeMextextendedtoMintwithMint}
Let the frame $({}^{(ext)}e_4', {}^{(ext)}e_3', {}^{(ext)}e_1', ^{(ext)}e_2')$ of Proposition \ref{prop:constructionsecondframeinMext}, defined on $\Mext$, and extended slightly into $\Mint$ in  Lemma \ref{lemma:extensionofframesofMextintoMint}. There exists a frame $(e_4'', e_3'', e_1'', e_2'')$ on $\Mext\cup\Mint(\ub\leq u_*^{(1)})$, as well as a pair of scalar functions $(r'', {J''}^{(0)})$, and a complex 1-form $\Jk''$,   such that
\begin{itemize}
\item[(a)] In $\Mext$, we have
\beaa
(e_4'', e_3'', e_1'', e_2'')=({}^{(ext)}\la{}^{(ext)}e_4',  {}^{(ext)}\la^{-1}{}^{(ext)}e_3', {}^{(ext)}e_1', ^{(ext)}e_2')
\eeaa
where  ${}^{(ext)}\la:=\frac{{}^{(ext)}\De}{|{}^{(ext)}q|^2}$, as well as $r''=\rext$, ${J''}^{(0)}=\cos({}^{(ext)}\th)$, and $\Jk''={}^{(ext)}\Jk$. 

\item[(b)] In $\Mint\cap\{\rint\leq r_0-1\}\cap\{\ub\leq u_*^{(1)}\}$, we have
\beaa
(e_4'', e_3'', e_1'', e_2'')= ({}^{(int)}e_4, {}^{(int)}e_3, {}^{(int)}e_1, ^{(int)}e_2),
\eeaa
as well as $r''=\rint$, ${J''}^{(0)}=\cos({}^{(int)}\th)$, and $\Jk''={}^{(int)}\Jk$. 
 
\item[(c)] In the matching region, we have
\beaa
\max_{0\leq k\leq k_{small}+127}\sup_{\mr_1}u^{1+\dec-2\de_0}\left|\dk^k(\Ga_g'', \Ga_b'')\right| &\les& \ep,
\eeaa
where  the Ricci coefficients and curvature components are the one associated to the frame $(e_4'', e_3'', e_1'', e_2'')$. 

\item[(d)] In the matching region, we also have
\beaa
\max_{0\leq k\leq k_{small}+128}\sup_{\mr_1}u^{1+\dec-2\de_0}\left|\dk^k(f, \fb, \log\la)\right| &\les& \ep,
\eeaa
where $(f, \fb, \la)$ denotes
\begin{itemize}
\item either the change of frame coefficients from $({}^{(ext)}\la{}^{(ext)}e_4',  {}^{(ext)}\la^{-1}{}^{(ext)}e_3', {}^{(ext)}e_1', ^{(ext)}e_2')$ to $(e_4'', e_3'', e_1'', e_2'')$,

\item or the one from $({}^{(int)}e_4, {}^{(int)}e_3, {}^{(int)}e_1, ^{(int)}e_2)$ to $(e_4'', e_3'', e_1'', e_2'')$.
\end{itemize}

\item[(e)] In the matching region, we have
\beaa
\max_{0\leq k\leq k_{small}+129}\sup_{\mr_1}u^{1+\dec-2\de_0}\Big|\dk^k(r'' -\rint,\, {J''}^{(0)} -\cos({}^{(int)}\th), \\
\Jk'' -{}^{(int)}\Jk)\Big| &\les& \ep.
\eeaa
\end{itemize}
\end{lemma}

\begin{proof}
See section \ref{sec:proooflemma:matchingof3rdframeMextextendedtoMintwithMint}.
\end{proof}

Next, we extend  the ingoing PG frame ${}^{(top)}\MM$ slightly inside $\Mext$ and $\Mint$. To this end, we first introduce a generalization of the notation $\,{}^{(ext)}\mathfrak{D}_k$ --- the decay norms in $\Mext$ introduced in section \ref{section:main-norms} ---  to deal with the extension in $\Mext$ for which we must take the weights in $r$ into account.

\begin{definition}
Given a subregion $\RR$ of $\Mext$ on which $r$ is arbitrary large, and given a null frame with corresponding $(\Ga_g, \Ga_b)$ on $\RR$, we then define the following decay norm
\bea
\bsplit
\mathfrak{D}_{k,\RR}'[\Ga_g, \Ga_b]  &:= \sup_{\RR}\Big(ru^{1+\dec'}+r^2u^{\frac{1}{2}+\dec'}\Big)|\dk^{\le k} \Ga_g|+\sup_{\RR}ru^{1+\dec'}|\dk^{\le k}\Ga_b|\\
&+\sup_{\RR}r^2 u^{1+\dec'} |\dk^{\le k-1} \nab_3\Ga_g|+\sup_{\RR}r^{\frac{7}{2}+\dec'}|\dk^{\le k} (A, B)|\\
&+\sup_{\RR}\Big(r^4u^{\frac{1}{2}+\dec'}+r^{\frac{9}{2}+\dec'}\Big)|\dk^{\le k-1}\nab_3 A|+\sup_{\RR}r^4u^{\frac{1}{2}+\dec'}|\dk^{\le k-1} \nab_3B|,
\end{split}
\eea
where $\dec'=\dec-2\de_0$.
\end{definition}

\begin{lemma}\lab{lemma:extensionofframesofMtopintoMextandMint}
We may extend the ingoing PG frame of ${}^{(top)}\MM$ into $\Mext(u\geq u_*^{(2)})$ and $\Mint(\ub\geq u_*^{(2)})$. Furthermore:
\begin{itemize}
\item we have
\beaa
\max_{0\leq k\leq k_{small}+126}\sup_{\Mint(\ub\geq u_*^{(2)})}u^{1+\dec-2\de_0}\left|\dk^k({}^{(top)}\Ga_g, {}^{(top)}\Ga_b)\right| &\les& \ep,
\eeaa
and 
\beaa
\mathfrak{D}_{k_{small}+127, \Mext(u\geq u_*^{(2)})}'[{}^{(top)}\Ga_g, {}^{(top)}\Ga_b] &\les& \ep,
\eeaa

\item we have
\beaa
\max_{0\leq k\leq k_{small}+127}\sup_{\Mint(\ub\geq u_*^{(2)})}u^{1+\dec-2\de_0}\left|\dk^k\left(f, \fb, \log\la\right)\right| &\les& \ep,
\eeaa
where $(f, \fb, \la)$ denotes the change of frame coefficients from the ingoing PG frame   of $\Mint$ to the ingoing PG frame of ${}^{(top)}\MM$,

\item we  have
\beaa
\max_{0\leq k\leq k_{small}+128}\sup_{\Mext(u\geq u_*^{(2)})}ru^{1+\dec-2\de_0}\left|\dk^k\left(f, \fb, \log\left(\frac{|{}^{(ext)}q|^2}{{}^{(ext)}\De}\la\right)\right)\right| &\les& \ep,
\eeaa
where $(f, \fb, \la)$ denotes the change of frame coefficients from the second frame of $\Mext$ constructed in Proposition \ref{prop:constructionsecondframeinMext}  to the ingoing PG frame of ${}^{(top)}\MM$,

\item we have
\beaa
&&\max_{0\leq k\leq k_{small}+127}\sup_{\Mint(\ub\geq u_*^{(2)})}u^{1+\dec-2\de_0}\\
&&\left|\dk^k\left({}^{(top)}r-\rint, \cos({}^{(top)}\th)-\cos({}^{(int)}\th), {}^{(top)}\Jk -{}^{(int)}\Jk\right)\right| \les \ep,
\eeaa
and
\beaa
\max_{0\leq k\leq k_{small}+129}\sup_{\Mext(u\geq u_*^{(2)})}u^{1+\dec-2\de_0}\Big(\left|\dk^k\left({}^{(top)}r-\rext\right)\right|\\
+r\left|\dk^k\left(\cos({}^{(top)}\th)-\cos({}^{(ext)}\th)\right)\right|+r^2\left|\dk^k\left({}^{(top)}\Jk -{}^{(ext)}\Jk\right)\right|\Big) &\les& \ep,
\eeaa

\item we have 
\beaa
\max_{0\leq k\leq k_{small}+126}\sup_{\Mext(u\geq u_*^{(2)})}r^3u^{\frac{1}{2}+\dec-2\de_0}\left|\dk^k{}^{(top)}\xi\right| &\les& \ep,
\eeaa
and 
\beaa
\max_{0\leq k\leq k_{small}+126}\sup_{\Mtop}r^3({}^{(top)}u)^{\frac{1}{2}+\dec-2\de_0}\left|\dk^k{}^{(top)}\xi\right| &\les& \ep.
\eeaa
\end{itemize}
\end{lemma}

\begin{proof}
See section \ref{sec:proofoflemma:extensionofframesofMtopintoMextandMint}.
\end{proof}

Finally, we glue the frame of Lemma \ref{lemma:matchingof3rdframeMextextendedtoMintwithMint} to the one of Lemma \ref{lemma:extensionofframesofMtopintoMextandMint} in the matching region 
\bea
\mr_2 &:=& (\Mint\cup\Mext)\cap\{u_*-2\leq \tau\leq u_*-1\},
\eea  
where $\tau$ is a smooth scalar function on $\Mint\cup\Mext$ such that 
\bea
\tau=u\,\,\,\textrm{in}\,\,\,\Mext\cap\{\rext\geq r_0+1\},\quad \tau=\ub+(u_*-1)-u_*^{(1)}\,\,\,\textrm{in}\,\,\, \Mint.
\eea

\begin{remark}\lab{rmk:bothframesusedforgluingglobalframedefinedinMatch3region}
Note that, in view of Remark \ref{rmk:controlofvalueustarjgluingregion}, and the definition of $\mr_2$ and $\tau$, we have 
\bea
\mr_2\cap\Mint\subset\{u_*^{(2)}\leq \ub\leq u_*^{(1)}\}
\eea
and we may choose $\tau$ in $\Mext\cap\{\rext\leq r_0+1\}$ such that we also have
\bea
\mr_2\cap\Mext\subset\{u\geq u_*^{(2)}\}.
\eea
In particular, the frames of Lemma \ref{lemma:matchingof3rdframeMextextendedtoMintwithMint} and Lemma \ref{lemma:extensionofframesofMtopintoMextandMint} are both defined in $\mr_2$.
\end{remark}

\begin{lemma}\lab{lemma:wefinallygetourglobalframegluingMtoptoframeMintandext}
Let the frame $(e_4'', e_3'', e_1'', e_2'')$ on $\Mext\cup\Mint(\ub\leq u_*^{(1)})$ of Lemma \ref{lemma:matchingof3rdframeMextextendedtoMintwithMint}, and let $({}^{(top)}e_4, {}^{(top)}e_3, {}^{(top)}e_1, ^{(top)}e_2)$ the ingoing PG frame of ${}^{(top)}\MM$ extended as in Lemma \ref{lemma:extensionofframesofMtopintoMextandMint}. There exists a global null frame $({}^{(glo)}e_4, {}^{(glo)}e_3, {}^{(glo)}e_1, ^{(glo)}e_2)$ defined on $\MM$, as well as a pair of scalar functions $({}^{(glo)}r, {}^{(glo)}J^{(0)})$, and a complex 1-form ${}^{(glo)}\Jk$,   such that:
\begin{itemize}
\item[(a)] In $(\Mint\cup\Mext)\cap\{\tau\leq u_*-2\}$, we have
\beaa
({}^{(glo)}e_4, {}^{(glo)}e_3, {}^{(glo)}e_1, ^{(glo)}e_2)=(e_4'', e_3'', e_1'', e_2''),
\eeaa
as well as ${}^{(glo)}r=r''$, ${}^{(glo)}J^{(0)}={J''}^{(0)}$, and ${}^{(glo)}\Jk=\Jk''$. 

\item[(b)] In ${}^{(top)}\MM\cup((\Mint\cup\Mext)\cap\{\tau\geq u_*-1\})$, we have
\beaa
({}^{(glo)}e_4, {}^{(glo)}e_3, {}^{(glo)}e_1, ^{(glo)}e_2) = ({}^{(top)}e_4, {}^{(top)}e_3, {}^{(top)}e_1, {}^{(top)}e_2),
\eeaa
 as well as ${}^{(glo)}r={}^{(top)}r$, ${}^{(glo)}J^{(0)}={}^{(top)}J^{(0)}$, and ${}^{(glo)}\Jk={}^{(top)}\Jk$. 
 
\item[(c)] In the matching region, we have
\beaa
\max_{0\leq k\leq k_{small}+126}\sup_{\mr_2\cap\Mint}u^{1+\dec-2\de_0}\left|\dk^k({}^{(glo)}\Ga_g, {}^{(glo)}\Ga_b)\right| &\les& \ep,
\eeaa
and 
\beaa
\mathfrak{D}_{k_{small}+126, \mr_2\cap\Mext}'[{}^{(glo)}\Ga_g, {}^{(glo)}\Ga_b] &\les& \ep.
\eeaa

\item[(d)] In the matching region, we  have
\beaa
\max_{0\leq k\leq k_{small}+127}\sup_{\mr_2\cap\Mext}ru^{1+\dec-2\de_0}\left|\dk^k(f, \fb, \log\la)\right|\\
+\max_{0\leq k\leq k_{small}+127}\sup_{\mr_2\cap\Mint}u^{1+\dec-2\de_0}\left|\dk^k(f, \fb, \log\la)\right| &\les& \ep,
\eeaa
where $(f, \fb, \la)$ denotes 
\begin{itemize}
\item either the change of frame coefficients from $(e_4'', e_3'', e_1'', e_2'')$ to $({}^{(glo)}e_4, {}^{(glo)}e_3, {}^{(glo)}e_1, ^{(glo)}e_2)$,

\item or the one from $({}^{(top)}e_4, {}^{(top)}e_3, {}^{(top)}e_1, ^{(top)}e_2)$ to $({}^{(glo)}e_4, {}^{(glo)}e_3, {}^{(glo)}e_1, ^{(glo)}e_2)$.
\end{itemize}

\item[(e)] In the matching region, we  have
\beaa
&&\max_{0\leq k\leq k_{small}+127}\sup_{\mr_2}u^{1+\dec-2\de_0}\Big(\left|\dk^k\left({}^{(glo)}r- r''\right)\right|\\
&&+r\left|\dk^k\left({}^{(glo)}J^{(0)}- {J''}^{(0)}\right)\right|+r^2\left|\dk^k\left({}^{(glo)}\Jk- \Jk''\right)\right|\Big) \les \ep
\eeaa
and
\beaa
&&\max_{0\leq k\leq k_{small}+127}\sup_{\mr_2}u^{1+\dec-2\de_0}\Big(\left|\dk^k\left({}^{(glo)}r-{}^{(top)}r\right)\right|\\
&&+r\left|\dk^k\left({}^{(glo)}J^{(0)}- \cos({}^{(top)}\th)\right)\right|+r^2\left|\dk^k\left({}^{(glo)}\Jk- {}^{(top)}\Jk\right)\right|\Big) \les \ep.
\eeaa

\item[(f)] In the part of the matching region inside $\Mext$, we  have
\beaa
\max_{0\leq k\leq k_{small}+125}\sup_{\mr_2\cap\Mext}r^{3+\dec-2\de_0}\left|\dk^k{}^{(glo)}\xi\right| &\les& \ep,\\
\max_{0\leq k\leq k_{small}+125}\sup_{\mr_2\cap\Mext}r^3u^{\frac{1}{2}+\dec-2\de_0}\left|\dk^{k-1}\nab_{{}^{(glo)}e_3}{}^{(glo)}\xi\right| &\les& \ep.
\eeaa
\end{itemize}
\end{lemma}

\begin{proof}
See section \ref{sec:proofoflemma:wefinallygetourglobalframegluingMtoptoframeMintandext}.
\end{proof}

We are now ready to prove Proposition \ref{prop:existenceandestimatesfortheglobalframe:bis}. The global frame of Proposition \ref{prop:existenceandestimatesfortheglobalframe:bis} is the one constructed in Lemma \ref{lemma:wefinallygetourglobalframegluingMtoptoframeMintandext}. First, note from the definition of $\mr_2$ and Remark \ref{rmk:bothframesusedforgluingglobalframedefinedinMatch3region} that $\mr_2\subset\mr$ where $\mr$ is the matching region of Definition \ref{def:cutofffunctionforthematchingregion:bis} appearing in the statement of  Proposition \ref{prop:existenceandestimatesfortheglobalframe:bis}. Next, note that property (c) of  Proposition \ref{prop:existenceandestimatesfortheglobalframe:bis} follows from property (b) of Lemma \ref{lemma:wefinallygetourglobalframegluingMtoptoframeMintandext}. Also, property (a) and (b) of  Proposition \ref{prop:existenceandestimatesfortheglobalframe:bis} follow from property (a) of Lemma \ref{lemma:wefinallygetourglobalframegluingMtoptoframeMintandext}, as well as properties (a) and (b) of Lemma \ref{lemma:matchingof3rdframeMextextendedtoMintwithMint}. Furthermore, properties (d), (e) and (f) of Proposition \ref{prop:existenceandestimatesfortheglobalframe:bis} follow from properties (c), (d) and (e) of Lemma \ref{lemma:wefinallygetourglobalframegluingMtoptoframeMintandext}, as well as property (d) of Lemma \ref{lemma:matchingof3rdframeMextextendedtoMintwithMint}. Finally, property (g) of Proposition \ref{prop:existenceandestimatesfortheglobalframe:bis} follows from property (f) of Lemma \ref{lemma:wefinallygetourglobalframegluingMtoptoframeMintandext} for $\Mext\cap\mr_2$, from Proposition \ref{prop:constructionsecondframeinMext} for $\Mext\setminus\mr_2$, and from Lemma \ref{lemma:extensionofframesofMtopintoMextandMint} for $\Mtop$.  This concludes the proof of Proposition \ref{prop:existenceandestimatesfortheglobalframe:bis}.

%%%%%%%%%%%%%%%%%%%%%%%%%%%%%%%%%%%%%%%%%%%

\subsection{Proof of Lemma \ref{lemma:extensionofframesofMextintoMint}}
\lab{sec:proofoflemma:extensionofframesofMextintoMint}

%%%%%%%%%%%%%%%%%%%%%%%%%%%%%%%%%%%%%%%%%%%

First, since the frame of $\Mext$ constructed in Proposition \ref{prop:constructionsecondframeinMext} and the outgoing PG frame of $\Mext$ agree in $\Mext(r\leq u^\frac{1}{2})$, it suffices in fact that show that the outgoing PG frame of $\Mext$ can be extended slightly into $\Mint$ and that its extension satisfies the properties of Lemma \ref{lemma:extensionofframesofMextintoMint}.  

To simplify  the notations, in this section, we denote
\begin{itemize}
\item by $(e_4, e_3, e_1, e_2)$ the ingoing PG frame of $\Mint$, 

\item by $(\ub, r, \th)$ and by $\Jk$ respectively triplet  of scalar functions and  the complex 1-form associated to the ingoing PG structure of $\Mint$, 

\item by $(e_4', e_3', e_1', e_2')$ the outgoing PG frame of $\Mext$ slightly extended into $\Mint$,  

\item by $(u, r', \th')$ and by $\Jk'$ respectively triplet  of scalar functions and  the complex 1-form associated to the outgoing PG structure of $\Mext$ slightly extended into $\Mint$,   

\item by $(f, \fb, \la)$ the change of frame from $(e_4, e_3, e_1, e_2)$ to  $(e_4', e_3', e_1', e_2')$.
\end{itemize}

Recall from section \ref{sec:initalizationadmissiblePGstructure} that we have, in view of the initialization of the ingoing PG frame of $\Mint$ on $\TT=\{r=r_0\}$ from the outgoing PG frame of $\Mext$, 
\bea\lab{eq:intializationofTTofingoingPGframetotransportslightlyinMint}
\ub=u, \quad r'=r, \quad \th'=\th, \quad f=\fb=0, \quad \la=\frac{|q|^2}{\De}\quad\textrm{on}\quad\{r=r_0\}.
\eea
Let 
\beaa
F:=f+i\dual f, \qquad \underline{F}:=\fb+i\dual \fb.
\eeaa
Then, we have in view of Corollary \ref{cor:transportequationine4forchangeofframecoeffinformFFbandlamba}
\beaa
\nab_{\la^{-1}e_4'}F+\frac{1}{2}\tr X F+2\om F &=& -2\Xi -\chih\c F +E_1(f, \Ga),\\
\la^{-1}\nab_4'(\log\la) &=& 2\om+f\c(\ze-\etab)+E_2(f, \Ga),\\
\nab_{\la^{-1}e_4'}\underline{F}+\frac{1}{2}\tr X\underline{F} &=& -2(\Hb+Z)   +2\DD'(\log\la)  +2\omb F   -  \fb^b\nab'f_b+E_3(f, \fb, \Ga),
\eeaa
where  $E_1(f, \Ga)$ and $E_2(f, \Ga)$ contain expressions of the type $O(\Ga f^2)$ with no derivatives, and $E_3(f, \fb, \Ga)$ contain expressions of the type $O(\Ga(f, \fb)^2)$ with no derivatives.

Since $\Xi$ and $\chih$ belong to $\Ga_g$, we have
\beaa
\nab_{\la^{-1}e_4'}F+\frac{1}{2}\tr X F+2\om F &=& \Ga_g +\Ga_g\c F+E_1(f, \Ga).
\eeaa
Integrating this transport equation for $F$ from $\{r=r_0\}$, where $F=0$ in view of \eqref{eq:intializationofTTofingoingPGframetotransportslightlyinMint}, to the region $\Mint\cap\{r\geq r_0-2\}\cap\{u\leq u_*\}$, using the bootstrap assumptions on decay and boundedness  for the ingoing PG structure of $\Mint$, and the above mentioned structure of the error term $E_1(f, \Ga)$, we easily deduce, since $F=f+i\dual f$, 
\beaa
\max_{k\leq k_{small}+129}\sup_{\Mint\cap\{r\geq r_0-2\}\cap\{u\leq u_*\}}u^{1+\dec-2\de_0}|\dk^kf| &\les& \ep.
\eeaa

Next, we estimate $\la$. We have 
\beaa
\la^{-1}\nab_4'\left(\log\left(\frac{\De}{|q|^2}\la\right)\right) &=& \la^{-1}\nab_4'(\log\la) +\frac{|q|^2}{\De}\la^{-1}\nab_4'\left(\frac{\De}{|q|^2}\right)
\eeaa
which together with the above equation for $\la^{-1}\nab_4'(\log\la)$ implies
\beaa
\la^{-1}\nab_4'\left(\log\left(\frac{\De}{|q|^2}\la\right)\right) &=& 2\om+f\c(\ze-\etab)+E_2(f, \Ga) +\frac{|q|^2}{\De}\la^{-1}\nab_4'\left(\frac{\De}{|q|^2}\right).
\eeaa
Also, we have by the transformation formula for $e_4'$
\beaa
\la^{-1}\nab_4'\left(\frac{\De}{|q|^2}\right) &=& \left(e_4+f\c\nab+\frac{1}{4}|f|^2e_3\right)\left(\frac{\De}{|q|^2}\right)\\
&=& e_4\left(\frac{\De}{|q|^2}\right)+\left(f\c\nab+\frac{1}{4}|f|^2e_3\right)\left(\frac{\De}{|q|^2}\right)
\eeaa
and hence
\beaa
\la^{-1}\nab_4'\left(\log\left(\frac{\De}{|q|^2}\la\right)\right) &=& 2\om+\frac{|q|^2}{\De}e_4\left(\frac{\De}{|q|^2}\right)
+f\c(\ze-\etab)+E_2(f, \Ga)\\
&& +\frac{|q|^2}{\De}\left(f\c\nab+\frac{1}{4}|f|^2e_3\right)\left(\frac{\De}{|q|^2}\right).
\eeaa
Note that we have for the frame of $\Mint$
\beaa
2\om+\frac{|q|^2}{\De}e_4\left(\frac{\De}{|q|^2}\right) &=& 2\omc+\De^{-1}\Ga_g=\Ga_g
\eeaa
since $\De^{-1}\sim r_0^{-2}$ in the region $r\geq r_0-2$. We deduce
\beaa
\la^{-1}\nab_4'\left(\log\left(\frac{\De}{|q|^2}\la\right)\right) = \Ga_g
+f\c(\ze-\etab)+E_2(f, \Ga) +\frac{|q|^2}{\De}\left(f\c\nab+\frac{1}{4}|f|^2e_3\right)\left(\frac{\De}{|q|^2}\right).
\eeaa
Integrating this transport equation  from $\{r=r_0\}$, where $\frac{\De}{|q|^2}\la=1$ in view of \eqref{eq:intializationofTTofingoingPGframetotransportslightlyinMint}, to the region $\Mint\cap\{r\geq r_0-2\}\cap\{u\leq u_*\}$, using the above control for $f$, the bootstrap assumptions on decay and boundedness  for the ingoing PG structure of $\Mint$, and the above mentioned structure of the error term $E_2(f, \Ga)$, we easily deduce 
\beaa
\max_{k\leq k_{small}+129}\sup_{\Mint\cap\{r\geq r_0-2\}\cap\{u\leq u_*\}}u^{1+\dec-2\de_0}\left|\dk^k\log\left(\frac{\De}{|q|^2}\la\right)\right| &\les& \ep.
\eeaa

Next, we estimate $\fb$. We have
\beaa
\DD'(\log\la) &=& \DD'\left(\log\left(\frac{\De}{|q|^2}\la\right)\right)+ \DD'\left(\log\left(\frac{|q|^2}{\De}\right)\right).
\eeaa
Using the  transformation formula for $e_a'$, we infer
\beaa
\DD'(\log\la) &=& \DD'\left(\log\left(\frac{\De}{|q|^2}\la\right)\right)+ \DD\left(\log\left(\frac{|q|^2}{\De}\right)\right)\\
&&+\left(\underline{F} f\c\nab+\underline{F}e_4+\left(\frac 1 2 F +\frac{1}{8}|f|^2\underline{F}\right)e_3\right)\left(\log\left(\frac{|q|^2}{\De}\right)\right).
\eeaa
Together with the above transport equation for $\underline{F}$, we deduce
\beaa
\nab_{\la^{-1}e_4'}\underline{F}+\frac{1}{2}\tr X\underline{F} &=& -2\left(\Hb+Z - \DD\left(\log\left(\frac{|q|^2}{\De}\right)\right)\right)  +2\DD'\left(\log\left(\frac{\De}{|q|^2}\la\right)\right)   \\
&&+2\omb F   -  \fb^b\nab'f_b+E_3(f, \fb, \Ga)\\
&&+2\left(\underline{F} f\c\nab+\underline{F}e_4+\left(\frac 1 2 F +\frac{1}{8}|f|^2\underline{F}\right)e_3\right)\left(\log\left(\frac{|q|^2}{\De}\right)\right).
\eeaa
Note that we have for the frame of $\Mint$, recalling that $\DD(r)=0$ and hence $\DD(\De)=0$, 
\beaa
\Hb+Z - \DD\left(\log\left(\frac{|q|^2}{\De}\right)\right) &=& -\frac{a\ov{q}}{|q|^2}\Jk+\Hbc+\frac{aq}{|q|^2}\Jk+\Zc - \frac{\DD(|q|^2)}{|q|^2}=\Ga_g
\eeaa
and hence
\beaa
\nab_{\la^{-1}e_4'}\underline{F}+\frac{1}{2}\tr X\underline{F} &=& \Ga_g +2\DD'\left(\log\left(\frac{\De}{|q|^2}\la\right)\right)   +2\omb F   -  \fb^b\nab'f_b+E_3(f, \fb, \Ga)\\
&&+2\left(\underline{F} f\c\nab+\underline{F}e_4+\left(\frac 1 2 F +\frac{1}{8}|f|^2\underline{F}\right)e_3\right)\left(\log\left(\frac{|q|^2}{\De}\right)\right).
\eeaa
Integrating this transport equation  from $\{r=r_0\}$, where $\underline{F}=0$ in view of \eqref{eq:intializationofTTofingoingPGframetotransportslightlyinMint}, to the region $\Mint\cap\{r\geq r_0-2\}\cap\{u\leq u_*\}$, using the above control for $f$ and $\la$, the bootstrap assumptions on decay and boundedness  for the ingoing PG structure of $\Mint$, and the above mentioned structure of the error term $E_3(f, \fb, \Ga)$, we easily deduce, since $\underline{F}=\fb+i\dual \fb$,  
\beaa
\max_{k\leq k_{small}+128}\sup_{\Mint\cap\{r\geq r_0-2\}\cap\{u\leq u_*\}}u^{1+\dec-2\de_0}\left|\dk^k\fb\right| &\les& \ep.
\eeaa
We have thus obtained the desired control of $(f, \fb, \la)$ in $\Mint\cap\{r\geq r_0-2\}\cap\{u\leq u_*\}$. 

Next, since $e_4'(r)=1$, $e_4'(\cos\th)=0$ and $\nab_4'\Jk'=-\frac{1}{q'}\Jk'$, we have
\beaa
e_4'(r'-r) &=& 1-e_4'(r)=1-\la\left(e_4+f\c\nab+\frac{1}{4}|f|^2e_3\right)r\\
&=& - \left(\la-\frac{|q|^2}{\De}\right)e_4(r)-\left(\frac{|q|^2}{\De}e_4(r)-1\right)- \frac{\la}{4}|f|^2e_3(r),\\
e_4'(\cos(\th')-\cos(\th)) &=& -e_4'(\cos\th)=-\la\left(e_4+f\c\nab+\frac{1}{4}|f|^2e_3\right)\cos\th\\
\nab_4'(q'\Jk'-q\Jk) &=& -e_4'(q\Jk)=-\la\left(e_4+f\c\nab+\frac{1}{4}|f|^2e_3\right)(q\Jk).
\eeaa
Integrating these transport equations  from $\{r=r_0\}$, where $r'=r$, $\cos(\th')=\cos(\th)$, $q'=q$ and $\Jk'=\Jk$  in view of \eqref{eq:intializationofTTofingoingPGframetotransportslightlyinMint}, to the region $\Mint\cap\{r\geq r_0-2\}\cap\{u\leq u_*\}$, using the above control for $f$ and $\la$, and the bootstrap assumptions on decay and boundedness  for the ingoing PG structure of $\Mint$, we easily deduce 
\beaa
\max_{k\leq k_{small}+129}\sup_{\Mint\cap\{r\geq r_0-2\}\cap\{u\leq u_*\}}u^{1+\dec-2\de_0}\left|\dk^k\left(r'-r, \, \cos(\th')-\cos(\th), \, \Jk'-\Jk\right)\right| &\les& \ep.
\eeaa

Finally, we consider the control of $(\Ga_g', \Ga_b')$. The change of frame formulas of Proposition \ref{Proposition:transformationRicci}, the above control of the change of frame coefficients $(f, \fb, \la)$,  the bootstrap assumptions on decay and boundedness  for the ingoing PG structure of $\Mint$, and the above control of $r'-r$, $\cos(\th')-\cos(\th)$ and $\Jk'-\Jk$ implies 
\beaa
\max_{k\leq k_{small}+127}\sup_{\Mint\cap\{r\geq r_0-2\}\cap\{u\leq u_*\}}u^{1+\dec-2\de_0}\left|\dk^k\left(\Gac', \Rc'\right)\right| &\les& \ep.
\eeaa
This concludes the proof of Lemma \ref{lemma:extensionofframesofMextintoMint}.

%%%%%%%%%%%%%%%%%%%%%%%%%%%%%%%%%%%%%%%%%%%

\subsection{Proof of Lemma \ref{lemma:matchingof3rdframeMextextendedtoMintwithMint}}
\lab{sec:proooflemma:matchingof3rdframeMextextendedtoMintwithMint}

%%%%%%%%%%%%%%%%%%%%%%%%%%%%%%%%%%%%%%%%%%%

To simplify  the notations, in this section, we denote
\begin{itemize}
\item by $(e_4, e_3, e_1, e_2)$ the ingoing PG frame of $\Mint$, 

\item by $(\ub, r, \th)$ and by $\Jk$ respectively triplet  of scalar functions and  the complex 1-form associated to the ingoing PG structure of $\Mint$, 

\item by $(e_4', e_3', e_1', e_2')$ the frame of Lemma \ref{lemma:extensionofframesofMextintoMint} which exists on $(\Mint\cap\{r\geq r_0-2\}\cap\{u\leq u_*\})\cup\Mext$,

\item by $(u, r', \th')$ and by $\Jk'$ respectively triplet  of scalar functions and  the complex 1-form associated to the outgoing PG structure of $\Mext$ extended to  $\Mint\cap\{r\geq r_0-2\}\cap\{u\leq u_*\}$.
\end{itemize}

We define the following null frame $(e_4'', e_3'', e_1'', e_2'')$ on $(\Mint\cap\{r\geq r_0-2\}\cap\{u\leq u_*\})\cup\Mext$
\beaa
e_4''=\frac{\De'}{|q'|^2}e_4', \quad e_3''=\frac{|q'|^2}{\De'}e_3', \quad e_a'=e_a, \,\, a=1,2.
\eeaa
We define the linearized quantities $(\Ga_g'', \Ga_b'')$ using $(r', \th')$ and $\Jk'$, with the normalization of $\Mint$. 
Let also $(f, \fb, \la)$ denote the  coefficients corresponding to the change of frame from  $(e_4, e_3, e_1, e_2)$ to $(e_4'', e_3'', e_1'', e_2'')$. In view of Lemma \ref{lemma:extensionofframesofMextintoMint} for $(e_4', e_3', e_1', e_2')$ and the above definition of $(e_4'', e_3'', e_1'', e_2'')$ and $(f, \fb, \la)$, we immediately obtain 
\beaa
\max_{k\leq k_{small}+127}\sup_{\Mint\cap\{r\geq r_0-2\}\cap\{u\leq u_*\}}u^{1+\dec-2\de_0}\left|\dk^k\left(\Ga_g'', \Ga_b''\right)\right| &\les& \ep.
\eeaa
and
\beaa
\max_{k\leq k_{small}+128}\sup_{\Mint\cap\{r\geq r_0-2\}\cap\{u\leq u_*\}}u^{1+\dec-2\de_0}\left|(f, \fb, \log(\la))\right| &\les& \ep.
\eeaa

Note that $\mr_1\subset \Mint\cap\{r\geq r_0-2\}\cap\{u\leq u_*\}$ so that all the above mentioned frames, scalars and complex 1-forms exist on $\mr_1$. Let $\psi$ a smooth cut-off function of $r$ such that $\psi=0$ for $r\leq r_0-1$ and $\psi=1$ for $r\geq r_0$. Then, we define the  null frame $(e_4''', e_3''', e_1''', e_2''')$ of $\Mext\cup(\Mint(\ub\leq u_*^{(1)})$ and the quantities $(r''', {J'''}^{(0)}, \Jk''')$ as follows
\begin{itemize}
\item In $\Mext$, we have
\beaa
(e_4''', e_3''', e_1''', e_2''')=(e_4'',  e_3'', e_1'', e_2''), \quad r'''=r', \quad {J'''}^{(0)}=\cos(\th'), \quad \Jk'''=\Jk'.
\eeaa

\item In $\Mint\cap\{\rint\leq r_0-1\}\cap\{\ub\leq u_*^{(1)}\}$, we have
\beaa
(e_4''', e_3''', e_1''', e_2''')= (e_4, e_3, e_1, e_2),\quad r'''=r, \quad {J'''}^{(0)}=\cos(\th), \quad \Jk'''=\Jk.
\eeaa
 
\item In the matching region $\mr_1$, $(e_4''', e_3''', e_1''', e_2''')$ is defined from $(e_4, e_3, e_1, e_2)$ using the change of frame coefficients $(f', \fb', \la')$ with 
\beaa
f'=\psi(r)f, \qquad \fb'=\psi(r)\fb, \qquad \la'=1-\psi(r)+\psi(r)\la,
\eeaa
where we recall that $(f, \fb, \la)$ denote the  coefficients corresponding to the change of frame from  $(e_4, e_3, e_1, e_2)$ to $(e_4'', e_3'', e_1'', e_2'')$.

\item In the matching region $\mr_1$, $r'''$, ${J'''}^{(0)}$ and $\Jk'''$ are defined by
\beaa
&& r'''=\psi(r)r'+(1-\psi(r))r, \quad {J'''}^{(0)}=\psi(r)\cos(\th')+(1-\psi(r))\cos\th, \\
&& \Jk'''=\psi(r)\Jk'+(1-\psi(r))\Jk.
\eeaa
\end{itemize}

In view of the above definitions, properties (a) and (b) of  Lemma \ref{lemma:matchingof3rdframeMextextendedtoMintwithMint} are immediate. Also, using the definition of $(f', \fb', \la')$ and the above control of $(f, \fb, \la)$, we have 
\beaa
\max_{k\leq k_{small}+128}\sup_{\mr_1}u^{1+\dec-2\de_0}\left|(f', \fb', \log(\la'))\right| &\les& \ep.
\eeaa
Also, if $(f'', \fb'', \la'')$ denotes the coefficients of the change of frame from $(e_4', e_3', e_1', e_2')$ to $(e_4'', e_3'', e_1'', e_2'')$, we easily obtain from the above control of $(f', \fb', \la')$ and  $(f, \fb, \la)$
\beaa
\max_{k\leq k_{small}+128}\sup_{\mr_1}u^{1+\dec-2\de_0}\left|(f'', \fb'', \log(\la''))\right| &\les& \ep,
\eeaa
which concludes the proof of property (d) of  Lemma \ref{lemma:matchingof3rdframeMextextendedtoMintwithMint}.

Finally, in view of the definition of $r'''$, ${J'''}^{(0)}$ and $\Jk'''$ in $\mr_1$, and the control of $r'-r$, $\cos(\th')-\cos\th$ and $\Jk'-\Jk$ of Lemma \ref{lemma:extensionofframesofMextintoMint}, we have
\beaa
\max_{k\leq k_{small}+129}\sup_{\mr_1}u^{1+\dec-2\de_0}\left|\dk^k\left(r'''-r, \, {J'''}^{(0)} -\cos(\th), \, \Jk'''-\Jk\right)\right| &\les& \ep.
\eeaa
Together with the change of frame formulas of Proposition \ref{Proposition:transformationRicci}, the above control of the change of frame coefficients $(f', \fb', \la')$,  and the bootstrap assumptions on decay and boundedness  for the ingoing PG structure of $\Mint$, we infer
\beaa
\max_{k\leq k_{small}+127}\sup_{\mr_1}u^{1+\dec-2\de_0}\left|\dk^k\left(\Ga_g''', \Ga_b'''\right)\right| &\les& \ep
\eeaa
which is property (c). This concludes the proof of Lemma \ref{lemma:matchingof3rdframeMextextendedtoMintwithMint}.

%%%%%%%%%%%%%%%%%%%%%%%%%%%%%%%%%%%%%%%%%%%

\subsection{Proof of Lemma \ref{lemma:extensionofframesofMtopintoMextandMint}}
\lab{sec:proofoflemma:extensionofframesofMtopintoMextandMint}

%%%%%%%%%%%%%%%%%%%%%%%%%%%%%%%%%%%%%%%%%%%

To simplify  the notations, in this section, we denote
\begin{itemize}
\item by $(e_4, e_3, e_1, e_2)$ the ingoing PG frame of $\Mint$, 

\item by $(\ub, r, \th)$ and by $\Jk$ respectively triplet  of scalar functions and  the complex 1-form associated to the ingoing PG structure of $\Mint$, 

\item by $(e_4', e_3', e_1', e_2')$ the second frame of $\Mext$ constructed in Proposition \ref{prop:constructionsecondframeinMext},

\item by $(u, r', \th')$ and by $\Jk'$ respectively triplet  of scalar functions and  the complex 1-form associated to the outgoing PG structure of $\Mext$,

\item by $(e_4'', e_3'', e_1'', e_2'')$ the frame of $\Mtop$,

\item by $(\ub'', r'', \th'')$ and by $\Jk''$ respectively triplet  of scalar functions and  the complex 1-form associated to the outgoing PG structure of ${}^{(top)}\MM$,

\item by $(f,\fb, \la)$ the coefficients of the change of frame between the frame of $\Mint$ and $\Mext$ in $\Mint\cap\{r\geq r_0-2\}\cap\{u\leq u_*\}$, 

\item by $(f',\fb', \la')$ the coefficients of the change of frame between the frame of $\Mint$ and ${}^{(top)}\MM$,

\item by $(f'',\fb'', \la'')$ the coefficients of the change of frame between the frame of $\Mext$ and ${}^{(top)}\MM$.
\end{itemize}

Recall from Proposition \ref{prop:constructionsecondframeinMext} that $(e_4', e_3', e_1', e_2')$ coincides with the PG frame of $\Mext$ on $\{u=u_*\}$. In particular, in view of the initialization of the PG frame of $\Mtop$ from the PG frame of $\Mext$ on $\{u=u_*\}$, we deduce that $(e_4'', e_3'', e_1'', e_2'')$ is initialized on $\{u=u_*\}$ by 
\beaa
e_4''=\frac{\De'}{|q'|^2}e_4', \qquad e_3''=\frac{|q'|^2}{\De'}e_3', \qquad e_a''=e_a', \,\, a=1,2.
\eeaa
In order to extend $(e_4'', e_3'', e_1'', e_2'')$, we assume that  the above initialization also holds on $\{u=u_*\}$ in $\Mint(\ub\geq u_*^{(2)})$, which is possible since the frame $(e_4', e_3', e_1', e_2')$  has been extend in Lemma \ref{lemma:extensionofframesofMextintoMint} to $\Mint\cap\{\rint\geq r_0-2\}\cap\{u\leq u_*\}$ and we have
\beaa
\{u=u_*\}\cap\Mint(\ub\geq u_*^{(2)})=\Mint\cap\{\rint\geq r_0-2\}\cap\{u\leq u_*\}.
\eeaa
Now, recall the following estimate derived in  Lemma \ref{lemma:extensionofframesofMextintoMint}
\beaa
\max_{0\leq k\leq k_{small}+128}\sup_{\Mint\cap\{\rint\geq r_0-2\}\cap\{u\leq u_*\}}u^{1+\dec-2\de_0}\left|\dk^k\left(f, \fb, \log\left(\frac{\De}{|q|^2}\la\right)\right)\right| \les \ep.
\eeaa
Together with the control of $r'-r$ and $\cos(\th')-\cos(\th)$ derived in the same lemma, the above initialization of  $(e_4'', e_3'', e_1'', e_2'')$, and the definition of $(f',\fb', \la')$, we infer
\beaa
\max_{0\leq k\leq k_{small}+128}\sup_{\{u=u_*\}\cap\Mint(\ub\geq u_*^{(2)})}u^{1+\dec-2\de_0}\left|\dk^k\left(f', \fb', \log(\la')\right)\right| \les \ep.
\eeaa

Next, let 
\beaa
F':=f'+i\dual f', \qquad \underline{F}':=\fb'+i\dual \fb'.
\eeaa
Then,  in view of the analog of Corollary \ref{cor:transportequationine4forchangeofframecoeffinformFFbandlamba} in the ingoing direction, we have
\beaa
\nab_{\la' e_3''}\underline{F}'+\frac{1}{2}\tr\Xb\,\underline{F}' &=&  -\chibh\c\underline{F}'+\underline{E}_1(\fb', \Ga),\\
\la'\nab_3''(\log\la') &=& \fb'\c(-\ze-\eta)+\underline{E}_2(\fb', \Ga),\\
\nab_{\la' e_3''}F+\frac{1}{2}\tr\Xb F' &=&   2\DD''(\log\la')     -  {f'}^b\nab''\fb_b'+\underline{E}_3(f', \fb', \Ga),
\eeaa
where  $\underline{E}_1(\fb', \Ga)$ and $\underline{E}_2(\fb', \Ga)$ contain expressions of the type $O(\Ga {\fb'}^2)$ with no derivatives, and $\underline{E}_3(f', \fb', \Ga)$ contain expressions of the type $O(\Ga(f', \fb')^2)$ with no derivatives. We propagate there transport equations from $\{u=u_*\}\cap\Mint(\ub\geq u_*^{(2)})$ to $\Mint(\ub\geq u_*^{(2)})$.  Using the above control of $(f', \fb', \log(\la'))$ on $\{u=u_*\}\cap\Mint(\ub\geq u_*^{(2)})$,   and the control of the foliation of $\Mint$, we infer, propagating first $\underline{F}'$, then $\la'$, and finally $F'$, 
\beaa
\max_{0\leq k\leq k_{small}+128}\sup_{\Mint(\ub\geq u_*^{(2)})}u^{1+\dec-2\de_0}\left|\dk^k\left(\fb', \log(\la')\right)\right| \les \ep
\eeaa
and 
\beaa
\max_{0\leq k\leq k_{small}+127}\sup_{\Mint(\ub\geq u_*^{(2)})}u^{1+\dec-2\de_0}\left|\dk^kf'\right| \les \ep.
\eeaa

Next, recall that $r''$, $\th''$ and $\Jk''$ are initialized on $\{u=u_*\}$ by 
\beaa
r''=r', \qquad \th''=\th', \qquad \Jk''=\Jk'.
\eeaa
In order to extend $r''$, $\th''$ and $\Jk''$, we assume that  the above initialization also holds on $\{u=u_*\}$ in $\Mint(\ub\geq u_*^{(2)})$. Now, recall the following estimate derived in  Lemma \ref{lemma:extensionofframesofMextintoMint}
\beaa
\max_{0\leq k\leq k_{small}+129}\sup_{\Mint\cap\{\rint\geq r_0-2\}\cap\{u=u_*\}}u^{1+\dec-2\de_0}\left|\dk^k\left(r'-r, \, \cos(\th')-\cos\th, \Jk'-\Jk\right)\right| \les \ep.
\eeaa
We infer
\beaa
\max_{0\leq k\leq k_{small}+129}\sup_{\Mint\cap\{\rint\geq r_0-2\}\cap\{u=u_*\}}u^{1+\dec-2\de_0}\left|\dk^k\left(r''-r, \, \cos(\th'')-\cos\th, \Jk''-\Jk\right)\right| \les \ep.
\eeaa
Since 
\beaa
&& e_3''(r''-r) = -1-e_3''(r),\qquad e_3''(\cos(\th'')-\cos(\th)) = -e_3''(\cos\th),\\
&&\nab_3''(q''\Jk''-q\Jk) = -\nab_3''(q\Jk),
\eeaa
together with the control of $\Mint$ and the above control of $(f', \fb', \la')$, we obtain  
\beaa
\max_{0\leq k\leq k_{small}+127}\sup_{\Mint(\ub\geq u_*^{(2)})}u^{1+\dec-2\de_0}\left|\dk^k\left(e_3''(r''-r), \, e_3''(\cos(\th'')-\cos\th), \nab_3''(q''\Jk''-q\Jk)\right)\right| \les \ep.
\eeaa
Integrating these transport equation from $\Mint\cap\{\rint\geq r_0-2\}\cap\{u=u_*\}$, we infer
\beaa
\max_{0\leq k\leq k_{small}+127}\sup_{\Mint(\ub\geq u_*^{(2)})}u^{1+\dec-2\de_0}\left|\dk^k\left(r''-r, \, \cos(\th'')-\cos\th, \Jk''-\Jk\right)\right| \les \ep.
\eeaa
Together with the change of frame formulas of Proposition \ref{Proposition:transformationRicci}, the above control of the change of frame coefficients $(f', \fb', \la')$,  and the bootstrap assumptions on decay and boundedness  for the ingoing PG structure of $\Mint$, we infer
\beaa
\max_{k\leq k_{small}+126}\sup_{\Mint(\ub\geq u_*^{(2)})}u^{1+\dec-2\de_0}\left|\dk^k\left(\Ga_g'', \Ga_b''\right)\right| &\les& \ep
\eeaa
as desired.

Next, we consider the extension of the ingoing PG structure of ${}^{(top)}\MM$ to $\Mext(u\geq u_*^{(2)})$. 
Recall from Proposition \ref{prop:constructionsecondframeinMext} that $(e_4', e_3', e_1', e_2')$ coincides with the PG frame of $\Mext$ on $\{u=u_*\}$. In particular, in view of the initialization of the PG structure of $\Mtop$ from the PG structure of $\Mext$ on $\{u=u_*\}$, we deduce that $(e_4'', e_3'', e_1'', e_2'')$ is initialized on $\{u=u_*\}$ by
\beaa
r''=r', \quad \th''=\th', \quad \Jk''=\Jk',\quad  f''=\fb''=0, \quad \la''=\frac{\De'}{|q'|^2}\quad\textrm{on}\quad\{u=u_*\}.
\eeaa
Let 
\beaa
F'':=f''+i\dual f'', \qquad \underline{F}'':=\fb''+i\dual \fb''.
\eeaa
Then, in view of the analog of Corollary \ref{cor:transportequationine4forchangeofframecoeffinformFFbandlamba} in the ingoing direction, we have
\beaa
\nab_{\la'' e_3''}\underline{F}''+\frac{1}{2}\tr\Xb'\,\underline{F}''+2\om'\underline{F}'' &=& -2\Xib' -\chibh'\c\underline{F}''+ \underline{E}_1(\fb'', \Ga'),\\
\la''\nab_3''(\log\la'') &=& 2\omb'+ \fb''\c(-\ze'-\eta')+\underline{E}_2(\fb'', \Ga'),\\
\nab_{\la'' e_3''}F''+\frac{1}{2}\tr\Xb' F'' &=&  -2(H'-Z')+ 2\DD''(\log\la'')    +2\om'\underline{F}'' -  {f''}^b\nab''\fb_b''\\
&&+\underline{E}_3(f'', \fb'', \Ga'),
\eeaa
where  $\underline{E}_1(\fb'', \Ga')$ and $\underline{E}_2(\fb'', \Ga')$ contain expressions of the type $O(\Ga {\fb''}^2)$ with no derivatives, and $\underline{E}_3(f'', \fb'', \Ga')$ contain expressions of the type $O(\Ga(f'', \fb'')^2)$ with no derivatives. 

Since $\Xib'$ and $\chibh'$ belong to $\Ga_b'$, we have
\beaa
\nab_{\la'' e_3''}\underline{F}''+\frac{1}{2}\tr\Xb'\,\underline{F}''+2\om'\underline{F}'' &=& \Ga_b'+\Ga_b'\c\underline{F}''+ \underline{E}_1(\fb'', \Ga').
\eeaa
Integrating this transport equation for $\underline{F}''$ from $\{u=u_*\}$, where $\underline{F}''=0$ in view of the above, to the region $\Mext(u\geq u_*^{(2)})$, noticing that the integration is on a region of size one along the integrable curves of $e_3''$ starting on $\{u=u_*\}$, using the bootstrap assumptions on decay and boundedness  for the outgoing PG structure of $\Mext$, and the above mentioned structure of the error term $\underline{E}_1(\fb'', \Ga')$, we easily deduce, since $\underline{F}''=\fb''+i\dual \fb''$, 
\beaa
\max_{k\leq k_{small}+129}\sup_{\Mext(u\geq u_*^{(2)})}ru^{1+\dec-2\de_0}|\dk^k\fb''| &\les& \ep.
\eeaa

Next, we estimate $\la''$. Similarly to the proof of Lemma \ref{lemma:extensionofframesofMextintoMint}, we may rewrite the transport equation for $\la''$ as 
\beaa
\la''\nab_3''\left(\log\left(\frac{|q'|^2}{\De'}\la''\right)\right) &=& \Ga_b'+ \frac{1}{2}\fb''\c(-\ze'-2\eta')+\underline{E}_2(\fb'', \Ga')\\
&& +\frac{\De'}{|q'|^2}\left(\fb''\c\nab+\frac{1}{4}|\fb''|^2e_3\right)\left(\frac{|q'|^2}{\De'}\right).
\eeaa
Integrating this transport equation  from $\{u=u_*\}$, where $\frac{|q'|^2}{\De'}\la''=1$  in view of the above, to the region $\Mext(u\geq u_*^{(2)})$,  noticing that the integration is on a region of size one along the integrable curves of $e_3''$ starting on $\{u=u_*\}$, 
using the above control for $f$, the bootstrap assumptions on decay and boundedness  for the outgoing PG structure of $\Mext$, and the above mentioned structure of the error term $\underline{E}_2(\fb'', \Ga')$, we easily deduce 
\beaa
\max_{k\leq k_{small}+129}\sup_{\Mext(u\geq u_*^{(2)})}ru^{1+\dec-2\de_0}\left|\dk^k\log\left(\frac{|q'|^2}{\De'}\la''\right)\right| &\les& \ep.
\eeaa

Next, we estimate $f''$. Similarly to the proof of Lemma \ref{lemma:extensionofframesofMextintoMint}, we may rewrite the transport equation for $f''$ as
\beaa 
\nab_{\la'' e_3''}F''+\frac{1}{2}\tr\Xb' F'' &=&  \Ga_b'+2\DD'\left(\log\left(\frac{|q'|^2}{\De'}\la''\right)\right)   +2\om'\underline{F}'' -  {f''}^b\nab''\fb_b'' +\underline{E}_3(f'', \fb'', \Ga')\\
&&+2\left(F''\underline{f}''\c\nab'+\left(\frac 1 2 \underline{F}'' +\frac{1}{8}|\fb''|^2F''\right)e_4'+F''e_3'\right)\left(\log\left(\frac{\De'}{|q'|^2}\right)\right).
\eeaa
Integrating this transport equation from $\{u=u_*\}$, where $F''=0$ in view of the above, to the region $\Mext(u\geq u_*^{(2)})$, noticing that the integration is on a region of size one along the integrable curves of $e_3''$ starting on $\{u=u_*\}$,  using the above control for $\fb''$ and $\la''$, the bootstrap assumptions on decay and boundedness  for the outgoing PG structure of $\Mext$, and the above mentioned structure of the error term $\underline{E}_3(f'', \fb'', \Ga')$, we easily deduce, since $F''=f''+i\dual f''$,  
\beaa
\max_{k\leq k_{small}+128}\sup_{\Mext(u\geq u_*^{(2)})}ru^{1+\dec-2\de_0}\left|\dk^kf''\right| &\les& \ep.
\eeaa
We have thus obtained the desired control of $(f'', \fb'', \la'')$ in $\Mext(u\geq u_*^{(2)})$. 

Next, we estimate $r''-r'$, $\th''-\th'$ and $\Jk''-\Jk'$. Similarly to the proof of Lemma \ref{lemma:extensionofframesofMextintoMint}, we derive transport equations in the $e_3''$ direction 
\beaa
e_3''(r''-r') &=& - \left(\la''-\frac{|q'|^2}{\De'}\right)e_3'(r')-\left(\frac{|q'|^2}{\De'}e_3'(r')-1\right)- \frac{\la''}{4}|\fb''|^2e_4'(r'),\\
e_3''(\cos(\th'')-\cos(\th')) &=& -\la''\left(e_3'+\fb''\c\nab+\frac{1}{4}|\fb''|^2e_4'\right)\cos(\th'),\\
\nab_3''(q'\Jk'-q\Jk) &=& -\la''\left(e_3'+\fb''\c\nab+\frac{1}{4}|\fb''|^2e_4'\right)(q'\Jk').
\eeaa
Integrating these transport equations  from $\{u=u_*\}$, where $r''=r'$, $\cos(\th'')=\cos(\th')$, $q''=q'$ and $\Jk''=\Jk'$, to the region $\Mext(u\geq u_*^{(2)})$, noticing that the integration is on a region of size one along the integrable curves of $e_3''$ starting on $\{u=u_*\}$, using the above control for $\fb''$ and $\la''$, and the bootstrap assumptions on decay and boundedness  for the outgoing PG structure of $\Mext$, we easily deduce 
\beaa
\max_{k\leq k_{small}+129}\sup_{\Mext(u\geq u_*^{(2)})}u^{1+\dec-2\de_0}\Big(\left|\dk^k\left(r''-r'\right)\right| +r\left|\dk^k\left(\cos(\th'')-\cos(\th')\right)\right|\\
+r^2\left|\dk^k\left(\Jk''-\Jk'\right)\right|\Big) &\les& \ep.
\eeaa

Next, we consider the control of $(\Ga_g'', \Ga_b'')$. The change of frame formulas of Proposition \ref{Proposition:transformationRicci}, the above control of the change of frame coefficients $(f'', \fb'', \la'')$,  the bootstrap assumptions on decay and boundedness  for the outgoing PG structure of $\Mext$, and the above control of $r''-r'$, $\cos(\th'')-\cos(\th')$ and $\Jk''-\Jk'$ implies 
\beaa
\mathfrak{D}_{k_{small}+127, \Mext(u\geq u_*^{(2)})}'[\Ga_g'', \Ga_b''] &\les& \ep.
\eeaa

Finally, we derive additional estimates for $\xi''$ in $\Mext(u\geq u_*^{(2)})$ and $\Mtop$. In view of the initialization of the  PG frame of $\Mtop$ from the PG frame of $\Mext$ on $\{u=u_*\}$, and since ${}^{(ext)}e_4$ is tangent to $\{u=u_*\}$, we infer $\xi''=(\frac{\De'}{|q'|^2})^2{}^{(ext)}\xi$ on $\{u=u_*\}$ and hence  
\beaa
\xi''=0\quad\textrm{on}\quad\{u=u_*\}. 
\eeaa
Now, the null structure equation for $\nab_3''\Xi''$ in the ingoing PG structure of $\Mtop$ yields
\beaa
\nab_3''\Xi'' &=& r^{-1}\dk^{\leq 1}\Ga_g''.
\eeaa
In view of the above control of $\Ga_g''$ on $\Mext(u\geq u_*^{(2)})$, as well as the bootstrap assumptions on decay and boundedness for the ingoing PG frame of $\Mtop$, we infer
\beaa
\max_{k\leq k_{small}+126}\left(\sup_{\Mext(u\geq u_*^{(2)})}r^3u^{\frac{1}{2}+\dec-2\de_0}|\dk^k\nab_3''\Xi''|+\sup_{\Mtop}r^3({}^{(top)}u)^{\frac{1}{2}+\dec-2\de_0}|\dk^k\nab_3''\Xi''| \right) \les \ep.
\eeaa
Integrating from $\{u=u_*\}$, where $\Xi''=0$ in view of the above, both towards $\Mtop$ and $\Mext(u\geq u_*^{(2)})$, and noticing that the integration is, in both cases, in a region of size one along the integral curves of $e_3''$ starting on $\{u=u_*\}$, we deduce 
\beaa
\max_{k\leq k_{small}+126}\left(\sup_{\Mext(u\geq u_*^{(2)})}r^3u^{\frac{1}{2}+\dec-2\de_0}|\dk^k\Xi''|+\sup_{\Mtop}r^3({}^{(top)}u)^{\frac{1}{2}+\dec-2\de_0}|\dk^k\Xi''| \right) &\les& \ep
\eeaa
as desired. This concludes the proof of Lemma \ref{lemma:extensionofframesofMtopintoMextandMint}.

%%%%%%%%%%%%%%%%%%%%%%%%%%%%%%%%%%%%%%%%%%%%%%%%%%

\subsection{Proof of Lemma \ref{lemma:wefinallygetourglobalframegluingMtoptoframeMintandext}}
\lab{sec:proofoflemma:wefinallygetourglobalframegluingMtoptoframeMintandext}

%%%%%%%%%%%%%%%%%%%%%%%%%%%%%%%%%%%%%%%%%%%%%%%%%%

To simplify  the notations, in this section, we denote
\begin{itemize}
\item by $(e_4, e_3, e_1, e_2)$ the frame on $\Mext\cup\Mint(\ub\leq u_*^{(1)})$ of Lemma \ref{lemma:matchingof3rdframeMextextendedtoMintwithMint},

\item by $(r, \th)$ and by $\Jk$ respectively the pair  of scalar functions and  the complex 1-form associated to $(e_4, e_3, e_1, e_2)$, 

\item by $(e_4', e_3', e_1', e_2')$ the ingoing PG frame of ${}^{(top)}\MM$, 

\item by $(r', \th')$ and by $\Jk'$ respectively the pair  of scalar functions and  the complex 1-form associated to the ingoing PG structure of ${}^{(top)}\MM$,

\item by $(f,\fb, \la)$ the coefficients of the change of frame between the frames $(e_4, e_3, e_1, e_2)$  and $(e_4', e_3', e_1', e_2')$.
\end{itemize}

Recall from Remark \ref{rmk:bothframesusedforgluingglobalframedefinedinMatch3region} that $(e_4, e_3, e_1, e_2)$ and $(e_4', e_3', e_1', e_2')$ are both defined in the matching region $\mr_2$, where we recall that 
\beaa
\mr_2 &=& (\Mint\cup\Mext)\cap\{u_*-2\leq \tau\leq u_*-1\},
\eeaa
where $\tau$ is a smooth scalar function on $\Mint\cup\Mext$ such that 
\beaa
\tau=u\,\,\,\textrm{in}\,\,\,\Mext\cap\{r\geq r_0+1\},\quad \tau=\ub+(u_*-1)-u_*^{(1)}\,\,\,\textrm{in}\,\,\, \Mint.
\eeaa

Let $\psi$ a smooth cut-off function of $\tau$ such that $\psi=0$ for $\tau\leq u_*-2$ and $\psi=1$ for $u\geq u_*-1$. Then, we define the global null frame $(e_4'', e_3'', e_1'', e_2'')$ of $\MM$ and the quantities $(r'', {J''}^{(0)}, \Jk'')$ as follows
\begin{itemize}
\item In $(\Mext\cup\Mint)\cap\{\tau\leq u_*-2\}$, we have
\beaa
(e_4'', e_3'', e_1'', e_2'')=(e_4,  e_3, e_1, e_2), \quad r''=r, \quad {J''}^{(0)}=\cos\th, \quad \Jk''=\Jk.
\eeaa

\item In ${}^{(top)}\MM\cup((\Mext\cup\Mint)\cap\{\tau\geq u_*-1\})$, we have
\beaa
(e_4'', e_3'', e_1'', e_2'')= (e_4', e_3', e_1', e_2'),\quad r''=r', \quad {J''}^{(0)}=\cos(\th'), \quad \Jk''=\Jk'.
\eeaa
 
\item In the matching region $\mr_2$, $(e_4'', e_3'', e_1'', e_2'')$ is defined from $(e_4, e_3, e_1, e_2)$ using the change of frame coefficients $(f', \fb', \la')$ with 
\beaa
f'=\psi(\tau)f, \qquad \fb'=\psi(\tau)\fb, \qquad \la'=1-\psi(\tau)+\psi(\tau)\la,
\eeaa
where we recall that $(f, \fb, \la)$ denote the  coefficients corresponding to the change of frame from  $(e_4, e_3, e_1, e_2)$ to $(e_4', e_3', e_1', e_2')$.

\item In the matching region $\mr_2$, $r''$, ${J''}^{(0)}$ and $\Jk''$ are defined by
\beaa
&& r''=\psi(r)r'+(1-\psi(r))r, \quad {J''}^{(0)}=\psi(r)\cos(\th')+(1-\psi(r))\cos\th, \\
&& \Jk''=\psi(r)\Jk'+(1-\psi(r))\Jk.
\eeaa
\end{itemize}

In view of the above definitions, properties (a) and (b) of  Lemma \ref{lemma:wefinallygetourglobalframegluingMtoptoframeMintandext} are immediate. Also, note that  in view of Lemma \ref{lemma:matchingof3rdframeMextextendedtoMintwithMint} and Lemma \ref{lemma:extensionofframesofMtopintoMextandMint}, $(f, \fb, \la)$ satisfies
\beaa
\max_{0\leq k\leq k_{small}+127}\sup_{\mr_2\cap\Mext}ru^{1+\dec-2\de_0}\left|\dk^k(f, \fb, \log\la)\right|\\
+\max_{0\leq k\leq k_{small}+127}\sup_{\mr_2\cap\Mint}u^{1+\dec-2\de_0}\left|\dk^k(f, \fb, \log\la)\right| &\les& \ep.
\eeaa
In view of the definition of $(f', \fb', \la')$ and the control of $(f, \fb, \la)$, we infer\footnote{Since  $\tau=u$ in $\mr_2\cap\{r\geq r_0+1\}$, and since 
\beaa
\dk(\psi(u)) &=& \psi'(u)(re_4(u), r\nab(u), e_3(u))=O(1)
\eeaa
note that there is no problem when $\dk$ falls on the cut-off $\psi(\tau)$.} 
\beaa
\max_{0\leq k\leq k_{small}+127}\sup_{\mr_2\cap\Mext}ru^{1+\dec-2\de_0}\left|\dk^k(f', \fb', \log(\la'))\right|\\
+\max_{0\leq k\leq k_{small}+127}\sup_{\mr_2\cap\Mint}u^{1+\dec-2\de_0}\left|\dk^k(f', \fb', \log(\la'))\right| &\les& \ep.
\eeaa
Also, if $(f'', \fb'', \la'')$ denotes the coefficients of the change of frame from $(e_4', e_3', e_1', e_2')$ to $(e_4'', e_3'', e_1'', e_2'')$, we easily obtain from the above control of $(f', \fb', \la')$ and  $(f, \fb, \la)$
\beaa
\max_{0\leq k\leq k_{small}+127}\sup_{\mr_2\cap\Mext}ru^{1+\dec-2\de_0}\left|\dk^k(f'', \fb'', \log(\la''))\right|\\
+\max_{0\leq k\leq k_{small}+127}\sup_{\mr_2\cap\Mint}u^{1+\dec-2\de_0}\left|\dk^k(f'', \fb'', \log(\la''))\right| &\les& \ep,
\eeaa
which concludes the proof of property (d) of  Lemma \ref{lemma:wefinallygetourglobalframegluingMtoptoframeMintandext}.

Finally, in view of the definition of $r''$, ${J''}^{(0)}$ and $\Jk''$ in $\mr_2$, and the control of $r'-r$, $\cos(\th')-\cos\th$ and $\Jk'-\Jk$ of Lemma \ref{lemma:extensionofframesofMtopintoMextandMint}, we have
\beaa
\max_{k\leq k_{small}+127}\sup_{\mr_2}u^{1+\dec-2\de_0}\Big(\left|\dk^k\left(r''-r'\right)\right| +r\left|\dk^k\left({J''}^{(0)} -\cos(\th')\right)\right|\\
+r^2\left|\dk^k\left(\Jk''-\Jk'\right)\right|\Big) &\les& \ep,
\eeaa
which is property (e). Together with the change of frame formulas of Proposition \ref{Proposition:transformationRicci}, the above control of the change of frame coefficients $(f', \fb', \la')$,  and the control of $(e_4, e_3, e_1, e_2)$ provided by  Lemma \ref{lemma:matchingof3rdframeMextextendedtoMintwithMint}, we infer
\beaa
\max_{k\leq k_{small}+126}\sup_{\mr_2\cap\Mint}u^{1+\dec-2\de_0}\left|\dk^k\left(\Ga_g'', \Ga_b''\right)\right|  &\les& \ep
\eeaa
and
\beaa
\mathfrak{D}_{k_{small}+126, \mr_2\cap\Mext}[\Ga_g'', \Ga_b''] &\les& \ep
\eeaa
which is property (c). Also, using the control of $\xi$ for $(e_4, e_3, e_1, e_2)$ provided by Proposition \ref{prop:constructionsecondframeinMext}, together with the change of frame formulas of Proposition \ref{Proposition:transformationRicci} and the above control of the change of frame coefficients $(f', \fb', \la')$, we obtain 
\beaa
\max_{k\leq k_{small}+126}\sup_{\mr_2\cap\Mext}r^3\left|\dk^k\xi''\right|  &\les& \ep,\\
\max_{k\leq k_{small}+126}\sup_{\mr_2\cap\Mext}r^3u^{\frac{1}{2}+\dec-2\de_0}\left|\dk^{k-1}\nab_3''\xi''\right|  &\les& \ep,
\eeaa
which is property (f).  This concludes the proof of Lemma \ref{lemma:wefinallygetourglobalframegluingMtoptoframeMintandext}.

%%%%%%%%%%%%%%%%%%%%%%%%%%%%%%%%%%%%%%%%

\chapter{Decay estimates on the last slice (Theorem M3)}
\lab{Chapter:decaySigmastar}

%%%%%%%%%%%%%%%%%%%%%%%%%%%%%%%%%%%%%%%%

The goal of this chapter is to prove Theorem M3, i.e. to improve our bootstrap assumptions on decay  for the integrable frame of $\Si_*$. This will be achieved in section \ref{sec::improvementofdecaybootassonSigmastar}, see Proposition \ref{prop:decayonSigamstarofallquantities}. We will then use these improvements to derive decay estimates for  the PG frame of $\Mext$ on $\Si_*$ in section \ref{sec:decayestimatesPGframeonSigmastar}, see Proposition \ref{prop:improvedesitmatesfortemporalframeofMextonSigmastar}.

In order to count the number of derivatives under control in this chapter, we introduce for convenience the following notation 
\bea\lab{eq:valueofkstarinchapter5forproofThmM3}
k_* &:=& k_{small}+80. 
\eea

%%%%%%%%%%%%%%%%%%%%%

\section{Geometric setting on $\Si_*$} 
\lab{sect.LastSlice-Intro}

%%%%%%%%%%%%%%%%%%%%%

We  recall the  properties  of      the  defining boundary  $\Si_*$ of our   GCM admissible   spacetime introduced in section \ref{sec:admissibleGMCPGdatasetonSigmastar}.  To start with  $\Si_*$ is  equipped, in view of section \ref{sec:admissibleGMCPGdatasetonSigmastar}, with a frame  $ (e_1, e_2, e_3, e_4)$  and a function  $r$ such that  $\big(\Si_*, r,  (e_1, e_2, e_3, e_4) \big)$  is a framed hypersurface, see Definition \ref{def:framedhypersurface}. Also, $\Si_*$ comes equipped with  function $u$  which verifies, for a constant\footnote{Recall that $c_*$ is fixed such that  $c_*=1+r(S_1)$ where $S_1$ is the only sphere of $\Si_*$ intersecting the curve of the south poles, see  section \ref{sec:initalizationadmissiblePGstructure}.}  $c_*$, 
\bea\lab{eq:relationbetweenrandudefininguonSigmastar}
u= c_*- r.
\eea 
Recall that we have also imposed the transversality conditions on $\Si_*$, see  \eqref{eq:tranversalityconditionforthefoliationonSi*},
   \bea\lab{eq:tranversalityconditionforthefoliationonSi*:chap5}
\xi=0, \qquad \om=0, \qquad \etab=-\ze, \qquad  e_4(r)=1, \qquad e_4(u)=0,
\eea
which allows to make sens of all the Ricci coefficients in the frame of $\Si_*$, i.e.,
 \beaa
 \chih, \quad  \ka=\trch, \quad  \eta,  \quad \ze,\quad  \etab,\quad   \xi,  \quad  \chibh,  \quad \kab= \trchb,  \quad  \xib,
  \eeaa  
 as well as make sense of all first order derivatives of $r$ and $u$ on $\Si_*$. 
 
For convenience, we introduce the following notations
\bea
y:= e_3(r), \qquad z:= e_3(u).
\eea
We thus have, using the transversality condition for $e_4(r)$ and $e_4(u)$, and the fact that $(e_1, e_2)$ is adapte to the $r$-foliation on $\Si_*$, 
 \bea
\nab(r)= \nab(u)=0, \quad e_4(u) =0, \quad e_4(r) =1, \quad e_3(r)=y,\quad e_3(u)=z.
\eea
 Also, recall that 
 \bea
 \nu= e_3+b_* e_4
 \eea
  denotes the  the vectorfield   tangent to $\Si_*$,  orthogonal  to the foliation and normalized by 
the condition $\g(\nu, e_4)=-2$.

%%%%%%%%%%%%%%%%%%%%%%%%%%%%%%%%%%%%%

 \subsection{Effective uniformization of  almost round $2$-spheres}
 \lab{subsection:effective-unifrmization}

%%%%%%%%%%%%%%%%%%%%%%%%%%%%%%%%%%%%% 

In this section, we  recall  some of the basic  results on effective uniformization from \cite{KS-GCM2} that will be useful in this chapter, as well as in chapter 8.
 \begin{definition}\lab{def:almostroundsphereS:chap5}
  A  $2$ dimensional, closed, Riemannian  surface $(S, g^S)$   is    said to be  almost round if  its Gauss curvature $K^S$ verifies, for a sufficiently small $\ep>0$,
     \bea
    \lab{eq:Almostroound-S} 
    \left|K^S-\frac{1}{( r^S)^2} \right|\le \frac{\ep}{( r^S)^2},
    \eea
    where $r^S$ is the area  radius of $S$. 
   \end{definition}

 The following  theorem is  Corollary 3.8 in  \cite{KS-GCM2}.
    \begin{theorem}[Effective uniformization]
    \lab{Thm:effectiveU1-Intro}
    Given an   almost round  sphere  $(S, g^S)$ as above  there exists, up to isometries\footnote{i.e.  all  the solutions are of the form $(\Phi\circ O, u\circ O)$ for $O\in O(3)$.} of $\SSS^2$,  a  unique   diffeomorphism   $\Phi:\SSS^2\to S$  and a  unique  conformal factor $u$ such that  
   \bea\lab{eq:centeredconditiononutobequotedlaterinchap5}
   \bsplit
    \Phi^\#( g^S)& =  (r^S)^2e^{2u}  \ga_{\SSS^2},\\
     \int_{\SSS^2} e^{2u} x^i &=0, \qquad i=1,2,3.
    \end{split}
    \eea    
    Moreover,  the  size of the conformal factor $u$ 
    is small with respect to the parameter $\ep$, i.e.
    \beaa
     \|  u\|_{L^\infty  (\SSS^2) } \les \ep.
    \eeaa 
    \end{theorem}

Theorem \ref{Thm:effectiveU1-Intro} is used to define  a canonical $\ell=1$ basis  on $S$, see  Definition 3.10 in \cite{KS-GCM2}.
 \begin{definition}[Basis of canonical $\ell=1$ modes  on $S$]
 \lab{definition:ell=1mpdesonS-intro}
  Let   $(\Phi, u)$  the unique, up to isometries of $\SSS^2$,  uniformization pair  given by Theorem \ref{Thm:effectiveU1-Intro}. 
We define the basis of canonical $\ell=1$ modes on $S$ by
\bea
J^S &:=& J^{\SSS^2}\circ \Phi^{-1}, 
\eea
where $J^{\SSS^2}$ denotes the $\ell=1$ spherical harmonics.
\end{definition}

The  main properties of  this basis are given in  Lemma 3.12 in \cite{KS-GCM2} which we review below.
 \begin{lemma}
 \lab{lemma:basicpropertiesofJforcanonicalell=1basis}
Let $J^S$ denote the basis of canonical $\ell=1$ modes on $S$ of Definition \ref{definition:ell=1mpdesonS-intro}. Then
 \bea
 \bsplit
 \lap_S J^{(p,S)}  &=-\frac{2}{(r^S)^2} J^{(p,S)}  +\frac{2}{(r^S)^2}\big(1- e^{-2v}\big)J^{(p,S)}, \\
 \int_S J^{(p,S)} J^{(q,S)} da_g &= \frac{4\pi}{3}(r^S)^2 \de_{pq}+ \int_S J^{(p,S)} J^{(q,S)} \big(1- e^{-2v} \big)   da_{g^S},\\
 \int_S J^{(p,S)} da_g &=0,
 \end{split}
 \eea
 with $\lap^S $ the Laplace-Beltrami of the metric $g^S$ and  $v:= u\circ\Phi^{-1}$. Moreover  we have
\bea
  \lab{eq:PropertiesofJ^S}
 \bsplit
 \lap_S J^{(p,S)}  &=\left(-\frac{2}{(r^S)^2} +O\left(\frac{\ep}{(r^S)^2}\right) \right)J^{(p,S)}, \\
 \int_S J^{(p,S)}J^{(q,S)} da_g &= \frac{4\pi}{3}(r^S)^2\de_{pq}+ O(\ep(r^S)^2),
 \end{split}
 \eea
 where  $\ep>0$ is the smallness constant appearing in  \eqref{eq:Almostroound-S}. 
  \end{lemma}

The following proposition is Proposition 4.15 in \cite{KS-GCM2}.
\begin{proposition}\lab{prop:asexpectedthebasiswidetildeJpisclosetocanonicalbasis}
Let $(S,g^S)$ an almost round sphere, i.e. verifying \eqref{eq:Almostroound-S}.   Consider an approximate  uniformization  pair  $( \widetilde{\Phi}, \widetilde{u})$, where $\widetilde{\Phi}:\SSS^2\to S$ is a smooth  diffeomorphism  and  $\widetilde{u}$  a smooth scalar function 
on $\SSS^2$ such that the following are verified, for $\de$ such that $0<\de\leq\ep$ and for $s\geq 2$, 
 \bea\lab{eq:approximateuniformizationuptoOofdelta}
  \left\|\widetilde{\Phi}^\#(g^S) -(r^S)^2e^{2\widetilde{u}}\ga_0\right\|_{H^s(\SSS^2)} &\leq& (r^S)^2\de,\qquad  \|\widetilde{u}\|_{H^2(\SSS^2)} \leq \ep.
  \eea
 Assume in addition  that the scalar functions $\widetilde{J}^{(p)}:=J^{(p,\SSS^2)}\circ\widetilde{\Phi}^{-1}$, $p\in\{0,+,-\}$, satisfy
  \bea\lab{eq:approximatebalancedconditionforwidetildeuOofdelta}
   \left|\int_{S}\widetilde{J}^{(p)}\right| \leq (r^S)^2\de.
  \eea
  Then we  can choose\footnote{Recall that  the pair   $(\Phi,u)$ is unique up to  isometries of $\SSS^2$.} the uniformization pair  $(\Phi, u)$  in Theorem \ref{Thm:effectiveU1-Intro} such that
  \bea\lab{eq:asexpectedthebasiswidetildeJpisclosetocanonicalbasis}
 (r^S)^{-1}\|\widetilde{u}\circ(\widetilde{\Phi})^{-1}-u\circ\Phi^{-1}\|_{\hk_{s}(S)}+\max_{p=0,+,-}(r^S)^{-1}\|\widetilde{J}^{(p)}-J^{(p,S)}\|_{\hk_{s+1}(S)} \les \de.
 \eea
\end{proposition}
 
 We recall  the following  definition from \cite{KS-GCM2} on  the calibration of   uniformization maps between almost round  spheres $S_1, S_2$ and  diffeomorphisms $\Psi: S_1\to S_2$.
   \begin{definition}
   \lab{def:calibration}
    On $\SSS^2$ we  fix\footnote{In particular, one can choose $N=(0,0,1)$ and $v=(1,0,0)$.} a point $N$ and a unit vector $v$ in the tangent space $T_N\SSS^2$. Given $\Psi:S_1\to S_2$,  we say that the  effective uniformization maps  $\Phi_1:\SSS^2 \to S_1 $, $\Phi_2:\SSS^2  \to S_2  $  are calibrated relative to $\Psi$  if  the map  $\Psih:=(\Phi_2)^{-1}\circ \Psi\circ \Phi_1:\SSS^2\to \SSS^2$ is such that:
\begin{enumerate}
\item  The map $\Psih$ fixes the point $N$, i.e. $\Psih(N)=N$.

\item    The tangent map    $\Psih_\#$ fixes the direction of $v$, i.e.    $\Psih_\#(v)=a_{1,2}v$ where $a_{1,2}>0$.

\item  The tangent map $\Psih_\#$ preserves the orientation of $T_N\SSS^2$.
 \end{enumerate}
 \end{definition}
 
The following theorem is Corollary 4.11 in  \cite{KS-GCM2}.
\begin{theorem}
 Let $\Psi:S_1\to S_2$   be  a given diffeomorphism and assume the following, for $\de$ such that $0<\de\leq\ep$, 
\begin{enumerate}
\item  The surfaces $S_1, S_2$ are close to each other, i.e.     for some $s\geq 0$, 
  \bea
   \lab{eq:unifor-twometricsB:higherregularityassumption}
  \left\| g^{S_1}-  \Psi^\#( g^{S_2}) \right\|_{\hk_{4+s}(S^1)} \le (r^{S_1})^3\de.
  \eea

  \item The maps $\Phi_1, \Phi_2$ are calibrated  according to Definition \ref{def:calibration}.
  \end{enumerate}
  Then
    \bea
  \lab{eq:unifor-twometrics1-calibr}
  \| \Psih-I\|_{L^\infty(\SSS^2)}+ \| \Psih-I\|_{H^1(\SSS^2)}&\les &\de,
 \eea
and  the  conformal factors $u_1, u_2$ verify
  \bea
  \lab{eq:unifor-twometrics2-calibr}
 \big \|u_1- \Psih^\# u_2\big \|_{L^\infty(\SSS^2)} \les \de.
  \eea
  \end{theorem}

%%%%%%%%%%%%%%%%%%%%%%%%%%%%%%%%%%%%%%% 
 
\subsection{The GCM conditions on  $\Si_*$}
\lab{sec:GCMconditionsonSigmastar}

%%%%%%%%%%%%%%%%%%%%%%%%%%%%%%%%%%%%%%%
 
 The  framed hypersurface   $\big(\Si_*, r,  (e_1, e_2, e_3, e_4) \big)$  introduced in section \ref{sec:admissibleGMCPGdatasetonSigmastar} 
   terminates in a future   boundary $S_*$    on which  the given function $r$ is constant, i.e. $S_*$ is a leaf of the $r$-foliation of $\Si_*$. 
     On $S_*$ there   exist  coordinates $(\th, \vphi)$    such that 
 \begin{enumerate}
 \item The induced metric $g$ on $S_*$  takes the form
 \bea\lab{eq:formofthemetriconSstarusinguniformization}
 g= r^2e^{2\phi}\Big( (d\th)^2+ \sin^2 \th (d\vphi)^2\Big).
 \eea
 
 \item The  functions 
 \bea
 J^{(0)} :=\cos\th, \qquad J^{(-)} :=\sin\th\sin\vphi, \qquad  J^{(+)} :=\sin\th\cos\vphi,
 \eea
 verify  the balanced  conditions 
 \bea\lab{eq:balancedconditionforJponSstar}
 \int_{S_*}  J^{(p)} =0, \qquad p=0,+,-.
 \eea
 \end{enumerate}
 
  Once $(\th, \vphi)$ are  chosen on $S_*$ we extend them to  $\Si_*$ by setting
  \bea
  \lab{eq:canonical-corrdsSi}
  \nu(\th)=\nu(\vphi)=0.
  \eea
  We  extend   the $\Jp$  functions to $\Si_*$ by setting
  \bea\lab{eq:canonical-ell=1modesonSi}
  \nu(\Jp)=0, \quad p=0,+,-.
  \eea
  We also impose the following transversality conditions on $\Si_*$ for $(\th, \vphi)$ and $\Jp$
  \bea\lab{eq:canonical-ell=1modesonSi:transcond}
  e_4(\th)=0, \qquad e_4(\vphi)=0, \qquad e_4(\Jp)=0, \quad p=0,+,-.
  \eea

  Recall the following, see Definition \ref{Definition:ell=1modesofascalarfunction}
  \begin{definition}
  \lab{definition:ell=1modes}
  Given a   scalar function $f$    on any  sphere $S$ of the $r$-foliation,  its  $\ell=1$ modes  are given by
  \beaa
  (f)_{\ell=1}=\left\{ \frac{1}{|S|} \int_S f \Jp, \, p=+, 0, - \right\}.
  \eeaa 
  \end{definition}
  
\begin{remark}
This definition is such that the $\ell=1$ mode have the same scaling in $r$ as the corresponding quantity. 
\end{remark}
  
  Endowed with  these  canonical coordinates and $\Jp$ basis,  $\Si_*$ is an admissible  GCM hypersurface, i.e.,
  \begin{enumerate}
  \item On $S_*$ we have
  \bea
  \ka=\frac 2 r, \qquad \kab =-\frac{2\Up}{r},
  \eea
  \bea
  (\div \b)_{\ell=1}=0, 
  \eea 
  as well as 
  \bea\lab{eq:S_*-GCM}
\int_{S_*}J^{(+)}\curl\b=0, \qquad \int_{S_*}J^{(-)}\curl\b=0.
\eea

  \item  On any sphere of the $r$-foliation of $\Si_*$, we have
  \bea
\lab{eq:Si_*-GCM1}
\bsplit
 \ka &=\frac{2}{r},\\
 \kab &=-\frac{2\Up}{r}+\underline{C}_0+\sum_{p=0, +,-}\underline{C}_pJ^{(p)},\\ 
 \mu &=\frac{2m}{r^3}+M_0+\sum_{p=0, +,-}M_p J^{(p)},
 \end{split}
 \eea
 where  $\underline{C}_0$, $\underline{C}_p$, $M_0$, $M_p$ are scalar functions on $\Si_*$ constant  on the leaves of the foliation. Also
 \bea
 \lab{eq:Si_*-GCM2}
(\div\eta)_{\ell=1} = (\div\xib)_{\ell=1} =0, \qquad   b_*\big|_{SP}=-1-\frac{2m}{r}, 
\eea
where $SP$ denotes the south poles of the spheres on $\Si_*$, i.e. $\th=\pi$. 

\item The mass $m$ is constant on $\Si_*$ and chosen to be the Hawking mass of $S_*$, i.e.,
\bea
\frac{2m}{r}=1+\frac{1}{16\pi} \int_{S_*} \ka\kab. 
\eea

\item The angular  momentum  is constant on $\Si_*$ and chosen as 
\bea
a &:=& \frac{r^3}{8\pi m}\int_{S_*} J^{(0)}\curl\b.
\eea

 \item Let $r_*$, $ u_*  $ denote the values of $r$  and $u$ on  $S_*$.  The function $r$ is monotonically decreasing on $\Si_*$  and the following dominance condition is verified
\bea\lab{eq:dominantconditiononronSigmastarchap5}
r_* &=& \de_*\ep_0^{-1}{u_*}^{1+\dec}
\eea 
which implies in particular on $\Si_*$
\beaa
r\ge r_* = \de_*\ep_0^{-1}{u_*}^{1+\dec} \ge \de_*\ep_0^{-1}u^{1+\dec}. 
\eeaa
 \end{enumerate}

%%%%%%%%%%%%%%%%%%%%%%%%%%%%%%%%%%%%%%% 
 
\subsection{Main equations in the frame of $\Si_*$}

%%%%%%%%%%%%%%%%%%%%%%%%%%%%%%%%%%%%%%%

Recall that the following transversality conditions hold on $\Si_*$, see \eqref{eq:tranversalityconditionforthefoliationonSi*:chap5}, 
\bea
\xi=0, \qquad \om=0, \qquad \etab=-\ze.
\eea

\begin{proposition}\label{prop-nullstrandBianchion.integrable}
We have on $\Si_*$, for the Ricci coefficients and curvature components associated to the GCM frame of $\Si_*$, the following\footnote{We only provide the  equations relevant for the control of the GCM frame of $\Si_*$.}:
\begin{enumerate}
\item The Ricci coefficients  verify the equations
\beaa
\nab_4\trch&=&-\frac 1 2  \trch^2-|\chih|^2,\\
\nab_4\chih +\trch\,  \chih &=& -\a,\\
\nab_4\trchb +\frac 1 2 \trch\trchb 
&=&   -  2   \div \ze + 2|\ze|^2+2\rho -\chih\c\chibh,\\
\nab_4\chibh +\frac 1 2 \big( \trchb \chih+\trch \chibh\big)
&=& -\nab\hot \ze +\ze\hot\ze,\\
\nab_3\trch +\frac 1 2 \trchb\trch
&=&         2   \div \eta+ 2 \omb \trch + 2 |\eta|^2+ 2\rho -\chibh\c\chih,
\eeaa
\beaa
\nab_3\trchb+\frac 1 2 \trchb^2 + 2\omb \trchb &=&2\div\xib   +  2 \xib\c(\eta-3\ze)-|\chibh|^2,\\
\nab_3\chih+\frac 1 2 \big( \trch \chibh+\trchb \chih\big)  -2 \omb \chih
&=& \nab\hot \eta +\eta\hot\eta,\\
\nab_3\chibh +\trchb\,  \chibh + 2\omb \chibh&=& \nab\hot \xib+   \xib\hot(\eta-3\ze)-\aa,
\eeaa
\beaa
\nab_4 \ze +\trch\ze &=& -2\chih\c\ze-\b,\\
\nab_3 \ze+2\nab\omb&=& -\chibh\c(\ze+\eta)-\frac{1}{2}\trchb(\ze+\eta)+ 2 \omb(\ze-\eta)\\
&&+\hch\c\xib+\frac{1}{2}\trch\,\xib  -\bb,
\\
\div\chih +\ze\c\chih &=& \frac{1}{2}\nab\trch+\frac{1}{2}\trch\ze - \b,\\
\div\chibh -\ze\c\chibh &=& \frac{1}{2}\nab\trchb-\frac{1}{2}\trchb\ze  +\bb,
\\
K &=& -\rho -\frac{1}{4}\trch\trchb +\frac{1}{2}\chih\c\chibh,
\\
\curl\ze&=&-\frac 1 2 \chih\wedge\chibh   +\dual \rho,\\
\curl \eta &=& \frac{1}{2}\chibh\wedge\chih +      \dual \rho,\\
\curl \xib &=&    \xib\wedge(\eta-\ze).
\eeaa

\item The curvature components   verify the equations
    \beaa
    \nab_3\a- \nab\hot \b&=&-\frac 1 2 \trchb\a+4\omb \a+
 (\ze+4\eta)\hot \b - 3 (\rho\chih +\rhod\dual\chih),\\
\nab_4\beta - \div\a &=&-2\trch\beta +\a\c  \ze,\\
     \nab_3\b - (\nab\rho+\dual\nab\rhod) &=&-\trchb \b+2 \omb\,\b+2\bb\c \chih+3 (\rho\eta+\rhod\dual \eta)+    \a\c\xib,\\
 \nab_4 \rho-\div \b&=&-\frac 3 2\trch \rho-\ze\c\b-\frac 1 2 \chibh \c\a,\\
   \nab_4 \rhod+\curl\b&=&-\frac 3 2\trch \rhod+\ze\c\dual \b+\frac 1 2 \chibh \c\dual \a, \\
     \nab_3 \rho+\div\bb&=&-\frac 3 2\trchb \rho  -(2\eta-\ze) \c\bb+2\xib\c\b-\frac{1}{2}\chih\c\aa,
 \\
   \nab_3 \rhod+\curl\bb&=&-\frac 3 2\trchb \rhod - (2\eta-\ze) \c\dual \bb-2\xib\c\dual\b-\frac 1 2 \chih\c\dual \aa.
\eeaa
\end{enumerate}
\end{proposition}

\begin{proof}
In view of the transversality conditions \eqref{eq:tranversalityconditionforthefoliationonSi*:chap5} and the fact that $(e_1, e_2)$ are tangent to the 2-spheres of the $r$-foliation of $\Si_*$, we have on $\Si_*$
\beaa
 \atrch=\atrchb=0, \qquad \xi=0, \quad \omb=0, \quad \etab+\ze=0.
 \eeaa
The proof follows then by plugging these identities in Propositions \ref{prop-nullstr} and \ref{prop:bianchi}.
\end{proof}

\begin{definition}
The mass aspect function $\mu$ is defined on $\Si$ by
\beaa
\mu &:=& -\div\ze -\rho +\frac{1}{2}\chih\c\chibh.
\eeaa
\end{definition}

\begin{lemma}
\lab{Lemma:eqts-nabeta,xib}
The following relations hold true for $y=e_3(r) $ and $z=e_3(u)$:
\bea
\nab   y =-\xib+\big(\ze-\eta) y, \qquad  \nab z = (\ze-\eta ) z.
\eea
Also, we have
\bea
b_*=-y-z.
\eea
\end{lemma}

\begin{proof}
We  make use of $\,[e_a, e_3]=(\ze_a-\eta_a ) e_3 -\xib_a e_4$ for $a=1,2$, which we apply to $r$. Using $\nab(r)=0$ and $e_4(r)=1$, we infer
\beaa
\nab(e_3(r)) &=& (\ze-\eta ) e_3(r)-\xib. 
\eeaa
Similarly, using $\nab(u)=0$ and $e_4(u)=0$, we have
\beaa
\nab(e_3(u)) &=& (\ze-\eta ) e_3(u).
\eeaa
Since $y=e_3(r) $ and $z=e_3(u)$, this concludes the proof of the first identities.

Next, since $u+r=c_*$ on $\Si_*$ and  $\nu$ is tangent to $\Si_*$, we have
\beaa
0 &=& \nu(u+r)=e_3(u)+e_3(r)+b_*e_4(u)+b_*e_4(r)=y+z+b_*
\eeaa
and hence $b_*=-y-z$ as stated. This concludes the proof of the lemma.
\end{proof}

%%%%%%%%%%%%%%%%%%%%%%%%%%%%%%%%%%%%%%%

\subsection{Linearized quantities and main quantitative assumptions}

%%%%%%%%%%%%%%%%%%%%%%%%%%%%%%%%%%%%%%%

Recall  the  notation  $\ka=\trch, \, \kab=\trchb$ and  $y=e_3(r), \, z= e_3(u)$. Also, recall
   the following linearized quantities introduced in Definition \ref{def:renormalizedquantitiesGCMfoliationofSigma*}: 
  \beaa
 \ba{cccccc}
\kac &:=&  \ds\trch-\frac{2}{r}, \qquad &
\kabc &:=& \ds\trchb+\frac{2\Up}{r},\\
\\
\ombc &:=& \ds\omb-\frac{m}{r^2},\qquad &\rhoc &:=& \ds \rho +\frac{2m}{r^3},\\
\\
\yc &:=& y +\Up,\qquad &  \zc &:=& z -2, \\
\\
 \widecheck{b_*}&:=& \ds b_*+1+\frac{2m}{r},\qquad &\widecheck{\mu} &:=& \ds \mu -\frac{2m}{r^3},\\
\ea
\eeaa
where $\Up = 1-\frac{2m}{r}$.

We denote by $\Ga_g, \Ga_b$ the sets of  linearized quantities\footnote{Note that these   were denoted by  $\Ga^*_g, \Ga^*_b$
 in Definition \ref{Definition:linearizedquantitiesSi*} } below.
 \begin{itemize}
\item The set $\Ga_g$
\bea
\Ga_{g}:&=\Big\{ \kac, \quad  \chih,  \quad   \ze,  \quad  \kabc, \quad  r\a,\quad  r\b, \quad  r \rhoc, \quad  r \rhod, \quad  r \muc\Big\}.
\eea
\item The set $\Ga_b$ 
\bea
\lab{eq:Definition-Ga_b}
\Ga_b&=&\Big\{ \eta,\quad  \chibh, \quad \widecheck{\omb}, \quad \xib,\quad r\bb,\quad \aa,  \quad  r^{-1}\yc , \quad r^{-1}\zc, \quad r^{-1} \widecheck{b_*}\Big\}.
\eea
\end{itemize}

\begin{definition}
\label{definition:norms-ProofThm.M3}
  We make  use of  the following   norms on  $S=S(u)\subset\Si_*$,
  \bea
  \bsplit
  \| f\|_{\infty} (u)&:=\| f\|_{L^\infty(S(u))}, \qquad \qquad  \| f\|_{2} (u):=\| f\|_{L^2(S(u))}, \\
  \|f\|_{\infty,k}(u)&:= \sum_{i=0}^k \|\dk_*^i f\|_{\infty }(u),  \qquad 
\|f\|_{2,k}(u):=\sum_{i=0}^k \|\dk_*^i f\|_{2}(u).
\end{split}
  \eea
  \end{definition}

Throughout this  chapter we rely on the following assumptions.

{\bf Ref 1.}  According to our bootstrap assumptions BA-D on decay, and BA-B on $r$-weighted sup norms, we have on $\Si_*$
\begin{enumerate}
\item For $0\le k\le k_{small}$,
\bea
\lab{Ref1-smallk:chap5}
\bsplit
\| \Ga_g\|_{\infty, k}&\leq\ep r^{-2} u^{-\frac{1}{2}-\dec},\\
 \|\nab_\nu\Ga_g\|_{\infty, k-1}  &\leq \ep  r^{-2} u^{-1-\dec}, \\   
\|\Ga_b\|_{\infty, k}&\leq\ep  r^{-1} u^{-1-\dec}.
\end{split}
\eea

\item For $0\le k\le k_{large}$,
\bea
\lab{Ref1-largek:chap5}
\| \Ga_g\|_{\infty, k}\leq\ep r^{-2},\qquad \|\Ga_b\|_{\infty, k}\leq\ep r^{-1}.
\eea
\end{enumerate}

{\bf Ref  2.}  Recall from \eqref{eq:valueofkstarinchapter5forproofThmM3} that $k_*=k_{small}+80$ in this chapter. The following estimates on $\Si_*$ are obtained in Theorems  M1 and M2:
\begin{enumerate}
\item The quantity $\qf$ satisfies, in view of Theorem M1, on $\Si_*$, for all  $0\le k\le k_{*}$,
 \bea
\bsplit
\| \qf\|_{\infty, k} &\les\ep_0 r^{-1} u^{-\frac{1}{2}-\dee},\\
\| \nab_3 \qf\|_{\infty, k-1} &\les\ep_0 r^{-1} u^{-1-\dee},
\end{split}
\eea
as well as 
\bea
\lab{improved-qf:onSi*}
\int_{\Sigma_*(u, u_*) } u^{2+2\dec}    |\nab_3\dk^{ k} \qf|^2    &\les \ep^2_0.
\eea

\item According to Theorem M2, $\aa$ satisfies,   for all $0\leq k\le k_{*}$,
\bea
\lab{improved-aa:onSi*}
\int_{\Sigma_*(u, u_*) } u^{2+2\dec}    |\dk^{ k}_*\aa|^2    &\les \ep^2_0
\eea

 \item According to Theorem  M1, $\a$ satisfies on $\Si_*$, for all  $0\le k\le k_{*}$,
\bea
\lab{improved-a:onSi*}
\begin{split}
\| \a\|_{\infty, k} &\les\ep_0r^{-\frac{7}{2}-\dee},\\
\|\nab_3\a\|_{\infty, k-1} &\les\ep_0r^{-\frac{9}{2}-\dee}.
\end{split}
\eea
\end{enumerate}

We conclude this section with an interpolation lemma.
\begin{lemma}\lab{lemma:interpolation}
We have for\footnote{Recall that $k_*= k_{small}+80$ throughout this chapter, see \eqref{eq:valueofkstarinchapter5forproofThmM3}.} $k\leq k_*$
\bea
\bsplit
\| \Ga_g\|_{\infty, k}&\les\ep r^{-2} u^{-\frac{1}{2}-\frac{\dec}{2}},\\
 \|\nab_\nu\Ga_g\|_{\infty, k-1}  &\les \ep  r^{-2} u^{-1-\frac{\dec}{2}}, \\   
\|\Ga_b\|_{\infty, k}&\les\ep  r^{-1} u^{-1-\frac{\dec}{2}}.
\end{split}
\eea
\end{lemma}

\begin{proof}
Since $k_{small}<k_*<k_{large}$, we have, interpolating between \eqref{Ref1-smallk:chap5} and \eqref{Ref1-largek:chap5}, for $k\leq k_*$,
\beaa
\bsplit
\| \Ga_g\|_{\infty, k}&\les\ep r^{-2}\left(u^{-\frac{1}{2}-\dec}\right)^{1-\frac{k_*-k_{small}}{k_{large}-k_{small}}},\\
 \|\nab_\nu\Ga_g\|_{\infty, k-1}  &\les \ep  r^{-2} \left(u^{-1-\dec}\right)^{1-\frac{k_*-k_{small}}{k_{large}-k_{small}}}, \\   
\|\Ga_b\|_{\infty, k}&\les\ep  r^{-1}\left(u^{-1-\dec}\right)^{1-\frac{k_*-k_{small}}{k_{large}-k_{small}}}.
\end{split}
\eeaa
Now, since $k_*$ satisfies\footnote{Note that we have in view of \eqref{eq:constraintsonthemainsmallconstantsepanddelta} and \eqref{eq:choiceksmallmaintheorem}
\beaa
\dec(k_{large}-k_{small})\geq \frac{1}{2}\dec k_{large}-\dec \gg 1,
\eeaa
and we may thus assume $\dec(k_{large}-k_{small})\geq 240$ so that, in view of  $k_*= k_{small}+80$, we have indeed $3(k_*-k_{small})=240\leq \dec(k_{large}-k_{small})$.}
\beaa
k_* \leq k_{small}+\frac{\dec}{3}(k_{large}-k_{small}),
\eeaa
we infer
\beaa
(2+2\dec)\frac{k_*-k_{small}}{k_{large}-k_{small}}\leq\dec
\eeaa
and hence
\beaa
&&\left(-\frac{1}{2}-\dec\right)\left(1-\frac{k_*-k_{small}}{k_{large}-k_{small}}\right)\\
 &=& -\frac{1}{2}-\frac{\dec}{2}-\frac{1}{2}\left(\dec -(1+2\dec)\frac{k_*-k_{small}}{k_{large}-k_{small}}\right) \leq  -\frac{1}{2}-\frac{\dec}{2}
\eeaa
as well as
\beaa
&&\left(-1-\dec\right)\left(1-\frac{k_*-k_{small}}{k_{large}-k_{small}}\right)\\
 &=& -1-\frac{\dec}{2}-\frac{1}{2}\left(\dec -(2+2\dec)\frac{k_*-k_{small}}{k_{large}-k_{small}}\right)\leq  -1-\frac{\dec}{2}.
\eeaa
This yields
\bea
\bsplit
\| \Ga_g\|_{\infty, k}&\les\ep r^{-2} u^{-\frac{1}{2}-\frac{\dec}{2}},\\
 \|\nab_\nu\Ga_g\|_{\infty, k-1}  &\les \ep  r^{-2} u^{-1-\frac{\dec}{2}}, \\   
\|\Ga_b\|_{\infty, k}&\les\ep  r^{-1} u^{-1-\frac{\dec}{2}}.
\end{split}
\eea
as stated.
\end{proof}

%%%%%%%%%%%%%%%%%%%%%

\subsection{Hodge operators}

%%%%%%%%%%%%%%%%%%%%%

We recall the following Hodge operators acting on  2-surfaces $S$ (see \cite{Ch-Kl} chapter 2 and \cite{KS}):

\begin{definition}\lab{def:HodgeoperatoronS:chap5} 
We  define   the following  Hodge operators:
\begin{enumerate}
\item The   operator $\ddd_1$ takes any 1-form $f$ into the pair of functions $(\div f, \curl f)$.

\item The   operator $\ddd_2$ takes any symmetric traceless 2-tensor $f$ into the 1-form $\div f$.

\item The   operator $\dds_1$ takes any pair of scalars $(h, \dual h)$ into the 1-form $-\nab h+\dual\nab\dual h$.

\item The   operator $\dds_2$ takes any 1-form $f$ into the symmetric traceless 2-tensor $-\frac{1}{2}\nab\hot f$.
\end{enumerate}
\end{definition}

One  can easily check that $  \dds_k $ is the formal adjoint on $L^2(S)$ of $\ddd_k$ for $k=1,2$. Moreover,
\bea
\begin{split}
\label{eq:dcalident}
\dds_1\, \ddd_1&=-\lap_1+K,\qquad \ddd_1\, \dds_1=-\lap_0,\\
\dds_2\, \ddd_2 &=-\f12\lap_2+K,\qquad \ddd_2\, \dds_2=-\f12(\lap_1+K).
\end{split}
\eea

Using \eqref{eq:dcalident} one can  prove the following (see  Chapter 2 in \cite{Ch-Kl}).
\begin{proposition} 
\label{prop:2D-hodge}
Let $(S,g)$ be
a compact manifold with Gauss curvature $K$. We have,

{\bf i.)}\quad The following identity holds for 1-forms  $f$
\beaa
\int_S\big(|\nab   f |^2+K|   f |^2\big)=\int_S|\ddd_1   f   |^2.
\eeaa

{\bf ii.)}\quad The following identity holds for symmetric traceless 2-tensors $f$
\beaa
\int_S\big(|\nab    f  |^2+2K| f  |^2\big)=2\int_S |\ddd_2   f  |^2.
\eeaa

{\bf iii.)}\quad The following identity holds for  scalars $f$
\beaa
  \int_S|\nab  f |^2=\int_S|\dds_1\, f|^2.
\eeaa

{\bf iv.)}\quad  The following identity holds for 1-forms $f$
\beaa
\int_S \big(|\nab  f   |^2-K|  f  |^2\big)=2\int_S|\dds_2   f   |^2.
\eeaa
\label{prop:hodgeident}
\end{proposition}

%%%%%%%%%%%%%%%%%%%%%%%%%%%%%%%%%%%%%%%

\subsection{Basic equations for the linearized quantities}

%%%%%%%%%%%%%%%%%%%%%%%%%%%%%%%%%%%%%%%

 \begin{proposition}
 \lab{Prop.NullStr+Bianchi-lastslice}
The following  equations hold true on $\Si_*$:
 \begin{enumerate}
\item  The linearized null structure equations are given by
 \beaa
 \bsplit
\nab_4\kac  &= \Ga_g\c\Ga_g,\\
 \nab_4\chih+ \frac{2}{r} \chih&=-\a,\\
 \nab_4 \ze+\frac{2}{r} \ze &=-\b+\Ga_g\c \Ga_g,\\
 \nab_4\kabc+\frac{1}{r}\kabc&= - 2 \div \ze + 2 \rhoc+\Ga_b\c \Ga_g,\\
 \nab_4 \chibh+\frac{1}{r} \chibh &= \frac{\Up}{r}\chih- \nab\hot \ze   +\Ga_b\c \Ga_g,
 \end{split}
 \eeaa
 and 
 \beaa
 \bsplit
 \nab_3\kac  &=         2   \div \eta + 2\rhoc -\frac{1}{r}\kabc + \frac{4}{r} \ombc +\frac{2}{r^2}\widecheck{y} +\Ga_b\c\Ga_b,\\
    \nab_3 \kabc -\frac{2\Up}{r}  \kabc &=2\div \xib +\frac{4\Up}{r} \ombc -\frac{2m}{r^2} \kabc  -\left(\frac{2}{r^2} -\frac{8m}{r^3} \right)\yc  +\Ga_b\c \Ga_b, \\
 \nab_3 \chibh -\frac{2\Up}{r}\chibh &=-\aa -\frac{2m}{r^2}\chibh +\nab\hot \xib +\Ga_b\c \Ga_b,\\
 \nab_3 \ze -\frac{\Up}{r} \ze &=-\bb-2 \nab \ombc +\frac{\Up}{r} (\eta+\ze)+\frac{1}{r}\xib+ \frac{2m}{r^2}(\ze-\eta) +\Ga_b\c \Ga_b,\\
 \nab_3 \chih -\frac{\Up}{r} \chih&=\nab\hot \eta -\frac{1}{r} \chibh+\frac{2m}{r^2}\chih+\Ga_b\c \Ga_b.
 \end{split}
 \eeaa
 Also,
\beaa
\div\chih &=& \frac{1}{r}\ze - \b+\Ga_g\c \Ga_g,\\
\div\chibh &=& \frac{1}{2}\nab\kabc +\frac{\Up}{r}\ze  +\bb+\Ga_b\c \Ga_g,
\eeaa
\beaa
\curl\ze&=& \rhod+\Ga_b\c \Ga_g,\\
\curl\eta &=& \rhod+\Ga_b\c \Ga_g,\\
\curl \xib &=&   \Ga_b\c \Ga_b,
\eeaa
and
\beaa
\widecheck{K}&=&-\frac{1}{2r} \kabc  -\rhoc +\Ga_b\c \Ga_g,\\
\widecheck{\mu} &=& -\div\ze -\rhoc +\Ga_b\c\Ga_g.
\eeaa

\item   The linearized Bianchi identities are given by
\beaa
\bsplit
   \nab_3\a -\frac{\Up}{r} \a&= \nab\hot \b +\frac{4m}{r^2}\a +\frac{6m}{r^3}\chih  +\Ga_b\c(\a,\b)+r^{-1}\Ga_g\c\Ga_g,\\
   \nab_4\beta +\frac{4}{r}\beta &=-\div\a  +r^{-1}\Ga_g\c\Ga_g,\\
      \nab_3\b -\frac{2\Up}{r}\b &= (\nab\rho+\dual\nab\rhod) +\frac{2m}{r^2}\b -\frac{6m}{r^3}\eta+r^{-1}\Ga_b\c\Ga_g,
      \end{split}
\eeaa
\beaa
\bsplit
       \nab_4 \rhoc+\frac{3}{r}\rhoc &=\div \b +r^{-1}\Ga_b\c\Ga_g,\\
           \nab_3 \rhoc  -\frac{3\Up}{r} \rhoc  &=  -\div\bb  +\frac{3m}{r^3}\kabc -\frac{6m}{r^4} \yc   -\frac{1}{2}\chih\c\aa+r^{-1}\Ga_b\c\Ga_b,\\
 \nab_4 \rhod+\frac{3}{r} \rhod&=-\curl\b+r^{-1}\Ga_b\c\Ga_g, \\
       \nab_3 \rhod -\frac{3\Up}{r} \rhod &=-\curl\bb - \frac 1 2 \chih\c\dual \aa+r^{-1}\Ga_b\c\Ga_b.
\end{split}
\eeaa
\end{enumerate}
 \end{proposition} 
 
 \begin{proof}
The proof  follow immediately from Proposition \ref{prop-nullstrandBianchion.integrable}, the definition of the linearized quantities, the definition of $\Ga_g$ and $\Ga_b$, the fact that $y=e_3(r)$, and the GCM condition $\kac=0$ on $\Si_*$.
 \end{proof}

%%%%%%%%%%%%%%%%%%%%%%%%

\subsection{Commutation lemmas}

%%%%%%%%%%%%%%%%%%%%%%%%

We start with the following lemma.
 \begin{lemma}
   \lab{lemma:comm_Si_*}
   For any tensor on $S$, the following  commutation formulas hold true
   \beaa
        \,[\nab_3, \nab_a] f &=&-\frac 1 2 \trchb \nab_a f +(\eta_a-\ze_a) \nab_3 f-\chibh_{ab}\nab_b f  +\xib_a \nab_4 f+(\underline{F}[f])_a,\\
         \,[\nab_4, \nab_a] f &=&-\frac 1 2 \trch \nab_a f+(\etab_a+\ze_a) \nab_4 f-\chih_{ab}\nab_b f+(F[f])_a.
       \eeaa
   where the tensors $F[f]$ and $\underline{F}[f]$ have the following schematic form
   \beaa
  F[f]=\big(\b, \chi\c\etab, \chib\c\xi\big)\c f, \qquad \underline{F}[f]= \big(\bb, \chib\c\eta, \chi\c\xib\big)\c f.
   \eeaa      
 \end{lemma}
 
\begin{proof}
See     Lemma 7.3.3 in \cite{Ch-Kl}.
\end{proof}

 \begin{lemma}
 \lab{Lemma:Commutation-Si_*}
 The following  commutation formulas hold true  for any  tensor $f$ on $S\subset\Si_*$:
 \begin{enumerate}
 \item We have
 \beaa
 \,[ \nab_3, \nab] f &=& \frac{\Up}{r}\nab f +\Ga_b \c\nab_3 f + r^{-1} \Ga_b \c \dk^{\leq 1} f,\\
  \,[ \nab_4, \nab] f &=&-\frac{1}{r}  \nab f +r^{-1} \Ga_g \c \dk^{\leq 1} f.
 \eeaa
  
 \item  We have
 \beaa
  \,[ \nab_3, \De] f&=& \frac{2\Up}{r} \lap  f+ r^{-1}\dkb^{\le 1}\Big( \Ga_b \c\nab_3 f + r^{-1} \Ga_b \c \dk f \Big),\\
   \,[ \nab_4, \De] f&=& -\frac{2}{r}\lap   f+ r^{-1} \dkb^{\le 1}\Big( \Ga_g \c  \dk f\Big).
 \eeaa
 
 \item We have 
 \beaa
\, [ \nab_{\nu}, \nab]f &=&  -\frac{2}{r}\nab f+\Ga_b \c \nab_\nu f+ r^{-1} \Ga_b \c \dk^{\leq 1} f.
 \eeaa

 \item We have 
 \beaa
 \, [ \nab_{\nu}, \De]f &=& - \frac{4}{r} \De f+ 
  r^{-1}\dkb^{\le 1}\Big( \Ga_b \c\nab_\nu f + r^{-1} \Ga_b \c \dk f \Big).
 \eeaa
 \end{enumerate}
 \end{lemma} 

 \begin{proof}
 The proof follows from Lemma \ref{lemma:comm_Si_*},  the definitions of $\Ga_g, \Ga_b$, the transversality conditions $\xi=0$ and $\etab=-\ze$ on $\Si_*$, and $\nu=e_3+b_*e_4$.
\end{proof}
 
\begin{corollary}
\lab{Corr:commutation-onSi_*}
The following  commutation formulas hold true  for any  tensor $f$ on $S\subset\Si_*$:
   \bea
   \bsplit
 \,[ \nab_3,  r\nab] f &= r \Ga_b \c\nab_3 f + \Ga_b \c \dk^{\leq 1}f,\\
  \,[ \nab_4,  r \nab] f &=\Ga_g \c \dk^{\leq 1}f,\\
  \,[ \nab_3,  r^2\De] f &=r \dkb^{\le 1}\Big( \Ga_b \c\nab_3 f + r^{-1} \Ga_b \c \dk f\Big),\\
   \,[ \nab_4,  r^2 \De] f &=  r \dkb^{\le 1}\Big( \Ga_g \c  \dk f\Big),
   \end{split}
 \eea
and 
 \bea
 \bsplit
\, [ \nab_{\nu},  r \nab]f &=  r \Ga_b \c \nab_{\nu} f + \Ga_b \c \dk^{\leq 1} f,\\
 \, [ \nab_{\nu},  r^2 \De]f &=   r \dkb^{\le 1}\Big( \Ga_b \c\nab_{\nu}  f \Big) + \dkb^{\le 1}\Big( \Ga_b \c \dk f\Big).
 \end{split}
 \eea
 \end{corollary}

 \begin{proof}
 The proof follows from Lemma \ref{Lemma:Commutation-Si_*},  the transversality condition $e_4(r)=1$ on $\Si_*$, and the fact that $e_3(r)=-\Up+r\Ga_b$ and $\nu(r)=-2+r\Ga_b$.
\end{proof}

%%%%%%%%%%%%%%%%%%%%%%%%%%

 \subsection{Additional equations}

%%%%%%%%%%%%%%%%%%%%%%%%%%

 \begin{proposition}
\lab{prop:additional.eqtsM4}
We have, schematically,
\beaa
2 \nab \ombc -\frac{1}{r}\xib &=& - \nab_3 \ze   -\bb +\frac{1}{r}\eta +r^{-1}\Ga_g+\Ga_b\c \Ga_b,\\
2\ddd_2\dds_2\eta &=&  -\nab_3\nab\kac    -\frac{2}{r}\nab_3 \ze  -\frac{2}{r}\bb  + r^{-2} \dkb^{\le 1} \Ga_g+ r^{-1} \dkb^{\le 1 } (\Ga_b\c \Ga_b),\\
 2\ddd_2\dds_2\xib  &=& -\nab_3\nab\kabc   -\frac{2}{r}\nab_3 \ze    -\frac{2}{r}\bb     +r^{-2}\dkb^{\leq 1}\Ga_g+r^{-1}\dkb^{\leq 1}(\Ga_b\c \Ga_b).
 \eeaa
\end{proposition}

\begin{proof}
See section \ref{sec:proofofprop:additional.eqtsM4} in the  appendix.  
\end{proof} 

The following   corollary   of Proposition \ref{prop:additional.eqtsM4}  will  be  very useful.
\begin{proposition}
\lab{Prop:nu*ofGCM:0}
 The following identities hold true   on $\Si_*$:
    \bea
   \lab{eq:D^5eta-prop:0}
  \bsplit
 2\dds_2\dds_1 \ddd_1\ddd_2\dds_2\eta &=  -\dds_2\dds_1 \ddd_1\nab_3\nab\kac    +\frac{2}{r}\nab_3\dds_2\dds_1\muc -\frac{4}{r}\dds_2\dds_1\div\bb        \\
& + r^{-5} \dkb^{\le 4} \Ga_g+ r^{-4} \dkb^{\le 4} (\Ga_b\c \Ga_b),
  \end{split} 
  \eea  
\bea
\lab{eq:D^5xib-prop:0}
\bsplit
 2\dds_2\dds_1 \ddd_1\ddd_2\dds_2\xib  &=  \nab_3\left(\dds_2\ddd_2+\frac{2}{r^2}\right)\dds_2\dds_1\kabc   +\frac{2}{r}\nab_3\dds_2\dds_1\muc   -\frac{4}{r}\dds_2\dds_1 \div\bb  \\
&   + r^{-5} \dkb^{\le 4} \Ga_g+ r^{-4} \dkb^{\le 4} (\Ga_b\c \Ga_b),
\end{split}
\eea
where by convention,  for any scalar $f$,
\beaa
\dds_1 f:= \dds_1(f, 0)=-\nab f.
\eeaa
\end{proposition}

\begin{proof}
See section \ref{sec:proofofProp:nu*ofGCM:0}  in the  appendix.  
\end{proof} 

 \begin{corollary}
\lab{Corr:nu*ofGCM}
 The following identities hold true   on $\Si_*$.
   \bea
   \lab{eq:D^5eta-prop}
  \bsplit
  2\dds_2\dds_1 \ddd_1\ddd_2\dds_2\eta &=  -\dds_2\dds_1 \ddd_1\nab_\nu\nab\kac    +\frac{2}{r}\nab_\nu\dds_2\dds_1\muc -\frac{4}{r}\dds_2\dds_1\div\bb        \\
& + r^{-5} \dkb^{\le 5} \Ga_g+ r^{-4} \dkb^{\le 4} (\Ga_b\c \Ga_b),
  \end{split} 
  \eea  
\bea
\lab{eq:D^5xib-prop}
\bsplit
 2\dds_2\dds_1 \ddd_1\ddd_2\dds_2\xib  &=  \nab_\nu\left(\dds_2\ddd_2+\frac{2}{r^2}\right)\dds_2\dds_1\kabc   +\frac{2}{r}\nab_\nu\dds_2\dds_1\muc   -\frac{4}{r}\dds_2\dds_1 \div\bb  \\
&   + r^{-5} \dkb^{\le 5} \Ga_g+ r^{-4} \dkb^{\le 4} (\Ga_b\c \Ga_b).
\end{split}
\eea
 \end{corollary}

\begin{proof}
We have in view of Proposition \ref{Prop.NullStr+Bianchi-lastslice}
\beaa
\nab_4\kac\in \Ga_g\c\Ga_g, \qquad \nab_4\kabc\in r^{-1}\dkb^{\leq 1}\Ga_g, \qquad \nab_4\muc\in r^{-2}\dkb^{\leq 1}\Ga_g, 
\eeaa
which together with  $b_*=-1-\frac{2m}{r}+r\Ga_b$ yields
\beaa
 &&\dds_2\dds_1 \ddd_1 b_*\nab_4\nab\kac \in r^{-5} \dkb^{\le 5} \Ga_g,\qquad \frac{b_*}{r}\nab_4\dds_2\dds_1\muc \in r^{-5} \dkb^{\le 5} \Ga_g,\\
 && b_*\nab_4\left(\dds_2\ddd_2+\frac{2}{r^2}\right)\dds_2\dds_1\kabc  \in r^{-5} \dkb^{\le 5} \Ga_g.
\eeaa
Since $\nu=e_3+b_*e_4$, the proof of \eqref{eq:D^5eta-prop} and \eqref{eq:D^5xib-prop} follows then  immediately  from Proposition \ref{Prop:nu*ofGCM:0}. 
\end{proof}

%%%%%%%%%%%%%%%%%%%%%%%%%%%%%%%%%%%%%%%

\subsection{Additional   renormalized  equations on $\Si_*$} 

%%%%%%%%%%%%%%%%%%%%%%%%%%%%%%%%%%%%%%%

\begin{lemma}
\lab{Lemma:transport.alongSi_*1} 
We have along  $\Si_*$
\bea
\bsplit
\nab_\nu\left( \lap\kabc+\frac{2\Up}{r}\div \ze\right)  &= O(r^{-1})\lap\kabc+2\lap\div \xib  +O(r^{-2})\lap\yc  +O(r^{-2})\div\ze \\
& +O(r^{-1})\div\bb +O(r^{-2})\div\eta+O(r^{-2})\div\xib   \\
& +2\left(1+O(r^{-1})\right)\lap\div \ze  - 2\left(1+O(r^{-1})\right)\lap\rhoc\\
& +O(r^{-1})\div\b  +r^{-2}\dkb^{\le 2 }(\Ga_b\c \Ga_b),\\
         \nab_\nu\div\b  &=  O(r^{-1})\div \b+\lap\rho +(1+O(r^{-1}))\div\div\a\\
         & +O(r^{-3})\div\eta+r^{-2}\dkb^{\leq 1}(\Ga_b\c\Ga_g),\\
    \nab_\nu\curl\b  &=  \frac{8}{r}(1+ O(r^{-1}))\curl \b -\lap\rhod +(1+O(r^{-1}))\curl\div\a \\
    &+O(r^{-3})\rhod+r^{-2}\dkb^{\leq 1}(\Ga_b\c\Ga_g),\\  
   \nab_\nu\left(\rhoc - \frac{1}{2}\chih\c\chibh \right)    &=  -\div\bb -(1+O(r^{-1}))\div\b +O(r^{-1})\rhoc   +O(r^{-3})\kabc\\
   & +O(r^{-4}) \yc   +r^{-1}\dkb^{\le  1}( \Ga_b \c \Ga_b),
  \end{split}
\eea
where the notation $O(r^a)$, for $a\in\mathbb{R}$, denotes an explicit function of $r$ which is bounded 
by $r^a$ as $r\to+\infty$. 
\end{lemma}

\begin{proof}
See section \ref{sec:proofofLemma:transport.alongSi_*1} in the  appendix.  
\end{proof}

%%%%%%%%%%%%%%%%%%%%%%%%%%

\subsection{Equations  involving   $\qf$}

%%%%%%%%%%%%%%%%%%%%%%%%%%

\begin{proposition}
\lab{prop:identitiesinqf}
Let $O(r^a)$ denote, for $a\in\RRR$, a function of $(r, \cos\th)$ bounded by $r^a$ as $r\to +\infty$. In the frame of $\Si_*$, the following holds:
\begin{enumerate}
\item We have
\bea
\lab{eq:alternateformulaforqfinvolvingtwoangularderrivativesofrho:M4}
\Re(\qf) &=&  r^4\dds_2\dds_1(-\rho, \rhod)  +O(r^{-2})+\dkb^{\leq 2}\Ga_b+r^{2}\dkb^{\le 2}\big(\Ga_b \c  \Ga_g\big).
\eea
 
\item  We have
\bea
\label{eq:Le-Teuk-Star1-M4}
\nn\Re(\nab_3(r\qf)) &=& r^5\dds_2\dds_1\big(\div\bb, -\curl\bb\big) +O(r)\dkb^{\leq 3}\aa\\
&& +O(r^{-2})+r\dkb^{\leq 3}\Ga_g+r^{3}\dkb^{\le 3}(\Ga_b \c  \Ga_g).
 \eea
\end{enumerate}
\end{proposition}

\begin{proof}
See section \ref{section:appendix-proofofprop:identitiesinqf}
 in the appendix.
\end{proof}

%%%%%%%%%%%%%%%%%%%%%%%%%%

\subsection{Hodge  elliptic systems}     
\lab{section:Hodge-elliptic systemsChapter5} 

%%%%%%%%%%%%%%%%%%%%%%%%%%

For a tensor $f$ on $S$, we define the following norms for any integer $k\geq 0$
\bea
\|f\|_{\hk_k(S)} &:=& \sum_{j=0}^k\|\dkb^jf\|_{L^2(S)}.
\eea
We record below the following   well known coercive properties of the operators $\ddd_1, \ddd_2, \dds_1$, see Chapter 2 in \cite{Ch-Kl}.

 \begin{lemma} 
\label{prop:2D-Hodge1}
For any sphere  $S=S(u) \subset \Si_* $   we have, for all $k\le k_{large}$:
\begin{enumerate}
\item  If   $f$ is a 1-form
\bea
 \|  f\|_{\hk_{k+1}   (S)}    \les r   \|\ddd_1   f  \|_{\hk_k(S)}.
\eea

\item If $f$ is a symmetric traceless 2-tensor
\bea
 \|  v\|_{\hk_{k+1}   (S)}    \les r   \|\ddd_2   v  \|_{\hk_k(S)}.
\eea

\item  If  $(h, \dual h)$ is a pair of scalars 
\bea
  \| (h-\ov{h}, \dual h-\ov{\dual h})\|_{\hk_{k+1} (S)} \les r  \|\dds_1\, (h, \dual h)\|_{\hk_{k}(S)}.
\eea
\end{enumerate}
\end{lemma}

\begin{proof}
Recall from Proposition \ref{Prop.NullStr+Bianchi-lastslice} the linearized Gauss equation
\beaa
K &=& -\frac{1}{2r}\kabc -\rhoc  +\Ga_b\c\Ga_g.
\eeaa
In view of the control \eqref{Ref1-largek:chap5} for $\Ga_g$ and $\Ga_b$, we have 
\beaa
\left\|\dkb^{\leq k_{large}}\left(K-\frac{1}{r^2}\right)\right\|_{L^\infty(S)} &\les& \frac{\ep}{r^3}.
\eeaa 
Together with Proposition \ref{prop:2D-hodge}, and a Poincar\'e inequality for $\dds_1$, this immediately yields the case $k=0$. The case $k\geq 1$ then follows by  the above control of $K$ and elliptic regularity.
\end{proof}

 The operator $\dds_2$ is not coercive but satisfies  instead  the following estimates.
 
 \begin{lemma} 
\label{prop:2D-Hodge2}
On a fixed sphere  $S=S(u) \subset \Si_* $,   we have for any 1-form $f$ and  all $k\le k_{large}$
\bea
  \|  f\|_{\hk_{k+1}  (S)}&\les&  r \|\dds_2\, f\|_{\hk_k(S)}+r^2 \big| (\ddd_1f)_{\ell=1}\big|.
  \eea
  Note also the straightforward  inequality
  \bea
  \lab{estimate:badmodeless}
  |(\ddd_1f)_{\ell=1}| &\les & r^{-1} \|\ddd_1f\|_{L^2(S)}.
  \eea
  \end{lemma}
  
\begin{proof}
The case $k=0$  is proved  in Lemma 2.19 of \cite{KS-GCM1}. The higher derivative estimates follow by elliptic regularity and the above control of $K$.
\end{proof}
 
 \begin{corollary} 
\label{prop:2D-Hodge3}
On a fixed sphere  $S=S(u) \subset \Si_* $,   we have for any pair of scalars  $(h, \dual h)$  and  all $k\le k_{large}-1$
\beaa
  \|(h-\ov{h}, \dual h-\ov{\dual h})\|_{\hk_{k+2}  (S)} &\les&  r^2\|\dds_2\,\dds_1(h, \dual h)\|_{\hk_k(S)}+r^3\big| (\Delta h)_{\ell=1}\big|+r^3\big| (\Delta\dual h)_{\ell=1}\big|.
  \eeaa
  \end{corollary}

\begin{proof}
Applying Lemma \ref{prop:2D-Hodge2} with $f=\dds_1(h, \dual h)$ yields the control of $\dds_1(h, \dual h)$. The control of $(h-\ov{h}, \dual h-\ov{\dual h})$ follows then from applying Lemma \ref{prop:2D-Hodge1}.
\end{proof}

%%%%%%%%%%%%%%%%%%%%%%%%%%%%

\section{Preliminary  estimates on $\Si_*$}
\lab{sec:preliminaryestimatesonSistar:chap5}

%%%%%%%%%%%%%%%%%%%%%%%%%%%%

 %%%%%%%%%%%%%%%%%%%%%%%%%%%
 
 \subsection{Behavior of $r$ on $\Si_*$}
 
 %%%%%%%%%%%%%%%%%%%%%%%%%%

The following lemma shows that $r$ is comparable to $r_*=r(S_*)$ on $\Si_*$.
\begin{lemma}
The function  $r$  satisfies on $\Si_*$ 
\bea\lab{eq:behaviorofronSigma*:almost constant}
 r_*\leq r\leq r_*\left(1+O(\ep_0)\right).
\eea
\end{lemma}

\begin{proof}
Recall that $u(S_*)=u_*$, $u(S_1)=1$ and $1\leq u\leq u_*$ along $\Si_*$. Since $r$ is a decreasing function of $u$ on $\Si_*$, we infer
\beaa
r_*\leq r\leq r(S_1)\quad \textrm{on}\quad\Si_*.
\eeaa
To conclude, it thus suffices to prove that $r(S_1)\leq r_*(1+O(\ep_0))$. To this end, recall first \eqref{eq:relationbetweenrandudefininguonSigmastar}, i.e.,  $u+r=c_*$ along $\Si_*$, so that 
\beaa
r(S_*)+u_*=c_*=r(S_1)+1
\eeaa
and hence
\beaa
r(S_1)=r_*+u_*-1=r_*\left(1+\frac{u_*-1}{r_*}\right)\leq r_*\left(1+\frac{u_*}{r_*}\right).
\eeaa
Together with the dominant condition on $r$ \eqref{eq:dominantconditiononronSigmastarchap5} on $\Si_*$, i.e., $r_*=\de_*\ep_0^{-1}u_*^{1+\dec}$, we infer $r(S_1)\leq r_*(1+O(\ep_0))$ which concludes the proof of the lemma.
\end{proof}

%%%%%%%%%%%%%%%%%%%%%%%%%%

\subsection{Transport lemmas along $\Si_*$}

%%%%%%%%%%%%%%%%%%%%%%%%%%

\begin{lemma}
\lab{Lemma:nuSof-integrals-again}
For any scalar function $h$ on $\Si_*$, we have  the formula
\bea
\nu\left(\int_Sh\right) &=&  z\int_S\frac{1}{z}\Big(\nu(h)+(\kab+b_*\ka) h     \Big). 
\eea
In particular\footnote{Recall that $r$ denotes  the area radius of GCM spheres along $\Si_*$.}  
 \bea
 \nu(r)  = \frac{rz}{2}\,\ov{z^{-1}(\kab + b_*\ka)}
 \eea
 where for a scalar function $f$, $\ov{f}$ denotes the average of $f$ with respect to the spheres of  the foliation of $\Si_*$.
\end{lemma}

\begin{proof}
Recall that on $\Si_*$ the coordinates $\th, \vphi$ where defined s.t. $\nu(\th)=\nu(\vphi)=0$. Moreover, in view of the transversality condition $e_4(u)=0$ on $\Si_*$, we have $\nu(u)=e_3(u)$. Thus $\nu(u)=z$ and  hence $\nu=z\pr_u$ in the $(u, \th, \vphi)$ coordinates system. We infer
\beaa
\nu\left(\int_Sh\right) &=& z\pr_u\left(\int_Sh\right)=z\int_S\Big(\pr_uh+g^{ab}g(\D_a\pr_u, e_b)\Big)=z\int_S\frac{1}{z}\Big(\nu(h)+g^{ab}g(\D_a\nu, e_b)\Big)\\
&=& z\int_S\frac{1}{z}\Big(\nu(h)+g^{ab}g(\D_a(e_3+b_*e_4), e_b)\Big)=z\int_S\frac{1}{z}\Big(\nu(h)+(\kab+b_*\ka) h     \Big)
\eeaa
which yields the first identity. The second identity follows then by choosing $h=1$ in the first identity.
\end{proof}

\begin{corollary}
\lab{Corr:nuSof integrals}
For any scalar function $h$ on $\Si_*$, we have  
 \bea
\nu\left(\int_Sh\right) = \int_S\nu (h)   - \frac{4}{r} \int_S h +r^3\Ga_b\nu(h)+r^2\Ga_b h
\eea
 and 
 \bea
 \nu(r)=  - 2 + r \Ga_b. 
 \eea
\end{corollary}

\begin{proof}
Since we have
 \beaa
  \ka =\frac 2 r +\Ga_g,\qquad  \kab=-\frac{2\Up}{r}+\Ga_g, \qquad b_*=  -1-\frac{2m}{r} +r \Ga_b,
  \eeaa
we infer
\bea
\lab{eq:kab+b_*ka}
\kab + b_*\ka&=& -\frac{4}{r}  + \Ga_b.
\eea
The proof follows then easily from Lemma \ref{Lemma:nuSof-integrals-again}, \eqref{eq:kab+b_*ka} and the fact that $z=2 +r\Ga_b$.
 \end{proof}

We now control transport equation in $\nu$ along $\Si_*$.
\begin{lemma}
\lab{lemma:integrate-transportSi_*}
Let $n$ and $m$ two integers, and let $f$ and $h$ two scalar functions on $\Si_*$. Assume that $f$ satisfies along $\Si_*$
\bea
\nu(f) -\frac{n}{r}f &=& h+\Ga_b f.
\eea
Then, we have for all $1\leq u\leq u_*$
\bea
\|r^m f\|_{L^\infty(S(u))} \les r_*^m\|f\|_{L^\infty(S_*)}+\int_u^{u_*}\Big(\|r^mh\|_{L^\infty(S(u'))}+\|r^m\Ga_bf\|_{L^\infty(S(u'))}\Big)du'
\eea
where the inequality is uniform in $u$.
\end{lemma}

\begin{proof}
We rewrite the transport equation for $f$ as
\beaa
\nu(r^{\frac{n}{2}}f) &=& r^{\frac{n}{2}}\left(\nu(f)+\frac{n}{2}\frac{\nu(r)}{r}f\right)=r^{\frac{n}{2}}\left(\frac{n}{r}f + h+\Ga_bf+\frac{n}{2}\frac{-2+r\Ga_b}{r}f\right)\\
&=& r^{\frac{n}{2}}\left(h+\Ga_b f\right)
\eeaa
where we have used $\nu(r)=-2+r\Ga_b$. Integrating from $S_*$, we deduce for all $1\leq u\leq u_*$
\beaa
\|r^{\frac{n}{2}} f\|_{L^\infty(S(u))} \les r_*^{\frac{n}{2}}\|f\|_{L^\infty(S_*)}+\int_u^{u_*}\Big(\|r^{\frac{n}{2}}h\|_{L^\infty(S(u'))}+\|r^{\frac{n}{2}}\Ga_bf\|_{L^\infty(S(u'))}\Big)du'
\eeaa
where the inequality is uniform in $u$. Multiplying this estimate by $r_*^{m-\frac{n}{2}}$, and since $r$ is comparable to $r_*$ on $\Si_*$ in view of \eqref{eq:behaviorofronSigma*:almost constant}, we infer the stated estimate.
\end{proof}

 \begin{corollary}
 \lab{cor:integrate-transportSi_*}
 Let $n$ and $m$ two integers, let $s\geq 0$ a positive real number, and let $f$ and $h$ two scalar functions on $\Si_*$. Assume that $f$ satisfies along $\Si_*$
\bea
\nu(f) -\frac{n}{r}f &=& h+\Ga_bf.
\eea
Assume also that there exists a constant $C>0$ such that 
\bea
\sup_{1\leq u\leq u_*}u^s\left(r_*^m\|f\|_{L^\infty(S_*)}+\int_u^{u_*}\|r^mh\|_{L^\infty(S(u'))}\right) &\leq& C.
\eea
Then, we have
\bea
\sup_{\Si_*}r^m u^s|f| &\les& C.
\eea
\end{corollary}
 
 \begin{proof}
 In view of Lemma \ref{lemma:integrate-transportSi_*}, there exists a constant $C_0>0$, uniform in $1\leq u\leq u_*$, such that  for all $1\leq u\leq u_*$
\beaa
\|r^m f\|_{L^\infty(S(u))} &\leq& C_0r_*^m\|f\|_{L^\infty(S_*)}+C_0\int_u^{u_*}\Big(\|r^mh\|_{L^\infty(S(u'))}+\|r^m\Ga_bf\|_{L^\infty(S(u'))}\Big)du'
\eeaa
and hence, in view of the assumptions of the corollary, and the assumption on $\Ga_b$, we infer
\beaa
u^s\|r^mf\|_{L^\infty(S(u))} &\leq& C_0C+u^s\ep\int_u^{u_*}\frac{\|r^m f\|_{L^\infty(S(u'))}}{{u'}^{1+\dec}}du'
\eeaa 
 and the proof easily follows by bootstrap from $u=u_*$ and the fact that $s\geq 0$. 
 \end{proof}

%%%%%%%%%%%%%%%%%%%%%%%%%%%%%%%%%%%
 
\subsection{Control of $\phi$ and the $\ell=1$ basis  $\Jp$ on $S_*$} 
 
%%%%%%%%%%%%%%%%%%%%%%%%%%%%%%%%%%%

In this section, we will rely on the following control of the Gauss curvature 
\bea\lab{eq:controlofkstarangularderivativesofKonSstar}
\sup_{S_*}\left|\dkb^{\leq k_*}\left(K-\frac{1}{r^2}\right)\right| &\les& \frac{\ep}{r^3u^{\frac{1}{2}+\frac{\dec}{2}}},
\eea 
which follows immediately from the fact that $\widecheck{K}\in r^{-1}\Ga_g$ in view of the linearized Gauss equation of Proposition \ref{Prop.NullStr+Bianchi-lastslice}, and  the control of Lemma \ref{lemma:interpolation} for $\Ga_g$.

Recall that $\phi$ is  the conformal factor of the metric $g$ on $S_*$, see \eqref{eq:formofthemetriconSstarusinguniformization}. We have the following lemma.
 
\begin{lemma}\lab{lemma:controloftheconformalfactorphi}
We have on $S_*$ 
\bea
\lab{eq:controlofphionSstar}
\|\dkb^{\leq k_*}\phi\|_{L^\infty(S_*)} &\les& \frac{\ep}{ru^{\frac{1}{2}+\frac{\dec}{2}}}.
\eea
\end{lemma}

\begin{proof}
In view of
\begin{itemize}
\item the control \eqref{eq:controlofkstarangularderivativesofKonSstar} of $K$ on $S_*$,

\item the fact that $g$ is conformal to the metric of $\SSS^2$ with conformal factor $r^2e^{2\phi}$, see \eqref{eq:formofthemetriconSstarusinguniformization},

\item and the balanced condition \eqref{eq:balancedconditionforJponSstar} for the $\ell=1$ modes $J^{(p)}$, 
\end{itemize}
we are in position to apply Corollary 3.8 in \cite{KS-GCM2} (restated here as Theorem \ref{Thm:effectiveU1-Intro}) which yields\footnote{Note that the vanishing integrals in \eqref{eq:centeredconditiononutobequotedlaterinchap5} are equivalent to the balanced condition \eqref{eq:balancedconditionforJponSstar}.}
\beaa
\|\phi\|_{\hk_{k_*+2}(S_*)}\les \frac{\ep}{ru^{\frac{1}{2}+\frac{\dec}{2}}}.
\eeaa
Together with Sobolev, this concludes the proof of the lemma.
\end{proof}
 
\begin{lemma}\lab{lemma:statementeq:DeJp.S_*}
We have on $S_*$
\bea
\lab{eq:DeJp.S_*}
\bsplit
\int_{S_*}J^{(p)} &= 0,\\
\int_{S_*}J^{(p)}J^{(q)} &= \frac{4\pi}{3}r^2\de_{pq}+O\left(\ep ru^{-\frac{1}{2}-\dec}\right),\\
\left\|\dkb^{\leq k_*}\left(\Delta \Jp+\frac{2}{r^2}\Jp\right)\right\|_{L^\infty(S_*)} & =  O\left(\ep r^{-3} u^{-\frac{1}{2}-\frac{\dec}{2}}\right).
\end{split}
\eea
\end{lemma}

\begin{proof}
The first identity in \eqref{eq:DeJp.S_*} is the balanced condition \eqref{eq:balancedconditionforJponSstar}, so we focus on the two other identities. In view of Lemma 3.12 in \cite{KS-GCM2} (restated here as Lemma \ref{lemma:basicpropertiesofJforcanonicalell=1basis}), we have\footnote{In view of \eqref{eq:formofthemetriconSstarusinguniformization}--\eqref{eq:balancedconditionforJponSstar}, $J^{(p)}$ is a canonical basis of $\ell=1$ mode of $S_*$ in the terminology of Definition 3.10 in \cite{KS-GCM2} (restated here as Definition \ref{definition:ell=1mpdesonS-intro}), so Lemma 3.12 in \cite{KS-GCM2} applies.}
\beaa
\int_{S_*}J^{(p)}J^{(q)} &=& \frac{4\pi}{3}r^2\de_{pq}+\int_{S_*}J^{(p)}J^{(q)}(1-e^{-2\phi}),\\
\Delta \Jp+\frac{2}{r^2}\Jp &=& \frac{2}{r^2}(1-e^{-2\phi})\Jp,
\eeaa
and the last two identities in \eqref{eq:DeJp.S_*} follow from the control \eqref{eq:controlofphionSstar} of $\phi$.
\end{proof}

%%%%%%%%%%%%%%%%%%%%%%%%%%%%%%%%
 
\subsection{Properties of the $\ell=1$ basis  $\Jp$ on $\Si_*$} 
 
%%%%%%%%%%%%%%%%%%%%%%%%%%%%%%%%

We are ready to  derive  the  basic properties of  the $\ell=1$ basis on $\Si_*$. 

\begin{lemma}
\lab{Le:Si*-ell=1modes}
The     functions $\Jp$ verify the following properties
\begin{enumerate}
\item We have on $\Si_*$
\bea
\lab{eq:basicestimatesforJp-onSi_*} 
\bsplit
\int_{S}J^{(p)} &= O\left(\ep ru^{-\dec}\right),\\
\int_{S}J^{(p)}J^{(q)} &= \frac{4\pi}{3}r^2\de_{pq}+O\left(\ep ru^{-\dec}\right).
\end{split}
\eea

\item We have on $\Si_*$
\bea
\nab_\nu \Big[\big( r^2\De +2)\Jp\Big]=O\big(\dkb^{\le 1} \Ga_b\big). 
\eea

\item For  any  $k\le k_*-1$, we have on $\Si_*$
\beaa
\left|\dk_*^k\left( \De +\frac{2}{r^2}\right)    \Jp\right|    \les  \ep r^{-3} u^{-\frac{\dec}{2}}.
\eeaa

\item We have for any $k\le k_*-3$ on $\Si_*$
\beaa
\left|\dk_*^k\dds_2\dds_1 \Jp\right|    \les  \ep r^{-3} u^{-\frac{\dec}{2}},
\eeaa
where by $\dds_1\Jp$, we mean either $\dds_1(\Jp,0)$ or $\dds_1(0,\Jp)$. 
\end{enumerate}
\end{lemma}

\begin{proof}
We proceed in steps as follows.

{\bf Step 1.} Since $\nu(\Jp)=0$ and $\Jp=O(1)$, we have, in view of Corollary \ref{Corr:nuSof integrals}
\beaa
\nu\left(\int_S J^{(p)}\right) &=& - \frac{4}{r} \int_S J^{(p)} +r^2\Ga_b,\\
\nu\left(\int_SJ^{(p)}J^{(q)}\right) &=& - \frac{4}{r} \int_S J^{(p)}J^{(q)} +r^2\Ga_b.
\eeaa
Since $\nu(r)=-2+r\Ga_b$, we infer
\beaa
\nu\left(r^{-2}\int_S J^{(p)}\right) &=& \Ga_b,\\
\nu\left(r^{-2}\int_SJ^{(p)}J^{(q)}-\frac{4\pi}{3}\de_{pq}\right) &=& \Ga_b.
\eeaa
Applying Corollary \ref{cor:integrate-transportSi_*} to the above transport equations with the choices $s=\dec>0$, $n=0$ and $m=1$, and using the control on $S_*$ provided by \eqref{eq:DeJp.S_*} and the control of $\Ga_b$ on $\Si_*$, we infer \eqref{eq:basicestimatesforJp-onSi_*}.

{\bf Step 2.} Using the  commutation formulas of Corollary \ref{Corr:commutation-onSi_*} we   have, since $\nu(\Jp)=0$,
\beaa
\nab_\nu\Big[\big( r^2\De +2\big)\Jp\Big]&=& \big[ \nab_{\nu} , r^2 \De \big] \Jp= r \dkb^{\le 1}\Big( \Ga_b \c\nab_{\nu}  \Jp \Big) + \dkb^{\le 1}\Big( \Ga_b \c \dk \Jp\Big)\\
&= & \dkb^{\le 1}\Big( \Ga_b \c \dk \Jp \Big)=           \dkb^{\le 1} \Ga_b 
\eeaa
as stated.

 {\bf Step 3.} Commuting the identity in Step 2 with $\dkb^k$, and using Corollary \ref{Corr:commutation-onSi_*}, we have 
\beaa
\nab_\nu\Big[\dkb^k\big( r^2\De +2\big)\Jp\Big] &=&         \dkb^{\le k+1} \Ga_b. 
\eeaa
Applying Corollary \ref{cor:integrate-transportSi_*} to the above transport equations with the choices $s=\frac{\dec}{2}>0$, $n=0$ and $m=1$, and using the control on $S_*$ provided by \eqref{eq:DeJp.S_*} and the control of $\Ga_b$ on $\Si_*$, we infer, for $k\leq k_*-1$, 
\beaa
\left|\dkb^k\left( \De +\frac{2}{r^2}\right)    \Jp\right|    \les  \ep r^{-3} u^{-\frac{\dec}{2}}.
\eeaa 
Together with the transport equation of Step 2 and the control of $\Ga_b$, we infer, for  any  $k\le k_*-1$, 
\beaa
\left|\dk_*^k\left( \De +\frac{2}{r^2}\right)    \Jp\right|    \les  \ep r^{-3} u^{-\frac{\dec}{2}}
\eeaa
as stated.

{\bf Step 4.} We introduce the notation $F:= \dds_2\dds_1(\Jp)$, where by $\dds_1\Jp$, we mean either $\dds_1(\Jp,0)$ or $\dds_1(0,\Jp)$. We have 
 \beaa
 \ddd_1\ddd_2 F&=& \ddd_1\ddd_2\dds_2\dds_1(\Jp).
 \eeaa
 Using  the identity $2\ddd_2\dds_2=\dds_1\ddd_1-2K$, we deduce\footnote{In fact $\De\Jp$ should be replaced by either $(\De\Jp, 0)$ or $(0, \De\Jp)$ depending whether we consider $\dds_1(\Jp,0)$ or $\dds_1(0,\Jp)$.} 
 \beaa
  2 \ddd_1\ddd_2 F&=& \ddd_1(\dds_1\ddd_1-2K)\dds_1(\Jp)\\
  &=&  \ddd_1\dds_1 \ddd_1\dds_1 \Jp - 2  K \ddd_1 \dds_1(\Jp)  +(\nab K, \dual\nab K)\c\dds_1(\Jp)\\
  &=& \De^2\Jp  +2 K \De \Jp  +(\nab K, \dual\nab K)\c\dds_1(\Jp).
 \eeaa
 Therefore,  since  $\widecheck{K} =K- r^{-2}\in  r^{-1} \Ga_g $, we have
 \beaa
  2 \ddd_1\ddd_2 F&=& \De^2\Jp  +\frac{2}{r^2} \De \Jp +2\widecheck{K}\De \Jp +(\nab \widecheck{K}, \dual\nab \widecheck{K})\c\dds_1(\Jp)\\
    &=&   \De\left(  \De  + \frac{2}{r^2}  \right)\Jp    + r^{-3} \dkb^{\le 1} \Ga_g.
 \eeaa
 In view of the coercivity of  the  operators $\ddd_1$ and $\ddd_2$, see Lemma \ref{prop:2D-Hodge1}, and in view of the definition of $F$, we deduce for $k\geq 2$
 \beaa
 \| \dds_2\dds_1 \Jp\|_{\hk_k(S)} &\les& 
 \left \|  \left(\De  + \frac{2}{r^2} \right)\Jp  \right\|_{\hk_k(S)} +  r^{-1}\|\Ga_g\|_{\hk_k(S)}.
 \eeaa
 In view of the estimate  proved in Step 3 and the control of $\Ga_b$, we infer
 \beaa
  \|   \dds_2\dds_1 \Jp \|_{\hk_{k_*-1}(S)} &\les&  \ep r^{-2} u^{ -\frac{\dec}{2}}.
 \eeaa
 Using Sobolev, we deduce, for $k\le k_*-3$,
\beaa
\left|\dkb^k\dds_2\dds_1 \Jp\right|    \les  \ep r^{-3} u^{-\frac{\dec}{2}}.
\eeaa
Together with the above transport equation along $\nu$ and the control of $\Ga_b$, we infer, for  any  $k\le k_*-3$, 
\beaa
\left|\dk_*^k\dds_2\dds_1 \Jp\right|    \les  \ep r^{-3} u^{-\frac{\dec}{2}}.
\eeaa
  as stated. This concludes the proof of  Lemma \ref{Le:Si*-ell=1modes}.
\end{proof}

We state below two corollaries of Lemma \ref{Le:Si*-ell=1modes}.
 \begin{corollary} 
\label{prop:2D-Hodge4}
On a fixed sphere  $S=S(u) \subset \Si_* $,   we have for any pair of scalars  $(h, \dual h)$  and  all $k\le k_{large}-1$
 \beaa
  \|(h, \dual h)\|_{\hk_{k+2}  (S)} &\les&  r^2\|\dds_2\,\dds_1(h, \dual h)\|_{\hk_k(S)}+r\big| (h)_{\ell=1}\big|+r\big| (\dual h)_{\ell=1}\big|+r|\ov{h}|+r|\ov{\dual h}|.
  \eeaa
  \end{corollary}
  
  \begin{proof}
  In view of Corollary \ref{prop:2D-Hodge3}, we have
\beaa
  \|(h-\ov{h}, \dual h-\ov{\dual h})\|_{\hk_{k+2}  (S)} \les  r^2\|\dds_2\,\dds_1(h, \dual h)\|_{\hk_k(S)}+r^3\big| (\Delta h)_{\ell=1}\big|+r^3\big| (\Delta\dual h)_{\ell=1}\big|.
  \eeaa  
  Since
  \beaa
  (\Delta h)_{\ell=1} &=& (\Delta (h-\ov{h}))_{\ell=1}=\frac{2}{r^2}(h-\ov{h})_{\ell=1}+\left(\left(\Delta+\frac{2}{r^2}\right)(h-\ov{h})\right)_{\ell=1},
  \eeaa
  and similarly for $\dual h$, we infer
\beaa
  &&\|(h-\ov{h}, \dual h-\ov{\dual h})\|_{\hk_{k+2}  (S)}\\ 
  &\les&  r^2\|\dds_2\,\dds_1(h, \dual h)\|_{\hk_k(S)}+r\big| (h-\ov{h})_{\ell=1}\big|+r\big| (\dual h-\ov{\dual h})_{\ell=1}\big|\\
  &&+r^3\left|\left(\left(\Delta+\frac{2}{r^2}\right)(h-\ov{h})\right)_{\ell=1}\right|+r^3\left|\left(\left(\Delta+\frac{2}{r^2}\right)(\dual h-\ov{\dual h})\right)_{\ell=1}\right|.
  \eeaa
  Now, using integration by parts and Lemma \ref{Le:Si*-ell=1modes}, we have, for $p=0,+,-$,
  \beaa
  \left|\int_S\left(\Delta+\frac{2}{r^2}\right)(h-\ov{h})\Jp\right| &=& \left|\int_S(h-\ov{h})\left(\Delta+\frac{2}{r^2}\right)\Jp\right|\\
  &\les& \|h-\ov{h}\|_{L^2(S)}r\left\|\left(\Delta+\frac{2}{r^2}\right)\Jp\right\|_{L^\infty(S)}\\
  &\les& \frac{\ep}{r^2}\|h-\ov{h}\|_{L^2(S)}
  \eeaa
  and hence
 \beaa
 \left|\left(\left(\Delta+\frac{2}{r^2}\right)(h-\ov{h})\right)_{\ell=1}\right| &\les&  \frac{\ep}{r^4}\|h-\ov{h}\|_{L^2(S)}.
 \eeaa   
 We infer
  \beaa
  \|(h-\ov{h}, \dual h-\ov{\dual h})\|_{\hk_{k+2}  (S)} &\les&  r^2\|\dds_2\,\dds_1(h, \dual h)\|_{\hk_k(S)}+r\big| (h-\ov{h})_{\ell=1}\big|+r\big| (\dual h-\ov{\dual h})_{\ell=1}\big|\\
  &&+\frac{\ep}{r}\|h-\ov{h}\|_{L^2(S)}+\frac{\ep}{r}\|\dual h-\ov{\dual h}\|_{L^2(S)}.
  \eeaa
  For $\ep>0$ small enough, we may absorb the last terms on the RHS and deduce
  \beaa
  \|(h-\ov{h}, \dual h-\ov{\dual h})\|_{\hk_{k+2}  (S)} &\les&  r^2\|\dds_2\,\dds_1(h, \dual h)\|_{\hk_k(S)}+r\big| (h-\ov{h})_{\ell=1}\big|+r\big| (\dual h-\ov{\dual h})_{\ell=1}\big|.
  \eeaa 
    In particular, we infer
  \beaa
  \|(h, \dual h)\|_{\hk_{k+2}  (S)} &\les&  r^2\|\dds_2\,\dds_1(h, \dual h)\|_{\hk_k(S)}+r\big| (h)_{\ell=1}\big|+r\big| (\dual h)_{\ell=1}\big|+r|\ov{h}|+r|\ov{\dual h}|
  \eeaa     
  as desired.
  \end{proof}

\begin{corollary}\lab{cor:Cb0CbpM0MareGagandrm1Gag}
The scalar functions $\underline{C}_0$, $\underline{C}_p$, $M_0$, $M_p$ on $\Si_*$ verify
\bea
\underline{C}_0\in \Ga_g, \qquad \underline{C}_p\in \Ga_g, \qquad M_0\in r^{-1}\Ga_g, \qquad M_p\in r^{-1}\Ga_g.
\eea
\end{corollary}

\begin{proof}
Recall that we have on $\Si_*$ in view of our GCM conditions
\beaa
\kabc &=& \underline{C}_0+\sum_p\underline{C}_p\Jp
\eeaa
and hence
\beaa
\underline{C}_0+\sum_p\underline{C}_p\Jp &=& \Ga_g.
\eeaa
Integrating this identity on $S$, as well as multiplying it by $J^{(q)}$ and integrating is also on $S$, we obtain, after dividing by $|S|$,  
\beaa
\underline{C}_0 &=& \Ga_g +O(r^{-2})\sum_p\left(\int_S\Jp\right)\underline{C}_p,\\
\underline{C}_q &=& \Ga_g +O(r^{-2})\left(\int_SJ^{(q)}\right)\underline{C}_0+O(r^{-2})\sum_{p}\left(\int_S\Jp J^{(q)} -\frac{4\pi}{3}r^2\de_{pq}\right)\underline{C}_p, \,\, q=0,+,-.
\eeaa
Together with Lemma \ref{Le:Si*-ell=1modes}, we deduce 
\beaa
\underline{C}_0 &=& \Ga_g +O(\ep)\sum_p\underline{C}_p,\\
\underline{C}_q &=& \Ga_g +O(\ep)\underline{C}_0+O(\ep)\sum_p\underline{C}_p, \quad q=0,+,-,
\eeaa
which yields the desired result for $\underline{C}_0$ and $\underline{C}_p$.

Next, recall that we have on $\Si_*$ in view of our GCM conditions
\beaa
\muc &=& M_0+\sum_pM_p\Jp
\eeaa
and hence
\beaa
M_0+\sum_pM_p\Jp &=& r^{-1}\Ga_g.
\eeaa
The proof for $M_0$ and $M_p$ follows then the same line as the one for $\underline{C}_0$ and $\underline{C}_p$. This concludes the proof of the corollary.
\end{proof}

Finally, we state below a corollary of Corollary \ref{cor:Cb0CbpM0MareGagandrm1Gag} and Corollary \ref{Corr:nu*ofGCM}.
 \begin{corollary}
\lab{Corr:nu*ofGCM:plugGCM}
 The following identities hold true   on $\Si_*$.
 \bea
   \lab{eq:D^5eta-xib-prop}
   \bsplit
  2\dds_2\dds_1 \ddd_1\ddd_2\dds_2\eta &=      -\frac{4}{r}\dds_2\dds_1\div\bb       + r^{-5} \dkb^{\le 5} \Ga_g+ r^{-4} \dkb^{\le 4} (\Ga_b\c \Ga_b)\\
  &+\sum_{p=0, +,-}r^{-1}\nab_\nu\Big(r^{-1}\Ga_g\dds_2\dds_1(J^{(p)})\Big),\\
 2\dds_2\dds_1 \ddd_1\ddd_2\dds_2\xib  &=     -\frac{4}{r}\dds_2\dds_1 \div\bb + r^{-5} \dkb^{\le 5} \Ga_g+ r^{-4} \dkb^{\le 4} (\Ga_b\c \Ga_b) \\
 & +\sum_{p=0, +,-}\nab_\nu\left(r^{-2}\Ga_g\dkb^{\leq 2}\dds_2\dds_1(J^{(p)})\right)   +\sum_{p=0, +,-}r^{-1}\nab_\nu\Big(r^{-1}\Ga_g\dds_2\dds_1(J^{(p)})\Big).
  \end{split}
  \eea  
 \end{corollary}

\begin{proof}
In view of the definition of the linearized quantities $\kac$, $\kabc$ and $\muc$, we may rewrite  the GCM conditions \eqref{eq:Si_*-GCM1} as follows
  \beaa
 \kac =0,\qquad \kabc =\underline{C}_0+\sum_{p=0, +,-}\underline{C}_pJ^{(p)},\qquad \muc =M_0+\sum_{p=0, +,-}M_p J^{(p)}.
 \eeaa
Since the scalar functions  $\underline{C}_0$, $\underline{C}_p$, $M_0$, and $M_p$ are constant  on the leaves of the $r$-foliation of $\Si_*$, and since $\nu$ is tangent to $\Si_*$, we infer
\beaa
\dds_2\dds_1 \ddd_1\nab_\nu\nab\kac &=& 0,\\
\nab_\nu\left(\dds_2\ddd_2+\frac{2}{r^2}\right)\dds_2\dds_1\kabc &=& \nab_\nu\left(\dds_2\ddd_2+\frac{2}{r^2}\right)\dds_2\dds_1\left(\underline{C}_0+\sum_{p=0, +,-}\underline{C}_pJ^{(p)}\right)\\
&=& \sum_{p=0, +,-}\nab_\nu\left(\underline{C}_p\left(\dds_2\ddd_2+\frac{2}{r^2}\right)\dds_2\dds_1(J^{(p)})\right),\\
\nab_\nu\dds_2\dds_1\muc &=& \nab_\nu\dds_2\dds_1\left(M_0+\sum_{p=0, +,-}M_p J^{(p)}\right)\\
&=& \sum_{p=0, +,-}\nab_\nu\Big(M_p\dds_2\dds_1(J^{(p)})\Big).
\eeaa
Since we have, in view of Corollary \ref{cor:Cb0CbpM0MareGagandrm1Gag}, $\underline{C}_p\in \Ga_g$ and $M_p\in r^{-1}\Ga_g$, we infer
\beaa
\dds_2\dds_1 \ddd_1\nab_\nu\nab\kac &=& 0,\\
\nab_\nu\left(\dds_2\ddd_2+\frac{2}{r^2}\right)\dds_2\dds_1\kabc &=& \sum_{p=0, +,-}\nab_\nu\left(r^{-2}\Ga_g\dkb^{\leq 2}\dds_2\dds_1(J^{(p)})\right),\\
\nab_\nu\dds_2\dds_1\muc &=& \sum_{p=0, +,-}\nab_\nu\Big(r^{-1}\Ga_g\dds_2\dds_1(J^{(p)})\Big).
\eeaa
Plugging these identities in \eqref{eq:D^5eta-prop} and \eqref{eq:D^5xib-prop}, we obtain
   \beaa
  \bsplit
  2\dds_2\dds_1 \ddd_1\ddd_2\dds_2\eta &=      -\frac{4}{r}\dds_2\dds_1\div\bb       + r^{-5} \dkb^{\le 5} \Ga_g+ r^{-4} \dkb^{\le 4} (\Ga_b\c \Ga_b)\\
  &+\sum_{p=0, +,-}r^{-1}\nab_\nu\Big(r^{-1}\Ga_g\dds_2\dds_1(J^{(p)})\Big),\\
 2\dds_2\dds_1 \ddd_1\ddd_2\dds_2\xib  &=     -\frac{4}{r}\dds_2\dds_1 \div\bb + r^{-5} \dkb^{\le 5} \Ga_g+ r^{-4} \dkb^{\le 4} (\Ga_b\c \Ga_b) \\
 & +\sum_{p=0, +,-}\nab_\nu\left(r^{-2}\Ga_g\dkb^{\leq 2}\dds_2\dds_1(J^{(p)})\right)   +\sum_{p=0, +,-}r^{-1}\nab_\nu\Big(r^{-1}\Ga_g\dds_2\dds_1(J^{(p)})\Big),
\end{split}
\eeaa
as stated.
\end{proof}

%%%%%%%%%%%%%%%%%%%%%%%%%%%%%%%%%%%%%%

\subsection{Propagation equations along $\Si_*$ for some $\ell=1$ modes}

%%%%%%%%%%%%%%%%%%%%%%%%%%%%%%%%%%%%%%

We have the following corollary of Lemma \ref{Lemma:transport.alongSi_*1}.
\begin{corollary}\lab{corofLemma:transport.alongSi_*1}
We have along $\Si_*$, for $p=0,+,-$,
\bea
\nn&&\nu\left(\int_S\left( \lap\kabc+\frac{2\Up}{r}\div \ze\right)\Jp\right)\\
\nn&=& O(r^{-3})\int_S\kabc\Jp   +O(r^{-2})\int_S\div\ze\Jp +O(r^{-1})\int_S\div\bb\Jp     +O(r^{-2})\int_S\rhoc\Jp \\
 &&+O(r^{-1})\int_S\div\b\Jp +r\left|\left(\Delta+\frac{2}{r^2}\right)\Jp\right|\dkb^{\leq 1}\Ga_b+\dkb^{\le 2 }(\Ga_b\c \Ga_b),
\eea
\bea
\nu\left(\int_S\div\b\Jp\right) &=& O(r^{-1})\int_S\div\b\Jp+O(r^{-2})\int_S\rhoc\Jp\\
\nn&&+r\left(\left|\left(\Delta+\frac{2}{r^2}\right)\Jp\right|+\left|\dds_2\dds_1\Jp\right|\right)\Ga_g+\dkb^{\leq 1}(\Ga_b\c\Ga_g),
\eea
\bea
\nn\nu\left(\int_S\curl\b\Jp\right) &=&  \frac{4}{r}(1+ O(r^{-1}))\int_S\curl\b\Jp+\frac{2}{r^2}(1+ O(r^{-1}))\int_S\rhod\Jp\\
\nn&&+r\left(\left|\left(\Delta+\frac{2}{r^2}\right)\Jp\right|+\left|\dds_2\dds_1\Jp\right|\right)\Ga_g\\
&&+\dkb^{\leq 1}(\Ga_b\c\Ga_g),
\eea
and
\bea
\nn\nu\left(\int_S\left(\rhoc - \frac{1}{2}\chih\c\chibh \right)\Jp\right) &=&  -\int_S\div\bb\Jp -(1+O(r^{-1}))\int_S\div\b\Jp \\
\nn&&  +O(r^{-1})\int_S\left(\rhoc - \frac{1}{2}\chih\c\chibh \right)\Jp +O(r^{-3})\int_S\kabc\Jp\\
 \nn  && +O(r^{-2})\int_S\div\ze\Jp+r\left|\left(\Delta+\frac{2}{r^2}\right)\Jp\right|\Ga_b \\
   &&+r\dkb^{\le  1}( \Ga_b \c \Ga_b),
\eea
where by $\dds_1\Jp$, we mean $\dds_1(\Jp,0)$ or $\dds_1(0,\Jp)$.
\end{corollary}

\begin{proof}
See section \ref{sec:corofLemma:transport.alongSi_*1}.
\end{proof}

%%%%%%%%%%%%%%%%%%%%%%%%%%%%%%%

\section{Control of the flux of some quantities on $\Sigma_*$}
\lab{sec:controlofthefluxofsomequantitiesonSigmastar}

%%%%%%%%%%%%%%%%%%%%%%%%%%%%%%%

The goal of this section is to establish the following.
\begin{proposition}
\lab{Prop.Flux-bb-vthb-eta-xib}
The following estimate holds true for all  $k\le k_{*}- 7$
\bea
\lab{Estimate:Flux-bb-vthb-eta-xib}
\int_{\Si_*}u^{2+2\dec}    \big|\dk_*^ k\Ga_b\big|^2  &\les& \ep_0^2.
\eea
We also have,  for $k\le k_{*} - 10$, 
\bea
\lab{Estimate:Flux-bb-vthb-eta-xib-weak}
\|\Ga_b \|_{\infty, k} &\les \ep_0 r^{-1} u^{-1-\dec}.
\eea
\end{proposition}

\begin{proof}
Note that \eqref{Estimate:Flux-bb-vthb-eta-xib-weak} follows immediately from \eqref{Estimate:Flux-bb-vthb-eta-xib} using the trace theorem and Sobolev. We thus concentrate our attention on deriving \eqref{Estimate:Flux-bb-vthb-eta-xib}.

{\bf Step 1.} We  first   prove the corresponding   estimates for   $\bb$  away from  its   $\ell=1$ mode.  More precisely we prove the following. 
\begin{lemma}
\lab{Le:Flux-dds_2bb}
The following estimates holds true  for all $k\le k_{*}-3$
\bea
\lab{Estim:Flux-dds_2bb}
\int_{\Si_*}  r^4 u^{2+2\dec}   \big| \dds_2 \big( \dk_*^ {k}   \bb  \big)\big|^2   &\les \ep_0^2.
\eea
\end{lemma}

\begin{proof}
Recall the identity \eqref{eq:Le-Teuk-Star1-M4}
\beaa
\Re(\nab_3(r\qf)) = r^5\dds_2\dds_1\big(\div\bb, -\curl\bb\big) +O(r)\dkb^{\leq 3}\aa +O(r^{-2})+r\dkb^{\leq 3}\Ga_g+r^{3}\dkb^{\le 3}(\Ga_b \c  \Ga_g).
\eeaa
Together with the control of $\Ga_b$ and $\Ga_g$ provided by {\bf Ref 1} and Lemma \ref{lemma:interpolation}, we infer, for $k\leq k_*-3$,
\beaa
r^5|\dk_*^k\dds_2\dds_1\big(\div\bb, -\curl\bb\big)| &\les& |\dk_*^{\leq k_*-3}\nab_3(r\qf)|+r|\dk_*^{\leq k_*}\aa|+\frac{1}{r^2}+\frac{\ep}{ru^{\frac{1}{2}+\frac{\dec}{2}}}+\frac{\ep^2}{u^{\frac{3}{2}+\frac{3\dec}{2}}}\\
&\les& |\dk_*^{\leq k_*-3}\nab_3(r\qf)|+r|\dk_*^{\leq k_*}\aa|+\frac{\ep_0}{u^{\frac{3}{2}+\frac{3\dec}{2}}}
\eeaa
where we used in the last inequality the dominance condition \eqref{eq:dominantconditiononronSigmastarchap5} on $r$ on $\Si_*$. Also, using the fact that $e_3(r)=-\Up+O(r)$ and {\bf Ref 2} for $\qf$, we have
\beaa
 |\dk_*^{\leq k_*-3}\nab_3(r\qf)| &\les&  r|\dk_*^{\leq k_*-3}\nab_3\qf|+ |\dk_*^{\leq k_*-3}\qf| \les r|\dk_*^{\leq k_*-3}\nab_3\qf|+\frac{\ep}{ru^{\frac{1}{2}+\dec}}\\
  &\les& r|\dk_*^{\leq k_*-3}\nab_3\qf|+\frac{\ep_0}{u^{\frac{3}{2}+\frac{3\dec}{2}}}
\eeaa
where we used again the dominance condition \eqref{eq:dominantconditiononronSigmastarchap5} on $r$ on $\Si_*$. We infer
\beaa
r^5|\dk_*^k\dds_2\dds_1\big(\div\bb, -\curl\bb\big)| &\les& r|\dk_*^{\leq k_*-3}\nab_3\qf|+r|\dk_*^{\leq k_*}\aa|+\frac{\ep_0}{u^{\frac{3}{2}+\frac{3\dec}{2}}}
\eeaa
Dividing by $r$, squaring, and integrating on $\Si_*$, we deduce
\beaa
&&\int_{\Si_*}r^8u^{2+2\dec}|\dk_*^k\dds_2\dds_1\big(\div\bb, -\curl\bb\big)|^2\\ 
&\les& \int_{\Si_*}u^{2+2\dec}|\dk_*^{\leq k_*-3}\nab_3\qf|^2+\int_{\Si_*}u^{2+2\dec}|\dk_*^{\leq k_*}\aa|^2+\int_{\Si_*}\frac{\ep_0^2}{r^2u^{1+\dec}}
\eeaa
and hence, since $\dec>0$ and using the control of $\aa$ and $\nab_3\qf$ provided by  {\bf Ref 2}, we deduce,  for $k\leq k_*-3$,
\beaa
\int_{\Si_*}r^8u^{2+2\dec}|\dk_*^k\dds_2\dds_1\big(\div\bb, -\curl\bb\big)|^2 &\les& \ep_0^2.
\eeaa
Taking into account the commutator Lemma \ref{Lemma:Commutation-Si_*},  this yields, for $k\le k_{*}-3$,
\beaa
\int_{\Si_*}r^8u^{2+2\dec}|\dds_1(\div, -\curl)\dds_2\,\dk_*^k\bb|^2 &\les& \ep_0^2.
\eeaa
Using the Hodge estimates of Lemma \ref{prop:2D-Hodge1} for $\dds_1$ and $\ddd_1$, we infer, for $k\le k_{*}-3$,
\beaa
\int_{\Si_*}r^4u^{2+2\dec}|\dds_2\,\dk_*^k\bb|^2 &\les& \ep_0^2
\eeaa
as stated.
This concludes the proof of Lemma  \ref{Le:Flux-dds_2bb}.
\end{proof}

{\bf Step 2.}    We next derive  the following  non sharp, preliminary,   estimate for the $\ell=1$ mode of $\ddd_1\bb$
\bea
\lab{eq:estimate(bb*)(ell=1)} 
\Big|\big(\ddd_1\nab_\nu^k\bb\big)_{\ell=1} \Big| &\les& \frac{\ep_0}{r^3u^{\frac{3}{2}+\frac{3\dec}{2}}}, \qquad k\le k_*-3.
\eea
To this end, we use the following consequence of the  Codazzi equation for $\chibh$
\beaa
\ddd_2\chibh &=&  \bb +r^{-1}\dkb^{\leq 1}\Ga_g+\Ga_b\c \Ga_g.
\eeaa
Differentiating w.r.t. $\nu^k\ddd_1$ for $k\le k_*-3$, we infer
\beaa
\nu^k\ddd_1\ddd_2\chibh &=&  \nu^k\ddd_1\bb +r^{-2}\dk_*^{\leq k_*-1}\Ga_g+r^{-1}\dk_*^{\leq k_*-2}(\Ga_b\c \Ga_g).
\eeaa
Taking into account the commutator Lemma \ref{Lemma:Commutation-Si_*},  this yields, for $k\le k_{*}-3$,
\beaa
\ddd_1\ddd_2\nab_\nu^k\chibh &=&  \ddd_1\nab_\nu^k\bb +r^{-2}\dk_*^{\leq k_*-1}\Ga_g+r^{-1}\dk_*^{\leq k_*-2}(\Ga_b\c \Ga_g).
\eeaa
Projecting on the $\ell=1$ modes, this yields for $k\leq k_*-3$
\beaa
(\ddd_1\ddd_2\nab_\nu^k\chibh)_{\ell=1} &=&  (\ddd_1\nab_\nu^k\bb)_{\ell=1} +r^{-2}\dk_*^{\leq k_*-1}\Ga_g+r^{-1}\dk_*^{\leq k_*-2}(\Ga_b\c \Ga_g).
\eeaa
Together with the control of $\Ga_b$ and $\Ga_g$ provided by {\bf Ref 1} and Lemma \ref{lemma:interpolation}, we infer, for $k\leq k_*-3$,
\beaa
|(\ddd_1\nab_\nu^k\bb)_{\ell=1}| &\les& |(\ddd_1\ddd_2\nab_\nu^k\chibh)_{\ell=1}|+\frac{\ep}{r^4u^{\frac{1}{2}+\frac{\dec}{2}}}
\eeaa
which together with  the dominance condition \eqref{eq:dominantconditiononronSigmastarchap5} on $r$ on $\Si_*$ implies, for $k\leq k_*-3$,
\beaa
|(\ddd_1\nab_\nu^k\bb)_{\ell=1}| &\les& |(\ddd_1\ddd_2\nab_\nu^k\chibh)_{\ell=1}|+\frac{\ep_0}{r^3u^{\frac{3}{2}+\frac{3\dec}{2}}}.
\eeaa

Next, we estimate $(\ddd_1\ddd_2\nab_\nu^k\chibh)_{\ell=1}$. We have
\beaa
(\ddd_1\ddd_2\nab_\nu^k\chibh)_{\ell=1,p}=\frac{1}{|S|}\int_S\ddd_1\ddd_2\nab_\nu^k\chibh\Jp =\frac{1}{|S|}\int_S\nab_\nu^k\chibh\c\dds_2\dds_1\Jp
\eeaa
and hence
\beaa
|(\ddd_1\ddd_2\nab_\nu^k\chibh)_{\ell=1}| &\les& |\dds_2\dds_1\Jp||\dk_*^k\Ga_b|.
\eeaa
Together with the control of $\Ga_b$  provided by  Lemma \ref{lemma:interpolation}, and the control of $\dds_2\dds_1\Jp$ provided by Lemma \ref{Le:Si*-ell=1modes}
we infer, for $k\leq k_*-3$,
\beaa
|(\ddd_1\ddd_2\nab_\nu^k\chibh)_{\ell=1}| &\les& \frac{\ep}{r^{3}u^{\frac{\dec}{2}}}\,\frac{\ep}{ru^{1+\frac{\dec}{2}}}\les \frac{\ep_0}{r^4u^{1+\dec}}\les \frac{\ep_0}{r^3u^{2+2\dec}}.
\eeaa
Plugging in the above, we deduce, for $k\leq k_*-3$,
\beaa
|(\ddd_1\nab_\nu^k\bb)_{\ell=1}| &\les& |(\ddd_1\ddd_2\nab_\nu^k\chibh)_{\ell=1}|+\frac{\ep_0}{r^3u^{\frac{3}{2}+\frac{3\dec}{2}}}\les \frac{\ep_0}{r^3u^{\frac{3}{2}+\frac{3\dec}{2}}}
\eeaa
which concludes the proof of \eqref{eq:estimate(bb*)(ell=1)}.

{\bf Step 3.} Using  the estimate \eqref{eq:estimate(bb*)(ell=1)}  of Step 2 and Lemma  \ref{prop:2D-Hodge2}, we have, for $k\le k_*-3$,
\beaa
  \|\dk_*^k  \bb\|_{L^2 (S)}&\les&  r\|\dds_2\, \dk_*^k  \bb \|_{L^2(S)}+  r^2  \big| (\ddd_1\nab_\nu^{\leq k}\bb)_{\ell=1}\big|\\
  &\les&   r  \|\dds_2\, \dk_*^k  \bb \|_{L^2(S)}+\frac{\ep_0}{ru^{\frac{3}{2}+\frac{3\dec}{2}}}.
  \eeaa
  Multiplying by $ru^{1+\dec}$, and squaring, we obtain
\beaa
  r^2u^{2+2\dec}\|\dk_*^k  \bb\|_{L^2 (S)}^2  &\les&   r^4u^{2+2\dec}  \|\dds_2\, \dk_*^k  \bb \|_{L^2(S)}^2+\frac{\ep_0^2}{u^{1+\dec}}.
  \eeaa  
  Integrating in $u$, we infer
\beaa
\int_{\Si_*}r^2u^{2+2\dec}|\dk_*^k  \bb|^2 &\les& \ep_0^2+\int_{\Si_*} r^4u^{2+2\dec}|\dds_2\, \dk_*^k  \bb|^2.
\eeaa
Together with the estimate \eqref{Estim:Flux-dds_2bb}  of Step 1, we infer
    \bea
  \lab{eq:flux-estimate-bb}
  \int_{\Si_*}r^2u^{2+2\dec}|\dk_*^k  \bb|^2 &\les& \ep_0^2, \qquad k\le k_*-3,
\eea
which is the desired estimate for $\bb$.

{\bf Step 4.} We   now  prove the desired  estimate for  $\chibh$,  i.e.
\bea\lab{eq:thefluxesimateforthetabaronSi*tobeprovedbelow}
\int_{\Si_*}  u^{2+2\dec}   \big|\dk_*^{k} \chibh|^2  \les \ep_0^2, \qquad k\le k_{*}-3.
\eea

\begin{proof}[Proof of \eqref{eq:thefluxesimateforthetabaronSi*tobeprovedbelow}]
One starts with the 
 the following consequence of the  Codazzi equation for $\chibh$
\beaa
\ddd_2\chibh &=&  \bb +r^{-1}\dkb^{\leq 1}\Ga_g+\Ga_b\c \Ga_g.
\eeaa
Differentiating w.r.t. $\dds_2$, we infer 
\beaa
  \dds_2 \ddd_2 \chibh &=&  \dds_2 \bb  +  r^{-2}\dkb^{\le 2} \Ga_g +r^{-1}\dkb^{\leq 1}\left( \Ga_b\c \Ga_g\right).
\eeaa
Taking into account the commutator Lemma \ref{Lemma:Commutation-Si_*},  this yields, for $k\le k_{*}-3$,
\beaa
  \dds_2 \ddd_2 \dk_*^k\chibh &=&  \dds_2(\dk_*^k\bb) +  r^{-2}\dkb^{\le k_*-1} \Ga_g +r^{-1}\dkb^{\leq k_*-2}\left( \Ga_b\c \Ga_g\right).
\eeaa
Together with the control of $\Ga_b$ and $\Ga_g$ provided by {\bf Ref 1} and Lemma \ref{lemma:interpolation}, we infer, for $k\leq k_*-3$,
\beaa
  |\dds_2 \ddd_2 \dk_*^k\chibh| &\les&  |\dds_2(\dk_*^k\bb)|  +  \frac{\ep}{r^4u^{\frac{1}{2}+\frac{\dec}{2}}}\\
  &\les&  |\dds_2(\dk_*^k\bb)|  +  \frac{\ep_0}{r^3u^{\frac{3}{2}+\frac{3\dec}{2}}},
\eeaa
where we used again the dominance condition \eqref{eq:dominantconditiononronSigmastarchap5} on $r$ on $\Si_*$. Squaring and multiplying by $r^4u^{2+2\dec}$, we infer
\beaa
  r^4u^{2+2\dec}|\dds_2 \ddd_2 \dk_*^k\chibh|^2   &\les&  r^4u^{2+2\dec}|\dds_2(\dk_*^k\bb)|^2  +  \frac{\ep_0^2}{r^2u^{1+\dec}},
\eeaa
which yields after integration on $\Si_*$ implies, for $k\leq k_*-3$,
\beaa
\int_{\Si_*}  r^4u^{2+2\dec}|\dds_2 \ddd_2 \dk_*^k\chibh|^2   &\les&  \int_{\Si_*}r^4u^{2+2\dec}|\dds_2(\dk_*^k\bb)|^2  +\int_{\Si_*}  \frac{\ep_0^2}{r^2u^{1+\dec}}\\
&\les& \ep_0^2+\int_{\Si_*}r^4u^{2+2\dec}|\dds_2(\dk_*^k\bb)|^2.
\eeaa
Hence, together with \eqref{Estim:Flux-dds_2bb}, we infer
\beaa
\int_{\Si_*}  r^4u^{2+2\dec}|\dds_2 \ddd_2 \dk_*^k\chibh|^2   &\les& \ep_0^2.
\eeaa
In view of the coercivity of $ \dds_2 \ddd_2 $, see \eqref{eq:dcalident},  we  deduce
\bea
\int_{\Si_*}  r^4u^{2+2\dec}|\dk_*^k\chibh|^2   &\les& \ep_0^2, \qquad k\le  k_{*}-3, 
\eea
which is the stated estimate \eqref{eq:thefluxesimateforthetabaronSi*tobeprovedbelow}.
 \end{proof}

{\bf Step 5.}  Next, we establish the estimates for $\eta$ and $\xib$  in Proposition \ref{Prop.Flux-bb-vthb-eta-xib}. To this end, we first   estimate $\dds_2\eta$ and $\dds_2\xib$.

\begin{lemma}\lab{lemma:firstfluxestimatesonSi*-etaxib}
We have for $k\le k_*-6$
\bea
\lab{fluxestimates-foreta-xib}
\int_{\Si_*}r^2u^{2+2\dec}\Big(|\dds_2\,\dk_*^k\eta|^2+|\dds_2\,\dk_*^k\xib|^2\Big) &\les& \ep_0^2.
\eea
\end{lemma}

\begin{proof}
Recall from  Corollary \ref{Corr:nu*ofGCM:plugGCM} that the following identities hold true   on $\Si_*$
 \beaa
   \bsplit
  2\dds_2\dds_1 \ddd_1\ddd_2\dds_2\eta &=      -\frac{4}{r}\dds_2\dds_1\div\bb       + r^{-5} \dkb^{\le 5} \Ga_g+ r^{-4} \dkb^{\le 4} (\Ga_b\c \Ga_b)\\
  &+\sum_{p=0, +,-}r^{-1}\nab_\nu\Big(r^{-1}\Ga_g\dds_2\dds_1(J^{(p)})\Big),\\
 2\dds_2\dds_1 \ddd_1\ddd_2\dds_2\xib  &=     -\frac{4}{r}\dds_2\dds_1 \div\bb + r^{-5} \dkb^{\le 5} \Ga_g+ r^{-4} \dkb^{\le 4} (\Ga_b\c \Ga_b) \\
 & +\sum_{p=0, +,-}\nab_\nu\left(r^{-2}\Ga_g\dkb^{\leq 2}\dds_2\dds_1(J^{(p)})\right)   +\sum_{p=0, +,-}r^{-1}\nab_\nu\Big(r^{-1}\Ga_g\dds_2\dds_1(J^{(p)})\Big).
 \end{split}
  \eeaa  
 Taking into account the commutator Lemma \ref{Lemma:Commutation-Si_*},  this yields, for $k\le k_{*}-6$,
 \beaa
   \bsplit
  2\dds_2\dds_1 \ddd_1\ddd_2\dds_2\,\dk_*^k\eta &=      r^{-4}\dk_*^{\leq k_*-3}\bb       + r^{-5}\dk_*^{\le k_*} \Ga_g+ r^{-4}\dk_*^{\leq k_*}(\Ga_b\c \Ga_b)\\
  &+\sum_{p=0, +,-}r^{-1}\dk_*^{\leq k_*-3}\Big(r^{-1}\Ga_g\dds_2\dds_1(J^{(p)})\Big),\\
 2\dds_2\dds_1 \ddd_1\ddd_2\dds_2\,\dk_*^k\xib  &=     r^{-4}\dk_*^{\leq k_*-3}\bb       + r^{-5}\dk_*^{\le k_*} \Ga_g+ r^{-4}\dk_*^{\leq k_*}(\Ga_b\c \Ga_b)\\
  &+\sum_{p=0, +,-}r^{-1}\dk_*^{\leq k_*-3}\Big(r^{-1}\Ga_g\dds_2\dds_1(J^{(p)})\Big).
 \end{split}
  \eeaa  
Together with the control of $\Ga_b$ and $\Ga_g$ provided by {\bf Ref 1} and Lemma \ref{lemma:interpolation}, and the control of $\dds_2\dds_1(J^{(p)})$ provided by Lemma \ref{Le:Si*-ell=1modes}, we infer, for $k\leq k_*-6$,
\beaa
|\dds_2\dds_1 \ddd_1\ddd_2\dds_2\,\dk_*^k\eta|+|\dds_2\dds_1 \ddd_1\ddd_2\dds_2\,\dk_*^k\xib| &\les& r^{-4}|\dk_*^{\leq k_*-3}\bb|+\frac{\ep}{r^7u^{\frac{1}{2}+\frac{\dec}{2}}}+\frac{\ep^2}{r^6u^{2+\frac{3\dec}{2}}}\\
&\les& r^{-4}|\dk_*^{\leq k_*-3}\bb|+\frac{\ep_0}{r^6u^{\frac{3}{2}+\frac{3\dec}{2}}}
\eeaa
where we used again the dominance condition \eqref{eq:dominantconditiononronSigmastarchap5} on $r$ on $\Si_*$. Squaring and multiplying by $r^{10}u^{2+2\dec}$, we infer
\beaa
r^{10}u^{2+2\dec}\Big(|\dds_2\dds_1 \ddd_1\ddd_2\dds_2\,\dk_*^k\eta|^2+|\dds_2\dds_1 \ddd_1\ddd_2\dds_2\,\dk_*^k\xib|^2\Big) &\les& r^2u^{2+2\dec}|\dk_*^{\leq k_*-3}\bb|^2+\frac{\ep_0^2}{r^2u^{1+\dec}}
\eeaa
which yields after integration on $\Si_*$ implies, for $k\leq k_*-6$,
\beaa
&&\int_{\Si_*}r^{10}u^{2+2\dec}\Big(|\dds_2\dds_1 \ddd_1\ddd_2\dds_2\,\dk_*^k\eta|^2+|\dds_2\dds_1 \ddd_1\ddd_2\dds_2\,\dk_*^k\xib|^2\Big)\\
 &\les& \int_{\Si_*}r^2u^{2+2\dec}|\dk_*^{\leq k_*-3}\bb|^2+\int_{\Si_*}\frac{\ep_0^2}{r^2u^{1+\dec}}\\
&\les& \ep_0^2+\int_{\Si_*}r^2u^{2+2\dec}|\dk_*^{\leq k_*-3}\bb|^2.
\eeaa
Hence, together with \eqref{eq:flux-estimate-bb}, we deduce, for $k\leq k_*-6$,
\beaa
\int_{\Si_*}r^{10}u^{2+2\dec}\Big(|\dds_2\dds_1 \ddd_1\ddd_2\dds_2\,\dk_*^k\eta|^2+|\dds_2\dds_1 \ddd_1\ddd_2\dds_2\,\dk_*^k\xib|^2\Big) &\les& \ep_0^2.
\eeaa
Since, $\dds_1\ddd_1=\ddd_2\dds_2+2K$, we have
\beaa
\dds_2\dds_1 \ddd_1\ddd_2 &=& \dds_2(\ddd_2\dds_2+2K)\ddd_2=\dds_2\ddd_2\left(\dds_2\ddd_2+\frac{2}{r^2}\right)+2\dds_2\widecheck{K}\ddd_2\\
&=& \dds_2\ddd_2\left(\dds_2\ddd_2+\frac{2}{r^2}\right)+2r^{-3}\dkb\Ga_g\dkb
\eeaa
and in view of the control of $\Ga_g$ and the coercivity of $\dds_2 \ddd_2$, see \eqref{eq:dcalident},  we  obtain, for $k\leq k_*-6$,
\beaa
\int_{\Si_*}r^2u^{2+2\dec}\Big(|\dds_2\,\dk_*^k\eta|^2+|\dds_2\,\dk_*^k\xib|^2\Big) &\les& \ep_0^2.
\eeaa
which is the stated estimate \eqref{fluxestimates-foreta-xib}. This completes the proof of Lemma \ref{lemma:firstfluxestimatesonSi*-etaxib}.
\end{proof}

{\bf Step 6.} In this step, we derive the desired estimates for $\eta$ and $\xib$, i.e. we show
\bea
\lab{Estimates:Flux-eta-xib}
 \int_{\Si_*}u^{2+2\dec}\Big(|\dk_*^k\eta|^2+|\dk_*^k\xib|^2\Big) &\les& \ep_0^2, \qquad k\leq k_*-6.
\eea
To this end, we apply Lemma  \ref{prop:2D-Hodge2} which yields
\beaa
\|\dk_*^k\eta\|_{L^2 (S)}+\|\dk_*^k\xib\|_{L^2 (S)} &\les& r\|\dds_2\dk_*^k\eta\|_{L^2 (S)}+r\|\dds_2\dk_*^k\xib\|_{L^2 (S)}\\
&&+r^2\left| (\ddd_1\nab_\nu^k\eta   )_{\ell=1}\right|+ r^2\left| (\ddd_1\nab_\nu^k \xib   )_{\ell=1} \right|
\eeaa
Squaring, multiplying by $u^{2+2\dec}$ and integrating in $u$, we infer
\beaa
 \int_{\Si_*}u^{2+2\dec}\Big(|\dk_*^k\eta|^2+|\dk_*^k\xib|^2\Big) &\les& \int_{\Si_*}r^2u^{2+2\dec}\Big(|\dds_2\,\dk_*^k\eta|^2+|\dds_2\,\dk_*^k\xib|^2\Big) \\
 &&+\int_{\Si_*}r^2u^{2+2\dec}\left(\left| (\ddd_1\nab_\nu^k\eta   )_{\ell=1}\right|^2+ \left| (\ddd_1\nab_\nu^k \xib   )_{\ell=1} \right|^2\right)
\eeaa
which together with the estimate \eqref{fluxestimates-foreta-xib} of Step 5 implies, for $k\leq k_*-6$, 
\beaa
 \int_{\Si_*}u^{2+2\dec}\Big(|\dk_*^k\eta|^2+|\dk_*^k\xib|^2\Big) &\les& \ep_0^2+\int_{\Si_*}r^2u^{2+2\dec}\left(\left| (\ddd_1\nab_\nu^k\eta   )_{\ell=1}\right|^2+ \left| (\ddd_1\nab_\nu^k \xib   )_{\ell=1} \right|^2\right).
\eeaa 
Now,  assume that we have 
\bea
\lab{Estimates:BadModed-eta-xib-onSi_*}
\left| (\ddd_1\nab_\nu^k\eta   )_{\ell=1}\right|+ \left| (\ddd_1\nab_\nu^k \xib   )_{\ell=1} \right|&\les&  \frac{\ep_0}{r^2u^{\frac{3}{2}+\frac{3\dec}{2}}}, \quad  k\le k_*.
\eea
Then, we infer, for $k\leq k_*-6$,
\beaa
 \int_{\Si_*}u^{2+2\dec}\Big(|\dk_*^k\eta|^2+|\dk_*^k\xib|^2\Big) &\les& \ep_0^2+\int_{\Si_*}\frac{\ep_0^2}{r^2u^{1+\dec}}\les\ep_0^2
\eeaa 
which is the desired estimate \eqref{Estimates:Flux-eta-xib}. 

In view of the above, to complete the proof of \eqref{Estimates:Flux-eta-xib}, it suffices to prove \eqref{Estimates:BadModed-eta-xib-onSi_*}. In view of the commutator Lemma \ref{Lemma:Commutation-Si_*},  we have
\beaa
\ddd_1\nab_\nu^k\eta &=& \nu^k\ddd_1\eta+r^{-2}\dkb_*^{\leq k}\Ga_b+\dkb_*^{\leq k}(\Ga_b\c\Ga_b),\\
\ddd_1\nab_\nu^k\xib &=& \nu^k\ddd_1\xib+r^{-2}\dkb_*^{\leq k}\Ga_b+\dkb_*^{\leq k}(\Ga_b\c\Ga_b).
\eeaa
Using $\ddd_1=(\div, \curl)$ and following consequence of the null structure equations
\beaa
\curl\eta=r^{-1}\Ga_g+\Ga_b\c\Ga_g, \qquad \curl\xib=\Ga_b\c\Ga_b,
\eeaa
we infer
\beaa
\ddd_1\nab_\nu^k\eta &=& \nu^k\div\eta+r^{-1}\dkb_*^{\leq k}\Ga_g+\dkb_*^{\leq k}(\Ga_b\c\Ga_b),\\
\ddd_1\nab_\nu^k\xib &=& \nu^k\div\xib+r^{-2}\dkb_*^{\leq k}\Ga_b+\dkb_*^{\leq k}(\Ga_b\c\Ga_b).
\eeaa
which together with the control of $\Ga_b$  provided by {\bf Ref 1} and Lemma \ref{lemma:interpolation} implies, for $k\leq k_*$,
\beaa
&&|(\ddd_1\nab_\nu^k\eta)_{\ell=1}|+|(\ddd_1\nab_\nu^k\xib)_{\ell=1}|\\
 &\les& |(\nu^k\div\eta)_{\ell=1}|+|(\nu^k\div\xib)_{\ell=1}|+\frac{\ep}{r^3u^{\frac{1}{2}+\frac{\dec}{2}}}+\frac{\ep^2}{r^2u^{2+\frac{3\dec}{2}}}\\
 &\les& |(\nu^k\div\eta)_{\ell=1}|+|(\nu^k\div\xib)_{\ell=1}|+\frac{\ep_0}{r^2u^{\frac{3}{2}+\frac{3\dec}{2}}}
\eeaa
where we used the dominance condition \eqref{eq:dominantconditiononronSigmastarchap5} on $r$ on $\Si_*$. Also, since $\nu(\Jp)=0$ and $\nu(r)=  - 2 + r \Ga_b$, we have in view of Corollary \ref{Corr:nuSof integrals}
\beaa
\nu^k\left(\frac{1}{|S|}\int_S\div\eta\Jp\right) &=& \frac{1}{|S|}\int_S\nu^k(\div\eta)\Jp +r^{-2}\dk_*^{\leq k_*-1}\Ga_b+\dk_*^{\leq k_*}(\Ga_b\c\Ga_b)\\
\nu^k\left(\frac{1}{|S|}\int_S\div\xib\Jp\right) &=& \frac{1}{|S|}\int_S\nu^k(\div\xib)\Jp +r^{-2}\dk_*^{\leq k_*-1}\Ga_b+\dk_*^{\leq k_*}(\Ga_b\c\Ga_b)
\eeaa
Now, recall the GCM conditions
\beaa
( \div\eta   )_{\ell=1}=( \div\xib   )_{\ell=1}=0.
\eeaa
Since $\nu$ is tangent to $\Si_*$, we infer
\beaa
\nu^k\left(\frac{1}{|S|}\int_S\div\eta\Jp\right) = 0, \qquad \nu^k\left(\frac{1}{|S|}\int_S\div\xib\Jp\right) = 0, \qquad p=0,+,-,
\eeaa
and plugging in the above
\beaa
(\nu^k\div\eta)_{\ell=1} = r^{-2}\dk_*^{\leq k_*-1}\Ga_b+\dk_*^{\leq k_*}(\Ga_b\c\Ga_b), \quad (\nu^k\div\xib)_{\ell=1} = r^{-2}\dk_*^{\leq k_*-1}\Ga_b+\dk_*^{\leq k_*}(\Ga_b\c\Ga_b).
\eeaa
Together with the control of $\Ga_b$  provided by {\bf Ref 1} and Lemma \ref{lemma:interpolation}, we infer, for $k\leq k_*$,
\beaa
|(\nu^k\div\eta)_{\ell=1}|+|(\nu^k\div\xib)_{\ell=1}| &\les& \frac{\ep}{r^3u^{1+\frac{\dec}{2}}}+\frac{\ep^2}{r^2u^{2+\frac{3\dec}{2}}}\les \frac{\ep_0}{r^2u^{2+\frac{3\dec}{2}}}
\eeaa
where we used the dominance condition \eqref{eq:dominantconditiononronSigmastarchap5} on $r$ on $\Si_*$. Recalling the above estimate, for $k\leq k_*$, 
\beaa
|(\ddd_1\nab_\nu^k\eta)_{\ell=1}|+|(\ddd_1\nab_\nu^k\xib)_{\ell=1}| &\les& |(\nu^k\div\eta)_{\ell=1}|+|(\nu^k\div\xib)_{\ell=1}|+\frac{\ep_0}{r^2u^{\frac{3}{2}+\frac{3\dec}{2}}},
\eeaa
we deduce, for $k\leq k_*$,
\beaa
|(\ddd_1\nab_\nu^k\eta)_{\ell=1}|+|(\ddd_1\nab_\nu^k\xib)_{\ell=1}| &\les& \frac{\ep_0}{r^2u^{\frac{3}{2}+\frac{3\dec}{2}}}
\eeaa
which is the desired estimate \eqref{Estimates:BadModed-eta-xib-onSi_*}. This concludes the proof of \eqref{Estimates:Flux-eta-xib}.

{\bf Step 7.} In this step, we derive the desired estimates for $y=e_3 (r)$, $z= e_3(u)$ and $b_*$, i.e. we show
\bea\lab{eq:fluxesitmateforxyb}
\int_{\Si_*}r^{-2}u^{2+2\dec}\left(\left|\dk_*^k(\widecheck{y})\right|^2+\left|\dk_*^k(\widecheck{z})\right|^2+\left|\dk_*^k(\widecheck{b_*})\right|^2\right) &\les& \ep_0^2, \quad k\leq k_*-7.
\eea
To this end,  we make use of the equations, see Lemma \ref{Lemma:eqts-nabeta,xib},
\beaa
\nab \widecheck{y} =-\xib+\big(\ze-\eta) y= -\xib +\eta+\Ga_g , \qquad \nab \widecheck{z} = (\ze-\eta ) z=-2  \eta+\Ga_g.
\eeaa
In view of the commutator Lemma \ref{Lemma:Commutation-Si_*},  we infer
\beaa
\nab\dk_*^k \widecheck{y} = -\dk_*^k\xib +\dk_*^k\eta+\dk_*^{\leq k}\Ga_g+r\dk_*^{\leq k}(\Ga_b\c\Ga_b),\qquad \nab\dk_*^k \widecheck{z} = -2\dk_*^k\eta+\dk_*^{\leq k}\Ga_g+r\dk_*^{\leq k}(\Ga_b\c\Ga_b).
\eeaa
Using the control of $\Ga_b$ and $\Ga_g$ provided by {\bf Ref 1} and Lemma \ref{lemma:interpolation}, and the dominance condition \eqref{eq:dominantconditiononronSigmastarchap5} on $r$ on $\Si_*$, we obtain 
\beaa
|\nab\dk_*^k \widecheck{y}|+|\nab\dk_*^k \widecheck{z}| &\les& |\dk_*^k\xib| +|\dk_*^k\eta|+\frac{\ep}{r^2u^{\frac{1}{2}+\frac{\dec}{2}}}+\frac{\ep^2}{ru^{2+\frac{3\dec}{2}}}\\
&\les& |\dk_*^k\xib| +|\dk_*^k\eta|+\frac{\ep_0}{ru^{\frac{3}{2}+\frac{3\dec}{2}}}.
\eeaa 
Squaring, multiplying by $u^{2+2\dec}$ and integrating on $\Si_*$, we infer, for $k\leq k_*-6$,
\beaa
\int_{\Si_*}u^{2+2\dec}\Big(|\nab\dk_*^k \widecheck{y}|^2+|\nab\dk_*^k \widecheck{z}|^2\Big) &\les& \ep_0^2+\int_{\Si_*}u^{2+2\dec}\Big(|\dk_*^k\eta|^2+|\dk_*^k\xib|^2\Big).
\eeaa
Together with \eqref{Estimates:Flux-eta-xib}, this yields, for $k\leq k_*-6$,
\beaa
\int_{\Si_*}u^{2+2\dec}\Big(|\nab\dk_*^k \widecheck{y}|^2+|\nab\dk_*^k \widecheck{z}|^2\Big) &\les& \ep_0^2.
\eeaa
Since $b_*=-y-z$ in view of Lemma \ref{Lemma:eqts-nabeta,xib}, we deduce, for $k\leq k_*-6$,
\bea\lab{eq:intermediatefluxesitmateaveragefreexyb}
\int_{\Si_*}u^{2+2\dec}\Big(|\nab\dk_*^k \widecheck{y}|^2+|\nab\dk_*^k \widecheck{z}|^2+|\nab\dk_*^k \widecheck{b_*}|^2\Big) &\les& \ep_0^2.
\eea

In view of \eqref{eq:intermediatefluxesitmateaveragefreexyb}, it remains to estimate the averages $\ov{\nu^k(\widecheck{y})}$, $\ov{\nu^k(\widecheck{z})}$ and $\ov{\nu^k(\widecheck{b_*})}$. We start with $b_*$. In view of the definition of $\widecheck{b_*}$ and \eqref{eq:Si_*-GCM2}, we have
\beaa
\widecheck{b_*}\big|_{SP}=0, 
\eeaa
where $SP$ denotes the south poles of the spheres on $\Si_*$, i.e. $\th=\pi$. Since $\nu(\th)=0$ along $\Si_*$, we infer 
\beaa
\nu^k(\widecheck{b_*})\big|_{SP}=0, \qquad k\leq k_{large}.
\eeaa
In particular, decomposing $\nu^k(\widecheck{b_*})=\ov{\nu^k(\widecheck{b_*})}+(\nu^k(\widecheck{b_*})-\ov{\nu^k(\widecheck{b_*})})$, we infer
\beaa
|\ov{\nu^k(\widecheck{b_*})}| \leq \left|(\nu^k(\widecheck{b_*})-\ov{\nu^k(\widecheck{b_*})})\big|_{SP}\right| \leq \|\nu^k(\widecheck{b_*})-\ov{\nu^k(\widecheck{b_*})}\|_{L^\infty(S)}\quad\textrm{for any}\quad S\subset\Si_*.
\eeaa
Using Poincar\'e and Sobolev, we infer
\beaa
\left|\ov{\nu^k(\widecheck{b_*})}\right| &\les& \|\nab\dkb^{\leq 1}\nu^k(\widecheck{b_*})\|_{L^2(S)}\les \|\nab\dk_*^{\leq k+1}\widecheck{b_*}\|_{L^2(S)}\quad\textrm{for any}\quad S\subset\Si_*.
\eeaa
Squaring, multiplying by $r^{-2}u^{2+2\dec}$ and integrating on $\Si_*$, we infer
\beaa
\int_{\Si_*}r^{-2}u^{2+2\dec}\left|\ov{\nu^k(\widecheck{b_*})}\right|^2 &\les& \int_{\Si_*}u^{2+2\dec}\left|\nab\dk_*^{\leq k+1}\widecheck{b_*}\right|^2.
\eeaa
Together with \eqref{eq:intermediatefluxesitmateaveragefreexyb}, we deduce, for $k\leq k_*-7$,
\bea\lab{eq:intermediatefluxesitmateaverageb}
\int_{\Si_*}r^{-2}u^{2+2\dec}\left|\ov{\nu^k(\widecheck{b_*})}\right|^2 &\les& \ep_0^2.
\eea

Next, we estimate  the average $\ov{\nu^k(\widecheck{y})}$. Recall from Lemma \ref{Lemma:nuSof-integrals-again} the following identity
 \beaa
 \nu(r)  = \frac{rz}{2}\,\ov{z^{-1}(\kab + b_*\ka)}
 \eeaa
Using the transversality condition $e_4(r)=1$ on $\Si_*$, the fact that $\nu=e_3+b_*e_4$, and the GCM condition $\ka=2/r$, we infer
\beaa
y+b_*  &=& \frac{rz}{2}\,\ov{z^{-1}\left(\kab +\frac{2}{r} b_*\right)}=\frac{r(2+\widecheck{z})}{2}\,\ov{\frac{1}{2+\widecheck{z}}\left(-\frac{2\Up}{r} +\frac{2}{r} b_*\right)}+r\Ga_g\\
&=& \left(1+\frac{1}{2}\widecheck{z}\right)\,\ov{\left(1-\frac{1}{2}\widecheck{z}\right)\left(-\Up + b_*\right)}+r\Ga_g+r^2\Ga_b\c\Ga_b
 \eeaa
and hence
\beaa
\widecheck{y}  &=&  -\left(\widecheck{z} - \ov{\widecheck{z}}\right)-\left(\widecheck{b_*}-\ov{\widecheck{b_*}}\right)+r\Ga_g+r^2\Ga_b\c\Ga_b.
 \eeaa
Taking the average, we infer
\beaa
\ov{\widecheck{y}}  &=&  r\Ga_g+r^2\Ga_b\c\Ga_b.
 \eeaa
 Together with Corollary \ref{Corr:nuSof integrals}, we deduce
\beaa
\ov{\nu^k(\widecheck{y})}  &=&  r\dk_*^{\leq k}\Ga_g+r^2\dk_*^{\leq k}(\Ga_b\c\Ga_b).
 \eeaa  
 Using the control of $\Ga_b$ and $\Ga_g$ provided by {\bf Ref 1} and Lemma \ref{lemma:interpolation}, and the dominance condition \eqref{eq:dominantconditiononronSigmastarchap5} on $r$ on $\Si_*$, we deduce, for $k\leq k_*$, 
\beaa
|\ov{\nu^k(\widecheck{y})}| &\les& \frac{\ep}{ru^{\frac{1}{2}+\frac{\dec}{2}}}+\frac{\ep^2}{u^{2+\frac{3\dec}{2}}}\les \frac{\ep_0}{u^{\frac{3}{2}+\frac{3\dec}{2}}}.
\eeaa
Squaring, multiplying by $r^{-2}u^{2+2\dec}$ and integrating on $\Si_*$, we infer, for $k\leq k_*$, 
\beaa
\int_{\Si_*}r^{-2}u^{2+2\dec}\left|\ov{\nu^k(\widecheck{y})}\right|^2 &\les& \ep_0^2.
\eeaa
Together with \eqref{eq:intermediatefluxesitmateaverageb} and the fact that $z=-y-b_*$ in view of Lemma \ref{Lemma:eqts-nabeta,xib}, we deduce, for $k\leq k_*-7$,
\beaa
\int_{\Si_*}r^{-2}u^{2+2\dec}\left(\left|\ov{\nu^k(\widecheck{y})}\right|^2+\left|\ov{\nu^k(\widecheck{z})}\right|^2+\left|\ov{\nu^k(\widecheck{b_*})}\right|^2\right) &\les& \ep_0^2.
\eeaa
Together with \eqref{eq:intermediatefluxesitmateaveragefreexyb}, and using a Poincar\'e inequality, we infer, for $k\leq k_*-7$,
\beaa
\int_{\Si_*}r^{-2}u^{2+2\dec}\left(\left|\dk_*^k(\widecheck{y})\right|^2+\left|\dk_*^k(\widecheck{z})\right|^2+\left|\dk_*^k(\widecheck{b_*})\right|^2\right) &\les& \ep_0^2
\eeaa
which is the desired estimate \eqref{eq:fluxesitmateforxyb}.

{\bf Step 8.}  In this final step, we derive the desired estimate for $\ombc$,  i.e. we show
\bea\lab{eq:fluxesitmateforombc}
\int_{\Si_*}u^{2+2\dec}\left|\dk_*^k(\ombc)\right|^2 &\les& \ep_0^2, \quad k\leq k_*-7.
\eea
We   start with the equations
\beaa
\nab_4\kac  &=& \Ga_g\c\Ga_g,\\
 \nab_3\kac  &=&         2   \div \eta   + \frac{4}{r} \ombc +\frac{2}{r^2}\widecheck{y} +r^{-1}\Ga_g+\Ga_b\c\Ga_b,
 \eeaa
which together with the fact that $\nu=e_3+b_*e_4$ yield
\beaa
\nab_\nu\kac  &=&         2   \div \eta   + \frac{4}{r} \ombc +\frac{2}{r^2}\widecheck{y} +r^{-1}\Ga_g+\Ga_b\c\Ga_b.
\eeaa
Since our GCM assumption for $\ka$ on $\Si_*$ implies $\kac=0$, and since $\nu$ is tangent to $\Si_*$, we infer
\beaa
0  &=&         2   \div \eta   + \frac{4}{r} \ombc +\frac{2}{r^2}\widecheck{y} +r^{-1}\Ga_g+\Ga_b\c\Ga_b.
\eeaa
and hence
\beaa
\ombc &=& -\frac{r}{2}\div \eta -\frac{1}{2r}\widecheck{y} +\Ga_g+r\Ga_b\c\Ga_b
\eeaa
 Using the control of $\Ga_b$ and $\Ga_g$ provided by {\bf Ref 1} and Lemma \ref{lemma:interpolation}, and the dominance condition \eqref{eq:dominantconditiononronSigmastarchap5} on $r$ on $\Si_*$, we deduce, for $k\leq k_*$, 
\beaa
|\dk_*^k\ombc| &\les& |\dk_*^{k+1}\eta|+r^{-1}|\dk_*^k\widecheck{y}|+\frac{\ep}{r^2u^{\frac{1}{2}+\frac{\dec}{2}}}+\frac{\ep^2}{ru^{2+\frac{3\dec}{2}}}\\
&\les& |\dk_*^{k+1}\eta|+r^{-1}|\dk_*^k\widecheck{y}|+\frac{\ep_0}{ru^{\frac{3}{2}+\frac{3\dec}{2}}}.
\eeaa
Squaring, multiplying by $u^{2+2\dec}$, and integrating on $\Si_*$, we infer
\beaa
\int_{\Si_*}u^{2+2\dec}|\dk_*^k\ombc|^2 &\les& \int_{\Si_*}u^{2+2\dec}|\dk_*^{k+1}\eta|^2+\int_{\Si_*}r^{-2}u^{2+2\dec}|\dk_*^k\widecheck{y}|^2+\ep_0^2.
\eeaa
Together with \eqref{Estimates:Flux-eta-xib} for $\eta$ and \eqref{eq:fluxesitmateforxyb} for $\widecheck{y}$, we infer, for $k\leq k_*-7$, 
\beaa
\int_{\Si_*}u^{2+2\dec}|\dk_*^k\ombc|^2 &\les& \ep_0^2
\eeaa
which is the desired estimate for $\ombc$. Together with the estimates of Step 1 to Step 7, we deduce, for $k\leq k_*-7$, 
\beaa
\int_{\Si_*}u^{2+2\dec}|\dk_*^k\Ga_b|^2 &\les& \ep_0^2
\eeaa
which is the desired estimate \eqref{Estimate:Flux-bb-vthb-eta-xib}. This concludes the proof of Proposition \ref{Prop.Flux-bb-vthb-eta-xib}.
\end{proof}

As a corollary  of  the  above  we derive the following improved version of Lemma \ref{Le:Si*-ell=1modes}.
\begin{corollary}
\lab{corr:Si*-ell=1modes-improved}
The     functions $\Jp$ verify the following properties
\begin{enumerate}
\item We have on $\Si_*$
\beaa
\bsplit
\int_{S}J^{(p)} &= O\left(\ep ru^{-\frac{1}{2}-\frac{\dec}{2}}\right),\\
\int_{S}J^{(p)}J^{(q)} &= \frac{4\pi}{3}r^2\de_{pq}+O\left(\ep ru^{-\frac{1}{2}-\frac{\dec}{2}}\right).
\end{split}
\eeaa

\item For  any  $k\le k_*-10$, we have on $\Si_*$
\beaa
\left|\dk_*^k\left( \De +\frac{2}{r^2}\right)    \Jp\right|    \les  \ep r^{-3} u^{-\frac{1}{2}-\frac{\dec}{2}}.
\eeaa

\item We have for any $k\le k_*-10$ on $\Si_*$
\beaa
\left|\dk_*^k\dds_2\dds_1 \Jp\right|    \les  \ep r^{-3} u^{-\frac{1}{2}-\frac{\dec}{2}},
\eeaa
where by $\dds_1\Jp$, we mean either $\dds_1(\Jp,0)$ or $\dds_1(0,\Jp)$. 
\end{enumerate}
\end{corollary}

\begin{proof}
The proof follows exactly the same lines  as the proof of Lemma \ref{Le:Si*-ell=1modes} by replacing the pointwise estimate of $\Ga_b$ by the improved flux  estimates of Proposition \ref{Prop.Flux-bb-vthb-eta-xib} for $\Ga_b$.
\end{proof}

%%%%%%%%%%%%%%%%%%%%%%%%%%%%%%%%%

\section{Estimates for $\ell=0$ and $\ell=1$ modes on $\Si_*$}

%%%%%%%%%%%%%%%%%%%%%%%%%%%%%%%%%

%%%%%%%%%%%%%%%%%%%%%%%%%%%%%%

\subsection{Estimates for some $\ell=1$ modes on $S_*$}
\lab{section:Estimates-ell=1mpdesS_*}

%%%%%%%%%%%%%%%%%%%%%%%%%%%%%%

We start with estimates for the $\ell=1$ modes on $S_*$. Recall that on $S_*$ we have in particular the following, see  section \ref{sec:GCMconditionsonSigmastar},
\beaa
\ka=\frac 2 r, \quad \kab=-\frac{2\Up}{r}, \quad (\div \b)_{\ell=1}=0, \quad   (\curl\b)_{\ell=1,\pm}=0,  \quad  (\curl\b)_{\ell=1,0}=\frac{2am}{r^5}.
  \eeaa

\begin{lemma}
\lab{Lemma: K-ell=1}
The Gauss curvature $K$ of $S_*$ verifies
\bea
 |(\widehat{K})_{\ell=1}| &\les&  \frac{\ep_0}{r^3u^{2+2\dec}}.
\eea
\end{lemma}

\begin{proof}
 Recall that the metric on $S_*$ is given by
\beaa
g_{S_*}=e^{2\phi} r^2\Big((d\th)^2+(\sin\th)^2(d\vphi)^2\Big).
\eeaa
The  Gauss curvature thus satisfies 
\beaa
K&=&-\Delta\phi+\frac{e^{-2\phi}}{r^2} 
\eeaa
and hence
\beaa
K-\frac{1}{r^2}=-\left(\Delta+\frac{2}{r^2}\right)\phi+\frac{h(\phi)}{r^2}, \qquad h(\phi)=e^{-2\phi}-1+2\phi.
\eeaa
Integrating by parts,
\beaa
\int_{S_*}\left(K-\frac{1}{r^2}\right) \Jp &=& -\int_{S_*}\left(\Delta+\frac{2}{r^2}\right)\phi \Jp+\int_{S_*}\frac{h(\phi)}{r^2}\Jp\\
&=& -\int_{S_*}\phi \left(\Delta+\frac{2}{r^2}\right)\Jp+\int_{S_*}\frac{h(\phi)}{r^2}\Jp
\eeaa
Using  the control of $\phi$ of Lemma \ref{lemma:controloftheconformalfactorphi}, i.e. $\phi=O(\ep r^{-1} u^{-\frac{1}{2} -\frac{\dec}{2}})$  and the estimate for  $(\Delta+2/r^2) \Jp$ in Corollary \ref{corr:Si*-ell=1modes-improved}, we easily deduce that
\beaa
\left|\int_{S_*}\left(K-\frac{1}{r^2}\right) \Jp \right| &\les&\frac{\ep}{r^2u^{1+\dec}}.
\eeaa
Using the dominance condition \eqref{eq:dominantconditiononronSigmastarchap5} on $r$ on $\Si_*$, we infer
\beaa
\left|\int_{S_*}\left(K-\frac{1}{r^2}\right) \Jp \right| &\les&\frac{\ep_0}{ru^{2+2\dec}},
\eeaa
i.e.  recalling  the definition   of $\ell=1$ modes,
\beaa
  |(\widehat{K})_{\ell=1}| &\les&  \frac{\ep_0}{r^3u^{2+2\dec}}
\eeaa
as stated.
\end{proof}

\begin{lemma}
\lab{Le:ell=1modesonS_*}
The following holds on $S_*$
\bea
\left|\left(\rhoc-\frac{1}{2}\chih\c\chibh\right)_{\ell=1} \right| &\les&  \ep_0 r^{-3}u^{-2-2\dec}.
\eea
\end{lemma}

\begin{proof} 
Recall that our GCM conditions on $S_*$ imply in particular $\kac=\kabc=0$ on $S_*$. This implies that the Gauss  equation on $S_*$     takes the form
\beaa
\widecheck{K}  &=& -\rhoc +\frac{1}{2}\chih\c\chibh
\eeaa
and hence
\beaa
\left(\rhoc-\frac{1}{2}\chih\c\chibh\right)_{\ell=1} &=& -\left(\widecheck{K}\right)_{\ell=1}.
\eeaa
Hence, in view of Lemma \ref{Lemma: K-ell=1}, we obtain 
\beaa
\left|\left(\rhoc-\frac{1}{2}\chih\c\chibh\right)_{\ell=1} \right| &\les&  \ep_0 r^{-3}u^{-2-2\dec}
\eeaa
as stated.
\end{proof}

%%%%%%%%%%%%%%%%%%%%%%%%%%%%%%%%%%%%

\subsection{Estimates for the $\ell=1$ modes on $ \Si_*$}
\lab{section:estimates-rll=1-Si_*}

%%%%%%%%%%%%%%%%%%%%%%%%%%%%%%%%%%%%

\begin{proposition}[Control of $\ell=1$ modes]\lab{prop:controlofell=1modesonSigmastar}
\lab{prop:control.ell=1modes-Si}
The following  estimates hold on $\Si_*$.
\bea
\bsplit
\left|(\div\b)_{\ell=1}\right|+\left|(\curl\b)_{\ell=1,\pm}\right|+\left|(\curl\b)_{\ell=1,0}-\frac{2am}{r^5}\right|&\les \ep_0 r^{-5} u^{-1-\dec}, \\
\left|(\nu\div\b)_{\ell=1}\right|+\left|(\nu\curl\b)_{\ell=1}\right| &\les \ep_0 r^{-5} u^{-1-\dec}, \\
\left|(\div\ze)_{\ell=1}\right|+\left|(\curl\ze )_{\ell=1,\pm}\right|+\left|(\curl\ze)_{\ell=1,0}-\frac{2am}{r^4}\right|&\les\ep_0  r^{-4} u^{-1-\dec}, \\
\left|(\kabc)_{\ell=1}\right| & \les \ep_0 r^{-2} u^{-2-2\dec}, 
\end{split}
\eea
\bea
\bsplit
\left|\left(\rhoc -\frac{1}{2}\chih\c\chibh \right)_{\ell=1}\right|&\les  \ep_0 r^{-3} u^{-2 -2\dec},\\
\left|\left(\dual \rho -\frac 1 2 \chih\wedge\chibh\right)_{\ell=1,\pm}\right|+\left|\left(\dual \rho -\frac 1 2 \chih\wedge\chibh\right)_{\ell=1,0}-\frac{2am}{r^4}\right|&\les \ep_0 r^{-3}u^{-2-2\dec},\\
\left|(\muc)_{\ell=1}\right|&\les \ep_0 r^{-3} u^{-2-2\dec}.
\end{split}
\eea
\end{proposition}

\begin{remark}
Note that the results are consistent with strong peeling.
\end{remark}

\begin{proof}
We  make the following local bootstrap   assumption on $\Si_*$
\bea
\lab{eq:BAlocal-Si_*}
\bsplit
|(\div\b)_{\ell=1}|+|(\curl\b)_{\ell=1,\pm}|+\left|(\curl\b)_{\ell=1,0}-\frac{2am}{r^5}\right|&\leq \ep r^{-5} u^{-1-\dec} ,\\
 \left|(\kabc)_{\ell=1}\right| &\leq \ep r^{-2} u^{-2-\dec}. 
\end{split}
\eea

{\bf Step 1.} We start with the control of $(\ddd_1\ze)_{\ell=1}$. Recall  the following consequence of the  Codazzi equation for $\chih$
\beaa
\ddd_2\chih &=& \frac{1}{r}\ze - \b+\Ga_g\c \Ga_g.
\eeaa
Differentiating w.r.t. $\ddd_1$, we infer
\beaa
\ddd_1\ddd_2\chih &=& \frac{1}{r}\ddd_1\ze - \ddd_1\b+r^{-1}\dkb^{\leq 1}(\Ga_g\c \Ga_g).
\eeaa
Projecting on the $\ell=1$ modes, this yields 
\beaa
(\ddd_1\ddd_2\chih)_{\ell=1} &=&  \frac{1}{r}(\ddd_1\ze)_{\ell=1}-(\ddd_1\b)_{\ell=1} +r^{-1}\dkb^{\leq 1}(\Ga_g\c \Ga_g).
\eeaa
Next, we estimate $(\ddd_1\ddd_2\chih)_{\ell=1}$. We have
\beaa
(\ddd_1\ddd_2\chih)_{\ell=1,p}=\frac{1}{|S|}\int_S\ddd_1\ddd_2\chih\Jp =\frac{1}{|S|}\int_S\chih\c\dds_2\dds_1\Jp
\eeaa
and hence
\beaa
|(\ddd_1\ddd_2\chih)_{\ell=1}| &\les& |\dds_2\dds_1\Jp||\Ga_g|.
\eeaa
We deduce
\beaa
 (\ddd_1\ze)_{\ell=1} &=& r(\ddd_1\b)_{\ell=1} + r|\dds_2\dds_1\Jp|\Ga_g+\dkb^{\leq 1}(\Ga_g\c \Ga_g)
\eeaa
and thus
\beaa
|(\div\ze)_{\ell=1}|+|(\curl\ze)_{\ell=1,\pm}| &\les& r|(\div\b)_{\ell=1}|+r|(\curl\b)_{\ell=1,\pm}| + r|\dds_2\dds_1\Jp|\Ga_g\\
&&+\dkb^{\leq 1}(\Ga_g\c \Ga_g),\\
\left|(\curl\ze)_{\ell=1,0}-\frac{2am}{r^4}\right| &\les& r\left|(\curl\b)_{\ell=1,0}-\frac{2am}{r^5}\right| + r|\dds_2\dds_1\Jp|\Ga_g+\dkb^{\leq 1}(\Ga_g\c \Ga_g).
\eeaa
Using the control of $\Ga_g$ provided by {\bf Ref 1}, the control of $\dds_2\dds_1\Jp$ provided by Corollary  \ref{corr:Si*-ell=1modes-improved}, and the local bootstrap assumption \eqref{eq:BAlocal-Si_*} on $(\ddd_1\b)_{\ell=1}$, we deduce
\bea
\lab{eq:control.ell=1modes-Si-Step1}
\left|(\div\ze)_{\ell=1}\right|+\left|(\curl\ze )_{\ell=1,\pm}\right|+\left|(\curl\ze)_{\ell=1,0}-\frac{2am}{r^4}\right|&\les& \frac{\ep}{r^4u^{1+\dec}}. 
\eea

{\bf Step 2.}   Next, we consider the control of $(\div\bb)_{\ell=1}$. Recall  the following consequence of the  Codazzi equation for $\chibh$
\beaa
\div\chibh &=& \frac{1}{2}\nab\kabc +\frac{\Up}{r}\ze  +\bb+\Ga_b\c \Ga_g.
\eeaa
Differentiating w.r.t. $\div$, we infer
\beaa
\div\ddd_2\chibh &=& \frac{1}{2}\Delta\kabc+\frac{\Up}{r}\div\ze + \div\bb+r^{-1}\dkb^{\leq 1}(\Ga_b\c \Ga_g).
\eeaa
Projecting on the $\ell=1$ modes, this yields 
\beaa
(\div\ddd_2\chibh)_{\ell=1} &=&  \frac{1}{2}(\Delta\kabc)_{\ell=1}+\frac{\Up}{r}(\div\ze)_{\ell=1} +(\div\bb)_{\ell=1} +r^{-1}\dkb^{\leq 1}(\Ga_b\c \Ga_g).
\eeaa
As in Step 1, we have
\beaa
|(\div\ddd_2\chibh)_{\ell=1}| &\les& |\dds_2\dds_1\Jp||\Ga_b|.
\eeaa
Also, we have
\beaa
(\Delta\kabc)_{\ell=1,p}=\frac{1}{|S|}\int_S\Delta\kabc \Jp =-\frac{2}{r^2}\frac{1}{|S|}\int_S\kabc\Jp+\frac{1}{|S|}\int_S\kabc\left(\Delta+\frac{2}{r^2}\right)\Jp
\eeaa
and hence
\beaa
|(\Delta\kabc)_{\ell=1}| &\les& r^{-2}|(\kabc)_{\ell=1}|+\left|\left(\Delta+\frac{2}{r^2}\right)\Jp\right||\Ga_g|.
\eeaa
We deduce
\beaa
|(\div\bb)_{\ell=1}| &\les& r^{-2}|(\kabc)_{\ell=1}|+r^{-1}|(\div\ze)_{\ell=1}|+\left|\left(\Delta+\frac{2}{r^2}\right)\Jp\right||\Ga_g|\\
&&+ |\dds_2\dds_1\Jp||\Ga_b|+r^{-1}|\dkb^{\leq 1}(\Ga_b\c \Ga_g)|.
\eeaa
Together with  the local bootstrap assumption \eqref{eq:BAlocal-Si_*} on $(\kabc)_{\ell=1}$, the control of $(\div\ze)_{\ell=1}$ in \eqref{eq:control.ell=1modes-Si-Step1}, the control of $\dds_2\dds_1\Jp$ and $(\De+2/r^2)\Jp$ provided by Corollary  \ref{corr:Si*-ell=1modes-improved}, and the  control of $\Ga_g$ provided by {\bf Ref 1}, we obtain 
\beaa
|(\div\bb)_{\ell=1}| &\les& \frac{\ep}{r^4u^{2+\dec}}+\frac{\ep}{r^5u^{1+\dec}}+\frac{\ep}{r^3u^{\frac{1}{2}+\frac{\dec}{2}}}|\dkb^{\leq 1}\Ga_b|.
\eeaa
Using the dominance condition \eqref{eq:dominantconditiononronSigmastarchap5} on $r$ on $\Si_*$, we infer
\beaa
|(\div\bb)_{\ell=1}| &\les& \frac{\ep_0}{r^3u^{3+2\dec}}+\frac{\ep_0}{r^2u^{\frac{3}{2}+\frac{3\dec}{2}}}|\dkb^{\leq 1}\Ga_b|.
\eeaa
By integration in $u$, and using Sobolev,  we deduce
\beaa
\int_u^{u_*}  r^3  \big    |  (\div\bb)_{\ell=1}\big| &\les& \frac{\ep_0}{u^{2+2\dec}} + \ep_0\int_u^{u_*}u^{-\frac{3}{2}-\frac{3\dec}{2}}    \| \dkb^{\le 3}\Ga_b\|_{L^2(S)}
\\
 &\les& \frac{\ep_0}{u^{2+2\dec}} +  \frac{\ep_0}{u^{2+2\dec}}\left(\int_{\Si_*}u^{2+2\dec}|\dkb^{\le 3}\Ga_b|^2\right)^{\frac{1}{2}}.
\eeaa
Hence, in view of Proposition \ref{Prop.Flux-bb-vthb-eta-xib}, we obtain
\bea
\lab{eq:control.ell=1modes-Si-Step2}
\int_u^{u_*}  r^3  \big    |  (\div\bb)_{\ell=1}\big| &\les& \frac{\ep_0}{u^{2+2\dec}}. 
\eea

{\bf Step 3.}  We provide  the estimate for   $(\rhoc-\frac{1}{2}\chih\c\chibh)_{\ell=1}$. Recall from Corollary \ref{corofLemma:transport.alongSi_*1} that we have along $\Si_*$, for $p=0,+,-$,
\beaa
\nn\nu\left(\int_S\left(\rhoc - \frac{1}{2}\chih\c\chibh \right)\Jp\right) &=&  O(r^{-1})\int_S\left(\rhoc - \frac{1}{2}\chih\c\chibh \right)\Jp +h_1,
\eeaa
where the scalar function $h_1$ is given by 
\beaa
h_1 &=& -\int_S\div\bb\Jp -(1+O(r^{-1}))\int_S\div\b\Jp  +O(r^{-3})\int_S\kabc\Jp\\
   && +O(r^{-2})\int_S\div\ze\Jp+r\left|\left(\Delta+\frac{2}{r^2}\right)\Jp\right|\Ga_b +r\dkb^{\le  1}( \Ga_b \c \Ga_b).
\eeaa
In view of  the local bootstrap assumption \eqref{eq:BAlocal-Si_*} on $(\kabc)_{\ell=1}$ and $(\div\b)_{\ell=1}$, the control of $(\div\ze)_{\ell=1}$ in \eqref{eq:control.ell=1modes-Si-Step1}, the control of $(\De+2/r^2)\Jp$ provided by Corollary  \ref{corr:Si*-ell=1modes-improved}, and the  control of $\Ga_g$ provided by {\bf Ref 1}, we obtain 
\beaa
|h_1| &\les& r^2|(\div\bb)_{\ell=1}|+\frac{\ep}{r^3u^{1+\dec}}+r|\dkb^{\le  1}( \Ga_b \c \Ga_b)|.
\eeaa
Using the dominance condition \eqref{eq:dominantconditiononronSigmastarchap5} on $r$ on $\Si_*$, we infer
\beaa
|h_1| &\les& r^2|(\div\bb)_{\ell=1}|+r|\dkb^{\le  1}( \Ga_b \c \Ga_b)|+\frac{\ep_0}{ru^{3+3\dec}}.
\eeaa
By integration in $u$, and using Sobolev,  we deduce
\beaa
\int_u^{u_*}r|h_1| &\les& \frac{\ep_0}{u^{2+3\dec}}+\int_u^{u_*}r^3|(\div\bb)_{\ell=1}|+\int_u^{u_*}\|\dkb^{\le  3}\Ga_b\|_{L^2(S)}^2\\
&\les& \frac{\ep_0}{u^{2+3\dec}}+\int_u^{u_*}r^3|(\div\bb)_{\ell=1}|+\frac{\ep_0}{u^{2+2\dec}}\left(\int_{\Si_*}u^{2+2\dec}|\dkb^{\le  3}\Ga_b|^2\right)^{\frac{1}{2}}.
\eeaa
Together with the control of $(\div\bb)_{\ell=1}$ in \eqref{eq:control.ell=1modes-Si-Step2}  and the control of $\Ga_b$ in Proposition \ref{Prop.Flux-bb-vthb-eta-xib}, we obtain
\beaa
\int_u^{u_*}r|h_1| &\les& \frac{\ep_0}{u^{2+2\dec}}.
\eeaa
Since 
\beaa
\nn\nu\left(\int_S\left(\rhoc - \frac{1}{2}\chih\c\chibh \right)\Jp\right) &=&  O(r^{-1})\int_S\left(\rhoc - \frac{1}{2}\chih\c\chibh \right)\Jp +h_1,
\eeaa
we may thus apply Corollary \ref{cor:integrate-transportSi_*} which implies
\beaa
ru^{2+2\dec}\left|\int_S\left(\rhoc - \frac{1}{2}\chih\c\chibh \right)\Jp\right| &\les& r_*u_*^{2+2\dec}\left|\int_{S_*}\left(\rhoc - \frac{1}{2}\chih\c\chibh \right)\Jp\right|+\ep_0. 
\eeaa
Together with the control of $(\rhoc-\frac{1}{2}\chih\c\chibh)_{\ell=1}$ on $S_*$ provided by Lemma \ref{Le:ell=1modesonS_*}, we infer
\beaa
ru^{2+2\dec}\left|\int_S\left(\rhoc - \frac{1}{2}\chih\c\chibh \right)\Jp\right| &\les& \ep_0
\eeaa
and hence
\bea\lab{eq:control.ell=1modes-Si-Step3}
\left|\left(\rhoc - \frac{1}{2}\chih\c\chibh \right)_{\ell=1}\right| &\les& \frac{\ep_0}{r^3u^{2+2\dec}}.
\eea

{\bf Step 4.} We provide  the estimate for   $(\kabc)_{\ell=1}$. Recall from Corollary \ref{corofLemma:transport.alongSi_*1} that we have along $\Si_*$, for $p=0,+,-$,
\beaa
\nu\left(\int_S\left( \lap\kabc+\frac{2\Up}{r}\div \ze\right)\Jp\right) &=& h_2
\eeaa
where the scalar function $h_2$ is given by 
\beaa
h_2 &=& O(r^{-3})\int_S\kabc\Jp   +O(r^{-2})\int_S\div\ze\Jp +O(r^{-1})\int_S\div\bb\Jp    \\
&& +O(r^{-2})\int_S\left(\rhoc-\frac{1}{2}\chih\c\chibh\right)\Jp +O(r^{-1})\int_S\div\b\Jp\\
&& +r\left|\left(\Delta+\frac{2}{r^2}\right)\Jp\right|\dkb^{\leq 1}\Ga_b+\dkb^{\le 2 }(\Ga_b\c \Ga_b).
\eeaa
In view of  the local bootstrap assumption \eqref{eq:BAlocal-Si_*} on $(\kabc)_{\ell=1}$ and $(\div\b)_{\ell=1}$, the control of $(\div\ze)_{\ell=1}$ in \eqref{eq:control.ell=1modes-Si-Step1}, and the control of $(\De+2/r^2)\Jp$ provided by Corollary  \ref{corr:Si*-ell=1modes-improved}, we obtain 
\beaa
|h_2| &\les& r|(\div\bb)_{\ell=1}|+\frac{\ep}{r^3u^{2+\dec}}+\frac{\ep}{r^4u^{1+\dec}}
+\frac{\ep}{r^2u^{\frac{1}{2}+\frac{\dec}{2}}}|\Ga_b|+|\dkb^{\le  2}( \Ga_b \c \Ga_b)|\\
&\les&  r|(\div\bb)_{\ell=1}|+\frac{\ep}{r^3u^{2+\dec}}+\frac{\ep}{r^4u^{1+\dec}}+|\dkb^{\le  2}( \Ga_b \c \Ga_b)|.
\eeaa
Using the dominance condition \eqref{eq:dominantconditiononronSigmastarchap5} on $r$ on $\Si_*$, we infer
\beaa
|h_2| &\les& r|(\div\bb)_{\ell=1}|+\frac{\ep_0}{r^2u^{3+2\dec}}+|\dkb^{\le  2}( \Ga_b \c \Ga_b)|.
\eeaa
By integration in $u$, and using Sobolev,  we deduce
\beaa
\int_u^{u_*}r^2|h_2| &\les& \frac{\ep_0}{u^{2+2\dec}}+\int_u^{u_*}r^3|(\div\bb)_{\ell=1}|+\int_u^{u_*}\|\dkb^{\le  4}\Ga_b\|_{L^2(S)}^2\\
&\les& \frac{\ep_0}{u^{2+2\dec}}+\int_u^{u_*}r^3|(\div\bb)_{\ell=1}|+\frac{\ep_0}{u^{2+2\dec}}\left(\int_{\Si_*}u^{2+2\dec}|\dkb^{\le  4}\Ga_b|^2\right)^{\frac{1}{2}}.
\eeaa
Together with the control of $(\div\bb)_{\ell=1}$ in \eqref{eq:control.ell=1modes-Si-Step2}  and the control of $\Ga_b$ in Proposition \ref{Prop.Flux-bb-vthb-eta-xib}, we obtain
\beaa
\int_u^{u_*}r^2|h_2| &\les& \frac{\ep_0}{u^{2+2\dec}}.
\eeaa
Since 
\beaa
\nu\left(\int_S\left( \lap\kabc+\frac{2\Up}{r}\div \ze\right)\Jp\right) &=& h_2,
\eeaa
we may thus apply Corollary \ref{cor:integrate-transportSi_*} which implies
\beaa
r^2u^{2+2\dec}\left|\int_S\left( \lap\kabc+\frac{2\Up}{r}\div \ze\right)\Jp\right| &\les& r_*^2u_*^{2+2\dec}\left|\int_{S_*}\left( \lap\kabc+\frac{2\Up}{r}\div \ze\right)\Jp\right|+\ep_0. 
\eeaa
Since $\kabc=0$ on $S_*$, and using  the control of $(\div\ze)_{\ell=1}$ in \eqref{eq:control.ell=1modes-Si-Step1}, we infer
\beaa
r^2u^{2+2\dec}\left|\int_S\lap\kabc \Jp\right| &\les& \ep_0+\frac{\ep u^{1+\dec}}{r}
\eeaa
and hence
\beaa
u^{2+2\dec}\left|\int_S\kabc \Jp\right| &\les& \ep_0+\frac{\ep u^{1+\dec}}{r}+r^4u^{2+2\dec}\left|\left(\Delta+\frac{2}{r^2}\right)\Jp\right||\Ga_g|.
\eeaa
Together with  the  control of $\Ga_g$ provided by {\bf Ref 1}, the control of $(\De+2/r^2)\Jp$ provided by Corollary  \ref{corr:Si*-ell=1modes-improved}, and the dominance condition \eqref{eq:dominantconditiononronSigmastarchap5} on $r$ on $\Si_*$, we infer
\beaa
u^{2+2\dec}\left|\int_S\kabc \Jp\right| &\les& \ep_0+\frac{\ep u^{1+\dec}}{r}\les \ep_0
\eeaa
and hence
\bea\lab{eq:control.ell=1modes-Si-Step4}
\left|\left(\kabc \right)_{\ell=1}\right| &\les& \frac{\ep_0}{r^2u^{2+2\dec}}.
\eea

{\bf Step 5.} Next, we estimate $(\rhod-\frac{1}{2}\chih\wedge\chibh)_{\ell=1}$. Recall that we have
\beaa
\curl\ze&=&   \dual \rho -\frac 1 2 \chih\wedge\chibh.
\eeaa
We infer from the control of $(\curl\ze)_{\ell=1}$ in \eqref{eq:control.ell=1modes-Si-Step1}
\beaa
\left|\left(\dual \rho -\frac 1 2 \chih\wedge\chibh\right)_{\ell=1,\pm}\right|+\left|\left(\dual \rho -\frac 1 2 \chih\wedge\chibh\right)_{\ell=1,0}-\frac{2am}{r^4}\right|&\les& \frac{\ep}{r^4u^{1+\dec}}. 
\eeaa
Together with  the   dominance condition \eqref{eq:dominantconditiononronSigmastarchap5} on $r$ on $\Si_*$, we infer
\bea
\lab{eq:control.ell=1modes-Si-Step5}
\left|\left(\dual \rho -\frac 1 2 \chih\wedge\chibh\right)_{\ell=1,\pm}\right|+\left|\left(\dual \rho -\frac 1 2 \chih\wedge\chibh\right)_{\ell=1,0}-\frac{2am}{r^4}\right|&\les& \frac{\ep_0}{r^3u^{2+2\dec}}. 
\eea

{\bf Step 6.} We provide  the estimate for   $(\div\b)_{\ell=1}$. Recall from Corollary \ref{corofLemma:transport.alongSi_*1} that we have along $\Si_*$, for $p=0,+,-$,
 \beaa
\nu\left(\int_S\div\b\Jp\right) &=& O(r^{-1})\int_S\div\b\Jp+h_3
\eeaa
where the scalar function $h_3$ is given by 
\beaa
h_3 &=& O(r^{-2})\int_S\left(\rhoc-\frac{1}{2}\chih\c\chibh\right)\Jp+r\left(\left|\left(\Delta+\frac{2}{r^2}\right)\Jp\right|+\left|\dds_2\dds_1\Jp\right|\right)\Ga_g+\dkb^{\leq 1}(\Ga_b\c\Ga_g).
\eeaa
In view of  the control of $(\rhoc-\frac{1}{2}\chih\c\chibh)_{\ell=1}$ in \eqref{eq:control.ell=1modes-Si-Step3}, the control of $(\De+2/r^2)\Jp$ and $\dds_2\dds_1\Jp$ provided by Corollary  \ref{corr:Si*-ell=1modes-improved}, and  the  control of $\Ga_g$ provided by {\bf Ref 1}, we obtain 
\beaa
|h_3| &\les& \frac{\ep_0}{r^3u^{2+2\dec}}+\frac{\ep^2}{r^4u^{1+\frac{3\dec}{2}}}+\frac{\ep}{r^2u^{\frac{1}{2}+\dec}}|\dkb^{\le  1}\Ga_b|.
\eeaa
Using the dominance condition \eqref{eq:dominantconditiononronSigmastarchap5} on $r$ on $\Si_*$, we infer
\beaa
|h_3| &\les& \frac{\ep_0}{r^3u^{2+2\dec}}+\frac{\ep}{r^2u^{\frac{1}{2}+\dec}}|\dkb^{\le  1}\Ga_b|.
\eeaa
By integration in $u$, and using Sobolev,  we obtain
\beaa
\int_u^{u_*}r^3|h_3| &\les& \frac{\ep_0}{u^{1+\dec}}+\int_u^{u_*}\frac{\ep}{{u'}^{\frac{1}{2}+\dec}}\|\dkb^{\le  3}\Ga_b\|_{L^2(S)}\\
&\les& \frac{\ep_0}{u^{1+\dec}}+\frac{\ep_0}{u^{1+2\dec}}\left(\int_{\Si_*}u^{2+2\dec}|\dkb^{\le  3}\Ga_b|^2\right)^{\frac{1}{2}}.
\eeaa
Together with the control of  $\Ga_b$ in Proposition \ref{Prop.Flux-bb-vthb-eta-xib}, we infer
\beaa
\int_u^{u_*}r^3|h_3| &\les& \frac{\ep_0}{u^{1+\dec}}.
\eeaa
Since 
 \beaa
\nu\left(\int_S\div\b\Jp\right) &=& O(r^{-1})\int_S\div\b\Jp+h_3,
\eeaa
we may thus apply Corollary \ref{cor:integrate-transportSi_*} which implies, together with the fact that $(\div\b)_{\ell=1}=0$ on $S_*$, 
\beaa
r^3u^{1+\dec}\left|\int_S\div\b\Jp\right| &\les& \ep_0,
\eeaa
and hence
\bea\lab{eq:control.ell=1modes-Si-Step6}
\left|\left(\div\b \right)_{\ell=1}\right| &\les& \frac{\ep_0}{r^5u^{1+\dec}}.
\eea

{\bf Step 7.} We provide  the estimate for   $(\curl\b)_{\ell=1}$. Recall from Corollary \ref{corofLemma:transport.alongSi_*1} that we have along $\Si_*$, for $p=0,+,-$,
\beaa
\nn\nu\left(\int_S\curl\b\Jp\right) &=&  \frac{4}{r}(1+ O(r^{-1}))\int_S\curl\b\Jp+\frac{2}{r^2}(1+ O(r^{-1}))\int_S\left(\rhod  -\frac 1 2 \chih\wedge\chibh\right)\Jp\\
\nn&&+r\left(\left|\left(\Delta+\frac{2}{r^2}\right)\Jp\right|+\left|\dds_2\dds_1\Jp\right|\right)\Ga_g+\dkb^{\leq 1}(\Ga_b\c\Ga_g).
\eeaa
In the case $p=\pm$, since we have $(\curl\b)_{\ell=1,\pm}=0$ on $S_*$, using \eqref{eq:control.ell=1modes-Si-Step5} to control $(\rhod  -\frac 1 2 \chih\wedge\chibh)_{\ell=1,\pm}$, and arguing exactly as for the control of $(\div\b)_{\ell=1}$ in Step 6, we obtain the following analog of \eqref{eq:control.ell=1modes-Si-Step6}
 \bea\lab{eq:control.ell=1modes-Si-Step7:intermediate}
\left|\left(\curl\b \right)_{\ell=1,\pm}\right| &\les& \frac{\ep_0}{r^5u^{1+\dec}}.
\eea
 
Next, we focus on the case $p=0$. We rewrite the above transport equation in this particular case
\beaa
\nu\left(\int_S\curl\b J^{(0)}\right) &=&  \frac{4}{r}(1+ O(r^{-1}))\int_S\curl\b J^{(0)}+\frac{2}{r^2}(1+ O(r^{-1}))\int_S\left(\rhod  -\frac 1 2 \chih\wedge\chibh\right)J^{(0)}\\
\nn&&+r\left(\left|\left(\Delta+\frac{2}{r^2}\right)J^{(0)}\right|+\left|\dds_2\dds_1J^{(0)}\right|\right)\Ga_g+\dkb^{\leq 1}(\Ga_b\c\Ga_g).
\eeaa 
Since $\nu(r)=-2+r\Ga_b$, we have
\beaa
\nu\left(r^3\int_S\curl\b J^{(0)}\right) &=& r^3\nu\left(\int_S\curl\b J^{(0)}\right)+3r^2\nu(r)\int_S\curl\b J^{(0)}\\
&=& r^3\nu\left(\int_S\curl\b J^{(0)}\right) -6r^2\int_S\curl\b J^{(0)}+r^5\Ga_b(\curl\b)_{\ell=1,0},
\eeaa
and hence
\beaa
&&\nu\left(r^3\int_S\curl\b J^{(0)}\right)\\
 &=&  -\frac{2}{r}r^3(1+ O(r^{-1}))\int_S\curl\b J^{(0)}+2r(1+ O(r^{-1}))\int_S\left(\rhod  -\frac 1 2 \chih\wedge\chibh\right)J^{(0)}\\
\nn&&+r^5\Ga_b(\curl\b)_{\ell=1,0}+r^4\left(\left|\left(\Delta+\frac{2}{r^2}\right)J^{(0)}\right|+\left|\dds_2\dds_1J^{(0)}\right|\right)\Ga_g+r^3\dkb^{\leq 1}(\Ga_b\c\Ga_g)
\eeaa
which we rewrite 
\beaa
\nu\left(r^3\int_S\curl\b J^{(0)} -8\pi am\right) &=&  -\frac{2}{r}\left(r^3\int_S\curl\b J^{(0)}-8\pi am\right)+h_4
\eeaa
where the scalar function $h_4$ is given by 
\beaa
h_4 &=& O(r^3)\left(\left(\rhod  -\frac 1 2 \chih\wedge\chibh\right)_{\ell=1,0} -\frac{2am}{r^4}\right)+O(r^3)(\curl\b J^{(0)})_{\ell=1,0}\\
&&+O(r^2)\left(\rhod  -\frac 1 2 \chih\wedge\chibh\right)_{\ell=1,0}+r^5\Ga_b(\curl\b)_{\ell=1,0}\\
&&+r^4\left(\left|\left(\Delta+\frac{2}{r^2}\right)J^{(0)}\right|+\left|\dds_2\dds_1J^{(0)}\right|\right)\Ga_g+r^3\dkb^{\leq 1}(\Ga_b\c\Ga_g).
\eeaa
In view of  the control of $(\rhod-\frac{1}{2}\chih\wedge\chibh)_{\ell=1}$ in \eqref{eq:control.ell=1modes-Si-Step5}, 
the local bootstrap assumption \eqref{eq:BAlocal-Si_*} on $(\curl\b)_{\ell=1,0}$, the control of $(\De+2/r^2)\Jp$ and $\dds_2\dds_1\Jp$ provided by Corollary  \ref{corr:Si*-ell=1modes-improved}, and  the  control of  $\Ga_g$ provided by {\bf Ref 1}, we obtain 
\beaa
|h_4| &\les& \frac{1}{r^2}+\frac{\ep_0}{u^{2+2\dec}}+\frac{\ep}{ru^{1+\dec}}+\frac{r\ep}{u^{\frac{1}{2}+\dec}}|\dkb^{\le  1}\Ga_b|.
\eeaa
Using the dominance condition \eqref{eq:dominantconditiononronSigmastarchap5} on $r$ on $\Si_*$, we infer
\beaa
|h_4| &\les& \frac{\ep_0}{u^{2+2\dec}}+\frac{r\ep}{u^{\frac{1}{2}+\dec}}|\dkb^{\le  1}\Ga_b|.
\eeaa
By integration in $u$, and using Sobolev,  we obtain
\beaa
\int_u^{u_*}|h_4| &\les& \frac{\ep_0}{u^{1+\dec}}+\int_u^{u_*}\frac{\ep}{{u'}^{\frac{1}{2}+\dec}}\|\dkb^{\le  3}\Ga_b\|_{L^2(S)}\\
&\les& \frac{\ep_0}{u^{1+\dec}}+\frac{\ep_0}{u^{1+2\dec}}\left(\int_{\Si_*}u^{2+2\dec}|\dkb^{\le  3}\Ga_b|^2\right)^{\frac{1}{2}}.
\eeaa
Together with the control of  $\Ga_b$ in Proposition \ref{Prop.Flux-bb-vthb-eta-xib}, we infer
\beaa
\int_u^{u_*}|h_4| &\les& \frac{\ep_0}{u^{1+\dec}}.
\eeaa
Since 
\beaa
\nu\left(r^3\int_S\curl\b J^{(0)} -8\pi am\right) &=&  -\frac{2}{r}\left(r^3\int_S\curl\b J^{(0)}-8\pi am\right)+h_4,
\eeaa
we may thus apply Corollary \ref{cor:integrate-transportSi_*} which implies, together with the fact that there holds $(\curl\b)_{\ell=1,0}=\frac{2am}{r^5}$ on $S_*$, 
\beaa
r^3u^{1+\dec}\left|\int_S\curl\b J^{(0)}-\frac{8\pi am}{r^3}\right| &\les& \ep_0,
\eeaa
and hence, together with \eqref{eq:control.ell=1modes-Si-Step7:intermediate}, we have obtained 
\bea\lab{eq:control.ell=1modes-Si-Step7}
\left|\left(\curl\b \right)_{\ell=1,\pm}\right|+\left|\left(\curl\b \right)_{\ell=1,0}-\frac{2am}{r^5}\right| &\les& \frac{\ep_0}{r^5u^{1+\dec}}.
\eea

\begin{remark}\lab{rmk:improvementonlocalbootasseq:BAlocal-Si_*}
Note that \eqref{eq:control.ell=1modes-Si-Step4} for $(\kabc)_{\ell=1}$, \eqref{eq:control.ell=1modes-Si-Step6} for $(\div\b)_{\ell=1}$, and \eqref{eq:control.ell=1modes-Si-Step7} for $(\curl\ze)_{\ell=1}$,  improve the local bootstrap assumptions \eqref{eq:BAlocal-Si_*}.
\end{remark}

{\bf Step 8.} We have by the definition of the mass aspect function $\mu$
\beaa
\muc &=& -\div\ze -\left(\rhoc -\frac{1}{2}\chih\c\chibh\right)
\eeaa
and hence
\beaa
(\muc)_{\ell=1} &=& -(\div\ze)_{\ell=1} -\left(\rhoc -\frac{1}{2}\chih\c\chibh\right)_{\ell=1}.
\eeaa
Together with the estimates \eqref{eq:control.ell=1modes-Si-Step1} for $(\div\ze)_{\ell=1}$ and \eqref{eq:control.ell=1modes-Si-Step3} for $(\rhoc -\frac{1}{2}\chih\c\chibh)_{\ell=1}$, we infer
\beaa
|(\muc)_{\ell=1}| &\les& \frac{\ep_0}{r^3u^{2+2\dec}}+\frac{\ep}{r^4u^{1+\dec}}.
\eeaa
Using the dominance condition \eqref{eq:dominantconditiononronSigmastarchap5} on $r$ on $\Si_*$, we deduce
\bea\lab{eq:control.ell=1modes-Si-Step8}
|(\muc)_{\ell=1}| &\les& \frac{\ep_0}{r^3u^{2+2\dec}}.
\eea

{\bf Step 9.} It remains to derive estimates for $(\nu\div\b)_{\ell=1}$ and $(\nu\curl\b)_{\ell=1}$. In view of Lemma \ref{Lemma:transport.alongSi_*1}, we have along  $\Si_*$
\beaa
\bsplit
         \nu\div\b  &=  O(r^{-1})\div \b+\lap\rho +(1+O(r^{-1}))\div\div\a\\
         & +O(r^{-3})\div\eta+r^{-2}\dkb^{\leq 1}(\Ga_b\c\Ga_g),\\
    \nu\curl\b  &=  \frac{8}{r}(1+ O(r^{-1}))\curl \b -\lap\rhod +(1+O(r^{-1}))\curl\div\a \\
    &+O(r^{-3})\rhod+r^{-2}\dkb^{\leq 1}(\Ga_b\c\Ga_g),
  \end{split}
\eeaa
where the notation $O(r^a)$, for $a\in\mathbb{R}$, denotes an explicit function of $r$ which is bounded 
by $r^a$ as $r\to+\infty$. We infer
\beaa
\bsplit
         (\nu\div\b)_{\ell=1}  &=  O(r^{-1})(\div \b)_{\ell=1}+(\lap\rho)_{\ell=1} +(1+O(r^{-1}))(\div\div\a)_{\ell=1}\\
         & +O(r^{-3})(\div\eta)_{\ell=1}+r^{-2}\dkb^{\leq 1}(\Ga_b\c\Ga_g),\\
    (\nu\curl\b)_{\ell=1}  &=  \frac{8}{r}(1+ O(r^{-1}))(\curl \b)_{\ell=1} -(\lap\rhod)_{\ell=1} +(1+O(r^{-1}))(\curl\div\a)_{\ell=1} \\
    &+O(r^{-3})(\rhod)_{\ell=1}+r^{-2}\dkb^{\leq 1}(\Ga_b\c\Ga_g).
  \end{split}
\eeaa
Using the fact that $(\div\eta)_{\ell=1}=0$, and integrating by parts $\Delta$, $\div\div$ and $\curl\div$, we obtain 
\beaa
&& |(\nu\div\b)_{\ell=1}|+|(\nu\curl\b)_{\ell=1}|\\
&\les& r^{-1}|(\div \b)_{\ell=1}|+r^{-1}|(\curl\b)_{\ell=1}|+r^{-2}\left|\left(\rhoc-\frac{1}{2}\chih\c\chibh\right)_{\ell=1}\right|+r^{-2}\left|\left(\rhod-\frac{1}{2}\chih\wedge\chibh\right)_{\ell=1}\right|\\&&+r^{-1}\left(\left|\left(\De+\frac{2}{r^2}\right)\Jp\right|+\left|\dds_2\dds_1\Jp\right|\right)|\Ga_g|+r^{-2}|\dkb^{\leq 1}(\Ga_b\c\Ga_g)|.
\eeaa
Together with the above estimates, we infer
\beaa
 |(\nu\div\b)_{\ell=1}|+|(\nu\curl\b)_{\ell=1}| &\les& \frac{1}{r^6}+\frac{\ep_0}{r^5u^{\frac{3}{2}+2\dec}}.
\eeaa
Using the dominance condition \eqref{eq:dominantconditiononronSigmastarchap5} on $r$ on $\Si_*$, we deduce
\bea
 |(\nu\div\b)_{\ell=1}|+|(\nu\curl\b)_{\ell=1}| &\les& \frac{\ep_0}{r^5u^{1+\dec}}.
\eea
This concludes the proof of Proposition \ref{prop:controlofell=1modesonSigmastar}.
\end{proof}

%%%%%%%%%%%%%%%%%%%%%%%%%%%%%%%%

\subsection {Estimates for $\ell=0$ modes on $\Si_*$} 

%%%%%%%%%%%%%%%%%%%%%%%%%%%%%%%%

In this section, we control the average (i.e. the $\ell=0$ mode) of $\kabc$, $\rhoc$, $\rhod$ and $\muc$. 
Recall the definition of the Hawking mass
\beaa
\frac{2m_H}{r} &=& 1+\frac{1}{16\pi}\int_S\ka\kab.
\eeaa
In order to control $\ell=0$ modes on $\Si_*$, we will need in particular to compare the Hawking mass $m_H$ with the constant $m$. To this end, we will rely on the following lemma.

\begin{lemma}\lab{lemma:transportequationforHawkingmass}
We have
\bea
\ov{\rho}= -\frac{2m_H}{r^3}+\Ga_b\c\Ga_g
\eea
and
\bea
\nu(m_H) &=& r^2\dkb^{\leq 1}(\Ga_b\c\Ga_b).
\eea
\end{lemma}

\begin{proof}
We start with the identity for the average of $\rho$. Recall the Gauss equation
\beaa
K &=& -\rho -\frac{1}{4}\trch\trchb +\frac{1}{2}\chih\c\chibh.
\eeaa
Integrating on $S$, and using the definition of the Hawking mass $m_H$, we obtain
\beaa
\int_SK &=& -\int_S\rho -4\pi\left(\frac{2m_H}{r}-1\right) +\frac{1}{2}\int_S\chih\c\chibh.
\eeaa
Since from Gauss Bonnet we have
\beaa
\int_SK &=& 4\pi,
\eeaa
we infer
\beaa
\int_S\rho &=&  -\frac{8\pi m_H}{r} +\frac{1}{2}\int_S\chih\c\chibh
\eeaa
and hence
\beaa
\ov{\rho}= -\frac{2m_H}{r^3}+\Ga_b\c\Ga_g
\eeaa
as stated. Note that this implies 
\beaa
\rho+\frac{2m_H}{r^3} &=& \left(\ov{\rho}+\frac{2m_H}{r^3}\right)+\rho-\ov{\rho}= \left(\ov{\rho}+\frac{2m_H}{r^3}\right)+\left(\rho+\frac{2m}{r^3}\right)-\ov{\rho+\frac{2m}{r^3}}\\
&=& r^{-1}\Ga_g+\Ga_b\c\Ga_g=r^{-1}\Ga_g
\eeaa
so that $\rho+\frac{2m_H}{r^3}\in r^{-1}\Ga_g$.

Next, we focus on the identity for $\nu(m_H)$. From the null structure equations, we have on $\Si_*$
\beaa
e_3(\trch\trchb) &=& \trch\left(-\frac 1 2 \trchb^2 - 2\omb \trchb + 2\div\xib   +  2 \xib\c(\eta-3\ze)-|\chibh|^2\right)\\
&&+\trchb\left( -\frac 1 2 \trchb\trch  + 2 \omb \trch
+      2   \div \eta + 2|\eta|^2+ 2\rho -\chibh\c\chih\right)\\
&=& -\trch\trchb^2 +2\trchb\rho+2\trch\div\xib +2\trchb\div\eta\\
&& +\trch\left(  2 \xib\c(\eta-3\ze)-|\chibh|^2\right)+\trchb\left(  2|\eta|^2 -\chibh\c\chih\right)
\eeaa
and
\beaa
e_4(\trch\trchb) &=& \trch\left(-\frac 1 2 \trch\trchb    -  2   \div \ze + 2|\ze|^2+2\rho -\chih\c\chibh\right)+\trchb\left(-\frac 1 2\trch^2 -|\chih|^2\right)\\
&=& -\trch^2\trchb +2\trch\rho -2\trch\div\ze  +\trch\left( 2|\ze|^2 -\chih\c\chibh\right)-\trchb|\chih|^2.
\eeaa
Hence, we obtain
\beaa
\nu(\trch\trchb) &=& e_3(\trch\trchb)+b_* e_4(\trch\trchb)\\
&=&  -\trch\trchb(\trchb+b_*\trch) +2(\trchb+b_*\trch)\rho+2\trch\div\xib +2\trchb\div\eta\\
&& +\trch\left(  2 \xib\c(\eta-3\ze)-|\chibh|^2\right)+\trchb\left(  2|\eta|^2 -\chibh\c\chih\right)\\
&& +b_*\Big\{   -2\trch\div\ze  +\trch\left( 2|\ze|^2 -\chih\c\chibh\right)-\trchb|\chih|^2\Big\},
\eeaa
which we rewrite
\beaa
\nu(\trch\trchb)+\trch\trchb(\trchb+b_*\trch) &=&   2(\trchb+b_*\trch)\rho+\frac{4}{r}\div\xib -\frac{4\Up}{r}\div\eta\\
&& +\frac{4\left(1+\frac{2m}{r}\right)}{r}\div\ze+r^{-1}\dkb^{\leq 1}(\Ga_b\c\Ga_b).
\eeaa
Together with Lemma \ref{Lemma:nuSof-integrals-again}, we infer
\beaa
\nu\left(\int_S\trch\trchb\right) &=&  z\int_S\frac{1}{z}\Big(\nu(\trch\trchb)+(\kab+b_*\ka)\trch\trchb     \Big)\\
&=& z\int_S\frac{1}{z}\Bigg\{2(\trchb+b_*\trch)\rho+\frac{4}{r}\div\xib -\frac{4\Up}{r}\div\eta\\
&& +\frac{4\left(1+\frac{2m}{r}\right)}{r}\div\ze+r^{-1}\dkb^{\leq 1}(\Ga_b\c\Ga_b\Bigg\}.
\eeaa
Integrating by parts the divergences, we deduce
\beaa
\nu\left(\int_S\trch\trchb\right) &=& 2z\int_S\frac{1}{z}(\trchb+b_*\trch)\rho +r\dkb^{\leq 1}(\Ga_b\c\Ga_b).
\eeaa
Thus, in view of the definition of the Hawking mass, 
\beaa
\frac{2m_H}{r} &=& 1+\frac{1}{16\pi}\int_S\trch\trchb,
\eeaa
we infer
\beaa
\nu\left(\frac{2m_H}{r}\right) &=& \frac{z}{8\pi}\int_S\frac{1}{z}(\trchb+b_*\trch)\rho +r\dkb^{\leq 1}(\Ga_b\c\Ga_b).
\eeaa
On the other hand, we have 
\beaa
\nu\left(\frac{2m_H}{r}\right) &=& \frac{2\nu(m_H)}{r}-\frac{2m_H\nu(r)}{r^2}
\eeaa
and hence, using again Lemma \ref{Lemma:nuSof-integrals-again}, 
\beaa
\nu\left(\frac{2m_H}{r}\right) &=& \frac{2\nu(m_H)}{r}-\frac{2m_H}{r^2}\frac{z}{8\pi r}\int_S\frac{1}{z}(\trchb+b_*\trch)
\eeaa
which yields
\beaa
\nu(m_H) &=& \frac{rz}{16\pi}\int_S\frac{1}{z}(\trchb+b_*\trch)\left(\rho+\frac{2m_H}{r^3}\right) +r^2\dkb^{\leq 1}(\Ga_b\c\Ga_b).
\eeaa
Recalling from above that we have $\rho+\frac{2m_H}{r^3}\in r^{-1}\Ga_g$, we infer
\beaa
\nu(m_H) &=& -\frac{1}{4\pi}\int_S\left(\rho+\frac{2m_H}{r^3}\right) +r^2\dkb^{\leq 1}(\Ga_b\c\Ga_b).
\eeaa
Together with the above control of the average of $\rho+\frac{2m_H}{r^3}$, we deduce
\beaa
\nu(m_H) &=& r^2\dkb^{\leq 1}(\Ga_b\c\Ga_b)
\eeaa
as desired. This concludes the proof of Lemma \ref{lemma:transportequationforHawkingmass}.
\end{proof}

\begin{proposition}[Control of $\ell=0$ modes on $\Si_*$]\lab{prop:controlofell=0modesonSigmastar}
We have on $\Si_*$
\bea
\sup_{\Si_*}u^{1+2\dec}\left(|m_H-m|+r^2\left|\ov{\kabc}\right|+r^3\left|\ov{\rhoc}\right|+r^3\left|\ov{\muc}\right|+r^3\left|\ov{\rhod}\right|\right) &\les& \ep_0.
\eea
\end{proposition}

\begin{proof}
Recall from Lemma \ref{lemma:transportequationforHawkingmass} that we have
\beaa
\nu(m_H) &=& r^2\dkb^{\leq 1}(\Ga_b\c\Ga_b).
\eeaa
Together with the control of $\Ga_b$ provided by  {\bf Ref 1}, and since $m$ is a constant, we infer
\beaa
|\nu(m_H-m)| &\les& \frac{\ep^2}{u^{2+2\dec}}\les \frac{\ep_0}{u^{2+2\dec}}. 
\eeaa
Also, recall that by definition of $m$, we have $m=m_H$ on $S_*$.  Integrating from $S_*$, we deduce on $\Si_*$
\beaa
|m_H-m| &\les& \frac{\ep_0}{u^{1+2\dec}} 
\eeaa
as desired.

Next, recall  from Lemma \ref{lemma:transportequationforHawkingmass} that we have
\beaa
\ov{\rho}= -\frac{2m_H}{r^3}+\Ga_b\c\Ga_g.
\eeaa
Together with the above control for $m_H-m$, and the control of $\Ga_b$ and $\Ga_g$ provided by  {\bf Ref 1}, we deduce
\beaa
\left|\ov{\rhoc}\right|=\left|\ov{\rho}+\frac{2m}{r^3}\right| &\les& \frac{\ep_0}{r^3u^{1+2\dec}} 
\eeaa
as desired.

Next, taking the average of 
\beaa
\curl\ze &=& \rhod -\frac{1}{2}\chih\wedge\chibh,
\eeaa
we infer
\beaa
\ov{\rhod} &=& \frac{1}{2}\ov{\chih\wedge\chibh} =\Ga_b\c\Ga_g,
\eeaa
and the conclusion follows from the control of $\Ga_b$ and $\Ga_g$ provided by  {\bf Ref 1}. 

Next, using the definition of the Hawking mass and the GCM condition for $\ka$ on $\Si_*$, we have
\beaa
\frac{2m_H}{r} &=& 1+\frac{1}{16\pi}\int_S\ka\kab=1+\frac{1}{8\pi r}\int_S\kab=1+\frac{r}{2}\ov{\kab}
\eeaa
and hence
\beaa
\ov{\kab} &=& -\frac{2}{r}\left(1-\frac{2m_H}{r}\right)=-\frac{2\Up}{r}+\frac{4}{r^2}(m_H-m).
\eeaa
Together with the above control for $m_H-m$, we deduce
\beaa
\left|\ov{\kab}+\frac{2\Up}{r}\right| &\les& \frac{\ep_0}{r^2u^{1+2\dec}} 
\eeaa
as desired.

Finally, we consider $\ov{\mu}$. We have by definition of $\mu$, and in view of Gauss equation, 
\beaa
\mu &=& -\div\ze -\rho +\frac{1}{2}\chih\c\chibh\\
&=& -\div\ze+K+\frac{1}{4}\trch\trchb.
\eeaa 
Integrating, and using integration by parts, we obtain 
\beaa
\int_S\mu &=& \int_SK+\frac{1}{4}\int_S\trch\trchb.
\eeaa 
Using Gauss Bonnet and the definition of the Hawking mass, we deduce
\beaa
\int_S\mu &=& 4\pi +4\pi\left(\frac{2m_H}{r}-1\right)=\frac{8\pi m_H}{r}\\
&=& \frac{8\pi m}{r}+\frac{8\pi (m_H-m)}{r}.
\eeaa
Together with the above control for $m_H-m$, we deduce
\beaa
\left|\ov{\mu}-\frac{2m}{r^3}\right| &\les& \frac{\ep_0}{r^3u^{1+2\dec}} 
\eeaa
as desired. This concludes the proof of the proposition.
\end{proof}

%%%%%%%%%%%%%%%%%%%%%%%%%%%%%%%

\section{Proof of Theorem M3}
\lab{sec::improvementofdecaybootassonSigmastar}

%%%%%%%%%%%%%%%%%%%%%%%%%%%%%%%

We are now ready to prove Theorem M3.  
\begin{proposition}\lab{prop:decayonSigamstarofallquantities}
We have along $\Si_*$, for all $k\le k_*-12$,
\bea\lab{eq:improvementofdecaybootassonSigmastar}
\bsplit
\big|\dk_*^{\leq k }\Ga_b\big| &\les  \ep_0 r^{-1} u^{-1-\dec},\\
\big|\dk_*^{\leq k }\Ga_g\big| &\les  \ep_0 r^{-2} u^{-\frac{1}{2} -\dec},\\
\big|\dk_*^{\leq k-1 }\nab_\nu\Ga_g\big| &\les  \ep_0 r^{-2}u^{-1-\dec}.
\end{split}
\eea
Moreover, for all $k\le k_*-12$,
 \bea
 \bsplit
 \big|\dk_*^{\leq k }\kac\big| &\les  \ep_0 r^{-2} u^{-1-\dec},\\
 \big|\dk_*^{\leq k }\muc\big| &\les  \ep_0 r^{-3} u^{-1-\dec},\\
 \big|\dk_*^{\leq k-1}\nab_\nu\b| &\les   \ep_0 r^{-4}u^{-\frac{1}{2}-\dec}. 
 \end{split}
 \eea
\end{proposition}

\begin{remark}
Note that \eqref{eq:improvementofdecaybootassonSigmastar} yields the proof of Theorem M3. Indeed, in view of the definition of the decay norm $\,^{*}\mathfrak{D}_k$ in section \ref{section:main-normsSigmastar},  \eqref{eq:improvementofdecaybootassonSigmastar} can be rewritten as $\,^{*}\mathfrak{D}_k\les \ep_0$ for $k\le k_*-12$. Thus, since $k_*=k_{small}+80$ in view of \eqref{eq:valueofkstarinchapter5forproofThmM3},  Proposition \ref{prop:decayonSigamstarofallquantities} yields in particular $\,^{*}\mathfrak{D}_{k_{small}+60}\les \ep_0$ and thus concludes the proof of Theorem M3.
\end{remark}

\begin{proof}
Note that the estimates for $\Ga_b$ have already been established  in Proposition \ref{Prop.Flux-bb-vthb-eta-xib}. Note also that $\kac=0$ in view of our GCM conditions, and that the estimate for $\a$ has already been established in Theorem M2.  Thus, it only remains to control the following quantities
\beaa
\kabc,\quad  \chih, \quad  \ze,  \quad  \rhoc, \quad \rhod,\quad  \muc,\quad \b.
\eeaa
We control these quantities as follows, starting first with estimates for angular derivatives. 

{\bf Step 1.} We start with $\dkb^k\kabc$ and $\dkb^k\muc$. Recall from our GCM conditions that we have on $\Si_*$
\beaa
\kabc=\underline{C}_0+\sum_p\underline{C}_p\Jp,\qquad \muc=M_0+\sum_pM_p\Jp.
\eeaa
Differentiating w.r.t. $\dds_2\dds_1$, and recalling that $\underline{C}_0$, $\underline{C}_p$, $M_0$ and $M_p$ are constant on the spheres $S$, we infer
\beaa
\dds_2\dds_1\kabc=\sum_p\underline{C}_p\dds_2\dds_1\Jp,\qquad \dds_2\dds_1\muc=\sum_pM_p\dds_2\dds_1\Jp,
\eeaa
which yields, for $k\geq 0$, 
\beaa
\|\dds_2\dds_1\kabc\|_{\hk_k(S)} &\les& r\sum_p|\underline{C}_p|\|\dkb^{\leq k}\dds_2\dds_1\Jp\|_{L^\infty(S)},\\ 
\|\dds_2\dds_1\muc\|_{\hk_k(S)} &\les& r\sum_p|M_p|\|\dkb^{\leq k}\dds_2\dds_1\Jp\|_{L^\infty(S)}.
\eeaa
Together with Corollary \ref{prop:2D-Hodge4}, we infer
\beaa
\|\kabc\|_{\hk_{k+2}(S)} &\les& r^3\sum_p|\underline{C}_p|\|\dkb^{\leq k}\dds_2\dds_1\Jp\|_{L^\infty(S)}+r|(\kabc)_{\ell=1}|+r|\ov{\kabc}|,\\ 
\|\muc\|_{\hk_{k+2}(S)} &\les& r^3\sum_p|M_p|\|\dkb^{\leq k}\dds_2\dds_1\Jp\|_{L^\infty(S)}+r|(\muc)_{\ell=1}|+r|\ov{\muc}|.
\eeaa
In view of the control of the $\ell=1$ mode of $\kabc$ and $\muc$ in Proposition \ref{prop:control.ell=1modes-Si}, the control of the average  of $\kabc$ and $\muc$ in Proposition \ref{prop:controlofell=0modesonSigmastar}, the fact that $\underline{C}_p\in \Ga_g$ and $M_p\in r^{-1}\Ga_g$ in view of Corollary  \ref{cor:Cb0CbpM0MareGagandrm1Gag}, the control of $\Ga_g$ provided by {\bf Ref 1}, and the control of $\dds_2\dds_1\Jp$ provided by Corollary \ref{corr:Si*-ell=1modes-improved}, we deduce, for $k\leq k_*-10$, 
\beaa
\|\kabc\|_{\hk_{k+2}(S)} \les \frac{\ep_0}{ru^{1+\dec}},\qquad \|\muc\|_{\hk_{k+2}(S)} \les \frac{\ep_0}{r^2u^{1+\dec}}.
\eeaa
Together with Sobolev, this implies, for $k\leq k_*-10$, 
\bea
|\dkb^k\kabc|\les \frac{\ep_0}{r^2u^{1+\dec}}, \qquad |\dkb^k\muc|\les \frac{\ep_0}{r^3u^{1+\dec}}.
\eea

{\bf Step 2.} Next, we focus on $\dkb^k\rhoc$ and $\dkb^k\rhod$. Recall from Proposition \ref{prop:identitiesinqf} that we have
\beaa
\Re(\qf) &=&  r^4\dds_2\dds_1(-\rho, \rhod)  +O(r^{-2})+\dkb^{\leq 2}\Ga_b+r^{2}\dkb^{\le 2}\big(\Ga_b \c  \Ga_g\big).
\eeaa
 We infer, for $k\geq 0$, 
\beaa
\|\dds_2\dds_1(-\rhoc, \rhod)\|_{\hk_k(S)} &\les& r^{-3}\|\dkb^{\leq k}\qf\|_{L^\infty(S)}+r^{-5}+r^{-3}\|\dkb^{\leq k+2}\Ga_b\|_{L^\infty(S)}\\
&&+r^{-1}\|\dkb^{\leq k+2}(\Ga_b\c\Ga_g)\|_{L^\infty(S)}.
\eeaa
In view of {\bf Ref 2} for $\qf$, the control of $\Ga_b$  established  in Proposition \ref{Prop.Flux-bb-vthb-eta-xib}, and the control of $\Ga_g$ provided by {\bf Ref 1}, we obtain, for $k\leq k_*-12$,
\beaa
\|\dds_2\dds_1(-\rhoc, \rhod)\|_{\hk_k(S)} &\les& \frac{1}{r^5}+\frac{\ep_0}{r^4u^{\frac{1}{2}+\dec}}.
\eeaa
Together with the dominance condition \eqref{eq:dominantconditiononronSigmastarchap5} on $r$ on $\Si_*$, this yields, for $k\leq k_*-12$,
\beaa
\|\dds_2\dds_1(-\rhoc, \rhod)\|_{\hk_k(S)} &\les& \frac{\ep_0}{r^4u^{\frac{1}{2}+\dec}}.
\eeaa
In view Corollary \ref{prop:2D-Hodge4}, we deduce, for $k\leq k_*-12$, 
\beaa
\|(-\rhoc, \rhod)\|_{\hk_{k+2}(S)} &\les& \frac{\ep_0}{r^2u^{\frac{1}{2}+\dec}}+r|(\rhoc)_{\ell=1}|+r|(\rhod)_{\ell=1}|+r|\ov{\rhoc}|+r|\ov{\rhod}|.
\eeaa
In view of the control of the $\ell=1$ mode of $\rhoc$ and $\rhod$ in Proposition \ref{prop:control.ell=1modes-Si}, and the control of the average  of $\rhoc$ and $\rhod$ in Proposition \ref{prop:controlofell=0modesonSigmastar},  we infer, for $k\leq k_*-12$, 
\beaa
\|\rhoc\|_{\hk_{k+2}(S)} \les \frac{\ep_0}{r^2u^{\frac{1}{2}+\dec}},\qquad \|\rhod\|_{\hk_{k+2}(S)} \les \frac{\ep_0}{r^2u^{\frac{1}{2}+\dec}}.
\eeaa
Together with Sobolev, this implies, for $k\leq k_*-12$,
\bea
|\dkb^k\rhoc|\les \frac{\ep_0}{r^3u^{\frac{1}{2}+\dec}}, \qquad |\dkb^k\rhod|\les \frac{\ep_0}{r^3u^{\frac{1}{2}+\dec}}.
\eea

{\bf Step 3.} Next, we focus on $\dkb^k\ze$. From the definition of $\mu$, and the null structure equation for $\curl\ze$, we have
\beaa
\ddd_1\ze &=& \left(-\muc-\rhoc+\frac{1}{2}\chih\c\chibh,\rhod -\frac{1}{2}\chih\c\chibh \right)\\
&=& \left(-\muc-\rhoc,\rhod\right)+\Ga_b\c\Ga_g.
\eeaa
In view of Lemma \ref{prop:2D-Hodge1}, we infer, for $k\geq 0$, 
\beaa
 \|\ze\|_{\hk_{k+1}   (S)} &\les&  r\|\muc\|_{\hk_k(S)}+r\|\rhoc\|_{\hk_k(S)}+r\|\rhod\|_{\hk_k(S)}+r^2\|\dkb^{\leq k}(\Ga_b\c\Ga_g)\|_{L^\infty(S)}.
\eeaa
Together with the control of $\muc$ derived in Step 1, the control of $(\rhoc, \rhod)$ derived in Step 2, the control of $\Ga_b$  established  in Proposition \ref{Prop.Flux-bb-vthb-eta-xib}, and the control of $\Ga_g$ provided by {\bf Ref 1}, we obtain, for $k\leq k_*-10$,
\beaa
 \|\ze\|_{\hk_{k+1}   (S)} &\les&  \frac{\ep_0}{ru^{\frac{1}{2}+\dec}}.
\eeaa
Together with Sobolev, this implies, for $k\leq k_*-11$,
\bea
|\dkb^k\ze|\les \frac{\ep_0}{r^2u^{\frac{1}{2}+\dec}}.
\eea

{\bf Step 4.} Next, we focus on $\dkb^k\b$. Recall the following consequence of Bianchi
\beaa
\dds_2\b &=& \nab_3\a+O(r^{-1})\a+O(r^{-3})\Ga_g+r^{-1}\Ga_b\c\Ga_g.
\eeaa
In view of Lemma \ref{prop:2D-Hodge2}, we infer, for $k\geq 0$, 
\beaa
  \|\b\|_{\hk_{k+1}  (S)}&\les&  r^2\|\dkb^{\leq k}\nab_3\a\|_{L^\infty(S)}+r\|\dkb^{\leq k}\a\|_{L^\infty(S)}+r^{-1}\|\dkb^{\leq k}\Ga_g\|_{L^\infty(S)}\\
  &&+r\|\dkb^{\leq k}(\Ga_b\c\Ga_g)\|_{L^\infty(S)}+r^2 \big| (\ddd_1\b)_{\ell=1}\big|.
 \eeaa
 Together with the control of $\a$ and $\nab_3\a$ provided by {\bf Ref 2}, the control of $\Ga_b$  established  in Proposition \ref{Prop.Flux-bb-vthb-eta-xib}, the control of $\Ga_g$ provided by Lemma \ref{lemma:interpolation}, and the control of $(\ddd_1\b)_{\ell=1}$ in Proposition \ref{prop:control.ell=1modes-Si}, we obtain, for $k\leq k_*$,
\beaa
  \|\b\|_{\hk_{k+1}  (S)}&\les&  \frac{\ep_0}{r^2u^{\frac{1}{2}+\dec}}+\frac{\ep_0}{r^{\frac{5}{2}+\dec}}+\frac{\ep}{r^3u^{\frac{1}{2}+\frac{\dec}{2}}}.
 \eeaa
 In view of the dominance condition \eqref{eq:dominantconditiononronSigmastarchap5} on $r$ on $\Si_*$, this yields, for $k\leq k_*$,
\beaa
  \|\b\|_{\hk_{k+1}  (S)}&\les&  \frac{\ep_0}{r^2u^{\frac{1}{2}+\dec}}.
 \eeaa 
Together with Sobolev, this implies, for $k\leq k_*-1$,
\bea
|\dkb^k\b|\les \frac{\ep_0}{r^3u^{\frac{1}{2}+\dec}}.
\eea

{\bf Step 5.} Next, we focus on $\dkb^k\chih$. Recall the following consequence of Codazzi
\beaa
\ddd_2\chih &=& \frac{1}{r}\ze -\b+\Ga_g\c\Ga_g.
\eeaa
In view of Lemma \ref{prop:2D-Hodge1}, we infer, for $k\geq 0$, 
\beaa
 \|\chih\|_{\hk_{k+1}   (S)} &\les&  \|\ze\|_{\hk_k(S)}+r\|\b\|_{\hk_k(S)}+r^2\|\dkb^{\leq k}(\Ga_g\c\Ga_g)\|_{L^\infty(S)}.
\eeaa
Together with the control of $\ze$ derived in Step 3, the control of $\b$ derived in Step 4, the control of $\Ga_b$  established  in Proposition \ref{Prop.Flux-bb-vthb-eta-xib}, and the control of $\Ga_g$ provided by {\bf Ref 1}, we obtain, for $k\leq k_*-10$,
\beaa
 \|\chih\|_{\hk_{k+1}   (S)} &\les&  \frac{\ep_0}{ru^{\frac{1}{2}+\dec}}.
\eeaa
Together with Sobolev, this implies, for $k\leq k_*-11$,
\bea
|\dkb^k\chih|\les \frac{\ep_0}{r^2u^{\frac{1}{2}+\dec}}.
\eea

{\bf Step 6.} In view of Steps 1--5, and the fact that $\kac=0$ in view of our GCM conditions, and that the estimate for $\a$ has already been established in Theorem M2, we have obtained, for $k\leq k_*-12$,
\bea
|\dkb^k\Ga_g|\les \frac{\ep_0}{r^2u^{\frac{1}{2}+\dec}}, \qquad |\dkb^{\leq k }\kac| \les  \frac{\ep_0}{r^2u^{1+\dec}},\qquad |\dkb^{\leq k }\muc| \les  \frac{\ep_0}{r^3u^{1+\dec}}.
\eea

{\bf Step 7.}  Next, we estimate $\nab_\nu\Ga_g$. From the null structure equations and Bianchi identities, one observes that all quantities in $\Ga_g$ verify schematically 
\beaa
\nab_\nu\Ga_g &=& r^{-1}\dkb^{\leq 1}\Ga_b+ r^{-1} \Ga_g+\Ga_b\c \Ga_b.
\eeaa
Together with the control of $\Ga_b$  established  in Proposition \ref{Prop.Flux-bb-vthb-eta-xib}, and the control of $\Ga_g$ provided by {\bf Ref 1}, we infer, for $k\leq k_*-10$,
\beaa
\big|\dk_*^{\leq k-1 }\nab_\nu\Ga_g\big| &\les&  \frac{\ep_0}{r^2u^{1+\dec}}+\frac{1}{r^3}.
\eeaa
In view of the dominance condition \eqref{eq:dominantconditiononronSigmastarchap5} on $r$ on $\Si_*$, this yields, for $k\leq k_*-10$,
\beaa
\big|\dk_*^{\leq k-1 }\nab_\nu\Ga_g\big| &\les&  \frac{\ep_0}{r^2u^{1+\dec}}.
\eeaa
Together with the estimates of Step 6, and the control of $\Ga_b$  established  in Proposition \ref{Prop.Flux-bb-vthb-eta-xib}, we deduce, for $k\leq k_*-12$,
\bea
\bsplit
\big|\dk_*^{\leq k }\Ga_b\big| &\les  \ep_0 r^{-1} u^{-1-\dec},\qquad \big|\dk_*^{\leq k }\Ga_g\big| &\les  \ep_0 r^{-2} u^{-\frac{1}{2} -\dec},\\
\big|\dk_*^{\leq k-1 }\nab_\nu\Ga_g\big| &\les  \ep_0 r^{-2}u^{-1-\dec}, \qquad \big|\dk_*^{\leq k }\kac\big| &\les  \ep_0 r^{-2} u^{-1-\dec},\\
\big|\dk_*^{\leq k }\muc\big| &\les  \ep_0 r^{-3} u^{-1-\dec}.
\end{split}
\eea

{\bf Step 8.} We conclude the proof with an estimate for $\nab_\nu\b$. In view of the Bianchi identities for $\nab_3\b$ and $\nab_4\b$, and using the fact that $\nu=e_3+b_*e_4$, we have
\beaa
\nab_\nu\b &=& r^{-2}\dkb^{\leq 1}\Ga_g.
\eeaa
Together with the estimate for $\Ga_g$ derived in Step 7, we infer, for $k\leq k_*-12$,
\bea
\big|\dk_*^{\leq k-1}\nab_\nu\b| &\les&   \ep_0 r^{-4}u^{-\frac{1}{2}-\dec}. 
\eea
This  concludes the proof of Proposition \ref{prop:decayonSigamstarofallquantities}.
\end{proof} 

We conclude this section with the following non sharp corollary of Proposition \ref{prop:decayonSigamstarofallquantities}, Proposition \ref{prop:controlofell=1modesonSigmastar} and Corollary \ref{corofLemma:transport.alongSi_*1} that will be useful in section \ref{sec:proofofThmM7}.
\begin{corollary}\lab{corofLemma:transport.alongSi_*1:again}
We have along $\Si_*$
\beaa
\left|\nu((\div\b)_{\ell=1})\right|+\left|\nu((\curl\b)_{\ell=1,\pm})\right|+\left|\nu\left((\curl\b)_{\ell=1,0}-\frac{2am}{r^5}\right)\right| &\les& \frac{\ep_0}{r^5u^{1+\dec}},\\
\left|\nu\big((\kabc)_{\ell=1}\big)\right| &\les& \frac{\ep_0^2}{r^2u^{2+2\dec}}+\frac{\ep_0}{r^3u^{1+\dec}}.
\eeaa
\end{corollary}

\begin{proof}
Recall the following identities derived in Corollary \ref{corofLemma:transport.alongSi_*1} along $\Si_*$, for $p=0,+,-$,
\beaa
\nn&&\nu\left(\int_S\left( \lap\kabc+\frac{2\Up}{r}\div \ze\right)\Jp\right)\\
\nn&=& O(r^{-3})\int_S\kabc\Jp   +O(r^{-2})\int_S\div\ze\Jp +O(r^{-1})\int_S\div\bb\Jp     +O(r^{-2})\int_S\rhoc\Jp \\
 &&+O(r^{-1})\int_S\div\b\Jp +r\left|\left(\Delta+\frac{2}{r^2}\right)\Jp\right|\dkb^{\leq 1}\Ga_b+\dkb^{\le 2 }(\Ga_b\c \Ga_b),
\eeaa
\beaa
\nu\left(\int_S\div\b\Jp\right) &=& O(r^{-1})\int_S\div\b\Jp+O(r^{-2})\int_S\rhoc\Jp\\
\nn&&+r\left(\left|\left(\Delta+\frac{2}{r^2}\right)\Jp\right|+\left|\dds_2\dds_1\Jp\right|\right)\Ga_g+\dkb^{\leq 1}(\Ga_b\c\Ga_g),
\eeaa
and
\beaa
\nn\nu\left(\int_S\curl\b\Jp\right) &=&  \frac{4}{r}(1+ O(r^{-1}))\int_S\curl\b\Jp+\frac{2}{r^2}(1+ O(r^{-1}))\int_S\rhod\Jp\\
\nn&&+r\left(\left|\left(\Delta+\frac{2}{r^2}\right)\Jp\right|+\left|\dds_2\dds_1\Jp\right|\right)\Ga_g\\
&&+\dkb^{\leq 1}(\Ga_b\c\Ga_g).
\eeaa
The estimates then easily follow  from  the control of the $\ell=1$ modes in Proposition \ref{prop:controlofell=1modesonSigmastar}, the control of $\Ga_g$ and $\Ga_b$ provided by Proposition \ref{prop:decayonSigamstarofallquantities},  the control of $(\De+\frac{2}{r^2})\Jp$ and $\dds_2\dds_1\Jp$ provided by Corollary  \ref{corr:Si*-ell=1modes-improved} (with $\ep_0$ smallness constant instead of $\ep$ thanks to Proposition \ref{prop:decayonSigamstarofallquantities}), and the dominance condition for $r$ on $\Si_*$. 
\end{proof}

%%%%%%%%%%%%%%%%%%%%%%%%%%%%%%%%%%%%%%%%%%%%%%%%%%

\section{Control of $\Jp$ and $\Jk$ on $\Si_*$}
\lab{sec:controlofJpandJkonSigmastar}

%%%%%%%%%%%%%%%%%%%%%%%%%%%%%%%%%%%%%%%%%%%%%%%%%%

Recall that the  induced metric $g$ on $S_*$  takes the form
 \beaa
 g= r^2e^{2\phi}\Big( (d\th)^2+ \sin^2 \th (d\vphi)^2\Big).
 \eeaa
For the constructions in Definition \ref{def:definitionoff0fplusfminus} below, we will rely on a special orthonormal basis $(e_1, e_2)$ of the tangent space of $S_*$ given by
\bea\lab{eq:specialorthonormalbasisofSstar}
e_1=\frac{1}{re^\phi}\pr_\th, \qquad e_2=\frac{1}{r\sin\th e^\phi}\pr_\vphi, \quad\textrm{on}\quad S_*.
\eea

To control the regularity of the basis of of $\ell=1$ modes $\Jp$, $p=0,+,-$, we introduce the following 1-forms.  

\begin{definition}\lab{def:definitionoff0fplusfminus}
Let $f_0$, $f_+$ and $f_-$ the 1-forms defined on $S_*$ by:
\bea
\bsplit
&(f_0)_1 =0,  \qquad (f_0)_2 =\sin\th, \qquad (f_+)_1 =\cos\th\cos\vphi,  \qquad (f_+)_2 =-\sin\vphi,\\
& (f_-)_1 =\cos\th\sin\vphi,  \qquad (f_-)_2=\cos\vphi, \quad \textrm{on}\quad S_*,
\end{split}
\eea
in the orthonormal basis $(e_1, e_2)$ of $S_*$ given by \eqref{eq:specialorthonormalbasisofSstar}, and extended to $\Si_*$ by:
\bea
\nab_\nu f_0=0, \qquad \nab_\nu f_+=0, \qquad \nab_\nu f_-=0. 
\eea
\end{definition}

This allows us to renormalize $\nab(J^{(p)})$, $p=0,+,-$ on $\Si_*$ as follows.
\begin{definition}\lab{def:renormalizationforJpbasisell=1modes}
We introduce the notations
\beaa
\widecheck{\nab J^{(0)}}:=\nab J^{(0)}+\frac{1}{r}\dual f_0, \qquad \widecheck{\nab J^{(+)}}:=\nab J^{(+)}-\frac{1}{r}f_+, \qquad \widecheck{\nab J^{(-)}}:=\nab J^{(-)}-\frac{1}{r}f_-.
\eeaa
\end{definition}

We also introduce the following renormalization for angular derivatives of $f_0$ and $f_\pm$.
\begin{definition}\lab{def:renormalizationforf0fpfm}
We introduce the notations
\beaa
\widecheck{\curl(f_0)}:=\curl(f_0)-\frac{2}{r}\cos\th, \qquad \widecheck{\div(f_\pm)}:=\div(f_\pm)+\frac{2}{r}J^{(\pm)}. 
\eeaa
\end{definition} 

Finally, note that the complex horizontal 1-form $\Jk$ introduced in Definition \ref{def:definitionofJkonMext} verifies  
\bea\lab{eq:relationbetweenJkandf0onSigmastar}
\Jk = \frac{1}{|q|}\left(f_0+i\dual f_0\right)\quad \textrm{on}\quad\Si_*.
\eea
We also introduce the following two complex horizontal 1-forms $\Jk_\pm$ given by
\bea\lab{eq:relationbetweenJkpmandfpmonSigmastar}
\Jk_\pm := \frac{1}{|q|}\left(f_\pm+i\dual f_\pm\right)\quad \textrm{on}\quad\Si_*,
\eea
as well as the following renormalizations
\bea
\bsplit
\widecheck{\ov{\DD}\c\Jk} &:= \ov{\DD}\c\Jk-\frac{4i(r^2+a^2)\cos\th}{|q|^4},\\ 
\widecheck{\ov{\DD}\c\Jk_\pm}&:=\ov{\DD}\c\Jk_\pm+\frac{4r^2}{|q|^4}J^{(\pm)}+\frac{2ia^2\cos\th}{|q|^4}J^{(\mp)}.
\end{split}
\eea

The goal of this section is to prove the following proposition.
\begin{proposition}\lab{prop:estimatesforf0fplusfminusandJp}
We have on $\Si_*$, for all $k\leq k_*-12$
\beaa
\left|\dk_*^k\Big[\widecheck{\nab\Jp}\Big]\right| &\les& \frac{\ep_0}{r^2u^{\frac{1}{2}+\dec}},
\eeaa
and on $\Si_*$, for all $k\leq k_*-13$
\beaa
\left|\dk_*^k\left[\div(f_0),\, \widecheck{\curl(f_0)},\, \nab\hot f_0,\, \nab f_0 - \frac{1}{r}\cos\th\in\right]\right| \\
+ \left|\dk_*^k\left[\widecheck{\div(f_\pm)},\, \curl(f_\pm),\, \nab\hot f_\pm,\, \nab f_\pm + \frac{1}{r}J^{(\pm)}\de\right]\right| &\les& \frac{\ep_0}{r^2u^{\frac{1}{2}+\dec}},
\eeaa
as well as
\beaa
 \left|\dk_*^k\widecheck{\ov{\DD}\c\Jk}\right|+\left|\dk_*^k\widecheck{\ov{\DD}\c\Jk_\pm}\right|  &\les& \frac{\ep_0}{r^3u^{\frac{1}{2}+\dec}}.
\eeaa
\end{proposition}

%%%%%%%%%%%%%%%%%%%%%%%%%%%%%%%%%%%%%%%%%%%%%%%%%%

\subsection{Control on $S_*$}

%%%%%%%%%%%%%%%%%%%%%%%%%%%%%%%%%%%%%%%%%%%%%%%%%%

First, we derive the following corollary of Lemma \ref{lemma:controloftheconformalfactorphi} and Lemma \ref{lemma:statementeq:DeJp.S_*}.
\begin{corollary}\lab{cor:statementeq:DeJp.S_*:improved}
The following holds on $S_*$:
\bea
\lab{eq:controlofphionSstar:onceagain}
\|\dkb^{\leq k_*-12}\phi\|_{L^\infty(S_*)} &\les& \frac{\ep_0}{ru^{\frac{1}{2}+\dec}},
\eea
and
\bea
\lab{eq:DeJp.S_*:onceagain}
\bsplit
\int_{S_*}J^{(p)} &= 0,\\
\int_{S_*}J^{(p)}J^{(q)} &= \frac{4\pi}{3}r^2\de_{pq}+O\left(\ep_0 ru^{-\frac{1}{2}-\dec}\right),\\
\left\|\dkb^{\leq 2}\left(\Delta +\frac{2}{r^2}\right)\Jp\right\|_{L^\infty(S_*)} & =  O\left(\ep_0 r^{-3} u^{-\frac{1}{2}-\dec}\right).
\end{split}
\eea
\end{corollary}

\begin{proof}
In view of the Gauss equation, we have $\widecheck{K}\in r^{-1}\Ga_g$. Together with the estimate for $\Ga_g$  in Proposition \ref{prop:decayonSigamstarofallquantities}, we infer 
\beaa
\sup_{S_*}\left|\dkb^{\leq k_*-12}\left(K-\frac{1}{r^2}\right)\right| &\les& \frac{\ep_0}{r^3u^{\frac{1}{2}+\dec}}.
\eeaa
The proof follows then immediately from the one of Lemma \ref{lemma:controloftheconformalfactorphi} and Lemma \ref{lemma:statementeq:DeJp.S_*} upon replacing the estimate \eqref{eq:controlofkstarangularderivativesofKonSstar} for $\widecheck{K}$ with the above one.
\end{proof}
 
The following lemma provides identities for first order derivatives of $\Jp$, $p=0,+,-$. 
\begin{lemma}\lab{lemma:computationforfirstorderderivativesJpell=1basis}
We have on $S_*$
\bea
\widecheck{\nab J^{(0)}} = -\frac{1}{r}(e^{-\phi}-1)\dual f_0,\quad \widecheck{\nab J^{(+)}} = \frac{1}{r}(e^{-\phi}-1)f_+,\quad \widecheck{\nab J^{(-)}} = \frac{1}{r}(e^{-\phi}-1)f_-.
\eea
\end{lemma}

\begin{proof}
Let the orthonormal basis $(e_1, e_2)$ of $S_*$ given by \eqref{eq:specialorthonormalbasisofSstar}. We have
\beaa
e_1(J^{(0)}) &=& \frac{1}{re^\phi}\pr_\th(\cos\th)=-\frac{1}{re^\phi}\sin\th,\\
e_2(J^{(0)}) &=& \frac{1}{r\sin\th e^\phi}\pr_\vphi(\cos\th)=0,\\
e_1(J^{(+)}) &=& \frac{1}{re^\phi}\pr_\th(\sin\th\cos\vphi)=\frac{1}{re^\phi}\cos\th\cos\vphi,\\
e_2(J^{(+)}) &=& \frac{1}{r\sin\th e^\phi}\pr_\vphi(\sin\th\cos\vphi)=-\frac{1}{re^\phi}\sin\vphi,\\
e_1(J^{(-)}) &=& \frac{1}{re^\phi}\pr_\th(\sin\th\sin\vphi)=\frac{1}{re^\phi}\cos\th\sin\vphi,\\
e_2(J^{(-)}) &=& \frac{1}{r\sin\th e^\phi}\pr_\vphi(\sin\th\sin\vphi)=\frac{1}{re^\phi}\cos\vphi.
\eeaa
Together with the definition of $f_0$ and $f_\pm$, we infer
\beaa
\nab J^{(0)} = -\frac{1}{re^\phi}\dual f_0,\qquad \nab J^{(+)} = \frac{1}{re^\phi}f_+,\qquad \nab J^{(-)} = \frac{1}{re^\phi}f_-.
\eeaa
In view of Definition \ref{def:renormalizationforJpbasisell=1modes}, this concludes the proof of the lemma.
\end{proof}

The following lemma provides identities for first order derivatives of $f_0$, $f_+$ and $f_-$.
\begin{lemma}\lab{lemma:computationforfirstorderderivativesf0fpfm}
We have on $S_*$
\bea
\bsplit
\div(f_0) &=f_0\c\nab\phi,\qquad\qquad\qquad\qquad \curl(f_0) = \frac{2}{re^{\phi}}\cos\th-f_0\wedge\nab\phi,\\
\div(f_+) &=-\frac{2}{re^\phi}J^{(+)}+f_+\c\nab\phi,\qquad \curl(f_+) = f_+\wedge\nab\phi,\\
\div(f_-) &=-\frac{2}{re^\phi}J^{(-)}+f_-\c\nab\phi,\qquad \curl(f_-) = -f_-\wedge\nab\phi,
\end{split}
\eea
and 
\bea
\bsplit
\nab\hot f_0 &= \left(\ba{cc}
 f_0\c\nab\phi & f_0\wedge\nab\phi\\
f_0\wedge\nab\phi & - f_0\c\nab\phi
\ea\right),\\
\nab\hot f_+ &= \left(\ba{cc}
(f_+)_2\nab_2\phi- (f_+)_1\nab_1\phi & -(f_+)_1\nab_2\phi- (f_+)_2\nab_1\phi\\
-(f_+)_1\nab_2\phi- (f_+)_2\nab_1\phi & -(f_+)_2\nab_2\phi+ (f_+)_1\nab_1\phi
\ea\right),\\
\nab\hot f_- &= \left(\ba{cc}
(f_-)_2\nab_2\phi- (f_-)_1\nab_1\phi & -(f_-)_1\nab_2\phi- (f_-)_2\nab_1\phi\\
-(f_-)_1\nab_2\phi- (f_-)_2\nab_1\phi & -(f_-)_2\nab_2\phi+ (f_-)_1\nab_1\phi
\ea\right).
\end{split}
\eea
In particular, in view of Definition \ref{def:renormalizationforf0fpfm}, we have
\beaa
\widecheck{\curl(f_0)} &=& \frac{2}{r}\cos\th(e^{-\phi}-1)  -f_0\wedge\nab\phi,\\
\widecheck{\div(f_\pm)} &=& -\frac{2}{r}J^{(\pm)}(e^{-\phi}-1)+f_\pm\c\nab\phi.
\eeaa
\end{lemma} 

\begin{proof}
See section \ref{sec:proofoflemma:computationforfirstorderderivativesf0fpfm} in the Appendix.
\end{proof}

\begin{lemma}\lab{lemma:estimatesforJpf0fpfmonsphereSstar}
On $S_*$, there holds, for $k\leq k_*-12$, 
\beaa
\left|\dkb^k\Big(\widecheck{\nab J^{(0)}},\,\, \widecheck{\nab J^{(+)}},\,\, \widecheck{\nab J^{(-)}}\Big)\right| &\les& \frac{\ep_0}{r^2u^{\frac{1}{2}+\dec}}.
\eeaa
Also, we have on $S_*$, for $k\leq k_*-13$,
\beaa
\left|\dkb^k\Big(\div(f_0),\, \widecheck{\curl(f_0)},\, \nab\hot f_0\Big)\right| &\les& \frac{\ep_0}{r^2u^{\frac{1}{2}+\dec}},\\
\left|\dkb^k\Big(\widecheck{\div(f_\pm)},\, \curl(f_\pm),\, \nab\hot f_\pm\Big)\right| &\les& \frac{\ep_0}{r^2u^{\frac{1}{2}+\dec}}.
\eeaa
In particular, we have on $S_*$, for $k\leq k_*-13$,
\beaa
\left|\dkb^k\Big(\nab f_0 -\frac{1}{r}\cos\th\in, \,\, \nab f_\pm +\frac{1}{r}J^{(\pm)}\de\Big)\right| &\les& \frac{\ep_0}{r^2u^{\frac{1}{2}+\dec}}.
\eeaa
\end{lemma}

\begin{proof}
The proof follows immediately from the identities of Lemma \ref{lemma:computationforfirstorderderivativesJpell=1basis} and Lemma \ref{lemma:computationforfirstorderderivativesf0fpfm} together with the control of $\phi$ provided by \eqref{eq:controlofphionSstar:onceagain}. 
\end{proof}

%%%%%%%%%%%%%%%%%%%%%%%%%%%%%%%%%%%%%%%%%%%%%%%%%%

\subsection{Proof of Proposition \ref{prop:estimatesforf0fplusfminusandJp}}

%%%%%%%%%%%%%%%%%%%%%%%%%%%%%%%%%%%%%%%%%%%%%%%%%%

We start with the following lemma.
\begin{lemma}\lab{lemma:simpleidentitiesforscalrproductsoff0fpfm}
We have on $\Si_*$
\beaa
&& |f_+|^2 = (\cos\th)^2(\cos\vphi)^2+(\sin\vphi)^2,\quad |f_-|^2 = (\cos\th)^2(\sin\vphi)^2+(\cos\vphi)^2,\\
&& |f_0|^2 = (\sin\th)^2, \quad f_+\c f_0=-J^{(-)}, \quad f_-\c f_0=J^{(+)}, \quad f_+\c f_-=-(\sin\th)^2\cos\vphi\sin\vphi.
\eeaa
\end{lemma}

\begin{proof}
Since $\nab_\nu f_\pm=0$, $\nab_\nu f_0=0$, $\nu(J^{(\pm)})=0$,  and $\nu(\th)=\nu(\vphi)=0$ on $\Si_*$, it suffices to prove these identities on $S_*$ which follows immediately from the definition of $f_\pm$, $f_0$ and $J^{(\pm)}$ on $S_*$.
\end{proof}

The next lemma relates angular derivatives of $\Jk$ and $\Jk_\pm$ with the ones of $f_0$ and $f_\pm$.
\begin{lemma}\lab{lemma:relationangularderivativesJkJkpmandff0fpfm}
We have on $\Si_*$
\bea
\bsplit
\widecheck{\ov{\DD}\c\Jk} &= O(r^{-4})+ \frac{2}{|q|}\div(f_0)  -\frac{2a^2\cos\th}{|q|^3}f_0\c\widecheck{\nab J^{(0)}}\\
&+i\left(\frac{2}{|q|}\widecheck{\curl(f_0)} -\frac{2a^2\cos\th}{|q|^3}f_0\c\dual\widecheck{\nab J^{(0)}}\right),
\end{split}
\eea
\bea
\bsplit
\widecheck{\ov{\DD}\c\Jk_\pm} &= O(r^{-4}) + \frac{2}{|q|}\widecheck{\div(f_\pm)} -\frac{2a^2\cos\th}{|q|^3}f_\pm\c\widecheck{\nab J^{(0)}}\\
&+i\left(\frac{2}{|q|}\curl(f_\pm)  -\frac{2a^2\cos\th}{|q|^3}f_\pm\c\dual\widecheck{\nab J^{(0)}}\right),
\end{split}
\eea
where $O(r^a)$ denotes, for $a\in\RRR$, a function of $(r, \cos\th)$ bounded by $r^a$ as $r\to +\infty$. 
\end{lemma}

\begin{proof}
See section \ref{sec:proofoflemma:relationangularderivativesJkJkpmandff0fpfm} in the Appendix.
\end{proof}

The following lemma provides a transport equation for $\nab\Jp$, $\nab f_0$ and $\nab f_\pm$ along $\Si_*$. 

\begin{lemma}\lab{lemma:transporteqnuofnabofJpf0fpfm}
Assume the following transversality conditions  on $\Si_*$ 
\beaa
\nu(\th)=0, \qquad \nu(\vphi)=0, \qquad \nab_4f_0=0, \qquad \nab_4f_+=0, \qquad \nab_4f_-=0.
\eeaa
Then, we have on $\Si_*$
\beaa
&&\nab_\nu\Big[r\nab f_0 - \cos\th\in\Big] =  \Ga_b \c \dkb^{\leq 1}f_0, \qquad \nab_\nu\Big[r\nab f_\pm + J^{(\pm)}\de\Big] =  \Ga_b \c \dkb^{\leq 1}f_\pm,\\
&& \nab_\nu\Big[r\widecheck{\nab\Jp}\Big]=\Ga_b \c \dkb^{\leq 1}\Jp, \quad p=0,+,-.
\eeaa
\end{lemma}

\begin{proof}
Recall from Corollary \ref{Corr:commutation-onSi_*} the following  commutation formula 
\beaa
\, [ \nab_{\nu},  r \nab]f &=  r \Ga_b \c \nab_{\nu} f + \Ga_b \c \dk^{\leq 1}f.
 \eeaa
 Applying it to $f_0$, $f_\pm$, and using the fact that $\nab_\nu(f_0, f_\pm)=0$ and $\nab_4(f_0, f_\pm)=0$, we infer
\beaa
&&\nab_\nu(r\nab f_0) =  \Ga_b \c \dkb^{\leq 1}f_0, \qquad \nab_\nu(r\nab f_\pm) =  \Ga_b \c \dkb^{\leq 1}f_\pm,\\
&& \nab_\nu(r\nab\Jp) =  \Ga_b \c \dkb^{\leq 1}\Jp, \quad p=0,+,-, 
\eeaa 
which yields, since $\nu(\th)=\nu(\vphi)=0$, 
\beaa
&&\nab_\nu\Big[r\nab f_0 - \cos\th\in\Big] =  \Ga_b \c \dkb^{\leq 1}f_0, \qquad \nab_\nu\Big[r\nab f_\pm + J^{(\pm)}\de\Big] =  \Ga_b \c \dkb^{\leq 1}f_\pm,\\
&& \nab_\nu\Big[r\widecheck{\nab\Jp}\Big]=\Ga_b \c \dkb^{\leq 1}\Jp, \quad p=0,+,-,
\eeaa
as stated. This concludes the proof of the lemma.
\end{proof}

We are now ready to prove Proposition \ref{prop:estimatesforf0fplusfminusandJp}. 

\begin{proof}[Proof of Proposition \ref{prop:estimatesforf0fplusfminusandJp}]
Recall from Lemma \ref{lemma:transporteqnuofnabofJpf0fpfm} that we have on $\Si_*$
\beaa
&&\nab_\nu\Big[r\nab f_0 - \cos\th\in\Big] =  \Ga_b \c \dkb^{\leq 1}f_0, \qquad \nab_\nu\Big[r\nab f_\pm + J^{(\pm)}\de\Big] =  \Ga_b \c \dkb^{\leq 1}f_\pm,\\
&& \nab_\nu\Big[r\widecheck{\nab\Jp}\Big]=\Ga_b \c \dkb^{\leq 1}\Jp, \quad p=0,+,-.
\eeaa
Using the commutation formula  from Corollary \ref{Corr:commutation-onSi_*}   
\beaa
\, [ \nab_{\nu},  r \nab]f &=  r \Ga_b \c \nab_{\nu} f + \Ga_b \c \dk^{\leq 1}f,
 \eeaa
we infer
\beaa
&&\nab_\nu\dkb^k\Big[r\nab f_0 - \cos\th\in\Big] = \dkb^{\leq k}(\Ga_b \c \dkb^{\leq 1}f_0)+ \dkb^{\leq k}\left(\Ga_b\c\Big[r\nab f_0 - \cos\th\in\Big]\right), \\ 
&& \nab_\nu\dkb^k\Big[r\nab f_\pm + J^{(\pm)}\de\Big] =  \dkb^{\leq k}(\Ga_b \c \dkb^{\leq 1}f_\pm)+ \dkb^{\leq k}\left(\Ga_b\c\Big[r\nab f_\pm + J^{(\pm)}\de\Big]\right),\\
&& \nab_\nu\dkb^k\Big[r\widecheck{\nab\Jp}\Big]=\Ga_b \c \dkb^{\leq 1}\Jp, \quad p=0,+,-.
\eeaa
In view of  Corollary \ref{cor:integrate-transportSi_*}, we infer on $\Si_*$, for all $k\leq k_*-12$
\beaa
r\left|\dkb^k\Big[r\widecheck{\nab\Jp}\Big]\right| &\les& \left|r\dkb^k\Big[r\widecheck{\nab\Jp}\Big]\right|_{L^\infty(S_*)}+\int_u^{u_*}r|\dkb^{\leq k}\Ga_b|
\eeaa
and on $\Si_*$, for all $k\leq k_*-13$
\beaa
r\left|\dkb^k\Big[r\nab f_0 - \cos\th\in\Big]\right| &\les& \left|r\dkb^k\Big[r\nab f_0 - \cos\th\in\Big]\right|_{L^\infty(S_*)}+\int_u^{u_*}r|\dkb^{\leq k}\Ga_b|,\\
r\left|\dkb^k\Big[r\nab f_\pm + J^{(\pm)}\de\Big]\right| &\les& \left|r\dkb^k\Big[r\nab f_\pm + J^{(\pm)}\de\Big]\right|_{L^\infty(S_*)}+\int_u^{u_*}r|\dkb^{\leq k}\Ga_b|. 
\eeaa
Now, in view of Sobolev and Proposition \ref{Prop.Flux-bb-vthb-eta-xib}, we have on $\Si_*$, for $k\leq k_*-12$,
\beaa
\int_u^{u_*}r|\dkb^{\leq k}\Ga_b| &\les& \frac{1}{u^{\frac{1}{2}+\dec}}\left(\int_u^{u_*}r^2{u'}^{2+2\dec}|\dkb^{\leq k_*-12}\Ga_b|^2\right)^{\frac{1}{2}}\\
&\les& \frac{1}{u^{\frac{1}{2}+\dec}}\left(\int_u^{u_*}{u'}^{2+2\dec}\|\dkb^{\leq k_*-10}\Ga_b\|_{L^2(S)}^2\right)^{\frac{1}{2}}\\
&\les& \frac{1}{u^{\frac{1}{2}+\dec}}\left(\int_{\Si_*}u^{2+2\dec}|\dkb^{\leq k_*-10}\Ga_b|^2\right)^{\frac{1}{2}}\\
&\les& \frac{\ep_0}{u^{\frac{1}{2}+\dec}}.
\eeaa
Together with the control on $S_*$ provided by Lemma \ref{lemma:estimatesforJpf0fpfmonsphereSstar}, we infer on $\Si_*$, for all $k\leq k_*-12$
\beaa
\left|\dkb^k\Big[\widecheck{\nab\Jp}\Big]\right| &\les& \frac{\ep_0}{r^2u^{\frac{1}{2}+\dec}},
\eeaa
and on $\Si_*$, for all $k\leq k_*-13$
\beaa
\left|\dkb^k\left[\nab f_0 - \frac{1}{r}\cos\th\in\right]\right| + \left|\dkb^k\left[\nab f_\pm + \frac{1}{r}J^{(\pm)}\de\right]\right| &\les& \frac{\ep_0}{r^2u^{\frac{1}{2}+\dec}}. 
\eeaa
Together with Lemma \ref{lemma:transporteqnuofnabofJpf0fpfm} and the control for $\Ga_b$ in Proposition \ref{Prop.Flux-bb-vthb-eta-xib}, we deduce
on $\Si_*$, for all $k\leq k_*-12$
\beaa
\left|\dk_*^k\Big[\widecheck{\nab\Jp}\Big]\right| &\les& \frac{\ep_0}{r^2u^{\frac{1}{2}+\dec}},
\eeaa
and on $\Si_*$, for all $k\leq k_*-13$
\beaa
\left|\dk_*^k\left[\nab f_0 - \frac{1}{r}\cos\th\in\right]\right| + \left|\dk_*^k\left[\nab f_\pm + \frac{1}{r}J^{(\pm)}\de\right]\right| &\les& \frac{\ep_0}{r^2u^{\frac{1}{2}+\dec}}. 
\eeaa
This implies in particular, on $\Si_*$, for all $k\leq k_*-13$,
\beaa
\left|\dk_*^k\div(f_0)\right| + \left|\dk_*^k\widecheck{\curl(f_0)}\right| +\left|\dk_*^k\nab\hot f_0\right| \\
+\left|\dk_*^k\widecheck{\div(f_\pm)}\right|+\left|\dk_*^k\curl(f_\pm)\right|+\left|\dk_*^k\nab\hot f_\pm\right|  &\les& \frac{\ep_0}{r^2u^{\frac{1}{2}+\dec}}. 
\eeaa

Finally, recalling from Lemma \ref{lemma:relationangularderivativesJkJkpmandff0fpfm} the following identities on $\Si_*$
\beaa
\widecheck{\ov{\DD}\c\Jk} &=& O(r^{-4})+ \frac{2}{|q|}\div(f_0)  -\frac{2a^2\cos\th}{|q|^3}f_0\c\widecheck{\nab J^{(0)}}\\
&&+i\left(\frac{2}{|q|}\widecheck{\curl(f_0)} -\frac{2a^2\cos\th}{|q|^3}f_0\c\dual\widecheck{\nab J^{(0)}}\right),
\eeaa
\beaa
\widecheck{\ov{\DD}\c\Jk_\pm} &=& O(r^{-4}) + \frac{2}{|q|}\widecheck{\div(f_\pm)} -\frac{2a^2\cos\th}{|q|^3}f_\pm\c\widecheck{\nab J^{(0)}}\\
&&+i\left(\frac{2}{|q|}\curl(f_\pm)  -\frac{2a^2\cos\th}{|q|^3}f_\pm\c\dual\widecheck{\nab J^{(0)}}\right),
\eeaa
where $O(r^a)$ denotes, for $a\in\RRR$, a function of $(r, \cos\th)$ bounded by $r^a$ as $r\to +\infty$, we immediately infer from the above that there holds on $\Si_*$, for all $k\leq k_*-13$,
\beaa
 \left|\dk_*^k\widecheck{\ov{\DD}\c\Jk}\right|+\left|\dk_*^k\widecheck{\ov{\DD}\c\Jk_\pm}\right|  &\les& \frac{\ep_0}{r^3u^{\frac{1}{2}+\dec}}+\frac{1}{r^4}. 
\eeaa
In view of the dominance condition \eqref{eq:dominantconditiononronSigmastarchap5} on $r$ on $\Si_*$, this yields, for $k\leq k_*-13$,
\beaa
 \left|\dk_*^k\widecheck{\ov{\DD}\c\Jk}\right|+\left|\dk_*^k\widecheck{\ov{\DD}\c\Jk_\pm}\right|  &\les& \frac{\ep_0}{r^3u^{\frac{1}{2}+\dec}}
\eeaa
which concludes the proof of Proposition \ref{prop:estimatesforf0fplusfminusandJp}. 
\end{proof}

%%%%%%%%%%%%%%%%%%%%%%%%%%%%%%%%%%%%%%%%%%%%%%%%%%

\subsection{A additional estimate for $\b$ on $\Si_*$}

%%%%%%%%%%%%%%%%%%%%%%%%%%%%%%%%%%%%%%%%%%%%%%%%%%

This section is devoted to a decay estimate for $\b$, see Corollary \ref{cor:lastmissingimproveddecayestimateforbetaonSigmastar}. We start with the following two lemmas.
\begin{lemma}\lab{lemma:decompositionofHodgeoperatorappliedtorenormalizeationfobeta}
We have on $\Si_*$, for $k\leq k_*-13$, 
\beaa
\left|\dkb^k\nab\hot\left(\b-\frac{3am}{r^4}f_0\right)\right| &\les& |\dkb^k\nab\hot\b| +\frac{\ep_0}{r^6u^{\frac{1}{2}+\dec}},\\
\left|\div\left(\b-\frac{3am}{r^4}f_0\right) -\div\b\right| &\les& \frac{\ep_0}{r^6u^{\frac{1}{2}+\dec}},\\
\left|\curl\left(\b-\frac{3am}{r^4}f_0\right) -\left(\curl\b -\frac{6am J^{(0)}}{r^5}\right)\right| &\les& \frac{\ep_0}{r^6u^{\frac{1}{2}+\dec}}.
\eeaa
\end{lemma}

\begin{proof}
The proof follows immediately from the control of $f_0$ provided by Proposition \ref{prop:estimatesforf0fplusfminusandJp} and the definition of $J^{(0)}$ and $\widecheck{\curl (f_0)}$. 
\end{proof}

\begin{lemma}\lab{lemma:statementeq:DeJp.Sigmastar:improvedd}
The     functions $\Jp$ verify the following properties on $\Si_*$
\beaa
\bsplit
\int_{S}J^{(p)} &= O\left(\ep_0 ru^{-\frac{1}{2}-\dec}\right),\\
\int_{S}J^{(p)}J^{(q)} &= \frac{4\pi}{3}r^2\de_{pq}+O\left(\ep_0 ru^{-\frac{1}{2}-\dec}\right),
\end{split}
\eeaa
\beaa
\left|\left( \De +\frac{2}{r^2}\right)    \Jp\right|    \les  \ep_0 r^{-3} u^{-\frac{1}{2}-\dec},
\eeaa
and
\beaa
\left|\dds_2\dds_1 \Jp\right|    \les  \ep_0 r^{-3} u^{-\frac{1}{2}-\dec},
\eeaa
where by $\dds_1\Jp$, we mean either $\dds_1(\Jp,0)$ or $\dds_1(0,\Jp)$. 
\end{lemma}

\begin{proof}
The proof follows exactly the same lines as the proof of Corollary \ref{corr:Si*-ell=1modes-improved} by replacing the control on $S_*$ of Lemma \ref{lemma:statementeq:DeJp.S_*} by the improved control on $S_*$ provided by  Corollary \ref{cor:statementeq:DeJp.S_*:improved}.
\end{proof}

The following lemma controls the $\ell=1$ modes of $J^{(0)}$.
\begin{lemma}\lab{lemma:controlofell=1modeofcosthetaonSistar}
We have on $\Si_*$
\beaa
\left|\left(J^{(0)}\right)_{\ell=1,0} - \frac{1}{3}\right|+\left|\left(J^{(0)}\right)_{\ell=1,\pm}\right| &\les& \frac{\ep_0}{ru^{\frac{1}{2}+\dec}}.
\eeaa
\end{lemma}

\begin{proof}
We have, by definition of the $\ell=1$ modes on $\Si_*$, for $p=0,+,-$, 
\beaa
\left(J^{(0)}\right)_{\ell=1,p} &=& \frac{1}{|S|}\int_S J^{(0)}\Jp
\eeaa
and hence
\beaa
\left(J^{(0)}\right)_{\ell=1,p} -\frac{1}{3}\de_{p0} &=& \frac{1}{4\pi r^2}\left(\int_S J^{(0)}\Jp -\frac{4\pi}{3}r^2\de_{p0}\right).
\eeaa
The proof follows then from Lemma \ref{lemma:statementeq:DeJp.Sigmastar:improvedd}.
\end{proof}

\begin{corollary}\lab{cor:lastmissingimproveddecayestimateforbetaonSigmastar}
We have on $\Si_*$
\beaa
\left|\left[\ddd_1\left(\b-\frac{3am}{r^4}f_0\right)\right]_{\ell=1}\right| &\les& \frac{\ep_0}{r^5u^{1+\dec}}.
\eeaa
and, for $k\leq k_*-14$, 
\beaa
\sup_{\Si_*}r^{\frac{7}{2}+\dee}\left|\dkb^{\leq k_{\large}}\left(\b-\frac{3am\sin\th}{r^4}f_0\right)\right| &\les& \ep_0.
\eeaa
\end{corollary}

\begin{proof}
We start with the first estimate. In view of the definition of $\ddd_1$ and Lemma \ref{lemma:decompositionofHodgeoperatorappliedtorenormalizeationfobeta}, we have 
\beaa
\left|\left[\ddd_1\left(\b-\frac{3am}{r^4}f_0\right)\right]_{\ell=1}\right| &\les& \left|(\div\b)_{\ell=1}\right|+\left|(\curl\b)_{\ell=1} -\frac{6am (J^{(0)})_{\ell=1}}{r^5}\right|+\frac{\ep_0}{r^6u^{\frac{1}{2}+\dec}}.
\eeaa
Together with Lemma \ref{lemma:controlofell=1modeofcosthetaonSistar}, we infer
\beaa
\left|\left[\ddd_1\left(\b-\frac{3am}{r^4}f_0\right)\right]_{\ell=1}\right| &\les& \left|(\div\b)_{\ell=1}\right|+\left|(\curl\b)_{\ell=1,\pm}\right|
+\left|(\curl\b)_{\ell=1,0} -\frac{2am}{r^5}\right|\\
&&+\frac{\ep_0}{r^6u^{\frac{1}{2}+\dec}}.
\eeaa
In view of Proposition \ref{prop:controlofell=1modesonSigmastar}, we infer
\beaa
\left|\left[\ddd_1\left(\b-\frac{3am}{r^4}f_0\right)\right]_{\ell=1}\right| &\les& \frac{\ep_0}{r^5u^{1+\dec}}+\frac{\ep_0}{r^6u^{\frac{1}{2}+\dec}}
\eeaa
and hence, using the dominance of $r$ on $\Si_*$, 
\beaa
\left|\left[\ddd_1\left(\b-\frac{3am}{r^4}f_0\right)\right]_{\ell=1}\right| &\les& \frac{\ep_0}{r^5u^{1+\dec}}
\eeaa
as stated.

Next, we focus on the second estimate. In view of Lemma \ref{prop:2D-Hodge2}, we have
\beaa
  \left\|\b-\frac{3am}{r^4}f_0\right\|_{\hk_{k_*-12}  (S)}&\les&  r \left\|\dds_2\left(\b-\frac{3am}{r^4}f_0\right)\right\|_{\hk_{k_*-13}(S)}+r^2 \left| \left(\ddd_1\left(\b-\frac{3am}{r^4}f_0\right)\right)_{\ell=1}\right|.
  \eeaa
In view of the above estimate, we deduce
\beaa
 \left\|\dkb^{k_*-14}\left(\b-\frac{3am}{r^4}f_0\right)\right\|_{L^\infty(S)} &\les&  r\left\|\dkb^{\leq k_*-13}\nab\hot\b\right\|_{L^\infty(S)}+\frac{\ep_0}{r^4u^{1+\dec}}.
  \eeaa 
It remains to control $\nab\hot\b$. We have the following consequence of Bianchi
\beaa
\nab\hot \b &=& \nab_3\a -\frac{\Up}{r} \a  +r^{-3}\Ga_g  +\Ga_b\c(\a,\b)+r^{-1}\Ga_g\c\Ga_g.
\eeaa 
Together with the control of $\a$ and $\nab_3\a$ given by {\bf Ref 2} and the control for $\Ga_g$ and $\Ga_b$ in {\bf Ref 1}, we infer
\beaa
\left\|\dkb^{\leq k_*-13}\nab\hot\b\right\|_{L^\infty(S)} &\les& \frac{\ep_0}{r^{\frac{9}{2}+\dee}}
\eeaa
and hence
\beaa
 \left\|\dkb^{k_*-14}\left(\b-\frac{3am}{r^4}f_0\right)\right\|_{L^\infty(S)} &\les&  \frac{\ep_0}{r^{\frac{7}{2}+\dee}} +\frac{\ep_0}{r^4u^{1+\dec}}\les \frac{\ep_0}{r^{\frac{7}{2}+\dee}}
  \eeaa 
  as stated. This concludes the proof of the corollary.
\end{proof}

%%%%%%%%%%%%%%%%%%%%%%%%%%%%%%%%%%%%%%%%%%%%%%%%%%%%

\subsection{An estimate for high order derivatives of $J^{(p)}$ and $\Jk$}

%%%%%%%%%%%%%%%%%%%%%%%%%%%%%%%%%%%%%%%%%%%%%%%%%%%%

In this section, we derive the following proposition on the control of $k_{large}$ derivatives of $J^{(p)}$, $f_0$, $f_\pm$ and $\Jk$. 

\begin{proposition}\lab{prop:estimatesforf0fplusfminusandJp:klarge}
We have on $\Si_*$, for all $k\leq k_{large}$
\beaa
\left|\dk_*^k\Big[\widecheck{\nab\Jp}\Big]\right| &\les& \frac{\ep}{r},
\eeaa
and on $\Si_*$, for all $k\leq k_{large}-1$
\beaa
\left|\dk_*^k\left[\div(f_0),\, \widecheck{\curl(f_0)},\, \nab\hot f_0,\, \nab f_0 - \frac{1}{r}\cos\th\in\right]\right| \\
+ \left|\dk_*^k\left[\widecheck{\div(f_\pm)},\, \curl(f_\pm),\, \nab\hot f_\pm,\, \nab f_\pm + \frac{1}{r}J^{(\pm)}\de\right]\right| &\les& \frac{\ep}{r},
\eeaa
as well as
\beaa
 \left|\dk_*^k\widecheck{\ov{\DD}\c\Jk}\right|+\left|\dk_*^k\widecheck{\ov{\DD}\c\Jk_\pm}\right|  &\les& \frac{\ep}{r^2}.
\eeaa
\end{proposition}

\begin{proof}
In view of the Gauss equation, we have $\widecheck{K}\in r^{-1}\Ga_g$. Together with the estimate for $\Ga_g$  in {\bf Ref 1}, we infer 
\beaa
\sup_{S_*}\left|\dkb^{\leq k_{large}}\left(K-\frac{1}{r^2}\right)\right| &\les& \frac{\ep}{r^3}.
\eeaa
Arguing as in Corollary \ref{cor:statementeq:DeJp.S_*:improved}, we infer 
\beaa
\sup_{S_*}\left|\dkb^{\leq k_{large}}\phi\right| &\les& \frac{\ep}{r}.
\eeaa
Then, arguing as in Lemma \ref{lemma:estimatesforJpf0fpfmonsphereSstar}, we deduce, for $k\leq k_{large}$,  
\beaa
\sup_{S_*}\left|\dkb^k\Big(\widecheck{\nab J^{(0)}},\,\, \widecheck{\nab J^{(+)}},\,\, \widecheck{\nab J^{(-)}}\Big)\right| &\les& \frac{\ep}{r^2}.
\eeaa
and for $k\leq k_{large}-1$,
\beaa
\sup_{S_*}\left|\dkb^k\Big(\div(f_0),\, \widecheck{\curl(f_0)},\, \nab\hot f_0\Big)\right| &\les& \frac{\ep}{r^2},\\
\sup_{S_*}\left|\dkb^k\Big(\widecheck{\div(f_\pm)},\, \curl(f_\pm),\, \nab\hot f_\pm\Big)\right| &\les& \frac{\ep}{r^2},\\
\sup_{S_*}\left|\dkb^k\Big(\nab f_0 -\frac{1}{r}\cos\th\in, \,\, \nab f_\pm +\frac{1}{r}J^{(\pm)}\de\Big)\right| &\les& \frac{\ep}{r^2}.
\eeaa

Next, recall from the proof of Proposition \ref{prop:estimatesforf0fplusfminusandJp}
\beaa
&&\nab_\nu\dkb^k\Big[r\nab f_0 - \cos\th\in\Big] = \dkb^{\leq k}(\Ga_b \c \dkb^{\leq 1}f_0)+ \dkb^{\leq k}\left(\Ga_b\c\Big[r\nab f_0 - \cos\th\in\Big]\right), \\ 
&& \nab_\nu\dkb^k\Big[r\nab f_\pm + J^{(\pm)}\de\Big] =  \dkb^{\leq k}(\Ga_b \c \dkb^{\leq 1}f_\pm)+ \dkb^{\leq k}\left(\Ga_b\c\Big[r\nab f_\pm + J^{(\pm)}\de\Big]\right),\\
&& \nab_\nu\dkb^k\Big[r\widecheck{\nab\Jp}\Big]=\Ga_b \c \dkb^{\leq 1}\Jp, \quad p=0,+,-.
\eeaa
In view of  Corollary \ref{cor:integrate-transportSi_*}, we infer on $\Si_*$, for all $k\leq k_{large}$
\beaa
r\left|\dkb^k\Big[r\widecheck{\nab\Jp}\Big]\right| &\les& \left|r\dkb^k\Big[r\widecheck{\nab\Jp}\Big]\right|_{L^\infty(S_*)}+\int_u^{u_*}r|\dkb^{\leq k}\Ga_b|
\eeaa
and on $\Si_*$, for all $k\leq k_{large}-1$
\beaa
r\left|\dkb^k\Big[r\nab f_0 - \cos\th\in\Big]\right| &\les& \left|r\dkb^k\Big[r\nab f_0 - \cos\th\in\Big]\right|_{L^\infty(S_*)}+\int_u^{u_*}r|\dkb^{\leq k}\Ga_b|,\\
r\left|\dkb^k\Big[r\nab f_\pm + J^{(\pm)}\de\Big]\right| &\les& \left|r\dkb^k\Big[r\nab f_\pm + J^{(\pm)}\de\Big]\right|_{L^\infty(S_*)}+\int_u^{u_*}r|\dkb^{\leq k}\Ga_b|. 
\eeaa
In view of the control of $\Ga_b$ provided by {\bf Ref 1}, we have, for $k\leq k_{large}$, 
\beaa
\int_u^{u_*}r|\dkb^{\leq k}\Ga_b| &\les& \ep\int_1^{u_*}\frac{du}{u^{1+\dec}}\les\ep.
\eeaa
Together with the above control on $S_*$, we obtain 
 on $\Si_*$, for all $k\leq k_{large}$
\beaa
\left|\dkb^k\Big[r\widecheck{\nab\Jp}\Big]\right| &\les& \frac{\ep}{r}
\eeaa
and on $\Si_*$, for all $k\leq k_{large}-1$
\beaa
\left|\dkb^k\Big[r\nab f_0 - \cos\th\in\Big]\right| &\les& \frac{\ep}{r},\\
\left|\dkb^k\Big[r\nab f_\pm + J^{(\pm)}\de\Big]\right| &\les& \frac{\ep}{r}. 
\eeaa
This implies in particular, on $\Si_*$, for all $k\leq k_{large}-1$,
\beaa
\left|\dk_*^k\div(f_0)\right| + \left|\dk_*^k\widecheck{\curl(f_0)}\right| +\left|\dk_*^k\nab\hot f_0\right| \\
+\left|\dk_*^k\widecheck{\div(f_\pm)}\right|+\left|\dk_*^k\curl(f_\pm)\right|+\left|\dk_*^k\nab\hot f_\pm\right|  &\les& \frac{\ep_0}{r}. 
\eeaa

Finally, recalling from Lemma \ref{lemma:relationangularderivativesJkJkpmandff0fpfm} the following identities on $\Si_*$
\beaa
\widecheck{\ov{\DD}\c\Jk} &=& O(r^{-4})+ \frac{2}{|q|}\div(f_0)  -\frac{2a^2\cos\th}{|q|^3}f_0\c\widecheck{\nab J^{(0)}}\\
&&+i\left(\frac{2}{|q|}\widecheck{\curl(f_0)} -\frac{2a^2\cos\th}{|q|^3}f_0\c\dual\widecheck{\nab J^{(0)}}\right),
\eeaa
\beaa
\widecheck{\ov{\DD}\c\Jk_\pm} &=& O(r^{-4}) + \frac{2}{|q|}\widecheck{\div(f_\pm)} -\frac{2a^2\cos\th}{|q|^3}f_\pm\c\widecheck{\nab J^{(0)}}\\
&&+i\left(\frac{2}{|q|}\curl(f_\pm)  -\frac{2a^2\cos\th}{|q|^3}f_\pm\c\dual\widecheck{\nab J^{(0)}}\right),
\eeaa
where $O(r^a)$ denotes, for $a\in\RRR$, a function of $(r, \cos\th)$ bounded by $r^a$ as $r\to +\infty$, we immediately infer from the above that there holds on $\Si_*$, for all $k\leq k_{large}-1$,
\beaa
 \left|\dk_*^k\widecheck{\ov{\DD}\c\Jk}\right|+\left|\dk_*^k\widecheck{\ov{\DD}\c\Jk_\pm}\right|  &\les& \frac{\ep}{r^2} 
\eeaa
which concludes the proof of Proposition \ref{prop:estimatesforf0fplusfminusandJp:klarge}. 
\end{proof}

Proposition \ref{prop:estimatesforf0fplusfminusandJp} and Proposition \ref{prop:estimatesforf0fplusfminusandJp:klarge} motivate the following definition.

\begin{definition}\lab{def:Gammabtilde}
We denote by $\widetilde{\Ga}_b$ the set of linearized quantities below
\beaa
\widetilde{\Ga}_b &:=& \Ga_b\cup\Big\{\widecheck{\nab\Jp}\Big\},\\
\dk_*^k\widetilde{\Ga}_b &:=& \dk_*^k\Ga_b\cup\Big\{\dk_*^k\widecheck{\nab\Jp}\Big\}\cup\dk_*^{k-1}\widetilde{\Ga}_{b,1},\quad\textrm{for}\quad k\geq 1,\\
\widetilde{\Ga}_{b,1} &=& \Big\{\div(f_0),\, \widecheck{\curl(f_0)},\, \nab\hot f_0,\,\widecheck{\div(f_\pm)},\, \curl(f_\pm),\, \nab\hot f_\pm,\, r\widecheck{\ov{\DD}\c\Jk},\, r\widecheck{\ov{\DD}\c\Jk_\pm}\Big\}.
\eeaa
\end{definition}

\begin{corollary}\lab{cor:estimatesforf0fplusfminusandJp:Gammatildeb}
We have on $\Si_*$, for $k\leq k_*-12$,
\beaa
|\dk_*^k\widetilde{\Ga}_b| &\les& \frac{\ep_0}{ru^{1+\dec}},
\eeaa
and, for $k\leq k_{large}$
\beaa
|\dk_*^k\widetilde{\Ga}_b| &\les& \frac{\ep}{r}.
\eeaa
\end{corollary}

\begin{proof}
This is an immediate consequence of {\bf Ref 1} and Proposition \ref{Prop.Flux-bb-vthb-eta-xib} for $\Ga_b$, and of 
Proposition \ref{prop:estimatesforf0fplusfminusandJp} and Proposition \ref{prop:estimatesforf0fplusfminusandJp:klarge} for the rest of $\widetilde{\Ga}_b$.
\end{proof}

\begin{remark}
In view of Corollary \ref{cor:estimatesforf0fplusfminusandJp:Gammatildeb}, $\widetilde{\Ga}_b$ enjoys the same estimates as $\Ga_b$. Note that the estimates of angular derivatives of $\Jp$, $f_0$, $f_\pm$ and $\Jk$ are
\begin{itemize}
\item consistent with $\Ga_g$ for $k\leq k_*-12$ derivatives in view of Proposition \ref{prop:estimatesforf0fplusfminusandJp},

\item consistent with $\Ga_b$ for $k\leq k_{large}$ derivatives in view of Proposition \ref{prop:estimatesforf0fplusfminusandJp:klarge}.
\end{itemize}
We do not need the better decay properties and simply treat these angular derivatives as $\Ga_b$, which justifies their inclusion in $\widetilde{\Ga}_b$.
\end{remark}

%%%%%%%%%%%%%%%%%%%%%%%%%%%%%%%%%%%%%%%%%%%%%%%%%%%%

\section{Decay estimates for the PG frame on $\Si_*$}\lab{sec:decayestimatesPGframeonSigmastar}

%%%%%%%%%%%%%%%%%%%%%%%%%%%%%%%%%%%%%%%%%%%%%%%%%%%%

In this section, we use the decay estimates derived for the integrable frame of $\Si_*$ in Proposition \ref{prop:decayonSigamstarofallquantities}, and the estimates of section \ref{sec:controlofJpandJkonSigmastar},  
to derive decay estimates for the PG frame of $\Mext$. This is a prerequisite to the improvement of the bootstrap assumptions on decay on $\Mext$ of Chapter 6.

%%%%%%%%%%%%%%%%%%%%%%%%%%%%%%%%%%%%%%%%%%%%%%%%%%

\subsection{Initialization of the PG frame on $\Si_*$}

%%%%%%%%%%%%%%%%%%%%%%%%%%%%%%%%%%%%%%%%%%%%%%%%%%

Let $(e_3, e_4, e_1, e_2)$ denote the null frame of $\Si_*$, and let $(e_3', e_4', e_1', e_2')$ denote the PG frame of $\Mext$. Then, $(e_3', e_4', e_1', e_2')$ is initialized on $\Si_*$ by
 \bea\lab{eq:changofframedefiningthefrmaeofMext:chap5}
 \bsplit
  e_4'&=e_4 + f^b  e_b +\frac 1 4 |f|^2  e_3,\\
  e_a' &= \left(\de_a^b +\frac{1}{2}\fb_af^b\right) e_b +\frac 1 2  \fb_a  e_4 +\left(\frac 1 2 f_a +\frac{1}{8}|f|^2\fb_a\right)   e_3,\\
 e_3'&=\left(1+\frac{1}{2}f\c\fb  +\frac{1}{16} |f|^2  |\fb|^2\right) e_3 + \left(\fb^b+\frac 1 4 |\fb|^2f^b\right) e_b  + \frac 1 4 |\fb|^2 e_4,
 \end{split}
 \eea
 where 
\bea\lab{eq:defoffandfbforthedefintionoftheframeofMext:chapp5}
f=\frac{a}{r}f_0, \qquad \fb = -\frac{(\nu(r)-b_*)}{1-\frac{1}{4}b_* |f|^2}f,
\eea
with the 1-form $f_0$ being defined in Definition \ref{def:definitionoff0fplusfminus}. All quantities with primes denote in section \ref{sec:decayestimatesPGframeonSigmastar} the ones corresponding to the  PG frame of $\Mext$.  Furthermore, the coordinates  $(r',u', \th', \vphi')$ associated to the PG frame of $\Mext$ are initialized on $\Si_*$ as follows
\bea
r'=r, \qquad u'=u, \qquad \th'=\th, \qquad\vphi'=\vphi, \quad\textrm{on}\quad\Mext.
\eea 

Also, note that the complex horizontal 1-form $\Jk$ introduced in Definition \ref{def:definitionofJkonMext} verifies
 \bea
 \bsplit
\Jk &= \frac{1}{|q|}\left(f_0+i\dual f_0\right)\quad \textrm{on}\quad\Si_*,\\
\nab_4'\Jk &= -\frac{1}{q}\Jk\quad \textrm{on}\quad\Mext.
\end{split}
\eea
We also introduce the following two complex horizontal 1-forms $\Jk_\pm$ given by
\bea
\bsplit
\Jk_\pm &= \frac{1}{|q|}\left(f_\pm+i\dual f_\pm\right)\quad \textrm{on}\quad\Si_*,\\
\nab_4'\Jk_{\pm} &= -\frac{1}{q}\Jk_{\pm}\quad \textrm{on}\quad\Mext.
\end{split}
\eea
 
\begin{remark} 
Recall that  the complex 1-form $\Jk$ is needed to linearize $Z$, $H$, $\DD(\cos\th)$, and $\DD(u)$, and the  complex 1-form $\Jk_\pm$ is needed to linearize $\DD(J^{(\pm)})$. 
\end{remark}

Also, note that the transformation formulas involve all Ricci coefficients of the foliation of $\Si_*$. We thus need to prescribe transversality conditions for the Ricci coefficients not defined on $\Si_*$, i.e. $\xi$, $\om$ and $\etab$. We recall that we choose them to be compatible with an outgoing geodesic foliation initialized on $\Si_*$, i.e.
\bea
\xi=0, \qquad \om=0, \qquad \etab=-\ze.
\eea

Recall the definition of $\Ga_b'$ and $\Ga_g'$ in the PG frame of $\Mext$.

\begin{definition}
Recall Definition \ref{def:renormalizationofallnonsmallquantitiesinPGstructurebyKerrvalue} for the definition of the linearized quantities in an outgoing PG frame. The set of all linearized quantities is of the form $\Ga_g\cup \Ga_b$ with  $\Ga_g,  \Ga_b$
 defined as follows.
 \begin{enumerate}
\item 
 The set   $\Ga_g$   with
 \bea
 \Ga_g' &= \Bigg\{\trXc',\quad  \Xh',\quad \Zc',\quad \trXbc', \quad r\Pc', \quad  rB', \quad  rA'\Bigg\}.
 \eea
 
 \item  The set  $\Ga_b'=\Ga_{b,1}'\cup \Ga_{b, 2}'\cup \Ga_{b, 3}'$   with
 \bea
 \bsplit
 \Ga_{b,1}' &= \Bigg\{\Hc', \quad \Xbh', \quad \ombc', \quad \Xib',\quad  r\Bb', \quad \Ab'\Bigg\},\\
  \Ga_{b, 2}' &= \Bigg\{r^{-1}\widecheck{e_3'(r)}, \quad  \widecheck{\DD'(\cos\th)}, \quad e_3'(\cos\th), \quad \widecheck{\DD'u}, \quad r^{-1}\widecheck{e_3'(u)},\\
 & \qquad \qquad \qquad \qquad  \widecheck{\DD'(J^{(+)})}, \quad \widecheck{\DD'(J^{(-)})}, \quad \widecheck{e_3'(J^{(+)})}, \quad \widecheck{e_3'(J^{(-)})}\Bigg\},\\
  \Ga_{b,3}&=\Bigg\{ r\,\widecheck{\ov{\DD'}\c\Jk}, \quad r\,\DD'\hot\Jk, \quad r\,\widecheck{\nab_3'\Jk}, \quad
   r\,\widecheck{\ov{\DD'}\c\Jk_\pm}, \quad r\,\DD'\hot\Jk_\pm, \quad r\,\widecheck{\nab_3'\Jk_\pm}\Bigg\}.
 \end{split}
 \eea
\end{enumerate}
\end{definition}

The goal of this section is to prove the following proposition concerning the control of the PG frame on $\Si_*$.  
\begin{proposition}\lab{prop:improvedesitmatesfortemporalframeofMextonSigmastar}
We have on $\Si_*$, for $k\leq k_*-15$, 
\beaa 
\sup_{\Si_*}\Big(ru^{1+\dec}|\dk^k\Ga_b'|+r^2u^{\frac{1}{2}+\dec}|\dk^k\Ga_g'|+r^2u^{1+\dec}|\dk^{k-1}\nab_3\Ga_g'|\Big) &\les& \ep_0,
\eeaa
\beaa
\sup_{\Si_*}\Big(r^2u^{1+\dec}|\dk^{k}\trXc'|+r^3u^{1+\dec}\left|\dk^{k}\left(\ov{\DD'}\c\Zc'+2\ov{\Pc'}\right)\right|\Big) &\les& \ep_0,
\eeaa 
\beaa
\sup_{\Si_*}\Big(r^{\frac{7}{2}+\dee}|\dk^{k}B'|+r^4u^{\frac{1}{2}+\dec}|\dk^{k-1}\nab_3B'|\Big) &\les& \ep_0,
\eeaa 
\beaa
\sup_{\Si_*}r^5u^{1+\dec}\left|\left[\left(\ov{\DD'}\c -\frac{a}{2}\ov{\Jk}\c\nab_{e_4'} -\frac{a}{2}\ov{\Jk}\c\nab_{e_3'}\right)\c\left( B'   -\frac{3a}{2}\ov{\Pc'}\Jk - \frac{a}{4} \ov{\Jk}\c A'\right)\right]_{\ell=1}\right| &\les& \ep_0,
\eeaa
and
\beaa
\sup_{\Si_*}r^5u^{1+\dec}\left|\left[\ov{\DD}\c\Lieb_{\T'} B' \right]_{\ell=1}\right| &\les& \ep_0.
\eeaa
\end{proposition} 
 
The proof of Proposition \ref{prop:improvedesitmatesfortemporalframeofMextonSigmastar} 
is done in section \ref{sec:proofofprop:improvedesitmatesfortemporalframeofMextonSigmastar}, relying on the estimates of sections \ref{sec:firstdecayesitmatesPGframeonSigmastar:chap5},  \ref{sec:decayesitmatesPGframeonSigmastar:chap5} and   \ref{sec:additionaldecayesitmatesPGframeonSigmastar:chap5}.

 %%%%%%%%%%%%%%%%%%%%%%%%%%%%%%%%%%%%%%%%%%%%%%%%%%

\subsection{First decay estimates for the PG frame on $\Si_*$ }
\lab{sec:firstdecayesitmatesPGframeonSigmastar:chap5}

%%%%%%%%%%%%%%%%%%%%%%%%%%%%%%%%%%%%%%%%%%%%%%%%%%

\begin{lemma}
We have
\beaa
e_4'(r')=1, \qquad e_4'(u')=e_4'(\th')=e_4'(\vphi')=0, \qquad \nab'(r')=0, 
\eeaa
and 
\beaa
\sup_{\Si_*}ru^{1+\dec}\Big(\left|\dk^{\leq k_*-12}\widecheck{\DD'(\cos\th')}\right|+\left|\dk^{\leq k_*-12}\widecheck{\DD'(u')}\right|+\left|\dk^{\leq k_*-12}\widecheck{\DD'({J'}^{(\pm)})}\right|\Big)\\
+\sup_{\Si_*}ru^{1+\dec}\Big(\left|\dk^{\leq k_*-12}e_3'(\cos\th')\right|+\left|\dk^{\leq k_*-12}\widecheck{e_3'({J'}^{(\pm)})}\right|\Big)\\
+\sup_{\Si_*}u^{1+\dec}\Big(\left|\dk^{\leq k_*-12}\widecheck{e_3'(r')}\right|+\left|\dk^{\leq k_*-12}\widecheck{e_3'(u')}\right|\Big) &\les& \ep_0.
\eeaa 
\end{lemma}

\begin{proof}
The identities  $e_4'(r')=1$, and $e_4'(u')=e_4'(\th')=e_4'(\vphi')=0$ are true  in $\Mext$ (and hence on $\Si_*$) by definition. In particular, since $\nab(r')=\nab(u')=0$ on $\Si_*$, this implies\footnote{Since $r'=r$, $u'=u$, $\nab(r)=0$ and $\nab(u)=0$ on $\Si_*$, and since  $\nab$ is tangent to $\Si_*$, we have indeed  $\nab(r')=0$ and $\nab(u')=0$, $\nab(\cos\th')=\nab(\cos\th)$ and $\nab({J'}^{(\pm)})=\nab(J^{(\pm)})$ on $\Si_*$.}
\beaa
1 &=& e_4(r')  +\frac 1 4 |f|^2  e_3(r'),\\
0 &=& e_4(u')  +\frac 1 4 |f|^2  e_3(u'),\\
0 &=& e_4(\cos\th') + f\c\nab(\cos\th) +\frac 1 4 |f|^2  e_3(\cos\th'),\\
0 &=& e_4({J'}^{(\pm)}) + f\c\nab(J^{(\pm)}) +\frac 1 4 |f|^2  e_3({J'}^{(\pm)}).
\eeaa
Since $e_3=\nu-b_*e_4$, and since $\nu(r')=\nu(r)$, $\nu(u')=\nu(u)$, $\nu(\cos\th')=0$ and $\nu({J'}^{(\pm)})=0$, we infer, using also $\nu(u)=-\nu(r)$, $f=\frac{a}{r}f_0$, and $|f|^2=\frac{a^2(\sin\th)^2}{r^2}$, 
\beaa
e_4(r') &=& \frac{1-\frac{\nu(r)a^2(\sin\th)^2}{4r^2}}{1-\frac{b_*a^2(\sin\th)^2}{4r^2}},\\
e_4(u') &=& \frac{\frac{\nu(r)a^2(\sin\th)^2}{4r^2}}{1-\frac{b_*a^2(\sin\th)^2}{4r^2}},\\
e_4(\cos\th')  &=&  -\frac{\frac{a}{r}f_0\c\nab(\cos\th)}{1-\frac{b_*a^2(\sin\th)^2}{4r^2}},\\
e_4({J'}^{(\pm)}) &=& -\frac{\frac{a}{r}f_0\c\nab(J^{(\pm)})}{1-\frac{b_*a^2(\sin\th)^2}{4r^2}}.
\eeaa
Since we have
\beaa
&&\nu(r)=-2+\Ga_b, \qquad b_*=-1-\frac{2m}{r}+r\Ga_b, \qquad \nab(\cos\th)=-\frac{1}{r}\dual f_0+\widecheck{\nab J^{(0)}}, \\
 &&\nab J^{(\pm)}=\frac{1}{r}f_\pm +\widecheck{ \nab J^{(\pm)}},
\eeaa
we infer, in view of Definition \ref{def:Gammabtilde} for $\widetilde{\Ga}_b$,
\beaa
e_4(r') &=& 1+O(r^{-2})+r^{-1}\widetilde{\Ga}_b,\\
e_4(u') &=& O(r^{-2})+r^{-1}\widetilde{\Ga}_b,\\
e_4(\cos\th')  &=&  r^{-1}\widetilde{\Ga}_b,\\
e_4({J'}^{(\pm)}) &=& O(r^{-2})+r^{-1}\widetilde{\Ga}_b.
\eeaa

Also, in view of the change of the definition \eqref{eq:changofframedefiningthefrmaeofMext:chap5} of the frame of $\Mext$, and using again  $\nab(r')=0$ on $\Si_*$, we have
\beaa
 \nab'(r') &=& \frac 1 2  \fb  e_4(r') +\left(\frac 1 2 f +\frac{1}{8}|f|^2\fb\right)   e_3(r')\\
 &=&  \frac 1 2\left(e_4(r')+\frac{1}{4}|f|^2e_3(r')\right)\fb+\frac 1 2 f e_3(r').
\eeaa
Together with \eqref{eq:defoffandfbforthedefintionoftheframeofMext:chapp5}, we infer
\beaa
 \nab'(r) &=& \frac 1 2\left[-\frac{(\nu(r)-b_*)}{1-\frac{1}{4}b_* |f|^2}\left(e_4(r')+\frac{1}{4}|f|^2e_3(r')\right)+e_3(r')\right] f. 
\eeaa
In view of the above, we have
\beaa
e_4(r')+\frac{1}{4}|f|^2e_3(r') =1, \qquad e_3(r')=\nu(r)-b_*e_4(r')
\eeaa
and hence
\beaa
 \nab'(r) &=& \frac 1 2\left[-\frac{(\nu(r)-b_*)}{1-\frac{1}{4}b_* |f|^2}+\nu(r)-b_*e_4(r')\right] f. 
\eeaa
Using again the above, we have
\beaa
e_4(r') &=& \frac{1-\frac{\nu(r)|f|}{4}}{1-\frac{b_*|f|^2}{4}}
\eeaa
and hence
\beaa
 \nab'(r) &=& \frac 1 2\left[-\frac{(\nu(r)-b_*)}{1-\frac{1}{4}b_* |f|^2}+\nu(r)-b_*\frac{1-\frac{\nu(r)|f|}{4}}{1-\frac{b_*|f|^2}{4}}\right] f = 0
\eeaa
as stated. 

Next, we focus on deriving the stated estimates. First, using the above identities for $e_4(r')$, $e_4(u')$, $e_4(\cos\th')$ and $e_4({J'}^{(\pm)})$, we have
\beaa
e_3(r') &=& \nu(r')-b_*e_4(r')=\nu(r)-b_*\Big(1+O(r^{-2})+r^{-1}\Ga_b\Big)=-\Up+O(r^{-2})+r^{-1}\widetilde{\Ga}_b,\\
e_3(u') &=& \nu(u')-b_*e_4(u')=\nu(u)-b_*\Big(O(r^{-2})+r^{-1}\Ga_b\Big)=2+O(r^{-2})+r^{-1}\widetilde{\Ga}_b,\\ 
e_3(\cos\th') &=& \nu(\cos\th')-b_*e_4(\cos\th')=\nu(\cos\th)-b_*O(r^{-1})\widecheck{\nab J^{(0)}}=r^{-1}\widetilde{\Ga}_b,\\
e_3({J'}^{(\pm)}) &=& \nu({J'}^{(\pm)})-b_*e_4({J'}^{(\pm)})=\nu(J^{(\pm)})-b_*\Big(O(r^{-2})+O(r^{-1})\widecheck{\nab J^{(\pm)}}\Big)\\
&=& O(r^{-2})+r^{-1}\widetilde{\Ga}_b.
\eeaa
Together with \eqref{eq:defoffandfbforthedefintionoftheframeofMext:chapp5} and the above identities for the $e_4$ derivatives, we infer
\beaa
\widecheck{e_3'(r')} &=& O(r^{-2})+r^{-1}\widetilde{\Ga}_b,\\
\widecheck{e_3'(u')} &=& O(r^{-2})+r^{-1}\widetilde{\Ga}_b,\\ 
e_3'(\cos\th') &=& r^{-1}\widetilde{\Ga}_b,\\
e_3'({J'}^{(\pm)}) &=&  O(r^{-2})+r^{-1}\widetilde{\Ga}_b.
\eeaa
Together with the dominance condition \eqref{eq:dominantconditiononronSigmastarchap5} on $r$ on $\Si_*$, and  the estimates for $\widetilde{\Ga}_b$ of Corollary \ref{cor:estimatesforf0fplusfminusandJp:Gammatildeb}, we obtain
\beaa
\sup_{\Si_*}ru^{1+\dec}\Big(\left|\dk^{\leq k_*-12}e_3'(\cos\th')\right|+\left|\dk^{\leq k_*-12}\widecheck{e_3'({J'}^{(\pm)})}\right|\Big)\\
+\sup_{\Si_*}u^{1+\dec}\Big(\left|\dk^{\leq k_*-12}\widecheck{e_3'(r')}\right|+\left|\dk^{\leq k_*-12}\widecheck{e_3'(u')}\right|\Big) &\les& \ep_0.
\eeaa 

Finally, we control the angular derivatives. We have in view of \eqref{eq:defoffandfbforthedefintionoftheframeofMext:chapp5} and the above identities for the $e_4$ derivatives and $e_3$ derivatives 
\beaa
\nab'(u') &=& O(r^{-1})e_4(u')+\frac{1}{2}f\Big(1+O(r^{-2})\Big)e_3(u')\\
&=& \frac{a}{r}f_0+O(r^{-3})+r^{-2}\widetilde{\Ga}_b,\\
\nab'(\cos\th') &=& \Big(1+O(r^{-2})\Big)\nab(\cos\th)+O(r^{-1})e_3(\cos\th')+O(r^{-1})e_4(\cos\th')\\
&=& -\frac{1}{r}\dual f_0+\widetilde{\Ga}_b,\\
\nab'({J'}^{(\pm)}) &=& \Big(1+O(r^{-2})\Big)\nab(J^{(\pm)})+O(r^{-1})e_3({J'}^{(\pm)})+O(r^{-1})e_4({J'}^{(\pm)})\\
&=& \frac{1}{r}f_\pm+\widetilde{\Ga}_b+O(r^{-3}),
\eeaa
and hence
\beaa
\widecheck{\nab'(u')} &=& O(r^{-3})+r^{-2}\widetilde{\Ga}_b,\\
\widecheck{\nab'(\cos\th')} &=&  O(r^{-3})+\widetilde{\Ga}_b,\\
\widecheck{\nab'({J'}^{(\pm)})} &=&  O(r^{-3})+\widetilde{\Ga}_b.
\eeaa
Together with the dominance condition \eqref{eq:dominantconditiononronSigmastarchap5} on $r$ on $\Si_*$, and the estimates for  $\widetilde{\Ga}_b$ of Corollary \ref{cor:estimatesforf0fplusfminusandJp:Gammatildeb}, we obtain 
\beaa
\sup_{\Si_*}ru^{1+\dec}\Big(\left|\dk^{\leq k_*-12}\widecheck{\DD'(\cos\th')}\right|+\left|\dk^{\leq k_*-12}\widecheck{\DD'(u')}\right|+\left|\dk^{\leq k_*-12}\widecheck{\DD'({J'}^{(\pm)})}\right|\Big) &\les& \ep_0.
\eeaa
This concludes the proof of the lemma.
\end{proof}

\begin{lemma}\lab{lemma:controloffpfbplapchangefrmaePGMexttointegSigmastar:chap5}
Consider the change of frame coefficients $(f', \fb', \la')$ from the PG frame $(e_3', e_4', e_1', e_2')$ of $\Mext$ to the frame $(e_3, e_4, e_1, e_2)$ on $\Si_*$, i.e. 
 \bea\lab{eq:inverseofchangofframedefiningthefrmaeofMext:chap5}
 \bsplit
  e_4&=\la'\left(e'_4 + f_b'  e'_b +\frac 1 4 |f'|^2  e'_3\right),\\
  e_a&= \left(\de_a^b  +\frac{1}{2}\fb'_a f'^b \right) e'_b +\frac 1 2  \fb'_a  e'_4 +\left(\frac 1 2 f'_a +\frac{1}{8}|f'|^2\fb'_a\right)   e'_3,\\
 e_3&=(\la')^{-1}\left( \left(1+\frac{1}{2}f'\c\fb'  +\frac{1}{16} |f'|^2  |\fb'|^2\right) e'_3 + \left(\fb'^b+\frac 1 4 |\fb'|^2f'^b\right) e'_b  + \frac 1 4 |\fb'|^2 e'_4 \right).
\end{split}
\eea
Then, we have
 \beaa
\bsplit
\la' &= 1+O(r^{-2})+r^{-1}\Ga_b,\\
f'  &=-\frac{a}{r}\Big(1+O(r^{-2})+r^{-1}\Ga_b\Big)f_0,\\
\fb' &= -\frac{a\Up}{r}\Big(1+O(r^{-2})+r^{-1}\Ga_b\Big)f_0.
\end{split}
\eeaa
\end{lemma}

\begin{proof}
Recall that $(f, \fb, \la)$, with $\la=1$ and $(f, \fb)$ given by \eqref{eq:defoffandfbforthedefintionoftheframeofMext:chapp5}, are  the change of frame coefficients from the frame $(e_3, e_4, e_1, e_2)$ to the frame $(e_3', e_4', e_1', e_2')$, so that $(f', \fb', \la')$ correspond to the inverse transformation of the one of $(f, \fb, \la)$. Thus, according to  \eqref{relations:laffb-to-primes}, $(f', \fb', \la')$ is related to the transition coefficients $(f, \fb, \la)$ by
\beaa
\bsplit
\la' &= \la^{-1} \left(1+\frac{1}{2}f\c\fb  +\frac{1}{16} |f|^2  |\fb|^2\right),\\
f_a'  &= -\frac{\la}{1+\frac{1}{2}f\c\fb  +\frac{1}{16} |f|^2  |\fb|^2}\left(f_a +\frac{1}{4}|f|^2\fb_a\right),\\
\fb_a' &= -\la^{-1}\left(\fb_a+\frac 1 4 |\fb|^2f_a\right).
\end{split}
\eeaa
Since  $(f, \fb)$ are given by \eqref{eq:defoffandfbforthedefintionoftheframeofMext:chapp5}, and since $|f_0|^2=(\sin\th)^2$, $\nu(r)=-2+r\Ga_b$ and $b_*=-1-\frac{2m}{r}+r\Ga_b$,  we have
\beaa
f &=& \frac{a}{r}f_0,  \\
\fb &=& -\frac{(\nu(r)-b_*)}{1-\frac{1}{4}b_*\frac{a^2(\sin\th)^2}{r^2}}\frac{a}{r}f_0=\frac{a\Up}{r}\left(1+O(r^{-2})+r\Ga_b\right)f_0,\\
f\c\fb &=& -\frac{(\nu(r)-b_*)}{1-\frac{1}{4}b_*\frac{a^2(\sin\th)^2}{r^2}}\frac{a^2(\sin\th)^2}{r^2}=\frac{a^2(\sin\th)^2\Up}{r^2}\left(1+O(r^{-2})+r\Ga_b\right),\\
|f|^2 &=& \frac{a^2(\sin\th)^2}{r^2},\\ 
|\fb|^2 &=& \frac{(\nu(r)-b_*)^2}{\left(1-\frac{1}{4}b_*\frac{a^2(\sin\th)^2}{r^2}\right)^2}\frac{a^2(\sin\th)^2}{r^2}=\frac{a^2(\sin\th)^2\Up^2}{r^2}\left(1+O(r^{-2})+r\Ga_b\right),
\eeaa
and hence, using also $\la=1$, we infer
\beaa
\bsplit
\la' &= 1+O(r^{-2})+r^{-1}\Ga_b,\\
f'  &=-\frac{a}{r}\Big(1+O(r^{-2})+r^{-1}\Ga_b\Big)f_0,\\
\fb' &= -\frac{a\Up}{r}\Big(1+O(r^{-2})+r^{-1}\Ga_b\Big)f_0,
\end{split}
\eeaa
as stated.
\end{proof}

\begin{lemma}\lab{lemma:simplepropertyderivativeOofporwerr}
We have
\beaa
\nab\left(O(r^{-j})\right)=O(r^{-j-1})+O(r^{-j})\widetilde{\Ga}_b,\\
\nu\left(O(r^{-j})\right)=O(r^{-j-1})+O(r^{-j})\widetilde{\Ga}_b,
\eeaa
where the notation $\widetilde{\Ga}_b$ has been introduced in Definition \ref{def:Gammabtilde}.
\end{lemma}
 
\begin{proof}
Recall that by $O(r^{-j})$, we mean a function $h(r,\cos\th)$ such that 
\beaa
|(r\pr_r)^k(\pr_{\cos\th})^lh|\les \frac{C_{k, l}}{r^j}, \quad\textrm{as}\quad r\to+\infty.
\eeaa
The proof is then an immediate consequence of the fact that $\nu(r)=-2+r\Ga_b$ and 
$\nab(\cos\th)=-\frac{1}{r}\dual f_0+\widetilde{\Ga}_b$.
\end{proof}

 \begin{lemma}\lab{lemma:computationfirstorderderivativeinchangeofframefromintegrabeltoprincial}
 We have
  \beaa
&&\div(f') =r^{-1}\dk_*^{\leq 1}\widetilde{\Ga}_b+O\left(\frac{1}{r^3}\right), \qquad \curl(f')= -\frac{2a\cos\th}{r^2}+r^{-1}\dk_*^{\leq 1}\widetilde{\Ga}_b+O\left(\frac{1}{r^3}\right), \\ 
&&\nab\hot f'=r^{-1}\dk_*^{\leq 1}\widetilde{\Ga}_b+O\left(\frac{1}{r^3}\right),\qquad \div(\fb') =r^{-1}\dk_*^{\leq 1}\widetilde{\Ga}_b+O\left(\frac{1}{r^3}\right), \\
&& \curl(\fb')= -\frac{2a\cos\th\Up}{r^2}+r^{-1}\dk_*^{\leq 1}\widetilde{\Ga}_b+O\left(\frac{1}{r^3}\right), \qquad \nab\hot(\fb')=r^{-1}\dk_*^{\leq 1}\widetilde{\Ga}_b+O\left(\frac{1}{r^3}\right),\\
&& \nab\la' = O(r^{-3})+r^{-2}\dk_*^{\leq 1}\widetilde{\Ga}_b.
\eeaa  
Also, we have
\beaa
\nab_{\nu}f' &=& \frac{2}{r}f'+r^{-1}\Ga_b+O\left(\frac{1}{r^3}\right),\\
\nab_{\nu}\fb' &=& \frac{2}{r}\fb'+r^{-1}\dk_*^{\leq 1}\Ga_b+O\left(\frac{1}{r^3}\right),\\
\nu(\la') &=& O(r^{-3})+r^{-1}\dk_*^{\leq 1}\widetilde{\Ga}_b,
\eeaa
where the notation $\widetilde{\Ga}_b$ has been introduced in Definition \ref{def:Gammabtilde}.
 \end{lemma}
 
\begin{proof}
Recall that 
\beaa
\bsplit
\la' &= 1+O(r^{-2})+r^{-1}\Ga_b,\\
f'  &=-\frac{a}{r}\Big(1+O(r^{-2})+r^{-1}\Ga_b\Big)f_0,\\
\fb' &= -\frac{a\Up}{r}\Big(1+O(r^{-2})+r^{-1}\Ga_b\Big)f_0.
\end{split}
\eeaa
The proof follows then immediately from the definition of $\widetilde{\Ga}_b$, Lemma \ref{lemma:simplepropertyderivativeOofporwerr}, and the fact that $\nu(r)=-2+r\Ga_b$ and $\nab_\nu f_0=0$. 
\end{proof}

%%%%%%%%%%%%%%%%%%%%%%%%%%%%%%%%%%%%%%%%%%%%%%%%%%

\subsection{Decay estimates for the  PG frame on $\Si_*$}
\lab{sec:decayesitmatesPGframeonSigmastar:chap5}

%%%%%%%%%%%%%%%%%%%%%%%%%%%%%%%%%%%%%%%%%%%%%%%%%%
 
We start with the control of the Ricci coefficients of the PG frame of $\Mext$ on $\Si_*$. 
\begin{lemma}\lab{lemma:controlofRiccicoeffPGframeMextonSigmastar}
We have on $\Si_*$, for $k\leq k_*-13$, 
\bea
\bsplit
\left|\dk_*^k\Big(\chibh', \, \xib', \, \ombc', \, \widecheck{\eta}'\Big)\right| &\les  \frac{\ep_0}{ru^{1+\dec}},\\
\left|\dk_*^k\Big(\chih', \, \trchbc', \, \widecheck{\atrchb}' , \, \widecheck{\ze}'\Big)\right| &\les \frac{\ep_0}{r^2u^{\frac{1}{2}+\dec}},\\
\left|\dk_*^{k-1}\nab_\nu\Big(\chih', \, \trchbc', \, \widecheck{\atrchb}' , \, \widecheck{\ze}'\Big)\right| &\les \frac{\ep_0}{r^2u^{1+\dec}},\\
\left|\dk_*^k\Big(\trchc', \, \widecheck{\atrch}'\Big)\right| &\les \frac{\ep_0}{r^2u^{1+\dec}}.
\end{split}
\eea
\end{lemma}
 
\begin{proof} 
We consider the frame transformation  from the frame $(e_1', e_2', e_3', e_4')$ of $\Mext$ to the frame $(e_1, e_2, e_3, e_4)$ of $\Si_*$, with corresponding change of frame coefficients $(f', \fb', \la')$. Using the transformation formulas of Proposition \ref{Proposition:transformationRicci}, Lemma \ref{lemma:controloffpfbplapchangefrmaePGMexttointegSigmastar:chap5} on the control of $(f', \fb', \la')$, and the fact that $\atrch=\atrchb=0$, $\xi'=0$, $\om'=0$, and $\etab'=-\ze'$, we have  
\beaa
\trch &=& \trch'  +  \div(f') +O(r^{-3})+O(r^{-1})\widecheck{\eta}'+r^{-1}\Ga_g'+r^{-2}\Ga_b,\\
0 &=& \atrch'  +  \curl(f')  +O(r^{-3})+O(r^{-1})\widecheck{\eta}'+r^{-1}\Ga_g'+r^{-2}\Ga_b,\\
\chih &=& \chih'  +  \nab\hot f' +O(r^{-3})+O(r^{-1})\widecheck{\eta}'+r^{-1}\Ga_g'+r^{-2}\Ga_b,
\eeaa
\beaa
\trchb &=& \trchb' +\div(\fb')   +O(r^{-3})+O(r^{-1})\xib'+r^{-1}\Ga_g'+r^{-2}\Ga_b,\\
0 &=& \atrchb' +\curl(\fb') +O(r^{-3})+O(r^{-1})\xib'+r^{-1}\Ga_g'+r^{-2}\Ga_b,\\
\chibh &=& \chibh' +\nab\hot\fb'  +O(r^{-3})+r^{-1}\Ga_b'+r^{-2}\Ga_b,
\eeaa
and
\beaa
\ze &=& \ze' -\nab(\log\la')  -\frac{1}{4}\trchb' f'  +\frac{1}{4}\fb'\trch' +  \frac{1}{4}\fb'\div(f')  + \frac{1}{4}\dual\fb' \curl(f')\\
&&+O(r^{-1})(\ombc',\chibh')+r^{-1}\Ga_g' +O(r^{-3})+r^{-1}\Ga_b'+r^{-2}\Ga_b+r^{-1}\Ga_g.
\eeaa
Together with Lemma \ref{lemma:computationfirstorderderivativeinchangeofframefromintegrabeltoprincial} on the control of first order derivatives of $(f', \fb', \la')$, and Corollary \ref{corollary:Kerr.PG.asymptotics} on the asymptotic for $r$ large of the Kerr values of the PG frame, we infer 
\beaa
\trchc &=& \trchc'  +O(r^{-3})+O(r^{-1})\widecheck{\eta}'+r^{-1}\Ga_g'+r^{-1}\dk_*^{\leq 1}\widetilde{\Ga}_b,\\
0 &=& \widecheck{\atrch}'  +O(r^{-3})+O(r^{-1})\widecheck{\eta}'+r^{-1}\Ga_g'+r^{-1}\dk_*^{\leq 1}\widetilde{\Ga}_b,\\
\chih &=& \chih'  +O(r^{-3})+O(r^{-1})\widecheck{\eta}'+r^{-1}\Ga_g'+r^{-1}\dk_*^{\leq 1}\widetilde{\Ga}_b,
\eeaa
\beaa
\trchbc &=& \trchbc' +O(r^{-3})+O(r^{-1})\xib'+r^{-1}\Ga_g'+r^{-1}\dk_*^{\leq 1}\widetilde{\Ga}_b,\\
0 &=& \widecheck{\atrchb}' +O(r^{-3})+O(r^{-1})\xib'+r^{-1}\Ga_g'+r^{-1}\dk_*^{\leq 1}\widetilde{\Ga}_b,\\
\chibh &=& \chibh' +O(r^{-3})+r^{-1}\Ga_b'+r^{-1}\dk_*^{\leq 1}\widetilde{\Ga}_b,
\eeaa
and
\beaa
\ze &=& \widecheck{\ze}'  +O(r^{-3})+O(r^{-1})(\ombc',\chibh')+r^{-1}\Ga_g'+r^{-1}\dk_*^{\leq 1}\widetilde{\Ga}_b.
\eeaa

Next, we use again the transformation formulas of Proposition \ref{Proposition:transformationRicci}, summing the one of $\xib$ with the one for $\etab$ multiplied by $\la^2b_*$, summing the one of $\omb$ with the one for $\om$ multiplied by $-\la^2b_*$, and summing the one for $\eta$ with the one for $\xi$ multiplied by $\la^2b_*$. Proceeding as above, and using in addition the fact that $\nu=e_3+b_*e_4$, and the transversality conditions $\xi=0$, $\om=0$ and $\etab=-\ze$ for the frame of $\Si_*$, we infer
\beaa
\xib -b_*\ze &=& \xib' -b_*\ze' + \frac{1}{2}\nab_\nu\fb'  + \frac{1}{4}\trchb'\,\fb' +\frac{b_*}{4}\trchb' f'  +O(r^{-3})+r^{-1}\Gac'+r^{-2}\Ga_b,\\
\omb &=& \omb'+\frac{1}{2}\nu(\log\la')    +O(r^{-3})+r^{-1}\Gac'+r^{-2}\Ga_b,\\
\eta &=& \eta' +\frac{1}{2}\nab_\nu f'  +\frac{1}{4}\fb'\trch'  +\frac{b_*}{4}\trch' f' +O(r^{-3})+r^{-1}\Gac'+r^{-2}\Ga_b.
\eeaa
Together with Lemma \ref{lemma:computationfirstorderderivativeinchangeofframefromintegrabeltoprincial} on the control of first order derivatives of $(f', \fb', \la')$, and Corollary \ref{corollary:Kerr.PG.asymptotics} on the asymptotic for $r$ large of the Kerr values of the PG frame, we deduce
\beaa
\xib -b_*\ze &=& \xib' -\widecheck{\ze}'  +O(r^{-3})+r^{-1}\Gac'+r^{-1}\dk_*^{\leq 1}\widetilde{\Ga}_b,\\
\ombc &=& \ombc'   +O(r^{-3})+r^{-1}\Gac'+r^{-1}\dk_*^{\leq 1}\widetilde{\Ga}_b,\\
\eta &=& \widecheck{\eta}'  +O(r^{-3})+r^{-1}\Gac'+r^{-1}\dk_*^{\leq 1}\widetilde{\Ga}_b.
\eeaa
In view of the above, we infer
\beaa
\Big(\chibh', \, \xib', \, \ombc', \, \widecheck{\eta}'\Big) &=& O(r^{-3})+\Ga_g'+r^{-1}\Ga_b'+\dk_*^{\leq 1}\widetilde{\Ga}_b.
\eeaa
Together with the bootstrap assumptions for  $\Ga_g'$ and $\Ga_b'$, and Corollary  \ref{cor:estimatesforf0fplusfminusandJp:Gammatildeb} on the control of $\widetilde{\Ga}_b$, this yields on $\Si_*$, for $k\leq k_*-13$,
\beaa 
\left|\dk_*^k\Big(\chibh', \, \xib', \, \ombc', \, \widecheck{\eta}'\Big)\right| &\les& \frac{\ep_0}{ru^{1+\dec}}+\frac{\ep}{r^2}+\frac{1}{r^3},
\eeaa
and hence, together with the dominance condition \eqref{eq:dominantconditiononronSigmastarchap5} on $r$ on $\Si_*$, we obtain, for $k\leq k_*-13$, 
\beaa 
\left|\dk_*^k\Big(\chibh', \, \xib', \, \ombc', \, \widecheck{\eta}'\Big)\right| &\les&  \frac{\ep_0}{ru^{1+\dec}}.
\eeaa

Next, using again the above, we have
\beaa
\Big(\chih', \, \trchbc', \, \widecheck{\atrchb}' , \, \widecheck{\ze}'\Big) &=& O(r^{-3})+\Ga_g+O(r^{-1})\Big(\chibh', \, \xib', \, \ombc', \, \widecheck{\eta}'\Big)+r^{-1}\Ga_g'+r^{-1}\dk_*^{\leq 1}\widetilde{\Ga}_b,\\
\Big(\trchc', \, \widecheck{\atrch}'\Big) &=& O(r^{-3})+\trchc+O(r^{-1})\Big(\chibh', \, \xib', \, \ombc', \, \widecheck{\eta}'\Big)+r^{-1}\Ga_g'+r^{-1}\dk_*^{\leq 1}\widetilde{\Ga}_b.
\eeaa
Together with the above control for $\chibh'$, $\xib'$, $\ombc'$, and $\widecheck{\eta}'$, the control for $\Ga_g$ and $\trchc$ of Proposition \ref{prop:decayonSigamstarofallquantities}, the bootstrap assumptions for  $\Ga_g'$, and Corollary  \ref{cor:estimatesforf0fplusfminusandJp:Gammatildeb} on the control of $\widetilde{\Ga}_b$, this yields on $\Si_*$, for $k\leq k_*-13$,
\beaa
\left|\dk_*^k\Big(\chih', \, \trchbc', \, \widecheck{\atrchb}' , \, \widecheck{\ze}'\Big)\right| &\les& \frac{\ep_0}{r^2u^{\frac{1}{2}+\dec}}+\frac{1}{r^3},\\
\left|\dk_*^{k-1}\nab_\nu\Big(\chih', \, \trchbc', \, \widecheck{\atrchb}' , \, \widecheck{\ze}'\Big)\right| &\les& \frac{\ep_0}{r^2u^{1+\dec}}+\frac{1}{r^3},\\
\left|\dk_*^k\Big(\trchc', \, \widecheck{\atrch}'\Big)\right| &\les& \frac{\ep_0}{r^2u^{1+\dec}}+\frac{1}{r^3},
\eeaa
and hence, together with the dominance condition \eqref{eq:dominantconditiononronSigmastarchap5} on $r$ on $\Si_*$, we obtain, for $k\leq k_*-13$, 
\beaa
\left|\dk_*^k\Big(\chih', \, \trchbc', \, \widecheck{\atrchb}' , \, \widecheck{\ze}'\Big)\right| &\les& \frac{\ep_0}{r^2u^{\frac{1}{2}+\dec}},\\
\left|\dk_*^{k-1}\nab_\nu\Big(\chih', \, \trchbc', \, \widecheck{\atrchb}' , \, \widecheck{\ze}'\Big)\right| &\les& \frac{\ep_0}{r^2u^{1+\dec}},\\
\left|\dk_*^k\Big(\trchc', \, \widecheck{\atrch}'\Big)\right| &\les& \frac{\ep_0}{r^2u^{1+\dec}}.
\eeaa
This concludes the proof of the lemma.
\end{proof}

\begin{lemma}\lab{lemma:controlofcurvaturePGframeMextonSigmastar}
We have on $\Si_*$, for $k\leq k_*-14$, 
 \beaa 
r^{\frac{7}{2}+\dee}\left|\dk^k\b'\right|  +r^4u^{\frac{1}{2}+\dec}\left|\dk^{k-1}\nab_\nu\b'\right|+r^3u^{\frac{1}{2}+\dec}\left|\dk^k\Pc'\right|\\
+r^3u^{1+\dec}\left|\dk^{k-1}\nab_\nu\Pc'\right|+r^2u^{1+\dec}\left|\dk^k\bb'\right| &\les& \ep_0,
\eeaa
and 
\beaa
\left|\dk_*^k\left(\b' - \left(\b - \frac{3am}{r^4}f_0\right) - \frac{3a}{2r}\left(  \rhoc\,' f_0+\widecheck{\rhod}'\dual  f _0 \right) - \frac{a}{2r}f_0\c\a'\right)\right| &\les& \frac{\ep_0}{r^4u^{1+\dec}}.
\eeaa
\end{lemma}

\begin{proof}
We consider the frame transformation  from the frame $(e_1', e_2', e_3', e_4')$ of $\Mext$ to the frame $(e_1, e_2, e_3, e_4)$ of $\Si_*$, with corresponding change of frame coefficients $(f', \fb', \la')$. Using the transformation formulas of Proposition \ref{Proposition:transformationRicci}, and Lemma \ref{lemma:controloffpfbplapchangefrmaePGMexttointegSigmastar:chap5} on the control of $(f', \fb', \la')$, we have  
   \beaa
\b &=& \b' +\frac 3 2\left( -\frac{a}{r}f_0\left(-\frac{2m}{r^3}+\rhoc'\right) -\frac{a}{r}\dual  f_0 \widecheck{\rhod}'\right)-\frac{a}{2r}f_0\c\a'+O(r^{-2})\bb'+O(r^{-3})\aa'\\
&&+O(r^{-3})\Ga_g'+O(r^{-5})+r^{-3}\Ga_b,\\
  \bb &=& \bb' +O(r^{-1})\aa'  +O(r^{-4}) +r^{-2}\Ga_g'+r^{-3}\Ga_b,\\
 \rho &=& \rho' +O(r^{-4})+r^{-2}\Ga_g'+O(r^{-1})\bb'+O(r^{-2})\aa',\\
  \rhod &=& \rhod' +O(r^{-4})+r^{-2}\Ga_g'+O(r^{-1})\bb'+O(r^{-2})\aa'.
  \eeaa 
In particular, we have
\beaa
 \bb' &=&  O(r^{-1})\aa'  +O(r^{-4}) +r^{-2}\Ga_g'+r^{-1}\Ga_b,
\eeaa
which together with the control of Theorem M3 for $\aa'$,  the control  of Proposition \ref{Prop.Flux-bb-vthb-eta-xib} for $\Ga_b$, and the bootstrap assumptions for  $\Ga_g'$, this yields on $\Si_*$, for $k\leq k_*-10$,
\beaa
|\dk_*^k\bb'| &\les& \frac{\ep_0}{r^2u^{1+\dec}}+\frac{1}{r^4},
\eeaa
and hence, together with the dominance condition \eqref{eq:dominantconditiononronSigmastarchap5} on $r$ on $\Si_*$, we obtain, for $k\leq k_*-10$,
\beaa
|\dk_*^k\bb'| &\les& \frac{\ep_0}{r^2u^{1+\dec}}.
\eeaa

Next, we rewrite the above identities for $\rho$ and $\rhod$ as
\beaa
\rhoc' &=& O(r^{-4})+r^{-2}\Ga_g'+O(r^{-1})\bb'+O(r^{-2})\aa'+r^{-1}\Ga_g,\\
  \widecheck{\rhod}' &=& O(r^{-4})+r^{-2}\Ga_g'+O(r^{-1})\bb'+O(r^{-2})\aa'+r^{-1}\Ga_g.
\eeaa 
which together with the above control of $\bb'$, the control of Theorem M3 for $\aa'$, the control for $\Ga_g$  of Proposition \ref{prop:decayonSigamstarofallquantities}, and the bootstrap assumptions for  $\Ga_g'$, this yields on $\Si_*$, for $k\leq k_*-10$, 
\beaa
|\dk_*^k(\rhoc',  \widecheck{\rhod}')| &\les& \frac{\ep_0}{r^3u^{\frac{1}{2}+\dec}}+\frac{1}{r^4},\\
|\dk_*^{k-1}\nab_\nu(\rhoc',  \widecheck{\rhod}')| &\les& \frac{\ep_0}{r^3u^{1+\dec}}+\frac{1}{r^4},\\
\eeaa
and hence, together with the dominance condition \eqref{eq:dominantconditiononronSigmastarchap5} on $r$ on $\Si_*$, we obtain, for $k\leq k_*-12$,
\beaa
|\dk_*^k(\rhoc',  \widecheck{\rhod}')| &\les& \frac{\ep_0}{r^3u^{\frac{1}{2}+\dec}},\\
|\dk_*^{k-1}\nab_\nu(\rhoc',  \widecheck{\rhod}')| &\les& \frac{\ep_0}{r^3u^{1+\dec}}.
\eeaa

Next, we rewrite the above identity for $\b$ as
   \beaa
\b' &=& \left(\b - \frac{3am}{r^4}f_0\right) +\frac{3a}{2r}\left(  \rhoc\,' f_0+\widecheck{\rhod}'\dual  f _0 \right) +\frac{a}{2r}f_0\c\a'+O(r^{-2})\bb'+O(r^{-3})\aa'\\
&&+O(r^{-3})\Ga_g'+O(r^{-5})+r^{-3}\Ga_b,
\eeaa
which together with the above control of $\bb'$, the control of Theorem M3 for $\aa'$, the control for $\Ga_g$  of Proposition \ref{prop:decayonSigamstarofallquantities}, and the bootstrap assumptions for  $\Ga_g'$, this yields on $\Si_*$, for $k\leq k_*-12$, 
\beaa
\left|\dk_*^k\left(\b' - \left(\b - \frac{3am}{r^4}f_0\right) - \frac{3a}{2r}\left(  \rhoc\,' f_0+\widecheck{\rhod}'\dual  f _0 \right) - \frac{a}{2r}f_0\c\a'\right)\right| &\les& \frac{\ep_0}{r^4u^{1+\dec}}+\frac{1}{r^5},
\eeaa
and hence, together with the dominance condition \eqref{eq:dominantconditiononronSigmastarchap5} on $r$ on $\Si_*$, we obtain, for $k\leq k_*-12$,
\beaa
\left|\dk_*^k\left(\b' - \left(\b - \frac{3am}{r^4}f_0\right) - \frac{3a}{2r}\left(  \rhoc\,' f_0+\widecheck{\rhod}'\dual  f _0 \right) - \frac{a}{2r}f_0\c\a'\right)\right| &\les& \frac{\ep_0}{r^4u^{1+\dec}}.
\eeaa
In particular, we have, for $k\leq k_*-12$,
\beaa
\left|\dk_*^k\b'\right| &\les& \left|\dk_*^k\left(\b - \frac{3am}{r^4}f_0\right)\right| +r^{-1}\left|\dk_*^k\left( \rhoc\,',\widecheck{\rhod}', \a'\right)\right| + \frac{\ep_0}{r^4u^{1+\dec}},\\
\left|\dk_*^{k-1}\nab_\nu\b' \right| &\les& \left|\dk_*^{k-1}\nab_\nu\left(\b - \frac{3am}{r^4}f_0\right)\right| +r^{-1}\left|\dk_*^k\left( \rhoc\,',\widecheck{\rhod}', \a'\right)\right| + \frac{\ep_0}{r^4u^{1+\dec}}\\
&\les& \left|\dk_*^{k-1}\nab_\nu\b\right|+r^{-1}\left|\dk_*^k\left( \rhoc\,',\widecheck{\rhod}', \a'\right)\right| +\frac{1}{r^5} + \frac{\ep_0}{r^4u^{1+\dec}},
\eeaa
were we also used the fact that $\nab_\nu f_0=0$, $\nu(r)=-2+r\Ga_b$, and the estimates {\bf Ref 1} for $\Ga_b$. 
Together with the estimate of Corollary \ref{cor:lastmissingimproveddecayestimateforbetaonSigmastar} for $\b-\frac{3am\sin\th}{r^4}f_0$,  the estimates of Proposition \ref{prop:decayonSigamstarofallquantities} for $\nab_\nu\b$, the above estimates for $\rhoc\,'$ and $\widecheck{\rhod}'$, the control of Theorem M2 for $\a'$, 
and  the dominance condition \eqref{eq:dominantconditiononronSigmastarchap5} on $r$ on $\Si_*$, we infer, for $k\leq k_*-14$, 
\beaa
\left|\dk_*^k\b'\right| &\les&  \frac{\ep_0}{r^{\frac{7}{2}+\dee}},\\
\left|\dk_*^{k-1}\nab_\nu\b'\right| &\les&  \frac{\ep_0}{r^4u^{\frac{1}{2}+\dec}},
\eeaa
which concludes the proof of the lemma.
\end{proof}

%%%%%%%%%%%%%%%%%%%%%%%%%%%%%%%%%%%%%%%%%%%%%%%%%%

\subsection{Additional decay estimates  on $\Si_*$}
\lab{sec:additionaldecayesitmatesPGframeonSigmastar:chap5}

%%%%%%%%%%%%%%%%%%%%%%%%%%%%%%%%%%%%%%%%%%%%%%%%%%

In this section, we prove the remaining estimates of Proposition \ref{prop:improvedesitmatesfortemporalframeofMextonSigmastar}. We start with the following lemma. 
\begin{lemma}\lab{cor:controloftheell=1modeofrenormalizedBonSigmarstar}
We have
\beaa
\sup_{\Si_*}r^5u^{1+\dec}\left|\left[\left(\ov{\DD'}\c -\frac{a}{2}\ov{\Jk}\c\nab_{e_4'} -\frac{a}{2}\ov{\Jk}\c\nab_{e_3'}\right)\c\left( B'   -\frac{3a}{2}\ov{\Pc'}\Jk - \frac{a}{4} \ov{\Jk}\c A'\right)\right]_{\ell=1}\right| &\les& \ep_0.
\eeaa
\end{lemma}

\begin{proof}
Recall from Lemma \ref{lemma:controlofcurvaturePGframeMextonSigmastar} that we have on $\Si_*$, for $k\leq k_*-14$, 
\beaa
\left|\dk_*^k\left(\b' - \left(\b - \frac{3am}{r^4}f_0\right) - \frac{3a}{2r}\left(  \rhoc\,' f_0+\widecheck{\rhod}'\dual  f _0 \right) - \frac{a}{2r}f_0\c\a'\right)\right| &\les& \frac{\ep_0}{r^4u^{1+\dec}},
\eeaa
and hence, in particular, 
\beaa
\left|\ddd_1\left(\b'  - \frac{3a}{2r}\left(  \rhoc\,' f_0+\widecheck{\rhod}'\dual  f _0 \right) - \frac{a}{2r}f_0\c\a'\right)-\ddd_1\left(\b - \frac{3am}{r^4}f_0\right)\right| &\les& \frac{\ep_0}{r^5u^{1+\dec}}.
\eeaa
Taking the $\ell=1$ mode, we infer
\beaa
\left|\left[\ddd_1\left(\b'  - \frac{3a}{2r}\left(  \rhoc\,' f_0+\widecheck{\rhod}'\dual  f _0 \right) - \frac{a}{2r}f_0\c\a'\right)\right]_{\ell=1}\right| &\les& \left|\left[\ddd_1\left(\b - \frac{3am}{r^4}f_0\right)\right]_{\ell=1}\right|\\
&&+\frac{\ep_0}{r^5u^{1+\dec}}.
\eeaa
Together with Corollary \ref{cor:lastmissingimproveddecayestimateforbetaonSigmastar}, we deduce 
\beaa
\left|\left[\ddd_1\left(\b'  - \frac{3a}{2r}\left(  \rhoc\,' f_0+\widecheck{\rhod}'\dual  f _0 \right) - \frac{a}{2r}f_0\c\a'\right)\right]_{\ell=1}\right| &\les& \frac{\ep_0}{r^5u^{1+\dec}}.
\eeaa

Next, in view of the definition of $\Jk$ in terms of $f_0$, and since $|q|=r(1+O(r^{-2}))$, we have
\beaa
&&\left( B'   -\frac{3a}{2}\ov{\Pc'}\Jk - \frac{a}{4} \ov{\Jk}\c A'\right)_1 \\
&=& \left( \b'   -\frac{3a}{2r}\left(  \rhoc\,' f_0+\widecheck{\rhod}'\dual  f _0 \right)- \frac{a}{2r} \a'\c f_0\right)_1\\
&&+i\left( \b'   -\frac{3a}{2r}\left(  \rhoc\,' f_0+\widecheck{\rhod}'\dual  f _0 \right)- \frac{a}{2r} \a'\c f_0\right)_2 +O(r^{-2})(\Pc', A'),
\eeaa
and hence
\beaa
\left( B'   -\frac{3a}{2}\ov{\Pc}\Jk - \frac{a}{4} \ov{\Jk}\c A\right) &=& w+i\dual w+O(r^{-2})(\Pc', A'),
\eeaa
with 
\beaa
w &:=& \b'   -\frac{3a}{2r}\left(  \rhoc\,' f_0+\widecheck{\rhod}'\dual  f _0 \right)- \frac{a}{2r} \a'\c f_0.
\eeaa
Then, we have
\beaa
&&\ov{\DD}\c\left( B'   -\frac{3a}{2}\ov{\Pc'}\Jk - \frac{a}{4} \ov{\Jk}\c A'\right)\\
&=& \ov{\DD}\c(w+i\dual w)+O(r^{-3})\dkb^{\leq 1}(\Pc', A') = (\nab-i\dual\nab)\c(w+i\dual w)+O(r^{-3})\dkb^{\leq 1}(\Pc', A')\\
&=& 2\div(w)+2i\curl(w)+O(r^{-3})\dkb^{\leq 1}(\Pc', A')
\eeaa
and hence
\beaa
&&\sup_{\Si_*}r^5u^{1+\dec}\left|\left[\ov{\DD}\c\left( B'   -\frac{3a}{2r}\ov{\Pc'}\Jk- \frac{a}{4r} \ov{\Jk}\c A'\right)\right]_{\ell=1}\right|\\
&\les& \sup_{\Si_*}r^5u^{1+\dec}|(\ddd_1w)_{\ell=1}|+\sup_{\Si_*}r^2u^{1+\dec}\dkb^{\leq 1}(\Pc', A').
\eeaa
In view of the definition of $w$, we infer from the above, using also the improved estimates for $\Pc'$ of Lemma \ref{lemma:controlofcurvaturePGframeMextonSigmastar} and the control of Theorem M2 for $A'$,
\beaa
\sup_{\Si_*}r^5u^{1+\dec}\left|\left[\ov{\DD}\c\left( B'   -\frac{3a}{2}\ov{\Pc'}\Jk - \frac{a}{4} \ov{\Jk}\c A'\right)\right]_{\ell=1}\right| &\les& \ep_0.
\eeaa

Next, using the decomposition 
\beaa
\DD' -\frac{a}{2}\Jk\nab_{e_4'} -\frac{a}{2}\Jk\nab_{e_3'} &=& (1+O(r^{-2}))\DD+O(r^{-2})\nab_{e_4'}+O(r^{-2})\nab_{e_3'}+O(r^{-3}),
\eeaa 
we have
\beaa
&&\left|\left(\ov{\DD'}\c -\frac{a}{2}\ov{\Jk}\c\nab_{e_4'} -\frac{a}{2}\ov{\Jk}\c\nab_{e_3'}\right)\c\left( B'   -\frac{3a}{2}\ov{\Pc'}\Jk - \frac{a}{4} \ov{\Jk}\c A'\right) - \ov{\DD}\c\left( B'   -\frac{3a}{2}\ov{\Pc'}\Jk - \frac{a}{4} \ov{\Jk}\c A'\right)\right| \\
&\les& r^{-4}|\dk^{\leq 1}\Ga_g'|+r^{-2}|\nab_\nu B'|.
\eeaa
Together with the bootstrap assumptions on $\Ga_g'$ and the estimate for $\nab_\nu B'$ in Lemma \ref{lemma:controlofcurvaturePGframeMextonSigmastar}, we infer on $\Si_*$
\beaa
\left|\left(\ov{\DD'}\c -\frac{a}{2}\ov{\Jk}\c\nab_{e_4'} -\frac{a}{2}\ov{\Jk}\c\nab_{e_3'}\right)\c\left( B'   -\frac{3a}{2}\ov{\Pc'}\Jk - \frac{a}{4} \ov{\Jk}\c A'\right) - \ov{\DD}\c\left( B'   -\frac{3a}{2}\ov{\Pc'}\Jk - \frac{a}{4} \ov{\Jk}\c A'\right)\right| \les \frac{1}{r^6},
\eeaa
and hence, together with the dominance in $r$ condition on $\Si_*$, we infer
\beaa
&&\left|\left(\ov{\DD'}\c -\frac{a}{2}\ov{\Jk}\c\nab_{e_4'} -\frac{a}{2}\ov{\Jk}\c\nab_{e_3'}\right)\c\left( B'   -\frac{3a}{2}\ov{\Pc'}\Jk - \frac{a}{4} \ov{\Jk}\c A'\right) - \ov{\DD}\c\left( B'   -\frac{3a}{2}\ov{\Pc'}\Jk - \frac{a}{4} \ov{\Jk}\c A'\right)\right|\\
& \les& \frac{\ep_0}{r^5u^{1+\dec}}.
\eeaa
Together with the above, we infer 
\beaa
\sup_{\Si_*}r^5u^{1+\dec}\left|\left[\left(\ov{\DD'}\c -\frac{a}{2}\ov{\Jk}\c\nab_{e_4'} -\frac{a}{2}\ov{\Jk}\c\nab_{e_3'}\right)\c\left( B'   -\frac{3a}{2}\ov{\Pc'}\Jk - \frac{a}{4} \ov{\Jk}\c A'\right)\right]_{\ell=1}\right| &\les& \ep_0
\eeaa
as desired.
\end{proof}

\begin{lemma}
We have
\beaa
\sup_{\Si_*}r^5u^{1+\dec}\left|\left[\ov{\DD}\c\Lieb_{\T'} B' \right]_{\ell=1}\right| &\les& \ep_0.
\eeaa
\end{lemma}

\begin{proof}
We have
\beaa
(\Lieb_{\T'} B')_{ab} &=& (\nab_{\T'}B')_{ab}+\g(\D_{e_a'}\T, e_c')B'_{cb}+\g(\D_{e_b'}\T, e_c')B'_{ac}. 
\eeaa
Now, recall that $k_{ab}=\g(\D_{e_a'}\T, e_b')$ verifies 
\beaa
k_{ab}=O(r^{-3})+\Ga_b.
\eeaa 
We infer
\beaa
\Lieb_{\T'} B' &=& \nab_{\T'}B'+O(r^{-3})B'+\Ga_b B'.
\eeaa
This yields
\beaa
2\Lieb_{\T'} B' &=& \nab_{e_3'}B'+\nab_{e_4'}B' +O(r^{-2})\dk^{\leq 1}B'+\Ga_b B'.
\eeaa
Using the Bianchi identities, we infer
\beaa
2\Lieb_{\T'} B' &=& \DD'\ov{\Pc'}+\frac{2}{r}B'+\frac{1}{2}\ov{\DD}'\c A'-\frac{4}{r}B'\\
&& +O(r^{-2})\dk^{\leq 1}B'+O(r^{-2})A'+r^{-3}\Ga_b'+r^{-1}\Ga_g'\c\Ga_b'\\
&=& \DD'\ov{\Pc'}+\frac{2}{r}B'+\frac{1}{2}\ov{\DD}'\c A'-\frac{4}{r}B'+r^{-3}\dk^{\leq 1}(\Ga_b')+r^{-1}\Ga_g'\c\Ga_b'.
\eeaa
and hence, using also the decomposition 
\beaa
\DD' &=& (1+O(r^{-2}))\DD+O(r^{-1})\nab_{e_4'}+O(r^{-1})\nab_{e_3'}+O(r^{-3}),
\eeaa 
we obtain
\beaa
2\Lieb_{\T'} B' &=& \DD\ov{\Pc'}+\frac{2}{r}B'+\frac{1}{2}\ov{\DD}\c A'-\frac{4}{r}B' +r^{-3}\dk^{\leq 1}(\Ga_b')+r^{-1}\Ga_g'\c\Ga_b'\\
 &=& \DD\ov{\Pc'}-\frac{2}{r}B'+\frac{1}{2}\ov{\DD}\c A' +r^{-3}\dk^{\leq 1}(\Ga_b')+r^{-1}\Ga_g'\c\Ga_b'.
\eeaa
This yields
\beaa
2\Lieb_{\T'} B' &=& \DD\ov{\Pc'}-\frac{2}{r}\left( B'   -\frac{3a}{2r}\ov{\Pc'}\Jk - \frac{a}{4r} \ov{\Jk}\c A'\right)+\frac{1}{2}\ov{\DD}\c A'\\
&&+r^{-3}\dk^{\leq 1}(\Ga_b')+r^{-1}\Ga_g'\c\Ga_b'.
\eeaa
We infer
\beaa
\ov{\DD}\c\Lieb_{\T'} B' &=& \frac{1}{2}\ov{\DD}\c\DD\ov{\Pc'}-\frac{1}{r}\ov{\DD}\c\left( B'   -\frac{3a}{2r}\ov{\Pc'}\Jk - \frac{a}{4r} \ov{\Jk}\c A'\right)+\frac{1}{4}\ov{\DD}\c\ov{\DD}\c A'\\
&& +r^{-4}\dk^{\leq 2}(\Ga_b')+r^{-2}\dk^{\leq 1}(\Ga_g'\c\Ga_b')\\
&=& \Delta\ov{\Pc'}-\frac{1}{r}\ov{\DD}\c\left( B'   -\frac{3a}{2r}\ov{\Pc'}\Jk - \frac{a}{4r} \ov{\Jk}\c A'\right)+\frac{1}{4}\ov{\DD}\c\ov{\DD}\c A'\\
&&  +r^{-4}\dk^{\leq 2}(\Ga_b')+r^{-2}\dk^{\leq 1}(\Ga_g'\c\Ga_b').
\eeaa
Recalling 
\beaa
 \rhoc\,' &=& \rhoc + r^{-2}\Ga_b+O\left(\frac{1}{r^4}\right),\\
  \widecheck{\rhod}' &=&  \rhod + r^{-2}\Ga_b+O\left(\frac{1}{r^4}\right),
  \eeaa  
we obtain 
\beaa
\ov{\DD}\c\Lieb_{\T'} B' &=& \Delta(\rhoc-i\rhod)-\frac{1}{r}\ov{\DD}\c\left( B'   -\frac{3a}{2r}\ov{\Pc'}\Jk - \frac{a}{4r} \ov{\Jk}\c A'\right)+\frac{1}{4}\ov{\DD}\c\ov{\DD}\c A'\\
&&+r^{-4}\dk^{\leq 2}(\Ga_b')+r^{-2}\dk^{\leq 1}(\Ga_g'\c\Ga_b').
\eeaa
We deduce 
\beaa
\left[\ov{\DD}\c\Lieb_{\T'} B'\right]_{\ell=1} &=& \left[\Delta\rhoc\right]_{\ell=1} -i \left[\Delta\rhod\right]_{\ell=1}+\frac{1}{4}\left[\ov{\DD}\c\ov{\DD}\c A'\right]_{\ell=1}\\
&& -\frac{1}{r}\left[\ov{\DD}\c\left( B'   -\frac{3a}{2r}\ov{\Pc'}\Jk - \frac{a}{4r} \ov{\Jk}\c A'\right)\right]_{\ell=1}\\
&& +r^{-4}\dk^{\leq 2}(\Ga_b')+r^{-2}\dk^{\leq 1}(\Ga_g'\c\Ga_b').
\eeaa
Since, for a scalar $f$,  
\beaa
\dual(\ov{\DD}\c\ov{\DD}\c)f &=& 2\dual(\div\div+i\curl\div)f\\
&=& 2\dds_2\dds_1(f,0)+2i\dds_2\dds_1(0,f),
\eeaa
we infer, by integration by parts of $\ov{\DD}\c\ov{\DD}\c$, 
\beaa
\left[\ov{\DD}\c\ov{\DD}\c A'\right]_{\ell=1} &\les& |A'||\dds_2\dds_1\Jp|
\eeaa
and hence
\beaa
\left|\left[\ov{\DD}\c\Lieb_{\T'} B'\right]_{\ell=1}\right| &\les&  \frac{1}{r}\left|\left[\ov{\DD}\c\left( B'   -\frac{3a}{2r}\ov{\Pc'}\Jk - \frac{a}{4r} \ov{\Jk}\c A'\right)\right]_{\ell=1}\right|+\frac{2}{r^2}\left|\left[\rhoc\right]_{\ell=1}\right| +\frac{2}{r^2}\left|\left[\rhod\right]_{\ell=1}\right|\\
&&+\left|\left[\left(\Delta+\frac{2}{r^2}\right)\rhoc\right]_{\ell=1}\right| +\left|\left[\left(\Delta+\frac{2}{r^2}\right)\rhod\right]_{\ell=1}\right| +|A'||\dds_2\dds_1\Jp|\\
&&+r^{-4}\dk^{\leq 2}(\Ga_b')+r^{-2}\dk^{\leq 1}(\Ga_g'\c\Ga_b').
\eeaa
The conclusion then follows from Proposition \ref{prop:controlofell=1modesonSigmastar} for the control of the $\ell=1$ mode of $\rhoc$ and $\rhod$, the control of Theorem M2 for $A'$, the control of $\dds_2\dds_1\Jp$ and $(\Delta+\frac{2}{r^2})\Jp$ in Lemma \ref{lemma:statementeq:DeJp.Sigmastar:improvedd}, the above improved control of $\Ga_b'$, and the one of Corollary \ref{cor:controloftheell=1modeofrenormalizedBonSigmarstar} for 
\beaa
\left[\ov{\DD}\c\left( B'   -\frac{3a}{2r}\ov{\Pc'}\Jk - \frac{a}{4r} \ov{\Jk}\c A'\right)\right]_{\ell=1}.
\eeaa
This concludes the proof of the corollary.
\end{proof}

\begin{lemma}
We have, for $k\leq k_*-14$, 
\beaa
\sup_{\Si_*}r^3u^{1+\dec}\left|\dk_*^{k}\left(\ov{\DD'}\c\Zc'+2\ov{\Pc'}\right)\right| &\les& \ep_0.
\eeaa 
\end{lemma}

\begin{proof}
Using  the decomposition 
\beaa
\DD' &=& (1+O(r^{-2}))\DD+O(r^{-1})\nab_{e_4'}+O(r^{-1})\nab_{e_3'}+O(r^{-3}),
\eeaa 
we have, using also $e_3'=\nu+r^{-1}\dk$, 
\beaa
\ov{\DD'}\c\Zc' &=& \ov{\DD}\c\Zc' +O(r^{-1})\nab_\nu\Zc' + r^{-2}\dk^{\leq }\Ga_g'
\eeaa
and hence 
\beaa
\ov{\DD'}\c\Zc'+2\ov{\Pc'}\ &=& \ov{\DD}\c\Zc' +2\ov{\Pc'}+O(r^{-1})\nab_\nu\Zc' + r^{-2}\dk^{\leq 1}\Ga_g'.
\eeaa
Recalling from the above the transformation formulas
\beaa
\Zc'=Z+O(r^{-1})(\ombc', \chibh')+r^{-1}\Ga_g'+O(r^{-3})+r^{-1}\dk_*^{\leq 1}\widetilde{\Ga_b},\\
\Pc'=\rhoc+i\rhod+O(r^{-4})+r^{-2}\Ga_g'+O(r^{-1})\bb'+O(r^{-2})\aa'.
\eeaa
We infer
\beaa
\ov{\DD'}\c\Zc'+2\ov{\Pc'}\ &=& \ov{\DD}\c Z +2(\rhoc-i\rhod) + O(r^{-4})+O(r^{-2})\dkb^{\leq 1}(\ombc', \chibh')+O(r^{-1})\bb'\\
&&+O(r^{-2})\aa' +r^{-2}\dk_*^{\leq 2}\widetilde{\Ga_b}+O(r^{-1})\nab_\nu\Zc' + r^{-2}\dk^{\leq 1}\Ga_g'.
\eeaa
Together with the control of the Ricci coefficients of the PG frame on $\Si_*$ provided by Lemma \ref{lemma:controlofRiccicoeffPGframeMextonSigmastar}, the control of the curvature components of the PG frame on $\Si_*$ provided by Lemma \ref{lemma:controlofcurvaturePGframeMextonSigmastar}, the bootstrap assumptions on $\Ga_g'$, and the control of $\widetilde{\Ga_b}$ provided by  Corollary  \ref{cor:estimatesforf0fplusfminusandJp:Gammatildeb}, we obtain on $\Si_*$, for $k\leq k_*-14$, 
\beaa
\left|\dk_*^k\left(\ov{\DD'}\c\Zc'+2\ov{\Pc'}\right)\right| &\les& \left|\dk_*^k\left(\ov{\DD}\c Z+2(\rhoc-i\rhod)\right)\right| + \frac{\ep_0}{r^3u^{1+\dec}}+\frac{1}{r^4},
\eeaa
and hence, in view of the dominance in $r$ condition for $r$ on $\Si_*$, this yields, for $k\leq k_*-14$, 
\beaa
\left|\dk_*^k\left(\ov{\DD'}\c\Zc'+2\ov{\Pc'}\right)\right| &\les& \left|\dk_*^k\left(\ov{\DD}\c Z+2(\rhoc-i\rhod)\right)\right| + \frac{\ep_0}{r^3u^{1+\dec}}.
\eeaa

Next, we compute
\beaa
\ov{\DD}\c Z+2(\rhoc-i\rhod) &=& (\nab-i\dual\nab)\c(\ze+i\dual\ze)+2(\rhoc-i\rhod)\\
&=& 2\big(\div\ze+\rhoc\big)+2i(\curl\ze-\rhod).
\eeaa
Together with the definition of $\mu$ and the null structure equation for $\curl\ze$, we infer
\beaa
\ov{\DD}\c Z+2(\rhoc-i\rhod) &=& 2\left(-\muc+\frac{1}{2}\chih\c\chibh\right) -i\chih\wedge\chibh =-2\muc+\Ga_b\c\Ga_g.
\eeaa 
We infer on $\Si_*$, for $k\leq k_*-14$, 
\beaa
\left|\dk_*^k\left(\ov{\DD'}\c\Zc'+2\ov{\Pc'}\right)\right| &\les& \left|\dk_*^k\muc\right| +\left|\dk_*^k(\Ga_b\c\Ga_g)\right|+ \frac{\ep_0}{r^3u^{1+\dec}}.
\eeaa
Together with the control for $\muc$, $\Ga_b$ and $\Ga_g$ of Proposition \ref{prop:decayonSigamstarofallquantities}, we deduce on $\Si_*$, for $k\leq k_*-14$, 
\beaa
\left|\dk_*^k\left(\ov{\DD'}\c\Zc'+2\ov{\Pc'}\right)\right| &\les&  \frac{\ep_0}{r^3u^{1+\dec}}.
\eeaa
This concludes the proof of the lemma.
\end{proof}

\begin{lemma}
We have 
\beaa
\sup_{\Si_*}ru^{\frac{1}{2}+\dec}\Big(\left|(r^2\Delta'+2)J^{(p)}\right|+\left|r^2\dds_2'\dds_1'(J^{(p)},0)\right|+\left|r^2\dds_2'\dds_1'(0, J^{(p)})\right|\Big) &\les& \ep_0.
\eeaa\end{lemma}

\begin{proof}
Since $\nu(\Jp)=0$, and since $\Jp$ is extended to $\Mext$ by $e_4'(\Jp)=0$, we easily infer from \eqref{eq:changofframedefiningthefrmaeofMext:chap5} the following estimate on $\Si_*$
\beaa
&&\left|(r^2\Delta'+2)J^{(p)}\right|+\left|r^2\dds_2'\dds_1'(J^{(p)},0)\right|+\left|r^2\dds_2'\dds_1'(0, J^{(p)})\right|\\
&\les& \left|(r^2\Delta+2)J^{(p)}\right|+\left|r^2\dds_2\dds_1(J^{(p)},0)\right|+\left|r^2\dds_2\dds_1(0, J^{(p)})\right| +r^{-2}|\dk^{\leq 2}\Jp|.
\eeaa
In view of the control of $\dds_2\dds_1\Jp$ and $(\Delta+\frac{2}{r^2})\Jp$ in Lemma \ref{lemma:statementeq:DeJp.Sigmastar:improvedd}, we infer on $\Si_*$
\beaa
\left|(r^2\Delta'+2)J^{(p)}\right|+\left|r^2\dds_2'\dds_1'(J^{(p)},0)\right|+\left|r^2\dds_2'\dds_1'(0, J^{(p)})\right| &\les& \frac{\ep_0}{ru^{\frac{1}{2}+\dec}}+\frac{1}{r^2}.
\eeaa
Together with the dominance condition on $r$ on $\Si_*$, we infer
\beaa
\left|(r^2\Delta'+2)J^{(p)}\right|+\left|r^2\dds_2'\dds_1'(J^{(p)},0)\right|+\left|r^2\dds_2'\dds_1'(0, J^{(p)})\right| &\les& \frac{\ep_0}{ru^{\frac{1}{2}+\dec}}
\eeaa
as stated.
\end{proof}

%%%%%%%%%%%%%%%%%%%%%%%%%%%%%%%%%%%%%%%%%%%%%%%%%%

\subsection{Proof of Proposition \ref{prop:improvedesitmatesfortemporalframeofMextonSigmastar}}
\lab{sec:proofofprop:improvedesitmatesfortemporalframeofMextonSigmastar}

%%%%%%%%%%%%%%%%%%%%%%%%%%%%%%%%%%%%%%%%%%%%%%%%%%

The results of sections \ref{sec:firstdecayesitmatesPGframeonSigmastar:chap5},  \ref{sec:decayesitmatesPGframeonSigmastar:chap5} and   \ref{sec:additionaldecayesitmatesPGframeonSigmastar:chap5} imply the proof of Proposition \ref{prop:improvedesitmatesfortemporalframeofMextonSigmastar} with $e_3'$ replaced by $\nu$ and $\dk$ replaced by $\dk_*$. The extension from $\dk_*$ to $\dk$ follows immediately from the null structure equations and Bianchi identities expressing derivatives in the $e_4'$ direction in terms of angular derivatives. Finally, since 
\beaa
\nu &=& e_3+b_*e_4=e_3'+r^{-1}\dk,
\eeaa 
we may replace $\nu$ by $e_3'$ in the corresponding estimates. This concludes the proof of  Proposition \ref{prop:improvedesitmatesfortemporalframeofMextonSigmastar}.

%%%%%%%%%%%%%%%%%%%%%%%%%%%%%%%%%%%%%%%%

\chapter{Decay estimates on the region $\Mext$ (Theorem M4)}
\lab{Chapter:DecayMext}

%%%%%%%%%%%%%%%%%%%%%%%%%%%%%%%%%%%%%%%%

The goal of this chapter is to  prove Theorem M4 by extending the decay estimates  on $\Si_*$ derived in section \ref{sec:decayestimatesPGframeonSigmastar}  to  the full spacetime region  $\Mext$.    The main result is stated in Proposition \ref{prop:improvedesitmatesfortemporalframeofMexton-new}. The estimates on $\Si_*$,  derived 
in section \ref{sec:decayestimatesPGframeonSigmastar},   are summarized in Proposition \ref{prop:improvedesitmatesfortemporalframeofMextonSigmastar-new}. 

 The results proved in this chapter rely  only on the main bootstrap assumptions 
on $\Mext$, the estimates for the extreme curvature  component  $A$  derived  in Theorem M1, and the estimates on the last slice proved in section \ref{sec:decayestimatesPGframeonSigmastar}.

In order to count the number of derivatives under control in this chapter, we introduce for convenience the following notation 
\bea\lab{eq:valueofkstarinchapter6forproofThmM4}
k_* &:=& k_{small}+60. 
\eea

%%%%%%%%%%%%%

\section{Preliminaries}

%%%%%%%%%%%%%

%%%%%%%%%%%%%%%%%%%%%%%%%%%

\subsection{The PG structure of $\Mext$}

%%%%%%%%%%%%%%%%%%%%%%%%%%%

Throughout this  chapter  we work with the spacetime region $\Mext$,   terminating in the  GCM last slice  $\Si_*$, as discussed in section \ref{section:GCMadmissible-spacetimes}.     For the convenience of the reader  we recall below  some the main facts concerning $\Mext$.   
\begin{enumerate}
 \item The PG structure of $\Mext$, given by $\big\{ r,  (e_3, e_4), \HH)\big\}$ together with adapted PG coordinates $(u, \th, \vphi)$,   is such that 
 \bea
 \xi=0, \qquad \om=0, \qquad \etab+\ze=0,
 \eea
 and
 \bea
 e_4(r)=1, \qquad e_4(u)= e_4(\th)=e_4(\vphi)=0, \qquad  \nab(r)=0.
 \eea 
 
  \item In $\Mext$, we have $1\leq u\leq u_*$ and $r\ge r_0$, with $r_0$ sufficiently large.
  
  \item The timelike hypersurface $\TT=\{r=r_0\}$ is a boundary of $\Mext$.

\item The constants $(a,m)$  are the ones  associated to $S_*$ according to Definition \ref{define:am-onSi}, see section \ref{sec:definitionofamthetandvphiadmissible}.
  
  \item   $\Mext$ comes equipped,  see section  \ref{section:can.ell=1basis.Mext},  with a basis of $\ell=1$ modes  $J^{(p)}$ with  $p=0,+,-$ verifying
\bea
e_4(\Jp)=0.
\eea
  
\item    $\Mext$ comes equipped with    complex, anti selfadjoint\footnote{i.e. $ \dual \Jk=-i\Jk$, $\dual \Jk_\pm=-i \Jk_{\pm}$.}   1-form $\Jk$,  $\Jk_\pm$   satisfying, see  section \ref{section:auxilliaryformsMext}   and section \ref{section:main-normsextregion}  
 \bea
 \lab{eq:propertiesJk}
 \bsplit
 \nab_4\Jk &= -\frac{1}{q}\Jk, \qquad \qquad\qquad \quad \quad  \nab_4\Jk_\pm  = -\frac{1}{q}\Jk_\pm,\\
 \Re(\Jk_+)\c\Re(\Jk)&=-\frac{1}{|q|^2}J^{(-)}, \qquad \quad \Re(\Jk_-)\c\Re(\Jk)=\frac{1}{|q|^2}J^{(+)},
 \end{split}
 \eea
 as well as 
 \bea
 \lab{eq:propertiesJk-more}
 \bsplit
   \Jk\c\ov{\Jk} &= \frac{2(\sin\th)^2}{|q|^2},\\
 \Jk_+\c\ov{\Jk_+}&=\frac{2(\cos\th)^2(\cos\vphi)^2+2(\sin\vphi)^2}{|q|^2},\\
  \Jk_-\c\ov{\Jk_-}&=\frac{2(\cos\th)^2(\sin\vphi)^2+2(\cos\vphi)^2}{|q|^2}.
  \end{split}
 \eea
 \end{enumerate}

%%%%%%%%%%%%%%%%%%%%%%%%%%%%%%%%%%%%%%%
 
\subsection{Linearized  quantities and definition of $\Ga_g$ and $\Ga_b$}
\lab{sec:linearizedquantitiesanddefGabGag:chap6}

%%%%%%%%%%%%%%%%%%%%%%%%%%%%%%%%%%%%%%%

 We recall below the linearized  quantities  obtained   by subtracting their   $\mbox{Kerr}(a, m)$  values, see Definition \ref{def:renormalizationofallnonsmallquantitiesinPGstructurebyKerrvalue}:
\begin{enumerate}
\item Linearization of Ricci and curvature coefficients.
\bea
\bsplit
\trXc &:= \tr X-\frac{2}{q}, \qquad\,\qquad \,\,\,\,    \trXbc := \tr\Xb+\frac{2q\Delta}{|q|^4},\\
\\
\Zc &:= Z-\frac{a\ov{q}}{|q|^2}\Jk,\qquad \qquad \quad 
\Hc := H-\frac{aq}{|q|^2}\Jk,\\
\\
  \ombc& := \omb  - \frac 1 2 \pr_r\left(\frac{\De}{|q|^2} \right),\qquad \,\, \Pc := P+\frac{2m}{q^3}.\\
\end{split}
\eea

\item Linearization of derivatives of $r, q, u$.
\bea
\bsplit
\widecheck{\DD q} &:=\DD q+a\Jk, \qquad\quad \widecheck{\DD \ov{q}} :=\DD \ov{q}-a\Jk,\\
\widecheck{e_3(r)} &:= e_3(r)+\frac{\Delta}{|q|^2},\\
\widecheck{\DD u} &:= \DD u -a\Jk,\qquad \widecheck{e_3(u)} := e_3(u) -\frac{2(r^2+a^2)}{|q|^2}.
\end{split}
\eea

\item Linearization for $\Jk$ and $\Jk_\pm$.
\bea
\bsplit
\widecheck{\ov{\DD}\c\Jk}& := \ov{\DD}\c\Jk-\frac{4i(r^2+a^2)\cos\th}{|q|^4}, \qquad \widecheck{\nab_3\Jk}:=\nab_3\Jk -\frac{\De q}{|q|^4}\Jk,\\
 \widecheck{\ov{\DD}\c\Jk_\pm}&:=\ov{\DD}\c \Jk_\pm+\frac{4}{r^2} J^{(\pm)} \pm \frac{4ia^2\cos\th}{|q|^4}J^{(\mp)}, \\  
 \widecheck{\nab_3 \Jk_\pm } &:= \nab_3 \Jk_\pm - \frac{\De q}{|q|^4} \, \Jk_{\pm} \pm   \,\frac{2a}{|q|^2}  \Jk_{\mp}. 
  \end{split}
\eea

\item Linearization for $\Jp$.
\bea
\lab{eq:linearizationsforJp}
\bsplit
\widecheck{\DD J^{(0)}}  &:= \DD J^{(0)} -i\Jk, \qquad \qquad  \qquad \widecheck{\DD(J^{(\pm)})}:= \DD(J^{(\pm)})-\Jk_{\pm},
\\
\widecheck{e_3(J^{(+)})}&:=e_3(J^{(+)}) + \frac{2a}{|q|^2}J^{(-)}, \qquad
 \widecheck{e_3(J^{(-)})}:=e_3(J^{(-)}) - \frac{2a}{|q|^2}J^{(+)}.
 \end{split}
\eea
 \end{enumerate}

We  also recall the sets $\Ga_g, \Ga_b$,  see Definition \ref{definition.Ga_gGa_b}:
\begin{enumerate}
\item  The set   $\Ga_g$   with
 \bea
 \bsplit
 \Ga_g &= \Big\{\trXc,\quad  \Xh,\quad \Zc,\quad \trXbc , \quad r\Pc, \quad  rB, \quad  rA\Big\}.
 \end{split}
 \eea
 
 \item  The set  $\Ga_b=\Ga_{b,1}\cup \Ga_{b, 2}\cup \Ga_{b,3} \cup \Ga_{b, 4}$   with
 \bea
 \bsplit
 \Ga_{b,1}&= \Big\{\Hc, \quad \Xbh, \quad \ombc, \quad \Xib,\quad  r\Bb, \quad \Ab\Big\},\\
  \Ga_{b, 2}&= \Big\{r^{-1}\widecheck{e_3(r)}, \quad \widecheck{\DD q}, \quad \widecheck{\DD \ov{q}}, \quad \widecheck{\DD u},  \quad r^{-1}\widecheck{e_3(u)}\Big\},\\
  \Ga_{b, 3}&= 
   \Big\{ \widecheck{\DD(J^{(0)})}, \quad \widecheck{\DD(J^{(\pm)})}, \quad e_3(J^{(0)}), \quad    \widecheck{e_3(J^{(\pm)})}\Big\}, \\
   \Ga_{b,4}&=\Bigg\{ r\,\widecheck{\ov{\DD}\c\Jk}, \quad r\,\DD\hot\Jk, \quad r\,\widecheck{\nab_3\Jk}, \quad
   r\,\widecheck{\ov{\DD}\c\Jk_\pm}, \quad r\,\DD\hot\Jk_\pm, \quad r\,\widecheck{\nab_3\Jk_\pm} \Bigg\}. 
   \end{split}
 \eea
\end{enumerate}

%%%%%%%%%%%%%%%%%%%%%%%%%%%%%

\subsection{Main assumptions}

%%%%%%%%%%%%%%%%%%%%%%%%%%%%%%%

\begin{definition}
\label{definition:norms-ProofThm.M4}
  We make  use of  the following   norms on  $S=S(u,r)\subset\Mext$,
  \bea
  \bsplit
  \| f\|_{\infty} (u,r):&=\| f\|_{L^\infty\big(S(u,r)\big)}, \qquad \quad  \| f\|_{2} (u,r):=\| f\|_{L^2\big(S(u,r)\big)}, \\
  \|f\|_{\infty,k}(u, r)&:= \sum_{i=0}^k \|\dk^i f\|_{\infty }(u, r),  \qquad 
\|f\|_{2,k}(u, r):=\sum_{i=0}^k \|\dk^i f\|_{2}(u, r).
\end{split}
  \eea
  We  shall also make use of 
  \beaa
   \|f\|_{\infty, k}(\Mext)=\sup_{\Mext} |\dk^{\le k} f|
  \eeaa 
  \end{definition}
  
  \begin{remark}
  We note that the derivatives $\dkb=(r\nab) $ and $\dk=( r\nab, r\nab_4, \nab_3)$  are defined with respect to the outgoing  PG frame of $\Mext$, which is not adapted to  the spheres $S(u, r)$.
  \end{remark}

 \begin{definition}[Order of magnitude notation]
 \lab{def:ordermagnitude}
 Throughout this chapter,  we will be  using the  notation $O(r^{-p})$ to  denote:
\begin{enumerate} 
\item a scalar function depending only on $(r, \th)$ which is smooth and such that
\beaa
r^p|(r\pr_r, \pr_\th)^kO(r^{-p})|\les_k 1\quad\textrm{for }\quad k\geq 0\quad\textrm{and}\quad r\geq r_0,
\eeaa

\item a 1-form of the type $O(r^{-p+1})\Jk$ where $O(r^{-p+1})$ denotes a scalar function as above,

\item a symmetric traceless 2-tensor of the type $O(r^{-p+2})\Jk\hot\Jk$ where $O(r^{-p+2})$ denotes a scalar function as above.
\end{enumerate}
Often in the text we shall use the notation $r^{-p} U$, where $U$ is a small quantity,  instead of $ O(r^{-p} ) U$.
\end{definition}

We will also make use of  the following notation.
  \begin{definition}
  \lab{Notation:forGo_k}
 We introduce  the notation $ U\in r^{-p}\Go_k$ for horizontal tensors $U$  satisfying 
\bea
\left|\dk^{\leq k} U \right| &\les& \ep_0 r^{-k} u^{-1-\dec}.
\eea
\end{definition}

 For the benefit of the reader we state below  the main assumptions which will be used throughout  this  chapter:

{\bf Ref 1.} In view of our main bootstrap assumptions on decay and boundedness, the following estimates hold on $\Mext$:
\begin{enumerate}
\item  For $0\le k\le k_{small}$, we have 
\bea
\lab{Ref1-smallk}
\bsplit
\| \Ga_g\|_{\infty, k} &\les\ep \min\left\{ r^{-2} u^{-\frac{1}{2}-\dec}, \, r^{-1} u^{-1-\dec}   \right\},\\
\| \nab_3  \Ga_g\|_{\infty, k-1}  &\les \ep  r^{-2} u^{-1-\dec}, \\   
\|\Ga_b\|_{\infty, k} &\les\ep  r^{-1} u^{-1-\dec},
\end{split}
\eea
and
\bea
\lab{Ref1-smallk:B}
\begin{split}
\| B\|_{\infty, k} &\les\ep\min\left\{r^{-3} (u+2r)^{-\frac{1}{2}-\dec}, \,  r^{-2} (u+2r)^{-1-\dec}\right\},\\
\|\nab _3  B\|_{\infty, k-1} &\les\ep\min\left\{r^{-4} (u+2r)^{-\frac{1}{2}-\dec}, \,  r^{-3} (u+2r)^{-1-\dec}\right\}.
\end{split}
\eea

\item We also make the auxiliary bootstrap assumption\footnote{This auxiliary bootstrap assumption will be improved in Proposition \ref{lemma:firstcontroloftrXbcinMext}.} for $0\le k\le k_{small}$
\bea
\lab{eq:improvedRef1-trXc}
\| \trXc\|_{\infty, k} &\les& \ep r^{-2} u^{-1-\dec}.
\eea

\item   For $ k\le k_{large}$, we have
\bea
\label{Ref1-largek}
\|\Ga_g\|_{\infty, k}&\les&\ep  r^{-2}, \qquad  \| \Ga_b\|_{\infty,k} \les  \ep  r^{-1}.
\eea
\end{enumerate}

\begin{remark}
\lab{remark:interpolation}
Note that we can interpolate  between the estimates   \eqref{Ref1-smallk}  for $k\le k_{small} $ and  \eqref{Ref1-largek} for $k\leq k_{large}$ to derive on $\Mext$, for\footnote{Recall from \eqref{eq:valueofkstarinchapter6forproofThmM4} that $k_* = k_{small}+60$ throughout this chapter.} all $k\le k_*$, 
\bea
\lab{Ref1-k_*}
\bsplit
\| \Ga_g\|_{\infty, k} &\les\ep \min\left\{ r^{-2} u^{-\frac{1}{2}-\frac{\dec}{2}}, \, r^{-1} u^{-1-\frac{\dec}{2}}\right\},\\
 \|\nab_3  \Ga_g\|_{\infty, k-1} &\les \ep  r^{-2} u^{-1-\frac{\dec}{2}}, \\   
\|\Ga_b\|_{\infty, k} &\les\ep  r^{-1} u^{-1-\frac{\dec}{2}},
\end{split}
\eea
see Lemma \ref{lemma:interpolation} for the corresponding statement on $\Si_*$. 
\end{remark}

{\bf Ref  2.}   According to\footnote{The estimate for $\nab_3^2A$ follows immediately from the estimate in Theorem M1 for $A$, $\nab_3A$ and $\qf$, and the definition of $\qf$ which yields $\nab_3^2A=O(r^{-4})\qf+O(r^{-1})\nab_3A+O(r^{-2})A$.} Theorem M1,   we have on $\MMext$, for\footnote{Recall from \eqref{eq:valueofkstarinchapter6forproofThmM4} that $k_* = k_{small}+60$ throughout this chapter.} all  $0\le k\le k_{*}$, 
\bea
\begin{split}
\| A\|_{\infty, k} &\les\ep_0\min\left\{r^{-3} (u+2r)^{-\frac{1}{2}-\dee}, \,  \log(1+u)r^{-2} (u+2r)^{-1-\dee}\right\},\\
\|\nab_3 A\|_{\infty, k-1} &\les\ep_0\min\left\{r^{-4} (u+2r)^{-\frac{1}{2}-\dee},\,  \log(1+u)r^{-3} (u+2r)^{-1-\dee}\right\},\\
\|\nab_3^2A\|_{\infty, k-1} &\les\ep_0\min\left\{r^{-5}u^{-\frac{1}{2}-\dee},\,  \log(1+u)  r^{-4}u^{-1-\dee}\right\},
\end{split}
\eea
where we recall that $\dee>\dec$.

\begin{remark}
Let 
\bea
\de':=\frac 1 2\left(\dee -\dec\right),
\eea
where $\de'>0$ since $\dee>\dec$. In view of Definition \ref{Notation:forGo_k} and {\bf Ref  2}, we have
 \bea
A \in r^{-2-\de'}\Go_{k_*},\qquad \nab_3A \in r^{-3-\de'}\Go_{k_*}, \qquad \nab_3^2A \in r^{-4-\de'}\Go_{k_*}.
 \eea
\end{remark}

%%%%%%%%%%%%%%%%%%%%%%%%%%%

\subsection{Main  equations in $\Mext$}

%%%%%%%%%%%%%%%%%%%%%%%%%%%

We recall below a subset of the  null structure and  Bianchi equations holding for an outgoing PG structure,    see  Proposition  \ref{prop-nullstrandBianchi:complex:outgoing}.

 \begin{proposition}\lab{prop-nullstrandBianchi:complex:outgoing:again:chap6}
  In the outgoing PG structure of $\Mext$, we have
\lab{prop-nullstrandBianchi:complex:outgoing-again}
\beaa
\nab_4\tr X +\frac{1}{2}(\tr X)^2 &=& -\frac{1}{2}\Xh\c\ov{\Xh},\\
\nab_4\Xh+\Re(\tr X)\Xh &=& -A,
\\
\nab_4\tr\Xb +\frac{1}{2}\tr X\tr\Xb &=& -\DD\c\ov{Z}+Z\c\ov{Z}+2\ov{P}-\frac{1}{2}\Xh\c\ov{\Xbh},\\
\nab_4\widehat{\Xb} +\frac{1}{2}\tr X\, \widehat{\Xb}  &=& -\frac{1}{2}\DD\hot Z  +\frac{1}{2}Z\hot Z -\frac{1}{2}\ov{\tr\Xb} \widehat{X},
\\
\nab_4Z +\tr X Z &=&  -\widehat{X}\c\ov{Z} -B,\\
\nab_4H +\frac{1}{2}\ov{\tr X}(H+Z) &=&   -\frac{1}{2}\Xh\c(\ov{H}+\ov{Z}) -B,
\\
\nab_4\omb  -(2\eta+\ze)\c\ze &=&   \rho,
\eeaa
and
\beaa
\frac{1}{2}\ov{\DD}\c\Xh +\frac{1}{2}\Xh\c\ov{Z} &=& \frac{1}{2}\DD\ov{\tr X}+\frac{1}{2}\ov{\tr X}Z-i\Im(\tr X)H -B,\\
\frac{1}{2}\ov{\DD}\c\Xbh -\frac{1}{2}\Xbh\c\ov{Z} &=& \frac{1}{2}\DD\ov{\tr\Xb}-\frac{1}{2}\ov{\tr\Xb}Z-i\Im(\tr\Xb)(-Z+\Xib)+\Bb.
\eeaa

Also, we have
 \beaa
 \nab_3A -\frac{1}{2}\DD\hot B &=& -\frac{1}{2}\tr\Xb A+4\omb A +\frac{1}{2}(Z+4H)\hot B -3\ov{P}\Xh,\\
\nab_4B -\frac{1}{2}\ov{\DD}\c A &=& -2\ov{\tr X} B  +\frac{1}{2}A\c  \ov{Z},\\
\nab_3B-\DD\ov{P} &=& -\tr\Xb B+2\omb B+\ov{\Bb}\c \Xh+3\ov{P}H +\frac{1}{2}A\c\ov{\Xib},\\
\nab_4P -\frac{1}{2}\DD\c \ov{B} &=& -\frac{3}{2}\tr X P -\frac{1}{2}Z\c\ov{B}  -\frac{1}{4}\Xbh\c \ov{A}, \\
\nab_3P +\frac{1}{2}\ov{\DD}\c\Bb &=& -\frac{3}{2}\ov{\tr\Xb} P -\frac{1}{2}(\ov{2H-Z})\c\Bb +\Xib\c \ov{B} -\frac{1}{4}\ov{\Xh}\c\Ab, \\
\nab_4\Bb+\DD P &=& -\tr X\Bb+\ov{B}\c \Xbh+3P Z,\\
\nab_3\Bb +\frac{1}{2}\ov{\DD}\c\Ab &=& -2\ov{\tr\Xb}\,\Bb -2\omb\,\Bb -\frac{1}{2}\Ab\c (\ov{-2Z +H})-3P \,\Xib,\\
\nab_4\Ab +\frac{1}{2}\DD\hot\Bb &=& -\frac{1}{2}\ov{\tr X} \Ab +\frac{5}{2}Z\hot \Bb -3P\Xbh.
\eeaa
\end{proposition} 
    
 We also recall the following transport    equations in the $e_4$ direction  for derivatives of the outgoing PG coordinates $(r,u, \th)$,  see Proposition  \ref{prop:e_4(xyz)}.
\begin{proposition}
\lab{prop:e_4(xyz)-again}
The following equations hold true for the coordinates $(u,r,\th)$ associated to an outgoing  PG structure
\beaa
e_4(e_3(r)) &=&  -2\omb,
\eeaa
\beaa
\nab_4\DD u +\frac{1}{2}\tr X\DD u &=& -\frac{1}{2}\Xh\c\ov{\DD}u,\\
e_4(e_3(u)) &=& -\Re\Big((Z+H)\c\ov{\DD} u\Big),
\eeaa
\beaa
\nab_4\DD\cos\th +\frac{1}{2}\tr X\DD\cos\th &=& -\frac{1}{2}\Xh\c\ov{\DD}\cos\th,\\
e_4(e_3(\cos\th)) &=& -\Re\Big((Z+H)\c\ov{\DD} \cos\th\Big).
\eeaa
\end{proposition}

%%%%%%%%%%%%%%%%%%%%%%%%%%%%%%%%%%%%%%%%%%%%%

\subsection{Commutator formulas revisited}
\lab{subsection:commutationlemmas revisited}

%%%%%%%%%%%%%%%%%%%%%%%%%%%%%%%%%%%%%%%%%%%%%

We record below the main commutation formulas  which will be used  in this chapter.

%%%%%%%%%%%%%%%%

\subsubsection{Real case}

%%%%%%%%%%%%%%%%

The following commutation formulas    are an immediate adaptation of   those  in Lemma \ref{lemma:comm-gen} to the case of an outgoing PG structure\footnote{That is $ \xi=0 $, $\om=0$ and $\etab+\ze=0$.}.

\begin{lemma}\lab{lemma:commutationformulasintherealcase:chap6}
We have the following commutations formulas:
\begin{enumerate}
\item If   $f $ is a scalar, we have
\beaa
\,[\nab_4, \nab_b]  f&=& -\frac{1}{2}\Big(\trch\nab_bf+\atrch\dual\nab_bf\Big) - \chih_{bc} \nab_c f,\\
\, [\nab_4, \nab_3]  f&=&-2 (\ze+\eta)\c\nab f - 2\omb  \nab_4 f.
\eeaa

\item  If U is a horizontal tensor, we have
\beaa
\,[\nab_4, \nab_b] U&=& -\frac{1}{2}\Big(\trch\nab_b+\atrch\dual\nab_b\Big)U -\chih_{bc}\nab_c U \\
&&+O(1)\trch\ze U+O(r^{-2})\chih U +    O(1)\b U,\\
\, [\nab_4, \nab_3] U&=& -2(\ze+\eta )\c  \nab U -2\omb \nab_4 U+  O(1)\eta\ze U+    O(1)\rhod U. 
\eeaa
\end{enumerate}
\end{lemma}

%%%%%%%%%%%%%%%%%%%%%%%%%%%%%%%%%%%%

\subsubsection{Complex case}

%%%%%%%%%%%%%%%%%%%%%%%%%%%%%%%%%%%%%

The    following  commutation formulas  can be easily derived from the  ones above, see also section 6.1 in \cite{GKS1}.
\begin{lemma}
\lab{lemma:commutation-complexM6}
The following commutation  formulas hold true.
\begin{enumerate}
\item  For  a scalar complex function $F$, we have
 \bea
 \lab{eq:commutation-complexM6-scalar}
         \,[\nab_4, \DD] F &=&-\frac 1 2 \tr X \DD F +r^{-1}\Ga_g\c \dkb  F.
       \eea
      
\item For    an anti-self dual horizontal 1-form $U$,  we have 
\bea
\lab{eq:commutation-complexM6-1form}
\, [\nab_4, \DD\hot ]U&=& - \frac 1 2 \tr X( \DD\hot U -  Z \hot U)+r^{-1}\Ga_g\c \dkb^{\le 1}  U
\eea
and
 \bea
\lab{eq:commutation-complexM6-1form:div}
  \, [\nab_4, \DDov\c] U &=& -\frac 1 2\ov{\tr X}\big(\DDov\c U +\ov{Z}\c U)  + r^{-1}\Ga_g\c \dkb^{\leq 1} U.
 \eea

\item For  an anti-self dual symmetric traceless horizontal 2-form $U$,  we have 
 \bea
\lab{eq:commutation-complexM6-2form:div}
  \, [\nab_4, \DDov\c] U &=& -\frac 1 2\ov{\tr X}\big(\DDov\c U +2\ov{Z}\c U)  + r^{-1}\Ga_g\c\dkb^{\leq 1} U.
 \eea
 
 \item For an anti-self dual horizontal  $k$-tensor, we have
 \bea
 \lab{eq:commutation-complexM6-2form:general}
\, [\nab_4, \DD]U&=& - \frac 1 2 \tr X\DD U+O(1)Z U+r^{-1}\Ga_g\c \dkb^{\le 1}  U.
\eea
 \end{enumerate}      
\end{lemma}

\begin{remark}\lab{remark:commutation-complexM6}
We note  that   the terms denoted by $\Ga_g$ in     \eqref{eq:commutation-complexM6-scalar} only contain  $\Xh$.  Also,   the terms denoted by $\Ga_g$   in \eqref{eq:commutation-complexM6-1form}--\eqref{eq:commutation-complexM6-2form:general} contain only $\Xh$,  $\Zc$ and $B$. 
\end{remark}

\begin{corollary}
\lab{cor:commutation-complexM6}    
If $U$ denotes 
\begin{enumerate}
\item a scalar, we have
 \beaa
  \, [\nab_4, q\DD] U &=& \Ga_g \c \dkb U,
 \eeaa

\item an anti-self dual horizontal 1-form,  we have 
     \beaa
     \, [\nab_4, q\DD\hot ]U&=& O(r^{-2})U +\Ga_g\c \dkb^{\le 1}  U,
 \eeaa 
 
\item an anti-self dual horizontal 1-form, or  an anti-self dual symmetric traceless horizontal 2-form,  we have 
     \beaa
  \, [\nab_4, \ov{q}\,\DDov\c]U &=& O(r^{-2})U+ \Ga_g \c \dkb^{\leq 1} U,
 \eeaa    
 
\item an anti-self dual horizontal  $k$-tensor, we have
   \beaa
     \, [\nab_4, q\DD]U&=& O(r^{-2})U +\Ga_g\c \dkb^{\le 1}  U.
 \eeaa 
 \end{enumerate}   
\end{corollary}

\begin{proof}
This follows immediately from  Lemma \ref{lemma:commutation-complexM6} and the fact that $e_4(q)=1$, $\tr X=\frac{2}{q}+\Ga_g$ and $Z=O(r^{-2})+\Ga_g$. 
\end{proof}

%%%%%%%%%%%%%%%%%%%%%%%%%%%%%%%%%%%%%%%%%%%%%%%%%%%%

\subsection{Linearized null structure equations and Bianchi identities for outgoing PG structures}

%%%%%%%%%%%%%%%%%%%%%%%%%%%%%%%%%%%%%%%%%%%%%%%%%%%%

Recall the definition of the linearized quantities and of $\Ga_g$ and $\Ga_b$ in section \ref{sec:linearizedquantitiesanddefGabGag:chap6}. We use extensively the notation $O(r^{-p})$ made in Definition \ref{def:ordermagnitude} to denote lower order  linear terms.
The following lemma provides the linearized null structure equations and Bianchi identities in $\Mext$. 

\begin{lemma}
\lab{Lemma:linearized-nullstr}
The linearized null structure equations in the  $e_4$ direction are 
\beaa
\nab_4(\widecheck{\tr X}) +\frac{2}{q}\widecheck{\tr X} &=& \Ga_g\c\Ga_g,\\
\nab_4\Xh+\frac{2r}{|q|^2}\Xh &=& -A+\Ga_g\c\Ga_g,\\
\nab_4\Zc + \frac{2}{q}\Zc  &=&    - \frac{aq}{|q|^2}\ov{\Jk}\c\widehat{X} -B +O(r^{-2})\trXc+\Ga_g\c\Ga_g,\\
\nab_4\Hc+\frac{1}{\ov{q}}\Hc &=& -\frac{1}{\ov{q}}\Zc     -\frac{ar}{|q|^2}\ov{\Jk}\c\Xh -B +O(r^{-2})\trXc +\Ga_b\c\Ga_g,\\
\nab_4\widecheck{\tr\Xb} +\frac{1}{q}\trXbc &=& -\DD\c\ov{\Zc}+2\ov{\Pc }+O(r^{-2})\Zc+O(r^{-1})\trXc\\
&&+ O(r^{-1})\widecheck{\DD\c\ov{\Jk}} +O(r^{-3})\widecheck{\DD(\cos\th)}+\Ga_b\c\Ga_g,\\
\nab_4\Xbh +\frac{1}{q} \Xbh &=& -\frac{1}{2}\DD\hot\Zc +O(r^{-2})\Zc+O(r^{-1})\Xh+O(r^{-1})\DD\hot\Jk+O(r^{-3})\widecheck{\DD(\cos\th)}\\
&&+\Ga_b\c\Ga_g,\\
\nab_4(\ombc) &=& \Re(\Pc)+O(r^{-2})\Zc+O(r^{-2})\Hc+\Ga_b\c\Ga_g,\\
\nab_4\Xib  +\frac{1}{q}\Xib &=& O(r^{-1})\dkb^{\leq 1}(\ombc)+O(r^{-2})\Zc+O(r^{-2})\Hc+O(r^{-2})\trXbc\\
&&+O(r^{-3})\widecheck{\DD(\cos\th)}+\Ga_b\c\Big(\ombc,\Ga_g\Big).
\eeaa
The linearized Codazzi  for $\Xh$ takes the form
\beaa
\frac{1}{2}\ov{\DD}\c\Xh &=& \frac{1}{\ov{q}}\Zc -B +O(r^{-2})\Xh +O(r^{-2})\Hc +O(r^{-1})\dkb^{\leq 1}\trXc +O(r^{-2})\widecheck{\DD(\cos\th)}+\Ga_b\c\Ga_g.
\eeaa
The linearized Bianchi equations for $B, P, \Bb$  are 
 \beaa
\nab_4B +\frac{4}{\ov{q}} B &=&  \frac{1}{2}\ov{\DD}\c A +\frac{aq}{2|q|^2}\ov{\Jk}\c A+\Ga_g\c(B,A),\\
\nab_4\left(\Pc \right)-\frac{1}{2}\DD\c \ov{B} &=& -\frac{3}{q} \Pc -\frac{a\ov{q}}{2|q|^2}\Jk\c\ov{B}+O(r^{-3})\trXc+r^{-1}\Ga_g\c\Ga_g+\Ga_b\c A,\\
\nab_4\Bb+\DD\left(\Pc\right) &=& -\frac{2}{q}\Bb +O(r^{-2})\Pc+O(r^{-3})\Zc +O(r^{-4})\widecheck{\DD(\cos\th)}+r^{-1}\Ga_b\c\Ga_g.
\eeaa
Also
\beaa
\nab_3B -\DD\ov{\Pc} &=& \frac{2}{r}B+O(r^{-2}) B +O(r^{-2})\Pc+O(r^{-3})\Hc +O(r^{-4})\widecheck{\DD(\cos\th)}+r^{-1}\Ga_b\c\Ga_g.
\eeaa
\end{lemma}

\begin{proof}
The proof of the lemma relies on the null structure equations and Bianchi identities of Proposition \ref{prop-nullstrandBianchi:complex:outgoing:again:chap6}, the definition of the linearized quantities and of $\Ga_g$ and $\Ga_b$ in section \ref{sec:linearizedquantitiesanddefGabGag:chap6}, the notation $O(r^{-p})$ made in Definition \ref{def:ordermagnitude}, the fact that  $a$ and $m$ are constants, and the following identities 
\beaa
e_4(r)=1, \qquad e_4(\th)=0, \qquad \nab_4\Jk=-\frac{1}{q}\Jk, \qquad e_4(q)=1, \qquad e_4(\ov{q})=1, \qquad \nab(r)=0,
\eeaa
where we used in particular the fact that $q=r+ai\cos\th$ and $\ov{q}=r-ai\sin\th$. See section \ref{appendix: Proof of Lemma-ref{Lemma:linearized-nullstr}} for the proof.  
\end{proof}

%%%%%%%%%%%%%%%%%%%%%%%%%%%%%%%%%%%%%%%%%%%%%

\subsection{Other linearized equations  for outgoing PG structures}

%%%%%%%%%%%%%%%%%%%%%%%%%%%%%%%%%%%%%%%%%%%%%

\begin{lemma}
\lab{Lemma:otherlinearizedquant}
We have
\beaa
e_4\left(\widecheck{e_3(r)}\right) &=& -2\ombc,\\
\nab_4\widecheck{\DD u}+\frac{1}{q}\widecheck{\DD u} &=& O(r^{-1})\trXc+O(r^{-1})\Xh+\Ga_b\c\Ga_g,\\
e_4\left(\widecheck{e_3(u)}\right) &=& O(r^{-1})\Hc+O(r^{-1})\Zc+O(r^{-2})\widecheck{\DD u}+\Ga_b\c\Ga_b,
\\
\nab_4\widecheck{\DD\cos\th}+\frac{1}{q}\widecheck{\DD\cos\th} &=& \frac{i}{2}\ov{\Jk}\c\Xh+O(r^{-1})\trXc+\Ga_b\c\Ga_g,\\
e_4(e_3(\cos\th)) &=& O(r^{-1})\Hc+O(r^{-1})\Zc+O(r^{-2})\widecheck{\DD\cos\th}+\Ga_b\c\Ga_b.
\eeaa
\end{lemma}

\begin{proof}
One can proceed as in Lemma \ref{Lemma:linearized-nullstr}     starting with the corresponding equations in  Lemma \ref{prop:e_4(xyz)-again}. See section \ref{section:Proof-lemmaLemma:otherlinearizedquant}  for the details.
\end{proof}

\begin{lemma}
\lab{Lemma:linearizedJk}
The following equations hold  for   the tensors $\Jk, \Jk_\pm$:
\begin{enumerate}
\item We have
\beaa
\nab_4 (\DD\hot\Jk)+\frac{2}{q}\DD\hot\Jk &=& O(r^{-1})B +O(r^{-2})\trXc+O(r^{-2})\Xh\\
&&+O(r^{-2})\Zc +O(r^{-3})\widecheck{\DD(\cos\th)}+r^{-1}\Ga_b\c\Ga_g,\\
\nab_4\big(\widecheck{\ov{\DD}\c\Jk} \big)+\Re\left(\frac{2}{q}\right)\widecheck{\ov{\DD}\c\Jk} &=& O(r^{-1})B+O(r^{-2})\trXc+O(r^{-2})\Xh +O(r^{-2})\Zc\\
&&+O(r^{-3})\widecheck{\DD(\cos\th)}+r^{-1}\Ga_b\c\Ga_g,\\
\nab_4\big(\widecheck{\nab_3\Jk}\big) +\frac{1}{q}\widecheck{\nab_3\Jk} &=& O(r^{-3})\widecheck{e_3(r)}+O(r^{-3})e_3(\cos\th)+O(r^{-2})\ombc\\
&&+O(r^{-2})\Hc+O(r^{-2})\Zc+O(r^{-2})\widecheck{\nab\Jk}+O(r^{-1})\Pc\\
&&+r^{-1}\Ga_b\c\Ga_b.
\eeaa

\item We also have
\beaa
\nab_4 (\DD\hot\Jk_\pm)+\frac{2}{q}\DD\hot\Jk _\pm &=& O(r^{-1})B +O(r^{-2})\trXc+O(r^{-2})\Xh\\
&&+O(r^{-2})\Zc +O(r^{-3})\widecheck{\DD(\cos\th)}+r^{-1}\Ga_b\c\Ga_g,\\
\nab_4\big(\widecheck{\ov{\DD}\c\Jk_\pm} \big)+\Re\left(\frac{2}{q}\right)\widecheck{\ov{\DD}\c\Jk}_\pm &=& O(r^{-1})B+O(r^{-2})\trXc+O(r^{-2})\Xh +O(r^{-2})\Zc\\
&&+O(r^{-3})\widecheck{\DD(\cos\th)}+r^{-1}\Ga_b\c\Ga_g,\\
\nab_4\big(\widecheck{\nab_3\Jk_\pm}\big) +\frac{1}{q}\widecheck{\nab_3\Jk_\pm} &=& O(r^{-3})\widecheck{e_3(r)}+O(r^{-3})e_3(\cos\th)+O(r^{-2})\ombc\\
&&+O(r^{-2})\Hc+O(r^{-2})\Zc+O(r^{-2})\widecheck{\nab\Jk}+O(r^{-1})\Pc\\
&&+r^{-1}\Ga_b\c\Ga_b.
\eeaa
\end{enumerate}
\end{lemma} 

\begin{proof}
To prove the first part of the lemma,  we make use of  the transport equation $ \nab_4\Jk=-q^{-1}\Jk$,  the relations
\beaa
 \ov{\DD}\c\Jk = \frac{4i(r^2+a^2)\cos\th}{|q|^4}+\widecheck{ \ov{\DD}\c\Jk }, \quad  \Jk\c\ov{\Jk} = \frac{2(\sin\th)^2}{|q|^2}, \quad 
  \Jk=O(r^{-1}), 
  \eeaa
and the commutation formulas  of Lemma \ref{lemma:commutation-complexM6} (see also Remark \ref{remark:commutation-complexM6})  as follows
\beaa
&&\nab_4 \DD\hot\Jk = \DD\hot\nab_4\Jk+[\nab_4,\DD\hot]\Jk\\
&=& -\DD\hot\left(\frac{1}{q}\Jk\right) -\frac{1}{2}\tr X(\DD\hot\Jk -Z\hot\Jk) +O(r^{-1})B +O(r^{-2})\Xh\\
&=& -\frac{2}{q}\DD\hot\Jk +\frac{\DD(q)}{q^2}\hot\Jk +\frac{1}{q}\left(\frac{a\ov{q}}{|q|^2}\Jk+\Zc\right)\hot\Jk +O(r^{-1})B +O(r^{-2})\trXc+O(r^{-2})\Xh\\
&&+r^{-1}\Ga_b\c\Ga_g\\
&=& -\frac{2}{q}\DD\hot\Jk  +O(r^{-1})B +O(r^{-2})\trXc+O(r^{-2})\Xh+O(r^{-2})\Zc +O(r^{-3})\widecheck{\DD(\cos\th)}+r^{-1}\Ga_b\c\Ga_g.
\eeaa

Also, in the same vein,
\beaa
\nab_4\widecheck{\ov{\DD}\c\Jk} &=& \nab_4\left(\ov{\DD}\c\Jk-\frac{4i(r^2+a^2)\cos\th}{|q|^4}\right)\\
&=& -\ov{\DD}\c\left(\frac{1}{q}\Jk\right)+[\nab_4, \ov{\DD}\c]\Jk -\frac{8ir\cos\th}{|q|^4}+\frac{8i(r^2+a^2)\cos\th}{|q|^6}e_4(|q|^2)\\
&=& -\frac{1}{q}\ov{\DD}\c\Jk+\frac{\ov{\DD}(q)}{q^2}\c\Jk -\frac{1}{2}\ov{\tr X}(\ov{\DD}\c\Jk +\ov{Z}\c\Jk) +O(r^{-1})B +O(r^{-2})\Xh\\
&& -\frac{8ir\cos\th}{|q|^4}+\frac{16ir(r^2+a^2)\cos\th}{|q|^6}\\
&=& -\left(\frac{1}{q}+\frac{1}{\ov{q}}\right)\ov{\DD}\c\Jk+\frac{a}{q^2}\ov{\Jk}\c\Jk -\frac{1}{\ov{q}}\frac{aq}{|q|^2}\ov{\Jk}\c\Jk -\frac{8ir\cos\th}{|q|^4}+\frac{16ir(r^2+a^2)\cos\th}{|q|^6}\\
&&+O(r^{-1})B+O(r^{-2})\trXc +O(r^{-2})\Xh+O(r^{-2})\Zc+O(r^{-3})\widecheck{\DD(\ov{q})}+r^{-1}\Ga_b\c\Ga_g\\
&=& -\Re\left(\frac{2}{q}\right)\widecheck{\ov{\DD}\c\Jk}+O(r^{-1})B+O(r^{-2})\trXc+O(r^{-2})\Xh +O(r^{-2})\Zc\\
&&+O(r^{-3})\widecheck{\DD(\cos\th)}+r^{-1}\Ga_b\c\Ga_g.
\eeaa

Finally, using  the commutation  formula for $[\nab_4, \nab_3]$ in Lemma \ref{lemma:comm-gen}, we have
\beaa
\nab_4\widecheck{\nab_3\Jk} &=& \nab_4\left(\nab_3\Jk -\frac{\De q}{|q|^4}\Jk\right)= -\nab_3\left(\frac{1}{q}\Jk\right) +[\nab_4,\nab_3]\Jk-\pr_r\left(\frac{\De q}{|q|^4}\right)\Jk+\frac{\De}{|q|^4}\Jk\\
&=& -\frac{1}{q}\nab_3\Jk+\frac{e_3(q)}{q^2}\Jk -2\omb\nab_4\Jk-2(\eta+\ze)\c\nab\Jk-2(\ze\c\Jk)\eta+2(\eta\c\Jk)\ze\\
&&-2\dual\rho\dual\Jk-\pr_r\left(\frac{\De q}{|q|^4}\right)\Jk+\frac{\De}{|q|^4}\Jk.
\eeaa
Continuing
\beaa
 \nab_4\widecheck{\nab_3\Jk}&=& -\frac{1}{q}\widecheck{\nab_3\Jk}+O(r^{-3})\widecheck{e_3(r)}+O(r^{-3})e_3(\cos\th)+O(r^{-2})\ombc+O(r^{-2})\Hc\\
&&+O(r^{-2})\Zc+O(r^{-2})\widecheck{\nab\Jk}+O(r^{-1})\Pc+r^{-1}\Ga_b\c\Ga_b.
\eeaa
This concludes the proof of the first part  of the lemma. The proof of the second part  is similar and left to the reader.
\end{proof}

\begin{lemma}
\lab{lemma:linearizationforJ+-}
The following  equations hold true\footnote{Similar equations hold for $J^{(0)}=\cos \th$, see Lemma \ref{Lemma:otherlinearizedquant}.}
\bea
\bsplit
\nab_4\big(  \widecheck{\DD(J^{(\pm)})}  \big)+\frac{1}{q}  \widecheck{\DD(J^{(\pm)})}   &=  O(r^{-1} ) \widecheck{\tr X} +
O(r^{-1}) \Xh+\Ga_b\c\Ga_g,\\
 \nab_4 \big( \widecheck{\nab_3 J^{(\pm)}}\big)&=O(r^{-2} ) \widecheck{\DD(J^{(\pm)})}   +O(r^{-1}) \Zc+ O(r^{-1})  \Hc+\Ga_b\c\Ga_b.
 \end{split}
\eea
\end{lemma}

\begin{proof}
Recall that, see \eqref{eq:linearizationsforJp},
\beaa
\bsplit
\widecheck{\DD(J^{(\pm)})}= \DD(J^{(\pm)})-\Jk_{\pm}, \qquad 
\widecheck{e_3(J^{(\pm )})}&=e_3(J^{(\pm)}) \pm  \frac{2a}{|q|^2}J^{(\mp)}.
\end{split}
\eeaa
Starting  with the equation $\nab_4(J^{(\pm)}) =0$ we use  the commutator formulas of Lemma \ref{lemma:commutation-complexM6} (see also Remark \ref{remark:commutation-complexM6}) to derive
\beaa
\nab_4( \DD\Jp) &=& [\nab_4, \DD]\Jp=-\frac 1 2 \tr X\DD \Jp  +r^{-1}\Xh \c \dk  \Jp.
\eeaa
Hence
\beaa
\nab_4( \DD\Jp)+\frac 1 2 \tr X \DD \Jp &=& r^{-1}\Xh \c \dk  \Jp.
\eeaa
We further deduce, using $\nab_4\Jk_\pm +\frac{1}{q} \Jk_\pm=0$,
\beaa
\nab_4\big(  \widecheck{\DD(J^{(\pm)})}  \big)+\frac 1 2 \tr X  \widecheck{\DD(J^{(\pm)})}  &=&\nab_4\big(\DD(J^{(\pm)})-\Jk_{\pm}\big)+ \frac 1 2 \tr X \big(\DD(J^{(\pm)})-\Jk_{\pm}\big)\\
&=& r^{-1}\Xh \c \dk  \Jp-\left(\nab_4 \Jk_{\pm} +\frac 1 2 \tr X \Jk_\pm\right)\\
&=&  r^{-1}\Xh \c \dk  \Jp -\frac 1 2 \widecheck{\tr X}  \Jk_{\pm}.
\eeaa
Thus, since $J^{(\pm)}=O(1)$ and $\Jk_\pm=O(r^{-1}) $,
\beaa
\nab_4\big(  \widecheck{\DD(J^{(\pm)})}  \big)+\frac{1}{q} \widecheck{\DD(J^{(\pm)})} &=&  O(r^{-1} ) \widecheck{\tr X} +
O(r^{-1}) \Xh +\Ga_b\c\Ga_g
\eeaa
as stated.

Similarly
\beaa
\nab_4 \nab_3 J^{(\pm)} &=&[\nab_4, \nab_3] J^{(\pm)}=  -2(\ze+\eta )\c \nab J^{(\pm)}=-\Re\Big((\ov{Z+H}) \c \DD  J^{(\pm)}\Big)\\
&=& -\Re\Big((\ov{Z+H}) \c\big(\widecheck{\DD(J^{(\pm)})} +\Jk_{\pm} \big)\Big)
\eeaa
and thus
\beaa
 \nab_4\Big( \nab_3 J^{(\pm)} \pm \frac{2a}{|q|^2} J^{(\mp)}\Big)&=&-\Re\Big((\ov{Z+H}) \c\big(\widecheck{\DD(J^{(\pm)})}+\Jk_{\pm}\big) \Big) \mp\frac{4a r}{|q|^4}  J^{(\mp)}. 
\eeaa
Hence
\beaa
 \nab_4\Big( \widecheck{\nab_3 J^{(\pm)}}\Big) &=&-\Re\Big((\ov{Z+H}) \c \Jk_{\pm} \Big)  \mp\frac{4a r}{|q|^4}  J^{(\mp)} + O(r^{-2})  \widecheck{\DD(J^{(\pm)})}+\Ga_b\c\Ga_b.
 \eeaa
 Now, recalling that , 
 \beaa
 Z= \Zc+     \frac{a\ov{q}}{|q|^2}\Jk, \qquad  H=\Hc+  \frac{aq}{|q|^2}\Jk, \qquad      \Re(\Jk_\pm)\c\Re(\Jk)=\mp \frac{1}{|q|^2}J^{(\mp)},
 \eeaa
  we deduce
\beaa
\Re\Big(( \ov{Z+H}) \c \Jk_\pm \Big)&=&\Re\Big(\ov{ ( \Zc+\Hc) }\c \Jk_\pm\Big) +\left(\frac{a\ov{q} }{|q|^2} +\frac{a q}{|q|^2} \right)\Re\big(  \ov{\Jk}\c \Jk_\pm\big)\\
&=&\frac{4a r}{|q|^2 } \Re(\Jk) \c\Re(\Jk_\pm) +O(r^{-1}) \Zc+ O(r^{-1})  \Hc \\
&=& \mp\frac{4 ar}{|q|^4} J^{(\mp)} +O(r^{-1}) \Zc+ O(r^{-1})  \Hc.
\eeaa
Hence
\beaa
 \nab_4\Big( \widecheck{\nab_3 J^{(\pm)}}\Big) &=& O(r^{-2})  \widecheck{\DD(J^{(\pm)})}+ O(r^{-1}) \Zc+ O(r^{-1})  \Hc+\Ga_b\c\Ga_b
 \eeaa
 as stated.
\end{proof}

%%%%%%%%%%%%%%%%%%%%%%%%%%%%%%%%%%%%%%%%%%%%%

\subsection{The vectorfield $\T$ in $\Mext$}

%%%%%%%%%%%%%%%%%%%%%%%%%%%%%%%%%%%%%%%%%%%%%

We recall that  in  $\Mext$, the vectorfield $\T$  was defined by, see   section \ref{section:Definition-boldT},
\bea
\lab{definition:T-Mext}
\T &:=& \frac{1}{2}\left(e_3+\frac{\Delta}{|q|^2}e_4 -2a\Re(\Jk)^be_b\right).
\eea

\begin{lemma}
\lab{Lemma:$T(rthuJ}
The following hold true.
\begin{enumerate}
\item We have
\bea
\lab{eq:definitionT-Mext1}
\g(\T, \T) &=& -1+\frac{2mr}{|q|^2}.
\eea

\item We have
\bea\lab{eq:computationTuTrTth}
\bsplit 
&\T(u) = 1+\frac{1}{2}\left(\widecheck{e_3(u)} -2a\Re(\Jk)\c\widecheck{ \nab u }\right),\qquad \T(r) = \frac{1}{2}\widecheck{e_3(r)}.
\end{split}
\eea

\item We have
\bea
\lab{eq:definitionT-Mext2}
\bsplit
\T(\cos\th) &= \frac{1}{2}\left(e_3(\cos\th) - 2a\Re(\Jk)\c\widecheck{ \nab \cos\th }\right),\\
\T (J^{(\pm)} )&=\frac 1 2 \widecheck{e_3(J^{(\pm)})}  -a\Re(\Jk)\c \Re\big(  \widecheck{\DD(J^{(\pm)})} \big).
\end{split}
\eea
\end{enumerate}
In particular
\beaa
\T(r)\in r\Ga_b,  \quad \T(u)=1+ r \Ga_b, \quad  \T(\cos\th)\in \Ga_b, \quad \T (J^{(\pm)} )\in\Ga_b.
\eeaa
\end{lemma}

\begin{proof}
Equation \eqref{eq:definitionT-Mext1} follows easily from   $|\Re(\Jk)|^2= \frac{(\sin\th)^2}{|q|^2}$. 

The identities in \eqref{eq:computationTuTrTth}  follow easily from  the relations $e_4(r)=1$, $e_4(u)=0$, $\nab(r)=0$, and the definition of  the linearized quantities $\widecheck{e_3(r)}$,  $\widecheck{e_3(u)}$ and 
$\widecheck{\nab u}$.

To check the  first  identity in \eqref{eq:definitionT-Mext2},  we  make use of $e_4(\th)=0$, the definition of the linearized quantity  $ \widecheck{  \nab\cos \th}=  \nab  \cos \th - \Re( i\Jk)$, and  $\Re(\Jk)\c \Im(\Jk)=0$, which yields
\beaa
\T (\cos \th) &=& \frac{1}{2}  e_3(\cos \th)  -a\Re(\Jk)\c \nab\cos \th=\frac{1}{2}  e_3(\cos \th)  -a\Re(\Jk)\c\widecheck{ \nab\cos \th}.
\eeaa

To check the second identity in \eqref{eq:definitionT-Mext2},     we  make use of $e_4(J^{(\pm)})=0$, the  definition of  the linearized quantities
$\widecheck{e_3(J^{(\pm)})}=e_3(J^{(\pm )}) \pm  \frac{2a}{|q|^2}J^{(\mp)}$, $\widecheck{\DD(J^{(\pm)})}= \DD(J^{(\pm)})-\Jk_{\pm}$, and the identity
$\Re(\Jk_\pm)\c\Re(\Jk)=\mp\frac{1}{|q|^2}J^{(\mp)}$, see   \eqref{eq:propertiesJk}. Thus
\beaa
\T (J^{(\pm)} )&=& \frac{1}{2}\left(e_3+\frac{\Delta}{|q|^2}e_4 -2a\Re(\Jk)^be_b\right)  J^{(\pm)}=\frac 1 2 e_3  ( J^{(\pm)} )   -a\Re(\Jk)\c \nab  J^{(\pm)}\\
&=&\frac 1 2 \widecheck{e_3(J^{(\pm)})} \mp\frac{a}{|q|^2}J^{(\mp)}  -a\Re(\Jk)\c \Re\big(  \DD J^{(\pm)}\big)\\
&=&\frac 1 2 \widecheck{e_3(J^{(\pm)})} \mp\frac{a}{|q|^2}J^{(\mp)}  -a\Re(\Jk)\c \Re\big(  \widecheck{\DD(J^{(\pm)})} +\Jk_\pm\big)\\
&=&\frac 1 2 \widecheck{e_3(J^{(\pm)})}  -a\Re(\Jk)\c \Re\big(  \widecheck{\DD(J^{(\pm)})} \big)  -a\Re(\Jk)\c \Re\big( \Jk_\pm\big) \mp\frac{a}{|q|^2}J^{(\mp)}\\
&=&\frac 1 2 \widecheck{e_3(J^{(\pm)})}  -a\Re(\Jk)\c \Re\big(  \widecheck{\DD(J^{(\pm)})} \big) 
\eeaa
as stated.
\end{proof}

\begin{remark}
In Kerr, we have $\T=\pr_t$ where $\pr_t$ is the coordinate vectorfield corresponding to the Boyer-Lindquist coordinates.
\end{remark}

We   recall below Proposition  \ref{Proposition:deftensorT}
\begin{proposition}
\lab{Proposition:deftensorT-again}
We have $\piT_{44} =0$,  $ \piT_{4a}\in \Ga_g $ and all  other components of  $\piT$, relative to the  frame of $\Mext$,  are in $\Ga_b$.  Moreover
\beaa
g^{ab}\,\piT_{ab} &=& \Ga_g.
\eeaa
In addition
\beaa
\g(\D_a \T, e_4),\,  \g(\D_4 \T, e_a) \in \Ga_g, \qquad 
\g(\D_a \T, e_3),\,  \g(\D_3 \T, e_a) \in \Ga_b,
\eeaa
and
\beaa
 k_{ab}:=\g(\D_a\T,  e_b)= -\frac{2amr\cos\th}{|q|^4}\in_{ab}+\Ga_b.
 \eeaa
\end{proposition}

%%%%%%%%%%%%%%%%%%%%%%%%%%

\subsection{Commutation formulas with $\Lieb_\T$}

%%%%%%%%%%%%%%%%%%%%%%%%%%

We recall the following definition of projected Lie derivative  from section 2.4 of \cite{GKS1}. Given   vectorfields  $X$, $Y$,  the projected Lie  derivative  $\Lieb_XY$ is given by
 \beaa
 \Lieb_X Y :=\Lie_X Y+ \frac 1 2 \g(\Lie_XY, e_3) e_4+  \frac 1 2 \g(\Lie_XY, e_4) e_3.
 \eeaa
 Given  a horizontal covariant k-tensor $U$,  the horizontal  Lie derivative $\Lieb_X U $ is defined   to be the projection of $\Lie_X U $
  to the  horizontal space. Thus, for horizontal indices $A=a_1\ldots a_k$,
  \bea
  \lab{eq:projectedLie}
 ( \Lieb_X U) _{A}:&=&\nab_X U_{A}  +\D_{a_1} X^b U_{b \ldots a_k} +\cdots  +\D_{a_k} X^b U_{ a_1\ldots  b}.
  \eea
   We  recall below Lemma 2.35 in \cite{GKS1}.
 
  \begin{lemma}
  \lab{Lemma:commLieb_Tnab}
  The following commutation formulas hold true for a horizontal covariant k-tensor $U$ and  a vectorfield $X$
  \bea
  \bsplit
  \nab_b(\Lieb_X U_{A})-\Lieb_X(\nab_b U_{A})&=\sum_{j=1}^k  \GabbX_{a_jb  c }U_{a_1\ldots\,\,\,\,\ldots a_k}^{\,\,\,\,\,\,\,\,\,\,\,c},\\
   \nab_4(\Lieb_X U_A)-\Lieb_X(\nab_4 U_A) +\nab_{\Lieb_X e_4} U_{A}&=\sum_{j=1}^k  \GabbX_{a_j4  c }U_{a_1\ldots\,\,\,\,\ldots a_k}^{\,\,\,\,\,\,\,\,\,\,\,c},\\
   \nab_3(\Lieb_X U_A)-\Lieb_X(\nab_3 U_A) +\nab_{\Lieb_X e_3} U_{A}&=\sum_{j=1}^k  \GabbX_{a_j3  c }U_{a_1\ldots\,\,\,\,\ldots a_k}^{\,\,\,\,\,\,\,\,\,\,\,c},
  \end{split}
  \eea
  with\footnote{Here,  $\piX_{ab}$ is treated as a horizontal  symmetric 2-tensor, and $\piX_{a4}$,  $\piX_{a3}$, as   horizontal 1-forms.}
  \bea
  \bsplit
  \GabbX_{abc }&=\frac 1 2 (\nab_a\piX_{bc }+\nab_b\piX_{a c }-\nab_c \piX_{ab}),\\
  \GabbX_{a4b} &= \frac 1 2 (\nab_a\piX_{4b }+\nab_4\piX_{a b  }-\nab_b  \piX_{a4}),\\
\GabbX_{a3b} &= \frac 1 2 (\nab_a\piX_{3b }+\nab_3\piX_{a b  }-\nab_b  \piX_{a3}).
  \end{split}
  \eea
  \end{lemma}

 We apply Lemma \ref{Lemma:commLieb_Tnab} to the case when $X$ is the vectorfield  $\T$.  
\begin{lemma}
\lab{Lemma:commLieb_Tnab-simb}
The following holds true:
\begin{enumerate}
\item We have
\bea
\lab{eq:commLieb_Te_4}
\Lieb_\T e_4 \in  \Ga_b.
\eea

\item   For any   horizontal $k$-tensor $U$, we have
\bea
\lab{eq:commLieb_Tnab}
 \nab_4 (\Lieb_\T U_{A})-\Lieb_\T(\nab_4 U_{A})&= r^{-1}\Ga_b\c\dkb U+r^{-1}\dk\Ga_b\c U.
\eea

\item For any   horizontal $k$-tensor $U$, we have
\bea
\lab{eq:commLieb_Tnab2}
 \nab(\Lieb_\T U_{A})-\Lieb_\T(\nab U_{A})&= r^{-1} \dkb  \Ga_b \c U.
\eea
\end{enumerate}
\end{lemma}

\begin{proof}
In view of the definition of $\T$ in \eqref{definition:T-Mext}, we have
\beaa
2[e_4, \T] &=& \left[e_4, e_3+\frac{\Delta}{|q|^2}e_4 -2a\Re(\Jk)^be_b\right]\\
&=&[e_4, e_3] +e_4\left(\frac{\Delta}{|q|^2}\right)e_4 -2a\nab_4\Re(\Jk)^be_b+2a\Re(\Jk)^b\nab_be_4\\
&=& -2(\ze+\eta)_be_b -2\omb e_4+\pr_r\left(\frac{\Delta}{|q|^2}\right)e_4 +2a\Re\left(\frac{1}{q}\Jk\right)^be_b+2a\Re(\Jk)^b(\chi_{bc}e_c-\ze_be_4)\\
&=& -2(\widecheck{\ze}+\widecheck{\eta})_be_b -2\ombc e_4 +2a\Re(\Jk)^b\left(\frac{1}{2}\trchc\de_{bc}e_c+\frac{1}{2}\widecheck{\atrch}\in_{bc}e_c+\chih_{bc}e_c-\widecheck{\ze}_be_4\right).
\eeaa
Hence, since  $\Lieb_\T( e_4 )$ is the  horizontal projection of $[\T, e_4]$, we obtain $\Lieb_T( e_4 )\in  \Ga_b$ 
which is \eqref{eq:commLieb_Te_4}. 

Next, in view of the form of $ \GabbX$ in Lemma \ref{Lemma:commLieb_Tnab} with the particular choice $X=\T$, and together with Proposition \ref{Proposition:deftensorT-again}, we have 
 \beaa
  \bsplit
  ^{(\T)} \Gabb_{abc }&=\frac 1 2 (\nab_a\,{}^{(\T)}\pi_{bc }+\nab_b\,{}^{(\T)}\pi_{a c }-\nab_c \,{}^{(\T)}\pi_{ab}) \in r^{-1} \dkb \Ga_b ,\\
 ^{(\T)} \Gabb_{a4b} &= \frac 1 2 (\nab_a\,{}^{(\T)}\pi_{4b }+\nab_4\,{}^{(\T)}\pi_{a b  }-\nab_b\,{}^{(\T)}\pi_{a4})\in  r^{-1} \dk \Ga_b.
  \end{split}
  \eeaa
Using Lemma \ref{Lemma:commLieb_Tnab} with $X=\T$, we infer
  \beaa
  \bsplit
  \nab_b(\Lieb_\T U_{A})-\Lieb_\T(\nab_b U_{A})&= r^{-1} \dkb \Ga_b\c U,\\
   \nab_4(\Lieb_\T U_A)-\Lieb_\T(\nab_4 U_A) +\nab_{\Lieb_X e_4} U_{A}&= r^{-1} \dk \Ga_b\c U.
  \end{split}
  \eeaa
The first identity yields \eqref{eq:commLieb_Tnab2} while the second identity, together with  \eqref{eq:commLieb_Te_4}, yields \eqref{eq:commLieb_Tnab}. This concludes the proof of the lemma. 
\end{proof}

%%%%%%%%%%%%%%%%%%%%%%%%%%%%%

\subsection{Relation between $\Lieb_\T$ and $\nab_3$}

%%%%%%%%%%%%%%%%%%%%%%%%%%%%%

\begin{lemma}
If $U$ is a horizontal $k$-tensor, we have
\bea
\lab{eq:Leib_T-nab_3U}
\Lieb_\T U_A  &=& \frac 1 2 \nab_3 U_A +\frac 1 2 \frac{\De}{|q|^2} \nab_4 U_A + O(r^{-1})\nab  U_A +  O(r^{-3} ) U+ \Ga_b \c U. 
\eea
In particular
\bea\lab{eq:Leib_T-nab_3U:immediateconsequences}
\bsplit
\Lieb_\T U &=  \frac{1}{2}\nab_3 U+O(r^{-1})\dk^{\leq 1}U+\Ga_b\c U,\\
\Lieb^2_\T U&=  \frac{1}{4}\nab_3^2U +O(r^{-1})\dk^{\leq 1}\nab_3U+O(r^{-2})\dk^{\leq 2}U+\dk^{\leq 1}(\Ga_b\c U).
\end{split}
\eea
\end{lemma}

\begin{proof}
Recall from Proposition \ref{Proposition:deftensorT-again} that we have
\beaa
 k_{ab}=\g(\D_a\T,  e_b)= -\frac{2amr\cos\th}{|q|^4}\in_{ab}+\Ga_b.
 \eeaa
Together with  \eqref{eq:projectedLie},   for $p=1,2 $ and $A = a_1 \ldots a_p$, we obtain 
\beaa
 ( \Lieb_\T U) _{A} &=&\nab_\T U_{A}  +\D_{a_1} \T^b U_{b \ldots a_p} +\cdots  +\D_{a_p} \T^b U_{ a_1\ldots  b}\\
 &=& \nab_\T U_{A}   -\frac{2 p amr\cos\th}{|q|^4} \dual U_A+\Ga_b \c U \\
 &=& \frac 1 2 \nab_3 U_A +\frac 1 2  \frac{\Delta}{|q|^2}\nab_4 U_A -  \Re(\Jk)^b\nab_b U_A +O(r^{-3})U+\Ga_b \c U 
 \\
 &=& \frac 1 2 \nab_3 U_A +\frac 1 2  \frac{\Delta}{|q|^2}\nab_4 U_A  + O(r^{-1})\nab  U_A +  O(r^{-3})U + \Ga_b \c U 
 \eeaa
 as stated. The two other identities easily follow from this one.
\end{proof}

In the same spirit, we have the following more precise decomposition.
\begin{corollary}\lab{cor:decompositionofnab3inLiebTnab4andlot}
We have
\bea
\lab{eq:Leib_T-nab_3U-Bis}
\nab_3  &=& 2\Lieb_\T - \frac{\De}{|q|^2}\nab_4+O(r^{-1})\nab+O(r^{-3})+\Ga_b.
\eea
\end{corollary}

\begin{proof}
Note from the proof of the previous lemma the identity 
\beaa
\Lieb_\T f &=& \nab_\T f+f\c k.
\eeaa
Since 
\beaa
2\T = e_3+\frac{\De}{|q|^2}e_4+O(r^{-1})\nab, \qquad k=O(r^{-3})+\Ga_b,
\eeaa
we infer
\beaa
\nab_3  &=& 2\Lieb_\T- \frac{\De}{|q|^2}\nab_4+O(r^{-1})\nab+O(r^{-3})+\Ga_b
\eeaa
as desired.
\end{proof}

%%%%%%%%%%%%%%%%%%%%%%%%%%%%%%%%%%%%%%%%%%%%%

\section{Properties of the  spheres  $S(u,r)$}

%%%%%%%%%%%%%%%%%%%%%%%%%%%%%%%%%%%%%%%%%%%%%

%%%%%%%%%%%%%%%%%%%%%%%%

\subsection{An orthonormal frame  of $S(u,r)$}
\lab{subsection:orthonormalbasisS(u,r)}

%%%%%%%%%%%%%%%%%%%%%%%%

The following lemma exhibits an orthonormal frame of $S(u,r)$.  
\begin{lemma}\lab{Lemma:orthonormalbasisS(u,r)}
Let $(e_1, e_2)$ an orthonormal  basis of the horizontal structure associated to the PG structure of $\Mext$. Then, there exists an orthonormal basis $(e_1', e_2')$ of the tangent space of $S(u,r)$ of the following form
\beaa
e_a' &=& \left(\de_a^b+\frac{1}{2}\fb_af^b\right) e_b+\frac{1}{2}\fb_a e_4+\left(\frac{1}{2}f_a+\frac{1}{8}|f|^2\fb_a\right)e_3,
\eeaa
where the 1-forms $f$ and $\fb$ are given by
\beaa
\bsplit
f  &= -\frac{4}{e_3(u)+ \sqrt{(e_3(u))^2 +4|\nab u|^2e_3(r)}}\nab u,\quad 
\fb = \frac{2e_3(r)}{\sqrt{(e_3(u))^2 +4|\nab u|^2e_3(r)}}\nab u.
\end{split}
\eeaa
Also,
\bea
\lab{eq:formulaforchangeofframecoeffandfbarforframetangenttosphereSur}
f=\left(-\frac{2|q|^2}{r^2+a^2+\Si}+r\Ga_b\right)\nab u, \qquad \fb=\left(-\frac{\De}{\Si}+r\Ga_b\right)\nab u.
\eea
\end{lemma}

\begin{proof}
In view of the definition of $\widecheck{e_3(u)}$, $\widecheck{\nab u}$ and $\widecheck{e_3(r)}$, and since $\widecheck{e_3(u)}\in r\Ga_b$, $\widecheck{e_3(r)}\in r\Ga_b$ and $\widecheck{\nab u}\in\Ga_b$, we have
\beaa
e_3(u)=2+O(r^{-2})+r\Ga_b, \qquad (e_3(u))^2 +4|\nab u|^2e_3(r)=4+O(r^{-2})+r\Ga_b.
\eeaa 
Thus, since $r\geq r_0$ in $\Mext$, we infer
\beaa
e_3(u)>0, \qquad (e_3(u))^2 +4|\nab u|^2e_3(r)>0, \quad\textrm{on}\quad\Mext.
\eeaa
Thus, we may apply Lemma \ref{Lemma:Transformation-principal-to-integrable frames} which yields the existence of an orthonormal basis $(e_1', e_2')$ of the tangent space of $S(u,r)$ of the form
\beaa
e_a' &=& \left(\de_a^b+\frac{1}{2}\fb_af^b\right) e_b+\frac{1}{2}\fb_a e_4+\left(\frac{1}{2}f_a+\frac{1}{8}|f|^2\fb_a\right)e_3,
\eeaa
with the 1-forms $f$ and $\fb$ given by
\beaa
\bsplit
f  &= -\frac{4}{e_3(u)+ \sqrt{(e_3(u))^2 +4|\nab u|^2e_3(r)}}\nab u,\quad 
\fb = \frac{2e_3(r)}{\sqrt{(e_3(u))^2 +4|\nab u|^2e_3(r)}}\nab u
\end{split}
\eeaa
as stated. Then, \eqref{eq:formulaforchangeofframecoeffandfbarforframetangenttosphereSur} follows from the above form of $f$ and $\fb$, the definition of $\widecheck{e_3(u)}$, $\widecheck{\nab u}$ and $\widecheck{e_3(r)}$, and the fact that $\widecheck{e_3(u)}\in r\Ga_b$, $\widecheck{e_3(r)}\in r\Ga_b$ and $\widecheck{\nab u}\in\Ga_b$.
\end{proof}

\begin{lemma}\lab{lemma:behavioroffandfbinchangeofframeorthonormalbasisSur}
Let $f$ and $\fb$ the horizontal 1-forms  given by \eqref{eq:formulaforchangeofframecoeffandfbarforframetangenttosphereSur}. Also, let the scalar function $\la$ given by \eqref{def:transition-functs:ffbla}, i.e. 
\beaa
\la &=& 1+\frac{1}{2}f\c\fb  +\frac{1}{16}|f|^2|\fb|^2.
\eeaa
Then, we have
\beaa
f=O(r^{-1}), \qquad \fb=O(r^{-1}), \qquad f=\fb+O(r^{-2})+\Ga_b, \qquad \la=1+O(r^{-2}),
\eeaa
as well as 
\beaa
\fb -\frac{\De}{|q|^2}\left(f+\frac{1}{4}|f|^2\fb\right) &=&  \Ga_b.
\eeaa
\end{lemma}

\begin{proof}
In view of \eqref{eq:formulaforchangeofframecoeffandfbarforframetangenttosphereSur}, we have
\beaa
f=\left(-\frac{2|q|^2}{r^2+a^2+\Si}+r\Ga_b\right)\nab u, \qquad \fb=\left(-\frac{\De}{\Si}+r\Ga_b\right)\nab u,
\eeaa
which, together with the formula for $\la$ in terms of $f$ and $\fb$, immediately implies the first three claimed identities.

Concerning the last one, we have
\beaa
\fb -\frac{\De}{|q|^2}\left(f+\frac{1}{4}|f|^2\fb\right) &=& \left(h+r\Ga_b\right)\nab u
\eeaa
where
\beaa
h &=& -\frac{\De}{\Si}+\frac{\De}{|q|^2}\left(\frac{2|q|^2}{r^2+a^2+\Si}+\frac{1}{4}\left(\frac{2|q|^2}{r^2+a^2+\Si}\right)^2|\nab u|^2\frac{\De}{\Si}\right)\\
&=& \frac{\De}{\Si(r^2+a^2+\Si)^2}\Big(-(r^2+a^2+\Si)^2+2(r^2+a^2+\Si)\Si+\De a^2(\sin\th)^2\Big)\\
&=& \frac{\De}{\Si(r^2+a^2+\Si)^2}\Big(\Si^2-(r^2+a^2)^2+\De a^2(\sin\th)^2\Big)=0
\eeaa
so that 
\beaa
\fb -\frac{\De}{|q|^2}\left(f+\frac{1}{4}|f|^2\fb\right) &=& r\Ga_b\nab u =\Ga_b
\eeaa
as desired.
\end{proof}

%%%%%%%%%%%%%%%%%%%%%%%%%%%%%%%%%%%%%%%%%%%%

\subsection{Comparison of horizontal derivatives and derivatives tangential to $S(u,r)$}
\lab{sec:comparisionhorizontalderivativesandtangentialtoSur}

%%%%%%%%%%%%%%%%%%%%%%%%%%%%%%%%%%%%%%%%%%%%

\begin{lemma}
\lab{lemma:Comparison.horizontalderivarives}
We have
\beaa
\g(\D_{e_a'}e_b', e_c') &=& \left(\de_a^d+\frac{1}{2}\fb_af^d\right)\g\left(\D_{e_d}e_b, e_c\right)+\frac{1}{2}\fb_a\g\left(\D_{ e_4}e_b, e_c\right)\\
&&+\left(\frac{1}{2}f_a+\frac{1}{8}|f|^2\fb_a\right)\g\left(\D_{e_3}e_b, e_c\right) +\frac{1}{2}\fb_cf_be_a'(\log\la)
  - \frac{1}{2}\fb_c\la^{-1}\chi_{ab}' \\
  &&+\frac{1}{2}\fb_b\la^{-1}\chi_{ac}'  + \frac{1}{2}\ze_a'\fb_cf_b - \frac{1}{4}\fb_c\la^{-1}\chi_{ad}' \fb_bf_d -\frac{1}{2}f_c\chib_{ab}+\frac{1}{2}f_b\chib_{ac}+\err[\g(\D_{e_a'}e_b', e_c')],
\eeaa
where $\err[\g(\D_{e_a'}e_b', e_c')]$ contains all the terms depending on $(f, \fb, \Ga)$, without derivative, and at least quadratic in $(f, \fb)$, and where the scalar function $\la$ is given by \eqref{def:transition-functs:ffbla}.\end{lemma}

\begin{proof}
See  appendix \ref{appendix:ProofLemma-lemma:Comparison.horizontalderivarives}.
\end{proof}

\begin{proposition}\lab{cor:formuallikningnabprimeandnabintransfoformula}
Let $V$ a horizontal k-tensor. Then, for horizontal indices  $B=b_1\cdots b_k$ and $B^c_{(j)} = b_1\cdots b_{j-1}c\,b_{j+1}\cdots b_k$ we have
\beaa
\nab_a' V_B&=& \nab_a V_B+\frac{1}{2}\fb_a f^c\nab_c V_B+\frac{1}{2}\fb_a\nab_4V_B+\left(\frac{1}{2}f_a+\frac{1}{8}|f|^2\fb_a\right)\nab_3V_B\\
&& -\sum_{j=1}^k\Bigg\{\frac{1}{2}\fb_cf_{b_j}\nab_a'(\log\la) - \frac{1}{2}\fb_c\la^{-1}\chi_{ab_j}' +\frac{1}{2}\fb_{b_j}\la^{-1}\chi_{ac}'  + \frac{1}{2}\ze_a'\fb_cf_{b_j}\\
&& - \frac{1}{4}\fb_c\la^{-1}\chi_{ad}' \fb_{b_j}f_d -\frac{1}{2}f_c\chib_{ab_j}+\frac{1}{2}f_{b_j}\chib_{ac}+\err[\g(\D_{e_a'}e_{b_j}', e_c')]\Bigg\}V_{B^c_{(j)}},
\eeaa
\beaa
\dual\nab_a' V_B &=& \dual\nab_a V_B+\frac{1}{2}\dual\fb_a f^c\nab_c V_B+\frac{1}{2}\dual\fb_a\nab_4V_B+\left(\frac{1}{2}\dual f_a+\frac{1}{8}|f|^2\dual\fb_a\right)\nab_3V_B\\
&& -\sum_{j=1}^k\Bigg\{\frac{1}{2}\fb_cf_{b_j}\dual\nab_a'(\log\la) - \frac{1}{2}\fb_c\la^{-1}\dual\chi_{ab_j}' +\frac{1}{2}\fb_{b_j}\la^{-1}\dual\chi_{ac}'  + \frac{1}{2}\dual\ze_a'\fb_cf_{b_j}\\
&& - \frac{1}{4}\fb_c\la^{-1}\dual\chi_{ad}' \fb_{b_j}f_d -\frac{1}{2}f_c\dual\chib_{ab_j}+\frac{1}{2}f_{b_j}\dual\chib_{ac}+\in_{ad}\err[\g(\D_{e_d'}e_{b_j}', e_c')]\Bigg\}V_{B^c_{(j)}},
\eeaa
and
\beaa
\DD_a' V_B &=& \DD_a V_B+\frac{1}{2}\underline{F}_a f^c\nab_c V_B+\frac{1}{2}\underline{F}_a\nab_4V_B+\left(\frac{1}{2}F_a+\frac{1}{8}|f|^2\underline{F}_a\right)\nab_3V_B +(E[V])_{aB},
\eeaa
with
\beaa
(E[V])_{aB} & =&-\sum_{j=1}^k\Bigg\{\frac{1}{2}\fb_cf_{b_j}\DD_a'(\log\la) - \frac{1}{2}\fb_c\la^{-1}X_{ab_j}' +\frac{1}{2}\fb_{b_j}\la^{-1}X_{ac}'  + \frac{1}{2}Z_a'\fb_cf_{b_j}\\
&& - \frac{1}{4}\fb_c\la^{-1}X_{ad}' \fb_{b_j}f_d -\frac{1}{2}f_c\Xb_{ab_j}+\frac{1}{2}f_{b_j}\Xb_{ac}   \Bigg\}V_{B^c_{(j)}}\\
&& +\bigg\{\err[\g(\D_{e_a'}e_{b_j}', e_c')]+i\in_{ad}\err[\g(\D_{e_d'}e_{b_j}', e_c')]\bigg\}  V_{B^c_{(j)}},
\eeaa
where the horizontal 1-forms $f$ and $\fb$ are given by \eqref{eq:formulaforchangeofframecoeffandfbarforframetangenttosphereSur}, and the horizontal complex 1-forms $F$ and $\underline{F}$ are given by
\beaa
F:=f+i\dual f, \qquad \underline{F}:=\fb+i\dual\fb.
\eeaa
\end{proposition}

\begin{proof}
We prove the formula for $\nab'V$.  For simplicity, we  do  it for a 1-tensor $V$. We have
\beaa
\nab_a' V_b &=& e_a'(V_b) -\g(\D_{e_a'}e_b', e_c')V_c.
\eeaa
Since
\beaa
e_a' &=& \left(\de_a^b+\frac{1}{2}\fb_af^b\right) e_b+\frac{1}{2}\fb_a e_4+\left(\frac{1}{2}f_a+\frac{1}{8}|f|^2\fb_a\right)e_3,
\eeaa
and
\beaa
\g(\D_{e_a'}e_b', e_c') &=& \left(\de_a^d+\frac{1}{2}\fb_af^d\right)\g\left(\D_{e_d}e_b, e_c\right)+\frac{1}{2}\fb_a\g\left(\D_{ e_4}e_b, e_c\right)\\
&&+\left(\frac{1}{2}f_a+\frac{1}{8}|f|^2\fb_a\right)\g\left(\D_{e_3}e_b, e_c\right)
  - \frac{1}{2}\fb_c\chi_{ab}' +\frac{1}{2}\fb_b\chi_{ac}'  + \frac{1}{2}\ze_a'\fb_cf_b\\
  && - \frac{1}{4}\fb_c\chi_{ad}' \fb_bf_d -\frac{1}{2}f_c\chib_{ab}+\frac{1}{2}f_b\chib_{ac}+\err[\g(\D_{e_a'}e_b', e_c')]
\eeaa
we infer
\beaa
\nab_a' V_b &=& \left[\left(\de_a^b+\frac{1}{2}\fb_af^d\right) e_d+\frac{1}{2}\fb_a e_4+\left(\frac{1}{2}f_a+\frac{1}{8}|f|^2\fb_a\right)e_3\right](V_b)\\
&& -\Bigg\{\left(\de_a^d+\frac{1}{2}\fb_af^d\right)\g\left(\D_{e_d}e_b, e_c\right)+\frac{1}{2}\fb_a\g\left(\D_{ e_4}e_b, e_c\right)\\
&&+\left(\frac{1}{2}f_a+\frac{1}{8}|f|^2\fb_a\right)\g\left(\D_{e_3}e_b, e_c\right)
  - \frac{1}{2}\fb_c\chi_{ab}' +\frac{1}{2}\fb_b\chi_{ac}'  + \frac{1}{2}\ze_a'\fb_cf_b\\
&& - \frac{1}{4}\fb_c\chi_{ad}' \fb_bf_d\ -\frac{1}{2}f_c\chib_{ab}+\frac{1}{2}f_b\chib_{ac}+\err[\g(\D_{e_a'}e_b', e_c')]\Bigg\}V_c
\eeaa
and hence
\beaa
\nab_a' V_b &=& \left(\de_a^b+\frac{1}{2}\fb_af^d\right)\nab_dV_b+\frac{1}{2}\fb_a\nab_4V_b+\left(\frac{1}{2}f_a+\frac{1}{8}|f|^2\fb_a\right)\nab_3V_b\\
&& -\Bigg\{ - \frac{1}{2}\fb_c\chi_{ab}' +\frac{1}{2}\fb_b\chi_{ac}'  + \frac{1}{2}\ze_a'\fb_cf_b - \frac{1}{4}\fb_c\chi_{ad}' \fb_bf_d\\
&& -\frac{1}{2}f_c\chib_{ab}+\frac{1}{2}f_b\chib_{ac}+\err[\g(\D_{e_a'}e_b', e_c')]\Bigg\}V_c
\eeaa
as stated. The formulas for $\dual\nab_a' V_B$ and $\DD_a' V_B$ are proved in the same manner. 
\end{proof}

\begin{corollary}\lab{cor:formuallikningnabprimeandnabintransfoformula:bis}
Let $f$ and $\fb$ the horizontal 1-forms  given by \eqref{eq:formulaforchangeofframecoeffandfbarforframetangenttosphereSur}, and let the scalar function $\la$ given by \eqref{def:transition-functs:ffbla}. Then, recalling the notation
of Definition \ref{def:ordermagnitude}, we have
\bea
\nab' &=& \Big(1+O(r^{-2})\Big)\nab+O(r^{-1})\Lieb_\T +O(r^{-3})+r^{-1}\Ga_b\dk^{\leq 1}.
\eea
Also, we have 
\bea
\nab' &=& (1+O(r^{-2}) )\nab+O(r^{-1})\nab_4+O(r^{-1})\nab_3+O(r^{-3})  +r^{-1}\Ga_b\dk^{\leq 1}.
\eea
\end{corollary}

\begin{proof}
Note that the   second identity is an immediate  consequence of the first in view of \eqref{eq:Leib_T-nab_3U-Bis}. Thus, we focus on proving the first one. Recall from Proposition \ref{cor:formuallikningnabprimeandnabintransfoformula} that we have for a horizontal tensor $V$  and indices $B=b_1\cdots b_k$
\beaa
\nab_a' V_B &=& \nab_a V_B+\frac{1}{2}\fb_a f^c\nab_c V_B+\frac{1}{2}\fb_a\nab_4V_B+\left(\frac{1}{2}f_a+\frac{1}{8}|f|^2\fb_a\right)\nab_3V_B\\
&& -\sum_{j=1}^k\Bigg\{\frac{1}{2}\fb_cf_{b_j}\nab_a'(\log\la) - \frac{1}{2}\fb_c\la^{-1}\chi_{ab_j}' +\frac{1}{2}\fb_{b_j}\la^{-1}\chi_{ac}'  + \frac{1}{2}\ze_a'\fb_cf_{b_j}\\ 
&& - \frac{1}{4}\fb_c\la^{-1}\chi_{ad}' \fb_{b_j}f_d -\frac{1}{2}f_c\chib_{ab_j}+\frac{1}{2}f_{b_j}\chib_{ac}+\err[\g(\D_{e_a'}e_{b_j}', e_c')]\Bigg\}  V_{B^c_{(j)}}.
\eeaa
Thus, using in particular $f=O(r^{-1})$ and $\fb=O(r^{-1})$ and $\la=1+O(r^{-2})$,  see  Lemma \ref{lemma:behavioroffandfbinchangeofframeorthonormalbasisSur},  
\beaa
\nab' &=& \Big(1+O(r^{-2})\Big)\nab+\frac{1}{2}\fb\nab_4+\left(\frac{1}{2}f+\frac{1}{8}|f|^2\fb\right)\nab_3+E+O(r^{-3}),
\eeaa
where $E$ is given by
\beaa
(EV)_{a B} =  -\sum_{j=1}^k\left\{ - \frac{1}{2}\fb_c\la^{-1}\chi_{ab_j}' +\frac{1}{2}\fb_{b_j}\la^{-1}\chi_{ac}'   -\frac{1}{2}f_c\chib_{ab_j}+\frac{1}{2}f_{b_j}\chib_{ac}\right\}V_{B^c_{(j)}}
\eeaa
We use  $f=O(r^{-1})$ and $\fb=O(r^{-1})$  and  the transformation formulas  for the Ricci coefficients 
  of Proposition \ref{Proposition:transformationRicci}         for $\chi'$   to deduce
 \beaa
 \la^{-1}\chi_{ab}' &=& \chi_{ab}+O(r^{-2})+\Ga_g
 \eeaa
 and hence, using also $f=\fb+O(r^{-2})+\Ga_b$ in view of Lemma \ref{lemma:behavioroffandfbinchangeofframeorthonormalbasisSur}, we infer
\beaa
(EV)_{a B} &=&  -\frac{1}{2}\sum_{j=1}^k\Big\{ - f_c\big(\chi_{ab_j}+\chib_{ab_j}\big) +f_{b_j}\big(\chi_{ac}  +\chib_{ac}\big)+O(r^{-3})+r^{-1}\Ga_b\Big\}V_{B^{(j)}}.
\eeaa
Since  $ \chi_{ab}+\chib_{ab} = O(r^{-2})+\Ga_b$,
\bea
\lab{eq:(EV)_{aB}}
(EV)_{aB} =  -\frac{1}{2}\sum_{j=1}^k\Big\{O(r^{-3})+r^{-1}\Ga_b\Big\}V_{B^{(j)}}
\eea
and hence
\beaa
\nab' &=& \Big(1+O(r^{-2})\Big)\nab+\frac{1}{2}\fb\nab_4+\left(\frac{1}{2}f+\frac{1}{8}|f|^2\fb\right)\nab_3+E+O(r^{-3})+r^{-1}\Ga_b.
\eeaa
Recalling the decomposition \eqref{eq:Leib_T-nab_3U-Bis}, i.e.
\beaa
\nab_3  &=& 2\Lieb_\T -\frac{\De}{|q|^2}\nab_4+O(r^{-1})\nab+O(r^{-3})+\Ga_b,
\eeaa
we infer, using also  $f=O(r^{-1})$ and $\fb=O(r^{-1})$ in view of Lemma \ref{lemma:behavioroffandfbinchangeofframeorthonormalbasisSur}, 
\beaa
\nab' &=& \Big(1+O(r^{-2})\Big)\nab+O(r^{-1})\Lieb_\T+\frac{1}{2}\left(\fb -\frac{\De}{|q|^2}\left(f+\frac{1}{4}|f|^2\fb\right)\right)\nab_4\\
&&+O(r^{-3})+r^{-1}\Ga_b.
\eeaa
Recalling also, from Lemma \ref{lemma:behavioroffandfbinchangeofframeorthonormalbasisSur},  
\beaa
\fb -\frac{\De}{|q|^2}\left(f+\frac{1}{4}|f|^2\fb\right) &=& \Ga_b,
\eeaa
we deduce
\beaa
\nab' &=& \Big(1+O(r^{-2})\Big)\nab+O(r^{-1})\Lieb_\T +O(r^{-3})+r^{-1}\Ga_b\dk^{\leq 1}
\eeaa
as stated.
\end{proof}

%%%%%%%%%%%%%%%%%%%%%%%%%%%%%%%%%%%%%%%%

\subsection{Derivatives in $e_4$  of integrals on 2-spheres $S(u,r)$}

%%%%%%%%%%%%%%%%%%%%%%%%%%%%%%%%%%%%%%%%

\begin{lemma}\lab{lemma:e4derivativeofintegralonS}
We have for a scalar function $h$
\beaa
e_4\left(\int_{S(r,u)}h\right) &=& \int_{S(r,u)}\Big(e_4(h)+\de^{ab}\g(\D_{e_a'}e_4, e_b')h\Big)
\eeaa 
where
\beaa
\de^{ab}\g(\D_{e_a'}e_4, e_b') &=& \left(1+f\c\fb+\frac{1}{2}|\fb|^2|f|^2\right)\trch+2f\c\chih\c\fb+|\fb|^2f\c\chih\c f\\
&&+f\c\eta+\frac{1}{2}(f\c\fb)(f\c\eta)+\frac{1}{4}|f|^2(\fb\c\eta)+\frac{1}{8}|f|^2|\fb|^2 f\c\eta\\
&&+ f\c\ze+\frac{1}{2}(f\c\fb)(f\c\ze)+\frac{1}{4}|f|^2(\fb\c\ze)+\frac{1}{8}(f\c\ze)|f|^2|\fb|^2\\
&&-\omb\left(|f|^2+\frac{1}{2}|f|^2(f\c\fb)+\frac{1}{16}|f|^4|\fb|^2\right).
\eeaa
\end{lemma}

\begin{proof}
Recall that $e_4(r)=1$ and $e_4(u)=0$. Consider coordinates $(x^1, x^2)$ on $S(u, r_0)$ and transport it by $e_4(x^1)=e_4(x^2)=0$. Then, $e_4=\pr_r$ and hence, for a scalar function $h$, we have
\beaa
e_4\left(\int_{S(r,u)}h\right) &=& \pr_r\left(\int h\sqrt{|g|}dx^1 dx^2\right)\\
&=& \int \left(\pr_r(h)+\frac{1}{\sqrt{|g|}}\pr_r(\sqrt{|g|})h\right)dx^1 dx^2.
\eeaa 
Now, we have
\beaa
\pr_r(g_{ab}) &=& \pr_r\left(g\left(\frac{\pr}{\pr x^a},\frac{\pr}{\pr x^b}\right)\right)= \g\left(\D_{\frac{\pr}{\pr r}}\frac{\pr}{\pr x^a},\frac{\pr}{\pr x^b}\right)+\g\left(\frac{\pr}{\pr x^a},\D_{\frac{\pr}{\pr r}}\frac{\pr}{\pr x^b}\right)\\
&=& \g\left(\D_{\frac{\pr}{\pr x^a}}\frac{\pr}{\pr r},\frac{\pr}{\pr x^b}\right)+\g\left(\frac{\pr}{\pr x^a},\D_{\frac{\pr}{\pr x^b}}\frac{\pr}{\pr r}\right)\\
&=& \g\left(\D_{\frac{\pr}{\pr x^a}}e_4,\frac{\pr}{\pr x^b}\right)+\g\left(\frac{\pr}{\pr x^a},\D_{\frac{\pr}{\pr x^b}}e_4\right).
\eeaa
Since
\beaa
\frac{1}{\sqrt{|g|}}\pr_r(\sqrt{|g|}) &=& \frac{1}{2}g^{ab}\pr_rg_{ab},
\eeaa
we infer
\beaa
\frac{1}{\sqrt{|g|}}\pr_r(\sqrt{|g|}) &=& g^{ab}\g\left(\D_{\frac{\pr}{\pr x^a}}e_4,\frac{\pr}{\pr x^b}\right)
\eeaa
and hence, for an orthonormal basis $(e_1', e_2')$ of $S(u,r)$, we have
\beaa
\frac{1}{\sqrt{|g|}}\pr_r(\sqrt{|g|}) &=& \de^{ab}\g(\D_{e_a'}e_4, e_b').
\eeaa
Thus, we deduce 
\beaa
e_4\left(\int_{S(r,u)}h\right) &=& \int_{S(r,u)}\Big(e_4(h)+\de^{ab}\g(\D_{e_a'}e_4, e_b')h\Big).
\eeaa 

Next, for simplicity, we write
\beaa
e_a' &=& j_a^be_b+k_ae_4+l_ae_3
\eeaa
where
\bea
\lab{definition-jkl}
j_a^b=\de_a^b+\frac{1}{2}\fb_a f^b, \qquad k_a=\frac{1}{2}\fb_a, \qquad l_a=\frac{1}{2}f_a+\frac{1}{8}|f|^2\fb_a.
\eea
We compute
\beaa
\g(\D_{e_a'}e_4, e_b') &=& \g\left(\D_{ j_a^ce_c+k_ae_4+l_ae_3}e_4, e_b'\right)= \g\left(\D_{ j_a^ce_c+l_ae_3}e_4, e_b'\right)\\
&=& \g\left(\D_{ j_a^ce_c+l_ae_3}e_4,  j_b^de_d+k_be_4+l_be_3\right)=\g\left(\D_{ j_a^ce_c+l_ae_3}e_4,  j_b^de_d+l_be_3\right)\\
&=& j_a^c j_b^d\g\left(\D_{ e_c}e_4, e_d\right)+l_a j_b^d\g\left(\D_{ e_3}e_4, e_d\right)+ j_a^cl_b\g\left(\D_{e_c}e_4,  e_3\right)+l_al_b\g\left(\D_{ e_3}e_4,  e_3\right)\\
&=& j_a^c j_b^d\chi_{cd}+2l_a j_b^d\eta_d+ 2j_a^cl_b\ze_c-4\omb l_al_b.
\eeaa
In view of \eqref{definition-jkl}
we infer
\beaa
&&\g(\D_{e_a'}e_4, e_b')\\
&& = \left(\de_a^c+\frac{1}{2}\fb_a f^c\right)\left(\de_b^d+\frac{1}{2}\fb_b f^d\right)\chi_{cd}+2\left(\frac{1}{2}f_a+\frac{1}{8}|f|^2\fb_a\right)\left(\de_b^d+\frac{1}{2}\fb_b f^d\right)\eta_d\\
&&+ 2\left(\de_a^c+\frac{1}{2}\fb_a f^c\right)\left(\frac{1}{2}f_b+\frac{1}{8}|f|^2\fb_b\right)\ze_c-4\omb\left(\frac{1}{2}f_a+\frac{1}{8}|f|^2\fb_a\right)\left(\frac{1}{2}f_b+\frac{1}{8}|f|^2\fb_b\right).
\eeaa
Taking the trace, this yields
\beaa
\de^{ab}\g(\D_{e_a'}e_4, e_b') &=& \left(1+f\c\fb+\frac{1}{2}|\fb|^2|f|^2\right)\trch+2f\c\chih\c\fb+|\fb|^2f\c\chih\c f\\
&&+f\c\eta+\frac{1}{2}(f\c\fb)(f\c\eta)+\frac{1}{4}|f|^2(\fb\c\eta)+\frac{1}{8}|f|^2|\fb|^2 f\c\eta\\
&&+ f\c\ze+\frac{1}{2}(f\c\fb)(f\c\ze)+\frac{1}{4}|f|^2(\fb\c\ze)+\frac{1}{8}(f\c\ze)|f|^2|\fb|^2\\
&&-\omb\left(|f|^2+\frac{1}{2}|f|^2(f\c\fb)+\frac{1}{16}|f|^4|\fb|^2\right)
\eeaa
which concludes the proof of the lemma.
\end{proof}

\begin{corollary}\lab{cor:e4derivativeofintegralonS}
We have for a scalar function $h$
\beaa
e_4\left(\int_{S(r,u)}h\right) &=& \int_{S(r,u)}\frac{e_4(\Si h)}{\Si}+O\left(\frac{\ep}{u^{1+\dec}}\right)h.
\eeaa 
\end{corollary}

\begin{proof}
Recall that we have
\beaa
\de^{ab}\g(\D_{e_a'}e_4, e_b') &=& \left(1+f\c\fb+\frac{1}{2}|\fb|^2|f|^2\right)\trch+2f\c\chih\c\fb+|\fb|^2f\c\chih\c f\\
&&+f\c\eta+\frac{1}{2}(f\c\fb)(f\c\eta)+\frac{1}{4}|f|^2(\fb\c\eta)+\frac{1}{8}|f|^2|\fb|^2 f\c\eta\\
&&+ f\c\ze+\frac{1}{2}(f\c\fb)(f\c\ze)+\frac{1}{4}|f|^2(\fb\c\ze)+\frac{1}{8}(f\c\ze)|f|^2|\fb|^2\\
&&-\omb\left(|f|^2+\frac{1}{2}|f|^2(f\c\fb)+\frac{1}{16}|f|^4|\fb|^2\right).
\eeaa
Hence  in view of the form of $f$ and $\fb$ in \eqref{eq:formulaforchangeofframecoeffandfbarforframetangenttosphereSur}, and in view of the control of the outgoing PG structure of $\Mext$ provided by  {\bf Ref 1}, including the  estimate \eqref{eq:improvedRef1-trXc} for $\trXc$, we infer
\beaa
\de^{ab}\g(\D_{e_a'}e_4, e_b') &=& \nu_0(r,\th)+O\left(\frac{\ep}{r^2u^{1+\dec}}\right),
\eeaa
where the function $\nu_0(r,\th)$ denotes the Kerr value. While it is in principle computable from the above formula, it is easier to compute it directly in Kerr. We have in Kerr in the $(\th, \vphi)$ coordinates 
\beaa
\pr_r\left(\int_S h\right)=\pr_r\left(\int_0^{2\pi}\int_0^\pi h\Si\sin\th d\th d\vphi\right)= \int_S\left(\pr_r(h)+\frac{\pr_r\Si}{\Si}h\right)
\eeaa
and hence
\beaa
\nu_0(r,\th) &=& \frac{\pr_r\Si}{\Si}.
\eeaa
We deduce in general
\beaa
e_4\left(\int_{S(r,u)}h\right) &=&  \int_{S(r,u)}\Big(e_4(h)+\de^{ab}\g(\D_{e_a'}e_4, e_b')h\Big)\\
&=& \int_{S(r,u)}\Big(e_4(h)+\nu_0(r,\th)h\Big)+O\left(\frac{\ep}{u^{1+\dec}}\right)h\\
&=& \int_{S(r,u)}\left(e_4(h)+\frac{\pr_r\Si}{\Si}h\right)+O\left(\frac{\ep}{u^{1+\dec}}\right)h\\
&=& \int_{S(r,u)}\frac{e_4(\Si h)}{\Si}+O\left(\frac{\ep}{u^{1+\dec}}\right)h,
\eeaa 
where we have used in particular the fact that $\Si=\Si(r,\th)$ and $e_4(r)=1$, $e_4(\th)=0$, so that $\pr_r\Si=e_4(\Si)$. This concludes the proof of the corollary.
\end{proof}

%%%%%%%%%%%%%%%%%%%%%%%%%%%%%%%%%%%%%%%%%%%%%

\subsection{Definition of $\ell=1$ modes on $S(u,r)$}

%%%%%%%%%%%%%%%%%%%%%%%%%%%%%%%%%%%%%%%%%%%%%

Recall the definition  of the basis of $\ell=1$ modes $J^{(p)}$, $p=0,+,-$ in $\Mext$, see section \ref{section:can.ell=1basis.Mext}. Relative to the PG coordinates $(\th, \vphi)$ of $\Mext$, we have
\beaa
J^{(0)}=\cos\th, \qquad J^{(+)}=\sin\th\cos\vphi, \qquad J^{(-)}=\sin\th\sin\vphi.
\eeaa

The $\ell=1$ modes of a scalar function on $S(u,r)$ are defined as follows.
\begin{definition}
Given a scalar function $f$ on a sphere $S=S(u, r)$,  we define   the $\ell=1$ modes of $f$  to be the triplet of numbers
\beaa
(f)_{\ell=1} &=& \left(\frac{1}{r^2}\int_{S}fJ^{(0)}, \frac{1}{r^2}\int_{S}fJ^{(+)}, \frac{1}{r^2}\int_{S}fJ^{(-)}\right).
\eeaa
\end{definition}

\begin{lemma}\lab{lemma:propertiesintegralofJpandJpJqonMext}
We have on $\Mext$
\beaa
\bsplit
\int_{S}J^{(p)} &= O\left(1+\ep ru^{-\frac{1}{2}-\dec}\right),\\
\int_{S}J^{(p)}J^{(q)} &= \frac{4\pi}{3}r^2\de_{pq}+O\left(1+\ep ru^{-\frac{1}{2}-\dec}\right).
\end{split}
\eeaa
\end{lemma}

\begin{proof}
Let $h$ a scalar function such that $e_4(h)=0$. Then, we have in view of Corollary \ref{cor:e4derivativeofintegralonS}
\beaa
e_4\left(r^{-2}\int_Sh\right) &=& \int_S\frac{e_4(r^{-2}\Si h)}{\Si} +O\left(\frac{\ep}{r^2u^{1+\dec}}\right)h\\
&=& \int_S\frac{e_4(r^{-2}\Si) }{\Si}h+O\left(\frac{\ep}{r^2u^{1+\dec}}\right)h.
\eeaa
Also, since $e_4(r)=1$ and $e_4(\th)=0$, we have
\beaa
e_4\left(\frac{\Si}{r^2}\right) &=& \pr_r\left(\frac{\Si}{r^2}\right)=\pr_r\sqrt{\left(1+\frac{a^2}{r^2}\right)^2+\frac{a^2(\sin\th)^2\De}{r^4}}=O(r^{-3})
\eeaa
and hence
\beaa
e_4\left(r^{-2}\int_Sh\right) &=& O\left(\frac{1}{r^3}+\frac{\ep}{r^2u^{1+\dec}}\right)h.
\eeaa
Applying this identity with $h=\Jp$ and $h=\Jp J^{(q)}$, we infer
\beaa
e_4\left(r^{-2}\int_S\Jp\right) &=& O\left(\frac{1}{r^3}+\frac{\ep}{r^2u^{1+\dec}}\right),\\
e_4\left(r^{-2}\int_S\Jp J^{(q)}\right) &=& O\left(\frac{1}{r^3}+\frac{\ep}{r^2u^{1+\dec}}\right).
\eeaa
Integrating from $\Si_*$, and together with the control on $\Si_*$ of Lemma \ref{lemma:statementeq:DeJp.Sigmastar:improvedd}, we infer 
\beaa
r^{-2}\int_S\Jp  &=& O\left(\frac{1}{r^2}+\frac{\ep}{ru^{\frac{1}{2}+\dec}}\right),\\
r^{-2}\int_S\Jp J^{(q)} &=& \frac{4\pi}{3}\de_{pq}+O\left(\frac{1}{r^2}+\frac{\ep}{ru^{\frac{1}{2}+\dec}}\right),
\eeaa
as stated.
\end{proof}

\begin{proposition}\lab{prop:controlofDDprimehotDDprimeJonMext}
Let $\nab'$ denote the covariant derivative on $S(u,r)$. Then, we have on $\Mext$, for $p=0,+,-$,
\beaa
|\DD'\hot\DD' \Jp|+|r^2\De' \Jp+2|\les \frac{\ep}{r^3u^{\frac{1}{2}+\dec}}+\frac{1}{r^4}.
\eeaa
\end{proposition}

\begin{proof}
See section \ref{sec:roofofprop:controlofDDprimehotDDprimeJonMext}.
\end{proof}

%%%%%%%%%%%%%%%%%%%%%%%%%%%%%%%%%%%%%

\subsection{Elliptic estimates on $S(u, r)$}

%%%%%%%%%%%%%%%%%%%%%%%%%%%%%%%%%%%%

We denote by $(e_1', e_2')$   the orthonormal  frame   of   $S=S(u, r) \subset\MMext$    defined in  section
\ref{subsection:orthonormalbasisS(u,r)}.  We first estimate the Gauss curvature of the spheres $S$.
\begin{lemma}
Let $K$ denote the Gauss curvature of the sphere $S=S(u, r) \subset\MMext$. Then, $K$ satisfies 
\bea
\sup_{\Mext}r^2\left|K-\frac{1}{r^2}\right| &\les& \frac{1}{r_0^2}+\frac{\ep}{r_0}.
\eea
\end{lemma}

\begin{remark}\lab{rmk:spheresofMextarealmostround:chap6}
In view of the above control of $K$, $S(u,r)$ is, for $r_0$ large enough, an almost round sphere in the sense of Definition \ref{def:almostroundsphereS:chap5}. 
\end{remark}

\begin{proof}
By Gauss equation, we have
\beaa
K &=& -\rho'-\frac{1}{4}\trch'\trchb'+\frac{1}{2}\chih'\c\chibh',
\eeaa
where $\trch'$, $\trchb'$, $\chih'$, $\chibh'$ and $\rho'$ correspond to the null frame $(e_4', e_3', e_1', e_2')$ adapted to $S(u,r)$, with $(e_1', e_2')$   the orthonormal  frame   of   $S$    defined in  section
\ref{subsection:orthonormalbasisS(u,r)}. The change of frame formulas of Proposition 
\ref{Proposition:transformationRicci}, the control of $(f, \fb)$ in \eqref{eq:formulaforchangeofframecoeffandfbarforframetangenttosphereSur}, and the computation in the Kerr case in Lemma \ref{Lemma:Transformation-principal-to-integrable-Kerr} imply
\beaa
\la^{-1}\trch' &=& \frac{2}{r}+O(r^{-3})+\dk^{\leq 1}\Ga_g,\\
\la\trchb' &=& -\frac{2\left(1-\frac{2m}{r}\right)}{r}+O(r^{-3})+\dk^{\leq 1}\Ga_g,\\
\la^{-1}\chih' &=& O(r^{-3})+\dk^{\leq 1}\Ga_g,\\
\la\chibh' &=& O(r^{-3})+\dk^{\leq 1}\Ga_b,\\
\rho' &=& -\frac{2m}{r^3}+O(r^{-5})+r^{-1}\Ga_g.
\eeaa
Plugging in Gauss equation, this yields 
\beaa
K &=& \frac{1}{r^2}+O(r^{-4})+r^{-1}\dk^{\leq 1}\Ga_g.
\eeaa
Together with the control of $\Ga_g$, and the fact that $r\geq r_0$ in $\Mext$, this concludes the proof of the lemma.
\end{proof}

We denote by $\ddd_1', \ddd_2',  \dds_1', \dds_2'$ the standard Hodge operators on $S$, see Definition \ref{def:HodgeoperatoronS:chap5}. Since the spheres $S$ are almost round, see Remark \ref{rmk:spheresofMextarealmostround:chap6}, and in view of the properties of the basis of $\ell=1$ modes $\Jp$, see Lemma \ref{lemma:propertiesintegralofJpandJpJqonMext} and Proposition \ref{prop:controlofDDprimehotDDprimeJonMext}, the Hodge elliptic estimates of Lemma \ref{prop:2D-Hodge1} and Lemma \ref{prop:2D-Hodge2} apply. We recall these results in the proposition below.

\begin{proposition} 
\label{prop:2D-hodge-Mext}
For any sphere $S=S(u, r) \subset\MMext$ we have for $k\leq k_{large}$:
\begin{enumerate}
\item  If $f$ is a 1-form
\bea
\|(\dkb' )^{\le k+1}f\|_{L^2(S)}\les r\|(\dkb' )^{\le k}\ddd_1'f\|_{L^2(S)}.
\eea

\item  If $f$ is a symmetric traceless 2-tensor
\bea
\|(\dkb' )^{\le k+1}f\|_{L^2(S)}\les r\|(\dkb' )^{\le k}\ddd_2'f\|_{L^2(S)}.
\eea

\item  If $(h, \dual h)$ is a pair of scalars 
\bea
\|(\dkb' )^{\le k+1}(h-\ov{h}, \dual h-\ov{\dual h})\|_{L^2(S)}\les r\|(\dkb' )^{\le k}\dds_1'(h, \dual h)\|_{L^2(S)}.
\eea

\item If $f$ is a 1-form 
\bea
  \| (\dkb' )^{\le k+1} f\|_{L^2  (S)}&\les&  r \|  (\dkb' )^{\le k } \dds'_2\, f\|_{L^2(S)}+r^2 \big| (\ddd'_1f)_{\ell=1}\big|.
  \eea
\end{enumerate}
 \end{proposition}

We have the following corollary of Proposition \ref{prop:2D-hodge-Mext}.
\begin{corollary}
\lab{Lemma:Hodge-estimateDDs_2} 
For any sphere $S=S(u, r) \subset\MMext$ we have for $k\leq k_{large}$:
\begin{enumerate}
\item If $U$ is an anti-selfdual  1-form
\bea
\|(\dkb' )^{\le k+1}  U \|_{L^2 (S)} &\les & r  \|(\dkb')^{\le k} ( \DD' \hot U)\|_{L^2 (S)}+ r^2  \Big| \big(\ov{ \DD'} \c  U\big)_{\ell=1}  \Big|.
\eea

\item If $U$ is an anti-selfdual  1-form
\bea
\|(\dkb' )^{\le k+1}  U \|_{L^2 (S)} &\les & r^2  \|(\dkb')^{\le k-1}\DD'(\ov{\DD'}\c U)\|_{L^2 (S)}.
\eea

\item If $U$ is an anti-selfdual  symmetric traceless 2-tensor
\bea
\|(\dkb' )^{\le k+1}  U \|_{L^2 (S)} &\les & r\|(\dkb')^{\le k}(\ov{\DD'}\c U)\|_{L^2 (S)}.
\eea
\end{enumerate}
\end{corollary}

\begin{proof}
We start with the first identity. Since $U$ is an anti-selfdual 1-form, $f=\Re(U)$ is a real 1-form and 
\beaa
U=f+i\dual f.
\eeaa
In particular, we have
\beaa
\DD' \hot U = 2(\nab'\hot f)+2i\dual(\nab'\hot f), \qquad \ov{\DD'}\c U=2\div'(f) +2i\curl'(f).
\eeaa
Thus, since 
\beaa
\dds_2'=-\frac{1}{2}\nab'\hot, \qquad \ddd_1'=(\div', \curl'),
\eeaa
the first identity follows immediately from the last estimate of Proposition \ref{prop:2D-hodge-Mext}.

Then, we consider the second identity. Since $U$ is an anti-selfdual 1-form, $f=\Re(U)$ is a real 1-form and 
\beaa
U=f+i\dual f.
\eeaa
In particular, we have
\beaa
\DD'(\ov{\DD'}\c U) &=& 2\Big(\nab'\div'(f)-\dual\nab'\curl'(f)\Big)+2i\dual\Big(\nab'\div'(f)-\dual\nab'\curl'(f)\Big),
\eeaa
and hence, we infer
\beaa
\DD'(\ov{\DD'}\c U) &=& -2\Big(\dds_1'\ddd_1'(f)+i\dual\big(\dds_1'\ddd_1'(f)\big)\Big).
\eeaa
Thus, the second identity follows immediately from the first and the third estimate of Proposition \ref{prop:2D-hodge-Mext}.

Finally, we consider the third identity. Since $U$ is an anti-selfdual symmetric traceless 2-tensor, $f=\Re(U)$ is a real symmetric traceless 2-tensor and 
\beaa
U=f+i\dual f.
\eeaa
In particular, we have
\beaa
\ov{\DD'}\c U &=& -2\Big(\div'(f)+i\div'(\dual f)\Big).
\eeaa
Thus, the third identity follows immediately from the first  estimate of Proposition \ref{prop:2D-hodge-Mext}. This concludes the proof of Corollary \ref{Lemma:Hodge-estimateDDs_2}. 
\end{proof}

%%%%%%%%%%%%%%%%%%%%%%%%%%%%%%%%%%%%%%%%%%%%%

\section{Renormalized quantities for outgoing PG structures}

%%%%%%%%%%%%%%%%%%%%%%%%%%%%%%%%%%%%%%%%%%%%%

%%%%%%%%%%%%%%%%%%%%%%%%%%%%%%%%%%%%%%%%%%%%%

\subsection{Renormalization of $\Hc$, $\widecheck{\cos\th}$ and $\ov{\DD}\c\Zc$}

%%%%%%%%%%%%%%%%%%%%%%%%%%%%%%%%%%%%%%%%%%%%%

We introduce the following renormalized quantities:
\bea
\bsplit
 [\Hc]_{ren}    :&=  \frac{1}{\ov{q}}\left(\ov{q}\Hc -q\Zc +\frac{1}{3}\left(-\ov{q}^2+|q|^2\right)B+\frac{a}{2}(q-\ov{q})\ov{\Jk}\c\Xh\right),\\
[\widecheck{\DD\cos\th}] _{ren}:&= \frac{1}{q}\left(q\widecheck{\DD\cos\th} +\frac{i}{2}|q|^2\ov{\Jk}\c\Xh\right),\\
[\Mc]_{ren} :&=\frac{1}{\ov{q} q^2}\Bigg[\ov{q}\,\ov{\DD}\c\left(q^2\Zc  +\left(-\frac{a}{2}q^2 -\frac{a}{2}|q|^2\right)\ov{\Jk}\c\Xh\right) +2\ov{q}^3\,\ov{\Pc} - 2aq^2\ov{\Jk}\c\Zc\\
&+\left(-\frac{1}{3}q^2\ov{q}^2 -\frac{1}{3}q\ov{q}^3+\frac{2}{3}\ov{q}^4\right)\ov{\DD}\c B +a\left( q^2\ov{q} +\frac{2}{3}q\ov{q}^2 - \frac{13}{6}\ov{q}^3\right)\ov{\Jk}\c B\\
&+a^2(q^2+|q|^2)\ov{\Jk}\c\Xh\c\ov{\Jk}\Bigg].
\end{split}
\eea

\begin{remark}
Note that, in the particular case  $a=0$,   $[\Mc]_{ren}$ is given by
\beaa
\ov{\DD}\c\Zc+2\ov{\Pc}=2(\div(\ze)+\rhoc)+2i(\curl(\ze)-\rhod),
\eeaa
where we have used the fact that $\widecheck{\ze}=\ze$ and $\widecheck{\rhod}=\rhod$ when $a=0$. In particular, the real part of $[\Mc]_{ren}$ coincides in that case, modulo a factor of $-2$, with the  linearized mass aspect function $\muc=-  \div\ze -\rhoc$. Thus, while there is no quantity denoted by $M$ in our work, the abuse of notation $\Mc$ should be thought as a complexified version of the linearized mass aspect function, and $[\Mc]_{ren}$ as its  corresponding renormalized version.
\end{remark}

\begin{proposition}
\lab{Prop: eqforrenormalized.qiantities}
We have
\bea
\bsplit
\nab_4\left(\ov{q}[\Hc]_{ren} \right)
&=  O(r^{-1})\trXc+O(1)\dkb^{\leq 1}A+r\Ga_b\c\Ga_g,\\
\nab_4\left(q [\widecheck{\DD\cos\th}]_{ren}  \right) &= O(1)\trXc+O(r)A+r\Ga_b\c\Ga_g,\\
\nab_4\left(\ov{q} q^2 [\Mc]_{ren}\right)&=   O(1)\dkb^{\leq 1}\trXc+O(r)\dkb^{\leq 2}A+r^2\dkb^{\leq 1}(\Ga_g\c\Ga_g)+r^3\Ga_b\c A.
\end{split}
\eea
\end{proposition}

\begin{proof}
See section \ref{appendix:ProofofProp{Prop: eqforrenormalized.qiantities}}.
\end{proof}

%%%%%%%%%%%%%%%%%%%%%%%%%%%%%%%%%%%%%%%%%%%%%

\subsection{Renormalization of the $\ell=1$ modes of  $\ov{\DD}\c B$}

%%%%%%%%%%%%%%%%%%%%%%%%%%%%%%%%%%%%%%%%%%%%%

\begin{definition}\lab{def:BrenovDren}
We introduce the following renormalized quantities
\bea
\bsplit
[B]_{ren}&:=B - \frac{3a}{2}\ov{\Pc}\Jk -\frac{a}{4}\ov{\Jk}\c A,\\
[\ov{\DD}\c]_{ren}&:= \ov{\DD}\c -\frac{a}{2}\ov{\Jk}\c\nab_4 -\frac{a}{2}\ov{\Jk}\c\nab_3.
\end{split}
\eea
\end{definition}

The goal of the section is to prove the following proposition. 
\begin{proposition}
\lab{Proposition:Identity.ell=1Div B}
Let $[B]_{ren}$ and $[\ov{\DD}\c]_{ren}$ given by Definition \ref{def:BrenovDren}. Then, the following identities hold true
\bea
\lab{eq:Proposition.Identity.ell=1Div B1}
\nn&&\nab_4 \left(\int_{S(r,u)} \frac{rJ^{(0)}}{\Si}   [\ov{\DD}\c]_{ren}    \left( r^4[B]_{ren} \right)    \right) \\
\nn&=&  O(1)\dk^{\leq 1}\Xh+O(r)\dk^{\leq 2}B+O(r^2)\dk^{\leq 1}\nab_3B+O(r)\dk^{\leq 2}\Pc\\
\nn&& +O(1)\dk^{\leq 1}\trXc
+O(r^{2})\dk^{\leq 1}\nab_3A+O(r)\dk^{\leq 2}A+r^4\dk^{\leq 1}\big(\Ga_g\c(B,A)\big)\\
&& +r^4\dk^{\leq 1}\big(\Ga_b\c\nab_3 A\big)
+r^2\dk^{\leq 2}\big(\Ga_g\c\Ga_g\big)
+O\left(\frac{r^2\ep}{u^{\frac 1 2+\dec}}\right)\dk^{\leq 1}\Big(A, B, r^{-1}\Pc\Big),
\eea
and
\bea
\lab{eq:Proposition.Identity.ell=1Div B2}
\nn&&\nab_4 \left(\int_{S(r,u)} \frac{rJ^{(\pm)}}{\Si}   [\ov{\DD}\c]_{ren}    \left( r^4[B]_{ren} \right)    \right)  \mp \frac{a}{r^2}  \int_{S(r,u)} \frac{rJ^{(\mp)}}{\Si}   [\ov{\DD}\c]_{ren}    \left( r^4[B]_{ren}    \right) \\
\nn&=&  O(1)\dk^{\leq 1}\Xh+O(r)\dk^{\leq 2}B+O(r^2)\dk^{\leq 1}\nab_3B+O(r)\dk^{\leq 2}\Pc\\
\nn&& +O(1)\dk^{\leq 1}\trXc
+O(r^{2})\dk^{\leq 1}\nab_3A+O(r)\dk^{\leq 2}A+r^4\dk^{\leq 1}\big(\Ga_g\c(B,A)\big)\\
&& +r^4\dk^{\leq 1}\big(\Ga_b\c\nab_3 A\big)
+r^2\dk^{\leq 2}\big(\Ga_g\c\Ga_g\big)
+O\left(\frac{r^2\ep}{u^{\frac 1 2 +\dec}}\right)\dk^{\leq 1}\Big(A, B, r^{-1}\Pc\Big).
\eea
\end{proposition}

The proof of Proposition \ref{Proposition:Identity.ell=1Div B} is based on  the following  identity.
\begin{lemma}
\lab{Le:mainpointwise-ell=1B}
The following identity holds true 
\beaa
\Big(\nab_4-a\Re(\Jk)^b\nab_b\Big)       \Big(    r [\ov{\DD}\c]_{ren}  \big(r^4 [B]_{ren} \big)\Big)
&=& \frac{r^5}{2}\ov{\DD}'\c\ov{\DD}'\c\Big( A  -a(\Jk\hot B)\Big) +\err,
\eeaa
where the  $\DD'$  is taken with respect  to the integral  frame  $(e_1', e_2')$   adapted to $S(u, r)$, see section \ref{subsection:orthonormalbasisS(u,r)}, and where the error term is given by 
\beaa
\err&=&    O(1)\dk^{\leq 1}\Xh+O(r)\dk^{\leq 2}B+O(r^2)\dk^{\leq 1}\nab_3B+O(r)\dk^{\leq 2}\Pc +O(1)\dk^{\leq 1}\trXc\\
&&+O(r^{2})\dk^{\leq 1}\nab_3A+O(r)\dk^{\leq 2}A+r^4\dk^{\leq 1}\big(\Ga_g\c(B,A)\big) +r^4\dk^{\leq 1}\big(\Ga_b\c\nab_3 A\big)\\
&&+r^2\dk^{\leq 2}\big(\Ga_g\c\Ga_g\big).
\eeaa
\end{lemma}

\begin{proof}
See section \ref{section:ProofLemma-mainpointwise-ell=1B}.
\end{proof}

We will also use the following lemma.
\begin{lemma}\lab{lemma:computationofRealpartofJknabJp}
We have
\beaa
\Re(\Jk)^b\nab_b(J^{(0)}) &=& r^{-1}\Ga_b,\\
\Re(\Jk)^b\nab_b(J^{(+)}) &=& -\frac{1}{r^2}J^{(-)}+O(r^{-4})+r^{-1}\Ga_b,\\
\Re(\Jk)^b\nab_b(J^{(-)}) &=& \frac{1}{r^2}J^{(+)}+O(r^{-4})+r^{-1}\Ga_b.
\eeaa
\end{lemma}

\begin{proof}
Since $\Im(\Jk)=\dual\Re(\Jk)$ and $J^{(0)}=\cos\th$, we have
\beaa
\Re(\Jk)^b\nab_b(J^{(0)}) &=& \Re(\Jk)\c\Re(\DD\cos\th)=\Re(\Jk)\c\Re(i\Jk+\widecheck{\DD\cos\th})= \Re(\Jk)\c\Im(\Jk)+r^{-1}\Ga_b\\
&=& \Re(\Jk)\c\dual\Re(\Jk)+r^{-1}\Ga_b= r^{-1}\Ga_b.
\eeaa
Also, recalling Definition \ref{def:renormalizationofallnonsmallquantitiesinPGstructurebyKerrvalue},
\beaa
\DD(J^{(\pm)}) &=& \Jk_\pm+\widecheck{\DD J^{(\pm)}}, \qquad \widecheck{\DD J^{(\pm)}}\in\Ga_b,
\eeaa
where the complex 1-forms $\Jk_+$ and $\Jk_-$ satisfy, see \eqref{eq:propertiesJk},
\beaa
\dual\Jk_\pm = -i\Jk_\pm,\qquad \Re(\Jk_+)\c\Re(\Jk)=-\frac{1}{|q|^2}J^{(-)}, \qquad \Re(\Jk_-)\c\Re(\Jk)=\frac{1}{|q|^2}J^{(+)}.
\eeaa
We infer
\beaa
\Re(\Jk)^b\nab_b(J^{(\pm)}) &=& \Re(\Jk)\c\Re(\DD(J^{(\pm)}))=\Re(\Jk)\c\Re(\Jk_\pm+\widecheck{\DD J^{(\pm)}})\\
&= &\Re(\Jk)\c\Re(\Jk_\pm)+r^{-1}\Ga_b,
\eeaa
and hence
\beaa
\Re(\Jk)^b\nab_b(J^{(+)}) &=& -\frac{1}{|q|^2}J^{(-)}+r^{-1}\Ga_b= -\frac{1}{r^2}J^{(-)}+O(r^{-4})+r^{-1}\Ga_b,\\
\Re(\Jk)^b\nab_b(J^{(-)}) &=& \frac{1}{|q|^2}J^{(+)}+r^{-1}\Ga_b=\frac{1}{r^2}J^{(+)}+O(r^{-4})+r^{-1}\Ga_b.
\eeaa
This concludes the proof of Lemma \ref{lemma:computationofRealpartofJknabJp}.
\end{proof}

We are now ready to prove Proposition \ref{Proposition:Identity.ell=1Div B}. 
 
\begin{proof}[Proof of Proposition \ref{Proposition:Identity.ell=1Div B}]
According to  Corollary \ref{cor:e4derivativeofintegralonS} we have for a scalar function $h$
\beaa
e_4\left(\int_{S(r,u)}h\right) &=& \int_{S(r,u)}\frac{e_4(\Si h)}{\Si}+O\left(\frac{\ep}{u^{1+\dec}}\right)h.
\eeaa 
We apply this identity with the choice $h = \frac{\Jp}{\Si}  r [\ov{\DD}\c]_{ren}  \big(r^4 [B]_{ren} \big)$. 
 We obtain, using also $e_4(J^{(p)})=0$,  $J^{(p)}=O(1)$, and $\Si\geq r^2$, 
 and recalling  the definition of  $[\ov{\DD}\c]_{ren}$,  $ [B]_{ren}$, see Definition \ref{def:BrenovDren}, for $p=0,+,-$,
\beaa
&&e_4\left(\int_{S(u,r)}    \frac{\Jp}{\Si}   \Big(    r [\ov{\DD}\c]_{ren}  \big(r^4 [B]_{ren} \big)\Big)\right)\\
&=& \int_{S(u,r)}       \frac{\Jp}{\Si} \nab_4  \Big(   r [\ov{\DD}\c]_{ren}  \big(r^4 [B]_{ren} \big)\Big)
+O\left(\frac{\ep}{  u^{1+\dec}}\right) \frac{\Jp}{\Si} r [\ov{\DD}\c]_{ren}  \big(r^4 [B]_{ren} \big)\\
&=&\int_{S(u,r)}       \frac{\Jp}{\Si} \nab_4  \Big(   r [\ov{\DD}\c]_{ren}  \big(r^4 [B]_{ren} \big)\Big) +O\left(\frac{r^2\ep}{u^{1+\dec}}\right)\dk^{\leq 1}\Big(B, r^{-1}\Pc, r^{-1} A\Big)\\
&=& \int_{S(u,r)}       \frac{\Jp}{\Si} \Big(  \nab_4   -a\Re(\Jk)^b\nab_b\Big)    \Big(   r [\ov{\DD}\c]_{ren}  \big(r^4 [B]_{ren} \big)\Big) \\
&&+\int_{S(u,r)}       \frac{\Jp}{\Si}  a\Re(\Jk)^b\nab_b   \Big(   r [\ov{\DD}\c]_{ren}  \big(r^4 [B]_{ren} \big)\Big) +O\left(\frac{r^2\ep}{u^{1+\dec}}\right)\dk^{\leq 1}\Big(B, r^{-1}\Pc, r^{-1} A\Big).
\eeaa

Making use  of Lemma \ref{Le:mainpointwise-ell=1B}, we infer, for $p=0,+,-$,
\beaa
&&e_4\left(\int_{S(u,r)}    \frac{\Jp}{\Si}   \Big(    r [\ov{\DD}\c]_{ren}  \big(r^4 [B]_{ren} \big)\Big)\right)\\
&=& \int_{S(r,u)}\frac{J^{(p)}}{\Si}\frac{r^5}{2}\ov{\DD}'\c\ov{\DD}'\c\Big( A  -a(\Jk\hot B)\Big) + a\int_{S(u,r)}       \frac{\Jp}{\Si}  \Re(\Jk)^b\nab_b   \Big(   r [\ov{\DD}\c]_{ren}  \big(r^4 [B]_{ren} \big)\Big)\\
&& +\int_{S(r,u)}\frac{J^{(p)}}{\Si}\err +O\left(\frac{r^2\ep}{u^{1+\dec}}\right)\dk^{\leq 1}\Big(B, r^{-1}\Pc, r^{-1} A\Big).
\eeaa
Since $\Si\geq r^2$ and $\Jp=O(1)$, we obtain,  with error terms $\err$ of the same form as the ones in Lemma \ref{Le:mainpointwise-ell=1B},
\beaa
&&e_4\left(\int_{S(u,r)}    \frac{\Jp}{\Si}   \Big(    r [\ov{\DD}\c]_{ren}  \big(r^4 [B]_{ren} \big)\Big)\right)\\
 &=& \int_{S(r,u)}\  J^{(p)} \frac{r^3}{2} \ov{\DD}'\c\ov{\DD}'\c\Big( A  -a(\Jk\hot B)\Big) + a\int_{S(u,r)}       \frac{\Jp}{\Si}  \Re(\Jk)^b\nab_b   \Big(   r [\ov{\DD}\c]_{ren}  \big(r^4 [B]_{ren} \big)\Big) \\
 &&+O\left(\frac{r^2\ep}{u^{1+\dec}}\right)\dk^{\leq 1}\Big(B, r^{-1}\Pc, r^{-1} A\Big) +\err.
\eeaa
Integrating by parts twice, we have 
\beaa
\int_{S(r,u)}J^{(p)}\frac{r^3}{2}\ov{\DD}'\c\ov{\DD}'\c\Big( A  -a(\Jk\hot B)\Big) &=& \frac{1}{2}\int_{S(r,u)}\frac{r^3}{2}\Big( A  -a(\Jk\hot B)\Big)\ov{\DD}'\hot\ov{\DD}' J^{(p)}\\
&=& O(r^{5})\Big( A  -a(\Jk\hot B)\Big)\ov{\DD'\hot\DD' J^{(p)}},
\eeaa
where we used the fact that $\Jp$ is real valued. Also, according to Proposition \ref{prop:controlofDDprimehotDDprimeJonMext}, we have, for $p=0,+,-$, 
\beaa
|\DD'\hot\DD'\Jp|\les \frac{\ep}{r^3u^{\frac{1}{2}+\dec}}+\frac{1}{r^4}.
\eeaa
We infer
\beaa
\int_{S(r,u)}J^{(p)}\frac{r^3}{2}\ov{\DD}'\c\ov{\DD}'\c\Big( A  -a(\Jk\hot B)\Big) &=& \left(O\left(\frac{\ep r^2}{u^{\frac{1}{2}+\dec}}\right)+O(r)\right)\Big( A  -a(\Jk\hot B)\Big).
\eeaa
Therefore, with error term $\err$ of the same form as the ones in Lemma \ref{Le:mainpointwise-ell=1B},
\beaa
&&e_4\left(\int_{S(u,r)}    \frac{\Jp}{\Si}   \Big(    r [\ov{\DD}\c]_{ren}  \big(r^4 [B]_{ren} \big)\Big)\right)\\
&=& a\int_{S(u,r)}       \frac{\Jp}{\Si}  \Re(\Jk)^b\nab_b   \Big(   r [\ov{\DD}\c]_{ren}  \big(r^4 [B]_{ren} \big)\Big)+O\left(\frac{r^2\ep}{u^{\frac{1}{2}+\dec}}\right)\dk^{\leq 1}\Big(A, B, r^{-1}\Pc\Big)+\err.
\eeaa
Next, in view of Corollary  \ref{cor:formuallikningnabprimeandnabintransfoformula:bis}, we have the rough decomposition 
\beaa
\nab' &=& (1+O(r^{-2}) )\nab+O(r^{-1})\nab_4+O(r^{-1})\nab_3+O(r^{-3})+r^{-1}\Ga_b\dk^{\leq 1}
\eeaa
or
\beaa
\nab &=& \nab'+O(r^{-1})\nab_3+O(r^{-2})\dk^{\leq 1}+r^{-1}\Ga_b\dk^{\leq 1}.
\eeaa
Hence, with error terms $\err$ of the same form as the ones in Lemma \ref{Le:mainpointwise-ell=1B},
\beaa
&&e_4\left(\int_{S(u,r)}    \frac{\Jp}{\Si}   \Big(    r [\ov{\DD}\c]_{ren}  \big(r^4 [B]_{ren} \big)\Big)\right)\\
&=& a\int_{S(u,r)}       \frac{\Jp}{\Si}  \Re(\Jk)^b\nab'_b   \Big(   r [\ov{\DD}\c]_{ren}  \big(r^4 [B]_{ren} \big)\Big)+
\err+O\left(\frac{r^2\ep}{u^{\frac{1}{2}+\dec}}\right)\dk^{\leq 1}\Big(A, B, r^{-1}\Pc\Big).
\eeaa
Integrating by parts
\beaa
\int_{S(u,r)}       \frac{\Jp}{\Si}  \Re(\Jk)^b\nab'_b   \Big(   r [\ov{\DD}\c]_{ren}  \big(r^4 [B]_{ren} \big)\Big) = - \int_{S(u,r)} 
 \div' \left(\frac{\Jp}{\Si}  \Re(\Jk) \right)   r [\ov{\DD}\c]_{ren}  \big(r^4 [B]_{ren}\big).
\eeaa
Using again the decomposition
\beaa
\nab' &=& \nab+O(r^{-1})\nab_3+O(r^{-2})\dk^{\leq 1}+r^{-1}\Ga_b\dk^{\leq 1},
\eeaa
recalling that $ \nab_3\Jk= \widecheck{\nab_3\Jk}+ \frac{\De q}{|q|^4}\Jk = r^{-1} \Ga_b + O(r^{-1})\Jk$, as well as  
 $ \widecheck{e_3(J^{(+)})}, \, \widecheck{e_3(J^{(-)})}\in  \Ga_b$, we deduce,
\beaa
&&e_4\left(\int_{S(u,r)}    \frac{\Jp}{\Si}   \Big(    r [\ov{\DD}\c]_{ren}  \big(r^4 [B]_{ren} \big)\Big)\right)\\
&=& - a\int_{S(u,r)}    \div\left(   \frac{\Jp}{\Si}  \Re(\Jk)\right)    \Big(   r [\ov{\DD}\c]_{ren}  \big(r^4 [B]_{ren} \big)\Big)+\err +O\left(\frac{r^2\ep}{u^{\frac{1}{2}+\dec}}\right)\dk^{\leq 1}\Big(A, B, r^{-1}\Pc\Big).
\eeaa
Since  $\div\left(\Re(\Jk)\right) = r^{-1}\Ga_b$, we infer, with error term $\err$ of the same form as the ones in Lemma \ref{Le:mainpointwise-ell=1B},
\beaa
&&e_4\left(\int_{S(u,r)}    \frac{\Jp}{\Si}   \Big(    r [\ov{\DD}\c]_{ren}  \big(r^4 [B]_{ren} \big)\Big)\right)\\
&=& - a\int_{S(u,r)}    \frac{1}{r^2}  \Re(\Jk)^b\nab_b(J^{(p)}) \Big(   r [\ov{\DD}\c]_{ren}  \big(r^4 [B]_{ren} \big)\Big)+\err +O\left(\frac{r^2\ep}{u^{\frac{1}{2}+\dec}}\right)\dk^{\leq 1}\Big(A, B, r^{-1}\Pc\Big).
\eeaa
Now, recall that we have according to Lemma \ref{lemma:computationofRealpartofJknabJp}
\beaa
\Re(\Jk)^b\nab_b(J^{(0)}) &=& r^{-1}\Ga_b,\\
\Re(\Jk)^b\nab_b(J^{(+)}) &=& -\frac{1}{r^2}J^{(-)}+O(r^{-4})+r^{-1}\Ga_b,\\
\Re(\Jk)^b\nab_b(J^{(-)}) &=& \frac{1}{r^2}J^{(+)}+O(r^{-4})+r^{-1}\Ga_b.
\eeaa
Since $\Si=r^2+O(1)$, we obtain, with error term $\err$ of the same form as the ones in Lemma \ref{Le:mainpointwise-ell=1B},
\beaa
e_4\left(\int_{S(u,r)}    \frac{J^{(0)}}{\Si}   \Big(    r [\ov{\DD}\c]_{ren}  \big(r^4 [B]_{ren} \big)\Big)\right) &=& \err +O\left(\frac{r^2\ep}{u^{\frac{1}{2}+\dec}}\right)\dk^{\leq 1}\Big(A, B, r^{-1}\Pc\Big)
\eeaa
and 
\beaa
e_4\left(\int_{S(u,r)}    \frac{J^{(\pm)}}{\Si}   \Big(    r [\ov{\DD}\c]_{ren}  \big(r^4 [B]_{ren} \big)\Big)\right) &=& \pm \frac{a}{r^2}\int_{S(u,r)} \frac{rJ^{(\mp)}}{\Si} [\ov{\DD}\c]_{ren}  \big(r^4 [B]_{ren} \big)\\
&&+\err +O\left(\frac{r^2\ep}{u^{\frac{1}{2}+\dec}}\right)\dk^{\leq 1}\Big(A, B, r^{-1}\Pc\Big),
\eeaa
as desired. This concludes the proof of Proposition \ref{Proposition:Identity.ell=1Div B}. 
\end{proof}

%%%%%%%%%%%%%%%%%%%%%%%%%%%%%%%%%%%%%%%%%%%%%%%%%%

\section{Main Estimates in $\Mext$}

%%%%%%%%%%%%%%%%%%%%%%%%%%%%%%%%%%%%%%%%%%%%%%%%%%

%%%%%%%%%%%%%%%%%%%%%%%%%%

\subsection{Transport lemmas}

%%%%%%%%%%%%%%%%%%%%%%%%%%

The transport lemmas derived in this section will be used repeatedly in the proof of Theorem M4. Recall the norms of  Definition \ref{definition:norms-ProofThm.M4}, i.e. 
    \beaa
  \| f\|_{\infty} (u,r):=\| f\|_{L^\infty\big(S(u,r)\big)}, \qquad \|f\|_{\infty,k}(u, r):= \sum_{i=0}^k \|\dk^i f\|_{\infty }(u, r).
  \eeaa
Recall also that  the weighted derivatives $\dkb=(r\nab) $ and $\dk=( r\nab, r\nab_4, \nab_3)$  are defined with respect to the outgoing  PG frame of $\Mext$.

\begin{lemma}
\lab{Lemma:TransportMext}
Let $U$ and $F$  anti-selfdual $k$-tensors. Assume  that  $U$ verifies one of the following equations, for a real constant $c$,
\bea
\lab{eq:maintransportMext}
\nab_4 U+\frac{c}{q} U=F
\eea
or
\bea
\lab{eq:maintransportMext2}
\nab_4 U+ \Re\left(\frac{c}{q}\right)  U=F.
\eea
In both cases we derive,   for  any $r_0\leq r\le r_* $  at fixed $u$, with  $1\leq u\leq u_*$, 
 \bea
 \lab{eq:Lemma:TransportMex1t}
\bsplit
    r^{c}  \|   U   \|_\infty (u, r)&\les    r_*  ^{c} \|  U \|_\infty (u,  r_*  )+ \int_r^{r_*}      \la ^{c} \| F\|_\infty(u, \la) d\la.
    \end{split}
\eea
\end{lemma} 

\begin{proof}
Assume first that $U$ satisfies \eqref{eq:maintransportMext}. Since $e_4(q)=1$, we can  rewrite the equation in the form
\beaa
\nab_4(q^{c}  U) = q^{c} F.
\eeaa
The desired inequality follows then immediately  by integration  in $r$, using the fact that $e_4(r)=1$, $e_4(u)=0$  and $r\leq |q|\leq 2r$.

Next, assume that $U$ satisfies \eqref{eq:maintransportMext2}. Since $e_4(q)=1$ and $e_4(\ov{q})=1$, we can  rewrite the equation in the form
\beaa
\nab_4\big( |q|^c  U\big)&=& |q|^c  \left( - \Re\Big(\frac{c}{q}\Big)  U+F\right)+ c|q|^{c-1} \nab_4 (|q|)  U\\
&=&  |q|^c  \left( - \Re\Big(\frac{c}{q}\Big)  U+F\right)+ c|q|^{c-1}\frac{1}{2|q|} \big(q+\ov{q}\big)  U\\
&=&  |q|^c  \left( - \Re\Big(\frac{c}{q}\Big)  U+F\right)+ c |q|^c  \frac{1}{2|q|^2 } \big(q+\ov{q}\big)  U\\
&=&   |q|^c  \left( - \Re\left(\frac{c}{q}\right)  U+F\right)+ c |q|^c  \Re\Big(\frac{1}{q} \Big) U =   |q|^c F,
\eeaa
i.e.
\beaa
\nab_4\big( |q|^c  U\big)&=&   |q|^c F.
\eeaa
We can then proceed in the same manner as  in the first case.
\end{proof}

\begin{proposition}
\lab{Prop:transportrp-f-Decay-knorms}
Solutions  $U$ of the equations \eqref{eq:maintransportMext}  or  \eqref{eq:maintransportMext2}  verify the following estimate, for all  $k\le k_{large}$,  $r_0\leq r\le r_*$, and $1\leq u\leq u_*$, 
 \bea\lab{eq:integration-knorms}
   r^{c}  \|U\|_{\infty,k} (u, r) \les    r_*  ^{c} \|  U \|_{\infty,k} (u,  r_*  )+ \int_r^{r_*}      \la ^{c}\|F\|_{\infty,k}(u, \la) d\la.
\eea
\end{proposition}

\begin{proof}
The proof being similar for \eqref{eq:maintransportMext} and  \eqref{eq:maintransportMext2}, we only treat the case where $U$ verifies \eqref{eq:maintransportMext}, i.e. 
\beaa
\nab_4 U+\frac{c}{q} U=F.
\eeaa

Recall from Corollary \ref{cor:commutation-complexM6} that we have
 \beaa
  \, [\nab_4, q\DD] U &=&  O(r^{-2})U+ \Ga_g \c \dkb^{\leq 1} U.
 \eeaa        
Also, recall from Lemma \ref{Lemma:commLieb_Tnab-simb} that we have
\beaa
 [\nab_4, \Lieb_\T]U &=& r^{-1}\Ga_b\c\dkb U+r^{-1}\dk\Ga_b\c U.
\eeaa
We infer, using also  $\T(q)\in r\Ga_b$ in view of Lemma \ref{Lemma:$T(rthuJ}, and $\DD(q)=O(r^{-1})+\Ga_b$, 
\bea\lab{eq:intermediarycommutedeqqddLieTfortranportlemma}
\nn\nab_4((q\DD)^j\Lieb_\T^lU)+\frac{c}{q}(q\DD)^j\Lieb_\T^lU &=& (q\DD)^j\Lieb_\T^lF+O(r^{-2})(q\DD,\Lieb_\T)^{\leq j+l}U\\
&&+r^{-1}\sum_{p=0}^{j+l}\dk^{\leq p}(\Ga_b)(q\DD,\Lieb_\T)^{j+l-p}U.
\eea

Applying Lemma \ref{Lemma:TransportMext} to \eqref{eq:intermediarycommutedeqqddLieTfortranportlemma}, and using the control of $\Ga_b$, we obtain, for $k\leq k_{large}$,  
 \beaa
 &&  \sum_{j+l\leq k} r^{c}  \|(q\DD)^j\Lieb_\T^lU\|_\infty (u, r)\\
 &\les&    r_*  ^{c} \|  U \|_{\infty,k} (u,  r_*  )+ \int_r^{r_*}      \la ^{c}\left( \|F\|_{\infty,k}(u, \la)+O(\la^{-2})\sum_{j+l\leq k} r^{c}  \|(q\DD)^j\Lieb_\T^lU\|_\infty (u, \la)\right) d\la.
\eeaa
Together with Gronwall lemma, we infer, for $k\leq k_{large}$, 
 \bea\lab{eq:intermediarycommutedeqqddLieTfortranportlemma:1}
   \sum_{j+l\leq k} r^{c}  \|(q\DD)^j\Lieb_\T^lU\|_\infty (u, r) &\les&    r_*  ^{c} \|  U \|_{\infty,k} (u,  r_*  )+ \int_r^{r_*}      \la ^{c}\|F\|_{\infty,k}(u, \la) d\la.
\eea

Next, multiplying \eqref{eq:intermediarycommutedeqqddLieTfortranportlemma} with $r$ and differentiating it w.r.t. $(r\nab_4)^p$, and using also the control of $\Ga_b$, we have, for $p\geq 1$ and $j+l+p\leq k_{large}$, 
\beaa
|(r\nab_4)^p(q\DD)^j\Lieb_\T^lU)| &\les& r|\dk^{\leq k-1}F|+ |(r\nab_4)^{p-1}(q\DD, \Lieb_\T)^{j+l}U)|.
\eeaa
Together with \eqref{eq:intermediarycommutedeqqddLieTfortranportlemma:1}, we deduce by iteration, for $k\leq k_{large}$, 
 \beaa
   \sum_{j+l+p\leq k} r^{c}  \|(r\nab_4)^p(q\DD)^j\Lieb_\T^lU\|_\infty (u, r) &\les&    r_*  ^{c} \|  U \|_{\infty,k} (u,  r_*  )+ r^{c+1}\|F\|_{\infty,k-1}(u, r)\\
   &&+ \int_r^{r_*}      \la ^{c}\|F\|_{\infty,k}(u, \la) d\la.
\eeaa
Using
\beaa
r^{c+1}\|F\|_{\infty,k-1}(u, r) &\les & \int_r^{r_*}      \la ^{c}\|(\la\nab_4)^{\leq 1}F\|_{\infty,k-1}(u, \la) d\la\\
&\les & \int_r^{r_*}      \la ^{c}\|F\|_{\infty,k}(u, \la) d\la,
\eeaa
we obtain, for $k\leq k_{large}$, 
 \beaa
   \sum_{j+l+p\leq k} r^{c}  \|(r\nab_4)^p(q\DD)^j\Lieb_\T^lU\|_\infty (u, r)  &\les&    r_*  ^{c} \|  U \|_{\infty,k} (u,  r_*  )+ \int_r^{r_*}      \la ^{c}\|F\|_{\infty,k}(u, \la) d\la.
\eeaa

Since $q=r+O(1)$, and since $e_3$ is spanned by $(\T, e_4, e_1, e_2)$ in view of the definition of $\T$, see \eqref{definition:T-Mext}, we infer,  for $k\leq k_{large}$, 
 \beaa
   r^{c}  \|U\|_{\infty,k} (u, r) &\les&    r_*  ^{c} \|  U \|_{\infty,k} (u,  r_*  )+ \int_r^{r_*}      \la ^{c}\|F\|_{\infty,k}(u, \la) d\la
\eeaa
as desired. This concludes the proof of Proposition \ref{Prop:transportrp-f-Decay-knorms}.
\end{proof}

%%%%%%%%%%%%%%%%%%%%%%%%%%%%%%%%%%%%%%%

\subsection{Estimates for the outgoing PG structure of $\Mext$ on $\Si_*$} 

%%%%%%%%%%%%%%%%%%%%%%%%%%%%%%%%%%%%%%%

In this section, we recall the main estimates derived in  section \ref{sec:decayestimatesPGframeonSigmastar} on $\Si_*$   with respect to  the outgoing  PG structure  of $\Mext$. More precisely, we restate below Proposition \ref{prop:improvedesitmatesfortemporalframeofMextonSigmastar}. Note that in the statement of that proposition, the PG frame is denoted by prime, while the  frame adapted to $\Si_*$, which is used in the proof,  is  unprimed. Since, in this chapter,  we only deal with the outgoing PG frame of $\Mext$, we therefore drop the primes in the statement of Proposition \ref{prop:improvedesitmatesfortemporalframeofMextonSigmastar} which thus takes the following form.

\begin{proposition}\lab{prop:improvedesitmatesfortemporalframeofMextonSigmastar-new}
We have on $\Si_*$, for\footnote{Recall from \eqref{eq:valueofkstarinchapter6forproofThmM4} that $k_*=k_{small}+60$ in this chapter.} $k\leq k_*$, 
\beaa 
\sup_{\Si_*}\Big(ru^{1+\dec}|\dk^k\Ga_b|+r^2u^{\frac{1}{2}+\dec}|\dk^k\Ga_g|+r^2u^{1+\dec}|\dk^{k-1}\nab_3\Ga_g|\Big) &\les& \ep_0,
\eeaa
\beaa
\sup_{\Si_*}\left(r^2u^{1+\dec}|\dk^{k}\trXc|+r^3u^{1+\dec}\left|\dk^{k}\left(\ov{\DD}\c\Zc+2\ov{\Pc}\right)\right|+r^4u^{\frac{1}{2}+\dec}|\dk^{k-1}\nab_3B|\right) &\les& \ep_0,
\eeaa 
and
\beaa
\sup_{\Si_*}r^5u^{1+\dec}\left(\Big |  \left( [\ov{\DD}\c]_{ren} [B]_{ren} \right)_{\ell=1}\Big| +\Big|\left[\ov{\DD}\c\Lieb_{\T} B \right]_{\ell=1}\Big|\right) &\les& \ep_0.
\eeaa
\end{proposition}

%%%%%%%%%%%%%%%%%%%%%%%%%%%%
 
\subsection{Strategy of the proof of Theorem M4}
\lab{section:MainResultsandStrategy}

%%%%%%%%%%%%%%%%%%%%%%%%%%%%

Our  goal in this chapter  is to extend  the results of Proposition \ref{prop:improvedesitmatesfortemporalframeofMextonSigmastar-new} to $\Mext$, i.e. to prove the following proposition which implies Theorem M4. 

\begin{proposition}\lab{prop:improvedesitmatesfortemporalframeofMexton-new}
We have on $\Mext$, for\footnote{Recall from \eqref{eq:valueofkstarinchapter6forproofThmM4} that $k_*=k_{small}+60$ in this chapter.} $k\leq k_{*}- 8$,
\beaa 
\sup_{\Mext} \Big(ru^{1+\dec}|\dk^k\Ga_b|+\big(r^2u^{\frac{1}{2}+\dec}+ru^{1+\dec}\big)|\dk^k\Ga_g|+r^2u^{1+\dec}|\dk^{k-1}\nab_3\Ga_g|\Big) &\les& \ep_0,
\eeaa
and
\beaa
\sup_{\Mext}r^2u^{1+\dec}|\dk^{k}\trXc| +\sup_{\Mext}r^4u^{\frac{1}{2}+\dec}|\dk^{k-1}\nab_3B| &\les& \ep_0.
\eeaa 
\end{proposition}

We now describe the strategy of the proof of Proposition \ref{prop:improvedesitmatesfortemporalframeofMexton-new}. To start with, we need to  distinguish between  the estimates  for the $\Ga_g$ quantities  which involve  $O(r^{-2}u^{-1/2-\dec})$  decay,  and those  which   involve $O(r^{-1}u^{-1-\dec})$ decay. The  first are  
 relatively easy to derive using   our main linearized  equations,  see Lemma \ref{Lemma:linearized-nullstr},     the corresponding estimates on the last slice, Proposition \ref{Prop:transportrp-f-Decay-knorms} and the assumptions {\bf Ref 1} and {\bf Ref 2}.   In what follows, we  describe the main steps in deriving the  much more subtle  $O(r^{-1}u^{-1-\dec})$
 estimates for the $\Ga_g$ quantities:
 
 \begin{enumerate}
\item  First derive an estimate for $ \widecheck{\tr X}$ using  the Raychadhouri equation it verifies, Proposition \ref{Prop:transportrp-f-Decay-knorms}, and its estimate on the last slice.  The resulting estimates improves the  stronger  assumption made in {\bf Ref 1}, see \eqref{eq:improvedRef1-trXc}, i.e. we obtain \eqref{eq:improvedRef1-trXc} with $\ep$ being replaced by $\ep_0$.

 \item  We observe that we  are not able to estimate directly the other  primary quantities  $\Xh, B, \Zc, \Pc$. Consider for example the equation verified by  $\Xh$
 \beaa
\nab_4\Xh+\Re(\tr X)\Xh &=& -A.
\eeaa
This works well, with the help of   Proposition \ref{Prop:transportrp-f-Decay-knorms},  to derive  an estimate of the form  $ O\big(\ep_0 r^{-2}  u^{-1/2-\dec}\big) $  but fails to provide an
$ O\big(\ep_0 r^{-1}  u^{-1-\dec}\big) $  estimate. Indeed according to    {\bf Ref 2},  we  only have  $A= O\big(\ep_0 r^{-2-\de'}  u^{-1-\dec} \big)$, for a small constant $\de'=\frac 1 2 \big(\de_{extra}-\dec\big)>0$. On the other hand  we can commute with
 $\Lieb_\T$ and derive an estimate of  the form  $\Lieb_\T \Xh = O\big(\ep_0 r^{-2}  u^{-1-\dec}\big) $ by making use of the fact that, according to {\bf Ref 2}, 
  $\Lieb_\T A= O\big(\ep_0 r^{-3-\de'}  u^{-1-\dec} \big)$.
  
\item  We  encounter  a  similar issue  with estimates for $B$.  To start with we can not use its natural transport equation
  \beaa
\nab_4B+2\ov{\tr X} B &=& \frac{1}{2}\ov{\DD}\c A+\frac{1}{2}(\ov{2Z+\Hb})\c A.
\eeaa
   for the simple fact that it is seriously overshooting in $r$. We   
   look  instead  at   another Bianchi equation,
    \beaa
\DD\hot B +(Z+4H)\hot B &=& \nab_3A+\left(\frac{1}{2}\tr\Xb-4\omb\right)A+3\ov{P}\Xh.
\eeaa
\item  We need first  to commute   with $\Lieb_\T$  to take advantage  of  the previously  derived  information for $\Lieb_\T\Xh$. This provides 
an estimate for  $\DD\hot \Lieb_\T B$ from which we would  have to recover  $\Lieb_\T B$. 

\item  The problem however is that $\DD\hot$ is not an operator on $S(t, u)$.  It is for this reason that we have to appeal to Proposition  \ref{cor:formuallikningnabprimeandnabintransfoformula} and Corollary 
\ref{cor:formuallikningnabprimeandnabintransfoformula:bis} to derive instead an estimate  for  $\DD'\hot \Lieb_\T B$. 

\item In the process however we generate another 
 $\Lieb_\T $ derivative for $B$, i.e. we first need  an estimates for $\Lieb_\T^2 B$. Fortunately this can be derived\footnote{Note that the corresponding  transport equation for $\Lieb_\T B$ is still overshooting in $r$, but the one for   $\Lieb_\T^2  B$   is not.}   by commuting  the overshooting equation
  $\nab_4B+2\ov{\tr X} B= \frac{1}{2}\ov{\DD}\c A+\ldots$  twice with $\Lieb_\T$ and thus derive   a transport equation for  $\Lieb_\T^2$ which is no longer overshooting. Thus we can  first estimate  $\Lieb_\T^2 B $  from which we derive also an estimate for  $\DD'\hot \Lieb_\T B$. Note that this step requires the full force of the  assumptions {\bf Ref 2} for $A$.

  \item It remains to  derive estimates for  $ \Lieb_\T B$ from  the one of   $\DD'\hot \Lieb_\T B$.  In view of Corollary \ref{Lemma:Hodge-estimateDDs_2}, we need first to estimate
  $  \big(\ov{ \DD'} \c  \Lieb_\T B\big)_{\ell=1} $.   Starting with the equation for $\nab_4 \Lieb_\T B +2\ov{\tr X}\Lieb_\T  B= \frac{1}{2}\ov{\DD}\c\Lieb_\T  A+\ldots$ we  first   derive a transport equation\footnote{Note that this is simpler to derive than the renormalized version of  the quantity  $  \frac{1}{\Si}\ov{q}\,\ov{\DD}\c(\ov{q}^4 B )J^{(p)}$  in   Proposition \ref{Proposition:Identity.ell=1Div B}.}  
   for   the quantity  $\frac{1}{\Si}\ov{q}\,\ov{\DD}\c(\ov{q}^4\Lieb_\T B )J^{(p)}$ which can therefore be estimated, using our information on  $\Si_*$.  
  Using the transformation  formulas   between the prime and unprimed frames we  then estimate  $  \big(\ov{ \DD'} \c  \Lieb_\T B\big)_{\ell=1} $ and thus, by  Corollary \ref{Lemma:Hodge-estimateDDs_2},  the derivatives   $(\dkb' )^k\Lieb_\T B$. Transforming back to the outgoing PG frame of $\Mext$, we derive the desired estimates for $\Lieb_\T B$.

  \item The estimates for $\Lieb_\T\Zc$ and $\Lieb_\T\Pc$ are then similar to the ones for $\Lieb_\T\Xh$. Our strategy  therefore  is to estimate  first $\Lieb_\T\Xh$, $\Lieb_\T^2B$, $\Lieb_\T B$, $\Lieb_\T\Zc$ and $\Lieb_\T\Pc$,  and then  derive   estimates for the primary quantities using elliptic theory on the spheres  $S=S(r, u)$.

\item  The estimates  for the quantities $\Xh, B, \Zc, \Pc$, though more subtle, follow a similar pattern. The  main difference is that they require the transport equations for  the renormalized quantities   $[H]_{ren}$, $[\widecheck{\DD\cos\th}] _{ren}$, $[\Mc]_{ren}$, as well as  
$\int_{S(r,u)} \frac{rJ^{(0)}}{\Si}   [\ov{\DD}\c]_{ren}     r^4[B]_{ren}$, derived in Proposition \ref{Prop: eqforrenormalized.qiantities} and Proposition \ref{Proposition:Identity.ell=1Div B}.
 
 \item   We first  derive conditional   decay  estimates, see Proposition \ref{prop:firstcontrolofBPcXhZcHcDDcosth},   in which the term  $(\ov{\DD}\c B)_{\ell=1}$ appear on the right hand side.
 The estimate are then made unconditional in Proposition \ref{Prop:estimatesforBPXhZcHc}. It is important to note that the estimates for  $\Xh, B, \Zc, \Pc$  in
 Proposition \ref{Prop:estimatesforBPXhZcHc} not only improve the corresponding assumptions in {\bf Ref 1}  by replacing $\ep$ with $\ep_0$; they are also better  in powers of $r$, that is they
  gain $r^{-\de'}$ for $\de'>0$. These improvements are needed later  to  derive  the correct estimates for the remaining  quantities in $\Ga_b$.

 \item  All remaining estimates, i.e. the $O(r^{-2}u^{-1/2-\dec})$  decay estimates for $\Ga_g$,  and the decay estimates for        $ \Ga_b$,    are derived in Propositions \ref{Proposition:estimatesforBPcXhZc-halfdecayinu} and Proposition \ref{Proposition:remaining.estimatesMext}.
 \end{enumerate}

%%%%%%%%%%%%%%%%%%%%%%%%%%%%%%%%%%%%%%%%%%%%%%%%%%

\subsection{Estimates for $\trXc$}

%%%%%%%%%%%%%%%%%%%%%%%%%%%%%%%%%%%%%%%%%%%%%%%%%%

\begin{proposition}\lab{lemma:firstcontroloftrXbcinMext}
We have on $\Mext$
\beaa
\max_{k\leq k_*}\sup_{\Mext}r^2u^{1+\dec}|\dk^k\trXc| &\les& \ep_0.
\eeaa
\end{proposition}

\begin{proof}
We  apply Proposition \ref{Prop:transportrp-f-Decay-knorms} to the equation  
\beaa
\nab_4\trXc +\frac{2}{q}\trXc &=& \Ga_g\c\Ga_g
\eeaa
and derive, for all $r_0\leq r\leq r_*$, 
\beaa
r^2\|\trXc  \|_{\infty, k}(u, r) &\les & r_*^2\| \trXc \|_{\infty, k}(u, r_*) +\int_r^{r_*} \la^2\| \Ga_g\c \Ga_g \|_{\infty, k}(u, \la) d\la.
\eeaa 
Thus,  in view of the  estimates  Proposition \ref{prop:improvedesitmatesfortemporalframeofMextonSigmastar-new}  for $\trXc$
on $\Si_*$, bootstrap assumptions  {\bf Ref 1} for   $\Ga_g$   and  interpolation  Remark \ref{remark:interpolation},  we derive
\beaa
 \sup_{r_{0} \le r\le r_*}  r^2\|\trXc  \|_{\infty, k}(u, r) &\les &    \ep_0u^{-1-\dec}, \qquad  k\le k_*, \qquad 1\leq u\leq u_*,
\eeaa
which concludes the proof of the proposition.
\end{proof}

%%%%%%%%%%%%%%%%%%%%%%%%%%%%%%%%%%

\subsection{Estimates for renormalized quantities in $\Mext$}

%%%%%%%%%%%%%%%%%%%%%%%%%%%%%%%%%%

We extend the estimates     for the renormalized quantities   $[H]_{ren}$, $[\widecheck{\DD\cos\th}] _{ren}$ and $[\Mc]_{ren}$  from $\Si_*$, as recorded
in Proposition   \ref{prop:improvedesitmatesfortemporalframeofMextonSigmastar-new}, to all of $\Mext$.   We do this  with the help of the transport equations derived in 
Proposition \ref{Prop: eqforrenormalized.qiantities}. We also make use of the estimates for $A$ in {\bf Ref 2}, as well as the above estimate for $\widecheck{\tr X}$. 
 
\begin{lemma}\lab{lemma:firstcontrolofrenormalizedquantitiesinMext}
We have on $\Mext$
\beaa
\max_{k\leq k_*}\sup_{\Mext}ru^{1+\dec}\left|\dk^k\big([H]_{ren}  \big) \right| &\les& \ep_0,
\eeaa
\beaa
\max_{k\leq k_*}\sup_{\Mext}ru^{1+\dec}\left|\dk^k\left([\widecheck{\DD\cos\th}] _{ren}   \right)\right| &\les& \ep_0,
\eeaa
and
\beaa
\max_{k\leq k_*-1}\sup_{\Mext}r^3u^{1+\dec}\Big | \dk^k\big([\Mc]_{ren} \big) \Big|  &\les& \ep_0.
\eeaa
\end{lemma}

\begin{proof}
 According to Proposition \ref{Prop: eqforrenormalized.qiantities}, we have
 \beaa
\bsplit
\nab_4\left(\ov{q}[\Hc]_{ren} \right) &=  O(r^{-1})\trXc+O(1)\dkb^{\leq 1}A+r\Ga_b\c\Ga_g,\\
\nab_4\left(q [\widecheck{\DD\cos\th}]_{ren}  \right) &= O(1)\trXc+O(r)A+r\Ga_b\c\Ga_g,\\
\nab_4\left(\ov{q} q^2 [\Mc]_{ren}\right)&=   O(1)\dkb^{\leq 1}\trXc+O(r)\dkb^{\leq 2}A+r^2\dkb^{\leq 1}(\Ga_g\c\Ga_g)+r^3\Ga_b\c A.
\end{split}
\eeaa
Applying Proposition \ref{Prop:transportrp-f-Decay-knorms} to the first identity,  and using the estimate   on $\Si_*$ for $[\Hc]_{ren}$,  we derive,  for all $k\le k_*$,  $r_0\leq r\le r_*$ and $1\leq u\leq u_*$,
\beaa
r \big\| [\Hc]_{ren} \big \|_{\infty, k}(u, r) &\les& r_* \big\| [\Hc]_{ren} \big \|_{\infty, k}(u, r_*)+
 \int_r^{r_*}   \| F\|_{\infty, k}(u, \la) d\la\\
 &\les&\ep_0 u^{-1-\dec}  +\int_r^{r_*}   \| F\|_{\infty, k}(u, \la) d\la
\eeaa
 where  $ F:= O(r^{-1})\trXc+O(1)\dkb^{\leq 1}A+r\Ga_b\c\Ga_g$.
 Using the available  estimates for $\trXc$ and  $A$, we obtain,  for all $k\le k_*$,  $r_0\leq r\le r_*$ and $1\leq u\leq u_*$,
 \beaa
\int_r^{r_*}   \| F\|_{\infty, k}(u, \la) d\la  &\les& \ep_0 r^{-1} u^{-1-\dec}.
 \eeaa
 Hence, we infer, for all $k\le k_*$,  $r_0\leq r\le r_*$ and $1\leq u\leq u_*$,
 \beaa
  r\big\| [\Hc]_{ren} \big \|_{\infty, k}(u, r)  &\les&\ep_0 u^{-1-\dec} 
 \eeaa
 as stated. The two other estimates are derived in the same manner.
 \end{proof}

%%%%%%%%%%%%%%%%%%%%%%%%%%%%%%%%%%%%%%%%%%%%%%%%%%

\subsection{Estimates for  some $\Lieb_\T$ derivatives in $\Mext$}

%%%%%%%%%%%%%%%%%%%%%%%%%%%%%%%%%%%%%%%%%%%%%%%%%%

\begin{proposition}\lab{prop:improveestiatesforsomeLiebTderivativesinMext}
We have on\footnote{Recall that $r\ge r_0$ on $\Mext$ for a sufficiently large $r_0$,  and that the small constant $\de'>0$, appearing in the estimate for $\Lieb_\T B$, is given by   $\de'=\frac 1 2 \big(\dee-\dec\big)$.} 
 $\Mext$, for  $r_0$ sufficiently large,
\beaa
\max_{k\leq k_*-1}\sup_{\Mext}r^2u^{1+\dec}|\dk^k\Lieb_\T\Xh| &\les& \ep_0,\\
\max_{k\leq k_*-2}\sup_{\Mext}r^4u^{1+\dec}|\dk^k\Lieb_\T^2B| &\les& \ep_0,\\
\max_{k\leq k_*-3}\sup_{\Mext}r^{3+\de'}u^{1+\dec}|\dk^k\Lieb_\T B| &\les& \ep_0,\\
\max_{k\leq k_*-3}\sup_{\Mext}r^2u^{1+\dec}|\dk^k\Lieb_\T\Zc| &\les& \ep_0,\\
\max_{k\leq k_*-4}\sup_{\Mext}r^3u^{1+\dec}|\dk^k\Lieb_\T\Pc| &\les& \ep_0.
\eeaa
\end{proposition}

\begin{proof}
The proof contains several steps. 

{\bf Step 1.} First, we estimate $\Lie_\T\Xh$. Starting with
\beaa
\nab_4\Xh+\Re(\tr X)\Xh &=& -A,
\eeaa
 and commuting  with $\Lieb_\T$ we derive, in view of Lemma \ref{Lemma:commLieb_Tnab-simb},
\beaa
\nab_4\Lieb_\T\Xh+\Re(\tr X)\Lieb_\T\Xh &=& -\Lieb_\T A+\dk^{\leq 1}(\Ga_g)\c\dk^{\leq 1}(\Ga_g)
\eeaa
and hence
\beaa
\nab_4\Lieb_\T\Xh+\Re\left( \frac{2}{q}\right) \Lieb_\T\Xh &=& -\Lieb_\T A+\dk^{\leq 1}(\Ga_g)\c\dk^{\leq 1}(\Ga_g).
\eeaa
Hence, applying Proposition \ref{Prop:transportrp-f-Decay-knorms}, we have, for all $r_0\leq r\leq r_*$ and $1\leq u\leq u_*$, 
\beaa
r^{2}\|\Lieb_\T \Xh \|_{\infty, k}(u, r)&\les  & r_*^{2}\|\Lieb_\T \Xh \|_{\infty, k}(u, r_*)  +\int_r^{r_*} \la^{2}\| F\|_{\infty, k}(u, \la) d\la
\eeaa
where
\beaa
F&=&  -\Lieb_\T A+\dk^{\leq 1}(\Ga_g)\c\dk^{\leq 1}(\Ga_g).
\eeaa
 According  to  the improved estimates for $A$ of Theorem M1, we have, recalling that the small constant $\de'>0$ is given by $\de'=1/2 (\dee-\dec)$,
 \beaa
\max_{k\leq k_*}\sup_{\Mext}r^{2+ \de'}u^{1+\dec}|\dk^kA| +\max_{k\leq k_*-1}\sup_{\Mext}r^{3+ \de'}u^{1+\dec}|\dk^k\nab_3 A| &\les& \ep_0.
\eeaa
Since  $\Lieb_\T A = \nab_3A+r^{-1}\dk^{\leq 1}A$, and together with the estimates for $\Ga_g$ and the form of $F$, we infer
\beaa
\max_{k\leq k_*-1}\sup_{\Mext}r^{3+\de'}u^{1+\dec}|\dk^k F| &\les& \ep_0.
\eeaa
Taking into account the  control on $\Si_*$, we deduce, for all $k\le k_*-1$, $r_0\leq r\leq r_*$ and $1\leq u\leq u_*$,  
\beaa
r^{2}\|\Lieb_\T \Xh \|_{\infty, k}(u, r)&\les  & \ep_0 u^{-1-\dec}
\eeaa
as stated.

{\bf Step 2.} Next, we estimate $\Lieb^2_\T B $. Recall that we have
\beaa
\nab_4B+\frac{4}{\ov{q}}B &=& \frac{1}{2}\ov{\DD}\c A+\frac{aq}{2|q|^2}\ov{\Jk}\c A+\Ga_g\c(B,A).
\eeaa
We commute with $\Lieb_\T$ and obtain, using Lemma \ref{Lemma:commLieb_Tnab-simb},  
\beaa
\nab_4\Lieb_\T B+\frac{4}{\ov{q}}\Lieb_\T B &=& \frac{1}{2}\ov{\DD}\c \Lieb_\T A+\frac{aq}{2|q|^2}\ov{\Jk}\c\Lieb_\T A+\dk^{\leq 1}(\Ga_g)\dk^{\leq 1}(A,B).
\eeaa
Commuting with  $\Lieb_\T$ again, we  derive
\beaa
\nab_4\Lieb_\T^2 B+\frac{4}{\ov{q}}\Lieb_\T^2B &=& \frac{1}{2}\ov{\DD}\c \Lieb_\T^2 A+\frac{aq}{2|q|^2}\ov{\Jk}\c\Lieb_\T^2 A+\dk^{\leq 2}(\Ga_g)\dk^{\leq 2}(A,B).
\eeaa
Hence, applying Proposition \ref{Prop:transportrp-f-Decay-knorms}, we infer, for all $r_0\leq r\leq r_*$ and $1\leq u\leq u_*$,
\beaa
r^{4}\|\Lieb^2_\T  B\|_{\infty, k}(u, r)&\les  & r_*^{4}\|\Lieb^2 _\T B \|_{\infty, k}(u, r_*) +\int_r^{r_*} \la^{4}\| F\|_{\infty, k}(u, \la) d\la
\eeaa
where
\beaa
F&=&  \frac{1}{2}\ov{\DD}\c \Lieb_\T^2 A+O(r^{-2})\Lieb_\T^2 A+\dk^{\leq 2}(\Ga_g)\dk^{\leq 2}(A,B).
\eeaa
Making use of the   improved estimate for $A$, $\nab_3A$ and $\nab_3^2 A$ in {\bf Ref 2}, 
we have, using also the relation between $\Lieb_\T^2$ and $\nab_3^2$ in \eqref{eq:Leib_T-nab_3U:immediateconsequences}, 
\beaa
\max_{k\leq k_*-1}\sup_{\Mext}r^{4+\de'}u^{1+\dec}|\dk^k\Lieb_\T^2 A| &\les& \ep_0.
\eeaa
Therefore, together with the definition of $F$, and the control for $B$ and $\Ga_g$ provided by {\bf Ref 1}, we infer, for $k\leq k_*-1$, and for all $r_0\leq r\leq r_*$ and $1\leq u\leq u_*$, 
\beaa
\int_r^{r_*} \la^{4}\| F\|_{\infty, k}(u, \la) d\la  &\les& \ep_0 u^{-1-\dec}.
\eeaa
Using also the estimates for $ \Lieb_\T^2 B$  on $\Si_*$, we obtain, for all $k\le k_*-2$, $r_0\leq r\leq r_*$ and $1\leq u\leq u_*$, 
\beaa
 r^{4}\|\Lieb^2_\T  B\|_{\infty, k}(u, r)&\les  &\ep_0 u^{-1-\dec}
\eeaa
as stated.

{\bf Step 3.} Next, we estimate   $\DD\hot\Lie_\T B$. We start with the equation 
\beaa
\DD\hot B +(Z+4H)\hot B &=& \nab_3A+\left(\frac{1}{2}\tr\Xb-4\omb\right)A+3\ov{P}\Xh.
\eeaa
Commuting with $\Lieb_\T$,  and  using Lemma \ref{Lemma:commLieb_Tnab-simb}   for the commutator, 
we infer
\beaa
\DD\hot\Lieb_\T B+(Z+4H)\hot\Lieb_\T B &=& \Lieb_\T\nab_3A+\left(\frac{1}{2}\tr\Xb-4\omb\right)\Lieb_\T A+3\ov{P}\Lieb_\T\Xh\\
&&+r^{-2}\dk^{\leq 1}(\Ga_b)\dk^{\leq 1}(\Ga_g)+\dk^{\leq 1}(\Ga_b)A.
\eeaa
Using  the relation between $\Lieb_\T$ and $\nab_3$ in \eqref{eq:Leib_T-nab_3U:immediateconsequences}, as well as the fact that $Z= O(r^{-2} )+ \Ga_g$ and $H=O(r^{-2}) +\Ga_b$,   we deduce
\beaa
\DD\hot\Lieb_\T B &=& O(r^{-2} )\Lieb_\T B +\nab_3^2A+O(r^{-1})\dk^{\leq 1}\nab_3A+O(r^{-2})\dk^{\leq 1}A+O(r^{-3})\Lieb_\T\Xh\\
&& +r^{-2}\dk^{\leq 1}(\Ga_b)\dk^{\leq 1}(\Ga_g)+\dk^{\leq 1}(\Ga_b)\dk^{\leq 1}(A, B).
\eeaa
In view of  Corollary  \ref{cor:formuallikningnabprimeandnabintransfoformula:bis} we have
\beaa
\DD' &=& \Big(1+O(r^{-2})\Big)\DD+O(r^{-1})\Lieb_\T +O(r^{-3})+r^{-1}\Ga_b\dk^{\leq 1}\\
&=&\DD + O(r^{-3})\dkb^{\leq 1} +O(r^{-1})\Lieb_\T+r^{-1}\Ga_b\dk^{\leq 1},
\eeaa
and thus
\beaa
\DD' \hot\Lieb_\T B &=& \DD\hot\Lieb_\T B+ O(r^{-3}) \dkb^{\leq 1}\Lieb_\T B +O(r^{-1})\Lieb_\T^2 B +r^{-1}\Ga_b\dk^{\leq 1}\Lieb_\T B. 
\eeaa
Hence,
\beaa
\DD' \hot\Lieb_\T B &=& O(r^{-1}) \Lieb^2_\T B+  O(r^{-2}) \dkb^{\leq 1}\Lieb_\T B  +\nab_3^2A+O(r^{-1})\dk^{\leq 1}\nab_3A+O(r^{-2})\dk^{\leq 1}A\\
&&+O(r^{-3})\Lieb_\T\Xh +r^{-2}\dk^{\leq 1}(\Ga_b)\dk^{\leq 1}(\Ga_g)+\dk^{\leq 1}(\Ga_b)\dk^{\leq 1}(A, B).
\eeaa
Using again  Corollary  \ref{cor:formuallikningnabprimeandnabintransfoformula:bis} we can express   $\dkb \Lieb_\T B$ in terms of $\dkb' \Lieb_\T B $ and derive
\beaa
\DD' \hot\Lieb_\T B &=& O(r^{-1}) \Lieb^2_\T B+  O(r^{-2}) (\dkb')^{\leq 1}\Lieb_\T B  +\nab_3^2A+O(r^{-1})\dk^{\leq 1}\nab_3A+O(r^{-2})\dk^{\leq 1}A\\
&&+O(r^{-3})\Lieb_\T\Xh +r^{-2}\dk^{\leq 1}(\Ga_b)\dk^{\leq 1}(\Ga_g)+\dk^{\leq 1}(\Ga_b)\dk^{\leq 1}(A, B).
\eeaa
Differentiating this identity with respect to $( \dkb' )^k$, using the  estimates for $\nab^2_3A$, $\nab_3A$ and $A$  in {\bf Ref 2},  the improved estimates for $\Lieb_\T\Xh$  in Step 1, the improved estimate for $\Lieb^2_\T B$ 
 in Step 2,  and the bootstrap assumptions  for $\Ga_b$, $\Ga_g$ and $B$, we deduce on any sphere $S=S(u, r)$ of $\Mext$, for $k\le k_*-2$,
\bea
\lab{eq:improveestiatesforsomeLiebTderivativesinMext2}
\bsplit
\big\|( \dkb' )^k \DD' \hot\Lieb_\T B\big\|_{L^2 (S)} &\les   O(r^{-2}) \big\|(\dkb')^{\le k+1}  \Lieb_\T B\big\|_{L^2 (S)}  + \ep_0r^{-3-\de'} u^{-1-\dec}.
\end{split}
\eea

{\bf Step 4.} Next, we derive a first estimate for  $\Lieb_\T B$. To this end, we apply  the first elliptic estimate of Corollary \ref{Lemma:Hodge-estimateDDs_2}  to derive
\beaa
\|(\dkb')^{k+1}  \Lieb_\T B\big\|_{L^2 (S)} &\les & r  \big\| (\dkb')^{\le k} \DD' \hot\Lieb_\T B\big\|_{L^2 (S)}+ r^2  \Big| \big(\ov{ \DD'} \c \Lieb_\T B\big)_{\ell=1}  \Big|.
\eeaa
In view of \eqref{eq:improveestiatesforsomeLiebTderivativesinMext2}, we deduce on any sphere $S=S(u, r)$ of $\Mext$, for $k\le k_*-2$,
\beaa
\| (\dkb' )^{k+1}\Lieb_\T B\big\|_{L^2 (S)} &\les &   O(r^{-1}) \big\|(\dkb')^ {\le k+ 1}  \Lieb_\T B\big\|_{L^2 (S)}  + r^2  \Big| \big(\ov{ \DD'} \c \Lieb_\T B\big)_{\ell=1}  \Big|\\
&& + \ep_0r^{-2-\de'} u^{-1-\dec}.
\eeaa
Thus, recalling that $r \ge r_0$ on $\Mext$,  and provided that $r_0$ is  sufficiently large, we  can absorb the first term on the RHS and deduce, for $k\le k_*-2$,
\bea
\| (\dkb' )^{k+1}\Lieb_\T B\big\|_{L^2 (S)}&\les &  r^2  \Big| \big(\ov{ \DD'} \c \Lieb_\T B\big)_{\ell=1}  \Big|+ \ep_0r^{-2-\de'} u^{-1-\dec}.
\eea
Using again
\beaa
\DD &=& \DD'+O(r^{-3})\dkb^{\leq 1}+O(r^{-1})\Lieb_\T+r^{-1}\Ga_b\dk^{\leq 1},
\eeaa
and  the above improved estimate for $\Lieb_\T^2 B$, we further  deduce, for $k\le k_*-2$,
\beaa
\|\dkb^{k+1}  \Lieb_\T B\big\|_{L^2 (S)}&\les &  r^2  \Big| \big(\ov{ \DD} \c \Lieb_\T B\big)_{\ell=1}  \Big|+ \ep_0r^{-2-\de'} u^{-1-\dec}.
\eeaa
By Sobolev, we finally obtain on any sphere $S=S(u, r)$ of $\Mext$, for $k\le k_*-3$,
\bea
\lab{eq:improveestiatesforsomeLiebTderivativesinMext3} 
\|\dkb^{k}  \Lieb_\T B\big\|_{L^\infty (S)} &\les &  r  \Big| \big(\ov{ \DD} \c \Lieb_\T B\big)_{\ell=1}  \Big|+ \ep_0 r^{-3-\de'} u^{-1-\dec}.
\eea

{\bf Step 5.} Next, we estimate $(\ov{\DD}\c\Lieb_\T B)_{\ell=1}$. Starting from 
\beaa
\nab_4B+\frac{4}{\ov{q}}B &=& \frac{1}{2}\ov{\DD}\c A+\frac{aq}{2|q|^2}\ov{\Jk}\c A+\Ga_g\c(B,A),
\eeaa
we commute with $\Lieb_\T$ and obtain, using Lemma \ref{Lemma:commLieb_Tnab-simb},  the following  more precise transport equation for $\nab_4\Lieb_\T B$ compared the one used in Step 2
\beaa
\nab_4\Lieb_\T B+\frac{4}{\ov{q}}\Lieb_\T B &=& \frac{1}{2}\ov{\DD}\c \Lieb_\T A+\frac{aq}{2|q|^2}\ov{\Jk}\c\Lieb_\T A+(r^{-1}\dk^{\leq 1}(\Ga_b),\Lieb_\T\Ga_g)\dk^{\leq 1}(A,B)\\
&&+\Ga_g\Lieb_\T(A,B).
\eeaa
Since $\nab_3\Ga_g=r^{-1}\dk^{\leq 1}\Ga_b$, and in view of the link between $\Lieb_\T$ and $\nab_3$, we infer
\beaa
\nab_4\Lieb_\T B+\frac{4}{\ov{q}}\Lieb_\T B &=& \frac{1}{2}\ov{\DD}\c \Lieb_\T A+\frac{aq}{2|q|^2}\ov{\Jk}\c\Lieb_\T A+r^{-1}\dk^{\leq 1}(\Ga_b)\dk^{\leq 1}(A,B)+\Ga_g\nab_3(A,B)\\
&=& \frac{1}{2}\ov{\DD}\c \Lieb_\T A+O(r^{-2})\Lieb_\T A+r^{-1}\dk^{\leq 1}(\Ga_b)\dk^{\leq 1}(A,B)+\Ga_g\nab_3(A,B)\\
&=& \frac{1}{2}\ov{\DD}\c \Lieb_\T A+O(r^{-2} )\nab_3  A  +O(r^{-3}) \dk^{\le 1} A +r^{-1}\dk^{\leq 1}(\Ga_b)\dk^{\leq 1}(A,B)\\
&&+\Ga_g\nab_3(A,B)
\eeaa
and hence
\beaa
\nab_4(\ov{q}^4\Lieb_\T B) &=& \frac{\ov{q}^4}{2}\ov{\DD}\c \Lieb_\T A+O(r^{2} )\nab_3  A  +O(r) \dk^{\le 1} A +r^3\dk^{\leq 1}(\Ga_b)\dk^{\leq 1}(A,B)\\
&&+r^4\Ga_g\nab_3(A,B).
\eeaa
We commute with $\ov{q}\,\ov{\DD}\c$ relying on the following commutator estimate, see Corollary  \ref{cor:commutation-complexM6},
\beaa
 \, [\nab_4, \ov{q}\,\ov{\DD}\c] &=& O(r^{-2})+ \Ga_g \c \dkb^{\leq 1}.
 \eeaa
 Using also $\DD(q) = O(r^{-1})+\Ga_b$,  $q=r+O(1)$, and the link between $\Lieb_\T$ and $\nab_3$, we deduce
 \beaa
\nab_4(\ov{q}\,\ov{\DD}\c(\ov{q}^4\Lieb_\T B)) &=& O(r^2)\Lieb_\T B+\frac{r^5}{2}\,\ov{\DD}\c\ov{\DD}\c \Lieb_\T A +O(r^2)\dkb^{\leq 2}\nab_3A+O(r)\dk^{\leq 3}A\\
&&+r^3\dkb^{\leq 1}(\dk^{\leq 1}(\Ga_b)\dk^{\leq 1}(A,B))+r^4\dk^{\leq 1}(\Ga_g\nab_3(A,B)).
\eeaa
Since $e_4(J^{(p)})=0$ on $\Mext$, and using $J^{(p)}=O(1)$, we derive, for $p=0,+,-$, 
\beaa
\nab_4(\ov{q}\,\ov{\DD}\c(\ov{q}^4\Lieb_\T B)\Jp) &=& O(r^2)\Lieb_\T B+\frac{r^5}{2}\,\ov{\DD}\c\ov{\DD}\c \Lieb_\T A \Jp+O(r^2)\dkb^{\leq 2}\nab_3A+O(r)\dk^{\leq 3}A\\
&&+r^3\dkb^{\leq 1}(\dk^{\leq 1}(\Ga_b)\dk^{\leq 1}(A,B))+r^4\dk^{\leq 1}(\Ga_g\nab_3(A,B)).
\eeaa
Next we integrate on $S$ with the help of  Corollary \ref{cor:e4derivativeofintegralonS} according to which, for a   scalar function $h$ on $S=S(u,r)$,
\beaa
e_4\left(\int_{S(r,u)}h \right) &=& \int_{S(r,u)}\frac{e_4(\Si h)}{\Si}+O\left(\frac{\ep}{u^{1+\dec}}\right)h.
\eeaa 
Applying this to 
\bea
h:= \frac{1}{\Si}\ov{q}\,\ov{\DD}\c(\ov{q}^4\Lieb_\T B)J^{(p)},
\eea
we deduce, since  $\Sigma=r^2+O(1)$,
\beaa
e_4\left(\int_{S(r,u)}h \right) &=& \int_S\frac{1}{\Si}\nab_4(\ov{q}\,\ov{\DD}\c(\ov{q}^4\Lieb_\T B)J^{(p)})+O\left(\frac{\ep}{u^{1+\dec}}\right)\frac{1}{\Si}\ov{q}\,\ov{\DD}\c(\ov{q}^4\Lieb_\T B)J^{(p)}
\\
&=&O(r^2)\Lieb_\T B+ \frac{r^3}{2}\left(\int_S\ov{\DD}\c\ov{\DD}\c \Lieb_\T A J^{(p)}\right)+O(r^2)\dkb^{\leq 2}\nab_3A+O(r)\dk^{\leq 3}A\\
&&+r^3\dkb^{\leq 1}(\dk^{\leq 1}(\Ga_b)\dk^{\leq 1}(A,B))+r^4\dk^{\leq 1}(\Ga_g\nab_3(A,B))\\
&&+O\left(\frac{r^2\ep}{u^{1+\dec}}\right)\dkb^{\leq 1}\Lieb_\T B.
\eeaa
Now, using again $ \DD =\DD'+O(r^{-3})\dkb^{\leq 1}+O(r^{-1})\Lieb_\T+r^{-1}\Ga_b\dk^{\leq 1}$, we have
\beaa
\ov{\DD}\c\ov{\DD}\c \Lieb_\T A &=& \ov{\DD'}\c\ov{\DD'}\c \Lieb_\T A+O(r^{-2})\dk^{\leq 1}\nab_3^2A+O(r^{-3})\dk^{\leq 2}\nab_3A+O(r^{-4})\dk^{\leq 3}A\\
&&+r^{-1}\dk^{\leq 1}(\Ga_g\c\dk^{\leq 1}\nab_3A)+r^{-2}\dk^{\leq 3}(\Ga_g\c A)
\eeaa 
and hence
\beaa
e_4\left(\int_{S(r,u)}h \right) &=&O(r^2)\Lieb_\T B+ \frac{r^3}{2}\left(\int_S\ov{\DD'}\c\ov{\DD'}\c \Lieb_\T A J^{(p)}\right)+O(r^3)\dk^{\leq 1}\nab_3^2A\\
&&+O(r^2)\dk^{\leq 2}\nab_3A+O(r)\dk^{\leq 3}A+r^3\dk^{\leq 3}(\Ga_b(A,B))\\
&&+r^4\dk^{\leq 1}(\Ga_g\dk^{\leq 1}\nab_3(A,B))+O\left(\frac{r^2\ep}{u^{1+\dec}}\right)\dkb^{\leq 1}\Lieb_\T B.
\eeaa
Integrating by parts and making use of, see Proposition \ref{prop:controlofDDprimehotDDprimeJonMext},
\beaa
\left|\DD'\hot\DD'J^{(p)}\right| &\les& \frac{\ep}{r^3u^{\frac{1}{2}+\dec}}+\frac{1}{r^4},
\eeaa
 we deduce
\beaa
\int_S \ov{\DD'}\c\ov{\DD'}\c  \Lieb_\T A J^{(p)} &=& \int_S\Lieb_\T A\,  \ov{\DD'\hot\DD' J^{(p)}}= O\left(\frac{\ep}{ru^{\frac{1}{2}+\dec}}\right)\Lieb_\T A+O(r^{-2})\Lieb_\T A.
\eeaa
Therefore
\beaa
e_4\left(\int_{S(r,u)}h \right) &=&O(r^2)\Lieb_\T B+O(r^3)\dk^{\leq 1}\nab_3^2A+O(r^2)\dk^{\leq 2}\nab_3A+O(r)\dk^{\leq 3}A\\
&&+r^3\dk^{\leq 3}(\Ga_b(A,B))+r^4\dk^{\leq 1}(\Ga_g\dk^{\leq 1}\nab_3(A,B))\\
&&+O\left(\frac{r^2\ep}{u^{1+\dec}}\right)\dkb^{\leq 1}\Lieb_\T B+O\left(\frac{r^2\ep}{u^{\frac{1}{2}+\dec}}\right)\Lieb_\T A.
\eeaa
Together with the improved  control on  $A$, $\nab_3 A$, $\nab^2_3 A$  in {\bf Ref 2},  as well as the bootstrap assumptions in {\bf Ref 1},  we obtain
\beaa
e_4\left(\int_{S(r,u)}h \right) &=&O(r^2)\Lieb_\T B+O\left(\frac{\ep_0}{r^{1+\de'}u^{1+\dec}}\right).
\eeaa
Now, recall that $h$ is given by
\beaa
h= \frac{1}{\Si}\ov{q}\,\ov{\DD}\c(\ov{q}^4\Lieb_\T B)J^{(p)},
\eeaa
so that, together with $\Si=r^2+O(1)$ and $q=r+O(1)$, we have 
\bea\lab{eq:intermediaryrelationbetweenintShandovDDcBell=1}
\int_Sh &=& r^5\big(\ov{ \DD} \c \Lieb_\T B\big)_{\ell=1} +O(r^3)\dkb^{\leq 1}\Lieb_\T B.
\eea
Together with \eqref{eq:improveestiatesforsomeLiebTderivativesinMext3}, we infer
\beaa
r^2\Lieb_\T B &=& O(r^{-2})\int_Sh+ O(\ep_0 r^{-1-\de'} u^{-1-\dec} ).
\eeaa
Plugging the second identity in the above, we obtain
\beaa
e_4\left(\int_{S(r,u)}h \right) &=& O(r^{-2})\int_Sh+O\left(\frac{\ep_0}{r^{1+\de'}u^{1+\dec}}\right).
\eeaa
Using Gronwall lemma, we deduce
\beaa
\left|\int_{S(r,u)}h\right| &\les& \left|\int_{S(r_*,u)}h\right|+\frac{\ep_0}{u^{1+\dec}}.
\eeaa
Together with \eqref{eq:intermediaryrelationbetweenintShandovDDcBell=1}, and the control of $\big(\ov{ \DD} \c \Lieb_\T B\big)_{\ell=1}$ and $\Lieb_\T B$ on $\Si_*$, this yields, on any sphere $S=S(u, r)$ of $\Mext$,
\bea
\lab{eq:improveestiatesforsomeLiebTderivativesinMext4} 
 \Big|(\ov{\DD}\c\Lieb_\T B)_{\ell=1}\Big| &\les& \ep_0 r^{-5} u^{-1-\dec}+ r^{-2} \|\dkb^{\leq 1}\Lieb_\T B\|_{L^\infty(S)}.
\eea

{\bf Step 6.} We are now in position to conclude the estimate for $\Lieb_\T B$. By combining  \eqref{eq:improveestiatesforsomeLiebTderivativesinMext4}  
with the estimate  \eqref{eq:improveestiatesforsomeLiebTderivativesinMext3}  derived in Step 4, we have, on any sphere $S=S(u, r)$ of $\Mext$, for $k\le k_*-3$,
\beaa
\|\dkb^{k}  \Lieb_\T B\big\|_{L^\infty (S)}&\les &  r  \Big| \big(\ov{ \DD} \c \Lieb_\T B\big)_{\ell=1}  \Big|+ \ep_0r^{-3-\de'} u^{-1-\dec}\\
&\les&  r^{-1}\|\dkb^{\leq 1}\Lieb_\T B\|_{L^\infty(S)}+  \ep_0r^{-3-\de'} u^{-1-\dec}.
\eeaa
Since $r\geq r_0$ on $\Mext$, for $r_0$ sufficiently large, we may absorb the first term on the RHS. We deduce, for every sphere $S= S(u, r)$ in $\Mext$,
\beaa
\|\dkb^{\le k_*-3}  \Lieb_\T B\big\|_{L^\infty (S)}&\les & \ep_0r^{-3-\de'} u^{-1-\dec}.
\eeaa
Together with the estimates for $\Lieb^2_\T B$ of Step 2, we infer
\beaa
\|(\Lieb_\T, \dkb)^{\le k_*-3}  \Lieb_\T B\big\|_{L^\infty (S)}&\les & \ep_0r^{-3-\de'} u^{-1-\dec}.
\eeaa
Finally, recalling the transport equation for $\nab_4\Lieb_\T B$ which yields
\beaa
\nab_4\Lieb_\T B &=& -\frac{4}{\ov{q}}\Lieb_\T B+O(r^{-1})\dk^{\leq 1}\nab_3  A  +O(r^{-2}) \dk^{\le 2} A +r^{-1}\dk^{\leq 1}(\Ga_b)\dk^{\leq 1}(A,B)\\
&&+\Ga_g\nab_3(A,B),
\eeaa
and together with the control of $A$ and $\nab_3A$ provided by {\bf Ref 2}, and the bootstrap assumptions in {\bf Ref 1}, we obtain 
\beaa
\|(re_4, \Lieb_\T, \dkb)^{\le k_*-3}  \Lieb_\T B\big\|_{L^\infty (S)}&\les & \ep_0r^{-3-\de'} u^{-1-\dec}.
\eeaa
Together with the fact that $e_3$ is spanned by $\T$, $e_4$ and $(e_1, e_2)$, we infer, for all $k\le k_*-3$, $r_0\leq r\leq r_*$ and $1\leq u\leq u_*$,  
\bea
r^{3+\de'}\|\Lieb_\T B \|_{\infty, k}(u, r)&\les  & \ep_0 u^{-1-\dec}
\eea
as stated.

{\bf Step 7.}  Next, we estimate $\Lieb_\T  \Pc$ with the help of the equation
\beaa
\nab_4\Pc+\frac{3}{q}\Pc &=& \frac{1}{2}\DD\c\ov{B} -\frac{a\ov{q}}{2|q|^2}\Jk\c\ov{B}+O(r^{-3})\trXc +\Ga_b\c A+r^{-1}\Ga_g\c\Ga_g.
\eeaa
We commute with $\Lieb_\T$ using, see  Lemma \ref{Lemma:commLieb_Tnab-simb},
\beaa
\,[\nab_4, \Lieb_\T] &=& r^{-1}\dk^{\leq 1}(\Ga_b)\dk^{\leq 1},\qquad 
\, [\Lieb_\T, \nab]   = r^{-1}  \dk^{\leq 1}(\Ga_b)\dk^{\leq 1},
\eeaa
and use also the fact that $\nab_3(\Ga_g)=r^{-1}\dk^{\leq 1}\Ga_b$ and  $\T(q)\in r\Ga_b$. We obtain
\beaa
\nab_4\Lieb_\T\Pc+\frac{3}{q}\Lieb_\T\Pc &=& \frac{1}{2}\DD\c\Lieb_\T\ov{B} -\frac{a\ov{q}}{2|q|^2}\Jk\c\Lieb_\T\ov{B}+O(r^{-3})\dk^{\leq 1}\trXc \\
&&+\dk^{\leq 1}(\Ga_b\c (A,B))+r^{-2}\dk^{\leq 1}(\Ga_b\c\Ga_g).
\eeaa
We now make use of Proposition \ref{Prop:transportrp-f-Decay-knorms} and deduce, for all  $r_0\leq r\leq r_*$ and $1\leq u\leq u_*$,  
\beaa
r^{3}\|\Lieb_\T \Pc \|_{\infty, k}(u, r) &\les & r_*^{3}\| \Lieb_\T \Pc \|_{\infty, k}(u, r_*) +\int_r^{r_*} \la^{c}\| F\|_{\infty, k}(u, \la) d\la 
\eeaa
with
\beaa
F&=&  \frac{1}{2}\DD\c\Lieb_\T\ov{B} -\frac{a\ov{q}}{2|q|^2}\Jk\c\Lieb_\T\ov{B}+O(r^{-3})\dk^{\leq 1}\trXc +\dk^{\leq 1}(\Ga_b\c (A,B))+r^{-2}\dk^{\leq 1}(\Ga_b\c\Ga_g).
\eeaa
Hence, using the estimates  for $ \Lieb_\T B$  derived in Step 6, the estimate for $\trXc$ of Proposition \ref{lemma:firstcontroloftrXbcinMext}, the bootstrap assumptions in {\bf Ref 1}, and the estimate for  $ \Lieb_\T \Pc$ on $\Si_*$, we deduce, for all $k\le k_*-4$, $r_0\leq r\leq r_*$ and $1\leq u\leq u_*$,  
\bea
\big|\Lieb_\T\Pc\big|_{\infty,k}(u,r) &\les& \ep_0 r^{-3} u^{-1-\dec}
\eea
as stated.

{\bf Step 8.}  Finally, we estimate $\Lieb_\T\Zc$. Recall
\beaa
\nab_4\Zc + \frac{2}{q}\Zc  &=&    - \frac{aq}{|q|^2}\ov{\Jk}\c\widehat{X} -B +O(r^{-2})\trXc+\Ga_g\c\Ga_g.
\eeaa
We commute with $\Lieb_\T$. Using again   $\,[\nab_4, \Lieb_\T] = r^{-1}\dk^{\leq 1}(\Ga_b)\dk^{\leq 1},
$
and 
\beaa
\T(q)=\T(r)+ia\T(\cos\th)\in r\Ga_b, \qquad  \Lieb_\T\Jk=\widecheck{\nab_3\Jk}+O(r^{-1})\widecheck{\nab\Jk}+r^{-1}\dk(\Ga_b)\Jk\in r^{-1}\dk^{\leq 1}\Ga_b,
\eeaa 
we infer
\beaa
\nab_4\Lieb_\T\Zc + \frac{2}{q}\Lieb_\T\Zc  &=&    - \frac{aq}{|q|^2}\ov{\Jk}\c\Lieb_\T\widehat{X} -\Lieb_\T B +O(r^{-2})\dk^{\leq 1}\trXc+\dk^{\leq 1}(\Ga_g\c\Ga_g)\\
&=&  O(r^{-2})\Lieb_\T\widehat{X} -\Lieb_\T B +O(r^{-2})\dk^{\leq 1}\trXc+\dk^{\leq 1}(\Ga_g\c\Ga_g).
\eeaa
Proceeding as before with the help of  Proposition \ref{Prop:transportrp-f-Decay-knorms}, using  the estimates  for $ \Lieb_\T\Xh$  derived in Step 1, the estimates  for $ \Lieb_\T B$  derived in Step 6, the estimate for $\trXc$ of Proposition \ref{lemma:firstcontroloftrXbcinMext}, the bootstrap assumptions in {\bf Ref 1}, and the estimate for  $ \Lieb_\T \Zc$ on $\Si_*$, we deduce,  for all $k\le k_*-3$, $r_0\leq r\leq r_*$ and $1\leq u\leq u_*$,  
\bea
\big|\Lieb_\T\Zc\big|_{\infty,k}(u,r) &\les& \ep_0 r^{-2} u^{-1-\dec}
\eea
as stated. This concludes the proof of Proposition \ref{prop:improveestiatesforsomeLiebTderivativesinMext}.
\end{proof}

%%%%%%%%%%%%%%%%%%%%%%%%%%%%%%%%%%%%%%%%%%%%%%%%%%%%%

\section{Improved decay estimates for $B$, $\Pc$, $\Xh$, $\Zc$, $\Hc$, and $\widecheck{\DD\cos\th}$}

%%%%%%%%%%%%%%%%%%%%%%%%%%%%%%%%%%%%%%%%%%%%%%%%%%%%%

%%%%%%%%%%%%%%%%%%%%%%%%%%%%%%%%%%%%%%%%%%%%%%%%%%

\subsection{Conditional  control of $B$, $\Pc$, $\Xh$, $\Zc$, $\Hc$, and $\widecheck{\DD\cos\th}$}

%%%%%%%%%%%%%%%%%%%%%%%%%%%%%%%%%%%%%%%%%%%%%%%%%%

The goal of this section is to prove the following proposition.
\begin{proposition}\lab{prop:firstcontrolofBPcXhZcHcDDcosth}
We have on $\Mext$, with  $\de'=\frac 1 2( \dee-\dec)>0$,  for all $k\le   k_*-5$,
\bea
\lab{eq-prop:firstcontrolofBPcXhZcHcDDcosth1}
 r\big(|\dk^{\leq k}B|+|\dk^{\leq k}\Pc|\big)+|\dk^{\leq k}\Xh|+|\dk^{\leq k}\Zc|
 &\les& r^2|(\ov{\DD}\c B)_{\ell=1}|+ \ep_0 r^{-1-\de'}u^{-1-\dec}
\eea
and
\bea
\lab{eq-prop:firstcontrolofBPcXhZcHcDDcosth2}
 |\dk^{\leq k }\Hc|+|\dk^{\leq k}\widecheck{\DD\cos\th}| &\les& r^2|(\ov{\DD}\c B)_{\ell=1}|+\ep_0 r^{-1} u^{-1-\dec}.
\eea
In addition we prove the following preliminary estimate\footnote{A stronger estimate, i.e. with $\ep$ replaced by $\ep_0$, will be proved  later.}   for $B$
\bea
\lab{eq:Prliminary-forB}
|\dk ^{\leq {k_*-2}}B| &\les& r |(\ov{\DD}\c B)_{\ell=1}|+      \ep r^{-3-\de'} u^{-1/2-\dec}.
\eea
\end{proposition}

\begin{proof}
We will repeatedly make use of   the notation introduced in Definition \ref{Notation:forGo_k}, i.e.  $U\in r^{-p}\Go_k$  if $U$ satisfies $|\dk^{\leq k} U| \les \ep_0 r^{-p} u^{-1-\dec}$.

{\bf Step 1.} In view of {\bf Ref 2} and  Proposition \ref{lemma:firstcontroloftrXbcinMext}, we have on $\Mext$
\bea\lab{eq:improvmentwithGoodtermsinMext:1}
A \in r^{-2-\de'}\Go_{k_*},\quad \nab_3A \in r^{-3-\de'}\Go_{k_*-1},\quad \trXc\in r^{-2}\Go_{k_*-1},
\eea
Also, in view of Lemma \ref{lemma:firstcontrolofrenormalizedquantitiesinMext}, we have on $\Mext$
\bea
\lab{eq:improvmentwithGoodtermsinMext:11}
[\Hc]_{ren},  \, \,   [\widecheck{\DD\cos\th}] _{ren}    \in  r^{-1}\Go_{k_*},\qquad [\Mc]_{ren}\in   r^{-3}\Go_{k_*-1},
\eea 
 where, recall,
\beaa
\bsplit
\ov{q} [\Hc]_{ren}    &=       \ov{q}\Hc -q\Zc +\frac{1}{3}\left(-\ov{q}^2+|q|^2\right)B+\frac{a}{2}(q-\ov{q})\ov{\Jk}\c\Xh,\\
q [\widecheck{\DD\cos\th}] _{ren} &= q\widecheck{\DD\cos\th} +\frac{i}{2}|q|^2\ov{\Jk}\c\Xh,\\
\ov{q} q^2 [\Mc]_{ren} &=\ov{q}\,\ov{\DD}\c\left(q^2\Zc  +\left(-\frac{a}{2}q^2 -\frac{a}{2}|q|^2\right)\ov{\Jk}\c\Xh\right) +2\ov{q}^3\,\ov{\Pc} - 2aq^2\ov{\Jk}\c\Zc\\
&+\left(-\frac{1}{3}q^2\ov{q}^2 -\frac{1}{3}q\ov{q}^3+\frac{2}{3}\ov{q}^4\right)\ov{\DD}\c B +a\left( q^2\ov{q} +\frac{2}{3}q\ov{q}^2 - \frac{13}{6}\ov{q}^3\right)\ov{\Jk}\c B\\
&+a^2(q^2+|q|^2)\ov{\Jk}\c\Xh\c\ov{\Jk}.
\end{split}
\eeaa

Since $q=r+O(1)$, $\ov{q}=r+O(1)$, $-\ov{q}^2+|q|^2=O(r)$, and
\beaa
-\frac{1}{3}q^2\ov{q}^2 -\frac{1}{3}q\ov{q}^3+\frac{2}{3}\ov{q}^4 =O(r^3),
\eeaa
we infer
\bea\lab{eq:improvmentwithGoodtermsinMext:2}
\bsplit
\Hc &= \Zc +O(r^{-1})\Zc+ O(1)B+O(r^{-2})\Xh+r^{-1}\Go_{k_*},\\
\widecheck{\DD\cos\th}  &= -\frac{i}{2}r\,\ov{\Jk}\c\Xh+O(r^{-1})\Xh+r^{-1}\Go_{k_*},\\
r\ov{\DD}\c\Zc   &= -2r\,\ov{\Pc}+O(r^{-1})\Zc+O(r^{-1})\dkb^{\leq 1}\Xh+O(1)\Pc+O(1)\dkb^{\leq 1}B\\
& +r^{-2}\Go_{k_*-1}.
\end{split}
 \eea

{\bf Step 2.} Recall the linearized Codazzi for $\Xh$
\beaa
\frac{1}{2}\ov{\DD}\c\Xh &=& \frac{1}{\ov{q}}\Zc -B +O(r^{-2})\Xh +O(r^{-2})\Hc +O(r^{-1})\dkb^{\leq 1}\trXc\\
&& +O(r^{-2})\widecheck{\DD(\cos\th)}+\Ga_b\c\Ga_g.
\eeaa
In view of  the assumptions {\bf Ref 1}  and   the estimate for $\trXc$ in    \eqref{eq:improvmentwithGoodtermsinMext:1}, we obtain
\beaa
\ov{\DD}\c\Xh&=& \frac{2}{r}\Zc - 2B +O(r^{-2})\Zc+O(r^{-2})\Xh +O(r^{-2})\Hc +O(r^{-2})\widecheck{\DD(\cos\th)}\\
&& +r^{-3}\Go_{k_*-1}.
\eeaa
Eliminating   $\Hc$ and  $\widecheck{\DD(\cos\th)}$ on the RHS  with the help of  the two first equations of \eqref{eq:improvmentwithGoodtermsinMext:2}, we infer
\bea\lab{eq:improvmentwithGoodtermsinMext:3}
\ov{\DD}\c\Xh = \frac{2}{r}\Zc - 2B +O(r^{-2})\Zc+O(r^{-2})\Xh +O(r^{-2})B +r^{-3}\Go_{k_*-1}.
\eea

{\bf Step 3.} Starting  with the  Bianchi identity 
\beaa
 \nab_3A -\frac{1}{2}\DD\hot B &=& -\frac{1}{2}\tr\Xb A+4\omb A +\frac{1}{2}(Z+4H)\hot B -3\ov{P}\Xh,
\eeaa
we have
\beaa
 \DD\hot B &=&  2\nab_3A +O(r^{-1}) A +O(r^{-2}) B +O(r^{-3})\Xh+r^{-1}\Ga_g\c\Ga_g+\Ga_b\c (A,B)
\eeaa
and hence, using the estimates for $\nab_3 A$ and $A$ in  \eqref{eq:improvmentwithGoodtermsinMext:1}, as well as {\bf Ref 1} for the nonlinear  terms, we deduce
\bea\lab{eq:improvmentwithGoodtermsinMext:4}
 \DD\hot B &=&  O(r^{-2}) B+O(r^{-3})\Xh+r^{-3-\de'}\Go_{k_*-1}. 
\eea

{\bf Step 4.} Starting with the  linearized Bianchi identity 
\beaa
\nab_3B -\DD\ov{\Pc} &=& \frac{2}{r}B+O(r^{-2}) B +O(r^{-2})\Pc+O(r^{-3})\Hc +O(r^{-4})\widecheck{\DD(\cos\th)}\\
&&+r^{-1}\Ga_b\c\Ga_g,
\eeaa
we  deduce
\bea
\lab{eq:Bianchi-Pc-step4-1}
\bsplit
 \DD\ov{\Pc} &= \nab_3B- \frac{2}{r}B+O(r^{-2}) B +O(r^{-2})\Pc+O(r^{-3})\Hc +O(r^{-4})\widecheck{\DD(\cos\th)} \\&+r^{-1}\Ga_b\c\Ga_g.
 \end{split}
\eea
Using the identity 
\beaa
\Lieb_\T &=& \frac{1}{2}\nab_3 +\frac{1}{2}\frac{\De}{|q|^2}\nab_4 +O(r^{-2})\dkb^{\leq 1} +\Ga_b,
\eeaa
we write
\beaa
\nab_3B &=& 2\Lieb_\T B - \frac{\De}{|q|^2}\nab_4B+O(r^{-2})\dkb^{\leq 1}B +r^{-1}\Ga_b\c\Ga_g.
\eeaa
Also, using the linearized Bianchi identity
\beaa
\nab_4B +\frac{4}{\ov{q}} B &=&  \frac{1}{2}\ov{\DD}\c A +\frac{aq}{2|q|^2}\ov{\Jk}\c A+\Ga_g\c(B,A),
\eeaa
we have
\beaa
\nab_4B &=& -\frac{4}{r}B+O(r^{-2})B+O(r^{-1})\dkb^{\leq 1}A+r^{-1}\Ga_g\c\Ga_g.
\eeaa
Hence
\beaa
\nab_3B &=& \frac{4}{r}B+ 2\Lieb_\T B  +O(r^{-2})\dkb^{\leq 1}B +O(r^{-1})\dkb^{\leq 1}A+r^{-1}\Ga_b\c\Ga_g.
\eeaa
Back to \eqref{eq:Bianchi-Pc-step4-1} we substitute $\nab_3 B$  to  deduce
\beaa
 \DD\ov{\Pc} &=& \frac{2}{r}B  + 2\Lieb_\T B  +O(r^{-2})\dkb^{\leq 1} B +O(r^{-2})\Pc+O(r^{-3})\Hc +O(r^{-4})\widecheck{\DD(\cos\th)}\\
 &&  +O(r^{-1})\dkb^{\leq 1}A+r^{-1}\Ga_b\c\Ga_g.
\eeaa
In view of  \eqref{eq:improvmentwithGoodtermsinMext:1} for $A$ and  the control of $\Lieb_\T B$ in  Proposition \ref{prop:improveestiatesforsomeLiebTderivativesinMext},  we deduce
\beaa
 \DD\ov{\Pc} &=& \frac{2}{r}B    +O(r^{-2})\dkb^{\leq 1} B +O(r^{-2})\Pc+O(r^{-3})\Hc +O(r^{-4})\widecheck{\DD(\cos\th)} \\
 &&+r^{-3-\de'}\Go_{k_*-3}.
\eeaa
Using  the first two  equations of \eqref{eq:improvmentwithGoodtermsinMext:2} to eliminate 
$\Hc$ and  $\widecheck{\DD(\cos\th)}$
\bea\lab{eq:improvmentwithGoodtermsinMext:5}
\bsplit
 \DD\ov{\Pc} &= \frac{2}{r}B    +O(r^{-2})\dkb^{\leq 1} B +O(r^{-2})\Pc+O(r^{-3})\Zc +O(r^{-4})\Xh\\
 &+r^{-3-\de'}\Go_{k_*-3}.
 \end{split}
\eea

{\bf Step 5.} We start with  equation \eqref{eq:improvmentwithGoodtermsinMext:4}  of Step  3.
\beaa
 \DD\hot B &=&  O(r^{-2}) B+O(r^{-3})\Xh+r^{-3-\de'}\Go_{k_*-1}. 
\eeaa
We appeal again to  Corollary \ref{cor:formuallikningnabprimeandnabintransfoformula:bis} to write
\bea
\lab{eq:nabtonab'}
\nab' &=& \Big(1+O(r^{-2})\Big)\nab+O(r^{-1})\Lieb_\T +O(r^{-3})+r^{-1}\Ga_b\dk^{\leq 1}.
\eea
Together with the control of $\Lieb_\T B$ provided by Proposition \ref{prop:improveestiatesforsomeLiebTderivativesinMext}, we infer
\beaa
 \DD'\hot B &=&  O(r^{-2})(\dkb')^{\leq 1} B+O(r^{-3})\Xh+r^{-3-\de'}\Go_{k_*-3}. 
\eeaa
 We appeal to the first elliptic estimate of Corollary \ref{Lemma:Hodge-estimateDDs_2}  to derive for $k\le k_*-3$,
\beaa
\|(\dkb')'^{\le k+1} B\|_{L^2(S(u,r))} &\les& r^2|(\ov{\DD'}\c B)_{\ell=1}|+
r^{-1}\|(\dkb')^{\leq k+1} B\|_{L^2(S(u,r))}\\
&&+r^{-2}\|(\dkb')^{\leq k}\Xh\|_{L^2(S(u,r))}+\ep_0 r^{-1- \de'}u^{-1-\dec}.
\eeaa
Using again \eqref{eq:nabtonab'} and  the control of $\Lieb_\T B$  as above  we deduce
\beaa
\|\dkb^{\leq k+1 }B\|_{L^2(S(u,r))} &\les& r^2|(\ov{\DD}\c B)_{\ell=1}|+
r^{-1}\|\dkb^{\leq k+1 }B\|_{L^2(S(u,r))}+r^{-2}\|\dkb^{\leq k }\Xh\|_{L^2(S(u,r))}\\
&&+\ep_0 r^{-1-\de'}u^{-1-\dec}.
\eeaa
Since $r\geq r_0$ on $\Mext$, we infer, for $r_0$ large enough,  for all $k\le k_*-3$,
\bea\lab{eq:improvmentwithGoodtermsinMext:6}
\bsplit
\|\dkb^{\leq k+1 }B\|_{L^2(S(u,r))} &\les r^2|(\ov{\DD}\c B)_{\ell=1}| +r^{-2}\|\dkb^{\leq k }\Xh\|_{L^2(S(u,r))}
+\ep_0 r^{-1-\de'}u^{-1-\dec}.
\end{split}
\eea

{\bf Step 6.}  We start with  the last equation of \eqref{eq:improvmentwithGoodtermsinMext:2}
\beaa
r\ov{\DD}\c\Zc   &=& -2r\,\ov{\Pc}+O(r^{-1})\Zc+O(r^{-1})\dkb^{\leq 1}\Xh+O(1)\Pc+O(1)\dkb^{\leq 1}B +r^{-2}\Go_{k_*-1}.
 \eeaa
Differentiating w.r.t. $r\DD$, we infer
\beaa
r^2\DD\ov{\DD}\c\Zc  & =& -2r^2\DD\,\ov{\Pc}+O(r^{-1})\dkb^{\leq 1}\Zc+O(r^{-1})\dkb^{\leq 2}\Xh+O(1)\dkb^{\leq 1}\Pc+O(1)\dkb^{\leq 2}B \\
&&+r^{-2}\Go_{k_*-2}.
 \eeaa
On the other hand, recall \eqref{eq:improvmentwithGoodtermsinMext:5}, 
\beaa
 \DD\ov{\Pc} = \frac{2}{r}B    +O(r^{-2})\dkb^{\leq 1} B +O(r^{-2})\Pc+O(r^{-3})\Zc +O(r^{-4})\Xh+r^{-3-\de'}\Go_{k_*-3}.
\eeaa
Plugging in the previous identity, we infer
\beaa
r^2\DD\ov{\DD}\c\Zc   &=& O(r)\dkb^{\leq 2}B +O(r^{-1})\dkb^{\leq 1}\Zc+O(r^{-1})\dkb^{\leq 2}\Xh+O(1)\dkb^{\leq 1}\Pc \\
&&+r^{-1-\de'}\Go_{k_*-3}.
 \eeaa
Using again formula  \eqref{eq:nabtonab'}  to pass to the prime frame, as well as the  control of $\Lieb_\T\Zc$ provided by Proposition \ref{prop:improveestiatesforsomeLiebTderivativesinMext}, we deduce
\beaa
r^2\DD'\ov{\DD'}\c\Zc   &=& O(r)\dkb^{\leq 2}B +O(r^{-1})(\dkb')^{\leq 2}\Zc+O(r^{-1})\dkb^{\leq 2}\Xh+O(1)\dkb^{\leq 1}\Pc \\&&+r^{-1-\de'}\Go_{k_*-3}.
 \eeaa
Using the second elliptic estimate of Corollary \ref{Lemma:Hodge-estimateDDs_2},  we deduce for all $k\le k_*-3$,
\beaa
\|(\dkb') ^{\leq k+2}\Zc\|_{L^2(S(u,r))} &\les& 
r \| \dkb^{\leq k+2}  B\|_{L^2(S(u,r))}+r^{-1}\|( \dkb')^{\leq k+2} \Zc\|_{L^2(S(u,r))}\\
&&+r^{-1}\|\dkb^{\leq k+2 }  \Xh\|_{L^2(S(u,r))}+\|\dkb^{\leq k+1} \Pc\|_{L^2(S(u,r))}+\ep_0 r ^{-\de'}u^{-1-\dec}.
\eeaa
Using again formula  \eqref{eq:nabtonab'},   to pass back to the un-prime frame, 
as well as the control of $\Lieb_\T\Zc$ provided by Proposition \ref{prop:improveestiatesforsomeLiebTderivativesinMext},
 we deduce for all $k\le k_*-3$,
 \beaa
\|\dkb ^{\leq k+2}\Zc\|_{L^2(S(u,r))} &\les& 
r\|\dkb^{\leq k+2} B\|_{L^2(S(u,r))}+r^{-1}\| \dkb ^{\leq k+2}\Zc\|_{L^2(S(u,r))}\\
&&+r^{-1}\|\dkb^{\leq k+2}\Xh\|_{L^2(S(u,r))}+\|\dkb^{\leq k+1}\Pc\|_{L^2(S(u,r))}+\ep_0 r ^{-\de'}u^{-1-\dec}.
\eeaa
 Since $r\geq r_0$ on $\Mext$, we infer, for $r_0$ large enough, and for $k\le k_*-3$,
\bea\lab{eq:improvmentwithGoodtermsinMext:7}
\bsplit
\|\dkb^{\leq k+2}\Zc\|_{L^2(S(u,r))} &\les 
r\|\dkb^{\leq k+2}B\|_{L^2(S(u,r))}+r^{-1}\|\dkb^{\leq k+2}\Xh\|_{L^2(S(u,r))}\\
&+\|\dkb^{\leq k+1}\Pc\|_{L^2(S(u,r))}
+\ep_0r^{-\de'}u^{-1-\dec}.
\end{split}
\eea
 
{\bf Step  7. } 
So far we have established for $k\le k_*-3$, see   \eqref{eq:improvmentwithGoodtermsinMext:6}    and  \eqref{eq:improvmentwithGoodtermsinMext:7} , 
\beaa
\bsplit
\|\dkb^{\leq k+1 }B\|_{L^2(S(u,r))} &\les r^2|(\ov{\DD}\c B)_{\ell=1}| +r^{-2}\|\dkb^{\leq k }\Xh\|_{L^2(S(u,r))}
+\ep_0 r^{-1-\de'}u^{-1-\dec}.\\
\|\dkb^{\leq k+2}\Zc\|_{L^2(S(u,r))} &\les 
r\|\dkb^{\leq k+2}B\|_{L^2(S(u,r))}+r^{-1}\|\dkb^{\leq k+2}\Xh\|_{L^2(S(u,r))}\\
&+\|\dkb^{\leq k+1}\Pc\|_{L^2(S(u,r))}
+\ep_0r^{-\de'}u^{-1-\dec}.
\end{split}
\eeaa
Combining them we  derive, for $k\le k_*-4$, 
\beaa
\|\dkb^{\leq k+2}\Zc\|_{L^2(S(u,r))} &\les &
 r^3|(\ov{\DD}\c B)_{\ell=1}|  +r^{-1}\|\dkb^{\leq k+2}\Xh\|_{L^2(S(u,r))}+\|\dkb^{\leq k+1}\Pc\|_{L^2(S(u,r))}
\\
&&+\ep_0r^{-\de'}u^{-1-\dec}.
\eeaa
Therefore, for any $S=S(u, r)$ in $\Mext$
\bea
\lab{eq:improvmentwithGoodtermsinMext:6-7}
\bsplit
  r\|\dkb^{\leq k_*-2}B\|_{L^2(S)}+\|\dkb^{\leq k_*- 2}\Zc\|_{L^2(S) }
 \les& r^3|(\ov{\DD}\c B)_{\ell=1}|+r^{-1}\|\dkb^{\leq k_*-1}\Xh\|_{L^2(S)}\\
 &+\|\dkb^{\leq k_*-2}\Pc\|_{L^2(S)}+\ep_0 r^{-\de'} u^{-1-\dec}.
 \end{split}
\eea
It thus remain to eliminate  the terms in $\Pc, \Xh$ on the RHS.

{\bf Step 8.} Using the last equation of \eqref{eq:improvmentwithGoodtermsinMext:2}, i.e. 
\beaa
r\ov{\DD}\c\Zc   &=& -2r\,\ov{\Pc}+O(r^{-1})\Zc+O(r^{-1})\dkb^{\leq 1}\Xh+O(1)\Pc+O(1)\dkb^{\leq 1}B +r^{-2}\Go_{k_*-1},
 \eeaa
 and $r\ge  r_0$ sufficiently large, we have   for all $k\le k_*-1$     
\beaa
\nn\|\dkb^{\leq k-1}\Pc\|_{L^2(S(u,r))} &\les& 
r^{-1}\|\dkb^{\leq k}\Zc\|_{L^2(S(u,r))}+r^{-2}\|\dkb^{\leq k}\Xh\|_{L^2(S(u,r))}\\
&&+r^{-1}\|\dkb^{\leq k }B\|_{L^2(S(u,r))}
+\ep_0 r^{-2} u^{-1-\dec}.
\eeaa
Also, recall \eqref{eq:improvmentwithGoodtermsinMext:5} 
\beaa
 \DD\ov{\Pc} = \frac{2}{r}B    +O(r^{-2})\dkb^{\leq 1} B +O(r^{-2})\Pc+O(r^{-3})\Zc +O(r^{-4})\Xh+r^{-3-\de'}\Go_{k_*-3}.
\eeaa
Together with the previous estimate, we infer for all $k\le k_*-3$,
\beaa
\bsplit
\|\dkb^{\leq k}\Pc\|_{L^2(S(u,r))} &\les
r^{-1}\|\dkb^{\leq k}\Zc\|_{L^2(S(u,r))}+\|\dkb^{\leq k}B\|_{L^2(S(u,r))}\\
&+r^{-2}\|\dkb^{\leq k}\Xh\|_{L^2(S(u,r))}
+   \ep_0 r^{-1-\de'} u^{-1-\dec}.     
\end{split}
\eeaa
Thus, for all $S=S(u, r)\subset\Mext$,
\bea\lab{eq:improvmentwithGoodtermsinMext:8}
\bsplit
\|\dkb^{\leq k_*-3}\Pc\|_{L^2(S)} &\les
r^{-1}\|\dkb^{\leq k_*-3}\Zc\|_{L^2(S)}+\|\dkb^{\leq k_*-3}B\|_{L^2(S)}\\
&+r^{-2}\|\dkb^{\leq k_*-3}\Xh\|_{L^2(S(u,r))}
+  \ep_0 r^{-1-\de'} u^{-1-\dec}.
\end{split}
\eea

{\bf Step 9.} So far  we have established, see \eqref{eq:improvmentwithGoodtermsinMext:6-7} and \eqref{eq:improvmentwithGoodtermsinMext:8},
\beaa
\bsplit
  r\|\dkb^{\leq k_*-3}B\|_{L^2(S)}+\|\dkb^{\leq k_*- 3}\Zc\|_{L^2(S) }
 \les& r^3|(\ov{\DD}\c B)_{\ell=1}|+r^{-1}\|\dkb^{\leq k_*-2}\Xh\|_{L^2(S)}\\
 &+\|\dkb^{\leq k_*-3}\Pc\|_{L^2(S)}+\ep_0 r^{-\de'} u^{-1-\dec}\\
  \|\dkb^{\leq k_*-3}\Pc\|_{L^2(S)} &\les
r^{-1}\|\dkb^{\leq k_*-3}\Zc\|_{L^2(S)}+\|\dkb^{\leq k_*-3}B\|_{L^2(S)}\\
&+r^{-2}\|\dkb^{\leq k_*-3}\Xh\|_{L^2(S(u,r))}
+  \ep_0 r^{-1-\de'} u^{-1-\dec}.
 \end{split}
\eeaa
Combining we deduce
\beaa
&&r\|\dkb^{\leq k_*-3}B\|_{L^2(S)}+\|\dkb^{\leq k_*-3}\Zc\|_{L^2(S)}\\
&&\les r^3|(\ov{\DD}\c B)_{\ell=1}|+r^{-1}\|\dkb^{\leq k_*-2}\Xh\|_{L^2(S)}+\|\dkb^{\leq k_*-3}\Pc\|_{L^2(S)}+\ep_0 r^{-\de'} u^{-1-\dec}\\
&&\les r^3|(\ov{\DD}\c B)_{\ell=1}|+r^{-1}\|\dkb^{\leq k_*-2}\Xh\|_{L^2(S)}+\ep_0 r^{-\de'} u^{-1-\dec}\\
&&+ r^{-1}\|\dkb^{\leq k_*-3}\Zc\|_{L^2(S)}+\|\dkb^{\leq k_*-3}B\|_{L^2(S)}+r^{-2}\|\dkb^{\leq k_*-3}\Xh\|_{L^2(S)}+  \ep_0 r^{-1} u^{-1-\dec}.  
\eeaa
Hence, for $r\ge r_0$ sufficiently large, we derive
\beaa
r\|\dkb^{\leq k_*-3}B\|_{L^2(S)}+\|\dkb^{\leq k_*-3}\Zc\|_{L^2(S)}  &\les& r^3|(\ov{\DD}\c B)_{\ell=1}|+r^{-1}\|\dkb^{\leq k_*-2}\Xh\|_{L^2(S)}+\ep_0 r^{-\de'} u^{-1-\dec}.
\eeaa
Back to \eqref{eq:improvmentwithGoodtermsinMext:8}, we deduce
 \beaa
 \|\dkb^{\leq k_*-3}\Pc\|_{L^2(S))} &\les&r^2|(\ov{\DD}\c B)_{\ell=1}|+r^{-2}\|\dkb^{\leq k_*-2}\Xh\|_{L^2(S))}+\ep_0 r^{-1-\de'} u^{-1-\dec}.
\eeaa
Combining the two estimates, we infer
\bea
\lab{eq:combinedB-P-Z}
\nn && r\Big(\|\dkb^{\leq k_*-3}B\|_{L^2(S)} +|\dkb^{\leq k_*-3}P\|_{L^2(S)} \Big)     +\|\dkb^{\leq k_*-3}\Zc\|_{L^2(S)}\\
 &\les& r^3|(\ov{\DD}\c B)_{\ell=1}|+r^{-1}\|\dkb^{\leq k_*-2}\Xh\|_{L^2(S)}+\ep_0 r^{-\de'} u^{-1-\dec}.
\eea
It thus remains to estimate $\Xh$.

{\bf Step 10.} Recall  \eqref{eq:improvmentwithGoodtermsinMext:3}
\beaa
\ov{\DD}\c\Xh = \frac{2}{r}\Zc - 2B +O(r^{-2})\Zc+O(r^{-2})\Xh +O(r^{-2})B +r^{-3}\Go_{k_*-1}.
\eeaa
 Using  formula  \eqref{eq:nabtonab'} to pass to the prime frame, as well as the control of $\Lieb_\T\Xh$ provided by Proposition \ref{prop:improveestiatesforsomeLiebTderivativesinMext}, we deduce
\beaa
r\ov{\DD'}\c\Xh = 2\Zc - 2rB +O(r^{-1})\Zc+O(r^{-1})\Xh +O(r^{-1})B +r^{-2}\Go_{k_*-2}.
\eeaa
Using the third elliptic estimate of Corollary \ref{Lemma:Hodge-estimateDDs_2}, we infer
\beaa
\|(\dkb)'^{\leq k_*-2}\Xh\|_{L^2(S(u,r))} &\les& 
\|(\dkb)^{\leq k_*-3}\Zc\|_{L^2(S(u,r))}+r\|(\dkb)^{\leq k_*-3}B\|_{L^2(S(u,r))}\\
&&+r^{-1}\|(\dkb)^{\leq k_*-3}\Xh\|_{L^2(S(u,r))}+\ep_0 r^{-1} u^{-1-\dec}.
\eeaa
Passing back to the un-primed frame, and using again the  control of $\Lieb_\T\Xh$ provided by Proposition \ref{prop:improveestiatesforsomeLiebTderivativesinMext}, we deduce
\beaa
\|(\dkb)^{\leq k_*-2}\Xh\|_{L^2(S(u,r))} &\les& 
\|(\dkb)^{\leq k_*-3}\Zc\|_{L^2(S(u,r))}+r\|(\dkb)^{\leq k_*-3}B\|_{L^2(S(u,r))}\\
&&+r^{-1}\|(\dkb)^{\leq k_*-3}\Xh\|_{L^2(S(u,r))} +\ep_0     r^{-1} u^{-1-\dec}.
\eeaa
Since $r\geq r_0$ on $\Mext$, we infer, for $r_0$ large enough, 
\bea\lab{eq:improvmentwithGoodtermsinMext:9}
\bsplit
\|\dkb^{\leq k_*-2}\Xh\|_{L^2(S(u,r))} &\les
\|\dkb^{\leq k_*-3}\Zc\|_{L^2(S(u,r))}+r\|\dkb^{\leq k_*-3}B\|_{L^2(S(u,r))}\\
&+\ep_0     r^{-1} u^{-1-\dec}.
\end{split}
\eea

{\bf Step 11.}
We  combine \eqref{eq:combinedB-P-Z} with \eqref{eq:improvmentwithGoodtermsinMext:9}. Thus for any $S=S(u, r)\subset\Mext$, we have
 \beaa
&& r\Big(\|\dkb^{\leq k_*-3}B\|_{L^2(S)} +|\dkb^{\leq k_*-3}P\|_{L^2(S)} \Big)     +\|\dkb^{\leq k_*-3}\Zc\|_{L^2(S)}\\
&\les& r^3|(\ov{\DD}\c B)_{\ell=1}|+r^{-1}\|\dkb^{\leq k_*-2}\Xh\|_{L^2(S)}+\ep_0 r^{-\de'} u^{-1-\dec}\\
 &\les & r^3|(\ov{\DD}\c B)_{\ell=1}| + r^{-1}  \|\dkb^{\leq k_*-3}\Zc\|_{L^2(S)}+\|\dkb^{\leq k_*-3}B\|_{L^2(S)}+\ep_0 r^{-\de'} u^{-1-\dec}.
\eeaa
 Thus,  absorbing the terms in $B$ and $\Zc$ on the RHS, we obtain 
  \beaa
 r\Big(\|\dkb^{\leq k_*-3}B\|_{L^2(S)} +|\dkb^{\leq k_*-3}P\|_{L^2(S)} \Big)     +\|\dkb^{\leq k_*-3}\Zc\|_{L^2(S)} &\les&
  r^3|(\ov{\DD}\c B)_{\ell=1}|\\
  &&+\ep_0 r^{-\de'} u^{-1-\dec}.
 \eeaa
 and
 \beaa
 \|\dkb^{\leq k_*-2}\Xh\|_{L^2(S(u,r))} &\les &  r^3|(\ov{\DD}\c B)_{\ell=1}|+\ep_0 r^{-\de'} u^{-1-\dec}.
 \eeaa
 We infer
 \beaa
 && r\Big(\|\dkb^{\leq k_*-3}B\|_{L^2(S)} +|\dkb^{\leq k_*-3}P\|_{L^2(S)} \Big)     +\|\dkb^{\leq k_*-3}\Zc\|_{L^2(S)} +\|\dkb^{\leq k_*-2}\Xh\|_{L^2(S(u,r))}\\
 &\les&  r^3|(\ov{\DD}\c B)_{\ell=1}|+\ep_0 r^{-\de'} u^{-1-\dec}.
 \eeaa
 By  Sobolev, we  deduce
\beaa
 && r\big(|\dkb^{\leq k_*-5}B|+|\dkb^{\leq k_*-5}\Pc|\big)+|\dkb^{\leq k_*-5}\Xh|+|\dkb^{\leq k_*-5}\Zc|\\
 & \les & r^2|(\ov{\DD}\c B)_{\ell=1}|+\ep_0 r^{-1-\de'} u^{-1-\dec}.
\eeaa

It remains to estimate the derivatives in $\nab_4, \nab_3$. In view of    the equations
\beaa
\nab_4\Xh &=&  -A+O(r^{-1})\Xh+\Ga_g\c\Ga_g,\\
\nab_4\Zc   &=&  -B+O(r^{-1})\Zc + O(r^{-2})\Xh  +O(r^{-2})\trXc+\Ga_g\c\Ga_g,\\
\nab_4B  &=&  O(r^{-1})\dkb^{\leq 1}A +O(r^{-1})B+\Ga_g\c(B,A),\\
\nab_4\left(\Pc \right) &=& O(r^{-1})\dkb^{\leq 1}B+O(r^{-1})\Pc+O(r^{-3})\trXc+r^{-1}\Ga_g\c\Ga_g+\Ga_b\c A,
\eeaa
and the control of $A$ and $\trXc$ in \eqref{eq:improvmentwithGoodtermsinMext:1}, we infer from the previous estimate
\beaa
 && r\big(|(r\nab_4, \dkb)^{\leq k_*-5}B|+|(r\nab_4, \dkb)^{\leq k_*-5}\Pc|\big)+|(r\nab_4, \dkb)^{\leq k_*-5}\Xh|+|(r\nab_4, \dkb)^{\leq k_*-5}\Zc|\\
&\les& r^2|(\ov{\DD}\c B)_{\ell=1}|+\ep_0 r^{-1-\de'} u^{-1-\dec}.
\eeaa
Together with the  control of $\Lieb_\T B$, $\Lieb_\T\Pc$, $\Lieb_\T\Xh$ and $\Lieb_\T\Zc$ provided by Proposition \ref{prop:improveestiatesforsomeLiebTderivativesinMext}, we deduce
\beaa
 r\big(|(r\nab_4, \dkb, \Lieb_\T)^{\leq k_*-5}B|+|(r\nab_4, \dkb, \Lieb_\T)^{\leq k_*-5}\Pc|\big)\\
 +|(r\nab_4, \dkb, \Lieb_\T)^{\leq k_*-5}\Xh|+|(r\nab_4, \dkb, \Lieb_\T)^{\leq k_*-5}\Zc| &\les& r^2|(\ov{\DD}\c B)_{\ell=1}|+\ep_0 r^{-1-\de'} u^{-1-\dec}.
\eeaa
Finally,      since  $\nab_3  = 2\Lieb_\T -\Big(1+O(r^{-2})\Big)\nab_4+O(r^{-1})\nab+O(r^{-3})+\Ga_b$, see     Corollary \ref{cor:decompositionofnab3inLiebTnab4andlot}, 
we infer that
\beaa
 r\big(|\dk^{\leq k_*-5}B|+|\dk^{\leq k_*-5}\Pc|\big)+|\dk^{\leq k_*-5}\Xh|+|\dk^{\leq k_*-5}\Zc| &\les& r^2|(\ov{\DD}\c B)_{\ell=1}|+\ep_0 r^{-1-\de'} u^{-1-\dec}
\eeaa
 which is precisely the estimate \eqref{eq-prop:firstcontrolofBPcXhZcHcDDcosth1}  of Proposition  \ref{prop:firstcontrolofBPcXhZcHcDDcosth}.
 
 {\bf Step 12.}  We are now ready  to  prove the  estimate \eqref{eq-prop:firstcontrolofBPcXhZcHcDDcosth2}. Indeed, combining  the first two equations of \eqref{eq:improvmentwithGoodtermsinMext:2} with the estimate  \eqref{eq-prop:firstcontrolofBPcXhZcHcDDcosth1} proved in Step 11, we obtain 
\beaa
 |\dk^{\leq k_*-5}\Hc|+|\dk^{\leq k_*-5}\widecheck{\DD\cos\th}| &\les& r^2|(\ov{\DD}\c B)_{\ell=1}|+\ep_0 r^{-1} u^{-1-\dec}
\eeaa
 which is the estimate \eqref{eq-prop:firstcontrolofBPcXhZcHcDDcosth2}  of Proposition  \ref{prop:firstcontrolofBPcXhZcHcDDcosth}.
  
  {\bf Step 13.} It remains to prove the auxiliary estimate \eqref{eq:Prliminary-forB}, i.e.
  \beaa
|\dkb ^{\leq {k_*-2}}B| &\les& r |(\ov{\DD}\c B)_{\ell=1}|+      \ep r^{-3-\de'} u^{-1/2-\dec}.
\eeaa
To do that, we start with the equation
\beaa
\DD\hot B +(Z+4H)\hot B &=& \nab_3A+\left(\frac{1}{2}\tr\Xb-4\omb\right)A+3\ov{P}\Xh.
\eeaa
from which we deduce, making use of {\bf Ref 1-2}, for $k\leq k_*-1$, 
\beaa
\big\|\dk^k\DD\hot B\big\|_{L^2(S)}&\les & r^{-2} \|\dk^{\leq k}B\|_{L^2(S)} +\ep_0  r^{-3-\de'} u^{-\frac{1}{2}-\dec}+ r^{-3}\big\|\dk^{\leq k}\Xh\big\|_{L^2(S)} \\
&\les& r^{-2} \|\dk^{\leq k}B\|_{L^2(S)} +\ep_0  r^{-3-\de'} u^{-\frac{1}{2}-\dec}+ \ep r^{-4}u^{-\frac{1}{2}-\dec}
\eeaa
and hence, for $k\leq k_*-1$,
\beaa
\big\|\dk^k\DD\hot B\big\|_{L^2(S)}&\les &  r^{-2} \|\dk^{\leq k}B\|_{L^2(S)} +\ep   r^{-3-\de'} u^{-\frac{1}{2}-\dec}.
\eeaa
Using once more  \eqref{eq:nabtonab'}    to pass the  prime frame,   and the  estimates for $\Lieb_\T B$ provided by Proposition \ref{prop:improveestiatesforsomeLiebTderivativesinMext},  we deduce, for $k\leq k_*-1$,
\beaa
\big\|(\dk')^k\DD'\hot B\big\|_{L^2(S)} &\les& r^{-2} \|\dk^{\leq k}B\|_{L^2(S)} +\ep  r^{-3-\de'} u^{-\frac{1}{2}-\dec}.
\eeaa
Applying the first estimate of Corollary \ref{Lemma:Hodge-estimateDDs_2},  we deduce, for all $k\le k_*-1$,
\beaa
\|(\dkb' )^{\le k+1}  B \|_{L^2 (S)} &\les & r  \|(\dkb' )^{\le k} \DD' \hot B\|_{L^2 (S)}+ r^2  \Big| \big(\ov{ \DD'} \c  B\big)_{\ell=1}  \Big|\\
&\les & r^{-1}   \|\dk^{\leq k}B\|_{L^2(S)}  + r^2  \Big| \big(\ov{ \DD'} \c  B\big)_{\ell=1}  \Big| +\ep  r^{-2-\de'} u^{-\frac{1}{2}-\dec}.
\eeaa
Passing back to the un-prime frame and absorbing the  $B$ term on the RHS, we obtain, for all $k\le k_*-1$,
\beaa
\|\dkb^{\le k+ 1}  B \|_{L^2 (S)} &\les&   r^2  \Big| \big(\ov{ \DD} \c  B\big)_{\ell=1}  \Big| +\ep  r^{-2-\de'} u^{-\frac{1}{2}-\dec}.
\eeaa
Thus, by Sobolev,
 \beaa
\|\dkb^{\le k_*-2}  B \|_{L^\infty (S)} &\les&   r   \Big| \big(\ov{ \DD} \c  B\big)_{\ell=1}  \Big| +\ep  r^{-3-\de'} u^{-\frac{1}{2}-\dec}.
\eeaa
The estimates for $\nab_4, \nab_3$ derivatives can then be recovered as in Step 11.  Thus, we finally obtain 
 \beaa
|\dk ^{\leq {k_*-2}}B| &\les& r |(\ov{\DD}\c B)_{\ell=1}|+      \ep r^{-3-\de'} u^{-\frac{1}{2}-\dec}
\eeaa
as stated in \eqref{eq:Prliminary-forB}. This concludes the proof of Proposition  \ref{prop:firstcontrolofBPcXhZcHcDDcosth}.
\end{proof}

%%%%%%%%%%%%%%%%%%%%%%%%%%%%%%%%%%%%%%%%%%%%%%%%%%%%%%%%%

\subsection{$O(u^{-1-\dec})$ type decay estimates for $B$, $\Pc$, $\Xh$, $\Zc$, $\Hc$, $\widecheck{\DD\cos\th}$}

%%%%%%%%%%%%%%%%%%%%%%%%%%%%%%%%%%%%%%%%%%%%%%%%%%%%%%%%%

The goal of this section is to prove the following proposition.
\begin{proposition}
\lab{Prop:estimatesforBPXhZcHc}
We have on $\Mext$
\bea
\lab{eq:Prop.estimatesforBPXhZcHc1}
\bsplit
 r\big(|\dk^{\leq k_*-5}B|+|\dk^{\leq k_*-5}\Pc|\big)+|\dk^{\leq k_*-5}\Xh|+|\dk^{\leq k_*-5}\Zc| &\les\ep_0 r^{-1-\de'} u^{-1-\dec}.,\\
|\dk^{\leq k_*-5}\Hc|+|\dk^{\leq k_*-5}\widecheck{\DD\cos\th}| &\les\ep_0 r^{-1} u^{-1-\dec},
\\
|(\ov{\DD}\c B)_{\ell=1}| &\le \ep_0 r^{-4-\de'} u^{-1-\dec}.
\end{split}
\eea
Also, we have on $\Mext$
\bea
\lab{eq:Prop.estimatesforBPXhZcHc2}
r|\dk^{\leq k_*-6}\nab_3B|+r|\dk^{\leq k_*-6}\nab_3\Pc|+|\dk^{\leq k_*-6}\nab_3\Xh|+|\dk^{\leq k_*-6}\nab_3\Zc| \les\ep_0 r^{-2} u^{-1-\dec}.
\eea
\end{proposition}

\begin{remark}\lab{rmk:takingadvantageofthegaininrtoavoidlogloss}
The above $r^{-\de'}$ gain for $B$, $\Pc$, $\Xh$ and $\Zc$ in \eqref{eq:Prop.estimatesforBPXhZcHc1} will be crucial to  derive $O(u^{-1-\dec})$ type decay estimates without log-loss in $r$ in particular for $\trXbc$, $\Bb$, $\widecheck{\DD u}$, $\widecheck{\D\vphi}$, $\DD\hot\Jk$ and $\widecheck{\ov{\DD}\c\Jk}$. 
\end{remark}

\begin{proof}
Recall from Proposition \ref{prop:firstcontrolofBPcXhZcHcDDcosth} that we have on $\Mext$
\bea
\lab{eq:result-prop.firstcontrolofBPcXhZcHcDDcosth}
\bsplit
 r\big(|\dk^{\leq k_*-5}B|+|\dk^{\leq k_*-5}\Pc|\big)+|\dk^{\leq k_*-5}\Xh|+|\dk^{\leq k_*-5}\Zc| &\les r^2|(\ov{\DD}\c B)_{\ell=1}|\\
 &+\ep_0 r^{-1-\de'} u^{-1-\dec},
\\
 |\dk^{\leq k_*-5}\Hc|+|\dk^{\leq k_*-5}\widecheck{\DD\cos\th}| &\les r^2|(\ov{\DD}\c B)_{\ell=1}|\\
 &+\ep_0 r^{-1} u^{-1-\dec}.
 \end{split}
\eea
Thus, to conclude the control of $B$, $\Pc$, $\Xh$, $\Zc$, $\Hc$ and $\widecheck{\DD\cos\th}$ in \eqref{eq:Prop.estimatesforBPXhZcHc1}, we need to control $(\ov{\DD}\c B)_{\ell=1}$, which is  the focus of Steps 1 and 2 below. Then, \eqref{eq:Prop.estimatesforBPXhZcHc2} is proved in Step 3.

{\bf Step 1.}   We  derive  transport  equations in $e_4$ for the following scalars $h^{(p)}$, for $p=0,+,-$,  
\bea
\bsplit
h^{(p)}&=\left(\int_{S(r,u)} \frac{rJ^{(p)}}{\Si}   [\ov{\DD}\c]_{ren}    \left( r^4[B]_{ren} \right)    \right)  \\
&=\int_{S(r,u)} \frac{rJ^{(p)}}{\Si}\left(\ov{\DD}\c -\frac{a}{2}\ov{\Jk}\c\nab_4 -\frac{a}{2}\ov{\Jk}\c\nab_3\right)\left\{r^4\left(B - \frac{3a}{2}\ov{\Pc}\Jk -\frac{a}{4}\ov{\Jk}\c A\right)\right\}.
\end{split}
\eea
To this end, we rely on  the crucial  identities  \eqref{eq:Proposition.Identity.ell=1Div B1} and  \eqref{eq:Proposition.Identity.ell=1Div B2}  of Proposition  \ref{Proposition:Identity.ell=1Div B},  which we rewrite below in the following form
\bea
\lab{eq:Proposition.Identity.ell=1Div B1-bis} 
\bsplit
e_4(h^{(0)}) &=  O(1)\dk^{\leq 1}\Xh+O(r)\dk^{\leq 2}B+O(r^2)\dk^{\leq 1}\nab_3B+O(r)\dk^{\leq 2}\Pc\\
& +O(1)\dk^{\leq 1}\trXc
+O(r^{2})\dk^{\leq 1}\nab_3A+O(r)\dk^{\leq 2}A+r^4\dk^{\leq 1}\big(\Ga_g\c(B,A)\big)\\
& +r^4\dk^{\leq 1}\big(\Ga_b\c\nab_3 A\big)
+r^2\dk^{\leq 2}\big(\Ga_g\c\Ga_g\big)
+O\left(\frac{r^2\ep}{u^{\frac 1 2 +\dec}}\right)\dk^{\leq 1}\Big(A, B, r^{-1}\Pc\Big),
\end{split}
\eea
and
\bea
\lab{eq:Proposition.Identity.ell=1Div B2-bis} 
\bsplit
e_4(h^{(\pm)}) \mp\frac{a}{r^2}   h^{(\mp)}   &=O(1)\dk^{\leq 1}\Xh+O(r)\dk^{\leq 2}B+O(r^2)\dk^{\leq 1}\nab_3B+O(r)\dk^{\leq 2}\Pc\\
& +O(1)\dk^{\leq 1}\trXc
+O(r^{2})\dk^{\leq 1}\nab_3A+O(r)\dk^{\leq 2}A\\
&+r^4\dk^{\leq 1}\big(\Ga_g\c(B,A)\big) +r^4\dk^{\leq 1}\big(\Ga_b\c\nab_3 A\big)
+r^2\dk^{\leq 2}\big(\Ga_g\c\Ga_g\big)\\
&+O\left(\frac{r^2\ep}{u^{\frac 1 2 +\dec}}\right)\dk^{\leq 1}\Big(A, B, r^{-1}\Pc\Big),
\end{split}
\eea

In view of the definition of $(\ov{\DD}\c B)_{\ell=1}$ and $h^{(p)}$, and since $\Si=r^2+O(1)$, we have
\bea\lab{eq:relationbetweenh0hplushminusandell=1modeofovDDcB}
(h^{(0)}, h^{(+)}, h^{(-)}) &=& r^5(\ov{\DD}\c B)_{\ell=1}+O(r^3)\dk^{\leq 1}(B, \Pc, A).
\eea
Together with the first equation in  \eqref{eq:result-prop.firstcontrolofBPcXhZcHcDDcosth}
we deduce
\beaa
&& r\big(|\dk^{\leq k_*-5}B|+|\dk^{\leq k_*-5}\Pc|\big)+|\dk^{\leq k_*-5}\Xh|+|\dk^{\leq k_*-5}\Zc|\\
&\les& r^{-3}|(h^{(0)}, h^{(+)}, h^{(-)})|+O(1)\dk^{\leq 1}(B, \Pc, A)+\ep_0 r^{-1-\de'} u^{-1-\dec}.
\eeaa
Since  $r\geq r_0$ on $\Mext$, we infer, for $r_0$ large enough,
\bea\lab{eq:controlofBPcXhZcfromh0hplushminus}
\bsplit
 &  r\big(|\dk^{\leq k_*-5}B|+|\dk^{\leq k_*-5}\Pc|\big)+|\dk^{\leq k_*-5}\Xh|+|\dk^{\leq k_*-5}\Zc||\\
&\les r^{-3}|(h^{(0)}, h^{(+)}, h^{(-)})|+O(1)\dk^{\leq 1}A+\ep_0 r^{-1-\de'} u^{-1-\dec}.
\end{split}
\eea
Similarly, from the second equation in \eqref{eq:result-prop.firstcontrolofBPcXhZcHcDDcosth}
\bea\lab{eq:controlofBPcXhZcfromh0hplushminus:bis}
\bsplit
|\dk^{\leq k_*-5}\Hc|+|\dk^{\leq k_*-5}\widecheck{\DD\cos\th}|
&\les r^{-3}|(h^{(0)}, h^{(+)}, h^{(-)})|+O(1)\dk^{\leq 1}A\\
&+\ep_0 r^{-1} u^{-1-\dec}.
\end{split}
\eea

Together with \eqref{eq:controlofBPcXhZcfromh0hplushminus} and \eqref{eq:controlofBPcXhZcfromh0hplushminus:bis},  we  deduce\footnote{To bound the terms $O(r^2)\dk^{\leq 1}\nab_3B$ on the RHS, we also use the following consequence of Bianchi 
\beaa
\nab_3B = O(r^{-1})\dkb^{\leq 1}\Pc+O(r^{-1})B+O(r^{-3})\Hc+O(r^{-4})\widecheck{\DD\cos\th}+r^{-1}\Ga_b\c\Ga_g.
\eeaa}  from \eqref{eq:Proposition.Identity.ell=1Div B1-bis} and  \eqref{eq:Proposition.Identity.ell=1Div B2-bis}
\beaa
 e_4\left(h^{(0)}\right) &=&  O(r^{-3})|(h^{(0)}, h^{(+)}, h^{(-)})| +O(1)\dk^{\leq 1}\trXc\\
&&+O(r^{2})\dk^{\leq 1}\nab_3A+O(r)\dk^{\leq 2}A+r^4\dk^{\leq 1}\big(\Ga_g\c(B,A)\big) +r^4\dk^{\leq 1}\big(\Ga_b\c\nab_3 A\big)\\
&&+r^2\dk^{\leq 2}\big(\Ga_g\c\Ga_g\big)+O\left(\frac{r^2\ep}{u^{\frac{1}{2}+\dec}}\right)\dk^{\leq 1}\Big(A, B, r^{-1}\Pc\Big)+O\big(\ep_0 r^{-1-\de'} u^{-1-\dec} \big),
\eeaa 
\beaa
 e_4\left(h^{(+)}\right) &=& \frac{a}{r^2}h^{(-)}+ O(r^{-3})|(h^{(0)}, h^{(+)}, h^{(-)})| +O(1)\dk^{\leq 1}\trXc\\
&&+O(r^{2})\dk^{\leq 1}\nab_3A+O(r)\dk^{\leq 2}A+r^4\dk^{\leq 1}\big(\Ga_g\c(B,A)\big) +r^4\dk^{\leq 1}\big(\Ga_b\c\nab_3 A\big)\\
&&+r^2\dk^{\leq 2}\big(\Ga_g\c\Ga_g\big)+O\left(\frac{r^2\ep}{u^{\frac{1}{2}+\dec}}\right)\dk^{\leq 1}\Big(A, B, r^{-1}\Pc\Big)+O\big(\ep_0 r^{-1-\de'} u^{-1-\dec} \big),
\eeaa 
and 
\beaa
 e_4\left(h^{(-)}\right) &=& -\frac{a}{r^2}h^{(+)}+ O(r^{-3})|(h^{(0)}, h^{(+)}, h^{(-)})| +O(1)\dk^{\leq 1}\trXc\\
&&+O(r^{2})\dk^{\leq 1}\nab_3A+O(r)\dk^{\leq 2}A+r^4\dk^{\leq 1}\big(\Ga_g\c(B,A)\big) +r^4\dk^{\leq 1}\big(\Ga_b\c\nab_3 A\big)\\
&& +r^2\dk^{\leq 2}\big(\Ga_g\c\Ga_g\big)+O\left(\frac{r^2\ep}{u^{\frac{1}{2}+\dec}}\right)\dk^{\leq 1}\Big(A, B, r^{-1}\Pc\Big)+O\big(\ep_0 r^{-1-\de'} u^{-1-\dec} \big).
\eeaa 
Using the control of $A$ in {\bf Ref 2},    the bootstrap assumptions {\bf Ref 1}, and    the control of $\trXc$ derived in Proposition \ref{lemma:firstcontroloftrXbcinMext}, we infer 
\beaa
 e_4\left(h^{(0)}\right) =  O(r^{-3})|(h^{(0)}, h^{(+)}, h^{(-)})| +O\left(\frac{r^2\ep}{u^{\frac{1}{2}+\dec}}\right)\dk^{\leq 1}B +O\big(\ep_0 r^{-1-\de'} u^{-1-\dec} \big),
\eeaa 
\beaa
 e_4\left(h^{(+)}\right) = \frac{a}{r^2}h^{(-)}+ O(r^{-3})|(h^{(0)}, h^{(+)}, h^{(-)})| +O\left(\frac{r^2\ep}{u^{\frac{1}{2}+\dec}}\right)\dk^{\leq 1}B+O\big(\ep_0 r^{-1-\de'} u^{-1-\dec} \big),
\eeaa 
and 
\beaa
 e_4\left(h^{(-)}\right) = -\frac{a}{r^2}h^{(+)}+ O(r^{-3})|(h^{(0)}, h^{(+)}, h^{(-)})| +O\left(\frac{r^2\ep}{u^{\frac{1}{2}+\dec}}\right)\dk^{\leq 1}B+O\big(\ep_0 r^{-1-\de'} u^{-1-\dec} \big).
\eeaa 
In order to estimate  the  term involving  $\dk^{\leq 1}B$ on the RHS, we rely on the  estimate \eqref{eq:Prliminary-forB} in Proposition \ref{prop:firstcontrolofBPcXhZcHcDDcosth}, i.e. 
\beaa
|\dk ^{\leq {k_*-2}}B| &\les& r |(\ov{\DD}\c B)_{\ell=1}|+      \ep r^{-3-\de'} u^{-1/2-\dec}.
\eeaa
Together with \eqref{eq:relationbetweenh0hplushminusandell=1modeofovDDcB}, we deduce 
\bea
\lab{eq:identitiesfor.h_123}
 e_4\left(h^{(p)}\right) =  O(r^{-2})|(h^{(0)}, h^{(+)}, h^{(-)})|  +O\big(\ep_0 r^{-1-\de'} u^{-1-\dec} \big),\quad p=0,+,-.
\eea

{\bf Step 2.} Using the estimate
\beaa
\sup_{\Si_*}r^5u^{1+\dec}\Big |  \left(  [\DD]_{ren} \c [B]_{ren} \right)_{\ell=1}\Big| &\les& \ep_0,
\eeaa
of Proposition  \ref{prop:improvedesitmatesfortemporalframeofMextonSigmastar-new}, and recalling that
\beaa
h^{(p)}&=& \int_{S(r,u)} \frac{rJ^{(p)}}{\Si}   [\ov{\DD}\c]_{ren}    \left( r^4[B]_{ren} \right), 
\eeaa
we easily deduce on $\Si_*$, using also $\Si=r^2+O(1)$ and the definition of $\ell=1$ modes, 
\beaa
|(h^{(0)}, h^{(+)}, h^{(-)})| &\les& \frac{\ep_0}{u^{1+\dec}}.
\eeaa
Thus, integrating the transport equations in $e_4$ of \eqref{eq:identitiesfor.h_123} from the last  slice  $\Si_*$,  we obtain on $\Mext$
\beaa
|(h^{(0)}, h^{(+)}, h^{(-)})| &\les& \frac{\ep_0}{u^{1+\dec}}.
\eeaa
Plugging  these in  \eqref{eq:controlofBPcXhZcfromh0hplushminus} and \eqref{eq:controlofBPcXhZcfromh0hplushminus:bis},  we deduce
\beaa
 r\big(|\dk^{\leq k_*-5}B|+|\dk^{\leq k_*-5}\Pc|\big)+|\dk^{\leq k_*-5}\Xh|+|\dk^{\leq k_*-5}\Zc| &\les& O(1)\dk^{\leq 1}A+O\big(\ep_0 r^{-1-\de'} u^{-1-\dec} \big),
\\
|\dk^{\leq k_*-5}\Hc|+|\dk^{\leq k_*-5}\widecheck{\DD\cos\th}| &\les& O(1)\dk^{\leq 1}A+O\big(\ep_0 r^{-1} u^{-1-\dec} \big).
\eeaa
Also, from \eqref{eq:relationbetweenh0hplushminusandell=1modeofovDDcB},
\beaa
|(\ov{\DD}\c B)_{\ell=1}| &\les& r^{-2}|\dk^{\leq 1}(B, \Pc, A)|+\frac{\ep_0}{r^5u^{1+\dec}}.
\eeaa
Thanks to the control of $A$ in {\bf Ref 2}, we  thus have obtained  on $\Mext$
\beaa
 r\big(|\dk^{\leq k_*-5}B|+|\dk^{\leq k_*-5}\Pc|\big)+|\dk^{\leq k_*-5}\Xh|+|\dk^{\leq k_*-5}\Zc| &\les&\ep_0 r^{-1-\de'} u^{-1-\dec},
\\
|\dk^{\leq k_*-5}\Hc|+|\dk^{\leq k_*-5}\widecheck{\DD\cos\th}| &\les& \ep_0 r^{-1} u^{-1-\dec} ,
\\
|(\ov{\DD}\c B)_{\ell=1}| &\les&\ep_0 r^{-4-\de'} u^{-1-\dec},
\eeaa
as stated in \eqref{eq:Prop.estimatesforBPXhZcHc1}.

{\bf Step 3.} It remains to prove  the estimate \eqref{eq:Prop.estimatesforBPXhZcHc2}
 for $\nab_3B$, $\nab_3\Pc$, $\nab_3\Xh$ and $\nab_3\Zc$. This follows easily by relying   on the formula,  see   Corollary \ref{cor:decompositionofnab3inLiebTnab4andlot},
\beaa
\nab_3  &=& 2\Lieb_\T +O(r^{-1})\dk^{\leq 1}+\Ga_b,
\eeaa
 the control of $\Lieb_\T(B, \Pc, \Xh, \Zc)$ provided by Proposition \ref{prop:improveestiatesforsomeLiebTderivativesinMext}, and the estimates  \eqref{eq:Prop.estimatesforBPXhZcHc1} derived in Step 2 above. This concludes the proof of Proposition \ref{Prop:estimatesforBPXhZcHc}.
\end{proof}

%%%%%%%%%%%%%%%%%%%%%%%%%%%%%%%%%%%%%%%%%%%%%%%%

\subsection{$O(u^{-\frac{1}{2}-\dec})$ type decay estimates for $B$, $\Pc$, $\Xh$ and $\Zc$}

%%%%%%%%%%%%%%%%%%%%%%%%%%%%%%%%%%%%%%%%%%%%%%%%

The goal of this section is to prove the following proposition.
\begin{proposition}
\lab{Proposition:estimatesforBPcXhZc-halfdecayinu}
We have on $\Mext$
\beaa
 r|\dk^{\leq k_*-6}\Pc|+|\dk^{\leq k_*-6}\Xh|+|\dk^{\leq k_*-6}\Zc| &\les&    \ep_0 r^{-2} u^{-1/2-\dec}    
\\
|\dk^{\leq k_*-6}B| &\les&    \ep_0 r^{-3-\de'} u^{-1/2-\dec}.    
\eeaa
\end{proposition}

\begin{proof}
As mentioned in  section \ref{section:MainResultsandStrategy}, the proof of these  estimates is much easier than the ones  derived so far. We sketch the main steps below. 
 
{\bf Step 1. } To get the  estimate for $\Xh$ we apply Proposition \ref{Prop:transportrp-f-Decay-knorms}
to the equation 
\beaa
\nab_4\Xh+\frac{2r}{|q|^2}\Xh &=& -A+\Ga_g\c\Ga_g.
\eeaa
Thus, using  the assumptions  {\bf Ref 1-2},  and the estimates for $\Xh$ on $\Si_*$ in Proposition \ref{prop:improvedesitmatesfortemporalframeofMextonSigmastar-new}, we derive, for all $r_0\leq r\leq r_*$ and $1\leq u\leq u_*$,
\beaa
r^2 \big|\dk^{\le k_*} \Xh (r, u)  \big|&\les& r_*^2  \big|\dk^{\le k_*}   \Xh(r_*, u)  \big|+ \ep_0  r^{-\de'} u^{-1/2-\dec}  \les \ep_0 u^{-1/2-\dec} 
\eeaa
as stated.

{\bf Step 2.} To get the  estimate for $B$ we  proceed exactly as in the proof for the auxiliary estimate \eqref{eq:Prliminary-forB} of Proposition \ref{prop:firstcontrolofBPcXhZcHcDDcosth}. More precisely we start
with
\beaa
\DD\hot B +(Z+4H)\hot B &=& \nab_3A+\left(\frac{1}{2}\tr\Xb-4\omb\right)A+3\ov{P}\Xh.
\eeaa
from which we deduce, making use of {\bf Ref 1-2}, and the above estimate for $\Xh$, for $k\leq k_*$, 
\beaa
\big\|\dkb^k\DD\hot B\big\|_{L^2(S)}&\les & r^{-2} \|\dkb^{\leq k}B\|_{L^2(S)} +\ep_0  r^{-3-\de'} u^{-1/2-\dec}+ r^{-3}\big\| \dkb^{\leq k}\Xh\big\|_{L^2(S)} \\
&\les& r^{-2} \|\dkb^{\leq k}B\|_{L^2(S)} +\ep_0  r^{-3-\de'} u^{-1/2-\dec}+ \ep_0 r^{-4}u^{-1/2-\dec},
\eeaa
which we write in the form
\beaa
\big\|\dkb^k\DD\hot B\big\|_{L^2(S)}&\les &  r^{-2} \|\dkb^{\leq k}B\|_{L^2(S)} +\ep_0   r^{-3-\de'} u^{-1/2-\dec}.
\eeaa
We  then  proceed   exactly as for the proof of \eqref{eq:Prliminary-forB} and deduce, for all $r_0\leq r\leq r_*$ and $1\leq u\leq u_*$,
\beaa
|\dk^{\leq k_*-4}B| &\les&    \ep_0 r^{-3-\de'} u^{-1/2-\dec}
\eeaa
as stated.

{\bf Step 3.} To get the estimate  for $\Zc$,  we apply Proposition \ref{Prop:transportrp-f-Decay-knorms} to the equation
\beaa
\nab_4\Zc + \frac{2}{q}\Zc  &=&    -B+ O(r^{-2})\Xh +O(r^{-2})\trXc+\Ga_g\c\Ga_g.
\eeaa
Using the estimate already derived for $\widehat{X}, B, \trXc$, assumption {\bf Ref 1},  and the estimate on the last slice for  $\Zc$,  we easily derive, for all $r_0\leq r\leq r_*$ and $1\leq u\leq u_*$,
\beaa
r^2 \big|\dk^{\le k_*-4} \Zc  (r, u)  \big| &\les& r_*^2\big|\dk^{\le k_*-4}   \Zc(r_*, u)  \big|+ \ep_0  r^{-\de'} u^{-1/2-\dec}\les \ep_0 u^{-1/2-\dec} 
\eeaa
as stated.

{\bf Step 4.} To get the estimate  for $\Pc$, we apply Proposition \ref{Prop:transportrp-f-Decay-knorms} to the equation
\beaa
\nab_4\left(\Pc \right)  +\frac{3}{q} \Pc&=& O(r^{-1})\dkb^{\leq 1}B+O(r^{-3})\trXc+r^{-1}\Ga_g\c\Ga_g+\Ga_b\c A.
\eeaa
Using the estimates already derived for $B$ and $\trXc$, assumption {\bf Ref 1},  and the estimate for $\Pc$ on the last slice,  we easily deduce, for all $r_0\leq r\leq r_*$ and $1\leq u\leq u_*$,
\beaa
r^3 \big|\dk^{\le k_*-5} \Pc  (r, u)  \big| &\les& r_*^3  \big|\dk^{\le k_*-5}   \Pc(r_*, u)  \big|+ \ep_0  r^{-\de'} u^{-1/2-\dec} \les \ep_0 u^{-1/2-\dec}  
\eeaa
as stated. This concludes the proof of Proposition \ref{Proposition:estimatesforBPcXhZc-halfdecayinu}.
\end{proof}

%%%%%%%%%%%%%%%%%%%%%%%%%%%%%%%%%%%%%%%%%%%%%%%%%%%%%%%%

\section{End of the proof of Proposition \ref{prop:improvedesitmatesfortemporalframeofMexton-new}}

%%%%%%%%%%%%%%%%%%%%%%%%%%%%%%%%%%%%%%%%%%%%%%%%%%%%%%%%%

\begin{remark} 
\lab{remark:Decay.Mext-sofar}
We summarize the estimates proved so far:
\begin{enumerate}
\item According to Proposition \ref{lemma:firstcontroloftrXbcinMext}, we have
\beaa
\big|\dk^{\le k_*} \trXc\big| &\les& \ep_0 r^{-2} u^{-1-\dec}.
\eeaa

\item According to Proposition  \ref{Prop:estimatesforBPXhZcHc},  we have
\beaa
r\big(|\dk^{\leq k_*-5}B|+|\dk^{\leq k_*-5}\Pc|\big)+|\dk^{\leq k_*-5}\Xh|+|\dk^{\leq k_*-5}\Zc| &\les& \ep_0 r^{-1-\de'} u^{-1-\dec},\\
|\dk^{\leq k_*-5}\Hc|+|\dk^{\leq k_*-5}\widecheck{\DD\cos\th}| &\les&\ep_0 r^{-1} u^{-1-\dec},
\eeaa
and
\beaa
r|\dk^{\leq k_*-6}\nab_3B|+r|\dk^{\leq k_*-6}\nab_3\Pc|+|\dk^{\leq k_*-6}\nab_3\Xh|+|\dk^{\leq k_*-6}\nab_3\Zc| &\les&\ep_0 r^{-2} u^{-1-\dec}.
\eeaa

\item  According to Proposition \ref{Proposition:estimatesforBPcXhZc-halfdecayinu}, we have
\beaa
 r|\dk^{\leq k_*-6}\Pc|+|\dk^{\leq k_*-6}\Xh|+|\dk^{\leq k_*-6}\Zc| &\les&    \ep_0 r^{-2} u^{-1/2-\dec},\\
|\dk^{\leq k_*-6}B| &\les&    \ep_0 r^{-3-\de'} u^{-1/2-\dec}.    
\eeaa
\end{enumerate}
This provides the desired $\Ga_g$ estimates  of Theorem M4  for the quantities 
\beaa
\trXc, \quad \Xh, \quad \Zc, \quad  r B, \quad  r \Pc,
\eeaa
as well as the desired $\Ga_b$ estimates for the quantities
\beaa
\Hc, \quad \widecheck{\DD\cos\th}.
\eeaa
\end{remark}

We now recover the remaining components of the outgoing PG structure of $\Mext$. They are stated in the following proposition.
\begin{proposition}
\lab{Proposition:remaining.estimatesMext}
The following estimates  hold true on $\Mext$
\bea
\lab{eq:Proposition.remaining.estimatesMext1}
\bsplit
 |\dk^{\leq k_*-7}\trXbc|+|\dk^{\leq k_*-6}\Xbh|+|\dk^{\leq k_*-7}\Xib| +|\dk^{\leq k_*-6}\ombc|+|\dk^{\leq k_*-7}\Ab| &\les         \ep_0 r^{-1} u^{-1-\dec},\\
 |\dk^{\leq k_*-6}\Bb|+|\dk^{\leq k_*-8}\nab_3\trXbc|   &\les         \ep_0 r^{-2} u^{-1-\dec},\\
 |\dk^{\leq k_*-7}\trXbc|  &\les \ep_0 r^{-2} u^{-1/2-\dec},
 \end{split}
\eea
\bea
\lab{eq:Proposition.remaining.estimatesMext2}
\begin{split}
|\dk^{\leq k_*-6}\widecheck{e_3(r)}|+|\dk^{\leq k_*-6}\widecheck{e_3(u)}|    &\les         \ep_0 u^{-1-\dec},\\
|\dk^{\leq k_*-6}e_3(\cos\th)|+|\dk^{\leq k_*-6}\widecheck{\DD u}|  &\les         \ep_0 r^{-1} u^{-1-\dec},
\end{split}
\eea
\bea
\lab{eq:Proposition.remaining.estimatesMext3}
|\dk^{\leq k_*-6}\widecheck{\ov{\DD}\c\Jk}|+|\dk^{\leq k_*-6}\DD\hot\Jk|+|\dk^{\leq k_*-6}\widecheck{\nab_3\Jk}|  &\les&         \ep_0 r^{-2} u^{-1-\dec}, 
\eea
\bea
\lab{eq:Proposition.remaining.estimatesMext4}
|\dk^{\leq k_*-6}\widecheck{\ov{\DD}\c\Jk_\pm}|+|\dk^{\leq k_*-6}\DD\hot\Jk_\pm|+|\dk^{\leq k_*-6}\widecheck{\nab_3\Jk_\pm}| &\les&\ep_0 r^{-2}  u^{-1-\dec},
\eea
and
\bea
\lab{eq:Proposition.remaining.estimatesMext5}
|\dk^{\leq k_*-6}\widecheck{\DD(J^{(\pm)})}|+|\dk^{\leq k_*-6}\widecheck{e_3(J^{(\pm)})}|&\les&  \ep_0 r^{-1} u^{-1-\dec}.
\eea
\end{proposition}

\begin{proof}
We proceed in steps as follows.

{\bf Step 1.} We  first  derive  the estimates for $\DD\hot\Jk$,     $\widecheck{\ov{\DD}\c\Jk} $,  $\DD\hot\Jk_\pm $ and $\widecheck{\ov{\DD}\c\Jk_\pm}$   in \eqref{eq:Proposition.remaining.estimatesMext3}, \eqref{eq:Proposition.remaining.estimatesMext4},  with the help of  their  transport equations, see Lemma  \ref{Lemma:linearizedJk}, and their estimates on the last slice $\Si_*$, see Proposition \ref{prop:improvedesitmatesfortemporalframeofMextonSigmastar-new}. Recall from Lemma  \ref{Lemma:linearizedJk} that we have 
\beaa
\nab_4 (\DD\hot\Jk)+\frac{2}{q}\DD\hot\Jk &=& O(r^{-1})B +O(r^{-2})\trXc+O(r^{-2})\Xh +O(r^{-2})\Zc\\
&& +O(r^{-3})\widecheck{\DD(\cos\th)}+r^{-1}\Ga_b\c\Ga_g,\\
\nab_4\big(\widecheck{\ov{\DD}\c\Jk} \big)+\Re\left(\frac{2}{q}\right)\widecheck{\ov{\DD}\c\Jk} &=& O(r^{-1})B+O(r^{-2})\trXc+O(r^{-2})\Xh +O(r^{-2})\Zc\\
&&+O(r^{-3})\widecheck{\DD(\cos\th)}+r^{-1}\Ga_b\c\Ga_g,\\
\nab_4 (\DD\hot\Jk_\pm)+\frac{2}{q}\DD\hot\Jk_\pm  &=& O(r^{-1})B +O(r^{-2})\trXc+O(r^{-2})\Xh +O(r^{-2})\Zc\\
&& +O(r^{-3})\widecheck{\DD(\cos\th)}+r^{-1}\Ga_b\c\Ga_g,\\
\nab_4\big(\widecheck{\ov{\DD}\c\Jk_\pm} \big)+\Re\left(\frac{2}{q}\right)\widecheck{\ov{\DD}\c\Jk_\pm} &=& O(r^{-1})B+O(r^{-2})\trXc+O(r^{-2})\Xh +O(r^{-2})\Zc\\
&&+O(r^{-3})\widecheck{\DD(\cos\th)}+r^{-1}\Ga_b\c\Ga_g.
\eeaa
Let $F$ denote the schematic right hand side of these transport equations, i.e.
\beaa
F= O(r^{-1})B +O(r^{-2})\trXc+O(r^{-2})\Xh +O(r^{-2})\Zc+O(r^{-3})\widecheck{\DD(\cos\th)}+r^{-1}\Ga_b\c\Ga_g.
\eeaa
In view of the estimates already derived  for $B,\, \trXc, \Zc, \widecheck{\DD(\cos\th)}$, see Remark \ref{remark:Decay.Mext-sofar}, and using {\bf Ref 1} for the nonlinear terms, we have on $\Mext$
\beaa
\big\|\dk^{\le k_*- 6} F\big\| _{\infty}(u, r) &\les& \ep_0 r^{-3-\de'} u^{-1-\dec}+ \ep_0 r^{-4} u^{-1-\dec}. 
\eeaa
Thus,  applying Proposition \ref{Prop:transportrp-f-Decay-knorms}  and making use of the estimates  on the last slice $\Si_*$,        we deduce
\beaa
 r^2\Big| \dk^{\le k_*-6}\big( \DD\hot\Jk\big)\Big|+ r^2 \Big| \dk^{\le k_*-6} \big(\widecheck{\ov{\DD}\c\Jk} \big)\Big|&\les& \ep_0 u^{-1-\dec}, \\
  r^2\Big| \dk^{\le k_*-6}\big( \DD\hot\Jk_\pm\big)\Big|+ r^2 \Big| \dk^{\le k_*-6} \big(\widecheck{\ov{\DD}\c\Jk_\pm} \big)\Big|&\les& \ep_0 u^{-1-\dec}, 
\eeaa
as stated in \eqref{eq:Proposition.remaining.estimatesMext3},  \eqref{eq:Proposition.remaining.estimatesMext4}.

{\bf Step 2.}   Next, we estimate $\trXbc$ with the help of the following equation, see Lemma \ref{Lemma:linearized-nullstr},  
\beaa
\nab_4\widecheck{\tr\Xb} +\frac{1}{q}\trXbc &=& F,
\eeaa
where
\beaa
 F&=& -\DD\c\ov{\Zc}+2\ov{\Pc }+O(r^{-2})\Zc+O(r^{-1})\trXc+ O(r^{-1})\widecheck{\DD\c\ov{\Jk}} +O(r^{-3})\widecheck{\DD(\cos\th)}+\Ga_b\c\Ga_g.
\eeaa
In view of the estimates already derived for $\Zc$, $\trXc$, $\widecheck{\DD(\cos\th)}$, $\Pc$, see Remark \ref{remark:Decay.Mext-sofar},  and the estimate for  $\widecheck{\DD\c\ov{\Jk}}$   obtained  in Step 1 above, we have on $\Mext$
\beaa
\big|\dk^{\le k_*-7} F \big|&\les& \min\{ r^{-3} u^{-1/2-\dec}, r^{-2-\de'} u^{-1-\dec} \}.
\eeaa
Thus, applying Proposition \ref{Prop:transportrp-f-Decay-knorms},  we infer on $\Mext$
\beaa
 r\|  \dk^{\le k_*-7}   \trXbc\|_{\infty}(u, r)  &\les    r_* \|\dk^{\le k_*-7}\trXbc\|_{\infty} (u, r_*) + \ep_0  \min\left\{ r^{-1} u^{-1/2-\dec}, r^{-\de'} u^{-1-\dec}\right\}. 
\eeaa
Hence,  making use  of the estimates for $\trXbc$  on $\Si_*$,  we derive
\beaa
\big| \dk^{\le k_*-7}   \trXbc \big|&\les&    \ep_0 \min\left\{r^{-2} u^{-1/2-\dec},   r^{-1} u^{-1-\dec}\right\}
\eeaa
as stated in \eqref{eq:Proposition.remaining.estimatesMext1}.

{\bf Step 3.} Next, we derive the desired estimate for $\Xbh$ using the following equation 
\beaa
\nab_4\Xbh +\frac{1}{q} \Xbh &=& -\frac{1}{2}\DD\hot\Zc +O(r^{-2})\Zc+O(r^{-1})\Xh+O(r^{-1})\DD\hot\Jk+O(r^{-3})\widecheck{\DD(\cos\th)}\\
&&+\Ga_b\c\Ga_g.
\eeaa
Denoting the right hand side by $F$ and  using the estimates already derived  for $\Zc$, $\trXc$, $\widecheck{\DD(\cos\th)}$, $\Pc$, see Remark \ref{remark:Decay.Mext-sofar},  and the estimate for  $\DD\hot\Jk$   obtained  in Step 1 above, we have on $\Mext$
\beaa
\big|\dk^{\le k_*-6} F \big|&\les&\ep_0 r^{-2-\de'} u^{-1-\dec}.
\eeaa
Thus, applying Proposition \ref{Prop:transportrp-f-Decay-knorms}, we infer on $\Mext$ 
\beaa
 r\big\| \dk^{\le k_*-6}   \Xbh \big\|_{\infty}(u, r) &\les&  r_* \big\| \dk^{\le k_*-6}   \Xbh \big\|_{\infty}(u, r_*)  +\ep_0 u^{-1-\dec}.
\eeaa
Hence, using the estimates on the last slice $\Si_*$ for $\Xbh$, we obtain
\beaa
\big| \dk^{\le k_*-6}   \Xbh \big| &\les& \ep_0  r^{-1} u^{-1-\dec}
\eeaa
as stated in   \eqref{eq:Proposition.remaining.estimatesMext1}.

{\bf Step  4.} Next, we   estimate $\ombc$ using the following equation, see Lemma \ref{Lemma:linearized-nullstr},  
\beaa
\nab_4(\ombc) &=& \Re(\Pc)+O(r^{-2})\Zc+O(r^{-2})\Hc+\Ga_b\c\Ga_g.
\eeaa
Denoting the right hand side by $F$ and  using the estimates already derived  for $\Zc$, $\Pc$, $\Hc$, see Remark \ref{remark:Decay.Mext-sofar}, we have on $\Mext$
\beaa
\big|\dk^{\le k_*- 6} F\big|&\les& \ep_0 r^{-2-\de'} u^{-1-\dec}.
\eeaa
Hence, applying Proposition \ref{Prop:transportrp-f-Decay-knorms},  and making use of the estimates  on the last slice $\Si_*$ for $\ombc$, we deduce
\bea
\lab{eq:Proposition.remaining.estimatesMext1-strong}
\big| \dk^{\le k_*-6}\ombc\big|(r, u) &\les  \ep_0 r_*^{-1} u^{-1-\dec}+ \ep_0 r^{-1-\de'} u^{-1-\dec},
\eea
which implies the estimate for $\ombc$ in \eqref{eq:Proposition.remaining.estimatesMext1}.

{\bf Step 5.} Next, we estimate  $\Xib$ using the following equation, see Lemma \ref{Lemma:linearized-nullstr},
\beaa
\nab_4\Xib  +\frac{1}{q}\Xib &=& O(r^{-1})\dkb^{\leq 1}(\ombc)+O(r^{-2})\Zc+O(r^{-2})\Hc+O(r^{-2})\trXbc\\
&&+O(r^{-3})\widecheck{\DD(\cos\th)}+\Ga_b\c\Big(\ombc,\Ga_g\Big).
\eeaa
Denoting the right hand side by $F$ and using the estimates already  derived  for  $\Zc$, $\Hc$, $\widecheck{\DD(\cos\th)}$, see Remark \ref{remark:Decay.Mext-sofar}, as well as  the estimate \eqref{eq:Proposition.remaining.estimatesMext1-strong} for $\ombc$, and the estimate for $\trXbc$ of Step 2 above, we derive on $\Mext$
\beaa
\big|\dk^{\le k_*-7} F \big|&\les&    \ep_0 r_*^{-2} u^{-1-\dec}+ \ep_0 r^{-2-\de'} u^{-1-\dec}.
\eeaa
Applying Proposition \ref{Prop:transportrp-f-Decay-knorms}, and making use of the estimate for $\Xib$ on the last slice, we deduce
\beaa
 r \big| \dk^{\le k_*-7}\Xib \big|(r, u) &\les& r _*\big| \dk^{\le k_*-6}\Xib\big|(r_*, u) +  \frac{(r^*-r)\ep_0}{r_*u^{1+\dec}}+\frac{\ep_0}{u^{1+\dec}}\\
  &\les&\ep_0 u^{-1-\dec}.
\eeaa
Consequently, we obtain
\bea
\big| \dk^{\le k_*-7}\Xib \big|(r, u) &\les& \ep_0  r^{-1} u^{-1-\dec}
\eea
as stated.

{\bf Step 6.} To estimate $\Bb$, we make use of the following equation
\beaa
\nab_4\Bb+\frac{2}{q}\Bb    &=& -    \DD\left(\Pc\right) +O(r^{-2})\Pc+O(r^{-3})\Zc +O(r^{-4})\widecheck{\DD(\cos\th)}+r^{-1}\Ga_b\c\Ga_g.
\eeaa
Denoting the right hand side by $F$ and making use of the estimates already derived for $\Pc$, $\Zc$,  $\widecheck{\DD(\cos\th)}$,   see Remark \ref{remark:Decay.Mext-sofar}, we have on $\Mext$
\beaa
\big|\dk^{\le k_*-6} F \big|&\les&    \ep_0 r^{-3-\de'} u^{-1-\dec}.
\eeaa
Thus, integrating with the help of   Proposition \ref{Prop:transportrp-f-Decay-knorms}, and making use of the estimate for $\Bb$ on the last slice, we deduce
\beaa
\big| \dk^{\le k_*-6}\Bb \big| &\les& \ep_0  r^{-2} u^{-1-\dec}
\eeaa
as stated in \eqref{eq:Proposition.remaining.estimatesMext1}.

{\bf Step 7.} To estimate $\Ab$, we make use of the following equation, see Proposition \ref{prop-nullstrandBianchi:complex:outgoing:again:chap6}, 
\beaa
\nab_4\Ab +\frac{1}{2}\DD\hot\Bb &=& -\frac{1}{2}\ov{\tr X} \Ab +\frac{5}{2}Z\hot \Bb -3P\Xbh,
\eeaa
which yields
\beaa
\nab_4\Ab+\frac{1}{\ov{q}}\Ab &=& O(r^{-1})\dkb^{\leq 1}\Bb+O(r^{-3})\Xbh+\Ga_g\c\Ga_b.
\eeaa
Denoting the right hand side by $F$ and making use of the estimates already derived for $\Bb$ and $\Xbh$, we have on $\Mext$
\beaa
\big|\dk^{\le k_*-7} F \big|&\les&    \ep_0 r^{-3} u^{-1-\dec}.
\eeaa
Thus, integrating with the help of   Proposition \ref{Prop:transportrp-f-Decay-knorms}, and making use of the estimate for $\Ab$ on the last slice, we deduce
\beaa
\big| \dk^{\le k_*-7}\Ab \big| &\les& \ep_0  r^{-1} u^{-1-\dec}
\eeaa
as stated in \eqref{eq:Proposition.remaining.estimatesMext1}.

{\bf Step 8.}  Next, we estimate   $\widecheck{e_3(r)}$ using the following equation, see Lemma  \ref{Lemma:otherlinearizedquant},
\beaa
e_4\left(\widecheck{e_3(r)}\right) &=& -2\ombc,
\eeaa
Integrating, and using the  estimate  \eqref{eq:Proposition.remaining.estimatesMext1-strong}  for $ \ombc$ derived above, we deduce
\beaa
\Big| \dk^{\le k_*-6} \left(\widecheck{e_3(r)}\right) (r, u) \Big|&\les & \Big| \dk^{\le k_*-6} \left(\widecheck{e_3(r)}\right) \Big|(r_*, u) + \frac{(r^*-r)\ep_0}{r_*u^{1+\dec}}+\frac{\ep_0}{u^{1+\dec}}\\
&\les & \frac{\ep_0}{u^{1+\dec}}.
\eeaa
Hence
\beaa
\Big| \dk^{\le k_*-6} \left(\widecheck{e_3(r)}\right)  \Big|&\les & \ep_0 u^{-1-\dec}
\eeaa
as stated in \eqref{eq:Proposition.remaining.estimatesMext2}.

{\bf Step 9.}  Next, we derive   estimates for $\widecheck{\DD\cos\th}$,   $\widecheck{e_3(u)}$, and $e_3(\cos \th)$  by relying on the following   equations,  see  Lemma \ref{Lemma:otherlinearizedquant},
\beaa
\nab_4\widecheck{\DD u}+\frac{1}{q}\widecheck{\DD u} &=& O(r^{-1})\trXc+O(r^{-1})\Xh+\Ga_b\c\Ga_g,\\
e_4\left(\widecheck{e_3(u)}\right) &=& O(r^{-1})\Hc+O(r^{-1})\Zc+O(r^{-2})\widecheck{\DD u}+\Ga_b\c\Ga_b,\\
e_4(e_3(\cos\th)) &=& O(r^{-1})\Hc+O(r^{-1})\Zc+O(r^{-2})\widecheck{\DD\cos\th}+\Ga_b\c\Ga_b.
\eeaa
 Note that the right side $F_1$ of the first equation verifies
\beaa
\big|\dk^{\le k_*-6}  F_1\big|&\les&\ep_0 r^{-2-\de'} u^{-1-\dec}.
\eeaa
Proceeding exactly as before, we infer 
\beaa
 r\| \dk^{\le k_*-6} \widecheck{\DD u}\|_{L^\infty}(u, r)  &\les&  r_*\| \dk^{\le k_*-6} \widecheck{\DD u}\|_{L^\infty}(u, r_*)  + \ep_0r^{-\de'} u^{-1-\dec}\les\ep_0 u^{-1-\dec}.
\eeaa
Thus,
\beaa
\big|\dk^{\le k_*-6} \widecheck{\DD u}\big| &\les& \ep_0  r^{-1}u^{-1-\dec}
\eeaa
as stated in \eqref{eq:Proposition.remaining.estimatesMext2}.

Also, the  right hand side $F_2$ of the second equation verifies
\beaa
\big|\dk^{\le k_*-6}  F_2\big|&\les&\ep_0 r^{-2} u^{-1-\dec}.
\eeaa
Hence, by integration, in the same manner, we obtain 
\beaa
\|\dk^{\le k_*-6} \widecheck{e_3(u)}\|_{L^\infty}(u, r)&\les& \|\dk^{\le k_*-6}\widecheck{e_3(u)}\|_{L^\infty}(u, r_*)+\ep_0 r^{-1}u^{-1-\dec}.
\eeaa
Thus, according to the estimate on the last slice for  $\widecheck{e_3(u)}$, see Proposition \ref{prop:improvedesitmatesfortemporalframeofMextonSigmastar-new}, we deduce
\beaa
\big|\dk^{\le k_*-6} \widecheck{e_3( u)}\big| &\les& \ep_0  u^{-1-\dec}
\eeaa
as stated in \eqref{eq:Proposition.remaining.estimatesMext2}. The estimate \eqref{eq:Proposition.remaining.estimatesMext2} for $e_3(\cos\th)$ follows exactly   in the same manner. This ends the proof of the estimates   \eqref{eq:Proposition.remaining.estimatesMext2}.

{\bf Step 10.} Next, we estimate $\nab_3\trXbc$. In view of Proposition \ref{prop-nullstr:complex}, we have
\beaa
\nab_3\tr\Xb &=& -\frac{1}{2}(\tr\Xb)^2-2\omb \tr\Xb+O(r^{-1})\dkb^{\leq 1}\Xib+\Ga_b\c\Ga_b
\eeaa
which implies
\beaa
\nab_3\trXbc &=& O(r^{-1})\dkb^{\leq 1}\Xib+O(r^{-1})\trXbc+O(r^{-1})\ombc+O(r^{-2})\widecheck{e_3(r)}+O(r^{-2})e_3(\cos\th)\\
&& +\Ga_b\c\Ga_b.
\eeaa
We deduce from the above estimates for $\Xib$, $\trXbc$, $\ombc$, $\widecheck{e_3(r)}$ and $e_3(\cos\th)$ that 
\beaa
\big|\dk^{\le k_*-8}\nab_3\trXbc\big| &\les& \ep_0  r^{-2}u^{-1-\dec}.
\eeaa
This ends the proof of the estimates   \eqref{eq:Proposition.remaining.estimatesMext1}.

{\bf Step 11.} Next, we estimate  $\widecheck{\nab_3\Jk}$ and   $\widecheck{\nab_3\Jk_\pm}$. According to Lemma \ref{Lemma:linearizedJk}, we have
\beaa
\nab_4\big(\widecheck{\nab_3\Jk}\big) +\frac{1}{q}\widecheck{\nab_3\Jk} &=& O(r^{-3})\widecheck{e_3(r)}+O(r^{-3})e_3(\cos\th)+O(r^{-2})\ombc+O(r^{-2})\Hc\\
&&+O(r^{-2})\Zc+O(r^{-2})\widecheck{\nab\Jk}+O(r^{-1})\Pc+r^{-1}\Ga_b\c\Ga_g,\\
\nab_4\big(\widecheck{\nab_3\Jk_\pm}\big) +\frac{1}{q}\widecheck{\nab_3\Jk_\pm} &=& O(r^{-3})\widecheck{e_3(r)}+O(r^{-3})e_3(\cos\th)+O(r^{-2})\ombc+O(r^{-2})\Hc\\
&&+O(r^{-2})\Zc+O(r^{-2})\widecheck{\nab\Jk}+O(r^{-1})\Pc+r^{-1}\Ga_b\c\Ga_b.
\eeaa
Denoting the right hand sides by $F$, we easily check, using the previously derived estimates,
\beaa
\big|\dk^{\le k_*-6}  F\big|&\les&\ep_0 r^{-3} u^{-1-\dec}.
\eeaa
Proceeding exactly as before, we infer
\beaa
 r\| \dk^{\le k_*-6} \widecheck{\nab_3(\Jk, \Jk_\pm )}\|_{L^\infty}(u, r)  &\les&  r_*\| \dk^{\le k_*-6} \widecheck{\nab_3\Jk}\|_{L^\infty}(u, r_*)  + \ep_0r^{-1} u^{-1-\dec}\les\ep_0  r^{-1}u^{-1-\dec}.
\eeaa
Hence
\beaa
\big|\dk^{\le k_*-6} \widecheck{\nab_3\Jk}\big| &\les& \ep_0  r^{-2} u^{-1-\dec},  \qquad \big|\dk^{\le k_*-6} \widecheck{\nab_3\Jk_\pm}\big| \les  \ep_0  r^{-2} u^{-1-\dec},
\eeaa
as stated. This ends the proof of the estimates \eqref{eq:Proposition.remaining.estimatesMext3} and  \eqref{eq:Proposition.remaining.estimatesMext4}.

{\bf Step 12.}  It only remains to prove the estimates  \eqref{eq:Proposition.remaining.estimatesMext5}  for $J^{(\pm)}$. We make use of Lemma \ref{lemma:linearizationforJ+-} according to which
\beaa
\bsplit
\nab_4\big(  \widecheck{\DD(J^{(\pm)})}  \big)+\frac{1}{q}\widecheck{\DD(J^{(\pm)})} &=  O(r^{-1} ) \widecheck{\tr X} +
O(r^{-1}) \Xh+\Ga_b\c\Ga_g,\\
 \nab_4 \big( \widecheck{\nab_3 J^{(\pm)}}\big)&=O(r^{-2} ) \widecheck{\DD(J^{(\pm)})}   +O(r^{-1}) \Zc+ O(r^{-1})  \Hc+\Ga_b\c\Ga_b.
 \end{split}
\eeaa
From the first equation we easily derive,  using the estimates for $\trXc$ and $\Xh$ derived before, 
\beaa
r \big\| \dk^{\le k_*-6}   \widecheck{\DD J^{(\pm)} }        \big\| _{L^\infty}(u, r) &\les& r_* \big\| \dk^{\le k_*-6}   \widecheck{\DD J^{(\pm)} }        \big\| _{L^\infty}(u, r_*) + \ep_0 r^{-\de'} u^{-1 -\dec}\\
&\les& \ep_0 u^{-1-\dec}.
\eeaa
Thus, we obtain 
\beaa
\big|\dk^{\le k_*-6}   \widecheck{\DD J^{(\pm)} }        \big| &\les& \ep_0  r^{-1} u^{-1-\dec}
\eeaa
as stated. Finally integrating the  equation for  $\nab_4 \big( \widecheck{\nab_3 J^{(+)}}\big)$ we derive
\beaa
\big|\dk^{\le k_*-6} \widecheck{\nab_3J^{(\pm)}}\big| \les  \ep_0  r^{-1} u^{-1-\dec}.
\eeaa
This ends the proof of \eqref{eq:Proposition.remaining.estimatesMext5}  and concludes the proof of Proposition \ref{Proposition:remaining.estimatesMext}.
\end{proof}

The estimates in Remark \ref{remark:Decay.Mext-sofar}, together with the ones of Proposition \ref{Proposition:remaining.estimatesMext}, conclude the proof of  Proposition \ref{prop:improvedesitmatesfortemporalframeofMexton-new}. Also, recalling from \eqref{eq:valueofkstarinchapter6forproofThmM4} that $k_*=k_{small}+60$ in this chapter, this concludes the proof of        Theorem M4 as stated in section \ref{sec:mainintermediateresults:chap3}.

%%%%%%%%%%%%%%%%%%%%%%%%%%%%%%%%%%%%%%%%

\chapter{Decay estimates on $\Mint$ and $\Mtop$ (Theorem M5)}
\lab{Chapter:DecayMintMtop}

%%%%%%%%%%%%%%%%%%%%%%%%%%%%%%%%%%%%%%%%

The goal of this chapter is to prove Theorem M5, i.e. to derive decay estimates on $\Mint$ and $\Mtop$.

%%%%%%%%%%%%%%%%%%%%%%%%%%%%%%%%%%%%%%%%%

\section{Linearized equations for ingoing PG structures}
\lab{sec:linearizedequationsforingoingPGstructures}

%%%%%%%%%%%%%%%%%%%%%%%%%%%%%%%%%%%%%%%%%

Ingoing PG structures have been introduced in section \ref{sec:Principalingoinggeodesicstructures}. In particular, recall that such structures verify the following identities, see section  \ref{sec:basicdefinitionforingoingPGstruct}, 
 \beaa
&&  \xib=0, \quad \omb=0, \quad \eta=\ze, \quad  e_3(r)=-1, \quad \nab(r)=0,\quad e_3(\ub)=e_3(\th)= e_3(\vphi)=0,\\
&&\nab_3\Jk = \frac{1}{\ov{q}}\Jk, \quad  \nab_3\Jk _\pm= \frac{1}{\ov{q}}\Jk_\pm,\quad e_3(\Jp)=0,\,\, p=0,+,-.
\eeaa

In this section we provide the linearized equations for ingoing PG structures that will be used to derive decay estimates for the ingoing PG structures of $\Mint$ and $\Mtop$. Recall that  the definition of the linearized quantities for  ingoing PG structures  can be found in Definition \ref{def:renormalizationofallnonsmallquantitiesinPGstructurebyKerrvalue:ingoingcase}. 

\begin{remark}\lab{rmk:howtogetingoingequationfromoutgoingones:chap7}
The equations for ingoing PG structures stated in this section can be easily deduced from their analog for outgoing PG structures by performing the following substitutions
\beaa
&& u\to \ub, \quad r\to r, \quad \th\to \th, \quad \vphi\to \vphi, \quad e_4\to e_3, \quad e_3\to e_4, \quad e_a\to e_a, \\
&& \a\to \aa, \quad \b\to -\bb, \quad \rho\to \rho, \quad \rhod\to -\rhod, \quad \bb\to -\b,\quad \aa\to \a,\\
&& \xi\to \xib, \quad \om\to \omb, \quad \chi\to \chib, \quad \eta\to\etab, \quad \etab\to \eta, \quad\ze\to -\ze, \quad \chib\to \chi,\quad \omb\to\om, \quad \xib\to\xi,\\
&&\Jp\to \Jp, \quad \Jk\to \Jk, \quad \Jk_\pm\to \Jk_\pm.
\eeaa
\end{remark}

In view of Remark \ref{rmk:howtogetingoingequationfromoutgoingones:chap7}, the following lemma can be easily obtained from its analog for the outgoing case in Lemma \ref{Lemma:linearized-nullstr}.
\begin{lemma}
\lab{Lemma:linearized-nullstr:ingoingcase}
The linearized null structure equations in the  $e_3$ direction are 
\beaa
\nab_3(\trXbc) -\frac{2}{\ov{q}}\trXbc &=& \Ga_b\c\Ga_b,\\
\nab_3\Xbh -\frac{2r}{|q|^2}\Xbh &=& -\Ab+\Ga_b\c\Ga_b,\\
\nab_3\Zc - \frac{2}{\ov{q}}\Zc  &=&  -\Bb+  O(r^{-2})\Xbh  +O(r^{-2})\trXbc+\Ga_b\c\Ga_g,\\
\nab_3\Hbc-\frac{1}{\ov{q}}\Hbc &=& \Bb+O(r^{-1})\Zc     +O(r^{-2})\Xbh  +O(r^{-2})\trXbc +\Ga_b\c\Ga_g,\\
\nab_3\trXc -\frac{1}{\ov{q}}\trXc &=& \DD\c\ov{\Zc}+2\Pc +O(r^{-2})\Zc+O(r^{-1})\trXbc\\
&&+ O(r^{-1})\widecheck{\DD\c\ov{\Jk}} +O(r^{-3})\widecheck{\DD(\cos\th)}+\Ga_b\c\Ga_g,\\
\nab_3\Xh -\frac{1}{\ov{q}} \Xh &=& \frac{1}{2}\DD\hot\Zc +O(r^{-2})\Zc+O(r^{-1})\Xbh+O(r^{-1})\DD\hot\Jk+O(r^{-3})\widecheck{\DD(\cos\th)}\\
&&+\Ga_b\c\Ga_g,\\
\nab_3(\omc) &=& \Re(\Pc)+O(r^{-2})\Zc+O(r^{-2})\Hbc+\Ga_b\c\Ga_g,\\
\nab_3\Xi  -\frac{1}{\ov{q}}\Xi &=& O(r^{-1})\dkb^{\leq 1}(\omc)+O(r^{-2})\Zc+O(r^{-2})\Hbc+O(r^{-2})\trXc\\
&&+O(r^{-3})\widecheck{\DD(\cos\th)}+\Ga_b\c\Big(\omc,\Ga_g\Big).
\eeaa
The  linearized Bianchi equations for $B, P, \Bb$  are 
 \beaa
\nab_3\Bb -\frac{4}{\ov{q}}\Bb &=&  \frac{1}{2}\ov{\DD}\c\Ab +O(r^{-2})\Ab+\Ga_b\c(\Bb,\Ab),\\
\nab_3\left(\Pc \right)+\frac{1}{2}\ov{\DD}\c\Bb &=& \frac{3}{\ov{q}} \Pc +O(r^{-2})\Bb+O(r^{-3})\trXbc+r^{-1}\Ga_b\c\Ga_g+\Ga_g\c \Ab,\\
\nab_3B+\DD\ov{\Pc} &=& \frac{2}{\ov{q}}B +O(r^{-2})\Pc+O(r^{-3})\Zc +O(r^{-4})\widecheck{\DD(\cos\th)}+r^{-1}\Ga_b\c\Ga_g.
\eeaa
\end{lemma}

In view of Remark \ref{rmk:howtogetingoingequationfromoutgoingones:chap7}, the following lemma can be easily obtained from its analog for the outgoing case in Lemma \ref{Lemma:otherlinearizedquant}.
\begin{lemma}
\lab{Lemma:otherlinearizedquant:ingoingcase}
We have
\beaa
e_3\left(\widecheck{e_4(r)}\right) &=& -2\omc,\\
\nab_3\widecheck{\DD\ub}-\frac{1}{\ov{q}}\widecheck{\DD\ub} &=& O(r^{-1})\trXbc+O(r^{-1})\Xbh+\Ga_b\c\Ga_b,\\
e_3\left(\widecheck{e_4(\ub)}\right) &=& O(r^{-1})\Hbc+O(r^{-1})\Zc+O(r^{-2})\widecheck{\DD\ub}+\Ga_g\c\Ga_g,
\\
\nab_3\widecheck{\DD\cos\th}-\frac{1}{\ov{q}}\widecheck{\DD\cos\th} &=& \frac{i}{2}\ov{\Jk}\c\Xbh+O(r^{-1})\trXbc+\Ga_b\c\Ga_b,\\
e_3(e_4(\cos\th)) &=& O(r^{-1})\Hbc+O(r^{-1})\Zc+O(r^{-2})\widecheck{\DD\cos\th}+\Ga_g\c\Ga_g.
\eeaa
\end{lemma}

In view of Remark \ref{rmk:howtogetingoingequationfromoutgoingones:chap7}, the following lemma can be easily obtained from its analog for the outgoing case in Lemma \ref{Lemma:linearizedJk}.
\begin{lemma}
\lab{Lemma:linearizedJk:ingoingcase}
The following equations hold  for   the tensors $\Jk, \Jk_\pm$.
\begin{enumerate}
\item We have
\beaa
\nab_3(\DD\hot\Jk)-\frac{2}{\ov{q}}\DD\hot\Jk &=& O(r^{-1})\Bb +O(r^{-2})\trXbc+O(r^{-2})\Xbh\\
&&+O(r^{-2})\Zc +O(r^{-3})\widecheck{\DD(\cos\th)},\\
\nab_3\big(\widecheck{\ov{\DD}\c\Jk} \big)-\Re\left(\frac{2}{\ov{q}}\right)\widecheck{\ov{\DD}\c\Jk} &=& O(r^{-1})\Bb+O(r^{-2})\trXbc+O(r^{-2})\Xbh +O(r^{-2})\Zc\\
&&+O(r^{-3})\widecheck{\DD(\cos\th)},\\
\nab_3\big(\widecheck{\nab_4\Jk}\big) -\frac{1}{\ov{q}}\widecheck{\nab_4\Jk} &=& O(r^{-3})\widecheck{e_4(r)}+O(r^{-3})e_4(\cos\th)+O(r^{-2})\omc\\
&&+O(r^{-2})\Hbc+O(r^{-2})\Zc+O(r^{-2})\widecheck{\nab\Jk}+O(r^{-1})\Pc.
\eeaa
\item We also have
\beaa
\nab_3(\DD\hot\Jk_\pm)-\frac{2}{\ov{q}}\DD\hot\Jk _\pm &=& O(r^{-1})\Bb +O(r^{-2})\trXbc+O(r^{-2})\Xbh\\
&&+O(r^{-2})\Zc +O(r^{-3})\widecheck{\DD(\cos\th)},\\
\nab_3\big(\widecheck{\ov{\DD}\c\Jk_\pm} \big)-\Re\left(\frac{2}{\ov{q}}\right)\widecheck{\ov{\DD}\c\Jk}_\pm &=& O(r^{-1})\Bb+O(r^{-2})\trXbc+O(r^{-2})\Xbh +O(r^{-2})\Zc\\
&&+O(r^{-3})\widecheck{\DD(\cos\th)},\\
\nab_3\big(\widecheck{\nab_4\Jk_\pm}\big) -\frac{1}{\ov{q}}\widecheck{\nab_4\Jk_\pm} &=& O(r^{-3})\widecheck{e_4(r)}+O(r^{-3})e_4(\cos\th)+O(r^{-2})\omc\\
&&+O(r^{-2})\Hbc+O(r^{-2})\Zc+O(r^{-2})\widecheck{\nab\Jk}+O(r^{-1})\Pc.
\eeaa
\end{enumerate}
\end{lemma} 

In view of Remark \ref{rmk:howtogetingoingequationfromoutgoingones:chap7}, the following lemma can be easily obtained from its analog for the outgoing case in Lemma \ref{lemma:linearizationforJ+-}.
\begin{lemma}
\lab{lemma:linearizationforJ+-:ingoingcase}
The following  equations hold true\footnote{Similar equations hold for $J^{(0)}=\cos \th$, see Lemma \ref{Lemma:otherlinearizedquant:ingoingcase}.}.
\bea
\bsplit
\nab_3\big(  \widecheck{\DD(J^{(\pm)})}  \big) -\frac{1}{\ov{q}} \widecheck{\DD(J^{(\pm)})}  &=  O(r^{-1} )\trXbc +O(r^{-1}) \Xbh +\Ga_b\c\Ga_g,\\
 \nab_3\big( \widecheck{\nab_4 J^{(\pm)}}\big)&=O(r^{-2} ) \widecheck{\DD(J^{(\pm)})}   +O(r^{-2}) \Zc+ O(r^{-2})  \Hbc+\Ga_g\c\Ga_g.
 \end{split}
\eea
\end{lemma}

%%%%%%%%%%%%%%%%%%%%%%%%%%%%%%%%%%%%%%%%%

\section{Decay estimates for the PG structure of $\Mint$ on $\TT$}
\lab{section:DecayPGMint-onTT}

%%%%%%%%%%%%%%%%%%%%%%%%%%%%%%%%%%%%%%%%%

To simplify  the notations, in this section, we denote
\begin{itemize}
\item by $(e_4, e_3, e_1, e_2)$ the outgoing PG frame of $\Mext$, with all quantities associated to the outgoing PG structure of $\Mext$ being unprimed, 

\item by $(e_4', e_3', e_1', e_2')$ the ingoing PG frame of $\Mint$, with all quantities associated to the ingoing PG structure of $\Mint$ being primed. 
\end{itemize}

Recall that $\Mext\cap\Mint=\TT=\{r=r_0\}$. In view of the above notations, and the initialization of the ingoing PG structure of $\Mint$ from the outgoing PG structure of $\Mext$ on $\TT$, see section \ref{sec:initalizationadmissiblePGstructure}, we have  
\bea
\ub=u, \quad r'=r, \quad {J'}^{(p)}=J^{(p)},\,\, p=0,+,-,  \quad \Jk'=\Jk, \quad \Jk_\pm'=\Jk_\pm\quad\textrm{on}\quad\TT,
\eea
and
\bea
e_4'=\la e_4, \qquad e_3'=\la^{-1}e_3, \qquad e_a'=e_a, \,\, a=1,2, \quad\textrm{on}\quad\TT,
\eea
where $\la$ is given by
\bea
\la=\frac{\Delta}{|q|^2}.
\eea

In order to derive decay estimates for the ingoing PG structure of $\Mint$ on $\TT$, we will  rely on the following lemma.

\begin{lemma}\lab{lemma:relationsPGMintandMextonTT}
We have on $\TT$
\beaa
&& A'=\la^2A, \qquad B'=\la B,\qquad \Pc'=\Pc, \qquad \Bb'=\la^{-1}\Bb, \qquad \Ab'=\la^{-2}\Ab, \\
&& \Xib'=0, \qquad \omb'=0, \qquad H'=Z',\\ 
&& \trXc'=\la\trXc, \qquad \Xh'=\la\Xh, \qquad \trXbc'=\la^{-1}\trXbc, \qquad \Xbh'=\la^{-1}\Xbh,\\
&& e_3'(r')=-1, \qquad \nab'(r')=0, \qquad e_3'(\ub)=0, \qquad e_3'({J'}^{(p)})=0,\,\,p=0,+,-,\\
&& \nab_3'\Jk'=\frac{1}{\ov{q'}}\Jk', \qquad \nab_3'\Jk_\pm'=\frac{1}{\ov{q'}}\Jk_\pm',\\
&& \widecheck{\nab'(\ub)}= \widecheck{\nab(u)}, \qquad \widecheck{\nab'({J'}^{(p)})}= \widecheck{\nab(J^{(p)})}, \,\, p=0,+,-,\\
&& \widecheck{\ov{\DD}'\c\Jk'}=\widecheck{\ov{\DD}\c\Jk}, \qquad \DD'\hot\Jk'=\DD\hot\Jk, \qquad \widecheck{\ov{\DD}'\c\Jk_\pm'}=\widecheck{\ov{\DD}\c\Jk_\pm}, \qquad \DD'\hot\Jk_\pm'=\DD\hot\Jk_\pm,
\eeaa
\beaa
\Zc' &=& \Zc+\frac{1}{q}\widecheck{\DD(q)}+\frac{1}{\ov{q}}\widecheck{\DD(\ov{q})},\\
 \Hbc' &=&  - \Zc - \frac{1}{e_3(r)}\Xib,\\
\Xi' &=& \frac{\la^2}{e_3(r)}\left(\Zc +\frac{1}{q}\widecheck{\DD(q)}+\frac{1}{\ov{q}}\widecheck{\DD(\ov{q})} - \Hc\right),\\
\omc' &=& \frac{\la}{e_3(r)}\ombc+\frac{1}{2e_3(r)}\pr_r\left(\frac{\De}{|q|^2}\right)\widecheck{e_3(r)}  -\frac{4a^2\la\cos\th}{e_3(r)|q|^2}e_3(\cos\th),
\eeaa
and
\beaa
\widecheck{e_4'(r')} &=& -\frac{\la}{e_3(r)}\widecheck{e_3(r)},\\
\widecheck{e_4'(\ub)} &=& -\frac{\la}{e_3(r)}\widecheck{e_3(u)} - \frac{2(r^2+a^2)}{e_3(r)|q|^2}\widecheck{e_3(r)},\\
e_4'({J'}^{(0)}) &=& -\frac{\la}{e_3(r)}e_3(J^{(0)}),\\
\widecheck{e_4'({J'}^{(\pm)})}  &=& -\frac{\la}{e_3(r)}\widecheck{e_3(J^{(\pm)})} \pm\frac{2a}{|q|^2e_3(r)}J^{(\mp)}\widecheck{e_3(r)},\\
\widecheck{\nab_{e_4'}\Jk'} &=& -\frac{\la}{e_3(r)}\widecheck{\nab_3\Jk},\\
\widecheck{\nab_{e_4'}\Jk_\pm'}   &=& -\frac{\la}{e_3(r)}\widecheck{\nab_{e_3}\Jk_\pm} \pm\frac{2a}{|q|e_3(r)}\widecheck{e_3(r)}\Jk_\mp,
\eeaa
where the definition of the linearized quantities for the outgoing PG structure of $\Mext$ can be found in Definition \ref{def:renormalizationofallnonsmallquantitiesinPGstructurebyKerrvalue}, while  definition of the linearized quantities for the ingoing PG structure of $\Mint$ can be found in Definition \ref{def:renormalizationofallnonsmallquantitiesinPGstructurebyKerrvalue:ingoingcase}. 
\end{lemma}

\begin{proof}
The identities for $\Xib'$, $\omb'$, $H'-Z'$, $e_3'(r')$, $e_3'(\ub)$, $e_3'({J'}^{(p)})$, $\nab_3'\Jk'$ and $\nab_3'\Jk_\pm'$ come from the ingoing PG structure assumption on $\Mint$. Also, the identities for $A'$, $B'$, $\Bb'$, $\Ab'$, $\Xh'$ and $\Xbh'$ follow immediately from the change of frame formulas of Proposition \ref{Proposition:transformationRicci} with coefficients $(f=0, \fb=0, \la)$ and the fact that $(e_1, e_2)$ are tangent to $\TT$. Also, the identities for $\Pc'$, $\trXc'$ and $\trXbc'$, follow immediately from the change of frame formulas of Proposition \ref{Proposition:transformationRicci} with coefficients $(f=0, \fb=0, \la)$, the explicit choice for $\la$, the fact that $q'=q$ on $\TT$, and the fact that $(e_1, e_2)$ are tangent to $\TT$. Also, the identities for $\nab'(r')$, $\widecheck{\nab'(\ub)}$, $\widecheck{\nab'({J'}^{(p)})}$, $\widecheck{\ov{\DD}'\c\Jk'}$, $\DD'\hot\Jk'$, $ \widecheck{\ov{\DD}'\c\Jk_\pm'}$, and $\DD'\hot\Jk_\pm'$, follow immediately from the fact that we have, on $\TT$,  $\nab'=\nab$, $r'=r$, $\ub=u$, ${J'}^{(p)}=J^{(p)}$, $\Jk'=\Jk$, $\Jk_\pm'=\Jk_\pm$, together with the fact that $\nab$ is tangent to $\TT$.

It remains to derive the identities for $\Zc'$, $\Xi'$, $\om'$, $\Hbc' $,  $\widecheck{e_4'(r)}$, $\widecheck{e_4'(\ub)}$, $e_4'({J'}^{(0)})$, $\widecheck{e_4'({J'}^{(\pm)})}$, $\widecheck{\nab_4'\Jk'}$ and $\widecheck{\nab_4'\Jk_{\pm}'}$. We start with $Z'$. In view of the change of frame formulas of Proposition \ref{Proposition:transformationRicci} with coefficients $(f=0, \fb=0, \la)$, and the fact that $(e_1, e_2)$ are tangent to $\TT$, we have
\beaa
Z' &=& Z -\DD'(\log\la). 
\eeaa
Since $\DD'=\DD$, and using the explicit form of $\la$, as well as $\nab(r)=0$, we infer
\beaa
Z' &=& Z -\DD\left(\log\left(\frac{\De}{|q|^2}\right)\right)= Z +\frac{1}{|q|^2}\DD(|q|^2)= Z +\frac{1}{q}\DD(q)+\frac{1}{\ov{q}}\DD(\ov{q})
\eeaa
which yields, together with the fact that $\Jk'=\Jk$ and $q'=q$ on $\TT$, in view of the linearization of the various quantities, and taking the different linearization for $Z'$ (ingoing PG structure) and $Z$ (outgoing PG structure) into account,
\beaa
\Zc' &=& \Zc +\frac{1}{q}\widecheck{\DD(q)}+\frac{1}{\ov{q}}\widecheck{\DD(\ov{q})}
\eeaa  
as desired. 

Next, since $e_4'=\la e_4$, $e_3'=\la^{-1}e_3$, and $e_a'=e_a$, $a=1,2$, on $\TT$, and since $e_3-e_3(r)e_4$ is tangent to $\TT$, we have 
\beaa
\g(\D_{e_3-e_3(r)e_4}e_3', e_a') &=& \g(\D_{e_3-e_3(r)e_4}(\la^{-1}e_3), e_a) ,\\
\g(\D_{e_3-e_3(r)e_4}e_4', e_a') &=& \g(\D_{e_3-e_3(r)e_4}(\la e_4), e_a) ,\\
\g(\D_{e_3-e_3(r)e_4}e_4', e_3') &=& \g(\D_{e_3-e_3(r)e_4}(\la e_4), \la^{-1}e_3).
\eeaa
Hence, using again that $e_4'=\la e_4$ and $e_3'=\la^{-1}e_3$ on $\TT$, we deduce 
\beaa
\g(\D_{\la e_3'-e_3(r)\la^{-1}e_4'}e_3', e_a') &=& \la^{-1}\g(\D_{e_3-e_3(r)e_4}e_3, e_a) ,\\
\g(\D_{\la e_3'-e_3(r)\la^{-1}e_4'}e_4', e_a') &=& \la\g(\D_{e_3-e_3(r)e_4}e_4, e_a) ,\\
\g(\D_{\la e_3'-e_3(r)\la^{-1}e_4'}e_4', e_3') &=& -2(e_3-e_3(r)e_4)\log(\la)+\g(\D_{e_3-e_3(r)e_4}e_4, e_3),
\eeaa
and hence
\beaa
\la\xib' -\la^{-1}e_3(r)\etab' &=& \la^{-1}(\xib -e_3(r)\etab),\\
\la\eta' - e_3(r)\la^{-1}\xi' &=& \la(\eta -e_3(r)\xi),\\
-\la\omb' -e_3(r)\la^{-1}\om' &=& -\omb-e_3(r)\om -2(e_3-e_3(r)e_4)\log(\la).
\eeaa
Since $\xib'=0$, $\xi=0$, $\omb'=0$, $\om=0$, $\eta'=\ze'$ and $\etab=-\ze$, we infer
\beaa
 \etab' &=&  - \ze - \frac{1}{e_3(r)}\xib,\\
\xi' &=& \frac{\la^2}{e_3(r)}(\ze' - \eta),\\
 \om' &=& \frac{\la}{e_3(r)}\omb  +\frac{2\la}{e_3(r)}(e_3-e_3(r)e_4)\log(\la).
\eeaa
In view of the lineraizations for ingoing and outgoing PG structures,  the fact that $\Jk'=\Jk$ and $q'=q$ on $\TT$, and the above identity for $\Zc'$, we obtain
\beaa
 \Hbc' &=&  - \Zc - \frac{1}{e_3(r)}\Xib,\\
\Xi' &=& \frac{\la^2}{e_3(r)}\left(\Zc +\frac{1}{q}\widecheck{\DD(q)}+\frac{1}{\ov{q}}\widecheck{\DD(\ov{q})} - \Hc\right),\\
\omc' &=& \frac{\la}{e_3(r)}\ombc+\frac{1}{2}\pr_r\left(\frac{\De}{|q|^2}\right)\left(1+\frac{\la}{e_3(r)}\right)  +\frac{2\la}{e_3(r)}(e_3-e_3(r)e_4)\log(\la).
\eeaa
These are the desired identities for $\Hbc'$ and $\Xi'$. For $\omc'$, we note that, in view of the formula for $\la$, we have $e_3(r)=-\la+\widecheck{e_3(r)}$. Also, $(e_3-e_3(r)e_4)(\De)=0$ and $e_4(\th)=0$. Hence
\beaa
&&\frac{1}{2}\pr_r\left(\frac{\De}{|q|^2}\right)\left(1+\frac{\la}{e_3(r)}\right)  +\frac{2\la}{e_3(r)}(e_3-e_3(r)e_4)\log(\la)\\
&=& \frac{1}{2e_3(r)}\pr_r\left(\frac{\De}{|q|^2}\right)\widecheck{e_3(r)}  -\frac{4a^2\la\cos\th}{e_3(r)|q|^2}e_3(\cos\th)
\eeaa 
which yields 
\beaa
\omc' &=& \frac{\la}{e_3(r)}\ombc+\frac{1}{2e_3(r)}\pr_r\left(\frac{\De}{|q|^2}\right)\widecheck{e_3(r)}  -\frac{4a^2\la\cos\th}{e_3(r)|q|^2}e_3(\cos\th)
\eeaa
as desired. 

It remains to derive the identities for $\widecheck{e_4'(r)}$, $\widecheck{e_4'(\ub)}$, $e_4'({J'}^{(0)})$, $\widecheck{e_4'({J'}^{(\pm)})}$, $\widecheck{\nab_4'\Jk'}$ and $\widecheck{\nab_4'\Jk_{\pm}'}$. Since we have $r'=r$, $\ub=u$, ${J'}^{(p)}=J^{(p)}$, $\Jk'=\Jk$ and $\Jk_\pm'=\Jk$ on $\TT$, and since  $e_3-e_3(r)e_4$ is tangent to $\TT$, we have 
\beaa
(e_3-e_3(r)e_4)r' &=& (e_3-e_3(r)e_4)r,\\
(e_3-e_3(r)e_4)\ub &=& (e_3-e_3(r)e_4)u,\\
(e_3-e_3(r)e_4){J'}^{(p)} &=& (e_3-e_3(r)e_4)\Jp, \quad p=0,+,-,\\
\nab_{e_3-e_3(r)e_4}\Jk' &=& \nab_{e_3-e_3(r)e_4}\Jk,\\
\nab_{e_3-e_3(r)e_4}\Jk_\pm' &=& \nab_{e_3-e_3(r)e_4}\Jk_\pm.
\eeaa
Using the fact that $e_4'=\la e_4$ and $e_3'=\la^{-1}e_3$ on $\TT$, and since  
\beaa
&& e_3'(r')=-1, \quad e_3'(\ub)=0, \quad e_3'({J'}^{(p)})=0,\,\,p=0,+,-,\quad \nab_3'\Jk'=\frac{1}{\ov{q'}}\Jk', \quad \nab_3'\Jk_\pm'=\frac{1}{\ov{q'}}\Jk_\pm',\\
&& e_4(r)=1, \quad e_4(u)=0, \quad e_4(J^{(p)})=0,\,\,p=0,+,-,\quad \nab_4\Jk=-\frac{1}{q}\Jk, \quad \nab_4\Jk_\pm=-\frac{1}{q}\Jk_\pm,
\eeaa
we infer
\beaa
e_4'(r') &=& -\frac{\la^2}{e_3(r)},\\
e_4'(\ub) &=& -\frac{\la}{e_3(r)}e_3(u),\\
e_4'({J'}^{(p)}) &=& -\frac{\la}{e_3(r)}e_3(\Jp), \quad p=0,+,-,\\
\nab_{e_4'}\Jk' &=& -\frac{\la}{e_3(r)}\left(\nab_{e_3}\Jk+\frac{e_3(r)}{q}\Jk-\frac{\la}{\ov{q'}}\Jk'\right),\\
\nab_{e_4'}\Jk_\pm' &=& -\frac{\la}{e_3(r)}\left(\nab_{e_3}\Jk_\pm+\frac{e_3(r)}{q}\Jk_\pm-\frac{\la}{\ov{q'}}\Jk_\pm'\right).
\eeaa
Since $r'=r$, $q'=q$, ${J'}^{(p)}=\Jp$, $\Jk'=\Jk$ and $\Jk_\pm'=\Jk_\pm$ on $\TT$, and since $e_3(r)=-\la+\widecheck{e_3(r)}$, we deduce, in view of the linerizations for ingoing and outgoing PG structures,
\beaa
\widecheck{e_4'(r')} &=& -\frac{\la}{e_3(r)}\widecheck{e_3(r)},\\
\widecheck{e_4'(\ub)} &=& -\frac{\la}{e_3(r)}\widecheck{e_3(u)} - \frac{2(r^2+a^2)}{e_3(r)|q|^2}\widecheck{e_3(r)},\\
e_4'({J'}^{(0)}) &=& -\frac{\la}{e_3(r)}e_3(J^{(0)}),\\
\widecheck{e_4'({J'}^{(\pm)})}  &=& -\frac{\la}{e_3(r)}\widecheck{e_3(J^{(\pm)})} \pm\frac{2a}{|q|^2e_3(r)}J^{(\mp)}\widecheck{e_3(r)},\\
\widecheck{\nab_{e_4'}\Jk'} &=& -\frac{\la}{e_3(r)}\widecheck{\nab_3\Jk},\\
\widecheck{\nab_{e_4'}\Jk_\pm'}   &=& -\frac{\la}{e_3(r)}\widecheck{\nab_{e_3}\Jk_\pm} \pm\frac{2a}{|q|e_3(r)}\widecheck{e_3(r)}\Jk_\mp,
\eeaa
as desired. This concludes the proof of the lemma.
\end{proof}

We are now ready  to derive decay estimates for the ingoing PG structure of $\Mint$ on $\TT$.

\begin{lemma}\lab{lemma:controlofPGstructureMintonTT}
The following decay estimates hold on $\TT$ for the ingoing PG structure of $\Mint$ 
\bea
\sup_{\TT}r'\ub^{1+\dec}|\dk^{\leq k_{small}+40}(\Ga_g', \Ga_b')| &\les& \ep_0.
\eea
\end{lemma}

\begin{proof}
In view of the control of outgoing PG structure of $\Mext$ established in Theorem M4 and the fact that $\TT\subset\Mext$, we have
\beaa
\sup_{\TT}ru^{1+\dec}|\dk^{\leq k_{small}+40}(\Ga_g, \Ga_b)| &\les& \ep_0.
\eeaa

Note that the tangential derivatives to $\TT=\{r=r_0\}$ are generated by $\nab$ and $\nab_3-e_3(r)\nab_4$. We introduce the following notation for $r$-weighted tangential derivatives to $\TT$
\beaa
\widetilde{\dk} &:=& \big(\nab_3-e_3(r)\nab_4, r\nab\big).
\eeaa
From the above estimate for $(\Ga_g, \Ga_b)$, together the identities of Lemma \ref{lemma:relationsPGMintandMextonTT} on $\TT$ and the fact that $r'=r$ and $\ub=u$ on $\TT$, we infer
\beaa
\sup_{\TT}r'\ub^{1+\dec}|\widetilde{\dk}^{\leq k_{small}+40}(\Ga_g', \Ga_b')| &\les& \ep_0.
\eeaa

Finally, since $\dk$ is generated by $e_3'$ and $\widetilde{\dk}$, the previous estimate and the control of $e_3'$ derivatives provided by the null structure equations and Bianchi identities of the ingoing PG structure of $\Mint$ immediately imply 
\beaa
\sup_{\TT}r'\ub^{1+\dec}|\dk^{\leq k_{small}+40}(\Ga_g', \Ga_b')| &\les& \ep_0
\eeaa
as stated.
\end{proof}

%%%%%%%%%%%%%%%%%%%%%%%%%%%%%%%%%%%%%%%%%

\section{Decay estimates in $\Mint$}
\lab{sec:decayestimatesMintmainstepsorder}

%%%%%%%%%%%%%%%%%%%%%%%%%%%%%%%%%%%%%%%%%

In this section, all quantities appearing correspond to the ingoing PG structure of $\Mint$. We are now ready to prove the part of Theorem M5 concerning $\Mint$, i.e. to derive decay estimates for the ingoing PG structure of $\Mint$. To this end, recall first that $\Ab$ has already been estimated in Theorem M2 and satisfies 
\bea\lab{eq:controlofalphabarneededforThmM5andcomingfromThmM3}
\sup_{\Mint}\ub^{1+\dec}|\dk^{\leq k_{small}+80}\Ab| &\les& \ep_0.
\eea

Relying on the estimates of the  ingoing PG structure of $\Mint$ on $\TT$ derived in Lemma \ref{lemma:controlofPGstructureMintonTT}, we propagate these estimates to $\Mint$ thanks to the linearized transport equations in the $e_3$ direction of section \ref{sec:linearizedequationsforingoingPGstructures} for ingoing PG structures. Recalling that $\Ab$ has already been estimated, see \eqref{eq:controlofalphabarneededforThmM5andcomingfromThmM3}, the other 
 quantities are recovered in the following order:
\begin{enumerate}
\item We recover $\trXbc$, with a control of $k_{small}+40$ derivatives, from
\beaa
\nab_3(\trXbc) -\frac{2}{\ov{q}}\trXbc &=& \Ga_b\c\Ga_b.
\eeaa

\item We recover $\Xbh$, with a control of $k_{small}+40$ derivatives, from
\beaa
\nab_3\Xbh -\frac{2r}{|q|^2}\Xbh &=& -\Ab+\Ga_b\c\Ga_b.
\eeaa

\item We recover $\Bb$, with a control of $k_{small}+40$ derivatives, from 
\beaa
\nab_3\Bb -\frac{4}{q}\Bb &=&  \frac{1}{2}\ov{\DD}\c\Ab +O(r^{-2})\Ab+\Ga_b\c(\Bb,\Ab).
\eeaa

\item We recover $\Zc$, with a control of $k_{small}+40$ derivatives, from
\beaa
\nab_3\Zc - \frac{2}{\ov{q}}\Zc  &=&  -\Bb+  O(r^{-2})\Xbh  +O(r^{-2})\trXbc+\Ga_b\c\Ga_g.
\eeaa

\item We recover $\Hb$, with a control of $k_{small}+40$ derivatives, from 
\beaa
\nab_3\Hbc-\frac{1}{q}\Hbc &=& \Bb+O(r^{-1})\Zc     +O(r^{-2})\Xbh  +O(r^{-2})\trXbc +\Ga_b\c\Ga_g.
\eeaa

\item We recover $\widecheck{\DD\cos\th}$,  with a control of $k_{small}+40$ derivatives, from
\beaa
\nab_3\widecheck{\DD\cos\th}-\frac{1}{\ov{q}}\widecheck{\DD\cos\th} &=& \frac{i}{2}\ov{\Jk}\c\Xbh+O(r^{-1})\trXbc+\Ga_b\c\Ga_g.
\eeaa

\item We recover $\DD\hot\Jk$,  with a control of $k_{small}+40$ derivatives, from
\beaa
\nab_3(\DD\hot\Jk)-\frac{2}{\ov{q}}\DD\hot\Jk &=& O(r^{-1})\Bb +O(r^{-2})\trXbc+O(r^{-2})\Xbh\\
&&+O(r^{-2})\Zc +O(r^{-3})\widecheck{\DD(\cos\th)}.
\eeaa

\item We recover $\widecheck{\ov{\DD}\c\Jk}$,  with a control of $k_{small}+40$ derivatives, from
\beaa
\nab_3\big(\widecheck{\ov{\DD}\c\Jk} \big)-\Re\left(\frac{2}{\ov{q}}\right)\widecheck{\ov{\DD}\c\Jk} &=& O(r^{-1})\Bb+O(r^{-2})\trXbc+O(r^{-2})\Xbh +O(r^{-2})\Zc\\
&&+O(r^{-3})\widecheck{\DD(\cos\th)}.
\eeaa

\item We recover $e_4(\cos\th)$,  with a control of $k_{small}+40$ derivatives, from
\beaa
e_3(e_4(\cos\th)) &=& O(r^{-1})\Hbc+O(r^{-1})\Zc+O(r^{-2})\widecheck{\DD\cos\th}+\Ga_b\c\Ga_b.
\eeaa

\item We recover $\Pc$, with a control of $k_{small}+39$ derivatives, from
\beaa
\nab_3\left(\Pc \right) - \frac{3}{q} \Pc &=& -\frac{1}{2}\ov{\DD}\c\Bb  +O(r^{-2})\Bb+O(r^{-3})\trXbc+r^{-1}\Ga_b\c\Ga_g+\Ga_g\c \Ab.
\eeaa

\item We recover $\trXc$, with a control of $k_{small}+39$ derivatives, from 
\beaa
\nab_3\trXc -\frac{1}{\ov{q}}\trXc &=& \DD\c\ov{\Zc}+2\ov{\Pc }+O(r^{-2})\Zc+O(r^{-1})\trXbc\\
&&+ O(r^{-1})\widecheck{\DD\c\ov{\Jk}} +O(r^{-3})\widecheck{\DD(\cos\th)}+\Ga_b\c\Ga_g.
\eeaa

\item We recover $\Xh$, with a control of $k_{small}+39$ derivatives, from 
\beaa
\nab_3\Xh -\frac{1}{\ov{q}} \Xh &=& \frac{1}{2}\DD\hot\Zc +O(r^{-2})\Zc+O(r^{-1})\Xbh+O(r^{-1})\DD\hot\Jk+O(r^{-3})\widecheck{\DD(\cos\th)}\\
&&+\Ga_b\c\Ga_g.
\eeaa

\item We recover $\check{\om}$, with a control of $k_{small}+39$ derivatives, from 
\beaa
\nab_3(\omc) &=& \Re(\Pc)+O(r^{-2})\Zc+O(r^{-2})\Hbc+\Ga_b\c\Ga_g.
\eeaa

\item We recover $\widecheck{e_4(r)}$, with a control of $k_{small}+39$ derivatives, from 
\beaa
e_3\left(\widecheck{e_4(r)}\right) &=& -2\omc.
\eeaa

\item We recover $\widecheck{\nab_4\Jk}$, with a control of $k_{small}+39$ derivatives, from 
\beaa
\nab_3\big(\widecheck{\nab_4\Jk}\big) -\frac{1}{\ov{q}}\widecheck{\nab_4\Jk} &=& O(r^{-3})\widecheck{e_4(r)}+O(r^{-3})e_4(\cos\th)+O(r^{-2})\omc\\
&&+O(r^{-2})\Hbc+O(r^{-2})\Zc+O(r^{-2})\widecheck{\nab\Jk}+O(r^{-1})\Pc.
\eeaa

\item We recover $B$, with a control of $k_{small}+38$ derivatives, from 
\beaa
\nab_3B - \frac{2}{\ov{q}}B &=& -\DD\ov{\Pc} +O(r^{-2})\Pc+O(r^{-3})\Zc +O(r^{-4})\widecheck{\DD(\cos\th)}+r^{-1}\Ga_b\c\Ga_g.
\eeaa

\item We recover $\Xi$, with a control of $k_{small}+38$ derivatives, from 
\beaa
\nab_3\Xi  -\frac{1}{\ov{q}}\Xi &=& O(r^{-1})\dkb^{\leq 1}(\omc)+O(r^{-2})\Zc+O(r^{-2})\Hbc+O(r^{-2})\trXc\\
&&+O(r^{-3})\widecheck{\DD(\cos\th)}+\Ga_b\c\Big(\omc,\Ga_g\Big).
\eeaa

\item We recover $A$, with a control of $k_{small}+37$ derivatives, from 
\beaa
\nab_3A -\frac{1}{\ov{q}}A &=& \frac{1}{2}\DD\hot B +O(r^{-2}) B +O(r^{-3})\Xh.
\eeaa

\item We recover $\widecheck{\DD\ub}$, with a control of $k_{small}+40$ derivatives, from
\beaa
\nab_3\widecheck{\DD\ub}-\frac{1}{\ov{q}}\widecheck{\DD\ub} &=& O(r^{-1})\trXbc+O(r^{-1})\Xbh+\Ga_b\c\Ga_g.
\eeaa

\item We recover $\widecheck{e_4(\ub)}$, with a control of $k_{small}+40$ derivatives, from
\beaa
e_3\left(\widecheck{e_4(\ub)}\right) &=& O(r^{-1})\Hbc+O(r^{-1})\Zc+O(r^{-2})\widecheck{\DD\ub}+\Ga_b\c\Ga_b.
\eeaa

\item We recover $\DD\hot\Jk_\pm$, with a control of $k_{small}+40$ derivatives, from
\beaa
\nab_3(\DD\hot\Jk_\pm)-\frac{2}{\ov{q}}\DD\hot\Jk _\pm &=& O(r^{-1})\Bb +O(r^{-2})\trXbc+O(r^{-2})\Xbh\\
&&+O(r^{-2})\Zc +O(r^{-3})\widecheck{\DD(\cos\th)}.
\eeaa

\item We recover $\widecheck{\ov{\DD}\c\Jk_\pm}$, with a control of $k_{small}+40$ derivatives, from
\beaa
\nab_3\big(\widecheck{\ov{\DD}\c\Jk_\pm} \big)-\Re\left(\frac{2}{\ov{q}}\right)\widecheck{\ov{\DD}\c\Jk}_\pm &=& O(r^{-1})\Bb+O(r^{-2})\trXbc+O(r^{-2})\Xbh +O(r^{-2})\Zc\\
&&+O(r^{-3})\widecheck{\DD(\cos\th)}.
\eeaa

\item We recover $\widecheck{\nab_4\Jk_\pm}$, with a control of $k_{small}+39$ derivatives, from
\beaa
\nab_3\big(\widecheck{\nab_4\Jk_\pm}\big) -\frac{1}{\ov{q}}\widecheck{\nab_4\Jk_\pm} &=& O(r^{-3})\widecheck{e_4(r)}+O(r^{-3})e_4(\cos\th)+O(r^{-2})\omc\\
&&+O(r^{-2})\Hbc+O(r^{-2})\Zc+O(r^{-2})\widecheck{\nab\Jk}+O(r^{-1})\Pc.
\eeaa

\item We recover $\widecheck{\DD(J^{(\pm)})}$, with a control of $k_{small}+40$ derivatives, from
\beaa
\nab_3\big(  \widecheck{\DD(J^{(\pm)})}  \big) -\frac{1}{\ov{q}} \widecheck{\DD(J^{(\pm)})}  &=&  O(r^{-1} )\trXbc +O(r^{-1}) \Xbh +\Ga_b\c\Ga_g.
\eeaa

\item We recover $\widecheck{\nab_4 J^{(\pm)}}$, with a control of $k_{small}+40$ derivatives, from
\beaa
\nab_3\big( \widecheck{\nab_4 J^{(\pm)}}\big)&=O(r^{-2} ) \widecheck{\DD(J^{(\pm)})}   +O(r^{-2}) \Zc+ O(r^{-2})  \Hbc+\Ga_b\c\Ga_b.
\eeaa
\end{enumerate}

As the estimates are significantly simpler to derive\footnote{Note that $r$ is bounded on $\Mint$ and that all quantities behave the same in $\Mint$.} 
and in the same spirit as the corresponding ones in Theorem M4, we leave the details to the reader. This concludes the proof of Theorem M5 for the part of $\Mint$.

%%%%%%%%%%%%%%%%%%%%%%%%%%%%%%%%%%%%%%%%%

\section{Decay estimates for the PG structure of $\Mtop$ on $\{u=u_*\}$}

%%%%%%%%%%%%%%%%%%%%%%%%%%%%%%%%%%%%%%%%%

To simplify  the notations, in this section, we denote
\begin{itemize}
\item by $(e_4, e_3, e_1, e_2)$ the outgoing PG frame of $\Mext$, with all quantities associated to the outgoing PG structure of $\Mext$ being unprimed, 

\item by $(e_4', e_3', e_1', e_2')$ the ingoing PG frame of $\Mtop$, with all quantities associated to the ingoing PG structure of $\Mtop$ being primed. 
\end{itemize}

\begin{remark}\lab{rmk:noJpmandnoJkpminMtop:chap7}
Note that in $\Mtop$, we do not need to define $\vphi'$, ${J'}^{(\pm)}$ and $\Jk'_\pm$. In particular, recall from Remark \ref{rmk:noJpmandnoJkpminMtop:chap3} that the quantities $\Ga_g'$, $\Ga_b'$ in $\Mtop$ correspond to the ones in Definition \ref{definition.Ga_gGa_b:ingoingcase} where all linearized quantities based on ${J'}^{(\pm)}$ and $\Jk'_\pm$ have been removed. 
\end{remark}

Recall that $\Mext\cap\Mtop=\{u=u_*\}$. In view of the above notations, and the initialization of the ingoing PG structure of $\Mtop$ from the outgoing PG structure of $\Mext$ on $\{u=u_*\}$, see section \ref{sec:initalizationadmissiblePGstructure}, we have  
\bea
r'=r, \quad {J'}^{(0)}=J^{(0)},\,\, p=0,+,-,  \quad \Jk'=\Jk,
\eea
\bea
\ub=u+2\int_{r_0}^r\frac{{\tilde{r}}^2+a^2}{{\tilde{r}}^2-2m\tilde{r}+a^2}d\tilde{r}, 
\eea
and
\bea
 \bsplit
   e_4' &=\la e_4,\\
  e_a' &= e_a +\frac 1 2  \fb_a e_4,\\
  e_3' &= \la^{-1}\left(e_3 + \fb^b e_b  + \frac 1 4 |\fb|^2 e_4\right),
 \end{split}
 \eea
where
\bea
\la=\frac{\Delta}{|q|^2}, \qquad \fb=h \widecheck{e_3(r)}\nab(u),
\eea
with the scalar function $h$ given by
\bea
h &=& \frac{4}{e_3(u)+\sqrt{(e_3(u))^2+4|\nab u|^2\widecheck{e_3(r)}}}.
\eea
Note in particular that $h=1+O(r^{-2})$, which together with the fact that $\widecheck{e_3(r)}\in r\Ga_b$ and $\nab(u)=O(r^{-1})$ implies
\bea
\fb\in \Ga_b.
\eea

In order to derive decay estimates for the ingoing PG structure of $\Mtop$ on $\{u=u_*\}$, we will  rely on the following lemma.

\begin{lemma}\lab{lemma:necessaryidentitiesonuequalustarforthecontrolofMtopnuequlustar}
We have on $\{u=u_*\}$
\beaa
&& A'=\la^2A, \qquad B'=\la B+\Ga_b\c A,\qquad \Pc'=\Pc+r^{-1}\Ga_b\c\Ga_g, \\
&&\Bb'=\la^{-1}\Bb+r^{-3}\Ga_b+r^{-1}\Ga_b\c\Ga_g, \qquad \Ab'=\la^{-2}\Ab+r^{-1}\Ga_b\c\Ga_b, \\
&& \Xib'=0, \qquad \omb'=0, \qquad H'=Z',\\ 
&& \Xi'=0, \qquad \om'=-\frac{1}{2}\pr_r\left(\frac{\De}{|q|^2}\right), \qquad \Hbc'=-\Zc+r^{-1}\dk^{\leq 1}\Ga_b,\\
&& e_3'(r')=-1, \qquad \nab'(r')=0, \qquad e_3'(\ub)=0, \qquad e_3'({J'}^{(0)})=0,\qquad \nab_3'\Jk'=\frac{1}{\ov{q'}}\Jk',\\
&&\widecheck{e_4'(r')} =0, \qquad \widecheck{e_4'(\ub)} =0,\qquad e_4'({J'}^{(0)}) =0,\qquad \widecheck{\nab_{e_4'}\Jk'} =0,
\eeaa
\beaa
&& \widecheck{\nab'(\ub)}=\widecheck{\nab(u)}+\Ga_b, \qquad \widecheck{\nab'({J'}^{(0)})}=\widecheck{\nab(J^{(0)})}+r^{-1}\Ga_b, \\ 
&& \DD'\hot\Jk'=\DD\hot\Jk+ r^{-2}\Ga_b, \qquad \ov{\DD}'\c\Jk'=\ov{\DD}\c\Jk+ r^{-2}\Ga_b,
\eeaa
\beaa
\trXc'=\la\trXc+r^{-1}\Ga_b, \qquad \Xh'=\la\Xh+r^{-1}\Ga_b,\\
\trXbc'=\la^{-1}\trXbc+r^{-1}\dkb^{\leq 1}\Ga_b, \qquad \Xbh'=\la^{-1}\Xbh+r^{-1}\dkb^{\leq 1}\Ga_b,
\eeaa
and
\beaa
\Zc' &=& \Zc +\frac{1}{q}\widecheck{\DD(q)}+\frac{1}{\ov{q}}\widecheck{\DD(\ov{q})}+r^{-1}\Ga_b,
\eeaa  
where the definition of the linearized quantities for the outgoing PG structure of $\Mext$ can be found in Definition \ref{def:renormalizationofallnonsmallquantitiesinPGstructurebyKerrvalue}, while  definition of the linearized quantities for the ingoing PG structure of $\Mtop$ can be found in Definition \ref{def:renormalizationofallnonsmallquantitiesinPGstructurebyKerrvalue:ingoingcase}. 
\end{lemma}

\begin{proof}
The identities for $\Xib'$, $\omb'$, $H'-Z'$, $e_3'(r')$, $e_3'(\ub)$, $e_3'({J'}^{(0)})$ and $\nab_3'\Jk'$  come from the ingoing PG structure assumption on $\Mtop$. Also, the identities for $A'$, $B'$, $\Bb'$, $\Ab'$, and $\Xi'$ follow immediately from the change of frame formulas of Proposition \ref{Proposition:transformationRicci} with coefficients $(f=0, \fb, \la)$, the fact that $\fb\in\Ga_b$, and the fact that $e_4$ is tangent to $\{u=u_*\}$. Also, the identities for $\Pc'$, $\omc'$ and $\Hbc'$ follow immediately from the change of frame formulas of Proposition \ref{Proposition:transformationRicci} with coefficients $(f=0, \fb, \la)$, the explicit choice for $\la$, the fact that $\fb\in\Ga_g$, the fact that $q'=q$ on $\{u=u_*\}$, and the fact that $e_4$ is tangent to $\{u=u_*\}$. Also, the identities for $\widecheck{e_4'(r')}$, $\widecheck{e_4'(\ub)}$, $e_4'({J'}^{(0)})$ and $\widecheck{\nab_{e_4'}\Jk'}$ follow immediately from the fact that we have, on $\{u=u_*\}$,  $e_4'=\la e_4$, $r'=r$, ${J'}^{(0)}=J^{(0)}$, $\Jk'=\Jk$, and
\beaa
\ub=u+2\int_{r_0}^r\frac{{\tilde{r}}^2+a^2}{{\tilde{r}}^2-2m\tilde{r}+a^2}d\tilde{r},
\eeaa
together with the fact that $e_4$ is tangent to $\{u=u_*\}$.

It remains to derive the identities for $\nab'(r')$, $\Zc'$, $\trXc'$, $\Xh'$, $\trXbc'$, $\Xbh'$, $\widecheck{\nab'(\ub)}$,   $\widecheck{\nab'({J'}^{(0)})}$, $\widecheck{\ov{\DD}'\c\Jk'}$ and $\DD'\hot\Jk'$. Since, $\nab - \frac{1}{e_3(u)}\nab(u) \nab_3$ is tangent to $\{u=u_*\}$ and since
\beaa
\nab - \frac{1}{e_3(u)}\nab(u) \nab_3 &=&  \nab' -\frac 1 2  \fb \la^{-1}\nab_{e_4'} - \frac{1}{e_3(u)}\nab(u)\nab_{\la e_3'  -   \fb^ae_a' +\frac 1 4 |\fb|^2\la^{-1} e_4'},
\eeaa
we have on $\{u=u_*\}$
\beaa
\left(\nab' -\frac 1 2  \fb \la^{-1}\nab_{e_4'} - \frac{1}{e_3(u)}\nab(u)\nab_{\la e_3'  -   \fb^ae_a' +\frac 1 4 |\fb|^2\la^{-1} e_4'}\right)r' &=& \left(\nab - \frac{1}{e_3(u)}\nab(u) \nab_3\right)r.
\eeaa
Since $\nab(r)=0$, $e_3'(r')=-1$, and, in view of the above, $e_4'(r')=\la$, we infer, using also $e_3(r)=-\la+\widecheck{e_3(r)}$, 
\beaa
\nab'(r')  + \frac{\fb\c\nab'(r')}{e_3(u)}\nab(u)  &=& \frac 1 2  \fb  + \frac{|\fb|^2}{4e_3(u)}\nab(u) - \frac{\widecheck{e_3(r)}}{e_3(u)}\nab(u).
\eeaa
Now, recall that we have chosen $\fb=h \widecheck{e_3(r)}\nab(u)$. We infer
\beaa
\nab'(r')  + \frac{h \widecheck{e_3(r)}}{e_3(u)}(\nab(u)\c\nab'(r'))\nab(u)  &=& \Big(\widecheck{e_3(r)}|\nab(u)|^2h^2 +2e_3(u)h  - 4\Big)\frac{\widecheck{e_3(r)}}{4e_3(u)}\nab(u).
\eeaa
Also, recall that we have made the following choice for $h$
\beaa
h &=& \frac{4}{e_3(u)+\sqrt{(e_3(u))^2+4|\nab u|^2\widecheck{e_3(r)}}}
\eeaa
and note that such a choice satisfies 
\beaa
\widecheck{e_3(r)}|\nab(u)|^2h^2 +2e_3(u)h  - 4 &=& 0.
\eeaa
We deduce 
\beaa
\nab'(r')  + \frac{h \widecheck{e_3(r)}}{e_3(u)}(\nab(u)\c\nab'(r'))\nab(u)  &=& 0.
\eeaa
Since $\widecheck{e_3(r)}\in r\Ga_b$, and since $h=1+O(r_0^{-2})$, $e_3(u)=2+O(r_0^{-2})$ and $|\nab u|^2=O(r_0^{-2})$ on $\{u=u_*\}$, we infer $|\nab'(r)|\les\ep|\nab'(r)|$ and hence
\beaa
\nab'(r') &=& 0
\eeaa
as desired.

Next, using $\fb\in\Ga_b$, note that we have
\beaa
\nab - \frac{1}{e_3(u)}\nab(u) \nab_3 &=&  \nab' -\frac 1 2  \fb \la^{-1}\nab_{e_4'} - \frac{1}{e_3(u)}\nab(u)\nab_{\la e_3'  -   \fb^ae_a' +\frac 1 4 |\fb|^2\la^{-1} e_4'}\\
&=&  \nab'  - \frac{\la}{e_3(u)}\nab(u)\nab_{e_3'}+r^{-1}\Ga_b\dk.
\eeaa
Together with the fact that that we have, on $\{u=u_*\}$,  ${J'}^{(0)}=J^{(0)}$, $\Jk'=\Jk$, and
\beaa
\ub=u+2\int_{r_0}^r\frac{{\tilde{r}}^2+a^2}{{\tilde{r}}^2-2m\tilde{r}+a^2}d\tilde{r},
\eeaa
and together with the fact that $\nab - \frac{1}{e_3(u)}\nab(u) \nab_3$ is tangent to $\{u=u_*\}$, we infer
\beaa
\left(\nab'  - \frac{\la}{e_3(u)}\nab(u)\nab_{e_3'}+r^{-1}\Ga_b\dk\right)\ub &=& \left(\nab - \frac{1}{e_3(u)}\nab(u) \nab_3\right)\left(u+2\int_{r_0}^r\frac{{\tilde{r}}^2+a^2}{{\tilde{r}}^2-2m\tilde{r}+a^2}d\tilde{r}\right),\\
\left(\nab'  - \frac{\la}{e_3(u)}\nab(u)\nab_{e_3'}+r^{-1}\Ga_b\dk\right){J'}^{(0)} &=& \left(\nab - \frac{1}{e_3(u)}\nab(u) \nab_3\right)J^{(0)},\\
\left(\nab'  - \frac{\la}{e_3(u)}\nab(u)\nab_{e_3'}+r^{-1}\Ga_b\dk\right)\Jk' &=& \left(\nab - \frac{1}{e_3(u)}\nab(u) \nab_3\right)\Jk+r^{-2}\Ga_b.
\eeaa
Using the fact that $q'=q$ on $\{u=u_*\}$, and since  
\beaa
e_3'(r')=-1, \quad e_3'(\ub)=0, \quad e_3'({J'}^{(0)})=0,\quad \nab_3'\Jk'=\frac{1}{\ov{q'}}\Jk',\quad \nab(r)=0,
\eeaa
we infer
\beaa
\nab'(\ub) &=& \nab(u) - \frac{1}{e_3(u)}\nab(u)\left(e_3(u)+2e_3(r)\frac{r^2+a^2}{\De}\right)+\Ga_b=\nab(u) +\Ga_b,\\
\nab'({J'}^{(0)}) &=& \nab(J^{(0)})+r^{-1}\Ga_b,\\
\nab'\Jk' - \frac{\la}{e_3(u)\ov{q}}\nab(u)\Jk' &=& \left(\nab - \frac{1}{e_3(u)}\nab(u) \nab_3\right)\Jk+r^{-2}\Ga_b,
\eeaa
and hence, since $q'=q$ and $\Jk'=\Jk$ on $\{u=u_*\}$, we obtain 
\beaa
&& \widecheck{\nab'(\ub)}=\widecheck{\nab(u)}+\Ga_b, \qquad \widecheck{\nab'({J'}^{(0)})}=\widecheck{\nab(J^{(0)})}+r^{-1}\Ga_b, \\ 
&& \DD'\hot\Jk'=\DD\hot\Jk+ r^{-2}\Ga_b, \qquad \ov{\DD}'\c\Jk'=\ov{\DD}\c\Jk+ r^{-2}\Ga_b,
\eeaa
as stated. 

It remains to derive the identities for  $\Zc'$, $\trXc'$, $\Xh'$, $\trXbc'$ and $\Xbh'$. In view of 
\beaa
e_a  &=&   e_a' -\frac 1 2  \fb_a \la^{-1} e_4',\\
e_3 &=&   \la e_3'  -   \fb^ae_a' +\frac 1 4 |\fb|^2\la^{-1} e_4',
 \eeaa
 we have on $\{u=u_*\}$, since $\fb\in\Ga_b$, 
\beaa
e_a'  - \frac{\la}{e_3(u)}e_a(u)e_3' &=& e_a - \frac{1}{e_3(u)}e_a(u)e_3 +r^{-1}\Ga_b\dk.
\eeaa
Since $e_a - \frac{1}{e_3(u)}e_a(u)e_3$ is tangent to $\{u=u_*\}$, and since 
\bea
 \bsplit
   e_4' &=\la e_4,\\
  e_a' &= e_a +\frac 1 2  \fb_a e_4,\\
  e_3' &= \la^{-1}\left(e_3 + \fb^b e_b  + \frac 1 4 |\fb|^2 e_4\right),
 \end{split}
 \eea
we infer on $\{u=u_*\}$
\beaa
\g\left(\D_{e_a'  - \frac{\la}{e_3(u)}e_a(u)e_3'}e_4', e_3'\right) &=& \g\left(\D_{e_a - \frac{1}{e_3(u)}e_a(u)e_3}(\la e_4), \la^{-1}\left(e_3 + \fb^b e_b  + \frac 1 4 |\fb|^2 e_4\right)\right)+r^{-1}\Ga_b,\\
\g\left(\D_{e_a'  - \frac{\la}{e_3(u)}e_a(u)e_3'}e_4', e_b'\right) &=& \g\left(\D_{e_a - \frac{1}{e_3(u)}e_a(u)e_3}(\la e_4), e_b +\frac 1 2  \fb_b e_4\right)+r^{-1}\Ga_b,\\
\g\left(\D_{e_a'  - \frac{\la}{e_3(u)}e_a(u)e_3'}e_3', e_b'\right) &=& \g\left(\D_{e_a - \frac{1}{e_3(u)}e_a(u)e_3}\left(\la^{-1}\left(e_3 + \fb^c e_c  + \frac 1 4 |\fb|^2 e_4\right)\right), e_b +\frac 1 2  \fb_b e_4\right)\\&&+r^{-1}\Ga_b,
\eeaa
and hence, using again $\fb\in \Ga_b$, 
\beaa
\g\left(\D_{e_a'  - \frac{\la}{e_3(u)}e_a(u)e_3'}e_4', e_3'\right) &=& \la^{-1}\g\left(\D_{e_a - \frac{1}{e_3(u)}e_a(u)e_3}(\la e_4), e_3\right)+r^{-1}\Ga_b,\\
\g\left(\D_{e_a'  - \frac{\la}{e_3(u)}e_a(u)e_3'}e_4', e_b'\right) &=& \la\g\left(\D_{e_a - \frac{1}{e_3(u)}e_a(u)e_3}e_4, e_b\right)+r^{-1}\Ga_b,\\
\g\left(\D_{e_a'  - \frac{\la}{e_3(u)}e_a(u)e_3'}e_3', e_b'\right) &=& \la^{-1}\g\left(\D_{e_a - \frac{1}{e_3(u)}e_a(u)e_3}e_3, e_b\right)+r^{-1}\dk^{\leq 1}\Ga_b.
\eeaa
We infer on $\{u=u_*\}$
\beaa
2\ze_a' +\frac{4\la}{e_3(u)}e_a(u)\omb' &=& -2\left(e_a - \frac{1}{e_3(u)}e_a(u)e_3\right)\log\la+2\ze_a+\frac{4}{e_3(u)}e_a(u)\omb+r^{-1}\Ga_b,\\
\chi_{ab}' - \frac{2\la}{e_3(u)}e_a(u)\eta_b' &=& \la\chi_{ab} - \frac{2\la}{e_3(u)}e_a(u)\eta_b +r^{-1}\Ga_b,\\
\chib_{ab}' - \frac{2\la}{e_3(u)}e_a(u)\xib_b' &=& \la^{-1}\chib_{ab} - \frac{2\la^{-1}}{e_3(u)}e_a(u)\xib_b +r^{-1}\dk^{\leq 1}\Ga_b.
\eeaa
Since $\xib'=0$, $\omb'=0$ and $\eta'=\ze'$, and since $\xib\in\Ga_b$ and $\ombc\in\Ga_b$, we deduce on $\{u=u_*\}$
\beaa
2\ze_a'  &=& -2\left(e_a - \frac{1}{e_3(u)}e_a(u)e_3\right)\log\la+2\ze_a+\frac{2}{e_3(u)}e_a(u)\pr_r\left(\frac{\De}{|q|^2}\right)+r^{-1}\Ga_b,\\
\chi_{ab}' - \frac{2\la}{e_3(u)}e_a(u)\ze_b' &=& \la\chi_{ab} - \frac{2\la}{e_3(u)}e_a(u)\eta_b +r^{-1}\Ga_b,\\
\chib_{ab}'  &=& \la^{-1}\chib_{ab}  +r^{-1}\dk^{\leq 1}\Ga_b.
\eeaa
In particular,  since $q'=q$ on $\{u=u_*\}$, we have
\beaa
\trXbc'=\la^{-1}\trXbc+r^{-1}\dk^{\leq 1}\Ga_b, \qquad \Xbh'=\la^{-1}\Xbh+r^{-1}\dk^{\leq 1}\Ga_b,
\eeaa
as stated. Also, since we have
\beaa
-2\left( - \frac{1}{e_3(u)}e_a(u)e_3\right)\log\la+\frac{2}{e_3(u)}e_a(u)\pr_r\left(\frac{\De}{|q|^2}\right) &=& r^{-1}\Ga_b,
\eeaa
we infer from the first identity 
\beaa
Z' &=& Z -\DD(\log\la)+r^{-1}\Ga_b. 
\eeaa
Using the explicit for of $\la$, as well as $\nab(r)=0$, we infer on $\{u=u_*\}$
\beaa
Z' &=& Z -\DD\left(\log\left(\frac{\De}{|q|^2}\right)\right)+r^{-1}\Ga_b= Z +\frac{1}{q}\DD(q)+\frac{1}{\ov{q}}\DD(\ov{q})+r^{-1}\Ga_b,
\eeaa
which yields, together with the fact that $\Jk'=\Jk$ and $q'=q$ on $\{u=u_*\}$, in view of the linearization of the various quantities, and taking the different linearization for $Z'$ (ingoing PG structure) and $Z$ (outgoing PG structure) into account,
\beaa
\Zc' &=& \Zc +\frac{1}{q}\widecheck{\DD(q)}+\frac{1}{\ov{q}}\widecheck{\DD(\ov{q})}+r^{-1}\Ga_b
\eeaa  
as desired. 

Finally, coming back to 
\beaa
\chi_{ab}' - \frac{2\la}{e_3(u)}e_a(u)\ze_b' &=& \la\chi_{ab} - \frac{2\la}{e_3(u)}e_a(u)\eta_b +r^{-1}\Ga_b
\eeaa
and since we have in view of the above on $\{u=u_*\}$
\beaa
Z' - H &=& Z -\DD(\log\la) -H+r^{-1}\Ga_b = Z +\frac{1}{q}\DD(q)+\frac{1}{\ov{q}}\DD(\ov{q}) -H+r^{-1}\Ga_b\\
&=& \Zc +\frac{1}{q}\widecheck{\DD(q)}+\frac{1}{\ov{q}}\widecheck{\DD(\ov{q})} -\Hc+r^{-1}\Ga_b = \Ga_b,
\eeaa
we deduce 
\beaa
\trXc'=\la\trXc+r^{-1}\Ga_b, \qquad \Xh'=\la\Xh+r^{-1}\Ga_b,
\eeaa
as desired. This concludes the proof of the lemma.
\end{proof}

We are now ready  to derive decay estimates for the ingoing PG structure of $\Mtop$ on $\{u=u_*\}$.

\begin{lemma}\lab{lemma:controlofPGstructureMtoptonu=u*}
The following decay estimates hold on $\{u=u_*\}$ for the ingoing PG structure of $\Mtop$ 
\bea
\sup_{\{u=u_*\}}\Big(ru^{1+\dec}+r^2u^{\frac{1}{2}+\dec}\Big)|\dk^{\leq k_{small}+39}\Ga_g'|+\sup_{\{u=u_*\}}ru^{1+\dec}|\dk^{\leq k_{small}+39}\Ga_b'| \les \ep_0.
\eea
\end{lemma}

\begin{proof}
In view of the control of outgoing PG structure of $\Mext$ established in Theorem M4 and the fact that $\{u=u_*\}\subset\Mext$, we have
\beaa
\sup_{\{u=u_*\}}\Big(ru^{1+\dec}+r^2u^{\frac{1}{2}+\dec}\Big)|\dk^{\leq k_{small}+40}\Ga_g|+\sup_{\{u=u_*\}}ru^{1+\dec}|\dk^{\leq k_{small}+40}\Ga_b| &\les& \ep_0.
\eeaa

Note that the tangential derivatives to $\{u=u_*\}$ are generated by $e_4$ and $\nab -\frac{1}{e_3(u)}\nab(u)\nab_3$. We introduce the following notation for $r$-weighted tangential derivatives to $\{u=u_*\}$
\beaa
\widetilde{\dk} &:=& \left(r\left(\nab -\frac{1}{e_3(u)}\nab(u)\nab_3\right), r\nab_4\right).
\eeaa
From the above estimate for $(\Ga_g, \Ga_b)$, together the identities of the previous lemma on $\{u=u_*\}$, we infer
\beaa
\sup_{\{u=u_*\}}\Big(ru^{1+\dec}+r^2u^{\frac{1}{2}+\dec}\Big)|\widetilde{\dk}^{\leq k_{small}+39}\Ga_g'|+\sup_{\{u=u_*\}}ru^{1+\dec}|\widetilde{\dk}^{\leq k_{small}+39}\Ga_b'| &\les& \ep_0.
\eeaa

Finally, since $\dk$ is generated by $e_3'$ and $\widetilde{\dk}$, the previous estimate and the control of $e_3'$ derivatives provided by the null structure equations and Bianchi identities of the ingoing PG structure of $\Mtop$ immediately imply 
\beaa
\sup_{\{u=u_*\}}\Big(ru^{1+\dec}+r^2u^{\frac{1}{2}+\dec}\Big)|\dk^{\leq k_{small}+39}\Ga_g'|+\sup_{\{u=u_*\}}ru^{1+\dec}|\dk^{\leq k_{small}+39}\Ga_b'| &\les& \ep_0
\eeaa
as stated.
\end{proof}

%%%%%%%%%%%%%%%%%%%%%%%%%%%%%%%%%%%%%%%%%

\section{Decay estimates in ${}^{(top)}\MM$}

%%%%%%%%%%%%%%%%%%%%%%%%%%%%%%%%%%%%%%%%%

We are now ready to prove the part of Theorem M5 concerning $\Mtop$, i.e. to derive decay estimates for the ingoing PG structure of $\Mtop$. To this end, recall first that $\Ab$ has already been estimated in Theorem M2 on $\Mtop$. Relying on the estimates of the  ingoing PG structure of $\Mtop$ on $\{u=u_*\}$ derived in Lemma \ref{lemma:controlofPGstructureMtoptonu=u*}, we propagate these estimates to $\Mtop$ thanks to the linearized transport equations in the $e_3$ direction of section \ref{sec:linearizedequationsforingoingPGstructures} for ingoing PG structures. Recalling that $\Ab$ has already been estimated in Theorem M2, the other 
 quantities are recovered following the same scheme\footnote{Note that the steps 21--25 in section \ref{sec:decayestimatesMintmainstepsorder} are not needed in the case of $\Mtop$ in view of Remark \ref{rmk:noJpmandnoJkpminMtop:chap7}.} as the one for $\Mint$ outlined in section \ref{sec:decayestimatesMintmainstepsorder}. 

As the estimates are significantly simpler to derive\footnote{Along a level set of $\ub$ in $\Mtop$, denoting $r_+(\ub)$ the maximal value of ${}^{(top)}r$, i.e. the one on $\{u=u_*\}$, and $r_-(\ub)$ the minimal value of ${}^{(top)}r$, i.e. the one on ${}^{(top)}\Si$, we have $0<r_+(\ub)-r_-(\ub)\les 1$, see \eqref{eq:upperboundrpubminusrmubonMtop:chap3}. In particular, the integration along $e_3$ is always finite in $\Mtop$.} and in the same spirit as the corresponding ones in Theorem M4, we leave the details to the reader. This concludes the proof of Theorem M5 for $\Mtop$. Together with the proof of Theorem M5 for $\Mint$ in section \ref{sec:decayestimatesMintmainstepsorder}, this concludes the proof of Theorem M5.

%%%%%%%%%%%%%%%%%%%%%%%%%%%%%%%%%%%%%%%%

\chapter{Initialization and extension (Theorems M0, M6 and M7)}
\lab{chapter:initializationandextension}

%%%%%%%%%%%%%%%%%%%%%%%%%%%%%%%%%%%%%%%%

The goal of this chapter is to prove Theorems M0, M6, and M7. To this end, we first review our GCM procedure in section \ref{sec:GCMpapersreview}, and construct an auxiliary outgoing geodesic foliation in the part $\Lext$ of the initial data layer in section \ref{section:geodesicfoliation8}.  Theorems M0, M6, and M7 are then proved respectively in sections \ref{sec:proofofTheoremM0}, \ref{sec:proofofTheoremM6} and \ref{sec:proofofThmM7}.

%%%%%%%%%%%%%%%%%%%

\section{GCM procedure}
\lab{sec:GCMpapersreview}

%%%%%%%%%%%%%%%%%%%%

In this section, we review the main results on the existence of GCM spheres in \cite{KS-GCM1} \cite{KS-GCM2} and on the existence of GCM hypersurfaces in \cite{Shen}. These results will be used repeatedly in the proof of Theorems M0, M6 and M7.

%%%%%%%%%%%%%%%%%%%

\subsection{Background spacetime} 

%%%%%%%%%%%%%%%%%%%

Our GCM results  hold true in a vacuum    spacetime region, denoted $\RR$, foliated by two functions $(u,s)$ such that:
\begin{enumerate}
\item On $\RR$, $(u, s)$  is a geodesic foliation  of lapse $\vsi$, i.e. 
\begin{itemize} 
\item[-] $ u$ is an optical function, $L=-\g^{\a\b}\pr_\b u \pr_\a u$, and $L(\vsi)=1$,
\item[-]  $e_4=\vsi L$ and $ L(s)=1$,
\item[-]   $e_3$ is the correspond null companion  to $ e_4$, perpendicular to the surfaces $S(u, s)$ induced by  the level surfaces of $(u, s)$, and such that $\g(e_3, e_4)=-2$.
\end{itemize}
In particular, it follows  from the above that 
\beaa
\atrch=\atrchb=0, \qquad \om=\xi=0,  \qquad  \etab = -\ze, \qquad  \vsi=\frac{2}{e_3(u)}.
\eeaa

\item We define the following renormalized quantities
\beaa
\bsplit
\widecheck{\trch} &:= \trch -\frac{2}{r}, \qquad \widecheck{\trchb} := \trchb +\frac{2\Up}{r},\qquad\,\,\, \widecheck{\omb} := \omb -\frac{m}{r^2},\\
\widecheck{K} &:= K -\frac{1}{r^2},\qquad\quad\,\,\,\widecheck{\rho} := \rho +\frac{2m}{r^3},\qquad \quad\,\,\,
\widecheck{\mu} := \mu -\frac{2m}{r^3}, \\
\widecheck{\Omb} &:=\Omb+\Up, \qquad\qquad\, \widecheck{\varsigma} := \varsigma-1,
\end{split}
\eeaa
where 
\beaa
\Omb:=e_3(s), \qquad \Up :=1-\frac{2m}{r},
\eeaa
and group them in the sets $\Ga_g$ and $\Ga_b$ defined as follows
\bea
\lab{definition:Ga_gGa_b:backgroundfoliationGCM}
\bsplit
\Ga_g &:= \Bigg\{\widecheck{\trch},\,\, \chih, \,\, \ze, \,\, \widecheck{\trchb},\,\,  r\widecheck{\mu} ,\,\,  r\widecheck{\rho}, \,\, r\dual\rho, \,\, r\b, \,\, r\a, \,\, r\widecheck{K}, \,\, r^{-1} \big(e_4(r)-1\big), \,\, r^{-1}e_4(m)\Bigg\},\\
\Ga_b &:= \Bigg\{\eta, \,\,\chibh, \,\, \ombc, \,\, \xib,\,\,  r\bb, \,\, \aa, \,\, r^{-1}\widecheck{\Omb}, \,\,r^{-1}\widecheck{\varsigma}, \,\, r^{-1}(e_3(r)+\Up\big), \,\, r^{-1}e_3(m)  \Bigg\}.
\end{split}
\eea

\item Let $(\ug, \sg)$ two real numbers. Let $\ovS := S(\ug, \sg),$
$\rg$ the area radius of $\ovS$, and $\mg$ the Hawking mass of $\ovS$.

\item $\RR$ is covered by two coordinates charts $\RR=\RR_N\cup \RR_S$ such that:
\begin{enumerate}
\item The North coordinate chart   $\RR_N$ is given by the coordinates
$(u, s, y_{N}^1, y_{N}^2)$ with    $(y^1_{N})^2+(y^2_{N})^2<2$. 

\item The South coordinate chart  $\RR_S$ is  given by the coordinates
$(u, s, y_{S}^1, y_{S}^2)$  with $(y^1_{S})^2+(y^2_{S})^2<2$. 

 \item   The two coordinate charts   intersect in  the  open equatorial region
 $\RR_{Eq}:=\RR_N\cap \RR_S$ in which both coordinate systems are defined.
 
 \item  In $\RR_{Eq}$, the transition functions  between the two coordinate  systems are given by  the smooth  functions $ \varphi_{SN}$ and $\varphi_{NS}= \varphi_{SN}^{-1} $. 
 \end{enumerate}
 
 \item The metric coefficients for the two coordinate systems   are given by
 \beaa
\g &=& - 2\vsi du ds + \vsi^2\Omb  du^2 +g^{N}_{ab}\big( dy_N^a- \vsi \undB_{N}^a du\big) \big( dy_N^b-\vsi \undB_N^b du\big),\\
\g &=& - 2\vsi du ds + \vsi^2\Omb  du^2 +g^{S}_{ab}\big( dy_S^a- \vsi \undB_{S}^a du\big) \big( dy_S^b-\vsi \undB_S^b du\big),
\eeaa
where
\beaa
\Omb=e_3(s), \qquad \undB_N^a =\frac{1}{2} e_3(y_N^a), \qquad \undB_S^a =\frac{1}{2} e_3(y_S^a).
\eeaa

\item We restrict the region $\RR$  such that,  for $\epg$ sufficiently small, $\rg$ the  area  radius  of $S(\ovu,\ovs)$ sufficiently large, i.e \, $\epg\ll m_0\ll \rg$,
\bea
\lab{definition:RR(dg,epg)}
\RR:=\left\{|u-\ug|\leq\epg,\quad |s-\sg|\leq  \epg \right\}.
\eea
\item We assume that on $\RR$ the following assumptions are verified, for an integer $s_{max}$ sufficiently large
\begin{enumerate}
\item[\bf A1.]
For  $k\le s_{max}$
\bea\lab{eq:assumtioninRRforGagandGabofbackgroundfoliation}
\bsplit
\| \Ga_g\|_{k, \infty}&\leq  \epg  r^{-2},\qquad \| \Ga_b\|_{k, \infty}&\leq  \epg  r^{-1}.
\end{split}
\eea

\item[\bf A2.]  The Hawking mass $m=m(u,s)$ of  $S(u, s)$ verifies 
\bea\lab{eq:assumtionsonthegivenusfoliationforGCMprocedure:Hawkingmass} 
\sup_{\RR}\left|\frac{m}{m_0}-1\right| &\leq& \epg.
\eea

\item[\bf A3.] 
In the  region of their respective validity\footnote{That is  the quantities on the left verify the  same estimates as those for $\Ga_b$, respectively $\Ga_g$.}   we have
\bea
 \undB_N^a,\,\, \undB_S^a \in r^{-1}\Ga_b, 
 \eea
 and
 \bea
 r^{-2} \widecheck{g}^{N}_{ab},  \,\, r^{-2} \widecheck{g}^{S}_{ab} \in r\Ga_g,
 \eea
 where
 \beaa
 \widecheck{g} ^{N}\!_{ab} &=&   g^N_{ab}   -   \frac{4r^2}{1+(y^1_{N})^2+(y^2_{N})^2) }\de_{ab},\qquad 
 \widecheck{g}^{S}\!_{ab} =   g^S_{ab}   -   \frac{4r^2}{(1+(y^1_{S})^2+(y^2_{S})^2) } \de_{ab}.
 \eeaa
\end{enumerate}

{\bf A4.} We assume  the existence of a   smooth family of  scalar functions $\Jp:\RR\to\RRR$, for $p=0,+,-$,   verifying the following properties
\begin{enumerate}
\item On the sphere $\ovS$ of the background foliation, there holds
 \bea
 \lab{eq:Jpsphericalharmonics}
\bsplit
  \Big((\rg)^2\lapzero+2\Big) \Jp  &= O(\epg),\qquad p=0,+,-,\\
\frac{1}{|\ovS|} \int_{\ovS}  \Jp J^{(q)} &=  \frac{1}{3}\de_{pq} +O(\epg),\qquad p,q=0,+,-,\\
\frac{1}{|\ovS|}  \int_{\ovS}   \Jp   &=O(\epg),\qquad p=0,+,-.
\end{split}
\eea

\item We extend $\Jp$ from $\ovS$ to $\RR$ by $\pr_s\Jp=\pr_u\Jp=0$, i.e.
\bea\lab{eq:extensionofJpfromovStoRR}
\Jp(u,s,y^1,y^2)=\Jp(\ug, \sg, y^1, y^2).
\eea  
\end{enumerate}
\end{enumerate}

%%%%%%%%%%%%%%%%%%%%%%%%%%%%%%%%%%%%%%

\subsection{Deformations of surfaces}

%%%%%%%%%%%%%%%%%%%%%%%%%%%%%%%%%%%%%%

We review the results in \cite{KS-GCM1} on deformations of surfaces that will be useful in this chapter.

 \begin{definition}
 \label{definition:Deformations}
 We say that    $\S$ is a  deformation of $ \ovS$ if there exist  smooth  scalar functions $U, S$ defined on $\ovS$ and a  map 
  a map $\Psi:\ovS\to \S $  verifying, on either coordinate  chart  $(y^1, y^2) $ of $\ovS$,  
   \bea
 \Psi(\ovu, \ovs,  y^1, y^2)=\left( \ovu+ U(y^1, y^2 ), \, \ovs+S(y^1, y^2 ), y^1, y^2  \right).
 \eea
 \end{definition}

\begin{definition}
\lab{definition:framechange.toadaptedframes}
Given a deformation $\Psi:\ovS\to \S$ we  say that 
 a new frame   $(e_3', e_4',  e_1', e_2')$ on $\S$, obtained from the standard frame $(e_3, e_4, e_1, e_2)$  via the general frame   transformation  \eqref{General-frametransformation},  is  $\S$-adapted  if   the horizontal  vectorfields $e'_1, e'_2$   are tangent to $\S$.
\end{definition}

The following result combines  Lemma 5.8,  Corollary 5.9 and  Corollary 5.17  in \cite{KS-GCM1}.
\begin{proposition}
\lab{proposition:58-59GCM1}
 Let $\ovS \subset \RR$.    Let  $\Psi:\ovS\to \S $  be   a  deformation generated by the  functions $(U, S)$ as in Definition \ref{definition:Deformations}. Assume the bound
 \bea
 \label{assumption-UV-dg}
   \| (U, S)\|_{L^\infty(\ovS)} +r  \| \nabzero(U, S)\|_{L^\infty(\ovS)}  +r^2  \| \nabzero^2(U, S)\|_{L^\infty(\ovS)} &\les&  \dg.
 \eea
  Then
   \begin{enumerate}
 \item   We have
 \bea
 \lab{equation:difference-ofgas}
\sum_{a,b=1}^2  \big| g^{\S, \#}_{ab} -\ovg_{ab} \big|&\les r\dg.
 \eea 

\item   We have
  \bea
 \frac{r^\S}{\ovr}= 1 + O(r ^{-1}  \dg )
 \eea
 where $r^\S$ is the area radius of $\S$ and $\ovr$ that of $\ovS$.
 \item  Let $m=m(u,s)$  the Hawking mass of  the surfaces $S(u, s)$ and      $m^S$       the Hawking mass of $S$. We have
   \beaa
\sup_{\S}|m-m^\S| \les \dg.
 \eeaa
 
 \item   For an arbitrary scalar function $f$ on $\RR$,
\beaa
\left|\int_\S f -\int_{\ovS} f\right| &\les&  \dg  \,  \rg \,\left(    \sup_{\RR}|f|+\rg\sup_{\RR}\big(|\pr_uf|+|\pr_sf|\big)\right).
\eeaa
\end{enumerate}
\end{proposition}

The following results combine Lemma   7.3  with  Corollary  7.7 in  \cite{KS-GCM1}.

\begin{lemma}\lab{lemma:consequencedeformationsurfaceusedinTheoremM0}
Let $\Psi:\ovS\to \S $ be a deformation defined by  $(U, S)$, as in Definition \ref{definition:Deformations}  with         $(f, \fb)$  the   transition function of  the frame transformation 
 from the  frame of $\RR$     to that adapted to $\S$, as in Definition \ref{definition:framechange.toadaptedframes}.
 There exists a small enough   constant $\de_1$ such that for given  $f, \fb$ on $\RR$   satisfying
\beaa
\|f\|_{\hk_{s_{max}}(\S)}+(r^\S)^{-1}\|\fb\|_{\hk_{s_{max}}(\S)} &\leq& \de_1,
\eeaa
the following holds
\begin{enumerate}
\item We have
\beaa
(\ovr)^{-1} \|U\|_{\hk_{s_{max}+1}(\ovS)} + (\ovr)^{-2} \|S\|_{\hk_{s_{max}+1}(\ovS)}  &\les& \de_1.
\eeaa 
In particular, we have
\beaa
\sup_{\S}|u-\ug| \les \de_1, \qquad \sup_{\S}|s-\sg| \les \rg\de_1.
\eeaa

\item  We have,
\beaa
\sum_{a,b=1}^2  \big| g^{\S, \#}_{ab} -\ovg_{ab} \big|&\les r^2\de_1.
 \eeaa

\item We have 
  \beaa
\left| \frac{r^\S}{\rg}-1\right|+\sup_{\S}\left| \frac{r^\S}{r}-1\right|\les \de_1. 
 \eeaa

\item The following estimate holds true for an arbitrary scalar function $h$ on $\RR$,
\beaa
\left|h^\# - h\right| &\les& \de_1   \sup_{\RR}   |\dk h|.
\eeaa
\item The following estimate holds true for an arbitrary scalar function $h$ on $\RR$,
\beaa
\left|\int_\S h -\int_{\ovS} h\right| &\les& \de_1 (\rg)^2\left(\sup_{\RR}\big(|f|+|\pr_uf|\big)+\rg\sup_{\RR}|\pr_sf|\right).
\eeaa

 \item  If $V\in \hk_s(\S)$ and $V^\#$ is its pull-back by $\Psi$, we have for all $0\leq s\leq s_{max}$,
 \bea
 \lab{eq:Prop:comparison1:forTheoremM0}
 \|V\|_{\hk_s(\S)}= \|V^\#\|_{\hk_s(\ovS,\, g^{\S,\#})} = \| V^\#\|_{\hk_s(\ovS, \ovg)}\big(1+O(\de_1)\big).
 \eea

 \item  For any  tensor $h$  on $\RR$
 \bea
 \lab{eq:Prop:comparison2:forTheoremM0}
 \|h\|_{\hk_s(\S)} \les  r \sup_{\RR}\big|\dk^{\leq s}h\big|, \qquad 0\leq s \leq s_{max}.
 \eea

 \item We have
\bea
\sum_{a,b,c=1,2}\Big\|(\Ga^{\S, \#})_{ab}^c-(\ovGa)_{ab}^c\Big\|_{\hk_{s_{max}-1}(\ovS)} &\les r^2\de_1.
\eea

\item We also have, for $\mg$ the Hawking mass of $\ovS$,
\beaa
|m^\S -\mg|\les \de_1+(\epg)^2.
\eeaa
\end{enumerate}
\end{lemma}

%%%%%%%%%%%%%%%%%%%%%%%%%%%%%%%%%%

\subsection{Existence of intrinsic GCM  spheres}
\lab{section:RecallresultsGCM}

%%%%%%%%%%%%%%%%%%%%%%%%%%%%%%%%%%

We review in this section the results of \cite{KS-GCM2} useful for this chapter. We start  with the following definition of canonical $\ell=1$ modes on a deformed sphere $\S$.
   \begin{definition}
    \lab{definition:canmodesS}   
    Given a  deformation map $\Psi:\ovS\to \S$   and a fixed effective uniformization map
       $(\ovPhi, \ovphi)$ for $\ovS$  we let $(\Phi, \phi) $ be the  unique effective  uniformization map of  $\S$ calibrated 
       with $(\ovPhi, \ovphi)$ relative to the map $\Psi$, in the sense of Definition \ref{def:calibration}. With this choice, we define the canonical $\ell=1$  modes of $\S$  by the formula
       \bea
    \lab{def:JpSS}
    \JpSS &=&J^{(p, \SSS^2)} \circ \Phi^{-1} 
    \eea        
  with $J^{(p, \SSS^2)}$  denoting the $\ell=1$ spherical harmonics of $\SSS^2$.
  \end{definition}

Consider as before a vacuum   spacetime region  $\RR$ verifying the assumptions {\bf A1-A4}. In addition we make the following stronger assumptions on {\bf A1} and {\bf A4}.

{\bf A1-Strong.}  For $k\le s_{max}$, and for a small enough constant $\de_1>0$, with $\de_1\geq\epg$,
\bea
\lab{eq:assumtioninRRforGagandGabofbackgroundfoliation'}
\Big\|(\Ga_g, \Ga_b) \Big\|_{k, \infty}  \les \de_1 r^{-2}, \qquad \Big\|\nab_3\Ga_g\Big\|_{k, \infty}  \les \de_1 r^{-3}.
\eea

{\bf  A4-Strong.}  We assume  the existence of a   smooth family of  scalar functions $\Jp:\RR\to\RRR$, for $p=0,+,-$,   verifying the following properties
\begin{enumerate}
\item On the sphere $\ovS$ of the background foliation, there holds
the following stronger version of \eqref{eq:Jpsphericalharmonics}
 \bea
 \lab{eq:Jpsphericalharmonics-strong}
\bsplit
  \Big((\rg)^2\lapzero+2\Big) \Jp  &= O(\epg r^{-1} ),\qquad p=0,+,-,\\
\frac{1}{|S|} \int_{S}  \Jp J^{(q)} &=  \frac{1}{3}\de_{pq} +O(\epg r^{-1} ),\qquad p,q=0,+,-,\\
\frac{1}{|S|}  \int_{S}   \Jp   &=O(\epg r^{-1}),\qquad p=0,+,-.
\end{split}
\eea

\item On  $\ovS$  we have
\bea\lab{eq:JpisclosedtocanonicaloneofovS:stronger}
\max_{p=0,-,+}\|\Jp-J^{(p,\ovS)}\|_{\hk_{s_{max}+1}(\ovS)}\les \dg,
\eea
where $J^{(p,\ovS)}$ denotes the canonical basis of $\ell=1$ modes on $\ovS$ corresponding to the effective uniformization map $(\ovPhi, \ovphi)$ for $\ovS$ appearing in Definition  \ref{definition:canmodesS}.

\item We extend $\Jp$ from $\ovS$ to $\RR$ as in \eqref{eq:extensionofJpfromovStoRR}, i.e. by $\pr_s\Jp=\pr_u\Jp=0$.
\end{enumerate}

We  state below the results  of Corollary 7.2. in \cite{KS-GCM2}.

\begin{corollary}\lab{Lemma:ComparisonJ-strong}
 Let $\Jp$ satisfying {\bf  A4-Strong}, and let $J^{(p,\ovS)}$ denotes the canonical basis of $\ell=1$ modes on $\ovS$ corresponding to the effective uniformization map $(\ovPhi, \ovphi)$ for $\ovS$ appearing in Definition  \ref{definition:canmodesS}.  Let a deformation $\Psi:\ovS\to \S$  with the corresponding deformation parameters $(U, S)$  satisfying
\bea
\lab{eq:ThmGCMS4}
 \|( U, S)\|_{\hk_{s_{max}+1}(\ovS)}  &\les&  r \dg, \qquad  s_{max}\geq 2. 
 \eea
    Let    $\JpSS$ be the  corresponding  canonical  basis  of $\ell=1$ modes of $\S$   calibrated  according to  Definition \ref{definition:canmodesS}. Then, the following estimate holds true
   \bea\lab{eq:comparsonsJs}
  \max_{p=0,-,+} \sup_{\S}|\Jp - \JpS| \les r^{-1}\dg.
   \eea
   \end{corollary}
   
We state below Theorem  7.3. in \cite{KS-GCM2}, which is the main result in that paper,  concerning the construction of  intrinsic GCM  spheres.  Recall, see Definition \ref{Definition:ell=1modesofascalarfunction}, that the  $\ell=1$  modes   of a scalar function $f$ 
are defined to be the triplet\footnote{Note that  this definition differs from the one in \cite{KS-GCM2}  by a factor of $r^{-2}$.}
\beaa
(f)_{\ell=1} &=& \left(\frac{1}{r^2}\int_{S}fJ^{(0)}, \frac{1}{r^2}\int_{S}fJ^{(+)}, \frac{1}{r^2}\int_{S}fJ^{(-)}\right).
\eeaa
Theorem  7.3. in \cite{KS-GCM2} holds for spacetime  regions $\RR$  verifying, in addition to  {\bf A1--A4}, 
\bea\lab{eq:thedecompositionofkakabandmuintermsofkadotkabdotandmudotplusmodes}
\bsplit
\ka&=\frac{2}{r}+\dot{\ka},\\
\kab&=-\frac{2\Up}{r} +  \Cb_0+\sum_p \Cbp \Jp+\dot{\kab},\\
\mu&= \frac{2m}{r^3} + M_0+\sum _p\Mp \Jp+\dot{\mu},
\end{split}
\eea
where the scalar functions $\Cb_0=\Cb_0(u,s)$, $\Cbp=\Cbp(u,s)$, $M_0=M_0(u,s)$ and $\Mp=\Mp(u,s)$, defined on the spacetime region $\RR$, depend only on the coordinates $(u,s)$, and where $\dot{\ka}$, $\dot{\kab}$ and $\dot{\mu}$ satisfy the following estimates
\bea
\lab{eq:GCM-improved estimate2-again}
\sup_{\RR}\big|\dk^{\leq s_{max}}(\kadot, \kabdot)|\les r^{-2}\dg,\qquad\qquad 
 \sup_{\RR}\big|\dk^{\leq s_{max}}\mudot| \les r^{-3}\dg.
\eea

\begin{theorem}[Existence of intrinsic GCM spheres]
\lab{theorem:ExistenceGCMS2} 
Assume that the spacetime region  $\RR$ verifies the assumptions {\bf A1-Strong, A2, A3}, {\bf A4-Strong}, as well as \eqref{eq:thedecompositionofkakabandmuintermsofkadotkabdotandmudotplusmodes} \eqref{eq:GCM-improved estimate2-again}. We further  assume that, relative to the   $\ell=1$ modes of  the background foliation,
\bea
\lab{Assumptions:theorem-ExistenceGCMS2}
(\div \b )_{\ell=1}=O( \dg r^{-5}), \qquad        (\widecheck{\trch})_{\ell=1}=O(\dg r^{-3}),   \qquad (\widecheck{\trchb})_{\ell=1}=O(\dg r^{-3}).
\eea
Then,  there  exist unique constants  $\MpS$,   $p\in\{-,0, +\}$, such that
 \bea
\lab{def:GCMC2}
\bsplit
\ka^\S&=\frac{2}{r^\S},\\
\kab^\S &=-\frac{2}{r^\S}\Up^\S,\\
\mu^\S&= \frac{2m^\S}{(r^\S)^3} +\sum _p\MpS \JpSS, 
\end{split}
\eea
and
\bea
\lab{def:GCMC2-b}
\int_\S \div^\S \b^\S  \JpSS  =0,
\eea
where $ \JpSS  $ is a  canonical  $\ell=1$ basis for $\S$ calibrated,  relative  by $\Psi $, with the canonical    $\ell=1$ basis  of $\ovS$.  Moreover the deformation  verifies the properties:
\begin{enumerate}
\item The volume radius $r^\S$  verifies
\beaa
\left|\frac{r^\S}{\rg}-1\right|\les  r^{-1} \dg.
\eeaa
\item  The parameter  functions $U, S$  of the deformation verify
\beaa
 \|( U, S)\|_{\hk_{s_{max}+1}(\ovS)}  &\les&  r \dg.
 \eeaa

\item The Hawking mass  $m^\S$  of $\S$ verifies the estimate
\beaa
 \big|m^\S-\ovm\big|&\les &\dg. 
 \eeaa
 \item The well defined\footnote{Note  that  while  the Ricci coefficients $\ka^\S, \kab^\S,  \chih^\S, \chibh^\S, \ze^\S$ as well as all curvature  components  and   mass aspect function $\mu^\S$     are well defined on $\S$, this in not the case  of $\eta^\S, \etab^\S, \xi^\S, \xib^\S, \om^\S, \omb^\S$ which require  the  derivatives of the frame in the $e_3^\S$ and $e_4^\S$ directions.} Ricci and curvature coefficients of $\S$  verify,
   \beaa
\bsplit
\| \Ga^\S_g\|_{\hk_{s_{max} }(\S) }&\les  \epg  r^{-1},\qquad 
\| \Ga^\S_b\|_{\hk_{s_{max} }(\S) }\les  \epg.
\end{split}
\eeaa
\end{enumerate}
\end{theorem}

The following    corollary is   Corollary  7.7 in \cite{KS-GCM2}.
  \begin{corollary}
  \lab{Corr:ExistenceGCMS2}
  Under the same assumptions as in Theorem \ref{theorem:ExistenceGCMS2}      we have,  in addition to   \eqref{def:GCMC2} and   \eqref{def:GCMC2-b}, 
  \begin{itemize}
  \item   either, for any choice of a canonical ${\ell=1}$ basis of $\S$, 
  \beaa
  (\curl^\S \b^\S)_{\ell=1}=0,
  \eeaa

  \item   or there exists a canonical basis of  $\ell=1$ modes of $\S$ such that
   \bea
   \lab{angular-momentum}
   \int_\S  \curl^\S \b^\S\, J^{(\pm, \S)}=0, \qquad  \int_\S  \curl^\S \b^\S\, J^{(0, \S)}\neq 0.
   \eea 
  \end{itemize} 
  We then define   the angular   parameter $a^\S$ on $\S$ by the  formula 
   \bea
 a^\S:=\frac{(r^\S)^3}{8\pi m^\S}   \int_\S  \curl^\S  \b^\S J^{(0, \S)}.  
   \eea  
With this definition, we have $a^\S=0$ in the first case, while $a^\S\neq 0$ in the second case.
  \end{corollary}

Finally we state below the results of Proposition 8.1 in \cite{KS-GCM2}. 
\begin{proposition}
\lab{corr:GCM-rigidity2}
  Let  a fixed  spacetime region 
 $\RR$    verifying assumptions {\bf A1--A4} and \eqref{eq:thedecompositionofkakabandmuintermsofkadotkabdotandmudotplusmodes} \eqref{eq:GCM-improved estimate2-again}, as well as, for any background sphere $S$ of $\RR$,
 \bea\lab{eq:thisisthesmalladditionGCMconditionforGCMrigidity2background}
 |(\div\b)_{\ell=1}|\les r^{-4}\dg, \qquad |(\kabc)_{\ell=1}|\les r^{-2} \dg.
 \eea 
Assume that $\S$ is a deformed  sphere in $\RR$ which verifies   the   GCM conditions
\bea\lab{eq:GMCconditionsrepeatedforthenthtime:andonemodetime}
\bsplit
\ka^\S&=\frac{2}{r^\S},\\
\kab^\S &=-\frac{2}{r^\S}\Up^\S+  \Cb^\S_0+\sum_p \CbpS \widetilde{J}^{(p)},\\
\mu^\S&= \frac{2m^\S}{(r^\S)^3} +   M^\S_0+\sum _p\MpS \widetilde{J}^{(p)},
\end{split}
\eea
for some basis\footnote{$\widetilde{J}^{(p)}$ is not assumed to be a canonical basis of $\ell=1$ modes on $\S$.} of $\ell=1$ modes $\widetilde{J}^{(p)}$ on $\S$, such that for a small enough constant $\de_1>0$,
\begin{itemize}
\item  
 The transition coefficients  $(f, \fb, \la)$ from the background frame  of $\RR$     to that of $\S$  verifies, for some $4\leq s\leq s_{max}$, the bound
\bea\lab{eq:boundforffbasumptioncorrigidity2}
\|f\|_{\hk_{s}(\S)}+(r^\S)^{-1}\|(\fb, \ovla)\|_{\hk_{s}(\S)} &\leq& \de_1,
\eea
\item  
 The difference between the basis of $\ell=1$ modes $\widetilde{J}^{(p)}$ on $\S$ and the basis of $\ell=1$ modes of the background foliation $\Jp$ verifies 
\bea\lab{eq:boundfordifferenceofell=1modescorrigidity2}
r^{-1}\|\widetilde{J}^{(p)}-\Jp\|_{\hk_{s}(\S)} &\leq& \de_1.
\eea
\end{itemize}

Assume in addition that we have, with respect to the basis of $\ell=1$ modes $\widetilde{J}^{(p)}$ on $\S$, 
\bea\lab{eq:thisisthesmalladditionGCMconditionforGCMrigidity2deformedsphere}
 |(\div^\S\b^\S)_{\ell=1}|\les r^{-4}\dg, \qquad |(\kabc^\S)_{\ell=1}|\les  r^{-2}\dg.
\eea
Then  $(f, \fb, \la=1+\ovla)$ verify the  estimates
 \beaa
\|(f, \fb, \check{\ovla}^\S)\|_{\hk_{s+1}(\S) } &\les& r\dg+r(\epg)^2+r\de_1\left(\frac{1}{\rg}+\epg+ \de_1  \right)
\eeaa
and 
\beaa
r|\ov{\ovla}^\S| &\les& r\dg+r(\epg)^2+r\de_1\left(\frac{1}{\rg}+\epg+ \de_1  \right)+\sup_\S|r-r^\S|.
\eeaa
\end{proposition}

%%%%%%%%%%%%%%%%%%%%%%%%%%%%%%%%%%%%%%

\subsection{Existence of GCM  hypersurfaces}

%%%%%%%%%%%%%%%%%%%%%%%%%%%%%%%%%%%%%%

In this section, we review the results on the construction  of GCM hypersurfaces (GCMH) from \cite{Shen}. We state below Theorem 4.1 in \cite{Shen} on the construction of GCM hypersurfaces.

\begin{theorem}[Construction of GCMH]
\lab{theorem:constuctionGCMH}
Assume that the spacetime region  $\RR$ verifies the assumptions {\bf A1--A4}.
We further  assume that, relative to the   $\ell=1$ modes of  the background foliation,
\bea
\sup_{\RR} r\left|\dk^{\leq s_{max}}e_3(\Jp)\right| &\les& \dg.
\eea
In addition, we assume on $\RR$
\bea
|(\div\eta)_{\ell=1}|\les \dg, \quad |(\div\xib)_{\ell=1}|\les \dg,
\eea
as well as the existence of a constant $m^{(0)}$ such that we have on $\RR$
\bea
\left|(e_3(u)+e_3(s))\Big|_{SP} -1-\frac{2m^{(0)}}{r}\right| &\les& \dg,
\eea
where $SP$ denotes the South Pole, i.e. $y^1_S=y^2_S=0$ w.r.t. the South coordinates $(u,s,y^1_S, y^2_S)$. 

Let $\S_0$ be a fixed sphere included in the region $\RR$, let a pair of triplets $\La_0, \Lab_0\in\RRR^3$  such that 
\bea
|\La_0|+|\Lab_0|\les r^{-2}\dg, 
\eea
and let $\Jp[\S_0]$ a basis of $\ell=1$ modes on $\S_0$, such that we have on $\S_0$
\bea
\bsplit
\ka^{\S_0}&=\frac{2}{r^{\S_0}},\\
\kab^{\S_0} &=-\frac{2\Up^{\S_0}}{r^{\S_0}}+  \Cb^{\S_0}_0+\sum_p \Cb^{(\S_0,p)}\Jp[\S_0],\\
\mu^{\S_0}&= \frac{2m^\S}{(r^{\S_0})^3} +   M^{\S_0}_0+\sum _pM^{(p,\S_0)}\Jp[\S_0],
\end{split}
\eea 
as well as 
\bea
(\div f_0)_{\ell=1}=\La_0, \qquad (\div\fb_0)_{\ell=1}=\Lab_0,
\eea
with $(f_0, \fb_0)$ corresponding to the coefficients from the background frame to the frame adapted to $\S_0$, and the $\ell=1$ modes being taken w.r.t. the basis $\Jp[\S_0]$, and where 
\beaa
\|\Jp[\S_0]-\Jp\|_{\hk_{s_{max}+1}(\S_0)} &\les& r\dg. 
\eeaa
Then, there exists a unique, local, smooth, space like hypersurface $\Si_0$ passing through $\S_0$, a scalar function $u^\S$ defined on $\Si_0$, whose level sets a topological spheres  denoted by $\S$, a smooth collection of constants $\La^\S, \Lab^\S$ and a triplet of functions $\Jp[\S]$ defined on $\Si_0$ verifying 
\beaa
\La^{\S_0}=\La_0, \qquad \Lab^{\S_0}=\Lab_0, \qquad \Jp[\S]\Big|_{\S_0}=\Jp[\S_0],
\eeaa
such that the following conditions are verified: 
\begin{enumerate}
\item The following GCM conditions hold on $\Si_0$
\bea
\lab{def:GCMC:onGCMH}
\bsplit
\ka^\S&=\frac{2}{r^\S},\\
\kab^\S &=-\frac{2\Up^\S}{r^\S}+  \Cb^\S_0+\sum_p \CbpS\Jp[\S],\\
\mu^\S&= \frac{2m^\S}{(r^\S)^3} +   M^\S_0+\sum _p\MpS\Jp[\S].
\end{split}
\eea 

\item There exists a constant $c_0$ such that
\bea
u^\S+r^\S=c_0\quad\textrm{along}\quad\Si_0.
\eea 

\item Let $\nu^\S$ the unique vectorfield tangent to the hypersurfaces $\Si_0$, normal to $\S$, and normalized by $\g(\nu^\S, e_4^\S)=-2$, and let $b^\S$ be the unique scalar function on $\Si_0$ such that $\nu^\S$ is given by
\bea
\nu^\S=e_3^\S+b^\S e_4^\S.
\eea 
Then, the following normalization condition holds true at the south pole SP of every sphere $\S$, i.e. at $y^1_S=y^2_S=0$, 
\bea
b^\S\Big|_{SP}=-1-\frac{2m^{(0)}}{r^\S}.
\eea

\item The triplet of functions $\Jp[\S]$ verifies on $\Si_0$
\bea
\nu^\S(\Jp[\S])=0, \quad p=0,+,-.
\eea

\item The following transversality conditions are verified on $\Si_0$
\bea\lab{eq:transversalityconditionTheoremGCMHDawei}
\xi^\S=0, \qquad \om^\S=0, \qquad \etab^\S+\ze^\S=0, \qquad e_4^\S(u^\S)=0, \qquad e_4^\S(r^\S)=1. 
\eea

\item In view of \eqref{eq:transversalityconditionTheoremGCMHDawei}, the Ricci coefficients $\eta^\S$ and $\xib^\S$ are well defined on $\Si_0$. They verify on $\Si_0$
\bea
(\div^\S\eta^\S)_{\ell=1}=0, \qquad (\div^\S\xib^\S)_{\ell=1}=0,
\eea
where the $\ell=1$ modes are taken w.r.t. $\Jp[\S]$.

\item The transition coefficients $(f, \fb, \la)$ from the background foliation to that of $\Si_0$ verify
\bea
\|\dk^{\leq s_{max}+1}(f, \fb, \la-1)\|_{L^2(\S)} \les \dg.
\eea 
\end{enumerate}
\end{theorem}

We state below Corollary 4.2 in \cite{Shen}.

\begin{corollary}\lab{cor:constuctionGCMH}
Assume that the spacetime region  $\RR$ verifies the assumptions {\bf A1--A4} and the small GCM conditions \eqref{eq:GCM-improved estimate2-again}. Assume given a GCM hypersurface $\Si_0\subset\RR$ foliated by hypersurfaces $\S$ such that 
\bea
\bsplit
\ka^\S&=\frac{2}{r^\S},\\
\kab^\S &=-\frac{2}{r^\S}\Up^\S+  \Cb^\S_0+\sum_p \CbpS\Jp[\S],\\
\mu^\S&= \frac{2m^\S}{(r^\S)^3} +   M^\S_0+\sum _p\MpS\Jp[\S],
\end{split}
\eea 
and 
\bea
(\div^\S\eta^\S)_{\ell=1}=0, \qquad (\div^\S\xib^\S)_{\ell=1}=0,
\eea
where the triplet of functions $\Jp[\S]$ verifies on $\Si_0$ 
\bea
\nu^\S(\Jp[\S])=0, \quad p=0,+,-.
\eea

\begin{enumerate}
\item If we assume in addition that for a given sphere $\S_0$ on $\Si_0$ the transition coefficients $(f, \fb, \la)$ from the background foliation to $\S_0$ verify 
\bea
\|(f, \fb, \la-1)\|_{\hk_{s_{max}+1}(\S_0)} &\les& \dg,
\eea
then 
\bea
\|\dk^{\leq s_{max}+1}(f, \fb, \la-1)\|_{L^2(\S_0)} &\les& \dg.
\eea

\item  If we assume in addition that for a given sphere $\S_0$ on $\Si_0$ the transition coefficients $(f, \fb, \la)$ from the background foliation to $\S_0$ verify 
\bea
\|f\|_{\hk_{s_{max}+1}(\S_0)}+(r^{\S_0})^{-1}\|(\fb, \la-1)\|_{\hk_{s_{max}+1}(\S_0)} &\les& \dg,
\eea
then 
\bea
\nn\|\dk^{\leq s_{max}+1}f\|_{L^2(\S_0)}+(r^{\S_0})^{-1}\|\dk^{\leq s_{max}+1}(\fb, \la-1)\|_{L^2(\S_0)}\\
+\|\dk^{\leq s_{max}}\nu^\S(\fb, \la-1)\|_{L^2(\S_0)} &\les& \dg.
\eea
\end{enumerate}
\end{corollary}

%%%%%%%%%%%%%%%%%%%%%%%%%%%%%%%%%%%%%%%%

\section{An auxiliary geodesic foliation in $\Lext$}
\lab{section:geodesicfoliation8}

%%%%%%%%%%%%%%%%%%%%%%%%%%%%%%%%%%%%%%%%

Recall from section \ref{sec:defintionoftheinitialdatalayer} that the initial data layer $\LL_0$ is given by 
$\LL_0=\Lint\cup\Lext$, with $\Lint$ and $\Lext$ covered by PG structures. The goal of this section is to construct and control an auxiliary outgoing geodesic foliation in $\Lext$ 
that will be used in the proof of Theorem M0 and Theorem M6. To this end:
\begin{itemize}
\item We recall basic properties of the outgoing PG structure of $\Lext$ in section \ref{sec:auxiliaryoutgoing:preliminaries}.

\item We construct an auxiliary outgoing geodesic foliation in $\Lext$ in section \ref{sec:auxiliaryoutgoingutandst:construction}.

\item We use the transformation formulas to compare the Ricci coefficients and curvature components of the auxiliary outgoing geodesic foliation to the Ricci coefficients and curvature components of the outgoing PG structure of $\Lext$ in section \ref{sec:auxiliaryoutgoingRiccicoeeffandcurvature:identties}.

\item Finally, we control the auxiliary outgoing geodesic foliation of $\Lext$ in section \ref{sec:auxiliaryoutgoing:finalcontrol}, see Proposition \ref{proposition:geodesicfoliationLextt}.
\end{itemize}

%%%%%%%%%%%%%%%%%%%%%

\subsection{Preliminaries}
\lab{sec:auxiliaryoutgoing:preliminaries}

%%%%%%%%%%%%%%%%%%%%%%

In section \ref{section:geodesicfoliation8},  we concentrate  to  the region $\Lext$. To ease notations, throughout section \ref{section:geodesicfoliation8}, we denote:
\begin{itemize}
\item $(a, m)$ instead of $(a_0, m_0)$,

\item by $(u, r, \th, \vphi)$ the PG coordinates of $\Lext$,

\item by $E=(e_1, e_2, e_3, e_4)$ the outgoing PG frame of $\Lext$,

\item by $\chi$, $\chib$, $\ze$, $\eta$, $\etab$, $\xi$, $\xib$, $\om$, $\omb$, and $\a$, $\aa$, $\b$, $\bb$, $\rho$, $\rhod$ the   Ricci  and curvature coefficients of the outgoing PG structure of $\Lext$.
\end{itemize}
Recall that the PG structure of $\Lext$ verifies the following identities
\beaa
e_4(r)=1, \qquad e_4(u)= e_4(\th)=e_4(\vphi)=0, \qquad  e_1(r) = e_2(r)=0,
\eeaa
as well as 
\beaa
\xi=0, \qquad \om=0, \qquad \etab=-\ze.
\eeaa
Moreover,   $\Lext$ is also endowed with a  complex  $1$-form $\Jk$ verifying 
\beaa
\nab_{e_4} ( r \Jk)=0, \qquad  \dual \Jk =- i \Jk, \qquad\big|\Re(\Jk)\big|^2=\frac{(\sin\th)^2}{|q|^2}.
\eeaa

Also, we define the linearized quantities for the outgoing PG structure of $\Lext$ as in  Definition \ref{def:renormalizationofallnonsmallquantitiesinPGstructurebyKerrvalue}, and the corresponding quantities $\Ga_g, \Ga_b$ as in Definition \ref{definition.Ga_gGa_b}. Since $u$ is bounded in $\Lext $, we simply write $\Ga_g= r^{-1} \Ga_b$. Our initial data control, see \eqref{def:initialdatalayerassumptions}, implies in particular the bounds
 \bea\lab{def:initialdatalayerassumptions:chap8}
\, ^{(ext)} \Ik_{k_{large}+10} &\leq& \ep_0, \qquad   ^{(ext)} \Ik_{3} \leq  \ep^2_0,
\eea
 where, see section \ref{sec:initialdatalayernorm}, 
\beaa
\, ^{(ext)} \Ik_k &=& \sup_{\Lext} \Big\{   r|\dk^{\leq k}\Ga_b|+r^{\frac{7}{2}+\frac{\dt}{2}}\big(|\dk^{\leq k}A|+|\dk^{\leq k}B|\big)\Big\}.
\eeaa
 We make the additional assumption
 \bea\lab{eq:additionalassumption-8} 
 \widecheck{ e_3(u)} \in  \Ga_b.
 \eea
 
 \begin{remark}
 To justify the above additional assumption in $\Lext$, we recall the following equation, see Lemma \ref{Lemma:otherlinearizedquant},
  \bea
  e_4\left(\widecheck{e_3(u)}\right) &=& O(r^{-1})\Hc+O(r^{-1})\Zc+O(r^{-2})\widecheck{\DD u}+\Ga_b\c\Ga_b\in r^{-1} \Ga_b.
\eea
 Thus, the assumption \eqref{eq:additionalassumption-8}   follows from integrating backwards the transport equation  $e_4\left(\widecheck{e_3(u)}\right)\in r^{-1} \Ga_b$,  and from using  fact that  $\widecheck{e_3(u)}\to 0$ as $r\to \infty$.
 \end{remark}

 Finally, note that we have, in view of Definitions \ref{def:renormalizationofallnonsmallquantitiesinPGstructurebyKerrvalue} and \ref{definition.Ga_gGa_b},
 \beaa
&&  e_3(u)=  2+\frac{2a^2\sin^2\th}{r^2}  + \Ga_b +O(r^{-4}), \quad 
   \nab u= a \Re(\Jk) + \Ga_b,\\
&& e_3(r)= -\Up+ r\Ga_b+O(r^{-2} ),  \quad  \nab(\cos \th) =\nab(J^{(0)} )=-\Im (\Jk) +\Ga_b, \quad e_3(\cos \th)=\Ga_b,
  \eeaa
  where $\Up=1-\frac{2m}{r}$.

%%%%%%%%%%%%%%%%%%%%%%%%%%%%%%%%%%
  
 \subsection{Construction and asymptotic of the geodesic foliation}
 \lab{sec:auxiliaryoutgoingutandst:construction}

%%%%%%%%%%%%%%%%%%%%%%%%%%%%%%%%%%

We look for an optical function $\ut$ such that $\ut\sim u$  as $r\to\infty$. Its existence is provided by the following lemma.
 \begin{lemma}
 \lab{lemma:optical-Lext}
There exists a   unique  optical function $\ut$   defined  in   $\Lext$ and  verifying
\bea
\ut&=& u  -\frac{a^2(\sin\th)^2}{2r}+  h, \qquad  h=  \Ga_b +O(r^{-3}).
\eea
\end{lemma} 

\begin{proof}
Let 
\beaa
\ut_0:= u  -\frac{a^2(\sin\th)^2}{2r}. 
\eeaa
We calculate
 \beaa
  \g^{\a\b}\pr_\a u\pr_\b u&=&-   e_3( u) e_4 (u)+|\nab u|^2= |\nab u|^2 = \big|a \Re(\Jk) + \Ga_b\big|^2= \frac{a^2\sin^2\th}{|q|^2} +    r^{-1} \Ga_b \\
  &=&  \frac{a^2\sin^2\th}{r^2} +O(r^{-4})+ r^{-1} \Ga_b.
  \eeaa
Hence, since $e_4(u)=0$, and in view of the definition of $\ut_0$, we infer
\beaa
 &&\g^{\a\b}\pr_\a \ut_0 \pr_\b \ut_0\\
 &=&  \g^{\a\b}\pr_\a u\pr_\b u - 2 \g^{\a\b}\pr_\a u\pr_\b\left(\frac{a^2(\sin\th)^2}{2r}\right)+\g^{\a\b}\pr_\a\left(\frac{a^2(\sin\th)^2}{2r}\right)\pr_\b\left(\frac{a^2(\sin\th)^2}{2r}\right)\\
&=&  \frac{a^2\sin^2\th}{r^2} +O(r^{-4})+ r^{-1} \Ga_b+ e_3(u)e_4\left(\frac{a^2(\sin\th)^2}{2r} \right) -2\nab(u)\c\nab\left(\frac{a^2(\sin\th)^2}{2r} \right)\\
&& -e_4\left( \frac{a^2(\sin\th)^2}{2r} \right)e_3\left( \frac{a^2(\sin\th)^2}{2r} \right)+\left|\nab\left( \frac{a^2(\sin\th)^2}{2r} \right)\right|^2.
\eeaa
Note that
\beaa
&&e_3(u)e_4\left(\frac{a^2(\sin\th)^2}{2r} \right) -2\nab(u)\c\nab\left(\frac{a^2(\sin\th)^2}{2r} \right)\\ && -e_4\left( \frac{a^2(\sin\th)^2}{2r} \right)e_3\left( \frac{a^2(\sin\th)^2}{2r} \right)+\left|\nab\left( \frac{a^2(\sin\th)^2}{2r} \right)\right|^2\\
&=& -\Big(2 +\Ga_b +O(r^{-2})\Big)\frac{a^2(\sin\th)^2}{2r^2} +4\Big(a \Re(\Jk) + \Ga_b\Big)\c\frac{a^2\cos\th}{2r}\Big(-\Im (\Jk) +\Ga_b\Big)\\
&& +O(r^{-4})+r^{-2}\Ga_b\\
&=&  - \frac{a^2(\sin\th)^2}{r^2} + O(r^{-4} ) + r^{-2} \Ga_b
\eeaa
where we used the fact that $\Re(\Jk)\c\Im(\Jk)=0$ since $\Im(\Jk)=\dual\Re(\Jk)$. Hence, we infer
  \beaa
   \g^{\a\b}\pr_\a \ut_0 \pr_\b \ut_0&=&  O(r^{-4} ) + r^{-1} \Ga_b.
  \eeaa
   Thus, since $\ut=\ut_0+h$, we deduce
    \beaa
    \g^{\a\b}\pr_\a \ut \pr_\b \ut &=&   \g^{\a\b}\pr_\a (\ut_0+h)  \pr_\b( \ut_0+h)  =   2 \g^{\a\b} \pr_\a \ut_0 \pr_\b  h +             \g^{\a\b}\pr_\a  h \pr_\b  h
   +   O(r^{-4} ) + r^{-1} \Ga_b\\
   &=&- e_3(\ut_0)e_4(h) - e_4(\ut_0)e_3(h) +2\nab(\ut_0)\c\nab(h) +    \g^{\a\b}\pr_\a  h \pr_\b  h \\
   &&+   O(r^{-4} ) + r^{-1} \Ga_b.
    \eeaa
    Now, we have
    \beaa
    e_3(\ut_0) &=& e_3(u) +\frac{a^2\sin^2\th}{  2 r^2} e_3(r)+r^{-1}\Ga_b=2+O(r^{-2})  + \Ga_b,\\
     e_4(\ut_0) &=& e_4(u) +\frac{a^2\sin^2\th}{  2 r^2} = \frac{a^2\sin^2\th}{  2 r^2},\\
     \nab(\ut_0) &=& \nab(u)+O(r^{-2})+r^{-1}\Ga_b=a\Re(\Jk)+O(r^{-2})+\Ga_b,
    \eeaa
    and hence
    \beaa
     \g^{\a\b}\pr_\a \ut \pr_\b \ut &=&  -\Big(2+ O(r^{-2} )+ \Ga_b\Big) e_4 (h)  +O(r^{-2} ) e_3 (h) \\
     &&+2\Big(a\Re(\Jk)+O(r^{-2})+\Ga_b\Big)\c\nab(h)   +    \g^{\a\b}\pr_\a  h \pr_\b  h +  O(r^{-4} ) + r^{-1} \Ga_b.
    \eeaa
Since $\ut$ is an optical function, we are thus led to solve the following equation  for $h$
   \bea\lab{eq:reducedEikonal:equationh}
     \nn\Big(1+ O(r^{-2} )+ \Ga_b\Big) e_4 (h)  +O(r^{-2} ) e_3 (h) \\
     -\Big(a\Re(\Jk)+O(r^{-2})+\Ga_b\Big)\c\nab(h)   +    \g^{\a\b}\pr_\a  h \pr_\b  h &=&  O(r^{-4} ) + r^{-1} \Ga_b,
    \eea
    with the initialization 
    \bea\lab{eq:reducedEikonal:initilizationh}
    \lim_{r\to +\infty}h &=& 0.
    \eea
\eqref{eq:reducedEikonal:equationh} is a nonlinear transport equation for $h$. If $h=\Ga_b +O(r^{-3})$, $h$ satisfies in particular $|\dk^{\leq 1}h|\les \ep_0r^{-1}+r^{-3}$, and the uniqueness follows immediately from backward integration of \eqref{eq:reducedEikonal:equationh} and the initialization \eqref{eq:reducedEikonal:initilizationh}. Also, commuting  \eqref{eq:reducedEikonal:equationh}, using \eqref{def:initialdatalayerassumptions:chap8} to estimate the RHS of \eqref{eq:reducedEikonal:equationh}, integrating backward, and using the initialization \eqref{eq:reducedEikonal:initilizationh}, we easily obtain the following a priori bounds for $h$
 \beaa
    \big|\dk^{\leq k_{large}+10} h \big|  \les \ep_0 r^{-1}  +r^{-3}, \qquad \big|\dk^{\leq 3} h \big|  \les \ep_0 ^2 r^{-1}  +r^{-3}. 
\eeaa
 These bounds yield $h=\Ga_b +O(r^{-3})$, and can be used to prove the existence of $h$. Thus, there exists a  unique   solution $h$  of  \eqref{eq:reducedEikonal:equationh}   defined  on $\Lext$, such that $h\to 0$ as $r\to \infty$ and $h=\Ga_b +O(r^{-3})$. This concludes the proof of the lemma.
       \end{proof}

\begin{definition}
\lab{definition:geodesifolaition}
Let  $\ut$ the outgoing optical function of Lemma \ref{lemma:optical-Lext}. Then, we define the following:
\begin{itemize}
\item Let $\et_4:=-\g^{\a\b} \pr_\a \ut \pr_\b$ the null  outgoing  geodesic vectorfield   associated  to $\ut$. 

\item  Let $\st$  be the associated affine parameter, 
i.e. $\et_4(\st)=1$,  with $\st $ normalized    such that $\st= r$ as $r\to \infty$.  

\item We define  the region $\Lextt\subset \Lext$ to be the region\footnote{Recall that $0\le u\le 3 $ and $r\ge r_0$ in $\Lext$, so that we have indeed $\Lextt\subset \Lext$ in view of the asymptotic of $\ut$ provided by  Lemma \ref{lemma:optical-Lext} and the asymptotic for $\st$ provided by Proposition \ref{Prop:transitionparam(E-Et)}.}  
      \bea 
      \Lextt := \left\{ 0\le \ut\le 2, \quad \st\ge \frac{\de_*}{2} \ep_0^{-1}\right\}, 
      \eea
      with $\de_*>0$ the small constant  introduced in section \ref{sec:discussionofsmallnessconstantforthemaintheorem}. See Figure \ref{figL0tilde} below where $\Lextt$ is sketched in red inside the initial data layer $\LL_0$.

\item We denote by $\rt$ the area radius of the spheres $S(\ut, \st)$. 

\item The foliation induced by $(\ut, \st)$ is called the   outgoing geodesic foliation of  $\Lext$  normalized at infinity.  We denote by $\Et=(\et_4, \et_3, \et_1, \et_2)$ the corresponding  null   frame. 

\item We denote by  $(f, \fb, \la)$ the transition coefficients  of the frame   transformation which  takes  the   outgoing   PG frame   of $\Lext$, denoted $E$,  into $\Et$.
\end{itemize}
\end{definition}

\begin{figure}[h!]
\centering
\includegraphics[scale=0.9]{Kerr_2.pdf}
\caption{The initial data layer $\LL_0$ and the region $\Lextt$ (in red)}
\label{figL0tilde}
\end{figure}

\begin{proposition}\lab{Prop:transitionparam(E-Et)} 
The following holds true in $\Lextt$:
\begin{enumerate}
\item
The transition functions $(f, \fb, \la)$  are given by the formulas
\bea
\bsplit
\la &= 1 +\frac{3a^2(\sin\th)^2}{4r^2}+O(r^{-3}) + r^{-1} \dk^{\le 1} \Ga_b,\\
f &=  - a\Re(\Jk) +\frac{a^2\cos\th}{r}\Im(\Jk) +O(r^{-3}) +  \dk^{\le 1} \Ga_b,\\
\fb &= - a\Up\Re(\Jk)-\frac{a^2\cos\th}{r}\Im(\Jk)+O(r^{-3})+  \dk^{\le 2} \Ga_b.
\end{split}
\eea
Using complex notations, $F=f+i\dual f$ and $\Fb=\fb+i\dual\fb$ satisfy
\beaa
F&=&-a\left(1 + i\frac{a \cos\th}{r} \right)\Jk+O(r^{-3}) + \dk^{\le 1} \Ga_b =-\frac{a q}{r} \Jk +O(r^{-3}  ) + \dk^{\le 1} \Ga_b,\\
\Fb&=& -a\left(\Up - i\frac{a \cos\th}{r} \right)\Jk+O(r^{-3}) + \dk^{\le 1} \Ga_b= - \frac{a\Up\ov{q} }{r} \Jk +O(r^{-3}  ) +  \dk^{\le 2} \Ga_b.
\eeaa

\item The function $\st$   behaves as follows
\bea
\st= r+ \frac{ a^2 \sin^2 \th}{2r}+ O(r^{-2} )+ \dk^{\leq 1}\Ga_b.
\eea
\end{enumerate}
\end{proposition}

 \begin{proof}
 The proof proceeds in several steps.
 
 {\bf Step 1.} We start with the control of $f$. To this end, we rely on the  formula 
 \beaa
 \et_4 &=& \la\left(e_4+f^ae_a+\frac{1}{4}|f|^2e_3\right).
 \eeaa 
 On the other hand, in view of the definition of $\et_4$, we  have
\beaa
\widetilde{e}_4 &=& -\g^{\a\b}\pr_\a(\ut )\pr_\b= \frac{1}{2}e_3(\ut)e_4+\frac{1}{2}e_4(\ut)e_3-e_a(\ut)e_a.
\eeaa
Hence
\beaa
\la=\frac 1 2 e_3(\ut), \qquad   f =- \frac{2}{\et_3(\ut)}\nab( \ut).
\eeaa
We calculate, using  $e_3(r) =-\Up+O(r^{-2})+r \Ga_b$,  $e_3(u)= 2+\frac{2a^2\sin^2\th}{r^2}+ O(r^{-4} )+    \Ga_b$, $e_3(\cos \th)=\Ga_b$, $\ut=u  -\frac{a^2(\sin\th)^2}{2r}+ h$ and $h=O(r^{-3})+\Ga_b$, 
\beaa
e_3(\ut)&=&  e_3\left(u  -\frac{a^2(\sin\th)^2}{2r}+ h \right)= 2 +\frac{2a^2\sin^2\th}{r^2} +O(r^{-4}) +\Ga_b\\
&& +\frac{a^2(\sin\th)^2}{2r^2} e_3(r)  +\frac{a^2}{r} \cos \th e_3(\cos \th)+ e_3(h)
\\
&=& 2+\frac{ 3a^2\sin^2\th}{2r^2}+ O(r^{-3} )+  \dk^{\le 1}  \Ga_b.
\eeaa

Also, since  $ \nab u= a \Re(\Jk)+\Ga_b $ and $\nab(\cos\th)=-\Im(\Jk) +\Ga_b$, we have
\beaa
\nab(\ut)&=&  \nab\left(u  -\frac{a^2(\sin\th)^2}{2r}+ h \right)= a \Re(\Jk) +\Ga_b+ \frac{a^2\cos\th}{r}\nab(\cos\th)+\nab(h)\\
&=& a \Re(\Jk) -\frac{ a^2 \cos\th}{r} \Im (\Jk) + O( r^{-4})+\dk^{\le 1} \Ga_b. 
\eeaa
We deduce,
\beaa
\la&=& \frac 1 2  e_3(\ut)= 1+\frac{ 3a^2\sin^2\th}{4r^2}+ O(r^{-3} )+  \dk^{\le 1}\Ga_b,\\
f &=& - \frac{2}{\et_3(\ut)}\nab( \ut) =- a \Re(\Jk) +\frac{ a^2 \cos\th}{r} \Im (\Jk) + O(r^{-3})+\dk^{\le 1}\Ga_b,
\eeaa
which is the desired estimate for $f$, but not  for $\la$.

{\bf Step 2.} Next, we derive the desired estimate for $\la$. To this end,  we need to  improve the estimate for $\la$ of Step 1. This improvement  will be needed to get the correct asymptotic for $\st$. 
\begin{lemma}\lab{lemma:improvedestimateforlaingeodesicfoliationauxilliaryLext0auxlemma}
We have
\beaa
\la =1+\frac{ 3a^2\sin^2\th}{4r^2}+ O(r^{-3} )+ r^{-1}\dk^{\leq 1}\Ga_b.
\eeaa
\end{lemma}

\begin{proof}
  Recall  the transport equation for $e_4(\log \la) $ in Corollary
  \ref{cor:transportequationine4forchangeofframecoeff:simplecasefirst}, which, in view of that fact that $\om=0$ and $\etab=-\ze$, takes the following form
  \beaa
  \la^{-1}\nab_{\et_4}(\log\la) &=& 2  f\c\ze +E_2(f, \Ga),\\
  E_2(f, \Ga) &=& - \frac{1}{2}|f|^2\omb  -\frac{1}{4}\trchb|f|^2+O(f^3\Ga+f^2\chibh).
  \eeaa
 Using in particular the form of $f$ in Step 1, note that
   \beaa
 E_2(f, \Ga) &=&   -\frac{1}{4}\trchb|f|^2+O(r^{-4}) + r^{-2}\dk^{\leq 1}\Ga_b. 
 \eeaa
 We deduce
 \beaa
 \et_4 (\log \la) &=& 2 f \c \ze  -\frac{1}{4}\trchb|f|^2+O(r^{-4})+r^{-2}\dk^{\leq 1}\Ga_b\\ 
 &=&  \Re(F\c \ov{Z})-\frac  1  4   \trchb|f|^2+O(r^{-4})+r^{-2}\dk^{\leq 1}\Ga_b.
 \eeaa
 Next, recall that
 \beaa
 Z=\frac{a\ov{q}}{|q|^2}\Jk+\Ga_g , \qquad \trchb=-\frac{2\Up}{r} +O(r^{-3} )+\Ga_g.
 \eeaa
  Thus, together with the form of $f$ in Step 1, we infer
  \beaa
   F\c \ov{Z}&=& \left(- \frac{aq}{r}\Jk  + O(r^{-3} ) + \dk^{\leq 1}\Ga_b\right)\c\left(\frac{aq}{|q|^2}\ov{\Jk}+\Ga_g\right)\\
   &=& - \frac{aq}{r}\Jk \c  \left(\frac{a q }{|q|^2}\ov{\Jk}\right)+O( r^{-4} ) +r^{-2}\dk^{\leq 1}\Ga_b = -\frac{2a^2\sin^2\th}{r^3} +O( r^{-4} ) +r^{-2} \dk^{\leq 1}\Ga_b,
  \eeaa
  \beaa
  |f|^2&=& \left| - a \Re(\Jk) +\frac{ a^2 \cos\th}{r} \Im (\Jk) + O( r^{-3} )+\dk^{\leq 1}\Ga_b\right|^2 =\frac{a^2\sin^2\th}{r^2} +O(r^{-3} )+ r^{-1} \dk^{\leq 1}\Ga_b,
  \eeaa
  and
  \beaa
  \frac  1  4   \trchb|f|^2&=& \left(- \frac{\Up}{2r}  + O(r^{-3} )+  \Ga_g\right)\left(   \frac{a^2\sin^2\th}{r^2}+O(r^{-3} )+ r^{-1} \dk^{\leq 1}\Ga_b \right)\\
  &=&-  \frac{a^2\sin^2\th}{2r^3}+ O(r^{-4} )+  r^{-2}\dk^{\leq 1}\Ga_b.
  \eeaa
 We deduce
 \beaa
 \et_4 (\log \la ) &=&  \Re(F\c \ov{Z})-\frac  1  4   \trchb|f|^2+O(r^{-4})+ r^{-2} \dk^{\leq 1}\Ga_b\\
 \\&=&  -\frac{2a^2\sin^2\th}{r^3}+   \frac{a^2\sin^2\th}{2r^3}+   O( r^{-4} )+ r^{-2}\dk^{\leq 1}\Ga_b\\
 &=&-\frac  3 2 \frac{a^2\sin^2\th}{r^3}+O(r^{-4} )+ r^{-2} \dk^{\leq 1}\Ga_b.
 \eeaa
 Integrating backwards from $r\to +\infty$, and noticing that $\log(\la)=0$ at $r=+\infty$ in view of the control for $\la$ derived in Step 1,  we obtain  
\beaa
\log \la &=&\frac{3}{ 4 }  \frac{a^2\sin^2\th}{r^2}+ O(r^{-3} )+ r^{-1}\dk^{\leq 1}\Ga_b
\eeaa
i.e.
\beaa
\la &=&1 +\frac{3}{ 4 }  \frac{a^2\sin^2\th}{r^2}+ O(r^{-3} )+ r^{-1}\dk^{\leq 1}\Ga_b
\eeaa
as stated. This concludes the proof of Lemma \ref{lemma:improvedestimateforlaingeodesicfoliationauxilliaryLext0auxlemma}.
\end{proof}

{\bf Step 3.} We  look for $\st$ in the form  
\bea
\st= r + \st_0, \qquad \lim_{r\to +\infty}\st_0=0.
\eea
To have $\et_4(\st)=1$ we thus need to solve the following   transport equation for $\st_0$
\bea
\et_4( \st_0) &=&1 -   \et_4(r).
\eea
  Since $e_a(r)=0$  and   $\big|\Re(\Jk)\big|^2=\frac{(\sin\th)^2}{|q|^2}= \frac{(\sin\th)^2}{r^2}+O(r^{-4})$, we have
  \beaa
\et_4(r) &=& \la\left( e_4 + f^a  e_a +\frac 1 4 |f|^2 e_3\right) r = \la  +\frac 1 4\la |f|^2 e_3(r) \\
&=&\la  +\frac 1 4\la  |f|^2\Big( -\Up  +r\Ga_b  +O(r^{-2}) \Big)= \left( \la- \frac 1 4 |f|^2 \right)+O( r^{-3} )+ r^{-1}\dk^{\leq 1}\Ga_b.
\eeaa
Using the improved asymptotic for $\la$  and  the asymptotic of  $f$ derived above
\beaa
\la-\frac 1 4  |f|^2  &=&1 +\frac{3}{ 4 }  \frac{a^2\sin^2\th}{r^2}+ O(r^{-3} )+ r^{-1}\dk^{\leq 1}\Ga_b\\
&&- \frac 1  4 \left|- a \Re(\Jk) +\frac{ a^2 \cos\th}{r} \Im (\Jk) + O( r^{-3} )+\dk^{\leq 1}\Ga_b\right|^2\\
&=& 1 +\frac{3}{ 4 }  \frac{a^2\sin^2\th}{r^2}-\frac 1 4 \frac{a^2 \sin^2\th }{r^2} + O( r^{-3} )+ r^{-1}\dk^{\leq 1}\Ga_b\\
&=&1+\frac{ a^2 \sin^2\th}{r^2} + O( r^{-3} )+ r^{-1} \dk^{\leq 1}\Ga_b.
\eeaa
Hence
\beaa
\et_4(r) &=& \left( \la- \frac 1 4 |f|^2 \right)+O( r^{-3} )+ r^{-1}\dk^{\leq 1}\Ga_b= 1+\frac{ a^2 \sin^2\th}{r^2} + O( r^{-3} )+ r^{-1}\dk^{\leq 1}\Ga_b
\eeaa
and we deduce
\beaa
\et_4(\st_0) &=&1 -   \et_4(r)=-  \frac{ a^2 \sin^2\th}{r^2} + O( r^{-3} )+ r^{-1}\dk^{\leq 1}\Ga_b.
\eeaa
Integrating backwards  from $r\to +\infty$, and since $\st_0=0$ at $r=+\infty$, we infer 
\beaa
\st_0&=& \frac{ a^2 \sin^2\th}{r} + O(r^{-2} )+ \dk^{\leq 1}\Ga_b.
\eeaa
Therefore 
 \beaa
 \st &=& r+\st_0= r+  \frac{ a^2 \sin^2\th}{r} + O(r^{-2} )+ \dk^{\leq 1}\Ga_b
 \eeaa
as  stated.

{\bf Step 4.} In this last step, we control $\fb$.  To  this end,  we recall  the  formula
\beaa
\widetilde{e}_a &=& \left(\de_a^b +\frac{1}{2}\fb_af^b\right) e_b +\frac 1 2  \fb_a  e_4 +\left(\frac 1 2 f_a +\frac{1}{8}|f|^2\fb_a\right)   e_3.
\eeaa
Since  $e_a(r)=0$,  $\st= r+  \frac{ a^2 \sin^2 \th}{2r}+ O(r^{-2}) +\dk^{\leq 1}\Ga_b$ and $f= O(r^{-1})$,  we derive
\beaa
\nab(\st)&=& \nab\left(   \frac{ a^2 \sin^2 \th}{2r}\right) +O(r^{-3} ) + r^{-1}  \dk^{\le 2} \Ga_b=-\frac{a^2}{r} \cos \th \nab(\cos \th)  +O(r^{-3} ) + r^{-1} \dk^{\le 2} \Ga_b\\
&=& -\frac{a^2}{r} \cos \th  \big(-\Im(\Jk) +\Ga_b\big)  +O(r^{-3} )+ r^{-1} \dk^{\le 2} \Ga_b\\
&=& \frac{a^2}{r} \cos \th\Im(\Jk)  +O(r^{-3}) + r^{-1} \dk^{\le 2} \Ga_b,\\
e_3(\st) &=& \left( 1-\frac{a^2\sin^2\th}{2r^2}\right)\Big( -\Up+r\Ga_b \Big)+O(r^{-3} )+ \dk^{\le 2} \Ga_b \\
&=& -\Up +\frac{a^2\sin^2\th}{r^2}   +O(r^{-3} )+ r\dk^{\le 2} \Ga_b,\\
 e_4(\st)&=& 1- \frac{a^2\sin^2\th}{2r^2}  +O(r^{-3} )+ r^{-1}  \dk^{\le 2} \Ga_b.
\eeaa
Therefore, since $\widetilde{e}_a(\st)=0$,
\beaa
0 &=& \nab(\st)+\frac{1}{2}(f\c\nab(\st))\fb+\frac 1 2  \fb  e_4(\st) +\left(\frac 1 2 f +\frac{1}{8}|f|^2\fb\right)   e_3(\st)\\
&=&   \nab(\st) +  \frac 1 2  \fb  e_4(\st) +\frac 1 2 f   e_3(\st) + O(r^{-3} )+ r^{-1}  \dk^{\le 2} \Ga_b\\
&=& \frac{a^2}{r} \cos \th\Im(\Jk) + O(r^{-3} ) + r^{-1} \dk^{\le 2} \Ga_b +\frac 1 2 \fb \left(1- \frac{a^2\sin^2\th}{2r^2}+O(r^{-3} ) + \dk^{\le 2}\Ga_b\right) \\
&&+\frac 1 2 f\left(-\Up +\frac{a^2\sin^2\th}{r^2} + O(r^{-3} )+ r \dk^{\le 2}\Ga_b \right) + O(r^{-3} )+ r^{-1}  \dk^{\le 2} \Ga_b.
\eeaa
We deduce,
\beaa
0&=& \frac{a^2}{r} \cos \th\Im(\Jk)+ \frac 1 2 \fb -\frac{1}{2}\Up f +O(r^{-3} )+\dk^{\le 2} \Ga_b
\eeaa
i.e.
\beaa
\fb &=& \Up f - \frac{2a^2\cos\th}{r}\Im (\Jk)+O(r^{-3}) +  \dk^{\le 2} \Ga_b\\
 &=&\Up \left(   - a\Re(\Jk) +\frac{a^2\cos\th}{r}\Im(\Jk) +O(r^{-3})    +  \dk^{\le 2} \Ga_b \right) - \frac{2a^2\cos\th}{r}\Im (\Jk)+O(r^{-3}) + \dk^{\le 2} \Ga_b\\
&=& - a\Up \Re(\Jk)  - \frac{a^2\cos\th}{r}\Im (\Jk)+O(r^{-3}) +  \dk^{\le 2} \Ga_b
\eeaa
as stated.
This ends the proof of Proposition \ref{Prop:transitionparam(E-Et)}.
\end{proof}

%%%%%%%%%%%%%%%%%%%%%%%%%%%%%%%%%%%%%%%%

\subsection{Ricci and curvature coefficients in the geodesic frame of $\Lext$}
\lab{sec:auxiliaryoutgoingRiccicoeeffandcurvature:identties}

%%%%%%%%%%%%%%%%%%%%%%%%%%%%%%%%%%%%%%%%

\begin{lemma}
\lab{Lemma:geodesic-RicciCurv}
 The Ricci and curvature coefficients  relative to the outgoing geodesic foliation of $\Lext$  verify following:
 \begin{enumerate}
 
\item  The curvature coefficients  satisfy 
\bea
\bsplit
\at&=\a+O(r^{-5} ) + r^{-3} \dk^{\le 2} \Ga_b, \qquad \quad\,\, \bt= \b+  \frac{3am} {r^3} \Re(\Jk)+O(r^{-5} )+ r^{-3} \dk^{\le 2}\Ga_b,\\
\rhot&=-\frac{2m}{r^3} + O(r^{-5}) + r^{-3} \dk^{\le 2} \Ga_b,\quad  \rhodt=\frac{6 am\cos \th }{r^4}  + O( r^{-5})+ r^{-3}  \dk^{\le 2}\Ga_b,\\
\bbt&=\bb+ O(r^{-5})+  r^{-1} \dk^{\le 2} \Ga_b, \qquad \quad\,\,  \aat=\aa + r^{-1} \dk^{\le 2}\Ga_b+ O(r^{-5} ).
\end{split}
\eea

\item The Ricci coefficients satisfy\footnote{Recall that $\xit=0$, $\omt=0$ and $\etat=\zet=-\etabt$ as  the foliation  is outgoing geodesic.}
\bea
\lab{eq:RiccicgeodframeinKerr}
\bsplit
\trcht& = \frac 2 r  +O(r^{-3} )+r^{-1} \dk^{\le 3} \Ga_b, \qquad \qquad \,
 \chiht=   O(r^{-3})  + r^{-1}   \dk^{\le 3} \Ga_b, \\
\trchbt&= -\frac{2 \Up }{r} + O(r^{-3} ) + r^{-1}  \dk^{\le 3}  \Ga_b , \qquad \,\,\,\,\chibht=     O(r^{-3})+ r^{-1}   \dk^{\le 3} \Ga_b, \\
 \zet&= O(r^{-3} )+ r^{-1}  \dk^{\le 3}  \Ga_b,  \qquad\qquad\qquad\,  \xibt=O\big( r^{-3} \big)+ r^{-1}   \dk^{\le 3} \Ga_b,\\
   \ombt&=\frac{m}{r^2} +O\big(r^{-3}\big) + r^{-1}   \dk^{\le 3} \Ga_b.
\end{split}
\eea

\item  Let $\rt$ be the area radius of $S(\ut, \st)$, i.e.  $ 4\pi (\rt)^2=| S(\ut, \st)|$. Then, we have
\bea
\rt &=&   r+\frac{a^2(\sin\th)^2}{2r}+O(r^{-2})+r\dk^{\leq 2}\Ga_b.
\eea

\item Let
 \beaa
 \widecheck{\trcht}:= \trcht-\frac{2}{\rt}, \qquad   \widecheck{\trchbt}:= \trchbt+ \frac{2(1-\frac{2m}{\rt }) }{\rt}.
 \eeaa
 Then
 \bea
 \lab{eq:ell=1modesfortrcht-trchbt}
 \widecheck{\trcht} =O(r^{-4})+ r^{-1} \dk^{\le 3}\Ga_b,  \qquad  \widecheck{\trchbt} =O(r^{-4})+ r^{-1} \dk^{\le 3}\Ga_b. 
 \eea
 
 \item Let the basis of $\ell=1$ modes $\Jt^{(p)} $  of  $S(\ut, \st)$ be given by $\nab_{\et_4}  \Jt^{(p)} =0$ and $\Jt^{(p)}= J^{(p)}$ as $r\to \infty$. Then, we have
 \bea
\Jt^{(p)}-J^{(p)} &=&  O(r^{-1})+\dk^{\leq 1}\Ga_b, \qquad p=0,+,-.
\eea
 
\item We have 
\bea
\divt \bt= O(r^{-6} ) + r^{-3} \dk^{\le 3} \Ga_b, \,\,\,\curlt \bt=\frac{6a_0m_0}{ \rt^5}\Jt^{(0)} +O(r^{-6})  + r^{-3} \dk^{\le 3} \Ga_b.
\eea
 \end{enumerate}
\end{lemma}

\begin{proof}
The proof proceeds in several steps. 

{\bf Step 1.} Recall from Proposition \ref{Prop:transitionparam(E-Et)} that the transition coefficients $(f, \fb, \la)$  from the frame $E$  to the frame $\Et$ are given by the formulas
\beaa
\bsplit
\la &= 1 +\frac{3a^2(\sin\th)^2}{4r^2}+O(r^{-3}) + r^{-1} \dk^{\le 1} \Ga_b,\\
f &=  - a\Re(\Jk) +\frac{a^2\cos\th}{r}\Im(\Jk) +O(r^{-3}) +  \dk^{\le 1} \Ga_b,\\
\fb &= - a\Up\Re(\Jk)-\frac{a^2\cos\th}{r}\Im(\Jk)+O(r^{-3})+  \dk^{\le 2} \Ga_b.
\end{split}
\eeaa
To derive the formulas for the curvature  and Ricci coefficients, it suffices to consider the simplified formulas 
\beaa
\bsplit
\la &= 1 + O(r^{-2}) + r^{-1} \dk^{\le 1} \Ga_b,\\
f &=  - a\Re(\Jk) +O(r^{-2}) +  \dk^{\le 1} \Ga_b,\\
\fb &= - a\Up\Re(\Jk)+O(r^{-2})+  \dk^{\le 2} \Ga_b.
\end{split}
\eeaa
Together with the transformation  formulas of  Proposition \ref{Proposition:transformationRicci}, this easily yields  the formulas for the curvature  and Ricci coefficients. Note that, modulo the terms in $\Ga_b$,  the  formulas are exactly the same as those in Lemma \ref{Lemma:Transformation-principal-to-integrable-Kerr} regarding the Ricci and curvature coefficients for the  integrable frame.  
  
{\bf Step 2.} Next, we derive more precise transformation formulas for  $\trcht, \trchbt$.   In view of  the  transformation formulas  of  Proposition \ref{Proposition:transformationRicci},  and the fact that $ \atrch, \,\omb  = O(r^{-2} ) +\Ga_b$ and  $f, \fb=O(r^{-1})$,  we have 
 \beaa
\bsplit
\la^{-1}\trcht &= \trch  +  \divt f + f\c\eta + f\c\ze+\err(\trch,\trcht ),\\
\err(\trch,\trcht ) &=    -\frac{1}{4}|f|^2\trchb + O(r^{-4}) + r^{-2} \dk^{\leq 2}\Ga_b,
\\
\la\trchbt  &= \trchb +\divt \fb +\fb\c\etab  -  \fb\c\ze +\err(\trchb, \trchbt ),\\
\err(\trchb, \trchbt ) &= \frac{1}{2}(f\c\fb)\trchb   -\frac 1 4 |\fb|^2 \trch + O(r^{-4}) + r^{-1} \dk^{\leq 2}\Ga_b,
\end{split}
\eeaa 
where $(f, \fb, \la) $ are the transition coefficients  from the frame $E$  to the frame $\Et$. In view of 
\beaa
\bsplit
f &=  - a\Re(\Jk)  +O(r^{-2}) +  \dk^{\le 1} \Ga_b,\\
\fb &= - a\Re(\Jk)+O(r^{-2})+  \dk^{\le 2} \Ga_b,\\
\trch &= \frac{2}{r} -\frac{2a^2(\cos\th)^2}{r^3}+O(r^{-5})+\Ga_g,\\
\trchb &= -\frac{2\Up}{r} -\frac{2a^2}{r^3} +\frac{4a^2(\cos\th)^2}{r^3}+O(r^{-4})+\Ga_g,\\
\ze &= \Re\left(\frac{a\ov{q}}{|q|^2}\Jk\right)+\Ga_g=\frac{a}{r}\Re(\Jk)+O(r^{-3})+\Ga_g,\\
\eta &= \Re\left(\frac{aq}{|q|^2}\Jk\right)+\Ga_b=\frac{a}{r}\Re(\Jk)+O(r^{-3})+\Ga_b,\\
\end{split}
\eeaa
and $\etab=-\ze$, we deduce
 \beaa
\bsplit
\la^{-1}\trcht &= \frac{2}{r} -\frac{2a^2(\cos\th)^2}{r^3} - \frac{3a^2(\sin\th)^2}{2r^3} +  \divt f      + O(r^{-4}) + r^{-1} \dk^{\leq 2}\Ga_b,
\\
\la\trchbt  &= -\frac{2\Up}{r} -\frac{2a^2}{r^3} +\frac{4a^2(\cos\th)^2}{r^3} +\frac{a^2(\sin\th)^2}{2r^3}  +\divt \fb   + O(r^{-4}) + r^{-1} \dk^{\leq 3}\Ga_b.
\end{split}
\eeaa 

Next, we compute $ \divt f$ and $\divt\fb$. Arguing similarly to section \ref{sec:comparisionhorizontalderivativesandtangentialtoSur}, see in particular \ref{cor:formuallikningnabprimeandnabintransfoformula:bis}, and using $f, \fb=O(r^{-1})$ and $\la=1+O(r^{-2})$, we obtain 
\beaa
\divt f&=& \div  f + \frac 1 2  f\c (\nab_3 f) + \frac 1 2  \fb \c( \nab_4  f) +O(r^{-4} )   +  r^{-1}  \dk^{\le  2}  \Ga_b,\\
\divt \fb&=& \div  \fb + \frac 1 2  f\c (\nab_3 \fb) + \frac 1 2  \fb \c( \nab_4  \fb) +O(r^{-4} )   +  r^{-1}  \dk^{\le  3}  \Ga_b.  
\eeaa
Using again
\beaa
\bsplit
f &=  - a\Re(\Jk) +\frac{a^2\cos\th}{r}\Im(\Jk) +O(r^{-3}) +  \dk^{\le 1} \Ga_b,\\
\fb &= - a\Up\Re(\Jk)-\frac{a^2\cos\th}{r}\Im(\Jk)+O(r^{-3})+  \dk^{\le 2} \Ga_b,
\end{split}
\eeaa
we infer
\beaa
\divt f&=& \div\left(- a\Re(\Jk) +\frac{a^2\cos\th}{r}\Im(\Jk)\right) + \frac{a^2}{2}\Re(\Jk)\c\nab_3\Re(\Jk) + \frac{a^2}{2}\Re(\Jk)\c\nab_4\Re(\Jk)\\
&& +O(r^{-4} )   +  r^{-1}  \dk^{\le  2}  \Ga_b,\\
\divt \fb&=& \div\left( - a\Up\Re(\Jk)-\frac{a^2\cos\th}{r}\Im(\Jk)\right) + \frac{a^2}{2}\Re(\Jk)\c\nab_3\Re(\Jk) + \frac{a^2}{2}\Re(\Jk)\c\nab_4\Re(\Jk)\\
&&+O(r^{-4} )   +  r^{-1}  \dk^{\le  3}  \Ga_b.  
\eeaa
Note that 
\beaa
2\Re(\Jk)\c\nab_3\Re(\Jk) + 2\Re(\Jk)\c\nab_4\Re(\Jk) &=& \nab_{e_3+e_4}(|\Re(\Jk)|^2)\\
&=& O(r^{-3})\Big(e_3(r)+e_4(r)\Big)+O(r^{-3})e_3(\cos\th)\\
&=& O(r^{-4})+r^{-2}\Ga_b
\eeaa
where we used $e_4(\th)=0$, $e_3(\cos\th)\in\Ga_b$, $e_4(r)=1$ and $e_3(r)=-1+O(r^{-1})+r\Ga_b$. This yields
\beaa
\divt f&=& \div\left(- a\Re(\Jk) +\frac{a^2\cos\th}{r}\Im(\Jk)\right)  +O(r^{-4} )   +  r^{-1}  \dk^{\le  2}  \Ga_b,\\
\divt \fb&=& \div\left( - a\Up\Re(\Jk)-\frac{a^2\cos\th}{r}\Im(\Jk)\right) +O(r^{-4} )   +  r^{-1}  \dk^{\le  3}  \Ga_b.  
\eeaa
Next, we use 
\beaa
\ov{\DD}\c\Jk=\frac{4i(r^2+a^2)\cos\th}{|q|^4}+r^{-1}\Ga_b
\eeaa
which yields
\beaa
\div(\Re(\Jk)) = r^{-1}\Ga_b, \qquad \div(\Im(\Jk))=\frac{2\cos\th}{r^2}+O(r^{-4})+r^{-1}\Ga_b,
\eeaa
and hence
\beaa
\divt f&=& \frac{2a^2(\cos\th)^2}{r^3} +\frac{a^2}{r}\nab(\cos\th)\c\Im(\Jk)  +O(r^{-4} )   +  r^{-1}  \dk^{\le  2}  \Ga_b,\\
\divt \fb&=& -\frac{2a^2(\cos\th)^2}{r^3} -\frac{a^2}{r}\nab(\cos\th)\c\Im(\Jk) +O(r^{-4} )   +  r^{-1}  \dk^{\le  3}  \Ga_b.  
\eeaa
Since $\nab(\cos\th) = -\Im(\Jk)+\Ga_b$ and $|\Im(\Jk)|^2=\frac{a^2(\sin\th)^2}{|q|^2}$, we obtain
\beaa
\divt f&=& \frac{2a^2(\cos\th)^2}{r^3} -\frac{a^2(\sin\th)^2}{|q|^2} +O(r^{-4} )   +  r^{-1}  \dk^{\le  2}  \Ga_b,\\
\divt \fb&=& -\frac{2a^2(\cos\th)^2}{r^3} +\frac{a^2(\sin\th)^2}{|q|^2}  +O(r^{-4} )   +  r^{-1}  \dk^{\le  3}  \Ga_b.  
\eeaa
Plugging in  
 \beaa
\bsplit
\la^{-1}\trcht &= \frac{2}{r} -\frac{2a^2(\cos\th)^2}{r^3} - \frac{3a^2(\sin\th)^2}{2r^3} +  \divt f      + O(r^{-4}) + r^{-1} \dk^{\leq 2}\Ga_b,
\\
\la\trchbt  &= -\frac{2\Up}{r} -\frac{2a^2}{r^3} +\frac{4a^2(\cos\th)^2}{r^3} +\frac{a^2(\sin\th)^2}{2r^3}  +\divt \fb   + O(r^{-4}) + r^{-1} \dk^{\leq 2}\Ga_b,
\end{split}
\eeaa 
we infer
\beaa
\bsplit
\la^{-1}\trcht &= \frac{2}{r}  - \frac{5a^2(\sin\th)^2}{2r^3}     + O(r^{-4}) + r^{-1} \dk^{\leq 2}\Ga_b,
\\
\la\trchbt  &= -\frac{2\Up}{r} -\frac{2a^2}{r^3} +\frac{2a^2(\cos\th)^2}{r^3} +\frac{3a^2(\sin\th)^2}{2r^3}    + O(r^{-4}) + r^{-1} \dk^{\leq 3}\Ga_b.
\end{split}
\eeaa 
As 
\beaa
\la &=& 1 +\frac{3a^2(\sin\th)^2}{4r^2}+O(r^{-3}) + r^{-1} \dk^{\le 1} \Ga_b,
\eeaa
we deduce
\beaa
\bsplit
\trcht &= \frac{2}{r} -\frac{a^2(\sin\th)^2}{r^3}   + O(r^{-4}) + r^{-1} \dk^{\leq 2}\Ga_b,
\\
\trchbt  &= -\frac{2\Up}{r}  +\frac{a^2(\sin\th)^2}{r^3}     + O(r^{-4}) + r^{-1} \dk^{\leq 3}\Ga_b,
\end{split}
\eeaa 
which is the precise form that will be used in Step 4 and Step 5.  

{\bf Step 3.} Next, we derive a first, non sharp, asymptotic for $r-\rt$. Given coordinates\footnote{In practice, as in sections \ref{sec:outgoingPGcoordinatesinKerr:chap2} and  \ref{sec:outgoingPGcoordinatesinKerr:additionalcoordsystems:chap2}, we cover the spheres $S(u,r)$ with the coordinates systems $(x^1, x^2)=(\th, \vphi)$ and $(x^1, x^2)=(J^{(+)}, J^{(-)})$.} $(x^1, x^2)$ on the spheres $S(u,r)$ with $e_4(x^A)=0$, we consider the coordinates $(\xt^1, \xt^2)$ on $S(\ut, \st)$ given by 
\beaa
\xt^1=x^1, \qquad \xt^2=x^2.
\eeaa
We introduce the vectorfields
\beaa
\Xt_A &=& \pr_{x^A}+\frac{1}{2}\g(\et_4, \pr_{x^A})\et_3+\frac{1}{2}\g(\et_3, \pr_{x^A})\et_4.
\eeaa
which are tangent to the spheres $S(\ut, \st)$.  We have
\beaa
\Xt_A(\xt^B) &=& \de_A^B +\frac{1}{2}\g(\et_4, \pr_{x^A})\et_3(\xt^B)+\frac{1}{2}\g(\et_3, \pr_{x^A})\et_4(\xt^B).
\eeaa
Since $\xt^A=x^A$, $e_4(x^A)=0$, $e_3(x^A)=O(r^{-2})+\Ga_b$, $e_B(x^A)=O(r^{-1})$, as well as\footnote{Note that $\pr_{x^A}=Y^B_Ae_B+z^3_Ae_3+z^4_Ae_4$, with $Y^B_A=O(r)$ and $z^3_Ae_3(u)=-Y^B_Ae_B(u)$, $z^4_A=-z^3_Ae_3(r)$, so that we have $z^3_A, \, z^4_A=O(1)$, and hence $\g(e_3, \pr_{x^A})= -2z_A^4=O(1)$ and $\g(e_4, \pr_{x^A})= -2z_A^3=O(1)$.} $\g(e_3, \pr_{x^a})=O(1)$ and $\g(e_4, \pr_{x^a})=O(1)$, and since 
\beaa
\et_4 &=& \la\left(e_4+f^ce_c+\frac{1}{4}|f|^2e_3\right)=\la\Big(e_4+f^ce_c+O(r^{-2})e_3\Big),\\
\et_3 &=& \la^{-1}\left(\left(1+\frac{1}{2}f\c\fb  +\frac{1}{16} |f|^2  |\fb|^2\right) e_3 + \left(\fb^b+\frac 1 4 |\fb|^2f^b\right) e_b  + \frac 1 4 |\fb|^2 e_4\right)\\
&=& \la^{-1}\Big(\left(1+O(r^{-2})\right) e_3 + \left(\fb^b+O(r^{-3})\right) e_b  + O(r^{-2}) e_4\Big),
\eeaa 
we infer
\beaa
\Xt_A(\xt^B) &=& \de_A^B +\frac{1}{2}\g(\et_4, \pr_{x^A})\et_3(\xt^B)+\frac{1}{2}\g(\et_3, \pr_{x^A})\et_4(\xt^B)\\
&=& \de_A^B +\frac{1}{2}\g\Big(e_4+f^ce_c+O(r^{-2})e_3, \pr_{x^A}\Big)\Big(\left(1+O(r^{-2})\right) e_3  + \left(\fb^b+O(r^{-3})\right) e_b\Big)(x^B)\\
&&+\frac{1}{2}\g\Big(\left(1+O(r^{-2})\right) e_3  + \left(\fb^b+O(r^{-3})\right) e_b+O(r^{-2})e_4, \pr_{x^A}\Big)\Big(
f^ce_c+O(r^{-2})e_3\Big)(x^B)\\
&=& \de_A^B+O(r^{-2})+\Ga_b.
\eeaa
Also, 
\beaa
&&\g(\Xt_A, \Xt_B)\\ 
&=& \g\left(\pr_{x^A}+\frac{1}{2}\g(\et_4, \pr_{x^A})\et_3+\frac{1}{2}\g(\et_3, \pr_{x^A})\et_4, \pr_{x^B}+\frac{1}{2}\g(\et_4, \pr_{x^B})\et_3+\frac{1}{2}\g(\et_3, \pr_{x^B})\et_4\right)\\
&=& \g(\pr_{x^A}, \pr_{x^B})+\frac{1}{2}\g(\et_4, \pr_{x^A})\g(\et_3, \pr_{x^B})+\frac{1}{2}\g(\et_4, \pr_{x^B})\g(\et_3, \pr_{x^A})\\
&=& g_{AB}+\frac{1}{2}\g\Big(e_4+f^ce_c+O(r^{-2})e_3, \pr_{x^A}\Big)\g\Big(\left(1+O(r^{-2})\right) e_3  + \left(\fb^b+O(r^{-3})\right) e_b+O(r^{-2})e_4, \pr_{x^B}\Big)\\
&& +\frac{1}{2}\g\Big(e_4+f^ce_c+O(r^{-2})e_3, \pr_{x^B}\Big)\g\Big(\left(1+O(r^{-2})\right) e_3  + \left(\fb^b+O(r^{-3})\right) e_b+O(r^{-2})e_4, \pr_{x^A}\Big)\\
&=& g_{AB}+O(1),
\eeaa
where $g_{AB}$ denotes the induced metric on $S(u,r)$. Denoting the  induced metric on $S(\ut,\st)$ by $\widetilde{g}_{AB}$, we easily infer from the above computation of $\Xt_A(\xt^B)$ and $\g(\Xt_A, \Xt_B)$, and from the fact that the vectorfields $\Xt_1$ and $\Xt_2$ are tangent to the spheres $S(\ut, \st)$, that 
\beaa
\widetilde{g}_{AB} &=& g_{AB}+O(1)+r^2\Ga_b.
\eeaa
In particular, we infer
\beaa
|S(\ut, \st)| &=& \int \sqrt{|\widetilde{g}_{AB}|}d\xt^Ad\xt^B\\
&=& \int \sqrt{|g_{AB}|}dx^Adx^B+O(1)+r^2\Ga_b\\
&=& |S(u,r)|+O(1)+r^2\Ga_b.
\eeaa
Together with the  bound $|S(u, r)| =  4\pi r^2+O(1)$ which is a non sharp consequence of Lemma \ref{lemma:asymptoticarearadiusSur},  we deduce
\beaa
|S(\ut, \st)| &=& 4\pi r^2+O(1)+r^2\Ga_b,
\eeaa
 and hence, since $|S(\ut, \st)| = 4\pi (\rt)^2$, we infer the following non sharp bound for $\rt-r$
\beaa
\rt &=& r+O(r^{-1})+r\Ga_b.
\eeaa

{\bf Step 4.} We now improve the bound for $\rt-r$ of Step 3. Recall that, as $(\ut, \st)$ is an outgoing geodesic foliation, we have
\beaa
\et_4\left(\trcht -\frac{2}{\rt}\right) &=& -\frac{1}{2}(\trcht)^2-|\chiht|^2+\frac{2}{(\rt)^2}\frac{\rt}{2}\ov{\trcht}^{S(\ut, \st)}.
\eeaa
In view of the above decomposition of $\chiht$, we infer
\beaa
\et_4\left(\trcht -\frac{2}{\rt}\right) &=& -\frac{1}{2}(\trcht)^2+\frac{1}{\rt}\ov{\trcht}^{S(\ut, \st)}+O(r^{-6})+r^{-4}\dk^{\leq 2}\Ga_b.
\eeaa
We infer
\beaa
\ov{\et_4\left(\trcht -\frac{2}{\rt}\right)}^{S(\ut, \st)} &=& -\frac{1}{2}\ov{(\trcht)^2}^{S(\ut, \st)} +\frac{1}{\rt}\ov{\trcht}^{S(\ut, \st)}+O(r^{-6})+r^{-4}\dk^{\leq 2}\Ga_b.
\eeaa
Since, for any scalar function $h$, we have
\beaa
\et_4\left(\ov{h}^{S(\ut, \st)}\right) &=& \frac{1}{|S(\ut, \st)|}\int_{S(\ut, \st)}\Big(\et_4(h)+\trcht h\Big) -\frac{2\et_4(\rt)}{\rt}\ov{h}^{S(\ut, \st)} \\
&=& \ov{e_4(h)}^{S(\ut, \st)}+\ov{\trcht h}^{S(\ut, \st)} -\ov{\trcht}^{S(\ut, \st)}\ov{h}^{S(\ut, \st)},
\eeaa
we infer
\beaa
\et_4\left(\ov{\trcht}^{S(\ut, \st)} -\frac{2}{\rt}\right) &=& -\frac{1}{2}\ov{(\trcht)^2}^{S(\ut, \st)} +\frac{1}{\rt}\ov{\trcht}^{S(\ut, \st)}+O(r^{-6})+r^{-4}\dk^{\leq 2}\Ga_b\\
&& +\ov{\trcht\left(\trcht -\frac{2}{\rt}\right)}^{S(\ut, \st)} -\ov{\trcht}^{S(\ut, \st)}\ov{\left(\trcht -\frac{2}{\rt}\right)}^{S(\ut, \st)}\\
&=& -\frac{1}{\rt}\left(\ov{\trcht}^{S(\ut, \st)} -\frac{2}{\rt}\right)+\frac{1}{2}\ov{\left(\trcht -\frac{2}{\rt}\right)^2}^{S(\ut, \st)}\\
&&-\left(\ov{\trcht}^{S(\ut, \st)} -\frac{2}{\rt}\right)^2+O(r^{-6})+r^{-4}\dk^{\leq 2}\Ga_b
\eeaa
and hence
\beaa
\et_4\left(\ov{\rt\trcht}^{S(\ut, \st)}-2\right) &=& \frac{1}{2\rt}\ov{\Big(\rt\trcht -2\Big)^2}^{S(\ut, \st)} -\frac{1}{2\rt}\left(\ov{\rt\trcht-2}^{S(\ut, \st)}\right)^2+O(r^{-5})+r^{-3}\dk^{\leq 2}\Ga_b.
\eeaa
Now, we have, in view of the decomposition of Step 3 for $\rt-r$, and  in view of the decomposition of Step 1 for $\trcht$,
\beaa
\rt\trcht-2 &=& O(r^{-1})(\rt-r)+r\trcht-2=O(r^{-2})+\dk^{\leq 2}\Ga_b
\eeaa
and hence
\beaa
\et_4\left(\ov{\rt\trcht}^{S(\ut, \st)}-2\right) &=& O(r^{-5})+r^{-2}\dk^{\leq 2}\Ga_b, \qquad\lim_{r\to +\infty}\left(\ov{\rt\trcht}^{S(\ut, \st)}-2\right)=0.
\eeaa
Integrating backwards from infinity, we infer
\beaa
\ov{\rt\trcht}^{S(\ut, \st)}-2 &=& O(r^{-4})+r^{-1}\dk^{\leq 2}\Ga_b.
\eeaa
Plugging the precise asymptotic for $\trcht$ of Step 2, we infer
\beaa
\ov{\rt\left(\frac{2}{r} -\frac{a^2(\sin\th)^2}{r^3}\right)}^{S(\ut, \st)}-2 &=& O(r^{-3})+\dk^{\leq 2}\Ga_b
\eeaa
which we rewrite
\beaa
\ov{\left(\frac{\rt}{r+\frac{a^2(\sin\th)^2}{2r}}\right)}^{S(\ut, \st)}-1 &=& O(r^{-3})+\dk^{\leq 2}\Ga_b.
\eeaa
On the other hand, we have
\beaa
&&\widetilde{\nab}\left(r+\frac{a^2(\sin\th)^2}{2r}\right)\\ 
&=& \left(\left(\de_a^b +\frac{1}{2}\fb_af^b\right) e_b +\frac 1 2  \fb_a  e_4 +\left(\frac 1 2 f_a +\frac{1}{8}|f|^2\fb_a\right)   e_3\right)\left(r+\frac{a^2(\sin\th)^2}{2r}\right)\\
&=& \frac 1 2  \fb_a   +\left(\frac 1 2 f_a +\frac{1}{8}|f|^2\fb_a\right)e_3(r)+e_a\left(\frac{a^2(\sin\th)^2}{2r}\right)+O(r^{-3}).
\eeaa
In view of the form of $f$, $\fb$, and the control of $e_3(r)$, we infer
\beaa
\widetilde{\nab}\left(r+\frac{a^2(\sin\th)^2}{2r}\right) &=& \frac{1}{2}\left( - a\Up\Re(\Jk)-\frac{a^2\cos\th}{r}\Im(\Jk)\right)+\frac{1}{2}\left(  - a\Re(\Jk) +\frac{a^2\cos\th}{r}\Im(\Jk) \right)(-\Up)\\
&&-\frac{a^2\cos\th}{r}\nab(\cos\th)+O(r^{-3})+  \dk^{\le 2}\Ga_b\\
&=& -\frac{a^2\cos\th}{r}\Big(\nab(\cos\th)+\Im(\Jk)\Big) +O(r^{-3})+  \dk^{\le 2}\Ga_b.
\eeaa
Since $\nab(\cos\th)=-\Im(\Jk)+\Ga_b$, we obtain
\beaa
\widetilde{\nab}\left(r+\frac{a^2(\sin\th)^2}{2r}\right) &=& O(r^{-3})+  \dk^{\le 2}\Ga_b.
\eeaa
Together with 
\beaa
\ov{\left(\frac{\rt}{r+\frac{a^2(\sin\th)^2}{2r}}\right)}^{S(\ut, \st)}-1 &=& O(r^{-3})+\dk^{\leq 2}\Ga_b,
\eeaa
we infer 
\beaa
\rt &=& r+\frac{a^2(\sin\th)^2}{2r}+O(r^{-2})+r\dk^{\leq 2}\Ga_b
\eeaa
which is the stated control of $\rt$. 

{\bf Step 5.} Next, recall the precise asymptotic for $\trcht$ and $\trchbt$ derived in Step 2
\beaa
\bsplit
\trcht &= \frac{2}{r} -\frac{a^2(\sin\th)^2}{r^3}   + O(r^{-4}) + r^{-1} \dk^{\leq 2}\Ga_b,
\\
\trchbt  &= -\frac{2\Up}{r}  +\frac{a^2(\sin\th)^2}{r^3}     + O(r^{-4}) + r^{-1} \dk^{\leq 3}\Ga_b.
\end{split}
\eeaa 
In view of the control for $\rt-r$ of Step 4, we deduce
\beaa
\bsplit
\trcht &= \frac{2}{\rt}  + O(r^{-4}) + r^{-1} \dk^{\leq 2}\Ga_b,
\\
\trchbt  &= -\frac{2\widetilde{\Up}}{\rt}    + O(r^{-4}) + r^{-1} \dk^{\leq 3}\Ga_b.
\end{split}
\eeaa 
In view of  the definition of $\widecheck{\trcht}$ and $\widecheck{\trchbt}$, we infer
\beaa
\bsplit
\widecheck{\trcht} &=  O(r^{-4}) + r^{-1} \dk^{\leq 2}\Ga_b,
\\
\widecheck{\trchbt} &=    O(r^{-4}) + r^{-1} \dk^{\leq 3}\Ga_b,
\end{split}
\eeaa 
as stated. 

{\bf Step 6.}  Next, we control $\Jt^{(p)}-J^{(p)}$ for $p=0,+,-$. Since $\et_4(\Jt^{(p)}) =0$, we have
\beaa
{\et_4}(\Jt^{(p)}-J^{(p)}) &=& -{\et_4}(J^{(0)}) =-\la\left(e_4+f\c\nab+\frac{1}{4}|f|^2e_3\right)J^{(p)}\\
&=& -\la\left(f\c\nab(J^{(p)})+\frac{1}{4}|f|^2e_3(J^{(p)})\right).
\eeaa
Using the control of $\la$ and $f$, we infer
\beaa
{\et_4}(\Jt^{(p)}-J^{(p)}) &=&  O(r^{-2})+r^{-1}\dk^{\leq 1}\Ga_b.
\eeaa
Integrating backwards from $r=+\infty$ where $\Jt^{(p)}= J^{(p)}$, we infer
\beaa
\Jt^{(p)}-J^{(p)} &=&  O(r^{-1})+\dk^{\leq 1}\Ga_b, \qquad p=0,+,-,
\eeaa
as desired.

{\bf Step 7.} Finally, we derive precise asymptotic for  $\divt \bt$ and $\curlt \bt$. Recall that we have obtained in Step 1
\beaa
\bt= \b+  \frac{3am} {r^3} \Re(\Jk)+O(r^{-5} )+ r^{-3} \dk^{\le 2}\Ga_b.
\eeaa
Since $\b\in r^{-2}\Ga_b$, we infer
\beaa
\bt= \frac{3am} {r^3} \Re(\Jk)+O(r^{-5} )+ r^{-2} \dk^{\le 2}\Ga_b
\eeaa
and hence
\beaa
\divt\bt &=& \divt\left(\frac{3am} {r^3} \Re(\Jk)\right)+O(r^{-6} )+ r^{-3} \dk^{\le 3}\Ga_b,\\
\curlt\bt &=& \curlt\left(\frac{3am} {r^3} \Re(\Jk)\right)+O(r^{-6} )+ r^{-3} \dk^{\le 3}\Ga_b.
\eeaa
One easily infers 
\beaa
\divt\bt &=& \frac{3am} {r^3} \div\left(\Re(\Jk)\right)+O(r^{-6} )+ r^{-3} \dk^{\le 3}\Ga_b,\\
\curlt\bt &=& \frac{3am} {r^3} \curl\left(\Re(\Jk)\right)+O(r^{-6} )+ r^{-3} \dk^{\le 3}\Ga_b.
\eeaa
Next, we use 
\beaa
\ov{\DD}\c\Jk=\frac{4i(r^2+a^2)\cos\th}{|q|^4}+r^{-1}\Ga_b
\eeaa
which yields
\beaa
\div(\Re(\Jk)) = r^{-1}\Ga_b, \qquad \curl(\Re(\Jk))=\frac{2\cos\th}{r^2}+O(r^{-4})+r^{-1}\Ga_b.
\eeaa
We infer, together with the control for $\rt-r$, 
\beaa
\divt\bt &=& O(r^{-6} )+ r^{-3} \dk^{\le 3}\Ga_b,\\
\curlt\bt &=& \frac{6am\cos\th}{\rt^5}+O(r^{-6} )+ r^{-3} \dk^{\le 3}\Ga_b.
\eeaa
This is the desired control for $\div\bt$. For $\curlt\bt$, we use in addition the following control derived in Step 6 
\beaa
\Jt^{(0)}-J^{(0)} &=&  O(r^{-1})+\dk^{\leq 1}\Ga_b.
\eeaa
Plugging in the above, and since $J^{(0)}=\cos\th$, we infer
\beaa
\curlt\bt &=& \frac{6am \Jt^{(0)}}{\rt^5}+O(r^{-6} )+ r^{-3} \dk^{\le 3}\Ga_b
\eeaa
as desired. This ends the proof of Lemma \ref{Lemma:geodesic-RicciCurv}.
\end{proof}

%%%%%%%%%%%%%%%%%%%%%%%%%%%%%%%%%%%%%

\subsection{Control of the geodesic foliation of $\Lext$}
\lab{sec:auxiliaryoutgoing:finalcontrol}

%%%%%%%%%%%%%%%%%%%%%%%%%%%%%%%%%%%%%

We recall in what follows that  the frame $\Et$, the scalar functions $(\ut, \st)$ and the scalars $\Jt^{(p)}$, $p=0,+,-$,  correspond to the outgoing geodesic foliation  of 
$\Lextt$, while $E$, $r$, $u$, $J^{(p)}$ and $\Jk$ correspond to the outgoing PG structure of $\Lext$. 
We also recall that  $(f, \fb, \la)$  denote the transition  coefficients  from  the frame $E$ to $\Et$. In addition, we  define
\bea
f_0:= r\Re(\Jk).
\eea

We  define  the following norms on $\Lextt$.
\beaa
\Ikt_{k} &:=& \Ikt_k ' +\Ikt_{k+1}''
\eeaa
where, for $ k=0$, 
\beaa
 \Ikt'_0 &:=&  \sup_{\Lextt}  \left[ r^{\frac{7}{2}  +\dt}\left( |\at| + \left|\bt - \frac{3a_0m_0}{\rt^4}f_0\right|\right)+r^3\ \left|\rhot+\frac{2m_0}{\rt^3}\right| +   r^2 |\bbt|+r|\aat|    \right]\\
  &+&\sup_{\Lextt}r^{\frac{9}{2} +\dt}  \left(  \big|\divt\b \big| +\left|\curlt \bt-\frac{6a_0 m_0}{\rt^5} \Jt^{(0)}\right| \right)\\
&+& \sup_{\Lextt}  r^2\left(|\chiht|+\left|\trcht-\frac{2}{\rt}\right|+|\zet|+ \left|\trchbt+\frac{2\left(1-\frac{2m_0}{\rt}\right)}{\rt}\right|\right)\\
&+& \sup_{\Lextt} r \left(|\chibht|+\left|\ombt-\frac{m_0}{\rt^2}\right|+|\xibt|+\left|\nabt \Jt^{(0)} +\frac{1}{\rt}\dual f_0 \right|+\left|\nabt f_0 -\frac{\Jt^{(0)} }{\rt}\in\right|\right)\\
&+&\sup_{\Lextt}\left(\left|\et_3(\rt)+1-\frac{2m_0}{\rt}\right| +|\et_3(\ut )-2| +|\rt-r|\right),
\eeaa
\beaa
 \Ikt_0''  &:=&  \sup_{\Lextt} r\left(  \left|f+\frac{a_0}{r}f_0\right|+\left|\fb+\frac{a_0\left(1-\frac{2m_0}{r}\right)}{r}f_0\right| +|\log(\la)|   \right)\\
&&+\sup_{\Lextt}\max_{p=0,+,-}\left| \Jt^{(p)}-J^{(p)}\right|+\sup_{\Lextt}r\left|\Jk-\frac{1}{|q|}(f_0+\dual f_0)\right|.
\eeaa
The higher derivative norms  $\Ikt', \Ikt''$ are then defined by  replacing each component with $\dk^{\le k}$ of it.

The following proposition provides the control of the norm $\Ikt_k$ under the assumptions \eqref{def:initialdatalayerassumptions:chap8}, i.e.  
\beaa
\, ^{(ext)} \Ik_{k_{large}+10} \leq \ep_0, \qquad   ^{(ext)} \Ik_{3} \leq  \ep^2_0.
\eeaa

\begin{proposition}
\lab{proposition:geodesicfoliationLextt}
The following holds true:
\begin{enumerate}
\item Under the assumption $\, ^{(ext)} \Ik_{k_{large}+10} \leq \ep_0$, we have
\bea
\Ikt_{k_{large + 7}} \les \ep_0.
\eea

\item If in addition the assumption  $\, ^{(ext)} \Ik_{3} \leq \ep^2_0$ also holds true,  then
\bea\lab{eq:strongerboundfordivbcurlbtrchtrchbatrlikeep0minus1inIDL}
\nn \sup_{\Lextt\cap\{\rt\sim \ep_0^{-1}\}  } r^5\left( \left|\divt\bt\right|+\left|\curlt\bt -\frac{6a_0m_0}{\rt^5}\Jt^{(0)} \right|\right)  \\
 +\sup_{\Lextt\cap\{\rt\sim \ep_0^{-1}\}  } r^3\left( \left|\widecheck{\trcht}\right|+\left|\widecheck{\trchbt}\right|\right)  &\les&  \ep_0.
\eea
\end{enumerate}
\end{proposition}

\begin{remark}\lab{rmk:whyonearthdoweneedep0squareandnotep0forinitaldata:thereturn}
According to \eqref{def:initialdatalayerassumptions}, we have the stronger bound $\Ik_{k_{large}+10}\leq \ep_0^2$. As it turns out, we only need the weaker bounds  $\Ik_{k_{large}+10}\leq \ep_0$ and $\,^{(ext)} \Ik_{3}\leq \ep_0^2$ as emphasized in Remark \ref{rmk:whyonearthdoweneedep0squareandnotep0forinitaldata}. In particular, $\,^{(ext)} \Ik_{3}\leq \ep_0^2$ is only used to ensure \eqref{eq:strongerboundfordivbcurlbtrchtrchbatrlikeep0minus1inIDL} which will be used in the proof of Theorem M0 and Theorem M6.
\end{remark}

\begin{proof}
We start with the first estimate. In view of the definition of $\Ikt_k$, Proposition \ref{Prop:transitionparam(E-Et)} and Lemma \ref{Lemma:geodesic-RicciCurv} imply
\beaa
\Ikt_k &\les& \sup_{\Lextt} r \left(\left|\dk^{\leq k}\left(\nab J^{(0)} +\frac{1}{r}\dual f_0\right)\right|+\left|\dk^{\leq k}\left(\nab f_0 -\frac{J^{(0)} }{r}\in\right)\right|\right)\\
&&+  \sup_{\Lextt}\Big(\left|\dk^{\leq k+1}\left(f_0-r\Re(\Jk)\right)\right|+r^{-1}+r|\dk^{\leq k+3}\Ga_b|\Big).
\eeaa
Together with the fact that $f_0= r\Re(\Jk)$, $r|\dk^{\leq k+3}\Ga_b|\les \Jk_{k+3}$, and $\, ^{(ext)} \Ik_{k_{large}+10} \leq \ep_0$, we infer, for $k\leq k_{large}+7$, 
\beaa
\Ikt_{k} &\les& \sup_{\Lextt} r \left(\left|\dk^{\leq k}\left(\nab J^{(0)} +\dual \Re(\Jk)\right)\right|+\left|\dk^{\leq k}\left(r\nab\Re(\Jk) -\frac{J^{(0)} }{r}\in\right)\right|\right)+\ep_0.
\eeaa
where we also used the fact that $r\gtrsim \ep_0^{-1}$ on $\Lextt$, since by definition $\st\gtrsim \ep_0^{-1}$ on $\Lextt$, and since $\st\sim\rt\sim r$ in view of Proposition \ref{Prop:transitionparam(E-Et)} and Lemma \ref{Lemma:geodesic-RicciCurv}. Since we have, in view of  Definitions \ref{def:renormalizationofallnonsmallquantitiesinPGstructurebyKerrvalue}  and \ref{definition.Ga_gGa_b},
\beaa
\nab J^{(0)} = - \dual\Re(\Jk) +\Ga_b,\qquad \nab\Re(\Jk)=\frac{J^{(0)}}{r^2}\in +O(r^{-4})+r^{-1}\Ga_b,
  \eeaa
  we deduce 
  \beaa
\Ikt_{k_{large}+7} &\les&  \ep_0
\eeaa
which is the first stated estimate. 

For the second estimate, we rely on the following identities in Lemma \ref{Lemma:geodesic-RicciCurv}   
\beaa
 \widecheck{\trcht} =O(r^{-4})+ r^{-1} \dk^{\le 3}\Ga_b,  \qquad  \widecheck{\trchbt} =O(r^{-4})+ r^{-1} \dk^{\le 3}\Ga_b, 
 \eeaa
\beaa
\divt \bt= O(r^{-6} ) + r^{-3} \dk^{\le 3} \Ga_b, \,\,\,\curlt \bt=\frac{6a_0m_0}{ \rt^5}\Jt^{(0)} +O(r^{-6})  + r^{-3} \dk^{\le 3} \Ga_b,
\eeaa
which yields
\beaa
 r^5\left( \left|\divt\bt\right|+\left|\curlt\bt -\frac{6a_0m_0}{\rt^5}\Jt^{(0)} \right|\right)   + r^3\left( \left|\widecheck{\trcht}\right|+\left|\widecheck{\trchbt}\right|\right)  &\les&  r^{-1}+r^2|\dk^{\leq 3}\Ga_b|.
\eeaa
We infer
\beaa
\nn \sup_{\Lextt\cap\{\rt\sim \ep_0^{-1}\}  } r^5\left( \left|\divt\bt\right|+\left|\curlt\bt -\frac{6a_0m_0}{\rt^5}\Jt^{(0)} \right|\right)  \\
 +\sup_{\Lextt\cap\{\rt\sim \ep_0^{-1}\}  } r^3\left( \left|\widecheck{\trcht}\right|+\left|\widecheck{\trchbt}\right|\right)  &\les&  \ep_0+\ep_0^{-1}\, ^{(ext)} \Ik_{3}
\eeaa
and the second stated estimate follows from the assumption $\, ^{(ext)} \Ik_{3} \leq \ep^2_0$. This concludes the proof of the proposition.
\end{proof}

%%%%%%%%%%%%%%%%%%%%%%%%%%%%%%%%%%%%%%%%

\section{Proof of Theorem M0}
\lab{sec:proofofTheoremM0}

%%%%%%%%%%%%%%%%%%%%%%%%%%%%%%%%%%%%%%%%

We recall the reader, see also Remark \ref{rmk:whyisTheoremM0onlyinchapter8},  that  Theorem M0  comes first in  the sequence of  steps,  Theorems M0--M8,  and therefore  its proof can only rely  on  the assumption \eqref{def:initialdatalayerassumptions} of our main theorem\footnote{In fact, under \eqref{def:initialdatalayerassumptions}, we have the stronger bound $\Ik_{k_{large}+10}\leq \ep_0^2$. As it turns out, we only need the weaker bounds  $\Ik_{k_{large}+10}\leq \ep_0$ and $\,^{(ext)} \Ik_{3}\leq \ep_0^2$ as emphasized in Remark \ref{rmk:whyonearthdoweneedep0squareandnotep0forinitaldata:thereturn}.}

\beaa
\Ik_{k_{large}+10}\leq \ep_0,\qquad \,^{(ext)} \Ik_{3}\leq \ep_0^2,
\eeaa
 as well as the fact that the space $\MM$    which we  consider is a   GCM admissible spacetime
 verifying our bootstrap assumptions ${\bf BA_\ep}$ made in section \ref{section:Bootstrap assumptions}.
 
{\bf Notation} \textit{Since we will be using  various frames in the proof  it is important to recall the main definitions. The PG structure 
of $\Mext$  is denoted by the usual symbols  $\{ E, r, u, \Jk, J\} $ where $E=\{ e_1, e_2, e_3, e_4\}$.  The PG  quantities
 of   the initial layer $\LL_0$ are denoted   with  $0$ indices, i.e.   $\{ E_0, r_0, u_0, \Jk, J_0\}$. The   quantities   related to the outgoing geodesic frame of $\Lextt$  are denoted by tildes, i.e.  $\{ \Et, \st, \rt  , \ut , \Jt\} $. We  shall also make use of  a second outgoing geodesic foliation  in $ \Lext $ defined starting  with the sphere  $S_1$, of the PG structure of $\Mext$ on $\Si_*$, whose related quantities will be denoted  by  primes.  At various stages of the proof,    when  only  two   foliations are needed,  we will  redefine notations accordingly.}

The proof of Theorem M0 proceeds in 24 steps which we summarize below for convenience:
\begin{enumerate}
\item In Steps 1--7, we propagate from $S_*$             along $\Si_*$ the $\ell=1$ modes of $\div\b$, $\curl\b$, $\rhoc$ and $\kabc$ to arrive at the estimate \eqref{eq:localbootassell1modedivbeta:proofThmM0:improved} on $S_1$.
 This  sequence of steps makes  use of the GCM assumptions on $S_*$ and the results of section \ref{section:Estimates-ell=1mpdesS_*}.

\item In Steps 8--16, we derive the control of $m-m_0$, see  \eqref{eq:finalesitateforffbchecklaonS1inproofofThmM0:reallyfinal} in Step 13, and 
 $a-a_0$, see \eqref{eq:controlofainthecasea0equal0ThM0} in Step 15 and \eqref{eq:controlofaminusa0andJminusJprimeforproofThM0} in Step 16. We also provide estimates for the transition coefficients  between the  auxiliary outgoing geodesic frame of $\Lextt$, introduced in section  \ref{section:geodesicfoliation8},  and the  frame of  $\Si_*$ induced  on the sphere $S_1=\Si_*\cap\{u=1\}$. In particular, we show that the sphere $S_1$ of $\Si_*$ is contained in $\Lextt\subset\Lext $.

\item In Steps 17--19, we control the change of frame coefficients  from the  outgoing PG frame of $\Lext$ to the outgoing PG  frame of $\Mext$ on the sphere $S_1$ in \eqref{eq:controlofffblabetweenPGLextandPGMextonS1}.

\item In steps 20--22, we propagate the control on  the change of frame coefficients  from the  outgoing PG frame of $\Lext$ to the outgoing PG  frame of $\Mext$  from $S_1$ to $\{u=1\}$ in \eqref{eq:controlofffblambdabetwenPGframeLextandMextonu=1}. 

\item In step 23, we control on  the change of frame coefficients  from the  ingoing PG frame of $\Lint$ to the ingoing PG  frame of $\Mint$ on $\{\ub=1\}$ in \eqref{eq:controlofffblambdabetwenPGframeLextandMextonu=1:ingoingcase}. 

\item Finally, we conclude the proof of Theorem M0 in step 24 by using the control of the change of frame coefficients, the control of the initial data layer, and the change of frame formulas to infer the control of the curvature components of $\Mext$ on $\{u=1\}$ and of $\Mint$ on $\{\ub=1\}$. 
\end{enumerate}

%%%%%%%%%%%%%%%

\subsection{Steps 1--7} 

%%%%%%%%%%%%%%%

We state the main result of Step 1--7 in the following
\begin{proposition}
\lab{Proposition:Step1-7.ThmM0}
The following estimates\footnote{Recall the definition of $a$ and $m$ in section \ref{sec:definitionofamthetandvphiadmissible}.} are true on $\Si_*$
\bea
\lab{eq:localbootassell1modedivbeta:proofThmM0:improved}
\bsplit
&\sup_{\Si_*}\left(r^5|(\div\b)_{\ell=1}|+r^5|(\curl\b)_{\ell=1,\pm}|+r^5\left|(\curl\b)_{\ell=1,0}-\frac{2am}{r^5}\right|\right)\les \ep_0,\\
&\sup_{\Si_*}\Big(r^3 |(\rhoc)_{\ell=1}|+r^2 |(\kabc)_{\ell=1}|\Big) \les \ep_0,
\end{split}
\eea
  where the quantities and the definition of the $\ell=1$ modes correspond to the frame of $\Si_*$. 
\end{proposition}

\begin{proof}
The proof,  which uses the  GCM conditions  on $S_*$  and the results    of section \ref{section:Estimates-ell=1mpdesS_*}       is outlined in Steps 1--7 below. 
\end{proof}

{\bf Step 1.}  We start by recalling Lemma  \ref{Le:Si*-ell=1modes} on the control of the $\ell=1$ basis on $\Si_*$.
\begin{lemma}\lab{Le:Si*-ell=1modes:proofofThmM0}
The     functions $\Jp$ verify the following properties
\begin{enumerate}
\item We have on $\Si_*$
\bea
\lab{eq:basicestimatesforJp-onSi_*:chap8} 
\bsplit
\int_{S}J^{(p)} &= O\left(\ep ru^{-\dec}\right),\\
\int_{S}J^{(p)}J^{(q)} &= \frac{4\pi}{3}r^2\de_{pq}+O\left(\ep ru^{-\dec}\right).
\end{split}
\eea

\item We have on $\Si_*$
\bea
\nab_\nu \Big[\big( r^2\De +2)\Jp\Big]=O\big(\dkb^{\le 1} \Ga_b\big). 
\eea

\item For  any  $k\le k_{small}$, we have on $\Si_*$
\beaa
\left|\dk_*^k\left( \De +\frac{2}{r^2}\right)    \Jp\right|    \les  \ep r^{-3} u^{-\frac{\dec}{2}}.
\eeaa

\item We have for any $k\le k_{small}$ on $\Si_*$
\beaa
\left|\dk_*^k\dds_2\dds_1 \Jp\right|    \les  \ep r^{-3} u^{-\frac{\dec}{2}},
\eeaa
where by $\dds_1\Jp$, we mean either $\dds_1(\Jp,0)$ or $\dds_1(0,\Jp)$. 
\end{enumerate}
\end{lemma}

Also, we recall  Lemma \ref{Le:ell=1modesonS_*} on the control of $(\rhoc-\frac{1}{2}\chih\c\chibh)_{\ell=1}$ on $S_*$. 
\begin{lemma}\lab{Le:ell=1modesonS_*:proofofThmM0}
The following holds on $S_*$
\bea
\left|\left(\rhoc-\frac{1}{2}\chih\c\chibh\right)_{\ell=1} \right| &\les&  \ep_0 r^{-3}u^{-2-2\dec}.
\eea
\end{lemma}

{\bf Step 2.} On $\Si_*$ we assume the following local bootstrap assumptions
\bea
\lab{eq:localbootassell1modedivbeta:proofThmM0}
\bsplit
&\sup_{\Si_*}\Big(r^5|(\div\b)_{\ell=1}|+r^5|(\curl\b)_{\ell=1,\pm}|+r^5\left|(\curl\b)_{\ell=1,0}-\frac{2am}{r^5}\right|\Big)\le \ep,\\
&\sup_{\Si_*}\Big(r^3u^{1+\dec}|(\rhoc)_{\ell=1}|+r^2u^{1+\dec}|(\kac)_{\ell=1}|\Big) \le \ep,
\end{split}
\eea
 which  will be improved in Steps 2--7.

We start with the control of $(\div\ze)_{\ell=1}$. Recall  the following consequence of the  Codazzi equation for $\chih$
\beaa
\ddd_2\chih &=& \frac{1}{r}\ze - \b+\Ga_g\c \Ga_g.
\eeaa
Differentiating w.r.t. $\div$, we infer
\beaa
\div\ddd_2\chih &=& \frac{1}{r}\div\ze - \div\b+r^{-1}\dkb^{\leq 1}(\Ga_g\c \Ga_g).
\eeaa
Projecting on the $\ell=1$ modes, this yields 
\beaa
(\div\ddd_2\chih)_{\ell=1} &=&  \frac{1}{r}(\div\ze)_{\ell=1}-(\div\b)_{\ell=1} +r^{-1}\dkb^{\leq 1}(\Ga_g\c \Ga_g).
\eeaa
Next, we estimate $(\ddd_1\ddd_2\chih)_{\ell=1}$. We have
\beaa
(\div\ddd_2\chih)_{\ell=1,p}=\frac{1}{|S|}\int_S\div\ddd_2\chih\Jp =\frac{1}{|S|}\int_S\chih\c\dds_2\nab\Jp
\eeaa
and hence
\beaa
|(\div\ddd_2\chih)_{\ell=1}| &\les& |\dds_2\dds_1\Jp||\Ga_g|.
\eeaa
We deduce
\beaa
 (\div\ze)_{\ell=1} &=& r(\div\b)_{\ell=1} + r|\dds_2\dds_1\Jp\Ga_g|+\dkb^{\leq 1}(\Ga_g\c \Ga_g).
\eeaa
Together with the local bootstrap assumption \eqref{eq:localbootassell1modedivbeta:proofThmM0} for $(\div\b)_{\ell=1}$, the control of Lemma \ref{Le:Si*-ell=1modes:proofofThmM0} for $\dds_2\dds_1\Jp$, and the bootstrap assumptions  for $\Ga_g$, we infer on $\Si_*$
\bea
|(\div\ze)_{\ell=1}| &\les& \frac{\ep}{r^4}.
\eea

{\bf Step 3.} Next, we control of $(\div\bb)_{\ell=1}$. Recall  the following consequence of the  Codazzi equation for $\chibh$
\beaa
\div\chibh &=& \frac{1}{2}\nab\kabc +\frac{\Up}{r}\ze  +\bb+\Ga_b\c \Ga_g.
\eeaa
Differentiating w.r.t. $\div$, we infer
\beaa
\div\ddd_2\chibh &=& \frac{1}{2}\Delta\kabc+\frac{\Up}{r}\div\ze + \div\bb+r^{-1}\dkb^{\leq 1}(\Ga_b\c \Ga_g).
\eeaa
Projecting on the $\ell=1$ modes, this yields 
\beaa
(\div\ddd_2\chibh)_{\ell=1} &=&  \frac{1}{2}(\Delta\kabc)_{\ell=1}+\frac{\Up}{r}(\div\ze)_{\ell=1} +(\div\bb)_{\ell=1} +r^{-1}\dkb^{\leq 1}(\Ga_b\c \Ga_g).
\eeaa
As in Step 2, we have
\beaa
|(\div\ddd_2\chibh)_{\ell=1}| &\les& |\dds_2\dds_1\Jp||\Ga_b|.
\eeaa
Also, we have
\beaa
(\Delta\kabc)_{\ell=1,p}=\frac{1}{|S|}\int_S\Delta\kabc \Jp =-\frac{2}{r^2}\frac{1}{|S|}\int_S\kabc\Jp+\frac{1}{|S|}\int_S\kabc\left(\Delta+\frac{2}{r^2}\right)\Jp
\eeaa
and hence
\beaa
|(\Delta\kabc)_{\ell=1}| &\les& r^{-2}|(\kabc)_{\ell=1}|+\left|\left(\Delta+\frac{2}{r^2}\right)\Jp\right||\Ga_g|.
\eeaa
We deduce
\beaa
|(\div\bb)_{\ell=1}| &\les& r^{-2}|(\kabc)_{\ell=1}|+r^{-1}|(\div\ze)_{\ell=1}|+\left|\left(\Delta+\frac{2}{r^2}\right)\Jp\right||\Ga_g|\\
&&+ |\dds_2\dds_1\Jp||\Ga_b|+r^{-1}|\dkb^{\leq 1}(\Ga_b\c \Ga_g)|.
\eeaa
Together with the local bootstrap assumption \eqref{eq:localbootassell1modedivbeta:proofThmM0} for $(\kabc)_{\ell=1}$, the control for $(\div\ze)_{\ell=1}$ of Step 2, the control of Lemma \ref{Le:Si*-ell=1modes:proofofThmM0} for $\dds_2\dds_1\Jp$ and $(\Delta+\frac{2}{r^2})\Jp$, and the bootstrap assumptions  for $\Ga_b$ and $\Ga_g$, we infer on $\Si_*$
\beaa
|(\div\bb)_{\ell=1}| &\les& \frac{\ep}{r^4u^{1+\dec}}+\frac{\ep}{r^5}.
\eeaa
Using the dominance condition \eqref{eq:behaviorofronS-star} on $r$ on $\Si_*$, we infer on $\Si_*$
\bea
|(\div\bb)_{\ell=1}| &\les& \frac{\ep}{r^4u^{1+\dec}}.
\eea

{\bf Step 4.} Recall from Corollary \ref{corofLemma:transport.alongSi_*1} the following transport along $\Si_*$, for $p=0,+,-$,
\beaa
\nn&&\nu\left(\int_S\left( \lap\kabc+\frac{2\Up}{r}\div \ze\right)\Jp\right)\\
\nn&=& O(r^{-3})\int_S\kabc\Jp   +O(r^{-2})\int_S\div\ze\Jp +O(r^{-1})\int_S\div\bb\Jp     +O(r^{-2})\int_S\rhoc\Jp \\
 &&+O(r^{-1})\int_S\div\b\Jp +r\left|\left(\Delta+\frac{2}{r^2}\right)\Jp\right|\dkb^{\leq 1}\Ga_b+\dkb^{\le 2 }(\Ga_b\c \Ga_b),
\eeaa
and
\beaa
\nn\nu\left(\int_S\left(\rhoc - \frac{1}{2}\chih\c\chibh \right)\Jp\right) &=&  -\int_S\div\bb\Jp -(1+O(r^{-1}))\int_S\div\b\Jp \\
\nn&&  +O(r^{-1})\int_S\left(\rhoc - \frac{1}{2}\chih\c\chibh \right)\Jp +O(r^{-3})\int_S\kabc\Jp\\
 \nn  && +O(r^{-2})\int_S\div\ze\Jp+r\left|\left(\Delta+\frac{2}{r^2}\right)\Jp\right|\Ga_b \\
   &&+r\dkb^{\le  1}( \Ga_b \c \Ga_b),
\eeaa
where by $\dds_1\Jp$, we mean $\dds_1(\Jp,0)$ or $\dds_1(0,\Jp)$. Together with the local bootstrap assumption \eqref{eq:localbootassell1modedivbeta:proofThmM0} for $(\div\b)_{\ell=1}$, $(\kabc)_{\ell=1}$, and $(\rhoc)_{\ell=1}$, 
the control for $(\div\ze)_{\ell=1}$ of Step 2, the control for $(\div\bb)_{\ell=1}$ of Step 3,  the control of Lemma \ref{Le:Si*-ell=1modes:proofofThmM0} for $(\Delta+\frac{2}{r^2})\Jp$, and the bootstrap assumptions  for $\Ga_b$ and $\Ga_g$, we infer on $\Si_*$, for $p=0,+,-$,
\beaa
\left|\nu\left(\int_S\left( \lap\kabc+\frac{2\Up}{r}\div \ze\right)\Jp\right)\right| &\les& \frac{\ep^2}{r^2u^{2+2\dec}}+\frac{\ep}{r^3u^{1+\dec}}+\frac{\ep}{r^4},\\
\left|\nu\left(\int_S\left(\rhoc - \frac{1}{2}\chih\c\chibh \right)\Jp\right)\right| &\les& \frac{\ep^2}{ru^{2+2\dec}}+\frac{\ep}{r^2u^{1+\dec}}+\frac{\ep}{r^3}.
\eeaa
Using the dominance condition \eqref{eq:behaviorofronS-star} on $r$ on $\Si_*$, we infer on $\Si_*$, for $p=0,+,-$,
\beaa
\left|\nu\left(\int_S\left( \lap\kabc+\frac{2\Up}{r}\div \ze\right)\Jp\right)\right| &\les& \frac{\ep_0}{r^2u^{2+2\dec}},\\
\left|\nu\left(\int_S\left(\rhoc - \frac{1}{2}\chih\c\chibh \right)\Jp\right)\right| &\les& \frac{\ep_0}{ru^{2+2\dec}}.
\eeaa
Integrating from $S_*$, and using the fact that $\nu(u)=2+O(\ep)$, we deduce on $\Si_*$, for $p=0,+,-$,
\beaa
\left|\int_S\left( \lap\kabc+\frac{2\Up}{r}\div \ze\right)\Jp\right| &\les& \left|\int_{S_*}\left( \lap\kabc+\frac{2\Up}{r}\div \ze\right)\Jp\right|+\frac{\ep_0}{r^2u^{1+\dec}},\\
\left|\int_S\left(\rhoc - \frac{1}{2}\chih\c\chibh \right)\Jp\right| &\les& \left|\int_{S_*}\left(\rhoc - \frac{1}{2}\chih\c\chibh \right)\Jp\right|+\frac{\ep_0}{ru^{1+\dec}}.
\eeaa
Now, since $\kabc=0$ on $S_*$ in view of our GCM conditions,  and in view of the control $(\rhoc - \frac{1}{2}\chih\c\chibh)_{\ell=1}$ on $S_*$ provided by Lemma \ref{Le:ell=1modesonS_*:proofofThmM0}, and the control for $(\div\ze)_{\ell=1}$ of Step 2, we have
\beaa
\left|\int_{S_*}\left( \lap\kabc+\frac{2\Up}{r}\div \ze\right)\Jp\right| &\les& r|(\div\ze)_{\ell=1}|\les\frac{\ep}{r^3},\\
 \left|\int_{S_*}\left(\rhoc - \frac{1}{2}\chih\c\chibh \right)\Jp\right| &\les& \frac{\ep_0}{ru^{1+\dec}},
\eeaa
and hence, we obtain on $\Si_*$, for $p=0,+,-$,
\beaa
\left|\int_S\left( \lap\kabc+\frac{2\Up}{r}\div \ze\right)\Jp\right| &\les& \frac{\ep_0}{r^2u^{1+\dec}}+\frac{\ep}{r^3}\les \frac{\ep_0}{r^2u^{1+\dec}}
\eeaa
where we used in the last inequality the dominance condition \eqref{eq:behaviorofronS-star} on $r$ on $\Si_*$, and 
\beaa
\left|\int_S\left(\rhoc - \frac{1}{2}\chih\c\chibh \right)\Jp\right| &\les& \frac{\ep_0}{ru^{1+\dec}}.
\eeaa
We infer on $\Si_*$, for $p=0,+,-$,
\beaa
\left|\int_S\kabc \Jp\right| &\les& r^2\left|\int_S\left( \lap+\frac{2}{r^2}\right)\kabc\Jp\right|+r\left|\int_S\div \ze \Jp\right| +\frac{\ep_0}{u^{1+\dec}},\\
\left|\int_S\rhoc\Jp\right| &\les& r^2\Ga_b\c\Ga_g+\frac{\ep_0}{ru^{1+\dec}}.
\eeaa
Together with the control for $(\div\ze)_{\ell=1}$ of Step 2, the control of Lemma \ref{Le:Si*-ell=1modes:proofofThmM0} for $(\Delta+\frac{2}{r^2})\Jp$, and the bootstrap assumptions  for $\Ga_b$ and $\Ga_g$, we infer on $\Si_*$, 
\beaa
r|(\rhoc)_{\ell=1}|+|(\kabc)_{\ell=1}| &\les& \frac{\ep_0}{r^2u^{1+\dec}}+\frac{\ep}{r^3}.
\eeaa
Using the dominance condition \eqref{eq:behaviorofronS-star} on $r$ on $\Si_*$, we deduce
\bea\lab{eq:controlofell=1modeforthorcandkabconSistarinproofoThmM0Step4}
\sup_{\Si_*}\Big(r^3u^{1+\dec}|(\rhoc)_{\ell=1}|+r^2u^{1+\dec}|(\kabc)_{\ell=1}|\Big) &\les& \ep_0.
\eea

{\bf Step 5.} Recall from Corollary \ref{corofLemma:transport.alongSi_*1} the following transport equation along $\Si_*$, for $p=0,+,-$,
\beaa
\nu\left(\int_S\div\b\Jp\right) &=& O(r^{-1})\int_S\div\b\Jp+O(r^{-2})\int_S\rhoc\Jp\\
\nn&&+r\left(\left|\left(\Delta+\frac{2}{r^2}\right)\Jp\right|+\left|\dds_2\dds_1\Jp\right|\right)\Ga_g+\dkb^{\leq 1}(\Ga_b\c\Ga_g).
\eeaa
Together with the local bootstrap assumption \eqref{eq:localbootassell1modedivbeta:proofThmM0} for $(\div\b)_{\ell=1}$,  
the control for $(\rhoc)_{\ell=1}$ of Step 4,  the control of Lemma \ref{Le:Si*-ell=1modes:proofofThmM0} for $(\Delta+\frac{2}{r^2})\Jp$ and $\dds_2\dds_1\Jp$, and the bootstrap assumptions  for $\Ga_b$ and $\Ga_g$, we infer on $\Si_*$, for $p=0,+,-$,
\beaa
\left|\nu\left(\int_S\div\b\Jp\right)\right| &\les& \frac{\ep_0+\ep^2}{r^3u^{1+\dec}}+\frac{\ep}{r^4},
\eeaa
and using the dominance condition \eqref{eq:behaviorofronS-star} on $r$ on $\Si_*$, we deduce
\beaa
\left|\nu\left(\int_S\div\b\Jp\right)\right| &\les& \frac{\ep_0}{r^3u^{1+\dec}}.
\eeaa
Integrating from $S_*$, using the fact that $\nu(u)=2+O(\ep)$, and the GCM condition $(\div\b)_{\ell=1}=0$ on $S_*$, 
we deduce on $\Si_*$, for $p=0,+,-$,
\bea
\sup_{\Si_*}r^5|(\div\b)_{\ell=1}| &\les& \ep_0.
\eea

{\bf Step 6.} Next, we control $(\rhod)_{\ell=1}$. To this end, we first control $(\curl\ze)_{\ell=1}$. Recall  the following consequence of the  Codazzi equation for $\chih$
\beaa
\ddd_2\chih &=& \frac{1}{r}\ze - \b+\Ga_g\c \Ga_g.
\eeaa
Differentiating w.r.t. $\curl$, we infer
\beaa
\curl\ddd_2\chih &=& \frac{1}{r}\curl\ze - \curl\b+r^{-1}\dkb^{\leq 1}(\Ga_g\c \Ga_g).
\eeaa
Proceeding as for the control of $(\div\ze)_{\ell=1}$ in Step 2, we infer
\beaa
 (\curl\ze)_{\ell=1} &=& r(\curl\b)_{\ell=1} + r|\dds_2\dds_1\Jp|\Ga_g+\dkb^{\leq 1}(\Ga_g\c \Ga_g).
\eeaa
Together with the local bootstrap assumption \eqref{eq:localbootassell1modedivbeta:proofThmM0} for $(\curl\b)_{\ell=1}$, the control of Lemma \ref{Le:Si*-ell=1modes:proofofThmM0} for $\dds_2\dds_1\Jp$, and the bootstrap assumptions  for $\Ga_g$, we infer on $\Si_*$
\beaa
|(\curl\ze)_{\ell=1,\pm}|+\left|(\curl\ze)_{\ell=1,0}-\frac{2am}{r^4}\right| &\les& \frac{\ep}{r^4}.
\eeaa
Now, note that  we have from  the null structure equations 
\beaa
\rhod=\curl\ze+\Ga_b\c\Ga_g.
\eeaa
Together with the above control of $(\curl\ze)_{\ell=1}$ and the bootstrap assumptions  for $\Ga_g$, we infer on $\Si_*$
\beaa
|(\rhod)_{\ell=1,\pm}|+\left|(\rhod)_{\ell=1,0}-\frac{2am}{r^4}\right| &\les& \frac{\ep}{r^4}+\frac{\ep^2}{r^3u^{1+\dec}}
\eeaa
and using the dominance condition \eqref{eq:behaviorofronS-star} on $r$ on $\Si_*$, we deduce
\bea
|(\rhod)_{\ell=1,\pm}|+\left|(\rhod)_{\ell=1,0}-\frac{2am}{r^4}\right| &\les& \frac{\ep_0}{r^3u^{1+\dec}}.
\eea

{\bf Step 7.} Recall from Corollary \ref{corofLemma:transport.alongSi_*1} the following transport along $\Si_*$, for $p=0,+,-$,
\beaa
\nn\nu\left(\int_S\curl\b\Jp\right) &=&  \frac{4}{r}(1+ O(r^{-1}))\int_S\curl\b\Jp+\frac{2}{r^2}(1+ O(r^{-1}))\int_S\rhod\Jp\\
\nn&&+r\left(\left|\left(\Delta+\frac{2}{r^2}\right)\Jp\right|+\left|\dds_2\dds_1\Jp\right|\right)\Ga_g+\dkb^{\leq 1}(\Ga_b\c\Ga_g).
\eeaa
In the case $p=\pm$, we have $(\curl\b)_{\ell=1, \pm}=0$ on $S_*$, so using Step 6 to control $(\rhod)_{\ell=1,\pm}$, and arguing exactly as for the control of $(\div\b)_{\ell=1}$ in Step 5, we obtain the corresponding estimate, i.e.
\bea
\sup_{\Si_*}r^5|(\curl\b)_{\ell=1,\pm}| &\les& \ep_0.
\eea 

Next, we focus on the case $p=0$. We rewrite the above transport equation in this particular case
\beaa
\nu\left(\int_S\curl\b J^{(0)}\right) &=&  \frac{4}{r}(1+ O(r^{-1}))\int_S\curl\b J^{(0)}+\frac{2}{r^2}(1+ O(r^{-1}))\int_S\rhod J^{(0)}\\
\nn&&+r\left(\left|\left(\Delta+\frac{2}{r^2}\right)J^{(0)}\right|+\left|\dds_2\dds_1J^{(0)}\right|\right)\Ga_g+\dkb^{\leq 1}(\Ga_b\c\Ga_g).
\eeaa 
Since $\nu(r)=-2+r\Ga_b$, we have
\beaa
\nu\left(r^3\int_S\curl\b J^{(0)}\right) &=& r^3\nu\left(\int_S\curl\b J^{(0)}\right)+3r^2\nu(r)\int_S\curl\b J^{(0)}\\
&=& r^3\nu\left(\int_S\curl\b J^{(0)}\right) -6r^2\int_S\curl\b J^{(0)}+r^5\Ga_b(\curl\b)_{\ell=1,0},
\eeaa
and hence
\beaa
&&\nu\left(r^3\int_S\curl\b J^{(0)}-8\pi am\right)\\
 &=&  -\frac{2}{r}r^3(1+ O(r^{-1}))\int_S\curl\b J^{(0)}+2r(1+ O(r^{-1}))\int_S\rhod J^{(0)}\\
\nn&&+r^5\Ga_b(\curl\b)_{\ell=1,0}+r^4\left(\left|\left(\Delta+\frac{2}{r^2}\right)J^{(0)}\right|+\left|\dds_2\dds_1J^{(0)}\right|\right)\Ga_g+r^3\dkb^{\leq 1}(\Ga_b\c\Ga_g).
\eeaa
In view of  the local bootstrap assumption \eqref{eq:localbootassell1modedivbeta:proofThmM0} for $(\curl\b)_{\ell=1,0}$,  the control of $(\rhod)_{\ell=1}$ of Step 6, the control of Lemma \ref{Le:Si*-ell=1modes:proofofThmM0} for $(\Delta+\frac{2}{r^2})\Jp$ and $\dds_2\dds_1\Jp$, and the bootstrap assumptions  for $\Ga_g$, we infer on $\Si_*$
\beaa
\left|\nu\left(r^3\int_S\curl\b J^{(0)} -8\pi am\right)\right| &\les& \frac{1}{r}+\frac{\ep_0+\ep^2}{u^{1+\dec}}.
\eeaa
Using the dominance condition \eqref{eq:behaviorofronS-star} on $r$ on $\Si_*$, we obtain
\beaa
\left|\nu\left(r^3\int_S\curl\b J^{(0)} -8\pi am\right)\right| &\les& \frac{\ep_0}{u^{1+\dec}}.
\eeaa
Integrating from $S_*$ where there holds $(\curl\b)_{\ell=1,0}=\frac{2am}{r^5}$ on $S_*$, we deduce on $\Si_*$
\bea
\left|(\curl\b)_{\ell=1,0}-\frac{2am}{r^5}\right| &\les& \frac{\ep_0}{r^5}.
\eea
Thus, in view of the control for $(\kabc)_{\ell=1}$ and $(\rhoc)_{\ell=1}$ of Step 4, the control for $(\div\b)_{\ell=1}$ of Step 5, and the control for $(\curl\b)_{\ell=1}$ of this step, we have finally obtained on $\Si_*$ the desired estimate \eqref{eq:localbootassell1modedivbeta:proofThmM0:improved}, i.e. 
\beaa
\bsplit
&\sup_{\Si_*}\left(r^5|(\div\b)_{\ell=1}|+r^5|(\curl\b)_{\ell=1,\pm}|+r^5\left|(\curl\b)_{\ell=1,0}-\frac{2am}{r^5}\right|\right)\les \ep_0,\\
&\sup_{\Si_*}\Big(r^3u^{1+\dec}|(\rhoc)_{\ell=1}|+r^2u^{1+\dec}|(\kac)_{\ell=1}|\Big) \les \ep_0,
\end{split}
\eeaa
thus improving  the local bootstrap assumption \eqref{eq:localbootassell1modedivbeta:proofThmM0}, and concluding the proof of Proposition \ref{Proposition:Step1-7.ThmM0}.

%%%%%%%%%%%%%%

\subsection{Steps 8--16}  

%%%%%%%%%%%%%%

 As a consequence of  \eqref{eq:localbootassell1modedivbeta:proofThmM0:improved}  we have,  in particular on $S_1=\Si_*\cap\{u=1\}$, the estimate
\bea\lab{eq:usefulesitateinthediscussionofThmM0:00}
\sup_{S_1}\left(r^2\left|(\kabc)_{\ell=1}\right| +r^5\left|(\div\b)_{\ell=1}\right|\right) &\les& \ep_0.
\eea
On $S_1$, we also have the GCM conditions 
\bea
\lab{eq:GCMonS1}
\ka=\frac{2}{r}, \quad \kab=-\frac{2\Up}{r}+\underline{C}_0+\sum_p\underline{C}_p\Jp, \quad \mu=\frac{2m}{r^3}+M_0+\sum_pM_p\Jp.
\eea

We introduce the following auxiliary construction.

\begin{definition}[The outgoing cone $\CC'_1$]
 Starting with the sphere  $S_1$  we define the outgoing geodesic null cone  $\CC'_1$ emanating from $S_1$ in the direction of $e_4$.  We denote: 
 \begin{itemize}
\item by $e_4'$ the geodesic extension of $e_4$, and by  $s'$ the corresponding affine parameter, i.e. $e_4'(s')=1$, normalized such that $s'=r$ on $S_1$,

\item by   $S'$ the spheres of constant $s'$ along $\CC'_1$ and by $r'$ the corresponding area radius,

\item by ${J'}^{(p)}$  the basis of $\ell=1$ modes verifying $e_4'({J'}^{(p)})=0$  with ${J'}^{(p)}=\Jp$ on $S_1$, for $p=0,+,-$. 
\end{itemize}
We restrict $\CC_1'$ to the region  $\{\de_*\ep_0^{-1}\leq r'\leq r(S_1)\}$. With this restriction, we will show below, 
see \eqref{remark:applyLemma1.4}, that $\CC'_1\subset\Lextt$. Finally, we denote by $(f, \fb, \la)$ the transition coefficients from the outgoing geodesic null frame of $\Lextt$ of section \ref{section:geodesicfoliation8} to the outgoing geodesic null frame of $\CC'_1$  initialized on  $S'_1=S_1$.
 \end{definition}

{\bf   Local Bootstrap Assumptions:} 
\begin{enumerate}
\item   Along $\CC'_1$, we assume 
\bea\lab{eq:localbootstrapssutionfortheproofofTheoremM0}
\sup_{S'\subset \CC'_1}\Big(\|f\|_{\hk_4(S')}+(r')^{-1}\|(\fb, \log\la)\|_{\hk_4(S')}\Big) \leq \ep.
\eea

\item  On $S'_1=S_1$, we assume
\bea\lab{eq:iterativeapssutiononS1fortheproofofTheoremM0}
\|f\|_{\hk_{k_{large}}(S'_1)}+r^{-1}\|(\fb, \log(\la))\|_{\hk_{k_{large}}(S'_1)} &\leq& \ep.
\eea

\item   In the case $a_0\neq 0$, we make the following assumption\footnote{Recall that $\dk_*$  refers to the properly normalized  tangential derivatives  along $\Si_*$, i.e. $\dk_*=(\dkb, \nab_\nu)$.}, on $S'_1$,  on the difference between the basis of $\ell=1$ modes $\Jp$ of $\Si_*$, and the basis of $\ell=1$ modes $\Jt^{(p)}$ of $\Lextt$
\bea\lab{eq:bootstrapassumptiondifferenceell=1basisofCC1andLextprime}
\max_{p=0,+,-}\left\|\dk_*^{\leq k_{large}}(\Jp-\Jt^{(p)})\right\|_{L^\infty(S'_1)} &\leq& \ep.
\eea
\end{enumerate}

\begin{remark}
\eqref{eq:iterativeapssutiononS1fortheproofofTheoremM0}  will be improved in Step 13, \eqref{eq:localbootstrapssutionfortheproofofTheoremM0} will be improve in Step 14, and  \eqref{eq:bootstrapassumptiondifferenceell=1basisofCC1andLextprime}  will be improved in Step 17.
\end{remark}

\begin{remark}
The discrepancy between  the  top number of derivatives in \eqref{eq:iterativeapssutiononS1fortheproofofTheoremM0} and  \eqref{eq:bootstrapassumptiondifferenceell=1basisofCC1andLextprime} is due to the fact that the former depends only on the  total number of derivatives allowed in $\Lextt$ while the  latter reflects the total number of derivatives allowed by the global bootstrap assumptions on $\Mext$. 
\end{remark} 

\begin{remark}\lab{rmk:anomalousbehaviorforfblaandnoproblemcurvatureesitmateproofThmM0}
The anomalous behavior for $\fb$ and $\la$ in \eqref{eq:localbootstrapssutionfortheproofofTheoremM0} \eqref{eq:iterativeapssutiononS1fortheproofofTheoremM0}, i.e. the fact that they display a $r$ loss compared to $f$, does not affect the desired estimates for the curvature components,  see \eqref{eq:gatheringtheconclusionsyieldingtheproofofThmM0:1}. This is due to the fact that, in the change of frame formulas for the curvature components, $\la$ and $\fb$ are multiplied by terms that decay faster in $r$. We refer to Remark \ref{rmk:anomalousbehaviorforfblaandLorentzboostproofThmM0} for a heuristic explanation of this anomalous behavior.
\end{remark}

\begin{remark}
Note that the control of $J_{\Lextt}^{(\pm)}$ provided by Proposition \ref{proposition:geodesicfoliationLextt} is invariant under the change\footnote{This holds true provided one also changes the quantities $J_{\Lext}^{(\pm)}$ and $\Jk_{\pm, \Lext}$, associated to the part $\Lext$ of the initial data layer, according to
\beaa
J_{\Lext}^{(\pm)} \to \cos(\vphi_0)J_{\Lext}^{(\pm)} \pm \sin(\vphi_0)J_{\Lext}^{(\mp)},\qquad  \Jk_{\pm, \Lext} \to \cos(\vphi_0)\Jk_{\pm, \Lext} \pm \sin(\vphi_0)\Jk_{\mp, \Lext}.
\eeaa
Note that these transformations leave invariant the linearized quantities in Definition \ref{def:renormalizationofallnonsmallquantitiesinPGstructurebyKerrvalue}, and hence $(\Ga_b, \Ga_g)$ in Definition \ref{definition.Ga_gGa_b}. This corresponds to the fact that Kerr is axially symmetric.}  
\beaa
J_{\Lextt}^{(\pm)} \longrightarrow \cos(\vphi_0)J_{\Lextt}^{(\pm)} \pm \sin(\vphi_0)J_{\Lextt}^{(\mp)}, \qquad\varphi_0\in[0,2\pi),
\eeaa
which corresponds to the invariance $\vphi_{\Lextt}\to \vphi_{\Lextt}-\vphi_0$, and hence to the fact that Kerr is axially symmetric. We may thus use this freedom to choose the particular pair $J_{\Lextt}^{(\pm)}$ such that 
\bea\lab{eq:fixtheaxialfreefomofchoiceofJpmbasis}
\int_{S_1}J_{\Lextt}^{(+)}J^{(-)} &=& 0.
\eea
This normalization will be needed to improve the bootstrap assumptions \eqref{eq:bootstrapassumptiondifferenceell=1basisofCC1andLextprime}.
\end{remark}

Assumption \eqref{eq:localbootstrapssutionfortheproofofTheoremM0} allows  to apply   Lemma 7.3 in \cite{KS-GCM1}  (recalled here in   Lemma \ref{lemma:consequencedeformationsurfaceusedinTheoremM0})  with $\de_1=\ep$ on each sphere\footnote{Note  that each sphere $S'$ of $\CC_1'$ intersects a unique sphere $\St$ of $\Lextt$ at its south pole. Hence,  we may consider $S'$ as a deformation of $\St$.}  $S'=S( s')$ of $\CC'_1$.   This allows us to deduce the following comparison estimate
\bea\lab{eq:anamolousdiffrenceofthesextminussextl}
\sup_{\CC'_1}\Big( \rt^{-1}|r'  -\rt | + |1- \ut| + \rt^{-1}|s' - \st |\Big) &\les& \ep.
\eea
In particular, since  $r'\geq \de_*\ep_0^{-1}$ on $\CC'_1$, we infer
\beaa
\sup_{\CC'_1}|\ut-1|\les \ep, \qquad \inf_{\CC'_1} \st\geq\frac{\de_*}{2}\ep_0^{-1}.
\eeaa
Since $\Lextt=\{0\leq \ut\leq 2, \,\st\geq\frac{\de_*}{2}\ep_0^{-1}\}$,  we deduce  
\bea\lab{remark:applyLemma1.4}
\CC'_1\subset \Lextt.
\eea
We note also that $r, \rt, r'$ are all comparable  along $\CC_1'$.

Also, since the sphere $S_1'=S_1$ satisfies the following:
\begin{itemize}
\item $S_1'$ is a sphere of $\Mext$ in $\Lextt$,

\item $S_1'$ is a sphere of the GCM hypersurface $\Si_*$,

\item the estimate \eqref{eq:usefulesitateinthediscussionofThmM0:00} holds on $S_1'$,

\item the estimate \eqref{eq:iterativeapssutiononS1fortheproofofTheoremM0} holds on $S_1'$,

\item the estimate \eqref{eq:bootstrapassumptiondifferenceell=1basisofCC1andLextprime} holds on $S_1'$,
\end{itemize}
we can invoke      Proposition 8.1 in \cite{KS-GCM2}   (restated  here in  Proposition \ref{corr:GCM-rigidity2}) with the choice $\epg=\dg=\ep_0$, $\de_1=\ep$, $s_{max}=k_{large}$, and with the background foliation being the outgoing geodesic foliation of the part $\Lextt$  of the initial data layer. We obtain 
\bea\lab{eq:usefulesitateinthediscussionofThmM0:00:first}
r^{-1}\| (f, \fb, \la-\ov{\la}^{S_1'}) \|_{\hk_{k_{large}+1}(S'_1)}   &\les& \ep_0,
\eea
and,  with    $\ov{\la}^{S'_1}$ denoting the average of $\la$ on $S'_1$,
\bea\lab{eq:usefulesitateinthediscussionofThmM0:00:first:ovla}
|\ov{\la}^{S'_1} -1| &\les& \ep_0+r^{-1}\sup_{S'_1}\left|r- \rt\right|.
\eea

\begin{remark}
In order to improve the bootstrap assumption \eqref{eq:iterativeapssutiononS1fortheproofofTheoremM0}   we need in particular to improve the estimate for $f$ in \eqref{eq:usefulesitateinthediscussionofThmM0:00:first} by $r^{-1}$. Obtaining this  improvement is the focus\footnote{In fact,  all intermediate estimates in Step 8 to 13 are only  needed   to derive this improvement on $f$.} of Step 8 to 13.
\end{remark}

\begin{remark}\lab{rmk:anomalousbehaviorforfblaandLorentzboostproofThmM0}
In view of  \eqref{eq:anamolousdiffrenceofthesextminussextl}, while $|1 - \ut| \les\ep$ on $S_1$, we only  have  $|s' - \st | \les r\ep$ on $S_1$. This, as well as the anomalous behavior of $\fb$ mentioned in Remark \ref{rmk:anomalousbehaviorforfblaandnoproblemcurvatureesitmateproofThmM0}, shows that  the sphere $S_1$ is a large deformation, along the outgoing direction, of spheres of the initial data layer $\Lextt$. This reflects the fact that $S_1$ (and $\Si_*$) captures the center of mass frame of the limiting Kerr solution, while the initial data layer foliation captures the center of mass frame of the initial Kerr solution. The behavior of $s' - \st$, as well as the one of $\fb$, is consistent with the presence of a  large Lorentz boost between these two center of mass frames.
\end{remark}

{\bf Local  Notation for Steps 8--18.}  \textit{In Steps 8--18 below, $(f, \fb, \la)$ denote  the transition coefficients from the frame  $\widetilde{E}=(\tilde{e}_1, \tilde{e}_2, \tilde{e}_3, \tilde{e}_4)$ of $\Lextt$ to the prime one of $\CC_1'$. Also, we will only make use of the prime frame along $\CC_1'$ and the geodesic frame of $\Lextt$. Since we are not making reference here to the PG structure of $\Lext$ and the one of $\Mext$, we drop the tilde on the quantities   associated to $\Lextt$.   }

 {\bf Step 8.}   We derive estimates for $\b' $  and $\dds_1'(-\rho', \rhod')$  on $S_1'$ with the help  of  transformation formulas of Proposition \ref{Proposition:transformationRicci}. To start with, we make use of 
\beaa
\b'&=& \la\left(\b +\frac 3 2\big(  f \rho+\dual  f  \rhod\big)+ \frac 1 2 \a\c\fb+\lot\right),
\eeaa
together with the  estimate \eqref{eq:usefulesitateinthediscussionofThmM0:00:first} for $f$ to estimate the linear terms $\rho f$ and $\rhod\dual f$, the estimate \eqref{eq:iterativeapssutiononS1fortheproofofTheoremM0} for $(f,\fb, \la)$ to estimate the other terms, and the control provided by Proposition \ref{proposition:geodesicfoliationLextt} for the  part $\Lextt$ of the initial data layer to control\footnote{In order to control $\hk_j(S_1')$ norms by sup norms, we use, here and in the remainder of the proof, Lemma 7.3 in \cite{KS-GCM1} restated  here in Lemma \ref{lemma:consequencedeformationsurfaceusedinTheoremM0}.}  the $\hk_j(S_{1}')$ norm of the Ricci coefficients and curvature components of the initial data foliation of $\Lextt$ in terms of their sup norm, we have
\bea
\sup_{k\le  k_{large}} r^2\|\dkb'^k\b'\|_{L^2(S_{1}')} &\les& \ep_0.
\eea

Also, we have
\beaa
\rho' &=& \rho + \fb\c\b - f\c\bb +\frac{3}{2}\rho(f\c\fb) -\frac{3}{2}\rhod (f\wedge\fb) +\lot,\\
 \rhod' &=& \rhod  -\fb\c\dual\b - f\c\dual\bb +\frac{3}{2}\rhod(f\c\fb) +\frac{3}{2}\rho (f\wedge\fb) +\lot
\eeaa
Differentiating the two equations w.r.t. $e_a'$, and using the decomposition of $e_a'$, we infer
\beaa
e_a'(\rho') &=& \left(\left(\de_a^b +\frac{1}{2}\fb_af^b\right) e_b +\frac 1 2  \fb_a  e_4 +\left(\frac 1 2 f_a +\frac{1}{8}|f|^2\fb_a\right)   e_3\right)\rho\\
&& + e_a'\left(\fb\c\b - f\c\bb +\frac{3}{2}\rho(f\c\fb) -\frac{3}{2}\rhod (f\wedge\fb) +\lot\right),\\
e_a'(\rhod') &=& \left(\left(\de_a^b +\frac{1}{2}\fb_af^b\right) e_b +\frac 1 2  \fb_a  e_4 +\left(\frac 1 2 f_a +\frac{1}{8}|f|^2\fb_a\right)   e_3\right)\rhod\\
&&  +e_a'\left(-\fb\c\dual\b - f\c\dual\bb +\frac{3}{2}\rhod(f\c\fb) +\frac{3}{2}\rho (f\wedge\fb) +\lot\right),
\eeaa
and hence
\beaa
e_a'(\rho') &=& e_a(\rho)+\frac 1 2  \fb_a  e_4(\rho)+\frac 1 2 f_a  e_3(\rho)\\
&& + e_a'\left(\fb\c\b - f\c\bb +\frac{3}{2}\rho(f\c\fb) -\frac{3}{2}\rhod (f\wedge\fb)\right)+\lot,\\
e_a'(\rhod') &=& e_a(\rho) +\frac 1 2  \fb_a  e_4(\rho) + \frac 1 2 f_a e_3(\rhod)\\
&&  +e_a'\left(-\fb\c\dual\b - f\c\dual\bb +\frac{3}{2}\rhod(f\c\fb) +\frac{3}{2}\rho (f\wedge\fb)\right)+\lot,
\eeaa
and hence
\beaa
\nab'(\rho') &=& \nab(\rho)+\frac 1 2  \fb  e_4(\rho)+\frac 1 2 f  e_3(\rho)\\
&& + \nab'\left(\fb\c\b - f\c\bb +\frac{3}{2}\rho(f\c\fb) -\frac{3}{2}\rhod (f\wedge\fb)\right)+\lot,\\
\dual\nab'(\rhod') &=& \dual\nab(\rho) +\frac 1 2  \dual\fb  e_4(\rho) + \frac 1 2\dual f e_3(\rhod)\\
&&  +\dual\nab'\left(-\fb\c\dual\b - f\c\dual\bb +\frac{3}{2}\rhod(f\c\fb) +\frac{3}{2}\rho (f\wedge\fb)\right)+\lot,
\eeaa
Together with the estimate \eqref{eq:usefulesitateinthediscussionofThmM0:00:first} for $f$ and $\fb$ to estimate the linear terms $\fb e_4(\rho)$, $\dual\fb e_4(\rho)$, $f e_3(\rho)$ and $\dual f e_3(\rho)$,  and the bootstrap  estimate \eqref{eq:iterativeapssutiononS1fortheproofofTheoremM0} for $(f,\fb, \la)$ to estimate the other terms, we deduce\footnote{We also make use of a  standard  elliptic estimate on $S_1$ and the $r$  dominance condition $r_*\sim   \ep_0^{-1} $, see \eqref{eq:behaviorofronS-star},   on $\Si_*$.},
\bea
\max_{k\leq k_{large}-1}r^3\|{\dkb'}^k\dds_1'(-\rho', \rhod')\|_{L^2(S_{1}')} &\les&\ep_0. 
\eea

{\bf Step 9.} Recall the definition of the mass aspect function $\mu'$
\beaa
\mu' &=& -\div\,'\ze' -\rho' +\frac{1}{2}\chih'\c\chibh'.
\eeaa
and the the null structure equation
\beaa
\curl'\ze' &=& \rhod' -\frac{1}{2}\chih'\wedge\chibh'.
\eeaa
Together with the GCM conditions for $\mu$ on $\Si_*$, this yields on $S_{1}'$, recalling that $S_{1}'\subset\Si_*$,
\beaa
\ddd_1'\ze' &=& (-\mu', 0)+(-\rho', \rhod')+\frac{1}{2}\left(\chih'\c\chibh', -\chih'\wedge\chibh'\right)\\
&=& -\left(M_{0}'+\sum_pM_{p}'{{J'}^{(p)}}, 0\right)+(-\rho', \rhod')+\frac{1}{2}\left(\chih'\c\chibh', -\chih'\wedge\chibh'\right)
\eeaa
and hence, since $M_{0}'$ and $M_{p}'$ are constant on $S_{1}'$, we infer
\beaa
\dds_2'\dds_1'\ddd_1'\ze' &=& -\sum_pM_{p}'\dds_2'\dds_1'({{J'}^{(p)}}, 0)+\dds_2'\dds_1'(-\rho', \rhod')+\frac{1}{2}\dds_2'\dds_1'\left(\chih'\c\chibh', -\chih'\wedge\chibh'\right).
\eeaa
In view of the identity $\dds_1\,'\ddd_1\,'=\ddd_2\,'\dds_2\,'+2K'$, we infer
\beaa
(\dds_2\,'\ddd_2\,'+2K')\dds_2\,'\ze' &=& \dds_2'\dds_1'(-\rho', \rhod') -\sum_pM_{p}'\dds_2'\dds_1'({{J'}^{(p)}}, 0)\\
&&+\frac{1}{2}\dds_2'\dds_1'\left(\chih'\c\chibh', -\chih'\wedge\chibh'\right)+\nab'(K')\hot\ze'.
\eeaa
Using  the estimate for $(\rho', \rhod')$ of Step 8, the fact that $M_{p}'\in r^{-1}\Ga_g$ in view of Corollary \ref{cor:Cb0CbpM0MareGagandrm1Gag}, the control of $\dds_2'\dds_1'({{J'}^{(p)}}, 0)$ provided by Lemma \ref{Le:Si*-ell=1modes:proofofThmM0}, and   an elliptic estimate for $\dds_2\,'\ddd_2\,'+2K'$, we obtain 
\bea
\max_{k\leq k_{large}-1}r^2\|{\dkb'}^k\dds_2\,'\ze'\|_{L^2(S_{1}')} &\les&\ep_0. 
\eea 
Note that the quadratic terms involving $\chih'\c\chibh'$, $\chih'\wedge\chibh'$ and $\nab'(K')\hot\ze'$ are  estimated using the transformation formulas\footnote{In fact, in view of the Gauss equation $K'=-\rho'-\frac{1}{4}\ka'\kab'+\frac{1}{2}\chih'\chibh'$, the GCM conditions for $\ka'$ and $\kab'$, and the control of $\rho'$ in Step 8, we only need the transformation formulas for $\chih'$ and $\chibh'$. These formulas involve at most one angular derivative of $f$ and $\fb$, and no transversal derivative.}, the estimates  \eqref{eq:iterativeapssutiononS1fortheproofofTheoremM0}  for $(f,\fb, \la)$, and   the control provided by Proposition \ref{proposition:geodesicfoliationLextt}  for the curvature components and the Ricci coefficients of the  part $\Lextt$ of the initial data layer.

{\bf Step 10.} Recall Codazzi for $\chih'$
\beaa
\ddd_2\,'\chih' +\ze'\c\chih' &=& \frac{1}{2}\nab'\ka'+\frac{1}{2}\ka'\ze' - \b'.
\eeaa
We differentiate w.r.t. $\dds_2\,'$ and use the GCM condition $\ka'=2/r'$ which holds on $\Si_*$ and $S_{1}'\subset\Si_*$ to deduce 
\beaa
\dds_2\,'\ddd_2\,'\chih' &=& -\dds_2\,'\b'  +\frac{1}{r'}\dds_2\,'\ze'  - \dds_2\,'(\ze'\c\chih').
\eeaa
Together with the estimate of Step 8 for $\b'$, the estimate of Step 9 for $\dds_2\,'\ze'$,  dealing with the quadratic terms as above, and using an elliptic estimate, we infer, 
\beaa
\max_{k\leq k_{large}}r\|{\dkb'}^k\chih'\|_{L^2(S_{1}')} &\les&\ep_0. 
\eeaa

Next, recall from Proposition \ref{Proposition:transformationRicci} the transformation formula
\beaa
\chih' &=& \la\Bigg(\chih  +  \nab'\hot f + f\hot\eta + f\hot\ze +\frac{1}{4}\fb\hot f\ka    -\omb f\hot f \\
&&  +\frac 1 4  (f\hot\fb) \la^{-1}\ka'  +\frac 1 2  \fb\hot (f\c\la^{-1}\chih')+\lot\Bigg)
\eeaa
where we have used the fact that   $\atrch'=\atrch=0$ and $\atrchb=\atrchb'=0$ since both frame are outgoing geodesic and hence integrable, and also that $\om=0$ and $\xi=0$. Together with the above estimate for $\chih'$,   the estimate \eqref{eq:iterativeapssutiononS1fortheproofofTheoremM0} for $(f,\fb, \la)$, and the estimates  \eqref{def:initialdatalayerassumptions} for the Ricci coefficients of  $\Lextt\subset\Lext $, we infer
\beaa
\max_{k\leq k_{large}}r\|{\dkb'}^k\dds_2\,' f\|_{L^2(S_{1}')} &\les& \ep_0+\ep^2\les\ep_0.
\eeaa
Together with  the  Hodge elliptic  estimates of  Lemma \ref{prop:2D-Hodge2}, chapter 5, we infer
\bea
\max_{k\leq k_{large}+1}\|{\dkb'}^kf\|_{L^2(S_{1}')} &\les& \ep_0+r^2 \big| (\ddd_1'f)_{\ell=1}\big|.
\eea

{\bf Step 11.} Next, recall from Proposition \ref{Proposition:transformationRicci} the transformation formula
\beaa
\la^{-1}\atrch' &=& \atrch  +  \curl'f + f\wedge\eta + f\wedge\ze + \fb\wedge\xi+\frac{1}{4}\left(\fb\wedge f\trch +(f\c\fb)\atrch\right)\\
&& +\om f\wedge\fb  -\frac{1}{4}|f|^2\atrchb -  \frac 1 4 ( f\c\fb) \la^{-1}\atrch' +\frac 1 4   \la^{-1}(f\wedge \fb)\trch'+\lot
\eeaa
Since   $\atrch'=\atrch=0$ and $\atrchb=\atrchb'=0$ as both frame are outgoing geodesic and hence integrable, and also since $\om=0$ and $\xi=0$, we infer 
\beaa
\curl'f &=& - f\wedge\eta - f\wedge\ze - \frac{1}{4}\fb\wedge f\ka     - \frac 1 4   \la^{-1}(f\wedge \fb)\ka'+\lot
\eeaa 
Together with  the estimate \eqref{eq:iterativeapssutiononS1fortheproofofTheoremM0} for $(f,\fb, \la)$, and  our   estimates  for the Ricci coefficients  of $\Lextt$, we deduce
\beaa
\max_{k\leq k_{large}+7}r\|{\dkb'}^k\curl' f\|_{L^2(S_{1}')} &\les& \ep_0.
\eeaa
Together with the estimate for $f$ of Step 10, and since $\ddd_1'=(\div',\curl')$, we infer
\bea
\lab{eqEstimateStep11-8}
\max_{k\leq k_{large}+8}\|{\dkb'}^kf\|_{L^2(S_{1}')} &\les& \ep_0+r^2 \big| (\div' f)_{\ell=1}\big|.
\eea

{\bf Step 12.} In view of Step 11,  it remains  to control $(\div' f)_{\ell=1}$.  We begin by making the following local   bootstrap assumptions  
\bea\lab{eq:localbootstrapssutionfortheproofofTheoremM0:Ricci}
\sup_{\CC_1'}\left(r^2|\dkb'^{\leq 5}(\kac', \chih', \ze')|+r^2|\dkb'^{\leq 4}\ze'|+r^3|\dkb'^{\leq 4}\b'|+\left|\frac{r'}{r}-1\right|\right) \leq \ep,
\eea
where $\kac'=\ka'-\frac{2}{r'}$, and where we recall that $\CC_1'$ denotes the portion of the past directed outgoing null cone initialized on the sphere $S_1'$ and restricted to $r'\geq \de_*\ep_0^{-1}$. Recall also that $\CC_1'\subset\Lextt$ and that $r$ denotes the area radius for the outgoing geodesic foliation  of $\Lextt$ while $r'$ denotes the area radius of the spheres $S'\subset \CC_1'$. The  local bootstrap assumptions   \eqref{eq:localbootstrapssutionfortheproofofTheoremM0:Ricci} will be improved in Step 14. The goal of this step is  to prove the following.
\begin{lemma}
\lab{lemma:estimatefodiv'f_{ell=1}}
The following estimate holds true
\bea
\sup_{\CC_1'} r^2\left|(\div'(f))_{\ell=1}\right| &\les& \ep_0.
\eea
\end{lemma}

\begin{proof}
We proceed in four steps.

{\bf Step 12a.} We start by deriving an estimate for $\nab'\log\la$ and $\fb$. In view of Corollary \ref{cor:transportequationine4forchangeofframecoeff:simplecasefirst},  since $\atrch=0$ and $\xi=\om=0$, we have 
\beaa
\la^{-1}e_4'(\log\la) &=& 2f\c\ze +E_2(f, \Ga)
\eeaa
and hence
\beaa
\nab_4'(r'\nab'\log\la) &=& r'\nab'\left(2\la f\c\ze +\la E_2(f, \Ga)\right) +[\nab_4', r'\nab'](\log\la)\\
&=& r'\nab'\left(2\la f\c\ze +\la E_2(f, \Ga)\right) +r'\left(-\frac{1}{2}(\kac'-\ov{\kac'})\nab' -\chih'\c\nab'\right)\log\la.
\eeaa
Together with the bootstrap assumptions \eqref{eq:localbootstrapssutionfortheproofofTheoremM0:Ricci} for $\kac'$ and $\chih'$, the bootstrap assumption \eqref{eq:localbootstrapssutionfortheproofofTheoremM0} for $(f,\la)$ along $\CC_1'$,  and the  estimates \eqref{def:initialdatalayerassumptions} for the Ricci coefficients of the part $\Lextt$ of the initial data layer, we infer along $\CC_1'$ 
\beaa
|\nab_4'(r'\nab'\log\la)| &\les& \frac{\ep_0+\ep^2}{{r'}^2}\les\frac{\ep_0}{{r'}^2}.
\eeaa
Integrating from $S_1'$ where $\nab'\log\la$ verifies \eqref{eq:usefulesitateinthediscussionofThmM0:00:first}, we deduce along $\CC_1'$
\beaa
r'|\nab'\la| &\les& \ep_0.
\eeaa
Also, using the transformation formula for $\ze'$,  we derive 
\beaa
|\fb| &\les& r\Big(|\nab'\la|+|\ze'|+|\chih'|+|\kac'|+|\Ga_g|\Big) +|{\dkb'}^{\leq 1}f|+\lot,
\eeaa
and hence, together with the above estimate for $\nab'\la$, the bootstrap assumptions \eqref{eq:localbootstrapssutionfortheproofofTheoremM0:Ricci} for $\ze'$, $\kac'$ and $\chih'$, the bootstrap assumption \eqref{eq:localbootstrapssutionfortheproofofTheoremM0} for $f$,  and the properties of the background foliation of $\Lextt$, we infer along $\CC_1'$, recalling that $r'\gtrsim \ep_0^{-1}$ along $\CC_1'$,  
\beaa
|\fb| &\les& \ep_0+\frac{\ep}{r}\les \ep_0
\eeaa
and hence, we have obtained along $\CC_1'$ 
\bea\lab{eq:improvedintermediarycontroloflaandfb}
r'|\nab'\la|+|\fb| &\les& \ep_0.
\eea

{\bf Step 12b.} We  derive  next an estimate for $(\div'\b')_{\ell=1}$ along $\CC_1'$.

 We start with the following identity for the outgoing geodesic foliation initialized on $S_1'$
\beaa
e_4'\left({r'}^3\int_{S'}\div'\b' {{J'}^{(p)}}\right) &=& {r'}^3\int_{S'}\Big(e_4'\div'\b' +\ka'\div'\b' \Big){{J'}^{(p)}}+3e_4'(r'){r'}^2\int_{S'}\div'\b' {{J'}^{(p)}}\\
&=& {r'}^3\int_{S'}\Big(\div'\nab_4'\b'+[\nab_4',\div']\b' +\ka'\div'\b' \Big){{J'}^{(p)}}\\
&&+\frac{3}{2}\ov{\ka'}{r'}^3\int_{S'}\div'\b' {{J'}^{(p)}}.
\eeaa
Using the Bianchi identity for $\nab_4'\b'$ and the structure of the commutator, we infer
\beaa
e_4'\left({r'}^3\int_{S'}\div'\b' {{J'}^{(p)}}\right) &=& {r'}^3\int_{S'}\Big(\div'\div'\a' -2\nab'\ka'\c\b'+\div'(\a'\c\ze')  -\frac{1}{2}\ka'\ze'\c\b'\\
&&+|\b'|^2 -\chih'\c\nab'\b'+\ze'\c\chih'\c\b'  \Big){{J'}^{(p)}}\\
&& -\frac{3}{2}(\kac'-\ov{\kac'}){r'}^3\int_{S'}\div'\b' {{J'}^{(p)}},
\eeaa
which yields, after integration by parts,
\beaa
e_4'\left({r'}^3\int_{S'}\div'\b' {{J'}^{(p)}}\right) &=& {r'}^3\int_{S'}\a' \dds_2\dds_1({{J'}^{(p)}},0)+{r'}^3\int_{S'}\Big(-2\nab'\ka'\c\b'+\div'(\a'\c\ze') \\
&& -\frac{1}{2}\ka'\ze'\c\b' +|\b'|^2 -\chih'\c\nab'\b'+\ze'\c\chih'\c\b'  \Big){{J'}^{(p)}}\\
&& -\frac{3}{2}(\kac'-\ov{\kac'}){r'}^3\int_{S'}\div'\b' {{J'}^{(p)}}.
\eeaa
Together with the bootstrap assumptions \eqref{eq:localbootstrapssutionfortheproofofTheoremM0:Ricci} and the control of $\dds_2\dds_1({{J'}^{(p)}},0)$ provided by Lemma \ref{Le:Si*-ell=1modes:proofofThmM0}, we obtain along $\CC_1'$,  
\beaa
\left|e_4'\left({r'}^3\int_{S'}\div'\b' {{J'}^{(p)}}\right)\right| &\les& \frac{\ep^2}{{r'}^\frac{3}{2}}\les \frac{\ep_0}{{r'}^\frac{3}{2}}.
\eeaa
Transporting along $\CC_1'$ from $S_{1}'$, and using  the control of $(\div'\b')_{\ell=1}$  in \eqref{eq:usefulesitateinthediscussionofThmM0:00} on $S_{1}'$,  we infer
\beaa
\sup_{\CC_1'}{r'}^5\left|(\div'\b')_{\ell=1}\right| &\les& \ep_0.
\eeaa
In particular, consider the sphere $S'(\de_*\ep_0^{-1})=\CC_1'\cap\{r'=\de_*\ep_0^{-1}\}$. Then
\bea
\lab{estimate:div'b'_{ell=1}-CC'}
{r'}^5\left|(\div'\b')_{\ell=1}\right| &\les& \ep_0\quad \textrm{on}\quad S'(\de_*\ep_0^{-1}).
\eea

{\bf Step 12c.} We next use  estimate \eqref{estimate:div'b'_{ell=1}-CC'} to derive  an estimate for  $(\div'(f))_{\ell=1}$  on $S'(\de_*\ep_0^{-1})$, see 
\eqref{estimate:div' f_{ell=1}-CC'}.

To this end, we invoke  again   the transformation formula
\beaa
\b'&=& \la\left(\b +\frac 3 2\big(  f \rho+\dual  f  \rhod\big)+ \frac 1 2 \a\c\fb+\lot\right)
\eeaa
 from which  we derive
\beaa
\div'\b'&=& \la\left(\div'\b +\frac 3 2\div'\big(  f \rho+\dual  f  \rhod\big)+ \frac 1 2 \div'(\a\c\fb)+\lot\right)\\
&&+\nab'\la\c\left(\b +\frac 3 2\big(  f \rho+\dual  f  \rhod\big)+ \frac 1 2 \a\c\fb+\lot\right)\\
&=&  \div\b +\frac{1}{2}f\c\nab_3\b+\frac{1}{2}\fb\c\nab_4\b -\frac{3m}{{r'}^3}\div'(f)+3m\left(\frac{1}{r^3}-\frac{1}{{r'}^3}\right)\div'(f)\\
&&+\frac 3 2\left(  \div'(f) \left(\rho+\frac{2m}{r^3}\right)+\curl'(f)  \rhod\right)\\
&&+\frac 3 2\left(  f\c\nab\rho+\frac{1}{2}f\c(f\nab_3+\fb\nab_4)\rho+\dual  f\c\nab\rhod+\frac{1}{2}\dual f\c(f\nab_3+\fb\nab_4)\rhod\right)\\
&&+ \frac 1 2 \div(\a\c\fb) +(\la-1)\left(\div\b +\frac 3 2\big(  \div'(f) \rho+\curl'(f)  \rhod\big)\right)\\
&& +\nab'\la\c\left(\b +\frac 3 2\big(  f \rho+\dual  f  \rhod\big)\right)+\lot
\eeaa
Together with   the bootstrap assumptions  \eqref{eq:localbootstrapssutionfortheproofofTheoremM0} for $(f, \fb, \la)$,  the bootstrap assumption \eqref{eq:localbootstrapssutionfortheproofofTheoremM0:Ricci} for       $r-r'$, and the control provided by Proposition \ref{proposition:geodesicfoliationLextt} for $\Lextt$, we infer
\beaa
r^5\left|\div'\b'  +\frac{3m}{{r'}^3}\div'(f) -\div\b -\frac{1}{2}f\c\nab_3\b-\frac{1}{2}\fb\c\nab_4\b-\nab'\la\c\b\right| &\les& r^{\frac{1}{2}}\ep\ep_0+\ep^2\\
&\les& r^{\frac{1}{2}}\ep_0^{\frac{5}{3}}+\ep_0
\eeaa
where we used the fact that $\ep=\ep_0^{\frac{2}{3}}$. Also, using the Bianchi identities\footnote{Concerning the term $\fb\c\nab_4\b$, note that Bianchi identities imply $\nab_4\b=-2\trch\beta+\lot$, thus, since $\b=\frac{3a_0m_0}{r^4}f_0+  O(  \ep_0  r^{-7/2})$, we obtain $\fb\c\nab_4\b=-\frac{12a_0m_0}{r^5}f_0\c\fb+O(r^{-\frac{9}{2}}\ep_0)\fb=O(r^{-5}+r^{-\frac{9}{2}}\ep_0)\fb$ so that  $\fb\c\nab_4\b=O(r^{-5}\ep_0+r^{-\frac{9}{2}}\ep_0^2)$ thanks to the estimate \eqref{eq:improvedintermediarycontroloflaandfb} for $\fb$.}  for $\nab_4\b$ and $\nab_3\b$, the control of $(\fb, \la)$ provided by \eqref{eq:improvedintermediarycontroloflaandfb},  we obtain along $\CC_1'$ 
\beaa
r^5\left|\div'\b'  +\frac{3m}{{r'}^3}\div'(f) -\div\b\right| &\les& r^{\frac{1}{2}}\ep_0^{\frac{5}{3}}+\ep_0+r^{-1}.
\eeaa
Thus, on the sphere $S'(\de_*\ep_0^{-1}) $, where $r'=\de_*\ep_0^{-1}$, we infer, see also  Remark \ref{remark:needforstrongerinitialdatalayernorminproofofThmM0:1} below, 
\beaa
\sup_{S'(\de_*\ep_0^{-1})}    r^5\left|\div'\b' +\frac{3m}{{r'}^3}\div'(f)\right| &\les& \ep_0.
\eeaa

\begin{remark}\lab{remark:needforstrongerinitialdatalayernorminproofofThmM0:1}
In the estimate above, we used in particular the following estimate of Proposition \ref{proposition:geodesicfoliationLextt} for the control of $\div\b$ which relies on the stronger bound on   the initial data layer norm $\, ^{(ext)}\Ik_{3}\leq \ep_0^2$, i.e.
\beaa
\sup_{\Lextt\cap\{r\sim \ep_0^{-1}\}}\Big(r^5|\div\b|\Big) \les \ep_0.
\eeaa
\end{remark}

This yields on $S'(\de_*\ep_0^{-1})$
\beaa
mr^2|(\div'(f))_{\ell=1}| &\les& \ep_0+r^5|(\div'\b')_{\ell=1}|.
\eeaa
Since  $|m-m_0|\leq\ep$ and\footnote{By the local  bootstrap assumption \eqref{eq:localbootstrapssutionfortheproofofTheoremM0:Ricci}.} $|r-r'|\leq \ep r$, together with  the  estimate \eqref{estimate:div'b'_{ell=1}-CC'}  for $r^5|(\div'\b')_{\ell=1}|$ on $S'(\de_*\ep_0^{-1})$, we infer on $S'(\de_*\ep_0^{-1})$
\beaa
m_0r^2|(\div'(f))_{\ell=1}| &\les& \ep_0
\eeaa
and hence
\bea
\lab{estimate:div' f_{ell=1}-CC'}
r^2|(\div'(f))_{\ell=1}| &\les& \ep_0\quad \textrm{on}\quad S'(\de_*\ep_0^{-1}).
\eea

{\bf Step 12d.} We  will  next  propagate forward   the information provided by estimate \eqref{estimate:div' f_{ell=1}-CC'}.

We use the following identity for the outgoing geodesic foliation initialized on $S_1'$
\beaa
e_4'\left(\int_{S'}\div'(f) {{J'}^{(p)}}\right) &=& \int_{S'}\Big(e_4'\div'(f) +\ka'\div'(f) \Big){{J'}^{(p)}}\\
&=& \int_{S'}\Big(\div'\nab_4'(f)+[\nab_4',\div']f +\ka'\div'(f) \Big){{J'}^{(p)}}.
\eeaa
Also recall, see  Corollary  \ref{cor:transportequationine4forchangeofframecoeff:simplecasefirst},  that $f$ satisfies the  transport equation\footnote{Recall that $\atrch=0$ and $\om=\xi=0$ for the outgoing geodesic foliation of  $\Lextt$. } along $\CC_1'$
\beaa
\nab_{\la^{-1}e_4'}f+\frac{1}{2}\ka f  &=& - f\c\chih+E_1(f, \Ga).
\eeaa
Plugging in the above, and using a commutator formula for $[\nab_4',\div']$, we infer
\beaa
&&e_4'\left(\int_{S'}\div'(f) {{J'}^{(p)}}\right)\\
 &=& \int_{S'}\Bigg[\div'\left(\la\left(-\frac{1}{2}\ka f   - f\c\chih+E_1(f, \Ga)\right)\right)\\
&& +\left(-\frac{1}{2}\ka'\div' -\frac{1}{2}\ka'\ze'\c +\dual\b'\c\dual -\chih'\c\nab' +\ze'\c\chih'\c\right)f +\ka'\div'(f) \Bigg]{{J'}^{(p)}}\\
&=&  \int_{S'}\Bigg[-\frac{1}{2}\nab'\ka'\c f -\frac{1}{2}\div'\left((\la\ka-\ka')f\right)   +\div'\left(\la\left(   - f\c\chih+E_1(f, \Ga)\right)\right)\\
&& +\left( -\frac{1}{2}\ka'\ze'\c +\dual\b'\c\dual -\chih'\c\nab' +\ze'\c\chih'\c\right)f  \Bigg]{{J'}^{(p)}}.
\eeaa
In view of the bootstrap assumption \eqref{eq:localbootstrapssutionfortheproofofTheoremM0} for $f$ and $\la$,   and the bootstrap assumptions \eqref{eq:localbootstrapssutionfortheproofofTheoremM0:Ricci}, we deduce
\beaa
\left|e_4'\left(\int_{S'}\div'(f) {{J'}^{(p)}}\right)\right| &\les& \frac{\ep_0+\ep^2}{r^2}+\ep\left|{\dkb'}^{\leq 1}\left(\ka-\la^{-1}\ka'\right)\right|+r\left|{\dkb'}^{\leq 1}\left(E_1(f, \Ga)\right)\right|.
\eeaa
In view of  the form of the error term  $E_1$ in Corollary \ref{cor:transportequationine4forchangeofframecoeff:simplecasefirst}, using  the transformation formula for $\ka'$, together with the bootstrap assumption \eqref{eq:localbootstrapssutionfortheproofofTheoremM0} for $f$, and the bootstrap assumptions  \eqref{eq:localbootstrapssutionfortheproofofTheoremM0:Ricci}, we obtain
\beaa
\left|e_4'\left(\int_{S'}\div'(f) {{J'}^{(p)}}\right)\right| &\les& \frac{\ep_0+\ep^2}{r^2}\les\frac{\ep_0}{r^2}.
\eeaa
Integrating forward from $r=\ep_0^{-1}$, and using  estimate  \eqref{estimate:div' f_{ell=1}-CC'} for $(\div'(f))_{\ell=1}$  on $S'(\de_*\ep_0^{-1})$, we obtain
\beaa
\sup_{\CC_1'} r^2\left|(\div'(f))_{\ell=1}\right| &\les& \ep_0.
\eeaa
This ends the proof of  Lemma  \ref{lemma:estimatefodiv'f_{ell=1}}. 
\end{proof}

\noindent{\bf Step 13.} Combining  the estimate  \eqref{eqEstimateStep11-8} for ${\dkb'}^kf$ of Step 11 and the  estimate for $(\div'(f))_{\ell=1}$ of Step 12, we obtain on $S_1=S_1'$,
\beaa
\max_{k\leq k_{large}+1}\|{\dkb'}^kf\|_{L^2(S_{1}')} &\les& \ep_0.
\eeaa
Together with \eqref{eq:usefulesitateinthediscussionofThmM0:00:first}, we infer
\beaa
\|f\|_{\hk_{k_{large}+1}(S_{1}')}+r^{-1}\|(\fb,  \la-\ov{\la}^{S_{1}'})\|_{\hk_{k_{large}+1}(S_{1}')} &\les& \ep_0.
\eeaa
In particular, the above estimate for $(f, \fb)$ allows to  reapply        Lemma 7.3 in \cite{KS-GCM1} (restated here in  Lemma    \ref{lemma:consequencedeformationsurfaceusedinTheoremM0}),  with $\de_1=\ep_0$ which yields
\beaa
\sup_{S_{1}'}\left|\frac{r'}{r}-1\right| &\les& \ep_0.
\eeaa
Together with \eqref{eq:usefulesitateinthediscussionofThmM0:00:first:ovla}, we infer
\beaa
\|f\|_{\hk_{k_{large}+1}(S_{1}')}+r^{-1}\|(\fb,  \log\la)\|_{\hk_{k_{large}+1}(S_{1}')} &\les& \ep_0.
\eeaa
This  improves  the iteration assumption \eqref{eq:iterativeapssutiononS1fortheproofofTheoremM0}.  

We next appeal to Corollary 4.2 in \cite{Shen} (restated here in Corollary \ref{cor:constuctionGCMH}) with $\dg=\ep_0$, with background foliation 
 given by $\Lextt$  and $s_{max}=k_{large}$ which allows us to make use of  the above estimate for $(f, \fb, \la)$ on $S_1'\subset\Si_* $ to derive
 \beaa
 \sup_{k\le k_{ large}+1} \Big( \|{\dk}^k f \|_{L^2(S_{1}')} +r^{-1} \|{\dk}^{k} (\fb, \log \la) \|_{L^2(S_{1}')} +\|\dk^{\le k-1}\nab_\nu' (\fb, \log \la) \|_{L^2(S_{1}')}\Big) &\les& \ep_0.
 \eeaa

The above control of $(f, \fb)$, together with  Lemma 7.3 in \cite{KS-GCM1} (restated here in Lemma     \ref{lemma:consequencedeformationsurfaceusedinTheoremM0})  for $\de_1=\ep_0$, and  
Corollary 7.7 in \cite{KS-GCM1}  (restated here in Lemma     \ref{lemma:consequencedeformationsurfaceusedinTheoremM0}) with $\epg=\ep_0$, implies
\beaa
 \sup_{S_{1}'}\left(\left|\frac{m_H'}{m_0}-1\right|+\left|\frac{r'}{r}-1\right| \right)\les \ep_0,
\eeaa
where $m_H'$ denotes the Hawking mass of $S'$. 

We appeal next to the argument used to derive  estimate  $m_H'-m$  in Proposition \ref{prop:controlofell=0modesonSigmastar}  which leads to 
\beaa
\sup_{\Si_*}u^{1+2\dec}|m_H'-m| &\les& \ep_0.
\eeaa
We have thus obtained on $S_{1}'$
\bea\lab{eq:finalesitateforffbchecklaonS1inproofofThmM0:reallyfinal}
\bsplit
\sup_{k\le k_{ large}+1} \Big( \|{\dk}^k f \|_{L^2(S_{1}')} +r^{-1} \|{\dk}^{k} (\fb, \log \la) \|_{L^2(S_{1}')} +\|\dk^{\le k-1}\nab_\nu' (\fb, \log \la) \|_{L^2(S_{1}')}\Big) &\les \ep_0.\\
\sup_{S_{1}'}\left(\left|\frac{m}{m_0}-1\right|+\left|\frac{r'}{r}-1\right| \right)  &\les \ep_0.
\end{split}
\eea 
In view of \eqref{eq:finalesitateforffbchecklaonS1inproofofThmM0:reallyfinal} and the transformation formulas from  the frame\footnote{Recall that $S_1'=S_1$ and  the  primed frame  on  $S_1'$ coincides with that induced by $\Si_*$.} of $\Lextt$ to that
of of $\Lextt$ we deduce, in addition to $\kac'=0$, $\xi'=0$, $\om'=0$ and $\etab'=\ze'$,
\bea\lab{eq:finalesitateforRiccicoeffSigmastaronS1inproofofThmM0:reallyfinal}
 \nn r\|{\dk}^{\le k_{large}}\chih'\|_{L^2(S_{1}')}+ r\|{\dk}^{\le k_{large}-1}\ze'\|_{L^2(S_{1}')} \\
 +\|{\dk}^{\le k_{large}}(\eta', \trchbc', \chibh', \ombc', \xib')\|_{L^2(S_{1}')} &\les& \ep_0.
\eea

\begin{remark}
The estimate for $\trchbc'$ in \eqref{eq:finalesitateforRiccicoeffSigmastaronS1inproofofThmM0:reallyfinal} displays a loss of $r^{-1}$ with respect to the expected behavior. It is due to the anomalous behavior for $\fb$ and $\la$ in \eqref{eq:finalesitateforffbchecklaonS1inproofofThmM0:reallyfinal}. A priori, the same  loss should also occur for $\ze'$ which would then create problems in Step 22. To avoid this issue, we estimate first $\b'$ using the transformation formulas and the control of $(f, \fb, \la)$ and $r'-r$ in \eqref{eq:finalesitateforffbchecklaonS1inproofofThmM0:reallyfinal} to  obtain
\beaa
r^2\|{\dk}^{\le k_{large}}\b'\|_{L^2(S_{1}')} &\les& \ep_0.
\eeaa
We can then use the Codazzi equation for $\chih'$ to estimate $\ze'$, Indeed 
together with the estimate for $\chih'$ in \eqref{eq:finalesitateforRiccicoeffSigmastaronS1inproofofThmM0:reallyfinal} and the fact that $\kac'=0$, we obtain the claimed estimate for $\ze'$ in \eqref{eq:finalesitateforRiccicoeffSigmastaronS1inproofofThmM0:reallyfinal}, where the loss of one derivative is due to the term $\div'\chih'$ in Codazzi.
\end{remark}

\noindent{\bf Step 14.} In this step, we  improve the bootstrap assumption \eqref{eq:localbootstrapssutionfortheproofofTheoremM0} for $(f, \fb, \la)$ and the bootstrap assumptions \eqref{eq:localbootstrapssutionfortheproofofTheoremM0:Ricci} on $(\kac', \chih', \ze')$, $\b'$ and $r'-r$. In view of Corollary \ref{cor:transportequationine4forchangeofframecoeff:simplecasefirst}, and since $\atrch=0$ and $\xi=\om=0$, we have 
\beaa
\nab_{\la^{-1}e_4'}f+\frac{1}{2}\ka f  &=&  - f\c\chih+E_1(f, \Ga),\\
\la^{-1}e_4'(\log\la) &=& 2f\c\ze +E_2(f, \Ga).
\eeaa
Since 
\beaa
\la^{-1}e_4' = e_4+f^ae_a+\frac{1}{4}|f|^2e_3, \qquad e_4(r)=\frac{r}{2}\ov{\ka}, \qquad e_4(e_3(r))=-2\omb,
\eeaa
we infer
\beaa
\nab_{\la^{-1}e_4'}(rf)  &=& -\frac{r}{2}(\kac -\ov{\kac}) f - rf\c\chih+rE_1(f, \Ga)+\frac{1}{4}|f|^2e_3(r)f.
\eeaa

Then, we proceed as follows for the estimates of $(f, \fb, \la)$, $(\kac', \chih', \ze')$, $\b'$ and $r'-r$:
\begin{enumerate}
\item Integrating the above  transport equations for  $f$ and  $\la$ from $S_1'$ where \eqref{eq:finalesitateforffbchecklaonS1inproofofThmM0:reallyfinal} holds, we obtain on $\CC_1'$
\bea
r'|\dkb'^{\leq 5}f|+|\dkb'^{\leq 5}\log\la| &\les& \ep_0.
\eea

\item We estimate $(\kac', \chih', \ze')$ and $\b'$ as follows:
\begin{enumerate}
\item one first controls 5 derivatives of $\a'$ relying on the corresponding change of frame formula in Proposition \ref{Proposition:transformationRicci}, the control of the foliation of $\Lextt$, and using in particular the fact that the change of frame formula for $\a'$ only involves   $(f, \la)$ but not $\fb$,

\item one then controls 5 derivatives $\chih'$ using the null structure equation for $\nab_4'\chih'$, the above control for 5 derivatives of $\a'$,  and integrating from $S_1'$ where $\chih'$ is under control from \eqref{eq:finalesitateforRiccicoeffSigmastaronS1inproofofThmM0:reallyfinal},

\item one then also controls 5 derivatives of $\kac'$ by integrating Raychadhuri from $S_1'$ where $\kac'=0$ in view of the GCM condition on $\Si_*$, 

\item then, using the transformation formula for $\b'$,  the control of the foliation of $\Lextt$, and the above control of $f$ and $\la$, we control 4 derivatives $\b'$,

\item then, we control 4 derivatives of $\ze'$ from the codazzi for $\chih'$, thanks to the control of 4 derivatives of $\b'$, and of 5 derivatives of $\chih'$ and $\kac'$.
\end{enumerate}
The above steps thus lead to the following control on  $\CC_1'$
\bea
\lab{eq:EstimatesStep14}
{r'}^2|\dkb'^{\leq 5}(\chih', \kac')|+{r'}^2|\dkb'^{\leq 4}\ze'|+{r'}^3|\dkb'^{\leq 4}\b'|+{r'}^{\frac{7}{2}}|\dkb'^{\leq 5}\a'| &\les& \ep_0.
\eea

\item  Using the transformation formulas for $\ze'$ and the control of 5 derivatives of $f$ and $\la$, and of 4 derivatives of $\ze'$, we obtain the control of 4 derivatives of $\fb$ on $\CC_1'$:
\bea
|\dkb'^{\leq 4}\fb| &\les& \ep_0.
\eea

\item Finally, the above control of $(f, \fb)$, together with Lemma 7.3 in \cite{KS-GCM1}   (recalled here in   Lemma \ref{lemma:consequencedeformationsurfaceusedinTheoremM0})   for $\de_1=\ep_0$ implies the following control of $r'-r$ on $\CC_1'$: 
\bea
|r'-r| &\les& \ep_0r.
\eea
\end{enumerate}
In view of the above, we have improved the bootstrap assumption \eqref{eq:localbootstrapssutionfortheproofofTheoremM0} for $(f, \fb, \la)$ and the bootstrap assumptions \eqref{eq:localbootstrapssutionfortheproofofTheoremM0:Ricci} on $(\kac', \chih', \ze')$, $\b'$ and $r'-r$.  

Proceeding as above for $\b'$, and using in addition the above estimate for $r'-r$, we get the following improved estimates for $\b'$
\bea
|\dkb'^{\leq 4}\b'| &\les& \frac{\ep_0}{{r'}^\frac{7}{2}}+\frac{1}{{r'}^4}.
\eea

\noindent{\bf Step 15.}  In Steps 15--16, we   estimate\footnote{Recall that $a$ is defined  in   section \ref{sec:definitionofamthetandvphiadmissible}.}  $a-a_0$.  Proceeding as in Step 12b  we obtain
\beaa
e_4'\left({r'}^3\int_{S'}\curl'\b' {{J'}^{(p)}}\right) &=& {r'}^3\int_{S'}\a' \dds_2\dds_1(0,{{J'}^{(p)}})+{r'}^3\int_{S'}\Big(-2\dual\nab'\ka'\c\b'+\curl'(\a'\c\ze') \\
&& -\frac{1}{2}\ka'\ze'\c\dual\b'  -\dual\chih'\c\nab'\b'+\ze'\c\dual\chih'\c\b'  \Big){{J'}^{(p)}}\\
&& -\frac{3}{2}(\kac'-\ov{\kac'}){r'}^3\int_{S'}\curl'\b' {{J'}^{(p)}}.
\eeaa

Together with the control of $(\kac', \chih', \ze', \b', \a')$ from Step 14 and the control of $\dds_2\dds_1(0,{{J'}^{(p)}})$ provided by Lemma \ref{Le:Si*-ell=1modes:proofofThmM0}, we obtain along $\CC_1'$
\beaa
\left|e_4'\left({r'}^3\int_{S'}\curl'\b' {{J'}^{(p)}} - 8\pi ma\de_{p0}\right)\right| &\les&  \frac{\ep_0}{{r'}^\frac{3}{2}}.
\eeaa
Transporting along $\CC_1'$ from $S_{1}'$, using  the control of $(\curl'\b')_{\ell=1}$ in \eqref{eq:localbootassell1modedivbeta:proofThmM0:improved} on $\Si_*$, and hence on $S_{1}'$, we infer
\beaa
\sup_{\CC_1'}{r'}^5\left(|(\curl'\b')_{\ell=1,\pm}|+\left|(\curl'\b')_{\ell=1,0}-\frac{2am}{{r'}^5}\right|\right) &\les& \ep_0.
\eeaa
In particular, consider the sphere $S'(\de_*\ep_0^{-1})=\CC_1'\cap\{r'=\de_*\ep_0^{-1}\}$. Then
\beaa
{r'}^5\left(|(\curl'\b')_{\ell=1,\pm}|+\left|(\curl'\b')_{\ell=1,0}-\frac{2am}{{r'}^5}\right|\right) &\les& \ep_0\quad \textrm{on}\quad S'(\de_*\ep_0^{-1}).
\eeaa
Also, using the change of frame formula for $\b'$ in Proposition \ref{Proposition:transformationRicci}, the control for $(f, \fb, \la)$ of Step 14, and the control of the the curvature components of $\Lextt$, we have
\beaa
\left|\curl'\b'-\curl\b\right| &\les& \frac{\ep_0}{{r'}^5}+\frac{\ep_0^2}{{r'}^{\frac{9}{2}}}.
\eeaa
Together with the above, we obtain on $S'(\de_*\ep_0^{-1}) $
\beaa
{r'}^3\left(\left|\int_{}{J'}^{(+)}\curl\b\right|+\left|\int_{}{J'}^{(-)}\curl\b\right| +\left|\int_{S'(\de_*\ep_0^{-1})}{J'}^{(0)}\curl\b -\frac{8\pi ma}{{r'}^3}\right|\right) &\les& \ep_0.
\eeaa
Using the  estimates for $m-m_0$ of Step 13, we deduce 
\beaa
{r'}^3\left(\left|\int_{S'(\de_*\ep_0^{-1})}{J'}^{(+)}\curl\b\right|+\left|\int_{S'(\de_*\ep_0^{-1})}{J'}^{(-)}\curl\b\right| +\left|\int_{S'(\de_*\ep_0^{-1})}{J'}^{(0)}\curl\b -\frac{8\pi m_0a}{{r'}^3}\right|\right) &\les& \ep_0.
\eeaa
Also,  making use of the estimate for $\curl \b$ in Proposition \ref{proposition:geodesicfoliationLextt}
\beaa
\curl\b &=& \frac{6a_0m_0}{r^5}J^{(0)}+O(r^{-5}\ep_0) \quad \mbox{on}  \,\,  S'(\de_*\ep_0^{-1}).
\eeaa
Using the  estimates for $r'-r$ of Step 14, this yields on $S'(\de_*\ep_0^{-1})$
\beaa
\curl\b &=& \frac{6a_0m_0}{{r'}^5}J^{(0)}+O({r'}^{-5}\ep_0).
\eeaa
Plugging in the above, and dividing by $m_0$, we deduce
\beaa
{r'}^{-2}\left(\left|\int_{S'(\de_*\ep_0^{-1})}{J'}^{(+)}a_0J^{(0)}\right|+\left|\int_{S'(\de_*\ep_0^{-1})}{J'}^{(-)}a_0J^{(0)}\right| +\left|\int_{S'(\de_*\ep_0^{-1})}{J'}^{(0)}a_0J^{(0)} -\frac{4\pi a}{3}{r'}^2\right|\right) &\les& \ep_0.
\eeaa

Now, recall that we have either $a_0=0$ or $|a_0|\gg \ep_0$. In particular, we have 
\bea\lab{eq:controlofainthecasea0equal0ThM0}
|a| &\les& \ep_0\quad\textrm{if}\quad a_0=0.
\eea
In the other case, we have, since $|a_0|\gg \ep_0$,
\bea\lab{eq:notyetthecontrolofainthecasea0notequal0ThM0}
\nn {r'}^{-2}\Bigg(\left|\int_{S'(\de_*\ep_0^{-1})}{J'}^{(+)}J^{(0)}\right|+\left|\int_{S'(\de_*\ep_0^{-1})}{J'}^{(-)}J^{(0)}\right|\\ +\left|\int_{S'(\de_*\ep_0^{-1})}{J'}^{(0)}J^{(0)} -\frac{4\pi a}{3a_0}{r'}^2\right|\Bigg) &\les& \ep_0.
\eea

\noindent{\bf Step 16.} In this step, we consider the case $a_0\neq 0$.

{\bf Step 16a.}  We introduce the following scalar function 
\bea
\widetilde{J}:=J^{(0)}-(1+c_0){J'}^{(0)}, \qquad c_0:=\frac{\int_{S'(\de_*\ep_0^{-1})}{J'}^{(0)}J^{(0)}}{\int_{S'(\de_*\ep_0^{-1})}{({J'}^{(0)})^2}}-1.
\eea
We summarize the results of this step in the lemma below.
\begin{lemma}
\lab{Lemma:estmatesforwidetildeJandc_0}
The following estimates hold true
\bea
\bsplit
\|\widetilde{J}\|_{L^\infty(S'(\de_*\ep_0^{-1}))} &\les \ep_0,\qquad 
|c_0| \les \ep_0.
\end{split}
\eea
\end{lemma}

\begin{proof}
In view of the bootstrap assumption \eqref{eq:bootstrapassumptiondifferenceell=1basisofCC1andLextprime} for ${J'}^{(p)}-\Jp$, we have
\bea\lab{eq:localbootstrapassumptiononc0:proofThmM0}
|c_0|\les \ep.
\eea
In view of the definition of $\widetilde{J}$ and $c_0$, we have
\beaa
\int_{S'(\de_*\ep_0^{-1})}{J'}^{(0)}\widetilde{J} &=& 0,
\eeaa
i.e. $(\widetilde{J})_{\ell=1,0}=0$ on $S'(\de_*\ep_0^{-1})$. 
 Using  an  elliptic estimate  for  $\De'+\frac{2}{{r'}^2}$,   see section \ref{section:Hodge-elliptic systemsChapter5},   we deduce
\beaa
r^{-1}\|\widetilde{J}\|_{\hk_{2}(S'(\de_*\ep_0^{-1}))} &\les& r\left\|\left(\De'+\frac{2}{{r'}^2}\right)\widetilde{J}\right\|_{L^2(S'(\de_*\ep_0^{-1}))} +|(\widetilde{J})_{\ell=1,\pm}|.
\eeaa
Since  ${J'}^{(p)}$ was extended  from $S_1'$ by $e_4'({J'}^{(p)})=0$, we have
\beaa
e_4'\left({r'}^{-2}\int_{S'}{J'}^{(p)} {J'}^{(q)} - \frac{4\pi}{3}\de_{pq}\right) &=& {r'}^{-2}\int_{S'}(\kac'-\ov{\kac'}){J'}^{(p)} {J'}^{(q)}=O\left(\frac{\ep_0}{{r'}^2}\right)
\eeaa
where we have used the estimate for $\kac'$ of Step 14, see \eqref{eq:EstimatesStep14}. Integrating from $S_1'\subset\Si_*$ where we have\footnote{Note  the change of notation,  the unprimed quantities  in  \eqref{eq:basicestimatesforJp-onSi_*:chap8}  are primed here.}  \eqref{eq:basicestimatesforJp-onSi_*:chap8}, we infer, for any $S'\subset\CC_1'$, 
\bea
\lab{eq:EstimateJ'-1.Step16}
\left|\int_{S'}{J'}^{(p)} {J'}^{(q)} - \frac{4\pi}{3}{r'}^2\de_{pq}\right| &\les& \ep r'.
\eea
Hence, since $r'\geq \de_*\ep_0^{-1}$ on $\CC_1'$, we deduce 
\beaa
\left|\int_{S'}{J'}^{(p)} {J'}^{(q)} - \frac{4\pi}{3}{r'}^2\de_{pq}\right| &\les& {r'}^2\ep_0, \qquad S'\subset\CC_1'.
\eeaa
We deduce, in particular, $|(\widetilde{J})_{\ell=1,\pm}| \les \ep_0$ and hence
\bea
\lab{eq:Estimatewidetilde{J}-1.Step16}
r^{-1}\|\widetilde{J}\|_{\hk_{2}(S'(\de_*\ep_0^{-1}))} &\les& r\left\|\left(\De'+\frac{2}{{r'}^2}\right)\widetilde{J}\right\|_{L^2(S'(\de_*\ep_0^{-1}))} +\ep_0.
\eea
Also, since we have extended ${J'}^{(p)}$ from $S_1'$ by $e_4'({J'}^{(p)})=0$,
\beaa
\bsplit
e_4'\left[{r'}^2\left(\De'+\frac{2}{{r'}^2}\right){J'}^{(p)}\right] &= [e_4', {r'}^2\De']{J'}^{(p)}= {\dkb'}^{\leq 1}(\kac'-\ov{\kac'}, \chih', \ze', r'\b'){\dkb'}^{\leq 2}{J'}^{(p)}\\
&= O(\ep_0 {r'}^{-2})
\end{split}
\eeaa
where we have used the estimate for $\kac'$, $\chih'$, $\ze'$ and $\b'$ of Step 14. Integrating from $S_1'\subset\Si_*$, and using the control on $\Si_*$ (and hence on $S_1'$) provided by Lemma \ref{Le:Si*-ell=1modes:proofofThmM0}, we infer along $\CC_1'$
\beaa
\left|\left(\De'+\frac{2}{{r'}^2}\right){J'}^{(p)}\right| &\les& \frac{\ep}{{r'}^3}.
\eeaa
In particular 
\beaa
\left|\left(\De'+\frac{2}{{r'}^2}\right){J'}^{(p)}\right| &\les& \frac{\ep_0}{{r'}^2}\qquad \mbox{on}\,\,\, S'(\de_*\ep_0^{-1}).
\eeaa

In view of the definition of $\widetilde{J}$ we infer  from \eqref{eq:Estimatewidetilde{J}-1.Step16}
\beaa
r^{-1}\|\widetilde{J}\|_{\hk_{2}(S'(\de_*\ep_0^{-1}))} &\les& r\left\|\left(\De'+\frac{2}{{r'}^2}\right)J^{(0)}\right\|_{L^2(S'(\de_*\ep_0^{-1}))} +\ep_0.
\eeaa
Together with the control of $r'-r$ of Step 14, this yields
\beaa
r^{-1}\|\widetilde{J}\|_{\hk_{2}(S'(\de_*\ep_0^{-1}))} &\les& r^{-1}\left\|\left(r^2\De'+2\right)J^{(0)}\right\|_{L^2(S'(\de_*\ep_0^{-1}))} +\ep_0.
\eeaa
Using the change of frame formula for $\De'$, and the control of $f$ and $\fb$ in Step 14, we infer
\beaa
r^{-1}\|\widetilde{J}\|_{\hk_{2}(S'(\de_*\ep_0^{-1}))} &\les& r^{-1}\left\|\left(r^2\De+2\right)J^{(0)}\right\|_{L^2(S'(\de_*\ep_0^{-1}))} +\ep_0.
\eeaa
Finally, in view of the control we have  for  the geodesic foliation of $\Lextt$, using also  Sobolev, we deduce
\beaa
\|\widetilde{J}\|_{L^\infty(S'(\de_*\ep_0^{-1}))} &\les& \ep_0
\eeaa
as stated in  Lemma \ref{Lemma:estmatesforwidetildeJandc_0}.

Next, we estimate $c_0$.  Using our control  for $\Lextt$, see Proposition \ref{proposition:geodesicfoliationLextt},
\beaa
\left|\int_{S}({J}^{(0)})^2  - \frac{4\pi}{3}{r}^2\right| &\les& \ep_0 r.
\eeaa
Using Lemma 7.3 in \cite{KS-GCM1}  (restated here in   Lemma \ref{lemma:consequencedeformationsurfaceusedinTheoremM0})   on the comparison between integrals on $S$ and on $S'$, and using also the control of $r'-r$ of Step 14, we obtain 
\beaa
\left|\int_{S'(\de_*\ep_0^{-1})}({J}^{(0)})^2  - \frac{4\pi}{3}{r'}^2\right| &\les& {r'}^2\ep_0.
\eeaa
 Also,  in view of  \eqref{eq:EstimateJ'-1.Step16},
\beaa
\left|\int_{S'(\de_*\ep_0^{-1})}({J'}^{(0)})^2  - \frac{4\pi}{3}{r'}^2\right| &\les& {r'}^2\ep_0,
\eeaa
and hence
\beaa
\left|\int_{S'(\de_*\ep_0^{-1})}({J}^{(0)})^2  - \int_{S'(\de_*\ep_0^{-1})}({J'}^{(0)})^2\right| &\les& {r'}^2\ep_0.
\eeaa
On the other hand, in view of the above control for $\widetilde{J}$ and its definition, we have
\beaa
\int_{S'(\de_*\ep_0^{-1})}({J}^{(0)})^2  &=& \int_{S'(\de_*\ep_0^{-1})}\Big((1+c_0){J'}^{(0)}+\widetilde{J}\Big)^2\\
&=&  (1+c_0)^2\int_{S'(\de_*\ep_0^{-1})}({J'}^{(0)})^2+2(1+c_0)\int_{S'(\de_*\ep_0^{-1})}{J'}^{(0)}\widetilde{J}+\int_{S'(\de_*\ep_0^{-1})}(\widetilde{J})^2.
\eeaa 
Together with our weak  control  for $c_0$ in \eqref{eq:localbootstrapassumptiononc0:proofThmM0}  and the above control for $\widetilde{J}$, we deduce
\beaa
2c_0\int_{S'(\de_*\ep_0^{-1})}({J'}^{(0)})^2 &=& \int_{S'(\de_*\ep_0^{-1})}({J}^{(0)})^2 - \int_{S'(\de_*\ep_0^{-1})}({J'}^{(0)})^2  +{r'}^2O(\ep_0+\ep^2).
\eeaa 
We obtain
\beaa
|c_0|{r'}^2 &\les& \left|\int_{S'(\de_*\ep_0^{-1})}({J}^{(0)})^2 - \int_{S'(\de_*\ep_0^{-1})}({J'}^{(0)})^2\right|  +{r'}^2\ep_0
\eeaa 
and hence, in view of the above, we infer
\beaa
|c_0| &\les& \ep_0.
\eeaa
This ends the proof of Lemma \ref{Lemma:estmatesforwidetildeJandc_0}.
\end{proof}

{\bf Step 16b.} The goal of this step is to prove the following lemma.
\begin{lemma}\lab{lemma:controlofaminusa0andJminusJprimeforproofThM0labellemma} 
The following  estimates hold true.
\bea\lab{eq:controlofaminusa0andJminusJprimeforproofThM0}
|a-a_0|+\max_{p=0,+,-}\sup_{S_1'}\left|{J}^{(p)}-{J'}^{(p)}\right| &\les& \ep_0\quad\textrm{if}\quad a_0\neq 0.
\eea
\end{lemma}

\begin{proof}
In view of the definition of $\widetilde{J}$, and the estimate  we have derived for it  and $c_0$ in   Lemma \ref{Lemma:estmatesforwidetildeJandc_0}, we have 
\bea
\sup_{S'(\de_*\ep_0^{-1})}|J^{(0)}-{J'}^{(0)}| &\les& \sup_{S'(\de_*\ep_0^{-1})}(|\widetilde{J}|+|c_0|)\les\ep_0.
\eea

Now, recalling from \eqref{eq:notyetthecontrolofainthecasea0notequal0ThM0} that we have
\beaa
{r'}^{-2}\left|\int_{S'(\de_*\ep_0^{-1})}{J'}^{(0)}J^{(0)} -\frac{4\pi a}{3a_0}{r'}^2\right| \les \ep_0,
\eeaa
we infer
\beaa
{r'}^{-2}\left|\int_{S'(\de_*\ep_0^{-1})}({J'}^{(0)})^2 -\frac{4\pi a}{3a_0}{r'}^2\right| \les \ep_0.
\eeaa
Together with the above control of $\int_{S'(\de_*\ep_0^{-1})}({J'}^{(0)})^2$, we deduce
\beaa
{r'}^{-2}\left|\frac{4\pi}{3}{r'}^2 - \frac{4\pi a}{3a_0}{r'}^2\right| \les \ep_0
\eeaa
and hence 
\beaa
|a-a_0| &\les&  \ep_0
\eeaa
as stated in  \eqref{eq:controlofaminusa0andJminusJprimeforproofThM0}. 

It remains to  prove the estimate for $|{J}^{(\pm)}-{J'}^{(\pm)}|$. To achieve  this, we introduce
 the following scalar functions 
\beaa
\widetilde{J}^{+}:=J^{(+)}-\Big((1+c_{++}){J'}^{(+)}+c_{+-}{J'}^{(-)}\Big), \quad \widetilde{J}^{-}:=J^{(-)}-\Big(c_{-+}{J'}^{(+)}+(1+c_{--}){J'}^{(-)}\Big),
\eeaa
where the constants $c_{++}$ and $c_{+-}$ are the solutions of the following 2 by 2 system
\beaa
c_{++}\int_{S'(\de_*\ep_0^{-1})}({J'}^{(+)})^2+c_{+-}\int_{S'(\de_*\ep_0^{-1})}{J'}^{(+)}{J'}^{(-)} &=& \int_{S'(\de_*\ep_0^{-1})}(J^{(+)}-{J'}^{(+)}){J'}^{(+)},\\
c_{++}\int_{S'(\de_*\ep_0^{-1})}{J'}^{(+)}{J'}^{(-)}+c_{+-}\int_{S'(\de_*\ep_0^{-1})}({J'}^{(-)})^2 &=& \int_{S'(\de_*\ep_0^{-1})}(J^{(+)}-{J'}^{(+)}){J'}^{(-)},
\eeaa
and where the constants $c_{-+}$ and $c_{--}$ are the solutions of the following 2 by 2 system
\beaa
c_{--}\int_{S'(\de_*\ep_0^{-1})}({J'}^{(-)})^2+c_{-+}\int_{S'(\de_*\ep_0^{-1})}{J'}^{(-)}{J'}^{(+)}&=&\int_{S'(\de_*\ep_0^{-1})}(J^{(-)}-{J'}^{(-)}){J'}^{(-)},\\
c_{--}\int_{S'(\de_*\ep_0^{-1})}{J'}^{(-)}{J'}^{(+)}+c_{-+}\int_{S'(\de_*\ep_0^{-1})}({J'}^{(+)})^2 &=& \int_{S'(\de_*\ep_0^{-1})}(J^{(-)}-{J'}^{(-)}){J'}^{(+)}.
\eeaa
In view of the bootstrap assumption \eqref{eq:bootstrapassumptiondifferenceell=1basisofCC1andLextprime} for ${J'}^{(p)}-\Jp$, we have
\bea\lab{eq:localbootstrapassumptiononc++c+-c--c-+:proofThmM0}
|c_{++}|+|c_{+-}|+|c_{-+}|+|c_{--}|\les \ep.
\eea
Also, in view of the definition of $\widetilde{J}^{\pm}$ and $c_{++}$, $c_{+-}$, $c_{-+}$ and $c_{--}$, we have
\beaa
\int_{S'(\de_*\ep_0^{-1})}{J'}^{(\pm)}\widetilde{J}^{+} = 0,\qquad \int_{S'(\de_*\ep_0^{-1})}{J'}^{(\pm)}\widetilde{J}^{-}=0, 
\eeaa
i.e. $(\widetilde{J}^+)_{\ell=1,\pm}=0$ and $(\widetilde{J}^-)_{\ell=1,\pm}=0$ on $S'(\de_*\ep_0^{-1})$. This yields, together with a Hodge elliptic estimate,
\beaa
r^{-1}\|\widetilde{J}^{\pm}\|_{\hk_{2}(S'(\de_*\ep_0^{-1}))} &\les& r\left\|\left(\De'+\frac{2}{{r'}^2}\right)\widetilde{J}^{\pm}\right\|_{L^2(S'(\de_*\ep_0^{-1}))} +|(\widetilde{J}^{\pm})_{\ell=1,0}|.
\eeaa
Arguing as above for $\widetilde{J}$, we have
\beaa
r\left\|\left(\De'+\frac{2}{{r'}^2}\right)\widetilde{J}^{\pm}\right\|_{L^2(S'(\de_*\ep_0^{-1}))} &\les& \sup_{S'(\de_*\ep_0^{-1})}\left|\left(\De'+\frac{2}{{r'}^2}\right)J^{(\pm)}\right|+\sup_{S'(\de_*\ep_0^{-1})}\left|\left(\De'+\frac{2}{{r'}^2}\right){J'}^{(\pm)}\right|\les\ep_0,
\eeaa
and hence
\beaa
r^{-1}\|\widetilde{J}^{\pm}\|_{\hk_{2}(S'(\de_*\ep_0^{-1}))} &\les& |(\widetilde{J}^{\pm})_{\ell=1,0}|+\ep_0.
\eeaa
Also, in view of the definition of $\widetilde{J}^\pm$, we have
\beaa
|(\widetilde{J}^{\pm})_{\ell=1,0}| &\les& r^{-2}\left|\int_{S'(\de_*\ep_0^{-1})}\widetilde{J}^\pm {J'}^{(0)}\right|\\
&\les& \sup_{{S'(\de_*\ep_0^{-1})}}|{J'}^{(0)}-J^{(0)}|+r^{-2}\left|\int_{S'(\de_*\ep_0^{-1})}{J'}^{(\pm)} {J'}^{(0)}\right|+r^{-2}\left|\int_{S'(\de_*\ep_0^{-1})}J^{(\pm)} J^{(0)}\right|.
\eeaa
In view of the above, we infer
\beaa
|(\widetilde{J}^{\pm})_{\ell=1,0}| &\les& \ep_0+r^{-2}\left|\int_{S'(\de_*\ep_0^{-1})}J^{(\pm)} J^{(0)}\right|.
\eeaa
Also, denoting $S_0$ the sphere of $\Lextt$ sharing the same south pole with $S'(\de_*\ep_0^{-1})$, we have in view of the control of $\Lextt$
\beaa
|(\widetilde{J}^{\pm})_{\ell=1,0}| &\les& \ep_0+r^{-2}\left|\int_{S'(\de_*\ep_0^{-1})}J^{(\pm)} J^{(0)}-\int_{S_0}J^{(\pm)} J^{(0)}\right|+r^{-2}\left|\int_{S_0}J^{(\pm)} J^{(0)}\right|\\
&\les& \ep_0+r^{-2}\left|\int_{S'(\de_*\ep_0^{-1})}J^{(\pm)} J^{(0)}-\int_{S_0}J^{(\pm)} J^{(0)}\right|.
\eeaa
Finally, in view of  the control of $(f, \fb)$ in Step 14, we may apply Lemma 7.3 in \cite{KS-GCM1} (restated here in   Lemma \ref{lemma:consequencedeformationsurfaceusedinTheoremM0})        with $\de_1=\ep_0$ which yields, together with $e_4(\Jp)=0$, 
\beaa
r^{-2}\left|\int_{S'(\de_*\ep_0^{-1})}J^{(\pm)} J^{(0)}-\int_{S_0}J^{(\pm)} J^{(0)}\right| &\les& \ep_0\max_{p=0,+,-}\sup_{\Lextt}|\dk^{\leq 1}\Jp|\les\ep_0
\eeaa
so that 
\beaa
|(\widetilde{J}^{\pm})_{\ell=1,0}| &\les& \ep_0
\eeaa
and hence
\beaa
r^{-1}\|\widetilde{J}^{\pm}\|_{\hk_{2}(S'(\de_*\ep_0^{-1}))} &\les& \ep_0.
\eeaa
Together with Sobolev, we deduce
\beaa
\sup_{S'(\de_*\ep_0^{-1})}|\widetilde{J}^{\pm}| &\les& \ep_0.
\eeaa

Next, using again Lemma 7.3 in \cite{KS-GCM1}   (restated here in   Lemma \ref{lemma:consequencedeformationsurfaceusedinTheoremM0})  with $\de_1=\ep_0$, we have 
\beaa
r^{-2}\left|\int_{S'(\de_*\ep_0^{-1})}(J^{(\pm)})^2-\int_{S_0}(J^{(\pm)})^2\right|+r^{-2}\left|\int_{S'(\de_*\ep_0^{-1})}J^{(+)}J^{(-)}-\int_{S_0}J^{(+)}J^{(-)}\right| &\les& \ep_0.
\eeaa
Together with the control of $\Lextt$ and the control of $r'-r$ of Step 14, we deduce
\beaa
\int_{S_{\de_*\ep_0^{-1}}'}(J^{(\pm)})^2=\frac{4\pi}{3}{r'}^2\Big(1+O(\ep_0)\Big), \qquad \int_{S_{\de_*\ep_0^{-1}}'}J^{(+)}J^{(-)}=O({r'}^2\ep_0).
\eeaa
Together with the above control of $\widetilde{J}^{\pm}$, we infer
\beaa
\frac{4\pi}{3}{r'}^2\Big(1+O(\ep_0)\Big) &=& \int_{S_{\de_*\ep_0^{-1}}'}(J^{(+)})^2= \int_{S_{\de_*\ep_0^{-1}}'}\left(\widetilde{J}^{+}+\Big((1+c_{++}){J'}^{(+)}+c_{+-}{J'}^{(-)}\Big)\right)^2\\
&=& \frac{4\pi}{3}{r'}^2\Big((1+c_{++})^2+(c_{+-})^2+O(\ep_0)\Big),
\eeaa
\beaa
O({r'}^2\ep_0) &=& \int_{S_{\de_*\ep_0^{-1}}'}J^{(+)}J^{(-)}\\
&=& \int_{S_{\de_*\ep_0^{-1}}'}\left(\widetilde{J}^{+}+\Big((1+c_{++}){J'}^{(+)}+c_{+-}{J'}^{(-)}\Big)\right)\left(\widetilde{J}^{-}+\Big(c_{-+}{J'}^{(+)}+(1+c_{--}){J'}^{(-)}\Big)\right)\\
&=& \frac{4\pi}{3}{r'}^2\Big((1+c_{++})c_{-+}++c_{+-}(1+c_{--})+O(\ep_0)\Big),
\eeaa
and 
\beaa
\frac{4\pi}{3}{r'}^2\Big(1+O(\ep_0)\Big) &=& \int_{S_{\de_*\ep_0^{-1}}'}(J^{(-)})^2= \int_{S_{\de_*\ep_0^{-1}}'}\left(\widetilde{J}^{-}+\Big(c_{-+}{J'}^{(+)}+(1+c_{--}){J'}^{(-)}\Big)\right)^2\\
&=& \frac{4\pi}{3}{r'}^2\Big((c_{-+})^2+(1+c_{--})^2+O(\ep_0)\Big),
\eeaa
which yields
\beaa
&& (1+c_{++})^2+(c_{+-})^2=1+O(\ep_0), \qquad (c_{-+})^2+(1+c_{--})^2=1+O(\ep_0), \\
&& (1+c_{++})c_{-+}+c_{+-}(1+c_{--})=O(\ep_0).
\eeaa
Together with the above control of $c_{++}$, $c_{+-}$, $c_{--}$ and $c_{--}$, we deduce
\beaa
1+c_{++}=\sqrt{1-(c_{+-})^2}+O(\ep_0),\quad 1+c_{--}=\sqrt{1-(c_{-+})^2}+O(\ep_0), \quad  c_{-+}=c_{+-}+O(\ep_0).
\eeaa
In particular,  there exists a real number $\varphi_0$ such that 
\bea
\bsplit
c_{+-}&=\sin(\varphi_0), \qquad c_{-+}=-\sin(\varphi_0)+O(\ep_0), \qquad |\varphi_0|\les\ep,\\
c_{++}&=\cos(\vphi_0)-1+O(\ep_0), \qquad c_{--}=\cos(\vphi_0)-1+O(\ep_0).
\end{split}
\eea
Together with the above definition, and the above control, of $\widetilde{J}^{\pm}$, we infer
\bea
\sup_{S_{\de_*\ep_0^{-1}}'}\left|J^{(\pm)}-\left(\cos(\vphi_0){J'}^{(\pm)}\pm \sin(\vphi_0){J'}^{(\mp)}\right)\right| &\les& \ep_0.
\eea

Also, we have
\beaa
e_4'\left({J}^{(0)}-{J'}^{(0)}\right) &=& e_4'\left({J}^{(0)}\right)=\la\left(e_4+f\c\nab'+\frac{1}{4}|f|^2e_3\right){J}^{(0)} = O(r^{-1}f)
\eeaa
and
\beaa
e_4'\left({J}^{(\pm)}-\left(\cos(\vphi_0){J'}^{(\pm)}\pm \sin(\vphi_0){J'}^{(\mp)}\right)\right) &=& e_4'\left({J}^{(\pm)}\right)\\
&=& \la\left(e_4+f\c\nab'+\frac{1}{4}|f|^2e_3\right){J}^{(\pm)} = O(r^{-1}f).
\eeaa
Together with the control of $f$ in Step 14, we obtain on $\CC_1'$
\beaa
\left|e_4'\left({J}^{(0)}-{J'}^{(0)}\right)\right|+\left|e_4'\left({J}^{(\pm)}-\left(\cos(\vphi_0){J'}^{(\pm)}\pm \sin(\vphi_0){J'}^{(\mp)}\right)\right)\right| &\les& \frac{\ep_0}{{r'}^2}
\eeaa
and hence, integrating forward from $S_{\de_*\ep_0^{-1}}'$, and using the above estimate for ${J}^{(0)}-{J'}^{(0)}$ and for ${J}^{(\pm)}-(\cos(\vphi_0){J'}^{(\pm)}\pm \sin(\vphi_0){J'}^{(\mp)})$ on $S_{\de_*\ep_0^{-1}}'$, we deduce on $\CC_1'$
\beaa
\left|{J}^{(0)}-{J'}^{(0)}\right|+\left|{J}^{(\pm)}-(\cos(\vphi_0){J'}^{(\pm)}\pm \sin(\vphi_0){J'}^{(\mp)})\right| &\les& \ep_0.
\eeaa
This estimate holds thus in particular on $S_1'$. Finally, recalling \eqref{eq:fixtheaxialfreefomofchoiceofJpmbasis}, we have 
\beaa
0 &=& \int_{S_1'}J^{(+)}{J'}^{(-)}=\int_{S_1'}\left(\cos(\vphi_0){J'}^{(+)}+ \sin(\vphi_0){J'}^{(-)}+O(\ep_0)\right){J'}^{(-)}\\
&=& \frac{4\pi}{3}{r'}^2\left(\sin(\vphi_0)+O(\ep_0)\right)
\eeaa
which implies $\sin(\vphi_0)=O(\ep_0)$.  We have thus obtain
\beaa
\max_{p=0,+,-}\sup_{S_1'}\left|{J}^{(p)}-{J'}^{(p)}\right| &\les& \ep_0\quad\textrm{if}\quad a_0\neq 0,
\eeaa
 which together with the previous estimate for $a-a_0$   establishes \eqref{eq:controlofaminusa0andJminusJprimeforproofThM0}. This  ends the proof of Lemma \ref{lemma:controlofaminusa0andJminusJprimeforproofThM0labellemma}.
\end{proof}

%%%%%%%%%%%%%%%%%

\subsection{Steps 17--24}

%%%%%%%%%%%%%%%%%

\noindent{\bf Step 17.}  
In the case $a_0\neq 0$, we derive  higher derivative  estimates for $f_0'-f_0$ and ${J}^{(p)}-{J'}^{(p)}$ on $\Si_*$.

In what follows we recall that $\phi$ is the  effective uniformization factor on  $S_*$ with  $(\th,\vphi)$ the corresponding coordinates and $\Jp$ the  corresponding  $\ell=1$ balanced  modes, see  \eqref{eq:balancedconditionforJponSstar}. We also recall the definition of the $1-$ forms $f_-, f_0, f_+$,   see   Definition \ref{def:definitionoff0fplusfminus}. Following our conventions above\footnote{Quantities  related to the PG frame are denoted by primes.} we  denote these by   $J'^{(p)} $ and  $f'_-, f'_0, f'_+$.

\begin{lemma}
\lab{Lemma:unifomfactor-ThmM0}
 The  effective   uniformization factor  $\phi$   of   $S_*$ verifies  
\bea
\|\phi\|_{\hk_{k+2}(S_*)}\les r'\left\|\Ga'_g\right\|_{\hk_{k}(S_*)}, \qquad  0\le k\le k_{large}.
\eea
\end{lemma}

\begin{proof}
According to Theorem   \ref{Thm:effectiveU1-Intro}, we have
\beaa
\|\phi\|_{\hk_{k+2}(S_*)}\les r'^2\left\|K'-\frac{1}{r'^2}\right\|_{\hk_{k}(S_*)}, \qquad 0\le  k\le k_{large}.
\eeaa
Also, in view  of the linearized Gauss equation of Proposition \ref{Prop.NullStr+Bianchi-lastslice}, we have $\widecheck{K}'\in r^{-1}\Ga_g'$ and hence
\beaa
\|\phi\|_{\hk_{k+2}(S_*)}\les r'\left\|\Ga_g'\right\|_{\hk_{k}(S_*)}
\eeaa
which concludes the proof of the lemma.
\end{proof}

We now proceed as follows: 
\begin{enumerate}
 \item Using   Lemmas \ref{lemma:computationforfirstorderderivativesJpell=1basis} and \ref{lemma:computationforfirstorderderivativesf0fpfm}, and simple elliptic estimates,  we deduce
\beaa
&&\|\phi\|_{\hk_{k+2}(S_*)}\les r'\|\Ga_g'\|_{\hk_k(S_*)}, \\
&&\max_{p=0,+,-}\|r'\widecheck{\nab'{J'}^{(p)}}\|_{\hk_{k+2}(S_*)}\les \|\phi\|_{\hk_{k+2}(S_*)},\\
&& \|r'\nab f_0' - {J'}^{(0)}\in\|_{\hk_{k+1}(S_*)}+\|r'\nab f_{\pm}' + {J'}^{(\pm)}\de\|_{\hk_{k+1}(S_*)}\les \|\phi\|_{\hk_{k+2}(S_*)}.
\eeaa
 Together with the bootstrap assumption for $\Ga_g'$ and $\Ga_b'$, we infer, for $k\leq k_{large}$, 
\beaa
&&\|\phi\|_{\hk_{k+2}(S_*)}\les \ep, \\
&&\max_{p=0,+,-}\|r'\widecheck{\nab'{J'}^{(p)}}\|_{\hk_{k+2}(S_*)}\les\ep,\\
&&\|r'\nab f_0' - {J'}^{(0)}\in\|_{\hk_{k+1}(S_*)}+\|r'\nab f_{\pm}' + {J'}^{(\pm)}\de\|_{\hk_{k+1}(S_*)}\les \ep.
\eeaa
By  Sobolev, we deduce
\beaa
\max_{p=0,+,-} r'\|{\dkb'}^{\leq k_{large}}r'\widecheck{\nab'{J'}^{(0)}}\|_{L^\infty(S_*)}+r'\|{\dkb'}^{\leq k_{large}-1}(r'\nab f_0' - {J'}^{(0)}\in)\|_{L^\infty(S_*)}\\
+r'\|{\dkb'}^{\leq k_{large}-1}(r'\nab f_{\pm}' + {J'}^{(\pm)}\de)\|_{L^\infty(S_*)}&\les& \ep.
\eeaa
\item On $\Si_*$ we have, using the transport equations  along $\Si_*$ of  Lemma \ref{lemma:transporteqnuofnabofJpf0fpfm}
\beaa
&&\nab_\nu\Big[r'\nab f_0' - {J'}^{(0)}\in\Big] =  \Ga_b' \c \dkb^{\leq 1}f_0',\qquad \nab_\nu\Big[r'\nab f_{\pm}' + {J'}^{(\pm)}\de\Big] =  \Ga_b' \c \dkb^{\leq 1}f_{\pm}',\\
&& \nab_\nu\Big[r'\widecheck{\nab'{J'}^{(p)}}\Big]=\Ga_b' \c \dkb^{\leq 1}{J'}^{(p)}, \quad p=0,+,-.
\eeaa
We integrate from $S_*$ and obtain, using the bootstrap assumption for  $\Ga_b'$,  
\beaa
\|{\dkb'}^{\leq k_{large}-1}(r'\nab f_0' - {J'}^{(0)}\in)\|_{L^\infty(\Si_*)}&\les& \|{\dkb'}^{\leq k_{large}-1}(r'\nab f_0' - {J'}^{(0)}\in)\|_{\hk_k(S_*)}+\frac{u_*}{r},\\
\|{\dkb'}^{\leq k_{large}-1}(r'\nab f_{\pm}' + {J'}^{(\pm)}\de)\|_{L^\infty(\Si_*)}&\les& \|{\dkb'}^{\leq k_{large}-1}(r'\nab f_{\pm}' + {J'}^{(\pm)}\de)\|_{\hk_k(S_*)}+\frac{u_*}{r},\\
\|{\dkb'}^{\leq k_{large}}(r'\widecheck{\nab'{J'}^{(p)}})\|_{L^\infty(\Si_*)} &\les& \|{\dkb'}^{\leq k_{large}}(r'\widecheck{\nab'{J'}^{(p)}})\|_{\hk_k(S_*)}+\frac{u_*}{r}, \quad p=0,+,-.
\eeaa
Together the above control on $S_*$ and the dominance condition for $r$ on $\Si_*$, we deduce 
\beaa
\|{\dkb'}^{\leq k_{large}-1}(r'\nab f_0' - {J'}^{(0)}\in)\|_{L^\infty(\Si_*)}&\les& \ep_0,\\
\|{\dkb'}^{\leq k_{large}-1}(r'\nab f_{\pm}' + {J'}^{(\pm)}\de)\|_{L^\infty(\Si_*)}&\les& \ep_0,\\
\|{\dkb'}^{\leq k_{large}}(r'\widecheck{\nab'{J'}^{(p)}})\|_{L^\infty(\Si_*)} &\les& \ep_0, \quad p=0,+,-.
\eeaa
\end{enumerate}

In view of the above, using Definition   \ref{def:renormalizationforJpbasisell=1modes}        of $\widecheck{\nab'{J'}^{(p)}}$, we infer on $S_1'$
\bea\lab{eq:controloff0primeonS1prime}
\left\|{\dkb'}^{\leq k_{large}-1}\left(\nab' f_0' - \frac{{J'}^{(0)}}{r'}\in\right)\right\|_{L^\infty(S_1')}  +\left\|{\dkb'}^{\leq k_{large}-1}\left(\nab' f_{\pm}' + \frac{{J'}^{(\pm)}}{r'}\de\right)\right\|_{L^\infty(S_1')} \\
\nn + \left\|{\dkb'}^{\leq k_{large}}\left(\nab'{J'}^{(0)}+\frac{1}{r'}\dual f_0'\right)\right\|_{L^\infty(S_1')}+ \left\|{\dkb'}^{\leq k_{large}}\left(\nab'{J'}^{(\pm)}-\frac{1}{r'}f_{\pm}'\right)\right\|_{L^\infty(S_1')} &\les& \frac{\ep_0}{r'}.
\eea

On the other hand, from the control of $\Lextt$, the change of frame formula for $\nab'$,  the control of the change of frame $(f, \fb, \la)$ from $\Lextt$ to $\Si_*$ of Step 13, and the control of $r'-r$ on $S_1'$ provided by \eqref{eq:finalesitateforffbchecklaonS1inproofofThmM0:reallyfinal}, we have
\beaa
\left\|{\dkb'}^{\leq k_{large}-1}\left(\nab' f_0 - \frac{{J}^{(0)}}{r'}\in\right)\right\|_{L^\infty(S_1')} +\left\|{\dkb'}^{\leq k_{large}-1}\left(\nab' f_{\pm} + \frac{{J}^{(\pm)}}{r'}\de\right)\right\|_{L^\infty(S_1')}\\
+ \left\|{\dkb'}^{\leq k_{large}}\left(\nab'{J}^{(0)}+\frac{1}{r'}\dual f_0\right)\right\|_{L^\infty(S_1')} + \left\|{\dkb'}^{\leq k_{large}}\left(\nab'{J}^{(\pm)}-\frac{1}{r'}f_{\pm}\right)\right\|_{L^\infty(S_1')}  &\les& \frac{\ep_0}{r},
\eeaa
and hence
\beaa
\left\|{\dkb'}^{\leq k_{large}-1}\left(\nab' (f_0'-f_0) - \frac{{J'}^{(0)}-J^{(0)}}{r'}\in\right)\right\|_{L^\infty(S_1')} \\
+\left\|{\dkb'}^{\leq k_{large}-1}\left(\nab' (f_{\pm}'-f_{\pm}) + \frac{{J'}^{(\pm)}-J^{(\pm)}}{r'}\de\right)\right\|_{L^\infty(S_1')}\\
+ \left\|{\dkb'}^{\leq k_{large}}\left(\nab'({J'}^{(0)}-J^{(0)})+\frac{1}{r'}\dual (f_0'-f_0)\right)\right\|_{L\infty(S_1')}\\
+ \left\|{\dkb'}^{\leq k_{large}}\left(\nab'({J'}^{(\pm)}-J^{(\pm)})-\frac{1}{r'}\dual (f_{\pm}'-f_{\pm})\right)\right\|_{L\infty(S_1')} &\les& \ep_0.
\eeaa
We deduce
\beaa
\left\|{\dkb'}^{\leq k_{large}}(f_0'-f_0)\right\|_{L^\infty(S_1')} + \left\|{\dkb'}^{\leq k_{large}+1}({J'}^{(0)}-J^{(0)})\right\|_{L^\infty(S_1')} \les \ep_0+\|{J'}^{(0)}-J^{(0)}\|_{L^\infty(S_1')}
\eeaa
and 
\beaa
\left\|{\dkb'}^{\leq k_{large}}(f_{\pm}'-f_{\pm})\right\|_{L^\infty(S_1')} + \left\|{\dkb'}^{\leq k_{large}+1}({J'}^{(\pm)}-J^{(\pm)})\right\|_{L^\infty(S_1')} \les \ep_0+\|{J'}^{(\pm)}-J^{(\pm)}\|_{L^\infty(S_1')}.
\eeaa
Together with \eqref{eq:controlofaminusa0andJminusJprimeforproofThM0}, we  obtain  in the case $a_0\neq 0$ 
\beaa
\max_{p=0,+,-}\left(\left\|{\dkb'}^{\leq k_{large}}(f_p'-f_p)\right\|_{L^\infty(S_1')} + \left\|{\dkb'}^{\leq k_{large}+1}({J'}^{(p)}-J^{(p)})\right\|_{L^\infty(S_1')}\right) &\les& \ep_0.
\eeaa
In view of the fact that $\nab_\nu f_p'=0$ for $p=0,+,-$ on $\Si_*$, and in view of the control of $\Lextt$, and hence of $f_p$  for $p=0,+,-$, we deduce for $a_0\neq 0$
\bea
\max_{p=0,+,-}\left(\left\|\dk_*^{\leq k_{large}}(f_p'-f_p)\right\|_{L^\infty(S_1')} + \left\|\dk_*^{\leq k_{large}+1}({J'}^{(p)}-J^{(p)})\right\|_{L^\infty(S_1')}\right) \les \ep_0.
\eea
In particular, this improves the bootstrap assumption \eqref{eq:bootstrapassumptiondifferenceell=1basisofCC1andLextprime} on ${J'}^{(p)}-J^{(p)}$ for $p=0,+,-$.

\noindent{\bf Step 18.} Next, we control $\nu(r')$ and $b_*$ on $S_1'$. First,  recall from \eqref{eq:finalesitateforRiccicoeffSigmastaronS1inproofofThmM0:reallyfinal} that we have, for $k\leq k_{large}+7$ for the frame of $\Si_*$, 
\beaa
r'|\dk_*^k(\eta', \xib', \ombc')| &\les& \ep_0. 
\eeaa
Together with Lemma \ref{Lemma:eqts-nabeta,xib}, and since  $b_*=-y'-z'$ in view of Lemma \ref{Lemma:eqts-nabeta,xib}, we infer, for $k\leq k_{large}+7$, 
\bea
r'|\nab'\dk_*^k(e_3(r)', e_3(u)', b_*)| &\les& \ep_0.
\eea

Then, proceeding as in Step 7 in the proof of Proposition \ref{Prop.Flux-bb-vthb-eta-xib} for the averages, we obtain on $S_1'$, for $k\leq  k_{large}+6$, 
\bea
\left|\dk_*^k\left(\nu(r')+2, b_*+1+\frac{2m}{r'}\right)\right| &\les& \ep_0.
\eea

\noindent{\bf Step 19.}  We consider the following change of frame coefficients: 
\begin{itemize}
\item $(f, \fb, \la)$ are the change of frame coefficients from the outgoing geodesic frame of $\Lextt$ to the  frame of $\Si_*$. They  satisfy, according to  \eqref{eq:finalesitateforffbchecklaonS1inproofofThmM0:reallyfinal},
\beaa
\sup_{k\le k_{ large}+1} \Big( \|{\dk}^k f \|_{L^2(S_{1}')} +r^{-1} \|{\dk}^{k} (\fb, \log \la) \|_{L^2(S_{1}')} +\|\dk^{\le k-1}\nab_\nu' (\fb, \log \la) \|_{L^2(S_{1}')}\Big) &\les \ep_0.
\eeaa

\item $(f', \fb', \la')$ are the change  of frame coefficients from the  outgoing PG  frame of $\Lext$ to the outgoing geodesic frame of $\Lextt$. $(f', \fb', \la')$ satisfies in view of Proposition \ref{proposition:geodesicfoliationLextt} 
\beaa
\sup_{\Lextt}r\left|\dk^{\leq k_{large}+8}\left(f'+\frac{a_0}{r}f_0, \fb'+\frac{a_0\Upsilon}{r}f_0, \log\la'\right)\right| &\les& \ep_0.
\eeaa

\item $(f'', \fb'', \la'')$ are the change of frame coefficients from the frame of $\Si_*$ to the outgoing PG frame of $\Mext$. $(f'', \fb'', \la'')$ satisfies by the initialization of the PG   structure on $\Si_*$, see section \ref{sec:initalizationadmissiblePGstructure}, 
\beaa
\la''=1,\qquad f''=\frac{a}{r'}f_0', \qquad \fb''=-\frac{(\nu(r')-b_*)}{1-\frac{1}{4}b_*|f''|^2}f''.
\eeaa
\end{itemize}

We now consider the change of frame   coefficients $(f''', \fb''', \la''')$  from the  outgoing PG frame of $\Lext$ to the outgoing PG  frame of $\Mext$.  In view of:
\begin{itemize}
\item the above estimates for $(f, \fb, \la)$ and $(f', \fb', \la')$,

\item the above formula for $(f'', \fb'', \la'')$,

\item the control for $r-r'$ and $m-m_0$ given by \eqref{eq:finalesitateforffbchecklaonS1inproofofThmM0:reallyfinal}, and the control for $\nu(r')$ and $b_*$ in Step 18, 

\item the control of $a$ in \eqref{eq:controlofainthecasea0equal0ThM0} in the case $a_0=0$, 

\item the control for $a-a_0$ in \eqref{eq:controlofaminusa0andJminusJprimeforproofThM0} and the control for $f_0'-f_0$ in Step 17 in the case $a_0\neq 0$,  
\end{itemize}
we infer the following estimates 
\beaa
\sup_{S_1'}r|\dk^{\leq k_{large}}f'''| +\sup_{S_1'}|\dk^{\leq k_{large}}(\fb''', \log\la''')| &\les& \ep_0+\frac{1}{r}.
\eeaa
Together with the dominance condition for $r$ on $\Si_*$, we infer
\bea\lab{eq:controlofffblabetweenPGLextandPGMextonS1}
\sup_{S_1'}r|\dk^{\leq k_{large}}f'''| +\sup_{S_1'}|\dk^{\leq k_{large}}(\fb''', \log\la''')| &\les& \ep_0.
\eea

\noindent{\bf Step 20.}  In Steps 20--22, $(e_1, e_2, e_3, e_4)$ denotes the outgoing PG frame of $\Lext$, and $(e_1', e_2', e_3', e_4')$ denotes the outgoing PG frame of $\Mext$. Also, from now on, $(f, \fb, \la)$ denotes\footnote{ Denoted earlier by triple prime.} the change of frame   coefficients   from the  outgoing PG frame of $\Lext$ to the outgoing PG  frame of $\Mext$.  In view of Step 19, we have on $S_1'$
\beaa
\sup_{S_1'}r|\dk^{\leq k_{large}}f| +\sup_{S_1'}|\dk^{\leq k_{large}}(\fb, \log\la)| &\les& \ep_0.
\eeaa

Let
\beaa
F=f+i\dual f.
\eeaa
Since 
\beaa
\Xi'=0, \qquad \om'=0, \qquad \Xi=0, \qquad \om=0, 
\eeaa
we  have, in view of Corollary \ref{cor:transportequationine4forchangeofframecoeffinformFFbandlamba}, 
\beaa
\nab_{\la^{-1}e_4'}(\ov{q}F) &=& E_4(f, \Ga),\\
\la^{-1}\nab_4'(\log\la) &=& f\c(\ze-\etab)+E_2(f, \Ga).
\eeaa
We integrate the above transport equations for $F$ and $\la$  in the order they appear from $S_1'$.  In view of the control for $(f, \la)$ on $S_1'$ derived in Step 19, and in view of the assumptions on the initial data layer norm, we infer 
\bea
\sup_{\{u'=1\}}\Big(r|\dk^{\leq k_{large}}f|+|\dk^{\leq k_{large}}\log(\la)|\Big) &\les& \ep_0,
\eea 
where $u'$ denotes from now on the scalar function of the PG structure of $\Mext$. In particular, we have by construction $S_1'=\Si_*\cap\{u'=1\}$, $e_4'(u')=0$ and $\{u'=1\}\subset\Lext$.

\noindent{\bf Step 21.} In this step, we estimate $r'-r$, as well as $A'$, $\trXc'$ and $\Xh'$. Moreover, in the case $a_0\neq 0$, we also estimate ${J'}^{(0)}-J^{(0)}$ and $\Jk'-\Jk$. First, since ${J'}^{(0)}$ is propagated from $\Si_*$ by $e_4'({J'}^{(0)})=0$, and using the change of frame formula between the PG frame of $\Lext$ and the PG frame of $\Mext$, we infer
\beaa
e_4'({J'}^{(0)}-J^{(0)}) &=& -\la\left(f\c\nab+\frac{1}{4}|f|^2e_3\right)J^{(0)}.
\eeaa
Together with the control of $f$ and $\la$ of Step 20, and the control of $\Lext$, we infer
\beaa
\sup_{\{u'=1\}}r^2|\dk^{\leq k_{large}}(e_4'({J'}^{(0)}-J^{(0)}))| &\les& \ep_0.
\eeaa
Integrating from $S_1'$, where ${J'}^{(0)}-J^{(0)}$ is under control in view of Step 17 in the case $a_0\neq 0$, we infer
\bea
\sup_{\{u'=1\}}|\dk^{\leq k_{large}}({J'}^{(0)}-J^{(0)})| &\les& \ep_0\quad\textrm{if}\quad a_0\neq 0.
\eea

Next, we control $\trXc'$. To this end, we also need to control $A'$ and $\Xh'$. First, note that the change of frame formula for $A'$, the control of the foliation of $\Lext$, the control of $f$ and $\la$ of Step 20, and the fact that the transformation formula for $A'$ does not depend on $\fb$ implies 
\bea
\sup_{\{u'=1\}}r^{\frac{7}{2}+\dt}|\dk^{\leq k_{large}}A'| &\les& \ep_0.
\eea
Then, using the control of the Ricci coefficients of the frame of $\Si_*$ obtained in  \eqref{eq:finalesitateforRiccicoeffSigmastaronS1inproofofThmM0:reallyfinal}, the control of $a$ in \eqref{eq:controlofainthecasea0equal0ThM0} in the case $a_0=0$, the control for $a-a_0$ in \eqref{eq:controlofaminusa0andJminusJprimeforproofThM0} and the control for $f_0'-f_0$ in  \eqref{eq:controloff0primeonS1prime} in the case $a_0\neq 0$, the control of $\nu(r')$ and $b_*$ of Step 18, and the change of frame formula between the frame of $\Si_*$ and the outgoing PG frame of $\Mext$ on $\Si_*$, we infer\footnote{Note that while the control of the  Ricci coefficients $\trchbc$ of $\Si_*$ in  \eqref{eq:finalesitateforRiccicoeffSigmastaronS1inproofofThmM0:reallyfinal} displays a loss of $r^{-1}$, this allows nevertheless to obtain the correct power of $r$ of $\trXc'$ and $\Xh'$ on $S_1'$.}
\beaa
\sup_{S_1'}r^2\Big(|\dk^{\leq k_{large}-1}\trXc'|+|\dk^{\leq k_{large}-1}\Xh'|\Big) &\les& \ep_0.
\eeaa
Then, propagating Raychadhuri for $\trXc'$ and the null structure equation for $\nab_4'\Xh'$ from $S_1'$ where $\trXc'$ and $\Xh'$  are under control in view of the above estimate, we infer, using the above control of $A'$, 
\bea
\sup_{\{u'=1\}}r^2\Big(|\dk^{\leq k_{large}-1}\trXc'|+|\dk^{\leq k_{large}-1}\Xh'|\Big) &\les& \ep_0.
\eea

Next, we notice
\beaa
\tr X' -\la\tr X &=& \frac{2}{q'}-\frac{2}{q}+(\la-1)\tr X+\trXc' -\trXc
\eeaa
so that 
\beaa
q'-q &=& \frac{qq'}{2}\left(-\la\Big(\la^{-1}\tr X' -\tr X\Big)+(\la-1)\tr X+\trXc' -\trXc\right)
\eeaa
Together with the control of Step 20 for $f$ and $\la$, the above control $\trXc'$, the above control for  ${J'}^{(0)}-J^{(0)}$ if $a_0\neq 0$, the control of $a$ in \eqref{eq:controlofainthecasea0equal0ThM0} in the case $a_0=0$, the control for $a-a_0$ in \eqref{eq:controlofaminusa0andJminusJprimeforproofThM0} in the case $a_0\neq 0$, the control of the foliation $\Lext$, and the fact that $q=r+ia_0J^{(0)}$ and $q'=r'+ia{J'}^{(0)}$, we infer
\beaa
\sup_{\{u'=1\}}\left|\dk^{\leq k_{large}-1}\left(\frac{r'}{r}-1 +\frac{qq'}{2r}\Big(\la^{-1}\tr X' - \tr X\Big)\right)\right| &\les& \ep_0.
\eeaa
Moreover, from the change of frame formulas for $\trch$ and $\atrch$
we have, schematically, 
\beaa
\la^{-1}\tr X' -\tr X &=& {r'}^{-1}{\dk'}f+\Ga\c f+f\c\fb\c\Ga+f\c\fb\c\tr X'+\lot
\eeaa
Together with the control of Step 20 for $f$, the above control for  $\trXc'$ and the control of the foliation $\Lext$, we deduce
\bea
\sup_{\{u'=1\}}\left|\dk^{\leq k_{large}-1}\left(\frac{r'}{r}-1\right)\right| &\les& \ep_0+\ep_0\sup_{\{u'=1\}}\left|\dk^{\leq k_{large}-1}\fb\right|.
\eea

Next, recall the definition of $\Jk'$ on $\Si_*$
\beaa
\Jk' &=& \frac{1}{|q'|}(f_0'+i\dual f_0')=\frac{1}{\sqrt{{r'}^2+a^2({J'}^{(0)})^2}}(f_0'+i\dual f_0').
\eeaa
Recall that we control $\Jk'-\Jk$ only in the case $a_0\neq 0$. Together with the control of $r-r'$ and $m-m_0$ given by \eqref{eq:finalesitateforffbchecklaonS1inproofofThmM0:reallyfinal}, the control of $a-a_0$ in Step 16 for $a_0\neq 0$,  the control for $f_0'-f_0$ and ${J'}^{(0)}-J^{(0)}$ in Step 17 for $a_0\neq 0$, and the control of $\Lext$, we infer
\beaa
r'\|\dk_*^{\leq k_{large}}(\Jk'-\Jk)\|_{L^\infty(S_1')} &\les& \ep_0+\frac{1}{r}. 
\eeaa
Together with the dominance condition for $r$ on $\Si_*$, we infer
\beaa
r'\|\dk_*^{\leq k_{large}}(\Jk'-\Jk)\|_{L^\infty(S_1')} &\les& \ep_0. 
\eeaa
Together with the identity
\beaa
q'\Jk'-q\Jk &=& q'(\Jk'-\Jk)+(q'-q)\Jk\\
&=& q'(\Jk'-\Jk)+\Big(r'-r+i(a{J'}^{(0)}-a_0J^{(0)})\Big)\Jk,
\eeaa
the control of $r'-r$ given by \eqref{eq:finalesitateforffbchecklaonS1inproofofThmM0:reallyfinal}, the control of ${J'}^{(0)}-J^{(0)}$ in Step 17 and the control of $a-a_0$ of Step 16, we obtain
\beaa
\|\dk_*^{\leq k_{large}}(q'\Jk'-q\Jk)\|_{L^\infty(S_1')} &\les& \ep_0. 
\eeaa
Also, recall that $\Jk'$ and $\Jk$ satisfy in $\{u'=1\}$
\beaa
\nab_4'\Jk' = -\frac{1}{q'}\Jk', \qquad \nab_4\Jk= -\frac{1}{q}\Jk,
\eeaa
and hence
\beaa
\nab_4'(q'\Jk')=0, \qquad \nab_4(q\Jk)=0.
\eeaa
We infer
\beaa
\nab_{\la^{-1}4'}(q'\Jk')=0, \qquad \nab_{\la^{-1}4'}(q\Jk)=\left(f\c\nab +\frac{1}{4}|f|^2e_3\right)(q\Jk).
\eeaa
Together with the control of $f$ and $\la$ of Step 20, and the control of $\Lext$, we obtain
\beaa
\sup_{\{u'=1\}}r^2|\dk^{\leq k_{large}}(\nab_4'(q'\Jk'-q\Jk))| &\les& \ep_0.
\eeaa
Integrating from $S_1'$ where $q'\Jk'-q\Jk$ is under control in view of the above, we infer
\beaa
\sup_{\{u'=1\}}|\dk^{\leq k_{large}}(q'\Jk'-q\Jk)| &\les& \ep_0.
\eeaa
Using again the above identity for $q'\Jk'-q\Jk$, as well as the above control of ${J'}^{(0)}-J^{(0)}$ and $r'-r$,  and the control of $a-a_0$ of Step 16, we deduce
\bea
\sup_{\{u'=1\}}r|\dk^{\leq k_{large}-1}(\Jk'-\Jk)| &\les& \ep_0+\ep_0\sup_{\{u'=1\}}\left|\dk^{\leq k_{large}-1}\fb\right|\quad\textrm{if}\quad a_0\neq 0.
\eea

\noindent{\bf Step 22.} In this step, we control $\fb$ on $\Mext$. To this end, we first control $B'$ and $Z'$.  The change of frame formula for $B'$, the control of the foliation of $\Lext$,  the control of $f$ and $\la$ of Step 20, and the fact that the terms involving $\fb$ in the transformation formula for $B'$ are at least quadratic,  implies 
\bea
\sup_{\{u'=1\}}r^{\frac{7}{2}+\dt}|\dk^{\leq k_{large}-1}B'| &\les& \ep_0+\ep_0\sup_{\{u'=1\}}\left|\dk^{\leq k_{large}-1}\fb\right|.
\eea
Then, propagating  the null structure equation for $\nab_4'\Zc'$ from $S_1'$ where it is under control in view of  \eqref{eq:finalesitateforRiccicoeffSigmastaronS1inproofofThmM0:reallyfinal} for the frame of $\Si_*$ and the change of frame formula, we infer, using the above control of $B'$, 
\bea
\sup_{\{u'=1\}}r^2|\dk^{\leq k_{large}-1}\Zc'| &\les& \ep_0+\ep_0\sup_{\{u'=1\}}\left|\dk^{\leq k_{large}-1}\fb\right|.
\eea

Next, the transformation formula for $Z'$, together with the control of the foliation of $\Lext$,  the control of $f$ and $\la$ of Step 20, and the control of $\trXc'$ and $\Xh'$ of Step 21,  yields
\beaa
\sup_{\{u'=1\}}r\left|\dk^{\leq k_{large}-1}\left(Z'-Z-\frac{1}{4}\tr X(\fb+i\dual\fb)\right)\right| &\les& \ep_0+\ep_0\sup_{\{u'=1\}}\left|\dk^{\leq k_{large}-1}\fb\right|.
\eeaa
Since we have
\beaa
Z'-Z &=& \frac{a\ov{q'}}{|q'|^2}\Jk' -\frac{a_0\ov{q}}{|q|^2}\Jk+\Zc'-\Zc,
\eeaa
we deduce, together with the control of the foliation of $\Lext$, the above control of $\Zc'$, the control of $r'-r$ of Step 21, the control of $a$ in \eqref{eq:controlofainthecasea0equal0ThM0} in the case $a_0=0$, the control for $a-a_0$ in \eqref{eq:controlofaminusa0andJminusJprimeforproofThM0} in the case $a_0\neq 0$,  and the control of  ${J'}^{(0)}-J^{(0)}$ and $\Jk'-\Jk$ of Step 21 in the case $a_0\neq 0$, 
\beaa
\sup_{\{u'=1\}}\left|\dk^{\leq k_{large}-1}\fb\right| &\les& \ep_0+\ep_0\sup_{\{u'=1\}}\left|\dk^{\leq k_{large}-1}\fb\right|
\eeaa
and hence
\beaa
\sup_{\{u'=1\}}\left|\dk^{\leq k_{large}-1}\fb\right| &\les& \ep_0.
\eeaa
Together with the control of $f$ and $\la$ of Step 20, we have finally obtained
\bea\lab{eq:controlofffblambdabetwenPGframeLextandMextonu=1}
\sup_{\{u'=1\}}\Big(r|\dk^{\leq k_{large}}f|+|\dk^{\leq k_{large}}\log(\la)|+|\dk^{\leq k_{large}-1}\fb|\Big) &\les& \ep_0.
\eea
Also, together with the estimates of Step 21 for $r'-r$, we obtain 
\bea
\sup_{\{u'=1\}}\left|\dk^{\leq k_{large}-1}\left(\frac{r'}{r}-1\right)\right|&\les& \ep_0.
\eea

\noindent{\bf Step 23.} Let  $(f', \fb', \la')$ denote the change of frame   coefficients   from the principal outgoing null frame of $\Lint$ to the principal outgoing null frame of $\Mint$. From
\begin{itemize}
\item the estimates of Step 22 on $\{u'=1\}$, 

\item the fact that $\Mint\cap\Mext=\{r'=r_0\}$,

\item the fact that $\{u=1\}\cap\{\ub'=1\}$ is included in $\Lext\cap\Lint$,

\item the initialization of the frame of $\Mint$ as an explicit renormalization of the frame of $\Mext$ on $\{r'=r_0\}$, 

\item the control in $\Lext\cap\Lint$ of the difference between the frame of $\Lint$ and an explicit renormalization of the frame of $\Lext$,
\end{itemize}
we easily infer, using also $\ub'=u'$ on $\{r=r_0\}$, 
\beaa
\sup_{\{r=r_0\}\cap\{\ub'=1\}}\Big(|\dk^{\leq k_{large}-1}(f', \fb', \log\la')|\Big) &\les& \ep_0.
\eeaa 

Next, we proceed as in Step 20, exchanging the role of $e_3$ and $e_4$, and we propagate along $e_3$ the above estimate to $\{\ub'=1\}$ for $\fb'$ and $\la'$. We also propagate the control of Step 21 for ${J'}^{(0)}-J^{(0)}$ on $\{u'=1\}$ in the case $a_0\neq 0$, and hence on its boundary $\{r'=r_0\}$ to $\{\ub'=1\}$. Also one propagates the control of Step 22 for $r-r'$ on $\{u'=1\}$, and hence on its boundary $\{r'=r_0\}$ to $\{\ub'=1\}$ using the transport equation\footnote{Note that we could not have used this transport equation in $\Lext$ in view of the lack of decay in $r$ for $\la-1$. This is why we avoided this transport equation in Step 21 and used instead the control of $\trXc'$. On the other hand, $r$ is bounded in $\Lint$ so that one can simply rely on the transport equation for $e_3'(r'-r)$ in $\Lint$.}
\beaa
e_3'(r'-r) &=& 1-\la'\left(e_3+{\fb'}^ae_a+\frac{1}{4}|\fb'|^2e_4\right)r\\
&=& -(\la'-1)+\fb'\c\nab(r)+ \frac{1}{4}|\fb'|^2e_4(r).
\eeaa
Finally, we propagate $f$ similarly to Step 22. We finally obtain   
\bea\lab{eq:controlofffblambdabetwenPGframeLextandMextonu=1:ingoingcase}
\sup_{\{\ub'=1\}}\Big(|\dk^{\leq k_{large}-1}(\fb', \log\la')|+|\dk^{\leq k_{large}-2}(r'-r, f')|\Big) &\les& \ep_0,
\eea 
and
\bea
\sup_{\{\ub'=1\}}\left|\dk^{\leq k_{large}-2}\left({J'}^{(0)}-J^{(0)}\right)\right| &\les& \ep_0, \quad\textrm{if}\quad a_0\neq 0.
\eea

\noindent{\bf Step 24.} Note that the desired estimate for $m-m_0$ has been obtained in Step 13. Also, note that the desired estimate for $a-a_0$ has been obtained in Step 15 in the case $a_0=0$, and in Step 16 in the case $a_0\neq 0$. To conclude the proof of Theorem M0, it remains to control $k_{large}-2$ derivatives, with suitable $r$-weights and $O(\ep_0)$ smallness constant, of $A'$, $B'$, $\Pc'$, $\Bb'$ and $\Ab'$ in $\{u'=1\}\cup\{\ub'=1\}$, i.e. 
\bea\lab{eq:gatheringtheconclusionsyieldingtheproofofThmM0:1}
\max_{0\leq k\leq k_{large-2}}\Bigg\{ \sup_{\{u'=1\}}\Bigg[r^{\frac{7}{2}  +\de_B}\left( |\dk^k\,{}^{(ext)}A'| + |\dk^k\,{}^{(ext)}B'|\right)  +r^3 \left|\dk^k\left(\,{}^{(ext)}P'+\frac{2m}{{q'}^3}\right)\right|\\
\nn+r^2|\dk^k\,{}^{(ext)}\Bb'|+r|\dk^k\,{}^{(ext)}\Ab'|\Bigg]\Bigg\} +\max_{0\leq k\leq k_{large-2}}\sup_{\BBb_1}\Bigg[ |\dk^k\,{}^{(int)}A'| \\
\nn+ |\dk^k\,{}^{(int)}B'|+ \left|\dk^k\left(\,{}^{(int)}P' +\frac{2m}{{q'}^3}\right)\right|+|\dk^k\,{}^{(int)}\Bb'|+|\dk^k\,{}^{(int)}\Ab'|\Bigg] &\les& \ep_0.
\eea
This follows from: 
\begin{itemize}
\item the control of $(f, \fb, \la)$ on $\{u'=1\}$ derived in Step 22, 

\item the control of $(f', \fb', \la')$ on $\{\ub'=1\}$ derived in Step 23,

\item the fact that $(f, \fb, \la)$ denote the change of frame   coefficients   from the PG frame of $\Lext$ to the PG  frame of $\Mext$, and the fact that $(f', \fb', \la')$ denote the change of frame   coefficients   from the principal outgoing null frame of $\Lint$ to the principal outgoing null frame of $\Mint$,

\item the change of frame formulas for the curvature components,

\item in the particular case of the estimate for $\Pc$, the fact that 
\beaa
P'-P &=& -\frac{2m}{{q'}^3}+\frac{2m_0}{q^3}+\Pc'-\Pc,
\eeaa
together with the control of $m-m_0$ derived in Step 13, the control of $r'-r$ in Steps 22 and 23, the control of $a$ in Step 15 in the case $a_0=0$, the control of $a-a_0$ in Step 16 in the case $a_0\neq 0$, and the control of ${J'}^{(0)}-J^{(0)}$ in Step 21 and 23 in the case $a_0\neq 0$,

\item the assumptions on the initial data layer norm.
\end{itemize}
The proof of Theorem M0 is now complete.

%%%%%%%%%%%%%%%%%%%%%%%%%%%%%%%%%%%%%%%%

\section{Proof of Theorem M6}
\lab{sec:proofofTheoremM6}

%%%%%%%%%%%%%%%%%%%%%%%%%%%%%%%%%%%%%%%%

The proof of Theorem M6 proceeds in 8 steps which we summarize below for convenience:
\begin{enumerate}
\item In Steps 1--3, we construct our last sphere $S_*$, and then our last slice $\Si_*$ inside the part $\Lextt$ of the initial data layer, by relying on the control of $\Lextt$ provided by Proposition \ref{proposition:geodesicfoliationLextt}, and the GCM constructions of \cite{KS-GCM2} and \cite{Shen} recalled in section \ref{sec:GCMpapersreview}.

\item In Steps 4--8, we construct from $\Si_*$ a GCM admissible spacetime $\MM$ and  
we control the change of frame coefficients between the frames of the initial data layer $\LL_0$, and the corresponding frames of $\MM$. In view of the control of $\LL_0$ and of the change of frames coefficients, we infer the desired control of $\MM$ thanks to the change of frame formulas.
\end{enumerate}

We now proceed with the proof of Theorem M6.

{\bf Step 1.} Let $r_{(0)}$ such that 
\bea\lab{eq:choicer0forproofThmM6}
r_{(0)} &:=& d_0\de_*\ep_0^{-1},
\eea
where the small constant $\de_*$ appears in \eqref{eq:behaviorofronS-star}, and where the constant $d_0$ satisfies 
\beaa
\frac{1}{2}\leq d_0\leq 2
\eeaa 
and will be suitably chosen in Step 3. Also,  let $\de_0>0$ sufficiently small. Consider the unique sphere $\ovS$ of the part $\Lextt$ of the initial data layer on $\{\ut=1+\de_0\}$ with area radius $r_{(0)}$. Then, denoting $S(\ut, \st)$ the spheres of the outgoing geodesic foliation of $\Lextt$, we have 
\beaa
\ovS=S(\ug, \sg), \qquad \ug=1+\de_0, \qquad |\sg-r_{(0)}|\les \ep_0,
\eeaa
where the control of $\sg - r_{(0)}$ follows from the assumptions on the control of $\Lextt$. Relying on the control of the initial data layer given by \eqref{def:initialdatalayerassumptions}, i.e. 
\beaa
\Ik_{k_{large}+10}\leq \ep_0, \qquad {}^{(ext)}\Ik_{3}\leq \ep_0^2,
\eeaa
we are in position to apply  Theorem 7.3 of \cite{KS-GCM2}   (restated here  as  Theorem \ref{theorem:ExistenceGCMS2})   and Corollary 7.7 of \cite{KS-GCM2} (restated here  as Corollary \ref{Corr:ExistenceGCMS2}), with the choices
\beaa
\dg=\epg=\ep_0, \quad s_{max}=k_{large}+7,
\eeaa
to produce a unique GCM sphere $S_*$, which is a deformation of $\ovS$, satisfying 
\bea\lab{eq:GCMconditionsinwidetildeStar:ThmM6}
 \bsplit
& \widecheck{\ka}^{S_*}=0, \qquad  \widecheck{\kab}^{S_*}=0, \qquad  \widecheck{\mu}^{S_*}=\sum_pM^{S_*}_pJ^{(p,S_*)}, \\
& (\div^{S_*}\b^{S_*})_{\ell=1}=0, \qquad  (\curl^{S_*}\b^{S_*})_{\ell=1,\pm}=0, \qquad  (\curl^{S_*}\b^{S_*})_{\ell=1,0}=\frac{2a^{S_*}m^{S_*}}{(r^{S_*})^5}, 
\end{split}
\eea
where
\begin{itemize}
\item $J^{(p, S_*)}$ denotes the canonical basis of $\ell=1$ mode on $S_*$ in the sense of Definition 3.10 of \cite{KS-GCM2} (recalled  here  in  Definition \ref{definition:ell=1mpdesonS-intro}),

\item the $\ell=1$ modes in \eqref{eq:GCMconditionsinwidetildeStar:ThmM6} are defined w.r.t. the basis of $\ell=1$ modes $J^{(p,S_*)}$, 

\item $m^{S_*}$ denotes the Hawking mass of $S_*$, $r^{S_*}$ denotes the area radius of $S_*$, and the identity for $(\curl^{S_*}\b^{S_*})_{\ell=1,0}$ in \eqref{eq:GCMconditionsinwidetildeStar:ThmM6} should be understood as providing the definition of $a^{S_*}$. 
\end{itemize}

\begin{remark}\lab{rem:thisistheonlyplacewhereweneedstrongercontrolofinitialdata}
In order to  apply   Theorem 7.3 of \cite{KS-GCM2}   (restated here  as  Theorem \ref{theorem:ExistenceGCMS2})   and Corollary 7.7 of \cite{KS-GCM2} (restated here  as Corollary \ref{Corr:ExistenceGCMS2})
to  the above setting, one needs to check that  foliation of $\Lextt$  satisfies the assumptions of the theorem, and  in particular, in the region $r\sim r_{(0)}$ of $\Lextt$,   
\beaa
r^5\left(|(\div\b)_{\ell=1}|+|(\curl\b)_{\ell=1,\pm}|+ \left|(\curl\b)_{\ell=1, 0}-\frac{2a_0m_0}{r^5}\right|\right)
 +r^3|(\kac)_{\ell=1}|+r^3|(\kabc)_{\ell=1}|  \les \dg,
\eeaa
as well as 
\beaa
r^2|\dk^k\Ga_b| &\les& \de_*,\qquad k\leq s_{max}.
\eeaa
Now,  in view of the above choice for $s_{max}$, $\dg$, $\epg$ and $r_{(0)}$, this follows from 
\beaa
r| \dk^{\leq k_{large}+7}\Ga_b |\les \ep_0
 \eeaa
 and 
 \beaa
\nn \sup_{\Lextt\cap\{r\sim \ep_0^{-1}\}  } r^5\left( \left|\div\b\right|+\left|\curl\b -\frac{6a_0m_0}{r^5}J^{(0)} \right|\right)  +\sup_{\Lextt\cap\{r\sim \ep_0^{-1}\}  } r^3\left( \left|\kac\right|+\left|\kab\right|\right)  &\les&  \ep_0
\eeaa
 and hence, in view of Proposition \ref{proposition:geodesicfoliationLextt}, from 
 \beaa
\Ik_{k_{large}+10}\leq \ep_0, \qquad \, ^{(ext)} \Ik_{3} \leq \ep^2_0.
\eeaa
\end{remark}

From now on, we denote for simplicity  
\bea
m:=m^{S_*}, \qquad a:=a^{S_*}.
\eea

{\bf Step 2.} Starting from $S_*$ constructed in Step 1, and relying on the control provided by Proposition \ref{proposition:geodesicfoliationLextt} for the foliation of $\Lextt$, we may then apply  Theorem 4.1  in  \cite{Shen} (restated here in Theorem \ref{theorem:constuctionGCMH}),  with $s_{max}=k_{large}+7$, which yields the existence of a smooth small piece of spacelike hypersurface $\Si_*$ passing through the sphere $S_*$, together with a scalar function $u$ defined on $\Si_*$, whose level surfaces are topological spheres denoted by $S$, so that 
\begin{itemize}
\item The following GCM conditions are verified on $\Si_*$
 \beaa
&& \kac=0, \qquad \kabc=\Cb_0+\sum_p\Cb_p\Jp, \qquad \muc=M_0+\sum_pM_p\Jp,\\
&& (\div\eta)_{\ell=1}=0, \qquad (\div\xib)_{\ell=1}=0.
\eeaa

\item $\Cb_0$, $\Cb_p$, $M_0$ and $M_p$ are constant on each leaf of the $u$-foliation of  $\Si_*$.

\item We have, for some constant $c_{\Si_*}$,  
\beaa
u+r=c_{\Si_*} , \quad \textrm{along} \quad \Si_*,
\eeaa
where $r$ denotes the area radius of the spheres $S$ of the $u$-foliation of  $\Si_*$. 

\item The following normalization condition  holds true  at the  south pole $SP$  of every sphere $S$,
 \beaa
b_*\Big|_{SP}=-1 -\frac{2m}{r},
 \eeaa
 where $b_*$ is such that we have
 \beaa
 \nu =e_3+b_*e_4,
 \eeaa
 with $\nu$  the unique vectorfield tangent to the hypersurface $\Si_*$, normal to $S$, and normalized by $\g(\nu, e_4)=-2$. 
 
\item The basis of $\ell=1$ modes $\Jp$ is given by $\Jp=J^{(p,S_*)}$ on $S_*$, and extended to $\Si_*$ by $\nu(\Jp)=0$. Also, the $\ell=1$ modes of  $\div\eta$ and $\div\xib$ above are computed with respect to  this basis.
 \end{itemize}

Furthermore, we have\footnote{We have in fact 
\beaa
\max_{k\leq k_{large}+8}\sup_{\Si_*}\Big(\|\dk^kf\|_{L^2(S)}+\|\dk^k\fb\|_{L^2(S)}+\|\dk^k\log(\la)\|_{L^2(S)}\Big) &\les& \ep_0, 
\eeaa
and then use the Sobolev embedding on the 2-spheres $S$ foliating $\Sigma_*$ to deduce \eqref{eq:GCMHoutcomeestimate1}.}
\bea\lab{eq:GCMHoutcomeestimate1}
\max_{k\leq k_{large}+6}\sup_{\Si_*}r\Big(|\dk^kf|+|\dk^k\fb|+|\dk^k\log(\la)|\Big) &\les& \ep_0, 
\eea
and 
\bea\lab{eq:GCMHoutcomeestimate2}
|m - m_0|+\sup_{\Si_*}|r-r_{(0)}| &\les& \ep_0,
\eea
where $(f, \fb, \la)$ are the transition function from the frame of $\Lextt$ to the frame of $\Si_*$. 

\begin{remark}\lab{rmk:calibrationofu}
To fix $u$, we need to pick a specific constant $c_{\Si_*}$ such that $u+r=c_{\Si_*}$ along $\Si_*$. We choose 
$c_{\Si_*}=1+r(S_1)$ where $S_1$ is the only sphere of $\Si_*$ intersecting the curve of the south poles\footnote{Note that this curve is transversal to $\Si_*$ and hence intersect $\Si_*$ at exactly one point given by 
\beaa
\Si_*\cap(\{\ut=1\}\cap\{\widetilde{\th}=\pi\}).
\eeaa} of the outgoing null cone $\{\ut=1\}$ of $\Lextt$.
\end{remark}

{\bf Step 3.} From now on, $u$ is calibrated according to\footnote{Indeed, provided $\delta_0>0$ has been chosen sufficiently small, the spacelike hypersurface $\Si_*$ of Step 2 intersects the curve of the south poles of the spheres foliating  the outgoing cone $\{\ut=1\}$ of the part $\Lextt$ of the initial data layer, which allows to calibrate $u$ as in Remark \ref{rmk:calibrationofu}.} Remark \ref{rmk:calibrationofu}, which also fixes the sphere  $S_1=\Si_*\cap\{u=1\}$. We can then compare $\ug=1+\de_0$ to $u(S_*)$ and obtain 
\beaa
|u(S_*)-1-\de_0|\les \ep_0\de_0,
\eeaa
so that
\beaa
1\leq u\leq u(S_*)\quad \textrm{ on }\Si_*\textrm{ where }1<u(S_*)<1+2\de_0.
\eeaa
Together with the estimate \eqref{eq:GCMHoutcomeestimate2}, and in view of the choice \eqref{eq:choicer0forproofThmM6} for $r_{(0)}$, we have 
\beaa
 r(S_*) &=& r_{(0)}+O(\ep_0)=d_0\de_*\ep_0^{-1}+O(\ep_0)\\
&=& \de_*\ep_0^{-1}(u(S_*))^{1+\dec}\left(d_0+O(\de_0)+O\left(\de_*^{-1}\ep_0^{2}\right)\right).
\eeaa
Thus, we may choose the constant $d_0$ in the range  $\frac{1}{2}\leq d_0\leq 2$ such that
\beaa
 r(S_*)=\de_*\ep_0^{-1}(u(S_*))^{1+\dec}
\eeaa
so that the condition \eqref{eq:behaviorofronS-star} for $r$ is satisfied. 

{\bf Step 4.} In view of Step 1 to Step 3, $\Si_*$ satisfies all the required properties for the future spacelike boundary of a GCM admissible spacetime, see section \ref{sec:admissibleGMCPGdatasetonSigmastar}. We now introduce 
\begin{itemize}
\item  the outgoing geodesic  frame $(\widetilde{e}_4, \widetilde{e}_3, \widetilde{e}_1, \widetilde{e}_2)$ of $\Lextt$, 

\item  the outgoing PG frame $((e_0)_4, (e_0)_3, (e_0)_1, (e_0)_2)$ of $\Lext$, 

\item  the outgoing PG frame  $(e_4', e_3', e_1', e_2')$ initialized on $\Si_*$ from the GCM frame $(e_4, e_3, e_1, e_2)$ by the change of frame with coefficients $(f'', \fb'', \la'')$ given by
\beaa
\la''=1,\qquad f''=\frac{a}{r}f_0, \qquad \fb''=-\frac{(\nu(r)-b_*)}{1-\frac{1}{4}b_*|f''|^2}f'',
\eeaa
where the 1-form $f_0$  is chosen on $\Si_*$ by
\beaa
(f_0)_1=0, \quad (f_0)_2=\sin(\th), \quad\textrm{on}\quad S_*, \qquad \nab_\nu f_0=0\quad\textrm{on}\quad\Si_*,
\eeaa
with $(e_1, e_2)$ specified on $S_*$ by \eqref{eq:canonical-e1ande2onSstar:0}.
\end{itemize}

We have the following change of frame coefficients:
\begin{itemize}
\item $(f, \fb, \la)$, introduced in Step 2, and corresponding to the change from the outgoing geodesic frame $(\widetilde{e}_4, \widetilde{e}_3, \widetilde{e}_1, \widetilde{e}_2)$ of $\Lextt$ to the GCM frame $(e_4, e_3, e_1, e_2)$ of $\Si_*$,

\item  $(f', \fb', \la')$, which we now introduce, corresponding to the change from the  outgoing geodesic frame $(\widetilde{e}_4, \widetilde{e}_3, \widetilde{e}_1, \widetilde{e}_2)$ of $\Lextt$   to   the outgoing PG frame $((e_0)_4, (e_0)_3, (e_0)_1, (e_0)_2)$ of $\Lext$,

\item  $(f'', \fb'', \la'')$, provided explicitly above, and corresponding to the change from the GCM frame $(e_4, e_3, e_1, e_2)$ of $\Si_*$ to the PG frame  $(e_4', e_3', e_1', e_2')$,

\item $(f''', \fb''', \la''')$, which we now introduce, corresponding to the change from  the outgoing PG frame $((e_0)_4, (e_0)_3, (e_0)_1, (e_0)_2)$ of the part $\Lext$ of the initial data layer to the outgoing PG frame  $(e_4', e_3', e_1', e_2')$ initialized on $\Si_*$.
\end{itemize}

In this step, our goal is to control the change of frame coefficients $(f''', \fb''', \la''')$. In view of the above, we have schematically 
\beaa
(f''', \fb''', \la''')=(f'', \fb'', \la'')\circ (f, \fb, \la)\circ (f', \fb', \la')^{-1}
\eeaa
where $(f', \fb', \la')^{-1}$ denote the coefficients corresponding to  the inverse transformation coefficients of the transformation with coefficients $(f', \fb', \la')$. We infer
\beaa
\sup_{\Si_*}\left|\dk_*^k(f''', \fb''', \la'''-1)\right| &\les& \sup_{\Si_*}\left|\dk_*^k(f, \fb, \la-1)\right| +\sup_{\Si_*}\left|\dk_*^k(f''-f', \fb''-\fb', \la''-\la')\right|. 
\eeaa
Together with \eqref{eq:GCMHoutcomeestimate1}, we infer, for $k\leq k_{large}+6$, 
\beaa
\sup_{\Si_*}r\left|\dk_*^k(f''', \fb''', \la'''-1)\right| &\les&  \ep_0+\sup_{\Si_*}r\left|\dk_*^k(f''-f', \fb''-\fb', \la''-\la')\right| 
\eeaa
Together with the explicit formulas above for  $(f'', \fb'', \la'')$, we obtain, for $k\leq k_{large}+6$, 
\beaa
\sup_{\Si_*}r\left|\dk_*^k(f''', \fb''', \la'''-1)\right| &\les&  \ep_0+\sup_{\Si_*}r\left|\dk_*^k\left(f'-\frac{a}{r}f_0, \fb' + \frac{(\nu(r)-b_*)}{1-\frac{1}{4}b_*\frac{a^2}{r^2}|f_0|^2}\frac{a}{r}f_0, \la' -1\right)\right| 
\eeaa
We deduce,  using also the control \eqref{eq:GCMHoutcomeestimate2} of $m-m_0$, for $k\leq k_{large}+6$, 
\beaa
&&\sup_{\Si_*}r\left|\dk_*^k(f''', \fb''', \la'''-1)\right|\\
\nn &\les&  \ep_0+\sup_{\Si_*}r\left|\dk_*^k\left(f'-\frac{a_0}{\rextl}{}^{(\idl)}f_0, \fb' -\frac{a_0\left(1-\frac{2m_0}{\rextl}\right)}{\rextl}{}^{(\idl)}f_0, \la' -1\right)\right| +|a-a_0|\\
\nn &&+\sup_{\Si_*}\left(|\dk_*^k(f_0-{}^{(\idl)}f_0)|+r^{-1}|\dk_*^k(r-\rextl)|+\left|\dk_*^k\left(b_*+1+\frac{2m}{r}\right)\right|+|\dk_*^k(\nu(r)+2)|\right).
\eeaa
Also, since the change of frame coefficients $(f', \fb', \la')$ correspond to the change from the  outgoing geodesic frame $(\widetilde{e}_4, \widetilde{e}_3, \widetilde{e}_1, \widetilde{e}_2)$ of $\Lextt$   to   the outgoing PG frame $((e_0)_4, (e_0)_3, (e_0)_1, (e_0)_2)$ of $\Lext$, we have by the control provided by Proposition \ref{proposition:geodesicfoliationLextt}, for $k\leq k_{large}+7$,
\beaa
\sup_{\Si_*}r\left|\dk_*^k\left(f'-\frac{a_0}{\rextl}{}^{(\idl)}f_0, \fb' -\frac{a_0\left(1-\frac{2m_0}{\rextl}\right)}{\rextl}{}^{(\idl)}f_0, \la' -1\right)\right| &\les& \ep_0.
\eeaa
We deduce,  for $k\leq k_{large}+6$, 
\bea\lab{eq:Step14ThmM7:boundonffblatripleprime:bis:ThmM6}
\nn\sup_{\Si_*}r\left|\dk_*^k(f''', \fb''', \la'''-1)\right| &\les&  \ep_0  +\sup_{\widetilde{\Si}_*}\Bigg(|\dk_*^k(af_0-a_0{}^{(\idl)}f_0)|+r^{-1}|\dk_*^k(r-\rextl)|\\
&&+\left|\dk_*^k\left(b_*+1+\frac{2m}{r}\right)\right|+|\dk_*^k(\nu(r)+2)|\Bigg).
\eea

{\bf Step 5.} In this step, we focus on the control of the terms on the RHS of \eqref{eq:Step14ThmM7:boundonffblatripleprime:bis:ThmM6}. To this end, we first estimate $r-\rextl$ on $\Si_*$. In view of the control of the part $\Lextt$  of the initial data layer, we have, for $k\leq k_{large}+7$,  
\beaa
\sup_{\Lextt}(\rt)^2|\widecheck{\kat}| &\les& \ep_0.
\eeaa
Together with the GCM condition $\kac=0$, we infer, for $k\leq k_{large}+7$,  
\beaa
\sup_{\Si_*}r^2\left|\dk_*^k\left(\kac - \widecheck{\kat}\right)\right| &\les&  \ep_0.
\eeaa
Now, we have 
\beaa
\kac-\widecheck{\kat} &=& \ka- \kat-\frac{2}{r}+\frac{2}{\rt}= \ka- \kat-\frac{2(\rt-r)}{r\rt}
\eeaa
so that
\beaa
r-\rt &=& \frac{r\rt}{2}\left(\ka- \kat -\left(\kac - \widecheck{\kat}\right)\right)
\eeaa
and hence, using the above estimate for $\kac - \widecheck{\kat}$, we have, for $k\leq k_{large}+7$,  
\beaa
\sup_{\Si_*}\left|\dk_*^k\left(r - \rt\right)\right| &\les& \ep_0+\sup_{\Si_*}r^2\left|\dk_*^k\left(\ka - \kat\right)\right|.
\eeaa
Using the change of frame formula for $\ka$, together with  the control \eqref{eq:GCMHoutcomeestimate1}  for $(f, \fb, \la)$ and the part $\Lextt$ of the initial data layer, we deduce, for $k\leq k_{large}+5$,
\bea\lab{eq:controlofwidetilderminusronSigmatildestar:ThmM7:ThmM6:rtilde}
\sup_{\Si_*}\left|\dk_*^k\left(r - \rt\right)\right| &\les& \ep_0.
\eea
Together with the control of $\rt-\rextl$ provided by Proposition \ref{proposition:geodesicfoliationLextt}, we infer, for $k\leq k_{large}+5$,
\bea\lab{eq:controlofwidetilderminusronSigmatildestar:ThmM7:ThmM6}
\sup_{\Si_*}\left|\dk_*^k\left(r - \rextl\right)\right| &\les& \ep_0.
\eea

Next, we control $b_*+1+\frac{2m}{r}$ and $\nu(r)+2$. First, note that we have
\beaa
\nu(r-\rt) &=& \nu(r) -e_3(\rt)-b_*e_4(\rt).
\eeaa
Together with the control of $r-\rt$ in \eqref{eq:controlofwidetilderminusronSigmatildestar:ThmM7:ThmM6:rtilde},  the fact that $\nu$ is tangent to $\Si_*$, the change of frame formulas,  the control \eqref{eq:GCMHoutcomeestimate1} for $(f, \fb, \la)$ and the control of the part $\Lextt$ of the initial data layer, we deduce for $k\leq k_{large}+4$,
 \beaa
\sup_{\Si_*}\left|\dk_*^k\left(\nu(r) +1-\frac{2m_0}{\rt}- b_*\right)\right| &\les&  \ep_0.
\eeaa
Together with the control of $r-\rt$ in \eqref{eq:controlofwidetilderminusronSigmatildestar:ThmM7:ThmM6:rtilde} and the  control \eqref{eq:GCMHoutcomeestimate2} of $m-m_0$, we infer, for $k\leq k_{large}+4$,
 \bea\lab{eq:letskeepthisboundforlaterinThmM6:infactjustbelow}
\sup_{\Si_*}\left|\dk_*^k\left(\nu(r) +1-\frac{2m}{r}- b_*\right)\right| &\les&  \ep_0.
\eea
Since $\nu(u+r)=0$ on $\Si_*$, we infer for $k\leq k_{large}+4$
 \beaa
\sup_{\Si_*}\left|\dk_*^k\left(b_*+\nu(u) -\left(1-\frac{2m}{r}\right)\right)\right| &\les&  \ep_0.
\eeaa
Now, as part of the construction of $\Si_*$, the following  transversality conditions on $\Si_*$ are assumed, see  \eqref{eq:transversalityconditionTheoremGCMHDawei} in Theorem \ref{theorem:constuctionGCMH}, 
\bea\lab{eq:transversalityconditionsonwidetildeSigmastar:ThmM7:ThmM6}
\xi=\om=0, \qquad \etab=-\ze, \qquad e_4(r)=1, \qquad e_4(u)=0.
\eea
We infer
\beaa
\nu(u)=e_3(u)+b_*e_4(u)=e_3(u)
\eeaa
and hence, for $k\leq k_{large}+4$,
 \beaa
\sup_{\Si_*}\left|\dk_*^k\left(b_*+e_3(u) -\left(1-\frac{2m}{r}\right)\right)\right| &\les&  \ep_0.
\eeaa
Also, using again the transversality conditions \eqref{eq:transversalityconditionsonwidetildeSigmastar:ThmM7:ThmM6}, we have
\beaa
\nab(e_3(u))=(\ze-\eta)e_3(u).
\eeaa
We deduce, for $k\leq k_{large}+3$,
 \beaa
\sup_{\Si_*}r\left|\dk_*^k\left(\nab(b_*)-(\ze-\eta)\left(b_*-\left(1-\frac{2m}{r}\right)\right)\right)\right| &\les&  \ep_0.
\eeaa
The control for $\ze$ and $\eta$ inferred from the transformation formula, the control of $(f, \fb, \la)$ and the control of the initial data layer implies, for $k\leq k_{large}+2$,
\beaa
\left\|\dk_*^k\nab(b_*)\right\|_{\hk_1(S)} &\les& \ep_0+r^{-1}\ep_0\left\|\dk_*^k(b_*)\right\|_{\hk_1(S)}.
\eeaa
Also, by our GCM condition on $\Si_*$ for $b_*$, we have
\beaa
b_* \Big|_{SP}=-1-\frac{2m}{r}\quad\textrm{on}\quad\Si_*,
\eeaa
and hence, since $\nu$ is tangent to $\Si_*$, we have 
\beaa
\nu^k\left(b_*+1+\frac{2m}{r}\right) \Big|_{SP}=0\quad\textrm{on}\quad\Si_*.
\eeaa 
Thus, introducing the scalar $h$ on $\Si_*$ given by
\beaa
h &:=& b_*+1+\frac{2m}{r},
\eeaa
we have obtained so far on  $\Si_*$, for $k\leq k_{large}+2$,
\beaa
\left\|\dk_*^k\nab(h)\right\|_{\hk_1(S)} &\les& \ep_0+r^{-1}\ep_0\left\|\dk_*^k(h)\right\|_{\hk_1(S)}
\eeaa
and for any $k$ 
\beaa
\nu^k\left(h\right) \Big|_{SP}=0\quad\textrm{on}\quad\Si_*.
\eeaa 
We decompose $\nu^k(h)=\ov{\nu^k(h)}+(\nu^k(h)-\ov{\nu^k(h)})$ and since $\ov{\nu^k(h)}$ is constant on $S$ we have
\beaa
|\ov{\nu^k(h)}| &=& |\ov{\nu^k(h)}_{|_{SP}}|\les  |(\nu^k(h))_{|_{SP}}|+\|\nu^k(h)-\ov{\nu^k(h)}\|_{L^\infty(S)}\\
&\les& r^{-1}\|\nu^k(h)-\ov{\nu^k(h)}\|_{\hk_2(S)}
\eeaa
where we have used Sobolev and the vanishing of $\nu^k(h)$ at the south pole. Together with Poincar\'e inequality, we infer
\beaa
|\ov{\nu^k(h)}| &\les& \left\|\nab\dk_*^{k}h\right\|_{\hk_1(S)}.
\eeaa
In view of the above, and using again Poincar\'e inequality, we deduce for $k\leq k_{large}+2$,
\beaa
r^{-1}\left\|\dk_*^{k}h\right\|_{\hk_2(S)} &\les& |\ov{\nu^k(h)}|+\left\|\dk_*^{k}\nab h\right\|_{\hk_1(S)} \les \left\|\dk_*^{k}\nab h\right\|_{\hk_1(S)}+\ep_0\\
&\les& r^{-1}\ep_0\left\|\dk_*^{k}h\right\|_{\hk_1(S)}+\ep_0.
\eeaa
For $\ep_0$ small enough, we infer,  for $k\leq k_{large}+2$,
\beaa
r^{-1}\left\|\dk_*^{k}h\right\|_{\hk_2(S)} &\les& \ep_0.
\eeaa
Using Sobolev, and recalling the definition of $h$, we infer,  for $k\leq k_{large}+2$,
\beaa
\sup_{\Si_*}\left|\dk_*^{k}\left(b_*+1+\frac{2m}{r}\right)\right| &\les& \ep_0.
\eeaa
Together with \eqref{eq:letskeepthisboundforlaterinThmM6:infactjustbelow}, we have obtain,  for $k\leq k_{large}+2$,
\bea
\sup_{\Si_*}\left|\dk_*^{k}\left(b_*+1+\frac{2m}{r}, \,\nu(r)+2\right)\right| &\les& \ep_0.
\eea
In view of \eqref{eq:Step14ThmM7:boundonffblatripleprime:bis:ThmM6} and \eqref{eq:controlofwidetilderminusronSigmatildestar:ThmM7:ThmM6}, this yields, for $k\leq k_{large}+2$,
\bea\lab{eq:Step14ThmM7:boundonffblatripleprime:ter:ThmM6}
\sup_{\Si_*}r\left|\dk_*^k(f''', \fb''', \la'''-1)\right| &\les&  \ep_0  +\sup_{\widetilde{\Si}_*}|\dk_*^k(af_0-a_0{}^{(\idl)}f_0)|.
\eea

{\bf Step 6.} In this step, we focus on the control of the terms on the RHS of \eqref{eq:Step14ThmM7:boundonffblatripleprime:ter:ThmM6}. To this end, we first control $a$. We have in view of Proposition \ref{proposition:geodesicfoliationLextt}
\beaa
\sup_{\Lextt(\rt\sim\ep_0^{-1})}(\rt)^4\left|\curlt\bt -\frac{6a_0m_0}{(\rt)^5}\widetilde{J}^{(0)}\right|\les \ep_0.
\eeaa
Since, in view of the transformation formulas, the control of $(f, \fb, \la)$ and the control provided by Proposition \ref{proposition:geodesicfoliationLextt} for the foliation of $\Lextt$, we have
\beaa
\sup_{S_*}\left|\curl\b-\curlt\bt\right| &\les& \frac{\ep_0}{r_{(0)}^5},
\eeaa
we infer
\beaa
\max_{p=0,+,-}(\rt)^5\left|\frac{1}{|S_*|}\int_{S_*}\curl\b J^{(p)} - \frac{1}{|S_*|}\int_{S_*}\frac{6a_0m_0}{(\rt)^5}\widetilde{J}^{(0)}J^{(p)}\right|  &\les& \ep_0.
\eeaa
Recalling from \eqref{eq:GCMconditionsinwidetildeStar:ThmM6} that the following  holds on the sphere $S_*$ of $\Si_*$
\beaa
\frac{1}{|S_*|}\int_{S_*}\curl\b J^{(0)}=\frac{2am}{r^5}, \qquad \frac{1}{|S_*|}\int_{S_*}\curl\b J^{(\pm)}=0,
\eeaa
we infer, using also \eqref{eq:controlofwidetilderminusronSigmatildestar:ThmM7:ThmM6:rtilde} to control $r -\rt$ on $S_*$, 
\beaa
\left|am - 3a_0m_0\frac{1}{|S_*|}\int_{S_*}\widetilde{J}^{(0)}J^{(0)}\right|
+|a_0|m_0\left|\frac{1}{|S_*|}\int_{S_*}\widetilde{J}^{(0)}J^{(\pm)}\right|  &\les& \ep_0.
\eeaa
Together with the control of $m - m_0$ in \eqref{eq:GCMHoutcomeestimate2}, and dividing by $m_0$, we obtain
\beaa 
\left|a - 3a_0\frac{1}{|S_*|}\int_{S_*}\widetilde{J}^{(0)}J^{(0)}\right|
+|a_0|\left|\frac{1}{|S_*|}\int_{S_*}\widetilde{J}^{(0)}J^{(\pm)}\right|  &\les& \ep_0.
\eeaa

Next, in view of Corollary 7.2 of \cite{KS-GCM2} (restated here as  Corollary \ref{Lemma:ComparisonJ-strong}),  there exists a canonical basis of $\ell=1$ modes on $S_*$  in the sense of  Definition 3.10 of \cite{KS-GCM2} (recalled  here  in  Definition \ref{definition:ell=1mpdesonS-intro}), which we denote by $J_0^{(p,S_*)}$, such that 
\beaa
\max_{p=0,+,-}\left\|J_0^{(p,S_*)} -\widetilde{J}^{(p)}\right\|_{\hk_{k_{large}+7}(S_*)} &\les& \ep_0. 
\eeaa
Also, recall that $J^{(p)}=J^{(p, S_*)}$ on $S_*$, where $J^{(p, S_*)}$ is in general another canonical basis of $\ell=1$ modes on $S_*$. In view of Definition \ref{definition:ell=1mpdesonS-intro}, note that the canonical basis of $\ell=1$ modes on $S_*$ are unique modulo isometries of $\SSS^2$, i.e. there exists  $O\in O(3)$ such that 
\bea\lab{eq:theuseoftheisometriyOofSSS2inthechoiceofthecanonicalbasiswidetildeS*:ThmM6}
J^{(p, S_*)}=\sum_{q=0,+,-}O_{pq}J_0^{(q,S_*)}, \qquad p=0,+,-.
\eea

\begin{remark}
In general, we have $O\neq I$ in \eqref{eq:theuseoftheisometriyOofSSS2inthechoiceofthecanonicalbasiswidetildeS*:ThmM6}. In fact, the role of $O$ corresponds in Step 1 to the application of Corollary \ref{Corr:ExistenceGCMS2} which ensures that the following holds on $S_*$ w.r.t. the canonical basis of $\ell=1$ modes $J^{(p, S_*)}$, see \eqref{eq:GCMconditionsinwidetildeStar:ThmM6}, 
\beaa
(\curl\b)_{\ell=1,\pm}=0.
\eeaa
This corresponds to fixing the axis of $S_*$. Note that this condition (and hence the axis of $S_*$) is preserved by multiplying the basis $J^{(p, S_*)}$ by $O=-I$ or by any  $O$ fixing $J^{(0, S_*)}$, so that we may assume in \eqref{eq:theuseoftheisometriyOofSSS2inthechoiceofthecanonicalbasiswidetildeS*:ThmM6} that $O$ satisfies 
\bea\lab{eq:conditiononOpossibletorestrictfreedom:ThmM6}
O_{00}\geq 0, \qquad O_{++}\geq 0, \qquad O_{+-}=0,\qquad O_{--}\geq 0.
\eea
\end{remark}

Since $J^{(p)}=J^{(p, S_*)}$ on $S_*$, we infer
\bea\lab{eq:firstcontroloftildeJminusJonwidetildeSstarwithrotation:ThmM6}
\max_{p=0,+,-}\left\|J^{(p)} -\sum_{q=0,+,-}O_{pq}\widetilde{J}^{(q)}\right\|_{\hk_{k_{large}+7}(S_*)} &\les& \ep_0,
\eea
where $O$ satisfies \eqref{eq:conditiononOpossibletorestrictfreedom:ThmM6}. Plugging \eqref{eq:firstcontroloftildeJminusJonwidetildeSstarwithrotation:ThmM6} in the above, we deduce
\beaa 
\left|a - 3a_0\sum_{q=0,+,-}O_{0q}\frac{1}{|S_*|}\int_{S_*}\widetilde{J}^{(0)}\widetilde{J}^{(q)}\right|
+|a_0|\left|\frac{1}{|S_*|}\sum_{q=0,+,-}O_{\pm q}\int_{S_*}\widetilde{J}^{(0)}\widetilde{J}^{(q)}\right|  &\les& \ep_0.
\eeaa
Now, recall that $\ovS=S(\ug, \sg)$ is the sphere of the foliation of $\Lextt$ which shares the same south pole a $S_*$. Relying on Corollary 5.9 in \cite{KS-GCM1} (see also Proposition \ref{proposition:58-59GCM1} here), we have, for $q=0,+,-$, 
\beaa
&&\left|\int_{S_*}\widetilde{J}^{(0)}\widetilde{J}^{(q)} -\int_{\ovS}\widetilde{J}^{(0)}\widetilde{J}^{(q)}\right|\\
 &\les& r\ep_0\left(\sup_{\widetilde{\RR}}|\dkb^{\leq 1}(\widetilde{J}^{(0)}\widetilde{J}^{(q)})|+r\sup_{\widetilde{\RR}}(|\nab_3(\widetilde{J}^{(0)}\widetilde{J}^{(q)})|+|\nab_4(\widetilde{J}^{(0)}\widetilde{J}^{(q)})|)\right).
\eeaa
Together with the control of $\widetilde{J}^{(p)}$ in $\Lextt$ provided by Proposition  \ref{proposition:geodesicfoliationLextt}, we deduce
\beaa
\left|\int_{S_*}\widetilde{J}^{(0)}\widetilde{J}^{(q)} -\int_{\ovS}\widetilde{J}^{(0)}\widetilde{J}^{(q)}\right| &\les& r\ep_0.
 \eeaa
Using the properties  of $\widetilde{J}^{(p)}$ on the sphere $\ovS$ of $\Lextt$ and the control of $r-\rt$ in  \eqref{eq:controlofwidetilderminusronSigmatildestar:ThmM7:ThmM6:rtilde}, this yields 
\beaa
\left|\frac{1}{|S_*|}\int_{S_*}\widetilde{J}^{(0)}\widetilde{J}^{(q)} - \frac{1}{3}\de_{0q}\right| &\les& \ep_0,
 \eeaa
and hence, plugging in the above, we obtain  
\bea
\left|a - a_0O_{00}\right|+|a_0|\left|O_{+0}\right| +|a_0|\left|O_{-0}\right|  &\les& \ep_0.
\eea
Now, recall that we have either $a_0=0$ or $|a_0|\gg\ep_0$. In particular, we have in view of the above estimate 
\bea\lab{eq:firstestimateforaincasea0equal0:ThmM6}
|a|\les \ep_0\quad\textrm{if}\quad a_0=0.
\eea
In the other case, we have, since $|a_0|\gg \ep_0$, 
\bea\lab{eq:startinginequalitywidetildeaandawithrotationO:ThmM6} 
\left|\frac{a}{a_0} - O_{00}\right|+\left|O_{+0}\right| +\left|O_{-0}\right|  &\les& \ep_0\quad\textrm{if}\quad a_0\neq 0.
\eea
This allows us to control, in the case $a_0=0$, the change of frame coefficients $(f''', \fb''', \la''')$ introduced in Step 4 from  the outgoing PG frame $((e_0)_4, (e_0)_3, (e_0)_1, (e_0)_2)$ of the part $\Lext$ of the initial data layer to the outgoing PG frame  $(e_4', e_3', e_1', e_2')$ initialized on $\Si_*$. Indeed, \eqref{eq:Step14ThmM7:boundonffblatripleprime:ter:ThmM6} and \eqref{eq:firstestimateforaincasea0equal0:ThmM6} yield, for $k\leq k_{large}+2$,
\bea\lab{eq:Step14ThmM7:boundonffblatripleprime:ter:ThmM6:firstcase}
\sup_{\Si_*}r\left|\dk_*^k(f''', \fb''', \la'''-1)\right| &\les&  \ep_0\quad\textrm{if}\quad a_0=0.
\eea

{\bf Step 7.} Next, we focus on controlling the RHS of \eqref{eq:Step14ThmM7:boundonffblatripleprime:ter:ThmM6} in the case $a_0\neq 0$.  Since $O\in O(3)$, we have $\sum_{p}O_{p0}^2=1$, and recalling also that $O_{00}\geq 0$ in view of \eqref{eq:conditiononOpossibletorestrictfreedom:ThmM6}, we infer from \eqref{eq:startinginequalitywidetildeaandawithrotationO:ThmM6}
\bea\lab{eq:controlofwidetildeaminusawhenaislarger:ThmM6}
|a-a_0|\les \ep_0\quad\textrm{if}\quad a_0\neq 0
\eea
and
\beaa
\left|O_{00} - 1\right|+\left|O_{+ 0}\right|+\left|O_{- 0}\right| &\les& \ep_0.
\eeaa
Also, since $O\in O(3)$, we also have
\beaa
0=\sum_pO_{p+}O_{p0}=O_{0+}+O(\ep_0), \qquad 0=\sum_pO_{p-}O_{p0}=O_{0-}+O(\ep_0),
\eeaa
and hence
\beaa
|O_{0+}|+|O_{0-}|\les \ep_0.
\eeaa
Together with the fact that $O_{+-}=0$ and $O_{--}\geq 0$ in view of \eqref{eq:conditiononOpossibletorestrictfreedom:ThmM6}, and since $\sum_{p}O_{p-}^2=1$, we infer
\beaa
\left|O_{--} - 1\right| &\les& \ep_0.
\eeaa
Finally $O_{++}\geq 0$ in view of \eqref{eq:conditiononOpossibletorestrictfreedom:ThmM6}, since we have obtained above that $|O_{0+}|\les \ep_0$, and since $\sum_{p}O_{p-}^2=1$ and $\sum_pO_{p+}O_{p-}=0$, we infer
\beaa
\left|O_{++} - 1\right|+\left|O_{-+}\right| &\les& \ep_0.
\eeaa
We have thus obtained
\beaa
|O-I| &\les& \ep_0,
\eeaa
which together with \eqref{eq:firstcontroloftildeJminusJonwidetildeSstarwithrotation:ThmM6} implies
\beaa
\max_{p=0,+,-}r^{-1}\left\|J^{(p)} - \widetilde{J}^{(q)}\right\|_{\hk_{k_{large}+7}(S_*)} &\les& \ep_0\quad\textrm{if}\quad a_0\neq 0.
\eeaa

Next, we control $\Jp- \widetilde{J}^{(p)}$ for $p=0,+,-$ on $\Si_*$. Recall that we have $\nu(\Jp)=0$ along $\Si_*$. We infer
\beaa
\nu\left(\Jp- \widetilde{J}^{(p)}\right) &=& -\nu\left(\widetilde{J}^{(p)}\right)=-e_3\left(\widetilde{J}^{(p)}\right)-b_*e_4\left(\widetilde{J}^{(p)}\right).
\eeaa
Using the change of frame formulas, and the control \eqref{eq:GCMHoutcomeestimate1} of the change of frame coefficients $(f, \fb, \la)$, and the control of $\widetilde{J}^{(p)}$, we easily obtain, for $k\leq k_{large}+6$,  
\beaa
\sup_{\Si_*}r\left|\dk^k_*\nu\left(\Jp- \widetilde{J}^{(p)}\right)\right| &\les& \ep_0.
\eeaa
Integrating along $\Si_*$  from $S_*$, and using the above control on $S_*$ and Sobolev, as well as the fact that $r\sim \ep_0^{-1}$ on $\Si_*$, we infer, for $k\leq k_{large}+5$,  
\bea
\max_{p=0,+,-}\sup_{\Si_*}\left|\dk_*^k\left(\Jp- \widetilde{J}^{(p)}\right)\right|  &\les& \ep_0\quad\textrm{if}\quad a_0\neq 0. 
\eea

Next, we control $f_0-{}^{(\idl)}f_0$ on $\Si_*$. First,  from the change of frame formulas, the control \eqref{eq:GCMHoutcomeestimate1}  for $(f, \fb, \la)$, and the control of the part $\Lextt$ of the initial data layer,  we have, for $k\leq k_{large}+7$,
\bea\lab{eq:controlofGabGaginTheoremM6asconsequencecontrolIDL}
\sup_{\Si_*}\left(r^2\left|\dk_*^k\Ga_g\right|+r\left|\dk_*^k\Ga_b\right|\right) &\les& \ep_0.
\eea
Proceeding as in Proposition \ref{prop:estimatesforf0fplusfminusandJp}, we infer, for $k\leq k_{large}+7$,
\beaa
\sup_{\Si_*}r^2\left|\dk_*^k\left(\nab J^{(0)} +\frac{1}{r}\dual f_0\right)\right| &\les& \ep_0.
\eeaa
Also, in view of the control of the part $\Lextt$ of the initial data layer,  we have, for $k\leq k_{large}+7$,
\beaa
\sup_{\Lextt}(\rt)^2\left|\dk^k\left(\widetilde{\nab} \widetilde{J}^{(0)} +\frac{1}{\rt}\dual ({}^{(\idl)}f_0)\right)\right| &\les& \ep_0.
\eeaa
Together with the control \eqref{eq:controlofwidetilderminusronSigmatildestar:ThmM7:ThmM6} for $r-\rextl$, the change of frame formulas, the control \eqref{eq:GCMHoutcomeestimate1}  for $(f, \fb, \la)$, and the control of the part $\Lextt$ of the initial data layer, we infer, for $k\leq k_{large}+7$,
\beaa
\sup_{\Si_*}r^2\left|\dk^k\left(\nab \widetilde{J}^{(0)} +\frac{1}{r}\dual ({}^{(\idl)}f_0)\right)\right| &\les& \ep_0.
\eeaa
Subtracting the two estimates, we infer, for $k\leq k_{large}+7$,
\beaa
\sup_{\Si_*}r^2\left|\dk_*^k\left(\nab (J^{(0)}-\widetilde{J}^{(0)}) +\frac{1}{r}\dual (f_0-{}^{(\idl)}f_0)\right)\right| &\les& \ep_0.
\eeaa
Together with the above control for $J^{(0)}- \widetilde{J}^{(0)}$, we deduce, for $k\leq k_{large}+7$,
\bea
\sup_{\Si_*}\left|\dk_*^k\left(f_0-{}^{(\idl)}f_0\right)\right| &\les& \ep_0\quad\textrm{if}\quad a_0\neq 0.
\eea

We are now ready to control, in the case $a_0\neq 0$, the change of frame coefficients $(f''', \fb''', \la''')$ introduced  in Step 4 and  corresponding to the change from    the outgoing PG frame $((e_0)_4, (e_0)_3, (e_0)_1, (e_0)_2)$ of $\Lext$ to the outgoing PG frame  $(e_4', e_3', e_1', e_2')$ initialized on $\Si_*$.  Indeed, \eqref{eq:Step14ThmM7:boundonffblatripleprime:ter:ThmM6}, and the above control of $a-a_0$, $J^{(0)}- \widetilde{J}^{(0)}$ and $f_0-{}^{(\idl)}f_0$ in the case $a_0\neq 0$ yields, for $k\leq k_{large}+5$,
\bea\lab{eq:Step14ThmM7:boundonffblatripleprime:ter:ThmM6:secondcase}
\sup_{\Si_*}r\left|\dk_*^k(f''', \fb''', \la'''-1)\right| &\les&  \ep_0\quad\textrm{if}\quad a_0\neq 0.
\eea

We conclude this step with the control of $\Jk-\Jk_{\idl}$ on $\Si_*$ in the case $a_0\neq 0$. Recall that $\Jk$ is given on $\Si_*$ by 
\beaa
\Jk &=& \frac{1}{|q|}\left(f_0+i\dual f_0\right)\quad\text{on}\quad\Si_*.
\eeaa
Together with the above estimates for $a-a_0$, $f_0-{}^{(\idl)}f_0$, $r-\rt$ and $J^{(0)}- \widetilde{J}^{(0)}$ in the case $a_0\neq 0$, we infer, for $k\leq k_{large}+7$, 
\beaa
\sup_{\Si_*}r\left|\dk_*^k\left(\Jk- \Jk_{\idl}\right)\right| &\les& \ep_0+\sup_{\Lextt}\rt\left|\dk^k\left(\Jk_{\idl}-\frac{1}{|q_{\idl}|}\left({}^{(\idl)}f_0+i\dual({}^{(\idl)}f_0)\right)\right)\right|
\eeaa
which together with the control provided by Proposition \ref{proposition:geodesicfoliationLextt} for the part $\Lextt$ of the initial data layer implies, for $k\leq k_{large}+7$, 
\bea
\sup_{\Si_*}r\left|\dk_*^k\left(\Jk- \Jk_{\idl}\right)\right| &\les& \ep_0\quad\textrm{if}\quad a_0\neq 0.
\eea

Finally, we have obtained in this step, for $k\leq k_{large}+2$, 
\bea\lab{eq:controlffblafinalonSigma*widetildeforThmM7:ThmM6}
 \left|a-a_0\right|+\sup_{\Si_*}\Bigg(r\left|\dk_*^k(f''', \fb''', \la'''-1)\right| +\left|\dk_*^k\left(J^{(0)}- J_{\LL_0}^{(0)}\right)\right|\\
\nn  +r\left|\dk_*^k\left(\Jk- \Jk_{\idl}\right)\right|\Bigg) &\les&  \ep_0\quad\textrm{if}\quad a_0\neq 0,
\eea
where we have also used the control for $\widetilde{J}^{(0)}-J_{\LL_0}^{(0)}$ provided by Proposition \ref{proposition:geodesicfoliationLextt}.

{\bf Step 8.} We now control the outgoing PG structure  initialized on $\Si_*$, and covering the region we denote by $\Mext$, which is included in the initial data layer. For convenience, we change our notation. From   now on:
\begin{itemize}
\item $(e_4, e_3, e_1, e_2)$ denotes the outgoing PG frame of the part $\Lext$ of the initial data layer,

\item $(e_4', e_3', e_1', e_2')$ denotes the outgoing PG frame   of $\Mext$ initialized on $\Si_*$,

\item $(f,\fb, \la)$ denote the transition coefficients from the outgoing PG frame $(e_4, e_3, e_1, e_2)$ to the outgoing PG frame $(e_4', e_3', e_1', e_2')$,

\item $(r, f_0, J^{(0)}, \Jk)$ and $(r', f_0', {J^{(0)}}', \Jk')$ correspond respectively to the outgoing PG structure of $\Lext$ and to the outgoing PG structure of $\Mext$.
\end{itemize}
In view of \eqref{eq:GCMHoutcomeestimate2}, \eqref{eq:controlofwidetilderminusronSigmatildestar:ThmM7:ThmM6}, \eqref{eq:firstestimateforaincasea0equal0:ThmM6}, \eqref{eq:Step14ThmM7:boundonffblatripleprime:ter:ThmM6:firstcase}, \eqref{eq:controlffblafinalonSigma*widetildeforThmM7:ThmM6}, using the above new notations, and noticing that the structure equations in the $e_4'$ direction for the outgoing PG structure initialized on $\Si_*$  allow to recover the $e_4'$ derivatives (which are transversal to $\Si_*$), we have, for $k\leq k_{large}+2$,   
\bea\lab{eq:controlffblafinalonSigma*widetildeforThmM7:ThmM6bis}
\nn \left|m-m_0\right|+\left|a-a_0\right|+\sup_{\Si_*}\Bigg(r\left|\dk_*^k(f, \fb, \la -1)\right| +\left|\dk_*^k\left(r'-r\right)\right|\\
+\left|\dk_*^k\left(a{{J'}^{(0)}}- a_0J^{(0)}\right)\right| +r\left|\dk_*^k\left(a\Jk'- a_0\Jk\right)\right|\Bigg) &\les&  \ep_0.
\eea

We introduce the notations
\beaa
F:=f+i\dual f, \qquad \underline{F}:=\fb+i\dual \fb.
\eeaa
Since $(e_4, e_3, e_1, e_2)$ and $(e_4', e_3', e_1', e_2')$ are outgoing PG frames, we have 
\beaa
\Xi=0, \qquad \om=0, \qquad \Hb+Z=0, \qquad \Xi'=0, \qquad \om'=0, \qquad \Hb'+Z'=0.
\eeaa 
In view of Corollary \ref{cor:transportequationine4forchangeofframecoeffinformFFbandlamba}, we have the following transport equations
\beaa
\nab_{\la^{-1}e_4'}\left(qF\right) &=&  E_4(f, \Ga),\\
\la^{-1}\nab_{e_4'}(\log\la) &=& 2f\c\ze+E_2(f, \Ga),\\
\nab_{\la^{-1}e_4'}\left[q\Big(\underline{F} -2q\DD'(\log\la)+e_3(r) F\Big)\right] &=&  -3q^2\DD'\left(f\c\ze\right) +E_5({\nab'}^{\leq 1}f,\fb, {\nab'}^{\leq 1}\la, \D^{\leq 1}\Ga),
\eeaa
where  $E_2$, $E_4$ and $E_5$ are given in  Corollary \ref{cor:transportequationine4forchangeofframecoeffinformFFbandlamba}. Integrating these transport equations from $\Si_*$ in the order they appear, using the control in \eqref{eq:controlffblafinalonSigma*widetildeforThmM7:ThmM6bis} for $(f, \fb, \la)$ on $\Si_*$, and together with the control  of the part $\Lext$ of the initial data layer, we obtain, for $k\leq k_{large}+2$,    
\begin{equation}\lab{eq:controlffblaonMextwidetildeThmM7ThmM6}
\sup_{\Mext}r\Big(|\dk^k(f,  \log(\la))|+|\dk^{k-1}\fb|\Big) \les \ep_0.
\end{equation}

Also, we have
\beaa
e_4'(r'-\la^{-1}r) &=& 1-\la^{-1}e_4(r)+\la^{-1}e_4'(\log(\la)).
\eeaa
Using the change of frame formula and the above transport equation for $\log(\la)$, we infer
\beaa
e_4(r'-\la^{-1}r) &=& 1-\left(e_4+f\c\nab+\frac{1}{4}|f|^2e_3\right)r+\frac{3}{2}f\c\ze+E_2(f, \Ga)\\ 
&=& -\frac{1}{4}|f|^2e_3(r)+\frac{3}{2}f\c\ze+E_2(f, \Ga).
\eeaa
Integrating from $\Si_*$ where $r'-r$ is under control in view of \eqref{eq:controlffblafinalonSigma*widetildeforThmM7:ThmM6bis}, and using the control \eqref{eq:controlffblaonMextwidetildeThmM7ThmM6} for $f$ and $\la$ as well as the control of $\Lext$, we infer, for $k\leq k_{large}+2$, 
\bea\lab{eq:controlffblaonMextwidetildeThmM7:rtildeThmM6}
\sup_{\Mext}\left|\dk^k(r'-r)\right| &\les& \ep_0. 
\eea
Also, we have
\beaa
e_4'(a{J'}^{(0)}-a_0J^{(0)}) &=& -e_4'(a_0J^{(0)})=-a_0\la\left(e_4+f\c\nab+\frac{1}{4}|f|^2e_3\right)J^{(0)}\\
&=& -a_0\la\left(f\c\nab+\frac{1}{4}|f|^2e_3\right)J^{(0)}
\eeaa
and
\beaa
\nab_4'(aq'\Jk'- a_0q\Jk) &=& -\nab_4'(a_0q\Jk)=-a_0\la\left(\nab_4+f\c\nab+\frac{1}{4}|f|^2\nab_3\right)(q\Jk)\\
&=& -a_0\la\left(f\c\nab+\frac{1}{4}|f|^2\nab_3\right)(q\Jk).
\eeaa
Integrating from $\Si_*$ where ${J'}^{(0)}-J^{(0)}$ and $\Jk' - \Jk$ are under control in view of \eqref{eq:controlffblafinalonSigma*widetildeforThmM7:ThmM6bis}, and using the control \eqref{eq:controlffblaonMextwidetildeThmM7ThmM6} for $f$ and $\la$ as well as the control of $\Lext$, we infer, for $k\leq k_{large}+2$, 
\bea\lab{eq:controlffblaonMextwidetildeThmM7:JandJktildeThmM6}
\sup_{\Mext}\left(\left|\dk^k(a{J'}^{(0)}-a_0J^{(0)})\right|+r\left|\dk^k(a\Jk' -a_0\Jk)\right|\right) &\les& \ep_0. 
\eea

Then, using the outgoing PG structure of $\Mext$, we initialize
\begin{itemize} 
\item the ingoing PG structure  of $\Mint$ on $\TT=\{r'=r_0\}$,

\item  the  ingoing PG structure  of $\,^{(top)}\MM$ on $\{u'=u_*\}$,
\end{itemize} 
as in section \ref{sec:initalizationadmissiblePGstructure}. Using the control of $(f, \fb, \la)$, $r'-r$, ${J'}^{(0)}-J^{(0)}$ and $\Jk' -\Jk$ induced on $\{r'=r_0\}$ and $\{u'=u_*\}$ by \eqref{eq:controlffblaonMextwidetildeThmM7ThmM6}, \eqref{eq:controlffblaonMextwidetildeThmM7:rtildeThmM6} and \eqref{eq:controlffblaonMextwidetildeThmM7:JandJktildeThmM6}, 
and using the analog in the $e_3'$ direction for ingoing PG structures of the above transport equation in the $e_4'$ direction for outgoing PG structures,  we obtain for $\Mint$ and $k\leq k_{large}+1$
\bea\lab{eq:controlffblaonMintwidetildeThmM7ThmM6}
\nn\sup_{\Mint}\Big(|\dk^k(\fb,  \log(\la))|+|\dk^{k-1}f|+|\dk^k(r'-r)|\\
+|\dk^k(a{J'}^{(0)}-a_0J^{(0)})|+|\dk^k(a\Jk'-a_0\Jk)| \Big)&\les& \ep_0,
\eea
and a similar estimate for $\,^{(top)}\MM$. 

Let now
\beaa
\MM &:=& \Mext\cup \Mint\cup \,^{(top)}\MM.
\eeaa
Then, in view of \eqref{eq:controlffblaonMextwidetildeThmM7ThmM6}-\eqref{eq:controlffblaonMintwidetildeThmM7ThmM6}, the control of $a-a_0$ and $m-m_0$ in \eqref{eq:controlffblafinalonSigma*widetildeforThmM7:ThmM6bis}, and using the transformation formulas of Proposition \ref{Proposition:transformationRicci}, and well as the definition of the linearized quantities based on $a$, $m$, $r$, $J^{(0)}=\cos\th$ and $\Jk$, we deduce
\beaa
 \Nk^{(Sup)}_{k_{large}}+ \Nk^{(Dec)}_{k_{small}} &\les&  \ep_0
\eeaa
which concludes the proof of Theorem M6.

%%%%%%%%%%%%%%%%%%%%%%%%%%%%%%%%%%%%%%%%

\section{Proof of Theorem M7}
\lab{sec:proofofThmM7}

%%%%%%%%%%%%%%%%%%%%%%%%%%%%%%%%%%%%%%%%

The proof of Theorem M7 proceeds in 18 steps which we summarize below for convenience:
\begin{enumerate}
\item In Steps 1--5, we use local existence to extend the spacetime $\MM$ a little bit, and then focus on the region in the future of $\Si_*$ in which we derive additional estimates. 

\item In Steps 6--7,  we construct a new last sphere $\widetilde{S}_*$, and a new last slice $\widetilde{\Si}_*$ inside  the region of the extended spacetime in the future of $\Si_*$ by relying on the control derived in Step 1--5 and the GCM constructions of \cite{KS-GCM2} and \cite{Shen} recalled in section \ref{sec:GCMpapersreview}.

\item In Steps 8--11, we show that the new last slice $\widetilde{\Si}_*$, which a priori exists only in a small neighborhood of the new last sphere $\widetilde{S}_*$, extends in fact all the way to the initial data layer. 

\item In Steps 12--13, we complete the proof of the fact that the new last slice $\widetilde{\Si}_*$ satisfies all the required properties for being the future spacelike boundary of a GCM admissible spacetime.

\item In Steps 14--18, we construct from $\widetilde{\Si}_*$ a new GCM admissible spacetime $\widetilde{\MM}$ and  we control the change of frame coefficients between the frames of the extended spacetime, and the corresponding frames of $\widetilde{\MM}$. In view of the control of the extended spacetime and of the change of frames coefficients, we infer the desired control of $\widetilde{\MM}$ thanks to the change of frame formulas. 
\end{enumerate}

We now proceed with the proof of Theorem M7. For convenience, we introduce the following notation 
\bea
k_* &:=& k_{small}+20.
\eea
Then, in view of the assumptions, we are given a GCM admissible spacetime $\MM=\MM(u_*) \in\aleph(u_*)$  verifying   the following improved  bounds  
\bea
\lab{Fullestimates:onMM}
 \Nk^{(Dec)}_{k_*}(\MM)\le  C \ep_0, 
 \eea
for a universal constant $C>0$  provided by Theorems M1-M5.

%%%%%%%%%%%%%

\subsection{Steps 1--5}

%%%%%%%%%%%%%
  
  {\bf Step 1.} We extend $\MM$ by a local existence argument, to a strictly   larger  spacetime  $\MM^{(extend)}$,  with a naturally  extended foliation and the following slightly increased bounds\footnote{The loss of three derivatives occurs due to the fact that local existence holds in $L^2$ based spaces while $\Nk^{(Dec)}_k$ is based on $L^\infty$.}
\beaa
 \Nk^{(Dec)}_{k_*- 3}(\MMextend) \le 2C  \ep_0,
\eeaa
but which may  not verify our admissibility criteria. 

{\bf Step 2.} We then invoke Theorem 4.1  in  \cite{Shen} (restated here in Theorem \ref{theorem:constuctionGCMH}) to extend $\Si_*$ in $\MMextend\setminus\MM$  as a smooth spacelike hypersurface  $\Si^{(extend)}_*$, together with a scalar function $u^{(extend)}$, satisfying the same GCM conditions than $\Si_*$. 

{\bf Step 3.} We consider the outgoing geodesic foliation $(u^{(extend)}, s^{(extend)})$ initialized on $\Si^{(extend)}_*$ to the future of $\Si^{(extend)}_*$ in $\MMextend$. Note in particular that we have from the definition  of $\Si_*$ and $\Si^{(extend)}_*$ 
\beaa
u^{(extend)}+s^{(extend)}=c_{\Si_*}\quad\textrm{on}\quad\Si^{(extend)}_*.
\eeaa
We define the following spacetime region to the future of $\Si^{(extend)}_*$
\beaa
\widetilde{\RR} &:=& \Big\{1\leq u^{(extend)}\leq u_*+\de_{ext}, \quad c_{\Sigma_*}\leq u^{(extend)}+s^{(extend)}\leq c_{\Sigma_*}+\Delta_{ext}\Big\},
\eeaa
where 
\beaa
\Delta_{ext} := \frac{d_0 r_*}{u_*}\de_{ext},\qquad r_*:=r(S_*), \qquad S_*:=\Si_*\cap\CC_*,
\eeaa
with $\de_{ext}>0$  chosen sufficiently small so that $\widetilde{\RR}\subset\MMextend$, and with $d_0$ a constant satisfying 
\beaa
\frac{1}{2}\leq d_0\leq 1
\eeaa 
which will  be suitably chosen in Step 12 below.  From now on, for convenience, we drop the index $(extend)$ and simply denote $u^{(extend)}$ and $s^{(extend)}$ by $u$ and $s$. 

{\bf Step 4.} On $\Si^{(extend)}_*$, the following GCM conditions hold by construction 
\bea\lab{eq:werepeattheGCMcondisitononextendedSigmastarforproofThM8}
\bsplit
&\kac=0, \qquad \kabc=\Cb_0+\sum_p\Cb_p\Jp, \qquad \muc=M_0+\sum_pM_p\Jp, \\
&(\div\eta)_{\ell=1}=0, \qquad (\div\xib)_{\ell=1}=0,
\end{split}
\eea
where the basis of $\ell=1$ modes satisfies in particular 
\bea\lab{eq:propagationogell=1basistoregionextendedSigma}
\nu(\Jp)=0, \quad p=0,+,-, \quad\textrm{along}\quad\Si^{(extend)}_*,
\eea 
and where the scalar functions $\Cb_0$, $\Cb_p$, $M_0$ and $M_p$ are constant on the leaves of the $u$-foliation of $\Si^{(extend)}_*$, i.e. they are functions of $u$ along $\Si^{(extend)}_*$. We propagate $\Jp$, $\Cb_0$, $\Cb_p$, $M_0$ and $M_p$ from $\Si^{(extend)}_*$ to the  spacetime region $\widetilde{\RR}$ along $e_4$ as follows
\bea\lab{eq:propagationogell=1basisandGCMconstansttoregionRRwidetilde}
\bsplit
& e_4(\Jp)=0, \qquad e_4(r^2\Cb_0)=0, \qquad e_4(r^2\Cb_p)=0, \\ 
& e_4(r^3M_0)=0, \qquad e_4(r^3M_p)=0 \quad\textrm{on}\quad\widetilde{\RR},
\end{split}
\eea
so that we have\footnote{More precisely, we have $\Cb_0=r^{-2}\widetilde{\Cb}_0$, $\Cb_p=r^{-2}\widetilde{\Cb}_p$, $M_0=r^{-3}\widetilde{M}_0$, and $M_p=r^{-3}\widetilde{M}_p$, , with $\widetilde{\Cb}_0$, $\widetilde{\Cb}_p$, $\widetilde{M}_0$ and $\widetilde{M}_p$ given by the restriction of $r^2\Cb_0$, $r^2\Cb_p$, $r^3M_0$ and $r^3M_p$ to $\Si^{(extend)}_*$ so that $\widetilde{\Cb}_0=\widetilde{\Cb}_0(u)$, $\widetilde{\Cb}_p=\widetilde{\Cb}_p(u)$, $\widetilde{M}_0=\widetilde{M}_0(u)$ and $\widetilde{M}_p=\widetilde{M}_p(u)$. Note also that $r=r(u,s)$.} $\Cb_0=\Cb_0(u,s)$, $\Cb_p=\Cb_p(u,s)$, $M_0=M_0(u,s)$,  and $M_p=M_p(u,s)$ in $\widetilde{\RR}$. In view of \eqref{eq:propagationogell=1basisandGCMconstansttoregionRRwidetilde}, we have in particular
\beaa
e_4\left(r^2\left(\kabc -\left(\Cb_0+\sum_p\Cb_p\Jp\right)\right)\right) &=& e_4(r^2\kabc),\\
e_4\left(r^3\left(\muc -\left(M_0+\sum_pM_p\Jp\right)\right)\right) &=& e_4(r^3\muc).
\eeaa
Propagating from $\Si^{(extend)}_*$ where  \eqref{eq:werepeattheGCMcondisitononextendedSigmastarforproofThM8} holds, and using the bounds of Step 1 on 
 $\MM^{(extend)}$, and hence on $\widetilde{\RR}$, for $\kac$, $\kabc$, $\muc$, $\eta$ and $\xib$, we obtain, for all $k\leq k_*-4$, 
\bea\lab{eq:controlofthesmallGCMquantitieskackabcmucforThmM7}
\nn\sup_{\widetilde{\RR}}\Bigg(r^2\left|\dk^k\left(\kac\right)\right|+r^2\left|\dk^{k}\left(\kabc -\left(\Cb_0+\sum_p\Cb_p\Jp\right)\right)\right|\\
+r^2\left|\dk^{k}\left(\muc -\left(M_0+\sum_pM_p\Jp\right)\right)\right|\Bigg) &\les& \frac{\ep_0}{r}\Delta_{ext}
\eea
and
\bea\lab{eq:controlofthesmallGCMquantitiesell=1modedivetaanddivxibforThmM7}
\sup_{\widetilde{\RR}}r^2\Big(|(\div\eta)_{\ell=1}|+|(\div\xib)_{\ell=1}|\Big)  &\les& \frac{\ep_0}{r}\Delta_{ext}.
\eea

Next, recall that $\nu=e_3+b_*e_4$ denotes the unique tangent vectorfield to $\Si_*$ which is orthogonal to the $u$-foliation and normalized by $\g(\nu, e_4)=-2$. In view of Corollary \ref{corofLemma:transport.alongSi_*1:again}, we have on $\Si^{(extend)}_*$
\beaa
\left|\nu((\div\b)_{\ell=1})\right|+\left|\nu((\curl\b)_{\ell=1,\pm})\right|+\left|\nu\left((\curl\b)_{\ell=1,0}-\frac{2am}{r^5}\right)\right| &\les& \frac{\ep_0}{r^5u^{1+\dec}},\\
\left|\nu((\kabc)_{\ell=1})\right| &\les& \frac{\ep_0}{r^3u^{1+\dec}}+\frac{\ep_0^2}{r^2u^{2+2\dec}}.
\eeaa
In particular, since $r(S_*)=\de_*\ep_0^{-1}(u(S_*))^{1+\dec}$ in view of \eqref{eq:behaviorofronS-star} and $u(S_*)=u_*$, we infer $r\sim \de_*\ep_0^{-1}u^{1+\dec}$ on $\Si^{(extend)}_*(u_*\leq u\leq u_*+\de_{ext})$ and hence
\beaa
\left|\nu((\div\b)_{\ell=1})\right|+\left|\nu((\curl\b)_{\ell=1,\pm})\right|+\left|\nu\left((\curl\b)_{\ell=1,0}-\frac{2am}{r^5}\right)\right| &\les& \frac{\ep_0}{r^5u^{1+\dec}},\\
\left|\nu((\kabc)_{\ell=1})\right| &\les& \frac{\ep_0}{r^3u^{1+\dec}}.
\eeaa
We integrate from $S_*$ where we have 
\beaa
(\div\b)_{\ell=1}=0, \qquad (\curl\b)_{\ell=1,\pm}=0,\qquad (\curl\b)_{\ell=1, 0}=\frac{2am}{r^5}, \qquad (\kabc)_{\ell=1}=0,  
\eeaa 
and obtain
\beaa
\sup_{\Si^{(extend)}_*(u_*\leq u\leq u_*+\de_{ext})}\Bigg[r^5\left(|(\div\b)_{\ell=1}|+|(\curl\b)_{\ell=1,\pm}|+ \left|(\curl\b)_{\ell=1, 0}-\frac{2am}{r^5}\right|\right)\\
+r^3|(\kabc)_{\ell=1}|\Bigg] &\les& \frac{\ep_0}{u_*}\de_{ext}.
\eeaa 
We now integrate in the $e_4$ direction from $\Si^{(extend)}_*(u_*\leq u\leq u_*+\de_{ext})$ where we have the above estimate as well as $\kac=0$. We obtain 
\bea\lab{eq:controlofthesmallGCMquantitiesell=1modedivbetacurlbetaandkabcforThmM7}
\nn\sup_{\widetilde{\RR}\cap\{u\geq u_*\}}\Bigg[r^5\left(|(\div\b)_{\ell=1}|+|(\curl\b)_{\ell=1,\pm}|+ \left|(\curl\b)_{\ell=1, 0}-\frac{2am}{r^5}\right|\right)\\
\nn +r^3|(\kac)_{\ell=1}|+r^3|(\kabc)_{\ell=1}|\Bigg]  &\les& \frac{\ep_0}{u_*}\de_{ext}+\frac{\ep_0}{r}\Delta_{ext}\\
&\les& \frac{\ep_0}{r}\Delta_{ext}.
\eea
 
Also, one has, since $u+r$ is constant on $\Si^{(extend)}_*$  and $s=r$ on $\Si^{(extend)}_*$
\beaa
0=\nu(u+s)=e_3(u)+b_*e_4(u)+e_3(s)+b_*e_4(s)=e_3(u)+e_3(s)+b_*,
\eeaa
where we used $e_4(s)=1$ and $e_4(u)=0$ for an outgoing geodesic foliation, and hence
\beaa
b_* &=& -e_3(u)-e_3(s)\textrm{ on }\Si_*.
\eeaa
Together with the GCM condition on $b_*$, we infer
\beaa
\big(e_3(u)+e_3(s)\big)\Big|_{SP} &=& 1+\frac{2m}{r}\textrm{ on }\Si_*,
\eeaa 
where $SP$ denotes the south poles of the spheres of the $u$-foliation, i.e. $\th=\pi$.  As above, propagating forward in $e_4$, using that $\th$ is extended from $\Si^{(extend)}_*$ by $e_4(\th)=0$ so that the integral curve of $e_4$ starting on SP stays on SP, and using the bounds of Step 1 on 
 $\MM^{(extend)}$, and hence on $\widetilde{\RR}$, for $\widecheck{e_3(u)}$ and $\widecheck{e_3(s)}$, we infer
\bea\lab{eq:controlofthesmallGCMquantitiesbsmallatSPforThmM7}
\sup_{\widetilde{\RR}}\left|\big(e_3(u)+e_3(s)\big)\Big|_{SP} -\left(1+\frac{2m}{r}\right)\right| &\les& \frac{\ep_0}{r}\Delta_{ext}.
\eea

Also, arguing as we did above on $\Si^{(extend)}_*(u_*\leq u\leq u_*+\de_{ext})$, we have $r\sim \de_*\ep_0^{-1}u^{1+\dec}$ on $\widetilde{\RR}\cap\{u\geq u_*\}$ and hence, for all $k\leq k_*-3$,
\bea\lab{eq:A1strongindeedholds:ThmM7}
\sup_{\widetilde{\RR}\cap\{u\geq u_*\}}r^2|\dk^k\Ga_b| &\les& \sup_{\widetilde{\RR}\cap\{u\geq u_*\}}\left(\frac{r\ep_0}{u^{1+\dec}}\right) \les \de_*,
\eea
which corresponds, for $\de_*>0$ small  enough, to assumption {\bf A1-strong} in \cite{KS-GCM2}, see \eqref{eq:assumtioninRRforGagandGabofbackgroundfoliation'}. 

Finally, we consider the control of the Hawking mass $m_H$ of the sphere $S(u,s)$ of the $(u,s)$ foliation in the region $\widetilde{\RR}\cap\{u\geq u_*\}$, where we recall that $m_H$ is given by the formula
\beaa
\frac{2m_H}{r} &=& 1+\frac{1}{16\pi}\int_S\ka\kab.
\eeaa
First, we have from Lemma \ref{lemma:transportequationforHawkingmass} on $\Si^{(extend)}_*$
\beaa
\nu(m_H) &=& r^2\dkb^{\leq 1}(\Ga_b\c\Ga_b)
\eeaa
and hence, since $m$ is a constant, and using the bounds of Step 1 on 
 $\MM^{(extend)}$, 
\beaa
\sup_{\Si^{(extend)}_*}u^{2+2\dec}|\nu(m_H-m)| &\les& \ep_0^2. 
\eeaa
Since $m_H=m$ on $S_*$ by definition of $m$, we infer, propagating the above transport equation from $S_*$, 
\beaa
\sup_{\Si^{(extend)}_*(u_*\leq u\leq u_*+\de_{ext})}|m_H-m| &\les& \frac{\ep_0^2\de_{ext}}{u_*^{2+2\dec}}. 
\eeaa
Next, recall from the proof of Lemma \ref{lemma:transportequationforHawkingmass} the following computation
\beaa
e_4(\ka\kab)+\ka^2\kab &=& 2\ka\rho -2\ka\div\ze  +\ka\left( 2|\ze|^2 -\chih\c\chibh\right)-\kab|\chih|^2\\
&=& 2\ka\rho -2\ka\div\ze+r^{-1}\Ga_b\c\Ga_g.
\eeaa
This yields, using a well-known identity for the $e_4$ derivative of the integral on $S$ of a scalar function in an outgoing geodesic foliation, 
\beaa
e_4\left(\int_S\ka\kab\right) &=& \int_S\Big(e_4(\ka\kab)+\ka^2\kab\Big)\\
&=& \int_S\Big(2\ka\rho -2\ka\div\ze+r^{-1}\Ga_b\c\Ga_g\Big)
\eeaa
and hence, using integration by parts,
\beaa
e_4\left(\int_S\ka\kab\right) &=& \int_S\Big(2\ka\rho +r^{-1}\Ga_b\c\dk^{\leq 1}\Ga_g\Big).
\eeaa
Together with the definition of $m_H$, we infer
\beaa
2e_4(m_H) &=& \frac{2m_H}{r}e_4(r)+\frac{r}{16\pi}e_4\left(\int_S\ka\kab\right)\\
&=& m_H \ov{\ka}+\frac{r}{8\pi}\int_S\Big(\ka\rho +r^{-1}\Ga_b\c\dk^{\leq 1}\Ga_g\Big).
\eeaa
Now, using in particular Gauss equation and Gauss-Bonnet formula, we have
\beaa
\int_S\ka\rho &=& \ov{\ka}\int_S\rho+\int_S(\ka-\ov{\ka})(\rho-\ov{\rho})=\ov{\ka}\int_S\rho+\int_Sr^{-1}\Ga_g\c\Ga_g\\
&=& \ov{\ka}\int_S\left(-K-\frac{1}{4}\ka\kab+\frac{1}{2}\chih\c\chibh\right)+\int_Sr^{-1}\Ga_g\c\Ga_g\\
&=&  \ov{\ka}\left(-4\pi-\frac{1}{4}\int_S\ka\kab\right) +\int_Sr^{-1}\Ga_g\c\Ga_g.
\eeaa
In view of the above, and using again the definition of $m_H$, we infer 
\beaa
e_4(m_H) &=& r^2\Ga_b\c\dk^{\leq 1}\Ga_g.
\eeaa
Hence, since $m$ is a constant, and using the bounds of Step 1 on 
 $\MM^{(extend)}$, we have
\beaa
\sup_{\widetilde{\RR}\cap\{u\geq u_*\}}ru^{\frac{3}{2}+2\dec}|e_4(m_H-m)| &\les& \ep_0^2. 
\eeaa
Integrating from $\Si^{(extend)}_*(u_*\leq u\leq u_*+\de_{ext})$ where we control $m_H-m$ in view of the above, we infer
\bea\lab{eq:controlofmHminusmonwidetildeRRmoreustar:ThmM7}
\sup_{\widetilde{\RR}\cap\{u\geq u_*\}}|m_H-m| &\les& \frac{\ep_0^2}{ru_*^{\frac{3}{2}+2\dec}}\De_{ext}+\frac{\ep_0^2\de_{ext}}{u_*^{2+2\dec}} \les \frac{\ep_0^2}{ru_*}\De_{ext}.
\eea

{\bf Step 5.} In this step, we control the basis of $\ell=1$ mode $\Jp$ in the spacetime region $\widetilde{\RR}$. Recall that $\Jp$ is chosen on $S_*$ to be a canonical basis of $\ell=1$ modes in the sense of  Definition 3.10 of \cite{KS-GCM2} (recalled  here  in  Definition \ref{definition:ell=1mpdesonS-intro}), i.e.  on $S_*$ there   exist  coordinates $(\th, \vphi)$    such that:
 \begin{enumerate}
 \item The induced metric $g$ on $S_*$  takes the form
 \bea\lab{eq:formofthemetriconSstarusinguniformization:againThM7}
 g= r^2e^{2\phi}\Big( (d\th)^2+ \sin^2 \th (d\vphi)^2\Big).
 \eea
 
 \item The  functions 
 \bea
 J^{(0)} :=\cos\th, \qquad J^{(-)} :=\sin\th\sin\vphi, \qquad  J^{(+)} :=\sin\th\cos\vphi,
 \eea
 verify  the balanced  conditions 
 \bea\lab{eq:balancedconditionforJponSstar:againThM7}
 \int_{S_*}  J^{(p)} =0, \qquad p=0,+,-.
 \eea
 \end{enumerate}
Recall also that we   extend $(\th, \vphi)$ and $\Jp$   to $\Si^{(extend)}_*$ by setting
  \bea\lab{eq:canonical-ell=1modesonSi:againThM7}
  \nu(\th)=0, \qquad \nu(\vphi)=0, \qquad \nu(\phi)=0, \qquad\nu(\Jp)=0, \quad p=0,+,-,
  \eea
  where we have also extended the conformal fact $\phi$. 
We then extend $(\th, \vphi)$, $\phi$ and $\Jp$   to $\widetilde{\RR}$ as follows
  \bea\lab{eq:canonical-ell=1modesonSi:againThM7:transcond}
  e_4(\th)=0, \qquad e_4(\vphi)=0, \qquad e_4(\phi)=0, \qquad e_4(\Jp)=0, \quad p=0,+,-.
  \eea
In what follows, to avoid  the singularities at $\th=0$ and $\th=\pi$ of the $(\th, \vphi)$ coordinates system on the spheres $S=S(u,s)$ of the outgoing geodesic $(u,s,)$-foliation of $\widetilde{\RR}$, we use instead  two regular coordinates charts based on $(\th, \vphi)$.
\begin{definition}
We define two coordinates charts on $S(u,s)$ as follows
\begin{enumerate}
\item The coordinates $(x^1_N, x^2_N)$ are defined for $0\leq\th<\pi$ by
\beaa
x^1_N:=\frac{\sin\th\cos\vphi}{1+\cos\th}, \qquad x^2_N:=\frac{\sin\th\sin\vphi}{1+\cos\th}.
\eeaa

\item The coordinates $(x^1_S, x^2_S)$ are defined for $0<\th\leq\pi$ by
\beaa
x^1_S:=\frac{\sin\th\cos\vphi}{1-\cos\th}, \qquad x^2_S:=\frac{\sin\th\sin\vphi}{1-\cos\th}.
\eeaa
\end{enumerate}
\end{definition}

\begin{lemma}
Let $g$ denote the metric induced by $\g$ on $S(u,s)$. Then,  on $S_*$, the metric $g$ takes the following form in  the $(x_N^1, x_N^2)$ coordinates system and in the $(x_S^1, x_S^2)$ coordinates system, for $(x^1, x^2)\in\RRR^2$, 
\bea\lab{eq:formofthemetriconSstarusinguniformization:againThM7:regular}
g &=& \frac{4r^2e^{2\phi}}{\big(1+(x^1)^2+(x^2)^2\big)^2}\Big[(dx^1)^2+(dx^2)^2\Big].
\eea
\end{lemma}

\begin{proof}
This follows immediately from \eqref{eq:formofthemetriconSstarusinguniformization:againThM7} and the definition of $(x^1_N, x^2_N)$ and $(x^1_S, x^2_S)$ in terms of $(\th, \vphi)$. 
\end{proof}

From now on, let  $(x^1, x^2)$ denote either $(x_N^1, x_N^2)$ or $(x_S^1, x_S^2)$. In view of  the definition of $(x^1, x^2)$, we have $\nu(x^1)=\nu(x^2)=0$ on $\Si^{(extend)}_*$. Since $\nu=e_3+b_*e_4$ and $\nu(u)=e_3(u)$, we infer $\pr_u=\frac{1}{e_3(u)}\nu$. We easily derive the following formula on $\Si^{(extend)}_*$ in the $(x^1, x^2)$ coordinates system of $S$ 
\beaa
\pr_u g_{ab} &=& 2\chi\left(\frac{\pr}{\pr_a}, \frac{\pr}{\pr_b}\right)+2b_*\chi\left(\frac{\pr}{\pr_a}, \frac{\pr}{\pr_b}\right), 
\eeaa 
and hence
\beaa
\pr_u g_{ab} &=& (\kab+b_*\ka)g_{ab}+2\chibh\left(\frac{\pr}{\pr_a}, \frac{\pr}{\pr_b}\right)+2b_*\chih\left(\frac{\pr}{\pr_a}, \frac{\pr}{\pr_b}\right),
\eeaa
which we rewrite as follows, recalling in particular that $\kac=0$ on $\Si^{(extend)}_*$, and using $\nu(x^a)=0$ and $\nu(\phi)=0$, 
\beaa
\pr_u\left(r^{-2}g_{ab}-\frac{4e^{2\phi}}{1+(x^1)^2+(x^2)^2}\de_{ab}\right) &=& \left(\kabc+\frac{2}{r}\left(1-\frac{1}{e_3(u)}\right)\widecheck{b_*}-\frac{2}{re_3(u)}\widecheck{e_3(u)}\right)r^{-2}g_{ab}\\
&&+2r^{-2}\chibh\left(\frac{\pr}{\pr_a}, \frac{\pr}{\pr_b}\right)+2r^{-2}b_*\chih\left(\frac{\pr}{\pr_a}, \frac{\pr}{\pr_b}\right).
\eeaa
Together with the control of Step 1 on  $\MM^{(extend)}$, and hence on $\Si^{(extend)}_*$, for $\kabc$, $\widecheck{b_*}$, $\widecheck{e_3(u)}$, $\chibh$ and $\chih$, we infer, for all $k\leq k_*-3$,
\beaa
\left|\dk^k\pr_u\left(r^{-2}g_{ab}-\frac{4e^{2\phi}}{1+(x^1)^2+(x^2)^2}\de_{ab}\right)\right| &\les&  \frac{\ep_0}{ru^{1+\dec}}r^{-2}|g|. 
\eeaa
Integrating from $S_*$, where \eqref{eq:formofthemetriconSstarusinguniformization:againThM7:regular} holds, we infer on $\Si^{(extend)}_*(u_*\leq u\leq u_*+\de_{ext})$, for all $k\leq k_*-3$,
\beaa
\left|\dk^k\left(r^{-2}g_{ab}-\frac{4e^{2\phi}}{1+(x^1)^2+(x^2)^2}\de_{ab}\right)\right| &\les& \frac{\ep_0}{ru_*^{1+\dec}}\de_{ext}.
\eeaa

Next, we estimate $r^{-2}g_{ab}$ in $\widetilde{\RR}\cap\{u\geq u_*\}$. In view of the definition of $(x^1, x^2)$, we have $e_4(x^1)=e_4(x^2)=0$. Since we also have $e_4(u)=0$ and $e_4(s)=1$, we infer $e_4=\pr_s$. We easily derive the following formula on $\widetilde{\RR}$ in the $(x^1, x^2)$ coordinates system of $S$ 
\beaa
\pr_s g_{ab} &=& 2\chi\left(\frac{\pr}{\pr_a}, \frac{\pr}{\pr_b}\right),
\eeaa 
and hence
\beaa
\pr_s g_{ab} &=& \ka g_{ab}+2\chih\left(\frac{\pr}{\pr_a}, \frac{\pr}{\pr_b}\right),
\eeaa
which we rewrite as follows, using $e_4(x^a)=0$ and $e_4(\phi)=0$, 
\beaa
\pr_s\left(r^{-2}g_{ab}-\frac{4e^{2\phi}}{1+(x^1)^2+(x^2)^2}\de_{ab}\right) &=& \left(\kac-\ov{\kac}\right)r^{-2}g_{ab}+2r^{-2}\chih\left(\frac{\pr}{\pr_a}, \frac{\pr}{\pr_b}\right).
\eeaa
Together with the control of Step 1 on  $\MM^{(extend)}$, and hence on $\widetilde{\RR}$, for $\kac$ and $\chih$, we infer, for all $k\leq k_*-3$,
\beaa
\left|\dk^k\pr_s\left(r^{-2}g_{ab}-\frac{4e^{2\phi}}{1+(x^1)^2+(x^2)^2}\de_{ab}\right)\right| &\les&  \frac{\ep_0}{r^2u^{\frac{1}{2}+\dec}}r^{-2}|g|. 
\eeaa
 Integrating from $\Si^{(extend)}_*(u_*\leq u\leq u_*+\de_{ext})$, and using the above control of $r^{-2}g_{ab}$ on $\Si^{(extend)}_*(u_*\leq u\leq u_*+\de_{ext})$, we infer in $\widetilde{\RR}\cap\{u\geq u_*\}$, for all $k\leq k_*-3$,
\bea\lab{eq:estimateformetricclosetounifforJnearcanonical1}
\left|\dk^k\left(r^{-2}g_{ab}-\frac{4e^{2\phi}}{1+(x^1)^2+(x^2)^2}\de_{ab}\right)\right| \les \frac{\ep_0}{ru_*^{1+\dec}}\de_{ext}+\frac{\ep_0}{r^2u_*^{\frac{1}{2}+\dec}}\De_{ext}\les \frac{\ep_0}{r^2}\Delta_{ext}.
\eea

Next, we estimate $\phi$ in $\widetilde{\RR}\cap\{u\geq u_*\}$. First, recall from Corollary \ref{cor:statementeq:DeJp.S_*:improved} that we have obtained on $S_*$
\beaa
\|\dkb^{\leq k_*}\phi\|_{L^\infty(S_*)} &\les& \frac{\ep_0}{ru_*^{\frac{1}{2}+\dec}}.
\eeaa
Since we have extended $\phi$ to $\Si^{(extend)}_*$ by $\nu(\phi)=0$ and then to $\widetilde{\RR}$ by $e_4(\phi)=0$, we easily infer  in $\widetilde{\RR}\cap\{u\geq u_*\}$, for all $k\leq k_*-3$,
\bea\lab{eq:estimateformetricclosetounifforJnearcanonical2}
|\dk^k\phi| &\les& \frac{\ep_0}{ru_*^{\frac{1}{2}+\dec}}.
\eea

Next, we estimate $\int_S\Jp$ for $p=0,+,-$ in $\widetilde{\RR}\cap\{u\geq u_*\}$. Recall that $\Jp$ is balanced on $S_*$, i.e.
\beaa
\int_{S_*}\Jp=0, \quad p=0,+,-.
\eeaa
Also, since $\nu(\Jp)=0$ on $\Si^{(extend)}_*$, we have, in view of Corollary \ref{Corr:nuSof integrals},
 \beaa
\nu\left(r^{-2}\int_S\Jp\right) =  \Ga_b, \quad\textrm{for}\quad S\subset\Si^{(extend)}_*,\quad p=0,+,-,
\eeaa
and since $\Jp$ is extended to $\widetilde{\RR}$ by $e_4(\Jp)=0$, we have
 \beaa
e_4\left(r^{-2}\int_S\Jp\right) =  \Ga_g, \quad\textrm{for}\quad S\subset\widetilde{\RR},\quad p=0,+,-.
\eeaa
We deduce in $\widetilde{\RR}\cap\{u\geq u_*\}$
\beaa
r^{-2}\left|\int_S\Jp\right| &\les& \frac{\ep_0}{ru_*^{1+\dec}}\de_{ext}+\frac{\ep_0}{r^2u_*^{\frac{1}{2}+\dec}}\De_{ext}\les \frac{\ep_0}{r^2}\Delta_{ext}
\eeaa
and hence
\bea\lab{eq:estimateformetricclosetounifforJnearcanonical3}
r^{-2}\left|\int_S\Jp\right| &\les&  \frac{\ep_0}{r^2}\Delta_{ext}\quad\textrm{for}\quad S\subset\widetilde{\RR}\cap\{u\geq u_*\},\quad p=0,+,-.
\eea

In view of \eqref{eq:estimateformetricclosetounifforJnearcanonical1}, \eqref{eq:estimateformetricclosetounifforJnearcanonical2} and \eqref{eq:estimateformetricclosetounifforJnearcanonical3}, we may apply  Proposition 4.15 in \cite{KS-GCM2}  (restated here as Proposition \ref{prop:asexpectedthebasiswidetildeJpisclosetocanonicalbasis})
which yields on any sphere $S$ of $\widetilde{\RR}\cap\{u\geq u_*\}$ the existence of a canonical basis $J^{(p,S)}$ of $\ell=1$ modes such that 
\bea\lab{eq:estimateformetricclosetounifforJnearcanonical4}
\max_{p=0,+,-}r^{-1}\|\Jp-J^{(p,S)}\|_{\hk_{k_*-2}(S)} \les  \frac{\ep_0}{r^2}\Delta_{ext},  \,\,\,\textrm{for}\,\,\, S\subset\widetilde{\RR}\cap\{u\geq u_*\}.
\eea
This, together with \eqref{eq:lastandeasiestparttocheckforA4strong:ThmM7} below,  corresponds to assumption {\bf A4-strong} which is stated at the beginning of section \ref{section:RecallresultsGCM}.

Next, we estimate $e_3(\Jp)$ in  the region $\widetilde{\RR}$. Since $\nu(\Jp)=0$ on $\Si^{(extend)}_*$, since $\nu=e_3+b_*e_4$ and  since $\Jp$ is extended from $\Si^{(extend)}_*$ by $e_4(\Jp)=0$, we have
\beaa
e_3(\Jp)=0\quad\textrm{on}\quad\Si^{(extend)}_*.
\eeaa
Then, we compute, using $e_4(\Jp)=0$, $\om=0$ and $\etab=-\ze$,  
\beaa
e_4(e_3(\Jp)) &=& [e_4, e_3]\Jp=\left(2\om e_3 -2\omb e_4+2(\etab-\eta)\c\nab\right)\Jp=-2(\eta+\ze)\c\nab\Jp.
\eeaa
Together with the control of Step 1 on  $\MM^{(extend)}$, and hence on $\widetilde{\RR}$, for $\eta$ and $\ze$, we infer on $\widetilde{\RR}$, for all $k\leq k_*-3$,
\beaa
|\dk^ke_4(e_3(\Jp))| &\les& \frac{\ep_0}{r^2u^{1+\dec}}.
\eeaa
Since $e_3(\Jp)=0$ on $\Si^{(extend)}_*$, we infer, for all $k\leq k_*-3$,
\bea\lab{eq:controlofthesmalle3JpforThmM7}
|\dk^ke_3(\Jp)| &\les& \frac{\ep_0}{r^2u^{1+\dec}}\De_{ext}\quad\textrm{on}\quad\widetilde{\RR}.
\eea

Finally, arguing as in Corollary \ref{corr:Si*-ell=1modes-improved}, the following holds on $\Si^{(extend)}_*$
 \beaa
\bsplit
  \Big(r^2\De+2\Big) \Jp  &= O(\ep_0r^{-1}u^{-\frac{1}{2}-\dec}),\qquad p=0,+,-,\\
\frac{1}{|S|} \int_{S}  \Jp J^{(q)} &=  \frac{1}{3}\de_{pq} +O(\ep_0r^{-1}u^{-\frac{1}{2}-\dec}),\qquad p,q=0,+,-,\\
\frac{1}{|S|}  \int_{S}   \Jp   &=O(\ep_0r^{-1}u^{-\frac{1}{2}-\dec}),\qquad p=0,+,-.
\end{split}
\eeaa
Since $\Jp$ is extended to $\widetilde{\RR}$ by $e_4(\Jp)=0$, we propagate from $\Si^{(extend)}_*$ and easily obtain on $\widetilde{\RR}$
 \bea\lab{eq:lastandeasiestparttocheckforA4strong:ThmM7}
\bsplit
  \Big(r^2\De+2\Big) \Jp  &= O(\ep_0r^{-1}u^{-\frac{1}{2}-\dec}),\qquad p=0,+,-,\\
\frac{1}{|S|} \int_{S}  \Jp J^{(q)} &=  \frac{1}{3}\de_{pq} +O(\ep_0r^{-1}u^{-\frac{1}{2}-\dec}),\qquad p,q=0,+,-,\\
\frac{1}{|S|}  \int_{S}   \Jp   &=O(\ep_0r^{-1}u^{-\frac{1}{2}-\dec}),\qquad p=0,+,-.
\end{split}
\eea

 %%%%%%%%%%%%%

\subsection{Steps 6--13}

%%%%%%%%%%%%%

{\bf Step 6.} We fix the following sphere of the $(u^{(extend)}, s^{(extend)})$ foliation in $\widetilde{\RR}\cap\{u\geq u_*\}$
\bea\lab{eq:checkthatthereisindeedalargeruinnewlastGCMS:5}
\ovS:=S(\ug, \sg), \qquad \ug:=u_*+\frac{\de_{ext}}{2}, \qquad \sg:=r_*+\frac{3d_0r_*}{4u_*}\de_{ext}.
\eea
Define
\beaa
\dg:=\frac{\ep_0}{r}\Delta_{ext}=\frac{d_0\ep_0\de_{ext}}{u_*}, \qquad \epg:= \ep_0,
\eeaa
and the small spacetime neighborhood of $\ovS$
\beaa
 \RR(\epg, \dg) :=\left\{|u-\ug|\leq\de_{\RR} ,\quad |s-\sg|\leq  \de_\RR \right\}, \qquad \de_\RR = \dg \big(\epg\big)^{-\frac{1}{2}}.
\eeaa
Note that $\RR(\epg, \dg)\subset \widetilde{\RR}$. In view of \eqref{eq:controlofthesmallGCMquantitieskackabcmucforThmM7}, \eqref{eq:controlofthesmallGCMquantitiesell=1modedivbetacurlbetaandkabcforThmM7}, \eqref{eq:A1strongindeedholds:ThmM7}, \eqref{eq:estimateformetricclosetounifforJnearcanonical4} and \eqref{eq:lastandeasiestparttocheckforA4strong:ThmM7}, we are in position to apply Theorem 7.3 and Corollary 7.7 of \cite{KS-GCM2} (restated here in Theorem \ref{theorem:ExistenceGCMS2} and Corollary \ref{Corr:ExistenceGCMS2}), with $s_{max}=k_{small}+k_*-4$, which yields the existence of a unique sphere $\widetilde{S}_*$, which is a deformation of $\ovS$, is included in $\RR(\epg, \dg)$, and is such that the following GCM conditions hold on it
 \bea\lab{eq:GCMconditionsinwidetildeStar:ThmM7}
 \bsplit
& \widecheck{\widetilde{\ka}}=0, \qquad  \widecheck{\widetilde{\kab}}=0, \qquad  \widecheck{\widetilde{\mu}}=\sum_p\widetilde{M}_pJ^{(p,\widetilde{S}_*)}, \\
& (\widetilde{\div}\widetilde{\b})_{\ell=1}=0, \qquad  (\widetilde{\curl}\widetilde{\b})_{\ell=1,\pm}=0, \qquad  (\widetilde{\curl}\widetilde{\b})_{\ell=1,0}=\frac{2\widetilde{a}\widetilde{m}}{\widetilde{r}^5}, 
\end{split}
\eea
where
\begin{itemize}
\item the tilde refer to the quantities and tangential operators on $\widetilde{S}_*$,

\item $J^{(p,\widetilde{S}_*)}$ denotes the canonical basis of $\ell=1$ mode on $\widetilde{S}_*$ in the sense of  Definition 3.10 of \cite{KS-GCM2} (recalled  here  in  Definition \ref{definition:ell=1mpdesonS-intro}),

\item the $\ell=1$ modes in \eqref{eq:GCMconditionsinwidetildeStar:ThmM7} are defined w.r.t. the basis of $\ell=1$ modes $J^{(p,\widetilde{S}_*)}$, 

\item $\widetilde{m}$ denotes the Hawking mass of $\widetilde{S}_*$, $\widetilde{r}$ denotes the area radius of $\widetilde{S}_*$, and the identity for $(\widetilde{\curl}\widetilde{\b})_{\ell=1,0}$ in \eqref{eq:GCMconditionsinwidetildeStar:ThmM7} should be understood as providing the definition of $\widetilde{a}$. 
\end{itemize}

{\bf Step 7.} Starting from $\widetilde{S}_*$ constructed in Step 6, and in view of \eqref{eq:controlofthesmallGCMquantitiesell=1modedivetaanddivxibforThmM7}, \eqref{eq:controlofthesmallGCMquantitiesbsmallatSPforThmM7} and \eqref{eq:controlofthesmalle3JpforThmM7}, we may apply Theorem 4.1 in  \cite{Shen} (restated here in Theorem \ref{theorem:constuctionGCMH}),  with $s_{max}=k_{small}+k_*-4$, which yields the existence of a smooth small piece of spacelike hypersurface $\widetilde{\Si}_*$ starting from $\widetilde{S}_*$ towards the initial data layer, together with a scalar function $\widetilde{u}$ defined on $\widetilde{\Si}_*$, whose level surfaces are topological spheres denoted by $\widetilde{S}$, so that 
\begin{itemize}
\item The following GCM conditions are verified on $\widetilde{\Si}_*$
 \beaa
&&\widecheck{\widetilde{\ka}}=0, \qquad \widecheck{\widetilde{\kab}}=\widetilde{\Cb}_0+\sum_p\widetilde{\Cb}_p\widetilde{J}^{(p)}, \qquad \widecheck{\widetilde{\mu}}=\widetilde{M}_0+\sum_p\widetilde{M}_p\widetilde{J}^{(p)},\\
&& (\widetilde{\div}\widetilde{\eta})_{\ell=1}=0, \qquad (\widetilde{\div}\widetilde{\xib})_{\ell=1}=0,
\eeaa
where the tilde refer to the quantities and tangential operators on $\widetilde{\Si}_*$. 

\item $\widetilde{\Cb}_0$, $\widetilde{\Cb}_p$, $\widetilde{M}_0$ and $\widetilde{M}_p$ are constant on each leaf of the $\widetilde{u}$-foliation of  $\widetilde{\Si}_*$.

\item We have, for some constant $c_{\widetilde{\Si}_*}$,  
\beaa
\widetilde{u}+\widetilde{r}=c_{\widetilde{\Si}_*} , \quad \textrm{along} \quad \widetilde{\Si}_*,
\eeaa
where $\widetilde{r}$ denotes the area radius of the spheres $\widetilde{S}$ of the $\widetilde{u}$-foliation of  $\widetilde{\Si}_*$. 

\item The following normalization condition  holds true  at the  south pole $SP$  of every sphere $\widetilde{S}$,
 \beaa
 \widetilde{b}\Big|_{SP}=-1 -\frac{2\widetilde{m}}{\widetilde{r}},
 \eeaa
 where $\widetilde{b}$ is such that we have
 \beaa
 \widetilde{\nu} =\widetilde{e}_3+ \widetilde{b} \widetilde{e}_4,
 \eeaa
 with $\widetilde{\nu}$  the unique vectorfield tangent to the hypersurface $\widetilde{\Si}_*$, normal to $\widetilde{S}$, and normalized by $\g(\widetilde{\nu}, \widetilde{e}_4)=-2$. 
 
\item The basis of $\ell=1$ modes $\widetilde{J}^{(p)}$ is given by $\widetilde{J}^{(p)}=J^{(p,\widetilde{S}_*)}$ on $\widetilde{S}_*$, and extended to $\widetilde{\Si}_*$ by $\widetilde{\nu}(\widetilde{J}^{(p)})=0$. Also, the $\ell=1$ modes of  $\widetilde{\div}\widetilde{\eta}$ and $\widetilde{\div}\widetilde{\xib}$ above are computed with respect to  this basis.
 
\item The transition functions $(f, \fb, \la)$ from the frame of $\MMextend$ to the frame of $\widetilde{\Si}_*$ satisfy on each sphere $\widetilde{S}\subset\widetilde{\Si}_*$
\beaa
\|(f, \fb, \log(\la))\|_{\hk_{k_*-3}(\widetilde{S})} &\les& \dg.
\eeaa
\end{itemize}

{\bf Step 8.} The spacelike GCM hypersurface $\widetilde{\Si}_*$ has been constructed in Step 7 in a small neighborhood of $\widetilde{S}_*$. We now focus on proving that it in fact extends all the way to the initial data layer. To this end, we denote by $u_1$ with
\beaa
1\leq u_1<\ug,
\eeaa
the minimal value of $u$ such that
\begin{itemize}
\item We have
\bea\lab{eq:aprioriestimatecrucialThmM7:0}
\widetilde{\Sigma}_*\cap\CC_u\neq\emptyset\textrm{ for any }u_1\leq u\leq \ug.
\eea

\item There exists a large constant $D\geq 1$ such that we have for any sphere $\widetilde{S}$ of $\widetilde{\Sigma}_*(u\geq u_1)$
\bea\lab{eq:aprioriestimatecrucialThmM7:1}
\|(f, \fb, \log(\la))\|_{\hk_{k_*- 3}(\widetilde{S})} &\leq& Du_*\dg.
\eea

\item For the same large constant $D\geq 1$ as above, we have along $\widetilde{\Sigma}_*(u\geq u_1)$
\bea\lab{eq:aprioriestimatecrucialThmM7:2}
|\psi(s)| &\leq& Du_*\dg,
\eea
where the function $\psi(s)$ is such that the curve 
\bea\lab{eq:checkthatthereisindeedalargeruinnewlastGCMS:6}
\Big(u=-s+c_{\widetilde{\Sigma}_*}+\psi(s),\, s, \,\th=\pi, \varphi\Big)\textrm{ with }\psi(\sg)=0,
\eea
coincides with the south poles of the sphere $\widetilde{S}$ of $\widetilde{\Sigma}_*$ and the constant $c_{\widetilde{\Sigma}_*}$ is fixed by the condition $\psi(\sg)=0$. 
\end{itemize}
The fact that $\psi(\sg)=0$ together with the bounds of Step 7 implies that \eqref{eq:aprioriestimatecrucialThmM7:0} \eqref{eq:aprioriestimatecrucialThmM7:1} \eqref{eq:aprioriestimatecrucialThmM7:2} hold for $u_1<\ug$ with $u_1$ close enough to $\ug$. By a continuity argument based on reapplying Theorem 4.1 in  \cite{Shen} (restated here in 
Theorem \ref{theorem:constuctionGCMH}), it suffices to show that we may improve the bounds \eqref{eq:aprioriestimatecrucialThmM7:1} \eqref{eq:aprioriestimatecrucialThmM7:2} independently of the value of $u_1$.

{\bf Step 9.} We now focus on improving the bounds \eqref{eq:aprioriestimatecrucialThmM7:1} \eqref{eq:aprioriestimatecrucialThmM7:2}. We first prove that $\widetilde{\Si}_*(u\geq u_1)$ is included in $\widetilde{\RR}$. Indeed, \eqref{eq:aprioriestimatecrucialThmM7:1} \eqref{eq:aprioriestimatecrucialThmM7:2} imply, using also the dominant condition on $r$ in $\widetilde{\RR}$, 
\beaa
\sup_{\widetilde{\Si}_*(u\geq u_1)}|u+s-c_{\widetilde{\Sigma}_*}| &\les& \sup_{\widetilde{\Si}_*(u\geq u_1)}\Big(|\psi|+r|f|+r|\fb|\Big)\\
&\les& Du_*\dg\\
&\les& \frac{Du_*}{r}\ep_0\Delta_{ext}\\
&\les& \ep_0 D\ep_0\Delta_{ext}\\
&\les& \ep_0\Delta_{ext}.
\eeaa
On the other hand, by construction, $\psi(\sg)=0$ and the south pole of $\ovS$ and $\widetilde{S}_*$ coincide, so that we have
\beaa
c_{\widetilde{\Sigma}_*} &=& \ug+\sg=u_*+r_*+\frac{\de_{ext}}{2}+\frac{3d_0r_*}{4u_*}\de_{ext}\\
&=& c_{\Sigma_*}+\frac{3}{4}\left(1+\frac{2u_*}{3d_0r_*}\right)\Delta_{ext}
\eeaa
and hence, using also the dominant condition on $r$ in $\widetilde{\RR}$, 
\beaa
\sup_{\widetilde{\Si}_*(u\geq u_1)}\left|u+s-c_{\Sigma_*}-\frac{3}{4}\Delta_{ext}\right| &\les& \left(\frac{u_*}{2d_0r_*}+\ep_0\right)\Delta_{ext}\\
&\les& \ep_0\Delta_{ext}.
\eeaa
In view of the definition of $\widetilde{\RR}$, we infer
\bea
\widetilde{\Si}_*(u\geq u_1) \subset \widetilde{\RR}
\eea
as claimed.

{\bf Step 10.} Since $\widetilde{\Si}_*(u\geq u_1) \subset \widetilde{\RR}$, the bounds \eqref{eq:controlofthesmallGCMquantitiesbsmallatSPforThmM7}, \eqref{eq:controlofthesmallGCMquantitieskackabcmucforThmM7} and \eqref{eq:controlofthesmallGCMquantitiesell=1modedivetaanddivxibforThmM7} 
 apply, and hence we have 
\beaa
\sup_{\widetilde{\RR}}\left|\left(e_3(u)+e_3(s)\right)\Big|_{SP} -\left(1+\frac{2m}{r}\right)\right| &\les& \frac{\ep_0}{r}\Delta_{ext}\les \dg,
\eeaa
and for all $k\leq k_*-4$
\beaa
\nn\sup_{\widetilde{\RR}}\Bigg(r^2\left|\dk^k\left(\kac\right)\right|+r^2\left|\dk^{k}\left(\kabc -\left(\Cb_0+\sum_p\Cb_p\Jp\right)\right)\right|\\
+r^2\left|\dk^{k}\left(\muc -\left(M_0+\sum_pM_p\Jp\right)\right)\right|\Bigg) &\les& \frac{\ep_0}{r}\Delta_{ext}\les\dg,
\eeaa 
as well as
\beaa
\sup_{\widetilde{\RR}}r^2\Big(|(\div\eta)_{\ell=1}|+|(\div\xib)_{\ell=1}|\Big)  &\les& \frac{\ep_0}{r}\Delta_{ext}\les\dg.
\eeaa
Together with the a priori estimates in the proof of Theorem 4.1 in  \cite{Shen}, this yields
\beaa
|\psi'(s)| &\les& \left|1+\frac{2\widetilde{m}}{\widetilde{r}}-\left(e_3(u)+e_3(s)\right)\Big|_{SP}\right|+|\la-1|\\
&\les& \left|\frac{\widetilde{m}}{\widetilde{r}}-\frac{m}{r}\right|+|\la-1|+\dg.
\eeaa

Now, we need to estimate $\widetilde{r}-r$ and $\widetilde{m}-m$. We claim
\bea\lab{eq:checkthatthereisindeedalargeruinnewlastGCMS:9}
\left|\widetilde{r}-r\right|+\left|\widetilde{m}-{m}\right| &\les&  Du_*\dg.
\eea
Indeed, in view of \eqref{eq:aprioriestimatecrucialThmM7:1}, and using Lemma 5.8 in \cite{KS-GCM1} (restated here in Proposition \ref{proposition:58-59GCM1}), we have
\beaa
\left|\widetilde{r}-r\right| &\les& \sup_{\widetilde{S}}r(|f|+|\fb|)\les Du_*\dg,
\eeaa
which is the stated estimate for $\widetilde{r}-r$ in \eqref{eq:checkthatthereisindeedalargeruinnewlastGCMS:9}. Next, recall that $\ovS=S(\ug, \sg)$ is the sphere of the foliation of $\widetilde{\RR}\cap\{u\geq u_*\}$ which shares the same south pole a $\widetilde{S}_*$. We denote $\mg$ the Hawking mass of $\ovS$ and recall that $\widetilde{m}$ denotes the Hawking mass of $\widetilde{S}_*$. Then, in view of \eqref{eq:aprioriestimatecrucialThmM7:1}, and using Corollary 5.17 in \cite{KS-GCM1} (restated here in Proposition \ref{proposition:58-59GCM1}), we have
\beaa
\left|\widetilde{m}-\mg\right| &\les& \sup_{\widetilde{S}}r(|f|+|\fb|)\les Du_*\dg.
\eeaa
Also, since $\ovS\subset\widetilde{\RR}\cap\{u\geq u_*\}$, and since $\mg$ denotes the  Hawking mass of $\ovS$, we have in view of \eqref{eq:controlofmHminusmonwidetildeRRmoreustar:ThmM7}
\beaa
|\mg -m| &\les&  \dg.
\eeaa
We deduce 
\beaa
\left|\widetilde{m}-m\right| &\les&  Du_*\dg
\eeaa
which concludes the proof of \eqref{eq:checkthatthereisindeedalargeruinnewlastGCMS:9}. 

We infer from \eqref{eq:checkthatthereisindeedalargeruinnewlastGCMS:9} and the estimate immediately above for $|\psi'(s)|$, using also the dominant condition on $r$ in $\widetilde{\RR}$, 
\beaa
|\psi'(s)| &\les& \frac{Du_*}{r}\dg+\dg\\
&\les& \left(1+\ep_0 D\right)\dg\\
&\les& \dg.
\eeaa
Integrating from $\sg$ where $\psi(\sg)=0$, we infer
\beaa
|\psi(s)| &\les& |s-\sg|\dg\\
&\les& u_*\dg
\eeaa
which improves \eqref{eq:aprioriestimatecrucialThmM7:2} for $D\geq 1$ large enough. 

Similarly, we obtain, using the a-priori estimates for GCM spheres in \cite{KS-GCM1},
\beaa
&&\|(f, \fb, \log(\la))\|_{\hk_{k_*-3}(\widetilde{S})}\\
 &\les& \max_{k\leq k_*-4}\sup_{\widetilde{\RR}}\Bigg(r^2\left|\dk^k\left(\kac\right)\right|+r^2\left|\dk^{k}\left(\kabc -\left(\Cb_0+\sum_p\Cb_p\Jp\right)\right)\right|\\
&&+r^2\left|\dk^{k}\left(\muc -\left(M_0+\sum_pM_p\Jp\right)\right)\right|\Bigg)+r\left(\left|(\widetilde{\div}f)_{\ell=1}\right|+\left|(\widetilde{\div}\fb)_{\ell=1}\right|\right),
\eeaa
and hence
\beaa
\|(f, \fb, \log(\la))\|_{\hk_{k_*-3}(\widetilde{S})} &\les& r^2\left(\left|(\widetilde{\div}f)_{\ell=1}\right|+\left|(\widetilde{\div}\fb)_{\ell=1}\right|\right)+\dg,
\eeaa
where the $\ell=1$ modes are taken w.r.t. the  basis $\widetilde{J}^{(p)}$. Also, using the a priori estimates in the proof of Theorem 4.1 in \cite{Shen}, we have
\beaa
\left|\tilde{\nu}\left((\widetilde{\div}f)_{\ell=1}\right)\right|+\left|\tilde{\nu}\left((\widetilde{\div}\fb)_{\ell=1}\right)\right|  &\les& r^{-2}\dg+\frac{1}{r}\left(\left|(\widetilde{\div}f)_{\ell=1}\right|+\left|(\widetilde{\div}\fb)_{\ell=1}\right|\right).
\eeaa
In view of \eqref{eq:aprioriestimatecrucialThmM7:1}, we infer
\beaa
\left|\tilde{\nu}\left((\widetilde{\div}f)_{\ell=1}\right)\right|+\left|\tilde{\nu}\left((\widetilde{\div}\fb)_{\ell=1}\right)\right| &\les& r^{-2}\dg+r^{-3}Du_*\dg
\eeaa
and integrating from $\widetilde{S}_*$, we infer, using also the dominant condition for $r$ in $\widetilde{\RR}$, 
\beaa
r^{2}\left(\left|(\widetilde{\div}f)_{\ell=1}\right|+\left|(\widetilde{\div}\fb)_{\ell=1}\right|\right) &\les& u_*\dg+\frac{D(u_*)^2}{r}\dg\\
&\les& \left(1+\ep_0 D\right)u_*\dg\\
&\les& u_*\dg.
\eeaa
This yields
\beaa
\|(f, \fb, \log(\la))\|_{\hk_{k_*-3}(\widetilde{S})} &\les& u_*\dg
\eeaa
which improves \eqref{eq:aprioriestimatecrucialThmM7:1} for $D\geq 1$ large enough. We thus conclude that $u_1=1$, i.e. $\widetilde{\Si}_*$  extends all the way to the initial data layer, $\widetilde{\Si}_*\subset\widetilde{\RR}$, and we have the bounds
\beaa
\|(f, \fb, \log(\la))\|_{\hk_{k_*-3}(\widetilde{S})} \les u_*\dg, \qquad |\psi(s)| \les u_*\dg.
\eeaa
In view of the definition of $\dg$, we infer in particular for any sphere $\widetilde{S}$ of $\widetilde{\Si}_*$
\bea\lab{eq:checkthatthereisindeedalargeruinnewlastGCMS:4}
\|(f, \fb, \log(\la))\|_{\hk_{k_*-3}(\widetilde{S})} \les \ep_0\de_{ext}, \qquad  |\psi(s)| \les \ep_0\de_{ext}.
\eea

{\bf Step 11.} As $\widetilde{\Si}_*$  extends all the way to the initial data layer, this allows us to calibrate $\tilde{u}$ along $\widetilde{\Si}_*$ by fixing the value $\widetilde{u}=1$ as in Remark \ref{rmk:calibrationofu}:
\bea
\widetilde{S}_1=\widetilde{\Si}_*\cap\{\widetilde{u}=1\}\textrm{ is such that }\widetilde{S}_1\cap\{\ut_{\idl}=1\}\cap\{\widetilde{\th}_{\idl}=\pi\}\neq \emptyset,  
\eea
i.e. $\widetilde{S}_1$  is the unique sphere of $\widetilde{\Si}_*$ intersecting the curve of the south poles of $\{\ut_{\idl}=1\}$ in the part $\Lextt$ of the initial data layer constructed in section \ref{section:geodesicfoliation8}.

Now that $\widetilde{u}$ is calibrated, we define
\bea
\tilde{u}_*:=\tilde{u}(\widetilde{S}_*).
\eea
For the proof of Theorem M7, we need in particular to show that $\tilde{u}_*>u_*$. First, note that, since $\tilde{u}+\tilde{r}$ is constant along $\widetilde{\Si}_*$, we have
\bea\lab{eq:checkthatthereisindeedalargeruinnewlastGCMS:8}
\widetilde{\Si}_*=\Big\{\widetilde{u}+\widetilde{r}=1+\widetilde{r}(\widetilde{S}_1)\Big\}.
\eea
Since $\widetilde{S}_*\subset\widetilde{\Si}_*$, and in view of \eqref{eq:checkthatthereisindeedalargeruinnewlastGCMS:8}, \eqref{eq:checkthatthereisindeedalargeruinnewlastGCMS:5}, \eqref{eq:checkthatthereisindeedalargeruinnewlastGCMS:6}, we infer, 
\beaa
\left|\widetilde{u}(\widetilde{S}_*)-\left(u_*+\frac{\de_{ext}}{2}\right)\right| &=& \left|\widetilde{u}(\widetilde{S}_*)-u(\ovS)\right|\\
&=& \left|1+\widetilde{r}(\widetilde{S}_1)-\widetilde{r}(\widetilde{S}_*)-\left(-s(\ovS)+c_{\widetilde{\Si}_*}\right)\right|.
\eeaa
Next, note from
\beaa
s=r\textrm{ on }\Si_*, \qquad e_4(r-s)=\frac{r}{2}\left(\ov{\ka}-\frac{2}{r}\right),
\eeaa
that we have
\bea\lab{eq:checkthatthereisindeedalargeruinnewlastGCMS:3}
\sup_{\widetilde{\RR}}|r-s| &\les& \frac{\ep_0}{r}\Delta_{ext}\les \ep_0\de_{ext}.
\eea
Together with  \eqref{eq:checkthatthereisindeedalargeruinnewlastGCMS:9}, this yields
\beaa
\left|\widetilde{u}(\widetilde{S}_*)-\left(u_*+\frac{\de_{ext}}{2}\right)\right| &\les& \Big|1+\widetilde{r}(\widetilde{S}_1)-c_{\widetilde{\Si}_*}\Big| +\ep_0\de_{ext}.
\eeaa
Since $c_{\widetilde{\Si}_*}$ in \eqref{eq:checkthatthereisindeedalargeruinnewlastGCMS:6} is a constant, we have in particular 
\beaa
c_{\widetilde{\Si}_*}=u(\widetilde{S}_1)+r(\widetilde{S}_1)-\psi(s(\widetilde{S}_1))
\eeaa
and thus
\beaa
\left|\widetilde{u}(\widetilde{S}_*)-\left(u_*+\frac{\de_{ext}}{2}\right)\right| &\les& \Big|1+\widetilde{r}(\widetilde{S}_1)-u(\widetilde{S}_1)-r(\widetilde{S}_1)+\psi(s(\widetilde{S}_1))\Big| +\ep_0\de_{ext}\\
&\les& \Big|1-u(\widetilde{S}_1)\Big|+\Big|\widetilde{r}(\widetilde{S}_1)-r(\widetilde{S}_1)\Big|+\Big|\psi(s(\widetilde{S}_1))\Big| +\ep_0\de_{ext}.
\eeaa
In view of \eqref{eq:checkthatthereisindeedalargeruinnewlastGCMS:4} and \eqref{eq:checkthatthereisindeedalargeruinnewlastGCMS:9}, we infer
\beaa
\left|\widetilde{u}(\widetilde{S}_*)-\left(u_*+\frac{\de_{ext}}{2}\right)\right| &\les& \Big|1-u(\widetilde{S}_1)\Big|+\ep_0\de_{ext}.
\eeaa
Also, since:
\begin{itemize}
\item we have by the calibration of $u$ 
\beaa
u=1\textrm{ on }S_1\cap\{\ut_{\LL_0}=1\}\cap\{\widetilde{\th}_{\idl}=\pi\},
\eeaa
i.e. $S_1$  is the unique sphere of $\Si_*$ intersecting the curve of the south poles of $\{\ut_{\idl}=1\}$ in the part $\Lextt$ of the initial data layer,

\item we have  the following control for the change of frame coefficients $(f_0, \fb_0, \la_0)$ between the outgoing geodesic frame of $\Lextt$ and the outgoing geodesic foliation $(u,s)$ initialized on  $\Si_*$
\beaa
|f_0|\les \frac{\ep_0}{r}, \qquad |\la_0-1|+|\fb_0|\les \ep_0,
\eeaa
see Step 14 in the proof of Theorem M0, and hence 
\beaa
\widetilde{e}_4^{\LL_0}(u-1)=\la_0\left(e_4+f_0\c\nab+\frac{1}{4}|f_0|^2e_3\right)u=\frac{\la_0}{4}|f_0|^2e_3(u)=O\left(\frac{\ep_0}{r^2}\right),
\eeaa 
\end{itemize}
we infer
\beaa
\sup_{\widetilde{\RR}\cap\{\ut_{\LL_0}=1\}\cap\{\widetilde{\th}_{\idl}=\pi\}}|u-1| \les  \Delta_{ext}\frac{\ep_0}{r^2}\les\ep_0\de_{ext}.
\eeaa
This yields $|1-u(\widetilde{S}_1)| \les\ep_0\de_{ext}$ and hence
\bea\lab{eq:checkthatthereisindeedalargeruinnewlastGCMS:10}
\left|\widetilde{u}(\widetilde{S}_*)-\left(u_*+\frac{\de_{ext}}{2}\right)\right| &\les& \ep_0\de_{ext}.
\eea
In particular, we deduce, for $\ep_0$ small enough, 
\bea\lab{eq:checkthatthereisindeedalargeruinnewlastGCMS:11}
\widetilde{u}(\widetilde{S}_*)>u_*
\eea
as desired.

{\bf Step 12.} We would like to check that the  condition \eqref{eq:behaviorofronS-star} for $r$ on $\widetilde{\Si}_*$ holds, i.e. we need to prove that there exists a choice of constant $d_0$ satisfying $\frac{1}{2}\leq d_0\leq 1$ such that  
\beaa
\widetilde{r}(\widetilde{S}_*) &=& \de_*\ep_0^{-1}(\widetilde{u}(\widetilde{S}_*))^{1+\dec}. 
\eeaa
To this end, note that we have in view of \eqref{eq:checkthatthereisindeedalargeruinnewlastGCMS:9}, \eqref{eq:checkthatthereisindeedalargeruinnewlastGCMS:3} and \eqref{eq:checkthatthereisindeedalargeruinnewlastGCMS:10}
\beaa
\widetilde{r}(\widetilde{S}_*) - \de_*\ep_0^{-1}(\widetilde{u}(\widetilde{S}_*))^{1+\dec} &=& s(\ovS) +O\left(\ep_0\de_{ext}\right) - \de_*\ep_0^{-1}\left(u_*+\frac{\de_{ext}}{2}+O\left(\ep_0\de_{ext}\right)\right)^{1+\dec}\\
&=& s(\ovS) - \de_*\ep_0^{-1}(u_*)^{1+\dec} -\frac{1+\dec}{2}\de_*\ep_0^{-1}(u_*)^{\dec}\de_{ext}\\
&&+\de_*\ep_0^{-1}(u_*)^{\dec}\de_{ext}O\left(\frac{\de_{ext}}{u_*}+\ep_0\right)+O\left(\ep_0\de_{ext}\right).
\eeaa
Together with \eqref{eq:checkthatthereisindeedalargeruinnewlastGCMS:5}, we infer
\beaa
\widetilde{r}(\widetilde{S}_*) - \de_*\ep_0^{-1}(\widetilde{u}(\widetilde{S}_*))^{1+\dec} &=& r_*+\frac{3d_0r_*}{4u_*}\de_{ext} - \de_*\ep_0^{-1}(u_*)^{1+\dec} -\frac{1+\dec}{2}\de_*\ep_0^{-1}(u_*)^{\dec}\de_{ext}\\
&&+\de_*\ep_0^{-1}(u_*)^{\dec}\de_{ext}O\left(\frac{\de_{ext}}{u_*}+\ep_0\right)+O\left(\ep_0\de_{ext}\right)\\
&=& r_*-\de_*\ep_0^{-1}(u_*)^{1+\dec} + \left(\frac{3d_0r_*}{4}-\frac{1+\dec}{2}\de_*\ep_0^{-1}(u_*)^{1+\dec}\right)\frac{\de_{ext}}{u_*}\\
&&+\de_*\ep_0^{-1}(u_*)^{\dec}\de_{ext}O\left(\frac{\de_{ext}}{u_*}+\ep_0\right)+O\left(\ep_0\de_{ext}\right).
\eeaa
Since we have by the condition \eqref{eq:behaviorofronS-star} of $r$ on $\Si_*$
\beaa
r_*= \de_*\ep_0^{-1}u_*^{1+\dec}, 
\eeaa
we deduce
\beaa
\widetilde{r}(\widetilde{S}_*) - \de_*\ep_0^{-1}(\widetilde{u}(\widetilde{S}_*))^{1+\dec} &= &  \left(\frac{3d_0}{4}-\frac{1+\dec}{2}\right)\frac{r_*\de_{ext}}{u_*}+\de_*\ep_0^{-1}(u_*)^{\dec}\de_{ext}O\left(\frac{\de_{ext}}{u_*}+\ep_0\right)+O\left(\ep_0\de_{ext}\right)\\
&=& \frac{3r_*\de_{ext}}{4u_*}\left(d_0 -\frac{2+2\dec}{3}+O\left(\ep_0+\frac{\de_{ext}}{u_*}\right)\right).
\eeaa
Thus, we may choose the constant $d_0$ such that $\frac{1}{2}\leq d_0\leq 1$ and 
\beaa
\widetilde{r}(\widetilde{S}_*) &=& \de_*\ep_0^{-1}(\widetilde{u}(\widetilde{S}_*))^{1+\dec} 
\eeaa
as desired.

{\bf Step 13.} We summarize the properties of $\widetilde{\Si}_*$ obtained so far:
\begin{itemize}
\item $\widetilde{\Si}_*$ is a spacelike hypersurface included in the spacetime region $\widetilde{\RR}$. 

\item The scalar function $\widetilde{u}$ is defined on $\widetilde{\Si}_*$ and its level sets  are topological 2-spheres denoted by $\widetilde{S}$.

\item The following GCM conditions holds on $\widetilde{\Si}_*$
 \beaa
&&\widecheck{\widetilde{\ka}}=0, \qquad \widecheck{\widetilde{\kab}}=\widetilde{\Cb}_0+\sum_p\widetilde{\Cb}_p\widetilde{J}^{(p)}, \qquad \widecheck{\widetilde{\mu}}=\widetilde{M}_0+\sum_p\widetilde{M}_p\widetilde{J}^{(p)},\\
&& (\widetilde{\div}\widetilde{\eta})_{\ell=1}=0, \qquad (\widetilde{\div}\widetilde{\xib})_{\ell=1}=0.
\eeaa

\item In addition, the following GCM conditions holds on the sphere $\widetilde{S}_*$ of $\widetilde{\Si}_*$
\beaa
 \widecheck{\widetilde{\kab}}=0, \qquad (\widetilde{\div}\widetilde{\b})_{\ell=1}=0, \qquad  (\widetilde{\curl}\widetilde{\b})_{\ell=1,\pm}=0, \qquad  (\widetilde{\curl}\widetilde{\b})_{\ell=1,0}=\frac{2\widetilde{a}\widetilde{m}}{\widetilde{r}^5}.
\eeaa

\item We have, for some constant $c_{\widetilde{\Si}_*}$,  
\beaa
\widetilde{u}+\widetilde{r}=c_{\widetilde{\Si}_*} , \quad \textrm{along} \quad \widetilde{\Si}_*.
\eeaa
\item The following normalization condition  holds true  at the  south pole $SP$  of every sphere $\widetilde{S}$,
 \beaa
 \widetilde{b}_*\Big|_{SP}=-1 -\frac{2\widetilde{m}}{\widetilde{r}},
 \eeaa
 where $\widetilde{b}_*$ is such that we have
 \beaa
 \widetilde{\nu} =\widetilde{e}_3+ \widetilde{b}_* \widetilde{e}_4,
 \eeaa
 with $\widetilde{\nu}$  the unique vectorfield tangent to the hypersurface $\widetilde{\Si}_*$, normal to $\widetilde{S}$, and normalized by $g(\widetilde{\nu}, \widetilde{e}_4)=-2$. 
 
 \item The  condition \eqref{eq:behaviorofronS-star} for $r$ on $\widetilde{S}_*$ holds, i.e. we have
\beaa
\widetilde{r}(\widetilde{S}_*) &=& \de_*\ep_0^{-1}(\widetilde{u}(\widetilde{S}_*))^{1+\dec}. 
\eeaa

\item $\tilde{u}$ is calibrated along $\widetilde{\Si}_*$ by fixing the value $\widetilde{u}=1$:
\bea
\widetilde{S}_1=\widetilde{\Si}_*\cap\{\widetilde{u}=1\}\textrm{ is such that }\widetilde{S}_1\cap\{u_{\LL_0}'=1\}\cap\{{}^{(ext)}\th_{\idl}'=\pi\}\neq \emptyset,  
\eea
i.e. $\widetilde{S}_1$  is the unique sphere of $\widetilde{\Si}_*$  intersecting the curve of the south poles of $\{u_{\idl}'=1\}$ in the part $\Lextt$ of the initial data layer.
\end{itemize}
Thus, $\widetilde{\Si}_*$ satisfies all the required properties for the future spacelike boundary of a GCM admissible spacetime, see section \ref{sec:admissibleGMCPGdatasetonSigmastar}. Furthermore, we have on  $\widetilde{\Si}_*$
\bea\lab{eq:conclusionofThmM7u'*>u*isok}
\widetilde{u}(\widetilde{S}_*)>u_*.
\eea

%%%%%%%%%%%%%%

\subsection{Steps 14--18}

%%%%%%%%%%%%%%

{\bf Step 14.} We introduce: 
\begin{itemize}
\item  the outgoing PG frame $(e_4', e_3', e_1', e_2')$ of $\Mext$ extended to  the   spacetime  $\MM^{(extend)}$, 

\item  the outgoing PG frame  $(\widetilde{e}_4', \widetilde{e}_3', \widetilde{e}_1', \widetilde{e}_2')$ initialized on $\widetilde{\Si}_*$ from the GCM frame $(\widetilde{e}_4, \widetilde{e}_3, \widetilde{e}_1, \widetilde{e}_2)$ by the change of frame with coefficients $(f'', \fb'', \la'')$ given by
\beaa
\la''=1,\qquad f''=\frac{\widetilde{a}}{\widetilde{r}}\widetilde{f}_0, \qquad \fb''=-\frac{(\widetilde{\nu}(\widetilde{r})-\widetilde{b}_*)}{1-\frac{1}{4}\widetilde{b}_*|f''|^2}f'',
\eeaa
where the 1-form $\widetilde{f}_0$  is chosen on $\widetilde{\Si}_*$ by
\beaa
(\widetilde{f}_0)_1=0, \quad (\widetilde{f}_0)_2=\sin(\widetilde{\th}), \quad\textrm{on}\quad \widetilde{S}_*, \qquad \widetilde{\nab}_{\widetilde{\nu}}\widetilde{f}_0=0\quad\textrm{on}\quad\widetilde{\Si}_*,
\eeaa
with $(e_1, e_2)$ specified on $\widetilde{S}_*$ by \eqref{eq:canonical-e1ande2onSstar:0}.
\end{itemize}

We have the following change of frame coefficients:
\begin{itemize}
\item $(f, \fb, \la)$, introduced in Step 7, and corresponding to the change from the outgoing geodesic frame $(e_4, e_3, e_1, e_2)$ of $\MMextend$ to the GCM frame $(\widetilde{e}_4, \widetilde{e}_3, \widetilde{e}_1, \widetilde{e}_2)$ of $\widetilde{\Si}_*$,

\item  $(f', \fb', \la')$, which we now introduce, corresponding to the change from the outgoing geodesic frame $(e_4, e_3, e_1, e_2)$ of $\MMextend$ to   the outgoing PG frame $(e_4', e_3', e_1', e_2')$ of $\Mext$ extended to  the   spacetime  $\MM^{(extend)}$,

\item  $(f'', \fb'', \la'')$, provided explicitly above, and corresponding to the change from the GCM frame $(\widetilde{e}_4, \widetilde{e}_3, \widetilde{e}_1, \widetilde{e}_2)$ of $\widetilde{\Si}_*$ to the outgoing PG frame $(\widetilde{e}_4', \widetilde{e}_3', \widetilde{e}_1', \widetilde{e}_2')$ initialized on $\widetilde{\Si}_*$,

\item $(f''', \fb''', \la''')$, which we now introduce, corresponding to the change from  the outgoing PG frame $(e_4', e_3', e_1', e_2')$ of $\Mext$ extended to  the   spacetime  $\MM^{(extend)}$ and the outgoing PG frame $(\widetilde{e}_4', \widetilde{e}_3', \widetilde{e}_1', \widetilde{e}_2')$  initialized on $\widetilde{\Si}_*$.
\end{itemize}

In this step, our goal is to control the change of frame coefficients $(f''', \fb''', \la''')$. In view of the above, we have schematically 
\beaa
(f''', \fb''', \la''')=(f'', \fb'', \la'')\circ (f, \fb, \la)\circ (f', \fb', \la')^{-1}
\eeaa
where $(f', \fb', \la')^{-1}$ denote the coefficients corresponding to  the inverse transformation coefficients of the transformation with coefficients $(f', \fb', \la')$. We infer
\beaa
\sup_{\widetilde{\Si}_*}\left|\widetilde{\dk}_*^k(f''', \fb''', \la'''-1)\right| &\les& \sup_{\widetilde{\Si}_*}\left|\widetilde{\dk}_*^k(f, \fb, \la-1)\right| +\sup_{\widetilde{\Si}_*}\left|\widetilde{\dk}_*^k(f''-f', \fb''-\fb', \la''-\la')\right|. 
\eeaa
Now, recall that  $(f, \fb, \la)$ satisfy in view of \eqref{eq:checkthatthereisindeedalargeruinnewlastGCMS:4} and Corollary 4.2 in \cite{Shen} (restated here in Corollary \ref{cor:constuctionGCMH})
\beaa
\sup_{\widetilde{S}\subset\widetilde{\Si}_*}\|\dk^{\leq k_*-3} (f, \fb, \log(\la))\|_{L^2(\widetilde{S})} \les \ep_0\de_{ext}.
\eeaa
Together with the Sobolev embedding on the spheres $\widetilde{S}$, we find
\bea\lab{eq:checkthatthereisindeedalargeruinnewlastGCMS:4:thereturn}
\sup_{\widetilde{\Si}_*}\widetilde{r}\,|\dk^{\leq k_*-5} (f, \fb, \log(\la))| \les \ep_0\de_{ext}.
\eea
We infer, for $\leq k_*-5$, 
\beaa
\sup_{\widetilde{\Si}_*}r\left|\widetilde{\dk}_*^k(f''', \fb''', \la'''-1)\right| &\les&  \ep_0\de_{ext}+\sup_{\widetilde{\Si}_*}r\left|\widetilde{\dk}_*^k(f''-f', \fb''-\fb', \la''-\la')\right|. 
\eeaa
Together with the explicit formulas above for  $(f'', \fb'', \la'')$, we obtain, for $\leq k_*-5$, 
\beaa
\sup_{\widetilde{\Si}_*}r\left|\widetilde{\dk}_*^k(f''', \fb''', \la'''-1)\right| &\les&  \ep_0\de_{ext}+\sup_{\widetilde{\Si}_*}r\left|\widetilde{\dk}_*^k\left(f'-\frac{\widetilde{a}}{\widetilde{r}}\widetilde{f}_0, \fb' + \frac{(\widetilde{\nu}(\widetilde{r})-\widetilde{b}_*)}{1-\frac{1}{4}\widetilde{b}_*\frac{(\widetilde{a})^2}{(\widetilde{r})^2}|\widetilde{f}_0|^2}\frac{\widetilde{a}}{\widetilde{r}}\widetilde{f}_0, \la' -1\right)\right|. 
\eeaa
We deduce,  for $\leq k_*-5$, 
\bea\lab{eq:Step14ThmM7:boundonffblatripleprime}
&&\sup_{\widetilde{\Si}_*}r\left|\widetilde{\dk}_*^k(f''', \fb''', \la'''-1)\right|\\
\nn &\les&  \ep_0\de_{ext}+\sup_{\widetilde{\Si}_*}r\left|\widetilde{\dk}_*^k\left(f'-\frac{a}{r}f_0, \fb' + \frac{e_3(s)}{1+\frac{1}{4}(e_3(u)+e_3(s))\frac{a^2}{r^2}|f_0|^2}\frac{a}{r}f_0, \la' -1\right)\right| \\
\nn &&+\sup_{\widetilde{\Si}_*}\left(|\widetilde{\dk}_*^k(\widetilde{a}\widetilde{f}_0-af_0)|+r^{-1}|\widetilde{\dk}_*^k(\widetilde{r}-r)|+|\widetilde{\dk}_*^k(\widetilde{b}_*+e_3(u)+e_3(s))|+|\widetilde{\dk}_*^k(\widetilde{\nu}(\widetilde{r})+e_3(u))|\right).
\eea

We now control the terms on the RHS of \eqref{eq:Step14ThmM7:boundonffblatripleprime}. Note first the have have on $\Si^{(extend)}_*$, in view of the initialization of the PG frame of $\Mext$, 
\beaa
\la'=1,\qquad f'=\frac{a}{r}f_0, \qquad \fb'=-\frac{(\nu(r)-b_*)}{1-\frac{1}{4}b_*\frac{a^2}{r^2}|f_0|^2}\frac{a}{r}f_0.
\eeaa
Also, we have $r=s$ on $\Si^{(extend)}_*$, and hence, using also $e_4(u)=0$ and $e_4(s)=1$, and the fact that $\nu$ is tangent to $\Si^{(extend)}_*$, 
\beaa
0&=& \nu(u+r)=\nu(u)+\nu(s)=e_3(u)+b_*e_4(u)+e_3(s)+b_*e_4(s)=e_3(u)+e_3(s)+b_*,\\
\nu(r) &=& \nu(s)=e_3(s)+b_*e_4(s)=e_3(s)+b_*,
\eeaa
and hence, we deduce 
\beaa
b_*=-e_3(u)-e_3(s), \qquad \nu(r)= -e_3(u)\quad\textrm{on}\quad \Si^{(extend)}_*.
\eeaa
In view of the above, this yields
\beaa
\la'-1=0,\qquad f'-\frac{a}{r}f_0=0, \qquad \fb'+\frac{e_3(s)}{1+\frac{1}{4}(e_3(u)+e_3(s))\frac{a^2}{r^2}|f_0|^2}\frac{a}{r}f_0=0\quad\textrm{on}\quad \Si^{(extend)}_*.
\eeaa
Next, using Corollary \ref{cor:transportequationine4forchangeofframecoeff:simplecasefirst}, we have the following transport equations for  $(f', \fb', \la')$
\beaa
\nab_{{\la'}^{-1}e_4'}f'+\frac{1}{2}\ka f'  &=&  - f'\c\chih+E_1(f', \Ga),\\
{\la'}^{-1}e_4'(\log\la') &=& 2f'\c\ze+E_2(f', \Ga),\\
\nab_{{\la'}^{-1}e_4'}\fb'+\frac{1}{2}\ka\fb'  &=&    2\nab'(\log\la')  +2\omb f'  -\fb'\c\chih -  \fb'\c\nab'f\\
&&+E_3(f', \fb', \Ga),
\eeaa
where $E_1(f', \Ga)$ and $E_2(f', \Ga)$ contain expressions of the type $O(\Ga {f'}^2)$ with no derivatives, and $E_3(f', \fb', \Ga)$ contain expressions of the type $O(\Ga(f', \fb')^2)$ with no derivatives. Together with  the control of Step 1 on $\MMextend$, we infer, using also the extension $\nab_4f_0=0$ of $f_0$ from $\Si^{(extend)}_*$ to $\widetilde{\RR}$, for $k\leq k_*-3$, 
\beaa
\sup_{\widetilde{\RR}}\left(r^{-1}\left|\dk^ke_4'\left(r\left(f'-\frac{a}{r}f_0\right)\right)\right|+\left|\dk^ke_4'\left(\log(\la')\right)\right|\right)\les \frac{1}{r^3}+\frac{1}{r}|\dk^{\leq k}\Ga_g| \les \frac{1}{r^3},
\eeaa
and, for $k\leq k_*-4$,
\beaa
&&\sup_{\widetilde{\RR}}r^{-1}\left|\dk^ke_4'\left(r\left(\fb'+\frac{e_3(s)}{1+\frac{1}{4}(e_3(u)+e_3(s))\frac{a^2}{r^2}|f_0|^2}\frac{a}{r}f_0 -2r\nab'(\log(\la'))\right)\right)\right|\\
&\les& \frac{1}{r^3}+\frac{1}{r}|\dk^{\leq k}\Ga_b|\les \frac{1}{r^3}+\frac{\ep_0}{r^2}.
\eeaa
Together with the above identities on $\Si^{(extend)}_*$, we integrate from $\Si^{(extend)}_*$ and obtain, using the dominant condition for $r$ on $\widetilde{\RR}$, for $k\leq k_*-3$,
\beaa
\sup_{\widetilde{\RR}}r\left|\dk^k\left(f'-\frac{a}{r}f_0, \la' -1\right)\right| &\les& \frac{\De_{ext}}{r^2}\les\frac{\ep_0\De_{ext}}{r}\les  \ep_0\de_{ext}
\eeaa
and for  $k\leq k_*-4$
\beaa
\sup_{\widetilde{\RR}}r\left|\dk^k\left(\fb' + \frac{e_3(s)}{1+\frac{1}{4}(e_3(u)+e_3(s))\frac{a^2}{r^2}|f_0|^2}\frac{a}{r}f_0\right)\right| &\les& \frac{\De_{ext}}{r^2}+\frac{\ep_0\De_{ext}}{r}\les\frac{\ep_0\De_{ext}}{r}\les  \ep_0\de_{ext}.
\eeaa
Together with \eqref{eq:Step14ThmM7:boundonffblatripleprime}, we infer, for $\leq k_*-5$, 
\bea\lab{eq:Step14ThmM7:boundonffblatripleprime:bis}
\sup_{\widetilde{\Si}_*}r\left|\widetilde{\dk}_*^k(f''', \fb''', \la'''-1)\right| &\les&  \ep_0\de_{ext}  +\sup_{\widetilde{\Si}_*}\Big(|\widetilde{\dk}_*^k(\widetilde{a}\widetilde{f}_0-af_0)|+r^{-1}|\widetilde{\dk}_*^k(\widetilde{r}-r)|\\
\nn&&+|\widetilde{\dk}_*^k(\widetilde{b}_*+e_3(u)+e_3(s))|+|\widetilde{\dk}_*^k(\widetilde{\nu}(\widetilde{r})+e_3(u))|\Big).
\eea

{\bf Step 15.} Next, we focus on the control of the terms on the RHS of \eqref{eq:Step14ThmM7:boundonffblatripleprime:bis}. To this end, we first estimate $\widetilde{r}-r$ on $\widetilde{\Si}_*$. Recall from \eqref{eq:controlofthesmallGCMquantitieskackabcmucforThmM7} that we have in particular, for $k\leq k_*-4$, 
\beaa
\sup_{\widetilde{\RR}}r^2\left|\dk^k\left(\kac\right)\right| &\les& \frac{\ep_0}{r}\Delta_{ext}\les \ep_0\de_{ext}.
\eeaa
Together with the GCM condition $\widecheck{\widetilde{\ka}}=0$, we infer
\beaa
\sup_{\widetilde{\Si}_*}r^2\left|\widetilde{\dk}_*^k\left(\widecheck{\widetilde{\ka}}-\kac\right)\right| &\les&  \ep_0\de_{ext}.
\eeaa
Now, we have 
\beaa
\widecheck{\widetilde{\ka}}-\kac &=& \widetilde{\ka}-\ka+\frac{2}{\widetilde{r}}-\frac{2}{r}= \widetilde{\ka}-\ka-\frac{2(\widetilde{r}-r)}{r\widetilde{r}}
\eeaa
so that
\beaa
\widetilde{r}-r &=& \frac{r\widetilde{r}}{2}\left(\widetilde{\ka}-\ka -\left(\widecheck{\widetilde{\ka}}-\kac\right)\right)
\eeaa
and hence, using the above estimate for $\widecheck{\widetilde{\ka}}-\kac$, we have, for $k\leq k_*-4$,  
\beaa
\sup_{\widetilde{\Si}_*}\left|\widetilde{\dk}_*^k\left(\widetilde{r}-r\right)\right| &\les&  \ep_0\de_{ext}+\sup_{\widetilde{\Si}_*}r^2\left|\widetilde{\dk}_*^k\left(\widetilde{\ka} -\ka\right)\right|.
\eeaa
Using the change of frame formula for $\widetilde{\ka}$, together with  the control \eqref{eq:checkthatthereisindeedalargeruinnewlastGCMS:4:thereturn} for $(f, \fb, \la)$ and the control of Step 1 on $\MMextend$, we deduce for $k\leq k_*-6$
\bea\lab{eq:controlofwidetilderminusronSigmatildestar:ThmM7}
\sup_{\widetilde{\Si}_*}\left|\widetilde{\dk}_*^k\left(\widetilde{r}-r\right)\right| &\les&  \ep_0\de_{ext}.
\eea

Next, we control $\widetilde{b}_*+e_3(u)+e_3(s)$. First, note that we have
\beaa
\widetilde{\nu}(\widetilde{r}-r)=\widetilde{\nu}(\widetilde{r})-\widetilde{e}_3(r)-\widetilde{b}_*\widetilde{e}_4(r).
\eeaa
Together with the control of $\widetilde{r}-r$ in \eqref{eq:controlofwidetilderminusronSigmatildestar:ThmM7}, the
 fact that $\widetilde{\nu}$ is tangent to $\widetilde{\Si}_*$, the change of frame formulas,  the control \eqref{eq:checkthatthereisindeedalargeruinnewlastGCMS:4:thereturn} for $(f, \fb, \la)$ and the control of Step 1 on $\MMextend$, we deduce for $k\leq k_*-7$,
 \beaa
\sup_{\widetilde{\Si}_*}\left|\widetilde{\dk}_*^k\left(\widetilde{\nu}(\widetilde{r})- e_3(r)-\widetilde{b}_*e_4(r)\right)\right| &\les&  \ep_0\de_{ext}.
\eeaa
Also, since $e_4(r)=1$ and $s=r$ on $\Si^{(extend)}_*$, and since $\nu$ is tangent to $\Si^{(extend)}_*$ with $\nu=e_4+b_*e_4$, we infer $e_3(r-s)=0$ on  $\Si^{(extend)}_*$. As in \eqref{eq:checkthatthereisindeedalargeruinnewlastGCMS:3}, we rely on the propagation equation $e_4(r-s)=\frac{r}{2}(\ov{\ka}-\frac{2}{r})$ and the control of $\kac$ in $\widetilde{\RR}$ provided by \eqref{eq:controlofthesmallGCMquantitieskackabcmucforThmM7} to infer,  for $k\leq k_*-4$,
\beaa
\sup_{\widetilde{\RR}}\left|\dk^k\left(e_3(r-s)\right)\right| &\les&\frac{\ep_0}{r}\De_{ext}\les  \ep_0\de_{ext}.
\eeaa
Together with \eqref{eq:checkthatthereisindeedalargeruinnewlastGCMS:3}, and the fact that $e_4(r-s)=\nab(r-s)=0$ on $\widetilde{\RR}$, we infer,  for $k\leq k_*-3$,
\bea\lab{eq:controlofrminussonwidetildeRRhigherorderderivatives:ThmM7}
\sup_{\widetilde{\RR}}\left|\dk^k\left(r-s\right)\right| &\les&   \ep_0\de_{ext}.
\eea
The control \eqref{eq:controlofrminussonwidetildeRRhigherorderderivatives:ThmM7} of $r-s$ allows us to replace  $r$ by $s$ in the above formula which yields, since $e_4(s)=1$,  for $k\leq k_*-7$,
 \beaa
\sup_{\widetilde{\Si}_*}\left|\widetilde{\dk}_*^k\left(\widetilde{\nu}(\widetilde{r})- e_3(s) -\widetilde{b}_*\right)\right| &\les&  \ep_0\de_{ext}.
\eeaa
Since $\widetilde{\nu}(\widetilde{u}+\widetilde{r})=0$ on $\widetilde{\Si}_*$, we infer for $k\leq k_*-7$
 \beaa
\sup_{\widetilde{\Si}_*}\left|\widetilde{\dk}_*^k\left(\widetilde{b}_*+\widetilde{\nu}(\widetilde{u})+e_3(s)\right)\right| &\les&  \ep_0\de_{ext}.
\eeaa
Now, as part of the construction of $\widetilde{\Si}_*$, the following  transversality conditions on $\widetilde{\Si}_*$ are assumed, see  Theorem 4.1 in  \cite{Shen} (restated here in 
Theorem \ref{theorem:constuctionGCMH}),
\bea\lab{eq:transversalityconditionsonwidetildeSigmastar:ThmM7}
\widetilde{\xi}=\widetilde{\om}=0, \qquad \widetilde{\etab}=-\widetilde{\ze}, \qquad \widetilde{e}_4(\widetilde{r})=1, \qquad \widetilde{e}_4(\widetilde{u})=0.
\eea
We infer
\beaa
\widetilde{\nu}(\widetilde{u}) &=& \widetilde{e}_3(\widetilde{u})+\widetilde{b}_*\widetilde{e}_4(\widetilde{u})=\widetilde{e}_3(\widetilde{u})
\eeaa
and hence, for $k\leq k_*-7$,
\beaa
\sup_{\widetilde{\Si}_*}\left|\widetilde{\dk}_*^k\left(\widetilde{b}_*+\widetilde{e}_3(\widetilde{u})+e_3(s)\right)\right| &\les&  \ep_0\de_{ext}.
\eeaa
Also, using again the transversality conditions \eqref{eq:transversalityconditionsonwidetildeSigmastar:ThmM7}, we have
\beaa
\widetilde{\nab}(\widetilde{e}_3(\widetilde{u}))=(\widetilde{\ze}-\widetilde{\eta})\widetilde{e}_3(\widetilde{u}).
\eeaa
Similarly, we have for the outgoing geodesic foliation of $\widetilde{\RR}$
\beaa
\nab(e_3(u))=(\ze-\eta)e_3(u).
\eeaa
Together with the change of frame formulas, the control \eqref{eq:checkthatthereisindeedalargeruinnewlastGCMS:4:thereturn} for $(f, \fb, \la)$ and the control of Step 1 on $\MMextend$, we deduce, for $k\leq k_*-6$,
\beaa
\sup_{\widetilde{\Si}_*}\left|\widetilde{\dk}_*^k\left(\widetilde{\nab}(e_3(u))-(\widetilde{\ze}-\widetilde{\eta})e_3(u)\right)\right| &\les& r^{-1}\ep_0\de_{ext}.
\eeaa
Subtracting to the above identity for $\widetilde{e}_3(\widetilde{u})$, we infer for $k\leq k_*-6$
\beaa
\sup_{\widetilde{\Si}_*}\left|\widetilde{\dk}_*^k\left(\widetilde{\nab}(\widetilde{e}_3(\widetilde{u})-e_3(u))-(\widetilde{\ze}-\widetilde{\eta})(\widetilde{e}_3(\widetilde{u})-e_3(u))\right)\right| &\les& r^{-1}\ep_0\de_{ext}.
\eeaa
Together with the above estimate for $\widetilde{b}_*+\widetilde{e}_3(\widetilde{u})+e_3(s)$, we infer, for $k\leq k_*-8$,
\beaa
\sup_{\widetilde{\Si}_*}r\left|\widetilde{\dk}_*^k\left(\widetilde{\nab}\left(\widetilde{b}_*+e_3(u)+e_3(s)\right) -(\widetilde{\ze}-\widetilde{\eta})\left(\widetilde{b}_*+e_3(u)+e_3(s)\right)\right)\right| &\les&  \ep_0\de_{ext}.
\eeaa
Using  the change of frame formulas, the control \eqref{eq:checkthatthereisindeedalargeruinnewlastGCMS:4:thereturn} for $(f, \fb, \la)$ and the control of Step 1 on $\MMextend$ to control $\widetilde{\ze}$ and $\widetilde{\eta}$, we infer, going also for convenience from $L^\infty$ to the weaker $L^2$ estimate,  for $k\leq k_*-9$,
\beaa
\left\|\widetilde{\dk}_*^{k}\widetilde{\nab}\left(\widetilde{b}_*+e_3(u)+e_3(s)\right)\right\|_{\hk_1(\widetilde{S})} &\les& r^{-1}\ep_0\left\|\widetilde{\dk}_*^{k}\left(\widetilde{b}_*+e_3(u)+e_3(s)\right)\right\|_{\hk_1(\widetilde{S})}
\eeaa
Also, recall that we have obtained in Step 4 
\beaa
\big(e_3(u)+e_3(s)\big)\Big|_{SP} &=& 1+\frac{2m}{r}\textrm{ on }\Si^{(extend)}_*,
\eeaa 
where $SP$ denotes the south poles of the spheres of the $u$-foliation, i.e. $\th=\pi$.  Since at the south pole we have $b_*=-1-\frac{2m}{r}$ by our GCM condition on $b_*$, we infer on $\Si^{(extend)}_*$, since $\nu=e_3+b_*e_4$ is tangent to $\Si^{(extend)}_*$, 
\beaa
\left(\left(e_3-\left(1+\frac{2m}{r}\right)e_4\right)^k\left(e_3(u)+e_3(s)-\left(1+\frac{2m}{r}\right)\right)\right)\Big|_{SP} &=& 0\textrm{ on }\Si_*.
\eeaa 
Arguing as for \eqref{eq:controlofthesmallGCMquantitiesbsmallatSPforThmM7}, we  propagate forward in $e_4$,  and using the bounds of Step 1 on 
 $\MM^{(extend)}$, we infer, for $k\leq k_*-3$,
\beaa
\sup_{\widetilde{\RR}}\left|\left(\left(e_3-\left(1+\frac{2m}{r}\right)e_4\right)^k\left(e_3(u)+e_3(s)-\left(1+\frac{2m}{r}\right)\right)\right)\Big|_{SP}\right| &\les& \frac{\ep_0}{r}\Delta_{ext}\les\ep_0\de_{ext}.
\eeaa
On the other have, by our GCM condition on $\widetilde{\Si}_*$ for $\widetilde{b}_*$, we have on  $\widetilde{\Si}_*$
\beaa
\widetilde{b}_* \Big|_{SP}=-1-\frac{2\widetilde{m}}{\widetilde{r}}\quad\textrm{on}\quad\Si^{(extend)}_*,
\eeaa
and hence, we have on $\widetilde{\Si}_*$
\beaa
\left(\left(\widetilde{e}_3-\left(1+\frac{2\widetilde{m}}{\widetilde{r}}\right)\widetilde{e}_4\right)^k\left(\widetilde{b_*}+\left(1+\frac{2\widetilde{m}}{\widetilde{r}}\right)\right)\right)\Big|_{SP} &=& 0.
\eeaa
Subtracting the identity to the previous estimate, using  the change of frame formulas, the control \eqref{eq:checkthatthereisindeedalargeruinnewlastGCMS:4:thereturn} for $(f, \fb, \la)$ and the control of Step 1 on $\MMextend$, we infer, for $k\leq k_*-6$,
\beaa
\sup_{\widetilde{\Si}_*}\left|\left(\widetilde{\nu}^k\left(\widetilde{b}_*+e_3(u)+e_3(s)\right)\right)\Big|_{SP}\right| &\les& \ep_0\de_{ext}+r^{-1}|\widetilde{m}-m|+r^{-2}\sup_{\widetilde{\Si}_*}|\widetilde{\nu}^{\leq k}(\widetilde{r}-r)|.
\eeaa
Together with \eqref{eq:checkthatthereisindeedalargeruinnewlastGCMS:9} for $\widetilde{m}-m$, and \eqref{eq:controlofwidetilderminusronSigmatildestar:ThmM7} for $\widetilde{r}-r$, we deduce, for $k\leq k_*-6$,
\beaa
\sup_{\widetilde{\Si}_*}\left|\left(\widetilde{\nu}^k\left(\widetilde{b}_*+e_3(u)+e_3(s)\right)\right)\Big|_{SP}\right| &\les& \ep_0\de_{ext}.
\eeaa
Thus, introducing the scalar $h$ on $\widetilde{\Si}_*$ given by
\beaa
h &:=& \widetilde{b}_*+e_3(u)+e_3(s),
\eeaa
we have obtained so far on  $\widetilde{\Si}_*$,  for $k\leq k_*-9$,
\beaa
\left\|\widetilde{\dk}_*^{k}\widetilde{\nab}h\right\|_{\hk_1(\widetilde{S})} &\les& r^{-1}\ep_0\left\|\widetilde{\dk}_*^{k}h\right\|_{\hk_1(\widetilde{S})},
\eeaa
and, for $k\leq k_*-6$, 
\beaa
\sup_{\widetilde{\Si}_*}\left|\left(\widetilde{\nu}^k(h)\right)\Big|_{SP}\right| &\les& \ep_0\de_{ext}.
\eeaa
We decompose $\widetilde{\nu}^k(h)=\ov{\widetilde{\nu}^k(h)}+(\widetilde{\nu}^k(h)-\ov{\widetilde{\nu}^k(h)})$ and since $\ov{\widetilde{\nu}^k(h)}$ is constant on $\widetilde{S}$ we have
\beaa
|\ov{\widetilde{\nu}^k(h)}| &=& |\ov{\widetilde{\nu}^k(h)}_{|_{SP}}|\les  |(\widetilde{\nu}^k(h))_{|_{SP}}|+\|\widetilde{\nu}^k(h)-\ov{\widetilde{\nu}^k(h)}\|_{L^\infty(\widetilde{S})}\\
&\les& |(\widetilde{\nu}^k(h))_{|_{SP}}|+r^{-1}\|\widetilde{\nu}^k(h)-\ov{\widetilde{\nu}^k(h)}\|_{\hk_2(\widetilde{S})}
\eeaa
where we have used Sobolev. Together with Poincar\'e, we infer
\beaa
|\ov{\widetilde{\nu}^k(h)}| &\les&  |(\widetilde{\nu}^k(h))_{|_{SP}}|+\left\|\widetilde{\nab}\widetilde{\dk}_*^{k}h\right\|_{\hk_1(\widetilde{S})}.
\eeaa
In view of the above, we deduce for $k\leq k_*-9$, 
\beaa
|\ov{\widetilde{\nu}^k(h)}| &\les& \left\|\widetilde{\dk}_*^{k}\widetilde{\nab}h\right\|_{\hk_1(\widetilde{S})}+\ep_0\de_{ext}
\eeaa
and then, using again Poincar\'e inequality, as well as the above, we obtain for $k\leq k_*-9$, 
\beaa
r^{-1}\left\|\widetilde{\dk}_*^{k}h\right\|_{\hk_2(\widetilde{S})} &\les& |\ov{\widetilde{\nu}^k(h)}|+\left\|\widetilde{\dk}_*^{k}\widetilde{\nab}h\right\|_{\hk_1(\widetilde{S})} \les \left\|\widetilde{\dk}_*^{k}\widetilde{\nab}h\right\|_{\hk_1(\widetilde{S})}+\ep_0\de_{ext}\\
&\les& r^{-1}\ep_0\left\|\widetilde{\dk}_*^{k}h\right\|_{\hk_1(\widetilde{S})}+\ep_0\de_{ext}.
\eeaa
For $\ep_0$ small enough, we infer, for $k\leq k_*-9$, 
\beaa
r^{-1}\left\|\widetilde{\dk}_*^{k}h\right\|_{\hk_2(\widetilde{S})} &\les& \ep_0\de_{ext}.
\eeaa
Using Sobolev, and recalling the definition of $h$, we infer, for $k\leq k_*-9$,
\beaa
\sup_{\widetilde{\Si}_*}\left|\widetilde{\dk}_*^{k}(\widetilde{b}_*+e_3(u)+e_3(s))\right| &\les& \ep_0\de_{ext}.
\eeaa

Next, we control $\widetilde{\nu}(\widetilde{r})+e_3(u)$. We have
 \beaa
\widetilde{\nu}(\widetilde{r})+e_3(u) &=& \widetilde{\nu}(\widetilde{r}-r)+\widetilde{\nu}(r)+e_3(u)\\
&=& \widetilde{\nu}(\widetilde{r}-r)+\widetilde{e}_3(r)+\widetilde{b}_*\widetilde{e}_4(r)+e_3(u)\\
&=& \widetilde{\nu}(\widetilde{r}-r)+e_3(r)+\widetilde{b}_*e_4(r)+e_3(u)+(\widetilde{e}_3(r)-e_3(r))+\widetilde{b}_*(\widetilde{e}_4(r)-e_4(r))\\
&=& \widetilde{\nu}(\widetilde{r}-r)+e_3(r-s)+e_3(s)+\widetilde{b}_*+\widetilde{b}_*e_4(r-s)+e_3(u)\\
&&+(\widetilde{e}_3(r)-e_3(r))+\widetilde{b}_*(\widetilde{e}_4(r)-e_4(r))\\
&=& \left(\widetilde{b}_*+e_3(s)+e_3(u)\right)+\widetilde{\nu}(\widetilde{r}-r)+e_3(r-s)+\widetilde{b}_*e_4(r-s)\\
&& +(\widetilde{e}_3(r)-e_3(r))+\widetilde{b}_*(\widetilde{e}_4(r)-e_4(r))
\eeaa
Together with 
\begin{itemize}
\item the above control of $\widetilde{b}_*+e_3(u)+e_3(s)$ on $\widetilde{\Si}_*$,

\item the above control of $\widetilde{r}-r$ on $\widetilde{\Si}_*$ and the fact that $\widetilde{\nu}$ is tangent to $\widetilde{\Si}_*$,

\item  the control \eqref{eq:controlofrminussonwidetildeRRhigherorderderivatives:ThmM7} of $r-s$ on $\widetilde{\RR}$, which yields the control of $e_3(r-s)$ and $e_4(r-s)$,

\item  the change of frame formulas, the control \eqref{eq:checkthatthereisindeedalargeruinnewlastGCMS:4:thereturn} for $(f, \fb, \la)$ and the control of Step 1 on $\MMextend$, which yields the control of $\widetilde{e}_3(r)-e_3(r)$ and $\widetilde{e}_4(r)-e_4(r)$,
\end{itemize}
we deduce,  for $k\leq k_*-9$,
\beaa
\sup_{\widetilde{\Si}_*}\left|\widetilde{\dk}_*^{k}(\widetilde{\nu}(\widetilde{r})+e_3(u))\right| &\les& \ep_0\de_{ext}.
\eeaa

The above estimates, together with \eqref{eq:Step14ThmM7:boundonffblatripleprime:bis}, yield for $\leq k_*-9$, 
\bea\lab{eq:Step14ThmM7:boundonffblatripleprime:ter}
\sup_{\widetilde{\Si}_*}r\left|\widetilde{\dk}_*^k(f''', \fb''', \la'''-1)\right| &\les&  \ep_0\de_{ext}  +\sup_{\widetilde{\Si}_*}|\widetilde{\dk}_*^k(\widetilde{a}\widetilde{f}_0-af_0)|.
\eea

{\bf Step 16.} Next, we focus on the control of the RHS of \eqref{eq:Step14ThmM7:boundonffblatripleprime:ter}. To this end, we first control $\widetilde{a}$. Recall that we have in view of \eqref{eq:controlofthesmallGCMquantitiesell=1modedivbetacurlbetaandkabcforThmM7}
\beaa
\nn\sup_{\widetilde{\RR}\cap\{u\geq u_*\}}r^5\left(\left|\frac{1}{|S|}\int_S\curl\b J^{(0)} -\frac{2am}{r^5}\right|+\frac{1}{|S|}\left|\int_S\curl\b J^{(\pm)}\right|\right)  &\les& \frac{\ep_0}{r}\Delta_{ext}\les\ep_0\de_{ext}.
\eeaa
Also, recall that $\ovS=S(\ug, \sg)$ is the sphere of the foliation of $\RR\cap\{u\geq u_*\}$ which shares the same south pole a $\widetilde{S}_*$. Relying on Corollary 5.9 in \cite{KS-GCM1} (see also Proposition \ref{proposition:58-59GCM1} here), we have, in view of the control \eqref{eq:controlofdeformationfunctionsUandS:ThmM7} of the deformation map $\Psi:\ovS\to \widetilde{S}_*$, for $p=0,+,-$, 
\beaa
&&\left|\int_{\widetilde{S}_*}\curl\b J^{(p)} -\int_{\ovS}\curl\b J^{(p)}\right|\\
 &\les& r\ep_0\de_{ext}\left(\sup_{\widetilde{\RR}}|\dkb^{\leq 1}(\curl\b J^{(p)})|+r\sup_{\widetilde{\RR}}(|\nab_3(\curl\b J^{(p)})|+|\nab_4(\curl\b J^{(p)})|)\right)\\
 &\les& \ep_0\de_{ext}\sup_{\widetilde{\RR}}\Big(|\dk^{\leq 2}\b|+r|\dk\nab_3\b|+r|\dk\b|(|\nab_3J^{(p)}|+|\nab_4J^{(p)}|\Big).
\eeaa
Together with the control of $\b$, the fact that $e_4(\Jp)=0$, and the control \eqref{eq:controlofthesmalle3JpforThmM7} for $e_3(\Jp)$, we deduce
\beaa
\max_{p=0,+,-}\left|\int_{\widetilde{S}_*}\curl\b J^{(p)} -\int_{\ovS}\curl\b J^{(p)}\right| &\les& \frac{\ep_0\de_{ext}}{r^3}.
\eeaa
Using the control of $\widetilde{r}_*-\rg$ in Proposition \ref{proposition:58-59GCM1}, we infer
\beaa
\max_{p=0,+,-}\widetilde{r}^5\left|\frac{1}{|\widetilde{S}_*|}\int_{\widetilde{S}_*}\curl\b J^{(p)} -\frac{1}{|\ovS|}\int_{\ovS}\curl\b J^{(p)}\right| &\les& \ep_0\de_{ext}.
\eeaa
Plugging in the above, using the fact that $\ovS\subset\widetilde{\RR}\cap\{u\geq u_*\}$, and using again the control of $\widetilde{r}_*-\rg$ in Proposition \ref{proposition:58-59GCM1}, we obtain 
\beaa
\widetilde{r}^5\left(\left|\frac{1}{|\widetilde{S}_*|}\int_{\widetilde{S}_*}\curl\b J^{(0)} -\frac{2am}{\widetilde{r}^5}\right|+\frac{1}{|\widetilde{S}_*|}\left|\int_{\widetilde{S}_*}\curl\b J^{(\pm)}\right|\right)  &\les& \ep_0\de_{ext}.
\eeaa
Together with \eqref{eq:checkthatthereisindeedalargeruinnewlastGCMS:9}, this yields
\beaa
\widetilde{r}^5\left(\left|\frac{1}{|\widetilde{S}_*|}\int_{\widetilde{S}_*}\curl\b J^{(0)} -\frac{2a\widetilde{m}}{\widetilde{r}^5}\right|+\frac{1}{|\widetilde{S}_*|}\left|\int_{\widetilde{S}_*}\curl\b J^{(\pm)}\right|\right)  &\les& \ep_0\de_{ext}.
\eeaa

Next, in view of the change of frame formulas, we have, using also the control of Step 1 on $\MMextend$,  
\beaa
\left|\widetilde{\curl}\widetilde{\b} - \curl\b\right| &\les& \left|\widetilde{\curl}\left(\la\left(\b+\frac{3}{2}(f\rho+\dual f\rhod)+\frac{1}{2}\a\c\fb+\cdots\right)\right) - \curl\b\right|\\
&\les& r^{-4}|\widetilde{\dkb}^{\leq 1}(f, \fb, \la-1)|
\eeaa
and hence, using the control \eqref{eq:checkthatthereisindeedalargeruinnewlastGCMS:4:thereturn} for $(f, \fb, \la)$, we deduce on $\widetilde{\Si}_*$
\beaa
\left|\widetilde{\curl}\widetilde{\b} - \curl\b\right| &\les& r^{-5}\ep_0\de_{ext}.
\eeaa
Plugging in the above, we infer
\beaa
\widetilde{r}^5\left(\left|\frac{1}{|\widetilde{S}_*|}\int_{\widetilde{S}_*}\widetilde{\curl}\widetilde{\b} J^{(0)} -\frac{2a\widetilde{m}}{\widetilde{r}^5}\right|+\frac{1}{|\widetilde{S}_*|}\left|\int_{\widetilde{S}_*}\widetilde{\curl}\widetilde{\b} J^{(\pm)}\right|\right)  &\les& \ep_0\de_{ext}.
\eeaa

Next, in view of Corollary 7.2 of \cite{KS-GCM2} (restated here as  Corollary \ref{Lemma:ComparisonJ-strong}),  there exists a canonical basis of $\ell=1$ modes on $\widetilde{S}_*$  in the sense of  Definition 3.10 of \cite{KS-GCM2} (recalled  here  in  Definition \ref{definition:ell=1mpdesonS-intro}), which we denote by $J_0^{(p,\widetilde{S}_*)}$, such that 
\beaa
\max_{p=0,+,-}\left\|J_0^{(p,\widetilde{S}_*)} -J^{(p)}\right\|_{\hk_{k_*-4}(\widetilde{S}_*)} &\les& \ep_0\de_{ext}. 
\eeaa
Also, recall that $\widetilde{J}^{(p)}=J^{(p, \widetilde{S}_*)}$ on $\widetilde{S}_*$, where $J^{(p, \widetilde{S}_*)}$ is in general another canonical basis of $\ell=1$ modes on $\widetilde{S}_*$. In view of Definition \ref{definition:ell=1mpdesonS-intro}, note that the canonical basis of $\ell=1$ modes on $\widetilde{S}_*$ are unique modulo isometries of $\SSS^2$, i.e. there exists  $O\in O(3)$ such that 
\bea\lab{eq:theuseoftheisometriyOofSSS2inthechoiceofthecanonicalbasiswidetildeS*:ThmM7}
J^{(p, \widetilde{S}_*)}=\sum_{q=0,+,-}O_{pq}J_0^{(q,\widetilde{S}_*)}, \qquad p=0,+,-.
\eea

\begin{remark}
In general, we have $O\neq I$ in \eqref{eq:theuseoftheisometriyOofSSS2inthechoiceofthecanonicalbasiswidetildeS*:ThmM7}. In fact, the role of $O$ corresponds in Step 6 to the application of Corollary \ref{Corr:ExistenceGCMS2} which ensures that the following holds on $\widetilde{S}_*$ w.r.t. the canonical basis of $\ell=1$ modes $J^{(p, \widetilde{S}_*)}$, see \eqref{eq:GCMconditionsinwidetildeStar:ThmM7}, 
\beaa
(\widetilde{\curl}\widetilde{\b})_{\ell=1,\pm}=0,
\eeaa
which corresponds to fixing the axis of $\widetilde{S}_*$. Note that this condition (and hence the axis of $\widetilde{S}_*$) is preserved by multiplying the basis $J^{(p, \widetilde{S}_*)}$ by $O=-I$ or by any  $O$ fixing $J^{(0, \widetilde{S}_*)}$, so that we may assume in \eqref{eq:theuseoftheisometriyOofSSS2inthechoiceofthecanonicalbasiswidetildeS*:ThmM7} that $O$ satisfies 
\bea\lab{eq:conditiononOpossibletorestrictfreedom:ThmM7}
O_{00}\geq 0, \qquad O_{++}\geq 0, \qquad O_{+-}=0,\qquad O_{--}\geq 0.
\eea
\end{remark}

Since $\widetilde{J}^{(p)}=J^{(p, \widetilde{S}_*)}$ on $\widetilde{S}_*$, we infer
\bea\lab{eq:firstcontroloftildeJminusJonwidetildeSstarwithrotation:ThmM7}
\max_{p=0,+,-}\left\|\widetilde{J}^{(p)} -\sum_{q=0,+,-}O_{pq}J^{(q)}\right\|_{\hk_{k_*-4}(\widetilde{S}_*)} &\les& \ep_0\de_{ext}.
\eea
where $O$ satisfies \eqref{eq:conditiononOpossibletorestrictfreedom:ThmM7}. Plugging \eqref{eq:firstcontroloftildeJminusJonwidetildeSstarwithrotation:ThmM7} in the above, we deduce
\beaa
\widetilde{r}^5\left(\left|\frac{1}{|\widetilde{S}_*|}\sum_{q=0,+,-}O_{0q}\int_{\widetilde{S}_*}\widetilde{\curl}\widetilde{\b}\widetilde{J}^{(q)} -\frac{2a\widetilde{m}}{\widetilde{r}^5}\right|+\frac{1}{|\widetilde{S}_*|}\left|\sum_{q=0,+,-}O_{\pm q}\int_{\widetilde{S}_*}\widetilde{\curl}\widetilde{\b} \widetilde{J}^{(p)}\right|\right)\les \ep_0\de_{ext}.
\eeaa
On the other hand, we have in view of \eqref{eq:GCMconditionsinwidetildeStar:ThmM7} and the fact that $\widetilde{J}^{(p)}=J^{(p, \widetilde{S}_*)}$ on $\widetilde{S}_*$
\beaa
\frac{1}{|\widetilde{S}_*|}\int_{\widetilde{S}_*}\widetilde{\curl}\widetilde{\b}\widetilde{J}^{(0)} \frac{2\widetilde{a}\widetilde{m}}{\widetilde{r}^5}, \qquad \frac{1}{|\widetilde{S}_*|}\int_{\widetilde{S}_*}\widetilde{\curl}\widetilde{\b}\widetilde{J}^{(\pm)}=0.
\eeaa
Plugging in the previous estimate, we obtain
\beaa
\widetilde{m}\Big(\left|O_{00}\widetilde{a} - a\right|+\left|\widetilde{a}\right|\left|O_{\pm 0}\right|\Big)  &\les& \ep_0\de_{ext},
\eeaa
and since $|\widetilde{m}-m|\les \de_{ext}\ep_0$ and $|m-m_0|\les\ep_0$, we may divide by $\widetilde{m}$ and hence
\bea\lab{eq:startinginequalitywidetildeaandawithrotationO:ThmM7}
\left|O_{00}\widetilde{a} - a\right|+\left|\widetilde{a}\right|\left|O_{+ 0}\right|+\left|\widetilde{a}\right|\left|O_{- 0}\right| &\les& \ep_0\de_{ext}.
\eea
Since $O\in O(3)$, we have $\sum_{p}O_{p0}^2=1$ which together with \eqref{eq:startinginequalitywidetildeaandawithrotationO:ThmM7} implies
\beaa
(\widetilde{a})^2 &=& (\widetilde{a})^2\left(\sum_{p}O_{p0}^2\right)= \left(a+O(\ep_0\de_{ext})\right)^2+O(\ep_0^2\de_{ext}^2) \les a^2+\ep_0^2\de_{ext}^2.
\eeaa
In particular, we infer the following estimate for $\widetilde{a}$
\bea\lab{eq:controlofwidetildeawhenaissmall:ThmM7}
|\widetilde{a}|\les \sqrt{\ep_0\de_{ext}}\quad\textrm{if}\quad|a|\leq \sqrt{\ep_0\de_{ext}}.
\eea
This allows us to control, in the case $|a|\leq \sqrt{\ep_0\de_{ext}}$, the change of frame coefficients $(f''', \fb''', \la''')$ introduced  in Step 14 and  corresponding to the change from  the outgoing PG frame $(e_4', e_3', e_1', e_2')$ of $\Mext$ extended to  the   spacetime  $\MM^{(extend)}$ and the  PG frame $(\widetilde{e}_4', \widetilde{e}_3', \widetilde{e}_1', \widetilde{e}_2')$. Indeed, \eqref{eq:controlofwidetildeawhenaissmall:ThmM7} and  \eqref{eq:Step14ThmM7:boundonffblatripleprime:ter} yield, for $\leq k_*-9$,
\bea\lab{eq:Step14ThmM7:boundonffblatripleprime:ter:firstcase}
\sup_{\widetilde{\Si}_*}r\left|\widetilde{\dk}_*^k(f''', \fb''', \la'''-1)\right| &\les&  \sqrt{\ep_0\de_{ext}}\quad\textrm{if}\quad |a|\leq \sqrt{\ep_0\de_{ext}}.
\eea

{\bf Step 17.} Next, we focus on controlling the RHS of \eqref{eq:Step14ThmM7:boundonffblatripleprime:ter} in the case $|a|>\sqrt{\ep_0\de_{ext}}$. In this case, we have from \eqref{eq:startinginequalitywidetildeaandawithrotationO:ThmM7} and the fact that $\sum_{p}O_{p0}^2=1$ 
\beaa
(\widetilde{a})^2 &=& (\widetilde{a})^2\left(\sum_{p}O_{p0}^2\right)= \left(a+O(\ep_0\de_{ext})\right)^2+O(\ep_0^2\de_{ext}^2) = a^2+O(\ep_0^{\frac{3}{2}}\de_{ext}^{\frac{3}{2}})
\eeaa
and hence $\widetilde{a}\geq \frac{1}{2}\sqrt{\ep_0\de_{ext}}$. Thus, dividing \eqref{eq:startinginequalitywidetildeaandawithrotationO:ThmM7} by $\widetilde{a}$, we obtain 
\beaa
\left|O_{00} - \frac{a}{\widetilde{a}}\right|+\left|O_{+ 0}\right|+\left|O_{- 0}\right| &\les& \sqrt{\ep_0\de_{ext}}.
\eeaa
Together with the fact that $\sum_{p}O_{p0}^2=1$, and recalling also that $O_{00}\geq 0$ in view of \eqref{eq:conditiononOpossibletorestrictfreedom:ThmM7}, we infer
\bea\lab{eq:controlofwidetildeaminusawhenaislarger:ThmM7}
|\widetilde{a}-a|\les \sqrt{\ep_0\de_{ext}}\quad\textrm{if}\quad|a|>\sqrt{\ep_0\de_{ext}}
\eea
and
\beaa
\left|O_{00} - 1\right|+\left|O_{+ 0}\right|+\left|O_{- 0}\right| &\les& \sqrt{\ep_0\de_{ext}}.
\eeaa
Also, since $O\in O(3)$, we also have
\beaa
0=\sum_pO_{p+}O_{p0}=O_{0+}+O(\sqrt{\ep_0\de_{ext}}), \qquad 0=\sum_pO_{p-}O_{p0}=O_{0-}+O(\sqrt{\ep_0\de_{ext}}),
\eeaa
and hence
\beaa
|O_{0+}|+|O_{0-}|\les \sqrt{\ep_0\de_{ext}}.
\eeaa
Together with the fact that $O_{+-}=0$ and $O_{--}\geq 0$ in view of \eqref{eq:conditiononOpossibletorestrictfreedom:ThmM7}, and since $\sum_{p}O_{p-}^2=1$, we infer
\beaa
\left|O_{--} - 1\right| &\les& \sqrt{\ep_0\de_{ext}}.
\eeaa
Finally $O_{++}\geq 0$ in view of \eqref{eq:conditiononOpossibletorestrictfreedom:ThmM7}, since we have obtained above that $|O_{0+}|\les \sqrt{\ep_0\de_{ext}}$, and since $\sum_{p}O_{p-}^2=1$ and $\sum_pO_{p+}O_{p-}=0$, we infer
\beaa
\left|O_{++} - 1\right|+\left|O_{-+}\right| &\les& \sqrt{\ep_0\de_{ext}}.
\eeaa
We have thus obtained
\beaa
|O-I| &\les& \sqrt{\ep_0\de_{ext}},
\eeaa
which together with \eqref{eq:firstcontroloftildeJminusJonwidetildeSstarwithrotation:ThmM7} implies
\beaa
\max_{p=0,+,-}r^{-1}\left\|\widetilde{J}^{(p)} -\Jp\right\|_{\hk_{k_*-4}(\widetilde{S}_*)} &\les& \sqrt{\ep_0\de_{ext}}\quad\textrm{if}\quad|a|>\sqrt{\ep_0\de_{ext}}.
\eeaa

Next, we control $\widetilde{J}^{(p)}-J^{(p)}$ for $p=0,+,-$ on $\widetilde{\Si}_*$. Recall that we have $\widetilde{\nu}(\widetilde{J}^{(p)})=0$ along $\widetilde{\Si}_*$. We infer
\beaa
\widetilde{\nu}\left(\widetilde{J}^{(p)}-\Jp\right) &=& -\widetilde{\nu}(\Jp)=-\widetilde{e}_4(\Jp)-\widetilde{b}_*\widetilde{e}_3(\Jp).
\eeaa
Using the change of frame formulas, and the fact that $e_4(\Jp)=0$, we obtain
\beaa
\widetilde{\nu}\left(\widetilde{J}^{(p)}-\Jp\right) &=& -\la\left(f\c\nab+\frac{1}{4}|f|^2e_3\right)\Jp\\
&& -\widetilde{b}_*\la^{-1}\left( \left(1+\frac{1}{2}f\c\fb  +\frac{1}{16} |f|^2  |\fb|^2\right) e_3 + \left(\fb+\frac 1 4 |\fb|^2\right)\c\nab\right)\Jp.
\eeaa
Together with the  control \eqref{eq:checkthatthereisindeedalargeruinnewlastGCMS:4:thereturn} for $(f, \fb, \la)$, and \eqref{eq:controlofthesmalle3JpforThmM7} for $e_3(\Jp)$, we infer 
\beaa
\sup_{\widetilde{\Si}_*}\left|\widetilde{\dk}_*^{k_*-5}\widetilde{\nu}\left(\widetilde{J}^{(p)}-\Jp\right)\right| &\les& \frac{\ep_0\de_{ext}}{ru_*}+\frac{\ep_0\de_{ext}}{r^2}.
\eeaa
Integrating along $\widetilde{\Si}_*$  from $\widetilde{S}_*$, using the above estimate on $\widetilde{S}_*$ and 
Sobolev, and using the dominant condition on $\widetilde{r}$ on $\widetilde{\Si}_*$, we infer
\bea
\max_{p=0,+,-}\sup_{\widetilde{\Si}_*}\left|\widetilde{\dk}_*^{k_*-6}\left(\widetilde{J}^{(p)}-\Jp\right)\right|  &\les& \sqrt{\ep_0\de_{ext}}\quad\textrm{if}\quad|a|>\sqrt{\ep_0\de_{ext}}. 
\eea

Next, we control $\widetilde{f}_0-f_0$ on $\widetilde{\Si}_*$. First, recall from Lemma \ref{lemma:computationforfirstorderderivativesJpell=1basis}  that the following identity holds on $S_*$
\beaa
\nab J^{(0)} = -\frac{1}{re^\phi}\dual f_0.
\eeaa
Next, recalling that we have extended $\phi$, $J^{(0)}$ and $f_0$ first to  $\Si^{(extend)}_*$ by $\nu(\phi)=0$, $\nab_\nu f_0=0$ and $\nu(J^{(0)})=0$, and then to $\widetilde{\RR}$ by $e_4(\phi)=0$, $\nab_4f_0=0$ and $e_4(J^{(0)})=0$, we have on $\Si^{(extend)}_*$
\beaa
\nab_\nu\left(r\left(\nab J^{(0)} +\frac{1}{re^\phi}\dual f_0\right)\right) &=& \Ga_b\dkb^{\leq 1}J^{(0)}
\eeaa
and on $\widetilde{\RR}$
\beaa
\nab_4\left(r\left(\nab J^{(0)} +\frac{1}{re^\phi}\dual f_0\right)\right) &=& \Ga_g\dkb^{\leq 1}J^{(0)}.
\eeaa
Integrating first from $S_*$, where the above identity holds, to $\Si^{(extend)}_*(u\geq u_*)$, and then from $\Si^{(extend)}_*(u\geq u_*)$ to $\widetilde{\RR}(u\geq u_*)$, we infer, using also the control of Step 1 on $\MMextend$,  for $k\leq k_*-3$,
\beaa
\sup_{\widetilde{\RR}(u\geq u_*)}\left|\dk^k\left(\nab J^{(0)} +\frac{1}{re^\phi}\dual f_0\right)\right| &\les& \frac{\ep_0}{r^3u_*^{\frac{1}{2}+\dec}}\De_{ext}+\frac{\ep_0}{r^2u_*^{1+\dec}}\de_{ext}\les \frac{\ep_0}{r^2u_*^{1+\dec}}\de_{ext}.
\eeaa
In particular, since the GCM sphere $\widetilde{S}_*$ constructed in Step 6 is included in $\widetilde{\RR}(u\geq u_*)$, we infer,  for $k\leq k_*-3$, 
\beaa
\sup_{\widetilde{S}_*}\left|\dk^k\left(\nab J^{(0)} +\frac{1}{re^\phi}\dual f_0\right)\right| &\les&  \frac{\ep_0}{r^2u_*^{1+\dec}}\de_{ext}.
\eeaa
Together with the control \eqref{eq:controlofwidetilderminusronSigmatildestar:ThmM7} for $\widetilde{r}-r$,  the change of frame formulas, the control \eqref{eq:checkthatthereisindeedalargeruinnewlastGCMS:4:thereturn} for $(f, \fb, \la)$ and the control of Step 1 on $\MMextend$, we infer, for $k\leq k_*-6$,
\beaa
\sup_{\widetilde{S}_*}\left|\widetilde{\dkb}^k\left(\widetilde{\nab} J^{(0)} +\frac{1}{\widetilde{r}e^\phi}\dual f_0\right)\right| &\les&  \frac{\ep_0}{r^2u_*^{1+\dec}}\de_{ext}.
\eeaa
Also, recall that $\ovS=S(\ug, \sg)$ is the sphere of the foliation of $\RR\cap\{u\geq u_*\}$ which shares the same south pole a $\widetilde{S}_*$.  As a byproduct of the construction of the GCM sphere\footnote{ based on  Theorem 7.3 and Corollary 7.7 of \cite{KS-GCM2}   (restated here  as  Theorem \ref{theorem:ExistenceGCMS2}   and Corollary  \ref{Corr:ExistenceGCMS2})  }  $\widetilde{S}_*$ in Step 6  we have, see  Corollary 7.2 of \cite{KS-GCM2} (stated here as Corollary  \ref{Lemma:ComparisonJ-strong}), the deformation map $\Psi:\ovS\to S$ is given by
\beaa
\Psi(\ug, \sg, x^1, x^2) &=& (\ug+U(x^1, x^2), \sg+S(x^1, x^2)),
\eeaa
where the scalar function $U$ and $S$ on $\ovS$ satisfy
\bea\lab{eq:controlofdeformationfunctionsUandS:ThmM7}
r^{-1}\|(U,S)\|_{\hk_{k_*-2}(\ovS)} &\les& \ep_0\de_{ext}. 
\eea 
The above control of the deformation $(U,S)$, together with the estimate  \eqref{eq:estimateformetricclosetounifforJnearcanonical1}, \eqref{eq:estimateformetricclosetounifforJnearcanonical2} and \eqref{eq:estimateformetricclosetounifforJnearcanonical3} that hold on any sphere $S$ of $\RR\cap\{u\geq u_*\}$, and hence in particular on $\ovS$, imply corresponding estimates with the metric $g$ on $\ovS$, and the scalar functions $\phi$ and $\Jp$ on $\ovS$ replaced with the corresponding analogs $(\Psi^{-1})^\#g$, $\phi\circ\Psi^{-1}$ and $\Jp\circ\Psi^{-1}$ on $\widetilde{S}_*$. We may thus apply   Proposition 4.15 in \cite{KS-GCM2}  (restated here as Proposition \ref{prop:asexpectedthebasiswidetildeJpisclosetocanonicalbasis}) which yields on $\widetilde{S}_*$ the estimate 
\beaa
r^{-1}\|\widetilde{\phi}-\phi\circ\Psi^{-1}\|_{\hk_{k_*-3}(\widetilde{S}_*)} \les  \frac{\ep_0}{r^2}\Delta_{ext}\les\frac{\ep_0}{r}\de_{ext}.
\eeaa
Together with the above form of $\Psi$, the above control of $(U, S)$, and the control of $\phi$ provided by \eqref{eq:estimateformetricclosetounifforJnearcanonical2}, we infer
\beaa
r^{-1}\|\widetilde{\phi}-\phi\|_{\hk_{k_*-3}(\widetilde{S}_*)} &\les& \frac{\ep_0}{r}\de_{ext}.
\eeaa 
Together with Sobolev, we deduce from the above, for $k\leq k_*-6$,
\beaa
\sup_{\widetilde{S}_*}\left|\widetilde{\dkb}^k\left(\widetilde{\nab} J^{(0)} +\frac{1}{\widetilde{r}e^{\widetilde{\phi}}}\dual f_0\right)\right| &\les&  \frac{\ep_0}{r^2}\de_{ext}.
\eeaa
On the other hand, in view of Lemma \ref{lemma:computationforfirstorderderivativesJpell=1basis} applied to the GCM sphere $\widetilde{S}_*$,  the following identity holds on $\widetilde{S}_*$
\beaa
\widetilde{\nab} \widetilde{J}^{(0)} = -\frac{1}{\widetilde{r}e^{\widetilde{\phi}}}\dual \widetilde{f}_0.
\eeaa
We infer, for $k\leq k_*-6$,
\beaa
\sup_{\widetilde{S}_*}\left|\widetilde{\dkb}^k\left(\widetilde{\nab}(\widetilde{J}^{(0)}-J^{(0)}) +\frac{1}{\widetilde{r}e^{\widetilde{\phi}}}\dual (\widetilde{f}_0-f_0)\right)\right| &\les&  \frac{\ep_0}{r^2}\de_{ext}.
\eeaa
Together with the above control of $\widetilde{J}^{(p)}-J^{(p)}$ on $\widetilde{\Si}_*$, we infer, for $k\leq k_*-7$,
\beaa
\sup_{\widetilde{S}_*}\left|\widetilde{\dkb}^k\left(\widetilde{f}_0-f_0\right)\right| &\les&  \sqrt{\ep_0\de_{ext}}.
\eeaa
Next, we control $\nab_3f_0$ on $\widetilde{\RR}$. Recall that we have $\nab_4f_0=0$ on $\widetilde{\RR}$ and $\nab_\nu f_0=0$ on $\Si^{(extend)}_*$. Since $\nu=e_3+b_*e_4$, we infer $\nab_3f_0=0$ on $\Si^{(extend)}_*$. Since $\nab_4f_0=0$ on $\widetilde{\RR}$, integrating from $\Si^{(extend)}_*$, we easily infer on $\widetilde{\RR}$, for $k\leq k_*-3$,
\beaa
\sup_{\widetilde{\RR}}\left|\dk^k\nab_3f_0\right| &\les& \frac{\ep_0}{r^2}\De_{ext}\les\frac{\ep_0}{r}\de_{ext}. 
\eeaa
Also, using the change of frame formulas and the fact that $\nab_4f_0=0$, we have on $\widetilde{\Si}_*$
\beaa
\widetilde{\nab}_{\widetilde{\nu}}f_0 &=& \widetilde{\nab}_{\widetilde{e}_3}f_0+\widetilde{b}_*\widetilde{\nab}_{\widetilde{e}_4}f_0\\
&=& \la^{-1}\left( \left(1+\frac{1}{2}f\c\fb  +\frac{1}{16} |f|^2  |\fb|^2\right) \nab_3 + \left(\fb+\frac 1 4 |\fb|^2f\right)\c\nab \right)f_0\\
&&+\widetilde{b}_*\la\left(f\c\nab +\frac 1 4 |f|^2  \nab_3\right)f_0.
\eeaa
Together with the above estimate for $\nab_3f_0$ and  the control \eqref{eq:checkthatthereisindeedalargeruinnewlastGCMS:4:thereturn} for $(f, \fb, \la)$, we infer, for $k\leq k_*-5$,
\beaa
\sup_{\widetilde{\Si}_*}\left|\dk^k\widetilde{\nab}_{\widetilde{\nu}}f_0\right| &\les& \frac{\ep_0}{r}\de_{ext}.
\eeaa
Since $\widetilde{\nab}_{\widetilde{\nu}}\widetilde{f}_0=0$, we deduce, for $k\leq k_*-5$,
\beaa
\sup_{\widetilde{\Si}_*}\left|\dk^k\widetilde{\nab}_{\widetilde{\nu}}\left(\widetilde{f}_0-f_0\right)\right| &\les& \frac{\ep_0}{r}\de_{ext}.
\eeaa
Integrating from $\widetilde{S}_*$ and using the above control on $\widetilde{S}_*$ for $\widetilde{f}_0-f_0$, as well as  the dominant condition for $r$ on $\widetilde{\RR}$, we deduce, for $k\leq k_*-7$,
\bea
\sup_{\widetilde{\Si}_*}\left|\dk^k\left(\widetilde{f}_0-f_0\right)\right| &\les& \sqrt{\ep_0\de_{ext}}\quad\textrm{if}\quad|a|>\sqrt{\ep_0\de_{ext}}.
\eea

We are now in position to control, in the case $|a|>\sqrt{\ep_0\de_{ext}}$, the change of frame coefficients $(f''', \fb''', \la''')$ introduced  in Step 14 and  corresponding to the change from  the outgoing PG frame $(e_4', e_3', e_1', e_2')$ of $\Mext$ extended to  the   spacetime  $\MM^{(extend)}$ and the  PG frame $(\widetilde{e}_4', \widetilde{e}_3', \widetilde{e}_1', \widetilde{e}_2')$. Recall \eqref{eq:Step14ThmM7:boundonffblatripleprime:ter} that holds for $\leq k_*-9$
\beaa
\sup_{\widetilde{\Si}_*}r\left|\widetilde{\dk}_*^k(f''', \fb''', \la'''-1)\right| &\les&  \ep_0\de_{ext}  +\sup_{\widetilde{\Si}_*}|\widetilde{\dk}_*^k(\widetilde{a}\widetilde{f}_0-af_0)|.
\eeaa
In view of the above estimates for $\widetilde{a}-a$ and $\widetilde{f}_0-f_0$, we infer, $\leq k_*-9$, 
\bea\lab{eq:Step14ThmM7:boundonffblatripleprime:ter:secondcase}
\sup_{\widetilde{\Si}_*}r\left|\widetilde{\dk}_*^k(f''', \fb''', \la'''-1)\right| &\les&  \sqrt{\ep_0\de_{ext}}\quad\textrm{if}\quad|a|>\sqrt{\ep_0\de_{ext}}.
\eea
 
We conclude this step with the control of $\widetilde{\Jk}-\Jk$ on $\widetilde{\Si}_*$ in the case $|a|>\sqrt{\ep_0\de_{ext}}$. Recall from Step 14 that  $(f', \fb', \la')$ denote the  change of frame coefficients introduced  in Step 14 and  corresponding to the change from the outgoing geodesic frame $(e_4, e_3, e_1, e_2)$ of $\MM^{(extend)}$ to the outgoing PG frame $(e_4', e_3', e_1', e_2')$ of $\Mext$ extended to  the   spacetime  $\MM^{(extend)}$. Recall from \eqref{eq:relationbetweenJkandf0onSigmastar} and \eqref{eq:basicdefinitionofJkonSstarSigmastarMextstar}  that $\Jk$ satisfies in $\MM^{(extend)}$, and hence in $\widetilde{\RR}$, the following identities
\beaa
\Jk=\frac{1}{|q|}\left(f_0+i\dual f_0\right)\quad\text{on}\quad\Si^{(extend)}_*, \qquad \nab_{e_4'}\Jk = -\frac{1}{q}\Jk\quad\text{on}\quad\widetilde{\RR}.
\eeaa
Using the change of frame formulas, we have on $\widetilde{\RR}$
\beaa
\nab_{e_4'}\Jk &=& \nab_{\la'\left(e_4+f'\c\nab+\frac{1}{4}|f'|^2e_3\right)}\Jk\\
&=& \nab_4\Jk +(\la'-1)\nab_4\Jk+\la' f'\c \nab\Jk+\frac{\la'}{4}|f'|^2\nab_3\Jk
\eeaa
and hence
\beaa
\nab_4\Jk -\frac{1}{q}\Jk &=& -(\la'-1)\nab_4\Jk-\la' f'\c \nab\Jk-\frac{\la'}{4}|f'|^2\nab_3\Jk.
\eeaa
Recall also that we have derived the following control for $(f', \la')$ in Step 14, for $k\leq k_*-3$,
\beaa
\sup_{\widetilde{\RR}}r\left|\dk^k\left(f'-\frac{a}{r}f_0, \la' -1\right)\right| &\les&  \ep_0\de_{ext}.
\eeaa
Together with the control of $\Jk$ in $\MM^{(extend)}$, we infer, for $k\leq k_*-3$,
\beaa
\sup_{\widetilde{\RR}}r^3\left|\dk^k\left(\nab_4\Jk -\frac{1}{q}\Jk\right)\right| &\les& \ep_0\de_{ext}. 
\eeaa
Also, recall that we have $\nab_4f_0=0$ in $\widetilde{\RR}$ and hence
\beaa
\nab_4\left(|q|\Jk-\left(f_0+i\dual f_0\right)\right) &=& |q|\left(\nab_4 +\frac{\nab_4(|q|)}{|q|}\right)\Jk=|q|\left(\nab_4 +\frac{r}{|q|^2}\right)\Jk\\
&=& |q|\left(\nab_4 +\frac{1}{q}\right)\Jk -\frac{ia\cos\th}{|q|}\Jk.
\eeaa
Together with the above estimate for $\nab_4\Jk -\frac{1}{q}\Jk$, we infer, for $k\leq k_*-3$,
\beaa
\sup_{\widetilde{\RR}}r^2\left|\dk^k\nab_4\left(|q|\Jk-\left(f_0+i\dual f_0\right)\right)\right| &\les& 1.
\eeaa
Integrating from $\Si^{(extend)}_*$ where $|q|\Jk-(f_0+i\dual f_0)=0$, we infer, for $k\leq k_*-3$,
\beaa
\sup_{\widetilde{\RR}}\left|\dk^k\left(|q|\Jk-\left(f_0+i\dual f_0\right)\right)\right| &\les& \frac{1}{r^2}\De_{ext}\les \frac{1}{r}\de_{ext}. 
\eeaa
Together with the dominant condition for $r$ on $\widetilde{\RR}$, we obtain, for $k\leq k_*-3$,
\beaa
\sup_{\widetilde{\RR}}r\left|\dk^k\left(\Jk-\frac{1}{|q|}\left(f_0+i\dual f_0\right)\right)\right| &\les& \ep_0\de_{ext}. 
\eeaa
On the other hand, we have on $\widetilde{\Si}_*$
\beaa
\widetilde{\Jk} &=& \frac{1}{|\widetilde{q}|}\left(\widetilde{f}_0+i\dual \widetilde{f}_0\right). 
\eeaa
Since $\widetilde{\Si}_*\subset\widetilde{\RR}$, this yields, for $k\leq k_*-3$,
\beaa
\sup_{\widetilde{\Si}_*}r\left|\widetilde{\dk}_*^k\left(\widetilde{\Jk}-\Jk\right)\right| &\les& \ep_0\de_{ext}+\sup_{\widetilde{\Si}_*}r\left|\widetilde{\dk}_*^k\left(\frac{1}{|\widetilde{q}|}\left(\widetilde{f}_0+i\dual \widetilde{f}_0\right) - \frac{1}{|q|}\left(f_0+i\dual f_0\right)\right)\right|\\
&\les& \ep_0\de_{ext}+\sup_{\widetilde{\Si}_*}\left(\left|\widetilde{\dk}_*^k\left(\widetilde{f}_0-f_0\right)\right|+r^{-1}\left|\widetilde{\dk}_*^k\left(\widetilde{r}-r\right)\right|+r^{-1}\left|\widetilde{\dk}_*^k\left(\widetilde{J}^{(0)} - J^{(0)}\right)\right|\right).
\eeaa
In view of the estimates for $\widetilde{r}-r$ in Step 15, and the above estimates for $\widetilde{f}_0-f_0$ and $\widetilde{J}^{(0)} - J^{(0)}$ in the case $|a|>\sqrt{\ep_0\de_{ext}}$, we infer,  for $\leq k_*-7$, 
\bea
\sup_{\widetilde{\Si}_*}r\left|\widetilde{\dk}_*^k\left(\widetilde{\Jk}-\Jk\right)\right| &\les& \sqrt{\ep_0\de_{ext}}\quad\textrm{if}\quad|a|>\sqrt{\ep_0\de_{ext}}.
\eea

Finally, we have obtained in this step,  for $\leq k_*-9$, 
\beaa
\sup_{\widetilde{\Si}_*}\Bigg(r\left|\widetilde{\dk}_*^k(f''', \fb''', \la'''-1)\right| +\left|\widetilde{a}-a\right|+\left|\widetilde{\dk}_*^k\left(\widetilde{J}^{(0)}- J^{(0)}\right)\right|\\
 +r\left|\widetilde{\dk}_*^k\left(\widetilde{\Jk}-\Jk\right)\right|\Bigg) &\les&  \sqrt{\ep_0\de_{ext}}\quad\textrm{if}\quad|a|>\sqrt{\ep_0\de_{ext}}.
\eeaa
Together with the control \eqref{eq:checkthatthereisindeedalargeruinnewlastGCMS:9} for $\widetilde{r}-{r}$ and $\widetilde{m}-{m}$, \eqref{eq:controlofwidetildeawhenaissmall:ThmM7} and \eqref{eq:Step14ThmM7:boundonffblatripleprime:ter:firstcase} for the case $|a|\leq\sqrt{\ep_0\de_{ext}}$, and 
possibly reducing the size of $\de_{ext}>0$ (which can be chosen arbitrarily small), we deduce,  for $\leq k_*-9$, 
\bea\lab{eq:controlffblafinalonSigma*widetildeforThmM7}
\nn\sup_{\widetilde{\Si}_*}\widetilde{u}^{1+\dec}\Bigg(r\left|\widetilde{\dk}_*^k(f''', \fb''', \la'''-1)\right| +\left|\widetilde{\dk}_*^k\left(\widetilde{r}-r\right)\right|+\left|\widetilde{m}-m\right|\\
+\left|\widetilde{a}-a\right|+\left|\widetilde{\dk}_*^k\left(\widetilde{a}\widetilde{J}^{(0)}- aJ^{(0)}\right)\right| +r\left|\widetilde{\dk}_*^k\left(\widetilde{a}\widetilde{\Jk}-a\Jk\right)\right|\Bigg) &\les&  \ep_0.
\eea

{\bf Step 18.} We now control the outgoing PG structure  initialized on $\widetilde{\Si}_*$. We denote by $\,^{(ext)}\widetilde{\MM}$ the region covered by this outgoing PG structure. For convenience, we change our notation. From   now on:
\begin{itemize}
\item $(e_4, e_3, e_1, e_2)$ denotes the outgoing PG frame   of $\Mext$ extended to the   spacetime  $\MM^{(extend)}$,

\item $(\widetilde{e}_4, \widetilde{e}_3, \widetilde{e}_1, \widetilde{e}_2)$ denotes the outgoing PG frame initialized on $\widetilde{\Si}_*$,

\item $(f,\fb, \la)$ denote the transition coefficients from the PG frame $(e_4, e_3, e_1, e_2)$ to the PG frame $(\widetilde{e}_4, \widetilde{e}_3, \widetilde{e}_1, \widetilde{e}_2)$.
\end{itemize}
In view of \eqref{eq:controlffblafinalonSigma*widetildeforThmM7}, using the above new notations for $(f,\fb, \la)$, and noticing that the structure equations in the $\widetilde{e}_4$ direction for the outgoing PG structure initialized on $\widetilde{\Si}_*$  allow to recover the $\widetilde{e}_4$ derivatives (which are transversal to $\widetilde{\Si}_*$), we have, for $\leq k_*-9$,   
\bea\lab{eq:controlffblafinalonSigma*widetildeforThmM7:bis}
\nn\sup_{\widetilde{\Si}_*}\widetilde{u}^{1+\dec}\Bigg(r\left|\dk^k(f, \fb, \la-1)\right| +\left|\dk^k\left(\widetilde{r}-r\right)\right|+\left|\widetilde{m}-m\right|\\
+\left|\widetilde{a}-a\right|+\left|\dk^k\left(\widetilde{a}\widetilde{J}^{(0)}- aJ^{(0)}\right)\right| +r\left|\dk^k\left(\widetilde{a}\widetilde{\Jk}-a\Jk\right)\right|\Bigg) &\les&  \ep_0.
\eea
Also $(e_4, e_3, e_1, e_2)$, as discussed in Step 1 to Step 3, satisfies
\bea\lab{eq:extendedspacetimebootstrapassumptionsforThmM7}
 \Nk^{(Dec)}_{k_*-3}(\MMextend) &\les&  \ep_0.
\eea

We introduce the notations
\beaa
F:=f+i\dual f, \qquad \underline{F}:=\fb+i\dual \fb.
\eeaa
Since $(e_4, e_3, e_1, e_2)$ and $(\widetilde{e}_4, \widetilde{e}_3, \widetilde{e}_1, \widetilde{e}_2)$ are outgoing PG frames, we have 
\beaa
\Xi=0, \qquad \om=0, \qquad \Hb+Z=0, \qquad \widetilde{\Xi}=0, \qquad \widetilde{\om}=0, \qquad \widetilde{\Hb}+\widetilde{Z}=0.
\eeaa 
In view of Corollary \ref{cor:transportequationine4forchangeofframecoeffinformFFbandlamba}, we have the following transport equations
\beaa
\nab_{\la^{-1}\widetilde{e}_4}\left(qF\right) &=&  E_4(f, \Ga),\\
\la^{-1}\nab_{\widetilde{e}_4}(\log\la) &=& 2f\c\ze+E_2(f, \Ga),\\
\nab_{\la^{-1}\widetilde{e}_4}\left[q\Big(\underline{F} -2q\widetilde{\DD}(\log\la)+e_3(r) F\Big)\right] &=&  -3q^2\widetilde{\DD}\left(f\c\ze\right) +E_5({\widetilde{\nab}}^{\leq 1}f,\fb, {\widetilde{\nab}}^{\leq 1}\la, \D^{\leq 1}\Ga),
\eeaa
where  $E_2$, $E_4$ and $E_5$ are given in  Corollary \ref{cor:transportequationine4forchangeofframecoeffinformFFbandlamba}. Integrating these transport equations from $\widetilde{\Si}_*$ in the order they appear, using the control in \eqref{eq:controlffblafinalonSigma*widetildeforThmM7:bis} for $(f, \fb, \la)$ on $\widetilde{\Si}_*$, and together with the control \eqref{eq:extendedspacetimebootstrapassumptionsforThmM7} for the Ricci coefficients of the foliation of $\MM^{(extend)}$, we obtain, for $\leq k_*-9$,    
\begin{equation}\lab{eq:controlffblaonMextwidetildeThmM7}
\sup_{\,^{(ext)}\widetilde{\MM}}\Big(\widetilde{r}\,\widetilde{u}^{\frac{1}{2}+\dec}+\widetilde{u}^{1+\dec}\Big)\Big(|\dk^k(f,  \log(\la))|+|\dk^{k-1}\fb|\Big) \les \ep_0.
\end{equation}

Also, we have
\beaa
\widetilde{e}_4(\widetilde{r}-\la^{-1}r) &=& 1-\la^{-1}\widetilde{e}_4(r)+\la^{-1}\widetilde{e}_4(\log(\la)).
\eeaa
Using the change of frame formula and the above transport equation for $\log(\la)$, we infer
\beaa
\widetilde{e}_4(\widetilde{r}-\la^{-1}r) &=& 1-\left(e_4+f\c\nab+\frac{1}{4}|f|^2e_3\right)r+\frac{3}{2}f\c\ze+E_2(f, \Ga)\\ 
&=& -\frac{1}{4}|f|^2e_3(r)+\frac{3}{2}f\c\ze+E_2(f, \Ga).
\eeaa
Integrating from $\widetilde{\Si}_*$ where $\widetilde{r}-r$ is under control in view of \eqref{eq:controlffblafinalonSigma*widetildeforThmM7:bis}, and using the control \eqref{eq:controlffblaonMextwidetildeThmM7} for $f$ and $\la$ as well as the control \eqref{eq:extendedspacetimebootstrapassumptionsforThmM7} for the foliation of $\MM^{(extend)}$, we infer, for $\leq k_*-9$, 
\bea\lab{eq:controlffblaonMextwidetildeThmM7:rtilde}
\sup_{\,^{(ext)}\widetilde{\MM}}\widetilde{u}^{1+\dec}\left|\dk^k(\widetilde{r}-r)\right| &\les& \ep_0. 
\eea
Also, we have
\beaa
\widetilde{e}_4(\widetilde{a}\widetilde{J}^{(0)}-aJ^{(0)}) &=& -\widetilde{e}_4(aJ^{(0)})=-a\la\left(e_4+f\c\nab+\frac{1}{4}|f|^2e_3\right)J^{(0)}\\
&=& -a\la\left(f\c\nab+\frac{1}{4}|f|^2e_3\right)J^{(0)}
\eeaa
and
\beaa
\nab_{\widetilde{e}_4}(\widetilde{a}\widetilde{q}\widetilde{\Jk}- aq\Jk) &=& -\nab_{\widetilde{e}_4}( aq\Jk)=-a\la\left(\nab_4+f\c\nab+\frac{1}{4}|f|^2\nab_3\right)(q\Jk)\\
&=& -a\la\left(f\c\nab+\frac{1}{4}|f|^2\nab_3\right)(q\Jk).
\eeaa
Integrating from $\widetilde{\Si}_*$ where $\widetilde{a}\widetilde{J}^{(0)}-aJ^{(0)}$ and $\widetilde{a}\widetilde{\Jk}- a\Jk$ are under control in view of \eqref{eq:controlffblafinalonSigma*widetildeforThmM7:bis}, and using the control \eqref{eq:controlffblaonMextwidetildeThmM7} for $f$ and $\la$ as well as the control \eqref{eq:extendedspacetimebootstrapassumptionsforThmM7} for the foliation of $\MM^{(extend)}$, we infer, for $\leq k_*-9$, 
\bea\lab{eq:controlffblaonMextwidetildeThmM7:JandJktilde}
\sup_{\,^{(ext)}\widetilde{\MM}}\widetilde{u}^{1+\dec}\left(\left|\dk^k(\widetilde{a}\widetilde{J}^{(0)}-aJ^{(0)})\right|+r\left|\dk^k(\widetilde{a}\widetilde{\Jk}-a\Jk)\right|\right) &\les& \ep_0. 
\eea

Then, using the outgoing PG structure of $\,^{(ext)}\widetilde{\MM}$, we initialize
\begin{itemize} 
\item the ingoing PG structure  of $\,^{(int)}\widetilde{\MM}$ on $\TT=\{\widetilde{r}=r_0\}$,

\item  the  ingoing PG structure  of $\,^{(top)}\widetilde{\MM}$ on $\{\widetilde{u}=\widetilde{u}_*\}$,
\end{itemize} 
as in section \ref{sec:initalizationadmissiblePGstructure}. Using the control of $(f, \fb, \la)$, $\widetilde{r}-r$, $\widetilde{a}\widetilde{J}^{(0)}-aJ^{(0)}$ and $\widetilde{a}\widetilde{\Jk}-a\Jk$ induced on $\{\widetilde{r}=r_0\}$ and $\{\widetilde{u}=\widetilde{u}_*\}$ by \eqref{eq:controlffblaonMextwidetildeThmM7}, \eqref{eq:controlffblaonMextwidetildeThmM7:rtilde} and \eqref{eq:controlffblaonMextwidetildeThmM7:JandJktilde}, 
and using the analog in the $\widetilde{e}_3$ direction for ingoing PG structures of the above transport equation in the $\widetilde{e}_4$ direction for outgoing PG structures,  we obtain for $\,^{(int)}\widetilde{\MM}$ and $\leq k_*-10$
\bea\lab{eq:controlffblaonMintwidetildeThmM7}
\nn\sup_{\,^{(int)}\widetilde{\MM}}\widetilde{\ub}^{1+\dec}\Big(|\dk^k(\fb,  \log(\la))|+|\dk^{k-1}f|+|\dk^k(\widetilde{r}-r)|\\
+|\dk^k(\widetilde{a}\widetilde{J}^{(0)}-aJ^{(0)})|+|\dk^k(\widetilde{a}\widetilde{\Jk}-a\Jk)| \Big)&\les& \ep_0,
\eea
and a similar estimate for $\,^{(top)}\widetilde{\MM}$. 

Let now
\beaa
\widetilde{\MM} &:=& \,^{(ext)}\widetilde{\MM}\cup \,^{(int)}\widetilde{\MM}\cup \,^{(top)}\widetilde{\MM}.
\eeaa
Then, in view of \eqref{eq:controlffblaonMextwidetildeThmM7}-\eqref{eq:controlffblaonMintwidetildeThmM7}, the control of $\widetilde{a}-a$ and $\widetilde{m}-m$ in \eqref{eq:controlffblafinalonSigma*widetildeforThmM7:bis}, and  \eqref{eq:extendedspacetimebootstrapassumptionsforThmM7}, and using the transformation formulas of Proposition \ref{Proposition:transformationRicci}, and well as the definition of the linearized quantities based on $a$, $m$, $r$, $aJ^{(0)}=a\cos\th$ and $a\Jk$, we deduce
\beaa
 \Nk^{(Dec)}_{k_*-12}(\widetilde{\MM}) &\les&  \ep_0.
\eeaa
In particular, since  $k_*= k_{small}+20$, we infer
\beaa
 \Nk^{(Dec)}_{k_{small}}(\widetilde{\MM}) &\les&  \ep_0
\eeaa
which concludes the proof of Theorem M7.

%%%%%%%%%%%%%%%%%%%%%%%%%%%%%%%%%%%%%%

%%%%%%%%%%%%%%%%%%%%%%%%%%%%%%%%%%%%%%

%%%%%%%%%%%%%%%%%%%%%%%%%%%%%%%%%%%%%%

\chapter{Top order estimates (Theorem M8)}

%%%%%%%%%%%%%%%%%%%%%%%%%%%%%%%%%%%%%%

The goal of this chapter is to prove Theorem M8, i.e. to improve our bootstrap assumptions on boundedness on $\Si_*$, and for the PG structures of $\Mext$, $\Mint$ and $\Mtop$. Now, while the PG structures we have studied so far are perfectly adequate  for  deriving  decay estimates,   they are  deficient in terms of loss of derivatives and thus inadequate   for deriving  boundedness estimates for the top derivatives of the Ricci coefficients. Hence, we cannot rely on PG structures in the proof of Theorem M8, and will instead rely on the PT structures introduced in section \ref{sec:introductionofPTstrctures}. Once boundedness estimates for top order derivatives are obtained for the PT structures, see Theorem \ref{theorem:Main-PT}, they will induce the improvement of our bootstrap assumptions on boundedness on $\Si_*$, and for the PG structures of $\Mext$, $\Mint$ and $\Mtop$, see section \ref{sec:endoftheproofofTheoremM8:chap9}. 

\begin{remark}
Recall that we in fact only prove part of Theorem M8 in this paper. Indeed, the estimates for the curvature components in the PT frame, see Theorem \ref{prop:rpweightedestimatesiterationassupmtionThM8}, will be done in the separate paper \cite{KS:Kerr-B}. 
\end{remark}

%%%%%%%%%%%%%%%%%%%%%%%%%%%

\section{Principal temporal structures in $\MM$}
\lab{section:introduction-PTframes}

%%%%%%%%%%%%%%%%%%%%%%%%%%%%%%

Before introducing the temporal structure of $\MM$, we recall the main definitions of outgoing and ingoing temporal structures, see section \ref{sec:introductionofPTstrctures}.

%%%%%%%%%%%%%%%%%%%%%%%%%%%

\subsection{Outgoing PT structures}

%%%%%%%%%%%%%%%%%%%%%%%%%%%%%%

We recall Definitions \ref{definition:outgoingPT} and \ref{def:initialoutgoingPTdataset} and Lemma \ref{lemma:constructionoutgoingPTframes} of section \ref{sec:outgoingPTstructures:chap2}. 
\begin{definition}
\lab{definition:outgoingPT:chap9}
An outgoing  PT structure    $\{ (e_3, e_4, \HH), r, \th, \Jk\}$     on $\MM$ consists of  a null pair $(e_3, e_4)$, the induced horizontal structure   $\HH$,   functions $(r, \th)$, and a horizontal  1-form $\Jk$ such that the following hold true:
\begin{enumerate}
\item   $e_4$ is geodesic.

\item We have
\bea
e_4(r)=1,\qquad    e_4(\th)=0, \qquad \nab_4 (q \Jk)=0, \qquad  q= r+a i \cos \th.
\eea

\item We have 
\bea
\Hb=-\frac{a\ov{q}}{|q|^2}\Jk.
\eea
\end{enumerate}
An extended outgoing PT structure possesses,  in addition,  a scalar function $u$ verifying $e_4(u)=0$. 
\end{definition}

\begin{definition}\lab{def:initialoutgoingPTdataset:chap9}
An outgoing PT  initial data set  consists of a hypersurface $\Si$ transversal to $e_4$ together with a null pair $(e_3, e_4)$, the induced horizontal structure $\HH$,  scalar   functions $(r, \th)$,  and a horizontal  1-form $\Jk$,  all defined on $\Si$.
\end{definition}

\begin{lemma}
\lab{lemma:constructionoutgoingPTframes:chap9}
Any outgoing PT  initial data set, as in Definition \ref{def:initialoutgoingPTdataset}, can be  locally extended to   an outgoing PT structure.
\end{lemma} 

%%%%%%%%%%%%%%%%%%%%%%%%%%%

\subsection{Ingoing PT structures}

%%%%%%%%%%%%%%%%%%%%%%%%%%%%%%

We recall Definitions \ref{definition:outgoingPT} and \ref{def:initialoutgoingPTdataset} and Lemma \ref{lemma:constructionoutgoingPTframes} of section \ref{sec:ingoingPTstructures:chap2}.  

\begin{definition}
\lab{definition:ingoingPT:chap9}
An ingoing  PT structure    $\{ (e_3, e_4, \HH), r, \th, \Jk\}$     on $\MM$ consists of  a null pair $(e_3, e_4)$, the induced horizontal structure   $\HH$,  functions $(r, \th)$, and a horizontal  1-form $\Jk$ such that the following hold true:
\begin{enumerate}
\item   $e_3$ is geodesic.

\item We have
\bea
e_3(r)=-1,\qquad    e_3(\th)=0, \qquad \nab_3 (\ov{q}\Jk)=0, \qquad  q= r+a i \cos \th.
\eea

\item We have 
\bea
H=\frac{aq}{|q|^2}\Jk.
\eea
\end{enumerate}
An extended  ingoing PT structure possesses,  in addition,  a function $\ub$ verifying $e_3(\ub)=0$. 
\end{definition}

\begin{definition}\lab{def:initialingoingPTdataset:chap9}
An ingoing PT  initial data set  consists of a hypersurface $\Si$ transversal to $e_3$ together with a null pair $(e_3, e_4)$, the induced horizontal structure $\HH$,  scalar   functions $(r, \th)$,  and a horizontal  1-form $\Jk$,  all defined on $\Si$.
\end{definition}

\begin{lemma}
\lab{lemma:constructioningoingPTframes:chap9}
Any ingoing PT  initial data set, as in Definition \ref{def:initialingoingPTdataset}, can be  locally extended to   an ingoing PT structure.
\end{lemma}

%%%%%%%%%%%%%%%%%%%%%%%%%%%%%

\subsection{Definition of the PT structures in $\MM$}
\lab{sec:defintionofthePTstructuresinMM}

%%%%%%%%%%%%%%%%%%%%%%%%%%%%%

Let $\MM$ denote  our admissible GCM spacetime  introduced in section \ref{section:GCMadmissible-spacetimes}. We decompose $\MM$ as follows
\beaa
\MM &=& \Mext\cup\Mint'\cup \,{}^{(top)}\MM',
\eeaa
where
\begin{itemize}
\item $\Mext$ is covered by an outgoing  PT structure initialized on $\Si_*$ as will be made precise below, 

\item $\Mint'$ is covered by an  ingoing  PT structure initialized on $\TT=\{r=r_0\}$ as will be made precise below,

\item $\Mtop'$ is covered by an  ingoing  PT structure initialized on $\{ u=u_*'\} $ for some $u_*'\in[u_*-2, u_*-1]$ as will be made precise below. 
\end{itemize}  

\begin{remark}
\lab{remark:outgoingPT-versiusPG}
Recall that $\Mext$ is covered by an outgoing PG structure. A priori, the outgoing PT structure introduced above covers a different region $\Mext'$. The fact that these two regions coincide, i.e. $\Mext'=\Mext$,  is due to the fact that $\Mext$ is defined purely in terms of $(u,r)$, and that the functions $(u,r)$ for the outgoing PT frame coincide with the ones  of the outgoing PG frame  in $\Mext$, see Lemma \ref{lemma:linkPGandTframeinMext}. 
\end{remark}

\begin{remark}
The constant $u_*'$ involved in the definition $\Mtop'$ and its associated ingoing PT structure  will be fixed  in section \ref{sec:choiceofuprimstar}. 
\end{remark}

%%%%%%%%%%%%%%%%%%%%%%%%%%%%%%%%%%%

\subsubsection{Initialization of the outgoing PT structure of $\Mext$}

%%%%%%%%%%%%%%%%%%%%%%%%%%%%%%%%%%%

Recall from section \ref{sec:admissibleGMCPGdatasetonSigmastar} that the GCM boundary $\Si_*$ comes together with a null frame,  as well as scalar functions $(u, r, \th)$ all defined on $\Si_*$. Also, recall  the definition of the 1-form $f_0$ on $\Si_*$, see Definition \ref{def:definitionoff0fplusfminus},
\beaa
(f_0)_1=0, \quad (f_0)_2=\sin\th, \quad \textrm{on}\quad S_*, \qquad \nab_\nu f_0=0\quad \textrm{on}\quad\Si_*,
\eeaa
where, on $S_*$, we consider the orthonormal basis $(e_1, e_2)$ of $S_*$ given by \eqref{eq:specialorthonormalbasisofSstar}.  

Next, recall that the outgoing PT structure of $\Mext$ is initialized on $\Si_*$. In view of Lemma \ref{lemma:constructionoutgoingPTframes:chap9}, it suffices to prescribe the corresponding outgoing PT  initial data set, as in Definition \ref{def:initialoutgoingPTdataset:chap9}, on $\Si_*$. We make the following choice of outgoing PT  initial data set $\{ (e_3, e_4, \HH), r, \th, \Jk\}$ on $\Si_*$:
\begin{enumerate}
\item The frame $(e_1, e_2, e_3, e_4)$, consisting of the null pair $(e_3, e_4)$ and the horizontal structure $\HH$, is obtained from the null frame attached to $\Si_*$ by the change of frame formula with frame coefficients $(f, \fb, \la)$ given by 
\bea
\la=1, \qquad f=\frac{a}{r}f_0, \qquad  \fb=\frac{a\Upsilon}{r}f_0.
\eea

\item The functions $(r, \th)$ coincide on $\Si_*$ with the ones of $\Si_*$. 

\item The complex horizontal 1-form $\Jk$ is given on $\Si_*$ by
\bea
\Jk &=& \frac{1}{r}(f_0+i\dual f_0). 
\eea
\end{enumerate}

\begin{remark}\lab{rmk:comparisionchangeofframePGandPTonSigmastar}
Let $(\widetilde{f}, \widetilde{\fb}, \widetilde{\la})$ denote the coefficients involved in the initialization of the outgoing PG frame of $\Mext$ on $\Si_*$, i.e.  corresponding to the change  from the frame of $\Si_*$ to the outgoing PG frame of $\Mext$. Recall from section \ref{section:InitializationofPG structure-Mext} that $(\widetilde{f}, \widetilde{\fb}, \widetilde{\la})$ is given on $\Si_*$ by
\beaa
\widetilde{\la}=1, \qquad \widetilde{f}=\frac{a}{r}f_0, \qquad  \widetilde{\fb}=-\frac{(\nu(r)-b_*)}{1-\frac{1}{4}b_*\frac{a^2(\sin\th)^2}{r^2}}\frac{a}{r}f_0.
\eeaa
Note in particular that $\widetilde{\la}=\la$ and $\widetilde{f}=f$, but $\widetilde{\fb}\neq \fb$ so that the outgoing PG frame and the outgoing PT frame of $\Mext$ are initialized differently on $\Si_*$. Indeed, the initialization used for the PG frame on $\Si_*$ would lead to a loss of one derivative and hence cannot be used for the PT frame. 
\end{remark}

%%%%%%%%%%%%%%%%%%%%%%%%%%%%%%%%%%%

\subsubsection{Initialization of the ingoing PT structure of $\Mint'$}

%%%%%%%%%%%%%%%%%%%%%%%%%%%%%%%%%%%

Recall that the ingoing PT structure of $\Mint'$ is initialized on the timelike hypersurface $\TT$. We denote by 
$\{ (e_3, e_4, \HH), r, \th, \Jk\}$ the outgoing PT initial data set induced by the outgoing PT structure of $\Mext$ on $\TT=\{r=r_0\}$. In view of Lemma \ref{lemma:constructioningoingPTframes:chap9}, it suffices to prescribe the  ingoing PT  initial data set, as in Definition \ref{def:initialingoingPTdataset:chap9}, corresponding to the PT structure of $\Mint'$ on $\TT$. We make the following choice of ingoing PT  initial data set $\{ (e_3', e_4', \HH'), r', \th', \Jk'\}$ on $\TT$:
\begin{enumerate}
\item The frame $(e_1', e_2', e_3', e_4')$, consisting of the null pair $(e_3', e_4')$ and the horizontal structure $\HH'$, is obtained from the null frame $(e_1, e_2, e_3, e_4)$ of the outgoing PT structure of $\Mext$ on $\TT$ by 
\bea
\lab{eq:initialization-frameMint'}
e_a'=e_a, \,\, a=1,2, \qquad e_4'=\frac{\Delta}{|q|^2}e_4, \qquad e_3'=\frac{|q|^2}{\Delta}e_3.
\eea

\item The functions $(r', \th')$ coincide on $\TT$ with the functions $(r, \th)$. 

\item The complex horizontal 1-form $\Jk'$ coincides on $\TT$ with the complex horizontal 1-form $\Jk$.
\end{enumerate}

\begin{definition}
\lab{def:initialization-frameMint'}
 We   define  $\Mint'\subset \MM\setminus \Mext $ to be  the region covered by the maximally extended\footnote{That is   the null ingoing geodesics generated by the PT structure end up either at $\AA$ or ${}^{(top)}\Si$.} ingoing PT structure initialized on $\TT=\{r=r_0\}$ by  the outgoing PT structure of $\Mext$ as above.
 \end{definition}

%%%%%%%%%%%%%%%%%%%%%%%%%%%%%%%%%%%

\subsubsection{Initialization of the ingoing PT structure of $\Mtop'$}

%%%%%%%%%%%%%%%%%%%%%%%%%%%%%%%%%%%

Recall that the ingoing PT structure of $\Mtop'$ is initialized on the hypersurface $\{u=u'_*\}$ of $\Mext$. We denote by $\{ (e_3, e_4, \HH), r, \th, \Jk\}$ the outgoing PT initial data set induced by the outgoing PT structure of $\Mext$ on $\{u=u'_*\}$. In view of Lemma \ref{lemma:constructioningoingPTframes:chap9}, it suffices to prescribe the  ingoing PT  initial data set, as in Definition \ref{def:initialingoingPTdataset:chap9}, corresponding to the PT structure of $\Mtop'$ on $\{u=u'_*\}$. We make the following choice of ingoing PT  initial data set $\{ (e_3'', e_4'', \HH''), r'', \th'', \Jk''\}$ on $\{u=u'_*\}$:
\begin{enumerate}
\item The frame $(e_1'', e_2'', e_3'', e_4'')$, consisting of the null pair $(e_3'', e_4'')$ and the horizontal structure $\HH''$, is obtained from the null frame $(e_1, e_2, e_3, e_4)$ of the outgoing PT structure of $\Mext$ on $\{u=u'_*\}$ by 
\bea
\lab{eq:initialization-frameMtop'}
e_a''=e_a, \,\, a=1,2, \qquad e_4''=\frac{\Delta}{|q|^2}e_4, \qquad e_3''=\frac{|q|^2}{\Delta}e_3.
\eea

\item The functions $(r'', \th'')$ coincide on $\{u=u'_*\}$ with the functions $(r, \th)$. 

\item The complex horizontal 1-form $\Jk''$ coincides on $\{u=u'_*\}$ with the complex horizontal 1-form $\Jk$.
\end{enumerate}

\begin{definition}
\lab{def:initialization-frameMtop'}
 We   define  $\Mtop'\subset \MM$ to be  the region covered by the maximally extended\footnote{That is   the null ingoing geodesics generated by the PT structure end up either at ${}^{(top)}\Si$ or $\AA$.} ingoing PT structure initialized on $\{u=u'_*\}$ by  the outgoing PT structure of $\Mext$ as above.
 \end{definition}

 \begin{remark}
 Since $u_*-2\leq u_*'\leq u_*-1$, note that $\Mtop'$ includes $\,{}^{(top)}\MM$, as well as a small portion of $\Mint$, $\Mint'$ and $\Mext$.
 \end{remark}

%%%%%%%%%%%%%%%%%%%%%%%%%%%%%%

\section{Outgoing PT structure of $\Mext$}

%%%%%%%%%%%%%%%%%%%%%%%%%%%%%%

%%%%%%%%%%%%%%%%%%%%%%%%%%%%%%%%%%%%%% 

\subsection{Null structure equations for the PT frame of $\Mext$} 

%%%%%%%%%%%%%%%%%%%%%%%%%%%%%%%%%%%%%%

We recall below Proposition \ref{Prop:NullStr-outgoingPTframe} for outgoing PT structures which applies in particular to the outgoing PT structure of $\Mext$.
\begin{proposition}
\lab{Prop:Mainequations-PTframe}
Consider an outgoing PT structure. Then,  the   equations in the $e_4$ direction for the Ricci coefficients  of the  outgoing PT frame take the form
\beaa
\nab_4\tr X +\frac{1}{2}(\tr X)^2 &=& -\frac{1}{2}\Xh\c\ov{\Xh},\\
\nab_4\Xh+\Re(\tr X)\Xh &=& -A,
\\
\nab_4\tr\Xb +\frac{1}{2}\tr X\tr\Xb  &=& -\DD\c\left(\frac{aq}{|q|^2}\ov{\Jk}\right) +  \frac{a^2}{|q|^2}|\Jk|^2 +2\ov{P} -\frac{1}{2}\Xh\c\ov{\Xbh},\\
\nab_4\widehat{\Xb} +\frac{1}{2}\tr X\, \widehat{\Xb}  &=& -\DD\hot\left(\frac{a\ov{q}}{|q|^2}\Jk\right)  +   \frac{a^2(\ov{q})^2}{|q|^4}\Jk\hot\Jk -\frac{1}{2}\ov{\tr\Xb} \widehat{X},
\\
\nab_4Z +\frac{1}{2}\tr X Z &=&  -\frac{1}{2}\tr X\frac{a\ov{q}}{|q|^2}\Jk-\frac{1}{2}\widehat{X}\c\left(\ov{Z}+\frac{aq}{|q|^2}\ov{\Jk}\right) -B,\\
\ \nab_4\Xib &=& -\nab_3\left(\frac{a\ov{q}}{|q|^2}\Jk\right) -\frac{1}{2}\ov{\tr\Xb}\left(\frac{a\ov{q}}{|q|^2}\Jk+H\right) -\frac{1}{2}\Xbh\c\left(\frac{aq}{|q|^2}\ov{\Jk}+\ov{H}\right) -\Bb,\\
\nab_4H &=&  -\frac{1}{2}\ov{\tr X}\left(H+\frac{a\ov{q}}{|q|^2}\Jk\right) -\frac{1}{2}\Xh\c\left(\ov{H} + \frac{aq}{|q|^2}\ov{\Jk}\right) -B,
\\
\nab_4\omb &=&   \left(\eta+\Re\left(\frac{a\ov{q}}{|q|^2}\Jk\right)\right)\c\ze +\eta\c\Re\left(\frac{a\ov{q}}{|q|^2}\Jk\right) +\rho.
\eeaa
\end{proposition}

%%%%%%%%%%%%%%%%%%%%%%%%%%%%%%%%%

\subsection{Other transport equations in the $e_4$ direction}

%%%%%%%%%%%%%%%%%%%%%%%%%%%%%%%%%

\begin{lemma}\lab{lemma:transportequationine4fornabthetaandrande3thetaandr}
We have
\beaa
\nab_4\DD\cos\th+\frac{1}{2}\tr X\DD\cos\th &=& -\frac{1}{2}\Xh\c\ov{\DD}\cos\th,\\
\nab_4\DD r+\frac{1}{2}\tr X\DD r &=& -\frac{1}{2}\Xh\c\ov{\DD}r+(Z+\Hb),\\
\nab_4\DD u+\frac{1}{2}\tr X\DD u &=& -\frac{1}{2}\Xh\c\ov{\DD}u,
\eeaa
and
\beaa
e_4(e_3(\cos\th)) &=&  -\Re\Big((H -\Hb)\c\ov{\DD}(\cos\th)\Big),\\
e_4(e_3(r)) &=& -2\omb   -\Re\Big((H -\Hb)\c\ov{\DD}r\Big),\\
e_4(e_3(u)) &=&   -\Re\Big((H -\Hb)\c\ov{\DD}u\Big).
\eeaa
\end{lemma}

\begin{proof}
Straightforward verification.
\end{proof}

%%%%%%%%%%%%%%%%%%%%%%%%%%%%%%%%%

\subsection{Linearized  quantities  for  the  outgoing PT frame}

%%%%%%%%%%%%%%%%%%%%%%%%%%%%%%%%%

Recall the definition of the  linearized quantities in the PT frame,  see Definition \ref{def:renormalizationofallnonsmallquantitiesinPTstructurebyKerrvalue}.

\begin{definition}\lab{def:renormalizationofallnonsmallquantitiesinPTstructurebyKerrvalue:chap9} 
We consider  the following renormalizations, for given constants $(a, m)$, 
\bea
\begin{split}
\trXc &:= \tr X-\frac{2}{q}, \qquad\qquad\,\,\,  \trXbc := \tr\Xb+\frac{2q\Delta}{|q|^4},\\
\Zc &:= Z -\frac{a\ov{q}}{|q|^2}\Jk,\qquad\qquad\,\,\,\, \Hc := H-\frac{aq}{|q|^2}\Jk,\\
\ombc &:= \omb-\frac 1 2 \pr_r\left(\frac{\De}{|q|^2} \right), \qquad \Pc := P+\frac{2m}{q^3},
\end{split}
\eea
as well as 
\bea
\begin{split}
\widecheck{e_3(r)} &:= e_3(r)+\frac{\Delta}{|q|^2},\qquad \widecheck{\DD(\cos\th)} := \DD(\cos(\th)) -i\Jk,\\
\widecheck{\DD u} &:=\DD u- a\Jk, \qquad\qquad\,\,\,\,  \widecheck{e_3(u)}:=e_3(u)-\frac{2(r^2+a^2)}{|q|^2},
\end{split}
\eea
and 
\bea
\widecheck{\ov{\DD}\c\Jk} := \ov{\DD}\c\Jk-\frac{4i(r^2+a^2)\cos\th}{|q|^4}, \qquad \widecheck{\nab_3\Jk}:=\nab_3\Jk -\frac{\De q}{|q|^4}\Jk.
\eea
\end{definition}

%%%%%%%%%%%%%%%%%%%%%%%%%%%%%%%%%%%%%%%%%%%%%

\subsection{Definition of the notations $\Ga_b$ and $\Ga_g$ for error terms}

%%%%%%%%%%%%%%%%%%%%%%%%%%%%%%%%%%%%%%%%%%%%%

\begin{definition}
\lab{definition.Ga_gGa_b:outgoingPTcase:chap9}
The set of all linearized quantities is of the form $\Ga_g\cup \Ga_b$ with  $\Ga_g,  \Ga_b$
 defined as follows.
 
\begin{enumerate}
\item 
 The set   $\Ga_g$   with
 \bea
 \bsplit
 \Ga_g &= \Bigg\{\trXc,\quad  \Xh,\quad \Zc,\quad \trXbc, \quad  r^{-1} \nab (r), \quad r\Pc, \quad  rB, \quad  rA\Bigg\}.
 \end{split}
 \eea
 
 \item  The set  $\Ga_b=\Ga_{b,1}\cup \Ga_{b, 2}\cup \Ga_{b, 3}$   with
 \bea
 \bsplit
 \Ga_{b,1}&= \Bigg\{\Hc, \quad \Xbh, \quad \ombc, \quad \Xib,\quad  r\Bb, \quad \Ab\Bigg\},\\
  \Ga_{b, 2}&= \Bigg\{r^{-1}\widecheck{e_3(r)}, \quad  \widecheck{\DD(\cos\th)}, \quad e_3(\cos\th), \quad \widecheck{\DD u}, \quad r^{-1}\widecheck{e_3(u)}\Bigg\},\\
  \Ga_{b,3}&=\Bigg\{ r\,\widecheck{\ov{\DD}\c\Jk}, \quad r\,\DD\hot\Jk, \quad r\,\widecheck{\nab_3\Jk}\Bigg\}.
 \end{split}
 \eea
 \end{enumerate}
\end{definition}

%%%%%%%%%%%%%%%%%%%%%%%%%%%%%%%%%

\subsection{Linearized  equations for outgoing PT structures}
\lab{sec:linearizedeqautionsforoutoingPTstructures:chap9}

%%%%%%%%%%%%%%%%%%%%%%%%%%%%%%%%%

Recall the convention  $O(r^{-p}) $ introduced earlier  in Definition \ref{def:ordermagnitude}.

\begin{definition}[Order of magnitude notation]
 \lab{def:ordermagnitude:chap9}
 Throughout this chapter,  we will be  using the  notation $O(r^{-p})$ to  denote:
\begin{enumerate} 
\item a scalar function depending only on $(r, \th)$ which is smooth and such that
\beaa
r^p|(r\pr_r, \pr_\th)^kO(r^{-p})|\les_k 1\quad\textrm{for }\quad k\geq 0\quad\textrm{and}\quad r\geq r_0,
\eeaa

\item a 1-form of the type $O(r^{-p+1})\Jk$ where $O(r^{-p+1})$ denotes a scalar function as above,

\item a symmetric traceless 2-tensor of the type $O(r^{-p+2})\Jk\hot\Jk$ where $O(r^{-p+2})$ denotes a scalar function as above.
\end{enumerate}
\end{definition}

The following proposition provides the linearized null structure equations for the outgoing PT frame of $\Mext$. 
\begin{proposition}\lab{Prop:linearizedPTstructure1}
In an outgoing PT frame, the  linearized null structure equations in the $e_4$ direction are 
\beaa
\nab_4\trXc +\frac{2}{q}\trXc &=& \Ga_g\c\Ga_g,\\
\nab_4\Xh+\Re\left(\frac{2}{q}\right)\Xh &=& -A+\Ga_g\c\Ga_g,\\
\nab_4\Zc + \frac{1}{q}\Zc &=&   - \frac{a\ov{q}}{|q|^2}\trXc \Jk    -\frac{aq}{|q|^2}\ov{\Jk}\c\widehat{X}  -B+\Ga_g\c\Ga_g,\\
\nab_4\Hc+ \frac{1}{\ov{q}}\Hc &=&  -\frac{ar}{|q|^2}\ov{\trXc}\Jk -\frac{ar}{|q|^2}\ov{\Jk}\c\Xh -B+\Ga_b\c\Ga_g,\\
 \nab_4\ombc  &=&   \Re\left(\Pc\right) +\frac{ar}{|q|^2}\Re\left(\Jk\c \ov{\Zc}\right)+\frac{2a}{|q|^2}\Re\left(q\ov{\Jk}\c \Hc\right)+\Ga_g\c\Ga_g,
\\
\nab_4\trXbc +\frac{1}{q}\trXbc  &=& 2\ov{\Pc} + \frac{q\Delta}{|q|^4}\trXc -\frac{aq}{|q|^2}\widecheck{\DD\c\ov{\Jk}} +\frac{a}{\ov{q}^2}\DD(r)\c\ov{\Jk} - \frac{ia^2}{\ov{q}^2}\widecheck{\DD(\cos\th)}\c\ov{\Jk} +\Ga_b\c\Ga_g,\\
\nab_4\Xbh +\frac{1}{q}\widehat{\Xb}  &=&  -\frac{a\ov{q}}{|q|^2}\DD\hot\Jk +\frac{a}{q^2}\DD(r)\hot\Jk+\frac{ia^2}{q^2}\widecheck{\DD(\cos\th)}\hot\Jk +\frac{q\Delta}{|q|^4} \widehat{X}+\Ga_b\c\Ga_g,
\\
 \nab_4\Xib &=&  \frac{\ov{q}\Delta}{|q|^4}\Hc - \frac{ar}{|q|^2}\ov{\trXbc}\Jk  -\frac{ar}{|q|^2}\ov{\Jk}\c\Xbh -\Bb -\frac{a\ov{q}}{|q|^2}\widecheck{\nab_3\Jk}\\
 &&+\frac{a}{q^2}\left(\widecheck{e_3(r)}+iae_3(\cos\th)\right)\Jk +\Ga_b\c\Ga_b.
\eeaa
\end{proposition}

\begin{proof}
See appendix \ref{appendix:ProofProp{Prop:linearizedPTstructure1}}. Compare also with the proof of  Lemma  \ref{Lemma:linearized-nullstr} in the outgoing PG frame.
\end{proof}

\begin{remark}
Note that the equations for $\nab_4\trXbc$, $\nab_4\Xbh$, and $\nab_4\Xib$ in Proposition \ref{Prop:linearizedPTstructure1} do  not loose derivative, unlike the corresponding ones   in  the outgoing PG frame  (compare with  Lemma \ref{Lemma:linearized-nullstr}).
\end{remark}

%%%%%%%%%%%%%%%%%%%%%%%%%%%%%%%%%

\subsection{Other linearized equations}
\lab{sec:otherlinearizedequationsoutgoingPTframeMext:chap9}

%%%%%%%%%%%%%%%%%%%%%%%%%%%%%%%%%

The following lemma follows immediately from Lemma \ref{lemma:transportequationine4fornabthetaandrande3thetaandr}, the definition of the linearized quantities, and the definition of $\Ga_g$ and $\Ga_b$. 
\begin{lemma}\lab{lemma:transportequationine4fornabthetaandrande3thetaandr:linearized}
We have
\beaa
\nab_4\widecheck{\DD\cos\th}+\frac{1}{q}\widecheck{\DD\cos\th} &=& O(r^{-1})\trXc+\frac{i}{2}\ov{\Jk}\c\Xh+\Ga_b\c\Ga_g,\\
\nab_4\DD r+\frac{1}{q}\DD r &=& \Zc+r\Ga_g\c\Ga_g,\\
\nab_4\widecheck{\DD u}+\frac{1}{q}\widecheck{\DD u} &=& O(r^{-1})\trXc -\frac{a}{2}\ov{\Jk}\c\Xh+\Ga_b\c\Ga_g,
\eeaa
and
\beaa
e_4(e_3(\cos\th)) &=&  -\Im\Big(\ov{\Jk}\c\Hc\Big)-\Re\left(\frac{2ar}{|q|^2}\Jk\c\ov{\widecheck{\DD(\cos\th)}}\right)+\Ga_b\c\Ga_b,\\
e_4(\widecheck{e_3(r)}) &=& -2\ombc  -\Re\left(\frac{2ar}{|q|^2}\Jk\c\ov{\DD}r\right)+r\Ga_b\c\Ga_g,\\
e_4(\widecheck{e_3(u)}) &=&   -\Re\Big(a\ov{\Jk}\c\Hc\Big) -\Re\left(\frac{2ar}{|q|^2}\Jk\c\ov{\widecheck{\DD u}}\right)+\Ga_b\c\Ga_b.
\eeaa
\end{lemma}
\begin{proof}
Straightforward verification.  Compare it  also with the proof of Lemma \ref{Lemma:otherlinearizedquant}.
\end{proof}

\begin{lemma}\lab{lemma:transportequationine4fornaJkandnab3Jk:linearized}
We have
\beaa
\nab_4 \DD\hot\Jk+\frac{2}{q}\DD\hot\Jk &=& O(r^{-1})B +O(r^{-2})\trXc+O(r^{-2})\Xh\\
&&+O(r^{-2})\Zc +O(r^{-3})\widecheck{\DD(\cos\th)},\\
\nab_4\widecheck{\ov{\DD}\c\Jk} +\Re\left(\frac{2}{q}\right)\widecheck{\ov{\DD}\c\Jk} &=& O(r^{-1})B+O(r^{-2})\trXc+O(r^{-2})\Xh +O(r^{-2})\Zc\\
&&+O(r^{-3})\widecheck{\DD(\cos\th)},\\
\nab_4\widecheck{\nab_3\Jk}+\frac{1}{q}\widecheck{\nab_3\Jk} &=& O(r^{-1})\Pc+ O(r^{-3})\widecheck{e_3(r)}+O(r^{-3})e_3(\cos\th)+O(r^{-2})\ombc\\
&&+O(r^{-2})\Hc+O(r^{-2})\widecheck{\nab\Jk}.
\eeaa
\end{lemma} 

\begin{proof}
Straightforward verification\footnote{ Note that the proof must be adapted from the case of a PG frame, i.e. $\Hb=-Z$, to the case of a PT frame, i.e. $\Hb=-\frac{a\ov{q}}{|q|^2}\Jk$. This introduces a slight change in the commutators $[\nab_4,\nab]$ and $[\nab_4, \nab_3]$.}.  Compare it  also with the proof of Lemma  \ref{Lemma:linearizedJk}. 
\end{proof}

\begin{remark}[Triangular structure of the main equations]
We can order the linearized quantities appearing in the equations of Proposition \ref{Prop:linearizedPTstructure1}, 
Lemma \ref{lemma:transportequationine4fornabthetaandrande3thetaandr:linearized} and Lemma \ref{lemma:transportequationine4fornaJkandnab3Jk:linearized} as follows
\beaa
\trXc, \, \Xh, \, \Zc,\, \Hc,\,  \widecheck{\DD\cos\th}, \,\ombc,\,  \, \DD r, \, \widecheck{\DD u},\,  e_3(\cos\th),\,\widecheck{e_3(r)},\, \widecheck{e_3(u)},\, \DD\hot\Jk,\, \widecheck{\ov{\DD}\c \Jk},\, \widecheck{\nab_3\Jk}, \trXbc, \, \Xbh,\, \Xib, 
\eeaa
and note that the  transport equation for each one of them  depends\footnote{At the linear level and excluding  curvature terms.}   only  on the previous  components of the sequence.   This triangular structure is essential  for estimating the terms one by one. This crucial   fact   will be used   in section \ref{sec:theoremM8recoverRicciawaytrapping} to estimate the Ricci coefficients of the PT frame of $\Mext$. 
\end{remark}

%%%%%%%%%%%%%%%%%%%%%%%%%%%%%%%%%%%%%%%%%%%%%

\subsection{Comparison between the PT and PG structures of $\Mext$}

%%%%%%%%%%%%%%%%%%%%%%%%%%%%%%%%%%%%%%%%%%%%%

The following lemma compares the outgoing PG and PT structures of $\Mext$. 
\begin{lemma}
\lab{lemma:linkPGandTframeinMext}
Let $\{  (e_3, e_4, \HH), r, \th, \Jk\} $ denote the extended  outgoing PG structure of $\Mext$, and let $\{ (e_3', e_4', \FF'), r' , \th', \Jk'\} $ denote the extended  outgoing PT structure of $\Mext$. Also, let $(f, \fb, \la)$ denote the transition coefficients from the PG frame to the PT one. Then, the following identities hold in $\Mext$:
 \begin{enumerate}
 \item We have  $f=0$  and $\la=1$. In particular, we have 
\beaa
e_4' = e_4.  
\eeaa

\item  We have
\beaa
u'=u, \qquad r'=r, \qquad \th'=\th, \qquad \Jk'=\Jk.
\eeaa

\item  Let $(f', \fb', \la')$ denote the transition coefficients from the PT frame to the PG one. Then, we have 
\beaa
\la'=1, \quad f'=0, \quad \fb'=-\fb.
\eeaa

\item We have, in view of  Definition \ref{def:renormalizationofallnonsmallquantitiesinPGstructurebyKerrvalue} for  linearized outgoing PG  quantities,
\beaa
\nab_4\fb &=& 2\widecheck{\ze}.
\eeaa

\item With the notation $\Fb'=\fb'+i\dual\fb'$, we have, in view of  Definition \ref{def:renormalizationofallnonsmallquantitiesinPTstructurebyKerrvalue:chap9} for linearized outgoing PT  quantities,
\beaa
\nab_4'\Fb'+\frac{1}{2}\tr X'\Fb' &=& -2\Zc' -\Fb'\c\chih'.
\eeaa
 \end{enumerate}
\end{lemma}

\begin{proof}
See section \ref{sec:proofoflemma:linkPGandTframeinMext}.
\end{proof}

%%%%%%%%%%%%%%%%%%%%%%%%%%%%%%%%%%%%%%%%%%%%%%%%%%

\subsection{The choice of the constant $u_*'$}
\lab{sec:choiceofuprimstar}

%%%%%%%%%%%%%%%%%%%%%%%%%%%%%%%%%%%%%%%%%%%%%%%%%%

Recall from section \ref{sec:defintionofthePTstructuresinMM} that the ingoing PT structure of $\Mtop'$ in initialized from the outgoing PT structure of $\Mext$ on the hypersurface $\{u=u_*'\}$ of $\Mext$ for $u_*'\in [u_*-2, u_*-1]$. We are now ready to make a specific choice of $u_*'$.  First, we introduce the following notation 
\bea\lab{eq:defintionofthenormsRcGacweighted}
 \bsplit
 \big|\Rc\big|_{w, k}^2 &:= r^{3+\dt}|\dk_*^{\le k}(A, B)|^2 +r^{3-\dt}|\dk_*^{\le k}\Pc |^2 +r^{1-\dt}|\dk_*^{\le k}\Bb|^2,\\
\big|\Gac\big|_{w, k}^2 &:= r^2|\dkb^{\leq k}\Ga_g|^2+|\dkb^{\leq k}\Ga_b|^2,
\end{split}
\eea
where $(A, B, \Pc, \Bb)$ denote  linearized curvature components w.r.t. the the  outgoing PT frame of $\Mext$,  
$\Ga_g$ and $\Ga_b$ are  defined w.r.t. the  outgoing PT frame of $\Mext$ as in Definition \ref{definition.Ga_gGa_b:outgoingPTcase:chap9}, $\dk_*$ denote weighted derivatives tangential to the hypersurface $\{u=u_*'\}$, and $\dkb$ denote weighted derivatives tangential to the sphere $\Si_*\cap\{u=u_*'\}$. 

With the notations in \eqref{eq:defintionofthenormsRcGacweighted}, we choose $u_*'$ such that we have 
\bea\lab{eq:choicofuprimestarbyLebesguepointarguement}
\nn &&\int_{\{u=u_*'\}}    \big|\Rc\big|_{w, k_{large}+7}^2 +\int_{\Si_*\cap\{u=u_*'\}} \big|\Gac\big|_{w, k_{large}+7}^2\\
 &=& \inf_{u*-2\leq u_1\leq u*-1}\left(\int_{\{u=u_1\}}    \big|\Rc\big|_{w, k_{large}+7}^2 +\int_{\Si_*\cap\{u=u_1\}} \big|\Gac\big|_{w, k_{large}+7}^2\right).
\eea

%%%%%%%%%%%%%%%%%%%%%%%%%%%%%%%%%%

\section{Ingoing  PT structures  of $\Mint'$ and $\Mtop'$} 

%%%%%%%%%%%%%%%%%%%%%%%%%%%%%%%%%%%
 
%%%%%%%%%%%%%%%%%%%%%%%%%%%%%%%%

\subsection{Linearized quantities in   an ingoing PT frame}

%%%%%%%%%%%%%%%%%%%%%%%%%%%%%%%%%

Given an extended  ingoing PT  structure   $\{ (e_3, e_4, \HH), r, \th,  \ub, \Jk\}$   hold true:
\begin{enumerate}
\item We have
 \beaa
 \xib=\omb=0, \quad e_3(r)=-1, \quad e_3(\ub)= e_3(\th)=0, \qquad \nab_3(\ov{q}\Jk)=0.
 \eeaa
 In addition, we have
 \beaa
 H=\frac{a q}{|q|^2}\Jk.
 \eeaa
 \item The  quantities
 \beaa
\Xh, \quad \Xbh, \quad \Xi, \quad A, \quad B, \quad \Bb, \quad \Ab, \quad \DD r, \quad e_4(\cos\th), \quad \DD\hot\Jk,
\eeaa
  vanish  in      Kerr  and therefore are   small in perturbations.
\end{enumerate}
We renormalize below all other quantities, not vanishing in Kerr\footnote{Since  $H=\frac{aq}{|q|^2}\Jk$,  $H$ does not need to be included in Definition \ref{def:linearizedPT-ingoingcase:chap9}.}, by subtracting their  $\mbox{Kerr}(a, m) $  values.

\begin{definition}\lab{def:linearizedPT-ingoingcase:chap9}
We consider the following renormalizations, for given constants $(a, m)$,
\bea
\bsplit
\trXc &:= \tr X-\frac{2\ov{q}\De}{|q|^4}, \qquad\,\qquad     \trXbc := \tr\Xb+\frac{2}{\ov{q}},\\
\Zc &:= Z-\frac{aq}{|q|^2}\Jk,\qquad \qquad \quad\,\,\,
\Hbc := \Hb+\frac{a\ov{q}}{|q|^2}\Jk,\\
\omc & := \om  + \frac{1}{2}\pr_r\left(\frac{\De}{|q|^2} \right),\qquad \quad\, \Pc := P+\frac{2m}{q^3},
\end{split}
\eea
as well as 
\bea
\bsplit
\widecheck{e_4(r)} &:= e_4(r)-\frac{\Delta}{|q|^2},\qquad  \widecheck{\DD(\cos\th)} := \DD(\cos(\th)) -i\Jk,\\
\widecheck{\DD\ub} &:= \DD\ub -a\Jk,\qquad\qquad\,\,\,\, \widecheck{e_4(\ub)} := e_4(\ub) -\frac{2(r^2+a^2)}{|q|^2},
\end{split}
\eea
and
\bea
\bsplit
\widecheck{\ov{\DD}\c\Jk}& := \ov{\DD}\c\Jk-\frac{4i(r^2+a^2)\cos\th}{|q|^4}, \qquad \widecheck{\nab_4\Jk}:=\nab_4\Jk +\frac{\De \ov{q}}{|q|^4}\Jk.
\end{split}
\eea
\end{definition}

%%%%%%%%%%%%%%%%%%%%%%%%%%%%%%%%%%%%%%%%%%%%%

\subsection{Definition of the notations $\Ga_b$ and $\Ga_g$ for error terms}

%%%%%%%%%%%%%%%%%%%%%%%%%%%%%%%%%%%%%%%%%%%%%

\begin{definition}
\lab{definition.Ga_gGa_b:outgoingPTcase}
The set of all linearized quantities is of the form $\Ga_g\cup \Ga_b$ with  $\Ga_g,  \Ga_b$
 defined as follows.
 
\begin{enumerate}
\item 
 The set   $\Ga_g=\Ga_{g,1}\cup \Ga_{g,2}$   with
 \bea
 \bsplit
 \Ga_{g,1}&= \Bigg\{\Xi, \quad \omc, \quad \trXc,\quad  \Xh,\quad \Zc,\quad \Hbc, \quad \trXbc, \quad r\Pc, \quad  rB, \quad  rA\Bigg\},\\
 \Ga_{g,2}&= \Bigg\{r^{-1} \nab (r), \quad r\widecheck{e_4(r)}, \quad r\widecheck{e_4(\ub)}, \quad re_4(\cos\th), \quad r^2\widecheck{\nab_4\Jk}\Bigg\}.
 \end{split}
 \eea
 
 \item  The set  $\Ga_b$   with
 \bea
 \bsplit
 \Ga_b &= \Bigg\{\Xbh,\quad  r\Bb, \quad \Ab, \quad \widecheck{\DD(\cos\th)}, \quad \widecheck{\DD\ub}, \quad r\,\widecheck{\ov{\DD}\c\Jk}, \quad r\,\DD\hot\Jk\Bigg\}.
 \end{split}
 \eea
 \end{enumerate}
\end{definition}

%%%%%%%%%%%%%%%%%%%%%%%%%%%%%%%%%
  
\subsection{Main linearized equations for ingoing PT structures}
 
%%%%%%%%%%%%%%%%%%%%%%%%%%%%%%%%%

\begin{proposition}
\lab{proposition:Maineqts-Mint'} 
In an ingoing PT frame, the  linearized null structure equations in the $e_3$ direction are
\beaa
\bsplit
\nab_3(\trXbc) -\frac{2}{\ov{q}}\trXbc  &=          \Ga_b\c\Ga_b, \\
\nab_3\Xbh -\frac{2r}{|q|^2}\Xbh &=  -\Ab +\Ga_b\c\Ga_b,\\
\nab_3\Zc - \frac{1}{\ov{q}}\Zc  &= -\Bb +  O(r^{-2}) \Xbh  +O(r^{-2})\trXbc+\Ga_b\c\Ga_g,\\
\nab_3\Hbc-\frac{1}{q}\Hbc &= \Bb     +O(r^{-2})\Xbh  +O(r^{-2}) \trXbc +\Ga_b\c\Ga_g,\\
\end{split}
\eeaa
\beaa
\bsplit
\nab_3\trXc -\frac{1}{\ov{q}}\trXc &= 2\Pc+O(r^{-1}) \trXbc+O(r^{-1})\widecheck{\DD\c\ov{\Jk}}+ O(r^{-3}) \DD r+O(r^{-3})\widecheck{\DD(\cos\th)}+\Ga_b\c \Ga_g,\\
\nab_3 \Xh-\frac{1}{\ov{q}}   \Xh&=O(r^{-1}) \DD\hot\Jk+O(r^{-4})\widecheck{ \DD(\cos \th)}+ O(r^{-1}) \Xbh+\Ga_b\c \Ga_g,\\
\nab_3\omc&=\Re(\Pc)+O(r^{-2})\Zc+O(r^{-2})\Hbc+\Ga_b\c\Ga_g,\\
 \nab_3\Xi &= O(r^{-1}) \Hbc+ O(r^{-2})\trXc+O(r^{-2}) \Xh +B+ O(r^{-1})\widecheck{\nab_4 \Jk}+O(r^{-3}) \widecheck{e_4(r)}\\
 &+O(r^{-3}) e_4(\cos\th) +\Ga_b\c \Ga_g.
\end{split}
\eeaa

Also, we have
\beaa
\nab_3\widecheck{\DD\cos\th}-\frac{1}{\ov{q}}\widecheck{\DD\cos\th} &=& O(r^{-1})\trXbc+O(r^{-1}) \Xbh+\Ga_b\c\Ga_b,\\
\nab_3\DD r-\frac{1}{\ov{q}}\DD r &=&  \Zc+r\Ga_b\c\Ga_g,\\
\nab_3\widecheck{\DD\ub}- \frac{1}{\ov{q}}\widecheck{\DD\ub} &=& O(r^{-1})\trXbc +O(r^{-1}) \Xbh+\Ga_b\c\Ga_b,
\eeaa
\beaa
e_3(e_4(\cos\th)) &=& O(r^{-1})\Hbc +O(r^{-2}) \widecheck{\DD(\cos\th)} +\Ga_g\c\Ga_g,\\
e_3(\widecheck{e_4(r)}) &=& -2\omc  +O(r^{-2}) \DD r +r\Ga_g\c\Ga_g,\\
e_3(\widecheck{e_4(\ub)}) &=&  O(r^{-1}) \Hbc  +O(r^{-2})  \widecheck{\DD u} +\Ga_g\c\Ga_g,
\eeaa
\beaa
\nab_3 \DD\hot\Jk- \frac{2}{\ov{q}}\DD\hot\Jk &=& O(r^{-1})\Bb +O(r^{-2})\trXbc+O(r^{-2}) \Xbh+ O(r^{-2})\Zc\\
&& +O(r^{-3})\widecheck{\DD(\cos\th)},\\
\nab_3\widecheck{\ov{\DD}\c\Jk} -\Re\left(\frac{2}{\ov{q}}\right)\widecheck{\ov{\DD}\c\Jk} &=& O(r^{-1}) \Bb  +O(r^{-2})\trXbc+O(r^{-2})\Xbh +O(r^{-2})\Zc\\
&&+O(r^{-3})\widecheck{\DD(\cos\th)},\\
\nab_3\widecheck{\nab_4\Jk}-\frac{1}{\ov{q}}\widecheck{\nab_4\Jk} &=& O(r^{-1}) \Pc+ O(r^{-3})\widecheck{e_4(r)}+O(r^{-3})e_4(\cos\th)+O(r^{-2})\omc\\
&&+O(r^{-2})\Hbc +O(r^{-2})\widecheck{\nab\Jk}.
\eeaa
\end{proposition}

\begin{proof}
Straightforward verification. Compare with the proof of Proposition \ref{Prop:linearizedPTstructure1}, and Lemmas \ref{lemma:transportequationine4fornabthetaandrande3thetaandr:linearized} and 
\ref{lemma:transportequationine4fornaJkandnab3Jk:linearized} for the corresponding equations for outgoing PT structures. Compare also with Lemma \ref{Lemma:linearized-nullstr:ingoingcase}, Lemma \ref{Lemma:otherlinearizedquant:ingoingcase} and Lemma \ref{Lemma:linearizedJk:ingoingcase} for the corresponding equations for ingoing PG structures. 
\end{proof}

In order to deal with the trapping in the region $\Mint'$, we will need additional equations provided in the proposition below.
\begin{proposition}
\lab{proposition:Maineqts-Mint':bis} 
In the ingoing PT structure of $\Mint'$, the linearized Codazzi of $\Xbh$ takes  the following schematic form\footnote{Note that $r\leq r_0$ in $\Mint'$ so that powers of $r$ do not matter. As a consequence, we may simply denote  linear terms by $\Ga_b$ and  nonlinear terms by $\Ga_b\c\Ga_b$.} 
\beaa
\ov{\DD}\c\Xbh &=& \DD\ov{\trXbc} +\Ga_b+\Ga_b\c\Ga_b.
\eeaa
Also, we can write the Bianchi identities in $\Mint'$ schematically  in the  form
 \beaa
 \bsplit
 \nab_3A -\frac{1}{2}\DD\hot B &=   \Ga_b+\Ga_b\c\Ga_b,\qquad\qquad \nab_4B -\frac{1}{2}\ov{\DD}\c A =  \Ga_b+\Ga_b\c\Ga_b,\\
\nab_3B-\DD\ov{\Pc} &=    \Ga_b+\Ga_b\c\Ga_b,\qquad\qquad \nab_4P -\frac{1}{2}\DD\c \ov{B} =  \Ga_b+\Ga_b\c\Ga_b, \\
\nab_3P +\frac{1}{2}\ov{\DD}\c\Bb &=  \Ga_b+\Ga_b\c\Ga_b, \qquad\qquad\,\,\,\quad \nab_4\Bb+\DD \Pc =  \Ga_b+\Ga_b\c\Ga_b,\\
\nab_3\Bb +\frac{1}{2}\ov{\DD}\c\Ab &=  \Ga_b+\Ga_b\c\Ga_b,\qquad\qquad \nab_4\Ab +\frac{1}{2}\DD\hot\Bb =   \Ga_b+\Ga_b\c\Ga_b.
\end{split}
\eeaa
\end{proposition}

%%%%%%%%%%%%%%%%%%%%%%%%%%%%%%%%%
  
\subsection{The scalar function $\tau$ on $\MM$}
\lab{sec:defintionofthescalarfunctiontau:chap9}
 
%%%%%%%%%%%%%%%%%%%%%%%%%%%%%%%%%
 
We introduce on $\MM$ a scalar function $\tau$ with the following properties. 
\begin{proposition}\lab{prop:propertiesoftauusefulfortheoremM8:chap9}
There exists a scalar function $\tau$ defined on $\MM$ such that:
\begin{enumerate}
\item We have on  $\MM$
\bea
\g(\D\tau, \D\tau) &\leq& -\frac{m^2}{8r^2},
\eea
so that the level sets of $\tau$ are spacelike and asymptotically null.

\item The future boundary ${}^{(top)}\Si$ of $\MM$ is given by 
\bea
{}^{(top)}\Si=\{\tau=u_*\}
\eea
and $\tau\leq u_*$ on $\MM$.

\item Denoting, on each level set of $\ub$ in $\,{}^{(top)}\MM'(r\geq r_0)$, by $r_+(\ub)$ the maximal value of $r$  and by $r_-(\ub)$ the minimal value of $r$, we have\footnote{Note that \eqref{eq:upperboundrpubminusrmubonMtop:chap9} depends on the choice of ${}^{(top)}\Si$ and hence on the choice of $\tau$.}
\bea\lab{eq:upperboundrpubminusrmubonMtop:chap9}
0\leq r_+(\ub)-r_-(\ub)\leq 2(2m+1).
\eea

\item In $\Mtop'(r\leq r_0)$, $\tau$ satisfies   
\bea
u_*-2(m+2)\leq \tau \leq u_*.
\eea

\item In $\MM(r\leq r_0)$, $\tau$ satisfies 
\bea
e_4(\tau) = \frac{2(r^2+a^2) -\frac{m^2}{r^2}\De}{|q|^2}+r^{-1}\Ga_g, \,\,\,  e_3(\tau) = \frac{m^2}{r^2},\,\,\, \nab(\tau) = a\Re(\Jk)+\Ga_b.
\eea
\end{enumerate}
\end{proposition}

\begin{proof}
See section \ref{sec:proofofprop:propertiesoftauusefulfortheoremM8:chap9}.
\end{proof}

\begin{remark}
In view of the third property and the fourth property of Proposition \ref{prop:propertiesoftauusefulfortheoremM8:chap9}, $\,{}^{(top)}\MM'$  is in fact a local existence type region.
\end{remark}

\begin{remark}
The fact that $r_+(\ub)-r_-(\ub)\leq 2(2m+1)$ in $\Mtop'(r\geq r_0)$ in view of Proposition \ref{prop:propertiesoftauusefulfortheoremM8:chap9} is crucial to recover the Ricci coefficients in $\Mtop'(r\geq r_0)$. In particular, it is  crucial to control $\Xbh$ from $\Ab$ as we have schematically 
\beaa
|\Xbh| &\les& |\Xbh_{|_{u=u_*'}}|+(r_+(\ub)-r_-(\ub))|\Ab|
\eeaa
so that we indeed need $r_+(\ub)-r_-(\ub)\les 1$.
\end{remark}

\begin{definition}
 Given a level hypersurface $\Si(\tau) $ of $\tau$, we denote\footnote{Note from Proposition \ref{prop:propertiesoftauusefulfortheoremM8:chap9} that 
 \beaa
 \g(N_\tau, N_\tau) = \g(\D\tau, \D\tau)\leq -\frac{m^2}{8r^2}<0.
 \eeaa}
 \bea
 N_\tau:=-\g^{\a\b}\pr_\b\tau \pr_\a,\qquad \Nhat_\tau:=    \frac{1}{\sqrt{|\g(N_\tau, N_\tau)|}}     N_\tau,
 \eea
so that $\Nhat_\tau$ is the future unit normal  to  $\Si(\tau)$.
\end{definition}

 We will consider  a  subregion  $\Mint'_*$ of $\Mint'$   defined as follows.
\begin{definition}\lab{def:u=introductionofregionMintprimestar:chap9}
Let $\tau_*$ the supremum\footnote{Note that $u_*-4(m+1)\leq \tau_*<u_*$ in view of Proposition \ref{prop:propertiesoftauusefulfortheoremM8:chap9}.} of the values of $\tau$ such that $\Si(\tau)\cap\Mint'\neq\emptyset$ and $\Si(\tau)$ does not intersect $\{\ub=u_*\}$. We denote by $\Mint_*'$ the subset of $\Mint'$ with $\tau\le \tau_*$, i.e. 
\bea
\Mint_*' &:=& \Mint'\cap\{\tau\leq \tau_*\}.
\eea
\end{definition}

\begin{remark}
Note that the region $\Mint'$ has as future boundary  the hypersurface $\{\ub=u_*\}$ which is not spacelike, while the region $\Mint'_*\subset\Mint'$ has as future boundary the hypersurface $\{\tau=\tau_*\}$ which is spacelike in view of Proposition \ref{prop:propertiesoftauusefulfortheoremM8:chap9}.
\end{remark}

%%%%%%%%%%%%%%%%%%%%%%%%%%%%%%%%%%%%%%%%%%%%%%%%%%

\section{Control of top order derivatives in the PT frame}

%%%%%%%%%%%%%%%%%%%%%%%%%%%%%%%%%%%%%%%%%%%%%%%%%%

In this section, we state our main result concerning the control of top order derivatives in the PT frame, see Theorem \ref{theorem:Main-PT}. As a consequence, this yields the control of the PG frame and concludes the proof of Theorem M8, see section \ref{sec:endoftheproofofTheoremM8:chap9}.

%%%%%%%%%%%%%%%%%%%%%%%%%%%%%%%%%%%%%%%%%%%%%%%%%%

\subsection{Main norms}
\lab{sec:mainnormsPTframe:chap9}

%%%%%%%%%%%%%%%%%%%%%%%%%%%%%%%%%%%%%%%%%%%%%%%%%%

We introduce here the main norms  appearing in the statement of our  main  PT-Theorem in section \ref{sec:statementofthemainPTTheorem}.

%%%%%%%%%%%%%%%%%%%

\subsubsection{Norms on  $\Si_*$}

%%%%%%%%%%%%%%%%%%%

\begin{definition}
\lab{definition:PT-normson-GacSi_*} 
We define the   following PT-Ricci  coefficients   norms on $\Si_*$
\bea
\bsplit
\Skstar_k^2&:= \int_{\Si_*}\Big(r^2\big|\dk^{\le k} (\Xh, \trXc,  \Zc, \trXbc)\big|^2 +  \big| \dk^{\le k}(\Xbh, \Hc, \ombc, \Xib )\big|^2\Big)\\
&+\int_{\Si_*} \Big( \big| \dk^{\le k}   \widecheck{\DD\cos\th}   \big|^2 + \big| \dk^{\le k}   \DD r    \big|^2 +    \big| \dk^{\le k}   e_3 (\cos \th)    \big|^2  +  r^{-2}   \big| \dk^{\le k}  \widecheck{ e_3 (r)  }  \big|^2   \Big)\\
&+\int_{\Si_*}r^2 \Big(\big|\dk^{\le k} \DD\hot\Jk|^2  +\big|\dk^{\le k} \widecheck{\ov{\DD}\c\Jk}\big|^2 +\big| \dk^{\le k}\widecheck{\nab_3\Jk}\big|^2 \Big)\\
&+\int_{\Si_*}\left(\left|\dk_*^{\leq k}\left(b_*+1+\frac{2m}{r}\right)\right|+\left|\dk_*^{\leq k}\left(\nu(r)+2\right)\right|\right),
\end{split}
\eea
where $\trXc$, $\Xh$, $\Zc$, $\Hc$, $\trXbc$,  $\Xbh$, $\ombc$, $\Xib$ are the linearized Ricci coefficients   of the outgoing PT frame  of $\Mext$, $\nu=e_3+b_*e_4$ is tangent to $\Si_*$ with $(e_3, e_4)$ denoting the null pair attached to $\Si_*$, and $\dk_*$ corresponds weighted derivatives tangent to $\Si_*$, i.e. $\dk_*=(\nab_\nu, \dkb)$ with $\dkb$ associated to $\Si_*$.  
\end{definition} 

\begin{definition}
\lab{definition:PT-normson-RcSi_*} 
We define the curvature norm on $\Si_*$ 
\bea
\Rkstar^2_k := \int_{\Si_*}\left[r^{4+\dt}\big(|\dk^{\leq k}A|^2+|\dk^{\leq k}B|^2\big)+r^{4}|\dk^{\leq k}\Pc|^2+r^{2}|\dk^{\leq k}\Bb|^2+|\dk^{\leq k}\Ab|^2\right],
\eea
where $A$, $B$, $\Pc$, $\Bb$, $\Ab$ denote the linearized curvature components relative to  the outgoing PT  frame of $\Mext$. 
\end{definition}

%%%%%%%%%%%%%%%%%%%

\subsubsection{Norms of $\Mext$}

%%%%%%%%%%%%%%%%%%%

\begin{definition}
We define the following norms  for the PT-Ricci coefficients  of   $\Mext$
\bea
\nn\Skext_k^2 &:=& \sup_{\la \ge r_0}  \int_{r=\la}  r^2\big|\dk^{\le k} (\Xh, \trXc,  \Zc  )\big|^2 + r^{2-\dt}|\dk^{\le k}\trXbc |^2+ \big| \dk^{\le k}(\Xbh, \Hc, \ombc)\big|^2+ r^{-\dt}|\dk^{\le k}\Xib |^2\\
\nn&&+\sup_{\la \ge r_0}  \int_{r=\la}  \Big( \big| \dk^{\le k}   \widecheck{\DD\cos\th}   \big|^2 + \big| \dk^{\le k}   \DD r    \big|^2 + \big| \dk^{\le k}   e_3 (\cos \th)    \big|^2  +  r^{-2}   \big| \dk^{\le k}  \widecheck{ e_3 (r)  }  \big|^2   \Big)\\
&&+\sup_{r\ge r_0}  \int_{ r=\la}  r^2 \Big(\big|\dk^{\le k} \DD\hot\Jk|^2  +\big|\dk^{\le k} \widecheck{\ov{\DD}\c\Jk}\big|^2 +\big| \dk^{\le k}\widecheck{\nab_3\Jk}\big|^2 \Big),
\eea
where $\trXc$, $\Xh$, $\Zc$, $\Hc$, $\trXbc$,  $\Xbh$, $\ombc$, $\Xib$ are the linearized Ricci coefficients   of the outgoing PT frame of $\Mext$.
\end{definition}

\begin{definition}
We define  the following norms for the PT curvature coefficients in $\Mext$
\bea
\Rkext^2_k :=\int_{\Mext} r^{3+\dt}|\dk^{\le k}(A, B)|^2 +r^{3-\dt}\big(|\dk^{\le k}\Pc |^2 +r^{-2} |\dk^{\le k}\Bb|^2 +r^{-4} |\dk^{\le k}\Ab|^2\big),
\eea
where $A$, $B$, $\Pc$, $\Bb$, $\Ab$ denote the linearized curvature components relative to  the outgoing PT  frame of $\Mext$. 
\end{definition}

%%%%%%%%%%%%%%%%%%%

\subsubsection{Norms of $\Mint'$} 

%%%%%%%%%%%%%%%%%%%

Recall that $r\leq r_0$ on $\Mint'$. We thus discard $r$-weights in the Ricci and curvature norms of $\Mint'$ introduced below.  
    
\begin{definition}
\lab{definit:norms-SkMint'}
We define the following norms  for the PT-Ricci coefficients  of   $\Mint'$
\beaa
\Skint_k^2 := \int_{\Mint'}\big|\dk^{\le k}\Gac\big|^2, 
\eeaa
 where $\Gac$ denotes the set of all  linearized Ricci and metric coefficients  with respect  to the ingoing PT frame of $\Mint'$, i.e.   
  \beaa
\Gac :=\Big\{  \trXbc,  \Xbh,  \Zc,  \Hbc, \widecheck{ \DD\cos \th},   \omc,  \DD r, \widecheck{\DD u},    e_4(\cos\th), \widecheck{e_4(r)},  \widecheck{e_4(\ub)},    \widecheck{ \ov{\DD}\c \Jk},  \DD\hot\Jk, \widecheck{\nab_4 \Jk},             \trXc, \Xh,  \Xi\Big\}.
  \eeaa
 \end{definition}

For the curvature norms in $\Mint'$, we rely in particular on the scalar function $\tau$ introduced in section \ref{sec:defintionofthescalarfunctiontau:chap9}. We also introduce the following vectorfield\footnote{In Kerr, we have $\Rhat=\frac{\De}{r^2+a^2}\pr_r$ in Boyer Lindquist coordinates.} in $\Mint'$
   \bea
   \lab{eq:ThatRhat-e_3e_4}
\Rhat := \frac 1 2 \left( \frac{|q|^2}{r^2+a^2} e_4-\frac{\De}{r^2+a^2}  e_3\right).
 \eea

\begin{definition}
\lab{definit:norms-RkMint'}
We define the following norms for the curvature coefficients in $\Mint'$
\beaa
\Rkint_k^2 &=&  \int_{\Mint'}\Big( \big| \nab_{\Rhat} \dk^{\le k-1}\Rc\big|^2+|\dk^{\le k-1}\Rc|^2\Big) +\sup_\tau\int_{\Mint'\cap\Si(\tau)}  |\dk^{\le k}  \Rc|^2,
\eeaa
 where $\Rc=\{A, B, \Pc, \Bb, \Ab\} $ is the set of all linearized  curvature coefficients  w.r.t.  the ingoing PT frame  of $\Mint'$.  The derivative $\nab_{\Rhat}$ is taken with respect to the vectorfield $\Rhat$ defined in \eqref{eq:ThatRhat-e_3e_4}.
\end{definition}

%%%%%%%%%%%%%%%%%%%%%%

\subsubsection{Norms of $\Mtop'$}
\lab{section:NormsMtop}

%%%%%%%%%%%%%%%%%%%%%%

In the norms on $\Mtop'$ introduced below, we separate $\Mtop'$ in $\Mtop'(r\leq r_0)$ and $\Mtop'(r\geq r_0)$, and we discard $r$-weights in the region $\Mtop'(r\leq r_0)$.
\begin{definition}
We define the following norms for the Ricci  coefficients in $\Mtop'$.
\bea
\lab{definit:norms-SkMtop'}
\Sktop_k^2:= \int_{\Mtop'(r\leq r_0)} \big| \dk^{\le k} \Gac |^2 +(\Sktop^{\ge r_0}_{k})^2
\eea
where $\Gac$ denotes the set of all  linearized Ricci and metric coefficients  with respect  to the ingoing PT frame of $\Mtop'$ as above, and where
\bea
\lab{definit:norms-SkMtop'>r_0}
\bsplit
(\Sktop^{\ge r_0}_{k})^2 &:= \sup_{\ub_1\geq u_*'}\int_{\Mtop'_{r_0, \ub_1}}  r^2\big|\dk^{\le k} (\Xi, \omc, \Xh, \trXc,  \Zc, \Hbc)\big|^2 +r^{2-\dt}|\dk^{\leq k}\trXbc|+  \big| \dk^{\le k}\Xbh\big|^2\\
&+\sup_{\ub_1\geq u_*'}\int_{\Mtop'_{r_0, \ub_1}}  \Big( \big| \dk^{\le k}   \widecheck{\DD\cos\th}   \big|^2 + \big| \dk^{\le k}   \DD r    \big|^2 \Big)\\
&+\sup_{\ub_1\geq u_*'}\int_{\Mtop'_{r_0, \ub_1}}  \Big(  r^4\big| \dk^{\le k}   e_4 (\cos \th)    \big|^2  +  r^4\big| \dk^{\le k}  \widecheck{ e_4 (r)  }  \big|^2   \Big)\\
&+\sup_{\ub_1\geq u_*'}\int_{\Mtop'_{r_0, \ub_1}} r^2 \Big(\big|\dk^{\le k} \DD\hot\Jk|^2  +\big|\dk^{\le k} \widecheck{\ov{\DD}\c\Jk}\big|^2\Big) +r^6\big| \dk^{\le k}\widecheck{\nab_4\Jk}\big|^2 \Big), 
\end{split}
\eea
with $\Xi$, $\omc$, $\trXc$, $\Xh$, $\Zc$, $\Hbc$, $\trXbc$,  $\Xbh$ the linearized Ricci coefficients   of the ingoing PT frame of $\Mtop'$, and with the notation
\bea
\Mtop'_{r_0,\ub_1} &:=& \Mtop'(r\geq r_0)\cap\{\ub_1\leq \ub\leq \ub_1+1\}.
\eea
\end{definition}

For the curvature norms in $\Mtop'$, we rely in particular on the scalar function $\tau$ introduced in section \ref{sec:defintionofthescalarfunctiontau:chap9}.  
\begin{definition}
We define the following norms for the curvature coefficients in $\Mtop'$.
\bea
\lab{definit:norms-RkMtop'}
\bsplit
\Rktop_k^2 &:= \int_{\Mtop(r\geq r_0)}\Big( r^{3+\dt}|\dk^{\le k}(A, B)|^2 +r^{3-\dt}|\dk^{\le k}\Pc |^2\Big)\\
&+\sup_{\ub_1\geq u_*'}\int_{\Mtop_{r_0, \ub_1}}\Big( r^{2} |\dk^{\le k}\Bb|^2 +|\dk^{\le k}\Ab|^2\Big)+\sup_\tau\int_{\Mtop'(r\leq r_0)\cap\Si(\tau)}  |\dk^{\le k}  \Rc|^2 ,
\end{split}
\eea
where $\Rc$ is the set of all linearized  curvature coefficients  w.r.t.  the ingoing PT frame  of $\Mtop'$ as above, and $A$, $B$, $\Pc$, $\Bb$, $\Ab$ denote the linearized curvature components relative to  the ingoing PT  frame of $\Mtop'$.   
\end{definition}

%%%%%%%%%%%%%%%%%%%%

\subsubsection{Global norms}

%%%%%%%%%%%%%%%%%%%%%

We   define the global  norms for the PT frames of $\MM$ as follows
\bea
\bsplit
\Sk_k &=\Skstar_k +\Skext_k+\Skint_k +\Sktop, \\
\Rk_k &=\Rkstar_k+\Rkext_k+\Rkint_k +\Rktop.
\end{split}
\eea

%%%%%%%%%%%%%%%%%%%%%

\subsubsection{Initial data norms} 

%%%%%%%%%%%%%%%%%%%%%

Given $u, \ub$  defined respectively relative to the outgoing PT frame of $\Mext$ and the ingoing PT frame of $\Mint'$, we  define
\bea
\BB_1:= \big\{ u=1\big\}, \qquad \BBb_1:= \big\{ \ub=1\big\}.
\eea

\begin{definition}
We define the following initial data  norms on $\BB_1\cup \BBb_1$
\bea
\bsplit
\IkPText_{k} &:= \sup_{S\subset\BB_1 }r^{\frac{5}{2}  +\de_B} \Big(  \big\| \dk^k\,{}^{(ext)}A\big\|_{L^2(S)}  + \big\| \dk^k\,{}^{(ext)}B\big\|_{L^2(S)}\Big)\\
&+ \sup_{S\subset\BB_1 }\Big(  r^2 \big\| \dk^k\,{}^{(ext)} \Pc  \big\|_{L^2(S)}+  r  \big\| \dk^k\,{}^{(ext)}\Bb\big\|_{L^2(S)} + \big\| \dk^k\,{}^{(ext)}\Ab\big\|_{L^2(S)}\Big),\\
\IkPTint_{k} &:= \sup_{S\subset\BBb_1 } \Big(  \big\| \dk^k\,{}^{(int)}A\big\|_{L^2(S)}  + \big\| \dk^k\,{}^{(int)}B\big\|_{L^2(S)}\Big)\\
&+\sup_{S\subset\BBb_1 }\Big(  \big\| \dk^k\,{}^{(int)}\Pc \big\|_{L^2(S)}+  r \big\| \dk^k\,{}^{(int)}\Bb\big\|_{L^2(S)}\ + \big\| \dk^k\,{}^{(int)}\Ab\big\|_{L^2(S)}\Big),
\end{split}
\eea
where the linearized curvature components are taken respectively w.r.t. the outgoing PT frame of $\Mext$ on $\BB_1$ and the ingoing PT frame of $\Mint'$ on $\BBb_1$. We also set
\bea
\IkPT_k&:=\IkPText_{k} +\IkPTint_{k}. 
\eea
\end{definition}

%%%%%%%%%%%%%%%%%%%%%%%%%%%%%%%%%

\subsection{Statement of the Main  PT-Theorem}
\lab{sec:statementofthemainPTTheorem}

%%%%%%%%%%%%%%%%%%%%%%%%%%%%%%%%%

We are now ready to state the main result of this chapter on the control of the PT structures of $\MM$. 
\begin{theorem}[Main PT-Theorem]
\lab{theorem:Main-PT}
Consider  a GCM admissible  spacetime  verifying  the initial data assumptions 
  of the  Main Theorem  in section \ref{section:MainTheorem}, i.e. 
\beaa
\Ik_{k_{large}+10}\leq \ep_0, \qquad {}^{(ext)}\Ik_3\leq \ep_0^2.
 \eeaa
 Then, relative to the global norms defined in section \ref{sec:mainnormsPTframe:chap9} for the PT frames of $\MM$, we have the following bounds
 \bea
 \Sk_k+\Rk_k &\les \ep_0, \qquad k\le k_{large}+7.
 \eea
\end{theorem}

%%%%%%%%%%%%%%%%%%%%%%%%%%%%%%%%%

\subsection{Proof of Theorem M8}
\lab{sec:endoftheproofofTheoremM8:chap9}

%%%%%%%%%%%%%%%%%%%%%%%%%%%%%%%%%

Using Theorem \ref{theorem:Main-PT}, we are ready to prove Theorem  M8,  stated in section \ref{sec:endoftheproofofthemaintheorem:chap3}. We proceed in several steps.

\noindent{\bf Step 1.} Let $(f, \fb, \la)$ denote the transition coefficients corresponding to the change  from the outgoing PT frame of $\Mext$ to the outgoing PG frame of $\Mext$. Also, we denote the quantities corresponding to the outgoing PT frame without primes, and the  quantities corresponding to the outgoing PG frame with primes. 
In view of Lemma \ref{lemma:linkPGandTframeinMext}, the following identities hold in $\Mext$:
 \begin{enumerate}
 \item We have  $f=0$  and $\la=1$. In particular, we have 
\beaa
e_4' = e_4.  
\eeaa

\item  We have
\beaa
r'=r, \qquad \th'=\th, \qquad q'=q, \qquad \Jk'=\Jk.
\eeaa

\item With the notation $\Fb=\fb+i\dual\fb$, we have
\beaa
\nab_4\Fb+\frac{1}{2}\tr X\Fb &=& -2\Zc -\Fb\c\chih.
\eeaa
 \end{enumerate}
Also, note from the initialization of the outgoing PG frame of $\Mext$ on $\Si_*$ of section \ref{sec:initalizationadmissiblePGstructure}, and the initialization of the outgoing PT frame of $\Mext$ on $\Si_*$ of section \ref{sec:defintionofthePTstructuresinMM}, that we have on $\Si_*$
\beaa
\fb &=& -\frac{(\nu(r)-b_*)}{1-\frac{1}{4}b_*\frac{a^2(\sin\th)^2}{r^2}}\frac{a}{r}f_0 - \frac{a\Upsilon}{r}f_0.
\eeaa
We infer on $\Si_*$
\bea
\fb &=& -\frac{a}{r}\left(\frac{(\nu(r)+2)-\left(b_*+1+\frac{2m}{r}\right)}{1-\frac{1}{4}b_*\frac{a^2(\sin\th)^2}{r^2}} +O(r^{-2})\right)f_0. 
\eea

\noindent{\bf Step 2.} We first control the $L^2(\Si_*)$ norm of $\Hc'$ on $\Si_*$ as this quantity is part of the  boundedness norm for the outgoing PG frame of $\Mext$, see \eqref{equation:defboudednessnormsMext:chap3}. In view of the transformation formula for $\eta'$ in Proposition \ref{Proposition:transformationRicci}, together with the fact that $f=0$ and $\la=1$, we have
\beaa
\eta' &=& \eta  +\frac{1}{4}\fb\trch -\frac{1}{4}\dual\fb\atrch   +\frac{1}{2}\fb\c\chih.
\eeaa
Since $r'=r$, $\th'=\th$, $q'=q$, and $\Jk'=\Jk$ in $\Mext$, we infer
\beaa
\etac' &=& \etac  +\frac{1}{4}\fb\trch -\frac{1}{4}\dual\fb\atrch   +\frac{1}{2}\fb\c\chih.
\eeaa
and hence
\beaa
\|\dk^{\leq k}\Hc'\|_{L^2(\Si_*)} &\les& \|\dk^{\leq k}\Hc\|_{L^2(\Si_*)}+\left\|\dk^{\leq k}\left(\fb\big(\trch, \atrch, \chih\big)\right)\right\|_{L^2(\Si_*)}.
\eeaa
Together with the control of the PT frames provided by Theorem \ref{theorem:Main-PT}, we infer
\beaa
\|\dk^{\leq k_{large}+7}\Hc'\|_{L^2(\Si_*)} &\les& \ep_0+\left\|r^{-1}\dk^{\leq k_{large}+7}\fb\right\|_{L^2(\Si_*)}.
\eeaa
Using the formula for $\fb$ on $\Si_*$ of Step 1, and using again the control of the PT frames provided by Theorem \ref{theorem:Main-PT}, we obtain 
\beaa
\|\dk^{\leq k_{large}+7}\Hc'\|_{L^2(\Si_*)} &\les& \ep_0+\left\|\dk^{\leq k_{large}+7}O(r^{-3})\right\|_{L^2(\Si_*)}.
\eeaa
In view of the dominance condition \eqref{eq:behaviorofronS-star} for $r$ on $\Si_*$, we deduce
\bea
\|\dk^{\leq k_{large}+7}\Hc'\|_{L^2(\Si_*)} &\les& \ep_0
\eea
which yields the desired behavior for $\Hc'$ on $\Si_*$. 

\noindent{\bf Step 3.} The remaining estimates for the PG frames of $\MM$ being all in sup norm, we first derive sup norms estimates for the PT frames of $\MM$. In view of  the control of the PT frames provided by Theorem \ref{theorem:Main-PT}, together with Sobolev and the trace theorem, we obtain, for $k\leq k_{large}+4$,  
\bea\lab{eq:applicationofSobolevtothePTframeinproofThmM8}
\bsplit
\sup_{\Mext} \Big\{ r^2| \dk^{\leq k}\Ga_g|  + r| \dk^{\leq k}\Ga_b|+r^{\frac{7}{2}+\frac{\dt}{2}}\big(|\dk^{\leq k}A|+|\dk^{\leq k}B|\big)\Big\}\\
+\sup_{\Mtop'} \Big\{ r^2| \dk^{\leq k}\Ga_g|  + r| \dk^{\leq k}\Ga_b|+r^{\frac{7}{2}+\frac{\dt}{2}}\big(|\dk^{\leq k}A|+|\dk^{\leq k}B|\big)\Big\}\\
+\sup_{\Mint'} \Big\{ | \dk^{\leq k}\Ga_g|  + | \dk^{\leq k}\Ga_b|\Big\} &\les \ep_0,
\end{split}
\eea
where in each case, $(\Ga_g, \Ga_b)$ is defined w.r.t. the linearized quantities in the PT frame of the corresponding region.

\begin{remark}
In view of the definition of $\Skext_k$ and $\Sktop_k$, \eqref{eq:applicationofSobolevtothePTframeinproofThmM8} holds a priori  only for $\Ga_g\setminus\{\trXbc\}$ and $\Ga_b\setminus\{\Xib\}$, while $\trXbc$ and $\Xib$  satisfy a priori the following weaker estimates (in terms of powers of $r$) in $\Mext$ and $\Mtop$,  for $k\leq k_{large}+4$,  
\beaa
\sup_{\Mext}\Big(r^{2-\frac{\dt}{2}}|\dk^{\leq k}\trXbc|+r^{1-\frac{\dt}{2}}|\dk^{\leq k}\Xib|\Big)+\sup_{\Mtop}r^{2-\frac{\dt}{2}}|\dk^{\leq k}\trXbc| &\les& \ep_0.
\eeaa
To recover the claimed estimates for $\trXbc$ and $\Xib$ of \eqref{eq:applicationofSobolevtothePTframeinproofThmM8} in $\Mext$, it suffices to integrate the transport equations for $\nab_4\trXc$ and $\nab_4\Xib$ of Proposition \ref{Prop:linearizedPTstructure1}, using the control provided by \eqref{eq:applicationofSobolevtothePTframeinproofThmM8} for $\Ga_g\setminus\{\trXbc\}$ and $\Ga_b\setminus\{\Xib\}$ for the RHS, and the control of $\trXbc$ and $\Xib$ on $\Si_*$ (which has no loss in $r$). Finally, integrating along the $e_3$ direction in $\Mtop$ the equation for $\nab_3\trXbc$ yields the desired claim \eqref{eq:applicationofSobolevtothePTframeinproofThmM8} also for $\trXbc$ in $\Mtop$  so that \eqref{eq:applicationofSobolevtothePTframeinproofThmM8}  indeed holds true.
\end{remark}

We will also need a sharper estimate for $\Zc$ on $\Mext$. Recall from Proposition \ref{Prop:linearizedPTstructure1} that $\Zc$, in the outgoing PT frame of $\Mext$, satisfies on $\Mext$
\beaa
\nab_4\Zc + \frac{1}{q}\Zc &=&   - \frac{a\ov{q}}{|q|^2}\trXc \Jk    -\frac{aq}{|q|^2}\ov{\Jk}\c\widehat{X}  -B+\Ga_g\c\Ga_g.
\eeaa
We infer from the above control, on $\Si_*$, 
\beaa
\left|\dk^{\leq k_{large}+4}(q\Zc)\right| &\les& \frac{\ep_0}{r^{\frac{5}{2}}}.
\eeaa
 Integrating from $\Si_*$, we infer the following improved bound on $\Mext$
\bea
\left|\dk^{\leq k_{large}+4}\Zc\right| &\les& \frac{\ep_0}{r_*^2}+\frac{\ep_0}{r^{\frac{5}{2}}}.
\eea

\noindent{\bf Step 4.} Next, we estimate $\fb$ on $\Mext$. First, in view of the identity for $\fb$ derived on $\Si_*$ in Step 1, together with the control of Step 3 for the PT frame of $\Mext$, and the dominant condition for $r$ on $\Si_*$, we have 
\beaa
\sup_{\Si_*}r|\dk^{\leq k_{large}+4}\fb|\les \ep_0.
\eeaa
Also, recall from Step 1 the following transport equation on $\Mext$
\beaa
\nab_4\Fb+\frac{1}{2}\tr X\Fb &=& -2\Zc -\Fb\c\chih,
\eeaa
where $\Fb=\fb+i\dual\fb$. Together with the control of the PT frames provided by Step 3, we infer
\beaa
\nab_4\Fb+\frac{1}{2}\tr X\Fb &=& -2\Zc -\Fb\c\chih,
\eeaa
In view of the bounds provided by Step 3 for the outgoing PT frame of $\Mext$, using in particular the improved bound for $\Zc$, we infer, on $\Mext$, 
\beaa
\left|\dk^{\leq k_{large}+4}\nab_4(q\Fb)\right| &=& \frac{\ep_0}{r_*}+\frac{\ep_0}{r^{\frac{3}{2}}} +\frac{\ep_0}{r^2}\left|\dk^{\leq k_{large}+4}(q\Fb)\right|.
\eeaa
Integrating from $\Si_*$, and using the above control for $\fb$ on $\Si_*$, we infer, for $\ep_0$ small enough, 
\beaa
\left|\dk^{\leq k_{large}+4}\Fb\right| &\les& \frac{r_*-r}{rr_*}\ep_0+\frac{\ep_0}{r^{\frac{3}{2}}}
\eeaa
and hence
\bea
\sup_{\Mext}r\left|\dk^{\leq k_{large}+4}\fb\right| &\les& \ep_0.
\eea

\noindent{\bf Step 5.} Next, let $(f', \fb', \la')$ denote the transition coefficients corresponding to the change  from the ingoing PT frame of $\Mtop'$ to the ingoing PG frame of $\Mtop$. Since $\Mtop\subset\Mtop'$, $(f', \fb', \la')$ are defined on $\Mtop$. In view of
\begin{enumerate}
\item The fact that $\{u=u_*\}=\Mtop\cap\Mext$.

\item The fact that $\la=1$, $f=0$, and $\fb$ is controlled on $\Mext$ in Step 4, where  $(f, \fb, \la)$ denote the transition coefficients corresponding to the change  from the outgoing PT frame of $\Mext$ to the outgoing PG frame of $\Mext$.

\item The initialization of the ingoing PG frame of $\Mtop$ from the ingoing PG frame of $\Mext$ in section \ref{sec:initalizationadmissiblePGstructure}.

\item The control in $\Mext\cap\Mtop'$ of the change of frame coefficients between the outgoing PT frame of $\Mext$ and the ingoing PT frame of $\Mtop'$ which results from
\begin{enumerate}
\item the  initialization of the ingoing PT frame of $\Mtop'$ from the outgoing PT frame of $\Mext$ on $\{u=u_*'\}$ in section \ref{sec:defintionofthePTstructuresinMM}, 

\item the transport equations for transition coefficients involving PT frames in section \ref{section:transport(f,fb,la):PTcase},
\end{enumerate}
 \end{enumerate}
we easily obtain
\beaa
\sup_{\{u=u_*\}}r\left|\dk^{\leq k_{large}+4}(f', \fb', \log(\la'))\right| &\les& \ep_0.
\eeaa

Next, using the analog of Corollary \ref{cor:transportequationine4forchangeofframecoeffinformFFbandlamba} for ingoing foliations, and the fact that both  the ingoing PT frame of $\Mtop'$ and the ingoing PG frame of $\Mtop$ verify $\xib=\omb=0$, we obtain the following transport equations 
\beaa
\nab_{{\la'}^{-1}e_3'}\Fb'+\frac{1}{2}\ov{\tr\Xb}\,\Fb' &=&  -\chibh\c\Fb'+\underline{E}_1(\fb', \Ga),\\
{\la'}^{-1}\nab_3'(\log\la') &=& \fb'\c(-\ze-\eta)+\underline{E}_2(\fb', \Ga),
\eeaa
where $\Fb'=\fb'+i\dual\fb'$, and where the Ricci coefficients appearing are the ones of the ingoing PT frame of $\Mtop'$. Integrating both transport equations from $\{u=u_*\}$ starting with the one for $\fb'$, using the control of the ingoing PT frame of $\Mtop'$ provided by Step 3, and the above control for $(\fb', \la')$ on $\{u=u_*\}$, we obtain 
\beaa
\sup_{\Mtop}r\left|\dk^{\leq k_{large}+4}(\fb', \log(\la'))\right| &\les& \ep_0.
\eeaa

Also, using again the analog of Corollary \ref{cor:transportequationine4forchangeofframecoeffinformFFbandlamba} for ingoing foliation, we have
\beaa
\nab_{{\la'}^{-1}e_3'}F'+\frac{1}{2}\tr\Xb F' &=& -2(H-Z)   +2\DD'(\log\la')  +2\om F'  +E_3({\nab'}^{\leq 1}\fb', f', \Ga, {\la'}^{-1}\chib'),
\eeaa
where $F'=f'+i\dual f'$, and where the Ricci coefficients appearing are the ones of the ingoing PT frame of $\Mtop'$. Now, we have $Z-H = \Zc$ in view of the fact that $H=\frac{aq}{|q|^2}\Jk$ for ingoing PT foliations, and hence 
\beaa
\nab_{{\la'}^{-1}e_3'}F'+\frac{1}{2}\tr\Xb F' &=& -2\Zc   +2\DD'(\log\la')  +2\om F'  +E_3({\nab'}^{\leq 1}\fb', f', \Ga, {\la'}^{-1}\chib').
\eeaa
Integrating the transport equations from $\{u=u_*\}$, using the control of the ingoing PT frame of $\Mtop'$ provided by Step 3, the above control of $(\fb', \la')$ on $\Mtop$, and the above control for $f'$ on $\{u=u_*\}$, we obtain 
\beaa
\sup_{\Mtop}r\left|\dk^{\leq k_{large}+3}f'\right| &\les& \ep_0.
\eeaa

\noindent{\bf Step 6.} We are now ready to conclude the proof of Theorem M8. First, let $(f'', \fb'', \la'')$ denote the transition coefficients corresponding to the change  from the ingoing PT frame of $\Mint'$ to the ingoing PG frame of $\Mint$. Starting from the timelike hypersurface $\Mext\cap\Mint=\{r=r_0\}\cap\{u\leq u_*\}$, and arguing as in Step 5, one easily obtains 
\beaa
\sup_{\Mint}\Big(\left|\dk^{\leq k_{large}+4}(\fb'', \log(\la'')\right|+\left|\dk^{\leq k_{large}+3}f''\right|\Big) &\les& \ep_0,
\eeaa

Next, relying on:
\begin{enumerate}
\item the fact that $\la=1$ and $f=0$ on $\Mext$ according to Step 1,

\item the control of $\fb$ on $\Mext$ in Step 4, 

\item the control of $(f', \fb', \la')$ on $\Mtop$ in Step 5, 

\item the above control of $(f'', \fb'', \la'')$ on $\Mint$,

\item the sup norm control of the PT structures of $\MM$ in Step 3 

\item the change of frame formulas of Proposition \ref{Proposition:transformationRicci},
\end{enumerate} 
we infer, for $k\leq k_{large}+2$, 
\beaa
&&\sup_{\Mext} \Big\{ r^2| \dk^{\leq k}\Ga_g'|  + r| \dk^{\leq k}\Ga_b'|+r^{\frac{7}{2}+\frac{\dt}{2}}\big(|\dk^{\leq k}A'|+|\dk^{\leq k}B'|\big)\Big\}\\
&&+\sup_{\Mtop} \Big\{ r^2| \dk^{\leq k}\Ga_g'|  + r| \dk^{\leq k}\Ga_b'|+r^{\frac{7}{2}+\frac{\dt}{2}}\big(|\dk^{\leq k}A'|+|\dk^{\leq k}B'|\big)\Big\}\\
&&+\sup_{\Mint} \Big\{ | \dk^{\leq k}\Ga_g'|  + | \dk^{\leq k}\Ga_b'|\Big\}\\ 
&\les& \ep_0+\sup_{\Mext}r\left|\dk^{\leq k+1}\Big(r'-r, \cos(\th')-\cos(\th), q'-q, r(\Jk'-\Jk)\Big)\right| \\
&&+\sup_{\Mtop}r\left|\dk^{\leq k+1}\Big(r'-r, \cos(\th')-\cos(\th), q'-q, r(\Jk'-\Jk)\Big)\right|\\
&&+\sup_{\Mint}\left|\dk^{\leq k+1}\Big(r'-r, \cos(\th')-\cos(\th), q'-q, \Jk'-\Jk\Big)\right| 
\eeaa
where in each case, $(\Ga_g', \Ga_b')$ is defined w.r.t. the linearized quantities in the PG frame of the corresponding region, and $r'-r$, $\cos(\th')-\cos(\th)$, $q'-q$ and $\Jk'-\Jk$ correspond to the difference between the un-primed quantity in the PT frame and the primed quantity in the PG frame of the corresponding region. In view of the definition of the combined sup norm $\Nk^{(Sup)}_k$ for the PG structures of $\Mext$, $\Mint$ and $\Mtop$, see section \ref{sec:definitionofconcatenatednorm}, and using also the estimate for $\Hc'$ on $\Si_*$ derived in Step 2, this yields
\beaa
 \Nk^{(Sup)}_{k_{large}+2} &\les& \ep_0+\sup_{\Mext}r\left|\dk^{\leq k_{large}+3}\Big(r'-r, \cos(\th')-\cos(\th), q'-q, r(\Jk'-\Jk)\Big)\right| \\
&&+\sup_{\Mtop}r\left|\dk^{\leq k_{large}+3}\Big(r'-r, \cos(\th')-\cos(\th), q'-q, r(\Jk'-\Jk)\Big)\right|\\
&&+\sup_{\Mint}\left|\dk^{\leq k_{large}+3}\Big(r'-r, \cos(\th')-\cos(\th), q'-q, \Jk'-\Jk\Big)\right|. 
\eeaa

Next, recall from Step 1 that we have on $\Mext$
\beaa
r'=r, \qquad \th'=\th, \qquad q'=q, \qquad \Jk'=\Jk.
\eeaa
Hence, we obtain 
\beaa
 \Nk^{(Sup)}_{k_{large}+2} &\les& \ep_0 +\sup_{\Mtop}r\left|\dk^{\leq k_{large}+3}\Big(r'-r, \cos(\th')-\cos(\th), q'-q, r(\Jk'-\Jk)\Big)\right|\\
&&+\sup_{\Mint}\left|\dk^{\leq k_{large}+3}\Big(r'-r, \cos(\th')-\cos(\th), q'-q, \Jk'-\Jk\Big)\right|. 
\eeaa
It thus remains to control the quantities $r'-r$, $\cos(\th')-\cos(\th)$, $q'-q$ and $\Jk'-\Jk$ on $\Mtop$ and $\Mint$. Denoting in each region by $(e_1', e_2', e_3', e_4')$ the corresponding ingoing PG frame, and by $(e_1, e_2, e_3, e_4)$ the corresponding ingoing PT frame, we have
\beaa
&& e_3(r)=-1, \qquad e_3(\th)=0, \qquad e_3(q)=-1, \qquad \nab_3(\ov{q}\Jk)=0,\\
&& e_3'(r')=-1, \qquad e_3'(\th')=0, \qquad e_3'(q')=-1, \qquad \nab_3'(\ov{q'}\Jk')=0.
\eeaa
Expressing $e_3'$ on the frame $(e_1, e_2, e_3, e_4)$ using the fame transformation formulas, using the control of  $(f', \fb', \la')$ on $\Mtop$ in Step 5, and the above control of $(f'', \fb'', \la'')$ on $\Mint$, we easily infer
\beaa
\sup_{\Mtop}r\left|\dk^{\leq k_{large}+3}\Big(e_3'(r'-r), e_3'(\cos(\th')-\cos(\th)), e_3'(q'-q), \nab_3'(q'(\Jk'-\Jk))\Big)\right|\\
+ \sup_{\Mint}\left|\dk^{\leq k_{large}+3}\Big(e_3'(r'-r), e_3'(\cos(\th')-\cos(\th)), e_3'(q'-q), \nab_3'(q'(\Jk'-\Jk))\Big)\right| &\les& \ep_0.
\eeaa
Together with the fact that we have on $\Mext$
\beaa
r'=r, \qquad \th'=\th, \qquad q'=q, \qquad \Jk'=\Jk,
\eeaa
we infer, integrating the transport equations in $e_3'$ from $\Mext\cap\Mtop=\{u=u_*\}$ to $\Mtop$ and from $\Mext\cap\Mint=\{r=r_0\}$ to $\Mint$, 
\beaa
\sup_{\Mtop}r\left|\dk^{\leq k_{large}+3}\Big(r'-r, \cos(\th')-\cos(\th), q'-q, r(\Jk'-\Jk)\Big)\right|\\
+\sup_{\Mint}\left|\dk^{\leq k_{large}+3}\Big(r'-r, \cos(\th')-\cos(\th), q'-q, \Jk'-\Jk\Big)\right| &\les& \ep_0.
\eeaa
Plugging in the above, we deduce 
\beaa
 \Nk^{(Sup)}_{k_{large}+2} &\les& \ep_0 
\eeaa
as desired. This concludes the proof of Theorem M8.

%%%%%%%%%%%%%%%%%%%%%%%%%%%%%%%%%

\subsection{Bootstrap assumptions for the  Main  PT-Theorem}

%%%%%%%%%%%%%%%%%%%%%%%%%%%%%%%%%

To prove Theorem \ref{theorem:Main-PT}, we make the following bootstrap assumptions.

{\bf  BA-PT.}  \quad   Relative to the global norms defined in section \ref{sec:mainnormsPTframe:chap9} for the PT frames of $\MM$, we have  
\bea\lab{eq:mainbootassforchapte9}
\Sk_k+\Rk_k &\leq& \ep, \qquad k\le k_{large}+7.
\eea

%%%%%%%%%%%%%%%%%%%%%%%%%%%%%%%%%

 \subsection{Control of the initial data in the PT frames} 

%%%%%%%%%%%%%%%%%%%%%%%%%%%%%%%%%

Recall that the control of the initial data for the PG frames of $\MM$ is provided by Theorem M0, see section \ref{sec:mainintermediateresults:chap3}. We will need a sharper\footnote{Sharper in terms of derivatives. Indeed, the conclusions of Theorem M0 hold for $k\le k_{large}-2$ while the ones of Theorem M0-PT hold for $k\le k_{large}+7$.} analog for the PT frames of $\MM$ which we state below.

\begin{theorem}[Theorem M0-PT]
\lab{Theorem:TheoremM0-PT} 
Assume that the initial data layer  $\LL_0$, as   defined  in section \ref{sec:defintionoftheinitialdatalayer},  satisfies
\beaa
\Ik_{k_{large}+10}\leq \ep_0, \qquad {}^{(ext)}\Ik_3\leq \ep_0^2.
\eeaa
Then,  under the bootstrap assumptions  {\bf BA-PT},  relative to the initial data norms defined in section \ref{sec:mainnormsPTframe:chap9} for the PT frames of $\MM$, we have
 \bea
 \IkPT_k &\les \ep_0, \qquad k\le k_{large}+7.
 \eea 
\end{theorem}

The proof of Theorem M0-PT is postponed to  section \ref{sec:proofof:Theorem:TheoremM0-PT}.

%%%%%%%%%%%%%%%%%%%%%%%%%%%%%%%%%

\subsection{Control of low derivatives of the PT frame}

%%%%%%%%%%%%%%%%%%%%%%%%%%%%%%%%%

The following lemma will allow us to initiate an iterative procedure in section \ref{sec:iterativeprocedureforThM8}.
\begin{lemma}\lab{lemma:proofofeq:initializationofiterationassumptioninproofThmM8}
Relative to the global norms defined in section \ref{sec:mainnormsPTframe:chap9} for the PT frames of $\MM$, we have the following bounds
\bea\lab{eq:initializationofiterationassumptioninproofThmM8}
\Sk_{k_{small}-1} + \Rk_{k_{small}-1} &\les&\ep_0.
 \eea
\end{lemma}

\begin{proof}
The PG structures of $\MM$ satisfy in view of Theorem M7,   stated in section \ref{sec:endoftheproofofthemaintheorem:chap3}, the following decay estimates, for $k\leq k_{small}$, 
\bea
\bsplit
\sup_{\Mext\cup\Mtop(r\geq r_0)} \Big\{ r^2u^{\frac{1}{2}+\dec}| \dk^{\leq k}\Ga_g|  + ru^{1+\dec}| \dk^{\leq k}\Ga_b|\Big\}\\
+\sup_{\Mtop(r\leq r_0)\cup\Mint} \ub^{1+\dec}\Big\{ | \dk^{\leq k}\Ga_g|  + | \dk^{\leq k}\Ga_b|\Big\} &\les \ep_0,
\end{split}
\eea
where in each case, $(\Ga_g, \Ga_b)$ is defined w.r.t. the linearized quantities in the PG frame of the corresponding region.

The proof is similar is spirit to the proof of Theorem M8 in section \ref{sec:endoftheproofofTheoremM8:chap9}. Here, instead of transferring the control of the PT frames to the PG frames of $\MM$, we instead transfer the above control of the PG frames to the PT frames. We proceed as follows:
\begin{enumerate}
\item We start with the control in $\Mext$. Let $(f, \fb, \la)$ denote the transition coefficients corresponding to the change  from the outgoing PG frame of $\Mext$ to the outgoing PT  frame of $\Mext$. Also, we denote the quantities corresponding to the outgoing PG frame without primes, and the  quantities corresponding to the outgoing PT frame with primes. 
In view of Lemma \ref{lemma:linkPGandTframeinMext}, the following identities hold in $\Mext$:
 \begin{enumerate}
 \item We have  $f=0$  and $\la=1$. In particular, we have 
\beaa
e_4' = e_4.  
\eeaa

\item  We have
\bea
r'=r, \qquad \th'=\th, \qquad q'=q, \qquad \Jk'=\Jk.
\eea

\item We have
\bea
\nab_4\fb &=& 2\zec.
\eea
 \end{enumerate}
Also, note from the initialization of the outgoing PG frame of $\Mext$ on $\Si_*$ of section \ref{sec:initalizationadmissiblePGstructure}, and the initialization of the outgoing PT frame of $\Mext$ on $\Si_*$ of section \ref{sec:defintionofthePTstructuresinMM}, that we have on $\Si_*$
\bea
\fb &=& \frac{a}{r}\left(\frac{(\nu(r)+2)-\left(b_*+1+\frac{2m}{r}\right)}{1-\frac{1}{4}b_*\frac{a^2(\sin\th)^2}{r^2}} +O(r^{-2})\right)f_0. 
\eea

\item Next, as in Step 4 of the proof of  Theorem M8 in section \ref{sec:endoftheproofofTheoremM8:chap9}, we control $\fb$ using the above transport equation for $\fb$ and the above initialization on $\Si_*$. Based on the above control of the outgoing PG structure of $\Mext$, we easily obtain 
\bea
\sup_{\Mext}ru^{\frac{1}{2}+\dec}|\dk^{\leq k_{small}}\fb| &\les& \ep_0.
\eea

\item Next, let $(f', \fb', \la')$ denote the transition coefficients corresponding to the change  from the ingoing PG frame of $\Mtop$ to the ingoing PT frame of $\Mtop'$. As in Step 5 of the proof of  Theorem M8 in section \ref{sec:endoftheproofofTheoremM8:chap9}, we control $(f', \fb', \la')$. We now rely on the analog for ingoing PT structures of the transport equations of Corollary \ref{cor:transportequationine4forchangeofframecoeffinformFFbandlamba:PTcase}. Note that these transport equations do not loose derivatives\footnote{Unlike their analogs for PG structures in Corollary \ref{cor:transportequationine4forchangeofframecoeffinformFFbandlamba}.} and we obtain 
\bea
\nn\sup_{\Mtop'(r\geq r_0)}ru^{\frac{1}{2}+\dec}|\dk^{\leq k_{small}}(f', \fb', \log(\la'))|\\
+\sup_{\Mtop'(r\leq r_0)}\ub^{\frac{1}{2}+\dec}|\dk^{\leq k_{small}}(f', \fb', \log(\la'))| &\les& \ep_0.
\eea

\item Next, let $(f'', \fb'', \la'')$ denote the transition coefficients corresponding to the change  from the ingoing PG frame of $\Mint$ to the ingoing PT frame of $\Mint'$. As in Step 6 of the proof of  Theorem M8 in section \ref{sec:endoftheproofofTheoremM8:chap9}, we control $(f'', \fb'', \la'')$, now relying on the analog for ingoing PT structures of the transport equations of Corollary \ref{cor:transportequationine4forchangeofframecoeffinformFFbandlamba:PTcase}, and  obtain 
\bea
\sup_{\Mint'}\ub^{\frac{1}{2}+\dec}|\dk^{\leq k_{small}}(f'', \fb'', \log(\la'')) &\les& \ep_0.
\eea

\item As in Step 6 of the proof of  Theorem M8 in section \ref{sec:endoftheproofofTheoremM8:chap9}, we then rely on the above control of the various change of frame coefficients, the above control of the PG structures of $\MM$, and the change of frame formulas of Proposition \ref{Proposition:transformationRicci} to infer, for $k\leq k_{small}-1$, 
\bea
\bsplit
\sup_{\Mext\cup\Mtop'(r\geq r_0)} \Big\{ r^2u^{\frac{1}{2}+\dec}| \dk^{\leq k}\Ga_g'|  + ru^{1+\dec}| \dk^{\leq k}\Ga_b'|\Big\}\\
+\sup_{\Mtop'(r\leq r_0)\cup\Mint'} \ub^{\frac{1}{2}+\dec}\Big\{ | \dk^{\leq k}\Ga_g'|  + | \dk^{\leq k}\Ga_b'|\Big\} &\les \ep_0,
\end{split}
\eea
where in each case, $(\Ga_g', \Ga_b')$ is defined w.r.t. the linearized quantities in the PT frame of the corresponding region. The weights in $r$, $u$ and $\ub$ are enough to take care of the spacetime integrations in 
 the global norms defined in section \ref{sec:mainnormsPTframe:chap9} for the PT frames of $\MM$, and we finally obtain  the following desired bounds for the PT frames of $\MM$
\bea
\Sk_{k_{small}-1} + \Rk_{k_{small}-1} &\les&\ep_0.
 \eea
\end{enumerate}
This concludes the proof of Lemma \ref{lemma:proofofeq:initializationofiterationassumptioninproofThmM8}.
\end{proof}

%%%%%%%%%%%%%%%%%%%%%%%%%%%%%%%%%%%%

\subsection{Iterative procedure for the proof of the Main PT-Theorem}
\lab{sec:iterativeprocedureforThM8}

%%%%%%%%%%%%%%%%%%%%%%%%%%%%%%%%%%%%%

First, recall our  bootstrap assumptions {\bf BA-PT} which are assumed throughout the proof of Theorem \ref{Theorem:TheoremM0-PT} 
\bea\lab{eq:bootstrapassumptiondiscussionThM8}
\Sk_{k_{large}+7} +\Rk_{k_{large}+7}&\leq& \ep.
\eea

 For $J$  in the range  $k_{small}-1\leq J\leq k_{large}+5$, we also make  the following iteration assumption 
\bea\lab{eq:iterationassumptiondiscussionThM8:bis}
\Sk_{J} +\Rk_{J}   &\les& \ep_0+L_*(J),
\eea
with $L_*(k)$ given by 
  \bea\lab{eq:defintionofLstarofk:chap9}
 L^2_*(k) := L^2_*(k, u_*') =  \int_{\{u=u_*'\}}    \big|\Rc\big|_{w, k}^2 +\int_{\Si_*\cap\{u=u_*'\}} \big|\Gac\big|_{w, k}^2
 \eea 
where we recall the notation \eqref{eq:defintionofthenormsRcGacweighted}
 \bea\lab{eq:defintionofthenormsRcGacweighted:bis}
 \bsplit
 \big|\Rc\big|_{w, k}^2 &:= r^{3+\dt}|\dk_*^{\le k}(A, B)|^2 +r^{3-\dt}|\dk_*^{\le k}\Pc |^2 +r^{1-\dt}|\dk_*^{\le k}\Bb|^2,\\
\big|\Gac\big|_{w, k}^2 &:= r^2|\dkb^{\leq k}\Ga_g|^2+|\dkb^{\leq k}\Ga_b|^2.
\end{split}
\eea
In \eqref{eq:defintionofthenormsRcGacweighted:bis}, $(A, B, \Pc, \Bb)$ denote  linearized curvature components w.r.t. the the  outgoing PT frame of $\Mext$,  
$\Ga_g$ and $\Ga_b$ are  defined w.r.t. the  outgoing PT frame of $\Mext$ as in Definition \ref{definition.Ga_gGa_b:outgoingPTcase:chap9}, $\dk_*$ denote weighted derivatives tangential to the hypersurface $\{u=u_*'\}$, and $\dkb$ denote weighted derivatives tangential to the sphere $\Si_*\cap\{u=u_*'\}$.

\begin{remark}\lab{rmk:theiterationassumptionholdsforklsmallminus2:chap9}
In view of \eqref{eq:initializationofiterationassumptioninproofThmM8}, i.e. 
\beaa
\Sk_{k_{small}-1} +\Rk_{k_{small}-1}   &\les& \ep_0,
\eeaa
\eqref{eq:iterationassumptiondiscussionThM8:bis} holds for $J=k_{small}-1$.
\end{remark}

We now state  the main sequence of estimates which  will allow us   to prove Theorem \ref{Theorem:TheoremM0-PT}  in the next section.

\begin{theorem}[Control of Curvature]
\lab{prop:rpweightedestimatesiterationassupmtionThM8}
Let $J$ such that $k_{small}-1\leq J\leq k_{large}+6$. Under the iteration assumption \eqref{eq:iterationassumptiondiscussionThM8:bis}, we have the following estimate in $\MM$ for the global PT  curvature norm of section \ref{sec:mainnormsPTframe:chap9}
\bea
 \Rk_{J+1} \les r_0^{-\dt}\Big( \Sktop^{\geq r_0}_{J+1} +\Skext_{J+1}\Big)+r_0^{10} \big(  \Rk_J+\Sk_J  \big)  +\ep_0
 \eea 
where the constant in $\les$ is independent of $r_0$. 
\end{theorem}

\begin{remark}
Theorem \ref{prop:rpweightedestimatesiterationassupmtionThM8} will be proved in a separate paper, see \cite{KS:Kerr-B}. The proof relies on energy, Morawetz and $r^p$ weighted estimates for the Bianchi system. Note also that the statement of  Theorem \ref{prop:rpweightedestimatesiterationassupmtionThM8} is the analog of Proposition 8.3.9 in \cite{KS} (which holds in the particular case of polarized perturbations of Schwarzschild). 
\end{remark}

\begin{proposition}\lab{prop:controlGaSistartiterationassupmtionThM8}
Let $J$ such that $k_{small}-1\leq J\leq k_{large}+6$. The following estimate holds true on $\Si_*$ for the Ricci and metric coefficients of the outgoing PT frame of $\Mext$
\bea
 \Skstar_{J+1}&\les&  \Rkstar_{J+1}+\ep_0,
 \eea
 where the constant in $\les$ is independent of $r_0$. 
\end{proposition}

\begin{proof}
 The proof is given in  section \ref{sec:theoremM8recoverRiccionSistar}, see Proposition \ref{Prop:controlGa-PTframe-Si_*}.
\end{proof} 

\begin{proposition}\lab{prop:controlGaextiterationassupmtionThM8}
Let $J$ such that $k_{small}-1\leq J\leq k_{large}+6$. The following estimates hold true for the Ricci and metric coefficients of the outgoing PT frame of $\Mext$
\bea
 \Skext_{J+1}&\les&  \Skstar_{J+1}+\Rkext_{J+1} +\ep_0,
 \eea
  where the constant in $\les$ is independent of $r_0$.
\end{proposition}

\begin{proof}
The proof is given in  section \ref{sec:theoremM8recoverRicciawaytrapping}, see Proposition \ref{Prop:EstmatesSkext}.
\end{proof}

\begin{proposition}\lab{prop:controlGaintiterationassupmtionThM8}
Let $J$ such that $k_{small}-1\leq J\leq k_{large}+6$.   The following estimates hold true for the Ricci and metric coefficients of the ingoing PT frame of $\Mint'$
\bea
 \Skint_{J+1} &\les&   \Skext_{J+1}+\Rkint_{J+1}+\ep_0.
 \eea
\end{proposition}

\begin{proof}
The proof is given in  section \ref{sec:controlofthePTRiccicoefficientsinMintprime}, see  Proposition \ref{Prop:MainestimatesMint}.
\end{proof}

\begin{proposition}
\lab{prop:improvementoftheiterationassupmtionThM8-Mtop}
Let $J$ such that $k_{small}-1\leq J\leq k_{large}+6$.  The following estimates hold true for the Ricci and metric coefficients of the ingoing PT frame of $\Mtop'$
\bea
 \Sktop_{J+1} &\les& \ep_0+   L_*(J+1)  +   \Rktop_{J+1}.
 \eea
 In particular, we have
\bea
 \Sktop^{\geq r_0}_{J+1} &\les& \ep_0+   L_*(J+1)  +   \Rktop_{J+1}
 \eea
  in which  case the constant in $\les$ is independent of $r_0$.
 \end{proposition}

\begin{proof}
The proof is given in section \ref{sec:controlRiccicoeffMtop}.
\end{proof} 

The following is a corollary of the above propositions.
\begin{corollary}
\lab{Corr:improvementoftheiterationassupmtionThM8-Mtop}
Let $J$ such that $k_{small}-1\leq J\leq k_{large}+6$. 
\begin{enumerate}
\item  We have
 \bea
\Rk_{J+1}+\Skstar_{J+1} &\les& \ep_0 +  r_0^{-\dt}  L_*(J+1)  +r_0^{10} \big(  \Rk_J+\Sk_J  \big),
 \eea
 where the constant in $\les$ is independent of $r_0$. 
 
\item   We have 
\bea
  \Sk_{J+1} &\les& \ep_0 +   L_*(J+1)  +r_0^{10} \big(  \Rk_J+\Sk_J  \big).
 \eea
\end{enumerate}
\end{corollary}

\begin{proof}
We start with the first estimate. In view of Theorem \ref{prop:rpweightedestimatesiterationassupmtionThM8}, and Propositions \ref{prop:controlGaSistartiterationassupmtionThM8}, \ref{prop:controlGaextiterationassupmtionThM8} and \ref{prop:improvementoftheiterationassupmtionThM8-Mtop}, we have
\beaa
 \Rk_{J+1}+\Skstar_{J+1} \les r_0^{-\dt}\Big( \Sktop^{\geq r_0}_{J+1} +\Skext_{J+1}\Big)+r_0^{10} \big(  \Rk_J+\Sk_J  \big)  +\ep_0
 \eeaa 
 and
 \beaa
 \Sktop^{\geq r_0}_{J+1}+\Skext_{J+1}&\les&  \Rk_{J+1}+\Skstar_{J+1} +\ep_0+L_*(J+1)
 \eeaa 
 where the constant in $\les$ is independent of $r_0$ in both inequalities. Plugging the second inequality in the first, we infer
\beaa
 \Rk_{J+1}+\Skstar_{J+1} &\les& r_0^{-\dt}\Big( \Rk_{J+1}+\Skstar_{J+1}\Big)+\ep_0 +  r_0^{-\dt}  L_*(J+1)  +r_0^{10} \big(  \Rk_J+\Sk_J  \big)
 \eeaa 
where the constant in $\les$ is independent of $r_0$. For $r_0$ large enough, we may absorb the first term on the RHS, and we obtain 
 \beaa
\Rk_{J+1}+\Skstar_{J+1} &\les \ep_0 +  r_0^{-\dt}  L_*(J+1)  +r_0^{10} \big(  \Rk_J+\Sk_J  \big),
 \eeaa
where the constant in $\les$ is independent of $r_0$, which is the first desired estimate. 

Next, we focus on the second estimate of the corollary. In view of Propositions \ref{prop:controlGaSistartiterationassupmtionThM8}, \ref{prop:controlGaextiterationassupmtionThM8}, \ref{prop:controlGaintiterationassupmtionThM8} and \ref{prop:improvementoftheiterationassupmtionThM8-Mtop}, we have
\beaa
  \Sk_{J+1} &\les& \Rk_{J+1}+\ep_0 +   L_*(J+1)  +r_0^{10} \big(  \Rk_J+\Sk_J  \big).
 \eeaa
Plugging the first estimate of the corollary proved above to control the term $\Rk_{J+1}$ in the RHS, we infer 
\beaa
  \Sk_{J+1} &\les& \ep_0 +   L_*(J+1)  +r_0^{10} \big(  \Rk_J+\Sk_J  \big)
 \eeaa
 which is the first desired estimate. This concludes the proof of Corollary \ref{Corr:improvementoftheiterationassupmtionThM8-Mtop}.
\end{proof}

%%%%%%%%%%%%%%%%%%%%%%%%%%%%%%%%%

\subsection{End of the proof of  the Main PT-Theorem}

%%%%%%%%%%%%%%%%%%%%%%%%%%%%%%%%%

{\bf Step 1.} As mentioned in Remark \ref{rmk:theiterationassumptionholdsforklsmallminus2:chap9},  the estimate \eqref{eq:initializationofiterationassumptioninproofThmM8} trivially implies the iteration assumption  \eqref{eq:iterationassumptiondiscussionThM8:bis} with $J=k_{small}-1$.  

 {\bf Step 2.}  We assume that the iteration assumption \eqref{eq:iterationassumptiondiscussionThM8:bis} holds for any fixed $J$ such that  $k_{small}-1\leq J\leq k_{large}+5 $, i.e. 
\beaa
\Sk_{J} +\Rk_{J}  &\les& \ep_0+L_*(J).
\eeaa 
Plugging the iteration assumption, for $J$ in the range  $k_{small}-1\leq J\leq k_{large}+5 $, in the estimates of 
Corollary \ref{Corr:improvementoftheiterationassupmtionThM8-Mtop}, we deduce
\beaa
\Rk_{J+1}+\Sk_{J+1} &\les& \ep_0 +    L_*(J+1)  +r_0^{10}(\ep_0+L_*(J)),
\eeaa
and hence
\beaa
\Rk_{J+1}+\Sk_{J+1} &\les& \ep_0 +    L_*(J+1),
\eeaa
 where the constant in $\les$ is depends on $r_0$. Thus, the iteration assumption \eqref{eq:iterationassumptiondiscussionThM8:bis} holds for $J+1$. This implies that 
 the iteration assumption \eqref{eq:iterationassumptiondiscussionThM8:bis} holds true for any 
  $J$ in the range  $k_{small}-1\leq J\leq k_{large}+6$. In particular, we have for $J=k_{large}+6$
\bea\lab{eq:iterationassumptionhodlstrueforkJequaitonklargeplus6:chap9}
\Rk_{k_{large}+6}+\Sk_{k_{large}+6} &\les& \ep_0 +    L_*(k_{large}+6).
\eea  
 
  {\bf Step 3.} In order to control the term $L_*(k_{large}+6)$ on the RHS of \eqref{eq:iterationassumptionhodlstrueforkJequaitonklargeplus6:chap9}, we rely on the following interpolation lemma.
  \begin{lemma}\lab{lemma:interpolationlemmatocontrolLstarofJiniteration:chap9}
 For any $k_{small}-1\leq k\leq k_{large}+7$, we have
 \bea
 L_*(k) &\les& (L_*(k_{small}-1))^{\frac{k_{large}+7-k}{k_{large}+7-(k_{small}-1)}} (L_*(k_{large}+7))^{\frac{k-(k_{small}-1)}{k_{large}+7-(k_{small}-1)}}.
 \eea
 \end{lemma}
 
 \begin{proof}
 We prove the lemma by iteration. Assume for $k_{small}-1< p\leq k_{large}+6$ the following iteration assumption 
 \beaa
 L_*(k) &\les& (L_*(k_{small}-1))^{\frac{p-k}{p-(k_{small}-1)}} (L_*(p))^{\frac{k-(k_{small}-1)}{p-(k_{small}-1)}},\quad k_{small}-1\leq k\leq p.
 \eeaa 
The iteration assumption holds trivially in the case $p=k_{small}$. 

Next, assume that the iteration assumption holds for some $p$ in the range $k_{small}-1< p\leq k_{large}+6$. We now consider whether the iteration assumption holds for $p+1$. Recall from \eqref{eq:defintionofLstarofk:chap9} that 
\beaa
 L^2_*(p) =  \int_{\{u=u_*'\}}    \big|\Rc\big|_{w, p}^2 +\int_{\Si_*\cap\{u=u_*'\}} \big|\Gac\big|_{w, p}^2
 \eeaa
In view of the definition \eqref{eq:defintionofthenormsRcGacweighted}, and noting that $\big|\Rc\big|_{w, k}$ only involves derivatives tangential to $\{u=u_*'\}$, while $ \big|\Gac\big|_{w, k}$ contains only derivatives tangential to $\Si_*\cap\{u=u_*'\}$, we may integrate by parts once, which yields
\beaa
 L^2_*(p) &\les& L_*(p-1)L_*(p+1).
 \eeaa
 Also, applying the iteration assumption with $k=p-1$, we have 
 \beaa
 L_*(p-1) &\les& (L_*(k_{small}-1))^{\frac{1}{p-(k_{small}-1)}} (L_*(p))^{\frac{p-1-(k_{small}-1)}{p-(k_{small}-1)}}\\
\eeaa
 and hence
\beaa
 L^2_*(p) &\les& (L_*(k_{small}-1))^{\frac{1}{p-(k_{small}-1)}} (L_*(p))^{\frac{p-1-(k_{small}-1)}{p-(k_{small}-1)}}L_*(p+1).
 \eeaa 
 or
\beaa
 L_*(p) &\les& (L_*(k_{small}-1))^{\frac{1}{p+1-(k_{small}-1)}} (L_*(p+1))^{\frac{p-(k_{small}-1)}{p+1-(k_{small}-1)}}.
 \eeaa 
 Plugging in the iteration assumption, we infer, for any $k_{small}-1\leq k\leq p$, 
 \beaa
 L_*(k) &\les& (L_*(k_{small}-1))^{\frac{p-k}{p-(k_{small}-1)}} (L_*(p))^{\frac{k-(k_{small}-1)}{p-(k_{small}-1)}}\\
 &\les&  (L_*(k_{small}-1))^{\frac{p+1-k}{p+1-(k_{small}-1)}} (L_*(p+1))^{\frac{k-(k_{small}-1)}{p+1-(k_{small}-1)}}\
 \eeaa 
 which, together with the fact that the case $k=p+1$ trivially holds, implies the iteration assumption for $p$ replaced by $p+1$. We deduce that the iteration assumption holds true for any $k_{small}-1< p\leq k_{large}+7$. In particular, for $p=k_{large}+7$, we infer that 
 \beaa
 L_*(k) &\les& (L_*(k_{small}-1))^{\frac{k_{large}+7-k}{k_{large}+7-(k_{small}-1)}} (L_*(k_{large}+7))^{\frac{k-(k_{small}-1)}{k_{large}+7-(k_{small}-1)}}
 \eeaa
for any $k_{small}-1\leq k\leq k_{large}+7$  as desired. This concludes the proof of Lemma \ref{lemma:interpolationlemmatocontrolLstarofJiniteration:chap9}.
 \end{proof}
 
  {\bf Step 4.} Applying  Lemma \ref{lemma:interpolationlemmatocontrolLstarofJiniteration:chap9} with $k=k_{large}+6$, we have
   \beaa
 L_*(k_{large}+6) &\les& (L_*(k_{small}-1))^{\frac{1}{k_{large}+7-(k_{small}-1)}} (L_*(k_{large}+7))^{1-\frac{1}{k_{large}+7-(k_{small}-1)}}.
 \eeaa
 Also, in view of \eqref{eq:initializationofiterationassumptioninproofThmM8}, we have
  \beaa
L_*(k_{small}-1) &\les& \ep_0 
 \eeaa
 and hence
   \beaa
 L_*(k_{large}+6) &\les& \ep_0^{\frac{1}{k_{large}+7-(k_{small}-1)}} (L_*(k_{large}+7))^{1-\frac{1}{k_{large}+7-(k_{small}-1)}}.
 \eeaa 
Plugging in  \eqref{eq:iterationassumptionhodlstrueforkJequaitonklargeplus6:chap9}, we deduce  
 \beaa
\Rk_{k_{large}+6}+\Sk_{k_{large}+6} &\les& \ep_0 +    L_*(k_{large}+6)\\
&\les& \ep_0+\ep_0^{\frac{1}{k_{large}+7-(k_{small}-1)}} (L_*(k_{large}+7))^{1-\frac{1}{k_{large}+7-(k_{small}-1)}}.
\eeaa  

Next, we use the first estimate of Corollary \ref{Corr:improvementoftheiterationassupmtionThM8-Mtop} with $J=k_{large}+6$, i.e.  
 \beaa
\Rk_{k_{large}+7}+\Skstar_{k_{large}+7} &\les& \ep_0 +  r_0^{-\dt}  L_*(k_{large}+7)  +r_0^{10} \big(  \Rk_{k_{large}+6}+\Sk_{k_{large}+6}   \big),
 \eeaa
 where the constant in $\les$ is independent of $r_0$. Plugging the above bound, we infer
 \beaa
&&\Rk_{k_{large}+7}+\Skstar_{k_{large}+7}\\ 
&\les&   r_0^{-\dt}  L_*(k_{large}+7)  +C_{r_0}\Big(\ep_0+\ep_0^{\frac{1}{k_{large}+7-(k_{small}-1)}} (L_*(k_{large}+7))^{1-\frac{1}{k_{large}+7-(k_{small}-1)}}\Big)
 \eeaa 
 where the constant in $\les$ is independent of $r_0$, and where $C_{r_0}$ depends on $r_0$. 

We now estimate $L_*(k_{large}+7)$. Recall from \eqref{eq:defintionofLstarofk:chap9} that 
\beaa
 L^2_*(k_{large}+7) =  \int_{\{u=u_*'\}}    \big|\Rc\big|_{w, k_{large}+7}^2 +\int_{\Si_*\cap\{u=u_*'\}} \big|\Gac\big|_{w, k_{large}+7}^2.
\eeaa 
Also, recall from \eqref{eq:choicofuprimestarbyLebesguepointarguement} that 
\beaa
\nn &&\int_{\{u=u_*'\}}    \big|\Rc\big|_{w, k_{large}+7}^2 +\int_{\Si_*\cap\{u=u_*'\}} \big|\Gac\big|_{w, k_{large}+7}^2\\
 &=& \inf_{u*-2\leq u_1\leq u*-1}\left(\int_{\{u=u_1\}}    \big|\Rc\big|_{w, k_{large}+7}^2 +\int_{\Si_*\cap\{u=u_1\}} \big|\Gac\big|_{w, k_{large}+7}^2\right).
\eeaa
We infer
\beaa
L^2_*(k_{large}+7) &\les& \int_{u_1=u_*-2}^{u_*-1}\left(\int_{\{u=u_1\}}    \big|\Rc\big|_{w, k_{large}+7}^2 +\int_{\Si_*\cap\{u=u_1\}} \big|\Gac\big|_{w, k_{large}+7}^2\right)du_1.
\eeaa
In view of the definition of $|\Rc|_{w, k}$ and $|\Gac|_{w, k}$ in \eqref{eq:defintionofthenormsRcGacweighted:bis}, and the definition of $\Rk_k$ and $\Skstar_k$ in section \ref{sec:mainnormsPTframe:chap9}, we deduce 
\beaa
L_*(k_{large}+7) &\les& \Rk_{k_{large}+7}+\Skstar_{k_{large}+7}
\eeaa
where the constant in $\les$ is independent of $r_0$. Plugging in the above, we obtain 
 \beaa
&&\Rk_{k_{large}+7}+\Skstar_{k_{large}+7}\\ 
&\les&   r_0^{-\dt}\Big(\Rk_{k_{large}+7}+\Skstar_{k_{large}+7}\Big) \\
&& +C_{r_0}\Big(\ep_0+\ep_0^{\frac{1}{k_{large}+7-(k_{small}-1)}} \Big(\Rk_{k_{large}+7}+\Skstar_{k_{large}+7}\Big)^{1-\frac{1}{k_{large}+7-(k_{small}-1)}}\Big)
 \eeaa 
where the constant in $\les$ is independent of $r_0$, and where $C_{r_0}$ depends on $r_0$. Taking $r_0$ large enough to absorb the first term on the RHS, and then $\ep_0$ small enough to absorb the last term on the RHS, we infer
 \beaa
\Rk_{k_{large}+7}+\Skstar_{k_{large}+7} &\les&   \ep_0.
 \eeaa 
 Together with the second estimate of Corollary \ref{Corr:improvementoftheiterationassupmtionThM8-Mtop} with $J=k_{large}+6$, we deduce 
\beaa
\Rk_{k_{large}+7}+\Sk_{k_{large}+7} &\les&   \ep_0
 \eeaa 
as desired. This concludes the proof of Theorem \ref{theorem:Main-PT}.

%%%%%%%%%%%%%%%%%%%%%%%%%%%%%%%

\section{Proof of  Theorem \ref{Theorem:TheoremM0-PT}}
\lab{sec:proofof:Theorem:TheoremM0-PT}

%%%%%%%%%%%%%%%%%%%%%%%%%%%%%%%

 As in the proof of Theorem M0, see section \ref{sec:proofofTheoremM0},  we divide the proof in  24 steps, denoted here by primes, i.e. 1'--24', which can thus be compared with  Steps 1--24 of section \ref{sec:proofofTheoremM0}.    
 
 \begin{remark}
 Steps  1'--19' differ little from Steps 1--19. The main difference occurs 
 with Steps 20'--24' where  the properties of the PT frames will  be important.  More generally, as the conclusions of Theorem M0 hold for $k\le k_{large}-2$ while the ones of Theorem \ref{Theorem:TheoremM0-PT} hold for $k\le k_{large}+7$, only estimates involving top regularity will differ. 
 \end{remark}

{\bf Steps  1'--7'.} Steps  1'--7' are identical  to Steps 1--7. More precisely, we propagate from $S_*$             along $\Si_*$, relative to the integrable  frame of $\Si_*$, the $\ell=1$ modes of $\div\b$, $\curl\b$, $\rhoc$ and $\kabc$, use the GCM assumptions on $S_*$,  and  arrive at the estimate  \eqref{eq:localbootassell1modedivbeta:proofThmM0:improved} of Proposition \ref{Proposition:Step1-7.ThmM0}, i.e.
\bea
\lab{eq:localbootassell1modedivbeta:proofThmM0:improved-M8}
\bsplit
&\sup_{\Si_*}\left(r^5|(\div\b)_{\ell=1}|+r^5|(\curl\b)_{\ell=1,\pm}|+r^5\left|(\curl\b)_{\ell=1,0}-\frac{2am}{r^5}\right|\right)\les \ep_0,\\
&\sup_{\Si_*}\Big(r^3 |(\rhoc)_{\ell=1}|+r^2 |(\kabc)_{\ell=1}|\Big) \les \ep_0.
\end{split}
\eea

{\bf Steps  8'--16'.} Next, as in   Steps 8--16,   $\CC_1'$ denotes the portion of the past directed outgoing, geodesic,  null cone initialized on the sphere $S_1'=S_1$ of $\Si_*$ and restricted to $r'\geq \de_*\ep_0^{-1}$.   Recall also that $\rt$ denotes the volume radius for the outgoing geodesic foliation  of $\Lextt$ constructed in section \ref{section:geodesicfoliation8}, while $r'$ is the are radius of the spheres $S'\subset \CC'_1$. In fact  all quantities associated to  the outgoing geodesic foliation of $\CC_1'$  are denoted by primes, while  the   quantities associated to the outgoing geodesic  foliation of $\Lextt$ are  denoted by tildes.

As in the proof of Theorem M0, we rely on  the estimates \eqref{eq:localbootassell1modedivbeta:proofThmM0:improved},  the GCM conditions on $S'_1=S_1$
and the following  bootstrap assumptions.

{\bf   Local Bootstrap Assumptions:} 
\begin{enumerate}
\item   Along $\CC'_1$, we have
\bea\lab{eq:localbootstrapssutionfortheproofofTheoremM0'}
\sup_{S'\subset \CC'_1}\Big(\|f\|_{\hk_4(S')}+(r')^{-1}\|(\fb, \log\la)\|_{\hk_4(S')}\Big) \leq \ep.
\eea

\item  On $S'_1=S_1$, we assume
\bea\lab{eq:iterativeapssutiononS1fortheproofofTheoremM0'}
\|f\|_{\hk_{k_{large}+7}(S'_1)}+r^{-1}\|(\fb, \log(\la))\|_{\hk_{k_{large}+7}(S'_1)} &\leq& \ep.
\eea

\item   In the case $a_0\neq 0$, we make the following assumption\footnote{Recall that $\dk_*$  refers to the properly normalized  tangential derivatives  along $\Si_*$.}, on $S'_1=S_1$,  on the difference between the basis of $\ell=1$ modes $\Jp$ of $\Si_*$, and the basis of $\ell=1$ modes $\Jt^{(p)}$ of $\Lextt$
\bea\lab{eq:bootstrapassumptiondifferenceell=1basisofCC1andLextprime'} 
\max_{p=0,+,-}\left\|\dk_*^{\leq k }(\Jp-\Jt^{(p)})\right\|_{L^2(S'_1)} &\leq& \ep, \qquad \mbox{for all} \quad k\le k_{large}+7.
\eea
\end{enumerate}

\begin{remark}\lab{rmk:comparisionnumberofderivativesinThM0abdThmM0'Steps8to16}
Note that \eqref{eq:iterativeapssutiononS1fortheproofofTheoremM0'} and  \eqref{eq:bootstrapassumptiondifferenceell=1basisofCC1andLextprime'} hold for $k\leq k_{large}+7$, while their analogs in Theorem M0 hold only for $k\leq k_{large}$. 
\end{remark}

The proof of Steps 8'--16' is then completely analogous to the one of Steps 8--16, with $k_{large}$ being replaced by $k_{large}+7$ in view of Remark \ref{rmk:comparisionnumberofderivativesinThM0abdThmM0'Steps8to16}. In particular, we obtain the following analog of \eqref{eq:finalesitateforffbchecklaonS1inproofofThmM0:reallyfinal}
\bea\lab{eq:finalesitateforffbchecklaonS1inproofofThmM0:reallyfinal'}
\bsplit
\sup_{k\le k_{ large}+8} \Big( \|{\dk}^k f \|_{L^2(S_{1}')} +r^{-1} \|{\dk}^{k} (\fb, \log \la) \|_{L^2(S_{1}')} +\|\dk^{\le k-1}\nab_\nu' (\fb, \log \la) \|_{L^2(S_{1}')}\Big) &\les \ep_0,\\
\sup_{S_{1}'}\left(\left|\frac{m}{m_0}-1\right|+\left|\frac{r'}{r}-1\right| \right)  &\les \ep_0,
\end{split}
\eea 
which improves the bootstrap assumption \eqref{eq:iterativeapssutiononS1fortheproofofTheoremM0'}. Also, we obtain the following analog of \eqref{eq:controlofainthecasea0equal0ThM0}
\bea\lab{eq:controlofainthecasea0equal0ThM0'}
|a| &\les& \ep_0\quad\textrm{if}\quad a_0=0,
\eea
and the following analog of \eqref{eq:controlofaminusa0andJminusJprimeforproofThM0}
\bea\lab{eq:controlofaminusa0andJminusJprimeforproofThM0'}
|a-a_0|+\max_{p=0,+,-}\sup_{S_1'}\left|{J}^{(p)}-{J'}^{(p)}\right| &\les& \ep_0\quad\textrm{if}\quad a_0\neq 0.
\eea

{\bf Step 17'.} The main difference with Step 17 is that we now rely on $L^2$ norms rather than sup norms. In particular, we rely on
\bea\lab{eq:bootstrapassumptionsonSigmastarforThmM0primeofThmM8'}
\|\dk_*^{\leq k}\Ga_b'\|_{L^2(\Si_*)} &\leq& \ep, \quad k\leq k_{large}+7,
\eea 
which follows immediately from the bootstrap assumptions {\bf  BA-PT}, see \eqref{eq:mainbootassforchapte9}.

As  in Step 17 of the proof of Theorem M0, using   Lemmas  \ref{Lemma:unifomfactor-ThmM0}, \ref{lemma:computationforfirstorderderivativesJpell=1basis} and \ref{lemma:computationforfirstorderderivativesf0fpfm} and simple elliptic estimates, we derive  first the estimates, for $k\le k_{large}+ 6$,
\beaa
&&\|\phi\|_{\hk_{k+2}(S_*)}\les \|\Ga_b'\|_{\hk_k(S_*)}, \\
&&\max_{p=0,+,-}\|r'\widecheck{\nab'{J'}^{(p)}}\|_{\hk_{k+2}(S_*)}\les \|\phi\|_{\hk_{k+2}(S_*)},\\
&& \|r'\nab f_0' - {J'}^{(0)}\in\|_{\hk_{k+1}(S_*)}+\|r'\nab f_{\pm}' + {J'}^{(\pm)}\de\|_{\hk_{k+1}(S_*)}\les \|\phi\|_{\hk_{k+2}(S_*)}.
\eeaa
Now, in view of   \eqref{eq:bootstrapassumptionsonSigmastarforThmM0primeofThmM8'} and the trace theorem, we have, for $k\le k_{large}+ 6$,
\beaa
\|\Ga_b'\|_{\hk_k(S_*)}\les \|\dk_*^{\leq k+1}\Ga_b'\|_{L^2(\Si_*)}\les\ep.
\eeaa
Hence, for  $k+1\le k_{large}+7$, 
\beaa
&&\|\phi\|_{\hk_{k+2}(S_*)}\les \ep, \\
&&\max_{p=0,+,-}\|r'\widecheck{\nab'{J'}^{(p)}}\in\|_{\hk_{k+2}(S_*)}\les\ep,\\
&&\|r'\nab f_0' - {J'}^{(0)}\in\|_{\hk_{k+1}(S_*)}+\|r'\nab f_{\pm}' + {J'}^{(\pm)}\de\|_{\hk_{k+1}(S_*)}\les \ep.
\eeaa
Next, on $\Si_*$, we have by Lemma \ref{lemma:transporteqnuofnabofJpf0fpfm},
\beaa
&&\nab_\nu\Big[r'\nab f_0' - {J'}^{(0)}\in\Big] =  \Ga_b' \c \dkb^{\leq 1}f_0',\qquad \nab_\nu\Big[r'\nab f_{\pm}' + {J'}^{(\pm)}\de\Big] =  \Ga_b' \c \dkb^{\leq 1}f_{\pm}',\\
&& \nab_\nu\Big[r'\widecheck{\nab'{J'}^{(p)}}\Big]=\Ga_b' \c \dkb^{\leq 1}{J'}^{(p)}, \quad p=0,+,-.
\eeaa
Integrating  from $S_*$,  we  obtain for $k\le k_{large}+7$
\beaa
\|r'\nab f_0' - {J'}^{(0)}\in\|_{L^\infty_u\hk_k(S')}&\les& \|r'\nab f_0' - {J'}^{(0)}\in\|_{\hk_k(S_*)}+\sqrt{u_*}\|\dk_*^k\Ga_b'\|_{L^2(\Si_*)},\\
\|r'\nab f_{\pm}' + {J'}^{(\pm)}\de\|_{L^\infty_u\hk_k(S')}&\les& \|r'\nab f_{\pm}' + {J'}^{(\pm)}\de\|_{\hk_k(S_*)}+\sqrt{u_*}\|\dk_*^k\Ga_b'\|_{L^2(\Si_*)},\\
\|r'\widecheck{\nab'{J'}^{(p)}}\|_{L^\infty_u\hk_k(S')} &\les& \|r'\widecheck{\nab'{J'}^{(p)}}\|_{\hk_k(S_*)}+\sqrt{u_*}\|\dk_*^k\Ga_b'\|_{L^2(\Si_*)}, \quad p=0,+,-.
\eeaa
Together with the above control on $S_*$, the control \eqref{eq:bootstrapassumptionsonSigmastarforThmM0primeofThmM8'} for  $\Ga_b'$, and the dominance condition for $r$ on $\Si_*$, we deduce, for all $k\le k_{large}+7,$
\beaa
\|r'\nab' f_0' - {J'}^{(0)}\in\|_{L^\infty_u\hk_{k}(S')}+\|r'\nab f_{\pm}' + {J'}^{(\pm)}\de\|_{L^\infty_u\hk_{k}(S')}\\ 
+ \max_{p=0,+,-}\|r'\widecheck{\nab'{J'}^{(p)}}\|_{L^\infty_u\hk_{k}(S')} &\les& \sqrt{u_*}\ep\les r\ep_0.
\eeaa

In view of the above, using the definition      of $\widecheck{\nab'{J'}^{(0)}}$, see Definition \ref{def:renormalizationforJpbasisell=1modes},  we infer on $S_1'$, for all $k\le k_{large}+7$,
\beaa
\nn\left\|\nab' f_0' - \frac{{J'}^{(0)}}{r'}\in\right\|_{\hk_{k}(S_1')} + \left\|\nab' f_{\pm}' + \frac{{J'}^{(p)}}{r'}\de\right\|_{\hk_{k}(S_1')}&\les & \ep_0,\\
+ \left\|\nab'{J'}^{(0)}+\frac{1}{r'}\dual f_0'\right\|_{\hk_{k}(S_1')}  + \left\|\nab'{J'}^{(\pm)}-\frac{1}{r'}f_{\pm}'\right\|_{\hk_{k}(S_1')} &\les& \ep_0.
\eeaa
On the other hand, from the control of $\Lextt$, the change of frame formula for $\nab'$, the control \eqref{eq:finalesitateforffbchecklaonS1inproofofThmM0:reallyfinal'} for the change of frame coefficients $(f, \fb, \la)$ from $\Lextt$ to $\Si_*$, and the control of $r'-r$ on $S_1'$ also provided by \eqref{eq:finalesitateforffbchecklaonS1inproofofThmM0:reallyfinal'}, we have, for all $k\le k_{large}+7$,
\beaa
\bsplit
\left\|\nab' f_0 - \frac{{J}^{(0)}}{r'}\in\right\|_{\hk_{k}(S_1')} + \left\|\nab' f_{\pm} + \frac{{J}^{(p)}}{r'}\de\right\|_{\hk_{k}(S_1')} &\les \ep_0,\\
 \left\|\nab'{J}^{(0)}+\frac{1}{r'}\dual f_0\right\|_{\hk_{k}(S_1')}  + \left\|\nab'{J}^{(\pm)}-\frac{1}{r'}f_{\pm}\right\|_{\hk_{k}(S_1')} &\les \ep_0,
 \end{split}
\eeaa
and hence,  for all $k\le k_{large}+7$,
\beaa
\left\|\nab' (f_0'-f_0) - \frac{{J'}^{(0)}-J^{(0)}}{r'}\in\right\|_{\hk_{k}(S_1')} +\left\|\nab' (f_{\pm} '-f_{\pm}) + \frac{{J'}^{(\pm)}-J^{(\pm)}}{r'}\de\right\|_{\hk_{k}(S_1')}&\les&\ep_0, \\
\left\|\nab'({J'}^{(0)}-J^{(0)})+\frac{1}{r'}\dual (f_0'-f_0)\right\|_{\hk_{k}(S_1')}
+ \left\|\nab'({J'}^{(\pm)}-J^{(\pm)})-\frac{1}{r'}\dual (f_{\pm}'-f_{\pm})\right\|_{\hk_{k}(S_1')} &\les& \ep_0.
\eeaa
We deduce,  for all $k\le k_{large}+7$,
\beaa
r^{-1}\left\|f_0'-f_0\right\|_{\hk_{k+1}(S_1')} + r^{-1}\left\|{J'}^{(0)}-J^{(0)}\right\|_{\hk_{k+1}(S_1')} \les \ep_0+\|{J'}^{(0)}-J^{(0)}\|_{L^2(S_1')}
\eeaa
and
\beaa
r^{-1}\left\|f_{\pm}'-f_{\pm}\right\|_{\hk_{k+1}(S_1')} + r^{-1}\left\|{J'}^{(\pm)}-J^{(\pm)}\right\|_{\hk_{k+1}(S_1')} \les \ep_0+r^{-1}\|{J'}^{(\pm)}-J^{(\pm)}\|_{L^2(S_1')}.
\eeaa
Together with \eqref{eq:controlofaminusa0andJminusJprimeforproofThM0'}, we  obtain in the case $a_0\neq 0$, for all $k\le k_{large}+7$,
\beaa
\max_{p=0,+,-}\left(r^{-1}\left\|f_p'-f_p\right\|_{\hk_{k+1}(S_1')} + r^{-1}\left\|{J'}^{(p)}-J^{(p)}\right\|_{\hk_{k+1}(S_1')}\right) \les \ep_0.
\eeaa
In view of the fact that $\nab_\nu f_p'=0$ for $p=0,+,-$ on $\Si_*$, and in view of the control of $\Lextt$, and hence of $f_p$ for $p=0,+,-$, we deduce for $a_0\neq 0$, for all $k\le k_{large}+8$, 
\bea
\max_{p=0,+,-}\left(r^{-1}\left\|\dk_*^{\leq k}(f_p'-f_p)\right\|_{L^2(S_1')} + r^{-1}\left\|\dk_*^{\leq k}({J'}^{(p)}-J^{(p)})\right\|_{L^2(S_1')}\right) \les \ep_0.
\eea
In particular, the above  improves\footnote{Both with respect to  the smallness constant  and to regularity.}  the bootstrap assumption \eqref{eq:bootstrapassumptiondifferenceell=1basisofCC1andLextprime'} on ${J'}^{(p)}-J^{(p)}$ for $p=0,+,-$.

\begin{remark}
Recall that Step 18 in the proof of Theorem M0 focuses on the control of $\nu(r')$ and $b_*$ which are involved in the change  from the  frame of $\Si_*$ to the PG frame of $\Mext$. Since the change  between the frame of $\Si_*$ and the PT frame of $\Mext$ involves neither  $\nu(r')$ nor $b_*$, we note that there is no need for Step 18'.
\end{remark}

\noindent{\bf Step 19'.} We consider the following change of frames coefficients:
\begin{itemize}
\item $(f, \fb, \la)$ are the change of frame coefficients from the outgoing geodesic frame of $\Lextt$ to the  frame of $\Si_*$. $(f, \fb, \la)$ satisfies according to  \eqref{eq:finalesitateforffbchecklaonS1inproofofThmM0:reallyfinal'}
\beaa
\sup_{k\le k_{ large}+8} \Big( \|{\dk}^k f \|_{L^2(S_{1}')} +r^{-1} \|{\dk}^{k} (\fb, \log \la) \|_{L^2(S_{1}')} +\|\dk^{\le k-1}\nab_\nu' (\fb, \log \la) \|_{L^2(S_{1}')}\Big) &\les \ep_0.
\eeaa

\item $(f', \fb', \la')$ are the change of  frame coefficients from the  outgoing PG frame of $\Lext$ to the outgoing geodesic frame of $\Lextt$. In view of 
 Proposition  \ref{proposition:geodesicfoliationLextt},   we have  
\beaa
\sup_{\Lextt}r\left|\dk^{\leq k_{large}+8}\left(f'+\frac{a_0}{r}f_0, \fb'+\frac{a_0\Upsilon}{r}f_0, \log\la'\right)\right| &\les& \ep_0.
\eeaa

\item $(f'', \fb'', \la'')$ are the change  of frame coefficients from the frame of $\Si_*$ to the outgoing PT frame of $\Mext$. We assume here that  $(f'', \fb'', \la'')$ satisfies by the initialization   of the PT structure  on $\Si_*$, see section \ref{sec:defintionofthePTstructuresinMM}, 
\beaa
\la''=1,\qquad f''=\frac{a}{r'}f_0', \qquad \fb''=\frac{a\Upsilon'}{r'}f_0'.
\eeaa
\end{itemize}

We now consider the change of frame   coefficients $(f''', \fb''', \la''')$  from the  outgoing  PG   frame  of $\Lext$ to the outgoing PT  frame of $\Mext$.  In view of:
\begin{itemize}
\item the above estimates for $(f, \fb, \la)$ and $(f', \fb', \la')$,

\item the above formula for $(f'', \fb'', \la'')$,

\item the control for $r-r'$ and $m-m_0$ given by \eqref{eq:finalesitateforffbchecklaonS1inproofofThmM0:reallyfinal'}, 

\item the control of $a$ in \eqref{eq:controlofainthecasea0equal0ThM0'} in the case $a_0=0$, 

\item the control for $a-a_0$ in \eqref{eq:controlofaminusa0andJminusJprimeforproofThM0'} and the control for $f_0'-f_0$ in Step 17' in the case $a_0\neq 0$,  
\end{itemize}
we infer the following estimates, for all $k\le k_{large}+8$, 
\beaa
\|\dk^{\leq k}f'''\|_{L^2(S_1')} +r^{-1}\|\dk^{\leq k}(\fb''', \log\la''')\|_{L^2(S_1')} &\les& \ep_0+\frac{1}{r}.
\eeaa
Together with the dominance condition for $r$ on $\Si_*$, we infer
\bea
\|\dk^{\leq  k}f'''\|_{L^2(S_1')} +r^{-1}\|\dk^{\leq k}(\fb''', \log\la''')\|_{L^2(S_1')} &\les& \ep_0, \quad\,\, k\le k_{large}+8.
\eea

\noindent{\bf Step 20'.}  In this step, $(e_1, e_2, e_3, e_4)$ denotes the outgoing PG frame of $\Lext$, and $(e_1', e_2', e_3', e_4')$ denotes the outgoing PT frame of $\Mext$. From now on, $(f, \fb, \la)$ denotes the change of frame   coefficients\footnote{Denoted  in the previous step by $(f''', \fb''', \la''')$.}    from the  outgoing PG frame of $\Lext$ to the outgoing PG  frame of $\Mext$.   In view of Step 19', we have on $S_1'$
\beaa
\|\dk^{\leq k}f\|_{L^2(S_1')} +r^{-1}\|\dk^{\leq k}(\fb, \log\la)\|_{L^2(S_1')} &\les& \ep_0, \quad k\le k_{large}+8.
\eeaa

Let
\beaa
F=f+i\dual f.
\eeaa
Since 
\beaa
\Xi'=0, \qquad \om'=0, \qquad \Xi=0, \qquad \om=0, \qquad \etab=-\ze,
\eeaa
we have, in view of Corollary \ref{cor:transportequationine4forchangeofframecoeffinformFFbandlamba:PTcase}, 
\beaa
\nab_{\la^{-1}e_4'}(\ov{q}F) &=& E_4(f, \Ga),\\
\la^{-1}\nab_4'(\log\la) &=& 2f\c\ze+E_2(f, \Ga).
\eeaa
We integrate the above transport equations for $F$ and $\la$   from $S_1'$.  In view of the control for $(f, \la)$ on $S_1'$ derived in Step 19', and in view of the assumptions on the initial data layer norm we infer,  for all $k\le k_{large}+8$, 
\bea
\lab{estimateforla-step21'}
\sup_{S'\subset\{u'=1\}}\Big(\|\dk'^{\leq k}f\|_{L^2(S')} +r^{-1}\|\dk'^{\leq k}\log\la\|_{L^2(S')}\Big) &\les& \ep_0,
\eea 
where\footnote{Note that $u'$ is associated both to the PT and to the PG structure of $\Mext$ in view of Lemma \ref{lemma:linkPGandTframeinMext}.} $u'$ denotes from now on the scalar function of the PT structure of $\Mext$. In particular, we have by construction $S_1'=\Si_*\cap\{u'=1\}$, $e_4'(u')=0$ and $\{u'=1\}\subset\Lext$.

\noindent{\bf Step 21'.} In this step, we estimate $r'-r$, ${J'}^{(0)}-J^{(0)}$ and $\Jk'-\Jk$, as well as $A'$, $\trXc'$ and $\Xh'$.  First, since ${J'}^{(0)}$ is propagated from $\Si_*$ by $e_4'({J'}^{(0)})=0$, and using the change of frame formula between the PG frame of $\Lext$ and the PT frame of $\Mext$, we infer
\beaa
e_4'({J'}^{(0)}-J^{(0)}) &=& -\la\left(f\c\nab+\frac{1}{4}|f|^2e_3\right)J^{(0)}.
\eeaa
Together with the control of $f$ and $\la$ of Step 20', and the control of $\Lext$, we infer,  for          $k \leq k_{large}+8$, 
\beaa
\sup_{S'\subset\{u'=1\}}r\|\dk^{\leq k}(e_4'({J'}^{(0)}-J^{(0)}))\|_{L^2(S')} &\les& \ep_0.
\eeaa
Integrating from $S_1'$ where ${J'}^{(0)}-J^{(0)}$ is under control in view of Step 17', we infer,  for    all    $k \leq k_{large}+8$, 
\bea
\sup_{S'\subset\{u'=1\}}r^{-1}\|\dk^{\leq k}({J'}^{(0)}-J^{(0)})\|_{L^2(S')} &\les& \ep_0.
\eea

Next, we control $\trXc'$. To this end, we also need to control $A'$ and $\Xh'$. First, note that the change of frame formula for $A'$, the control of the foliation of $\Lext$, the control of $f$ and $\la$ of Step 20', and the fact that the transformation formula for $A'$ does not depend on $\fb$ implies, since $k \leq k_{large}+8$,  
\bea
\sup_{S'\subset\{u'=1\}}r^{\frac{5}{2}+\dt}\|\dk^{\leq k}A'\|_{L^2(S')} &\les& \ep_0.
\eea
Then, using the control of $(f, \fb, \la)$ on $S_1'$ derived in Step 19', the control of the initial data layer, and the change of frame formulas between the PG frame of $\Lext$ and the PT frame of $\Mext$ on $S_1'$, we infer, for all $k\le k_{large}+7$,
\beaa
r\Big(\|\dk^{\leq k}\trXc'\|_{L^2(S_1')}+\|\dk^{\leq k}\Xh'\|_{L^2(S_1')}\Big) &\les& \ep_0.
\eeaa
Then, propagating Raychadhuri for $\trXc'$ and the null structure equation for $\nab_4'\Xh'$ from $S_1'$ where $\trXc'$ and $\Xh'$  are under control in view of the above estimate, we infer, using  also the above control of $A'$, 
\bea
\sup_{S'\subset\{u'=1\}}r\Big(\|\dk^{\leq k}\trXc'\|_{L^2(S')}+\|\dk^{\leq k}\Xh'\|_{L^2(S')}\Big) &\les& \ep_0, \qquad k\le k_{large}+7.
\eea

Next, we notice
\beaa
\tr X' -\la\tr X &=& \frac{2}{q'}-\frac{2}{q}+(\la-1)\tr X+\trXc' -\trXc
\eeaa
so that 
\beaa
q'-q &=& \frac{qq'}{2}\left(-\la\Big(\la^{-1}\tr X' -\tr X\Big)+(\la-1)\tr X+\trXc' -\trXc\right).
\eeaa
Together with the control of Step 20' for $f$ and $\la$, the above control for ${J'}^{(0)}-J^{(0)}$ and $\trXc'$, the control of $a-a_0$ in \eqref{eq:controlofainthecasea0equal0ThM0'} \eqref{eq:controlofaminusa0andJminusJprimeforproofThM0'}, the control of the foliation $\Lext$, and the fact that $q=r+ia_0J^{(0)}$ and $q'=r'+ia{J'}^{(0)}$ we infer,  for all $k\le k_{large}+7$.
\beaa
\sup_{S'\subset\{u'=1\}}r^{-1}\left\|\dk^{\leq k}\left(\frac{r'}{r}-1 +\frac{qq'}{2r}\Big(\la^{-1}\tr X' - \tr X\Big)\right)\right\|_{L^2(S')} &\les& \ep_0.
\eeaa
Moreover, from the change of frame formulas for $\trch$ and $\atrch$
we have, schematically, 
\beaa
\la^{-1}\tr X' -\tr X &=& {r'}^{-1}{\dk'}f+\Ga\c f+f\c\fb\c\Ga+f\c\fb\c\tr X'+\lot
\eeaa
Together with the control of Step 20' for $f$, the above control for  $\trXc'$ and the control of the foliation $\Lext$ we deduce,  for all $k\le k_{large}+7$,
\bea
\lab{eq:estimateform(r-r')-Step21'}
\sup_{S'\subset\{u'=1\}}r^{-1}\left\|\dk^{\leq k}\left(\frac{r'}{r}-1\right)\right\|_{L^2(S')} \les \ep_0+\ep_0\sup_{S'\subset\{u'=1\}}r^{-1}\left\|\dk^{\leq k}\fb\right\|_{L^2(S')}.
\eea

Next, recall the definition of $\Jk'$ on $\Si_*$, see Definition \ref{def:definitionofJkonMext},
\beaa
\Jk' &=& \frac{1}{|q'|}(f_0'+i\dual f_0')=\frac{1}{\sqrt{{r'}^2+a^2({J'}^{(0)})^2}}(f_0'+i\dual f_0').
\eeaa
Together with the control of $r-r'$ and $m-m_0$ given by \eqref{eq:finalesitateforffbchecklaonS1inproofofThmM0:reallyfinal'}, the control of $a-a_0$  in \eqref{eq:controlofainthecasea0equal0ThM0'} \eqref{eq:controlofaminusa0andJminusJprimeforproofThM0'},  the control for $f_0'-f_0$ and ${J'}^{(0)}-J^{(0)}$ in Step 17', and the control of $\Lext$, we infer
\beaa
\|\dk_*^{\leq k}(\Jk'-\Jk)\|_{L^2(S_1')} &\les& \ep_0+\frac{1}{r}. 
\eeaa
Together with the dominance condition for $r$ on $\Si_*$, we infer,  for all $k\le k_{large}+7$,
\beaa
\|\dk_*^{\leq k}(\Jk'-\Jk)\|_{L^2(S_1')} &\les& \ep_0. 
\eeaa
Together with the identity
\beaa
q'\Jk'-q\Jk &=& q'(\Jk'-\Jk)+(q'-q)\Jk\\
&=& q'(\Jk'-\Jk)+\Big(r'-r+i(a{J'}^{(0)}-a_0J^{(0)})\Big)\Jk,
\eeaa
the control of $r'-r$ given by \eqref{eq:finalesitateforffbchecklaonS1inproofofThmM0:reallyfinal'}, the control of ${J'}^{(0)}-J^{(0)}$ in Step 17' and the control of $a-a_0$  in \eqref{eq:controlofainthecasea0equal0ThM0'} \eqref{eq:controlofaminusa0andJminusJprimeforproofThM0'}, we obtain
\beaa
r^{-1}\|\dk_*^{\leq k} (q'\Jk'-q\Jk)\|_{L^2(S_1')} &\les& \ep_0. 
\eeaa
Also, recall that $\Jk'$ and $\Jk$ satisfy in $\{u'=1\}$
\beaa
\nab_4'\Jk' = -\frac{1}{q'}\Jk', \qquad \nab_4\Jk= -\frac{1}{q}\Jk,
\eeaa
and hence
\beaa
\nab_4'(q'\Jk')=0, \qquad \nab_4(q\Jk)=0.
\eeaa
We infer
\beaa
\nab_{\la^{-1}4'}(q'\Jk')=0, \qquad \nab_{\la^{-1}4'}(q\Jk)=\left(f\c\nab +\frac{1}{4}|f|^2e_3\right)(q\Jk).
\eeaa
Together with the control of $f$ and $\la$ of Step 20', (see estimate \eqref{estimateforla-step21'} ), and the control of $\Lext$, we obtain, for $k\leq k_{large}+8$,
\beaa
\sup_{S'\subset\{u'=1\}}r\|\dk^{\leq k}(\nab_4'(q'\Jk'-q\Jk))\|_{L^2(S')} &\les& \ep_0.
\eeaa
Integrating from $S_1'$ where $q'\Jk'-q\Jk$ is under control in view of the above we infer,   for all $k\le k_{large}+7$,
\beaa
\sup_{S'\subset\{u'=1\}}\|\dk^{\leq k}(q'\Jk'-q\Jk)\|_{L^2(S')} &\les& \ep_0.
\eeaa
Using again the above identity for $q'\Jk'-q\Jk$,  the above control of ${J'}^{(0)}-J^{(0)}$ and $r'-r$  and the control of $a-a_0$   in \eqref{eq:controlofainthecasea0equal0ThM0'} \eqref{eq:controlofaminusa0andJminusJprimeforproofThM0'}, we deduce  for all $k\le k_{large}+7$,
\bea
\lab{eq:EstimateforJk-Jk'.Step21'}
\sup_{S'\subset\{u'=1\}}\|\dk^{\leq k}(\Jk'-\Jk)\|_{L^2(S')} \les \ep_0+\ep_0\sup_{S'\subset\{u'=1\}}r^{-1}\left\|\dk^{\leq k}\fb\right\|_{L^2(S')}.
\eea

\noindent{\bf Step 22'.}  To estimate $\fb$ we make use of the last equation in Corollary \ref{cor:transportequationine4forchangeofframecoeffinformFFbandlamba:PTcase}, and the fact that $\Hb+\frac{a\ov{q}}{|q|^2}\Jk=-\widecheck{Z} $, see  Definition \ref{def:renormalizationofallnonsmallquantitiesinPGstructurebyKerrvalue}. For $\Fb=\fb+i\dual \fb$, we obtain
\beaa
\nab_{\la^{-1}e_4'}\underline{F} &=& -2\left(\frac{a\ov{q'}}{|q'|^2}\Jk'-\frac{a\ov{q}}{|q|^2}\Jk\right) +2\Zc -\frac{1}{2}\ov{\tr X}F  -F\c\chibh  +E_6(f, \fb, \Ga).
\eeaa
Together with the control of $a-a_0$  in \eqref{eq:controlofainthecasea0equal0ThM0'} \eqref{eq:controlofaminusa0andJminusJprimeforproofThM0'},  the control of ${J'}^{(0)}-J^{(0)}$, $r'-r$ and $\Jk'-\Jk$ of Step 21',  and the control on the 
  Ricci coefficients  of the PG frame of $\Lext$, we obtain for,  all $k\le k_{large}+7$, 
\beaa
\sup_{S'\subset\{u'=1\}}r\left\|\dk^{\leq k}\big( \la^{-1}\nab_{e_4'}\fb\big)\right\|_{L^2(S')} \les& \ep_0+\ep_0\sup_{S'\subset\{u'=1\}}r^{-1}\left\|\dk^{\leq k}\fb\right\|_{L^2(S')}.
\eeaa

Together with the above, and the control of $\la$ of Step 20', see estimate \eqref{estimateforla-step21'}, we infer
\beaa
\sup_{S'\subset\{u'=1\}}r\left\|\dk^{\leq  k}\nab_{e_4'}\fb\right\|_{L^2(S')} &\les& \ep_0+\ep_0\sup_{S'\subset\{u'=1\}}r^{-1}\left\|\dk^{\leq k}\fb\right\|_{L^2(S')}.
\eeaa
Integrating from $S_1'$ where $\fb$ is under control in view of Step 19', we infer
\beaa
\sup_{S'\subset\{u'=1\}}r^{-1}\left\|\dk^{\leq k}\fb\right\|_{L^2(S')} &\les& \ep_0+\ep_0\sup_{S'\subset\{u'=1\}}r^{-1}\left\|\dk^{\leq k}\fb\right\|_{L^2(S')}
\eeaa
and hence, for all $k\le k_{large}+7$,
\beaa
\sup_{S'\subset\{u'=1\}}r^{-1}\left\|\dk^{\leq k}\fb\right\|_{L^2(S')} &\les& \ep_0.
\eeaa
Together with the control of $f$ and $\la$ of Step 20', we have finally obtained,  for all $k\le k_{large}+7$,
\bea
\lab{eq:Step23'.1}
\sup_{S'\subset\{u'=1\}}\Big(\|\dk^{\leq k }f\|_{L^2(S')}+r^{-1}\|\dk^{\leq k }(\log(\la), \fb)\|_{L^2(S')}\Big) &\les& \ep_0.
\eea
Also, together with the estimates \eqref{eq:EstimateforJk-Jk'.Step21'}   and  \eqref{eq:estimateform(r-r')-Step21'} for $r'-r$ and $\Jk'-\Jk$   of Step 21', we obtain 
for all $k\le k_{large}+7$, 
\bea
\lab{eq:Step23'.2}
\sup_{S'\subset\{u'=1\}}\left(r^{-1}\left\|\dk^{\leq k}\left(\frac{r'}{r}-1\right)\right\|_{L^2(S')}+\|\dk^{\leq  k}(\Jk'-\Jk)\|_{L^2(S')}\right) &\les& \ep_0.
\eea

\noindent{\bf Step 23'.} Let  $(f', \fb', \la')$ denote the change of frame   coefficients   from the  ingoing PG  frame of $\Lint$ to the ingoing PT frame of $\Mint'$. From
\begin{itemize}
\item the estimates of Step 22' on $\{u'=1\}$, 

\item the fact that $\Mint'\cap\Mext=\{r'=r_0\}$,

\item the fact that $\{u'=1\}\cap\{\ub'=1\}$ is included in $\Lext\cap\Lint$,

\item the initialization of the frame of $\Mint'$ as an explicit renormalization of the frame of $\Mext$ on $\{r'=r_0\}$, 

\item the control in $\Lext\cap\Lint$ of the difference between the frame of $\Lint$ and an explicit renormalization of the frame of $\Lext$,
\end{itemize}
we easily infer, using also $\ub'=u'$ on $\{r=r_0\}$,  for all $k\le k_{large}+7$,
\beaa
\sup_{S'\subset\{r'=r_0\}\cap\{\ub'=1\}}\Big(\|\dk^{\leq k}(f', \log(\la'), \fb')\|_{L^2(S')}\Big) &\les& \ep_0.
\eeaa

Next, we proceed as in Step 20', exchanging the role of $e_3$ and $e_4$, and we propagate along $e_3$ the above estimate to $\{\ub'=1\}$ for $\fb'$ and $\la'$. We also propagate the control of Step 21' for ${J'}^{(0)}-J^{(0)}$ on $\{u'=1\}$, and hence on its boundary $\{r'=r_0\}$ to $\{\ub'=1\}$. Also one propagates the control of Step 22' for $r-r'$ on $\{u'=1\}$, and hence on its boundary $\{r'=r_0\}$ to $\{\ub'=1\}$ using the transport equation\footnote{Note that we could not have used this transport equation in $\Lext$ in view of the lack of decay in $r$ for $\la-1$. This is why we avoided this transport equation in Step 21' and used instead the control of $\trXc'$. On the other hand, $r$ is bounded in $\Lint$ so that one can simply rely on the transport equation for $e_3'(r'-r)$ in $\Lint$. This is in fact crucial: proceeding as in Step 21' would  only lead to the control $k_{large}+6$ derivatives of $r'-r$, and hence to a loss of at least one derivatives in Theorem \ref{Theorem:TheoremM0-PT}.}
\beaa
e_3'(r'-r) &=& 1-\la'\left(e_3+{\fb'}^ae_a+\frac{1}{4}|\fb'|^2e_4\right)r\\
&=& -(\la'-1)+\fb'\c\nab(r)+ \frac{1}{4}|\fb'|^2e_4(r).
\eeaa
Then, one propagates $\Jk'-\Jk$ similarly to Step 21'. Finally, we propagate $f'$ similarly to Step 22'. We finally obtain 
\bea
\sup_{S'\subset\{\ub'=1\}}\Big(\|\dk^{\leq k}(f', \fb', \log\la')\|_{L^2(S')}\Big) &\les& \ep_0, \qquad k\le k_{large}+7
\eea
and
\bea
\sup_{S'\subset\{\ub'=1\}}\left\|\dk^{\leq j+1}\left({J'}^{(0)}-J^{(0)}, r'-r, \Jk'-\Jk\right)\right\|_{L^2(S')} \les \ep_0, \qquad k\le k_{large}+7.
\eea

\noindent{\bf Step 24'.} Note that the desired estimates\footnote{Since they only require  a small number of derivatives their proof has in fact already been obtained in Theorem M0.}  for $m-m_0$ and $a-a_0$ have been obtained respectively in Step 13' and in\eqref{eq:controlofainthecasea0equal0ThM0'} \eqref{eq:controlofaminusa0andJminusJprimeforproofThM0'}. To conclude the proof of Theorem \ref{Theorem:TheoremM0-PT}, it remains to control $k\le k_{karge} +7$ derivatives, with suitable $r$-weights and $O(\ep_0)$ smallness constant, of $A$, $B$, $\Pc$, $\Bb$ and $\Ab$ in $\{u'=1\}\cup\{\ub'=1\}$. This follows from: 
\begin{itemize}
\item the control of $(f, \fb, \la)$ on $\{u'=1\}$ derived in Step 22', 

\item the control of $(f', \fb', \la')$ on $\{\ub'=1\}$ derived in Step 23',

\item the fact that $(f, \fb, \la)$ denote the change of frame   coefficients   from the outgoing PT frame of $\Lext$ to the outgoing PT   frame of $\Mext$, and the fact that $(f', \fb', \la')$ denote the change of frame   coefficients   from the ingoing PT frame of $\Lint$ to the ingoing PT frame of $\Mint'$,

\item the change of frame formulas for the curvature components,

\item in the particular case of the estimate for $\Pc$, the fact that 
\beaa
P'-P &=& -\frac{2m}{{q'}^3}+\frac{2m_0}{q^3}+\Pc'-\Pc,
\eeaa
together with the control of $m-m_0$ derived in Step 13',  the control of $a-a_0$ in \eqref{eq:controlofainthecasea0equal0ThM0'} \eqref{eq:controlofaminusa0andJminusJprimeforproofThM0'},  the control of ${J'}^{(0)}-J^{(0)}$ in Step 21' and 23', and the control of $r'-r$ in Steps 22' and 23',

\item the assumptions for the curvature components of the foliations of $\Lext$ and $\Lint$.
\end{itemize}

This concludes the proof of Theorem \ref{Theorem:TheoremM0-PT}.

%%%%%%%%%%%%%%%%%%%%%%%%%%%%%%%%%%%%%%%%%%%%%%%%%%

\section{Control of the PT-Ricci coefficients  on $\Si_*$}
\lab{sec:theoremM8recoverRiccionSistar}

%%%%%%%%%%%%%%%%%%%%%%%%%%%%%%%%%%%%%%%%%%%%%%%%%%

The goal of this section is to  provide the proof of Proposition \ref{prop:controlGaSistartiterationassupmtionThM8}.  For convenience, we restate the result below.
\begin{proposition}
\lab{Prop:controlGa-PTframe-Si_*}
The  Ricci and metric coefficients  of the outgoing PT frame of $\Mext$ verify  the following estimates on $\Si_*$, for all $k\le k_{large}+7$,
\bea
\Skstar_k\les \ep_0+\Rkstar_k.
\eea
\end{proposition}

Recall, see Definition \ref{def:definitionoff0fplusfminus}, the tangential 1-form $f_0$ on $\Si_*$ given by
\beaa
(f_0)_1=0, \quad (f_0)_2=\sin\th, \quad \textrm{on}\quad S_*, \qquad \nab_\nu f_0=0,
\eeaa
where, on $S_*$, we consider the orthonormal basis $(e_1, e_2)$ of $S_*$ given by \eqref{eq:specialorthonormalbasisofSstar}.  

Also, recall from section \ref{sec:defintionofthePTstructuresinMM} 
that we  initialize the PT frame of $\Mext$ from the integrable frame  on $\Si_*$, by  relying on the   change of frame formula 
with  the transition  coefficients 
\beaa
\la=1, \qquad f=\frac{a}{r}f_0, \qquad  \fb=\frac{a\Upsilon}{r}f_0.
\eeaa

In order to prove Proposition \ref{Prop:controlGa-PTframe-Si_*}, we proceed as follows:
\begin{enumerate}
\item  We first  derive analogous estimates for the integrable frame  of $\Si_*$ in section \ref{sec:controlofintegrableframeSigmastar:chap9}. 

\item  Then, we derive  estimates for  the 1-form $f_0$  on $\Si_*$   in section \ref{sec:controlofthe1formf0onSigramstar:cshap9}.

\item Finally, we conclude the proof of Proposition \ref{Prop:controlGa-PTframe-Si_*} in section \ref{sec:proofofProp:controlGa-PTframe-Si_*} by using  the transition coefficients 
\beaa
\la=1, \qquad f=\frac{a}{r}f_0, \qquad  \fb=\frac{a\Upsilon}{r}f_0,
\eeaa
the transformation formulas to  pass  from the integrable frame of $\Si_*$ to the PT frame of $\Mext$  on $\Si_*$, and the control of the integrable of $\Si_*$ and of the 1-form $f_0$. 
\end{enumerate}

%%%%%%%%%%%%%%%%%%%%%%%%%%%%%%%%%%%%%%%%%%%%%%%%%%

\subsection{Control of the  integrable frame of $\Si_*$}
\lab{sec:controlofintegrableframeSigmastar:chap9}

%%%%%%%%%%%%%%%%%%%%%%%%%%%%%%%%%%%%%%%%%%%%%%%%%%

Recall that on $\Si_*$,  the Ricci coefficients of the integrable frame of $\Si_*$   verify  the following GCM conditions
\beaa
&& \trchc=0, \qquad \trchbc=\sum_p\underline{C}_pJ^{(p)}, \qquad \muc =\sum_p M_pJ^{(p)}, \\
&& (\div\eta)_{\ell=1}=0, \qquad (\div\xib)_{\ell=1}=0,\qquad (b_*)_{|_{SP}}=1,
\eeaa
where $\underline{C}_p$, $M_p$, $p=0,+,-$, are functions of $r$ along $\Si_*$, and where $SP$ denotes the curve $\th=\pi$ of the south poles.

Additionally, the following holds on  $S_*$
\beaa
(\trchbc)_{\ell=1}=0, \qquad (\div\b)_{\ell=1}=0, \qquad (\curl\b)_{\ell=1, \pm}=0, \qquad (\curl\b)_{\ell=1,0}=\frac{2am}{r^5}. 
\eeaa
Also, remember that the integrable  frame of $\Si_*$ satisfies the following transversality conditions
\bea\lab{eq:transversalityconditionsonSigmastarcompatiblewithanoutgoingfoliation}
\xi=0, \qquad \om=0, \qquad \etab=-\ze, \qquad e_4(r)=1, \qquad e_4(u)=0,
\eea
where $r$ denotes the area radius of the $u$-foliation of $\Si_*$, and we have $u+r=c_{\Si_*}$ on $\Si_*$, where $c_{\Si_*}$ is a constant.

\begin{remark}
Note that the transversality conditions above are compatible with an outgoing geodesic foliation initialized on $\Si_*$. Note also that they must be specified to make sense of $\xib$ and $\eta$, and hence to make sense of the GCM conditions $(\div\eta)_{\ell=1}=0$ and $(\div\xib)_{\ell=1}=0$ on $\Si_*$. 
\end{remark}

To state the main result of this section we need to introduce  the following norms.
\begin{definition}
\lab{definition:Int-normson-GacSi_*} 
We define the   following Ricci  coefficients   norms on $\Si_*$ relative to the integrable frame of $\Si_*$
\bea
\bsplit
\Skstarr_k^2&:= \int_{\Si_*} r^2|\dk_*^{\le k+1}( \chih, \trchc, \ze, \trchbc)|^2 + |\dk^{\leq k+1}_*\chibh|+|\dkb^{\leq 1}\dk^{\leq k}_*(\eta, \ombc, \xib )|^2\\
&+\int_{\Si_*} r^{-2}|\dkb^{\leq 2}\mathfrak{d}_*^{\leq k}(\widecheck{\nu(r)}, \widecheck{\nu(u)}, \widecheck{b_*})|^2,
\end{split}
\eea
where $\trchc$, $\chih$, $\ze$, $\eta$, $\trchbc$, $\chibh$, $\ombc$, $\xib$ are the Ricci coefficients   of the integrable frame  of  $\Si_*$, and $\nu= e_3+b_* e_4$ is    tangent to $\Si_*$. 
\end{definition} 

\begin{definition}
\lab{definition:int-normson-RcSi_*} 
We define the curvature norm on $\Si_*$, relative to the integrable frame of $\Si_*$, 
\bea
\bsplit
\Rkstarr^2_k&:= \int_{\Si_*}\left[r^{4+\de}|\dk^{\leq k}(\a, \b)|^2+r^{4}|\dk^{\leq k}(\rhoc, \rhodc)|^2+r^{2}|\dk^{\leq k}\bb|^2+|\dk^{\leq k}\aa|^2\right],
\end{split}
\eea
where $\a$, $\b$, $\rho$, $\rhod$, $\bb$, $\aa$ denote the curvature components relative to  the integrable   frame of $\Si_*$.
\end{definition}

The following proposition, which provides the control of the Ricci coefficients of the integrable frame of $\Si_*$, is the main result of this section. 
\begin{proposition}
\lab{prop:controlofRiccicoefficientsonSigmastar}
The following estimates  hold true for the integrable frame of $\Si_*$
\bea
 \Skstarr_k&\les& \ep_0+\Rkstarr_k, \qquad k\le k_{large}+7.
\eea
\end{proposition}

\begin{proof}
We rely on $L^2(\Si_*)$-estimates  for  the main quantities, with the exception  of   $\ell=1$ modes and averages which are estimated  in $L^2_uL^\infty(S)$. Also, we rely on the following local bootstrap assumptions
\bea\lab{eq:localbootassforframeSistaronSistar:chap9}
 \Skstarr_k &\leq& \ep, \qquad k\le k_{large}+7.
\eea
Finally, we rely repeatedly on the results of sections \ref{sect.LastSlice-Intro} and \ref{sec:preliminaryestimatesonSistar:chap5}.

{\bf Step 1.}  The estimate for $\trchc$ follows immediately from our  GCM conditions according to which we have $\trchc=0$ on $\Si_*$. 

{\bf Step 2.} Since we have $\trch=\frac 2 r$ and $\trchb=-\frac{2\Up}{r}$ on $S_*$ according to our GCM conditions, we infer from  the Gauss equation  
\beaa
K &=& -\rho -\frac{1}{4}\trch\trchb +\frac{1}{2}\chih\c\chibh=-\rho+\frac{\Up}{r^2}+\frac{1}{2}\chih\c\chibh=\frac{1}{r^2}-\rhoc+\frac{1}{2}\chih\c\chibh\quad \textrm{on}\quad S_*. 
\eeaa
We deduce, using the trace theorem, 
\beaa
 r^2 \left\|K-\frac{1}{r^2}\right\| _{\hk_{k-1}(S_*)} \les  r^2 \| \rhoc \| _{\hk_{k-1} (S_*)}+ \|   \chih\c \chibh  \|_{\hk_{k-1}(S_*)}^2\les \ep^2 +\Rkstarr_{k}, \qquad k\le k_{large}+7.
\eeaa
Next, recall that $\phi$ denotes  the uniformization factor  of $S_*$, see \eqref{eq:formofthemetriconSstarusinguniformization:0} in section \ref{section:canonical-coord-BPG}. Using the  effective uniformization theorem, see Theorem   \ref{Thm:effectiveU1-Intro},
\beaa
\|\phi\|_{\hk_{k+1}(S_*)}\les r^2\left\|K-\frac{1}{r^2}\right\|_{\hk_{k-1}(S_*)} \les \ep_0 +\Rkstarr_{k}, \qquad k\le k_{large}+7,
\eeaa
we deduce
\bea
 \|\phi\|_{\hk_{k+1}(S_*)}\les  \ep_0 +\Rkstarr_{k}, \qquad  \mbox{for all}\,\, k\le k_{large} +7. 
\eea

{\bf Step 3.}  Using the formula $\Delta \Jp+\frac{2}{r^2}\Jp = \frac{2}{r^2}(1-e^{-2\phi})\Jp$, see Lemma  \ref{lemma:statementeq:DeJp.S_*}  and its proof, we infer from the control of $\phi$ provided by Step 2
\bea
r^2\left\|\left(\Delta \Jp+\frac{2}{r^2}\Jp\right)\right\|_{\hk_{k+1} (S_*)} & \les \ep_0 +\Rkstarr_{k}, \qquad  \mbox{for all}\,\, k\le k_{large} +7. 
\eea

{\bf Step 4.} Recall from Lemma  \ref{Le:Si*-ell=1modes} that we have
\beaa
\nab_\nu \Big[\big( r^2\De +2)\Jp\Big]=O\big(\dkb^{\le 1} \Ga_b\big). 
\eeaa
Integrating from $S_*$ along $\Si_*$ and using \eqref{eq:localbootassforframeSistaronSistar:chap9} as well as the control on $S_*$ of Step 3, we infer
\beaa
\sup_{S\subset\Si_*}r^2\left\|\left(\Delta \Jp+\frac{2}{r^2}\Jp\right)\right\|_{\hk_{k-1} (S)}   &\les&  \ep_0+\sqrt{u_*}\ep\les\sqrt{u_*}\ep, \qquad k\leq k_{large}+7.
\eeaa
Then,  proceeding  as in the proof of Lemma  \ref{Le:Si*-ell=1modes}, we have 
\beaa
 \left\|\dds_2\dds_1 \Jp\right\|_{\hk_k(S)}  &\les& \left\|\left(\Delta \Jp+\frac{2}{r^2}\Jp\right)\right\|_{\hk_{k} (S)}+r^{-1} \left\|\Ga_g\right\|_{\hk_k(S)}, 
\eeaa
where by $\dds_1\Jp$, we mean either $\dds_1(\Jp,0)$ or $\dds_1(0,\Jp)$. Together with the above estimate and \eqref{eq:localbootassforframeSistaronSistar:chap9} yields
\beaa
\sup_{S\subset\Si_*}r^2\left\|\dds_2\dds_1 \Jp\right\|_{\hk_{k-1} (S)}   \les\sqrt{u_*}\ep, \qquad k\leq k_{large}+7.
\eeaa
Together with the dominance condition for $r$ on $\Si_*$, we deduce
\beaa
\sup_{S\subset\Si_*}r\left(\left\|\left(\Delta \Jp+\frac{2}{r^2}\Jp\right)\right\|_{\hk_{k-1} (S)}+\left\|\dds_2\dds_1 \Jp\right\|_{\hk_{k-1} (S)}\right)   \les\frac{\sqrt{u_*}\ep}{r}\les\ep_0, \qquad k\leq k_{large}+7.
\eeaa
Using $\nu(\Jp)=0$, we recover $\nu$ derivatives as well, and hence, for $k\leq k_{large}+7$, 
\bea
\sup_{S\subset\Si_*}  r\left(\left\|\dk_*^{k-1}\left(\Delta \Jp+\frac{2}{r^2}\Jp\right)\right\|_{L^2(S)}+\left\|\dk_*^{k-1}\dds_2\dds_1 \Jp\right\|_{L^2(S)}\right)    \les  \ep_0.
\eea

{\bf Step 5.}  Next, differentiating Codazzi for $\chih$  by $\div$, and integrating against $\Jp$, we infer, see Steps 1  of the proof of Proposition \ref{prop:controlofell=1modesonSigmastar} in section \ref{section:estimates-rll=1-Si_*}, 
\beaa
(\div\ddd_2\chih)_{\ell=1} &=&  \frac{1}{r}(\div\ze)_{\ell=1}-(\div\b)_{\ell=1} +r^{-3}\int_S\dkb^{\leq 1}(\chih\c \ze)\Jp.
\eeaa
This then yields, after  after integration by parts, 
\beaa
(\div\ze)_{\ell=1} &=& r(\div\b)_{\ell=1}+r^{-2}\int_S\chih\c\dds_2\dds_1\Jp+r^{-3}\int_S\dkb^{\leq 1}(\chih\c\ze)\Jp. 
\eeaa
Together with \eqref{eq:localbootassforframeSistaronSistar:chap9} and the control for $\dds_2\,\dds_1\Jp$ of Step 4, we infer, using integration by parts on $S$ for angular derivatives to avoid  loss of derivatives, for $k\leq k_{large}+7$, 
\beaa
r^3\left( \int_{u=1} ^ {u_*}  \big|\nu^{\leq k+1}\big((\div \ze)_{\ell=1}\big)\big|^2\right)^{\frac{1}{2}}    &\les&  \ep^2 +\Rkstarr_{k}+r\left(\int_{u=1}^{u_*}\left|\int_S\nab_\nu^{k+1}\b\c\dkb\Jp\right|^2\right)^{\frac{1}{2}}.
\eeaa
Since $\nu=e_3+b_*e_4$, and using the Bianchi identities for $\nab_3\b$ and $\nab_4\b$, as well as the commutator Lemma \ref{Lemma:Commutation-Si_*}, we have
\beaa
\nab_\nu^{k+1}\b &=& \nab(\nab_\nu^k\rhoc)+\dual\nab(\nab_\nu^k\rhod) -b_*\div(\nab_\nu^k\a)+\lot
\eeaa 
Plugging in the above and integrating the angular derivatives by parts  to avoid a loss of derivatives, we infer
\beaa
r^3\left( \int_{u=1} ^ {u_*}  \big|\nu^{\leq k+1}\big((\div \ze)_{\ell=1}\big)\big|^2\right)^{\frac{1}{2}}    \les  \ep^2 +\Rkstarr_{k}, \quad k\leq k_{large}+7,
\eeaa
and hence
\bea
r^3\left( \int_{u=1} ^ {u_*}  \big|\nu^{\leq k+1}\big((\div \ze)_{\ell=1}\big)\big|^2\right)^{\frac{1}{2}}    \les  \ep_0 +\Rkstarr_{k}, \quad k\leq k_{large}+7.
\eea

{\bf Step 6.}  Next, recall that we have already established  an estimate for   $(\trchbc)_{\ell=1}$ on $\Si_*$ in  \eqref{eq:controlofell=1modeforthorcandkabconSistarinproofoThmM0Step4}, i.e. 
\beaa
\sup_{\Si_*}r^2u^{1+\dec}|(\kabc)_{\ell=1}| &\les& \ep_0.
\eeaa
To derive an estimate for higher order derivatives in $\nu$, we rely on the following equation of Corollary \ref{corofLemma:transport.alongSi_*1} 
\beaa
\nn&&\nu\left(\int_S\left( \lap\trchbc+\frac{2\Up}{r}\div \ze\right)\Jp\right)\\
\nn&=& O(r^{-3})\int_S\kabc\Jp   +O(r^{-2})\int_S\div\ze\Jp +O(r^{-1})\int_S\div\bb\Jp     +O(r^{-2})\int_S\rhoc\Jp \\
 &&+O(r^{-1})\int_S\div\b\Jp +r^{-1}\int_S\left(\Delta+\frac{2}{r^2}\right)\Jp\dkb^{\leq 1}\Ga_b+r^{-2}\int_S\dkb^{\le 2 }(\Ga_b\c \Ga_b)\Jp,
\eeaa
which yields
\beaa
\nn&&\nu\left(\int_S\left( \lap\kabc+\frac{2\Up}{r}\div \ze\right)\Jp\right)\\
\nn&=&   O(r^{-1})\int_S\div\bb\Jp  +  O(r^{-3})\int_S\Jp\dkb^{\leq 1}\Ga_g  +r^{-1}\int_S\left(\Delta+\frac{2}{r^2}\right)\Jp\dkb^{\leq 1}\Ga_b\\
&&+r^{-2}\int_S\dkb^{\le 2 }(\Ga_b\c \Ga_b)\Jp.
\eeaa
Differentiating in $\nu$, using integration by parts on $S$ for angular derivatives to avoid  loss of derivatives,  relying on the estimates of Step 4 and Step 5, and using \eqref{eq:localbootassforframeSistaronSistar:chap9}, we deduce, for $k\leq k_{large}+7$, 
\beaa
r^2\left( \int_{u=1}^{u_*} \big|\nu^{\leq k+1}\big((\trchbc)_{\ell=1}\big)\big|^2\right)^{\frac{1}{2}}  &\les& r^{-1}\ep+\ep^2+\Rkstarr_{k}+\left(\int_{u=1}^{u_*}\left|\int_S\nab_\nu^{k+1}\bb\c\dkb\Jp\right|^2\right)^{\frac{1}{2}}\\
&&+\left(\int_{u=1}^{u_*}\left|\int_S\nab_\nu^{k+1}\b\c\dkb\Jp\right|^2\right)^{\frac{1}{2}}.
\eeaa

Next, as in Step 5, we use the fact that $\nu=e_3+b_*e_4$ and the Bianchi identities, as well as the commutator Lemma \ref{Lemma:Commutation-Si_*}, to derive
\beaa
\nab_\nu^{k+1}\b &=& \nab(\nab_\nu^k\rhoc)+\dual\nab(\nab_\nu^k\rhod) -b_*\div(\nab_\nu^k\a)+\lot,\\
\nab_\nu^{k+1}\bb &=& -\div(\nab_\nu^k\aa)-b_*\nab(\nab_\nu^k\rhoc)+b_*\dual\nab(\nab_\nu^k\rhod) +\lot
\eeaa 
Plugging in the above and integrating the angular derivatives by parts to avoid a loss of derivatives, we infer
\beaa
r^2\left( \int_{u=1}^{u_*} \big|\nu^{\leq k+1}\big((\trchbc)_{\ell=1}\big)\big|^2\right)^{\frac{1}{2}}  &\les& r^{-1}\ep+\ep^2+\Rkstarr_{k}, \quad k\leq k_{large}+7.
\eeaa
Together with  the dominance condition for $r$ on $\Si_*$, we obtain 
\bea
r^2 \left( \int_{u=1}^{u_*} \big|\nu^{\leq k+1}\big((\trchbc)_{\ell=1}\big)\big|^2\right)^{\frac{1}{2}}    &\les&  \ep_0 +\Rkstarr_{k}, \quad k\leq k_{large}+7.
\eea

{\bf Step 7.} We are now ready to control $\trchbc$. 

{\bf Step 7a.} First, we control $\ov{\trchbc}$. Recall the definition of the Hawking mass $m_H$
\beaa
\frac{2m_H}{r} &=& 1+\frac{1}{16\pi}\int_S\trch \trchb.
\eeaa
Since we have $\trch=\frac{2}{r}$ by our GCM conditions on $\Si_*$, we infer on $\Si_*$
\beaa
\ov{\trchb} &=& -\frac{2\left(1-\frac{2m_H}{r}\right)}{r}.
\eeaa
where $\ov{\trchb}$ denotes the average of $\trchb$ on $S$. In particular, we infer
\beaa
\ov{\trchbc} &=& -\frac{4(m-m_H)}{r^2}.
\eeaa
Since we have by Lemma \ref{lemma:transportequationforHawkingmass} and by Proposition \ref{prop:controlofell=0modesonSigmastar}
\beaa
\nu(m_H) = \int_S\dkb^{\leq 1}(\Ga_b\c\Ga_b), \qquad \sup_{\Si_*}u^{1+2\dec}|m_H-m| \les \ep_0,
\eeaa
we infer, together with \eqref{eq:localbootassforframeSistaronSistar:chap9}, and  integrating the angular derivatives by parts  to avoid a loss of derivatives, 
\beaa
\left( \int_{u=1}^{u_*} \big|\nu^{\leq k+1}\big(m_H-m\big)\big|^2\right)^{\frac{1}{2}} +r^2 \left( \int_{u=1}^{u_*} \big|\nu^{\leq k+1}\big(\ov{\trchbc}\big)\big|^2\right)^{\frac{1}{2}}    &\les&  \ep_0, \quad k\leq k_{large}+7.
\eeaa

{\bf Step 7b.} Next, we control $\nu^k(\Cb_0)$ and $\nu^k(\Cbp)$. In view of our GCM conditions for $\trchbc$ and the fact that $\nu(\Jp)=0$, we have 
\beaa
\nu^k(\trchbc) &=& \nu^k(\Cb_0)+\sum_p\nu^k(\Cbp)\Jp=\ov{\nu^k(\Cb_0)}+\sum_p\ov{\nu^k(\Cbp)}\Jp+h_k,\\
h_k &:=& \Big(\nu^k(\Cb_0)-\ov{\nu^k(\Cb_0)}\Big)+\sum_p\Big(\nu^k(\Cbp)-\ov{\nu^k(\Cbp)}\Big)\Jp.
\eeaa
Using a Poincar\'e inequality on $S$, the fact that $\nab(\Cb_0)=0$ and $\nab(\Cbp)=0$, \eqref{eq:localbootassforframeSistaronSistar:chap9}, and the following commutator formula of Lemma \ref{Lemma:Commutation-Si_*}
  \beaa
\, [ \nab_{\nu}, \nab]f &=&  -\frac{2}{r}\nab f+\Ga_b \c \nab_\nu f+ r^{-1} \Ga_b \c \dk^{\leq 1} f,
 \eeaa
we easily obtain 
\beaa
r\|h_k\|_{L^2(\Si_*)}\les \ep^2\les \ep_0, \quad k\leq k_{large}+8.
\eeaa
Then, coming back to the above identity for $\nu^k(\trchbc)$, multiplying it respectively with 1 or $J^{(p)}$, and integrating on $S$, we easily obtain, using the properties of $\Jp$ in Lemma \ref{Le:Si*-ell=1modes}, for $k\leq k_{large}+8$, 
\beaa
r^2 \left( \int_{u=1}^{u_*} \big|\nu^k(\Cb_0, \Cbp)\big|^2\right)^{\frac{1}{2}}    &\les&  \ep_0+r^2 \left( \int_{u=1}^{u_*} \big|\nu^k\big(\ov{\trchbc}, \,(\trchbc)_{\ell=1}\big)\big|^2\right)^{\frac{1}{2}} , 
\eeaa
Together with the above estimate for $\ov{\trchbc}$, and the one of Step 6 for $(\trchbc)_{\ell=1}$, we infer
\beaa
r^2 \left( \int_{u=1}^{u_*} \big|\nu^{k+1}(\Cb_0, \Cbp)\big|^2\right)^{\frac{1}{2}}    &\les&  \ep_0+\Rkstarr_k, \quad  k\leq k_{large}+7.
\eeaa
  
{\bf Step 7c.} We finally control $\dk_*^k\trchbc$. First, the above identity for $\nu^k(\trchbc)$ and the above control of $\nu^k(\Cb_0)$ and $\nu^k(\Cbp)$ implies
\beaa
r\|\nu^{k+1}\trchbc\|_{L^2(\Si_*)}   \les \ep_0 +\Rkstarr_k, \qquad  k\le  k_{large}+7.
\eeaa
Next, to recover the other derivatives, we rely on Corollary  \ref{prop:2D-Hodge4} which yields 
\beaa
  \|\nu^k(\trchbc)\|_{\hk_{j+2}  (S)} &\les&  r^2\|\dds_2\,\dds_1\nu^k(\trchbc)\|_{\hk_j(S)}+r\big| (\nu^k(\trchbc))_{\ell=1}\big|+r|\ov{\nu^k(\trchbc)}|.
  \eeaa
Together with the above estimate for $\ov{\trchbc}$, and the one of Step 6 for $(\trchbc)_{\ell=1}$, we infer
\beaa
r\|\dkb^{\leq j+2}\nu^k(\trchbc)\|_{L^2(\Si_*)} \les \ep_0 +\Rkstarr_{k}+r^3\|\dkb^{\leq j}\dds_2\,\dds_1\nu^k(\trchbc)\|_{L^2(\Si_*)}, \,\,\, k+j\leq k_{large}+6.
\eeaa
Also, differentiating the above identity for  $\nu^k(\trchbc)$, and since $\nab(\Cb_0)=\nab(\Cbp)=0$, we have
\beaa
\dds_2\,\dds_1\nu^k(\trchbc) &=& [\dds_2\,\dds_1,\nu^k]\Cb_0+\sum_p[\dds_2\,\dds_1,\nu^k]\Cbp\Jp+\sum_p\nu^k(\Cbp)\dds_2\,\dds_1\Jp.
\eeaa
Using the above commutator formula for $[ \nab_{\nu}, \nab]$, \eqref{eq:localbootassforframeSistaronSistar:chap9}, and the control for $\dds_2\,\dds_1\Jp$ of Step 4, we infer
\beaa
r^3\|\dkb^{\leq j}\dds_2\,\dds_1\nu^k(\trchbc)\|_{L^2(\Si_*)} &\les& \ep^2\les\ep_0, \quad k+j\leq k_{large}+6,
\eeaa
and hence
\beaa
r\|\dkb^{\leq j+2}\nu^k(\trchbc)\|_{L^2(\Si_*)} &\les& \ep_0 +\Rkstarr_{k}, \quad k+j\leq k_{large}+6.
\eeaa
Together with the above estimate for $\nu^k(\trchbc)$, we deduce
\bea
r\|\dk_*^{\leq k+1}  \trchbc\|_{L^2(\Si_*)}   &\les& \ep_0 +\Rkstarr_{k}, \quad k\leq k_{large}+7.
\eea

{\bf Step 8.}  In view of the definition of $\mu$, we have
\beaa
\muc &=&  -\div\ze -\rhoc +\frac{1}{2}\chih\c\chibh.
\eeaa
Together with \eqref{eq:localbootassforframeSistaronSistar:chap9} and the  control of $\nu^{k }(\div\ze)_{\ell=1}$ in Step 5, we infer
\beaa
 r^3\left( \int_{u=1} ^ {u_*} \big|\nu^{\leq k+1}\big((\muc)_{\ell=1}\big)\big|^2\right)^{1/2}   & \les&  \ep_0 +\ep^2+\Rkstarr_{k}\\
 &&+ r\left( \int_{u=1} ^ {u_*} \left|\int_S\nu^{k+1}(\rhoc)\Jp\right|^2\right)^{1/2}, \quad k\leq k_{large}+7. 
\eeaa
Next, as in Step 5, we use the fact that $\nu=e_3+b_*e_4$ and the Bianchi identities, as well as the commutator Lemma \ref{Lemma:Commutation-Si_*}, to derive
\beaa
\nab_\nu^{k+1}\rhoc &=& -\div(\nab_\nu^k\bb)+b_*\div(\nab_\nu^k\b)+\lot
\eeaa 
Plugging in the above and integrating the angular derivatives by parts to avoid a loss of derivatives, we infer
\beaa
 r^3\left( \int_{u=1} ^ {u_*} \big|\nu^{\leq k+1}\big((\muc)_{\ell=1}\big)\big|^2\right)^{1/2}   & \les&  \ep_0+ \Rkstarr_{k}, \quad k\leq k_{large}+7.
\eeaa

{\bf Step 9.}  In view of the proof of Proposition \ref{prop:controlofell=0modesonSigmastar}, the average of $\muc$ is given by
\beaa
\ov{\muc} &=& \frac{2(m_H-m)}{r^3},
\eeaa
which together with the control of $m_H-m$ in Step 7a implies
\beaa
r^3 \left( \int_{u=1}^{u_*} \big|\nu^{\leq k+1}\big(\ov{\muc}\big)\big|^2\right)^{\frac{1}{2}}    &\les&  \ep_0, \quad k\leq k_{large}+7.
\eeaa
Also recall the GCM condition for $\mu$ on $\Si_*$
\beaa
\muc &=& M_0+\sum_pM_p\Jp. 
\eeaa
In view of the above formula for $\muc$ on $\Si_*$, the above control of $\ov{\muc}$, and the estimate for $(\muc)_{\ell=1}$ in Step 8, the control of $\muc$ is completely analogous to the one of $\trchbc$ in Step 7, and we infer the corresponding estimate
\bea
r^2\|\dk_*^{\leq k+1}\muc\|_{L^2(\Si_*)}   &\les& \ep_0 +\Rkstarr_{k}, \quad  k\leq k_{large}+7.
\eea

{\bf Step 10.} Next, recall that we have
\beaa
\div\ze &=&  -\muc -\rhoc +\frac{1}{2}\chih\c\chibh,\\
\curl\ze &=& \rhod -\frac{1}{2}\chih\wedge\chibh.
\eeaa
In view of \eqref{eq:localbootassforframeSistaronSistar:chap9} and the control of $\muc$ in Step 9, as well as the commutator Lemma \ref{Lemma:Commutation-Si_*}, we infer
\beaa
r^2\|\ddd_1(\dk_*^{\leq k}\ze)\|_{L^2(\Si_*)} &\les& \ep_0+\ep^2+ \Rkstarr_{k}, \quad k\leq k_{large}+7,
\eeaa
which together with the elliptic estimate of Lemma \ref{prop:2D-Hodge1} implies
\beaa
r\|\dkb^{\leq 1}\dk_*^{\leq k}\ze\|_{L^2(\Si_*)} &\les& \ep_0+ \Rkstarr_{k}, \quad k\leq k_{large}+7.
\eeaa

To control $\nab_\nu^{k+1}\ze$, we come back to the above system and differentiate it w.r.t. $\nab_\nu^{k+1}$. We obtain
\beaa
\ddd_1\left(\nab_\nu^{k+1}\ze\right) &=&  \nab_\nu^{k+1}(-\rhoc, \rhod) +h_{k+1}
\eeaa
where $h_{k+1}$ satisfies, in view of \eqref{eq:localbootassforframeSistaronSistar:chap9} and the control of $\muc$ in Step 9, as well as the commutator Lemma \ref{Lemma:Commutation-Si_*},
\beaa
r^2\|h_{k+1}\|_{L^2(\Si_*)} &\les& \ep_0+\ep^2\les \ep_0, \quad k\leq k_{large}+7.
\eeaa
Next, as in Step 5, we use the fact that $\nu=e_3+b_*e_4$ and the Bianchi identities, as well as the commutator Lemma \ref{Lemma:Commutation-Si_*}, to derive
\beaa
\nab_\nu^{k+1}(-\rhoc, \rhod) &=& \Big(\div(\nab_\nu^k\bb)-b_*\div(\nab_\nu^k\b), -\curl(\nab_\nu^k\bb)-b_*\curl(\nab_\nu^k\b)\Big)+\lot
\eeaa 
Plugging in the above, and using the above estimate for $h_{k+1}$, we infer
\beaa
r\|\nab_\nu^{k+1}\ze\|_{L^2(\Si_*)} &\les&  r\|(r\ddd_1)^{-1}\dkb \nab_\nu^{k}(\bb, \b)\|_{L^2(\Si_*)}+\ep_0, \quad k\leq k_{large}+7.
\eeaa
Using the ellipticity of $(r\ddd_1)^{-1}\dkb$ on the spheres $S$, this yields 
\beaa
r\|\nab_\nu^{k+1}\ze\|_{L^2(\Si_*)} &\les&  \ep_0+ \Rkstarr_{k}, \quad k\leq k_{large}+7.
\eeaa
Together with the above estimate for $\dkb^{\leq 1}\dk_*^{k}\ze$, we deduce
\bea
r\|\dk_*^{\leq k+1}\ze\|_{L^2(\Si_*)} &\les &  \ep_0+ \Rkstarr_{k}, \quad k\leq k_{large}+7.
\eea

{\bf Step 11.} Using  Codazzi for $\chih$,  \eqref{eq:localbootassforframeSistaronSistar:chap9}, the control of Step 1 for $\trchc$, and the control of Step 10 for $\ze$, as well as the commutator Lemma \ref{Lemma:Commutation-Si_*}, we have
\beaa
r^2\|\ddd_2(\dk_*^{\leq k}\chih)\|_{L^2(\Si_*)} &\les &  \ep_0+\ep^2+ \Rkstarr_{k}, \quad k\leq k_{large}+7,
\eeaa
which together with the elliptic estimate of Lemma \ref{prop:2D-Hodge1} implies
\bea
r\|\dkb^{\leq 1}(\dk_*^{\leq k}\chih)\|_{L^2(\Si_*)} &\les &  \ep_0+ \Rkstarr_{k}, \quad k\leq k_{large}+7.
\eea
Note that this is not yet the desired control for $\chih$ as we still need to recover $\nab_\nu^{k+1}\chih$ which will be done in Step 15.

{\bf Step 12.} Using  Codazzi for $\chibh$,  \eqref{eq:localbootassforframeSistaronSistar:chap9}, the control of Step 7 for $\trchbc$, and the control of Step 10 for $\ze$, as well as the commutator Lemma \ref{Lemma:Commutation-Si_*}, we have
\beaa
r\|\ddd_2(\dk_*^{\leq k}\chibh)\|_{L^2(\Si_*)} &\les &  \ep_0+\ep^2+ \Rkstarr_{k}, \quad k\leq k_{large}+7,
\eeaa
which together with the elliptic estimate of Lemma \ref{prop:2D-Hodge1} implies
\bea
\|\dkb^{\leq 1}(\dk_*^{\leq k}\chibh)\|_{L^2(\Si_*)} &\les &  \ep_0+ \Rkstarr_{k}, \quad k\leq k_{large}+7.
\eea
Note that this is not yet the desired control for $\chibh$ as we still need to recover $\nab_\nu^{k+1}\chibh$ which will be done in Step 15.

{\bf Step 13.} Recall the equation  for  $\ddd_2\dds_2\,\eta $ derived in  Proposition \ref{prop:additional.eqtsM4}  
\beaa
2\ddd_2\dds_2\eta &=&  -\nab_3\nab\kac    -\frac{2}{r}\nab_3 \ze  -\frac{2}{r}\bb  + r^{-2} \dkb^{\le 1} \Ga_g+ r^{-1} \dkb^{\le 1 } (\Ga_b\c \Ga_b).
 \eeaa
We infer
\beaa
\ddd_2\dds_2\big(\dk_*^{\leq k}\eta\big) &=&  \dkb^{\leq 1}\big(h_{k}\big), \quad k\leq k_{large}+7,
 \eeaa
where $h_{k}$ satisfies in view of  \eqref{eq:localbootassforframeSistaronSistar:chap9},  the estimate of Step 1 for $\trchc$ and  of Step 10 for $\ze$, and  the commutator Lemma \ref{Lemma:Commutation-Si_*}, 
\beaa
r^2\|h_{k}\|_{L^2(\Si_*)} &\les &  \ep_0+r^{-1}\ep+\ep^2+ \Rkstarr_{k}, \quad k\leq k_{large}+7.
\eeaa
Using  the ellipticity of $(r\ddd_2)^{-1}\dkb^{\leq 1}$ on the spheres $S$, and the elliptic estimate for $\dds_2$ of Lemma \ref{prop:2D-Hodge2}, we infer
\beaa
\|\dkb^{\leq 1}(\dk_*^{\leq k}\eta)\|_{L^2(\Si_*)} &\les &  r^2\|h_{k}\|_{L^2(\Si_*)}+r^2\big| (\ddd_1\nab_\nu^{\leq k}\eta)_{\ell=1}\big|, \quad k\leq k_{large}+7.
\eeaa
Together with the above estimate for $h_{k}$ and the dominance condition for $r$ on $\Si_*$, this yields
\beaa
\|\dkb^{\leq 1}(\dk_*^{\leq k}\eta)\|_{L^2(\Si_*)} &\les &  \ep_0+ \Rkstarr_{k}+r^2\big| (\ddd_1\nab_\nu^{\leq k}\eta)_{\ell=1}\big|, \quad k\leq k_{large}+7.
\eeaa
Since $\ddd_1=(\div, \curl)$, and using  the null structure equation $\curl\eta = \rhod -\frac{1}{2}\chih\wedge\chibh$, \eqref{eq:localbootassforframeSistaronSistar:chap9}, and  the commutator Lemma \ref{Lemma:Commutation-Si_*}, we infer
\beaa
\|\dkb^{\leq 1}(\dk_*^{\leq k}\eta)\|_{L^2(\Si_*)} &\les &  \ep_0+ \Rkstarr_{k}+r^2\big| (\div\nab_\nu^{\leq k}\eta)_{\ell=1}\big|, \quad k\leq k_{large}+7.
\eeaa
Together with the GCM condition $(\div\eta)_{\ell=1}=0$ on $\Si_*$, the fact that $\nu$ is tangent to $\Si_*$, \eqref{eq:localbootassforframeSistaronSistar:chap9}, and Corollary \ref{Corr:nuSof integrals}, we finally obtain 
\bea
\|\dkb^{\leq 1}(\dk_*^{\leq k}\eta)\|_{L^2(\Si_*)} &\les &  \ep_0+ \Rkstarr_{k}, \quad k\leq k_{large}+7.
\eea

{\bf Step 14.} Recall the equation  for  $\ddd_2\dds_2\,\xib$ derived in  Proposition \ref{prop:additional.eqtsM4}  
\beaa
 2\ddd_2\dds_2\xib  &=& -\nab_3\nab\kabc   -\frac{2}{r}\nab_3 \ze    -\frac{2}{r}\bb     +r^{-2}\dkb^{\leq 1}\Ga_g+r^{-1}\dkb^{\leq 1}(\Ga_b\c \Ga_b).
 \eeaa
We infer
\beaa
\ddd_2\dds_2\big(\dk_*^{\leq k}\xib\big) &=&  \dkb^{\leq 1}\big(h_{k}\big), \quad k\leq k_{large}+7, 
\eeaa
where $h_{k}$ satisfies in view of  \eqref{eq:localbootassforframeSistaronSistar:chap9},  the estimate of Step 7 for $\trchbc$ and  of Step 10 for $\ze$, and  the commutator Lemma \ref{Lemma:Commutation-Si_*}, 
\beaa
r^2\|h_{k}\|_{L^2(\Si_*)} &\les &  \ep_0+r^{-1}\ep+\ep^2+ \Rkstarr_{k}, \quad k\leq k_{large}+7.
\eeaa
Since $\xib$ satisfies $\curl\xib=\Ga_b\c\Ga_b$ and the GCM condition $(\div\xib)_{\ell=1}=0$ on $\Si_*$, we proceed exactly as in Step 13 and obtain 
\bea
\|\dkb^{\leq 1}(\dk_*^{\leq k}\xib)\|_{L^2(\Si_*)} &\les &  \ep_0+ \Rkstarr_{k}, \quad k\leq k_{large}+7.
\eea

{\bf Step 15.}  We may now  complete the estimates for $\chih$ and $\chibh$. Recall from Step 11 and Step 12 that we have
\beaa
r\|\dkb^{\leq 1}(\dk_*^{\leq k}\chih)\|_{L^2(\Si_*)}+\|\dkb^{\leq 1}(\dk_*^{\leq k}\chibh)\|_{L^2(\Si_*)} &\les &  \ep_0+ \Rkstarr_{k}, \quad k\leq k_{large}+7,
\eeaa
so that we still need to recover $\nab_\nu^{k+1}\chih$ and $\nab_\nu^{k+1}\chibh$. In view of the fact that $\nu=e_3+b_*e_4$, the null structure equations for $\nab_3\chih$, $\nab_3\chibh$, $\nab_4\chih$, and $\nab_4\chibh$, \eqref{eq:localbootassforframeSistaronSistar:chap9},  the estimate of  Step 10 for $\ze$, of Step 13 for $\eta$ and of Step 14 for $\xib$, we have 
\beaa
r\|\nab_\nu^{k+1}\chih\|_{L^2(\Si_*)}+\|\nab_\nu^{k+1}\chibh\|_{L^2(\Si_*)} &\les &  \ep_0+\ep^2+ \Rkstarr_{k}, \quad k\leq k_{large}+7.
\eeaa
Together with the above estimates, we deduce
\beaa
r\|\dk_*^{\leq k+1}\chih\|_{L^2(\Si_*)}+\|\dk_*^{\leq k+1}\chibh\|_{L^2(\Si_*)} &\les &  \ep_0+ \Rkstarr_{k}, \quad k\leq k_{large}+7.
\eeaa

{\bf Step 16.} We have 
\beaa
\nab_a(e_3(r)) &=& [e_a, e_3]r=(\ze_a-\eta_a)e_3(r)-\xib_ae_4(r)= (\ze_a-\eta_a)e_3(r)-\xib_a.
\eeaa
Together with the estimate of  Step 10 for $\ze$, of Step 13 for $\eta$ and of Step 14 for $\xib$, the commutator Lemma \ref{Lemma:Commutation-Si_*}, and a Poincar\'e inequality, we  derive
  \bea
r^{-1}\|\dkb^{\leq 2}\dkb\dk_*^{k-1}e_3(r)\big)\|_{L^2(\Si_*)} &\les& \ep_0+ \Rkstarr_k, \qquad  k\le  k_{large}+7,
\eea
and 
 \bea
r^{-1}\|\dkb^{\leq 2}\big(\nu^k(e_3(r))-\ov{\nu^k(e_3(r))}\big)\|_{L^2(\Si_*)} \les \ep_0+ \Rkstarr_k, \qquad  k\le  k_{large}+7.
\eea

{\bf Step 17.}  We have 
\beaa
\nab_a(e_3(u)) &=& [e_a, e_3]u=(\ze_a-\eta_a)e_3(u)-\xib_ae_4(u)= (\ze_a-\eta_a)e_3(u).
\eeaa
Together with the estimate of  Step 10 for $\ze$ and of Step 13 for $\eta$, the commutator Lemma \ref{Lemma:Commutation-Si_*}, and a Poincar\'e inequality, we  derive
 \bea
r^{-1}\|\dkb^{\leq 2}\dkb\dk_*^{k-1}e_3(u)\big)\|_{L^2(\Si_*)} &\les& \ep_0+ \Rkstarr_k, \qquad  k\le  k_{large}+7,
\eea
and 
 \bea
r^{-1}\|\dkb^{\leq 2}\big(\nu^k(e_3(u))-\ov{\nu^k(e_3(u))}\big)\|_{L^2(\Si_*)} \les \ep_0+ \Rkstarr_k, \qquad  k\le  k_{large}+7.
\eea

{\bf Step 18.}   Since $u+r=c_{\Si_*}$ on $\Si_*$, and since $\nu$ is tangent to $\Si_*$, we have $\nu(r+u)=0$ and hence 
\beaa
0 &=& \nu(u+r)=e_3(u+r)+b_*e_4(u+r)=e_3(u)+e_3(r)+b_*,
\eeaa
so that $b_*=-e_3(u)-e_3(r)$. In view of the above estimates for $e_3(u)-\ov{e_3(u)}$ and $e_3(r)-\ov{e_3(r)}$,  we deduce
 \bea
r^{-1}\|\dkb^{\leq 2}\dkb\dk_*^{k-1}b_*\big)\|_{L^2(\Si_*)} &\les& \ep_0+ \Rkstarr_k, \qquad  k\le  k_{large}+7,
\eea
and 
 \bea
r^{-1}\|\dkb^{\leq 2}\big(\nu^k(b_*)-\ov{\nu^k(b_*)}\big)\|_{L^2(\Si_*)} \les \ep_0+ \Rkstarr_k, \qquad  k\le  k_{large}+7.
\eea
In particular, by Sobolev,
\bea
\|   \nu^k(b_*)-\ov{\nu^k(b_*)} \|_{L_u^2L^\infty(S)} &\les& \ep_0+ \Rkstarr_k, \qquad  k\le  k_{large}+7.
\eea

{\bf Step 19.}
Since ${b_*}_{|_{SP}}=-1 -\frac{2m}{r} $, where  $SP$ corresponds to $\th=\pi$, and since $\nu(\th)=0$, we have ${\nu^{k}(b_*)}_{|_{SP}}=O(r^{-k-1})$ for $k\geq 1$, and hence, in view of the dominant condition for $r$ on $\Si_*$, we infer
\beaa
\big|{\nu^{k}(b_*)}_{|_{SP}}\big| &\les& \ep_0  u^{-1-\dec}, \quad 1\leq k\le k_{large}+7. 
\eeaa
Decomposing 
\beaa
\nu^{k}(b_*)=\ov{\nu^{k}(b_*)}+\Big(\nu^{k}(b_*)-\ov{\nu^{k}(b_*)}\Big),
\eeaa
restricting to $SP$, using the above bound, and using the fact that $\ov{\nu^{k}(b_*)}$ is constant on the spheres $S$, we infer
\beaa
|\ov{\nu^{k}(b_*)}| &\les& \ep_0  u^{-1-\dec}+\|   \nu^k(b_*)-\ov{\nu^k(b_*)} \|_{L^\infty(S)}, \quad 1\leq k\le k_{large}+7.
\eeaa
Together with the last bound of Step 18, we infer
\beaa
 \left( \int_{u=1}^{u_*}   \big|\ov{\nu^{k}(b_*)}\big|^2 \right)^{1/2}    \les  \ep_0 +\Rkstarr_{k}, \qquad 1\le  k\le k_{large}+7.
\eeaa
Plugging back in the bounds of Step 18, we deduce
 \beaa
r^{-1}\|\dkb^{\leq 2}\nu^k(b_*)\big)\|_{L^2(\Si_*)} \les \ep_0+ \Rkstarr_k, \qquad  1\leq k\le  k_{large}+7.
\eeaa
Together with the dominance in $r$, this yields
 \beaa
r^{-1}\|\dkb^{\leq 2}\nu^k(\widecheck{b_*})\big)\|_{L^2(\Si_*)} \les \ep_0+ \Rkstarr_k, \qquad  1\leq k\le  k_{large}+7.
\eeaa
Since the case $k=0$ follows from Theorem M3, we thus have
 \beaa
r^{-1}\|\dkb^{\leq 2}\nu^k(\widecheck{b_*})\big)\|_{L^2(\Si_*)} \les \ep_0+ \Rkstarr_k, \qquad  k\le  k_{large}+7.
\eeaa
Plugging back in the bounds of Step 18, we deduce
 \bea
r^{-1}\|\dkb^{\leq 2}\dk_*^{k}\widecheck{b_*}\|_{L^2(\Si_*)} &\les& \ep_0+ \Rkstarr_k, \qquad  k\le  k_{large}+7.
\eea

{\bf Step 20.} We make use of  Lemma \ref{Lemma:nuSof-integrals-again}, the definition of $\nu= e_4+b_* e_4$, and the GCM condition for $\trch$ on $\Si_*$,  to deduce
\beaa
 e_3(r)+b_* = \nu(r) &=& \frac{r e_3(u)}{2}\ov{\frac{1}{e_3(u)}(\trchb+b_*\trch)}\\
 &=&  \frac{r e_3(u)}{2}\ov{\frac{1}{e_3(u)}\left(-\frac{2\Up}{r}+\trchbc+\frac{2}{r}b_*\right)}\\
 &=& -2e_3(u)\ov{\frac{1}{e_3(u)}} +  \frac{r e_3(u)}{2}\ov{\frac{1}{e_3(u)}\left(\trchbc+\frac{2}{r}\widecheck{b_*}\right)}
\eeaa
and hence
\beaa
e_3(r) &=& -\Up+2\left(1 -e_3(u)\ov{\frac{1}{e_3(u)}}\right) -\widecheck{b_*} +  \frac{r e_3(u)}{2}\ov{\frac{1}{e_3(u)}\left(\trchbc+\frac{2}{r}\widecheck{b_*}\right)}.
\eeaa
Together with the control of $\trchbc$ in Step 7, the control of $\widecheck{b_*}$ in Step 19, and the dominant condition for $r$ on $\Si_*$, we infer, for $1\le  k\le k_{large}+7$,
\beaa
 \left( \int_{u=1}^{u_*}   \big|\ov{\nu^{k}(e_3(r))}\big|^2 \right)^{\frac{1}{2}} &\les& \left( \int_{u=1}^{u_*}   \left|\ov{\nu^{k}\left(1 -e_3(u)\ov{\frac{1}{e_3(u)}}\right)}\right|^2 \right)^{\frac{1}{2}} +\ep_0 +\Rkstarr_{k}.
\eeaa
Since $e_3(u)=2+\widecheck{e_3(u)}$, using \eqref{eq:localbootassforframeSistaronSistar:chap9}, we obtain, for $1\le  k\le k_{large}+7$, 
\beaa
 \left( \int_{u=1}^{u_*}   \big|\ov{\nu^{k}(e_3(r))}\big|^2 \right)^{\frac{1}{2}} &\les& \left( \int_{u=1}^{u_*}   \left|\ov{\nu^{k}\left(\widecheck{e_3(u)}-\ov{\widecheck{e_3(u)}}\right)}\right|^2 \right)^{\frac{1}{2}} +\ep_0 +\Rkstarr_{k}.
\eeaa
Together with Corollary \ref{Corr:nuSof integrals}, we infer, for $1\le  k\le k_{large}+7$, 
\beaa
 \left( \int_{u=1}^{u_*}   \big|\ov{\nu^{k}(e_3(r))}\big|^2 \right)^{\frac{1}{2}} &\les& \left( \int_{u=1}^{u_*}   \left\|\nu^{k}\left(\widecheck{e_3(u)}\right)-\ov{\nu^k\left(\widecheck{e_3(u)}\right)}\right\|^2_{L^\infty(S)} \right)^{\frac{1}{2}} +\ep_0 +\Rkstarr_{k},
\eeaa
and since $k\geq 1$ and $e_3(u)=2+\widecheck{e_3(u)}$, this yields, for $1\le  k\le k_{large}+7$, 
\beaa
 \left( \int_{u=1}^{u_*}   \big|\ov{\nu^{k}(e_3(r))}\big|^2 \right)^{\frac{1}{2}} &\les& \left( \int_{u=1}^{u_*}   \left\|\nu^{k}\left(e_3(u)\right)-\ov{\nu^k\left(e_3(u)\right)}\right\|^2_{L^\infty(S)} \right)^{\frac{1}{2}} +\ep_0 +\Rkstarr_{k}.
\eeaa
Together with the control of $\nu^{k}e_3(u)- \ov{\nu^ke_3(u)}$ in Step 17, we infer, for $1\le  k\le k_{large}+7$, 
\beaa
 \left( \int_{u=1}^{u_*}   \big|\ov{\nu^{k}(e_3(r))}\big|^2 \right)^{\frac{1}{2}} &\les& \ep_0 +\Rkstarr_{k}.
\eeaa
In view of the control of $e_3(r)$ in Step 16, this yields
 \beaa
r^{-1}\|\dkb^{\leq 2}\big(\nu^k(e_3(r))\big)\|_{L^2(\Si_*)} \les \ep_0+ \Rkstarr_k, \qquad 1\leq k\le  k_{large}+7.
\eeaa
Since the case $k=0$ is treated by Theorem M3, since $e_3(r)=-\Up+\widecheck{e_3(r)}$, and using the dominant condition on $r$, we obtain 
 \beaa
r^{-1}\|\dkb^{\leq 2}\big(\nu^k\widecheck{e_3(r)}\big)\|_{L^2(\Si_*)} \les \ep_0+\frac{1}{r}+ \Rkstarr_k\les \ep_0+ \Rkstarr_k, \qquad k\le  k_{large}+7.
\eeaa
Using again  the control of $e_3(r)$ in Step 16, we deduce 
  \beaa
r^{-1}\|\dkb^{\leq 2}\dk_*^{k}\widecheck{e_3(r)}\big)\|_{L^2(\Si_*)} &\les& \ep_0+ \Rkstarr_k, \qquad  k\le  k_{large}+7
\eeaa
Finally, since $\nu(r)=e_3(r)+b_*$, and since $\nu(u+r)=0$ along $\Si_*$, we infer, together with the control of $b_*$ in Step 19, 
  \bea
r^{-1}\left\|\dkb^{\leq 2}\dk_*^{k}\big(\widecheck{\nu(r)}, \, \widecheck{\nu(u)}\big)\right\|_{L^2(\Si_*)} &\les& \ep_0+ \Rkstarr_k, \qquad  k\le  k_{large}+7
\eea

{\bf Step 21.} It remains to control $\ombc$. Recall from Proposition \ref{Prop.NullStr+Bianchi-lastslice} the following linearized null structure equation
\beaa
 \nab_3 \ze -\frac{\Up}{r} \ze &=-\bb-2 \nab \ombc +\frac{\Up}{r} (\eta+\ze)+\frac{1}{r}\xib+ \frac{2m}{r^2}(\ze-\eta) +\Ga_b\c \Ga_b
\eeaa
The control of $\ze$ in Step 10, of $\eta$ in Step 13, of $\xib$ in Step 14, and \eqref{eq:localbootassforframeSistaronSistar:chap9} implies
  \beaa
\left\|\dkb\dk_*^{k}\ombc\right\|_{L^2(\Si_*)} &\les& \ep_0+ \Rkstarr_k, \qquad  k\le  k_{large}+7.
\eeaa
Also, using the null structure equations for $\nab_3\trchc$ and $\nab_4\trchc$, and the fact that $\nu=e_3+b_*e_4$, we have 
\beaa
\nab_{\nu}\trchc &=&  2   \div \eta + 2\rhoc -\frac{1}{r}\kabc + \frac{4}{r} \ombc +\frac{2}{r^2}\widecheck{e_3(r)} +\Ga_b\c\Ga_g.
\eeaa
Together with the control of $\trchc$ in Step 1,  of $\trchbc$ in Step 7, of $\eta$ in Step 13, of $\widecheck{e_3(r)}$ in Step 20, and \eqref{eq:localbootassforframeSistaronSistar:chap9}, we infer
  \beaa
\left\|\nu^{k}\ombc\right\|_{L^2(\Si_*)} &\les& \ep_0+ \Rkstarr_k, \qquad  k\le  k_{large}+7.
\eeaa
In view of the above estimate for $\dkb\dk_*^{k}\ombc$, this implies 
 \bea
\left\|\dkb^{\leq 1}\dk_*^{k}\ombc\right\|_{L^2(\Si_*)} &\les& \ep_0+ \Rkstarr_k, \qquad  k\le  k_{large}+7.
\eea

In view of Steps 1--21, and in view of the definition of the norm $ \Skstarr_k$, we have obtained 
\beaa
 \Skstarr_k&\les& \ep_0+\Rkstarr_k, \qquad k\le k_{large}+7,
\eeaa
which is an improvement of the  local bootstrap assumptions \eqref{eq:localbootassforframeSistaronSistar:chap9}. 
This concludes the proof of Proposition \ref{prop:controlofRiccicoefficientsonSigmastar}.
\end{proof}

%%%%%%%%%%%%%%%%%%%%%%%%%%%%%%%%%%%%%%%%%%%%%%%%%%

\subsection{Control of $J^{(0)}$, $f_0$ and $\Jk$ on $\Si_*$}
\lab{sec:controlofthe1formf0onSigramstar:cshap9}

%%%%%%%%%%%%%%%%%%%%%%%%%%%%%%%%%%%%%%%%%%%%%%%%%%

In view of section \ref{sec:admissibleGMCPGdatasetonSigmastar}, $\th$  initialized  on $S_*$ as in  section \ref{section:canonical-coord-BPG},  and propagated to $\Si_*$ by $\nu(\th)=0$. Recall also that 
\beaa
J^{(0)}=\cos\th
\eeaa
so that $\nu(J^{(0)})=0$ along $\Si_*$.

Moreover, the tangential 1-form $f_0$ on $\Si_*$ given by,  see Definition \ref{def:definitionoff0fplusfminus},
\beaa
(f_0)_1=0, \quad (f_0)_2=\sin\th, \quad \textrm{on}\quad S_*, \qquad \nab_\nu f_0=0,
\eeaa
where, on $S_*$, we consider the orthonormal basis $(e_1, e_2)$ of $S_*$ given by \eqref{eq:specialorthonormalbasisofSstar}.  

Finally, note that the complex horizontal 1-form $\Jk$ introduced in Definition \ref{def:definitionofJkonMext} verifies
\beaa
\Jk = \frac{1}{|q|}\left(f_0+i\dual f_0\right)\quad \textrm{on}\quad\Si_*.
\eeaa

Recall from section \ref{sec:defintionofthePTstructuresinMM} 
that we  initialize the PT frame of $\Mext$ from the integrable frame  on $\Si_*$, by  relying on the   change of frame formula 
with  the transition  coefficients 
\beaa
\la=1, \qquad f=\frac{a}{r}f_0, \qquad  \fb=\frac{a\Upsilon}{r}f_0.
\eeaa
Thus, in order to control the PT frame of $\Mext$ on $\Si_*$ in the next section, we first need to control the 1-form $f_0$. We also need to control $J^{(0)}=\cos\th$ and $\Jk$ as these quantities are involved in the definition of the linearized quantities corresponding to the PT frame of $\Mext$. The following lemma provides the control of $J^{(0)}$, $f_0$ and $\Jk$ on $\Si_*$.  
 
 \begin{lemma}
\lab{lemma:controlofthe1formf0onSigmastar}
The following estimates holds true on $\Si_*$, for all $k\le k_{large}+7 $,
\bea
\nn \|\dk_*^{\leq k}\widecheck{\nab J^{(0)}}\|_{L^2(\Si_*)}+\|\dk_*^{\leq k}\div(f_0)\|_{L^2(\Si_*)}\\
\nn +\|\dk_*^{\leq k}\widecheck{\curl(f_0)}\|_{L^2(\Si_*)}+\|\dk_*^{\leq k}\nab\hot f_0\|_{L^2(\Si_*)}\\
+\|r\dk_*^{\leq k}\widecheck{\ov{\DD}\c\Jk}\|_{L^2(\Si_*)}+\|r\dk_*^{\leq k}\DD\hot\Jk\|_{L^2(\Si_*)}  &\les& \ep_0 +      \Rkstarr_k,
\eea
where, see Definitions \ref{def:renormalizationforJpbasisell=1modes} and \ref{def:renormalizationforf0fpfm},   
\bea
\bsplit
\widecheck{\nab J^{(0)}}&:=\nab J^{(0)}+\frac{1}{r}\dual f_0, \qquad \widecheck{\curl(f_0)} := \curl(f_0)-\frac{2}{r}\cos\th,\\
\widecheck{\ov{\DD}\c\Jk} &:= \ov{\DD}\c\Jk-\frac{4i(r^2+a^2)\cos\th}{|q|^4}.
\end{split}
\eea
\end{lemma}
 
\begin{proof}
We first derive  estimates on $S_*$. We make use of 
Lemmas \ref{lemma:computationforfirstorderderivativesJpell=1basis} and \ref{lemma:computationforfirstorderderivativesf0fpfm} according to which we have on $S_*$
\beaa
\bsplit
\widecheck{\nab J^{(0)}} &= -\frac{1}{r}(e^{-\phi}-1)\dual f_0, \qquad \div(f_0) =f_0\c\nab\phi,\qquad \curl(f_0) = \frac{2}{re^{\phi}}\cos\th-f_0\wedge\nab\phi,\\ 
\nab\hot f_0 &= \left(\ba{cc}
 f_0\c\nab\phi & f_0\wedge\nab\phi\\
f_0\wedge\nab\phi & - f_0\c\nab\phi
\ea\right).
\end{split}
\eeaa
 We also use the following estimate on $S_*$
\beaa
 \|\phi\|_{\hk_{k+1}(S_*)}\les  \ep_0 +\Rkstarr_{k}, \qquad  \mbox{for all}\,\, k\le k_{large} +7,
\eeaa
derived in Step 2 of the proof of  Proposition \ref{prop:controlofRiccicoefficientsonSigmastar}. The above control of $\phi$ on $S_*$ and the above identities on $S_*$ implies immediately 
\beaa
r\left\|\Big(\widecheck{\nab J^{(0)}}, \,\div(f_0), \, \widecheck{\curl(f_0)}, \,\nab\hot f_0\Big)\right\|_{\hk_{k}(S_*)}\les  \ep_0 +\Rkstarr_{k}, \qquad  \mbox{for all}\,\, k\le k_{large} +7,
\eeaa
which is equivalent to 
\beaa
\left\|r\widecheck{\nab J^{(0)}}\right\|_{\hk_{k}(S_*)}+\left\|r\nab f_0-\cos\th\in\right\|_{\hk_{k}(S_*)}\les  \ep_0 +\Rkstarr_{k}, \qquad  \mbox{for all}\,\, k\le k_{large} +7.
\eeaa

 Next,  using   $\nu(J^{(0)})=0$ and $\nab_\nu f_0=0$, we have along $\Si_*$, see  Lemma  \ref{lemma:transporteqnuofnabofJpf0fpfm}, 
 \beaa
 \nab_\nu\Big[r\widecheck{\nab J^{(0)}}\Big]=\Ga_b \c \dkb^{\leq 1}J^{(0)}, \qquad \nab_\nu\Big[r\nab f_0-\cos\th\in\Big] &=& \Ga_b\dkb^{\leq 1}f_0.
\eeaa
 Commuting with $\dkb$ and integrating from $S_*$,  where $r\widecheck{\nab J^{(0)}}$ and $r\nab f_0-\cos\th\in$ are   under control in view of the above,  and making use   of Proposition \ref{prop:controlofRiccicoefficientsonSigmastar} to control $\Ga_b$, we deduce on $\Si_*$
 \beaa
 &&\sup_{S\subset\Si_*}\left(\left\|\dkb^{\leq k}(r\widecheck{\nab J^{(0)}})\right\|_{L^2(S)}+ \|\dkb^{\leq k}(r\nab f_0-\cos\th\in)\|_{L^2(S)}\right)\\
 &\les&\ep_0+ \Rkstarr_{k}+\int_1^{u_*}\|\dkb^{\leq k}(\Ga_b\dkb^{\leq 1}f_0)\|_{L^2(S)}\\
&\les&\ep_0+ \Rkstarr_{k}+\sqrt{u_*}\|\dkb^{\leq k}(\Ga_b\dkb^{\leq 1}f_0)\|_{L^2(\Si_*)}\\
&\les& \sqrt{u_*}\Big(\ep_0+ \Rkstarr_{k}\Big).
\eeaa
Thus, using the dominant condition for $r$ on $\Si_*$, we infer, for $k\le k_{large} +7$,
 \beaa
 \sup_{S\subset\Si_*}\left(\left\|\dkb^{\leq k}\widecheck{\nab J^{(0)}}\right\|_{L^2(S)}+ \left\|\dkb^{\leq k}\left(\nab f_0-\frac{\cos\th}{r}\in\right)\right\|_{L^2(S)}\right) &\les& u_*^{-\frac{1}{2}-\dec} \big(\ep_0 +      \Rkstarr_k\big).
 \eeaa
As a consequence, since $\de_{dec}>0$, we have
\beaa
 \left\|\dkb^{\leq k}\widecheck{\nab J^{(0)}}\right\|_{L^2(\Si_*)} +\left\|\dkb^{\leq k}\left(\nab f_0-\frac{\cos\th}{r}\in\right)\right\|_{L^2(\Si_*)} &\les& \ep_0 +      \Rkstarr_k, \quad k\le k_{large} +7.
 \eeaa
 Using again $\nab_\nu[r\widecheck{\nab J^{(0)}}]=\Ga_b \c \dkb^{\leq 1}J^{(0)}$ and  $\nab_\nu[r\nab f_0-\cos\th\in] = \Ga_b\dkb^{\leq 1}f_0$, together with the control of $\Ga_b$ provided by Proposition \ref{prop:controlofRiccicoefficientsonSigmastar}, we may recover the $\nab_\nu$ derivatives and obtain 
\beaa
\left\|\dk_*^{\leq k}\widecheck{\nab J^{(0)}}\right\|_{L^2(\Si_*)}+ \left\|\dk_*^{\leq k}\left(\nab f_0-\frac{\cos\th}{r}\in\right)\right\|_{L^2(\Si_*)} &\les& \ep_0 +      \Rkstarr_k, \quad k\le k_{large} +7,
 \eeaa 
 or equivalently, for  $k\le k_{large} +7$, 
 \beaa
\nn \|\dk_*^{\leq k}\widecheck{\nab J^{(0)}}\|_{L^2(\Si_*)}+\|\dk_*^{\leq k}\div(f_0)\|_{L^2(\Si_*)}\\
+\|\dk_*^{\leq k}\widecheck{\curl(f_0)}\|_{L^2(\Si_*)}+\|\dk_*^{\leq k}\nab\hot f_0\|_{L^2(\Si_*)}  &\les& \ep_0 +      \Rkstarr_k
\eeaa
as stated.

Finally, we control $\Jk$. Recall that we have on $\Si_*$
\beaa
\Jk = \frac{1}{|q|}\left(f_0+i\dual f_0\right)=\frac{1}{r}\left(f_0+i\dual f_0\right)+O(r^{-3})\left(f_0+i\dual f_0\right).
\eeaa
We infer on $\Si_*$
\beaa
\ov{\DD}\c\Jk &=& \frac{1}{r}(\nab-i\dual\nab)\c\left(f_0+i\dual f_0\right)+O(r^{-4}) = \frac{2}{r}\left(\div(f_0)+i\curl(f_0)\right)\\
&=& \frac{4i\cos\th}{r^2}+\frac{2}{r}\left(\div(f_0)+i\widecheck{\curl(f_0)}\right)+O(r^{-4})\\
&=& \frac{4i(r^2+a^2)\cos\th}{|q|^4}+\frac{2}{r}\left(\div(f_0)+i\widecheck{\curl(f_0)}\right)+O(r^{-4}),
\eeaa
and hence 
\beaa
\widecheck{\ov{\DD}\c\Jk} &=& \frac{2}{r}\left(\div(f_0)+i\widecheck{\curl(f_0)}\right)+O(r^{-4}),
\eeaa
as well as
\beaa
\DD\hot\Jk &=& \frac{1}{r}(\nab+i\dual\nab)\hot\left(f_0+i\dual f_0\right)+O(r^{-4})\\
&=&  \frac{2}{r}\left(\nab\hot f_0+i\dual(\nab\hot f_0)\right)+O(r^{-4}).
\eeaa
Together with the above control of $\div(f_0)$, $\widecheck{\curl(f_0)}$ and $\nab\hot f_0$, and the dominant condition for $r$ on $\Si_*$, we immediately infer, for all $k\le k_{large}+7 $,
\beaa
\|r\dk_*^{\leq k}\widecheck{\ov{\DD}\c\Jk}\|_{L^2(\Si_*)}+\|r\dk_*^{\leq k}\DD\hot\Jk\|_{L^2(\Si_*)}  \les \ep_0 +      \Rkstarr_k+r^{-2}\sqrt{u_*} &\les&  \ep_0 +      \Rkstarr_k.
\eeaa
This concludes the proof of Lemma \ref{lemma:controlofthe1formf0onSigmastar}.
\end{proof}

For convenience, we introduce the following notation. 
\begin{definition}\lab{def:Gammabtilde:chap9}
We denote by ${}^{(\Si_*)}\Ga_b$ and ${}^{(\Si_*)}\Ga_g$ the set of linearized quantities below
\beaa
\bsplit
{}^{(\Si_*)}\Ga_b &:= {}^{(\Si_*)}\Ga_{b,1}\cup {}^{(\Si_*)}\Ga_{b,2},\\
{}^{(\Si_*)}\Ga_{b,1} &:= \Big\{\eta, \,\, \chibh,\,  \ombc, \,\, \xib, r^{-1}\widecheck{\nu(r)},\,\, r^{-1}\widecheck{\nu(u)},\,\, r^{-1}\widecheck{b_*}\Big\},\\
{}^{(\Si_*)}\Ga_{b,2} &:= \Big\{ \widecheck{\nab J^{(0)}},\,\, \div(f_0),\,\, \widecheck{\curl(f_0)},\,\, \nab\hot f_0,\,\, r\widecheck{\ov{\DD}\c\Jk}, \,\, \DD\hot\Jk\Big\},\\
{}^{(\Si_*)}\Ga_g &:= \Big\{\chih, \,\,\trchc, \,\,\ze, \,\,\trchbc\Big\},
\end{split}
\eeaa
where $\trchc$, $\chih$, $\ze$, $\eta$, $\trchbc$, $\chibh$, $\ombc$, $\xib$ are the Ricci coefficients   of the integrable frame  of  $\Si_*$. 
\end{definition}

\begin{corollary}\lab{cor:estimatesforf0fplusfminusandJp:Gammatildeb:chap9}
We have on $\Si_*$, for $k\leq k_{large}+7$,
\beaa
\|\dk_*^{\leq k}({}^{(\Si_*)}\Ga_b)\|_{L^2(\Si_*)}+\|r\dk_*^{\leq k}({}^{(\Si_*)}\Ga_g)\|_{L^2(\Si_*)} &\les& \ep_0+\Rkstarr_k.
\eeaa
\end{corollary}

\begin{proof}
This is an immediate consequence of Proposition 
\ref{prop:controlofRiccicoefficientsonSigmastar} and Lemma \ref{lemma:controlofthe1formf0onSigmastar}.
\end{proof}

%%%%%%%%%%%%%%%%%%%%%%%%%%%%%%%%%%%%%%%%%%%%%%%%%%%%%

\subsection{Proof of Proposition \ref{Prop:controlGa-PTframe-Si_*}}
\lab{sec:proofofProp:controlGa-PTframe-Si_*}

%%%%%%%%%%%%%%%%%%%%%%%%%%%%%%%%%%%%%%%%%%%%%%%%%%%%%

Throughout this section, to avoid confusion between the integrable frame of $\Si_*$ and the outgoing PT frame of $\Mext$:
\begin{itemize}
\item the frame and quantities associated to the integrable frame of $\Si_*$ are denoted without prime,

\item the frame and quantities associated to the outgoing PT frame of $\Mext$ are denoted with prime.
\end{itemize}

Recall that we have in the outgoing PT frame on $\Mext$
\beaa
\xi'=0,\qquad \om'=0, \qquad \Hb'=-\frac{a\ov{q'}}{|q'|^2}\Jk', \qquad e_4'(r)=1, \qquad e_4'(u)=e_4'(\th)=0,
\eeaa  
while the integrable frame of $\Si_*$ satisfies the  transversality conditions \eqref{eq:transversalityconditionsonSigmastarcompatiblewithanoutgoingfoliation}, i.e.
\beaa
\xi=0, \qquad \om=0, \qquad \etab=-\ze, \qquad e_4(r)=1, \qquad e_4(u)=0.
\eeaa
Recall also that $\Jk'$ satisfies on $\Mext$
\beaa
 \nab_4'\Jk'=-\frac{1}{q'}\Jk'.
\eeaa
We denote: 
\begin{itemize}
\item by  $(f, \fb, \la)$  the transition coefficients from the  integrable frame  $E$ of $\Si_*$ to the outgoing PT frame $E'$ of $\Mext$,

\item by  $(f', \fb', \la')$  the transition coefficients of the reverse transformation, i.e.  the transition coefficients from  the frame $E'$ to the frame  $E$,
\end{itemize}
and recall from section \ref{sec:defintionofthePTstructuresinMM} 
that the  initialization of the PT frame of $\Mext$ given by the following choice for $(f, \fb, \la)$ on $\Si_*$
\beaa
\la=1, \qquad f=\frac{a}{r}f_0, \qquad  \fb=\frac{a\Upsilon}{r}f_0,
\eeaa
and that $r'$, $u'$, $\th'$ and $\Jk'$ are initialized on $\Si_*$ by
\beaa
r'=r, \qquad u'=u, \qquad \th'=\th, \qquad \Jk'=\Jk.
\eeaa

\begin{remark}\lab{rmk:equivalenceRkstarRkstarr}
We  note that, in view of the transformation formulas of Proposition \ref{Proposition:transformationRicci} in the particular case of the transformation from  the integrable frame  of $\Si_*$ to the outgoing PT frame of $\Mext$,  and in view of the form of the transition coefficients $(f, \fb, \la)$ on $\Si_*$ recalled above, the curvature norm $\Rkstarr_k$ of Definition \ref{definition:int-normson-RcSi_*} and the curvature  norm $\Rkstar_k$ of Definition \ref{definition:PT-normson-RcSi_*} are equivalent. 
\end{remark}

Our bootstrap assumption {\bf  BA-PT}, see \eqref{eq:mainbootassforchapte9}, implies in particular  for the Ricci and metric coefficients of the outgoing PT frame of $\Mext$ on $\Si_*$
\bea
\lab{Prop-bootstrap:finalcontrolofthePTframeonSigmastar}
\Skstar_k \leq \ep, \qquad k\leq k_{large}+7.
\eea

To control the PT frame on $\Si_*$, we will   rely on the reverse transformation, i.e. from the PT frame of $\Mext$ to the integrable frame of $\Si_*$. This is done in the following  lemma.

\begin{lemma}\lab{lemma:controloflapfpfbponSigmastar}
Consider the change of frame coefficients $(f', \fb', \la')$  from the outgoing PT frame of $\Mext$ to the integrable frame of $\Si_*$. Then, we have\footnote{See Definition \ref{def:ordermagnitude:chap9} for our notation  $O(r^{-p})$.} on $\Si_*$
\beaa
\la' =  1+O(r^{-2}),\qquad f'  = -\frac{a}{r}\Big(1+O(r^{-2})\Big)f_0,\qquad \fb' = -\frac{a\Upsilon}{r}\Big(1+O(r^{-2})\Big)f_0.
\eeaa
\end{lemma}

\begin{proof}
Recall that the change of frame coefficients from  the integrable frame of $\Si_*$ to the PT frame of $\Mext$ are given by
\beaa
\la=1, \qquad f=\frac{a}{r}f_0, \qquad  \fb=\frac{a\Upsilon}{r}f_0.
\eeaa
In view of   equation \eqref{relations:laffb-to-primes}  of Lemma \ref{Lemma:Generalframetransf}, we have
\beaa
\bsplit
\la' &= \la^{-1} \left(1+\frac{1}{2}f\c\fb  +\frac{1}{16} |f|^2  |\fb|^2\right),\\
f_a'  &= -\frac{\la}{1+\frac{1}{2}f\c\fb  +\frac{1}{16} |f|^2  |\fb|^2}\left(f_a +\frac{1}{4}|f|^2\fb_a\right),\\
\fb_a' &= -\la^{-1}\left(\fb_a+\frac 1 4 |\fb|^2f_a\right),
\end{split}
\eeaa
 from which the lemma easily follows.
\end{proof}

\begin{lemma}\lab{lemma:controloflapfpfbponSigmastar:derivatives}
Let
\beaa
\bsplit
\widecheck{\curl(f')} &:= \curl(f')+\frac{2a}{r^2}\cos\th, \qquad \widecheck{\curl(\fb')} := \curl(f')+\frac{2a\Upsilon}{r^2}\cos\th,\\
\widecheck{\nab_\nu f'} &:= \nab_\nu f' + \frac{2a}{r^2}f_0,\qquad\qquad\qquad\,  \widecheck{\nab_\nu\fb'} := \nab_\nu\fb' +\frac{2a}{r^2}f_0.
\end{split}
\eeaa
Then, we have, for all $k\le k_{large}+7$, 
\beaa
 \|r \dk_*^{\leq k+1}\log(\la')\|_{L^2(\Si_*)}+ \|r\dk_*^{\leq k}\div(f')\|_{L^2(\Si_*)}+ \|r\dk_*^{\leq k}\widecheck{\curl(f')}\|_{L^2(\Si_*)}\\
+ \|r\dk_*^{\leq k}\nab\hot f'\|_{L^2(\Si_*)}+\|r\dk_*^{\leq k}\widecheck{\nab_\nu f'}\|_{L^2(\Si_*)} +\|r\dk_*^{\leq k}\div(\fb')\|_{L^2(\Si_*)} \\
+\|r\dk_*^{\leq k}\widecheck{\curl(\fb')}\|_{L^2(S_*)}+ \|r\dk_*^{\leq k}\nab\hot \fb'\|_{L^2(\Si_*)} +\|r\dk_*^{\leq k}\widecheck{\nab_\nu\fb'}\|_{L^2(\Si_*)} &\les& \ep_0+\Rkstar_k.
\eeaa
\end{lemma}

\begin{proof}
Recall that we have
\beaa
\la' =  1+O(r^{-2}),\qquad f'  = -\frac{a}{r}\Big(1+O(r^{-2})\Big)f_0,\qquad \fb' = -\frac{a\Upsilon}{r}\Big(1+O(r^{-2})\Big)f_0.
\eeaa
We infer, using in particular $\nab(r)=0$, 
\beaa
\bsplit
&\nab\la' =O(r^{-3}), \quad \div(f') = -\frac{a}{r}\div(f_0)+O(r^{-4}), \quad \curl(f') = -\frac{a}{r}\curl(f_0)+O(r^{-4}), \\ 
&\nab\hot f' = -\frac{a}{r}\nab\hot f_0+O(r^{-4}), \quad  \div(\fb') = -\frac{a\Upsilon}{r}\div(f_0)+O(r^{-4}), \\
 &\curl(\fb') = -\frac{a\Upsilon}{r}\curl(f_0)+O(r^{-4}), \quad \nab\hot \fb' = -\frac{a\Upsilon}{r}\nab\hot f_0+O(r^{-4}).
\end{split}
\eeaa
Together with the control of $f_0$ in Lemma \ref{lemma:controlofthe1formf0onSigmastar}, the dominant condition for $r$ on $\Si_*$, the definition of $\widecheck{\curl(f')}$ and $\widecheck{\curl(\fb')}$, and Remark \ref{rmk:equivalenceRkstarRkstarr} on the equivalence between the norms $\Rkstarr_k$ and $\Rkstar_k$,  
we deduce, for all $k\leq k_{large}+7$, 
\beaa
\|r\dk_*^{\leq k}\nab(\la')\|_{L^2(\Si_*)}+\|r\dk_*^{\leq k}\div(f')\|_{L^2(\Si_*)}+\|r\dk_*^{\leq k}\widecheck{\curl(f')}\|_{L^2(\Si_*)}\\
+\|r\dk_*^{\leq k}\nab\hot f'\|_{L^2(\Si_*)}+\|r\dk_*^{\leq k}\div(\fb')\|_{L^2(\Si_*)}+\|r\dk_*^{\leq k}\widecheck{\curl(\fb')}\|_{L^2(S_*)}\\
+\|r\dk_*^{\leq k}\nab\hot \fb'\|_{L^2(\Si_*)} &\les& \ep_0+\Rkstar_k.
\eeaa

Also, using again 
\beaa
\la' =  1+O(r^{-2}),\qquad f'  = -\frac{a}{r}\Big(1+O(r^{-2})\Big)f_0,\qquad \fb' = -\frac{a\Upsilon}{r}\Big(1+O(r^{-2})\Big)f_0,
\eeaa
and since $\nab_\nu f_0=0$ and $\nu(\th)=0$ on $\Si_*$, and  $\nu(r)=-2+\widecheck{\nu(r)}$, we obtain
\beaa
&& \nu(\la') =  O(r^{-3})\left(-2+\widecheck{\nu(r)}\right),\qquad \nab_\nu f'  = \frac{a}{r^2}\Big(1+O(r^{-2})\Big)\left(-2+\widecheck{\nu(r)}\right)f_0,\\ 
&& \nab_\nu\fb' = \frac{a}{r^2}\Big(1+O(r^{-1})\Big)\left(-2+\widecheck{\nu(r)}\right)f_0,
\eeaa
and hence
\beaa
&&\nu(\la') =  O(r^{-3})\left(-2+\widecheck{\nu(r)}\right),\qquad \widecheck{\nab_\nu f'}  = \frac{a}{r^2}\widecheck{\nu(r)}f_0+O(r^{-4})\left(-2+\widecheck{\nu(r)}\right)f_0,\\
&& \widecheck{\nab_\nu\fb'} = \frac{a}{r^2}\widecheck{\nu(r)}f_0+O(r^{-3})\left(-2+\widecheck{\nu(r)}\right)f_0.
\eeaa
Together with the control of $\widecheck{\nu(r)}$ derived in Proposition \ref{prop:controlofRiccicoefficientsonSigmastar},  the control of $f_0$ in Lemma \ref{lemma:controlofthe1formf0onSigmastar},  the dominant condition for $r$ on $\Si_*$, and Remark \ref{rmk:equivalenceRkstarRkstarr} on the equivalence between the norms $\Rkstarr_k$ and $\Rkstar_k$, we infer,  for all $k\leq k_{large}+7$, 
\beaa
\| r \dk_*^{\leq k}\nu(\la')\|_{L^2(\Si_*)}+\|r\dk_*^{\leq k}\widecheck{\nab_\nu f'}\|_{L^2(\Si_*)}+\|r\dk_*^{\leq k}\widecheck{\nab_\nu\fb'}\|_{L^2(\Si_*)}
 &\les&\ep_0+ \Rkstar_k.
\eeaa
This concludes the proof of Lemma \ref{lemma:controloflapfpfbponSigmastar:derivatives}.
\end{proof}

%%%%%%%%%%%%%%%%%%%%%%%%%%%%%%%%%%%%%

\subsubsection{Proof of Proposition \ref{Prop:controlGa-PTframe-Si_*}}

%%%%%%%%%%%%%%%%%%%%%%%%%%%%%%%%%%%%%

We are now ready to prove Proposition \ref{Prop:controlGa-PTframe-Si_*}. 

{\bf Step 1.}    We start with  the following lemma.
\begin{lemma}\lab{lemma:controlofderivativescoordinatesPTonSistar:chap9} 
For all $k\le k_{large}+7$, we have 
\beaa
\|\dk_*^{\leq k}\DD'(r')\|_{L^2(\Si_*)} + \|\dk_*^{\leq k}\widecheck{\DD'(u')}\|_{L^2(\Si_*)} + \|\dk_*^{\leq k}\widecheck{\DD'(\cos(\th'))}\|_{L^2(\Si_*)} \\
+\|r^{-1}\dk_*^{\leq k}\widecheck{e_3'(r')}\|_{L^2(\Si_*)}+ \|r^{-1}\dk_*^{\leq k}\widecheck{e_3'(u')}\|_{L^2(\Si_*)} + \|\dk_*^{\leq k}e_3'(\cos(\th'))\|_{L^2(\Si_*)} &\les& \ep_0+\Rkstar_k.
\eeaa
\end{lemma}

\begin{proof}
Since $r'=r$ and $u'=u$ on $\Si_*$ by initialization, and since $\nab$ is tangent  to $\Si_*$, we infer $\nab(r')=\nab(u')=0$ on $\Si_*$. Using the change of frame transformation, and using the fact that $\la=1$, $\nab(r')=0$ and $\nab(u')=0$ on $\Si_*$, we have on $\Si_*$, 
\beaa
e_4'(r') &=& e_4(r')+\frac{1}{4}|f|^2e_3(r'),\\
e_4'(u') &=& e_4(u')+\frac{1}{4}|f|^2e_3(u'),\\
e_4'(\cos(\th')) &=& e_4(\cos(\th'))+f\c\nab\cos(\th)+\frac{1}{4}|f|^2e_3(\cos(\th')).
\eeaa
Also, since we have  $e_4'(r')=1$, $e_4'(\th')=0$, $e_4'(u')=0$ in $\Mext$, since $\nu=e_3+b_*e_4$, and since $\nu(\th')=0$ (since $\th'=\th$ on $\Si_*$, $\nu(\th)=0$ and $\nu$ is tangent to $\Si_*$), we infer
\beaa
e_4(r') &=& \frac{1-\frac{1}{4}|f|^2\nu(r)}{1-\frac{1}{4}|f|^2b_*},\\
e_4(u') &=&  -\frac{\frac{1}{4}|f|^2\nu(u)}{1-\frac{1}{4}|f|^2b_*},\\
e_4(\cos(\th')) &=& -\frac{f\c\nab\cos\th}{1-\frac{1}{4}|f|^2b_*}.
\eeaa
In view of the choice of $f$, we infer
\beaa
e_4(r') &=& 1+O(r^{-2})+O(r^{-2})(\widecheck{b_*}, \widecheck{\nu(r)}),\\
e_4(u') &=&  O(r^{-2})+O(r^{-2}) \widecheck{\nu(u)}+O(r^{-4})\widecheck{b_*},\\
e_4(\cos(\th')) &=& O(r^{-2})+O(r^{-1})\widecheck{\DD\cos\th}+O(r^{-4})\widecheck{b_*},
\eeaa
and hence, denoting for convenience ${}^{(\Si_*)}\Ga_b$ and ${}^{(\Si_*)}\Ga_g$ of Definition \ref{def:Gammabtilde:chap9} simply by $\Ga_b$ and $\Ga_g$,  we obtain 
\beaa
e_4(r')=1+O(r^{-2})+r^{-1}\Ga_b, \quad e_4(u')=O(r^{-2})+r^{-1}\Ga_b,  \quad e_4(\cos(\th'))=O(r^{-2})+r^{-1}\Ga_b.
\eeaa
Using again $\nu=e_3+b_*e_4$, and the fact that $\nu$ is tangent to $\Si_*$ so that we have $\nu(r'-r)=0$, $\nu(u'-u)=0$ and $\nu(\th')=0$, we infer that
\beaa
e_3(r')=-\Upsilon+O(r^{-2})+r\Ga_b, \quad e_3(u')=2+O(r^{-2})+r\Ga_b,  \quad e_3(\cos(\th'))=O(r^{-2})+r^{-1}\Ga_b.
\eeaa

Using the above identities on $\Si_*$ for $e_4(r')$, $e_4(u')$, $e_4(\cos(\th'))$ and $e_3(r')$, $e_3(u')$, $e_3(\cos(\th'))$, together with the change of frame transformation, the choice of $\la$, $f$ and $\fb$, and the fact that $\nab(r')=\nab(u')=0$ on $\Si_*$, we obtain
\beaa
\nab'(r') &=& \frac{1}{2}\left(\frac{a}{r}e_3(r')+\frac{a\Upsilon}{r}e_4(r')\right)f_0+O(r^{-3})= \Ga_b+O(r^{-3}),\\
\nab'(u') &=&  \frac{1}{2}\left(\frac{a}{r}e_3(u')+\frac{a\Upsilon}{r}e_4(u')\right)f_0+O(r^{-3})= \frac{a}{r}f_0+\Ga_b+O(r^{-3}),\\
\nab'(\cos(\th')) &=& \nab(\cos\th)+r^{-2}\Ga_b+O(r^{-3})=-\frac{1}{r}\dual f_0+\Ga_b+O(r^{-3}),
\eeaa
and 
\beaa
e_3'(r') &=& (1+O(r^{-2}))e_3(r')+O(r^{-2})e_4(r')= -\Upsilon+O(r^{-2})+r\Ga_b,\\
e_3'(u') &=& (1+O(r^{-2}))e_3(u')+O(r^{-2})e_4(u')= 2+O(r^{-2})+r\Ga_b,\\
e_3'(\cos(\th')) &=&  (1+O(r^{-2}))e_3(\cos(\th'))+\frac{a\Upsilon}{r}f_0\c\nab\cos(\th') +O(r^{-2})e_4(\cos(\th'))\\
&=& O(r^{-2})+r^{-1}\Ga_b.
\eeaa
Since we have by definition of $\mathfrak{J}$ on $\Si_*$
\beaa
\mathfrak{J} &=& \frac{1}{|q|}(f_0+i\dual f_0)=\frac{1}{r}(f_0+i\dual f_0)+O(r^{-3})(f_0+i\dual f_0),
\eeaa
we deduce 
\beaa
\bsplit
&\DD'(r), \, \widecheck{\DD'(u)}, \,  \widecheck{\DD'(\cos\th)}   =  O(r^{-3})+ \Ga_b, \qquad 
\widecheck{e_3'(r)}, \, \widecheck{e_3'(u)}, \,   = O(r^{-2})+r\Ga_b,\\
& e_3'(\cos\th)=O(r^{-2})+r^{-1}\Ga_b.
\end{split}
\eeaa
The proof of Lemma \ref{lemma:controlofderivativescoordinatesPTonSistar:chap9}  
 then follows from the control of $\Ga_b$ provided by Corollary \ref{cor:estimatesforf0fplusfminusandJp:Gammatildeb:chap9}, and the dominant condition for $r$ on $\Si_*$.
\end{proof}

{\bf Step 2.}  Next, we consider  the following lemma.
\begin{lemma}\lab{lemma:controlofderivativesJkPTonSistar:chap9} 
For all $k\le k_{large}+7$, we have 
\beaa
\|r\dk_*^{\leq k}\DD'\hot\Jk' \|_{L^2(\Si_*)}+\|r\dk_*^{\leq k}\widecheck{\ov{\DD'}\c\Jk'}\|_{L^2(\Si_*)}+\|r\dk_*^{\leq k}\widecheck{\nab_3'\Jk'} \|_{L^2(\Si_*)} &\les& \ep_0+\Rkstar_k.
\eeaa
\end{lemma}

\begin{proof}
The proof is in the same spirit as the one of Lemma \ref{lemma:controlofderivativescoordinatesPTonSistar:chap9}. First, since $\nab_4'\Jk'=-\frac{1}{q'}\Jk'$,  since $q'=q$ and $\Jk'=\Jk$ on $\Si_*$, and since $\nab\Jk'=\nab\Jk$ we have, using the transformation formulas
\beaa
-\frac{1}{q'}\Jk' &=& \nab_4'\Jk'=\nab_4\Jk'+f\c\nab\Jk+O(r^{-2})\nab_3\Jk'\\
&=& \nab_4\Jk'+f\c\nab\Jk+O(r^{-2})\nab_\nu\Jk+O(r^{-2})\nab_4\Jk'\\
&=& (1+O(r^{-2}))\nab_4\Jk'+O(r^{-3})+r^{-3}\Ga_b
\eeaa
where we have denoted for convenience ${}^{(\Si_*)}\Ga_b$ and ${}^{(\Si_*)}\Ga_g$ of Definition \ref{def:Gammabtilde:chap9} simply by $\Ga_b$ and $\Ga_g$. Hence, we have on $\Si_*$
\beaa
\nab_4\Jk' &=& -\frac{1}{q}\Jk+O(r^{-3})+r^{-3}\Ga_b.
\eeaa

Also, since $\nab_\nu\Jk'=\nab_\nu\Jk$ on $\Si_*$, and since $\nu=e_3+b_*e_4$, we have
\beaa
\nab_3\Jk' &=& \nab_\nu\Jk -b_*\nab_4\Jk'\\
&=& \frac{1}{r}\nab_\nu(r\Jk)  -\frac{\nu(r)}{r}\Jk +\frac{b_*}{q}\Jk+O(r^{-3})+r^{-3}\Ga_b\\
&=& \frac{1}{r}\Jk+O(r^{-3})+r^{-1}\Ga_b,
\eeaa
where we have used in particular $\nab_\nu(r\Jk)=\nab_\nu(f_0+i\dual f_0)=0$ on $\Si_*$. 

Using the above identities on $\Si_*$ for $\nab_4\Jk'$ and $\nab_3\Jk'$, together with the change of frame transformation, the choice of $\la$, $f$ and $\fb$, and the fact that $\nab(r')=\nab(u')=0$ on $\Si_*$, we obtain
\beaa
\DD'\Jk' &=& \DD\Jk+O(r^{-3})+r^{-2}\Ga_b,\\
\nab_3'\Jk' &=& \nab_3\Jk'+O(r^{-3})+r^{-2}\Ga_b=  \frac{1}{r}\Jk+O(r^{-3})+r^{-1}\Ga_b.
\eeaa
We infer
\beaa
\DD'\hot\Jk' &=& \DD\hot\Jk+O(r^{-3})+r^{-2}\Ga_b=r^{-1}\Ga_b,\\
\ov{\DD'}\c\Jk' &=& \ov{\DD}\c\Jk+O(r^{-3})+r^{-2}\Ga_b = \frac{4i(r^2+a^2)\cos\th}{|q|^4}+O(r^{-3})+r^{-1}\Ga_b,\\
\nab_3'\Jk' &=&  \frac{1}{r}\Jk+O(r^{-3})+r^{-1}\Ga_b = \frac{\De q}{|q|^4}\Jk+O(r^{-3})+r^{-1}\Ga_b.
\eeaa
and hence
\beaa
\DD'\hot\Jk' = r^{-1}\Ga_b,\qquad \widecheck{\ov{\DD'}\c\Jk'} = O(r^{-3})+r^{-1}\Ga_b,\qquad 
\widecheck{\nab_3'\Jk'} =  O(r^{-3})+r^{-1}\Ga_b.
\eeaa
The proof of Lemma \ref{lemma:controlofderivativesJkPTonSistar:chap9}   
 then follows from the control of $\Ga_b$ provided by Corollary \ref{cor:estimatesforf0fplusfminusandJp:Gammatildeb:chap9}, and the dominant condition for $r$ on $\Si_*$.
\end{proof}

{\bf Step 3.} We prove the  following Lemma
\begin{lemma}
\lab{corollary:PTfram-Si_*2}
We have, for all $k\leq k_{large}+7$, 
\beaa
\|\dk_*^{\leq k}(\Xbh', \Hc', \Xib', \ombc')\|_{L^2(\Si_*)} &\les& \ep_0+\Rkstar_k.
\eeaa
\end{lemma}

\begin{proof}
In this lemma, we use the following notation
\beaa
\Gac' &:=& \Big\{r\trXc', r\Xh',\, r\Zc',\,r\trXbc', \, \Hc', \,\Xbh',\, \omb', \, \Xib'\Big\},
\eeaa
i.e. $\Gac'$ contain all the linearized Ricci coefficients of the PT frame of $\Mext$. Also, we denote for convenience  ${}^{(\Si_*)}\Ga_b$ and ${}^{(\Si_*)}\Ga_g$ of Definition \ref{def:Gammabtilde:chap9} simply by $\Ga_b$ and $\Ga_g$. In particular, we have in view of \eqref{Prop-bootstrap:finalcontrolofthePTframeonSigmastar}  
\beaa
\|\dk_*^{\leq k}\Gac'\|_{L^2(\Si_*)} &\leq&  \ep, \qquad k\leq k_{large}+7,
\eeaa
and in view of  Corollary \ref{cor:estimatesforf0fplusfminusandJp:Gammatildeb:chap9} 
\beaa
\|r\dk_*^{\leq k}\Ga_g\|_{L^2(\Si_*)}+\|\dk_*^{\leq k}\Ga_b\|_{L^2(\Si_*)} &\les&  \ep_0+\Rkstar_k, \qquad k\leq k_{large}+7.
\eeaa

With these notations, recall from Lemma \ref{lemma:controlofRiccicoeffPGframeMextonSigmastar} the following consequence of  the change of frame formulas of  Proposition \ref{Proposition:transformationRicci}  
\beaa
\chibh &=& \chibh' +\nab\hot\fb'  +O(r^{-3})+r^{-1}\Gac'+r^{-2}\Ga_b,\\
\xib -b_*\ze &=& \xib' -b_*\ze' + \frac{1}{2}\nab_\nu\fb'  + \frac{1}{4}\trchb'\,\fb' +\frac{b_*}{4}\trchb' f'  +O(r^{-3})+r^{-1}\Gac'+r^{-2}\Ga_b,\\
\omb &=& \omb'+\frac{1}{2}\nu(\log\la')    +O(r^{-3})+r^{-1}\Gac'+r^{-2}\Ga_b,\\
\eta &=& \eta' +\frac{1}{2}\nab_\nu f'  +\frac{1}{4}\fb'\trch'  +\frac{b_*}{4}\trch' f' +O(r^{-3})+r^{-1}\Gac'+r^{-2}\Ga_b.
\eeaa
Together with the form of $(f', \fb', \la')$ in Lemma \ref{lemma:controloflapfpfbponSigmastar} and the asymptotic for large $r$ of the Kerr values of the PT frame, we infer
\beaa
\chibh' &=& \Ga_b +\nab\hot\fb'  +O(r^{-3})+r^{-1}\Gac',\\
\xib' &=&  \Ga_b + \frac{1}{2}\widecheck{\nab_\nu\fb'}   +O(r^{-3})+r^{-1}\Gac',\\
\ombc' &=& \Ga_b+\frac{1}{2}\nu(\log\la')    +O(r^{-3})+r^{-1}\Gac',\\
\etac' &=& \Ga_b +\frac{1}{2}\widecheck{\nab_\nu f'}  +O(r^{-3})+r^{-1}\Gac'.
\eeaa
Lemma \ref{corollary:PTfram-Si_*2} is then an immediate consequence of the above identities, the above control of $\Gac'$ and $\Ga_b$,  the control of $\nab\hot\fb' $, $\nab_\nu f'$, $\nab_\nu\fb'$ and $\nu(\log\la')$ provided by Lemma \ref{lemma:controloflapfpfbponSigmastar:derivatives}, and the dominant condition for $r$ on $\Si_*$. 
\end{proof}

{\bf Step 4.} We prove the following lemma.
\begin{lemma}
\lab{corollary:PTfram-Si_*1}
The  following  estimates hold true, for $k\leq k_{large}+7$, 
\bea
 \|r\dk_*^{\leq k}(\trXc', \Xh', \trXbc', \Zc')\|_{L^2(\Si_*)}  &\les  \ep_0 +\Rkstar_k.
\eea
\end{lemma}

\begin{proof}
We use again the notations $\Gac'$, $\Ga_b$ and $\Ga_g$ of Lemma \ref{corollary:PTfram-Si_*2}. Then, we have in view of Lemma \ref{lemma:controlofRiccicoeffPGframeMextonSigmastar} the following consequence of  the change of frame formulas of  Proposition \ref{Proposition:transformationRicci}  
\beaa
\trch &=& \trch'  +  \div(f') +O(r^{-3})+O(r^{-1})\widecheck{\eta}'+r^{-2}\Gac'+r^{-2}\Ga_b,\\
0 &=& \atrch'  +  \curl(f')  +O(r^{-3})+O(r^{-1})\widecheck{\eta}'+r^{-2}\Gac'+r^{-2}\Ga_b,\\
\chih &=& \chih'  +  \nab\hot f' +O(r^{-3})+O(r^{-1})\widecheck{\eta}'+r^{-2}\Gac'+r^{-2}\Ga_b,
\eeaa
\beaa
\trchb &=& \trchb' +\div(\fb')   +O(r^{-3})+O(r^{-1})\xib'+r^{-2}\Gac'+r^{-2}\Ga_b,\\
0 &=& \atrchb' +\curl(\fb') +O(r^{-3})+O(r^{-1})\xib'+r^{-2}\Gac'+r^{-2}\Ga_b,
\eeaa
and
\beaa
\ze &=& \ze' -\nab(\log\la')  -\frac{1}{4}\trchb' f'  +\frac{1}{4}\fb'\trch' +  \frac{1}{4}\fb'\div(f')  + \frac{1}{4}\dual\fb' \curl(f')\\
&&+O(r^{-1})(\ombc',\chibh')+r^{-2}\Gac' +O(r^{-3})+r^{-1}\Ga_g.
\eeaa
Together with the form of $(f', \fb', \la')$ in Lemma \ref{lemma:controloflapfpfbponSigmastar} and the asymptotic for large $r$ of the Kerr values of the PT frame, we infer
\beaa
\trchc' &=&  \Ga_g+  \div(f') +O(r^{-3})+O(r^{-1})\widecheck{\eta}'+r^{-2}\Gac',\\
 \widecheck{\atrch'} &=&  \widecheck{\curl(f')}  +O(r^{-3})+O(r^{-1})\widecheck{\eta}'+r^{-2}\Gac'+r^{-2}\Ga_b,\\
\chih' &=&  \Ga_g+ \nab\hot f' +O(r^{-3})+O(r^{-1})\widecheck{\eta}'+r^{-2}\Gac',
\eeaa
\beaa
\trchb' &=& \Ga_g +\div(\fb')   +O(r^{-3})+O(r^{-1})\xib'+r^{-2}\Gac'+r^{-2}\Ga_b,\\
\widecheck{ \atrchb'}  &=& \widecheck{\curl(\fb')} +O(r^{-3})+O(r^{-1})\xib'+r^{-2}\Gac'+r^{-2}\Ga_b,
\eeaa
and
\beaa
\zec' &=& \Ga_g -\nab(\log\la')  +  \frac{1}{4}\fb'\div(f')  + \frac{1}{4}\dual\fb' \widecheck{\curl(f')}\\
&&+O(r^{-1})(\ombc',\chibh')+r^{-2}\Gac' +O(r^{-3}).
\eeaa
Lemma \ref{corollary:PTfram-Si_*1} is then an immediate consequence of the above identities, the control of $(\chibh', \etac', \Xib', \ombc')$ in Lemma \ref{corollary:PTfram-Si_*2}, the control of $\Gac'$ and $\Ga_b$ recalled in Lemma \ref{corollary:PTfram-Si_*2},  the control of $\nab\hot\fb' $, $\nab_\nu f'$, $\nab_\nu\fb'$ and $\nu(\log\la')$ provided by Lemma \ref{lemma:controloflapfpfbponSigmastar:derivatives}, and the dominant condition for $r$ on $\Si_*$. 
 \end{proof}

The proof of Proposition \ref{Prop:controlGa-PTframe-Si_*} follows immediately from Lemmas \ref{lemma:controlofderivativescoordinatesPTonSistar:chap9}--\ref{corollary:PTfram-Si_*1} above.

%%%%%%%%%%%%%%%%%%%%%%%%%%%%%%%%%%%%%%%%%%%%%%%%%%

\section{Control of the PT-Ricci coefficients in $\Mext$}
\lab{sec:theoremM8recoverRicciawaytrapping}

%%%%%%%%%%%%%%%%%%%%%%%%%%%%%%%%%%%%%%%%%%%%%%%%%%

The goal of this section is to  provide the proof of Proposition \ref{prop:controlGaextiterationassupmtionThM8}. For convenience, we restate the result below.
\begin{proposition}
\lab{Prop:EstmatesSkext}
Relative to the PT frame of  $Mext$ we have
\bea
\Skext_k\les  \ep_0 +\Skstar_k+\Rkext_k, \qquad  k\le k_{large}+7.
\eea
\end{proposition} 

To prove Proposition \ref{Prop:EstmatesSkext}, we rely in particular on our bootstrap assumption {\bf  BA-PT}, see \eqref{eq:mainbootassforchapte9}, which implies   for the Ricci and metric coefficients of the outgoing PT frame of $\Mext$ 
\bea
\lab{Prop-bootstrap:finalcontrolofthePTframeonMext}
\Skext_k \leq \ep, \qquad k\leq k_{large}+7.
\eea
Also,  we make use of   the transport lemmas derived in the next section.

%%%%%%%%%%%%%%%%%%

\subsection{Transport lemmas}

%%%%%%%%%%%%%%%%%%%

 In what  follows  the weighted derivatives $\dkb=(r\nab) $ and $\dk=( r\nab, r\nab_4, \nab_3)$  are defined with respect to the outgoing  PT frame of $\Mext$. We revisit  Proposition \ref{Prop:transportrp-f-Decay-knorms}   and state  below its $L^2$ version using the following norms
 \beaa
  \| f\|_{2} (u,r):=\| f\|_{L^2\big(S(u,r)\big)}, \qquad \|f\|_{2,k}(u, r):= \sum_{i=0}^k \|\dk^i f\|_{2 }(u, r).
  \eeaa

\begin{proposition}
\lab{Prop:transportrp-L2-knorms-PT}
Let $U$ and $F$  anti-selfdual $k$-tensors. Assume  that  $U$ verifies one of the following equations, for a real constant $c$,
\bea
\lab{eq:maintransportMext-PT}
\nab_4 U+\frac{c}{q} U=F
\eea
or
\bea
\lab{eq:maintransportMext2-PT}
\nab_4 U+ \Re\left(\frac{c}{q}\right)  U=F.
\eea
In both cases we derive,   for  any $r_0\leq r\le r_* =r_*(u)$  at fixed $u$, with  $1\leq u\leq u_*$,  in $\Mext$
 \bea\lab{eq:integration-knorms=PT}
   r^{c-1}  \|U\|_{2,k} (u, r) &\les&    r_*  ^{c-1} \|  U \|_{2,k} (u,  r_*  )+ \int_r^{r_*}      \la ^{c-1}\|F\|_{2,k}(u, \la) d\la.
\eea
\end{proposition} 

\begin{proof}
The proof is completely analogous of the one of Proposition \ref{Prop:transportrp-f-Decay-knorms} where one simply has to replace $L^\infty(S(u,r))$ based norms by by $L^2(S(u,r))$ based norms.  
\end{proof}

\begin{corollary}
\lab{Corr:Prop:transportrp-L2-knorms-PT}
If $U$  verifies  \eqref{eq:maintransportMext-PT} or \eqref{eq:maintransportMext2-PT},  we   obtain,  
 for any $C> c- \frac 1 2 $, and for any $r_1\geq r_0$, 
\bea
\sup_{\la \ge r_1}  \int_{r=\la}  \la^ {2 (c-1)}  \big| \dk^{\le k}  U \big|^2 &\les &  \int_{\Si_*} r^{2(c-1) }  \big| \dk^{\le k}  U \big|^2 +\int_{\Mext(r\ge r_1)}  r^ {2C}   \big| \dk^{\le k}  F \big|^2. 
\eea
\end{corollary}

\begin{proof}
 Apply \eqref{eq:integration-knorms=PT}  for  $r_0\leq r_1\leq r\le r_* =r_*(u)$,    $1\leq u\leq u_*$ with $r_*(u)= c_*- u$ on $\Si_*$, take the  square and integrate in $u$ to derive
 \beaa
  \int_{r=\la}  \la^ {2 (c-1)}  \big| \dk^{\le k}  U \big|^2 &\les &  \int_{\Si_*} r^{2(c-1) }  \big| \dk^{\le k}  U \big|^2 +\int_{1 }^{u_*} du \left( \int_\la^{r_*}      {\la'}^{(c -1) }\|F\|_{2,k}(u, \la') d\la'\right)^2. 
 \eeaa
 By Cauchy-Schwartz, we have for any  $C> c- \frac 1 2 $
 \beaa
 \left(  \int_\la^{r_*}      {\la'}^{(c -1) }\|F\|_{2,k}(u, \la') d\la'\right)^2&\les&  \left( \int_\la^{r_*}      {\la'}^{2C}\|F\|_{2,k}(u, \la') d\la'\right)\left(\int_{\la}^{r_*}      {\la'}^{2 c-2-2C} d\la'\right) \\
 &\les &
   \int_\la^{r_*}      {\la'}^{2 C}\|F\|^2_{2,k}(u, \la') d\la'.
 \eeaa
 Hence, we infer for any $C> c- \frac 1 2 $, and for any $r_1\geq r_0$,
\beaa
\sup_{\la \ge r_1}  \int_{r=\la}  \la^ {2 (c-1)}  \big| \dk^{\le k}  U \big|^2 &\les &  \int_{\Si_*} r^{2(c-1) }  \big| \dk^{\le k}  U \big|^2 +\int_{\Mext(r\ge r_1)}  r^ {2C}   \big| \dk^{\le k}  F \big|^2
\eeaa
 as stated.
\end{proof}

%%%%%%%%%%%%%%%%%%%%%%%%%%%%%%%

\subsection{Proof of Proposition  \ref{Prop:EstmatesSkext}}
\lab{sec:proofofProp:EstmatesSkext}

%%%%%%%%%%%%%%%%%%%%%%%%%%%%%%%

We estimate the Ricci and metric coefficients of the outgoing PT structure of $\Mext$ in the following order 
\beaa
 \trXc,\,\,\, \Xh,\,\,\,  \DD \cos \th,\,\,\,  \Zc, \,\,\,  \DD r,\,\,\,  \Hc, \,\,\,    \widecheck{e_3(r)}, \,\,\,     \ombc, \,\,\, \DD\hot\Jk, \,\,\, \widecheck{\DD\c\ov{\Jk}},\,\,\,  e_3 (\cos \th),\,\,\,  \widecheck{e_3 \Jk},\,\,\,  \Xib,
\eeaa
making use   of the triangular structure of the linearized equations for  outgoing PT structures derived  in sections  \ref{sec:linearizedeqautionsforoutoingPTstructures:chap9} and \ref{sec:otherlinearizedequationsoutgoingPTframeMext:chap9}, the transport lemmas of the previous section, and the bootstrap assumption \eqref{Prop-bootstrap:finalcontrolofthePTframeonMext}.

{\bf Step 1.} Estimates for $\trXc$.

We  apply  Corollary  \ref{Corr:Prop:transportrp-L2-knorms-PT} to the following equation, see Proposition \ref{Prop:linearizedPTstructure1}, 
\beaa
\nab_4\trXc +\frac{2}{q}\trXc &=& \Ga_g\c\Ga_g,
\eeaa
with $c=2 $  and $C=2$ and derive, for $k\le k_{large}+7$, 
\beaa
 \sup_{\la \ge r_0}  \int_{r=\la}  \la^ {2 }  \big| \dk^{\le k}  \trXc \big|^2 &\les&   \int_{\Si_*} r^{2}  \big| \dk^{\le k}  \trXc \big|^2+ \int_{\Mext}  r^{4}\big |\dk ^{\le k} (\Ga_g\c \Ga_g) \big|^2\\
 &\les&  (\Skstar_k)^2+ \int_{\Mext}  r^{4}\big |\dk ^{\le k} (\Ga_g\c \Ga_g) \big|^2. 
\eeaa
Note that, in view of our bootstrap assumptions \eqref{Prop-bootstrap:finalcontrolofthePTframeonMext} on $\Ga_g$,  
\beaa
\big |\dk ^{\le k} (\Ga_g\c \Ga_g) \big|^2&\les&  \big |\dk ^{\le k/2} (\Ga_g)   \big|^2     \big|\dk ^{\le k} (\Ga_g)  \big|^2 \les  \ep^2 r^{-4}  \big|\dk ^{\le k} (\Ga_g)  \big|^2.
\eeaa
Hence
\beaa
 \int_{\Mext}  r^{4}\big |\dk ^{\le k} (\Ga_g\c \Ga_g) \big|^2 &\les &\ep^2 \int_{\Mext} \|\Ga_g\|_{2, k} ^2 \les\ep^2 \int_{r_0} ^{r_*}\int_{r=\la}  \big|\dk ^{\le k} (\Ga_g)  \big|^2\\
 &=&\ep^2 \int_{r_0} ^{r_*}\la^{-2}\int_{r=\la}  \la^2 \big|\dk ^{\le k} (\Ga_g)  \big|^2 \les\ep^2\sup_{\la\ge r_0}  \int_{r=\la}  \la^2 \big|\dk ^{\le k} (\Ga_g)  \big|^2\\
 &\les& \ep^2 \Sk_k^2\les \ep^4\les\ep_0^2. 
\eeaa
Thus,  we deduce
\bea
\sup_{\la \ge r_0} \left( \int_{r=\la}  \la^ {2 }  \big| \dk^{\le k}  \trXc \big|^2 \right)^{\frac{1}{2}}  &\les& \ep_0+ \Skstar_k, \qquad k\le k_{large}+7.
\eea

{\bf Step 2.}  Estimates for $\Xh$.  

We  apply  Corollary  \ref{Corr:Prop:transportrp-L2-knorms-PT} to the following equation, see Proposition \ref{Prop:linearizedPTstructure1}, 
\beaa
\nab_4\Xh+\Re\left(\frac{2}{q}\right)\Xh &=& -A+\Ga_g\c\Ga_g,
\eeaa
with $c=2 $  and $C= 3/2 +\de_B/ 2$ and derive, for $k\le k_{large}+7$, 
\beaa
 \sup_{\la \ge r_0} \left( \int_{r=\la}  \la^ {2 }  \big| \dk^{\le k}  \Xh \big|^2 \right) & \les &\ \int_{\Si_*} r^{2}  \big| \dk^{\le k}  \Xh \big|^2+ \int_{\Mext}  r^{ 3+\de_B}\big |\dk ^{\le k} A \big|^2 \\
 &&+ \int_{\Mext}  r^{ 3+\de_B}\big |\dk ^{\le k} (\Ga_g\c \Ga_g) \big|^2.
\eeaa
Recalling the definition of $\Skstar_k$ and $\Rkext_k$, and proceeding as in Step 1 with the quadratic term, we deduce
\bea
\sup_{\la \ge r_0} \left( \int_{r=\la}  \la^ {2 }  \big| \dk^{\le k}  \Xh \big|^2 \right)^{\frac{1}{2}} & \les &\ep_0+\Skstar_k+\Rkext_k, \qquad k\le k_{large}+7.
\eea

{ \bf Step 3.} Estimate for  $\DD\cos\th$.
 
We apply    Corollary  \ref{Corr:Prop:transportrp-L2-knorms-PT} to the following equation, see     Lemma  
\ref{lemma:transportequationine4fornabthetaandrande3thetaandr:linearized}, 
\beaa
\nab_4\widecheck{\DD\cos\th}+\frac{1}{q}\widecheck{\DD\cos\th} &=& O(r^{-1})\trXc+O(r^{-1})\Xh+\Ga_b\c\Ga_g,\\
\eeaa
 with $c=1$ and  $C=1$.  Using the estimate of Step 1 for $\trXc$, the estimate of Step 2 for $\Xh$,  and treating the quadratic terms as before,  we easily  derive
 \bea
\sup_{\la \ge r_0}\left( \int_{r=\la}   \big| \dk^{\le k}   \widecheck{\DD\cos\th}   \big|^2 \right)^{\frac{1}{2}} & \les &\ep_0+\Skstar_k+\Rkext_k, \qquad k\le k_{large}+7.
\eea

{\bf Step 4.}  Estimate for $\Zc$.

We apply\footnote{In order to obtain an  estimate that will be useful in Step 5, see \eqref{eq:intermediary-Step5}, we proceed differently compared to previous steps.}  Proposition \ref{Prop:transportrp-L2-knorms-PT}  to the following equation, see Proposition \ref{Prop:linearizedPTstructure1},    
\beaa
\nab_4\Zc + \frac{1}{q}\Zc &=& F,\\
F &=&   O(r^{-2})\trXc + O(r^{-2})\widehat{X}  -B+\Ga_g\c\Ga_g,
\eeaa
with $c=1$ and derive,  for $k\le k_{large}+7$, 
\beaa
\| \Zc\|_{2, k} (r, u)&\les&  \| \Zc\|_{2, k} (r_*, u)+\int_r^{r_*}  \| F\|(\la, u)\|_{2, k} (u, \la)  d\la, 
\eeaa
and hence,  for $k\le k_{large}+7$, 
\beaa
\| \Zc\|^2_{2, k} (r, u)&\les&  \| \Zc\|^2_{2, k} (r_*, u)+\left(\int_r^{r_*}  \| F\|(\la, u)\|_{2, k} (u, \la)  d\la \right)^2\\
&\les&  \| \Zc\|^2_{2, k} (r_*, u)+r^{-2-\dt}\left(\int_r^{r_*} \la^{3+\dt} \| F\|(\la, u)\|^2_{2, k} (u, \la)  d\la \right).
\eeaa
In view of the form of $F$, this yields,  for $k\le k_{large}+7$, 
\bea
\lab{eq:intermediary-Step5}
\bsplit
\| \Zc\|^2_{2, k} (r, u)&\les  \| \Zc\|^2_{2, k} (r_*, u)+  r^{-2-\de_B}  \int_r^{r_*}\Big( \la^{3+\de_B} \| B\|^2_{2, k}(u, \la)\\
&\qquad\qquad +\|(\trXc, \Xh)\|^2_{2, k}(u, \la)+\ep^2 \|\Ga_g\|^2_{2, k}(u, \la)\Big).
 \end{split}
\eea

 Multiplying by $r^2$, and integrating in $u$,  we derive,  for $k\le k_{large}+7$, 
\beaa
\bsplit
\sup_{\la \ge r_0}\left( \int_{r=\la}\la^2   \big| \dk^{\le k}\Zc\big|^2 \right)^{\frac{1}{2}} &\les  \int_{\Si_*}r^2|\dk^{\le k}\Zc|^2\\
&+   \int_{\Mext}\Big( r^{3+\dt}|\dk^{\leq k}B|^2+|\dk^{\leq k}(\trXc, \Xh)|^2+\ep^2|\dk^{\leq k}\Ga_g|^2\Big).
 \end{split}
\eeaa
Using the estimate of Step 1 for $\trXc$ and the estimate of Step 2 for $\Xh$,  we easily  derive 
\bea
\sup_{\la \ge r_0}\left( \int_{r=\la}\la^2   \big| \dk^{\le k}\Zc\big|^2 \right)^{\frac{1}{2}} &\les  \ep_0+\Skstar_k+\Rkext_k, \qquad k\le k_{large}+7.
\eea

{\bf Step 5.} Estimates for $\DD r$. 

We apply   Proposition \ref{Prop:transportrp-L2-knorms-PT}  to the following equation, see  Lemma \ref{lemma:transportequationine4fornabthetaandrande3thetaandr:linearized},    
\beaa
\nab_4\DD r+\frac{1}{q}\DD r &=& \Zc+r\Ga_g\c\Ga_g
\eeaa
with $c=1$ and deduce, for $k\le k_{large}+7$, 
\beaa
    \|\DD r\|_ {2,k} (u, r) \les     \|\DD r \| _{2,k} (u,  r_*  )+ \int_r^{r_*}      \|\Zc\|_{2,k}(u, \la) d\la+  \int_r^{r_*}  \la \|\Ga_g\c\Ga_g\|_{2,k}(u,\la)d\la.
\eeaa
Squaring and integrating in $u$, controlling the error term as before using \eqref{Prop-bootstrap:finalcontrolofthePTframeonMext},  we derive, for $k\le k_{large}+7$, 
\beaa
 \sup_{\la \ge r_0}\left( \int_{r=\la} \big| \dk^{\le k}\DD r\big|^2 \right) &\les& \Skstar^2_k +\sup_{r\geq r_0}\int_{u=1}^{u_*}  \left( \int_r^{r_*}      \|\Zc\|_{2,k}(u, \la) d\la\right)^2 + \ep_0^2. 
\eeaa
Now, recall \eqref{eq:intermediary-Step5}
\beaa
\bsplit
\| \Zc\|_{2, k} (r, u)&\les  \| \Zc\|_{2, k} (r_*, u)+  r^{-1-\frac{\de_B}{2}} \Bigg(\int_r^{r_*}\Big( \la^{3+\de_B} \| B\|^2_{2, k}(u, \la)\\
&\qquad\qquad +\|(\trXc, \Xh)\|^2_{2, k}(u, \la)+\ep^2 \|\Ga_g\|^2_{2, k}(u, \la)\Big)\Bigg)^{\frac{1}{2}}.
 \end{split}
\eeaa
Hence, we have
\beaa
\int_r^{r_*}   \| \Zc\|_{2, k} (\la,   u) d\la &\les& r_*   \| \Zc\|_{2, k} (r_*, u)+  \Bigg(\int_r^{r_*}\Big( \la^{3+\de_B} \| B\|^2_{2, k}(u, \la)\\
&&\qquad\qquad +\|(\trXc, \Xh)\|^2_{2, k}(u, \la)+\ep^2 \|\Ga_g\|^2_{2, k}(u, \la)\Big)\Bigg)^{\frac{1}{2}}
\eeaa
and thus, using the estimate of Step 1 for $\trXc$ and the estimate of Step 2 for $\Xh$, we obtain, for $k\le k_{large}+7$, 
\beaa
\int_{u=1} ^{u_*} \left(\int_r^{r_*}   \| \Zc\|_{2, k} (\la,   u) d\la \right)^2 &\les&\Skstar^2_k +\Rkext^2_k +\ep_0^2. 
\eeaa
Plugging in the above estimate for $\DD r$, we infer, for $k\le k_{large}+7$, 
\beaa
 \sup_{\la \ge r_0}\left( \int_{r=\la} \big| \dk^{\le k}\DD r\big|^2 \right) &\les& \Skstar^2_k +\sup_{r\geq r_0}\int_{u=1}^{u_*}  \left( \int_r^{r_*}      \|\Zc\|_{2,k}(u, \la) d\la\right)^2 + \ep_0^2\\
 &\les& \Skstar^2_k +\Rkext_k^2 +\ep_0^2,
\eeaa
i.e.
\bea
\sup_{\la \ge r_0}\left( \int_{r=\la} \big| \dk^{\le k}\DD r\big|^2 \right)^{\frac{1}{2}}  &\les & \ep_0 +\Skstar_k+ \Rkext_k, \qquad k\le k_{large}+7.
\eea

{\bf Step 6.}  Estimates for $\Hc$. 
 
 We  apply  Corollary  \ref{Corr:Prop:transportrp-L2-knorms-PT} to the following equation, see Proposition \ref{Prop:linearizedPTstructure1}, 
\beaa
\nab_4\Hc+ \frac{1}{\ov{q}}\Hc &=&  O(r^{-2})\ov{\trXc} +O(r^{-2})\Xh -B+\Ga_b\c\Ga_g.
 \eeaa
with $c=1$  and $C=1$ and derive, for $k\le k_{large}+7$, 
\beaa
 \sup_{\la \ge r_0} \left( \int_{r=\la}  \big| \dk^{\le k}  \Hc \big|^2 \right) & \les &\ \int_{\Si_*}  \big| \dk^{\le k}  \Hc \big|^2\\
 &&+ \int_{\Mext}  r^2\left|\dk ^{\le k}\Big(O(r^{-2})\ov{\trXc} +O(r^{-2})\Xh -B+\Ga_b\c\Ga_g\Big)\right|^2.
\eeaa
Using the estimate of Step 1 for $\trXc$ and the estimate of Step 2 for $\Xh$, and proceeding as in Step 1 with the quadratic term, we deduce
\bea
\sup_{\la \ge r_0} \left( \int_{r=\la}   \big| \dk^{\le k}  \Hc \big|^2 \right)^{\frac{1}{2}} & \les &\ep_0+\Skstar_k+\Rkext_k, \qquad k\le k_{large}+7.
\eea

 {\bf Step 7.} Estimate for $e_3(\cos \th) $. 
 
  We  apply    Proposition \ref{Prop:transportrp-L2-knorms-PT} to the following equation, see  Lemma \ref{lemma:transportequationine4fornabthetaandrande3thetaandr:linearized}, 
 \beaa
 e_4(e_3(\cos\th)) &=&  O(r^{-1})\Hc +O(r^{-2})\widecheck{\DD(\cos\th)}+\Ga_b\c\Ga_b,
 \eeaa
 with $c=0$  and deduce, for $k\le k_{large}+7$,  
 \beaa
   r^{-1}  \|e_3(\cos\th)\|_{2,k} (u, r) \les    r_*  ^{-1} \| e_3(\cos\th) \|_{2,k} (u,  r_*  )+ \int_r^{r_*}      \la ^{-1}\|F\|_{2,k}(u, \la) d\la,
\eeaa
 with 
 \beaa
 F &=& O(r^{-1})\Hc +O(r^{-2})\widecheck{\DD(\cos\th)}+\Ga_b\c\Ga_b.  
 \eeaa
  Multiplying by $r$, squaring, and integrating in $u$, we deduce, for $r\geq r_0$ and $k\leq k_{large}+7$, 
  \beaa
   \int_{u=1}^{u_*}  \|e_3(\cos\th)\|^2_{2,k} (u, r) &\les& \int_{\Si_*}\big| \dk^{\le k} e_3(\cos\th)  \big|^2+  \int_{u=1}^{u_*} \left( r\int_r^{r_*}      \la ^{-1}\|F\|_{2,k}(u, \la) d\la\right)^2\\
   &\les& (\Skstar_k)^2+  \int_{u=1}^{u_*} \left( r\int_r^{r_*}      \la ^{-1}\|F\|_{2,k}(u, \la) d\la\right)^2.
  \eeaa

  Also
  \beaa
 &&\int_r^{r_*}      \la ^{-1}\|F\|_{2,k}(u, \la) d\la \\
 &\les&  \int_r^{r_*}      \la ^{-2}\|(\Hc, \widecheck{\DD(\cos\th)})\|_{2,k}(u, \la) d\la + \ep   \int_r^{r_*}      \la ^{-2} \|\Ga_b\|_{2,k}(u, \la) d\la\\
 &\les & r^{-\frac{1}{2}} \left( \int_r^{r_*}  \la^{-2}   \|(\Hc, \widecheck{\DD(\cos\th)})\|^2 _{2,k}(u, \la) d\la \right)^{\frac{1}{2}} + \ep  r^{-\frac{1}{2}} \left( \int_r^{r_*}  \la^{-2}   \|\Ga_b\|^2 _{2,k}(u, \la) d\la \right)^{\frac{1}{2}} 
  \eeaa
  and hence
  \beaa
&&   \int_{u=1}^{u_*} \left(r\int_r^{r_*}      \la ^{-1}\|F\|_{2,k}(u, \la) d\la\right)^2\\
&\les& r\int_{u=1}^{u_*}   du    \left( \int_{r} ^{r_*}  \la^{-2}   \|(\Hc, \widecheck{\DD(\cos\th)})\|^2 _{2,k}(u, \la) d\la \right)\\
&&+ \ep^2 r\int_{u=1}^{u_*}   du    \left( \int_{r} ^{r_*}  \la^{-2}   \|\Ga_b\|^2 _{2,k}(u, \la) d\la \right)\\
&\les& r\left(\int_r^{r_*}\la^{-2}d\la\right)\sup_{r\geq r_0}\left[\int_{r=\la}   \big| \dk^{\le k} (\Hc, \widecheck{\DD(\cos\th)}) \big|^2+\ep^2\int_{r=\la}\big| \dk^{\le k}\Ga_b\big|^2\right]\\
&\les& \sup_{r\geq r_0}\left[\int_{r=\la}   \big| \dk^{\le k} (\Hc, \widecheck{\DD(\cos\th)}) \big|^2+\ep^2\int_{r=\la}\big| \dk^{\le k}\Ga_b\big|^2\right].
  \eeaa
  Together with \eqref{Prop-bootstrap:finalcontrolofthePTframeonMext}, the  estimate of Step 3 for $\widecheck{\DD(\cos\th)})$ and the estimate of Step 6 for $\Hc$,  we infer, for $k\le k_{large}+7$,  
  \beaa
\int_{u=1}^{u_*} \left(r\int_r^{r_*}      \la ^{-1}\|F\|_{2,k}(u, \la) d\la\right)^2 &\les& \Big(\ep_0+\Skstar_k+\Rkext_k\Big)^2.
  \eeaa
 Plugging in the above, we infer
 \bea
\sup_{\la \ge r_0} \left( \int_{r=\la}   \big| \dk^{\le k}e_3(\cos\th) \big|^2 \right)^{\frac{1}{2}} & \les &\ep_0+\Skstar_k+\Rkext_k, \qquad k\le k_{large}+7.
\eea

  {\bf Step 8.} Estimate for $\ombc$.
  
  We start with the following equation, see Proposition \ref{Prop:linearizedPTstructure1}, 
   \beaa
   \nab_4\ombc  &=& \Re\left(\Pc\right)+F,\\ 
   F&=&   O(r^{-2})\Zc+O(r^{-2}) \Hc+\Ga_g\c\Ga_g.
   \eeaa
By Proposition \ref{Prop:transportrp-L2-knorms-PT},   we deduce, for $k\leq k_{large}+7$, 
   \beaa
    r^{-1} \|\omb\|_{ 2, k} (r, u) &\les&   r^{-1}_* \|\omb\|_{ 2, k} (r_*, u) + \int_r^{r_*} \la^{-1} \| \Pc \|_{2, k} + \int_r^{r_*} \la^{-1} \|  F \|_{2, k}.
   \eeaa
   By Cauchy-Schwartz, we infer
    \bea
    \lab{eq:intermediate-Step8}
    \bsplit
    r^{-2} \|\omb\|^2_{ 2, k} (r, u) &\les  r^{-2}_* \|\omb\|^2_{ 2, k} (r_*, u)+ r^{-4+\de_B}  \left(\int_r^{r_*} \la^{3-\de_B}  \| \Pc \|^2_{2, k} d\la  \right) \\
    &+ r^{-3}   \left(\int_r^{r_*}\la^{-2}\|(\Hc, \Zc)\|^2_{2, k} d\la +\ep^2  \int_r^{r_*}\la^{-2}\| \Ga_g  \|^2_{2, k} d\la\right).
    \end{split}
    \eea
   Multiplying by $r^2$, integrating in $u$, using  \eqref{Prop-bootstrap:finalcontrolofthePTframeonMext}, the  estimate of Step 4 for $\Zc$ and the estimate of Step 6 for $\Hc$,  we easily derive
  \bea
\sup_{\la \ge r_0} \left( \int_{r=\la}   \big| \dk^{\le k}\ombc \big|^2 \right)^{\frac{1}{2}} & \les &\ep_0+\Skstar_k+\Rkext_k, \qquad k\le k_{large}+7.
\eea

{\bf Step 9.}  Estimates for  $\widecheck{e_3(r)}$.

In view of Lemma \ref{lemma:transportequationine4fornabthetaandrande3thetaandr:linearized},
we have
\beaa
e_4(\widecheck{e_3(r)}) &=& -2\ombc +F,\\
F &=& O(r^{-2})\DD r+r\Ga_b\c\Ga_g.
\eeaa
Proceeding  as  above we deduce
\beaa
 r^{-1} \|\widecheck{e_3(r)}\|_{2, k} &\les&    r^{-1}_* \|\widecheck{e_3(r)}\|_{ 2, k} (r,_* u) + \int_r^{r_*} \la^{-1} \|\ombc\|_{2, k}(u, \la)  d\la + \int_r^{r_*} \la^{-1} \|F\|_{2, k}(u, \la)  d\la 
\eeaa
and hence, squaring and integrating in $u$, 
\beaa
r^{-2} \int_{\la =r} \big|\dk^{\le k}( \widecheck{e_3(r)}) \big|^2&\les& (\Skstar_k)^2+ \int_{u=1}^{u_*}\left( \int_r^{r_*} \la^{-1} \|\ombc\|_{2, k}(u, \la)  d\la\right)^2\\
&&+ \int_{u=1}^{u_*} \left( \int_r^{r_*} \la^{-1} \|F\|_{2, k}(u, \la)  d\la\right)^2.
\eeaa
Using Cauchy Schwarz,  \eqref{Prop-bootstrap:finalcontrolofthePTframeonMext}, and the estimate of Step 5 for $\DD r$, we infer
\beaa
r^{-2} \int_{\la =r} \big|\dk^{\le k}( \widecheck{e_3(r)}) \big|^2&\les& (\Skstar_k)^2+\ep_0^2+ \int_{u=1}^{u_*}\left( \int_r^{r_*} \la^{-1} \|\ombc\|_{2, k}(u, \la)  d\la\right)^2.
\eeaa

The term in  $\ombc$ is the more  dangerous as it could lead to  a logarithmic divergence. We estimate it using 
the more precise estimate for $\ombc$ in \eqref{eq:intermediate-Step8}.
Thus,
\beaa
   \int_r^{r_*} \la^{-1} \|\ombc\|_{2, k}(u, \la)  d\la &\les &  \frac{r_*-r}{r_*}\|\ombc\|_{ 2, k} (r_*, u)+ r^{-1+\frac{\de_B}{2}}   \left(\int_{r}^{r_*} \la^{3-\de_B}  \| \Pc \|^2_{2, k} d\la  \right)^{1/2}\\
   &&+ r^{-\frac{1}{2}}\left(\int_r^{r_*}\la^{-2}\|(\Hc, \Zc)\|^2_{2, k} d\la +\ep^2  \int_r^{r_*}\la^{-2}\| \Ga_g  \|^2_{2, k} d\la\right)^{\frac{1}{2}}\\
   &\les &  \|\ombc\|_{ 2, k} (r_*, u)+  \left(\int_{r}^{r_*} \la^{3-\de_B}  \| \Pc \|^2_{2, k} d\la  \right)^{1/2}\\
   &&+ \left(\int_r^{r_*}\la^{-2}\|(\Hc, \Zc)\|^2_{2, k} d\la +\ep^2  \int_r^{r_*}\la^{-2}\| \Ga_g  \|^2_{2, k} d\la\right)^{\frac{1}{2}}.
 \eeaa
 Squaring, integrating in $u$, using  \eqref{Prop-bootstrap:finalcontrolofthePTframeonMext}, and the  estimate of Step 4 for $\Zc$ and the estimate of Step 6 for $\Hc$,  we easily derive
 \beaa
\int_{u=1}^{u_*} \left( \int_r^{r_*} \la^{-1} \|\ombc\|_{2, k}(u, \la)  d\la\right)^2\les  \Big(\Skstar_k+\Rkext_k+\ep_0\Big)^2, \qquad k\le k_{large}+7.
\eeaa
In vie of the above, we infer
  \bea
\sup_{\la \ge r_0} \la^{-2}\left( \int_{r=\la}   \big| \dk^{\le k}( \widecheck{e_3(r)}) \big|^2 \right)^{\frac{1}{2}} & \les &\ep_0+\Skstar_k+\Rkext_k, \quad k\le k_{large}+7.
\eea

{\bf Step 10.} Estimates for  $\DD\hot\Jk$ and $\widecheck{\DD\c\ov{\Jk}}$.

We make use of equations, see  Lemma \ref{lemma:transportequationine4fornaJkandnab3Jk:linearized}, 
\beaa
\nab_4 \DD\hot\Jk+\frac{2}{q}\DD\hot\Jk &=& O(r^{-1})B +O(r^{-2})\trXc+O(r^{-2})\Xh\\
&&+O(r^{-2})\Zc +O(r^{-3})\widecheck{\DD(\cos\th)},\\
\nab_4\widecheck{\ov{\DD}\c\Jk} +\Re\left(\frac{2}{q}\right)\widecheck{\ov{\DD}\c\Jk} &=& O(r^{-1})B+O(r^{-2})\trXc+O(r^{-2})\Xh +O(r^{-2})\Zc\\
&&+O(r^{-3})\widecheck{\DD(\cos\th)}.
\eeaa
 Using  Corollary  \ref{Corr:Prop:transportrp-L2-knorms-PT} with $c=2$ and  $C=2$, the estimates of Step1--4 for $\trXc$, $\Xh$, $\widecheck{\DD(\cos\th)}$ and $\Zc$,  we easily  derive
 \bea
\sup_{\la \ge r_0}\left( \int_{r=\la} \la^2\big|\dk^{\le k} (\DD\hot\Jk, \,  \widecheck{\ov{\DD}\c\Jk) }\big|^2 \right)^{\frac{1}{2}} \les \ep_0+\Skstar_k+\Rkext_k, \qquad k\le k_{large}+7.
\eea

{\bf Step 11.} Estimates for $\trXbc$.

We make use of the following equation, see  Proposition \ref{Prop:linearizedPTstructure1}, 
\beaa
\nab_4\trXbc +\frac{1}{q}\trXbc  &=& 2\ov{\Pc}+F,\\
F &=&   O(r^{-1})\trXc +O(r^{-1})\widecheck{\DD\c\ov{\Jk}} +O(r^{-3})\DD(r) +O(r^{-3})\widecheck{\DD(\cos\th)} +\Ga_b\c\Ga_g.
\eeaa
 Using Proposition \ref{Prop:transportrp-L2-knorms-PT}, we infer
 \beaa
 \| \trXbc\|_{2, k}(u, r) &\les & \| \trXbc\|_{2, k}(u, r_* )  +\int_r^{r_*} \|\Pc\|_{2, k} (u, \la) d\la+\int_r^{r_*} \| F\|_{2, k} (u, \la) d\la. 
 \eeaa
 Now, we have, for $k\leq k_{large}+7$, 
 \beaa
 \int_r^{r_*} \| \Pc \|_{2, k} (u, \la) d\la&\les & r^{-1+\frac{\de_B}{2} } \left( \int_r^{r_*} \la^{3-\de_B}  \| \Pc \|^2_{2, k} (u, \la) d\la\right)^{\frac{1}{2}},
  \\
\int_{u=1}^ {u_*}  \left( \int_r^{r_*} \| \Pc \|_{2, k} (u, \la) d\la\right)^2 &\les&  r^{-2+\de_B}(\Rkext_k)^2.
 \eeaa
 Also, we have
 \beaa
 &&\int_r^{r_*} \| F\|_{2, k} (u, \la) d\la\\
  &\les& r^{-1+\frac{\de_B}{2} }\left( \int_r^{r_*} \la^{3-\de_B}  \| F\|^2_{2, k} (u, \la) d\la\right)^{\frac{1}{2}}\\
 &\les& r^{-1+\frac{\de_B}{2} }\Bigg( \int_r^{r_*} \la^{1-\dt}\Big( \|(\trXc, \widecheck{\DD\c\ov{\Jk}})\|_{2, k} (u, \la)+\la^{-2} \|(\DD r, \widecheck{\DD(\cos\th)})\|_{2, k} (u, \la)\\
 &&+ \la^2\|\Ga_b\c\Ga_g\|_{2, k} (u, \la)\Big)^2\Bigg)^{\frac{1}{2}}.
 \eeaa
 Integrating in $u$, using  \eqref{Prop-bootstrap:finalcontrolofthePTframeonMext}, and the  estimates of Steps 1, 3,  5 and 10 respectively for  $\trXc$, $\widecheck{\DD(\cos\th)}$, $\DD r$ and  $\widecheck{\DD\c\ov{\Jk}}$,  we easily derive, for $k\leq k_{large}+7$, 
\beaa
\int_{u=1}^ {u_*}  \left( \int_r^{r_*} \| F \|_{2, k} (u, \la) d\la\right)^2 &\les&  r^{-2+\de_B}(\Skstar+\Rkext_k+\ep_0)^2.
\eeaa 
We deduce  
 \bea
\sup_{\la \ge r_0}\left( \int_{r=\la} \la^{2-\dt}\big|\dk^{\le k}\trXbc\big|^2 \right)^{\frac{1}{2}} \les \ep_0+\Skstar_k+\Rkext_k, \qquad k\le k_{large}+7.
\eea

 {\bf Step 12.} Estimates for $\Xbh$.
 
 We make use of  the following equation, see  Proposition \ref{Prop:linearizedPTstructure1}, 
 \beaa
 \nab_4\Xbh +\frac{1}{q}\widehat{\Xb}  &=&  O(r^{-1})\DD\hot\Jk +O(r^{-3})\DD r+O(r^{-3})\widecheck{\DD(\cos\th)} +O(r^{-1}) \widehat{X}+\Ga_b\c\Ga_g.
\eeaa
Using  Corollary  \ref{Corr:Prop:transportrp-L2-knorms-PT} with $c=1$ and  $C=1$,  \eqref{Prop-bootstrap:finalcontrolofthePTframeonMext}, and  the  estimates of Steps 2, 3, 5 and 10 respectively for  $\Xh$, $\widecheck{\DD(\cos\th)}$, $\DD r$ and  $\widecheck{\DD\c\ov{\Jk}}$,  we easily derive
 \bea
\sup_{\la \ge r_0}\left( \int_{r=\la}\big|\dk^{\le k}\Xbh\big|^2 \right)^{\frac{1}{2}} \les \ep_0+\Skstar_k+\Rkext_k, \qquad k\le k_{large}+7.
\eea

 {\bf Step 13.} Estimates for $\widecheck{\nab_3\Jk}$.
 
 We make use of the following equation, see Lemma \ref{lemma:transportequationine4fornaJkandnab3Jk:linearized},
\beaa
\nab_4\widecheck{\nab_3\Jk}+\frac{1}{q}\widecheck{\nab_3\Jk} &=& F,\\
F&=&O(r^{-1})\Pc+ O(r^{-3})\widecheck{e_3(r)}+O(r^{-3})e_3(\cos\th)+O(r^{-2})\ombc\\
&&+O(r^{-2})\Hc+O(r^{-2})\widecheck{\nab\Jk}.
\eeaa
Using Proposition \ref{Prop:transportrp-L2-knorms-PT}, we infer
 \beaa
 \| \widecheck{\nab_3\Jk}\|_{2, k}(u, r) &\les & \| \widecheck{\nab_3\Jk}\|_{2, k}(u, r_* )  +\int_r^{r_*} \| F\|_{2, k} (u, \la) d\la. 
 \eeaa
Now, we have
\beaa
&&\int_r^{r_*} \| F\|_{2, k} (u, \la) d\la\\ 
&\les& r^{-2+\frac{\dt}{2}}\left(\int_r^{r_*}r^{3-\dec}\|\Pc\|^2_{2, k} (u, \la) d\la\right)^{\frac{1}{2}}\\
&&+r^{-\frac{1}{2}}\left(\int_r^{r_*}r^{-2}\Big(\|(\ombc, \Hc, \widecheck{\nab\Jk})\|^2_{2, k} (u, \la)+r^{-2}\|(\widecheck{e_3(r)}, e_3(\cos\th))\|^2_{2, k} (u, \la)\Big) d\la\right)^{\frac{1}{2}}.
\eeaa
Squaring, integrating in $u$, using \eqref{Prop-bootstrap:finalcontrolofthePTframeonMext}, and  the  estimates of Steps 6--10 respectively for  $\Hc$, $e_3(\cos\th)$, $\ombc$, $\widecheck{e_3(r)}$ and $\widecheck{\nab\Jk}$,  we easily derive
\beaa
\int_{u=1}^ {u_*}  \left( \int_r^{r_*} \| F \|_{2, k} (u, \la) d\la\right)^2 &\les&  r^{-2}(\Skstar+\Rkext_k+\ep_0)^2.
\eeaa 
We deduce  
 \bea
\sup_{\la \ge r_0}\left( \int_{r=\la}\la^2\big|\dk^{\le k}\widecheck{\nab_3\Jk}\big|^2 \right)^{\frac{1}{2}} \les \ep_0+\Skstar_k+\Rkext_k, \qquad k\le k_{large}+7.
\eea

{\bf Step 14.} Estimates for $\Xib$.

We make use of  the following equation, see  Proposition \ref{Prop:linearizedPTstructure1}, 
\beaa
 \nab_4\Xib &=&  F,\\
 F &=& O(r^{-1})\Hc +O(r^{-2})\trXbc +O(r^{-2})\Xbh -\Bb +O(r^{-1})\widecheck{\nab_3\Jk}+O(r^{-3})\widecheck{e_3(r)}\\
 &&+O(r^{-3})e_3(\cos\th) +\Ga_b\c\Ga_b.
\eeaa
Using Proposition \ref{Prop:transportrp-L2-knorms-PT}, we infer
 \beaa
r^{-1} \| \Xib\|_{2, k}(u, r) &\les & r_*^{-1}\| \Xib\|_{2, k}(u, r_* )  +\int_r^{r_*} \la^{-1}\| F\|_{2, k} (u, \la) d\la. 
 \eeaa
Now, we have
\beaa
&&\int_r^{r_*} \| F\|_{2, k} (u, \la) d\la\\
&\les&  r^{-1+\frac{\dt}{2}}\left(\int_r^{r_*}\la^{1-\dt}\| \Bb\|^2_{2, k} (u, \la) d\la\right)^{\frac{1}{2}}\\
&& +r^{-\frac{1}{2}}\Bigg(\int_r^{r_*}\la^{-2}\Big(\|(\Hc, \widecheck{\nab_3\Jk})\|^2_{2, k} (u, \la)
+\la^{-2}\|(\trXbc, \Xbh )\|^2_{2, k} (u, \la)\\
&&\qquad\qquad\qquad +\la^{-4}\|(\widecheck{e_3(r)}, e_3(\cos\th) )\|^2_{2, k} (u, \la)\Big)\Bigg)^{\frac{1}{2}}.
\eeaa
Squaring, integrating in $u$, using \eqref{Prop-bootstrap:finalcontrolofthePTframeonMext}, and  the  estimates of Steps 6, 7, 9, 11, 12, 13 respectively for  $\Hc$, $e_3(\cos\th)$, $\widecheck{e_3(r)}$, $\trXbc$, $\Xbh$ and $\widecheck{\nab_3\Jk}$,  we easily derive
\beaa
\int_{u=1}^ {u_*}  \left( \int_r^{r_*} \| F \|_{2, k} (u, \la) d\la\right)^2 &\les&  r^{-2+\dt}(\Skstar+\Rkext_k+\ep_0)^2.
\eeaa 
We deduce 
 \bea
\sup_{\la \ge r_0}\left( \int_{r=\la} \la^{-\dt}\big|\dk^{\le k}\Xi\big|^2 \right)^{\frac{1}{2}} \les \ep_0+\Skstar_k+\Rkext_k, \qquad k\le k_{large}+7.
\eea 
 
Gathering the estimates derived in Steps 1--14, we infer, in view of the definition of $\Skext_k$,    
\beaa
\Skext_k\les  \ep_0 +\Skstar_k+\Rkext_k, \qquad  k\le k_{large}+7,
\eeaa
as stated. This concludes the proof of Proposition \ref{Prop:EstmatesSkext}.

%%%%%%%%%%%%%%%%%%%%%%%%%%%%%
 
  \section{Control of the PT-Ricci coefficients in $\Mint'$}
  \lab{sec:controlofthePTRiccicoefficientsinMintprime} 
  
%%%%%%%%%%%%%%%%%%%%%%%%%%%%%  

Since we will not need to refer to the old region $\Mint$, defined w.r.t. the PG frame, we drop  the prime of $\Mint'$ in this section.

Recall the following norms on $\Mint$ introduced in section \ref{sec:mainnormsPTframe:chap9}, for $k\le k_{large}+7$,
\beaa
\Skint_k^2 &=& \int_{\Mint}\big|\dk^{\le k}\Gac\big|^2, \\
\Rkint_k^2 &=&  \int_{\Mint}\Big( \big| \nab_{\Rhat} \dk^{\le k-1}\Rc\big|^2+|\dk^{\le k-1}\Rc|^2\Big) +\sup_\tau\int_{\Mint\cap\Si(\tau)}  |\dk^{\le k}  \Rc|^2,
\eeaa
where $\Gac$ denotes the set of all  linearized Ricci and metric coefficients  with respect  to the ingoing PT frame of $\Mint$, i.e.   
  \beaa
\Gac =\Big\{  \trXbc,  \Xbh,  \Zc,  \Hbc, \widecheck{ \DD\cos \th},   \omc,  \DD r, \widecheck{\DD u},    e_4(\cos\th), \widecheck{e_4(r)},  \widecheck{e_4(\ub)},    \widecheck{ \ov{\DD}\c \Jk},  \DD\hot\Jk, \widecheck{\nab_4 \Jk},             \trXc, \Xh,  \Xi\Big\},
  \eeaa
where $\Rc$ denotes the set of all  linearized curvature components with respect  to the ingoing PT frame of $\Mint$, i.e.
\beaa
\Rc=\{A, B, \Pc, \Bb, \Ab\}, 
\eeaa
and where the vectorfield $\Rhat$ in $\Mint$ is given by, see \eqref{eq:ThatRhat-e_3e_4},
   \beaa
\Rhat = \frac 1 2 \left( \frac{|q|^2}{r^2+a^2} e_4-\frac{\De}{r^2+a^2}  e_3\right).
 \eeaa

\begin{remark}
Recall  that $r\leq r_0$ in $\Mint$ so that $r$ is uniformly bounded in that region. In particular, all components have the same behavior in $\Mint$ which is reflected in   the definition of the norms $\Skint_k$ and $\Rkint_k$. 
\end{remark}

The goal of this section is to  provide the proof of  Proposition \ref{prop:controlGaintiterationassupmtionThM8}. For convenience, we restate the result below.
 \begin{proposition}[Control of $\Gac$ in $\Mint$]
 \lab{Prop:MainestimatesMint}
 Relative to the PT frame of $\Mint$, we have
 \beaa
 \Skint_k&\les & \ep_0 + \Skext_k+\Rkint_k, \qquad k \le  k_{large}+7.
 \eeaa
 \end{proposition} 
 
 To prove Proposition \ref{Prop:MainestimatesMint}, we rely in particular on our bootstrap assumption {\bf  BA-PT}, see \eqref{eq:mainbootassforchapte9}, which implies   for the Ricci and metric coefficients of the ingoing PT frame of $\Mint$
 \bea
 \lab{eq:Bootstrap-Mint-M8}
 \Skint_k&\leq & \ep, \qquad k\leq k_{large}+7.
 \eea

%%%%%%%%%%%%%%%%%%%%%%%%%

\subsection{Preliminaries}

%%%%%%%%%%%%%%%%%%%%%%%%%

 Recall that the following identities hold for an ingoing PT structure 
\beaa
\xib=0, \quad \omb=0, \quad  H=\frac{aq}{|q|^2} \Jk, \quad  e_3(r)=-1, \quad e_3(\ub)=e_3(\th)=0,\quad \nab_3\Jk = \frac{1}{\ov{q}}\Jk.
\eeaa

In view of Definition \ref{def:linearizedPT-ingoingcase:chap9} and the notation $\Gac$, we have
\beaa
e_4(r)=\frac{\Delta}{|q|^2}+\Gac , \qquad \nab(r)=\Gac, \qquad e_4(\ub)=\frac{2(r^2+a^2)}{|q|^2}+\Gac,\qquad \DD\ub=a\Jk+\Gac.
\eeaa
Also, note from the definition of $\Rhat$ that
\beaa
\Rhat(r)=\frac{\De}{r^2+a^2}+\Gac, \qquad e_4=\frac{2(r^2+a^2)}{|q|^2}\Rhat+\frac{\De}{|q|^2}  e_3.
\eeaa

%%%%%%%%%%%%%%%%%%%%%%%%%%%%%%%%%

\subsubsection{The region $\Mint_*$}

%%%%%%%%%%%%%%%%%%%%%%%%%%%%%%%%%

We recall the scalar function $\tau$ constructed in Proposition \ref{prop:propertiesoftauusefulfortheoremM8:chap9} whose level sets are uniformly spacelike in $\Mint$, and the region 
\bea
\Mint_*=\Mint\cap\{\tau\leq\tau_*\},
\eea 
see Definition \ref{def:u=introductionofregionMintprimestar:chap9}. In particular, the boundary of $\Mint_*$ is given by
\bea
\pr\Mint_*=\AA\cup\Sint_*\cup \TT\cup \{\ub=1\}
\eea
where $ \AA $  and 
\bea
\Sint_*=\Mint\cap\{\tau=\tau_*\}
\eea  
are  strictly spacelike, while $\TT=\{r=r_0\}$ and $\{\ub=1\}$ are timelike, with $\{\ub=1\}$ included in the initial data layer.

%%%%%%%%%%%%%%%%%%%%%%%%%%%%%%%%%

\subsubsection{A simple computation}

%%%%%%%%%%%%%%%%%%%%%%%%%%%%%%%%%

The following simple computation will be useful in Lemma \ref{lemma:general-transport-Mint}. 
\begin{lemma}
\lab{Lemma:div-e_3}
For any function $f$, we have
\beaa
\Div ( f e_3)&=& e_3(f)  -\frac{2r}{|q|^2} f +f\Gac.
\eeaa
\end{lemma}

\begin{proof}
 We have, with $\pi^{(3)}$ the deformation tensor of $e_3$,
\beaa
\Div ( e_3)=\frac 1 2 \tr \pi^{(3)}= \frac 1 2 \left(\de_{ab} \pi^{(3)}_{ab} -\pi^{(3)}_{34}\right)
\eeaa
Observe that
\beaa
\,\pi^{(3)}_{34}&=& \g(\D_3 e_3, e_4)+\g(\D_4 e_3, e_3)=0, \\
\,\pi^{(3)}_{ab}&=&\g(\D_a e_3, e_b)+\g(\D_b e_3, e_a)=  \chib_{ab}+\chib_{ba}=\trchb \de_{ab}+2\chibh_{ab}.
\eeaa
This yields
\beaa
\Div (e_3)= \trchb=-\frac{2r}{|q|^2}+\Gac.
\eeaa
Thus,
\beaa
\Div ( fe_3)= f \Div(e_3) + e_3(f) = -\frac{2r}{|q|^2} f +f\Gac+ e_3(f) 
\eeaa
 as stated.
 \end{proof}

%%%%%%%%%%%%%%%%%%%%%%%%%%%%%%%%%

\subsubsection{Commutation formulas}

%%%%%%%%%%%%%%%%%%%%%%%%%%%%%%%%%

\begin{lemma}
\lab{Le:Comm-Mint} 
Given  $U$  a  horizontal tensor  in $\Mint$, we have\footnote{Recall that $r\leq r_0$ in $\Mint$ so that weights in $r$ do not matter and are hence dropped.}
\beaa
\,[\nab_3, \nab] U &=& O(1)\nab U+\Gac\dk U+\Big(O(1)+\Gac+\Rc\Big)U,\\ 
\,[\nab_3, \nab_4] U &=& O(1)\nab U+ O(1)\nab_3U+\Gac\dk U +\Big(O(1)+\Gac+\Rc\Big)U.
\eeaa
\end{lemma} 

\begin{proof}
In view of Lemma \ref{lemma:comm-gen}, we have in $\Mint$, for a real horizontal tensor  in $\Mint$, 
\beaa
\,[\nab_3, \nab_b] U &=& -\chib_{bc}\nab_c U +(\eta_b-\ze_b)  \nab_3 U +\xib_b \nab_4 U+\Big(O(1)+\Gac+\Rc\Big)U,\\ 
\,[\nab_4, \nab_3] U &=& 2(\etab_b-\eta_b ) \nab_b U + 2 \om \nab_3 U -2\omb \nab_4 U +\Big(O(1)+\Gac+\Rc\Big)U.
\eeaa
Since, in view of the identities for ingoing PT structures, we have $\omb=0$, $\xib=0$, and $\ze-\eta=\zec$, we infer
\beaa
\,[\nab_3, \nab] U &=& O(1)\nab  +\Gac\dk U+\Big(O(1)+\Gac+\Rc\Big)U,\\ 
\,[\nab_4, \nab_3] U &=& O(1)\nab U+ O(1)\nab_3U+\Gac\dk U +\Big(O(1)+\Gac+\Rc\Big)U.
\eeaa
as stated.
\end{proof} 

\begin{corollary}
\lab{Cor:Comm-Mint} 
Given  $U$  a complex horizontal tensor  in $\Mint$ we have
\beaa
\,[\nab_3, \nab_\Rhat] U &=& O(1)\nab U+ O(1)\nab_3U+O(1)\nab_\Rhat U+\Gac\dk U +\Big(O(1)+\Gac+\Rc\Big)U.
\eeaa
\end{corollary} 

\begin{proof}
Recall that we have
  \beaa
 \Rhat= \frac 1 2 \left( \frac{|q|^2}{r^2+a^2} e_4-\frac{\De}{r^2+a^2}  e_3\right),
 \eeaa
so that 
\beaa
[\nab_3, \nab_\Rhat] &=&  \frac 1 2 \left( e_3\left(\frac{|q|^2}{r^2+a^2}\right) \nab_4-e_3\left(\frac{\De}{r^2+a^2}  \right)\nab_3\right)+ \frac 1 2  \frac{|q|^2}{r^2+a^2}[\nab_3, \nab_4]\\
&=& O(1)\nab_\Rhat+O(1)\nab_3+O(1)[\nab_3, \nab_4]
\eeaa
and the corollary  follows immediately from the formula for $[\nab_3, \nab_4]$ in Lemma \ref{Le:Comm-Mint}.
\end{proof}

 %%%%%%%%%%%%%%%%%%%%%%%%%%%%%
 
 \subsection{Commuted structure equations}
 
  %%%%%%%%%%%%%%%%%%%%%%%%%%%%%

In this section, we rely on the commutation formulas of Lemma \ref{Le:Comm-Mint} and Corollary \ref{Cor:Comm-Mint} to derive the equations for higher order derivatives of the equations for the Ricci coefficients and metric coefficients of the ingoing PT structure of $\Mint$. We start with the following proposition concerning the schematic form of the uncommuted equations in $\Mint$.  

\begin{proposition}\lab{prop:schamticformoflinearzizedstrucutreequationsofMint}
In the PT structure of $\Mint$,  the equations of Proposition \ref{proposition:Maineqts-Mint'} for the Ricci and metric coefficients take the following form\footnote{Recall that $r\leq r_0$ in $\Mint$ so that weights in $r$ do not matter and are hence dropped.}:
\begin{enumerate}
\item $\Xbh$ satisfies a transport equation of the following form 
\bea
\nab_3\Xbh &=& O(1)\Xh -\Ab+\Gac\c\Gac.
\eea

\item If $\Phi$ is among the quantities $\trXbc$, $\Xh$, $\widecheck{\DD\cos\th}$, $\DD r$, $e_4(\cos\th)$, $\widecheck{e_4(r)}$, then $\Phi$ satisfies a transport equation of the following form
\bea
\nab_3\Phi &=& O(1)\Phi +\Gac_0[\Phi]+\Gac\c\Gac,
\eea
where 
\bea
\bsplit
\Gac_0[\trXbc]&=0, \qquad\qquad\qquad\qquad\qquad \Gac_0[\Xh] =\Big(\DD\hot\Jk, \,\widecheck{ \DD(\cos \th)},\, \Xbh\Big),\\
\Gac_0[\widecheck{\DD\cos\th}] &=\Big(\trXbc, \,\Xbh\Big), \quad\,\,\,\,\qquad\qquad \Gac_0[\DD r] =\Zc,\\ 
\Gac_0[e_4(\cos\th)] &=\Big(\Hbc, \widecheck{\DD(\cos\th)}\Big), \qquad\quad\Gac_0[\widecheck{e_4(r)}] =\Big(\omc, \DD r\Big).
\end{split}
\eea

\item If $\Phi$ is among the quantities $\Zc$, $\Hbc$, $\trXc$, $\omc$, $\Xi$, $\DD\hot\Jk$, $\widecheck{\ov{\DD}\c\Jk}$, $\widecheck{\nab_4\Jk}$, then $\Phi$ satisfies a transport equation of the following form
\bea
\nab_3\Phi &=& O(1)\Phi +\Gac_0[\Phi]+\Rc_0[\Phi]+\Gac\c\Gac,
\eea
where $\Rc_0[\Phi]$ is a curvature component among the list
\bea
\Rc_0[\Phi] &=& \Big\{B,\,\,\, \Pc, \,\,\,\Bb\Big\},
\eea
and 
\bea
\bsplit
\Gac_0[\Zc]&=\Big(\Xbh,\, \trXbc\Big),\\
\Gac_0[\Hbc]&=\Big(\Xbh,\, \trXbc\Big),\\
\Gac_0[\trXc]&=\Big(\trXbc, \,\widecheck{\DD\c\ov{\Jk}},\, \DD r,\,\widecheck{\DD(\cos\th)}\Big),\\
\Gac_0[\omc]&=\Big(\Zc,\,\Hbc\Big),\\
\Gac_0[\Xi]&= \Big(\Hbc,\,\trXc,\,\Xh,\,\widecheck{\nab_4 \Jk},\, \widecheck{e_4(r)},\,e_4(\cos\th)\Big),\\
\Gac_0[\DD\hot\Jk]&=\Big(\trXbc,\,\Xbh,\,\Zc,\,\widecheck{\DD(\cos\th)}\Big),\\
\Gac_0[\widecheck{\ov{\DD}\c\Jk}]&=\Big(\trXbc,\,\Xbh,\,\Zc, \,\widecheck{\DD(\cos\th)}\Big),\\
\Gac_0[\widecheck{\nab_4\Jk}]&=\Big(\widecheck{e_4(r)},\,e_4(\cos\th),\,\omc,\,\Hbc,\,\widecheck{\nab\Jk}\Big). 
\end{split}
\eea
\end{enumerate}
In addition, the linearized Codazzi equation for $\Xbh$ in Proposition \ref{proposition:Maineqts-Mint':bis} 
\beaa
\ov{\DD}\c\Xbh &=& \DD\ov{\trXbc} +\Rc+\Gac+\Gac\c\Gac.
\eeaa
\end{proposition}

\begin{proof}
This follows immediately from Proposition \ref{proposition:Maineqts-Mint'}.
\end{proof}

We now commute the equations of Proposition \ref{prop:schamticformoflinearzizedstrucutreequationsofMint}.
\begin{proposition}\lab{prop:schamticformoflinearzizedstrucutreequationsofMint:bis}
Let $1\leq k\leq k_{large}$ an integer. In the PT structure of $\Mint$, the commutation of the equation of Proposition \ref{prop:schamticformoflinearzizedstrucutreequationsofMint} for the Ricci and metric coefficients take the following form:
\begin{enumerate}
\item $\Xbh$ satisfies a transport equation of the following form 
\bea
\nn\nab_3(\nab_{\Rhat}\dk^{k-1}\Xbh) &=& O(1)\nab_{\Rhat}\dk^{k-1}\Xbh+O(1)\nab\dk^{k-1}\Xbh\\
&& +\nab_{\Rhat}\dk^{k-1}\Rc+\dk^{\leq k-1}(\Rc, \Gac)+\dk^{\leq k}(\Gac\c\Gac).
\eea

\item If $\Phi$ is among the quantities $\trXbc$, $\Xh$, $\widecheck{\DD\cos\th}$, $\DD r$, $e_4(\cos\th)$, $\widecheck{e_4(r)}$, then $\Phi$ satisfies a transport equation of the following form
\bea
\nn\nab_3(\nab_{\Rhat}\dk^{k-1}\Phi) &=& O(1)\nab_{\Rhat}\dk^{k-1}\Phi +O(1)\nab\dk^{k-1}\Phi +\Gac_k[\Phi]\\
&&+\dk^{\leq k-1}(\Rc, \Gac)+\dk^{\leq k}(\Gac\c\Gac),
\eea
and 
\bea
\nn\nab_3(\nab\dk^{k-1}\Phi) &=& O(1)\nab_{\Rhat}\dk^{k-1}\Phi +O(1)\nab\dk^{k-1}\Phi +\Gac_k[\Phi]\\
&&+\dk^{\leq k-1}(\Rc, \Gac)+\dk^{\leq k}(\Gac\c\Gac),
\eea
where 
\bea
\bsplit
\Gac_k[\trXbc]&=0, \qquad\qquad\qquad\qquad\qquad\,\,\,\,\, \Gac_k[\Xh] =\dk^k\Big(\DD\hot\Jk, \,\widecheck{ \DD(\cos \th)},\, \Xbh\Big),\\
\Gac_k[\widecheck{\DD\cos\th}] &=\dk^k\Big(\trXbc, \,\Xbh\Big), \quad\,\,\,\,\qquad\qquad \Gac_k[\DD r] =\dk^k\Zc,\\ 
\Gac_k[e_4(\cos\th)] &=\dk^k\Big(\Hbc, \widecheck{\DD(\cos\th)}\Big), \qquad\quad\Gac_k[\widecheck{e_4(r)}] =\dk^k\Big(\omc, \DD r\Big).
\end{split}
\eea

\item If $\Phi$ is among the quantities $\Zc$, $\Hbc$, $\trXc$, $\omc$, $\Xi$, $\DD\hot\Jk$, $\widecheck{\ov{\DD}\c\Jk}$, $\widecheck{\nab_4\Jk}$, then $\Phi$ satisfies  transport equations of the following form
\bea
\nn\nab_3(\nab_{\Rhat}\dk^{k-1}\Phi) &=& O(1)\nab_{\Rhat}\dk^{k-1}\Phi+O(1)\nab\dk^{k-1}\Phi +\Gac_k[\Phi] \\
&&+\nab_{\Rhat}\dk^{k-1}\Rc+\dk^{\leq k-1}(\Rc, \Gac)+\dk^{\leq k}(\Gac\c\Gac),
\eea
and
\bea\lab{eq:slightlymodresubtlecommutationformulafortransporteqinMint'}
\nn\nab_3(\nab\dk^{k-1}\Phi+\dk^{k-1}\Rc) &=& O(1)\nab_{\Rhat}\dk^{k-1}\Phi+O(1)\nab\dk^{k-1}\Phi +\Gac_k[\Phi] \\
&&+\nab_{\Rhat}\dk^{k-1}\Rc+\dk^{\leq k-1}(\Rc, \Gac)+\dk^{\leq k}(\Gac\c\Gac),
\eea
where
\bea
\bsplit
\Gac_k[\Zc]&=\dk^k\Big(\Xbh,\, \trXbc\Big),\\
\Gac_k[\Hbc]&=\dk^k\Big(\Xbh,\, \trXbc\Big),\\
\Gac_k[\trXc]&=\dk^k\Big(\trXbc, \,\widecheck{\DD\c\ov{\Jk}},\, \DD r,\,\widecheck{\DD(\cos\th)}\Big),\\
\Gac_k[\omc]&=\dk^k\Big(\Zc,\,\Hbc\Big),\\
\Gac_k[\Xi]&= \dk^k\Big(\Hbc,\,\trXc,\,\Xh,\,\widecheck{\nab_4 \Jk},\, \widecheck{e_4(r)},\,e_4(\cos\th)\Big),\\
\Gac_k[\DD\hot\Jk]&=\dk^k\Big(\trXbc,\,\Xbh,\,\Zc,\,\widecheck{\DD(\cos\th)}\Big),\\
\Gac_k[\widecheck{\ov{\DD}\c\Jk}]&=\dk^k\Big(\trXbc,\,\Xbh,\,\Zc, \,\widecheck{\DD(\cos\th)}\Big),\\
\Gac_k[\widecheck{\nab_4\Jk}]&=\dk^k\Big(\widecheck{e_4(r)},\,e_4(\cos\th),\,\omc,\,\Hbc,\,\widecheck{\nab\Jk}\Big). 
\end{split}
\eea
\end{enumerate}
In addition, we have
\bea
\ov{\DD}\c(\dk^{k-1}\Xbh) &=& \dk^k\trXbc +\dk^{\leq k-1}(\Rc, \Gac)+\dk^{\leq k-1}(\Gac\c\Gac).
\eea
\end{proposition}

\begin{proof}
We focus on \eqref{eq:slightlymodresubtlecommutationformulafortransporteqinMint'} as the other commuted equations follow immediately from the corresponding uncommuted equation in Proposition \ref{prop:schamticformoflinearzizedstrucutreequationsofMint} and the commutation formulas of Lemma \ref{Le:Comm-Mint} and Corollary \ref{Cor:Comm-Mint}. 

Recall that the case of \eqref{eq:slightlymodresubtlecommutationformulafortransporteqinMint'}, $\Phi$ satisfies according to Proposition \ref{prop:schamticformoflinearzizedstrucutreequationsofMint}  the following uncommuted equation 
\beaa
\nab_3\Phi &=& O(1)\Phi +\Gac_0[\Phi]+\Rc_0[\Phi]+\Gac\c\Gac,
\eeaa
where $\Rc_0[\Phi]$ is a curvature component among the list $\{B,\, \Pc,\,\Bb\}$. Thus, using the commutation formulas of Lemma \ref{Le:Comm-Mint} and Corollary \ref{Cor:Comm-Mint}, we immediately infer
\beaa
\nn\nab_3(\nab\dk^{k-1}\Phi) &=& O(1)\nab_{\Rhat}\dk^{k-1}\Phi+O(1)\nab\dk^{k-1}\Phi+O(1)\nab_3\dk^{k-1}\Phi +\Gac_k[\Phi] \\
&&+\dk^{k-1}\nab\Rc_0[\Phi]+\dk^{\leq k-1}(\Rc, \Gac)+\dk^{\leq k}(\Gac\c\Gac).
\eeaa
Also, using 
\beaa
\nab_3(\dk^{k-1}\Phi )=\dk^{\leq k-1}\Gac+\dk^{\leq k-1}\Rc+\dk^{\leq k-1}(\Gac\c\Gac),
\eeaa
we may get rid of the term $O(1)\nab_3\dk^{k-1}\Phi $ on the RHS and obtain 
\beaa
\nn\nab_3(\nab\dk^{k-1}\Phi) &=& O(1)\nab_{\Rhat}\dk^{k-1}\Phi+O(1)\nab\dk^{k-1}\Phi +\Gac_k[\Phi] \\
&&+\dk^{k-1}\nab\Rc_0[\Phi]+\dk^{\leq k-1}(\Rc, \Gac)+\dk^{\leq k}(\Gac\c\Gac),
\eeaa
where the only term of the RHS which does not agree with \eqref{eq:slightlymodresubtlecommutationformulafortransporteqinMint'} is $\dk^{k-1}\nab\Rc_0[\Phi]$. 

Next, in view of Proposition \ref{proposition:Maineqts-Mint':bis},  we may rewrite some of the linearized Bianchi identities in the following schematic form
 \beaa
 \bsplit
 \DD\hot B &= \nab_3\Rc+  \Rc+\Gac+\Rc\c\Gac,\qquad\quad\, \DD\c \ov{B} =  \nab_4\Rc+  \Rc+\Gac+\Rc\c\Gac,\\
\DD\ov{\Pc} &= \nab_3\Rc+   \Rc+\Gac+\Rc\c\Gac,\qquad\qquad \DD \Pc =  \nab_4\Rc+  \Rc+\Gac+\Rc\c\Gac,\\
\ov{\DD}\c\Bb &=  \nab_3\Rc+  \Rc+\Gac+\Rc\c\Gac, \qquad\quad\, \DD\hot\Bb =  \nab_4\Rc+  \Rc+\Gac+\Rc\c\Gac.
\end{split}
\eeaa
Now, note that we have
\beaa
\nab\rhoc =\frac{1}{2}\Big(\Re(\DD\Pc)+\Re(\DD\ov{\Pc})\Big), \qquad \nab\rhodc =\frac{1}{2}\Big(\Im(\DD\Pc)-\Im(\DD\ov{\Pc})\Big)
\eeaa
so that $\DD\Pc$ and $\DD\ov{\Pc}$ generate any angular derivative of $\Pc$. Thus, in view of the above schematic linearized Bianchi identities, we have
\beaa
\nab\Pc &=& \nab_3\Rc+\nab_4\Rc+  \Rc+\Gac+\Rc\c\Gac.
\eeaa 
Also, note that for a complex anti-selfdual 1-form $F=f+i\dual f$, we have
\beaa
\div(f)=\frac{1}{2}\Re(\ov{\DD}\c F), \qquad \curl(f)=\frac{1}{2}\Im(\ov{\DD}\c F), \qquad \nab\hot f=\Re(\DD\hot F), 
\eeaa
and 
\beaa
\nab_af_b &=& \frac{1}{2}\div(f)\de_{ab}+\frac{1}{2}\curl(f)\in_{ab}+\frac{1}{2}(\nab\hot\b)_{ab},
\eeaa
so that $\ov{\DD}\c F$ and $\DD\hot F$ generate  any angular derivative of $F$. Thus, in view of the above schematic linearized Bianchi identities, we have
\beaa
\nab B &=& \nab_3\Rc+\nab_4\Rc+  \Rc+\Gac+\Rc\c\Gac, \\
\nab\Bb &=& \nab_3\Rc+\nab_4\Rc+  \Rc+\Gac+\Rc\c\Gac.
\eeaa 
Since $\Rc_0[\Phi]$ is a curvature component among the list $\{B,\, \Pc,\,\Bb\}$, we infer from the above identities for $\nab\Pc$, $\nab B$ and $\nab\Bb$ that we have
\beaa
\nab\Rc_0[\Phi] &=&  \nab_3\Rc+\nab_4\Rc+  \Rc+\Gac+\Rc\c\Gac.
\eeaa
Also, since $e_4$ is in the span of $\Rhat$ and $e_3$, we infer
\beaa
\nab\Rc_0[\Phi] &=&  \nab_3\Rc+\nab_{\Rhat}\Rc+  \Rc+\Gac+\Rc\c\Gac.
\eeaa

Next, recall from the above that we have
\beaa
\nn\nab_3(\nab\dk^{k-1}\Phi) &=& O(1)\nab_{\Rhat}\dk^{k-1}\Phi+O(1)\nab\dk^{k-1}\Phi +\Gac_k[\Phi] \\
&&+\dk^{k-1}\nab\Rc_0[\Phi]+\dk^{\leq k-1}(\Rc, \Gac)+\dk^{\leq k}(\Gac\c\Gac).
\eeaa
Plugging the above identity for $\nab\Rc_0[\Phi]$, we infer
\beaa
\nn\nab_3(\nab\dk^{k-1}\Phi) &=& O(1)\nab_{\Rhat}\dk^{k-1}\Phi+O(1)\nab\dk^{k-1}\Phi +\Gac_k[\Phi] \\
&&+\dk^{k-1}\nab_3\Rc+\dk^{k-1}\nab_{\Rhat}\Rc+\dk^{\leq k-1}(\Rc, \Gac)+\dk^{\leq k}(\Gac\c\Gac).
\eeaa
Using again the commutation formulas of Lemma \ref{Le:Comm-Mint} and Corollary \ref{Cor:Comm-Mint}, we infer 
\beaa
\nn\nab_3(\nab\dk^{k-1}\Phi) &=& O(1)\nab_{\Rhat}\dk^{k-1}\Phi+O(1)\nab\dk^{k-1}\Phi +\Gac_k[\Phi] \\
&&+\nab_3\dk^{k-1}\Rc+\nab_{\Rhat}\dk^{k-1}\Rc+\dk^{\leq k-1}(\Rc, \Gac)+\dk^{\leq k}(\Gac\c\Gac).
\eeaa
and hence
\beaa
\nn\nab_3(\nab\dk^{k-1}\Phi+\dk^{k-1}\Rc) &=& O(1)\nab_{\Rhat}\dk^{k-1}\Phi+O(1)\nab\dk^{k-1}\Phi +\Gac_k[\Phi] \\
&&+\nab_{\Rhat}\dk^{k-1}\Rc+\dk^{\leq k-1}(\Rc, \Gac)+\dk^{\leq k}(\Gac\c\Gac),
\eeaa
as stated in \eqref{eq:slightlymodresubtlecommutationformulafortransporteqinMint'}. This concludes the proof of Proposition \ref{prop:schamticformoflinearzizedstrucutreequationsofMint:bis}. 
\end{proof}

 %%%%%%%%%%%%%%%%%%%%%%%%%%%%%
 
 \subsection{Estimates for transport equations}
 
  %%%%%%%%%%%%%%%%%%%%%%%%%%%%%

We shall make use of the following divergence lemma in $\Mint_*$. 
\begin{lemma}[Divergence Lemma in $\Mint_*$]
Consider  a vectorfield $ X$  in $\Mint $ and    the region $\Mint_*=\Mint\cap\{\tau\leq\tau_*\}$. We have, 
\beaa
\int_{\Sint_*} \g(X, N) +\int_{\TT}\g(X, N)+\int_{\{\ub=1\}}\g(X, N) + \int_\AA \g(X, N) &=& \int_{\Mint_*} \Div(X),
\eeaa
where $N$ is the exterior unit normal to each  portion of the  boundary.
\end{lemma}

\begin{proof}
This follows immediately from the standard divergence lemma and the fact that  the region $\Mint_*$ has the boundary $\pr\Mint_*=\AA\cup\Sint_*\cup \TT\cup \{\ub=1\}$.  
\end{proof}

Next we  state an application of the divergence theorem to transport equations.
\begin{lemma}\lab{lemma:general-transport-Mint} 
Suppose $\Phi$ and  $F$ are two anti-selfdual horizontal tensors of the same type satisfying 
 \bea
 \lab{eq:TransportMint}
 \nab_3 \Phi=F.
 \eea
 Then, we have on $\Mint_*$, for any real number $p$,
\bea\label{eq:general-integrated-estimate-Mint}
\nn&&\int_{\Sint_*\cup\AA}|q|^{p}|\Phi|^2|\g(e_3, N)| + p\int_{\Mint_*}r|q|^{p-2}|\Phi|^2 \\
&\leq& \int_{\Mint_*}r^{-1}|q|^{p+2}|F|^2+\int_{\TT\cup\{\ub=1\}}|q|^{p-2}|\Phi|^2|\g(e_3, N)|.
\eea
 \end{lemma}
 
\begin{remark}
Recall that $r\leq r_0$ on $\Mint$ so that powers of $r$ are a priori irrelevant. They will be however useful in Proposition \ref{prop:maincontroloftransporteuqiaotningoingPTstructureMint} to produce an arbitrary large positive bulk term on the LHS, by choosing $p$ large enough in \eqref{eq:general-integrated-estimate-Mint}, to allow to absorb terms on the RHS. This is analogous to the use of  a Gronwall lemma for the transport equation \eqref{eq:TransportMint}.
\end{remark} 
 
 \begin{proof} 
 Let $p$ a real number. In view of Lemma \ref{Lemma:div-e_3}, we have
\beaa
\Div(  |q|^{p}|\Phi|^2 e_3)&=& e_3( |q|^{p}|\Phi|^2)-\frac{2r}{|q|^2} |q|^{p}|\Phi|^2+|q|^{p}|\Phi|^2\Gac.
\eeaa
Also, we have, since $e_3(|q|) =-\frac{r}{|q|}$, 
\beaa
e_3( |q|^{p}|\Phi|^2) &=& -pr|q|^{p-2}|\Phi|^2+2|q|^{p}\Re(\ov{\Phi}\c\nab_3\Phi)
\eeaa
so that for $\Phi$ satisfying \eqref{eq:TransportMint}, we have
\beaa
e_3( |q|^{p}|\Phi|^2) &=& -pr|q|^{p-2}|\Phi|^2+2|q|^{p}\Re(\ov{\Phi}\c F).
\eeaa
We deduce
\beaa
\Div(  |q|^{p}|\Phi|^2 e_3)&=&  -(p+2)r|q|^{p-2}|\Phi|^2+2|q|^{p}\Re(\ov{\Phi}\c F) +|q|^{p}|\Phi|^2\Gac
\eeaa
and hence, using in particular $\Gac=O(\ep)$, 
\beaa
\Div(  |q|^{p}|\Phi|^2 e_3)&\leq&  -\Big(p+1+O(\ep)\Big)r|q|^{p-2}|\Phi|^2+r^{-1}|q|^{p+2}|F|^2.
\eeaa
Integrating on $\Mint_*$, we infer, for $\ep>0$ small enough, 
\beaa
&&\int_{\Sint_*}|q|^{p}|\Phi|^2\g(e_3, N) + \int_\AA |q|^{p}|\Phi|^2\g(e_3, N)+ p\int_{\Mint_*}r|q|^{p-2}|\Phi|^2 \\
&\leq& \int_{\Mint_*}r^{-1}|q|^{p+2}|F|^2+\int_{\TT}|q|^{p}|\Phi|^2|\g(e_3, N)|+\int_{\{\ub=1\}}|q|^{p}|\Phi|^2|\g(e_3, N)|. 
\eeaa
Since both $\Sint_*$ and $\AA$ are spacelike, and taking the orientation into account, we have $\g(e_3, N)>0$ in both cases, and hence
\beaa
&&\int_{\Sint_*}|q|^{p}|\Phi|^2|\g(e_3, N)| + \int_\AA |q|^{p}|\Phi|^2|\g(e_3, N)|+ p\int_{\Mint_*}r^{-1}|q|^{p-2}|\Phi|^2 \\
&\leq& \int_{\Mint_*}r^{-1}|q|^{p+2}|F|^2+\int_{\TT}|q|^{p}|\Phi|^2|\g(e_3, N)|+\int_{\{\ub=1\}}|q|^{p}|\Phi|^2|\g(e_3, N)|
\eeaa
as stated.
\end{proof}

We now apply Lemma \ref{lemma:general-transport-Mint}  to the transport equations of Proposition \ref{prop:schamticformoflinearzizedstrucutreequationsofMint:bis}.
\begin{proposition}\lab{prop:maincontroloftransporteuqiaotningoingPTstructureMint}
If $\Phi$ is any quantity among $\Gac\setminus\{\Xbh\}$, i.e. $\Phi$ is any linearized Ricci or metric coefficient of the ingoing PT structure of $\Mint$ except $\Xbh$, then $\Phi$ satisfies, for $k\leq k_{large}+7$, 
\bea\label{eq:estimatesforPhiofthefirstypeinMint'}
\nn\int_{\Sint_*\cup\AA} |\dk^k\Phi|^2+ \int_{\Mint_*}|\dk^k\Phi|^2 &\les&  \Big(\ep_0+\Rkint_k+\Skint_{k-1}\Big)^2+ \int_{\Mint_*}|\Gac_k[\Phi]|^2\\
&&+\int_{\TT\cup\{\ub=1\}}|\dk^k\Phi|^2.
\eea
\end{proposition}

\begin{proof}
We start with the case where $\Phi$ is among the quantities $\trXbc$, $\Xh$, $\widecheck{\DD\cos\th}$, $\DD r$, $e_4(\cos\th)$, $\widecheck{e_4(r)}$. Then, according to Proposition \ref{prop:schamticformoflinearzizedstrucutreequationsofMint:bis}.
 $\Phi$ satisfies a transport equations of the following form
 \beaa
\nn\nab_3(\nab_{\Rhat}\dk^{k-1}\Phi) &=& O(1)\nab_{\Rhat}\dk^{k-1}\Phi +O(1)\nab\dk^{k-1}\Phi +\Gac_k[\Phi]\\
&&+\dk^{\leq k-1}(\Rc, \Gac)+\dk^{\leq k}(\Gac\c\Gac),
\eeaa
and 
\beaa
\nn\nab_3(\nab\dk^{k-1}\Phi) &=& O(1)\nab_{\Rhat}\dk^{k-1}\Phi +O(1)\nab\dk^{k-1}\Phi +\Gac_k[\Phi]\\
&&+\dk^{\leq k-1}(\Rc, \Gac)+\dk^{\leq k}(\Gac\c\Gac),
\eeaa
Then, applying \eqref{eq:general-integrated-estimate-Mint}, we obtain 
\beaa
\nn&&\int_{\Sint_*\cup\AA}|q|^{p}\Big(|\nab_\Rhat\dk^{k-1}\Phi|^2+|\nab\dk^{k-1}\Phi|^2\Big)|\g(e_3, N)| \\
&&+ p\int_{\Mint_*}r|q|^{p-2}\Big(|\nab_\Rhat\dk^{k-1}\Phi|^2+|\nab\dk^{k-1}\Phi|^2\Big) \\
&\leq& \int_{\Mint_*}O(1)r^{-1}|q|^{p+2}\Big(|\nab_\Rhat\dk^{k-1}\Phi|^2+|\nab\dk^{k-1}\Phi|^2\Big)\\
&&+ \int_{\Mint_*}r^{-1}|q|^{p+2}\Big|\Gac_k[\Phi]+\dk^{\leq k-1}(\Rc, \Gac)+\dk^{\leq k}(\Gac\c\Gac)\Big|^2\\
&&+\int_{\TT\cup\{\ub=1\}}|q|^{p-2}|\dk^k\Phi|^2|\g(e_3, N)|.
\eeaa
Note that $O(1)$ depends only on $k$ and the equation of $\Phi$, and hence not on $p$. Together with the fact that $r_+-\de_\HH\leq r\leq r_0$ in $\Mint$, there exists thus a constant $C$ independent of $p$ such that 
\beaa
\nn&&\int_{\Sint_*\cup\AA}|q|^{p}\Big(|\nab_\Rhat\dk^{k-1}\Phi|^2+|\nab\dk^{k-1}\Phi|^2\Big)|\g(e_3, N)| \\
&&+ p\int_{\Mint_*}r|q|^{p-2}\Big(|\nab_\Rhat\dk^{k-1}\Phi|^2+|\nab\dk^{k-1}\Phi|^2\Big) \\
&\leq& C\int_{\Mint_*}r|q|^{p-2}\Big(|\nab_\Rhat\dk^{k-1}\Phi|^2+|\nab\dk^{k-1}\Phi|^2\Big)\\
&&+ \int_{\Mint_*}r^{-1}|q|^{p+2}\Big|\Gac_k[\Phi]+\dk^{\leq k-1}(\Rc, \Gac)+\dk^{\leq k}(\Gac\c\Gac)\Big|^2\\
&&+\int_{\TT\cup\{\ub=1\}}|q|^{p-2}|\dk^k\Phi|^2|\g(e_3, N)|.
\eeaa
In particular, taking $p$ large enough, we may absorb the first term on the RHS by the LHS and obtain, using again  the fact that $r_+-\de_\HH\leq r\leq r_0$ in $\Mint$, 
\beaa
\nn&&\int_{\Sint_*\cup\AA}\Big(|\nab_\Rhat\dk^{k-1}\Phi|^2+|\nab\dk^{k-1}\Phi|^2\Big)|\g(e_3, N)| \\
&&+ \int_{\Mint_*}\Big(|\nab_\Rhat\dk^{k-1}\Phi|^2+|\nab\dk^{k-1}\Phi|^2\Big) \\
&\les& \int_{\Mint_*}\Big|\Gac_k[\Phi]+\dk^{\leq k-1}(\Rc, \Gac)+\dk^{\leq k}(\Gac\c\Gac)\Big|^2\\
&&+\int_{\TT\cup\{\ub=1\}}\Big(|\nab_\Rhat\dk^{k-1}\Phi|^2+|\nab\dk^{k-1}\Phi|^2\Big)|\g(e_3, N)|.
\eeaa
Also, using in particular Proposition \ref{prop:propertiesoftauusefulfortheoremM8:chap9} for the boundary $\Sint_*=\{\tau=\tau_*\}$, as well as the fact that $\AA=\{r=r_+-\de_\HH\}$ and $\TT\{r=r_0\}$, we infer
\beaa
|\g(e_3, N)|\les 1\quad\textrm{on}\quad\TT\cup\{\ub=1\}, \qquad |\g(e_3, N)|\gtrsim 1\quad\textrm{on}\quad\Sint_*\cup\AA
\eeaa
which yields
\beaa
\nn&&\int_{\Sint_*\cup\AA}\Big(|\nab_\Rhat\dk^{k-1}\Phi|^2+|\nab\dk^{k-1}\Phi|^2\Big) + \int_{\Mint_*}\Big(|\nab_\Rhat\dk^{k-1}\Phi|^2+|\nab\dk^{k-1}\Phi|^2\Big) \\
&\les& \int_{\Mint_*}\Big|\Gac_k[\Phi]+\dk^{\leq k-1}(\Rc, \Gac)+\dk^{\leq k}(\Gac\c\Gac)\Big|^2+\int_{\TT\cup\{\ub=1\}}|\dk^k\Phi|^2.
\eeaa

Next, we estimate the first term on the RHS. In view of the definition of $\Rkint_k$ and $\Skint_k$, and using the bootstrap assumption \eqref{eq:Bootstrap-Mint-M8}, we have, for $k\leq k_{large}+7$, 
\beaa
 \int_{\Mint_*}\Big|\dk^{\leq k-1}(\Rc, \Gac)+\dk^{\leq k}(\Gac\c\Gac)\Big|^2 &\les&  \Big(\ep^2+\Rkint_k+\Skint_{k-1}\Big)^2\\
 &\les& \Big(\ep_0+\Rkint_k+\Skint_{k-1}\Big)^2.
\eeaa
Plugging in the above, we deduce, for $k\leq k_{large}+7$, 
\beaa
\nn&&\int_{\Sint_*\cup\AA}\Big(|\nab_\Rhat\dk^{k-1}\Phi|^2+|\nab\dk^{k-1}\Phi|^2\Big) + \int_{\Mint_*}\Big(|\nab_\Rhat\dk^{k-1}\Phi|^2+|\nab\dk^{k-1}\Phi|^2\Big) \\
&\les&  \Big(\ep_0+\Rkint_k+\Skint_{k-1}\Big)^2+\int_{\Mint_*}|\Gac_k[\Phi]|^2+\int_{\TT\cup\{\ub=1\}}|\dk^k\Phi|^2.
\eeaa
Now, note that 
\beaa
|\dk^k\Phi| &\les& |\nab_\Rhat\dk^{k-1}\Phi|+|\nab\dk^{k-1}\Phi|+|\nab_3^k\Phi|+|\dk^{\leq k-1}\Phi|
\eeaa
and hence, for $k\leq k_{large}+7$, 
\beaa
\nn&&\int_{\Sint_*\cup\AA}|\dk^k\Phi|^2 + \int_{\Mint_*}|\dk^k\Phi|^2  \\
&\les&  \Big(\ep_0+\Rkint_k+\Skint_{k-1}\Big)^2+\int_{\Mint_*}|\Gac_k[\Phi]|^2+\int_{\TT\cup\{\ub=1\}}|\dk^k\Phi|^2\\
&&+\int_{\Sint_*\cup\AA}|\nab_3^k\Phi|^2 + \int_{\Mint_*}|\nab_3^k\Phi|^2
\eeaa
so that it still remains to control $\nab_3^k\Phi$. Using $\nab_3\Phi=\Gac+\Rc$, and iterating, we have
\beaa
\nab_3^k\Phi=\Gac+\dk^{\leq k-1}\Rc,
\eeaa
and hence, in view of the definition of $\Rkint_k$ and $\Skint_k$, and the trace theorem for $\Gac$, we have
\beaa
&&\int_{\Sint_*\cup\AA}|\nab_3^k\Phi|^2 + \int_{\Mint_*}|\nab_3^k\Phi|^2\\
 &\les& \int_{\Sint_*\cup\AA}|\Gac+\dk^{\leq k-1}\Rc|^2 + \int_{\Mint_*}|\Gac+\dk^{\leq k-1}\Rc|^2\\
 &\les& \Rkint_k^2+ \int_{\Mint_*}|\dk^{\leq 1}\Gac|^2\les  \Skint_1^2+\Rkint_k^2.
\eeaa
Since $\Skint_1\les \ep_0$ in view of \eqref{eq:initializationofiterationassumptioninproofThmM8}, we infer
\beaa
\int_{\Sint_*\cup\AA}|\nab_3^k\Phi|^2 + \int_{\Mint_*}|\nab_3^k\Phi|^2 &\les& \ep_0^2+\Rkint_k^2.
\eeaa
Plugging in the above estimate for $\dk^k\Phi$, we infer, for $k\leq k_{large}+7$, 
\beaa
\nn&&\int_{\Sint_*\cup\AA}|\dk^k\Phi|^2 + \int_{\Mint_*}|\dk^{k}\Phi|^2 \\
&\les&  \Big(\ep_0+\Rkint_k+\Skint_{k-1}\Big)^2+\int_{\Mint_*}|\Gac_k[\Phi]|^2+\int_{\TT\cup\{\ub=1\}}|\dk^k\Phi|^2
\eeaa
as stated.

It remains to treat the case where $\Phi$ is among the quantities $\Zc$, $\Hbc$, $\trXc$, $\omc$, $\Xi$, $\DD\hot\Jk$, $\widecheck{\ov{\DD}\c\Jk}$, $\widecheck{\nab_4\Jk}$. Then, according to Proposition \ref{prop:schamticformoflinearzizedstrucutreequationsofMint:bis}.
 $\Phi$ satisfies a transport equations of the following form
\beaa
\nn\nab_3(\nab_{\Rhat}\dk^{k-1}\Phi) &=& O(1)\nab_{\Rhat}\dk^{k-1}\Phi+O(1)\nab\dk^{k-1}\Phi +\Gac_k[\Phi] \\
&&+\nab_{\Rhat}\dk^{k-1}\Rc+\dk^{\leq k-1}(\Rc, \Gac)+\dk^{\leq k}(\Gac\c\Gac),
\eeaa
and
\beaa
\nn\nab_3(\nab\dk^{k-1}\Phi+\dk^{k-1}\Rc) &=& O(1)\nab_{\Rhat}\dk^{k-1}\Phi+O(1)\nab\dk^{k-1}\Phi +\Gac_k[\Phi] \\
&&+\nab_{\Rhat}\dk^{k-1}\Rc+\dk^{\leq k-1}(\Rc, \Gac)+\dk^{\leq k}(\Gac\c\Gac).
\eeaa
Proceeding exactly as in the first case, we obtain 
\beaa
\nn&&\int_{\Sint_*\cup\AA}|\dk^k\Phi+\dk^{k-1}\Rc|^2 + \int_{\Mint_*}|\dk^{k}\Phi+\dk^{k-1}\Rc|^2 \\
&\les&  \int_{\Mint_*}|\nab_{\Rhat}\dk^{k-1}\Rc|^2+\Big(\ep_0+\Rkint_k+\Skint_{k-1}\Big)^2+\int_{\Mint_*}|\Gac_k[\Phi]|^2\\
&&+\int_{\TT\cup\{\ub=1\}}|\dk^k\Phi|^2.
\eeaa
Since we have, in view of the definition of $\Rkint_k$, 
\beaa
&&\int_{\Sint_*\cup\AA}|\dk^{k-1}\Rc|^2 + \int_{\Mint_*}|\dk^{k-1}\Rc|^2 + \int_{\Mint_*}|\nab_{\Rhat}\dk^{k-1}\Rc|^2 \les (\Rkint_k)^2,
\eeaa
we infer, for $k\leq k_{large}+7$,  
\beaa
\nn&&\int_{\Sint_*\cup\AA}|\dk^k\Phi|^2 + \int_{\Mint_*}|\dk^{k}\Phi|^2 \\
&\les&  \Big(\ep_0+\Rkint_k+\Skint_{k-1}\Big)^2+\int_{\Mint_*}|\Gac_k[\Phi]|^2 +\int_{\TT\cup\{\ub=1\}}|\dk^k\Phi|^2
\eeaa
as stated. This concludes the proof of Proposition \ref{prop:maincontroloftransporteuqiaotningoingPTstructureMint}.
\end{proof}

%%%%%%%%%%%%%%%%%%%%%%%%%%%
 
  \subsection{Non-integrable Hodge estimates}

%%%%%%%%%%%%%%%%%%%%%%%%%%%
  
  Note that Proposition \ref{prop:maincontroloftransporteuqiaotningoingPTstructureMint} does not apply to $\Xbh$. In this section, we provide estimate for $\Xbh$ relying in particular on Codazzi. To this end, we start by considering 
 Hodge type systems  of the form
  \bea
  \ov{\DD}\c U &=& H
  \eea
   in $\Mint $ where $U$ is an anti-selfdual horizontal symmetric traceless 2-tensor, and $H$ is a anti-selfdual horizontal 1-form. We prove the following  non-integrable version of   elliptic-Hodge  estimates. 
      \begin{proposition}
      \lab{Prop:nonintegrableHodgeMint}
      Given $U$ is an anti-selfdual horizontal symmetric traceless 2-tensor,  we have the estimate 
     \Red{ \bea
      \lab{eq:nonintegrableHodgeMint}
\nn \int_{\Mint_* }  |\nab    U  |^2 &\les&   \int_{\Mint_* } \Big( |\ov{\DD} \c  U |^2  +|U||\dk U|  \Big)\\
 &&+\int_{\pr \Mint_* } \Big(|\nab_3 U|+|\nab_{\Rhat} U|+|\ov{\DD}\c  U|+|U|\Big)|U|.
      \eea}
      \end{proposition}
      
      \begin{proof}
       We rely on the following identity.    
\begin{lemma} 
\label{Le:2D-hodge-non-integrable}
Given a  horizontal structure, and given $U$ is an anti-selfdual horizontal symmetric traceless 2-tensor, the following point-wise identity  holds\footnote{Here  $ \Kh=  - \frac 14  \trch \trchb-\frac 1 4 \atrch \atrchb-  \rho + \Ga_g \cdot \Ga_b$  is a non-integrable version of the Gauss curvature appearing in \cite{GKS1}.}
\bea
\label{eq:hodgeident2-nonint}
\bsplit
|\nab    U  |^2+\Kh |U|^2 &=2 |\ov{\DD} \c  U  |^2 + \frac 1 2 \bigg(\Big(\atrch\nab_3+\atrchb \nab_4\Big) \dual U \bigg) \c U\\
 &+\nab_a \Big( \nab^a U \c U- \Red{(\ov{\DD}\c U)_b} U^{ab} \Big).
 \end{split}
\eea
\end{lemma}

\begin{proof}
See section 2 in \cite{GKS1}
\end{proof}

We rewrite Lemma \ref{Le:2D-hodge-non-integrable} schematically in $\Mint$.
\begin{corollary}
We have for  an anti-selfdual horizontal symmetric traceless 2-tensor $U$  in $\Mint$
\bea
\lab{eq:2D-hodge-non-integrable-CCC}
 |\nab    U  |^2&=& 2 |\ov{\DD} \c  U |^2  \Red{+O(1)|U||\dk U| } + \D_\a \Big(\D^a U \c U-  (\ov{\DD}\c  U )_b\c U^{\a b} \Big).
\eea
\end{corollary}

\begin{proof} 
In view of the following simple calculation for an anti-selfdual 1-form $V$
\beaa
\lab{eq:Div-div}
\D^\a V_\a&=& \nab^a V_a + (\eta +\etab) \c V,
\eeaa
this follows immediately  from Lemma   \ref{Le:2D-hodge-non-integrable}.  
\end{proof}
  
 We next   integrate  \eqref{eq:2D-hodge-non-integrable-CCC} in $\Mint_*$ and apply the divergence lemma to  derive
  \beaa
 \int_{\Mint_* }  |\nab    U  |^2 \les   \int_{\Mint_* } \Big( |\ov{\DD} \c  U |^2  \Red{+|U||\dk U| }  \Big)+\int_{\pr \Mint_* }  
   N_\a   \Big(\D^\a U \c U-  (\ov{\DD}\c  U )_b\c U^{\a b} \Big)
  \eeaa
    where  $\pr \Mint_*=\Sint_*\cup\AA\cup\{\ub=1\}\cup \TT $ and $N$ the corresponding normal. In order to decompose the first boundary term on the RHS, i.e. $N_\a\D^\a U\c U$, we introduce the vectorfields\footnote{Recall that $e_4(\ub)=\frac{2(r^2+a^2)}{|q|^2}+\Gac$ so that we may divide by $e_4(\ub)$.}  
\beaa
\ehat_b &=& e_b+ \left(e_b(r) -\frac{e_b(\ub)}{e_4(\ub)}e_4(r)\right)e_3 -\frac{e_b(\ub)}{e_4(\ub)}e_4,\quad b=1,2,
\eeaa
and note that with this definition, we have
\beaa
\ehat_1(\ub)=0, \qquad  \ehat_1(r)=0, \qquad \ehat_2(\ub)=0, \qquad \ehat_2(r)=0,
\eeaa
so that $(\ehat_1, \ehat_2)$ spans the tangent space of the spheres $S(\ub, r)$. We decompose
\beaa
N_\a   \D^\a U \c U &=& -\frac{1}{2}\g(N, e_4)\D_3 U\c U -\frac{1}{2}\g(N, e_3)\D_4 U\c U
+\g(N, e_b)\D_b U\c U\\
&=& O(1)\nab_3U\c U+O(1)\nab_4U\c U+\frac{1}{2}\g(N, e_b)\ehat_b(U\c U).
\eeaa
Now, since $\pr \Mint_*=\Sint_*\cup\AA\cup\{\ub=1\}\cup \TT $, and since $\Sint_*=\{\tau=\tau_*\}$, $\AA=\{r=r_+-\de_\HH\}$, $\TT=\{r=r_0\}$, and since $\tau=\ub+f(r)$ for some scalar function $f$ by construction, see the proof of  Proposition \ref{prop:propertiesoftauusefulfortheoremM8:chap9}, all parts of $\pr \Mint_*$ are foliated by spheres $S(\ub, r)$. In particular, we may integrate the term $g(N, e_b)\ehat_b(U\c U)$ by parts on $\pr \Mint_*$ and hence
 \beaa
 \int_{\pr \Mint_* }   N_\a   \D^\a U \c U &=& \int_{\pr \Mint_* }\Big(O(1)|\nab_3 U||U|+O(1)|\nab_4 U||U|+O(1)|U|^2\Big)\\
 &=& \int_{\pr \Mint_* }\Big(O(1)|\nab_3 U||U|+O(1)|\nab_{\Rhat} U||U|+O(1)|U|^2\Big).
 \eeaa
Plugging in the above, we infer 
\Red{ \beaa
 \int_{\Mint_* }  |\nab    U  |^2 &\les&   \int_{\Mint_* } \Big( |\ov{\DD} \c  U |^2 +|U||\dk U|   \Big)\\
 &&+\int_{\pr \Mint_* } \Big(|\nab_3 U|+|\nab_{\Rhat} U|+|\ov{\DD}\c  U|+|U|\Big)|U|
\eeaa}
as stated. This completes the  proof of Proposition \ref{Prop:nonintegrableHodgeMint}.
  \end{proof}
  
We are now ready to derive and estimate for $\Xbh$ on $\Mint_*$. 
\begin{proposition}\lab{prop:estimateforcoupledtransportCodazzisystemXbhinMint:chap9} 
$\Xbh$ satisfies, for $k\leq k_{large}+7$, 
\bea
\nn\int_{\Mint_*}|\dk^{k}\Xbh|^2 &\les& \int_{\Mint_* }|\dk^k\trXbc|^2+\int_{\Sint_*\cup\AA}\Big(|\dk^k\trXbc|^2+|\dk^{\leq k-1}(\Gac\setminus\{\Xbh\})|^2\Big)\\
&&+\int_{\TT\cup\{\ub=1\}}|\dk^{\leq k}\Gac|^2+\Big(\ep_0+\Rkint_k+\Skint_{k-1}\Big)^2.
\eea
\end{proposition} 
 
\begin{proof}
Recall from Proposition \ref{prop:schamticformoflinearzizedstrucutreequationsofMint:bis} that 
$\Xbh$ satisfies the following equation 
\beaa
\nn\nab_3(\nab_{\Rhat}\dk^{k-1}\Xbh) &=& O(1)\nab_{\Rhat}\dk^{k-1}\Xbh+O(1)\nab\dk^{k-1}\Xbh\\
&& +\nab_{\Rhat}\dk^{k-1}\Rc+\dk^{\leq k-1}(\Rc, \Gac)+\dk^{\leq k}(\Gac\c\Gac).
\eeaa
Then, applying \eqref{eq:general-integrated-estimate-Mint}, we obtain 
\beaa
\nn&&\int_{\Sint_*\cup\AA}|q|^{p}|\nab_{\Rhat}\dk^{k-1}\Xbh|^2|\g(e_3, N)| + p\int_{\Mint_*}r|q|^{p-2}|\nab_{\Rhat}\dk^{k-1}\Xbh|^2 \\
&\leq& \int_{\Mint_*}O(1)r^{-1}|q|^{p+2}|\nab_{\Rhat}\dk^{k-1}\Xbh|^2 +\int_{\Mint_*}O(1)r^{-1}|q|^{p+2}|\nab\dk^{k-1}\Xbh|^2\\
&&+\int_{\Mint_*}r^{-1}|q|^{p+2}\Big|\nab_{\Rhat}\dk^{k-1}\Rc+\dk^{\leq k-1}(\Rc, \Gac)+\dk^{\leq k}(\Gac\c\Gac)\Big|^2\\
&&+\int_{\TT\cup\{\ub=1\}}|q|^{p-2}|\nab_{\Rhat}\dk^{k-1}\Xbh|^2|\g(e_3, N)|.
\eeaa
Next, we choose $p$ large enough to absorb the first term on the RHS. Proceeding as in the proof of Proposition \ref{prop:maincontroloftransporteuqiaotningoingPTstructureMint}, we obtain 
\beaa
\nn&&\int_{\Sint_*\cup\AA}|\nab_{\Rhat}\dk^{k-1}\Xbh|^2 + \int_{\Mint_*}|\nab_{\Rhat}\dk^{k-1}\Xbh|^2 \\
&\les& \int_{\Mint_*}|\nab\dk^{k-1}\Xbh|^2+\int_{\TT\cup\{\ub=1\}}|\nab_{\Rhat}\dk^{k-1}\Xbh|^2\\
&&+\int_{\Mint_*}\Big|\nab_{\Rhat}\dk^{k-1}\Rc+\dk^{\leq k-1}(\Rc, \Gac)+\dk^{\leq k}(\Gac\c\Gac)\Big|^2.
\eeaa
In view of the definition of $\Rkint_k$ and $\Skint_k$, and using the bootstrap assumption \eqref{eq:Bootstrap-Mint-M8}, we infer, for $k\leq k_{large}+7$, 
\beaa
\nn&&\int_{\Sint_*\cup\AA}|\nab_{\Rhat}\dk^{k-1}\Xbh|^2 + \int_{\Mint_*}|\nab_{\Rhat}\dk^{k-1}\Xbh|^2 \\
&\les& \int_{\Mint_*}|\nab\dk^{k-1}\Xbh|^2+\int_{\TT\cup\{\ub=1\}}|\nab_{\Rhat}\dk^{k-1}\Xbh|^2+\Big(\ep_0+\Rkint_k+\Skint_{k-1}\Big)^2.
\eeaa

Next, we focus on controlling the first term on the RHS of the above estimate, i.e. the one involving $\nab\dk^{k-1}\Xbh$. We use Proposition \ref{Prop:nonintegrableHodgeMint} with $U=\dk^{k-1}\Xbh$ which yields
\beaa 
&& \int_{\Mint_* }  |\nab\dk^{k-1}\Xbh|^2\\ 
&\les&   \int_{\Mint_* } \Big( |\ov{\DD} \c\dk^{k-1}\Xbh |^2  +|\dk^{k-1}\Xbh||\dk^{k}\Xbh|  \Big)\\
 &&+\int_{\pr \Mint_* } \Big(|\nab_3\dk^{k-1}\Xbh|+|\nab_{\Rhat}\dk^{k-1}\Xbh|+|\ov{\DD}\c\dk^{k-1}\Xbh|+|\dk^{k-1}\Xbh|\Big)|\dk^{k-1}\Xbh|.
\eeaa
Also, recall from Proposition \ref{prop:schamticformoflinearzizedstrucutreequationsofMint:bis} that 
$\Xbh$ satisfies the following equation 
\beaa
\nab_3\dk^{k-1}\Xbh &=& \dk^{\leq k-1}(\Rc, \Gac)+\dk^{\leq k-1}(\Gac\c\Gac),\\
\ov{\DD}\c(\dk^{k-1}\Xbh) &=& \dk^k\trXbc +\dk^{\leq k-1}(\Rc, \Gac)+\dk^{\leq k-1}(\Gac\c\Gac).
\eeaa
Plugging in the above, we deduce 
\beaa 
&& \int_{\Mint_* }  |\nab\dk^{k-1}\Xbh|^2\\ 
&\les&   \int_{\Mint_* } \left( \Big|\dk^k\trXbc +\dk^{\leq k-1}(\Rc, \Gac)+\dk^{\leq k-1}(\Gac\c\Gac)\Big|^2 +|\dk^{k-1}\Xbh||\dk^{k}\Xbh|  \right)\\
 &&+\int_{\pr \Mint_* } \Big(|\dk^k\trXbc|+|\dk^{\leq k-1}(\Rc, \Gac)+\dk^{\leq k-1}(\Gac\c\Gac)|+|\nab_{\Rhat}\dk^{k-1}\Xbh| +|\dk^{k-1}\Xbh|\Big)|\dk^{k-1}\Xbh|.
\eeaa
In view of the definition of $\Rkint_k$ and $\Skint_k$, and using the bootstrap assumption \eqref{eq:Bootstrap-Mint-M8}, we infer, for $k\leq k_{large}+7$, 
\beaa 
&& \int_{\Mint_* }  |\nab\dk^{k-1}\Xbh|^2\\ 
&\les& \int_{\Mint_* } \left( |\dk^k\trXbc|^2 +|\dk^{k-1}\Xbh||\dk^{k}\Xbh|  \right)+\int_{\pr \Mint_* }\Big(|\dk^k\trXbc|^2+|\dk^{\leq k-1}(\Gac\setminus\{\Xbh\})|^2\Big)\\
 &&+\int_{\pr \Mint_* } \Big(|\nab_{\Rhat}\dk^{k-1}\Xbh| +|\dk^{k-1}\Xbh|\Big)|\dk^{k-1}\Xbh|+\Big(\ep_0+\Rkint_k+\Skint_{k-1}\Big)^2.
\eeaa
Also, using the trace theorem, we have
\beaa
\int_{\pr \Mint_* } |\dk^{k-1}\Xbh|^2 &\les&  \int_{\Mint_* }|\dk^{\leq k-1}\Xbh||\dk^{\leq k}\Xbh|. 
\eeaa
Together with 
\beaa
\int_{\Mint_* }|\dk^{\leq k-1}\Xbh||\dk^{\leq k}\Xbh| &\les& \Skint_{k-1}\left(\int_{\Mint_* }|\dk^{k}\Xbh|^2\right)^{\frac{1}{2}}+\Skint^2_{k-1},
\eeaa
we infer
\beaa 
&& \int_{\Mint_* }  |\nab\dk^{k-1}\Xbh|^2\\ 
&\les& \int_{\Mint_* }|\dk^k\trXbc|^2+\Skint_{k-1}\left(\int_{\Mint_* }|\dk^{k}\Xbh|^2\right)^{\frac{1}{2}}\\
&&+(\Skint_{k-1})^{\frac{1}{2}}\left(\int_{\Mint_* }|\dk^{k}\Xbh|^2\right)^{\frac{1}{4}}\left(\int_{\pr \Mint_* }|\nab_{\Rhat}\dk^{k-1}\Xbh|^2\right)^{\frac{1}{2}}\\
&&+\int_{\pr \Mint_* }\Big(|\dk^k\trXbc|^2+|\dk^{\leq k-1}(\Gac\setminus\{\Xbh\})|^2\Big)+\Big(\ep_0+\Rkint_k+\Skint_{k-1}\Big)^2
\eeaa
and hence, for $k\leq k_{large}+7$,
\beaa 
&& \int_{\Mint_* }  |\nab\dk^{k-1}\Xbh|^2\\ 
&\les& \int_{\Mint_* }|\dk^k\trXbc|^2+\Skint_{k-1}\left(\int_{\Mint_* }|\dk^{k}\Xbh|^2\right)^{\frac{1}{2}}\\
&&+(\Skint_{k-1})^{\frac{1}{2}}\left(\int_{\Mint_* }|\dk^{k}\Xbh|^2\right)^{\frac{1}{4}}\left(\int_{\Sint_*\cup\AA}|\nab_{\Rhat}\dk^{k-1}\Xbh|^2\right)^{\frac{1}{2}}\\
&&+\int_{\Sint_*\cup\AA}\Big(|\dk^k\trXbc|^2+|\dk^{\leq k-1}(\Gac\setminus\{\Xbh\})|^2\Big)+\Big(\ep_0+\Rkint_k+\Skint_{k-1}\Big)^2\\
&&+\int_{\TT\cup\{\ub=1\}}|\dk^{\leq k}\Gac|^2.
\eeaa

Next, recall the following above estimate, for $k\leq k_{large}+7$, 
\beaa
\nn&&\int_{\Sint_*\cup\AA}|\nab_{\Rhat}\dk^{k-1}\Xbh|^2 + \int_{\Mint_*}|\nab_{\Rhat}\dk^{k-1}\Xbh|^2 \\
&\les& \int_{\Mint_*}|\nab\dk^{k-1}\Xbh|^2+\int_{\TT\cup\{\ub=1\}}|\nab_{\Rhat}\dk^{k-1}\Xbh|^2+\Big(\ep_0+\Rkint_k+\Skint_{k-1}\Big)^2.
\eeaa
Plugging the above estimate for $\nab\dk^{k-1}\Xbh$, we obtain, for $k\leq k_{large}+7$,  
\beaa
\nn&&\int_{\Sint_*\cup\AA}|\nab_{\Rhat}\dk^{k-1}\Xbh|^2 + \int_{\Mint_*}|\nab_{\Rhat}\dk^{k-1}\Xbh|^2+ \int_{\Mint_*}|\nab\dk^{k-1}\Xbh|^2 \\
&\les& \int_{\Mint_* }|\dk^k\trXbc|^2+\Skint_{k-1}\left(\int_{\Mint_* }|\dk^{k}\Xbh|^2\right)^{\frac{1}{2}}\\
&&+(\Skint_{k-1})^{\frac{1}{2}}\left(\int_{\Mint_* }|\dk^{k}\Xbh|^2\right)^{\frac{1}{4}}\left(\int_{\Sint_*\cup\AA}|\nab_{\Rhat}\dk^{k-1}\Xbh|^2\right)^{\frac{1}{2}}\\
&&+\int_{\Sint_*\cup\AA}\Big(|\dk^k\trXbc|^2+|\dk^{\leq k-1}(\Gac\setminus\{\Xbh\})|^2\Big)\\
&&+\int_{\TT\cup\{\ub=1\}}|\dk^{\leq k}\Gac|^2+\Big(\ep_0+\Rkint_k+\Skint_{k-1}\Big)^2.
\eeaa
We use the boundary term for $\nab_{\Rhat}\dk^{k-1}\Xbh$ on the LHS to absorb the corresponding term on the RHS and infer
\beaa
\nn&& \int_{\Mint_*}|\nab_{\Rhat}\dk^{k-1}\Xbh|^2+ \int_{\Mint_*}|\nab\dk^{k-1}\Xbh|^2 \\
&\les& \int_{\Mint_* }|\dk^k\trXbc|^2+\Skint_{k-1}\left(\int_{\Mint_* }|\dk^{k}\Xbh|^2\right)^{\frac{1}{2}}\\
&&+\int_{\Sint_*\cup\AA}\Big(|\dk^k\trXbc|^2+|\dk^{\leq k-1}(\Gac\setminus\{\Xbh\})|^2\Big)\\
&&+\int_{\TT\cup\{\ub=1\}}|\dk^{\leq k}\Gac|^2+\Big(\ep_0+\Rkint_k+\Skint_{k-1}\Big)^2.
\eeaa
Using the fact that $\Xbh$ satisfies 
\beaa
\nab_3\dk^{k-1}\Xbh &=& \dk^{\leq k-1}\Gac+\dk^{\leq k-1}\Rc+\dk^{\leq k-1}(\Gac\c\Gac),
\eeaa
and using again the definition of $\Rkint_k$ and $\Skint_k$, and  the bootstrap assumption \eqref{eq:Bootstrap-Mint-M8}, we infer, for $k\leq k_{large}+7$, 
\beaa
\nn&& \int_{\Mint_*}|\nab_{\Rhat}\dk^{k-1}\Xbh|^2+ \int_{\Mint_*}|\nab\dk^{k-1}\Xbh|^2+ \int_{\Mint_*}|\nab_3\dk^{k-1}\Xbh|^2 \\
&\les& \int_{\Mint_* }|\dk^k\trXbc|^2+\Skint_{k-1}\left(\int_{\Mint_* }|\dk^{k}\Xbh|^2\right)^{\frac{1}{2}}\\
&&+\int_{\Sint_*\cup\AA}\Big(|\dk^k\trXbc|^2+|\dk^{\leq k-1}(\Gac\setminus\{\Xbh\})|^2\Big)\\
&&+\int_{\TT\cup\{\ub=1\}}|\dk^{\leq k}\Gac|^2+\Big(\ep_0+\Rkint_k+\Skint_{k-1}\Big)^2.
\eeaa
Since $\nab$, $\nab_{\Rhat}$ and $\nab_3$ span all derivatives, we deduce
\beaa
\int_{\Mint_*}|\dk^{k}\Xbh|^2 &\les& \int_{\Mint_* }|\dk^k\trXbc|^2+\Skint_{k-1}\left(\int_{\Mint_* }|\dk^{k}\Xbh|^2\right)^{\frac{1}{2}}\\
&&+\int_{\Sint_*\cup\AA}\Big(|\dk^k\trXbc|^2+|\dk^{\leq k-1}(\Gac\setminus\{\Xbh\})|^2\Big)\\
&&+\int_{\TT\cup\{\ub=1\}}|\dk^{\leq k}\Gac|^2+\Big(\ep_0+\Rkint_k+\Skint_{k-1}\Big)^2.
\eeaa
We may absorb the second term on the RHS by the LHS and obtain, for $k\leq k_{large}+7$, 
\beaa
\int_{\Mint_*}|\dk^{k}\Xbh|^2 &\les& \int_{\Mint_* }|\dk^k\trXbc|^2+\int_{\Sint_*\cup\AA}\Big(|\dk^k\trXbc|^2+|\dk^{\leq k-1}(\Gac\setminus\{\Xbh\})|^2\Big)\\
&&+\int_{\TT\cup\{\ub=1\}}|\dk^{\leq k}\Gac|^2+\Big(\ep_0+\Rkint_k+\Skint_{k-1}\Big)^2.
\eeaa
as stated. This concludes the proof of Proposition \ref{prop:estimateforcoupledtransportCodazzisystemXbhinMint:chap9}.
\end{proof}

%%%%%%%%%%%%%%%%%%%%%%%%%%%%%%%%%
 
  \subsection{Estimates for  the PT frame of $\Mint'$ on $\TT$}
  \label{section:PTframeMint-onTT}

%%%%%%%%%%%%%%%%%%%%%%%%%%%%%%%%%

The following lemma provides the control of the ingoing PT structure of $\Mint$ on $\TT$ from the outgoing PT structure of $\Mext$.

\begin{lemma}\lab{Lemma:EstimatesonTT-Mint}
The following estimates hold on  $\TT$  for the ingoing PT structure of $\Mint$  
\bea
\int_{\TT}\left|\dk^{\leq k}\Gacint\right|^2 \les \Big(\ep_0+\Skext_k+\Rkint_k\Big)^2, \qquad k\leq k_{large}+7.
\eea
\end{lemma}

\begin{proof}
The proof is the analog for the PT frames of $\Mext$ and $\Mint$ of the one of Lemma \ref{lemma:controlofPGstructureMintonTT} for the PG frames of $\Mext$ and $\Mint$. To simplify  the notations, in this proof, we denote:
\begin{itemize}
\item by $(e_4, e_3, e_1, e_2)$ the outgoing PT frame of $\Mext$, with all quantities associated to the outgoing PT structure of $\Mext$ being unprimed, 

\item by $(e_4', e_3', e_1', e_2')$ the ingoing PT frame of $\Mint$, with all quantities associated to the ingoing PT structure of $\Mint$ being primed. 
\end{itemize}

Recall that $\Mext\cap\Mint=\TT=\{r=r_0\}$. In view of the above notations, and the initialization of the ingoing PT structure of $\Mint$ from the outgoing PG structure of $\Mext$ on $\TT$, see section \ref{sec:defintionofthePTstructuresinMM}, we have  
\beaa
\ub=u, \qquad r'=r, \qquad \th'=\th,  \qquad \Jk'=\Jk\quad\quad\textrm{on}\quad\TT,
\eeaa
and
\beaa
e_4'=\la e_4, \qquad e_3'=\la^{-1}e_3, \qquad e_a'=e_a, \,\, a=1,2, \quad\textrm{on}\quad\TT,
\eeaa
where $\la$ is given by
\beaa
\la=\frac{\Delta}{|q|^2}.
\eeaa

Based on this initialization we derive in particular,   for any tangent vector $X$ on $\TT$, 
\beaa
\g(\D_Xe_4', e_a') &=& \g\left(\D_X\left(\la e_4\right), e_a\right)=\la\g\left(\D_Xe_4, e_a\right),\\
\g(\D_Xe_3', e_a') &=& \g\left(\D_X\left(\la^{-1}e_3\right), e_a\right)=\la^{-1}\g\left(\D_Xe_3, e_a\right),\\
\g(\D_Xe_4', e_3') &=& \g\left(\D_X\left(\la e_4\right), \la^{-1}e_3\right)=-2X(\log\la)+\g\left(\D_Xe_4, e_3\right).
\eeaa
Note that any such tangent   vector $X$  is a linear combination of  the following  three tangent directions to  $\TT$,
\beaa
X_1=e_1 - e_1(r)e_4, \qquad X_2=e_2 - e_2(r)e_4, \qquad X_3=e_3- e_3(r)e_4.
\eeaa
Also, note that we have
\beaa
  \nab(r)=\Gac, \qquad   e_3(r)=-\la+\Gac.
\eeaa
so that 
\beaa
X_1=e_1+\Gac e_4, \qquad X_2=e_2+\Gac e_4, \qquad X_3=e_3+\la e_4+\Gac e_4,
\eeaa
as well as 
\beaa
X_1=e_1'+\Gac e_4', \qquad X_2=e_2'+\Gac e_4', \qquad X_3=\la e_3'+e_4'+\Gac e_4'.
\eeaa
We infer
\beaa
\chi_{ba}' &=& \la\chi_{ba}+\Gac,\\
\chib_{ba}' &=& \la^{-1}\chib_{ba}+\Gac,\\
\ze_b' &=& -e_b(\log\la)+\ze_b+\Gac,
\eeaa
and
\beaa
\la\eta_a'+\xi_a' &=& \la\eta_a+\la^2\xi+\Gac,\\
\la\xib_a'+\etab_a' &=& \la^{-1}\xib_a+\etab_a+\Gac,\\
-\la\omb'+\om' &=&  -\frac{1}{2}e_3(\log\la)-\omb+\la\om+\Gac.
\eeaa
Together with
\begin{itemize}
\item the fact that, by the PT gauge choices, $\xib'=0$, $\omb'=0$, $H'=\frac{aq'}{|q'|^2}$, $\xi=0$, $\om=0$, $\Hb=-\frac{a\ov{q}}{|q|^2}\Jk$, 

\item the form of $\la$, 

\item the fact that $r'=r$, $\th'=\th$, and $\Jk'=\Jk$ on $\TT$, 

\item the definition of the linearized quantities for the PT frame in $\Mint$ and $\Mext$, 
\end{itemize}
we infer\footnote{See Lemma \ref{lemma:relationsPGMintandMextonTT} for more precise formulas concerning the analog situation for PG structures.} on $\TT$
\beaa
\trXc', \quad \Xh', \quad \trXbc',\quad \Xbh', \quad \Zc', \quad \Xi', \quad \Hbc', \quad \om'=\Gac.
\eeaa
In view of the definition of $\Skext_k$, we have
\beaa
\int_\TT|\dk^{\leq k}\Gac|^2 &\les& \Skext_k^2.
\eeaa
We thus infer
\beaa
\int_{\TT}\left|\widetilde{\dk}^{\leq k}\Big(\trXc', \, \Xh', \, \trXbc',\, \Xbh', \, \Zc', \, \Xi', \, \Hbc', \, \om'\Big)\right|^2 &\les& \Skext_k^2,
\eeaa
where $\widetilde{\dk}^{\leq k}$ denote tangential derivatives to $\TT$. Since $\dk$ is generated by $e_3'$ and $\widetilde{\dk}$, we obtain
\beaa
\int_{\TT}\left|\dk^{\leq k}\Big(\trXc', \, \Xh', \, \trXbc',\, \Xbh', \, \Zc', \, \Xi', \, \Hbc', \, \om'\Big)\right|^2 \les \Skext_k^2+\int_{\TT}\Big(|\dk^{k-1}\nab_3'\Gac'|^2+|\dk^{k-1}\Gac'|^2\Big).
\eeaa

Concerning the linearized metric coefficients, we consider 
\beaa
\ub=u, \qquad r'=r, \qquad \th'=\th,  \qquad \Jk'=\Jk\quad\quad\textrm{on}\quad\TT,
\eeaa
and apply to both sides of these identities the tangential vectorfields $X$ to $\TT$ as above. Then, proceeding 
similarly to the linearized Ricci coefficients\footnote{See Lemma \ref{lemma:relationsPGMintandMextonTT} for the analog situation for PG structures.}, we obtain on $\TT$ relations between the linearized metric coefficients on $\Mint$ and $\Mext$, which then, together with the above control for linearized Ricci coefficients, yields 
\beaa
\int_{\TT}\left|\dk^{\leq k}\Gac'\right|^2 \les \Skext_k^2+\int_{\TT}\Big(|\dk^{k-1}\nab_3'\Gac'|^2+|\dk^{k-1}\Gac'|^2\Big).
\eeaa

Then, since the PT structure equations take the form $\nab_3\Gac=\Gac+\Rc+\Gac\c\Gac$, and using the bootstrap assumption \eqref{eq:Bootstrap-Mint-M8}, we infer, for $k\leq k_{large}+7$, 
\beaa
\int_{\TT}\left|\dk^{\leq k}\Gac'\right|^2 \les \Skext_k^2+\int_{\TT}\Big(|\dk^{\leq k-1}\Rc'|^2+|\dk^{\leq k-1}\Gac'|^2\Big).
\eeaa
Together with the trace theorem and the definition of $\Rkint_k$, we obtain, for $k\leq k_{large}+7$, 
\beaa
\int_{\TT}\left|\dk^{\leq k}\Gac'\right|^2 \les \Skext_k^2+\Rkint_k^2+\int_{\TT}|\dk^{\leq k-1}\Gac'|^2.
\eeaa
Arguing by iteration on $k$, we deduce, for $k\leq k_{large}+7$, 
\beaa
\int_{\TT}\left|\dk^{\leq k}\Gac'\right|^2 \les \Skext_k^2+\Rkint_k^2+\int_{\TT}|\Gac'|^2.
\eeaa
Together with \eqref{eq:initializationofiterationassumptioninproofThmM8}, we infer, for $k\leq k_{large}+7$, 
\beaa
\int_{\TT}\left|\dk^{\leq k}\Gac'\right|^2 \les \Big(\ep_0+\Skext_k+\Rkint_k\Big)^2.
\eeaa
as desired.
\end{proof}

%%%%%%%%%%%%%%%%%%%%%%%%%%%%%%%%% 
    
  \subsection{Proof of Proposition \ref{Prop:MainestimatesMint}}
  \lab{sec:proofofProp:MainestimatesMint666}
  
%%%%%%%%%%%%%%%%%%%%%%%%%%%%%%%%%  

We are now ready to control $\Skint_k$ for $k\leq k_{large}+7$.

%%%%%%%%%%%%%%%%%%%%%%%%%%%%%%%%% 
    
  \subsubsection{Iteration assumption}
  
%%%%%%%%%%%%%%%%%%%%%%%%%%%%%%%%%  

Note first from \eqref{eq:initializationofiterationassumptioninproofThmM8} that we have
\bea\lab{eq:initializationofiterationassumptioninproofThmM8:Mintcase}
\Skint_{k_{small}-2}+\left(\int_{\Sint_*\cup\AA}|\dk^{\leq k_{small}-2}(\Gac\setminus\{\Xbh\})|^2\right)^{\frac{1}{2}}  &\les&\ep_0.
 \eea
This allows us to prove Proposition \ref{Prop:MainestimatesMint} by iteration. To this end, consider the following  iteration assumption
\bea\lab{eq:iterationassumptionforcontrolMintinThmM8}
\Skint_{k-1}+\left(\int_{\Sint_*\cup\AA}|\dk^{\leq k-1}(\Gac\setminus\{\Xbh\})|^2\right)^{\frac{1}{2}} \les \ep_0+\Skext_{k-1}+\Rkint_{k-1}. 
\eea
In view of \eqref{eq:initializationofiterationassumptioninproofThmM8:Mintcase}, \eqref{eq:iterationassumptionforcontrolMintinThmM8} holds for $k=k_{small}-1$. From now on, we assume \eqref{eq:iterationassumptionforcontrolMintinThmM8} for $k_{small}\leq k\leq k_{large}+7$. The proof of Proposition \ref{Prop:MainestimatesMint} will follow from proving \eqref{eq:iterationassumptionforcontrolMintinThmM8} for $k$ replaced by $k+1$.

%%%%%%%%%%%%%%%%%%%%%%%%%%%%%%%%%%%%%%% 
    
  \subsubsection{Control on $\Mint\setminus\Mint_*$ and on $1\leq\ub\leq 2$}
  
%%%%%%%%%%%%%%%%%%%%%%%%%%%%%%%%%%%%%%%%  

We start with controlling the solution in $\{1\leq\ub\leq 2\}$ and in $\Mint\setminus\Mint_*$. Note from the definition of $\Mint_*=\{\tau=\tau_*\}$ and $\tau_*$, see Definition \ref{def:u=introductionofregionMintprimestar:chap9}, and the properties of the scalar function $\tau$ constructed in Proposition \ref{prop:propertiesoftauusefulfortheoremM8:chap9} that we have
\bea\lab{eq:inclusionofMintminusMiintstartinregionRubaub2}
\Mint\setminus\Mint_*\subset\{u_*-4(m+1)\leq\ub\leq u_*\}.
\eea
We thus need in this step to control the solution on regions of the type 
\bea\lab{eq:regionRRub1ub2Mint:chap9}
\RR_{\ub_1,\ub_2} :=\big\{\ub_1\leq \ub\leq \ub_2\big\}, \,\,\, 1\leq\ub_1<\ub_1+1\leq\ub_2\leq u_*,\,\,\, \ub_2-\ub_1\leq 4(m+1).
\eea
 
We use the following lemma.
\begin{lemma}\lab{lemma:controlonMintminusMintstarbasiclemma:chap9}
Let $\Phi$ satisfying on $\Mint$ an equation of the type 
\bea
\nab_3(\dk^k\Phi) &=& O(1)\dk^k\Phi+\dk^{\leq k-1}\Gac+\dk^{\leq k}\Rc+\dk^{\leq k}(\Gac\c\Gac).
\eea
Then, $\Phi$ satisfies on regions $\RR_{\ub_1,\ub_2}$ as in \eqref{eq:regionRRub1ub2Mint:chap9}, for $k\leq k_{large}+7$, 
\bea
\int_{\RR_{\ub_1,\ub_2}}|\dk^k\Phi|^2 \les \int_{\TT}|\dk^k\Phi|^2+\Big(\ep_0+\Skint_{k-1}+\Rkint_k\Big)^2.
\eea
\end{lemma} 
 
\begin{proof}
Let $\ub_1\leq \ub\leq \ub_2$. We integrate the transport equation in $r$ along the level set of $\ub$ which yields, after applying Gronwall lemma and using the fact that $r\leq r_0$ in $\Mint$, 
\beaa
\sup_{r_+-\de_H\leq r\leq r_0}\int_{S(\ub, r)}|\dk^k\Phi|^2 &\les& \int_{S(\ub, r_0)}|\dk^k\Phi|^2+\int_{r_+-\de_\HH}^{r_0}\int_{S(r,\ub)}\Big|\dk^{\leq k-1}\Gac+\dk^{\leq k}\Rc+\dk^{\leq k}(\Gac\c\Gac)\Big|^2.
\eeaa
Since $r\leq r_0$ in $\Mint$, we infer
\beaa
\int_{r_+-\de_\HH}^{r_0}\int_{S(\ub, r)}|\dk^k\Phi|^2 &\les& \int_{S(\ub, r_0)}|\dk^k\Phi|^2+\int_{r_+-\de_\HH}^{r_0}\int_{S(r,\ub)}\Big|\dk^{\leq k-1}\Gac+\dk^{\leq k}\Rc+\dk^{\leq k}(\Gac\c\Gac)\Big|^2.
\eeaa
Integrating in $\ub$ for $\ub\in(\ub_1, \ub_2)$, we infer
\beaa
\int_{\RR_{\ub_1,\ub_2}}|\dk^k\Phi|^2 &\les& \int_{\TT}|\dk^k\Phi|^2+\int_{\RR_{\ub_1,\ub_2}}\Big|\dk^{\leq k-1}\Gac+\dk^{\leq k}\Rc+\dk^{\leq k}(\Gac\c\Gac)\Big|^2.
\eeaa
Then, using  the definition of $\Rkint_k$ and $\Skint_k$, and  the bootstrap assumption \eqref{eq:Bootstrap-Mint-M8}, we have, for $k\leq k_{large}+7$, 
\beaa
&&\int_{\RR_{\ub_1,\ub_2}}\Big|\dk^{\leq k-1}\Gac+\dk^{\leq k}\Rc+\dk^{\leq k}(\Gac\c\Gac)\Big|^2\\ 
&\les& \Skint_{k-1}^2+\ep_0^2+\Big(\max_{\RR_{\ub_1,\ub_2}}\tau - \max_{\RR_{\ub_1,\ub_2}}\tau\Big)\sup_{\tau}\int_{\Si(\tau)\cap\Mint}|\dk^{\leq k}\Rc|^2\\
&\les& \Big(\ep_0+\Skint_{k-1}+\Rkint\Big)^2
\eeaa
where we have also use the fact that $\ub_2-\ub_1\les 1$ by the definition of $\RR_{\ub_1,\ub_2}$ and that the variation of $\tau$ is controlled by the variation of $\ub$ on $\Mint$ in view of Proposition \ref{prop:propertiesoftauusefulfortheoremM8:chap9}. We thus infer, for $k\leq k_{large}+7$, 
\beaa
\int_{\RR_{\ub_1,\ub_2}}|\dk^k\Phi|^2 \les \int_{\TT}|\dk^k\Phi|^2+\Big(\ep_0+\Skint_{k-1}+\Rkint_k\Big)^2
\eeaa
as desired. This concludes the proof of Lemma \ref{lemma:controlonMintminusMintstarbasiclemma:chap9}.
\end{proof} 

The set of all linearized Ricci and metric coefficients $\Gac$ satisfies, in view of the equations in the ingoing PT structure, 
\beaa
\nab_3(\dk^k\Gac) &=& O(1)\dk^k\Phi+\dk^{\leq k-1}\Gac+\dk^{\leq k}\Rc+\dk^{\leq k}(\Gac\c\Gac).
\eeaa
We may thus apply Lemma \ref{lemma:controlonMintminusMintstarbasiclemma:chap9} which yields on regions $\RR_{\ub_1,\ub_2}$ as in \eqref{eq:regionRRub1ub2Mint:chap9}, for $k\leq k_{large}+7$,
 \beaa
\int_{\RR_{\ub_1,\ub_2}}|\dk^k\Gac|^2 &\les& \int_{\TT}|\dk^k\Gac|^2+\Big(\ep_0+\Skint_{k-1}+\Rkint_k\Big)^2.
\eeaa
 Together with the iteration assumption \eqref{eq:iterationassumptionforcontrolMintinThmM8} and the control on $\TT$ for $\Gac$ provided by Lemma \ref{Lemma:EstimatesonTT-Mint}, we infer, for $k\leq k_{large}+7$, 
 \bea\lab{eq:estimatestoappytotriangularstructureforregionRub1ub2}
\int_{\RR_{\ub_1,\ub_2}}|\dk^k\Gac|^2 &\les& \Big(\ep_0+\Skext_k+\Rkint_k\Big)^2.
\eea
In particular, since $\Mint\setminus\Mint_*$ is included in  $\RR_{\ub_1,\ub_2}$ with $\ub_1=u_*-4(m+1)$ and $\ub_2=u_*$ in view of \eqref{eq:inclusionofMintminusMiintstartinregionRubaub2}, \eqref{eq:estimatestoappytotriangularstructureforregionRub1ub2} yields
\bea\lab{eq:controlofGaconMintminusMintstrat}
\int_{\Mint\setminus\Mint_*}|\dk^k\Gac|^2 &\les& \Big(\ep_0+\Skext_k+\Rkint_k\Big)^2, \quad k\leq k_{large}+7.
\eea

Also, choosing $\ub_1=1$ and $\ub_2=2$, we have in view of \eqref{eq:estimatestoappytotriangularstructureforregionRub1ub2} 
 \bea\lab{eq:controlonubbetween1and2onMint}
\int_{1\leq\ub\leq 2}|\dk^k\Gac|^2 &\les& \Big(\ep_0+\Skext_k+\Rkint_k\Big)^2, \quad k\leq k_{large}+7.
\eea
Let $1\leq\ub_0\leq 2$ such that 
\beaa
\int_{\ub=\ub_0}|\dk^k\Gac|^2 &=& \inf_{1\leq\widetilde{\ub}\leq 2}\int_{\ub=\widetilde{\ub}}|\dk^k\Gac|^2.
\eeaa
Together with \eqref{eq:controlonubbetween1and2onMint}, this yields
\beaa
\int_{\ub=\ub_0}|\dk^k\Gac|^2  &\leq& \int_1^2\int_{\ub=\widetilde{\ub}}|\dk^k\Gac|^2 d\widetilde{\ub}\\
&\les& \int_{1\leq\ub\leq 2}|\dk^k\Gac|^2\\
&\les& \Big(\ep_0+\Skext_k+\Rkint_k\Big)^2
\eeaa
so that 
\bea\lab{eq:ub0between1and2realizesinfimumandthusundercontrolLebesguepoint}
\int_{\ub=\ub_0}|\dk^k\Gac|^2  &\les&  \Big(\ep_0+\Skext_k+\Rkint_k\Big)^2, \quad k\leq k_{large}+7. 
\eea
In view of \eqref{eq:ub0between1and2realizesinfimumandthusundercontrolLebesguepoint}, we should:
\begin{itemize}
\item reduce $\Mint_*$ to $\Mint_*\cap\{\ub\geq\ub_0\}$ so that the boundary terms  on $\{\ub=\ub_0\}$ in the  integration by parts of Propositions \ref{prop:maincontroloftransporteuqiaotningoingPTstructureMint} and  \ref{prop:estimateforcoupledtransportCodazzisystemXbhinMint:chap9} are under control,

\item control the remaining region $\{1\leq\ub\leq \ub_0\}$ thanks to \eqref{eq:controlonubbetween1and2onMint}.
\end{itemize} 
For simplicity, we  pretend that $\ub_0=1$, i.e. we assume from now by a slight abuse that we have \eqref{eq:ub0between1and2realizesinfimumandthusundercontrolLebesguepoint} with $\ub_0=1$ and hence
\bea\lab{eq:ub0between1and2realizesinfimumandthusundercontrolLebesguepoint:ub0=1}
\int_{\ub=1}|\dk^k\Gac|^2  &\les&  \Big(\ep_0+\Skext_k+\Rkint_k\Big)^2, \quad k\leq k_{large}+7. 
\eea

%%%%%%%%%%%%%%%%%%%%%%%%%%%%%%%%%%%% 
    
  \subsubsection{Control on $\Mint_*$}
  
%%%%%%%%%%%%%%%%%%%%%%%%%%%%%%%%%%%%  

Note that \eqref{eq:controlofGaconMintminusMintstrat} provides the desired control of $\Gac$ on $\Mint\setminus\Mint_*$. It thus remains to control $\Gac$ on $\Mint_*$. We start with the following corollary of Propositions \ref{prop:maincontroloftransporteuqiaotningoingPTstructureMint} and  \ref{prop:estimateforcoupledtransportCodazzisystemXbhinMint:chap9}.

\begin{corollary}\lab{cor:maincontroloftransporteuqiaotningoingPTstructureMint}
$\Phi$ is any quantity among $\Gac\setminus\{\Xbh\}$, i.e. $\Phi$ is any linearized Ricci or metric coefficient of the ingoing PT structure of $\Mint$, then $\Phi$ satisfies, for $k\leq k_{large}+7$, 
\bea\lab{eq:finalestimateinMintstarusingeverythingforGacminusXbh}
\nn\int_{\Sint_*\cup\AA} |\dk^k\Phi|^2+ \int_{\Mint_*}|\dk^k\Phi|^2 &\les&  \Big(\ep_0+\Rkint_k+\Skext_{k}\Big)^2\\
&&+ \int_{\Mint_*}|\Gac_k[\Phi]|^2
\eea

Moreover, the remaining component $\Xbh$ satisfies, for $k\leq k_{large}+7$, 
\bea\lab{eq:finalestimateinMintstarusingeverythingforXbh}
\nn\int_{\Mint_*}|\dk^{k}\Xbh|^2 &\les& \int_{\Mint_* }|\Gac_k[\Xbh]|^2+\int_{\Sint_*\cup\AA}|\Gac_k[\Xbh]|^2\\
&&+\Big(\ep_0+\Rkint_k+\Skext_k\Big)^2
\eea 
with the notation 
\bea
\Gac_k[\Xbh] &:=& \dk^k\trXc.
\eea
\end{corollary}

\begin{proof}
Consider first the case of $\Phi$ in $\Gac\setminus\{\Xbh\}$. Then, according to Proposition \ref{prop:maincontroloftransporteuqiaotningoingPTstructureMint}, $\Phi$ satisfies, for $k\leq k_{large}+7$, 
\beaa
\nn\int_{\Sint_*\cup\AA} |\dk^k\Phi|^2+ \int_{\Mint_*}|\dk^k\Phi|^2 &\les&  \Big(\ep_0+\Rkint_k+\Skint_{k-1}\Big)^2+ \int_{\Mint_*}|\Gac_k[\Phi]|^2\\
&&+\int_{\TT\cup\{\ub=1\}}|\dk^k\Gac|^2.
\eeaa
Together with the iteration assumption \eqref{eq:iterationassumptionforcontrolMintinThmM8}, the control on $\TT$ for $\Gac$ provided by Lemma \ref{Lemma:EstimatesonTT-Mint}, and the control of $\Gac$ on $\ub=1$ provided by \eqref{eq:ub0between1and2realizesinfimumandthusundercontrolLebesguepoint:ub0=1}, we infer, for $k\leq k_{large}+7$, 
\beaa
\nn\int_{\Sint_*\cup\AA} |\dk^k\Phi|^2+ \int_{\Mint_*}|\dk^k\Phi|^2 &\les&  \Big(\ep_0+\Rkint_k+\Skext_{k}\Big)^2+ \int_{\Mint_*}|\Gac_k[\Phi]|^2
\eeaa
as stated.

Next, we focus on the estimate for $\Xbh$. According to Proposition \ref{prop:estimateforcoupledtransportCodazzisystemXbhinMint:chap9}, $\Xbh$  satisfies, for $k\leq k_{large}+7$, 
\beaa
\nn\int_{\Mint_*}|\dk^{k}\Xbh|^2 &\les& \int_{\Mint_* }|\dk^k\trXbc|^2+\int_{\Sint_*\cup\AA}\Big(|\dk^k\trXbc|^2+|\dk^{\leq k-1}(\Gac\setminus\{\Xbh\})|^2\Big)\\
&&+\int_{\TT\cup\{\ub=1\}}|\dk^{\leq k}\Gac|^2+\Big(\ep_0+\Rkint_k+\Skint_{k-1}\Big)^2.
\eeaa
Together with the iteration assumption \eqref{eq:iterationassumptionforcontrolMintinThmM8}, the control on $\TT$ for $\Gac$ provided by Lemma \ref{Lemma:EstimatesonTT-Mint}, and the control of $\Gac$ on $\ub=1$ provided by \eqref{eq:ub0between1and2realizesinfimumandthusundercontrolLebesguepoint:ub0=1}, we infer, for $k\leq k_{large}+7$, 
\beaa
\nn\int_{\Mint_*}|\dk^{k}\Xbh|^2 &\les& \int_{\Mint_* }|\dk^k\trXbc|^2+\int_{\Sint_*\cup\AA}|\dk^k\trXbc|^2+\Big(\ep_0+\Rkint_k+\Skext_k\Big)^2.
\eeaa
Together with the notation $\Gac_k[\Xbh]=\dk^k\trXc$, we deduce, for $k\leq k_{large}+7$, 
\beaa
\nn\int_{\Mint_*}|\dk^{k}\Xbh|^2 &\les& \int_{\Mint_* }|\Gac_k[\Xbh]|^2+\int_{\Sint_*\cup\AA}|\Gac_k[\Xbh]|^2\\
&&+\Big(\ep_0+\Rkint_k+\Skext_k\Big)^2
\eeaa 
as stated. This concludes the proof of Corollary \ref{cor:maincontroloftransporteuqiaotningoingPTstructureMint}.
\end{proof}

Recall from Proposition \ref{prop:schamticformoflinearzizedstrucutreequationsofMint:bis} and Corollary \ref{cor:maincontroloftransporteuqiaotningoingPTstructureMint} that $\Gac_k[\Phi]$ given for each component of $\Gac$ by
 \bea\lab{eq:triangularstructureforthecontrolofMintinThmM8:chap9}
\bsplit
\Gac_k[\trXbc]&=0, \\
\Gac_k[\Xbh]&=\dk^k\trXbc, \qquad \Gac_k[\widecheck{\DD\cos\th}] =\dk^k\Big(\trXbc, \,\Xbh\Big),\\ 
\Gac_k[\Zc]&=\dk^k\Big(\Xbh,\, \trXbc\Big),\qquad \Gac_k[\Hbc]=\dk^k\Big(\Xbh,\, \trXbc\Big),\\
\Gac_k[\DD r] &=\dk^k\Zc,\qquad \Gac_k[e_4(\cos\th)] =\dk^k\Big(\Hbc, \widecheck{\DD(\cos\th)}\Big), \qquad \Gac_k[\omc]=\dk^k\Big(\Zc,\,\Hbc\Big),\\
\Gac_k[\DD\hot\Jk]&=\dk^k\Big(\trXbc,\,\Xbh,\,\Zc,\,\widecheck{\DD(\cos\th)}\Big),\,\,\,\, \Gac_k[\widecheck{\ov{\DD}\c\Jk}]=\dk^k\Big(\trXbc,\,\Xbh,\,\Zc, \,\widecheck{\DD(\cos\th)}\Big),\\
\Gac_k[\trXc]&=\dk^k\Big(\trXbc, \,\widecheck{\DD\c\ov{\Jk}},\, \DD r,\,\widecheck{\DD(\cos\th)}\Big),\quad  \Gac_k[\Xh] =\dk^k\Big(\DD\hot\Jk, \,\widecheck{ \DD(\cos \th)},\, \Xbh\Big),\\
\Gac_k[\widecheck{e_4(r)}] &=\dk^k\Big(\omc, \DD r\Big),\\
\Gac_k[\widecheck{\nab_4\Jk}]&=\dk^k\Big(\widecheck{e_4(r)},\,e_4(\cos\th),\,\omc,\,\Hbc,\,\widecheck{\nab\Jk}\Big),\\
\Gac_k[\Xi]&= \dk^k\Big(\Hbc,\,\trXc,\,\Xh,\,\widecheck{\nab_4 \Jk},\, \widecheck{e_4(r)},\,e_4(\cos\th)\Big).
\end{split}
\eea

We now proceed as follows to control $\Gac$ on $\Mint_*$. We use \eqref{eq:finalestimateinMintstarusingeverythingforGacminusXbh} to control all the components in $\Gac\setminus\{\Xbh\}$ and \eqref{eq:finalestimateinMintstarusingeverythingforXbh} to control $\Xbh$. The components of $\Gac$  are then recovered in the following order, suggested by the triangular structure in \eqref{eq:triangularstructureforthecontrolofMintinThmM8:chap9}, 
\beaa
\trXbc, \,\, \Xbh, \, \widecheck{\DD\cos\th},\,\, \Zc, \,\Hbc,\,\, \DD r, \, e_4(\cos\th),\,\, \omc,\,\DD\hot\Jk,\,\, \widecheck{\ov{\DD}\c\Jk}, \,\, \trXc,\,\,\Xh,\,\,\widecheck{e_4(r)},\,\,\widecheck{\nab_4\Jk},\,\,\Xi.
\eeaa 
The crucial point of this triangular structure is that at each step, when estimating a component $\Phi$ on $\Gac$, all the components in the term $\Gac_k[\Phi]$ on the RHS have been  already estimated. We easily obtain, for $k\leq k_{large}+7$, 
\beaa
\int_{\Sint_*\cup\AA} |\dk^k(\Gac\setminus\{\Xbh\})|^2+ \int_{\Mint_*}|\dk^k\Gac|^2 &\les&  \Big(\ep_0+\Rkint_k+\Skext_{k}\Big)^2.
\eeaa
Together with the control on $\Mint\setminus\Mint_*$ provided by \eqref{eq:controlofGaconMintminusMintstrat}, we infer
\bea
\int_{\Sint_*\cup\AA} |\dk^k(\Gac\setminus\{\Xbh\})|^2+ \int_{\Mint}|\dk^k\Gac|^2 &\les&  \Big(\ep_0+\Rkint_k+\Skext_{k}\Big)^2.
\eea
Together with iteration assumption  \eqref{eq:iterationassumptionforcontrolMintinThmM8}, and in view of the definition of $\Skint_k$, we deduce
\beaa
\Skint_k+\left(\int_{\Sint_*\cup\AA}|\dk^{\leq k}(\Gac\setminus\{\Xbh\})|^2\right)^{\frac{1}{2}} &\les&  \ep_0+\Rkint_k+\Skext_{k}.
\eeaa
This is the iteration assumption  \eqref{eq:iterationassumptionforcontrolMintinThmM8} with $k$ replaced by $k+1$. We have thus obtained
\beaa
\Skint_k &\les&  \ep_0+\Rkint_k+\Skext_{k}, \quad k\leq k_{large}+7.
\eeaa
This concludes the proof of Proposition \ref{Prop:MainestimatesMint}.

%%%%%%%%%%%%%%%%%%%%%%%%%%%%%%%%%%%%%%%%%%%%%%%%%%

\section{Control of the PT-Ricci coefficients in $\Mtop'$}
\lab{sec:controlRiccicoeffMtop}

%%%%%%%%%%%%%%%%%%%%%%%%%%%%%%%%%%%%%%%%%%%%%%%%%%

The goal of this section is to provide the proof of Proof of Proposition \ref{prop:improvementoftheiterationassupmtionThM8-Mtop} recalled in Proposition \ref{Prop:MainestimatesMtop:chap9} below. Since we will not need to refer to the old region $\Mtop$, defined w.r.t. the PG frame, we drop  the prime of $\Mtop'$ in this section.

%%%%%%%%%%%%%%%%%%%%%%%%%%%%%%%%%

\subsection{Preliminaries}
\lab{subsection::controlRiccicoeffMtop-Def}

%%%%%%%%%%%%%%%%%%%%%%%%%%%%%%%%%

Recall the following norm for the linearized Ricci and metric coefficients on $\Mtop$ introduced in section \ref{sec:mainnormsPTframe:chap9}, for $k\le k_{large}+7$,
\beaa
\Sktop_k^2= \int_{\Mtop(r\leq r_0)} \big| \dk^{\le k} \Gac |^2 +(\Sktop^{\ge r_0}_{k})^2
\eeaa
where $\Gac$ denotes the set of all  linearized Ricci and metric coefficients  with respect  to the ingoing PT frame of $\Mtop$, and where
\beaa
\bsplit
(\Sktop^{\ge r_0}_{k})^2 &:= \sup_{\ub_1\geq u_*'}\int_{\Mtop_{r_0, \ub_1}}  r^2\big|\dk^{\le k} (\Xi, \omc, \Xh, \trXc,  \Zc, \Hbc)\big|^2 +r^{2-\dt}|\dk^{\leq k}\trXbc|+  \big| \dk^{\le k}\Xbh\big|^2\\
&+\sup_{\ub_1\geq u_*'}\int_{\Mtop_{r_0, \ub_1}}  \Big( \big| \dk^{\le k}   \widecheck{\DD\cos\th}   \big|^2 + \big| \dk^{\le k}   \DD r    \big|^2 \Big)\\
&+\sup_{\ub_1\geq u_*'}\int_{\Mtop_{r_0, \ub_1}}  \Big(  r^4\big| \dk^{\le k}   e_4 (\cos \th)    \big|^2  +  r^4\big| \dk^{\le k}  \widecheck{ e_4 (r)  }  \big|^2   \Big)\\
&+\sup_{\ub_1\geq u_*'}\int_{\Mtop_{r_0, \ub_1}} r^2 \Big(\big|\dk^{\le k} \DD\hot\Jk|^2  +\big|\dk^{\le k} \widecheck{\ov{\DD}\c\Jk}\big|^2\Big) +r^6\big| \dk^{\le k}\widecheck{\nab_4\Jk}\big|^2 \Big), 
\end{split}
\eeaa
with $\Xi$, $\omc$, $\trXc$, $\Xh$, $\Zc$, $\Hbc$, $\trXbc$,  $\Xbh$ the linearized Ricci coefficients   of the ingoing PT frame of $\Mtop$, and with the notation
\beaa
\Mtop_{r_0,\ub_1} &:=& \Mtop(r\geq r_0)\cap\{\ub_1\leq \ub\leq \ub_1+1\}.
\eeaa

For the curvature norms in $\Mtop$, we rely in particular on the scalar function $\tau$ introduced in section \ref{sec:defintionofthescalarfunctiontau:chap9} and define
\beaa
\bsplit
\Rktop_k^2 &:=  \int_{\Mtop(r\geq r_0)}\Big( r^{3+\dt}|\dk^{\le k}(A, B)|^2 +r^{3-\dt}|\dk^{\le k}\Pc |^2\Big)\\
&+\sup_{\ub_1\geq u_*'}\int_{\Mtop_{r_0, \ub_1}}\Big( r^{2} |\dk^{\le k}\Bb|^2 +|\dk^{\le k}\Ab|^2\Big)+\sup_\tau\int_{\Mtop(r\leq r_0)\cap\Si(\tau)}  |\dk^{\le k}  \Rc|^2,
\end{split}
\eeaa
where $\Rc$ is the set of all linearized  curvature coefficients  w.r.t.  the ingoing PT frame  of $\Mtop$, and $A$, $B$, $\Pc$, $\Bb$, $\Ab$ denote the linearized curvature components relative to  the ingoing PT  frame of $\Mtop$.

Also, recall the definition of  $L_*(k)$, see \eqref{eq:defintionofLstarofk:chap9}, 
  \bea\lab{eq:defintionofLstarofk:chap9:ter} 
 L^2_*(k)  =  \int_{\{u=u_*'\}}    \big|\Rc\big|_{w, k}^2 +\int_{\Si_*\cap\{u=u_*'\}} \big|\Gac\big|_{w, k}^2
 \eea 
where, see \eqref{eq:defintionofthenormsRcGacweighted:bis}, 
 \bea\lab{eq:defintionofthenormsRcGacweighted:ter}
 \bsplit
 \big|\Rc\big|_{w, k}^2 &:= r^{3+\dt}|\dk_*^{\le k}(A, B)|^2 +r^{3-\dt}|\dk_*^{\le k}\Pc |^2 +r^{1-\dt}|\dk_*^{\le k}\Bb|^2,\\
\big|\Gac\big|_{w, k}^2 &:= r^2|\dkb^{\leq k}\Ga_g|^2+|\dkb^{\leq k}\Ga_b|^2.
\end{split}
\eea
In \eqref{eq:defintionofthenormsRcGacweighted:ter}, $(A, B, \Pc, \Bb)$ denote  linearized curvature components w.r.t. the the  outgoing PT frame of $\Mext$,  
$\Ga_g$ and $\Ga_b$ are  defined w.r.t. the  outgoing PT frame of $\Mext$ as in Definition \ref{definition.Ga_gGa_b:outgoingPTcase:chap9}, $\dk_*$ denote weighted derivatives tangential to the hypersurface $\{u=u_*'\}$, and $\dkb$ denote weighted derivatives tangential to the sphere $\Si_*\cap\{u=u_*'\}$.

The goal of  section \ref{sec:controlRiccicoeffMtop} is to provide the proof of Proof of Proposition \ref{prop:improvementoftheiterationassupmtionThM8-Mtop}. For convenience, we restate the result below.
 \begin{proposition}\lab{Prop:MainestimatesMtop:chap9}
 Relative to the ingoing PT frame of $\Mtop$, we have
\bea
 \Sktop_{k} &\les& \ep_0+   L_*(k)  +   \Rktop_{k}, \quad k\leq k_{large}+7.
 \eea
 In particular, we have
\bea
 \Sktop^{\geq r_0}_{k} &\les& \ep_0+   L_*(k)  +   \Rktop_{k}, \quad k\leq k_{large}+7,
 \eea
  in which  case the constant in $\les$ is independent of $r_0$.
 \end{proposition} 

 To prove Proposition \ref{Prop:MainestimatesMtop:chap9}, we rely in particular on our bootstrap assumption {\bf  BA-PT}, see \eqref{eq:mainbootassforchapte9}, which implies   for the Ricci and metric coefficients of the ingoing PT frame of $\Mtop$
 \bea
 \lab{eq:Bootstrap-Mtop-M8}
 \Sktop_k&\leq & \ep, \qquad k\leq k_{large}+7.
 \eea

%%%%%%%%%%%%%%%%%%%%%%%%%%%%%%%%%

\subsection{Control on the hypersurface $\{u=u_*'\}$} 
 
%%%%%%%%%%%%%%%%%%%%%%%%%%%%%%%%% 
 
Our goal in this section is to control the linearized Ricci and metric coefficients of the ingoing PT frame of $\Mtop$ on $\{u=u_*'\}$. We start with the following transport lemma along the hypersurface $\{u=u_*'\}$ of $\Mext$.

\begin{lemma}\lab{lemma:transportrp-L2-knorms-PT:foru=ustarprime}
Let $U$ and $F$  anti-selfdual $k$-tensors on $\{u=u_*'\}$. Assume  that  $U$ verifies one of the following equations, for a real constant $c$,
\bea
\lab{eq:maintransportMext-PT:bis}
\nab_4 U+\frac{c}{q} U=F
\eea
or
\bea
\lab{eq:maintransportMext2-PT:bis}
\nab_4 U+ \Re\left(\frac{c}{q}\right)  U=F.
\eea
In both cases we derive,   for  any $r_0\leq r\le r_* =r_*(u_*')$  
 \bea\lab{eq:integration-knorms=PT:bis}
   r^{c-1}  \|\widetilde{\dk}^{\leq k}U\|_{L^2(S(u_*',r))} \les    r_*  ^{c-1} \|\widetilde{\dk}^{\leq k}U\|_{L^2(S(u_*',r_*))}+ \int_r^{r_*}      \la ^{c-1} \|\widetilde{\dk}^{\leq k}F\|_{L^2(S(u_*',\la))} d\la,
\eea
where $\widetilde{\dk}$ denotes weighted derivatives tangential to $\{u=u_*'\}$, i.e.
\beaa
\widetilde{\dk} &:=& \left\{r\nab_4,\, r\left(\nab_1 - \frac{1}{e_3(u)}e_1(u)\nab_3\right),\, r\left(\nab_2 - \frac{1}{e_3(u)}e_2(u)\nab_3\right)\right\}.
\eeaa
\end{lemma}

\begin{proof}
This follows from a straightforward adaptation of Proposition \ref{Prop:transportrp-L2-knorms-PT} which is our main transport lemma in $\Mext$. 
\end{proof}

Next, we control the linearized Ricci and metric coefficients of the outgoing PT frame of $\Mext$ on the hypersurface $\{u=u_*'\}$. 
\begin{proposition}
\lab{proposition:controlofPTstructureMextonu=ustarprime}
Along the  hypersurface $\{u=u_*'\}$, the outgoing PT frame of $\Mext$ satisfies the following estimates, for $k\le k_{large}+7$,
\bea
\lab{eq:Estimatesin-MextPTonu=ustarprime}
\nn\sup_{S\subset\{u=u_*'\}}\Big(r^2\|\widetilde{\dk}^{\leq k}(\Ga_g\setminus\{\trXbc\})\|_{L^2(S)}^2+r^{2-\dt}\|\widetilde{\dk}^{\leq k}\trXbc\|_{L^2(S)}^2\Big)\\
+\sup_{S\subset\{u=u_*'\}}\Big(\|\widetilde{\dk}^{\leq k}(\Ga_b\setminus\{\Xib\})\|_{L^2(S)}^2+r^{-\dt}\|\widetilde{\dk}^{\leq k}\Xib\|_{L^2(S)}^2\Big) &\les& \ep_0^2+L_*^2(k)
\eea
where the constant in the definition of  $\les$ is independent of $r_0$.
\end{proposition}

\begin{proof}
We proceed as follows:
\begin{enumerate}
\item We rely on the equations for the linearized Ricci and metric coefficients of the PT structure of $\Mext$, see Proposition \ref{Prop:linearizedPTstructure1} and Lemmas \ref{lemma:transportequationine4fornabthetaandrande3thetaandr:linearized} and  \ref{lemma:transportequationine4fornaJkandnab3Jk:linearized}, and the control of transport equations along $e_4$ provided by Lemma \ref{lemma:transportrp-L2-knorms-PT:foru=ustarprime}. 

\item We integrate the transport equations in the order consistent with the triangular structure of the system, i.e., as in section \ref{sec:proofofProp:EstmatesSkext}, we estimate the linearized Ricci and metric coefficients of the outgoing PT structure of $\Mext$ in the following order 
\beaa
 \trXc,\,\,\, \Xh,\,\,\,  \DD \cos \th,\,\,\,  \Zc, \,\,\,  \DD r,\,\,\,  \Hc, \,\,\,    \widecheck{e_3(r)}, \,\,\,     \ombc, \,\,\, \DD\hot\Jk, \,\,\, \widecheck{\DD\c\ov{\Jk}},\,\,\,  e_3 (\cos \th),\,\,\,  \widecheck{e_3 \Jk},\,\,\,  \Xib.
\eeaa

\item We also make use  of the bootstrap assumption \eqref{Prop-bootstrap:finalcontrolofthePTframeonMext} for the outgoing PT frame of $\Mext$. 
\end{enumerate}
As a result, we obtain the control of  the linearized Ricci and metric coefficients of the outgoing PT structure of $\Mext$ by the following two contributions:
\begin{itemize}
\item their restriction to $S(u_*', r_*)=\Si_*\cap\{u=u_*'\}$,

\item $r^p$ weighed norms of curvature along $\{u=u_*\}$, 
\end{itemize}
which are precisely in the form of $L_*(k)$, see \eqref{eq:defintionofLstarofk:chap9:ter}. The proof is very similar to  the one in section \ref{sec:proofofProp:EstmatesSkext}, and in fact simpler. We leave the details to the reader.
\end{proof}

 We are now ready to state the main result of this section on the control of the linearized Ricci and metric coefficients of the ingoing PT frame of $\Mtop$ on $\{u=u_*'\}$. 
\begin{proposition}
\lab{proposition:Estimatesin-MtopPTonu=ustarprime}
Along the  hypersurface $\{u=u_*'\}$, the ingoing PT frame of $\Mtop$ satisfies the following estimates, for $k\le k_{large}+7$,
\bea
\lab{eq:Estimatesin-MtopPTonu=ustarprime}
\nn\sup_{S\subset\{u=u_*'\}}\Big(r^2\|\dk^{\leq k}(\Ga_g\setminus\{\trXbc\})\|_{L^2(S)}^2+r^{2-\dt}\|\dk^{\leq k}\trXbc\|_{L^2(S)}^2\Big)\\
+\sup_{S\subset\{u=u_*'\}}\|\dk^{\leq k}\Ga_b\|_{L^2(S)}^2 &\les& \ep_0^2+L_*^2(k)
\eea
where the constant in the definition of  $\les$ is independent of $r_0$.
\end{proposition}

 \begin{proof}
 The proof is the analog for the PT frames of $\Mext$ and $\Mtop$ of the one of Lemma \ref{lemma:controlofPGstructureMtoptonu=u*} for the PG frames of $\Mext$ and $\Mtop$. To simplify  the notations, in this proof, we denote:
\begin{itemize}
\item by $(e_4, e_3, e_1, e_2)$ the outgoing PT frame of $\Mext$, with all quantities associated to the outgoing PT structure of $\Mext$ being unprimed, 

\item by $(e_4', e_3', e_1', e_2')$ the ingoing PT frame of $\Mtop$, with all quantities associated to the ingoing PT structure of $\Mint$ being primed. 
\end{itemize}

In view of the above notations, and the initialization of the ingoing PT structure of $\Mtop$ from the outgoing PG structure of $\Mext$ on $\{u=u_*'\}$, see section \ref{sec:defintionofthePTstructuresinMM}, we have  
\beaa
\ub=u+2\int_{r_0}^r\frac{{\tilde{r}}^2+a^2}{{\tilde{r}}^2-2m\tilde{r}+a^2}d\tilde{r}, \qquad r'=r, \qquad \th'=\th,  \qquad \Jk'=\Jk\quad\quad\textrm{on}\quad\{u=u_*'\},
\eeaa
and
\beaa
e_4'=\la e_4, \qquad e_3'=\la^{-1}e_3, \qquad e_a'=e_a, \,\, a=1,2, \quad\textrm{on}\quad\{u=u_*'\},
\eeaa
where $\la$ is given by
\beaa
\la=\frac{\Delta}{|q|^2}.
\eeaa

Based on this initialization we derive in particular,   for any tangent vector $X$ on $\{u=u_*'\}$, 
\beaa
\g(\D_Xe_4', e_a') &=& \g\left(\D_X\left(\la e_4\right), e_a\right)=\la\g\left(\D_Xe_4, e_a\right),\\
\g(\D_Xe_3', e_a') &=& \g\left(\D_X\left(\la^{-1}e_3\right), e_a\right)=\la^{-1}\g\left(\D_Xe_3, e_a\right),\\
\g(\D_Xe_4', e_3') &=& \g\left(\D_X\left(\la e_4\right), \la^{-1}e_3\right)=-2X(\log\la)+\g\left(\D_Xe_4, e_3\right).
\eeaa
Note that any such tangent   vector $X$  is a linear combination of  the following  three tangent directions to  $\{u=u_*'\}$,
\beaa
X_1=e_1-\frac{e_1(u)}{e_3(u)}e_3, \qquad X_2=e_2-\frac{e_2(u)}{e_3(u)}e_3, \qquad X_4=e_4.
\eeaa
Also, note that we have
\beaa
  \nab(u)=a\Re(\Jk)+\widetilde{\Ga}_b, \qquad   e_3(u)=\frac{2(r^2+a^2)}{|q|^2}+r\widetilde{\Ga}_b,
\eeaa
where we have introduced the following notation\footnote{Recall that $\Xib$ and $\trXbc$ exhibit a $r^{-\frac{\dt}{2}}$ loss compared to the other components.}  
\beaa
\widetilde{\Ga}_b:=\Ga_b\setminus\{\Xib\}, \qquad \widetilde{\Ga}_g:=\Ga_g\setminus\{\trXbc\}.
\eeaa
We infer 
\beaa
X_b=e_b+\left(-\frac{a|q|^2}{2(r^2+a^2)}(\Re(\Jk))_b+\Ga_b\right) e_3, \,\, b=1,2, \qquad X_4= e_4,
\eeaa
as well as 
\beaa
X_b=e_b'+\left(-\frac{a|q|^2}{2(r^2+a^2)}(\Re(\Jk))_b+\widetilde{\Ga}_b\right)\la e_3', \,\, b=1,2, \qquad X_4=\la^{-1} e_4'.
\eeaa
We infer
\beaa
\chi_{ba}' +2\left(-\frac{a|q|^2}{(r^2+a^2)}(\Re(\Jk))_b+\widetilde{\Ga}_b\right)\la\eta_a' &=& \la\chi_{ba}+2\left(-\frac{a|q|^2}{2(r^2+a^2)}(\Re(\Jk))_b+\Ga_b\right)\la\eta_a,\\
\chib_{ba}' +2\left(-\frac{a|q|^2}{2(r^2+a^2)}(\Re(\Jk))_b+\widetilde{\Ga}_b\right)\la\xib_a' &=& \la^{-1}\chib_{ba}+2\left(-\frac{a|q|^2}{2(r^2+a^2)}(\Re(\Jk))_b+\widetilde{\Ga}_b\right)\la^{-1}\xib_a,\\
\ze_b' - 2\left(-\frac{a|q|^2}{(r^2+a^2)}(\Re(\Jk))_b+\widetilde{\Ga}_b\right)\la\omb' &=& -\left(e_b+\left(-\frac{a|q|^2}{2(r^2+a^2)}(\Re(\Jk))_b+\Ga_b\right) e_3\right)\log\la\\
&&+\ze_b - 2\left(-\frac{a|q|^2}{(r^2+a^2)}(\Re(\Jk))_b+\widetilde{\Ga}_b\right)\la\omb,
\eeaa
and
\beaa
\xi' &=& \la^2\xi,\\
\etab' &=& \etab,\\
\om' &=& -\frac{1}{2}e_4(\log\la)+\om.
\eeaa
Together with 
\begin{itemize}
\item the fact that, by the PT gauge choices, $\xib'=0$, $\omb'=0$, $H'=\frac{aq'}{|q'|^2}$, $\xi=0$, $\om=0$, $\Hb=-\frac{a\ov{q}}{|q|^2}\Jk$, 

\item the form of $\la$, 

\item the fact that $r'=r$, $\th'=\th$, and $\Jk'=\Jk$ on $\{u=u_*'\}$, 

\item the definition of the linearized quantities for the PT frame in $\Mtop$ and $\Mext$, 
\end{itemize}
we infer on $\{u=u_*'\}$
\beaa
\trXc', \,\Xh',\,\Zc' = \widetilde{\Ga}_g,\qquad \trXbc' = \la\trXbc+O(r^{-1})\Xib,\qquad \Xbh' = \widetilde{\Ga}_b,
\eeaa
and
\beaa
\Xi' = 0,\qquad \Hbc' = 0,\qquad \omc' = 0.
\eeaa
We deduce 
\beaa
&& \sup_{S\subset\{u=u_*'\}}\Big(r^2\|\widetilde{\dk}^{\leq k}(\Xi', \,\om',\, \trXc', \,\Xh',\,\Zc', \, \Hbc')\|_{L^2(S)}^2+r^{2-\dt}\|\widetilde{\dk}^{\leq k}\trXbc'\|_{L^2(S)}^2\Big)\\
&&+\sup_{S\subset\{u=u_*'\}}\|\widetilde{\dk}^{\leq k}\Xbh'\|_{L^2(S)}^2\\
&\les& \sup_{S\subset\{u=u_*'\}}\Big(r^2\|\widetilde{\dk}^{\leq k}\widetilde{\Ga}_g\|_{L^2(S)}^2+r^{2-\dt}\|\widetilde{\dk}^{\leq k}\trXbc\|_{L^2(S)}^2+\|\widetilde{\dk}^{\leq k}\widetilde{\Ga}_b\|_{L^2(S)}^+r^{-\dt}\|\widetilde{\dk}^{\leq k}\Xib\|_{L^2(S)}^2\Big).
\eeaa
Together with Proposition \ref{proposition:controlofPTstructureMextonu=ustarprime} on the control of the PT structure of $\Mext$ on $\{u=u_*'\}$, we infer, for $k\leq k_{large}+7$, 
\beaa
 \sup_{S\subset\{u=u_*'\}}\Big(r^2\|\widetilde{\dk}^{\leq k}(\Xi', \,\om',\, \trXc', \,\Xh',\,\Zc', \, \Hbc')\|_{L^2(S)}^2+r^{2-\dt}\|\widetilde{\dk}^{\leq k}\trXbc'\|_{L^2(S)}^2\Big)\\
+\sup_{S\subset\{u=u_*'\}}\|\widetilde{\dk}^{\leq k}\Xbh'\|_{L^2(S)}^2 &\les& \ep_0^2+L_*^2(k)
\eeaa
where the constant in the definition of  $\les$ is independent of $r_0$.

Concerning the linearized metric coefficients, we consider 
\beaa
\ub=u+2\int_{r_0}^r\frac{{\tilde{r}}^2+a^2}{{\tilde{r}}^2-2m\tilde{r}+a^2}d\tilde{r}, \qquad r'=r, \qquad \th'=\th,  \qquad \Jk'=\Jk\quad\quad\textrm{on}\quad\{u=u_*'\},
\eeaa
and apply to both sides of these identities the tangential vectorfields $X$ to $\{u=u_*'\}$ as above. Then, proceeding 
similarly to the linearized Ricci coefficients\footnote{See Lemma \ref{lemma:necessaryidentitiesonuequalustarforthecontrolofMtopnuequlustar} for the analog situation for PG structures.}, we obtain on $\{u=u_*'\}$ relations between the linearized metric coefficients on $\Mtop$ and $\Mext$, which then, together with the above control for linearized Ricci coefficients, yields 
\beaa
\nn\sup_{S\subset\{u=u_*'\}}\Big(r^2\|\widetilde{\dk}^{\leq k}(\Ga_g'\setminus\{\trXbc'\})\|_{L^2(S)}^2+r^{2-\dt}\|\widetilde{\dk}^{\leq k}\trXbc'\|_{L^2(S)}^2\Big)\\
+\sup_{S\subset\{u=u_*'\}}\|\widetilde{\dk}^{\leq k}\Ga_b'\|_{L^2(S)}^2 &\les& \ep_0^2+L_*^2(k)
\eeaa
where the constant in the definition of  $\les$ is independent of $r_0$.

Then, relying on the transport equations along $e_3$ of the ingoing PT structure of $\Mtop$, and on the fact that $e_3$ is transversal to $\{u=u_*'\}$ so that $\dk$ is spanned by $\nab_3$ and $\widetilde{\dk}$, we infer, for $k\le k_{large}+7$,
\beaa
\nn\sup_{S\subset\{u=u_*'\}}\Big(r^2\|\dk^{\leq k}(\Ga_g\setminus\{\trXbc\})\|_{L^2(S)}^2+r^{2-\dt}\|\dk^{\leq k}\trXbc\|_{L^2(S)}^2\Big)\\
+\sup_{S\subset\{u=u_*'\}}\|\dk^{\leq k}\Ga_b\|_{L^2(S)}^2 &\les& \ep_0^2+L_*^2(k)
\eeaa
where the constant in the definition of  $\les$ is independent of $r_0$. This concludes the proof of Proposition \ref{proposition:Estimatesin-MtopPTonu=ustarprime}.
\end{proof}

%%%%%%%%%%%%%%%%%%%%%%%%%%%%%%%%%%%%%%%%%%%%%%

\subsection{Proof of Proposition \ref{Prop:MainestimatesMtop:chap9}}

%%%%%%%%%%%%%%%%%%%%%%%%%%%%%%%%%%%%%%%%%%%%%%

We are now in position to prove Proposition \ref{Prop:MainestimatesMtop:chap9} on the control of the ingoing PT structure of $\Mtop$.

%%%%%%%%%%%%%%%%%%%%%%%%%%%%%%%%%%%%%%%%%%%%%%

\subsubsection{Propagation lemmas}

%%%%%%%%%%%%%%%%%%%%%%%%%%%%%%%%%%%%%%%%%%%%%%

Note that $\Mtop$ has by construction the following boundaries 
\bea
\pr\Mtop &=& \{u=u_*'\}\cup\{\ub=u_*'\}\cup{}^{(top)}\Si\cup(\AA\cap\{\ub\geq u_*'\}).
\eea
Also, recall from the construction of the scalar function $\tau$ in Proposition \ref{prop:propertiesoftauusefulfortheoremM8:chap9} that the following properties hold on $\Mtop$:
\begin{enumerate}
\item The future boundary ${}^{(top)}\Si$ of $\MM$  is given by 
\bea
{}^{(top)}\Si=\{\tau=u_*\}
\eea
and $\tau\leq u_*$ on $\MM$.

\item Denoting, on each level set of $\ub$ in $\,{}^{(top)}\MM(r\geq r_0)$, by $r_+(\ub)$ the maximal value of $r$  and by $r_-(\ub)$ the minimal value of $r$, we have\footnote{Note that \eqref{eq:upperboundrpubminusrmubonMtop:chap9:bis} depends on the choice of ${}^{(top)}\Si$ and hence on the choice of $\tau$.}
\bea\lab{eq:upperboundrpubminusrmubonMtop:chap9:bis}
0\leq r_+(\ub)-r_-(\ub)\leq 2(2m+1).
\eea

\item In $\Mtop(r\leq r_0)$, $\tau$ satisfies   
\bea\lab{eq:anotherpropertyoftauonMtopwhichisuseful}
u_*-2(m+2)\leq \tau \leq u_*.
\eea
\end{enumerate}

\begin{remark}\lab{rmk:Mtopisinfactalocalexistencetyperegion:chap9}
In view of \eqref{eq:upperboundrpubminusrmubonMtop:chap9:bis} and \eqref{eq:anotherpropertyoftauonMtopwhichisuseful}, $\,{}^{(top)}\MM$  is in fact a local existence type region.
\end{remark}

We introduce the following norms on $\Mtop$
 \beaa
  \| f\|_{2} (\ub,r):=\| f\|_{L^2\big(S(\ub,r)\big)}, \qquad \|f\|_{2,k}(\ub, r):= \sum_{i=0}^k \|\dk^i f\|_{2 }(\ub, r)
  \eeaa
 which allow us to state a first propagation lemma.

\begin{proposition}
\lab{Prop:transportrp-L2-knorms-PT:Mtop}
Let $U$ and $F$  anti-selfdual $k$-tensors. Assume  that  $U$ verifies one of the following equations, for a real constant $c$,
\bea
\lab{eq:maintransportMext-PT:Mtop}
\nab_3 U-\frac{c}{\ov{q}} U=F
\eea
or
\bea
\lab{eq:maintransportMext2-PT:Mtop}
\nab_3 U- \Re\left(\frac{c}{\ov{q}}\right)  U=F.
\eea
In both cases we derive,  for any $c'$, and  for  any $r_-(\ub)\leq r\leq r_1\leq r_+(\ub)$  at fixed $\ub$, with  $\ub\geq u_*'$,  in $\Mtop$
 \bea
   r^{c'}  \|U\|_{2,k} (\ub, r) &\les&    r_1^{c'} \|  U \|_{2,k} (\ub,  r_1)+ \int_r^{r_1}      \la ^{c'}\|F\|_{2,k}(\ub, \la) d\la,
\eea
where the constant in the definition of  $\les$ is independent of $r_0$. 
\end{proposition} 

\begin{proof}
Note first that we have  for  any $r_-(\ub)\leq r\leq r_1\leq r_+(\ub)$  
\beaa
   r^{c-1}  \|U\|_{2,k} (\ub, r) &\les&    r_1^{c-1} \|  U \|_{2,k} (\ub,  r_1)+ \int_r^{r_1}      \la ^{c-1}\|F\|_{2,k}(\ub, \la) d\la,
\eeaa
whose proof is completely analogous of the one of Proposition \ref{Prop:transportrp-L2-knorms-PT}. The stated estimate follows then by multiplying by $r^{c'-c+1}$ and by noticing that $r\leq r_1\leq r+2(2m+1)$ in view of 
\eqref{eq:upperboundrpubminusrmubonMtop:chap9:bis}.
\end{proof}

We infer the following corollary.
\begin{corollary}\lab{cor:transportrp-L2-knorms-PT:Mtop}
Let $U$ and $F$  anti-selfdual $k$-tensors. Assume  that  $U$ verifies \eqref{eq:maintransportMext-PT:Mtop} or 
\eqref{eq:maintransportMext2-PT:Mtop}. In both cases we derive the following:
\begin{enumerate}
\item In $\Mtop(r\geq r_0)$, we have, for any $c'$,  
\bea
\nn \sup_{\ub_1\geq u_*'}\int_{\Mtop_{r_0, \ub_1}}r^{2c'}|\dk^{\leq k}U|^2 &\les& \sup_{S\subset\{u=u_*'\}}r^{2c'}\|\dk^{\leq k}U\|^2_{L^2(S)}\\
 &&+ \sup_{\ub_1\geq u_*'}\int_{\Mtop_{r_0, \ub_1}}r^{2c'}|\dk^{\leq k}F|^2,
\eea
where the constant in the definition of  $\les$ is independent of $r_0$, and where we recall the notation 
\beaa
\Mtop_{r_0,\ub_1} &=& \Mtop(r\geq r_0)\cap\{\ub_1\leq \ub\leq \ub_1+1\}.
\eeaa

\item In $\Mtop(r\leq r_0)$, we have,  for any $c'$,   
\bea
\nn\int_{\Mtop(r\leq r_0)} \big| \dk^{\le k}U|^2 &\les& \sup_{S\subset\{u=u_*'\}}r^{2c'}\|\dk^{\leq k}U\|^2_{L^2(S)}+\sup_\tau\int_{\Mtop(r\leq r_0)\cap\Si(\tau)}\big| \dk^{\le k}F|^2\\
&&+\sup_{\ub_1\geq u_*'}\int_{\Mtop_{r_0, \ub_1}}r^{2c'}|\dk^{\leq k}F|^2.
\eea
\end{enumerate}
\end{corollary}

\begin{proof}
We start with the estimate for $U$ in $\Mtop(r\geq r_0)$. In view of Proposition \ref{Prop:transportrp-L2-knorms-PT:Mtop} with the choice $r_1=r_+(\ub)$, we have,  for  any $r_-(\ub)\leq r\leq r_+(\ub)$,
 \beaa
   r^{c'}  \|U\|_{2,k} (\ub, r) &\les&    (r_+(\ub))^{c'} \|  U \|_{2,k} (\ub,  r_+(\ub))+ \int_r^{r_+(\ub)}      \la ^{c'}\|F\|_{2,k}(\ub, \la) d\la.
\eeaa
Squaring, using Cauchy-Schwarz and the bound \eqref{eq:upperboundrpubminusrmubonMtop:chap9:bis}, as well as $r\geq \max(r_0, r_-(\ub))$, we infer
 \beaa
   r^{2c'}  \|U\|^2_{2,k} (\ub, r) &\les&    (r_+(\ub))^{2c'} \|  U \|^2_{2,k} (\ub,  r_+(\ub))+ \int_{\max(r_0, r_-(\ub))}^{r_+(\ub)}      \la ^{2c'}\|F\|^2_{2,k}(\ub, \la) d\la.
\eeaa
Since $S(\ub,  r_+(\ub))\subset\{u=u_*'\}$ by the definition of $r_+(\ub)$, we deduce
\beaa
   r^{2c'}  \|U\|^2_{2,k} (\ub, r) &\les&   \sup_{S\subset\{u=u_*'\}}r^{2c'}\|\dk^{\leq k}U\|^2_{L^2(S)}+ \int_{\max(r_0, r_-(\ub))}^{r_+(\ub)}      \la ^{2c'}\|F\|^2_{2,k}(\ub, \la) d\la.
\eeaa
Integrating in $r$ on $(\max(r_0, r_-(\ub)), r_+(\ub))$, and using again  the bound \eqref{eq:upperboundrpubminusrmubonMtop:chap9:bis}, this yields
\beaa
 \int_{\max(r_0, r_-(\ub))}^{r_+(\ub)}r^{2c'}  \|U\|^2_{2,k} (\ub, r)dr &\les&   \sup_{S\subset\{u=u_*'\}}r^{2c'}\|\dk^{\leq k}U\|^2_{L^2(S)}\\
 &&+ \int_{\max(r_0, r_-(\ub))}^{r_+(\ub)}r^{2c'}\|F\|^2_{2,k}(\ub, r) dr.
\eeaa
Integrating in $\ub$ on $(\ub_1, \ub_1+1)$, we infer, in view of the definition of $\Mtop_{r_0,\ub_1}$, 
 \beaa
 \int_{\Mtop_{r_0,\ub_1}}r^{2c'}|\dk^{\leq k}U|^2  &\les&   \sup_{S\subset\{u=u_*'\}}r^{2c'}\|\dk^{\leq k}U\|^2_{L^2(S)}+ \int_{\Mtop_{r_0,\ub_1}}r^{2c'}|\dk^{\leq k}F|^2.
\eeaa
Taking the supremum in $\ub_1\geq u_*'$, we deduce  
\beaa
\nn \sup_{\ub_1\geq u_*'}\int_{\Mtop_{r_0, \ub_1}}r^{2c'}|\dk^{\leq k}U|^2 &\les& \sup_{S\subset\{u=u_*'\}}r^{2c'}\|\dk^{\leq k}U\|^2_{L^2(S)}\\
 &&+ \sup_{\ub_1\geq u_*'}\int_{\Mtop_{r_0, \ub_1}}r^{2c'}|\dk^{\leq k}F|^2,
\eeaa
as stated. Note that the constant in the definition of  $\les$ is independent of $r_0$ since it is the case in Proposition \ref{Prop:transportrp-L2-knorms-PT:Mtop} and in the bound \eqref{eq:upperboundrpubminusrmubonMtop:chap9:bis}.

Next, we consider the second case, i.e. the estimate for $U$ in $\Mtop(r\leq r_0)$. We proceed as above, but instead of integrating in $r$ on $(\max(r_0, r_-(\ub)), r_+(\ub))$, we integrate in $r$ on $(r_-(\ub), r_+(\ub))$ and obtain 
\beaa
\nn \int_{\{\ub_1\leq \ub\leq \ub_1+1\}}r^{2c'}|\dk^{\leq k}U|^2 &\les& \sup_{S\subset\{u=u_*'\}}r^{2c'}\|\dk^{\leq k}U\|^2_{L^2(S)}+ \int_{\{\ub_1\leq \ub\leq \ub_1+1\}}r^{2c'}|\dk^{\leq k}F|^2.
\eeaa
Now, recall that $\ub\geq u_*'$ on $\Mtop$. Also, $\tau\leq u_*$ on $\Mtop$ and $\tau=u_*+f(r)$ for an explicit smooth function $f$, see the proof of Proposition \ref{prop:propertiesoftauusefulfortheoremM8:chap9}. Thus, we have $\ub\leq u_*-f(r)$ and hence
\beaa
\max_{\Mtop(r\leq r_0)}\ub &\leq& u_*+\max_{r_+-\de_\HH\leq r\leq r_0}f(r).
\eeaa
Since $u_*'\geq u_*-2$ by definition, we infer
\beaa
\max_{\Mtop(r\leq r_0)}\ub -u_*' &\leq&  u_*-u_*'+\max_{r_+-\de_\HH\leq r\leq r_0}f(r)\\
&\leq& 2+\max_{r_+-\de_\HH\leq r\leq r_0}f(r).
\eeaa
Thus, introducing the notation 
\beaa
\ub_{max} &:=& \max_{\Mtop(r\leq r_0)}\ub,
\eeaa
we have
\beaa
\ub_{max} -u_*' &\leq& 2+\max_{r_+-\de_\HH\leq r\leq r_0}f(r)\les 1.
\eeaa
In particular, we deduce 
\beaa
\int_{\Mtop(r\leq r_0)}r^{2c'}|\dk^{\leq k}U|^2 &\les& \sup_{u_*'\leq\ub\leq\ub_{max}}\int_{\{\ub_1\leq \ub\leq \ub_1+1\}}r^{2c'}|\dk^{\leq k}U|^2
\eeaa
and 
\beaa
\sup_{u_*'\leq\ub\leq\ub_{max}}\int_{\{\ub_1\leq \ub\leq \ub_1+1\}}r^{2c'}|\dk^{\leq k}F|^2 &\les& \int_{\Mtop(r\leq r_0)} \big| \dk^{\le k}F|^2\\
&&+\sup_{\ub_1\geq u_*'}\int_{\Mtop_{r_0, \ub_1}}r^{2c'}|\dk^{\leq k}F|^2.
\eeaa
Thus, coming back to the above estimate 
\beaa
\nn \int_{\{\ub_1\leq \ub\leq \ub_1+1\}}r^{2c'}|\dk^{\leq k}U|^2 &\les& \sup_{S\subset\{u=u_*'\}}r^{2c'}\|\dk^{\leq k}U\|^2_{L^2(S)}+ \int_{\{\ub_1\leq \ub\leq \ub_1+1\}}r^{2c'}|\dk^{\leq k}F|^2,
\eeaa
and taking the supremum on $u_*'\leq\ub\leq\ub_{max}$, we infer 
\beaa
\nn\int_{\Mtop(r\leq r_0)} \big| \dk^{\le k}U|^2 &\les& \sup_{S\subset\{u=u_*'\}}r^{2c'}\|\dk^{\leq k}U\|^2_{L^2(S)}+\int_{\Mtop(r\leq r_0)} \big| \dk^{\le k}F|^2\\
&&+\sup_{\ub_1\geq u_*'}\int_{\Mtop_{r_0, \ub_1}}r^{2c'}|\dk^{\leq k}F|^2
\eeaa
Finally, since we have $u_*-2(m+2)\leq \tau \leq u_*$ in $\Mtop(r\leq r_0)$, see \eqref{eq:anotherpropertyoftauonMtopwhichisuseful}, we have 
\beaa
\int_{\Mtop(r\leq r_0)} \big| \dk^{\le k}F|^2 &\leq& 2(m+2)\sup_\tau\int_{\Mtop(r\leq r_0)\cap\Si(\tau)}\big| \dk^{\le k}F|^2 
\eeaa
and hence
\beaa
\nn\int_{\Mtop(r\leq r_0)} \big| \dk^{\le k}U|^2 &\les& \sup_{S\subset\{u=u_*'\}}r^{2c'}\|\dk^{\leq k}U\|^2_{L^2(S)}+\sup_\tau\int_{\Mtop(r\leq r_0)\cap\Si(\tau)}\big| \dk^{\le k}F|^2\\
&&+\sup_{\ub_1\geq u_*'}\int_{\Mtop_{r_0, \ub_1}}r^{2c'}|\dk^{\leq k}F|^2
\eeaa
as stated. This concludes the proof of Corollary \ref{cor:transportrp-L2-knorms-PT:Mtop}.
\end{proof}

%%%%%%%%%%%%%%%%%%%%%%%%%%%%%%%%%%%%%%%%%%%%%%

\subsubsection{Control of $\Mtop(r\geq r_0)$}

%%%%%%%%%%%%%%%%%%%%%%%%%%%%%%%%%%%%%%%%%%%%%%

In the region $\Mtop(r\geq r_0)$, we proceed as follows: 
\begin{enumerate}
\item We rely on the transport equation in $e_3$ of the ingoing PT structure of $\Mtop$ for the linearized Ricci and metric coefficients, see Proposition \ref{proposition:Maineqts-Mint'}.  

\item We control the transport equations using the first estimate of Corollary \ref{cor:transportrp-L2-knorms-PT:Mtop} according to which, for $U$ satisfying a transport equation with RHS $F$, there holds, for any $c'$,  
\bea\lab{eq:rewritingthemaintransportestimateinMtoprgeqr0fortheproof:chap9}
\nn \sup_{\ub_1\geq u_*'}\int_{\Mtop_{r_0, \ub_1}}r^{2c'}|\dk^{\leq k}U|^2 &\les& \sup_{S\subset\{u=u_*'\}}r^{2c'}\|\dk^{\leq k}U\|^2_{L^2(S)}\\
 &&+ \sup_{\ub_1\geq u_*'}\int_{\Mtop_{r_0, \ub_1}}r^{2c'}|\dk^{\leq k}F|^2,
\eea
where the constant in the definition of  $\les$ is independent of $r_0$.

\item We control the first term of the RHS of \eqref{eq:rewritingthemaintransportestimateinMtoprgeqr0fortheproof:chap9}, i.e. the term involving the control of $U$ on $\{u=u_*'\}$, thanks to Proposition \ref{proposition:Estimatesin-MtopPTonu=ustarprime} on the control of the linearized Ricci and metric coefficients of the ingoing PT frame of $\Mtop$ on $\{u=u_*'\}$ by  $L_*(k)$. 

\item We control the second term on the RHS of \eqref{eq:rewritingthemaintransportestimateinMtoprgeqr0fortheproof:chap9}, i.e. the term involving the control of $F$ on $\Mtop_{r_0, \ub_1}$, using:
\begin{enumerate}
\item The bootstrap assumption \eqref{eq:Bootstrap-Mtop-M8} to control the quadratic terms.

\item The definition of the norm $\Rktop_k$ to control the curvature terms.

\item The triangular structure to control the linear terms involving the Ricci and metric coefficients. More precisely, we need the linear terms appearing on the RHS of the estimate of any given Ricci or metric coefficient to be already under control. As in section \ref{sec:proofofProp:MainestimatesMint666}, this is possible by estimating the quantities in the following order 
\beaa
\trXbc, \, \Xbh, \, \widecheck{\DD\cos\th},\, \Zc, \,\Hbc,\, \DD r, \, e_4(\cos\th),\, \omc,\,\DD\hot\Jk,\, \widecheck{\ov{\DD}\c\Jk}, \, \trXc,\,\Xh,\,\widecheck{e_4(r)},\,\widecheck{\nab_4\Jk},\,\Xi.
\eeaa 
\end{enumerate}
\end{enumerate}

The above scheme yields the desired estimate in $\Mtop(r\geq r_0)$, i.e. 
\bea
 \Sktop^{\geq r_0}_{k} &\les& \ep_0+   L_*(k)  +   \Rktop_{k}, \quad k\leq k_{large}+7,
 \eea
where the constant in $\les$ is independent of $r_0$. Since the proof is reminiscent\footnote{And in fact significantly simpler since $\Mtop$ is a local existence type region, see Remark \ref{rmk:Mtopisinfactalocalexistencetyperegion:chap9}.} of the strategy used both in $\Mext$, see section \ref{sec:proofofProp:EstmatesSkext}, and $\Mint$, see section \ref{sec:proofofProp:MainestimatesMint666}, we leave the details to the reader.

%%%%%%%%%%%%%%%%%%%%%%%%%%%%%%%%%%%%%%%%%%%%%%

\subsubsection{Control of $\Mtop(r\leq r_0)$}

%%%%%%%%%%%%%%%%%%%%%%%%%%%%%%%%%%%%%%%%%%%%%%

In the region $\Mtop(r\leq r_0)$,  the control on solutions of  transport equations given by  \eqref{eq:rewritingthemaintransportestimateinMtoprgeqr0fortheproof:chap9} is replaced by the second estimate of Corollary \ref{cor:transportrp-L2-knorms-PT:Mtop} according to which, for $U$ satisfying a transport equation with RHS $F$, there holds, for any $c'$,  
\beaa
\nn\int_{\Mtop(r\leq r_0)} \big| \dk^{\le k}U|^2 &\les& \sup_{S\subset\{u=u_*'\}}r^{2c'}\|\dk^{\leq k}U\|^2_{L^2(S)}+\sup_\tau\int_{\Mtop(r\leq r_0)\cap\Si(\tau)}\big| \dk^{\le k}F|^2\\
&&+\sup_{\ub_1\geq u_*'}\int_{\Mtop_{r_0, \ub_1}}r^{2c'}|\dk^{\leq k}F|^2.
\eeaa
Then, the proof in the region $\Mtop(r\leq r_0)$ follows the same steps as the one in the region $\Mtop(r\geq r_0)$ and we finally obtain 
\bea
 \Sktop_{k} &\les& \ep_0+   L_*(k)  +   \Rktop_{k}, \quad k\leq k_{large}+7.
 \eea
Again, since the proof is reminiscent\footnote{And in fact significantly simpler since $\Mtop$ is a local existence type region, see Remark \ref{rmk:Mtopisinfactalocalexistencetyperegion:chap9}.} of the strategy used in $\Mint$, see section \ref{sec:proofofProp:MainestimatesMint666}, we leave the details to the reader. This concludes the proof of Proposition \ref{Prop:MainestimatesMtop:chap9}.

%%%%%%%%

\appendix

%%%%%%%%

%%%%%%%%%%%%%%%%%%%%%%%%%%%%%%%%%

\chapter{Proof of results in Chapter   \ref{chapter:preliminaries}}
\lab{appendix-chapter2}

%%%%%%%%%%%%%%%%%%%%%%%%%%%%%%%%%

%%%%%%%%%%%%%%%%%%%%%%%%%%%%%%%%%%%%%%%%%%%%%%%%%

\section{Proof of Corollary \ref{cor:transportequationine4forchangeofframecoeff:simplecasefirst}}
\lab{proofof:cor:transportequationine4forchangeofframecoeff:simplecasefirst}

%%%%%%%%%%%%%%%%%%%%%%%%%%%%%%%%%%%%%%%%%%%%%%%%%

We start with the equation for $f$. Assuming that $\xi'=0$, we have
\beaa
0 &=& 2\xi'_a=\g\left(\D_{e_4'}e_4', e_a'\right) = \la\g\left(\D_{e_4'}(\la^{-1}e_4'), e_a'\right)\\
&=& \la\g\left(\D_{e_4'}(\la^{-1}e_4'), e_a+\frac{1}{2}\fb_a\la^{-1}e_4'+\frac{1}{2}f_ae_3\right)\\
&=& \la\g\left(\D_{e_4'}(\la^{-1}e_4'), e_a+\frac{1}{2}f_ae_3\right)\\
&=& \la^2\g\left(\D_{\la^{-1}e_4'}\left(e_4+f^be_b+\frac{1}{4}|f|^2e_3\right), e_a+\frac{1}{2}f_ae_3\right)
\eeaa
and hence 
\beaa
0 &=& \g\left(\D_{\la^{-1}e_4'}e_4, e_a+\frac{1}{2}f_ae_3\right)+\la^{-1}e_4'(f_a)+f^b\g\left(\D_{\la^{-1}e_4'}e_b, e_a+\frac{1}{2}f_ae_3\right)\\
&&+\frac{1}{4}|f|^2\g\left(\D_{\la^{-1}e_4'}e_3, e_a+\frac{1}{2}f_ae_3\right).
\eeaa

Next, we compute
\beaa
\g\left(\D_{\la^{-1}e_4'}e_4, e_a+\frac{1}{2}f_ae_3\right) &=& \g\left(\D_{e_4+f^be_b+\frac{1}{4}|f|^2e_3}e_4, e_a+\frac{1}{2}f_ae_3\right)\\
&=& 2\xi_a+2\om f_a+f^b\chi_{ba}+(f\c\ze)f_a+\frac{1}{2}|f|^2\eta_a +O(f^3\Ga),
\eeaa
\beaa
\g\left(\D_{\la^{-1}e_4'}e_b, e_a+\frac{1}{2}f_ae_3\right) &=& -\g\left(\D_{\la^{-1}e_4'}\left(e_a+\frac{1}{2}f_ae_3\right), e_b\right) \\
&=& -\g\left(\D_{\la^{-1}e_4'}\left(e_a+\frac{1}{2}f_ae_3\right), e_b+\frac{1}{2}f_be_3\right)\\
&& +\frac{1}{2}f_a\g\left(\D_{\la^{-1}e_4'}\left(e_a+\frac{1}{2}f_ae_3\right), e_3\right)\\
&=& -\g\left(\D_{\la^{-1}e_4'}\left(e_a' -\frac{1}{2}\fb_a \la^{-1}e_4'\right), e_b' -\frac{1}{2}\fb_b\la^{-1}e_4'\right)\\
&& -\frac{1}{2}f_b\g\left(\D_{\la^{-1}e_4'}e_3, e_a+\frac{1}{2}f_ae_3\right)\\
&=& -\g\left(\D_{\la^{-1}e_4'}e_a', e_b'\right)+\fb_a\xi_b'-\fb_b\xi_a' -\frac{1}{2}f_b\g\left(\D_{\la^{-1}e_4'}e_3, e_a\right)\\
&=& -\g\left(\D_{\la^{-1}e_4'}e_a', e_b'\right) -\etab_af_b -\frac{1}{2}f_bf_c\chib_{ca} -\frac{1}{4}|f|^2f_b\xib_a\\
&=& -\g\left(\D_{\la^{-1}e_4'}e_a', e_b'\right) -\etab_af_b +O(f^2\Ga),
\eeaa
and
\beaa
\g\left(\D_{\la^{-1}e_4'}e_3, e_a+\frac{1}{2}f_ae_3\right) &=& 2\etab_a+O(f\Ga).
\eeaa
Plugging in the above, we infer
\beaa
0 &=& 2\xi+2\om f +\frac{1}{2}(\trch f - \atrch\dual f)+\chih\c f+(f\c\ze)f+\frac{1}{2}|f|^2\eta +\nab_{\la^{-1}e_4'}f\\
&& -\frac{1}{2}|f|^2\etab +O(f^3\Ga)
\eeaa
and hence 
\beaa
\nab_{\la^{-1}e_4'}f+\frac{1}{2}(\trch f - \atrch\dual f)+ 2\om f &=& -2\xi -f\c\chih+E_1(f, \Ga)
\eeaa
where 
\beaa
E_1(f, \Ga) &=& -(f\c\ze)f-\frac{1}{2}|f|^2\eta +\frac{1}{2}|f|^2\etab +O(f^3\Ga)
\eeaa
as stated.

Next, we derive the equation for $\la$. Assuming that $\om'=0$ and $\xi'=0$, we have
\beaa
0 &=& 4\om' = \g\left(\D_{e_4'}e_4', e_3'\right)=-2e_4'(\log\la)+\la\g\left(\D_{\la^{-1}e_4'}(\la^{-1}e_4'), \la e_3'\right)\\
&=& -2e_4'(\log\la)+\la\g\left(\D_{\la^{-1}e_4'}(\la^{-1}e_4'), e_3+\fb^ae_a'-\frac{1}{4}|\fb|^2\la^{-1}e_4'\right)\\
&=& -2e_4'(\log\la)+\la\g\left(\D_{\la^{-1}e_4'}(\la^{-1}e_4'), e_3\right)+2\la^{-1}\fb^a\xi_a'\\
&=& -2e_4'(\log\la)+\la\g\left(\D_{\la^{-1}e_4'}(\la^{-1}e_4'), e_3\right)\\
&=& -2e_4'(\log\la)+\la\g\left(\D_{\la^{-1}e_4'}e_4, e_3\right)+\la f^a\g\left(\D_{\la^{-1}e_4'}e_a, e_3\right).
\eeaa
Next, we compute 
\beaa
\g\left(\D_{\la^{-1}e_4'}e_4, e_3\right) &=& 4\om+2f\c\ze - |f|^2\omb
\eeaa
and 
\beaa
\g\left(\D_{\la^{-1}e_4'}e_a, e_3\right) &=& -2\etab_a-f^b\chib_{ba}+O(f^2\Ga).
\eeaa
Plugging in the above, we infer
\beaa
0 &=& -2\la^{-1}e_4'(\log\la) +4\om+2f\c\ze - |f|^2\omb  -2f\c\etab -\frac{1}{2}\trchb|f|^2+O(f^3\Ga+f^2\chibh)
\eeaa
and hence
\beaa
\la^{-1}e_4'(\log\la) &=&  2\om+f\c(\ze-\etab)+E_2(f, \Ga)
\eeaa
where
\beaa
E_2(f, \Ga) &=& - \frac{1}{2}|f|^2\omb  -\frac{1}{4}\trchb|f|^2+O(f^3\Ga+f^2\chibh)
\eeaa
as stated.

Finally, we derive the equation for $\fb$. Summing the transformation formula for $\ze'$ and $\etab'$, and using $\etab'+\ze'=0$, we have 
\beaa
\nab_{\la^{-1}e_4'}\fb +\frac{1}{2}(\trch\fb +\atrch\dual\fb) &=& -2(\etab+\ze) +2\nab'(\log\la) +2\omb f\\
&& +2\err(\etab, \etab')+2\err(\ze, \ze').
\eeaa
Now, in view of the form of $\err(\etab, \etab')$ and $\err(\ze, \ze')$, and since $\xi'=0$, we have
\beaa
&& \err(\etab, \etab')+\err(\ze, \ze')\\ 
&=& \frac{1}{2}f\c\chibh + \frac{1}{2}(f\c\etab)\fb-\frac 1 4  (f\c\ze)\fb  -\frac 1 4 |\fb|^2\la^{-2}\xi'\\
&& -\frac{1}{2}\chibh\c f + \frac{1}{2}(f\c\ze)\fb -  \frac{1}{2}(f\c\etab)\fb +\frac{1}{4}\fb(f\c\eta) + \frac{1}{4}\fb(f\c\ze)  \\
&& + \frac{1}{4}\dual\fb(f\wedge\eta) + \frac{1}{4}\dual\fb(f\wedge\ze) +  \frac{1}{4}\fb\div'f  + \frac{1}{4}\dual\fb \curl'f +\frac{1}{2}\la^{-1}\fb\c\chih' \\
&&  -\frac{1}{16}(f\c\fb)\fb\la^{-1}\trch' +\frac{1}{16}  (\fb\wedge f) \fb\la^{-1}\atrch'  -  \frac{1}{16}\dual\fb ( f\c\fb) \la^{-1}\atrch'\\
&& +\frac{1}{16}\dual\fb \la^{-1}(f\wedge \fb)\trch' +\lot\\
&=& \frac{1}{4}(f\c\eta)\fb + \frac{1}{2}(f\c\ze)\fb  + \frac{1}{4}\dual\fb(f\wedge\eta) + \frac{1}{4}\dual\fb(f\wedge\ze) +  \frac{1}{4}\fb\div'f  \\
&&+ \frac{1}{4}\dual\fb \curl'f +\frac{1}{2}\la^{-1}\fb\c\chih' +O\Big((\la^{-1}\trch', \la^{-1}\atrch')(f, \fb)^3+(f, \fb)^3\Ga\Big)
\eeaa
and hence
\beaa
\nab_{\la^{-1}e_4'}\fb +\frac{1}{2}(\trch\fb +\atrch\dual\fb) = -2(\etab+\ze) +2\nab'(\log\la) +2\omb f+ E_3({\nab'}^{\leq 1}f, \fb, \Ga, \la^{-1}\chi'),
\eeaa
where 
\beaa
E_3({\nab'}^{\leq 1}f, \fb, \Ga, \la^{-1}\chi') &=&  \frac{1}{4}(f\c\eta)\fb + \frac{1}{2}(f\c\ze)\fb    + \frac{1}{4}\dual\fb(f\wedge\eta) + \frac{1}{4}\dual\fb(f\wedge\ze)\\
&& +  \frac{1}{4}\fb\div'f + \frac{1}{4}\dual\fb \curl'f +\frac{1}{2}\la^{-1}\fb\c\chih' \\
&&+O\Big((\la^{-1}\trch', \la^{-1}\atrch')(f, \fb)^3+(f, \fb)^3\Ga\Big)
\eeaa
as stated. This concludes the proof of Corollary \ref{cor:transportequationine4forchangeofframecoeff:simplecasefirst}.

%%%%%%%%%%%%%%%%%%%%%%%%%%%%%%%%%%%%%%%%%%%%%%%%%%%

\section{Proof of Corollary \ref{cor:transportequationine4forchangeofframecoeffinformFFbandlamba}}
\lab{proofof:cor:transportequationine4forchangeofframecoeffinformFFbandlamba}

%%%%%%%%%%%%%%%%%%%%%%%%%%%%%%%%%%%%%%%%%%%%%%%%%%%

The first three identities, i.e. 
\beaa
\nab_{\la^{-1}e_4'}F+\frac{1}{2}\ov{\tr X} F+2\om F &=& -2\Xi -\chih\c F+E_1(f, \Ga),\\
\la^{-1}\nab_4'(\log\la) &=& 2\om+f\c(\ze-\etab)+E_2(f, \Ga),\\
\nab_{\la^{-1}e_4'}\underline{F}+\frac{1}{2}\tr X\underline{F} &=& -2(\Hb+Z)   +2\DD'(\log\la)  +2\omb F  +E_3({\nab'}^{\leq 1}f, \fb, \Ga, \la^{-1}\chi'),
\eeaa
are an immediate consequence of Corollary \ref{cor:transportequationine4forchangeofframecoeff:simplecasefirst} and the notation for $F$ and $\underline{F}$. 

Next, we focus on proving the fourth identity. Since we assume $e_4(q)=1$, which also yields $e_4(\ov{q})=1$, we have 
\beaa
\la^{-1}e_4'(q) &=& \left(e_4+f\c\nab+\frac{1}{4}|f|^2e_3\right)q=1+\left(f\c\nab+\frac{1}{4}|f|^2e_3\right)q,\\
\la^{-1}e_4'(\ov{q}) &=& \left(e_4+f\c\nab+\frac{1}{4}|f|^2e_3\right)\ov{q}=1+\left(f\c\nab+\frac{1}{4}|f|^2e_3\right)\ov{q}.
\eeaa
We infer  
\beaa
\nab_{\la^{-1}e_4'}(\ov{q}F) &=& \ov{q}\nab_{\la^{-1}e_4'}F+\la^{-1}e_4'(\ov{q})F\\
&=& \ov{q}\left(-\frac{1}{2}\ov{\tr X} F-2\om F  -2\Xi -\chih\c F+E_1(f, \Ga)\right)+F+\left(f\c\nab+\frac{1}{4}|f|^2e_3\right)F\\
&=& -2\ov{q}\om F  -2\ov{q}\Xi -\frac{1}{2}\ov{q}\left(\tr X-\frac{2}{\ov{q}}\right)F -\ov{q}\chih\c F+\ov{q}E_1(f, \Ga)+f\c\nab(\ov{q})F\\
&&+\frac{1}{4}|f|^2e_3(\ov{q})F
\eeaa
and hence
\beaa
\nab_{\la^{-1}e_4'}(\ov{q}F) &=& -2\ov{q}\om F  -2\ov{q}\Xi +E_4(f, \Ga),
\eeaa
where $E_4(f, \Ga)$ is given by
\beaa
E_4(f, \Ga) &=& -\frac{1}{2}\ov{q}\left(\tr X-\frac{2}{\ov{q}}\right)F -\ov{q}\chih\c F+\ov{q}E_1(f, \Ga)+f\c\nab(\ov{q})F+\frac{1}{4}|f|^2e_3(\ov{q})F
\eeaa
as stated.

Finally, we consider the fifth and last identity. Using the first four identities and the above identities for $\la^{-1}e_4'(q)$ and $\la^{-1}e_4'(\ov{q})$, we have, using also $e_4(e_3(r))=-2\omb$, 
\beaa
&&\nab_{\la^{-1}e_4'}\left[q\Big(\underline{F}-2q\DD'(\log\la)\Big)+\ov{q}e_3(r)F\right] \\
&=& q\nab_{\la^{-1}e_4'}\underline{F} +q\Big(-2q\DD'\nab_{\la^{-1}e_4'}(\log\la) -2[\nab_{\la^{-1}e_4'},q\DD']\log\la\Big)\\
&&+\la^{-1}e_4'(q)\Big(\underline{F}-2q\DD'(\log\la)\Big) +e_3(r)\nab_{\la^{-1}e_4'}(\ov{q}F)+ \la^{-1}e_4'(e_3(r))\ov{q}F\\
&=& q\left(-\frac{1}{2}\tr X\underline{F}  -2(\Hb+Z)   +2\DD'(\log\la)  +2\omb F  +E_3({\nab'}^{\leq 1}f, \fb, \Ga, \la^{-1}\chi')\right)\\
&& -2q^2\DD'\Big(2\om+f\c(\ze-\etab)+E_2(f, \Ga)\Big) -2q[\nab_{\la^{-1}e_4'},q\DD']\log\la\\
&&+\underline{F}-2q\DD'(\log\la)+\left(f\c\nab(q)+\frac{1}{4}|f|^2e_3(q)\right)\Big(\underline{F}-2q\DD'(\log\la)\Big)\\
&&+e_3(r)\Big(-2\ov{q}\om F  -2\ov{q}\Xi +E_4(f, \Ga)\Big) -2\omb\ov{q}F.
\eeaa
We infer
\beaa
\nab_{\la^{-1}e_4'}\left[q\Big(\underline{F}-2q\DD'(\log\la)\Big)+\ov{q}e_3(r)F\right] &=&  -2q(\Hb+Z)  -2q^2\DD'\Big(2\om+f\c(\ze-\etab)\Big)\\
&& +e_3(r)\left(-2\ov{q}\om F  -2\ov{q}\Xi\right) +2\omb(q-\ov{q})F \\
&&+E_5({\nab'}^{\leq 1}f,\fb, {\nab'}^{\leq 1}\la, \D^{\leq 1}\Ga)
\eeaa
where $E_5({\nab'}^{\leq 1}f,\fb, {\nab'}^{\leq 1}\la, \D^{\leq 1}\Ga)$ is given by
\beaa
E_5({\nab'}^{\leq 1}f,\fb, {\nab'}^{\leq 1}\la, \D^{\leq 1}\Ga) &=&-\frac{q}{2}\left(\tr X-\frac{2}{q}\right)\underline{F}       +qE_3({\nab'}^{\leq 1}f, \fb, \Ga, \la^{-1}\chi')\\
&& -2q^2\DD'\left(E_2(f, \Ga)\right)  -2q[\nab_{\la^{-1}e_4'},q\DD']\log\la\\
&&+\left(f\c\nab(q)+\frac{1}{4}|f|^2e_3(q)\right)\Big(\underline{F}-2q\DD'(\log\la)\Big)\\
&&+e_3(r)E_4(f, \Ga)
\eeaa
as stated. This concludes the proof of Corollary \ref{cor:transportequationine4forchangeofframecoeffinformFFbandlamba}.

%%%%%%%%%%%%%%%%%%%%%%%%%%%%%%%%%%%%%%%%%%%

\section{Proof of Proposition \ref{prop:modesforwidecheck{trch, trchb}-Kerr}}
\lab{appendix-proofofprop:modesforwidecheck{trch, trchb}-Kerr}

%%%%%%%%%%%%%%%%%%%%%%%%%%%%%%%%%%%%%%%%%%%%

In the following lemma, we prove  refined  asymptotic   for $\trch'$, and $\trchb'$ which will be used in the proof of Proposition \ref{prop:modesforwidecheck{trch, trchb}-Kerr} below.
\begin{lemma}\lab{lemma:precise(trch'- trchb')-Kerr}
We have, for large $r$,
\bea
\lab{eq:trchtrchb-inKerr-refined}
\bsplit
\trch' &=  \frac{2}{r}  -\frac{2a^2(\cos\th)^2}{r^3}   -\frac{a^2(\sin\th)^2}{2r^3}  +O\left(\frac{1}{r^4}\right),
\\
\trchb' &= -\frac{2\left(1-\frac{2m}{r}\right)}{r}  -\frac{2a^2}{r^3}\left(1-2(\cos\th)^2\right) +\frac{3a^2(\sin\th)^2}{2r^3}+O\left(\frac{1}{r^4}\right).
\end{split}
\eea
\end{lemma}

\begin{proof}
We write, making use of  the transformation formulas in Proposition \ref{Proposition:transformationRicci},
\beaa
\bsplit
\la^{-1}\trch' &= \trch  +  \div'f + f\c\eta + f\c\ze+\err(\trch,\trch'),\\
\err(\trch,\trch') &= \fb\c\xi+\frac{1}{4}\fb\c\left(f\trch -\dual f\atrch\right) +\om (f\c\fb)  -\omb |f|^2 \\
& -\frac{1}{4}|f|^2\trchb -  \frac 1 4 ( f\c\fb) \la^{-1}\trch' +\frac 1 4  (\fb\wedge f) \la^{-1}\atrch'+\lot,
\\
\la\trchb' &= \trchb +\div'\fb +\fb\c\etab  -  \fb\c\ze +\err(\trchb, \trchb'),\\
\err(\trchb, \trchb') &= \frac{1}{2}(f\c\fb)\trchb+f\c\xib -|\fb|^2\om + (f\c\fb)\omb   -\frac 1 4 |\fb|^2\la^{-1}\trch'+\lot
\end{split}
\eeaa
Also, recall from Lemma \ref{Lemma:Transformation-principal-to-integrable-Kerr} the following asymptotic:
\begin{itemize}
\item We have $f_1=\fb_1=0$, and
 \beaa
f_2 =-\frac{a\sin\th}{r}+ O\big( \sin\th r^{-3} \big), \qquad \fb_2=- \frac{a\sin\th \Up}{r}+ O\big( \sin\th r^{-3} \big).
\eeaa

\item We have,
\beaa
\bsplit
 \div' f&=O(\sin^2\th r^{-5} ), \qquad   \nab'\hot f=O(\sin \th r^{-5} ),\\
  \div' \fb&=O(\sin^2 \th r^{-5} ), \qquad   \nab'\hot \fb=O(\sin \th r^{-5} ).
  \end{split}
 \eeaa 
\end{itemize}
Together with the asymptotic for the outgoing principal frame of Kerr, Corollary \ref{corollary:Kerr.PG.asymptotics}, we deduce
\beaa
\bsplit
\la^{-1}\trch' &= \trch  + f_2\eta_2 + f_2\ze_2+\frac{1}{4}\fb_2f_2\trch   
 -\frac{1}{4}(f_2)^2\trchb -  \frac 1 4 f_2\fb_2\trch' +O\left(\frac{1}{r^4}\right)\\
&= \trch  -\frac{3a^2(\sin\th)^2}{2r^3}  +O\left(\frac{1}{r^4}\right),
\\
\la\trchb' &= \trchb  +\fb_2\etab_2  -  \fb_2\ze_2 +\frac{1}{2}f_2\fb_2\trchb   
 -\frac 1 4 |\fb|^2\la^{-1}\trch' +O\left(\frac{1}{r^4}\right)\\
&= \trchb  +\frac{a^2(\sin\th)^2}{2r^3}+O\left(\frac{1}{r^4}\right).
\end{split}
\eeaa
Recall that, relative to the   principal frame of Kerr
\beaa
\trch &=& \frac{2r}{|q|^2} =\frac{2r}{r^2+a^2(\cos\th)^2}= \frac{2}{r}-\frac{2a^2(\cos\th)^2}{r^3}+O\left(\frac{1}{r^5}\right),\\
\trchb &=& -\frac{2r\Delta}{|q|^4}=-\frac{2r\left(r^2-2mr+a^2\right)}{(r^2+a^2(\cos\th)^2)^2}
= -\frac{2\Up}{r}-\frac{2a^2}{r^3}\left(1-2(\cos\th)^2\right)+O\left(\frac{1}{r^4}\right).
\eeaa
We infer that
\beaa
\bsplit
\la^{-1}\trch' &=  \frac{2}{r}  -\frac{2a^2(\cos\th)^2}{r^3}   -\frac{3a^2(\sin\th)^2}{2r^3}  +O\left(\frac{1}{r^4}\right),
\\
\la\trchb' &= -\frac{2\left(1-\frac{2m}{r}\right)}{r}  -\frac{2a^2}{r^3}\left(1-2(\cos\th)^2\right) +\frac{a^2(\sin\th)^2}{2r^3}+O\left(\frac{1}{r^4}\right).
\end{split}
\eeaa

To conclude, we need to derive an improved asymptotic for $\la$ compared to the one in Lemma \ref{Lemma:Transformation-principal-to-integrable-Kerr}. Recall from \eqref{def:transition-functs:ffbla} that 
\beaa
\la &=& 1+\frac{1}{2}f\c\fb  +\frac{1}{16}|f|^2|\fb|^2.
\eeaa
Together with the above asymptotic for $f$ and $\fb$, this  yields
 \beaa
 \la=1+\frac{a^2(\sin\th)^2}{2r^2}+O(r^{-3})
 \eeaa
and hence 
\beaa
\bsplit
\trch' &=  \frac{2}{r}  -\frac{2a^2(\cos\th)^2}{r^3}   -\frac{a^2(\sin\th)^2}{2r^3}  +O\left(\frac{1}{r^4}\right),
\\
\trchb' &= -\frac{2\left(1-\frac{2m}{r}\right)}{r}  -\frac{2a^2}{r^3}\left(1-2(\cos\th)^2\right) +\frac{3a^2(\sin\th)^2}{2r^3}+O\left(\frac{1}{r^4}\right),
\end{split}
\eeaa
as stated.
\end{proof}

In the following lemma, we derive the asymptotic of the area radius of $S(u,r)$ that will be used in the proof of  Proposition \ref{prop:modesforwidecheck{trch, trchb}-Kerr} below.
\begin{lemma}\lab{lemma:asymptoticarearadiusSur}
Let $r'$ denote the area radius of $S(u, r)$. $r'$ verifies
\beaa
r' &=& r\left(1+\frac{a^2}{3r^2}+O\left(\frac{1}{r^3}\right)\right).
\eeaa
\end{lemma}

\begin{proof}
Let $g$ denote the induced metric on $S(u,r)$. Then, $(\th, \vphi)$ forms a coordinates system on $S(u,r)$, and we have
\beaa
g_{\th\th}=\g_{\th\th}=|q|^2, \qquad g_{\th\vphi}=\g_{\th\vphi}=0, \qquad g_{\vphi\vphi}=\g_{\vphi\vphi}=\frac{\Si^2(\sin\th)^2}{|q|^2}.
\eeaa
We deduce
\beaa
\sqrt{g_{\th\th}g_{\vphi\vphi}-(g_{\th\vphi})^2} &=& \Si\sin\th.
\eeaa
Now, in view of the definition of $\Si$, we have
\beaa
\Si &=& \sqrt{(r^2+a^2)^2-a^2(\sin\th)^2(r^2-2mr+a^2)}\\
&=& r^2\sqrt{1+\frac{a^2(2-(\sin\th)^2)}{r^2}+O\left(\frac{1}{r^3}\right)}\\
&=& r^2\left(1+\frac{a^2(1+(\cos\th)^2)}{2r^2}+O\left(\frac{1}{r^3}\right)\right)
\eeaa
and hence
\beaa
\sqrt{g_{\th\th}g_{\vphi\vphi}-(g_{\th\vphi})^2} &=& r^2\sin\th\left(1+\frac{a^2(1+(\cos\th)^2)}{2r^2}+O\left(\frac{1}{r^3}\right)\right).
\eeaa
We deduce
\beaa
4\pi(r')^2 &=& |S(u,r)|=\int_0^{2\pi}\int_0^{\pi}\sqrt{g_{\th\th}g_{\vphi\vphi}-(g_{\th\vphi})^2}d\th d\vphi\\
&=& 2\pi r^2\int_0^\pi \sin\th\left(1+\frac{a^2(1+(\cos\th)^2)}{2r^2}+O\left(\frac{1}{r^3}\right)\right) d\th\\
&=& 4\pi r^2\left(1+\frac{2a^2}{3r^2}+O\left(\frac{1}{r^3}\right)\right)
\eeaa
and hence
\beaa
r' &=& r\left(1+\frac{a^2}{3r^2}+O\left(\frac{1}{r^3}\right)\right)
\eeaa
as stated.
\end{proof}

\begin{proof}[Proof of Proposition \ref{prop:modesforwidecheck{trch, trchb}-Kerr}]
We make use of the refined asymptotic calculation $ \trch', \trchb'$ in \eqref{eq:trchtrchb-inKerr-refined},  i.e.
\beaa
\bsplit
\trch' &=  \frac{2}{r}  -\frac{2a^2(\cos\th)^2}{r^3}   -\frac{a^2(\sin\th)^2}{2r^3}  +O\left(\frac{1}{r^4}\right),
\\
\trchb' &= -\frac{2\left(1-\frac{2m}{r}\right)}{r}  -\frac{2a^2}{r^3}\left(1-2(\cos\th)^2\right) +\frac{3a^2(\sin\th)^2}{2r^3}+O\left(\frac{1}{r^4}\right).
\end{split}
\eeaa
 Since the area radius $r'$ of $S(u, r)$ verifies, in view of Lemma \ref{lemma:asymptoticarearadiusSur},
\beaa
r' &=& r\left(1+\frac{a^2}{3r^2} +O\left(\frac{1}{r^3}\right)\right),
\eeaa
and since we have by definition 
\beaa
 \widecheck{\trch}'=\trch' -\frac{2}{r'}, \qquad   \widecheck{\trchb}'=\trchb' +\frac{2\left(1-\frac{2m}{r'}\right)}{r'},
\eeaa
we deduce
\bea
\lab{eq:widecheck{trch}'-widecheck{trch}'-Kerr}
\bsplit
 \widecheck{\trch}'&=  \frac{2a^2}{3r^3}  -\frac{2a^2(\cos\th)^2}{r^3}   -\frac{a^2(\sin\th)^2}{2r^3}  +O\left(\frac{1}{r^4}\right),
\\
 \widecheck{\trchb}'&= -\frac{2a^2}{3r^3} -\frac{2a^2}{r^3}\left(1-2(\cos\th)^2\right) +\frac{3a^2(\sin\th)^2}{2r^3}+O\left(\frac{1}{r^4}\right).
\end{split}
\eea
Since the $O(r^{-3})$ terms in the expansion of $\widecheck{\trch}'$ and $\widecheck{\trchb}'$ do not depend on $\vphi$, and in view of the form of $J^{(+)}$ and $J^{(-)}$, we immediately  deduce that
\beaa
\int_{S(u,r)}\widecheck{\trch}'J^{(p)} &=& O\left(\frac{1}{r^2}\right),\quad p=+,-,\\
\int_{S(u,r)}\widecheck{\trchb}'J^{(p)} &=& O\left(\frac{1}{r^2}\right),\quad p=+,-.
\eeaa

It remains to consider the case $p=0$. In view of \eqref{eq:widecheck{trch}'-widecheck{trch}'-Kerr} and  the form of $J^{(0)}$, this case   follows from  the  calculation
\beaa
\int_0^\pi \Big(1, (\cos\th)^2, (\sin\th)^2\Big)\cos\th\sin\th d\th=(0,0,0).
\eeaa
This concludes the proof of Proposition \ref{prop:modesforwidecheck{trch, trchb}-Kerr}.
\end{proof}

%%%%%%%%%%%%%%%%%%%%%%%%%%%%%%%%%%%%%%%%%%%%

\section{Proof of Proposition  \ref{Proposition:deftensorT}}
\lab{appendix:Proof-Proposition:deftensorT}

%%%%%%%%%%%%%%%%%%%%%%%%%%%%%%%%%%%%%%%%%%%%

We  decompose the proof in the  following steps.

\textit{Calculation of $\pi_{33}$.}   We have, using the fact that  $\omb= \ombc+ \frac 1 2  \pr_r\Big(\frac{\De}{|q|^2} \Big)$,
\beaa
2\g(\D_3\T, e_3) &=& \g\left(\D_3\left(e_3+\frac{\Delta}{|q|^2}e_4 - 2a\Re(\Jk)^be_b\right), e_3\right)\\
&=& -2e_3\left(\frac{\Delta}{|q|^2}\right) -\frac{4\Delta}{|q|^2}\omb +4a\Re(\Jk)\c\xib\\
&=& -2\pr_r\left(\frac{\Delta}{|q|^2}\right)e_3(r) -2\pr_{\cos\th}\left(\frac{\Delta}{|q|^2}\right)e_3(\cos\th) -\frac{4\Delta}{|q|^2}\omb +4a\Re(\Jk)\c\xib\\
&=& -2\pr_r\left(\frac{\Delta}{|q|^2}\right)\widecheck{e_3(r)} -2\pr_{\cos\th}\left(\frac{\Delta}{|q|^2}\right)e_3(\cos\th) -\frac{4\Delta}{|q|^2}\ombc +4a\Re(\Jk)\c\xib\\
&=& \Ga_b.
\eeaa
\textit{Calculation of $\pi_{44}$.} 
\beaa
2\g(\D_4\T, e_4) &=& \g\left(\D_4\left(e_3+\frac{\Delta}{|q|^2}e_4 -2a\Re(\Jk)^be_b, e_4\right)\right)= -4\om + 4a\Re(\Jk)\c\xi\\
&=& 0.
\eeaa

 \textit{Calculation of $\pi_{34}$.} Making use of
\beaa
\om=0, \qquad \omb= \frac{1}{2} \pr_r\left(\frac{\De}{|q|^2} \right)+\Ga_b, \qquad H+\Hb= \frac{a(q-\ov{q})} {|q|^2} \Jk+\Re(\widecheck{H} -\widecheck{Z}),\qquad \dual \Jk=-i \Jk,
\eeaa
we deduce
\beaa
2\g(\D_4\T, e_3)&=& \g\left(\D_3\left(e_3+\frac{\Delta}{|q|^2}e_4 - 2a\Re(\Jk)^be_b\right), e_4\right)= 4\omb +4a\Re(\Jk)\c\eta, \\
2\g(\D_3\T, e_4)&=& \g\left(\D_4\left(e_3+\frac{\Delta}{|q|^2}e_4 - 2a\Re(\Jk)^be_b\right), e_3\right)=-2e_4\left(\frac{\Delta}{|q|^2}\right) +4a\Re(\Jk)\c\etab,
\eeaa
and hence
\beaa
2\g(\D_4\T, e_3)+2\g(\D_3\T, e_4) 
&=& 4 \ombc+ 4a\Re(\Jk)\c(\eta+\etab)\\
&=&  4 \ombc+ 4a\Re(\Jk)\c\bigg(\Re\Big( \frac{a(q-\ov{q})} {|q|^2} \Jk\Big)  +\Re(\widecheck{H} -\widecheck{Z})\bigg)   \\
&=& 4a^2   \frac{2 a  \cos \th}  {|q|^2} \Re(\Jk)\c\Re \Big( i\Jk\Big) +\Ga_b \\
&=& -4a^2   \frac{2 a  \cos \th}  {|q|^2} \Re(\Jk)\c \dual \Re \Big( \Jk\Big) +\Ga_b=\Ga_b.
\eeaa

 \textit{Calculation of $\pi_{4a}$.}  Next, using
\beaa
\xi=0, \qquad \ze=\Re\left(\frac{a\ov{q} }{|q|^2 } \right)+ \Ga_g, \qquad  \, \tr X= \frac{2}{q} +\Ga_g,\qquad \nab_4\Jk=-q^{-1}\Jk,\qquad \dual \Jk=-i\Jk,
\eeaa
we deduce
\beaa
2\g(\D_4\T, e_c)&=& \g\left(\D_4\left(e_3+\frac{\Delta}{|q|^2}e_4 -2a\Re(\Jk)^be_b\right), e_c\right)=2\etab_c+\frac{2\Delta}{|q|^2}\xi_c -2a\nab_4\Re(\Jk)_c \\
&=& -2\Re\left(\frac{a\ov{q}}{|q|^2}\Jk\right)_c +2 a \Re\big(q^{-1}\Jk\big)_c+\Ga_g=\Ga_g,
\eeaa
\beaa
2\g(\D_c\T, e_4)&=&\g\left(\D_c\left(e_3+\frac{\Delta}{|q|^2}e_4 -2a\Re(\Jk)^be_b\right), e_4\right)\\
&=&-2\ze_c+2a\Re(\Jk)_b\chi_{cb}= -2\ze_c+a\Re(\Jk)_b\big(\trch\de_{bc}+\atrch\in_{bc}\big)+ r^{-1}\Ga_g\\
&=&- 2 \Re\left(\frac{2a\ov{q}}{|q|^2} \Jk\right)_c+  a \Big(\trch  \Re(\Jk)_c+ \atrch \dual \Re(\Jk)_c\Big)+\Ga_g\\
&=& - 2 \Re\left(\frac{2a\ov{q}}{|q|^2} \Jk\right)_c+  a \left(\frac{2 r  }{|q|^2}   \Re(\Jk)_c+ \frac{2a\cos\th}{|q|^2}  \dual \Re(\Jk)_c\right)+\Ga_g \\
&=&\frac{2 a^2\cos\th}{|q|^2} \Re\big(i\Jk+ \dual \Jk\big)_c +\Ga_g=\Ga_g. 
\eeaa
Thus both $\g(\D_4\T, e_c)$ and  $\g(\D_c\T, e_4) $ are $\Ga_g$  and so is $\piT_{c4} $.

  \textit{Calculation of $\pi_{3a}$.} 
  Since $ \nab_3\Jk=  \frac{\De q}{|q|^4}\Jk+  \widecheck{\nab_3\Jk}=\frac{\De q}{|q|^4}\Jk+ r^{-1} \Ga_b$ and   $\eta=\Re\Big( \frac{aq}{|q|^2} \Jk\Big)+\Ga_b$,
  \beaa
2\g(\D_3\T, e_c)&=& \g\left(\D_3\left(e_3+\frac{\Delta}{|q|^2}e_4  -2a\Re(\Jk)^be_b\right), e_c\right)\\
&=&\frac{2\Delta}{|q|^2}\eta_c -2a\nab_3\Re(\Jk)_c +\Ga_b     =   \frac{2\Delta}{|q|^2}\eta_c - 2 a \Re\left( \frac{\De q}{|q|^4}\Jk\right)_c  + r^{-1}\Ga_b  \\
&=&\frac{2\Delta}{|q|^2} \Re\left( \frac{aq}{|q|^2} \Jk\right)_c - 2 a \Re\left( \frac{\De q}{|q|^4}\Jk\right)_c  + \Ga_b=\Ga_b.
  \eeaa
  Also, since $\dual \Jk=-i\Jk$,
  \beaa
  2\g(\D_c\T, e_3) &=&\g\left(\D_c\left(e_3+\frac{\Delta}{|q|^2}e_4  -2a\Re(\Jk)^be_b\right), e_3\right)\\
&=&\frac{2\Delta}{|q|^2}\ze_c +2a\Re(\Jk)_b\chib_{cb}= \frac{2\Delta}{|q|^2}\ze_c +2a\Re(\Jk)_b\left(\frac 1 2 \de_{cb} \trchb +\frac 1 2 \in_{cb} \atrchb\right)+r^{-1} \Ga_b\\
&=&  \frac{2\Delta}{|q|^2}\ze_c +  a \big( \trchb \Re(\Jk)_c+\atrch \dual\Re( \Jk) _c\big) +r^{-1} \Ga_b\\
&=&  \frac{2\Delta}{|q|^2}\Re\left(\frac{a\ov{q} \Jk }{|q|^2} \right)_c+ a\left(-\frac{2r\De}{|q|^4}\Re(\Jk)_c-\frac{2a\De \cos \th}{|q|^4} \dual \Re(\Jk)_c \right) +r^{-1}\Ga_b\\
&=&\frac{2\Delta}{|q|^2}\Re\left(\frac{-  a^2i \cos \th \Jk }{|q|^2} \right)_c -\frac{2a^2\De \cos \th}{|q|^4} 
\dual \Re(\Jk)_c + r^{-1}\Ga_b\\
&=&-\frac{2a^2\De \cos \th}{|q|^4}  \Re\big( i \Jk+\dual  \Jk\big)_c+r^{-1}\Ga_b= r^{-1} \Ga_b.
\eeaa
  Thus  both $\g(\D_3\T, e_c)$ and  $\g(\D_c\T, e_3) $ are $\Ga_b$ and so is $\piT_{3c}$.

  \textit{Calculation of $\pi_{ab}$.}   We make us of the assumptions on  $ r\,\widecheck{\ov{\DD}\c\Jk}, \,\, r\,\DD\hot\Jk\in  r^{-1}\Ga_b$
  to deduce
\beaa
2\g(\D_c\T, e_d) &=& \g\left(\D_c\left(e_3+\frac{\Delta}{|q|^2}e_4  -2a\Re(\Jk)^be_b\right), e_d\right)=\chib_{cd}+\frac{\Delta}{|q|^2}\chi_{cd} -2a\nab_c\Re(\Jk)_d\\
&=&\frac 1 2 \de_{cd} \Big(\trchb + \frac{\Delta}{|q|^2}\trch\Big)+\frac 1 2 \in_{ cd} \Big(\atrchb + \frac{\Delta}{|q|^2}\atrch\Big) -2a\nab_c\Re(\Jk)_d +\Ga_b \\
&=&\frac 1 2\Big(\frac{-2r\De} {|q|^4} +  \frac{\Delta}{|q|^2} \big( \frac{2r}{|q|^2}\big)\Big)\de_{cd} +\frac 1 2 \in_{ cd} \Big( \frac{2a\Delta\cos\th}{|q|^4} + \frac{\Delta}{|q|^2} \frac{2a \cos\th}{|q|^2}\Big)- 2 a \nab_c\Re(\Jk_d) \\
&=& \frac 1 2 \in_{ cd} \Big( \frac{2a\Delta\cos\th}{|q|^4} + \frac{\Delta}{|q|^2} \frac{2a \cos\th}{|q|^2}\Big)- 2 a \nab_c\Re(\Jk_d) +\Ga_b\\
&=&-\frac{4amr\cos \th}{|q|^4}\in_{ab}  +\Ga_b.
\eeaa
Hence
\beaa
2\g(\D_c\T, e_d)  +2\g(\D_d\T, e_c)&=&- 2 a\big(  \nab_c\Re(\Jk_d) + \nab_d\Re(\Jk_c) \big)+\Ga_b\\
&=& -2a\Re(\nab\hot\Jk)_{cd} -2a\Re(\widecheck{\div\Jk})\de_{cd}+\Ga_b=\Ga_b.
\eeaa
Moreover, using the fact that $\chibh$, $\chih$ and $\nab\hot\Jk$, are traceless,
\beaa
2g^{cd}\,{}^{(\T)}\pi_{cd} &=& \trchc +\frac{\Delta}{|q|^2}\trchbc -2a\Re(\widecheck{\div\Jk})\\
&=& \Ga_g.
\eeaa
This concludes the proof of Proposition  \ref{Proposition:deftensorT}.

%%%%%%%%%%%%%%%%%%%%%%%%%%%%%%%%%%%%%%%%%%%%%

\chapter{Proof of results in Chapter   \ref{Chapter:decaySigmastar}}
\lab{appendix-chapter5}

%%%%%%%%%%%%%%%%%%%%%%%%%%%%%%%%%%%%%%%%%%%%%

The proofs in this section rely on the   the linearized null structure and  null Bianchi  identities of Proposition \ref{Prop.NullStr+Bianchi-lastslice}.

%%%%%%%%%%%%%%%%%%%%%%%%%%%%%%%%%%%%%%%

\section{Proof of Proposition \ref{prop:additional.eqtsM4}}
\lab{sec:proofofprop:additional.eqtsM4}

%%%%%%%%%%%%%%%%%%%%%%%%%%%%%%%%%%%%%%%

 We proceed as follows.
 
{\bf Step 1.}  Recall from Proposition \ref{Prop.NullStr+Bianchi-lastslice} the linearized the null structure equation for $\nab_3\ze $
\beaa
 \nab_3 \ze -\frac{\Up}{r} \ze &=& -\bb-2 \nab \ombc +\frac{\Up}{r} (\eta+\ze)+\frac{1}{r}\xib+ \frac{2m}{r^2}(\ze-\eta) +\Ga_b\c \Ga_b.
\eeaa
It can be rewritten in the form
\beaa
2 \nab \ombc -\frac{1}{r}\xib &=& - \nab_3 \ze   -\bb +\frac{1}{r}\eta +r^{-1}\Ga_g+\Ga_b\c \Ga_b
\eeaa
which is the first stated identity of the proposition.

{\bf Step 2.} We make use of the equations for $\ka$  and $\curl \eta$ in Proposition \ref{Prop.NullStr+Bianchi-lastslice}
\beaa
\nab_3\kac  &=&         2   \div \eta + 2\rhoc -\frac{1}{r}\kabc + \frac{4}{r} \ombc +\frac{2}{r^2}\widecheck{y} +\Ga_b\c\Ga_b,\\
\curl\eta &=& \rhod +\Ga_b\c \Ga_g,
\eeaa
which we rewrite in the form
\beaa
 2   \div \eta &=& \nab_3\kac    - \frac{4}{r} \ombc -\frac{2}{r^2}\widecheck{y}+r^{-1}\Ga_g +\Ga_b\c\Ga_b,\\
  2 \curl\eta &=&  r^{-1}\Ga_g +\Ga_b\c \Ga_g.
\eeaa
Recalling that $\ddd_1=(\div, \curl)$, we rewrite in the form
\beaa
2\ddd_1\eta&=&\left(\nab_3\kac   - \frac{4}{r} \ombc -\frac{2}{r^2}\widecheck{y}, 0 \right)+r^{-1}\Ga_g+\Ga_b\c \Ga_b.
\eeaa
Next, recall that for a pair of scalar functions $(f,h)$, $\dds_1(f, h)=-\nab f+\dual\nab h$.
Hence
\beaa
2\dds_1 \ddd_1\eta &=&-\nab\left(\nab_3\kac    - \frac{4}{r} \ombc -\frac{2}{r^2}\widecheck{y}\right)+r^{-2}\dkb^{\leq 1}\Ga_g+ r^{-1} \dkb^{\le 1 } (\Ga_b\c \Ga_b)\\
&=& -\nab_3\nab\kac -[\nab,\nab_3]\kac  +\frac{4}{r} \nab\ombc +\frac{2}{r^2}\nab\widecheck{y} + r^{-2} \dkb^{\le 1} \Ga_g+ r^{-1} \dkb^{\le 1 } (\Ga_b\c \Ga_b).
\eeaa
Now, in view of  Lemma \ref{Lemma:eqts-nabeta,xib},  we have
\beaa
\nab   \widecheck{y} &=& \nab y = -\xib+\big(\ze-\eta) y= -\xib -\Up\big(\ze-\eta) +r\Ga_b\c\Ga_b\\
&=&  -\xib +\eta+\Ga_g +r\Ga_b\c\Ga_b,
\eeaa
where we have also used $y=-\Up+r\Ga_b$. We infer
\beaa
2\dds_1 \ddd_1\eta &=& -\nab_3\nab\kac -[\nab,\nab_3]\kac  +\frac{4}{r} \nab\ombc -\frac{2}{r^2}\xib +\frac{2}{r^2}\eta + r^{-2} \dkb^{\le 1} \Ga_g+ r^{-1} \dkb^{\le 1 } (\Ga_b\c \Ga_b).
\eeaa

We make use  the commutation formula of Lemma \ref{lemma:comm_Si_*} in the particular case of a scalar 
  \beaa
-[\nab, \nab_3]\kac &=&[\nab_3, \nab] \kac =(\eta-\ze)\nab_3\kac  +\xib \nab_4 \kac+ r^{-2} \dkb^{\le 1} \Ga_g+ r^{-1} \dkb^{\le 1 } (\Ga_b\c \Ga_b).
\eeaa
Making use of the equations
\beaa
\nab_3\kac  &=&         2   \div \eta + 2\rhoc -\frac{1}{r}\kabc + \frac{4}{r} \ombc +\frac{2}{r^2}\widecheck{y} +\Ga_b\c\Ga_b,\\
\nab_4\kac  &=& \Ga_g\c\Ga_g,
\eeaa
we deduce,
\beaa 
(\eta-\ze)\nab_3\kac  +\xib \nab_4 \kac &=& r^{-1}  \dkb^{\le 1 } (\Ga_b\c \Ga_b). 
\eeaa
Hence,
\beaa
-[\nab, \nab_3]\kac &=&  r^{-2} \dkb^{\le 1} \Ga_g+ r^{-1}  \dkb^{\le 1 } (\Ga_b\c \Ga_b) 
\eeaa
and therefore,
\beaa
2\dds_1 \ddd_1\eta &=& -\nab_3\nab\kac   +\frac{4}{r} \nab\ombc -\frac{2}{r^2}\xib +\frac{2}{r^2}\eta + r^{-2} \dkb^{\le 1} \Ga_g+ r^{-1} \dkb^{\le 1 } (\Ga_b\c \Ga_b).
\eeaa
We now make use of Step 1 to substitute $\nab\omb$ and deduce
\beaa
2\dds_1 \ddd_1\eta &=& -\nab_3\nab\kac -\frac{2}{r^2}\xib +\frac{2}{r^2}\eta  +\frac{2}{r}\left(\frac{1}{r}\xib  - \nab_3 \ze   -\bb +\frac{1}{r}\eta\right) \\
&& + r^{-2} \dkb^{\le 1} \Ga_g+ r^{-1} \dkb^{\le 1 } (\Ga_b\c \Ga_b)\\
&=&  -\nab_3\nab\kac  +\frac{4}{r^2}\eta  -\frac{2}{r}\nab_3 \ze  -\frac{2}{r}\bb  + r^{-2} \dkb^{\le 1} \Ga_g+ r^{-1} \dkb^{\le 1 } (\Ga_b\c \Ga_b).
\eeaa
Recalling $\dds_1\ddd_1=\ddd_2\dds_2+2K$ and $\widecheck{K}=r^{-1}\Ga_g$, we obtain
\beaa
2\ddd_2\dds_2\eta &=&  -\nab_3\nab\kac    -\frac{2}{r}\nab_3 \ze  -\frac{2}{r}\bb  + r^{-2} \dkb^{\le 1} \Ga_g+ r^{-1} \dkb^{\le 1 } (\Ga_b\c \Ga_b)
\eeaa
as stated.

 {\bf Step 3.}  We start with the equation
 \beaa
    \nab_3 \kabc -\frac{2\Up}{r}  \kabc &=& 2\div \xib +\frac{4\Up}{r} \ombc -\frac{2m}{r^2} \kabc  -\left(\frac{2}{r^2} -\frac{8m}{r^3} \right)\yc  +\Ga_b\c \Ga_b, 
 \eeaa
 which we write in the form
 \beaa
2\div \xib  &=&     \nab_3 \kabc  -\frac{4}{r} \ombc  +\frac{2}{r^2}\yc  +r^{-1}\Ga_g+\Ga_b\c \Ga_b. 
 \eeaa 
  We also have
 \beaa
 \curl \xib&=& \Ga_b\c \Ga_b.
 \eeaa
 Hence,
 \beaa
 2\ddd_1\xib&=& \left( \nab_3 \kabc  -\frac{4}{r} \ombc  +\frac{2}{r^2}\yc, 0 \right)+r^{-1}\Ga_g+\Ga_b\c \Ga_b
 \eeaa
 and thus
 \beaa
 2\dds_1\ddd_1\xib&=& -\nab\left( \nab_3 \kabc  -\frac{4}{r} \ombc  +\frac{2}{r^2}\yc\right)+r^{-2}\dkb^{\leq 1}\Ga_g+r^{-1}\dkb^{\leq 1}(\Ga_b\c \Ga_b)\\
 &=& -\nab_3\nab\kabc -[\nab,\nab_3]\kabc   +\frac{4}{r} \nab\ombc  -\frac{2}{r^2}\nab\yc +r^{-2}\dkb^{\leq 1}\Ga_g+r^{-1}\dkb^{\leq 1}(\Ga_b\c \Ga_b).
 \eeaa 
 Plugging the identity  $\nab   \widecheck{y} =  -\xib +\eta+\Ga_g +r\Ga_b\c\Ga_b$ derived in Step 2, we infer
 \beaa
 2\dds_1\ddd_1\xib &=& -\nab_3\nab\kabc -[\nab,\nab_3]\kabc   +\frac{4}{r} \nab\ombc  +\frac{2}{r^2}\xib -\frac{2}{r^2}\eta +r^{-2}\dkb^{\leq 1}\Ga_g+r^{-1}\dkb^{\leq 1}(\Ga_b\c \Ga_b).
 \eeaa 

 Next, as in Step 2, we have
 \beaa
 -[\nab, \nab_3] \kabc &=& [\nab_3, \nab] \kabc = (\eta-\ze)\nab_3\kabc  +\xib \nab_4 \kabc+ r^{-2} \dkb^{\le 1} \Ga_g+ r^{-1} \dkb^{\le 1 } (\Ga_b\c \Ga_b).
 \eeaa
 Using the equations for $\nab_3\kabc$ and $\nab_4\kabc$, we infer
 \beaa
 -[\nab, \nab_3] \kabc &=& r^{-2} \dkb^{\le 1} \Ga_g+ r^{-1} \dkb^{\le 1 } (\Ga_b\c \Ga_b).
 \eeaa
 Therefore
 \beaa
 2\dds_1\ddd_1\xib &=& -\nab_3\nab\kabc   +\frac{4}{r} \nab\ombc  +\frac{2}{r^2}\xib -\frac{2}{r^2}\eta +r^{-2}\dkb^{\leq 1}\Ga_g+r^{-1}\dkb^{\leq 1}(\Ga_b\c \Ga_b).
 \eeaa 
 Making  use of Step 1 to substitute $\nab\omb$ we deduce
 \beaa
 2\dds_1\ddd_1\xib &=& -\nab_3\nab\kabc   +\frac{2}{r}\left(\frac{1}{r}\xib  - \nab_3 \ze   -\bb +\frac{1}{r}\eta\right)  +\frac{2}{r^2}\xib -\frac{2}{r^2}\eta \\
 &&+r^{-2}\dkb^{\leq 1}\Ga_g+r^{-1}\dkb^{\leq 1}(\Ga_b\c \Ga_b)\\
 &=& -\nab_3\nab\kabc   -\frac{2}{r}\nab_3 \ze    -\frac{2}{r}\bb   +\frac{4}{r^2}\xib  +r^{-2}\dkb^{\leq 1}\Ga_g+r^{-1}\dkb^{\leq 1}(\Ga_b\c \Ga_b).
 \eeaa  
 Using as in Step 2  that $\dds_1\ddd_1=\ddd_2\dds_2+2K$ and $\widecheck{K}=r^{-1}\Ga_g$, we obtain
 \beaa
 2\ddd_2\dds_2\xib  &=& -\nab_3\nab\kabc   -\frac{2}{r}\nab_3 \ze    -\frac{2}{r}\bb     +r^{-2}\dkb^{\leq 1}\Ga_g+r^{-1}\dkb^{\leq 1}(\Ga_b\c \Ga_b).
 \eeaa  
 This ends the proof of  Proposition \ref{prop:additional.eqtsM4}.

%%%%%%%%%%%%%%%%%%%%%%%%%%%%%%%%%%%%%%%

\section{Proof of Proposition \ref{Prop:nu*ofGCM:0}}
\lab{sec:proofofProp:nu*ofGCM:0} 

%%%%%%%%%%%%%%%%%%%%%%%%%%%%%%%%%%%%%%%

We start with the following identity  of Proposition   \ref{prop:additional.eqtsM4}   
\beaa
2\ddd_2\dds_2\eta &=&  -\nab_3\nab\kac    -\frac{2}{r}\nab_3 \ze  -\frac{2}{r}\bb  + r^{-2} \dkb^{\le 1} \Ga_g+ r^{-1} \dkb^{\le 1 } (\Ga_b\c \Ga_b)
\eeaa
and apply $\dds_2\dds_1 \ddd_1$ to derive
\beaa
2\dds_2\dds_1 \ddd_1\ddd_2\dds_2\eta &=&  -\dds_2\dds_1 \ddd_1\nab_3\nab\kac    -\frac{2}{r}\dds_2\dds_1 \ddd_1\nab_3 \ze  -\frac{2}{r}\dds_2\dds_1 \ddd_1\bb  + r^{-5} \dkb^{\le 4} \Ga_g\\
&&+ r^{-4} \dkb^{\le 4} (\Ga_b\c \Ga_b)\\
&=&  -\dds_2\dds_1 \ddd_1\nab_3\nab\kac    -\frac{2}{r}\nab_3\dds_2\dds_1 \ddd_1\ze  -\frac{2}{r}\dds_2\dds_1 \ddd_1\bb  \\
&&     -\frac{2}{r}[\dds_2\dds_1 \ddd_1,\nab_3]\ze + r^{-5} \dkb^{\le 4} \Ga_g+ r^{-4} \dkb^{\le 4} (\Ga_b\c \Ga_b).
\eeaa
Making use of    the commutation  formulas of Lemma        \ref{Lemma:Commutation-Si_*}, we have
\beaa
[\dds_2\dds_1 \ddd_1,\nab_3]\ze &=& \dds_2\dds_1 [\ddd_1,\nab_3]\ze+\dds_2[\dds_1,\nab_3] \ddd_1\ze+[\dds_2, \nab_3]\dds_1 \ddd_1\ze\\
&=& \frac{3\Up}{r}\dds_2\dds_1 \ddd_1\ze+r^{-2}\dkb^{\leq 2}(\Ga_b\nab_3\ze)+r^{-3}\dkb^{\leq 2}(\Ga_b\dk^{\leq 1}\ze)\\
&=& r^{-2}\dkb^{\leq 2}(\Ga_b\nab_3\ze)+r^{-2}\dkb^{\leq 2}(\Ga_b\nab_4\ze)+r^{-4}\dkb^{\leq 3}\Ga_g.
\eeaa
Since the null structure equations for $\nab_3\ze$ and $\nab_4 \ze$ on $\Si_*$ imply that $\nab_4\ze=r^{-1}\Ga_g$ and $\nab_3\ze=r^{-1}\dkb^{\leq 1}\Ga_b$, we infer
\beaa
[\dds_2\dds_1 \ddd_1,\nab_3]\ze &=& r^{-4}\dkb^{\leq 3}\Ga_g+r^{-3}\dkb^{\leq 3}(\Ga_b\c\Ga_b),
\eeaa
and hence
\beaa
2\dds_2\dds_1 \ddd_1\ddd_2\dds_2\eta &=&  -\dds_2\dds_1 \ddd_1\nab_3\nab\kac    -\frac{2}{r}\nab_3\dds_2\dds_1 \ddd_1\ze  -\frac{2}{r}\dds_2\dds_1 \ddd_1\bb  \\
&&      + r^{-5} \dkb^{\le 4} \Ga_g+ r^{-4} \dkb^{\le 4} (\Ga_b\c \Ga_b).
\eeaa

Next, recalling that  
\beaa
\mu=-\div \ze  -\rho+\frac 1 2\chih \c \chibh, \qquad   \curl \ze= \rhod-\frac{1}{2}\chih\wedge\chibh,
 \eeaa
  we have
\beaa
\ddd_1\ze &=& \left(-\mu -\rho+\frac 1 2 \chih \c \chibh, \rhod-\frac{1}{2}\chih\wedge\chibh\right)
\eeaa
and hence
\beaa
\nab_3\dds_2\dds_1\ddd_1\ze &=& \nab_3\dds_2\dds_1\left(-\mu -\rho+\frac 1 2 \chih \c \chibh, \rhod-\frac{1}{2}\chih\wedge\chibh\right)\\
&=& -\nab_3\dds_2\dds_1\muc+\dds_2\dds_1\nab_3\left( -\rhoc+\frac 1 2 \chih \c \chibh, \rhod -\frac{1}{2}\chih\wedge\chibh\right)\\
&&+[\nab_3,\dds_2\dds_1]\left( -\rhoc+\frac 1 2 \chih \c \chibh, \rhod -\frac{1}{2}\chih\wedge\chibh\right).
\eeaa
From Proposition \ref{Prop.NullStr+Bianchi-lastslice}, we have on $\Si_*$
\beaa
&&\nab_3\chih\in r^{-1}\dkb^{\leq 1}\Ga_b, \quad  \nab_3\chibh\in \dkb^{\leq 1}\Ga_b, \quad \nab_3\rhoc\in r^{-2}\dkb^{\leq 1}\Ga_b, \quad \nab_3\rhod\in r^{-2}\dkb^{\leq 1}\Ga_b,\\
&& \nab_4\chih\in r^{-1}\dkb^{\leq 1}\Ga_g, \quad  \nab_4\chibh\in r^{-1}\dkb^{\leq 1}\Ga_b, \quad \nab_4\rhoc\in r^{-2}\dkb^{\leq 1}\Ga_g, \quad \nab_4\rhod\in r^{-2}\dkb^{\leq 1}\Ga_g,
\eeaa
which, together with the commutation  formulas of Lemma        \ref{Lemma:Commutation-Si_*}, implies
\beaa
\nab_3\dds_2\dds_1\ddd_1\ze &=& -\nab_3\dds_2\dds_1\muc+\dds_2\dds_1\left( -\nab_3\rhoc+\frac 1 2 \chih \c\nab_3 \chibh, \nab_3\rhod -\frac{1}{2}\chih\wedge\nab_3\chibh\right)\\
&& + r^{-4} \dkb^{\le 2} \Ga_g+ r^{-3} \dkb^{\le 3} (\Ga_b\c \Ga_b).
\eeaa
Also, we have in view of Proposition \ref{Prop.NullStr+Bianchi-lastslice}
\beaa
-\nab_3\rhoc+\frac 1 2 \chih \c\nab_3 \chibh &=& -\left(-\div\bb -\frac{1}{2}\chih\c\aa\right) +\frac 1 2 \chih \c(-\aa)+ r^{-2}\Ga_g+ r^{-1} \dkb^{\le 1} (\Ga_b\c \Ga_b)\\
&=& \div\bb + r^{-2}\Ga_g+ r^{-1} \dkb^{\le 1} (\Ga_b\c \Ga_b),\\
 \nab_3\rhod -\frac{1}{2}\chih\wedge\nab_3\chibh &=& -\curl\bb -\frac{1}{2}\chih\wedge\aa -\frac{1}{2}\chih\wedge(-\aa)+ r^{-2}\Ga_g+ r^{-1} \dkb^{\le 1} (\Ga_b\c \Ga_b)\\
&=& -\curl\bb + r^{-2}\Ga_g+ r^{-1} \dkb^{\le 1} (\Ga_b\c \Ga_b)
\eeaa
so that 
\beaa
\nab_3\dds_2\dds_1\ddd_1\ze &=& -\nab_3\dds_2\dds_1\muc+\dds_2\dds_1\left( \div\bb, -\curl\bb\right) + r^{-4} \dkb^{\le 2} \Ga_g+ r^{-3} \dkb^{\le 3} (\Ga_b\c \Ga_b).
\eeaa
We deduce
\beaa
2\dds_2\dds_1 \ddd_1\ddd_2\dds_2\eta &=&  -\dds_2\dds_1 \ddd_1\nab_3\nab\kac    +\frac{2}{r}\nab_3\dds_2\dds_1\muc -\frac{2}{r}\dds_2\dds_1\left( \div\bb, -\curl\bb\right) -\frac{2}{r}\dds_2\dds_1 \ddd_1\bb  \\
&&      + r^{-5} \dkb^{\le 4} \Ga_g+ r^{-4} \dkb^{\le 4} (\Ga_b\c \Ga_b)
\eeaa
and hence
\beaa
2\dds_2\dds_1 \ddd_1\ddd_2\dds_2\eta &=&  -\dds_2\dds_1 \ddd_1\nab_3\nab\kac    +\frac{2}{r}\nab_3\dds_2\dds_1\muc -\frac{4}{r}\dds_2\dds_1\div\bb        \\
&& + r^{-5} \dkb^{\le 4} \Ga_g+ r^{-4} \dkb^{\le 4} (\Ga_b\c \Ga_b)
\eeaa
 as desired.

Next, we focus on the second identity. We consider  the following identity  of Proposition   \ref{prop:additional.eqtsM4}
\beaa
 2\ddd_2\dds_2\xib  &=& -\nab_3\nab\kabc   -\frac{2}{r}\nab_3 \ze    -\frac{2}{r}\bb     +r^{-2}\dkb^{\leq 1}\Ga_g+r^{-1}\dkb^{\leq 1}(\Ga_b\c \Ga_b)
\eeaa
and apply $\dds_2\dds_1 \ddd_1$ to derive
\beaa
 2\dds_2\dds_1 \ddd_1\ddd_2\dds_2\xib  &=& -\dds_2\dds_1 \ddd_1\nab_3\nab\kabc   -\frac{2}{r}\dds_2\dds_1 \ddd_1\nab_3 \ze    -\frac{2}{r}\dds_2\dds_1 \ddd_1\bb   + r^{-5} \dkb^{\le 4} \Ga_g\\
&&+ r^{-4} \dkb^{\le 4} (\Ga_b\c \Ga_b)\\
&=&  -\nab_3\dds_2\dds_1 \ddd_1\nab\kabc   -\frac{2}{r}\nab_3\dds_2\dds_1 \ddd_1 \ze    -\frac{2}{r}\dds_2\dds_1 \ddd_1\bb  \\
&& -[\dds_2\dds_1 \ddd_1,\nab_3]\nab\kabc -\frac{2}{r}[\dds_2\dds_1 \ddd_1,\nab_3 \ze] + r^{-5} \dkb^{\le 4} \Ga_g+ r^{-4} \dkb^{\le 4} (\Ga_b\c \Ga_b).
 \eeaa
 Recall from above that we have
\beaa
[\dds_2\dds_1 \ddd_1,\nab_3]\ze &=& r^{-4}\dkb^{\leq 3}\Ga_g+r^{-3}\dkb^{\leq 3}(\Ga_b\c\Ga_b),
\eeaa
so that 
\beaa
 2\dds_2\dds_1 \ddd_1\ddd_2\dds_2\xib  &=&  -\nab_3\dds_2\dds_1 \ddd_1\nab\kabc   -\frac{2}{r}\nab_3\dds_2\dds_1 \ddd_1 \ze    -\frac{2}{r}\dds_2\dds_1 \ddd_1\bb  \\
&& -[\dds_2\dds_1 \ddd_1,\nab_3]\nab\kabc  + r^{-5} \dkb^{\le 4} \Ga_g+ r^{-4} \dkb^{\le 4} (\Ga_b\c \Ga_b).
 \eeaa
Also, since the null structure equations for $\nab_3\kabc$ and $\nab_4 \kabc$ on $\Si_*$ imply  $\nab_4\kabc=r^{-1}\dkb^{\leq 1}\Ga_g$ and $\nab_3\kabc=r^{-1}\dkb^{\leq 1}\Ga_b$, and making use of    the commutation  formulas of Lemma        \ref{Lemma:Commutation-Si_*}, we infer, arguing as for the commutator $[\dds_2\dds_1 \ddd_1,\nab_3]\ze$, 
\beaa
[\dds_2\dds_1 \ddd_1,\nab_3]\kabc &=& r^{-4}\dkb^{\leq 3}\Ga_g+r^{-3}\dkb^{\leq 3}(\Ga_b\c\Ga_b),
\eeaa
and hence 
\beaa
 2\dds_2\dds_1 \ddd_1\ddd_2\dds_2\xib  &=&  -\nab_3\dds_2\dds_1 \ddd_1\nab\kabc   -\frac{2}{r}\nab_3\dds_2\dds_1 \ddd_1 \ze    -\frac{2}{r}\dds_2\dds_1 \ddd_1\bb  \\
&&   + r^{-5} \dkb^{\le 4} \Ga_g+ r^{-4} \dkb^{\le 4} (\Ga_b\c \Ga_b).
 \eeaa
Now, recall from above that we have
\beaa
\nab_3\dds_2\dds_1\ddd_1\ze &=& -\nab_3\dds_2\dds_1\muc+\dds_2\dds_1\left( \div\bb, -\curl\bb\right) + r^{-4} \dkb^{\le 2} \Ga_g+ r^{-3} \dkb^{\le 3} (\Ga_b\c \Ga_b).
\eeaa
We deduce
\beaa
 2\dds_2\dds_1 \ddd_1\ddd_2\dds_2\xib  &=&  -\nab_3\dds_2\dds_1 \ddd_1\nab\kabc   +\frac{2}{r}\nab_3\dds_2\dds_1\muc   -\frac{4}{r}\dds_2\dds_1 \div\bb  \\
&&   + r^{-5} \dkb^{\le 4} \Ga_g+ r^{-4} \dkb^{\le 4} (\Ga_b\c \Ga_b).
 \eeaa
Also, we have
\beaa
-\dds_2\dds_1 \ddd_1\nab\kabc &=& \dds_2\dds_1 \ddd_1\dds_1(\kabc,0)= \dds_2(\ddd_2 \dds_2+2K)\dds_1(\kabc,0)\\
&=& \dds_2\left(\ddd_2 \dds_2+\frac{2}{r^2}\right)\dds_1(\kabc,0)+2\dds_2\Big(\widecheck{K}\dds_1(\kabc,0)\Big)\\
&=& \left(\ddd_2 \dds_2+\frac{2}{r^2}\right)\dds_2\dds_1(\kabc,0)+2\dds_2\Big(\widecheck{K}\dds_1(\kabc,0)\Big)
\eeaa
and hence
\beaa
 2\dds_2\dds_1 \ddd_1\ddd_2\dds_2\xib  &=&  \nab_3\left(\ddd_2 \dds_2+\frac{2}{r^2}\right)\dds_2\dds_1\kabc   +\frac{2}{r}\nab_3\dds_2\dds_1\muc   -\frac{4}{r}\dds_2\dds_1 \div\bb  \\
 &&+\nab_3\left(2\dds_2\Big(\widecheck{K}\dds_1(\kabc,0)\Big)\right)   + r^{-5} \dkb^{\le 4} \Ga_g+ r^{-4} \dkb^{\le 4} (\Ga_b\c \Ga_b).
 \eeaa
Now, from Proposition \ref{Prop.NullStr+Bianchi-lastslice}, we have $\nab_4\kabc=r^{-1}\dkb^{\leq 1}\Ga_g$, $\nab_3\kabc=r^{-1}\dkb^{\leq 1}\Ga_b$, $\nab_4\widecheck{K}=r^{-2}\dkb^{\leq 1}\Ga_g$ and $\nab_3\widecheck{K}=r^{-2}\dkb^{\leq 1}\Ga_b$ on $\Si_*$, which, together with the commutation  formulas of Lemma        \ref{Lemma:Commutation-Si_*}, implies
\beaa
\nab_3\left( 2\dds_2\Big(\widecheck{K}\dds_1(\kabc,0)\Big)\right) &=&  r^{-4} \dkb^{\le 3} (\Ga_b\c \Ga_b)
\eeaa
and hence 
\beaa
 2\dds_2\dds_1 \ddd_1\ddd_2\dds_2\xib  &=&  \nab_3\left(\ddd_2 \dds_2+\frac{2}{r^2}\right)\dds_2\dds_1\kabc   +\frac{2}{r}\nab_3\dds_2\dds_1\muc   -\frac{4}{r}\dds_2\dds_1 \div\bb  \\
&&   + r^{-5} \dkb^{\le 4} \Ga_g+ r^{-4} \dkb^{\le 4} (\Ga_b\c \Ga_b)
 \eeaa
as stated. This concludes the proof of Proposition \ref{Prop:nu*ofGCM:0}.

%%%%%%%%%%%%%%%%%%%%%%%%%%%%%%%%%%%%%%%

\section{Proof of Lemma \ref{Lemma:transport.alongSi_*1}}
\lab{sec:proofofLemma:transport.alongSi_*1} 

%%%%%%%%%%%%%%%%%%%%%%%%%%%%%%%%%%%%%%%

Below, recall that the notation $O(r^a)$, for $a\in\mathbb{R}$, denotes an explicit function of $r$ which is bounded 
by $r^a$ as $r\to+\infty$.

{\bf First identity.} We start with the following  equation, see Proposition \ref{Prop.NullStr+Bianchi-lastslice},
\beaa
   \nab_3 \kabc -\frac{2\Up}{r}  \kabc &=& 2\div \xib +\frac{4\Up}{r} \ombc -\frac{2m}{r^2} \kabc  -\left(\frac{2}{r^2} -\frac{8m}{r^3} \right)\yc  +\Ga_b\c \Ga_b.
  \eeaa
Commuting the first equation with the Laplacian we deduce
\beaa
  \nab_3 \lap\kabc -\frac{2\Up}{r}  \lap\kabc &=& [\nab_3,\lap]\kabc+2\lap\div \xib +\frac{4\Up}{r} \lap\ombc -\frac{2m}{r^2} \lap\kabc  -\left(\frac{2}{r^2} -\frac{8m}{r^3} \right)\lap\yc\\
  &&  + r^{-2}\dkb^{\le 2 } \big(\Ga_b\c \Ga_b).
\eeaa
According to Lemma \ref{Lemma:Commutation-Si_*}, we have
\beaa
\,[ \nab_3, \De]\kabc &=& \frac{2\Up}{r} \lap\kabc+ r^{-1}\dkb^{\le 1}\Big( \Ga_b \c\nab_3\kabc + r^{-1} \Ga_b \c \dk\kabc \Big)
\eeaa
and, using again the equation for $\nab_3\kab$,
\beaa
\, [\nab_3, \lap]\kabc= \frac{2\Up}{r} \lap\kabc + r^{-2} \dkb^{\le 2 }(\Ga_b\c\Ga_b). 
\eeaa
Thus
\bea
\lab{eq:nab3kabc+}
\bsplit
 \nab_3 \lap\kabc -\frac{4\Up}{r}  \lap\kabc &= 2\lap\div \xib +\frac{4\Up}{r} \lap\ombc -\frac{2m}{r^2} \lap\kabc  -\left(\frac{2}{r^2} -\frac{8m}{r^3} \right)\lap\yc\\
  &  + r^{-2}\dkb^{\le 2 } \big(\Ga_b\c \Ga_b).
\end{split}
\eea
To remove  the term involving $\lap \omb$, we consider  the null structure equation for $\ze$ in the form
\beaa
 \nab_3 \ze -\frac{\Up}{r} \ze &=& -\bb-2 \nab \ombc +\frac{\Up}{r} (\eta+\ze)+\frac{1}{r}\xib+ \frac{2m}{r^2}(\ze-\eta) +\Ga_b\c \Ga_b.
 \eeaa
Differentiating w.r.t. $\div$,   and using the commutation 
\beaa
 \,[ \nab_3, \div] \ze&=& \frac{\Up}{r}\div \ze +\Ga_b \c\nab_3 \ze + r^{-1} \Ga_b \c \dk \ze= \frac{\Up}{r}\div \ze + r^{-1}\dkb^{\le 1} (\Ga_b\c \Ga_b),
\eeaa
we obtain
\beaa
 \nab_3\div \ze - \frac{2\Up}{r} \div\ze &=& -\div\bb-2 \lap\ombc +\frac{\Up}{r}\div(\eta+\ze)+\frac{1}{r}\div\xib+ \frac{2m}{r^2}\div(\ze-\eta) \\
 &&+r^{-1}\dkb^{\le 1} (\Ga_b\c \Ga_b)
 \eeaa
from which we easily deduce
\beaa
 \nab_3\left(\frac{2\Up}{r}\div \ze\right) &=& \frac{2\Up}{r}\nab_3\div \ze +\left(-\frac{2e_3(r)}{r^2}+\frac{8me_3(r)}{r^3}\right)\div\ze\\
 &=&  \frac{4\Up^2}{r^2} \div\ze  -\frac{2\Up}{r}\div\bb -\frac{4\Up}{r}\lap\ombc +\frac{2\Up^2}{r^2}\div(\eta+\ze)+\frac{2\Up}{r^2}\div\xib\\
 &&+ \frac{4m\Up}{r^3}\div(\ze-\eta) +\left(\frac{2\Up}{r^2}-\frac{8m\Up}{r^3}\right)\div\ze+r^{-2}\dkb^{\le 1} (\Ga_b\c \Ga_b)
 \eeaa
and hence
\beaa
\frac{4\Up}{r}\lap\ombc   &=& -\nab_3\left(\frac{2\Up}{r}\div \ze\right) +\frac{4\Up^2}{r^2} \div\ze  -\frac{2\Up}{r}\div\bb +\frac{2\Up^2}{r^2}\div(\eta+\ze)+\frac{2\Up}{r^2}\div\xib\\
 &&+ \frac{4m\Up}{r^3}\div(\ze-\eta) +\left(\frac{2\Up}{r^2}-\frac{8m\Up}{r^3}\right)\div\ze+r^{-2}\dkb^{\le 1} (\Ga_b\c \Ga_b)
 \eeaa
Combining this  with the previous identity \eqref{eq:nab3kabc+}, we deduce
\beaa
\bsplit
 \nab_3\left( \lap\kabc+\frac{2\Up}{r}\div \ze\right)  &= O(r^{-1})\lap\kabc+2\lap\div \xib  +O(r^{-2})\lap\yc  +O(r^{-2})\div\ze  +O(r^{-1})\div\bb\\
 & +O(r^{-2})\div\eta+O(r^{-2})\div\xib  + r^{-2}\dkb^{\le 2 } \big(\Ga_b\c \Ga_b).
\end{split}
\eeaa

Similarly, proceeding as above, starting with the following  equations, see Proposition \ref{Prop.NullStr+Bianchi-lastslice},
\beaa
  \nab_4\kabc+\frac{1}{r}\kabc&=& - 2 \div \ze + 2 \rhoc+\Ga_b\c \Ga_g,
  \eeaa
we deduce
\beaa
  \nab_4\lap\kabc+\frac{3}{r}\lap\kabc&=& - 2 \lap\div \ze + 2 \lap\rhoc+ r^{-2}\dkb^{\le 2 }(\Ga_b\c \Ga_g).
  \eeaa
Also, using
\beaa
 \nab_4 \ze+\frac{2}{r} \ze &=-\b+\Ga_g\c \Ga_g,
\eeaa
and differentiating w.r.t. $\div$, we have
\beaa
 \nab_4\div\ze+\frac{3}{r}\div\ze &=-\div\b+r^{-1} \dkb^{\le 1 } \big(\Ga_g\c \Ga_g\big)
\eeaa
and hence
\beaa
\nab_4\left( \frac{2\Up}{r}\div \ze\right)  &=& \frac{2\Up}{r}\nab_4\div\ze+\left(-\frac{2}{r^2}+\frac{8m}{r^3}\right)\div\ze\\
&=& -\frac{2\Up}{r}\div\b  -\frac{4}{r^2}\left(2-\frac{5m}{r}\right)\div\ze +r^{-2} \dkb^{\le 1 } \big(\Ga_g\c \Ga_g\big).
\eeaa
Combining, we obtain 
\beaa
\nab_4\left( \lap\kabc+\frac{2\Up}{r}\div \ze\right)  &=& O(r^{-1})\lap\kabc  - 2 \lap\div \ze + 2 \lap\rhoc +O(r^{-1})\div\b  +O(r^{-2})\div\ze\\
&& +r^{-2}\dkb^{\le 2 }(\Ga_b\c \Ga_g)
\eeaa
and, since $b_*=-1-\frac{2m}{r}+r\Ga_b$, 
\beaa
b_*\nab_4\left( \lap\kabc+\frac{2\Up}{r}\div \ze\right)  &=& O(r^{-1})\lap\kabc  +2\left(1+O(r^{-1})\right)\lap\div \ze - 2\left(1+O(r^{-1})\right)\lap\rhoc\\
&& +O(r^{-1})\div\b  +O(r^{-2})\div\ze +r^{-2}\dkb^{\le 2 }(\Ga_b\c \Ga_g).
\eeaa
Since
\beaa
\nab_\nu\left( \lap\kabc+\frac{2\Up}{r}\div \ze\right)  &=& \nab_3\left( \lap\kabc+\frac{2\Up}{r}\div \ze\right)+ b_*  \nab_4\left( \lap\kabc+\frac{2\Up}{r}\div \ze\right),
\eeaa
we infer from the above
\beaa
\nab_\nu\left( \lap\kabc+\frac{2\Up}{r}\div \ze\right)  &=& O(r^{-1})\lap\kabc+2\lap\div \xib  +O(r^{-2})\lap\yc  +O(r^{-2})\div\ze  +O(r^{-1})\div\bb\\
 && +O(r^{-2})\div\eta+O(r^{-2})\div\xib   +2\left(1+O(r^{-1})\right)\lap\div \ze \\
 && - 2\left(1+O(r^{-1})\right)\lap\rhoc +O(r^{-1})\div\b  +r^{-2}\dkb^{\le 2 }(\Ga_b\c \Ga_b)
\eeaa
as stated.

 {\bf Identities 2 and 3.} We have, see Proposition \ref{Prop.NullStr+Bianchi-lastslice},
 \beaa
      \nab_3\b -\frac{2\Up}{r}\b &=& (\nab\rho+\dual\nab\rhod) +\frac{2m}{r^2}\b -\frac{6m}{r^3}\eta+r^{-1}\Ga_b\c\Ga_g.
 \eeaa
Taking the divergence and the curl, we infer
\beaa
      \nab_3\div\b -\frac{2\Up}{r}\div\b &=&  [\nab_3,\div]\b+\lap\rho+\frac{2m}{r^2}\div\b -\frac{6m}{r^3}\div\eta+r^{-2}\dkb^{\leq 1}(\Ga_b\c\Ga_g),\\
    \nab_3\curl\b -\frac{2\Up}{r}\curl\b &=&  [\nab_3,\curl]\b-\lap\rhod+\frac{2m}{r^2}\curl\b -\frac{6m}{r^3}\curl\eta+r^{-2}\dkb^{\leq 1}(\Ga_b\c\Ga_g).
 \eeaa
According to Lemma \ref{Lemma:Commutation-Si_*}, using  the equation for $\nab_3\b$ and $\nab_4\b$, 
\beaa
 \,[\nab_3,\div]\b &=& \frac{\Up}{r} \div \b +\Ga_b\c  \nab_3\b +r^{-2}\dkb^{\le 1} ( \Ga_g \c \Ga_b)  \\
 &    =&    \frac{\Up}{r} \div \b      +r^{-2}\dkb^{\le 1} ( \Ga_g \c \Ga_b)  \\
\,  [\nab_3,\curl ] \b &=& \frac{\Up}{r} \curl \b +\Ga_b\c \nab_3\b +r^{-2}\dkb^{\le 1}( \Ga_g \c \Ga_b)
\\
&=& \frac{\Up}{r} \curl \b    +r^{-2}\dkb^{\le 1} ( \Ga_g \c \Ga_b).
\eeaa 
Hence
\beaa
      \nab_3\div\b  &=&  \left(\frac{3}{r}-\frac{4m}{r^2}  \right)\div \b+\lap\rho -\frac{6m}{r^3}\div\eta+r^{-2}\dkb^{\leq 1}(\Ga_b\c\Ga_g),\\
    \nab_3\curl\b  &=&  \left(\frac{3}{r}-\frac{4m}{r^2}  \right)\curl \b -\lap\rhod  -\frac{6m}{r^3}\curl\eta+r^{-2}\dkb^{\leq 1}(\Ga_b\c\Ga_g).
 \eeaa
Since  $\curl\eta = \rhod+\Ga_b\c \Ga_g$, we infer
\beaa
      \nab_3\div\b  &=&  \left(\frac{3}{r}-\frac{4m}{r^2}  \right)\div \b+\lap\rho -\frac{6m}{r^3}\div\eta+r^{-2}\dkb^{\leq 1}(\Ga_b\c\Ga_g),\\
    \nab_3\curl\b  &=&  \left(\frac{3}{r}-\frac{4m}{r^2}  \right)\curl \b -\lap\rhod  -\frac{6m}{r^3}\rhod+r^{-2}\dkb^{\leq 1}(\Ga_b\c\Ga_g).
 \eeaa 

Also, we have
\beaa
   \nab_4\beta +\frac{4}{r}\beta &=& -\div\a  +r^{-1}\Ga_g\c\Ga_g.
\eeaa
Taking the divergence and the curl, we infer
\beaa
   \nab_4\div\beta +\frac{4}{r}\div\beta &=& [\nab_4,\div]\beta  -\div\div\a  +r^{-2}\dkb^{\leq 1}(\Ga_g\c\Ga_g),\\
     \nab_4\curl\beta +\frac{4}{r}\curl\beta &=& [\nab_4,\curl]\beta  -\curl\div\a  +r^{-2}\dkb^{\leq 1}(\Ga_g\c\Ga_g).
\eeaa
According to Lemma \ref{Lemma:Commutation-Si_*}, using again the equation for $\nab_3\b$ and $\nab_4\b$, 
\beaa
 \,[\nab_4,\div]\b  &    =&    -\frac{1}{r} \div \b      +r^{-2}\dkb^{\le 1} ( \Ga_g \c \Ga_g),  \\
\,  [\nab_4,\curl ] \b &=& -\frac{1}{r} \curl \b     +r^{-2}\dkb^{\le 1} ( \Ga_g \c \Ga_g).
\eeaa 
Hence
\beaa
   \nab_4\div\beta  &=& -\frac{5}{r}\div\beta  -\div\div\a  +r^{-2}\dkb^{\leq 1}(\Ga_g\c\Ga_g),\\
     \nab_4\curl\beta  &=& -\frac{5}{r}\curl\beta  -\curl\div\a  +r^{-2}\dkb^{\leq 1}(\Ga_g\c\Ga_g).
\eeaa
Since $\nu=e_3+b_*e_4$ and $b_*=-1-\frac{2m}{r}+r\Ga_b$, we infer
\beaa
      \nab_\nu\div\b  &=&  O(r^{-1})\div \b+\lap\rho +(1+O(r^{-1}))\div\div\a +O(r^{-3})\div\eta+r^{-2}\dkb^{\leq 1}(\Ga_b\c\Ga_g),\\
    \nab_\nu\curl\b  &=& \frac{8}{r}(1+ O(r^{-1}))\curl \b -\lap\rhod +(1+O(r^{-1}))\curl\div\a +O(r^{-3})\rhod\\
    &&+r^{-2}\dkb^{\leq 1}(\Ga_b\c\Ga_g).
 \eeaa

{\bf Identity 4.} We have, see Proposition \ref{Prop.NullStr+Bianchi-lastslice},
\beaa
         \nab_3 \rhoc  -\frac{3\Up}{r} \rhoc  &=&  -\div\bb  +\frac{3m}{r^3}\kabc -\frac{6m}{r^4} \yc   -\frac{1}{2}\chih\c\aa+r^{-1}\Ga_b\c\Ga_b.
\eeaa
Also, using the equation for $\chih$ and $\chibh$, we have
\beaa
\nab_3(\chih\c\chibh) &=& \chih\c\left( \frac{2\Up}{r}\chibh -\aa -\frac{2m}{r^2}\chibh +\nab\hot \xib +\Ga_b\c \Ga_b\right)\\
&&+\left(\frac{\Up}{r} \chih+\nab\hot \eta -\frac{1}{r} \chibh+\frac{2m}{r^2}\chih+\Ga_b\c \Ga_b\right)\c\chibh\\
&=&-\chih \c\aa+  r^{-1}\dkb^{\le  1}( \Ga_b \c \Ga_b).
\eeaa
Hence, we deduce 
\beaa
   \nab_3\left(\rhoc - \frac{1}{2}\chih\c\chibh \right)    &=&  -\div\bb+\frac{3\Up}{r} \rhoc   +\frac{3m}{r^3}\kabc -\frac{6m}{r^4} \yc   +r^{-1}\dkb^{\le  1}( \Ga_b \c \Ga_b).
\eeaa
Also, we have
\beaa
     \nab_4 \rhoc+\frac{3}{r}\rhoc &=&\div \b +r^{-1}\Ga_b\c\Ga_g
\eeaa 
and
\beaa
\nab_4(\chih\c\chibh) &=& r^{-1}\dkb^{\le  1}( \Ga_b \c \Ga_g),
\eeaa
and hence
\beaa
   \nab_4\left(\rhoc - \frac{1}{2}\chih\c\chibh \right)    &=& \div \b -\frac{3}{r}\rhoc  +r^{-1}\dkb^{\le  1}( \Ga_b \c \Ga_g).
\eeaa
Since $\nu=e_3+b_*e_4$ and $b_*=-1-\frac{2m}{r}+r\Ga_b$, we infer
\beaa
   \nab_\nu\left(\rhoc - \frac{1}{2}\chih\c\chibh \right)    &=&  -\div\bb -(1+O(r^{-1}))\div\b +O(r^{-1})\rhoc   +O(r^{-3})\kabc +O(r^{-4}) \yc \\
   &&  +r^{-1}\dkb^{\le  1}( \Ga_b \c \Ga_b)
\eeaa
as stated. This concludes the proof of Lemma \ref{Lemma:transport.alongSi_*1}.

%%%%%%%%%%%%%%%%%%%%%%%%%%%%%%%%%%%

\section{Proof of Corollary \ref{corofLemma:transport.alongSi_*1}}
\lab{sec:corofLemma:transport.alongSi_*1}

%%%%%%%%%%%%%%%%%%%%%%%%%%%%%%%%%%%

Since $\nu(\Jp)=0$, we have in view of Corollary \ref{Corr:nuSof integrals}
for any scalar function $h$ on $\Si_*$ and any $S\subset\Si_*$
 \bea\lab{eq:transportequationforell=1modealongSigmastargeneral}
\nu\left(\int_Sh\Jp\right) = \int_S\nu (h)\Jp  -\frac{4}{r}\int_Sh\Jp +r^3\Ga_b\nu(h)+r^2\Ga_b h
\eea
where we also used $\Jp=O(1)$, and where we recall that  the notation $O(r^a)$, for $a\in\mathbb{R}$, denotes an explicit function of $r$ which is bounded by $r^a$ as $r\to+\infty$. 

Next, recall from Lemma \ref{Lemma:transport.alongSi_*1} that we have along  $\Si_*$
\beaa
\nab_\nu\left( \lap\kabc+\frac{2\Up}{r}\div \ze\right)  &=& O(r^{-1})\lap\kabc+2\lap\div \xib  +O(r^{-2})\lap\yc  +O(r^{-2})\div\ze \\
&& +O(r^{-1})\div\bb +O(r^{-2})\div\eta+O(r^{-2})\div\xib   \\
&& +2\left(1+O(r^{-1})\right)\lap\div \ze  - 2\left(1+O(r^{-1})\right)\lap\rhoc\\
&& +O(r^{-1})\div\b  +r^{-2}\dkb^{\le 2 }(\Ga_b\c \Ga_b).
\eeaa
Together with \eqref{eq:transportequationforell=1modealongSigmastargeneral}, and noticing that the terms $O(r^a)$ only depend on $r$ and are thus constant on $S$, we infer
\beaa
&&\nu\left(\int_S\left( \lap\kabc+\frac{2\Up}{r}\div \ze\right)\Jp\right)\\
 &=& O(r^{-1})\int_S\lap\kabc\Jp+2\int_S\lap\div\xib\Jp  +O(r^{-2})\int_S\lap\yc\Jp  +O(r^{-2})\int_S\div\ze\Jp \\
&& +O(r^{-1})\int_S\div\bb\Jp +O(r^{-2})\int_S\div\eta\Jp+O(r^{-2})\int_S\div\xib\Jp   \\
&& +2\left(1+O(r^{-1})\right)\int_S\lap\div\ze\Jp  - 2\left(1+O(r^{-1})\right)\int_S\lap\rhoc\Jp +O(r^{-1})\int_S\div\b\Jp\\
&&  +\dkb^{\le 2 }(\Ga_b\c \Ga_b).
\eeaa
Integrating by parts, we infer
\beaa
&&\nu\left(\int_S\left( \lap\kabc+\frac{2\Up}{r}\div \ze\right)\Jp\right)\\
 &=& O(r^{-3})\int_S\kabc\Jp+O(r^{-2})\int_S\div\xib\Jp  +O(r^{-2})\int_S\Delta\yc\Jp  +O(r^{-2})\int_S\div\ze\Jp \\
&& +O(r^{-1})\int_S\div\bb\Jp +O(r^{-2})\int_S\div\eta\Jp  +O(r^{-2})\int_S\rhoc\Jp +O(r^{-1})\int_S\div\b\Jp \\
&&+r\left|\left(\Delta+\frac{2}{r^2}\right)\Jp\right|\dkb^{\leq 1}\Ga_b+\dkb^{\le 2 }(\Ga_b\c \Ga_b).
\eeaa
Also, in view of Lemma \ref{Lemma:eqts-nabeta,xib}, we have
\beaa
\Delta\yc &=& \div(-\xib+\big(\ze-\eta) y)=-\div\xib+y(\div\ze-\div\eta)+(\ze-\eta)\c\nab y,
\eeaa
and hence, since $y=-\Up+r\Ga_b$, 
\beaa
\Delta\yc &=& -\div\xib -\Up(\div\ze-\div\eta)+\dkb(\Ga_b\c\Ga_b)
\eeaa
so that 
\beaa
&&\nu\left(\int_S\left( \lap\kabc+\frac{2\Up}{r}\div \ze\right)\Jp\right)\\
 &=& O(r^{-3})\int_S\kabc\Jp+O(r^{-2})\int_S\div\xib\Jp   +O(r^{-2})\int_S\div\ze\Jp  +O(r^{-1})\int_S\div\bb\Jp \\
 &&+O(r^{-2})\int_S\div\eta\Jp    +O(r^{-2})\int_S\rhoc\Jp +O(r^{-1})\int_S\div\b\Jp\\
&& +r\left|\left(\Delta+\frac{2}{r^2}\right)\Jp\right|\dkb^{\leq 1}\Ga_b+\dkb^{\le 2 }(\Ga_b\c \Ga_b).
\eeaa
Together with the GCM conditions $(\div\eta)_{\ell=1}=0$ and $(\div\xib)_{\ell=1}=0$, we infer
\beaa
&&\nu\left(\int_S\left( \lap\kabc+\frac{2\Up}{r}\div \ze\right)\Jp\right)\\
 &=& O(r^{-3})\int_S\kabc\Jp   +O(r^{-2})\int_S\div\ze\Jp +O(r^{-1})\int_S\div\bb\Jp     +O(r^{-2})\int_S\rhoc\Jp \\
 &&+O(r^{-1})\int_S\div\b\Jp +r\left|\left(\Delta+\frac{2}{r^2}\right)\Jp\right|\dkb^{\leq 1}\Ga_b+\dkb^{\le 2 }(\Ga_b\c \Ga_b)
\eeaa
as stated.

Next, recall from Lemma \ref{Lemma:transport.alongSi_*1} that we have along  $\Si_*$
\beaa
    \nab_\nu\div\b  &=&  O(r^{-1})\div \b+\lap\rho +(1+O(r^{-1}))\div\div\a\\
         && +O(r^{-3})\div\eta+r^{-2}\dkb^{\leq 1}(\Ga_b\c\Ga_g).
\eeaa
Together with \eqref{eq:transportequationforell=1modealongSigmastargeneral}, and noticing that the terms $O(r^a)$ only depend on $r$ and are thus constant on $S$, we infer
\beaa
\nu\left(\int_S\div\b\Jp\right) &=& O(r^{-1})\int_S\div\b\Jp +\int_S\lap\rhoc \Jp +(1+O(r^{-1}))\int_S\div\div\a \Jp\\
&&+O(r^{-3})\int_S\div\eta \Jp+\dkb^{\leq 1}(\Ga_b\c\Ga_g)
\eeaa
Together with the GCM condition $(\div\eta)_{\ell=1}=0$, we infer
\beaa
\nu\left(\int_S\div\b\Jp\right) &=& O(r^{-1})\int_S\div\b\Jp +\int_S\lap\rhoc \Jp +(1+O(r^{-1}))\int_S\div\div\a \Jp\\
&&+\dkb^{\leq 1}(\Ga_b\c\Ga_g).
\eeaa
Integrating by parts, we deduce the desired identity for $\div\b$
\beaa
\nu\left(\int_S\div\b\Jp\right) &=& O(r^{-1})\int_S\div\b\Jp+O(r^{-2})\int_S\rhoc\Jp\\
&&+r\left(\left|\left(\Delta+\frac{2}{r^2}\right)\Jp\right|+\left|\dds_2\dds_1\Jp\right|\right)\Ga_g+\dkb^{\leq 1}(\Ga_b\c\Ga_g),
\eeaa
where by $\dds_1\Jp$, we mean $\dds_1(\Jp,0)$. Similarly, starting from 
\beaa
 \nab_\nu\curl\b  &=&  \frac{8}{r}(1+ O(r^{-1}))\curl \b -\lap\rhod +(1+O(r^{-1}))\curl\div\a \\
    &&+O(r^{-3})\rhod+r^{-2}\dkb^{\leq 1}(\Ga_b\c\Ga_g),
\eeaa
 we deduce the desired identity for $\curl\b$
\beaa
\nu\left(\int_S\curl\b\Jp\right) &=&  \frac{4}{r}(1+ O(r^{-1}))\int_S\curl\b\Jp+\frac{2}{r^2}(1+ O(r^{-1}))\int_S\rhod\Jp\\
&&+r\left(\left|\left(\Delta+\frac{2}{r^2}\right)\Jp\right|+\left|\dds_2\dds_1\Jp\right|\right)\Ga_g+\dkb^{\leq 1}(\Ga_b\c\Ga_g),
\eeaa
where by $\dds_1\Jp$, we mean  $\dds_1(0,\Jp)$.

Finally, recall from Lemma \ref{Lemma:transport.alongSi_*1} that we have along  $\Si_*$
\beaa
   \nab_\nu\left(\rhoc - \frac{1}{2}\chih\c\chibh \right)    &=&  -\div\bb -(1+O(r^{-1}))\div\b +O(r^{-1})\rhoc   +O(r^{-3})\kabc\\
   && +O(r^{-4}) \yc   +r^{-1}\dkb^{\le  1}( \Ga_b \c \Ga_b).
\eeaa
Together with \eqref{eq:transportequationforell=1modealongSigmastargeneral}, and noticing that the terms $O(r^a)$ only depend on $r$ and are thus constant on $S$, we infer
\beaa
\nu\left(\int_S\left(\rhoc - \frac{1}{2}\chih\c\chibh \right)\Jp\right) &=&  -\int_S\div\bb\Jp -(1+O(r^{-1}))\int_S\div\b\Jp \\
&&  +O(r^{-1})\int_S\left(\rhoc - \frac{1}{2}\chih\c\chibh \right)\Jp +O(r^{-3})\int_S\kabc\Jp\\
   && +O(r^{-4})\int_S\yc\Jp +r\dkb^{\le  1}( \Ga_b \c \Ga_b)
\eeaa
Next, recalling from the above that 
\beaa
\Delta\yc &=& -\div\xib -\Up(\div\ze-\div\eta)+\dkb(\Ga_b\c\Ga_b)
\eeaa
we write, using also integration by parts,
\beaa
O(r^{-4})\int_S\yc\Jp &=& O(r^{-2})\int_S\Delta\yc\Jp+r\left|\left(\Delta+\frac{2}{r^2}\right)\Jp\right|\Ga_b\\
&=& O(r^{-2})\int_S\div\xib\Jp+O(r^{-2})\int_S\div\ze\Jp+O(r^{-2})\int_S\div\eta\Jp\\
&&+r\left|\left(\Delta+\frac{2}{r^2}\right)\Jp\right|\Ga_b +\dkb(\Ga_b\c\Ga_b).
\eeaa
Together with the GCM conditions $(\div\eta)_{\ell=1}=0$ and $(\div\xib)_{\ell=1}=0$, we infer
\beaa
O(r^{-4})\int_S\yc\Jp &=& O(r^{-2})\int_S\div\ze\Jp+r\left|\left(\Delta+\frac{2}{r^2}\right)\Jp\right|\Ga_b +\dkb(\Ga_b\c\Ga_b)
\eeaa
and hence
\beaa
\nu\left(\int_S\left(\rhoc - \frac{1}{2}\chih\c\chibh \right)\Jp\right) &=&  -\int_S\div\bb\Jp -(1+O(r^{-1}))\int_S\div\b\Jp \\
&&  +O(r^{-1})\int_S\left(\rhoc - \frac{1}{2}\chih\c\chibh \right)\Jp +O(r^{-3})\int_S\kabc\Jp\\
   && +O(r^{-2})\int_S\div\ze\Jp+r\left|\left(\Delta+\frac{2}{r^2}\right)\Jp\right|\Ga_b \\
   &&+r\dkb^{\le  1}( \Ga_b \c \Ga_b)
\eeaa
as stated. This concludes the proof of Corollary \ref{corofLemma:transport.alongSi_*1}.

%%%%%%%%%%%%%%%%%%%%%%%%%%%%%

\section{Proof of Proposition \ref{prop:identitiesinqf}}
\lab{section:appendix-proofofprop:identitiesinqf}

%%%%%%%%%%%%%%%%%%%%%%%%%%%%%

In the proof, we denote by
\begin{itemize}
\item $(e_1, e_2, e_3, e_4)$ the null frame of $\Si_*$,

\item $(e_1', e_2', e_3', e_4')$ the null frame of $\Mext$,

\item $(e_1'', e_2'', e_3'', e_4'')$ the second null frame of $\Mext$ of Proposition \ref{prop:constructionsecondframeinMext},

\item $(e_1''', e_2''', e_3''', e_4''')$ the global null frame of $\MM$ of Proposition \ref{prop:existenceandestimatesfortheglobalframe:bis}.
\end{itemize}
Also, we denote 
\begin{itemize}
\item by $(f, \fb, \la)$ the change coefficients from $(e_1, e_2, e_3, e_4)$ to $(e_1', e_2', e_3', e_4')$,

\item by $(f', \fb', \la')$ the change coefficients from $(e_1', e_2', e_3', e_4')$ to $(e_1'', e_2'', e_3'', e_4'')$.
\end{itemize}

Since Proposition \ref{prop:identitiesinqf} involves identities on $\Si_*$, it suffices to consider a neighborhood of $\Si_*$ where we have in view of Proposition \ref{prop:existenceandestimatesfortheglobalframe:bis}
\bea\lab{eq:linkbetween3primeand2primenullframe}
e_4'''=\la e_4'', \qquad e_3'''=\la^{-1}e_3'', \qquad e_a'''=e_a'', \,\,\, a=1,2, \qquad \la=\frac{\De}{|q|^2}.
\eea
Also, in view of the construction of $(e_1'', e_2'', e_3'', e_4'')$ in Proposition \ref{prop:constructionsecondframeinMext}, and in view of the definition of $(f', \fb', \la')$, we have $\fb'=0$, $\la'=1$,  and 
\bea\lab{eq:linkbetween2primeand1primenullframe}
e_4'' = e_4 + {f'}^a e_a' +\frac 1 4 |f'|^2  e_3',\qquad e_a'' = e_a' +\frac 1 2 {f'}_a e_3',\,\,\, a=1,2, \qquad e_3''=  e_3',
\eea
where $f'$ satisfies on $\Mext$, and hence in a neighborhood of $\Si_*$, for any $k\leq k_*$,
\bea
|\dk^kf'| \les \frac{\ep}{ru^{\frac{1}{2}+\frac{\dec}{2}}}, \qquad |\dk^{k-1}\nab_3f'| \les \frac{\ep}{ru^{1+\frac{\dec}{2}}}.
\eea
Moreover, in view of the initialization of the PG frame of $\Mext$ on $\Si_*$, see section \ref{sec:initalizationadmissiblePGstructure}, and in view of the definition of $(f, \fb, \la)$, we have $\la=1$,  and 
\bea\lab{eq:linkbetween1primeandunprimednullframe}
 \bsplit
   e_4' &=e_4 + f^b  e_b +\frac 1 4 |f|^2  e_3,\\
  e_a' &= \left(\de_a^b +\frac{1}{2}\fb_af^b\right) e_b +\frac 1 2  \fb_a  e_4 +\left(\frac 1 2 f_a +\frac{1}{8}|f|^2\fb_a\right)   e_3,\quad a=1,2,\\
 e_3' &= \left(1+\frac{1}{2}f\c\fb  +\frac{1}{16} |f|^2  |\fb|^2\right) e_3 + \left(\fb^b+\frac 1 4 |\fb|^2f^b\right) e_b  + \frac 1 4 |\fb|^2 e_4,
 \end{split}
 \eea
where  $f$ and $\fb$ are given respectively by 
\bea
f_1 =0, \qquad f_2 =\frac{a\sin\th}{r}, \quad\textrm{on}\quad S_*, \qquad \nab_\nu(rf)=0\quad\textrm{on}\quad\Si_*,
\eea
and
\bea
\fb = -\frac{(\nu(r)-b_*)}{1-\frac{1}{4}b_* |f|^2}f\quad\textrm{on}\quad\Si_*.
\eea

Finally, recall that the quantity $\qf$ is defined with respect to the frame $(e_1''', e_2''', e_3''', e_4''')$ as follows,  see section \ref{sec:defintionofquantityqf}, 
\bea
\qf &=& q \ov{q}^{3} \Big( (\nab_{e_3'''}-2\omb''')(\nab_{e_3'''}-4\omb''')A''' + C_1(\nab_{e_3'''}-4\omb''')A''' + C_2   A'''\Big)
\eea
where the scalar functions $C_1$ and $C_2$ are given by\footnote{See  section 8 in \cite{GKS1}.}
\bea\label{finalchoicefordefinition-C}
C_1&=& 2 \trchb''' -4 i  \atrchb''', \qquad C_2= \frac 1 2 {\trchb'''}^2 -2i \trchb'''\atrchb'''.
\eea

We now ready to prove the identities \eqref{eq:alternateformulaforqfinvolvingtwoangularderrivativesofrho:M4} and \eqref{eq:Le-Teuk-Star1-M4} starting with the first one.

%%%%%%%%%%%%%%%%%%%%%%%%%%%%%%%%%%%%%%%%%%%%%%%%

\subsection{Proof of \eqref{eq:alternateformulaforqfinvolvingtwoangularderrivativesofrho:M4}}

%%%%%%%%%%%%%%%%%%%%%%%%%%%%%%%%%%%%%%%%%%%%%%%%%

In view of  Proposition \ref{prop:bianchi:complex} applied to the frame $(e_1''', e_2''', e_3''', e_4''')$, we have in particular 
\beaa
(\nab_{e_3'''}-4\omb''')A''' &=& \frac{1}{2}\DD'''\hot B''' -\frac{1}{2}\tr\Xb'''A'''+\frac{1}{2}(Z'''+4H''')\hot B''' -3\ov{P'''}\Xh'''.
\eeaa
Also, recall that in the frame $(e_1''', e_2''', e_3''', e_4''')$, we have $\Hc'''\in \Ga_g$ and the normalization in ingoing, so that 
\beaa
\tr\Xb''' =-\frac{2}{q}+\Ga_g''', \quad Z'''=\frac{aq}{|q|^2}\Jk+\Ga_g''',\quad H'''=\frac{aq}{|q|^2}\Jk+\Ga_g''', \quad  P'''=-\frac{2m}{q^3}+r^{-1}\Ga_g'''.
\eeaa
Together with the fact that $A'''\in r^{-1}\Ga_g'''$ and $B'''\in r^{-1}\Ga_g'''$, we infer
\beaa
(\nab_{e_3'''}-4\omb''')A''' &=& \frac{1}{2}\DD'''\hot B''' +\frac{1}{q}A'''+\frac{5}{2}\frac{aq}{|q|^2}\Jk\hot B''' +\frac{6m}{\ov{q}^3}\Xh''' +r^{-1}\Ga_g'''\c\Ga_g'''.
\eeaa
Plugging in the definition of $\qf$, and using additionally 
\beaa
C_1=-\frac{4}{r}+O(r^{-2})+\Ga_g''', \qquad C_2=\frac{2}{r^2}+O(r^{-3})+r^{-1}\Ga_g''',\qquad q=r+O(1),
\eeaa
we infer
\beaa
\qf &=& q \ov{q}^{3} \Bigg[ (\nab_{e_3'''}-2\omb''')\left(\frac{1}{2}\DD'''\hot B''' +\frac{1}{q}A'''+\frac{5}{2}\frac{aq}{|q|^2}\Jk\hot B''' +\frac{6m}{\ov{q}^3}\Xh''' \right)\\
&& -\frac{4}{r}\left(\frac{1}{2}\DD'''\hot B''' +\frac{1}{q}A'''\right) +\frac{2}{r^2} A'''\Bigg]+\dk^{\leq 1}\Ga_g'''+r^{2}\dk^{\le 1}(\Ga_b''' \c  \Ga_g''').
\eeaa
where we have used in particular the fact that $\nab_{e_3'''}(\Ga_g''')=r^{-1}\dk\Ga_b'''$ and $\nab_{3'''}(r)=-1+r\Ga_b'''$. Also, using again $q=r+O(1)$, $\nab_{3'''}(r)=-1+r\Ga_b'''$, as well as
\beaa
&&\omb'''=O(r^{-2})+\Ga_b''', \qquad \nab_{e_3'''}\Jk=O(r^{-1})\Jk+r^{-1}\Ga_b,\\
&&[\nab_{3'''}, \DD'''\hot]B''' =\frac{1}{r}\DD'''\hot B'''+r^{-4}\dk^{\leq 1}\Ga_g'''+r^{-2}\dk^{\le 1}(\Ga_b''' \c  \Ga_g'''),
\eeaa
we obtain
\beaa
\qf &=& q \ov{q}^{3} \left[ \frac{1}{2}\DD'''\hot\nab_{e_3'''} B'''  +\frac{1}{q}\nab_{e_3'''}A'''+\frac{5}{2}\frac{aq}{|q|^2}\Jk\hot\nab_{e_3'''} B''' +\frac{6m}{\ov{q}^3}\nab_{e_3'''}\Xh'''\right] \\
&& +r^4\left[-\frac{3}{2r}\DD'''\hot B''' -\frac{1}{r^2}A''' \right]+\dk^{\leq 1}\Ga_g'''+r^{2}\dk^{\le 1}(\Ga_b''' \c  \Ga_g''').
\eeaa

Next, we use the null structure equation for $\nab_{e_3'''}\Xh'''$ of Proposition \ref{prop-nullstr:complex} and the Bianchi identities for $\nab_{e_3'''} B'''$ and $\nab_{e_3'''}A'''$ of Proposition \ref{prop:bianchi:complex} according to which we have
\beaa
\nab_{e_3'''}A''' &=& \frac{1}{2}\DD'''\hot B''' +\frac{1}{r}A'''+r^{-3}\Ga_g'''+r^{-1}\Ga_b'''\c\Ga_g''',\\
\nab_{e_3'''}B''' &=& \DD'''\ov{P'''}+\frac{2}{r} B'''+3\ov{P'''}H''' +r^{-3}\Ga_g'''+r^{-1}\Ga_b'''\c\Ga_g'''\\
&=& \DD'''\ov{P'''}+\frac{2}{r} B'''+O(r^{-5}) +r^{-3}\Ga_g'''+r^{-1}\Ga_b'''\c\Ga_g''',\\
\nab_{e_3'''}\Xh'''  &=&  -\frac{1}{r}\Xbh'''+r^{-1}\dk^{\leq 1}\Ga_g'''+\Ga_b''' \c  \Ga_g'''.
\eeaa
where we recall that $O(r^a)$ denotes, for $a\in\RRR$, a function of $(r, \cos\th)$ bounded by $r^a$ as $r\to +\infty$.  We infer
\beaa
\qf &=& r^4\left[ \frac{1}{2}\DD'''\hot\left(\DD'''\ov{P'''}+\frac{2}{r} B''' \right)  +\frac{1}{r}\left(\frac{1}{2}\DD'''\hot B''' +\frac{1}{r}A'''\right)-\frac{6m}{r^4}\Xbh'''\right] \\
&& +r^4\left[-\frac{3}{2r}\DD'''\hot B''' -\frac{1}{r^2}A''' \right]+O(r^{-2})+\dk^{\leq 2}\Ga_g'''+r^{2}\dk^{\le 2}(\Ga_b''' \c  \Ga_g''')
\eeaa
and hence
\bea\lab{eq:usefulforproofLe-Teuk-Star1-M4}
\qf &=& r^4\left[ \frac{1}{2}\DD'''\hot\DD'''\ov{P'''}  -\frac{6m}{r^4}\Xbh'''\right] +O(r^{-2})+\dk^{\leq 2}\Ga_g'''+r^{2}\dk^{\le 2}(\Ga_b''' \c  \Ga_g''').
\eea
This yields in particular
\beaa
\qf &=&  \frac{1}{2}r^4\DD'''\hot\DD'''\ov{P'''}   +O(r^{-2})+\dk^{\leq 2}\Ga_b'''+r^{2}\dk^{\le 2}(\Ga_b''' \c  \Ga_g''').
\eeaa

Next, we deduce an identity in the second frame of $\Mext$, i.e. in the frame $(e_1'', e_2'', e_3'', e_4'')$. In view of \eqref{eq:linkbetween3primeand2primenullframe}, the change of frame coefficients $(f'', \fb'', \la'')$ from $(e_1'', e_2'', e_3'', e_4'')$ to $(e_1''', e_2''', e_3''', e_4''')$ satisfy 
\beaa
f''=0, \qquad \fb''=0, \qquad \la''=1+O(r^{-1}).
\eeaa
Together with the transformation formulas of Proposition \ref{Proposition:transformationRicci}, we obtain
\beaa
P'''=P'',\qquad \DD'''\hot\DD'''=\DD''\hot\DD'',
\eeaa
and hence
\beaa
\qf &=& \frac{1}{2}r^4\DD''\hot\DD''\ov{P''}  +O(r^{-2})+\dk^{\leq 2}\Ga_b'''+r^{2}\dk^{\le 2}(\Ga_b''' \c  \Ga_g''').
\eeaa

Next, we deduce an identity in the frame of $\Mext$, i.e. in the frame $(e_1', e_2', e_3', e_4')$. Recall that the change of frame coefficients $(f', \fb', \la')$ from $(e_1', e_2', e_3', e_4')$ to $(e_1'', e_2'', e_3'', e_4'')$ satisfy in particular $\fb'=0$ and $\la'=1$, see \eqref{eq:linkbetween2primeand1primenullframe}. Together with the transformation formulas of Proposition \ref{Proposition:transformationRicci}, we obtain
\beaa
P'' = P'   +r^{-1}f'\c\Ga_b'  +\lot, \qquad \DD'' = \DD' +\frac 1 2 f' \nab_3'.
\eeaa
We deduce
 \beaa
\qf &=&  \frac{1}{2}r^4\DD'\hot\DD'\ov{P'}  +O(r^{-2})+\dk^{\leq 2}\Ga_b'''+r^{2}\dk^{\le 2}\big(\Ga_b''' \c  \Ga_g'''\big)+r^{2}\dk^{\le 2}\big(r^{-1}f'\c\Ga_b\big).
\eeaa 

Finally, we derive an identity in the frame of $\Si_*$, i.e. in the frame $(e_1, e_2, e_3, e_4)$. Recall that the change of frame coefficients $(f, \fb, \la)$ from $(e_1, e_2, e_3, e_4)$ to $(e_1', e_2', e_3', e_4')$ satisfy in particular $\la=1$, as well as $f=O(r^{-1})$ and $\fb=O(r^{-1})+\Ga_b$, see \eqref{eq:linkbetween1primeandunprimednullframe}. Together with the transformation formulas of Proposition \ref{Proposition:transformationRicci}, we obtain
\beaa
P' = P  +O(r^{-5})  +r^{-2}\Ga_b, \qquad \DD' = \DD +\frac 1 2 f\nab_3+\frac{1}{2}\fb\nab_4.
\eeaa
We deduce
 \beaa
\qf &=&  \frac{1}{2}r^4\DD\hot\DD\ov{P}  +O(r^{-2})+\dkb^{\leq 2}\Ga_b+\dk^{\leq 2}\Ga_b'''+r^{2}\dk^{\le 2}\big(\Ga_b''' \c  \Ga_g'''\big)+r^{2}\dk^{\le 2}\big(r^{-1}f'\c\Ga_b\big).
\eeaa
As $\Ga_b'''$ and $\Ga_g'''$, $r^{-1}f'$, satisfy for $k\leq k_*$ the same estimates as $\Ga_b$, respectively $\Ga_g$, we write, by a slight abuse of notations
\beaa
\qf &=&  \frac{1}{2}r^4\DD\hot\DD\ov{P}  +O(r^{-2})+\dkb^{\leq 2}\Ga_b+r^{2}\dkb^{\le 2}\big(\Ga_b \c  \Ga_g\big).
\eeaa
Since we have
\beaa
\Re(\DD\hot\DD\ov{P}) &=& 2\nab\hot\nab\rho+2\nab\hot\dual\nab\rhod = 2\dds_2\dds_1(-\rho, \rhod),
\eeaa
we infer
\beaa
\Re(\qf) &=&  r^4\dds_2\dds_1(-\rho, \rhod)  +O(r^{-2})+\dkb^{\leq 2}\Ga_b+r^{2}\dkb^{\le 2}\big(\Ga_b \c  \Ga_g\big)
\eeaa
which is the desired identity \eqref{eq:alternateformulaforqfinvolvingtwoangularderrivativesofrho:M4}.

%%%%%%%%%%%%%%%%%%%%%%%%%%%%%%%%%%%%%%%%%%%%%%%%

\subsection{Proof of \eqref{eq:Le-Teuk-Star1-M4}}

%%%%%%%%%%%%%%%%%%%%%%%%%%%%%%%%%%%%%%%%%%%%%%%%%

Recall \eqref{eq:usefulforproofLe-Teuk-Star1-M4}
\beaa
\qf &=& r^4\left[ \frac{1}{2}\DD'''\hot\DD'''\ov{P'''}  -\frac{6m}{r^4}\Xbh'''\right] +O(r^{-2})+\dk^{\leq 2}\Ga_g'''+r^{2}\dk^{\le 2}(\Ga_b''' \c  \Ga_g''').
\eeaa
We multiply by $r$ and differentiate w.r.t. $\nab_{e_3'''}$. Using $e_3'''(r)=-1+r\Ga_b'''$ and 
\beaa
[\nab_{e_3'''},\DD'''\hot\DD''']\ov{P'''} &=& r^{-4}\dk^{\leq 2}\Ga_g'''+r^{-2}\dk^{\leq 2}(\Ga_b'''\c\Ga_g'''),
\eeaa
we infer
\beaa
\nab_{e_3'''}(r\qf) = r^5\left[ \frac{1}{2}\DD'''\hot\DD'''\nab_{e_3'''}\ov{P'''}  -\frac{6m}{r^4}\nab_{e_3'''}\Xbh'''\right] +O(r^{-2})+r\dk^{\leq 3}\Ga_g'''+r^{3}\dk^{\le 3}(\Ga_b''' \c  \Ga_g''').
\eeaa
Next, we use the null structure equation for $\nab_{e_3'''}\Xbh'''$ of Proposition \ref{prop-nullstr:complex} and the Bianchi identity for $\nab_{e_3'''}P'''$ of Proposition \ref{prop:bianchi:complex} according to which we have
\beaa
\nab_{e_3'''}\Xbh''' &=& -\Ab'''+r^{-1}\dk^{\leq 1}\Ga_b''',\\
\nab_{e_3'''}P'''  &=& -\frac{1}{2}\ov{\DD'''}\c\Bb''' -\frac{6m}{\ov{q}q^3}  +r^{-2}\Ga_g''' +\Ga_b'''\c\Ga_g''',
\eeaa
and hence
\beaa
\nab_{e_3'''}(r\qf) &=& r^5\left[-\frac{1}{4}\DD'''\hot\DD'''\DD'''\c\ov{\Bb'''} -3m\DD'''\hot\DD'''\left(\frac{1}{q\ov{q}^3}\right)  +\frac{6m}{r^4}\Ab'''\right] \\
&&+O(r^{-2})+r\dk^{\leq 3}\Ga_g'''+r^{3}\dk^{\le 3}(\Ga_b''' \c  \Ga_g''').
\eeaa
Since $\DD'''(q)=O(r^{-1})+r\Ga_g'''$, we have
\beaa
\DD'''\hot\DD'''\left(\frac{1}{q\ov{q}^3}\right) &=& O(r^{-7})+r^{-5}\dk^{\leq 1}\Ga_g'''
\eeaa
and thus
\beaa
\nab_{e_3'''}(r\qf) &=& -\frac{1}{4}r^5\DD'''\hot\DD'''\DD'''\c\ov{\Bb'''} +O(r)\Ab'''+O(r^{-2})+r\dk^{\leq 3}\Ga_g'''+r^{3}\dk^{\le 3}(\Ga_b''' \c  \Ga_g''').
\eeaa

As in the above proof of \eqref{eq:alternateformulaforqfinvolvingtwoangularderrivativesofrho:M4}, we now come back to the frame of $\Si_*$. First, since $f''=0$, $\fb''=0$, and $\la''=1+O(r^{-1})$, we have, together with the transformation formulas of Proposition \ref{Proposition:transformationRicci}, 
\beaa
\Bb'''=(1+O(r^{-1}))\Bb'',\quad \Ab'''=(1+O(r^{-1}))\Ab'', \quad \DD'''=\DD'',\quad e_3'''=(1+O(r^{-1}))e_3'',
\eeaa
and hence
\beaa
\nab_{e_3''}(r\qf) &=& -\frac{1}{4}r^5\DD''\hot\DD''\DD''\c\ov{\Bb''}  +O(r)\Ab'' +O(r^{-2})+r\dk^{\leq 3}\Ga_g'''+r^{3}\dk^{\le 3}(\Ga_b''' \c  \Ga_g''').
\eeaa
Next, since $\fb'=0$ and $\la'=1$, we have, together with the transformation formulas of Proposition \ref{Proposition:transformationRicci}, 
\beaa
&& \Bb''=\Bb'+f'\c\Ga_b'+O(r^{-3})f',\qquad \Ab''=\Ab'+r^{-1}f'\c\Ga_b', \\ 
&& \DD'' = \DD' +\frac 1 2 f' \nab_3',\qquad e_3''=e_3'.
\eeaa
We deduce
\beaa
\nab_{e_3'}(r\qf) &=& -\frac{1}{4}r^5\DD'\hot\DD'\DD'\c\ov{\Bb'}  +O(r)\Ab' +O(r^{-2})+r\dk^{\leq 3}\Ga_g'''+\dk^{\leq 3}f'\\
&&+r^{3}\dk^{\le 3}(\Ga_b''' \c  \Ga_g''')+r^3\dk^{\le 2}\big(r^{-1}f'\c\Ga_b\big).
\eeaa
Finally, since $\la=1$, as well as $f=O(r^{-1})$ and $\fb=O(r^{-1})+\Ga_b$, together with the transformation formulas of Proposition \ref{Proposition:transformationRicci}, we have
\beaa
&& \Bb'=\Bb+O(r^{-1})\Ab+O(r^{-4})+r^{-3}\Ga_b,\qquad \Ab'=\Ab+r^{-2}\Ga_b, \\ 
&& \DD' = \DD +O(r^{-1})\nab_3+O(r^{-1})\nab_4,\qquad e_3'=(1+O(r^{-2}))e_3+O(r^{-1})\nab+O(r^{-2})\nab_4.
\eeaa
We infer
\beaa
\nab_3(r\qf) &=& -\frac{1}{4}r^5\DD\hot\DD\DD\c\ov{\Bb}  +O(r)\dkb^{\leq 3}\Ab +O(r^{-2})+r\dk^{\leq 3}\Ga_g'''+\dk^{\leq 3}f'+r\dkb^{\leq 3}\Ga_g\\
&&+r^{3}\dk^{\le 3}(\Ga_b''' \c  \Ga_g''')+r^3\dk^{\le 2}\big(r^{-1}f'\c\Ga_b\big)+r^{3}\dkb^{\le 3}(\Ga_b \c  \Ga_g)+O(r^{-1})\dk^{\leq 1}\qf.
\eeaa
As $\Ga_b'''$ and $\Ga_g'''$, $r^{-1}f'$, $r^{-1}\qf$, satisfy for $k\leq k_*$ the same estimates as $\Ga_b$, respectively $\Ga_g$, we write, by a slight abuse of notations
\beaa
\nab_3(r\qf) &=& -\frac{1}{4}r^5\DD\hot\DD\DD\c\ov{\Bb}  +O(r)\dkb^{\leq 3}\Ab +O(r^{-2})+r\dkb^{\leq 3}\Ga_g+r^{3}\dkb^{\le 3}(\Ga_b \c  \Ga_g).
\eeaa
Since we have
\beaa
\Re(\DD\hot\DD\DD\c\ov{\Bb}) &=& 4\nab\hot\nab\div\bb+4\nab\hot\dual\nab\curl\bb = 4\dds_2\dds_1(-\div\bb, \curl\bb),
\eeaa
we infer
\beaa
\Re(\nab_3(r\qf)) = r^5\dds_2\dds_1\big(\div\bb, -\curl\bb\big) +O(r)\dkb^{\leq 3}\aa +O(r^{-2})+r\dkb^{\leq 3}\Ga_g+r^{3}\dkb^{\le 3}(\Ga_b \c  \Ga_g).
\eeaa
which is the desired identity \eqref{eq:Le-Teuk-Star1-M4}. This concludes the proof of Proposition \ref{prop:identitiesinqf}.

%%%%%%%%%%%%%%%%%%%%%%%%%%%%%%%%%%%%%%%%%%%%%%%%

\section{Proof of Lemma \ref{lemma:computationforfirstorderderivativesf0fpfm}}
\lab{sec:proofoflemma:computationforfirstorderderivativesf0fpfm}

%%%%%%%%%%%%%%%%%%%%%%%%%%%%%%%%%%%%%%%%%%%%%%%%%

Let the orthonormal basis $(e_1, e_2)$ of $S_*$ given by \eqref{eq:specialorthonormalbasisofSstar}, i.e.
\beaa
e_1=\frac{1}{re^\phi}\pr_\th, \quad e_2=\frac{1}{r\sin\th e^\phi}\pr_\vphi.
\eeaa
In order to prove Lemma \ref{lemma:computationforfirstorderderivativesf0fpfm}, we will need the following simple lemma.
\begin{lemma}
On $S_*$, for  $(e_1, e_2)$  given by \eqref{eq:specialorthonormalbasisofSstar}, we have
\bea
\bsplit
(\La_1)_{12} &:= \g(\D_1e_2, e_1)= \frac{1}{r\sin\th e^{\phi}}\pr_\vphi(\phi),\\
(\La_2)_{12} &:= \g(\D_2e_2, e_1)= -\frac{1}{re^{\phi}}\Big(\cot\th+\pr_\th(\phi)\Big).
\end{split}
\eea
\end{lemma}
 
\begin{proof}
Recall that we have on $S_*$
\beaa
g_{S_*} =r^2e^{2\phi}\Big((d\th)^2+(\sin\th)^2(d\vphi)^2\Big), \qquad e_1=\frac{1}{re^{\phi}}\pr_\th, \qquad e_2=\frac{1}{r\sin\th e^{\phi}}\pr_\vphi.
\eeaa
This allows us to compute
\beaa
(\La_1)_{12} &=& \g(\D_1e_2, e_1)=\frac{1}{r^3\sin\th e^{3\phi}}\g(\D_{\pr_\th}\pr_\vphi, \pr_\th)=\frac{1}{r^3\sin\th e^{3\phi}}\g(\D_{\pr_\vphi}\pr_\th, \pr_\th)\\
&=& \frac{1}{r^3\sin\th e^{3\phi}}\frac{1}{2}\pr_\vphi(\g_{\th\th})=\frac{1}{r^3\sin\th e^{3\phi}}\frac{1}{2}\pr_\vphi(r^2e^{2\phi})\\
&=& \frac{1}{r\sin\th e^{\phi}}\pr_\vphi(\phi)
\eeaa
and 
\beaa
(\La_2)_{12} &=& \g(\D_2e_2, e_1)=\frac{1}{r^3(\sin\th)^2 e^{3\phi}}\g(\D_{\pr_\vphi}\pr_\vphi, \pr_\th)=-\frac{1}{r^3(\sin\th)^2 e^{3\phi}}\g(\pr_\vphi, \D_{\pr_\vphi}\pr_\th)\\
&=& -\frac{1}{r^3(\sin\th)^2 e^{3\phi}}\g(\pr_\vphi, \D_{\pr_\th}\pr_\vphi)=-\frac{1}{r^3(\sin\th)^2 e^{3\phi}}\frac{1}{2}\pr_\th(\g_{\vphi\vphi})\\
&=& -\frac{1}{r^3(\sin\th)^2 e^{3\phi}}\frac{1}{2}\pr_\th(r^2(\sin\th)^2e^{2\phi})=-\frac{1}{r(\sin\th)^2 e^{\phi}}\Big(\sin\th\cos\th +(\sin\th)^2\pr_\th(\phi)\Big)\\
&=& -\frac{1}{re^{\phi}}\Big(\cot\th+\pr_\th(\phi)\Big)
\eeaa
as stated.
\end{proof}

We are now ready to prove Lemma \ref{lemma:computationforfirstorderderivativesf0fpfm}. We start with the identities for $f_0$. We have
\beaa
\div(f_0) &=& \nab_1(f_0)_1+\nab_2(f_0)_2=e_1((f_0)_1)-\g(\D_1e_1, e_2)(f_0)_2 +e_2((f_0)_2)-\g(\D_2e_2, e_1)(f_0)_1\\
&=& -\g(\D_1e_1, e_2)\sin\th +e_2(\sin\th) = (\La_1)_{12}\sin\th +\frac{1}{r\sin\th e^{\phi}}\pr_\vphi(\sin\th) = \frac{1}{re^{\phi}}\pr_\vphi(\phi)\\
&=& f_0\c\nab\phi 
\eeaa
and 
\beaa
\curl(f_0) &=& \nab_1(f_0)_2-\nab_2(f_0)_1=e_1((f_0)_2) -\g(\D_1e_2, e_1)(f_0)_1 - e_2((f_0)_1) +\g(\D_2e_1, e_2)(f_0)_2\\
&=& \frac{1}{re^{\phi}}\pr_\th(\sin\th)-(\La_2)_{12}\sin\th= \frac{1}{re^{\phi}}\cos\th +\frac{1}{re^{\phi}}\Big(\cos\th+\sin\th\pr_\th(\phi)\Big)\\
&=& \frac{2}{re^{\phi}}\cos\th +\frac{1}{re^{\phi}}\sin\th\pr_\th(\phi)= \frac{2}{re^{\phi}}\cos\th-f_0\wedge\nab\phi
\eeaa
as stated. Also, we have
\beaa
(\nab\hot f_0)_{11} &=& \nab_1(f_0)_1- \nab_2(f_0)_2\\
&=& e_1((f_0)_1)-\g(\D_1e_1, e_2)(f_0)_2 -e_2((f_0)_2)+\g(\D_2e_2, e_1)(f_0)_1\\
&=& -\g(\D_1e_1, e_2)\sin\th -e_2(\sin\th) = (\La_1)_{12}\sin\th -\frac{1}{r\sin\th e^{\phi}}\pr_\vphi(\sin\th) = \frac{1}{re^{\phi}}\pr_\vphi(\phi)\\
&=& f_0\c\nab\phi, 
\eeaa
and 
\beaa
(\nab\hot f_0)_{12} &=& \nab_1(f_0)_2+\nab_2(f_0)_1=e_1((f_0)_2) -\g(\D_1e_2, e_1)(f_0)_1 + e_2((f_0)_1) -\g(\D_2e_1, e_2)(f_0)_2\\
&=& \frac{1}{re^{\phi}}\pr_\th(\sin\th)+(\La_2)_{12}\sin\th= \frac{1}{re^{\phi}}\cos\th -\frac{1}{re^{\phi}}\Big(\cos\th+\sin\th\pr_\th(\phi)\Big)\\
&=& -\frac{1}{re^{\phi}}\sin\th\pr_\th(\phi)=f_0\wedge\nab\phi
\eeaa
so that 
\beaa
\nab\hot f_0 &=& \left(\ba{cc}
 f_0\c\nab\phi & f_0\wedge\nab\phi\\
f_0\wedge\nab\phi & - f_0\c\nab\phi
\ea\right)
\eeaa
as stated.

Next, we consider the  identities for $f_+$. We have
\beaa
\div(f_+) &=& \nab_1(f_+)_1+\nab_2(f_+)_2=e_1((f_+)_1)-\g(\D_1e_1, e_2)(f_+)_2 +e_2((f_+)_2)-\g(\D_2e_2, e_1)(f_+)_1\\
&=& \frac{1}{re^\phi}\left(-\sin\th\cos\vphi +  \frac{1}{\sin\th}\pr_\vphi(\phi)(-\sin\vphi)-\frac{\cos\vphi}{\sin\th}+\Big(\cot\th+\pr_\th(\phi)\Big)\cos\th\cos\vphi\right)\\
&=& -\frac{2}{re^\phi}J^{(+)}+f_+\c\nab\phi
\eeaa
and 
\beaa
\curl(f_+) &=& \nab_1(f_+)_2-\nab_2(f_+)_1=e_1((f_+)_2) -\g(\D_1e_2, e_1)(f_+)_1 - e_2((f_+)_1) +\g(\D_2e_1, e_2)(f_+)_2\\
&=& \frac{1}{re^\phi}\left( - \frac{1}{\sin\th}\pr_\vphi(\phi)\cos\th\cos\vphi +\cot\th\sin\vphi +\Big(\cot\th+\pr_\th(\phi)\Big)(-\sin\vphi)\right)\\
&=& f_+\wedge\nab\phi,
\eeaa
as stated. Also, we have
\beaa
(\nab\hot f_+)_{11} &=& \nab_1(f_+)_1- \nab_2(f_+)_2\\
&=& e_1((f_+)_1)-\g(\D_1e_1, e_2)(f_+)_2 -e_2((f_+)_2)+\g(\D_2e_2, e_1)(f_+)_1\\
&=& \frac{1}{re^\phi}\left(-\sin\th\cos\vphi +  \frac{1}{\sin\th}\pr_\vphi(\phi)(-\sin\vphi)+\frac{\cos\vphi}{\sin\th}-\Big(\cot\th+\pr_\th(\phi)\Big)\cos\th\cos\vphi\right)\\
&=& (f_+)_2\nab_2\phi- (f_+)_1\nab_1\phi, 
\eeaa
and 
\beaa
(\nab\hot f_+)_{12} &=& \nab_1(f_+)_2+\nab_2(f_+)_1=e_1((f_+)_2) -\g(\D_1e_2, e_1)(f_+)_1 + e_2((f_+)_1) -\g(\D_2e_1, e_2)(f_+)_2\\
&=& \frac{1}{re^\phi}\left( - \frac{1}{\sin\th}\pr_\vphi(\phi)\cos\th\cos\vphi -\cot\th\sin\vphi -\Big(\cot\th+\pr_\th(\phi)\Big)(-\sin\vphi)\right)\\
&=& -(f_+)_1\nab_2\phi- (f_+)_2\nab_1\phi
\eeaa
so that 
\beaa
\nab\hot f_+ &=& \left(\ba{cc}
(f_+)_2\nab_2\phi- (f_+)_1\nab_1\phi & -(f_+)_1\nab_2\phi- (f_+)_2\nab_1\phi\\
-(f_+)_1\nab_2\phi- (f_+)_2\nab_1\phi & -(f_+)_2\nab_2\phi+ (f_+)_1\nab_1\phi
\ea\right)
\eeaa
as stated.

Finally, we consider the  identities for $f_-$. We have
\beaa
\div(f_-) &=& \nab_1(f_-)_1+\nab_2(f_-)_2=e_1((f_-)_1)-\g(\D_1e_1, e_2)(f_-)_2 +e_2((f_-)_2)-\g(\D_2e_2, e_1)(f_+)_1\\
&=& \frac{1}{re^\phi}\left(-\sin\th\sin\vphi +  \frac{1}{\sin\th}\pr_\vphi(\phi)\cos\vphi -\frac{\sin\vphi}{\sin\th}+\Big(\cot\th+\pr_\th(\phi)\Big)\cos\th\sin\vphi\right)\\
&=& -\frac{2}{re^\phi}J^{(-)}+f_-\c\nab\phi
\eeaa
and 
\beaa
\curl(f_-) &=& \nab_1(f_-)_2-\nab_2(f_-)_1=e_1((f_-)_2) -\g(\D_1e_2, e_1)(f_-)_1 - e_2((f_-)_1) +\g(\D_2e_1, e_2)(f_-)_2\\
&=& \frac{1}{re^\phi}\left( - \frac{1}{\sin\th}\pr_\vphi(\phi)\cos\th\sin\vphi -\cot\th\cos\vphi +\Big(\cot\th+\pr_\th(\phi)\Big)\cos\vphi\right)\\
&=& -f_-\wedge\nab\phi
\eeaa
as stated. Also, we have
\beaa
(\nab\hot f_-)_{11} &=& \nab_1(f_-)_1- \nab_2(f_-)_2\\
&=& e_1((f_-)_1)-\g(\D_1e_1, e_2)(f_-)_2 -e_2((f_-)_2)+\g(\D_2e_2, e_1)(f_-)_1\\
&=& \frac{1}{re^\phi}\left(-\sin\th\sin\vphi +  \frac{1}{\sin\th}\pr_\vphi(\phi)\cos\vphi +\frac{\sin\vphi}{\sin\th}-\Big(\cot\th+\pr_\th(\phi)\Big)\cos\th\sin\vphi\right)\\
&=& (f_-)_2\nab_2\phi- (f_-)_1\nab_1\phi, 
\eeaa
and 
\beaa
(\nab\hot f_-)_{12} &=& \nab_1(f_-)_2+\nab_2(f_-)_1=e_1((f_-)_2) -\g(\D_1e_2, e_1)(f_-)_1 + e_2((f_-)_1) -\g(\D_2e_1, e_2)(f_-)_2\\
&=& \frac{1}{re^\phi}\left( - \frac{1}{\sin\th}\pr_\vphi(\phi)\cos\th\sin\vphi +\cot\th\cos\vphi -\Big(\cot\th+\pr_\th(\phi)\Big)\cos\vphi\right)\\
&=& -(f_-)_1\nab_2\phi- (f_-)_2\nab_1\phi
\eeaa
so that 
\beaa
\nab\hot f_- &=& \left(\ba{cc}
(f_-)_2\nab_2\phi- (f_-)_1\nab_1\phi & -(f_-)_1\nab_2\phi- (f_-)_2\nab_1\phi\\
-(f_-)_1\nab_2\phi- (f_-)_2\nab_1\phi & -(f_-)_2\nab_2\phi+ (f_-)_1\nab_1\phi
\ea\right)
\eeaa
as stated. This concludes the proof of Lemma \ref{lemma:computationforfirstorderderivativesf0fpfm}.

%%%%%%%%%%%%%%%%%%%%%%%%%%%%%%%%%%%%%%%%%%%%%%%%

\section{Proof of Lemma \ref{lemma:relationangularderivativesJkJkpmandff0fpfm}}
\lab{sec:proofoflemma:relationangularderivativesJkJkpmandff0fpfm}

%%%%%%%%%%%%%%%%%%%%%%%%%%%%%%%%%%%%%%%%%%%%%%%%%

Let $f$ a 1-form and $F=f+i\dual f$. Then, we have
\beaa
\ov{\DD}\c F &=& 2\div(f)+2i\curl(f).
\eeaa
Also, since $|q|^2=r^2+a^2(\cos\th)^2$, $J^{(0)}=\cos\th$, and $\nab(r)=0$, we have
\beaa
\div\left(\frac{1}{|q|}f\right) &=& \frac{1}{|q|}\div(f) -\frac{\nab(|q|)}{|q|^2}\c f = \frac{1}{|q|}\div(f) -\frac{\nab(r^2+a^2(\cos\th)^2)}{2|q|^3}\c f\\
&=& \frac{1}{|q|}\div(f) -\frac{a^2\cos\th\nab(J^{(0)})}{|q|^3}\c f\\
&=& \frac{1}{|q|}\div(f) +\frac{a^2\cos\th}{r|q|^3}\dual f_0\c f -\frac{a^2\cos\th}{|q|^3}\widecheck{\nab J^{(0)}}\c f
\eeaa
and 
\beaa
\curl\left(\frac{1}{|q|}f\right) &=& \frac{1}{|q|}\curl(f) +\frac{\dual\nab(|q|)}{|q|^2}\c f = \frac{1}{|q|}\curl(f) -\frac{\dual\nab(r^2+a^2(\cos\th)^2)}{2|q|^3}\c f\\
&=& \frac{1}{|q|}\curl(f) -\frac{a^2\cos\th\dual\nab(J^{(0)})}{|q|^3}\c f\\
&=& \frac{1}{|q|}\curl(f) -\frac{a^2\cos\th}{r|q|^3}f_0\c f -\frac{a^2\cos\th}{|q|^3}\dual\widecheck{\nab J^{(0)}}\c f.
\eeaa
Hence, we have
\beaa
\ov{\DD}\c\left(\frac{1}{|q|}F\right) &=&  \frac{2}{|q|}\div(f) +\frac{2a^2\cos\th}{r|q|^3}\dual f_0\c f -\frac{2a^2\cos\th}{|q|^3}\widecheck{\nab J^{(0)}}\c f\\
&&+i\left(\frac{2}{|q|}\curl(f) -\frac{2a^2\cos\th}{r|q|^3}f_0\c f -\frac{2a^2\cos\th}{|q|^3}\dual\widecheck{\nab J^{(0)}}\c f\right).
\eeaa
Since we have by definition on $\Si_*$
\beaa
\Jk = \frac{1}{|q|}\left(f_0+i\dual f_0\right), \qquad \Jk_\pm = \frac{1}{|q|}\left(f_\pm+i\dual f_\pm\right),
\eeaa 
we infer 
\beaa
\ov{\DD}\c\Jk &=&  \frac{2}{|q|}\div(f_0) +\frac{2a^2\cos\th}{r|q|^3}\dual f_0\c f_0 -\frac{2a^2\cos\th}{|q|^3}f_0\c\widecheck{\nab J^{(0)}}\\
&&+i\left(\frac{2}{|q|}\curl(f_0) -\frac{2a^2\cos\th}{r|q|^3}f_0\c f_0 -\frac{2a^2\cos\th}{|q|^3}f_0\c\dual\widecheck{\nab J^{(0)}}\right)
\eeaa
and 
\beaa
\ov{\DD}\c\Jk_\pm &=&  \frac{2}{|q|}\div(f_\pm) +\frac{2a^2\cos\th}{r|q|^3}\dual f_0\c f_\pm -\frac{2a^2\cos\th}{|q|^3}f_\pm\c\widecheck{\nab J^{(0)}}\\
&&+i\left(\frac{2}{|q|}\curl(f_\pm) -\frac{2a^2\cos\th}{r|q|^3}f_0\c f_\pm -\frac{2a^2\cos\th}{|q|^3}f_\pm\c\dual\widecheck{\nab J^{(0)}}\right).
\eeaa
Since we have on $\Si_*$
\beaa
&& f_0\c f_0 = (\sin\th)^2, \quad f_+\c f_0=-J^{(-)}, \quad f_-\c f_0=J^{(+)}, \\
&& f_0\c\dual f_0 = 0, \quad f_+\c\dual f_0=\sin\th\cos\th\cos\vphi, \quad f_-\c\dual f_0=\sin\th\cos\th\sin\vphi,
\eeaa
we deduce
\beaa
\ov{\DD}\c\Jk &=&  \frac{2}{|q|}\div(f_0)  -\frac{2a^2\cos\th}{|q|^3}f_0\c\widecheck{\nab J^{(0)}}\\
&&+i\left(\frac{2}{|q|}\curl(f_0) -\frac{2a^2\cos\th(\sin\th)^2}{r|q|^3} -\frac{2a^2\cos\th}{|q|^3}f_0\c\dual\widecheck{\nab J^{(0)}}\right),
\eeaa
\beaa
\ov{\DD}\c\Jk_+ &=&  \frac{2}{|q|}\div(f_+) +\frac{2a^2\cos\th}{r|q|^3}\sin\th\cos\th\cos\vphi -\frac{2a^2\cos\th}{|q|^3}f_+\c\widecheck{\nab J^{(0)}}\\
&&+i\left(\frac{2}{|q|}\curl(f_+) +\frac{2a^2\cos\th\sin\th\sin\vphi}{r|q|^3} -\frac{2a^2\cos\th}{|q|^3}f_+\c\dual\widecheck{\nab J^{(0)}}\right),
\eeaa
and
\beaa
\ov{\DD}\c\Jk_- &=&  \frac{2}{|q|}\div(f_-) +\frac{2a^2\cos\th}{r|q|^3}\sin\th\cos\th\sin\vphi -\frac{2a^2\cos\th}{|q|^3}f_-\c\widecheck{\nab J^{(0)}}\\
&&+i\left(\frac{2}{|q|}\curl(f_-) -\frac{2a^2\cos\th\sin\th\cos\vphi}{r|q|^3} -\frac{2a^2\cos\th}{|q|^3}f_-\c\dual\widecheck{\nab J^{(0)}}\right).
\eeaa
In view of the definition of $\widecheck{\curl(f_0)}$, $\widecheck{\div(f_\pm)}$, $\widecheck{\ov{\DD}\c\Jk}$, and $\widecheck{\ov{\DD}\c\Jk_\pm}$
we obtain 
\beaa
\widecheck{\ov{\DD}\c\Jk} &=& O(r^{-4})+ \frac{2}{|q|}\div(f_0)  -\frac{2a^2\cos\th}{|q|^3}f_0\c\widecheck{\nab J^{(0)}}\\
&&+i\left(\frac{2}{|q|}\widecheck{\curl(f_0)} -\frac{2a^2\cos\th}{|q|^3}f_0\c\dual\widecheck{\nab J^{(0)}}\right),
\eeaa
\beaa
\widecheck{\ov{\DD}\c\Jk_\pm} &=& O(r^{-4}) + \frac{2}{|q|}\widecheck{\div(f_\pm)} -\frac{2a^2\cos\th}{|q|^3}f_\pm\c\widecheck{\nab J^{(0)}}\\
&&+i\left(\frac{2}{|q|}\curl(f_\pm)  -\frac{2a^2\cos\th}{|q|^3}f_\pm\c\dual\widecheck{\nab J^{(0)}}\right),
\eeaa
where $O(r^a)$ denotes, for $a\in\RRR$, a function of $(r, \cos\th)$ bounded by $r^a$ as $r\to +\infty$. 
This condoles the proof of Lemma \ref{lemma:relationangularderivativesJkJkpmandff0fpfm}.

%%%%%%%%%%%%%%%%%%%%%%%%%%%%%%%%%

\chapter{Proof of results in Chapter \ref{Chapter:DecayMext}}

%%%%%%%%%%%%%%%%%%%%%%%%%%%%%%%%%

%%%%%%%%%%%%%%%%%%%%%%%%%%%%%%%%%%%%%%%%%%%%%%%%

\section{Proof of Lemma \ref{Lemma:linearized-nullstr}}
\lab{appendix: Proof of Lemma-ref{Lemma:linearized-nullstr}}

%%%%%%%%%%%%%%%%%%%%%%%%%%%%%%%%%%%%%%%%%%%%%%%%

The proof relies on the null structure equations and Bianchi identities of Proposition \ref{prop-nullstrandBianchi:complex:outgoing:again:chap6}, the definition of the linearized quantities and of $\Ga_g$ and $\Ga_b$ in section \ref{sec:linearizedquantitiesanddefGabGag:chap6}, the notation $O(r^{-p})$ made in Definition \ref{def:ordermagnitude}, the fact that  $a$ and $m$ are constants, and the following identities 
\beaa
e_4(r)=1, \qquad e_4(\th)=0, \qquad \nab_4\Jk=-\frac{1}{q}\Jk, \qquad e_4(q)=1, \qquad e_4(\ov{q})=1, \qquad \nab(r)=0,
\eeaa
where we used in particular the fact that $q=r+ai\cos\th$ and $\ov{q}=r-ai\sin\th$.  

We start with the equation for $\trXc$. Recall
\beaa
\nab_4\tr X +\frac{1}{2}(\tr X)^2 &=& -\frac{1}{2}\Xh\c\ov{\Xh}.
\eeaa
Since $e_4(q)=1$, we infer
\beaa
\nab_4(\widecheck{\tr X}) &=& -\frac{1}{2}(\tr X)^2  -\frac{1}{2}\Xh\c\ov{\Xh} -\pr_r\left(\frac{2}{q}\right)\\
&=& -\frac{1}{2}(\tr X)^2 +\frac{1}{2}\left(\frac{2}{q}\right)^2 -\frac{1}{2}\Xh\c\ov{\Xh} 
\eeaa
and hence
\beaa
\nab_4(\widecheck{\tr X}) +\frac{2}{q} \widecheck{\tr X} &=& \Ga_g\c\Ga_g.
\eeaa

Next, recall
\beaa
\nab_4\Xh+\Re(\tr X)\Xh &=& -A
\eeaa
and  hence
\beaa
\nab_4\Xh+\frac{2r}{|q|^2}\Xh &=& -A+\Ga_g\c\Ga_g.
\eeaa

Next, recall
\beaa
\nab_4Z +\tr X Z &=&  -\widehat{X}\c\ov{Z} -B.
\eeaa
We infer
\beaa
\nab_4\Zc &=& -\tr X Z   -\widehat{X}\c\ov{Z} -B -\nab_4\left(\frac{a\ov{q}}{|q|^2}\Jk\right).
\eeaa
Now, since $e_4(q)=1$ and $\nab_4\Jk=-q^{-1}\Jk$, we have
\beaa
\nab_4\left(\frac{a\ov{q}}{|q|^2}\Jk\right) &=& \left(\frac{a}{|q|^2} -\frac{2a\ov{q}e_4(|q|)}{|q|^3}\right)\Jk -\frac{a\ov{q}}{q|q|^2}\Jk\\
&=& \frac{a}{|q|^4}\Big(|q|^2-\ov{q}(q+\ov{q})-\ov{q}^2\Big)\Jk = -\frac{2a\ov{q}^2}{|q|^4}\Jk
\eeaa
and thus
\beaa
\nab_4\Zc + \frac{2}{q}\Zc  &=&  -\frac{a\ov{q}}{|q|^2}\Jk\trXc  - \frac{aq}{|q|^2}\ov{\Jk}\c\widehat{X} -B +\Ga_g\c\Ga_g.
\eeaa
We infer
\beaa
\nab_4\Zc + \frac{2}{q}\Zc  &=&    - \frac{aq}{|q|^2}\ov{\Jk}\c\widehat{X} -B +O(r^{-2})\trXc+\Ga_g\c\Ga_g.
\eeaa

Next, recall
\beaa
\nab_4H +\frac{1}{2}\ov{\tr X}(H+Z) &=&   -\frac{1}{2}\Xh\c(\ov{H}+\ov{Z}) -B,
\eeaa
we infer
\beaa
\nab_4\Hc &=& -\frac{1}{2}\ov{\tr X}(H+Z)    -\frac{1}{2}\Xh\c(\ov{H}+\ov{Z}) -B  -\nab_4\left(\frac{aq}{|q|^2}\Jk\right).
\eeaa
Now, since $e_4(q)=1$ and $\nab_4\Jk=-q^{-1}\Jk$, we have
\beaa
\nab_4\left(\frac{aq}{|q|^2}\Jk\right) &=& \left(\frac{a}{|q|^2} -\frac{2aqe_4(|q|)}{|q|^3}\right)\Jk -\frac{a}{|q|^2}\Jk =  -\frac{aq(q+\ov{q})}{|q|^4}\Jk
\eeaa
and thus
\beaa
\nab_4\Hc+\frac{1}{\ov{q}}\Hc &=& -\frac{1}{\ov{q}}\Zc -\frac{ar}{|q|^2}\ov{\trXc}\Jk    -\frac{ar}{|q|^2}\ov{\Jk}\c\Xh -B  +\Ga_b\c\Ga_g.
\eeaa
We infer
\beaa
\nab_4\Hc+\frac{1}{\ov{q}}\Hc &=& -\frac{1}{\ov{q}}\Zc     -\frac{ar}{|q|^2}\ov{\Jk}\c\Xh -B +O(r^{-2})\trXc +\Ga_b\c\Ga_g.
\eeaa

Next, recall
\beaa
\nab_4\tr\Xb +\frac{1}{2}\tr X\tr\Xb &=& -\DD\c\ov{Z}+Z\c\ov{Z}+2\ov{P}-\frac{1}{2}\Xh\c\ov{\Xbh}.
\eeaa
We infer
\beaa
\nab_4\widecheck{\tr\Xb}  &=& \nab_4\left(\frac{2q\Delta}{|q|^4}\right)-\frac{1}{2}\tr X\tr\Xb -\DD\c\ov{Z}+Z\c\ov{Z}+2\ov{P}\\
&&-\frac{1}{2}\Xh\c\ov{\Xbh}\\
&=& \pr_r\left(\frac{2q\Delta}{|q|^4}\right)-\frac{1}{2}\left(\frac{2}{q}+\widecheck{\tr X}\right)\left(-\frac{2q\Delta}{|q|^4}+\widecheck{\tr\Xb} \right)\\
&& -\DD\c\ov{Z}+Z\c\ov{Z}+2\left(-\frac{2m}{(\ov{q})^3}+\ov{\Pc }\right)-\frac{1}{2}\Xh\c\ov{\Xbh}, 
\eeaa
and hence
\beaa
\nab_4\widecheck{\tr\Xb} +\frac{1}{q}\trXbc &=&  -\DD\c\ov{Z}+Z\c\ov{Z}+2\ov{\Pc }+\pr_r\left(\frac{2q\Delta}{|q|^4}\right)+\frac{2\Delta}{|q|^4} -\frac{4m}{(\ov{q})^3}\\
&&+O(r^{-1})\trXc+\Ga_b\c\Ga_g. 
\eeaa
Next, we compute
\beaa
&&  -\DD\c\ov{Z}+Z\c\ov{Z} \\
&=& -\DD\c\ov{\frac{a\ov{q}}{|q|^2}\Jk+\Zc} +\left(\frac{a\ov{q}}{|q|^2}\Jk+\Zc\right)\c\left(\ov{\frac{a\ov{q}}{|q|^2}\Jk+\Zc}\right) \\
&=& -\frac{aq}{|q|^2}\DD\c\ov{\Jk}-\DD\left(\frac{aq}{|q|^2}\right)\c\ov{\Jk}-\DD\c\ov{\Zc}+\frac{a^2}{|q|^2}\Jk\c\ov{\Jk}+\frac{a\ov{q}}{|q|^2}\Jk\c\ov{\Zc}+\frac{aq}{|q|^2}\ov{\Jk}\c\Zc+\Ga_g\c\Ga_g.
\eeaa
Since $\Jk\c\ov{\Jk}=\frac{2(\sin\th)^2}{|q|^2}$, we infer
\beaa
&&  -\DD\c\ov{Z}+Z\c\ov{Z} \\
&=& -\frac{aq}{|q|^2}\left(-\frac{4i(r^2+a^2)\cos\th}{|q|^4}+\widecheck{\DD\c\ov{\Jk}}\right)+\frac{aq^2}{|q|^4}\DD(\ov{q})\c\ov{\Jk}-\DD\c\ov{\Zc}+\frac{2a^2(\sin\th)^2}{|q|^4}\\
&&+\frac{a\ov{q}}{|q|^2}\Jk\c\ov{\Zc}+\frac{aq}{|q|^2}\ov{\Jk}\c\Zc+\Ga_g\c\Ga_g\\
&=& \frac{4aiq(r^2+a^2)\cos\th}{|q|^6}+\frac{2a^2(\sin\th)^2}{|q|^4} -\frac{aq}{|q|^2}\widecheck{\DD\c\ov{\Jk}} -\frac{ia^2q^2}{|q|^4}\left(i\Jk+\widecheck{\DD(\cos\th)}\right)\c\ov{\Jk}-\DD\c\ov{\Zc}\\
&&+\frac{a\ov{q}}{|q|^2}\Jk\c\ov{\Zc}+\frac{aq}{|q|^2}\ov{\Jk}\c\Zc+\Ga_g\c\Ga_g\\
&=& \frac{4aiq(r^2+a^2)\cos\th}{|q|^6}+\frac{2a^2(\sin\th)^2}{|q|^4} +\frac{a^2q^2}{|q|^4}\frac{2(\sin\th)^2}{|q|^2}\\
&&-\frac{aq}{|q|^2}\widecheck{\DD\c\ov{\Jk}} -\frac{ia^2q^2}{|q|^4}\widecheck{\DD(\cos\th)}\c\ov{\Jk}-\DD\c\ov{\Zc}+\frac{a\ov{q}}{|q|^2}\Jk\c\ov{\Zc}+\frac{aq}{|q|^2}\ov{\Jk}\c\Zc+\Ga_g\c\Ga_g.
\eeaa
This yields
\beaa
\nab_4\widecheck{\tr\Xb} +\frac{1}{q}\trXbc &=&  -\frac{aq}{|q|^2}\widecheck{\DD\c\ov{\Jk}} -\frac{ia^2q^2}{|q|^4}\widecheck{\DD(\cos\th)}\c\ov{\Jk}-\DD\c\ov{\Zc}+\frac{a\ov{q}}{|q|^2}\Jk\c\ov{\Zc}+\frac{aq}{|q|^2}\ov{\Jk}\c\Zc\\
&&+2\ov{\Pc }+O(r^{-1})\trXc+\Ga_b\c\Ga_g
\eeaa
and hence
\beaa
\nab_4\widecheck{\tr\Xb} +\frac{1}{q}\trXbc &=& -\DD\c\ov{\Zc}+2\ov{\Pc }+O(r^{-2})\Zc+O(r^{-1})\trXc\\
&&+ O(r^{-1})\widecheck{\DD\c\ov{\Jk}} +O(r^{-3})\widecheck{\DD(\cos\th)}+\Ga_b\c\Ga_g.
\eeaa

Next, recall that 
 \beaa
\nab_3Z +\frac{1}{2}\tr\Xb(Z+H)-2\omb(Z-H) &=& -2\DD\omb -\frac{1}{2}\widehat{\Xb}\c(\ov{Z}+\ov{H})+\frac{1}{2}\tr X\Xib -\Bb+\frac{1}{2}\ov{\Xib}\c\Xh,\\
\nab_3Z +\nab_4\Xib &=&  -\frac{1}{2}\ov{\tr\Xb}(Z+H) -\frac{1}{2}\Xbh\c(\ov{Z}+\ov{H}) -\Bb.
\eeaa
We infer
 \beaa
\nab_4\Xib  +\frac{1}{2}\tr X\Xib &=& 2\DD\omb  +i\Im(\tr\Xb)(Z+H)-2\omb(Z-H)+\Ga_b\c\Ga_g.
\eeaa
This yields
 \beaa
&&\nab_4\Xib  +\frac{1}{q}\Xib\\
 &=& \DD\left(\pr_r\left(\frac{\De}{|q|^2}\right)\right) +i\Im\left(-\frac{2q\Delta}{|q|^4}\right)\frac{a(q+\ov{q})}{|q|^2}\Jk -\pr_r\left(\frac{\De}{|q|^2}\right)\frac{a(\ov{q}-q)}{|q|^2}\Jk\\
&& +O(r^{-1})\dkb^{\leq 1}(\ombc)+O(r^{-2})\Zc+O(r^{-2})\Hc+O(r^{-2})\trXbc+\Ga_b\c\Big(\ombc,\Ga_g\Big).
\eeaa
Now, using
\beaa
\DD(r)=0,\quad \Im(q)=a\cos\th, \quad \DD(|q|^2)=2a^2\cos\th\DD(\cos\th),\quad \DD(\cos\th)=-i\Jk+\widecheck{\DD(\cos\th)},
\eeaa
we have by direct check
\beaa
 \DD\left(\pr_r\left(\frac{\De}{|q|^2}\right)\right) +i\Im\left(-\frac{2q\Delta}{|q|^4}\right)\frac{a(q+\ov{q})}{|q|^2}\Jk -\pr_r\left(\frac{\De}{|q|^2}\right)\frac{a(\ov{q}-q)}{|q|^2}\Jk &=&  O(r^{-3})\widecheck{\DD(\cos\th)}
\eeaa
and hence
 \beaa
\nab_4\Xib  +\frac{1}{q}\Xib &=& O(r^{-1})\dkb^{\leq 1}(\ombc)+O(r^{-2})\Zc+O(r^{-2})\Hc+O(r^{-2})\trXbc\\
&&+O(r^{-3})\widecheck{\DD(\cos\th)}+\Ga_b\c\Big(\ombc,\Ga_g\Big).
\eeaa

Next, recall that 
\beaa
\nab_4\Xbh +\frac{1}{2}\tr X \Xbh &=& -\frac{1}{2}\DD\hot Z+\frac{1}{2}Z\hot Z -\frac{1}{2}\ov{\tr\Xb}\Xh
\eeaa
and hence
\beaa
\nab_4\Xbh +\frac{1}{q} \Xbh &=& -\frac{1}{2}\DD\hot Z+\frac{1}{2}Z\hot Z +O(r^{-1})\Xh+\Ga_b\c\Ga_g.
\eeaa
Now, we have
\beaa
-\DD\hot Z+Z\hot Z &=& -\DD\hot\left(\frac{a\ov{q}}{|q|^2}\Jk+\Zc\right)+\left(\frac{a\ov{q}}{|q|^2}\Jk+\Zc\right)\hot\left(\frac{a\ov{q}}{|q|^2}\Jk+\Zc\right)\\
&=& -\frac{a\ov{q}}{|q|^2}\DD\hot\Jk -\DD\left(\frac{a\ov{q}}{|q|^2}\right)\hot\Jk -\DD\hot\Zc +\frac{a^2\ov{q}^2}{|q|^4}\Jk\hot\Jk+\frac{2a\ov{q}}{|q|^2}\Jk\hot\Zc+\Ga_g\c\Ga_g\\
&=&\frac{a^2\ov{q}^2}{|q|^4}\Jk\hot\Jk +\frac{a}{q^2}\DD(q)\hot\Jk -\DD\hot\Zc +O(r^{-2})\Zc+O(r^{-1})\DD\hot\Jk+\Ga_g\c\Ga_g\\
&=& \frac{a^2i\ov{q}^2}{|q|^4}\Big(\DD(\cos\th)-i\Jk\Big)\hot\Jk -\DD\hot\Zc +O(r^{-2})\Zc+O(r^{-1})\DD\hot\Jk+\Ga_g\c\Ga_g\\
&=&  -\DD\hot\Zc +O(r^{-2})\Zc+O(r^{-1})\DD\hot\Jk+O(r^{-3})\widecheck{\DD(\cos\th)}+\Ga_g\c\Ga_g.
\eeaa
We infer
\beaa
\nab_4\Xbh +\frac{1}{q} \Xbh &=& -\frac{1}{2}\DD\hot\Zc +O(r^{-2})\Zc+O(r^{-1})\Xh+O(r^{-1})\DD\hot\Jk+O(r^{-3})\widecheck{\DD(\cos\th)}+\Ga_b\c\Ga_g.
\eeaa

Next, recall that
\beaa
\nab_4\omb  -(2\eta+\ze)\c\ze &=&   \rho,
\eeaa
which we rewrite 
\beaa
\nab_4\omb  -\frac{1}{2}\Re((2H+Z)\c\ov{Z}) &=&   \Re(P).
\eeaa
We infer
\beaa
\nab_4(\ombc) &=&  \frac{1}{2}\Re((2H+Z)\c\ov{Z}) +   \Re(P) -\frac{1}{2}\nab_4\left(\pr_r\left(\frac{\De}{|q|^2}\right)\right)\\
&=& \frac{1}{2}\Re\left(\frac{a(2q+\ov{q})}{|q|^2}\Jk\c\frac{aq}{|q|^2}\ov{\Jk}\right) +   \Re\left(-\frac{2m}{q^3}\right) -\frac{1}{2}\pr_r^2\left(\frac{\De}{|q|^2}\right)\\
&&+\Re(\Pc)+O(r^{-2})\Zc+O(r^{-2})\Hc+\Ga_b\c\Ga_g.
\eeaa
Now, using $\Jk\c\ov{\Jk}=\frac{2(\sin\th)^2}{|q|^2}$, we have
\beaa
\frac{1}{2}\Re\left(\frac{a^2q(2q+\ov{q})}{|q|^4}\Jk\c\ov{\Jk}\right) +   \Re\left(-\frac{2m}{q^3}\right) -\frac{1}{2}\pr_r^2\left(\frac{\De}{|q|^2}\right) = 0
\eeaa
and hence
\beaa
\nab_4(\ombc) &=& \Re(\Pc)+O(r^{-2})\Zc+O(r^{-2})\Hc+\Ga_b\c\Ga_g.
\eeaa

Next, recall
\beaa
\frac{1}{2}\ov{\DD}\c\Xh +\frac{1}{2}\Xh\c\ov{Z} &=& \frac{1}{2}\DD\ov{\tr X}+\frac{1}{2}\ov{\tr X}Z-i\Im(\tr X)H -B.
\eeaa
We infer
\beaa
\frac{1}{2}\ov{\DD}\c\Xh +O(r^{-2})\Xh &=& \DD\left(\frac{1}{\ov{q}}\right)+\frac{1}{\ov{q}}\frac{a\ov{q}}{|q|^2}\Jk+\frac{1}{\ov{q}}\Zc -i\Im\left(\frac{2}{q}\right)\frac{aq}{|q|^2}\Jk+O(r^{-2})\Hc\\
&&+O(r^{-1})\dkb^{\leq 1}\trXc -B+\Ga_b\c\Ga_g.
\eeaa
Since we have
\beaa
\DD\left(\frac{1}{\ov{q}}\right)+\frac{1}{\ov{q}}\frac{a\ov{q}}{|q|^2}\Jk -i\Im\left(\frac{2}{q}\right)\frac{aq}{|q|^2}\Jk &=& -
\frac{1}{\ov{q}^2}\DD(\ov{q})+\frac{a}{|q|^2}\Jk +\frac{2ia^2\cos\th q}{|q|^4}\Jk\\
&=&  \frac{iaq^2}{|q|^4}\left(\DD(\cos\th)-i\Jk\right) = \frac{iaq^2}{|q|^4}\widecheck{\DD(\cos\th)},
\eeaa
we deduce
\beaa
\frac{1}{2}\ov{\DD}\c\Xh &=& \frac{1}{\ov{q}}\Zc -B +O(r^{-2})\Xh +O(r^{-2})\Hc +O(r^{-1})\dkb^{\leq 1}\trXc\\
&& +O(r^{-2})\widecheck{\DD(\cos\th)}+\Ga_b\c\Ga_g.
\eeaa

Next, recall
 \beaa
\nab_3B-\DD\ov{P} &=& -\tr\Xb B+2\omb B+\ov{\Bb}\c \Xh+3\ov{P}H +\frac{1}{2}A\c\ov{\Xib}.
\eeaa
We infer
 \beaa
\nab_3B -\DD\ov{\Pc} &=& \frac{2}{r}B+O(r^{-2}) B +O(r^{-2})\Pc+O(r^{-3})\Hc -\DD\left(\frac{2m}{(\ov{q})^3}\right) -\frac{6m}{(\ov{q})^3}\frac{aq}{|q|^2}\Jk+r^{-1}\Ga_b\c\Ga_g.
\eeaa
Now, we have
\beaa
-\DD\left(\frac{2m}{(\ov{q})^3}\right) -\frac{6m}{(\ov{q})^3}\frac{aq}{|q|^2}\Jk &=& \frac{6m}{\ov{q}^4}\DD(\ov{q}) -\frac{6am}{(\ov{q})^4}\Jk = -\frac{6iam}{\ov{q}^4}\DD(\cos\th) -\frac{6am}{(\ov{q})^4}\Jk\\
&=& -\frac{6iam}{\ov{q}^4}\Big(\DD(\cos\th) -i\Jk\Big)=-\frac{6iam}{\ov{q}^4}\widecheck{\DD(\cos\th)}
\eeaa
and hence
\beaa
\nab_3B -\DD\ov{\Pc} &=& \frac{2}{r}B+O(r^{-2}) B +O(r^{-2})\Pc+O(r^{-3})\Hc +O(r^{-4})\widecheck{\DD(\cos\th)}+r^{-1}\Ga_b\c\Ga_g.
\eeaa

Next, recall
\beaa
\nab_4B -\frac{1}{2}\ov{\DD}\c A &=& -2\ov{\tr X} B  +\frac{1}{2}A\c  \ov{Z},
\eeaa
and hence
\beaa
\nab_4B +\frac{4}{\ov{q}} B &=&  \frac{1}{2}\ov{\DD}\c A +\frac{aq}{2|q|^2}\ov{\Jk}\c A+\Ga_g\c(B,A).
\eeaa

Next, recall
\beaa
\nab_4P -\frac{1}{2}\DD\c \ov{B} &=& -\frac{3}{2}\tr X P -\frac{1}{2}Z\c\ov{B}  -\frac{1}{4}\Xbh\c \ov{A}.
\eeaa
We infer
\beaa
\nab_4\left(-\frac{2m}{q^3}+\Pc \right)-\frac{1}{2}\DD\c \ov{B} &=& \frac{6m\ov{q}}{|q|^2q^3} -\frac{3}{q} \Pc -\frac{a\ov{q}}{2|q|^2}\Jk\c\ov{B}+O(r^{-3})\trXc+r^{-1}\Ga_g\c\Ga_g+\Ga_b\c A.
\eeaa
Now, we have 
\beaa
\nab_4\left(-\frac{2m}{q^3} \right) = \pr_r\left(-\frac{2m}{(r+ia\cos\th)^4} \right) = \frac{6m}{(r+ia\cos\th)^4}=\frac{6m\ov{q}}{|q|^2q^3} 
\eeaa
and hence
\beaa
\nab_4\left(\Pc \right)-\frac{1}{2}\DD\c \ov{B} &=& -\frac{3}{q} \Pc -\frac{a\ov{q}}{2|q|^2}\Jk\c\ov{B}+O(r^{-3})\trXc+r^{-1}\Ga_g\c\Ga_g+\Ga_b\c A.
\eeaa

Finally, recall
 \beaa
\nab_4\Bb+\DD P &=& -\tr X\Bb+\ov{B}\c \Xbh+3P Z.
\eeaa
We infer
 \beaa
\nab_4\Bb+\DD\left(\Pc\right) &=& -\frac{2}{q}\Bb +O(r^{-2})\Pc+O(r^{-3})\Zc +\DD\left(\frac{2m}{q^3}\right) -\frac{6m}{q^3}\frac{a\ov{q}}{|q|^2}\Jk+r^{-1}\Ga_b\c\Ga_g.
\eeaa
Now, we have
\beaa
\DD\left(\frac{2m}{q^3}\right) -\frac{6m}{q^3}\frac{a\ov{q}}{|q|^2}\Jk &=& -\frac{6m}{q^4}\DD(q) -\frac{6am}{q^4}\Jk=-\frac{6iam}{q^4}\DD(\cos\th) -\frac{6am}{q^4}\Jk\\
&=& -\frac{6iam}{q^4}\Big(\DD(\cos\th) -i\Jk\Big) = -\frac{6iam}{q^4}\widecheck{\DD(\cos\th)}
\eeaa
and hence 
 \beaa
\nab_4\Bb+\DD\left(\Pc\right) &=& -\frac{2}{q}\Bb +O(r^{-2})\Pc+O(r^{-3})\Zc +O(r^{-4})\widecheck{\DD(\cos\th)}+r^{-1}\Ga_b\c\Ga_g.
\eeaa
This concludes the proof of Lemma \ref{Lemma:linearized-nullstr}.

%%%%%%%%%%%%%%%%%%%%%%%%%%%%%%%%

\section{Proof of  Lemma \ref{Lemma:otherlinearizedquant}}
\lab{section:Proof-lemmaLemma:otherlinearizedquant}

%%%%%%%%%%%%%%%%%%%%%%%%%%%%%%%%%

We have from Lemma    \ref{prop:e_4(xyz)-again}.
\beaa
e_4(e_3(r)) &=&  -2\omb,
\eeaa
and hence, since $e_4(r)=1$, $e_4(\th)=0$, 
\beaa
e_4\left(\widecheck{e_3(r)}\right) &=& e_4\left(e_3(r)+\frac{\Delta}{|q|^2}\right) =  -2\omb +\pr_r\left(\frac{\Delta}{|q|^2}\right)\\
&=&-2\ombc.
\eeaa

Also, we have from Proposition \ref{prop:e_4(xyz)-again}.
\beaa
\nab_4\DD u +\frac{1}{2}\tr X\DD u &=& -\frac{1}{2}\Xh\c\ov{\DD}u
\eeaa
and hence
\beaa
\nab_4\DD u +\frac{1}{q}\DD u &=& O(r^{-1})\trXc+O(r^{-1})\Xh+\Ga_b\c\Ga_g.
\eeaa
Since 
\beaa
\nab_4\DD u +\frac{1}{q}\DD u &=& \nab_4\left(a\Jk+\widecheck{\DD u}\right)+\frac{1}{q}\left(a\Jk+\widecheck{\DD u}\right) = \nab_4\widecheck{\DD u}+\frac{1}{q}\widecheck{\DD u},
\eeaa
where we have used $\nab_4\Jk=-q^{-1}\Jk$, we infer
\beaa
\nab_4\widecheck{\DD u}+\frac{1}{q}\widecheck{\DD u} &=& O(r^{-1})\trXc+O(r^{-1})\Xh+\Ga_b\c\Ga_g.
\eeaa

Also, recall that  we have
\beaa
e_4(e_3(u)) &=& -\Re\Big((Z+H)\c\ov{\DD} u\Big).
\eeaa
We infer
\beaa
e_4\left(\frac{2(r^2+a^2)}{|q|^2}+\widecheck{e_3(u)}\right) &=& -\Re\left(\frac{a(q+\ov{q})}{|q|^2}\Jk\c a\ov{\Jk}\right)+O(r^{-1})\Hc+O(r^{-1})\Zc+O(r^{-2})\widecheck{\DD u}\\
&&+\Ga_b\c\Ga_b.
\eeaa
Now, using $e_4(r)=1$, $e_4(\th)=0$ and $\Jk\c\ov{\Jk}=\frac{2(\sin\th)^2}{|q|^2}$, we have
\beaa
-e_4\left(\frac{2(r^2+a^2)}{|q|^2}\right)-\Re\left(\frac{a(q+\ov{q})}{|q|^2}\Jk\c a\ov{\Jk}\right) &=& -\pr_r\left(\frac{2(r^2+a^2)}{|q|^2}\right)-\Re\left(\frac{2ra^2}{|q|^2}\Jk\c\ov{\Jk}\right)\\
&=& -\frac{4r}{|q|^2}+\frac{4r(r^2+a^2)}{|q|^4} -\frac{4ra^2(\sin\th)^2}{|q|^4}=0
\eeaa
and hence
\beaa
e_4\left(\widecheck{e_3(u)}\right) &=& O(r^{-1})\Hc+O(r^{-1})\Zc+O(r^{-2})\widecheck{\DD u}+\Ga_b\c\Ga_b.
\eeaa

Also, recall that we have 
\beaa
\nab_4\DD\cos\th +\frac{1}{2}\tr X\DD\cos\th &=& -\frac{1}{2}\Xh\c\ov{\DD}\cos\th.
\eeaa
and hence
\beaa
\nab_4\DD\cos\th +\frac{1}{q}\DD\cos\th &=& \frac{i}{2}\ov{\Jk}\c\Xh+O(r^{-1})\trXc+\Ga_b\c\Ga_g.
\eeaa
Since 
\beaa
\nab_4\DD\cos\th +\frac{1}{q}\DD\cos\th &=& \nab_4\left(i\Jk+\widecheck{\DD\cos\th}\right)+\frac{1}{q}\left(i\Jk+\widecheck{\DD\cos\th}\right) = \nab_4\widecheck{\DD\cos\th}+\frac{1}{q}\widecheck{\DD\cos\th},
\eeaa
where we have used $\nab_4\Jk=-q^{-1}\Jk$, we infer
\beaa
\nab_4\widecheck{\DD\cos\th}+\frac{1}{q}\widecheck{\DD\cos\th} &=& \frac{i}{2}\ov{\Jk}\c\Xh+O(r^{-1})\trXc+\Ga_b\c\Ga_g.
\eeaa

Also, recall that  we have
\beaa
e_4(e_3(\cos\th)) &=& -\Re\Big((Z+H)\c\ov{\DD} \cos\th\Big)
\eeaa
and hence
\beaa
e_4(e_3(\cos\th)) &=& \Re\left(\frac{a(q+\ov{q})}{|q|^2}\Jk\c i\ov{\Jk}\right)
+O(r^{-1})\Hc+O(r^{-1})\Zc+O(r^{-2})\widecheck{\DD\cos\th}+\Ga_b\c\Ga_b\\
&=& \Re\left(i\frac{2ra}{|q|^2}\frac{2(\sin\th)^2}{|q|^2}\right)
+O(r^{-1})\Hc+O(r^{-1})\Zc+O(r^{-2})\widecheck{\DD\cos\th}+\Ga_b\c\Ga_b\\
&=& O(r^{-1})\Hc+O(r^{-1})\Zc+O(r^{-2})\widecheck{\DD\cos\th}+\Ga_b\c\Ga_b.
\eeaa
This concludes the proof of  Lemma \ref{Lemma:otherlinearizedquant}.

%%%%%%%%%%%%%%%%%%%%%%%%%%%%%%%%%%%%%%%

\section{Proof of Lemma \ref{lemma:Comparison.horizontalderivarives}}
\lab{appendix:ProofLemma-lemma:Comparison.horizontalderivarives}

%%%%%%%%%%%%%%%%%%%%%%%%%%%%%%%%%%%%%%%

In view of Lemma \ref{Lemma:Generalframetransf}, we have
\beaa
e_a' &=& e_a+\frac{1}{2}\fb_a\la^{-1}e_4'+\frac{1}{2}f_ae_3.
\eeaa
We compute
\beaa
&&\g(\D_{e_a'}e_b', e_c') = \g\left(\D_{e_a'}\left(e_b+\frac{1}{2}\fb_b\la^{-1}e_4'+\frac{1}{2}f_be_3\right),e_c'\right)\\
&=& \frac{1}{2}e_a'(f_b)\g(e_3, e_c')+\g(\D_{e_a'}e_b, e_c')+\frac{1}{2}\fb_b\la^{-1}\g(\D_{e_a'}e_4', e_c')+\frac{1}{2}f_b\g(\D_{e_a'}e_3, e_c')\\
&=& -\frac{1}{2}e_a'(f_b)\fb_c+\g(\D_{e_a'}e_b, e_c)+ \frac{1}{2}\fb_c\la^{-1}\g(\D_{e_a'}e_b,e_4')+\frac{1}{2}f_c\g(\D_{e_a'}e_b, e_3)\\
&&+\frac{1}{2}\fb_b\la^{-1}\chi_{ac}'+\frac{1}{2}f_b\g\left(\D_{e_a'}e_3, \left(\de_c^d+\frac{1}{2}\fb_cf^d\right) e_d+\frac{1}{2}\fb_c e_4\right)\\
&=& -\frac{1}{2}e_a'(f_b)\fb_c+\g\left(\D_{\left(\de_a^d+\frac{1}{2}\fb_af^d\right) e_d+\frac{1}{2}\fb_a e_4+\left(\frac{1}{2}f_a+\frac{1}{8}|f|^2\fb_a\right)e_3}e_b, e_c\right)\\
&&+ \frac{1}{2}\fb_c\la^{-1} e_a'(\g(e_b,e_4'))  - \frac{1}{2}\fb_c\la^{-1}\g(e_b, \D_{e_a'}e_4')\\
&&+\frac{1}{2}f_c\g\left(\D_{\left(\de_a^d+\frac{1}{2}\fb_af^d\right) e_d+\frac{1}{2}\fb_a e_4+\left(\frac{1}{2}f_a+\frac{1}{8}|f|^2\fb_a\right)e_3}e_b, e_3\right)+\frac{1}{2}\fb_b\la^{-1}\chi_{ac}'\\
&&+\frac{1}{2}f_b\g\left(\D_{\left(\de_a^p+\frac{1}{2}\fb_af^p\right) e_p+\frac{1}{2}\fb_a e_4+\left(\frac{1}{2}f_a+\frac{1}{8}|f|^2\fb_a\right)e_3}e_3, \left(\de_c^d+\frac{1}{2}\fb_cf^d\right) e_d+\frac{1}{2}\fb_c e_4\right).
\eeaa
We further deduce
\beaa
&&\g(\D_{e_a'}e_b', e_c') \\
&=& \left(\de_a^d+\frac{1}{2}\fb_af^d\right)\g\left(\D_{e_d}e_b, e_c\right)+\frac{1}{2}\fb_a\g\left(\D_{ e_4}e_b, e_c\right)+\left(\frac{1}{2}f_a+\frac{1}{8}|f|^2\fb_a\right)\g\left(\D_{e_3}e_b, e_c\right) +\frac{1}{2}\fb_cf_be_a'(\log\la)\\
&&  - \frac{1}{2}\fb_c\la^{-1}\g\Big(e_b, \chi_{ad}'e_d' -\ze_a'e_4'\Big) -\frac{1}{2}f_c\g\left(\D_{\left(\de_a^d+\frac{1}{2}\fb_af^d\right) e_d+\frac{1}{2}\fb_a e_4+\left(\frac{1}{2}f_a+\frac{1}{8}|f|^2\fb_a\right)e_3}e_3, e_b\right)\\
&&+\frac{1}{2}\fb_b\la^{-1}\chi_{ac}'+\frac{1}{2}f_b\left(\de_c^d+\frac{1}{2}\fb_cf^d\right) \left(\left(\de_a^p+\frac{1}{2}\fb_af^p\right)\chib_{pd}+\fb_a\etab_d+\left(f_a+\frac{1}{4}|f|^2\fb_a\right)\xib_d\right)\\
&&+\frac{1}{4}f_b\fb_c\left(-2\left(\de_a^p+\frac{1}{2}\fb_af^p\right)\ze_p -2\om\fb_a+\omb\left(2f_a+\frac{1}{2}|f|^2\fb_a\right)\right)
\eeaa
and hence
\beaa
&&\g(\D_{e_a'}e_b', e_c') \\
&= & \left(\de_a^d+\frac{1}{2}\fb_af^d\right)\g\left(\D_{e_d}e_b, e_c\right)+\frac{1}{2}\fb_a\g\left(\D_{ e_4}e_b, e_c\right)+\left(\frac{1}{2}f_a+\frac{1}{8}|f|^2\fb_a\right)\g\left(\D_{e_3}e_b, e_c\right)\\
&& +\frac{1}{2}\fb_cf_be_a'(\log\la) - \frac{1}{2}\fb_c\la^{-1}\chi_{ab}' +\frac{1}{2}\fb_b\la^{-1}\chi_{ac}'  + \frac{1}{2}\la^{-1}\ze_a'\fb_cf_b - \frac{1}{4}\fb_c\la^{-1}\chi_{ad}' \fb_bf_d\\
&&-\frac{1}{2}f_c\left(\de_a^d+\frac{1}{2}\fb_af^d\right)\chib_{db} -\frac{1}{2}f_c\fb_a\etab_b-\frac{1}{2}f_c\left(f_a+\frac{1}{4}|f|^2\fb_a\right)\xib_b\\
&&+\frac{1}{2}f_b\left(\de_c^d+\frac{1}{2}\fb_cf^d\right) \left(\left(\de_a^p+\frac{1}{2}\fb_af^p\right)\chib_{pd}+\fb_a\etab_d+\left(f_a+\frac{1}{4}|f|^2\fb_a\right)\xib_d\right)\\
&&+\frac{1}{4}f_b\fb_c\left(-2\left(\de_a^p+\frac{1}{2}\fb_af^p\right)\ze_p -2\om\fb_a+\omb\left(2f_a+\frac{1}{2}|f|^2\fb_a\right)\right).
\eeaa
We infer
\beaa
\g(\D_{e_a'}e_b', e_c') &=& \left(\de_a^d+\frac{1}{2}\fb_af^d\right)\g\left(\D_{e_d}e_b, e_c\right)+\frac{1}{2}\fb_a\g\left(\D_{ e_4}e_b, e_c\right)\\
&&+\left(\frac{1}{2}f_a+\frac{1}{8}|f|^2\fb_a\right)\g\left(\D_{e_3}e_b, e_c\right) +\frac{1}{2}\fb_cf_be_a'(\log\la)
  - \frac{1}{2}\fb_c\la^{-1}\chi_{ab}' \\
  &&+\frac{1}{2}\fb_b\la^{-1}\chi_{ac}'  + \frac{1}{2}\ze_a'\fb_cf_b - \frac{1}{4}\fb_c\la^{-1}\chi_{ad}' \fb_bf_d -\frac{1}{2}f_c\chib_{ab}+\frac{1}{2}f_b\chib_{ac}+\err[\g(\D_{e_a'}e_b', e_c')],
\eeaa
with 
\beaa
&& \err[\g(\D_{e_a'}e_b', e_c')]\\
&=& -\frac{1}{2}f_c\left(\de_a^d+\frac{1}{2}\fb_af^d\right)\chib_{db} -\frac{1}{2}f_c\fb_a\etab_b-\frac{1}{2}f_c\left(f_a+\frac{1}{4}|f|^2\fb_a\right)\xib_b\\
&&+\frac{1}{2}f_b\left(\de_c^d+\frac{1}{2}\fb_cf^d\right) \left(\left(\de_a^p+\frac{1}{2}\fb_af^p\right)\chib_{pd}+\fb_a\etab_d+\left(f_a+\frac{1}{4}|f|^2\fb_a\right)\xib_d\right)\\
&&+\frac{1}{4}f_b\fb_c\left(-2\left(\de_a^p+\frac{1}{2}\fb_af^p\right)\ze_p -2\om\fb_a+\omb\left(2f_a+\frac{1}{2}|f|^2\fb_a\right)\right).
\eeaa
where $\err[\g(\D_{e_a'}e_b', e_c')]$ contains all the terms depending on $(f, \fb, \Ga)$, without derivative, and at least quadratic in $(f, \fb)$. This concludes the proof of Lemma \ref{lemma:Comparison.horizontalderivarives}.

%%%%%%%%%%%%%%%%%%%%%%%%%%%%%%%%%%%%%%%%%

\section{Proof of Proposition \ref{prop:controlofDDprimehotDDprimeJonMext}}
\lab{sec:roofofprop:controlofDDprimehotDDprimeJonMext}

%%%%%%%%%%%%%%%%%%%%%%%%%%%%%%%%%%%%%%%%

In order to prove Proposition \ref{prop:controlofDDprimehotDDprimeJonMext}, we start with the following  lemma.
\begin{lemma}
For all $J=\Jp$, we have
\beaa
\nab_4(\DD\hot\DD J)+\frac{2}{q}\DD\hot\DD J &=& r^{-2}\dk^{\leq 1}\Ga_g
\eeaa
and
\beaa
\nab_4(\ov{\DD}\c\DD J)+\frac{2}{r}\ov{\DD}\c\DD J &=& O(r^{-5})+r^{-2}\dk^{\leq 1}\Ga_g.
\eeaa
\end{lemma}

\begin{proof}
We start with the first identity. We use  $e_4(J)=0$ and the commutation Lemma \ref{lemma:commutation-complexM6} to deduce
\beaa
\nab_4(\DD\hot\DD J) &=& [\nab_4, \DD\hot]\DD J+\DD\hot[\nab_4,\DD] J\\
&=& -\frac{1}{2}\tr X\Big(\DD\hot\DD J -Z\hot\DD J\Big)+r^{-1}\Ga_g\c\dk^{\leq 1}\DD J\\
&&+\DD\hot\left(-\frac{1}{2}\tr X\DD J+r^{-1}\Ga_g\c\dkb J\right)\\
&=&  -\tr X\DD\hot\DD J -\frac{1}{2}\Big(\DD\tr X -\tr X Z\Big)\hot\DD J+r^{-2}\dk^{\leq 1}\Ga_g\\
&=&  -\frac{2}{q}\DD\hot\DD J -\frac{1}{2}\left(\DD\left(\frac{2}{q}\right) -\frac{2}{q}\frac{a\ov{q}}{|q|^2}\Jk\right)\hot\DD J+r^{-2}\dk^{\leq 1}\Ga_g\\
&=&  -\frac{2}{q}\DD\hot\DD J -\frac{1}{2}\left(-\frac{2ia}{q^2}\DD(\cos\th) -\frac{2}{q}\frac{a\ov{q}}{|q|^2}\Jk\right)\hot\DD J+r^{-2}\dk^{\leq 1}\Ga_g\\
&=&  -\frac{2}{q}\DD\hot\DD J +\frac{ia}{q^2}\Big(\DD(\cos\th) -i\Jk\Big)\hot\DD J+r^{-2}\dk^{\leq 1}\Ga_g\\
&=&  -\frac{2}{q}\DD\hot\DD J +O(r^{-3})\widecheck{\DD(\cos\th)}+r^{-2}\dk^{\leq 1}\Ga_g\\
&=&  -\frac{2}{q}\DD\hot\DD J +r^{-2}\dk^{\leq 1}\Ga_g
\eeaa
and hence
\beaa
\nab_4(\DD\hot\DD J)+\frac{2}{q}\DD\hot\DD J &=& r^{-2}\dk^{\leq 1}\Ga_g
\eeaa
which is the first identity.

Similarly, to prove the second identity, we use  $e_4(J)=0$ and the commutation Lemma \ref{lemma:commutation-complexM6} to deduce
\beaa
\nab_4(\ov{\DD}\c\DD J) &=& [\nab_4, \ov{\DD}\c]\DD J+\ov{\DD}\c[\nab_4,\DD] J\\
&=& -\frac{1}{2}\ov{\tr X}\Big(\ov{\DD}\c\DD J +\ov{Z}\c\DD J\Big)+r^{-1}\Ga_g\c\dk^{\leq 1}\DD J\\
&&+\ov{\DD}\c\left(-\frac{1}{2}\tr X\DD J+r^{-1}\Ga_g\c\dkb J\right)\\
&=&  -\frac{1}{2}(\tr X+\ov{\tr X})\ov{\DD}\c\DD J -\frac{1}{2}\Big(\ov{\DD}\tr X +\ov{\tr X}\,\ov{Z}\Big)\c\DD J+r^{-2}\dk^{\leq 1}\Ga_g\\
&=&  -\frac{2r}{|q|^2}\ov{\DD}\c\DD J -\frac{1}{2}\left(\ov{\DD}\left(\frac{2}{q}\right) -\frac{2}{\ov{q}}\frac{aq}{|q|^2}\ov{\Jk}\right)\c\DD J+r^{-2}\dk^{\leq 1}\Ga_g\\
&=&  -\frac{2r}{|q|^2}\ov{\DD}\c\DD J -\frac{1}{2}\left(-\frac{2ia}{q^2}\ov{\DD}(\cos\th) -\frac{2}{\ov{q}}\frac{aq}{|q|^2}\ov{\Jk}\right)\c\DD J+r^{-2}\dk^{\leq 1}\Ga_g\\
&=&  -\frac{2r}{|q|^2}\ov{\DD}\c\DD J +a\left(\frac{i}{q^2}\ov{\DD(\cos\th) -i\Jk}+\left(\frac{1}{(\ov{q})^2}-\frac{1}{q^2}\right)\ov{\Jk}\right)\c\DD J+r^{-2}\dk^{\leq 1}\Ga_g\\
&=&  -\frac{2}{r}\ov{\DD}\c\DD J +O(r^{-3})\widecheck{\DD(\cos\th)}+O(r^{-5})+r^{-2}\dk^{\leq 1}\Ga_g\\
&=&  -\frac{2}{r}\ov{\DD}\c\DD J +O(r^{-5})+r^{-2}\dk^{\leq 1}\Ga_g\\
\eeaa
and hence
\beaa
\nab_4(\ov{\DD}\c\DD J) +\frac{2}{r}\ov{\DD}\c\DD J &=& O(r^{-5})+r^{-2}\dk^{\leq 1}\Ga_g
\eeaa
as desired.
\end{proof}

\begin{corollary}
For all $J=\Jp$, we have
\beaa
\sup_{\Mext}r^3u^{\frac{1}{2}+\dec}|\DD\hot\DD J|\les \ep
\eeaa
and 
\beaa
\left|\left(\De+\frac{2}{r^2}\right)J\right| &\les& \frac{1}{r^4}+\frac{\ep}{r^3u^{\frac{1}{2}+\dec}}\quad\textrm{on}\quad\Mext.
\eeaa
\end{corollary}

\begin{proof}
Recall that we have
\beaa
\nab_4(\DD\hot\DD J)+\frac{2}{q}\DD\hot\DD J &=& r^{-2}\dk^{\leq 1}\Ga_g
\eeaa
so that 
\beaa
\nab_4(q^2\DD\hot\DD J) &=& \dk^{\leq 1}\Ga_g.
\eeaa
Integrating from $\Si_*$, and together with the control on $\Si_*$ of Lemma \ref{lemma:statementeq:DeJp.Sigmastar:improvedd}, we infer
\beaa
\sup_{\Mext}r^3u^{\frac{1}{2}+\dec}|\DD\hot\DD J|\les \ep
\eeaa
as desired.

Also, recall that we have
\beaa
\nab_4(\ov{\DD}\c\DD J)+\frac{2}{r}\ov{\DD}\c\DD J &=& O(r^{-5})+r^{-2}\dk^{\leq 1}\Ga_g
\eeaa
so that 
\beaa
\nab_4(r^2\ov{\DD}\c\DD J+4) &=& O(r^{-3})+\dk^{\leq 1}\Ga_g.
\eeaa
Integrating from $\Si_*$, we infer on $\Mext$
\beaa
\left|r^2\ov{\DD}\c\DD J+4\right| &\les& \sup_{\Si_*}\left|r^2\ov{\DD}\c\DD J+4\right|+\frac{1}{r^2}+\frac{\ep}{ru^{\frac{1}{2}+\dec}}.
\eeaa
Now, we have, for a scalar function $h$ 
\beaa
\ov{\DD}\c\DD h &=& 2\De h+2i\in_{ab}\nab_a\nab_b h\\
&=& 2\De h+i\in_{ab}\Big(\D_a\D_b+\chi_{ab}e_3+\chib_{ab}e_4\Big)h\\
&=& 2\De h+i\Big(\atrch e_3(h)+\atrchb e_4(h)\Big). 
\eeaa
Using this formula with $h=\Jp$, and since $e_4(\Jp)=0$, $e_3(\Jp)=O(r^{-2})+\Ga_b$ and $\atrch=O(r^{-2})+\Ga_g$, we infer
\beaa
\ov{\DD}\c\DD J &=& 2\De J+\Big(O(r^{-2})+\Ga_g\Big)\Big(O(r^{-2})+\Ga_b\Big)\\
&=& 2\De J+O(r^{-4})+r^{-2}\Ga_b.  
\eeaa
Plugging in the above, and using the control of $\Ga_b$, we infer
\beaa
\left|r^2\De J+2\right| &\les& \sup_{\Si_*}\left|r^2\De J+2\right|+\frac{1}{r^2}+\frac{\ep}{ru^{\frac{1}{2}+\dec}}.
\eeaa
Together with the control on $\Si_*$ of Lemma \ref{lemma:statementeq:DeJp.Sigmastar:improvedd}, we deduce
\beaa
\left|r^2\De J+2\right| &\les& \frac{1}{r^2}+\frac{\ep}{ru^{\frac{1}{2}+\dec}}.
\eeaa
as desired.
\end{proof}

We are now ready to prove Proposition \ref{prop:controlofDDprimehotDDprimeJonMext}.
\begin{proof}[Proof of Proposition \ref{prop:controlofDDprimehotDDprimeJonMext}] 
Recall from Corollary \ref{cor:formuallikningnabprimeandnabintransfoformula:bis} that we have 
\beaa
\nab' &=& \Big(1+O(r^{-2})\Big)\nab+O(r^{-1})\Lieb_\T+O(r^{-1})\nab_4+O(r^{-3})+r^{-1}\Ga_b.
\eeaa
We infer
\beaa
\nab'\hot\nab' &=& \nab\hot\nab+O(r^{-4})\dkb^{\leq 2}+O(r^{-2})\dk^{\leq 1}\Lieb_\T+O(r^{-2})\dk^{\leq 1}\nab_4+O(r^{-4})+r^{-2}\dk^{\leq 1}\Ga_b,\\
\De' &=& \De+O(r^{-4})\dkb^{\leq 2}+O(r^{-2})\dk^{\leq 1}\Lieb_\T+O(r^{-2})\dk^{\leq 1}\nab_4+O(r^{-4})+r^{-2}\dk^{\leq 1}\Ga_b.
\eeaa
Since $\nab_4J=0$ and $\T(J)=\Ga_b$, see Lemma \ref{Lemma:$T(rthuJ},         we deduce
\beaa
\nab'\hot\nab'J &=& \nab\hot\nab J+O(r^{-4})+r^{-2}\dk^{\leq 1}\Ga_b,\\
\De'J &=& \De J+O(r^{-4})+r^{-2}\dk^{\leq 1}\Ga_b.
\eeaa
Together with the above estimate for $\nab\hot\nab J$ and $\De J$, and in view of the control of $\Ga_b$, we infer
\beaa
|\DD'\hot\DD' J|+|r^2\De' J+2|\les \frac{\ep}{r^3u^{\frac{1}{2}+\dec}}+\frac{1}{r^4}
\eeaa
as desired.
\end{proof}

%%%%%%%%%%%%%%%%%%%%%%%%%%%%%%%%%%%%

\section{Proof of Proposition \ref{Prop: eqforrenormalized.qiantities}}
\lab{appendix:ProofofProp{Prop: eqforrenormalized.qiantities}}

%%%%%%%%%%%%%%%%%%%%%%%%%%%%%%%%%%%%

Recall the definition of our renormalized quantities:
\beaa
\bsplit
 [\Hc]_{ren}    &=  \frac{1}{\ov{q}}\left(\ov{q}\Hc -q\Zc +\frac{1}{3}\left(-\ov{q}^2+|q|^2\right)B+\frac{a}{2}(q-\ov{q})\ov{\Jk}\c\Xh\right),\\
[\widecheck{\DD\cos\th}] _{ren} &= \frac{1}{q}\left(q\widecheck{\DD\cos\th} +\frac{i}{2}|q|^2\ov{\Jk}\c\Xh\right),\\
[\Mc]_{ren} &=\frac{1}{\ov{q} q^2}\Bigg[\ov{q}\,\ov{\DD}\c\left(q^2\Zc  +\left(-\frac{a}{2}q^2 -\frac{a}{2}|q|^2\right)\ov{\Jk}\c\Xh\right) +2\ov{q}^3\,\ov{\Pc} - 2aq^2\ov{\Jk}\c\Zc\\
&+\left(-\frac{1}{3}q^2\ov{q}^2 -\frac{1}{3}q\ov{q}^3+\frac{2}{3}\ov{q}^4\right)\ov{\DD}\c B +a\left( q^2\ov{q} +\frac{2}{3}q\ov{q}^2 - \frac{13}{6}\ov{q}^3\right)\ov{\Jk}\c B\\
&+a^2(q^2+|q|^2)\ov{\Jk}\c\Xh\c\ov{\Jk}\Bigg].
\end{split}
\eeaa
Our goal is to prove the following three identities
\beaa
\bsplit
\nab_4\left(\ov{q}[\Hc]_{ren} \right)
&=  O(r^{-1})\trXc+O(1)\dkb^{\leq 1}A+r\Ga_b\c\Ga_g,\\
\nab_4\left(q [\widecheck{\DD\cos\th}]_{ren}  \right) &= O(1)\trXc+O(r)A+r\Ga_b\c\Ga_g,\\
\nab_4\left(\ov{q} q^2 [\Mc]_{ren}\right)&=   O(1)\dkb^{\leq 1}\trXc+O(r)\dkb^{\leq 2}A+r^2\dkb^{\leq 1}(\Ga_g\c\Ga_g)+r^3\Ga_b\c A.
\end{split}
\eeaa

{\bf First identity.} To check
\beaa
\nab_4\left(\ov{q}[\Hc]_{ren} \right) &=&  O(r^{-1})\trXc+O(1)\dkb^{\leq 1}A+r\Ga_b\c\Ga_g,
\eeaa
we start with, see Lemma \ref{Lemma:linearized-nullstr},
\beaa
\nab_4\Hc+\frac{1}{\ov{q}}\Hc &=& -\frac{1}{\ov{q}}\Zc     -\frac{ar}{|q|^2}\ov{\Jk}\c\Xh -B +O(r^{-2})\trXc +\Ga_b\c\Ga_g,\\
\nab_4\Zc + \frac{2}{q}\Zc  &=&    - \frac{aq}{|q|^2}\ov{\Jk}\c\widehat{X} -B +O(r^{-2})\trXc+\Ga_g\c\Ga_g.
\eeaa
The  goal is  to  eliminate the presence of $B$ and $\Xh$  on the RHS.  We first combine the two as follows
\beaa
\nab_4(\ov{q}\Hc -q\Zc) &=& \ov{q}\nab_4\Hc+\Hc-q\nab_4\Zc-\Zc
\\
&& -q\left(- \frac{2}{q}\Zc      - \frac{aq}{|q|^2}\ov{\Jk}\c\widehat{X} -B \right) -\Zc +O(r^{-1})\trXc+r\Ga_b\c\Ga_g\\
&=& (q-\ov{q})B +\frac{a\big(q^2-r\ov{q}\big)}{|q|^2}\ov{\Jk}\c\widehat{X} +O(r^{-1})\trXc+r\Ga_b\c\Ga_g.
\eeaa
To eliminate $(q-\ov{q})B$  we make use of the equation
\beaa
\nab_4B +\frac{4}{\ov{q}} B &=&  O(r^{-1})\dkb^{\leq 1}A+\Ga_g\c(B,A).
\eeaa
Since  $-\ov{q}^2+|q|^2=O(r)$, we obtain
\beaa
&&\nab_4\left(\left(-\frac{1}{3}\ov{q}^2+\frac{1}{3}|q|^2\right)B\right)= \nab_4\left(-\frac{1}{3}\ov{q}^2+\frac{1}{3}|q|^2\right)B+\left(-\frac{1}{3}\ov{q}^2+\frac{1}{3}|q|^2\right)\nab_4B\\
&=& \left(-\frac{2}{3}\ov{q}+\frac{1}{3}(q+\ov{q})\right)B+\left(-\frac{1}{3}\ov{q}^2+\frac{1}{3}|q|^2\right)\left(-\frac{4}{\ov{q}} B +  O(r^{-1})\dkb^{\leq 1}A+\Ga_g\c(B,A)\right)\\
&=& (\ov{q}-q)B+O(1)\dkb^{\leq 1}A+r\Ga_g\c(B,A).
\eeaa
To eliminate $\frac{a\big(q^2-r\ov{q}\big)}{|q|^2}\ov{\Jk}\c\Xh  $ we make use of \beaa
\nab_4\Xh+\frac{2r}{|q|^2}\Xh &=& -A+\Ga_g\c\Ga_g.
\eeaa
We infer, using also $\nab_4\Jk=-q^{-1}\Jk$, $\Jk=O(r^{-1})$  and $q-\ov{q}=O(1)$,
\beaa
&&\nab_4\left(a\left(\frac{1}{2}q -\frac{1}{2}\ov{q}\right)\ov{\Jk}\c\Xh\right) \\
&=& a\nab_4\left(\frac{1}{2}q -\frac{1}{2}\ov{q}\right)\ov{\Jk}\c\Xh+a\left(\frac{1}{2}q -\frac{1}{2}\ov{q}\right)\nab_4\ov{\Jk}\c\Xh\\
&&+a\left(\frac{1}{2}q -\frac{1}{2}\ov{q}\right)\ov{\Jk}\c\nab_4\Xh\\
&=& -\frac{a}{\ov{q}}\left(\frac{1}{2}q -\frac{1}{2}\ov{q}\right)\ov{\Jk}\c\Xh+a\left(\frac{1}{2}q -\frac{1}{2}\ov{q}\right)\ov{\Jk}\c\left(-\frac{2r}{|q|^2}\Xh  -A+\Ga_g\c\Ga_g\right)\\
&=& a\left(-\frac{1}{\ov{q}}-\frac{2r}{|q|^2}\right)\left(\frac{1}{2}q -\frac{1}{2}\ov{q}\right)\ov{\Jk}\c\Xh+O(r^{-1})A+r^{-1}\Ga_g\c\Ga_g\\
&=& -\frac{a\big(q^2-r\ov{q}\big)}{|q|^2}\ov{\Jk}\c\widehat{X}+O(r^{-1})A+r^{-1}\Ga_g\c\Ga_g.
\eeaa
Combining the  above  three identities we deduce
\beaa
 \nab_4\left([\ov{q}\Hc ]_{red} \right) &=&\nab_4\left(\ov{q}\Hc -q\Zc +\left(-\frac{1}{3}\ov{q}^2+\frac{1}{3}|q|^2\right)B+a\left(\frac{1}{2}q -\frac{1}{2}\ov{q}\right)\ov{\Jk}\c\Xh\right)\\ 
&=&  O(r^{-1})\trXc+O(1)\dkb^{\leq 1}A+r\Ga_b\c\Ga_g
\eeaa
as desired.

{\bf Second identity.} To check
\beaa
\nab_4\left(q [\widecheck{\DD\cos\th}]_{ren}  \right) &= O(1)\trXc+O(r)A+r\Ga_b\c\Ga_g,
\eeaa
we  start with  the equation, see  Lemma \ref{Lemma:otherlinearizedquant}, 
\beaa
\nab_4\widecheck{\DD\cos\th}+\frac{1}{q}\widecheck{\DD\cos\th} &=& \frac{i}{2}\ov{\Jk}\c\Xh+O(r^{-1})\trXc+\Ga_b\c\Ga_g.
\eeaa
We infer
\beaa
\nab_4\Big(q\widecheck{\DD\cos\th}\Big) &=& \frac{i}{2}q\ov{\Jk}\c\Xh+O(1)\trXc+r\Ga_b\c\Ga_g.
\eeaa
To eliminate  the term in $\Xh$ we make use of
\beaa
\nab_4\Xh+\frac{2r}{|q|^2}\Xh &=& -A+\Ga_g\c\Ga_g.
\eeaa
We infer, using also $\nab_4\Jk=-q^{-1}\Jk$ and $\Jk=O(r^{-1})$,
\beaa
\nab_4\left(\frac{i}{2}|q|^2\ov{\Jk}\c\Xh\right) &=& \nab_4\left(\frac{i}{2}|q|^2\right)\ov{\Jk}\c\Xh+\left(\frac{i}{2}|q|^2\right)\nab_4\ov{\Jk}\c\Xh+\left(\frac{i}{2}|q|^2\right)\ov{\Jk}\c\nab_4\Xh\\
&=& \frac{i}{2}(q+\ov{q})\ov{\Jk}\c\Xh+\left(\frac{i}{2}|q|^2\right)\left(-\frac{1}{\ov{q}}\right)\ov{\Jk}\c\Xh\\
&&+\left(\frac{i}{2}|q|^2\right)\ov{\Jk}\c\left(-\frac{2r}{|q|^2}\Xh  -A+\Ga_g\c\Ga_g\right)\\
&=& -\frac{i}{2}q\ov{\Jk}\c\Xh  +O(r)A+r\Ga_g\c\Ga_g.
\eeaa
Summing the two identities above, we infer
\beaa
\nab_4\left(q [\widecheck{\DD\cos\th}]_{ren}  \right) =\nab_4\left(q\widecheck{\DD\cos\th} +\frac{i}{2}|q|^2\ov{\Jk}\c\Xh\right) &=& O(1)\trXc+O(r)A+r\Ga_b\c\Ga_g
\eeaa
as desired.

{\bf Third identity.}  To check 
\beaa
\nab_4\left(\ov{q} q^2 [\Mc]_{ren}\right)&=   O(1)\dkb^{\leq 1}\trXc+O(r)\dkb^{\leq 2}A+r^2\dkb^{\leq 1}(\Ga_g\c\Ga_g)+r^3\Ga_b\c A,
\eeaa
 we start with the equation 
\beaa
\nab_4\Zc + \frac{2}{q}\Zc  &=&    - \frac{aq}{|q|^2}\ov{\Jk}\c\widehat{X} -B +O(r^{-2})\trXc+\Ga_g\c\Ga_g
\eeaa
which  we rewrite in the form
\beaa
\nab_4(q^2\Zc)   &=&    - \frac{aq^3}{|q|^2}\ov{\Jk}\c\widehat{X} -q^2B +O(1)\trXc+r^2\Ga_g\c\Ga_g.
\eeaa
To eliminate the term in $\Xh$ on the right we make use of the equation
\beaa
\nab_4\Xh+\frac{2r}{|q|^2}\Xh &=& -A+\Ga_g\c\Ga_g.
\eeaa
Using also $\nab_4\Jk=-q^{-1}\Jk$ and $\Jk=O(r^{-1})$ we infer
\beaa
&&\nab_4\left(\left(-\frac{a}{2}q^2 -\frac{a}{2}|q|^2\right)\ov{\Jk}\c\Xh\right)\\
 &=& \nab_4\left(-\frac{a}{2}q^2 -\frac{a}{2}|q|^2\right)\ov{\Jk}\c\Xh+\left(-\frac{a}{2}q^2 -\frac{a}{2}|q|^2\right)\nab_4\ov{\Jk}\c\Xh\\
&&+\left(-\frac{a}{2}q^2 -\frac{a}{2}|q|^2\right)\ov{\Jk}\c\nab_4\Xh\\
&=& \left(-aq -\frac{a}{2}(q+\ov{q})\right)\ov{\Jk}\c\Xh+\left(-\frac{a}{2}q^2 -\frac{a}{2}|q|^2\right)\left(-\frac{1}{\ov{q}}\right)\ov{\Jk}\c\Xh\\
&&+\left(-\frac{a}{2}q^2 -\frac{a}{2}|q|^2\right)\ov{\Jk}\c\left(-\frac{2r}{|q|^2}\Xh  -A+\Ga_g\c\Ga_g\right)\\
&=& \frac{aq^3}{|q|^2}\ov{\Jk}\c\Xh + O(r)A+r\Ga_g\c\Ga_g.
\eeaa

Summing the two identities above, we infer
\beaa
\nab_4\left(q^2\Zc +\left(-\frac{a}{2}q^2 -\frac{a}{2}|q|^2\right)\ov{\Jk}\c\Xh\right)   &=&   -q^2B +O(1)\trXc+O(r)A+r^2\Ga_g\c\Ga_g.
\eeaa
Commuting with $\ov{q}\,\ov{\DD}\c$ we deduce
\beaa
&&\nab_4\left(\ov{q}\,\ov{\DD}\c\left(q^2\Zc +\left(-\frac{a}{2}q^2 -\frac{a}{2}|q|^2\right)\ov{\Jk}\c\Xh\right)\right) \\
&=& [\nab_4, \ov{q}\,\ov{\DD}\c]\left(q^2\Zc +\left(-\frac{a}{2}q^2 -\frac{a}{2}|q|^2\right)\ov{\Jk}\c\Xh\right)\\
&& -\ov{q}\,\ov{\DD}\c(q^2B) +O(1)\dkb^{\leq 1}\trXc+O(r)\dkb^{\leq 1}A+r^2\dkb^{\leq 1}(\Ga_g\c\Ga_g)\\
&=& \left(-\ov{q}\,\ov{\tr X}\,\ov{Z}+\Ga_g\c\dk^{\leq 1}\right)\c\left(q^2\Zc +\left(-\frac{a}{2}q^2 -\frac{a}{2}|q|^2\right)\ov{\Jk}\c\Xh\right)\\
&& -\ov{q}\,\ov{\DD}\c(q^2B) +O(1)\dkb^{\leq 1}\trXc+O(r)\dkb^{\leq 1}A+r^2\dkb^{\leq 1}(\Ga_g\c\Ga_g)
\eeaa
and hence
\beaa
&&\nab_4\left(\ov{q}\,\ov{\DD}\c\left(q^2\Zc +\left(-\frac{a}{2}q^2 -\frac{a}{2}|q|^2\right)\ov{\Jk}\c\Xh\right)\right) \\
&=& -\frac{2aq^3}{|q|^2}\ov{\Jk}\c\Zc +\frac{a^2q}{|q|^2}\left(q^2 +|q|^2\right)\ov{\Jk}\c\Xh\c\ov{\Jk} -\ov{q}\,\ov{\DD}\c(q^2B) \\
&&+O(1)\dkb^{\leq 1}\trXc+O(r)\dkb^{\leq 1}A+r^2\dkb^{\leq 1}(\Ga_g\c\Ga_g).
\eeaa
Next we  make use of 
\beaa
\nab_4\left(\Pc \right)-\frac{1}{2}\DD\c \ov{B} &=& -\frac{3}{q} \Pc -\frac{a\ov{q}}{2|q|^2}\Jk\c\ov{B}+O(r^{-3})\trXc+r^{-1}\Ga_g\c\Ga_g+\Ga_b\c A
\eeaa
from which we obtain 
\beaa
\nab_4(2q^3\Pc) &=& q^3\DD\c \ov{B} - aq^2\Jk\c\ov{B}+O(1)\trXc+r^2\Ga_g\c\Ga_g+r^3\Ga_b\c A.
\eeaa
Taking the complex conjugate, we infer
\beaa
\nab_4\left(2\ov{q}^3\,\ov{\Pc}\right) &=& \ov{q}^3\ov{\DD}\c B - a\ov{q}^2\ov{\Jk}\c B+O(1)\trXc+r^2\Ga_g\c\Ga_g+r^3\Ga_b\c A.
\eeaa
Adding to the previous identity, we deduce
\beaa
&&\nab_4\left(\ov{q}\,\ov{\DD}\c\left(q^2\Zc +\left(-\frac{a}{2}q^2 -\frac{a}{2}|q|^2\right)\ov{\Jk}\c\Xh\right)+2\ov{q}^3\,\ov{\Pc}\right) \\
&=& -\frac{2aq^3}{|q|^2}\ov{\Jk}\c\Zc +\frac{a^2q}{|q|^2}\left(q^2 +|q|^2\right)\ov{\Jk}\c\Xh\c\ov{\Jk} -\ov{q}\,\ov{\DD}\c(q^2B)+\ov{q}^3\ov{\DD}\c B - a\ov{q}^2\ov{\Jk}\c B\\
&&+O(1)\dkb^{\leq 1}\trXc+O(r)\dkb^{\leq 1}A+r^2\dkb^{\leq 1}(\Ga_g\c\Ga_g)+r^3\Ga_b\c A.
\eeaa
Since we have
\beaa
-\ov{q}\,\ov{\DD}\c(q^2B)+\ov{q}^3\ov{\DD}\c B &=&  \ov{q}(\ov{q}^2-q^2)\ov{\DD}\c B -2 |q|^2\ov{\DD}q\c B\\
&=&  \ov{q}(\ov{q}^2-q^2)\ov{\DD}\c B -2ai|q|^2\ov{\DD}(\cos\th)\c B\\
&=&  \ov{q}(\ov{q}^2-q^2)\ov{\DD}\c B -2a|q|^2\ov{\Jk}\c B +r\Ga_b\c\Ga_g,
\eeaa
we infer
\beaa
&&\nab_4\left(\ov{q}\,\ov{\DD}\c\left(q^2\Zc +\left(-\frac{a}{2}q^2 -\frac{a}{2}|q|^2\right)\ov{\Jk}\c\Xh\right)+2\ov{q}^3\,\ov{\Pc}\right) \\
&=& -\frac{2aq^3}{|q|^2}\ov{\Jk}\c\Zc +\ov{q}(\ov{q}^2-q^2)\ov{\DD}\c B -a\Big(2|q|^2+\ov{q}^2\Big)\ov{\Jk}\c B +\frac{a^2q}{|q|^2}\left(q^2 +|q|^2\right)\ov{\Jk}\c\Xh\c\ov{\Jk}\\
&&+O(1)\dkb^{\leq 1}\trXc+O(r)\dkb^{\leq 1}A+r^2\dkb^{\leq 1}(\Ga_g\c\Ga_g)+r^3\Ga_b\c A.
\eeaa
To eliminate the term  in $\Zc$ on the right 
we write
\beaa
\nab_4\Big(- 2aq^2\ov{\Jk}\c\Zc\Big) &=& a\nab_4(- 2q^2)\ov{\Jk}\c\Zc - 2aq^2\nab_4\ov{\Jk}\c\Zc - 2aq^2\ov{\Jk}\c\nab_4\Zc\\
&=& -4aq\ov{\Jk}\c\Zc +\frac{2aq^2}{\ov{q}}\ov{\Jk}\c\Zc\\
&&- 2aq^2\ov{\Jk}\c\left(- \frac{2}{q}\Zc     - \frac{aq}{|q|^2}\ov{\Jk}\c\widehat{X} -B +O(r^{-2})\trXc+\Ga_g\c\Ga_g\right)\\
&=& \frac{2aq^3}{|q|^2}\ov{\Jk}\c\Zc +\frac{2a^2q^3}{|q|^2}\ov{\Jk}\c\widehat{X}\c\ov{\Jk} +2aq^2\ov{\Jk}\c B +O(r^{-1})\trXc+r\Ga_g\c\Ga_g.
\eeaa
Summing with the previous identity, we infer
\bea
\lab{eq:Prop.eqforrenormalized.quantities3}
\bsplit
&\nab_4\left(\ov{q}\,\ov{\DD}\c\left(q^2\Zc +\left(-\frac{a}{2}q^2 -\frac{a}{2}|q|^2\right)\ov{\Jk}\c\Xh\right)+2\ov{q}^3\,\ov{\Pc} - 2aq^2\ov{\Jk}\c\Zc\right) \\
&= \ov{q}(\ov{q}^2-q^2)\ov{\DD}\c B +a\Big(-2|q|^2-\ov{q}^2 +2q^2\Big)\ov{\Jk}\c B  +\frac{a^2q}{|q|^2}\left(3q^2 +|q|^2\right)\ov{\Jk}\c\Xh\c\ov{\Jk}\\
&+O(1)\dkb^{\leq 1}\trXc+O(r)\dkb^{\leq 1}A+r^2\dkb^{\leq 1}(\Ga_g\c\Ga_g)+r^3\Ga_b\c A.
\end{split}
\eea
To eliminate the term in  $\ov{\DD}\c B$ we  make use of the equation
\beaa
\nab_4B +\frac{4}{\ov{q}} B &=&  O(r^{-1})\dkb^{\leq 1}A+\Ga_g\c(B,A).
\eeaa
Differentiating w.r.t. $\ov{\DD}\c$, we infer
\beaa
\nab_4( \ov{\DD}\c B )&=& [\nab_4, \ov{\DD}\c]B +  \ov{\DD}\c\left(-\frac{4}{\ov{q}} B +  O(r^{-1})\dkb^{\leq 1}A+\Ga_g\c(B,A)\right)\\
&=& -\frac{1}{2}\ov{\tr X}\Big(\ov{\DD}\c B +2\ov{Z}\c B\Big)+r^{-1}\Ga_g\c \dk^{\leq 1}B -\frac{4}{\ov{q}} \ov{\DD}\c B+\frac{4}{\ov{q}^2}\ov{\DD}(\ov{q})\c B\\
&&+O(r^{-2})\dkb^{\leq 2}A+r^{-2}\dk^{\leq 1}(\Ga_g\c\Ga_g)\\
&=& -\frac{5}{\ov{q}} \ov{\DD}\c B -\frac{2}{\ov{q}}\ov{Z}\c B -\frac{4ai}{\ov{q}^2}\ov{\DD}(\cos\th)\c B\\
&&+O(r^{-2})\dkb^{\leq 2}A+r^{-2}\dk^{\leq 1}(\Ga_g\c\Ga_g)
\eeaa
and hence
\beaa
\nab_4(\ov{\DD}\c B )&=& -\frac{5}{\ov{q}} \ov{\DD}\c B  -\frac{6a}{\ov{q}^2}\ov{\Jk}\c B +O(r^{-2})\dkb^{\leq 2}A+r^{-2}\dk^{\leq 1}(\Ga_g\c\Ga_g).
\eeaa
We then  compute
\beaa
&&\nab_4\left(\left(-\frac{1}{3}q^2\ov{q}^2 -\frac{1}{3}q\ov{q}^3+\frac{2}{3}\ov{q}^4\right)\ov{\DD}\c B\right)\\
 &=& \nab_4\left(-\frac{1}{3}q^2\ov{q}^2 -\frac{1}{3}q\ov{q}^3+\frac{2}{3}\ov{q}^4\right)\ov{\DD}\c B+\left(-\frac{1}{3}q^2\ov{q}^2 -\frac{1}{3}q\ov{q}^3+\frac{2}{3}\ov{q}^4\right)\nab_4\ov{\DD}\c B\\
&=&  \left(-\frac{2}{3}q\ov{q}^2 -\frac{2}{3}q^2\ov{q} -\frac{1}{3}\ov{q}^3 -q\ov{q}^2+\frac{8}{3}\ov{q}^3\right)\ov{\DD}\c B\\
&&+\left(-\frac{1}{3}q^2\ov{q}^2 -\frac{1}{3}q\ov{q}^3+\frac{2}{3}\ov{q}^4\right)\left( -\frac{5}{\ov{q}} \ov{\DD}\c B  -\frac{6a}{\ov{q}^2}\ov{\Jk}\c B +O(r^{-2})\dkb^{\leq 2}A+r^{-2}\dk^{\leq 1}(\Ga_g\c\Ga_g)\right)\\
&=& \ov{q}(q^2 -\ov{q}^2)\ov{\DD}\c B    -\frac{6a}{\ov{q}^2}\left(-\frac{1}{3}q^2\ov{q}^2 -\frac{1}{3}q\ov{q}^3+\frac{2}{3}\ov{q}^4\right)\ov{\Jk}\c B +O(r)\dkb^{\leq 2}A+r\dk^{\leq 1}(\Ga_g\c\Ga_g)
\eeaa
where we used in particular the fact that $-\frac{1}{3}q^2\ov{q}^2 -\frac{1}{3}q\ov{q}^3+\frac{2}{3}\ov{q}^4 = O(r^3).$
Together  with \eqref{eq:Prop.eqforrenormalized.quantities3}, we deduce, for the  intermediary quantity $[\Mc']_{red}$ defined by 
\beaa
[\Mc']_{red}:&=&\ov{q}\,\ov{\DD}\c\left(q^2\Zc +\left(-\frac{a}{2}q^2 -\frac{a}{2}|q|^2\right)\ov{\Jk}\c\Xh\right)+2\ov{q}^3\,\ov{\Pc} - 2aq^2\ov{\Jk}\c\Zc\\
&&+\left(-\frac{1}{3}q^2\ov{q}^2 -\frac{1}{3}q\ov{q}^3+\frac{2}{3}\ov{q}^4\right)\ov{\DD}\c B
\eeaa
the following identity
\bea
\lab{eq:Prop.eqforrenormalized.quantities4}
\bsplit
\nab_4\left([\Mc']_{red}\right) &=  a\Big( 4q^2 - 5\ov{q}^2\Big)\ov{\Jk}\c B  +\frac{a^2q}{|q|^2}\left(3q^2 +|q|^2\right)\ov{\Jk}\c\Xh\c\ov{\Jk}\\
&+O(1)\dkb^{\leq 1}\trXc+O(r)\dkb^{\leq 2}A+r^2\dkb^{\leq 1}(\Ga_g\c\Ga_g)+r^3\Ga_b\c A.
\end{split}
\eea
To eliminate  the term in $\Jk\c B$ on the RHS of \eqref{eq:Prop.eqforrenormalized.quantities4},  we compute
\beaa
&&\nab_4\left(-a\left(- q^2\ov{q} -\frac{2}{3}q\ov{q}^2+ \frac{13}{6}\ov{q}^3\right)\ov{\Jk}\c B\right)\\
 &=& -a\nab_4\left(- q^2\ov{q} -\frac{2}{3}q\ov{q}^2+ \frac{13}{6}\ov{q}^3\right)\ov{\Jk}\c B-a\left(- q^2\ov{q} -\frac{2}{3}q\ov{q}^2+ \frac{13}{6}\ov{q}^3\right)\nab_4\ov{\Jk}\c B\\
 &&-a\left(- q^2\ov{q} -\frac{2}{3}q\ov{q}^2+ \frac{13}{6}\ov{q}^3\right)\ov{\Jk}\c \nab_4B\\
&=& -a\left(-q^2   -\frac{10}{3}|q|^2 + \frac{35}{6}\ov{q}^2\right)\ov{\Jk}\c B-a\left(- q^2\ov{q} -\frac{2}{3}q\ov{q}^2+ \frac{13}{6}\ov{q}^3\right)\left(-\frac{1}{\ov{q}}\right)\ov{\Jk}\c B\\
&& -a\left(- q^2\ov{q} -\frac{2}{3}q\ov{q}^2+ \frac{13}{6}\ov{q}^3\right)\ov{\Jk}\c \left(-\frac{4}{\ov{q}} B +  O(r^{-1})\dkb^{\leq 1}A+\Ga_g\c(B,A)\right)\\
&=& -a(4q^2    -5\ov{q}^2)\ov{\Jk}\c B +O(r)\dkb^{\leq 1}A+r\Ga_g\c\Ga_g.
\eeaa
Summing with \eqref{eq:Prop.eqforrenormalized.quantities4}, we infer
\beaa
&&\nab_4\Bigg\{[\Mc']_{ren}  +a\left( q^2\ov{q} +\frac{2}{3}q\ov{q}^2 - \frac{13}{6}\ov{q}^3\right)\ov{\Jk}\c B\Bigg\} \\
&=&   \frac{a^2q}{|q|^2}\left(3q^2 +|q|^2\right)\ov{\Jk}\c\Xh\c\ov{\Jk} +O(1)\dkb^{\leq 1}\trXc+O(r)\dkb^{\leq 2}A+r^2\dkb^{\leq 1}(\Ga_g\c\Ga_g)+r^3\Ga_b\c A.
\eeaa
It remains to eliminate  the term in $\Xh$. To do this we write
\beaa
\nab_4(a^2(q^2+|q|^2)\ov{\Jk}\c\Xh\c\ov{\Jk}) &=& a^2\nab_4(q^2+|q|^2)\ov{\Jk}\c\Xh\c\ov{\Jk}+2a^2(q^2+|q|^2)\ov{\Jk}\c\Xh\c\nab_4\ov{\Jk}\\
&&+a^2(q^2+|q|^2)\ov{\Jk}\c\nab_4\Xh\c\ov{\Jk}\\
&=& a^2(3q+\ov{q})\ov{\Jk}\c\Xh\c\ov{\Jk}+2a^2(q^2+|q|^2)\ov{\Jk}\c\Xh\c\left(-\frac{1}{\ov{q}}\ov{\Jk}\right)\\
&&+a^2(q^2+|q|^2)\ov{\Jk}\c\left(-\frac{2r}{|q|^2}\Xh  -A+\Ga_g\c\Ga_g\right)\c\ov{\Jk}\\
&=& -\frac{a^2q}{|q|^2}(3q^2+|q|^2)\ov{\Jk}\c\Xh\c\ov{\Jk} +O(1)A+\Ga_g\c\Ga_g.
\eeaa
Summing with the previous identity above, we infer that
\beaa
 \nab_4\left( [\Mc]_{ren}\right) &=&\nab_4\Bigg\{ [\Mc']_{ren}  +a\left( q^2\ov{q} +\frac{2}{3}q\ov{q}^2 - \frac{13}{6}\ov{q}^3\right)\ov{\Jk}\c B+a^2(q^2+|q|^2)\ov{\Jk}\c\Xh\c\ov{\Jk}\Bigg\} \\
&=&   O(1)\dkb^{\leq 1}\trXc+O(r)\dkb^{\leq 2}A+r^2\dkb^{\leq 1}(\Ga_g\c\Ga_g)+r^3\Ga_b\c A
\eeaa
as stated. This concludes the proof of Proposition \ref{Prop: eqforrenormalized.qiantities}.

%%%%%%%%%%%%%%%%%%%%%%%%%%%%%%%% 
 
 \section{Proof of Lemma \ref{Le:mainpointwise-ell=1B}}
 \lab{section:ProofLemma-mainpointwise-ell=1B}
 
%%%%%%%%%%%%%%%%%%%%%%%%%%%%%%%% 
 
The goal  of this section is to prove  the following identity 
\beaa
\Big(\nab_4-a\Re(\Jk)^b\nab_b\Big)       \Big(    r [\ov{\DD}\c]_{ren}  \big(r^4 [B]_{ren} \big)\Big)
&=& \frac{r^5}{2}\ov{\DD}'\c\ov{\DD}'\c\Big( A  -a(\Jk\hot B)\Big) +\err,
\eeaa
where the  $\DD'$  is taken with respect  to the integral  frame  $(e_1', e_2')$   adapted to $S(u, r)$, see section \ref{subsection:orthonormalbasisS(u,r)}, and where the error term has the following structure
\beaa
\err&=&   + O(1)\dk^{\leq 1}\Xh+O(r)\dk^{\leq 2}B+O(r^2)\dk^{\leq 1}\nab_3B+O(r)\dk^{\leq 2}\Pc +O(1)\dk^{\leq 1}\trXc\\
&&+O(r^{2})\dk^{\leq 1}\nab_3A+O(r)\dk^{\leq 2}A+r^4\dk^{\leq 1}\big(\Ga_g\c(B,A)\big) +r^4\dk^{\leq 1}\big(\Ga_b\c\nab_3 A\big)\\
&&+r^2\dk^{\leq 2}\big(\Ga_g\c\Ga_g\big).
\eeaa

{\bf Step 1.} We start with the following lemma.
\begin{lemma}
\lab{Lemma:Le:mainpointwise-ell=1B1}
We have
\beaa
&&\Big(\nab_4-a\Re(\Jk)^b\nab_b\Big) [B]_{ren}   +\frac{4}{\ov{q}} [B]_{ren}\\
&=&   - \frac{3a}{4}(\ov{\DD}\c B)\Jk     -a\Re(\Jk)^b\nab_bB   +  \frac{1}{2}\ov{\DD}\c A  -\frac{a}{4}\ov{\Jk}\c \nab_4A -\frac{a}{4\ov{q}}\ov{\Jk}\c A\\
&& +O(r^{-3})B+O(r^{-3})\dkb^{\leq 1}\Pc+O(r^{-3})\dkb^{\leq 1}A +O(r^{-4})\trXc+\Ga_g\c(B,A) +r^{-2}\Ga_g\c\Ga_g.
\eeaa
\end{lemma}

\begin{proof}
We start with the linearized Bianchi  equations
\beaa
\nab_4B +\frac{4}{\ov{q}} B &=&  \frac{1}{2}\ov{\DD}\c A +\frac{aq}{2|q|^2}\ov{\Jk}\c A+\Ga_g\c(B,A),\\
\nab_4\left(\Pc \right)-\frac{1}{2}\DD\c \ov{B} &=& -\frac{3}{q} \Pc -\frac{a\ov{q}}{2|q|^2}\Jk\c\ov{B}+O(r^{-3})\trXc+r^{-1}\Ga_g\c\Ga_g+\Ga_b\c A.
\eeaa
We infer, using $\nab_4\Jk=-q^{-1}\Jk$, 
\beaa
\nab_4\big([B]_{ren}\big) &=& \nab_4B - \frac{3a}{2}\Jk\nab_4\ov{\Pc} - \frac{3a}{2}\ov{\Pc}\nab_4\Jk  -\frac{a}{4}\ov{\Jk}\c \nab_4A  -\frac{a}{4}A\c\nab_4\ov{\Jk}\\
&=& -\frac{4}{\ov{q}} B +  \frac{1}{2}\ov{\DD}\c A +\frac{aq}{2|q|^2}\ov{\Jk}\c A+\Ga_g\c(B,A) \\
&& - \frac{3a}{2}\Jk\left(\frac{1}{2}\ov{\DD}\c B  -\frac{3}{\ov{q}} \ov{\Pc} -\frac{aq}{2|q|^2}\ov{\Jk}\c B+O(r^{-3})\trXc+r^{-1}\Ga_g\c\Ga_g+\Ga_b\c A\right)\\
&& +\frac{3a}{2q}\ov{\Pc}\Jk  -\frac{a}{4}\ov{\Jk}\c \nab_4A  +\frac{a}{4\ov{q}}\ov{\Jk}\c A\\
&=& -\frac{4}{\ov{q}}[B]_{ren} - \frac{3a}{2}\Jk\left(\frac{1}{2}\ov{\DD}\c B  -\frac{aq}{2|q|^2}\ov{\Jk}\c B\right) +\frac{3a}{2}\left(\frac{1}{q}-\frac{1}{\ov{q}}\right)\ov{\Pc}\Jk\\
&& +  \frac{1}{2}\ov{\DD}\c A  -\frac{a}{4}\ov{\Jk}\c \nab_4A -\frac{a}{4\ov{q}}\ov{\Jk}\c A +O(r^{-4})\trXc+\Ga_g\c(B,A) +r^{-2}\Ga_g\c\Ga_g
\eeaa
and hence
\beaa
\nab_4 \big([B]_{ren}\big) +\frac{4}{\ov{q}} [B]_{ren}
&=&  - \frac{3a}{4}(\ov{\DD}\c B)\Jk     +  \frac{1}{2}\ov{\DD}\c A  -\frac{a}{4}\ov{\Jk}\c \nab_4A -\frac{a}{4\ov{q}}\ov{\Jk}\c A\\
&& +O(r^{-3})B+O(r^{-3})\Pc+O(r^{-4})\trXc+\Ga_g\c(B,A) +r^{-2}\Ga_g\c\Ga_g.
\eeaa

Also, we have
\beaa
a\Re(\Jk)^b\nab_b[B]_{ren}&=& a\Re(\Jk)^b\nab_bB - \frac{3a^2}{2}\Re(\Jk)^b\nab_b(\ov{\Pc}\Jk) +O(r^{-3})\dkb^{\leq 1}A
\eeaa
and hence
\beaa
&&\Big(\nab_4-a\Re(\Jk)^b\nab_b\Big)[B]_{ren}+\frac{4}{\ov{q}} [B]_{ren}\\ 
&=&   - \frac{3a}{4}(\ov{\DD}\c B)\Jk     -a\Re(\Jk)^b\nab_bB   +  \frac{1}{2}\ov{\DD}\c A  -\frac{a}{4}\ov{\Jk}\c \nab_4A -\frac{a}{4\ov{q}}\ov{\Jk}\c A\\
&& +O(r^{-3})B+O(r^{-3})\dkb^{\leq 1}\Pc+O(r^{-3})\dkb^{\leq 1}A +O(r^{-4})\trXc+\Ga_g\c(B,A) +r^{-2}\Ga_g\c\Ga_g
\eeaa
as desired. This concludes the proof of Lemma \ref{Lemma:Le:mainpointwise-ell=1B1}.
\end{proof}

{\bf Step 2.}  Next, we prove the following lemma.
\begin{lemma}
\lab{Lemma:Le:mainpointwise-ell=1B2}
The following holds true, with  $\nab'$ the covariant derivative on $S(u, r)$ and $\DD'$ the corresponding complex operator,
\beaa
  \ov{\DD}\c A -\frac{a}{2}\ov{\Jk}\c \nab_4A - \frac{a}{2\ov{q}}\ov{\Jk}\c A
&=& \ov{\DD}'\c A + \frac{a}{4}\ov{\Jk}\c\DD\hot B + O(r^{-3})B +O(r^{-4})\Xh\\
&&+O(r^{-2})\nab_3A+O(r^{-3})\dk^{\leq 1}A+\Ga_b\c\nab_3 A\\
&&+r^{-1}\Ga_b\c (A, B)+r^{-2}\Ga_g\c\Ga_g.
\eeaa
\end{lemma}

\begin{proof}
We make use of the  formula
\bea
\bsplit
\lab{eq:Step2[B]_{ren}1}
\DD_a' V_B &= \DD_a V_B+\frac{1}{2}\underline{F}_a f^c\nab_c V_B+\frac{1}{2}\underline{F}_a\nab_4V_B+\left(\frac{1}{2}F_a+\frac{1}{8}|f|^2\underline{F}_a\right)\nab_3V_B\\
&+\Big\{O(r^{-3})+ r^{-1} \Ga_b \Big\} V
\end{split}
\eea
which follows from   the last identity of   Proposition \ref{cor:formuallikningnabprimeandnabintransfoformula}, i.e.  
\beaa
\bsplit
\DD_a' V_B &= \DD_a V_B+\frac{1}{2}\underline{F}_a f^c\nab_c V_B+\frac{1}{2}\underline{F}_a\nab_4V_B+\left(\frac{1}{2}F_a+\frac{1}{8}|f|^2\underline{F}_a\right)\nab_3V_B\\
&+(E[V])_{aB},
\end{split}
\eeaa
and  formula  \eqref{eq:(EV)_{aB}}  for  $(E[V])_{aB}$   derived   in the proof of Corollary \ref{cor:formuallikningnabprimeandnabintransfoformula:bis}.
Here,  the  horizontal 1-forms $f$ and $\fb$ are given by \eqref{eq:formulaforchangeofframecoeffandfbarforframetangenttosphereSur} and 
\beaa
F=f+i\dual f, \qquad \underline{F}=\fb+i\dual\fb.
\eeaa
In view of \eqref{eq:Step2[B]_{ren}1}
we  deduce 
\beaa
\ov{\DD}'\c A &=& \ov{\DD}\c A+\frac{1}{2}\ov{\underline{F}}\c \nab_4A + \frac{1}{2}\ov{F}\c\nab_3A+O(r^{-3})   \dk^{\leq 1}A+r^{-1}\Ga_b\c A.
\eeaa

Also, in view of \eqref{eq:formulaforchangeofframecoeffandfbarforframetangenttosphereSur}, the horizontal 1-forms $f$ and $\fb$ satisfy 
\beaa
f=\Big(-1+O(r^{-1})+r\Ga_b\Big)\nab u, \qquad \fb=\Big(-1+O(r^{-1})+r\Ga_b\Big)\nab u
\eeaa
and hence
\beaa
F &=& \Big(-1+O(r^{-1})+r\Ga_b\Big)\DD u=-a\Jk+O(r^{-2})+\Ga_b, \\ 
\underline{F}   &=& \Big(-1+O(r^{-1})+r\Ga_b\Big)\DD u=-a\Jk+O(r^{-2})+\Ga_b.
\eeaa
We infer
\beaa
\ov{\DD}'\c A &=& \ov{\DD}\c A -\frac{a}{2}\ov{\Jk}\c \nab_4A - \frac{a}{2}\ov{\Jk}\c\nab_3A+O(r^{-2})\nab_3A+O(r^{-3})\dk^{\leq 1}A\\
&&+\Ga_b\c\nab_3 A+r^{-1}\Ga_b\c A.
\eeaa
Making use of  the Bianchi identity 
\beaa
 \nab_3A  &=&  \frac{1}{2}\DD\hot B -\frac{1}{2}\tr\Xb A+4\omb A +\frac{1}{2}(Z+4H)\hot B -3\ov{P}\Xh\\
 &=&  \frac{1}{2}\DD\hot B +\frac{1}{r}A + \frac{a(\ov{q}+4q)}{2|q|^2}\Jk\hot B +\frac{6m}{\ov{q}^3}\Xh +O(r^{-2})A+\Ga_b\c (A,B)+r^{-1}\Ga_g\c\Ga_g,
\eeaa
we deduce
\beaa
\ov{\DD}'\c A &=& \ov{\DD}\c A -\frac{a}{2}\ov{\Jk}\c \nab_4A - \frac{a}{2r}\ov{\Jk}\c A - \frac{a}{4}\ov{\Jk}\c\DD\hot B - \frac{a^2(\ov{q}+4q)}{4|q|^2}\ov{\Jk}\c(\Jk\hot B) - \frac{3am}{\ov{q}^3}\ov{\Jk}\c\Xh \\
&&+O(r^{-2})\nab_3A+O(r^{-3})\dk^{\leq 1}A+\Ga_b\c\nab_3 A+r^{-1}\Ga_b\c (A, B)+r^{-2}\Ga_g\c\Ga_g
\eeaa
and hence
\beaa
&&  \ov{\DD}\c A -\frac{a}{2}\ov{\Jk}\c \nab_4A - \frac{a}{2\ov{q}}\ov{\Jk}\c A= \ov{\DD}'\c A + \frac{a}{4}\ov{\Jk}\c\DD\hot B + O(r^{-3})B +O(r^{-4})\Xh\\
&& +O(r^{-2})\nab_3A+O(r^{-3})\dk^{\leq 1}A+\Ga_b\c\nab_3 A+r^{-1}\Ga_b\c (A, B)+r^{-2}\Ga_g\c\Ga_g
\eeaa
as stated. This concludes the proof of Lemma \ref{Lemma:Le:mainpointwise-ell=1B2}. 
\end{proof}

{\bf Step 3.} Next, we derive the following corollary.
\begin{corollary}
\lab{Corr:Step3[B]_{ren}}
We have
\beaa
&&\Big(\nab_4-a\Re(\Jk)^b\nab_b\Big)[B]_{red} +\frac{4}{r}[B]_{red}\\ 
&=& \frac{1}{2}\ov{\DD}'\c A  - \frac{3a}{4}(\ov{\DD}\c B)\Jk     -a\Re(\Jk)^b\nab_bB + \frac{a}{8}\ov{\Jk}\c\DD\hot B  -\frac{4ia\cos\th}{r^2}B\\
&&   + O(r^{-4})\Xh+O(r^{-3})B+O(r^{-3})\dk^{\leq 1}\Pc +O(r^{-4})\trXc\\
&&+O(r^{-2})\nab_3A+O(r^{-3})\dk^{\leq 1}A+\Ga_g\c(B,A)+\Ga_b\c\nab_3 A+r^{-2}\Ga_g\c\Ga_g.
\eeaa
\end{corollary}

\begin{proof}
Plugging the  identity  of Lemma  \ref{Lemma:Le:mainpointwise-ell=1B2}  in the identity of   Lemma  \ref{Lemma:Le:mainpointwise-ell=1B1}, we obtain 
\beaa
&&\Big(\nab_4-a\Re(\Jk)^b\nab_b\Big)[B]_{red} +\frac{4}{\ov{q}}[B]_{red}\\ 
&=& \frac{1}{2}\ov{\DD}'\c A  - \frac{3a}{4}(\ov{\DD}\c B)\Jk     -a\Re(\Jk)^b\nab_bB + \frac{a}{8}\ov{\Jk}\c\DD\hot B  \\
&& + O(r^{-4})\Xh+O(r^{-3})B+O(r^{-3})\dkb^{\leq 1}\Pc +O(r^{-4})\trXc\\
&&+O(r^{-2})\nab_3A+O(r^{-3})\dk^{\leq 1}A+\Ga_g\c(B,A)+\Ga_b\c\nab_3 A+r^{-2}\Ga_g\c\Ga_g.
\eeaa
Since 
\beaa
\frac{4}{\ov{q}}[B]_{red}\
&=& \frac{4}{r}[B]_{red}+\frac{4ia\cos\th}{r^2}B+O(r^{-3})\Pc+O(r^{-3})B+O(r^{-3})A,
\eeaa
we deduce
\beaa
&&\Big(\nab_4-a\Re(\Jk)^b\nab_b\Big)[B]_{red} +\frac{4}{r}[B]_{red}\\ 
&=& \frac{1}{2}\ov{\DD}'\c A  - \frac{3a}{4}(\ov{\DD}\c B)\Jk     -a\Re(\Jk)^b\nab_bB + \frac{a}{8}\ov{\Jk}\c\DD\hot B  -\frac{4ia\cos\th}{r^2}B\\
&&   + O(r^{-4})\Xh+O(r^{-3})B+O(r^{-3})\dk^{\leq 1}\Pc +O(r^{-4})\trXc\\
&&+O(r^{-2})\nab_3A+O(r^{-3})\dk^{\leq 1}A+\Ga_g\c(B,A)+\Ga_b\c\nab_3 A+r^{-2}\Ga_g\c\Ga_g
\eeaa
as desired. This concludes the proof of Corollary \ref{Corr:Step3[B]_{ren}}. 
\end{proof}

{\bf Step 4.} Next, we derive the following lemma.
\begin{lemma}
\lab{Le:Step4[B]_{ren}}
We have
\beaa
\ov{\DD}\c(\Jk\hot B) &=& -\frac{1}{4}\ov{\Jk}\c(\DD\hot B)+\frac{3}{2}(\ov{\DD}\c B)\Jk+2\Re(\Jk)^b\nab_bB+\frac{8i(r^2+a^2)\cos\th}{|q|^4}B\\
&&+r^{-1}\Ga_b\c B.
\eeaa
\end{lemma}

\begin{proof}
The proof follows immediately by combining  the following three identities
\bea
\lab{eq:LeStep4[B]_{ren}1}
\bsplit
\ov{\DD}\c(\Jk\hot B)_b &= (\ov{\DD}\c\Jk) B_b+\de_{cd}\ov{\DD}_c\Jk_b B_d - \ov{\DD}_b\Jk\c B+ \Jk_b\ov{\DD}\c B\\
&+\de_{cd}\Jk_d\ov{\DD}_c B_b - \Jk\c\ov{\DD}_b B,
\end{split}
\eea
\bea
\lab{eq:LeStep4[B]_{ren}2}
(\ov{\DD}\c\Jk) B_b+\de_{cd}\ov{\DD}_c\Jk_b B_d - \ov{\DD}_b\Jk\c B &=&  \frac{8i(r^2+a^2)\cos\th}{|q|^4}B_b  +r^{-1}\Ga_b\c B,
\eea
and
\bea
\lab{eq:LeStep4[B]_{ren}3}
\bsplit
 \Jk_b\ov{\DD}\c B+\de_{cd}\Jk_d\ov{\DD}_c B_b - \Jk\c\ov{\DD}_b B &= \left(-\frac{1}{4}\ov{\Jk}\c(\DD\hot B)+\frac{3}{2}(\ov{\DD}\c B)\Jk+2\Re(\Jk)^c\nab_cB\right)_b.
 \end{split}
\eea

Thus, from now on, we focus on the proof of \eqref{eq:LeStep4[B]_{ren}1}, \eqref{eq:LeStep4[B]_{ren}2} and \eqref{eq:LeStep4[B]_{ren}3}. To prove \eqref{eq:LeStep4[B]_{ren}1}, we write
\beaa
\ov{\DD}\c(\Jk\hot B)_b &=& \de_{cd}\ov{\DD}_c(\Jk\hot B)_{db}= \de_{cd}(\ov{\DD}_c\Jk\hot B)_{db}+ \de_{cd}(\Jk\hot\ov{\DD}_c B)_{db}\\
&=& \de_{cd}\Big(\ov{\DD}_c\Jk_d B_b+\ov{\DD}_c\Jk_b B_d -\de_{bd}\ov{\DD}_c\Jk\c B\Big)\\
&+& \de_{cd}\Big(\Jk_b\ov{\DD}_c B_d+\Jk_d\ov{\DD}_c B_b -\de_{bd}\Jk\c\ov{\DD}_c B\Big)\\
&=& (\ov{\DD}\c\Jk) B_b+\de_{cd}\ov{\DD}_c\Jk_b B_d - \ov{\DD}_b\Jk\c B+ \Jk_b\ov{\DD}\c B+\de_{cd}\Jk_d\ov{\DD}_c B_b - \Jk\c\ov{\DD}_b B,
\eeaa
as stated.

Next, to prove   \eqref{eq:LeStep4[B]_{ren}2}, we write, since   $\Im(\Jk)=\Re(\dual\Jk)$,
\beaa
\ov{\DD}\c\Jk &=& (\nab-i\dual\nab)\c(\Re(\Jk)+i\Im(\Jk))=2\div(\Re(\Jk))+2i\curl(\Re(\Jk)),\\
\DD\hot\Jk &=& (\nab+i\dual\nab)\hot(\Re(\Jk)+i\Im(\Jk)) = 2\nab\hot\Re(\Jk)+2i\dual\nab\Re(\Jk),
\eeaa
so that
\beaa
\nab\hot\Re(\Jk)\in r^{-1}\Ga_b, \quad \div(\Re(\Jk))\in r^{-1}\Ga_b, \quad \curl(\Re(\Jk))=\frac{2(r^2+a^2)\cos\th}{|q|^4}+r^{-1}\Ga_b.
\eeaa
This yields $\nab_a\Re(\Jk)_b=\frac{(r^2+a^2)\cos\th}{|q|^4}\in_{ab}+r^{-1}\Ga_b$
and hence
\beaa
\ov{\DD}_a\Jk_b &=& (\nab_a-i\in_{ac}\nab_c)(\Re(\Jk)_b+i\in_{bd}\Re(\Jk)_d)\\
&=& \nab_a\Re(\Jk)_b+\in_{ac}\in_{bd}\nab_c\Re(\Jk)_d+i\Big(\in_{bd}\nab_a\Re(\Jk)_d-\in_{ac}\nab_c\Re(\Jk)_b\Big)\\
&=& \frac{(r^2+a^2)\cos\th}{|q|^4}\Big[\in_{ab}+\in_{ac}\in_{bd}\in_{cd}+i\big(\in_{bd}\in_{ad}-\in_{ac}\in_{cb}\big)\Big]
+r^{-1}\Ga_b\\
&=& \frac{2(r^2+a^2)\cos\th}{|q|^4}\Big[\in_{ab}+i\de_{ab}\Big]+r^{-1}\Ga_b.
\eeaa
We infer
\beaa
&&(\ov{\DD}\c\Jk) B_b+\de_{cd}\ov{\DD}_c\Jk_b B_d - \ov{\DD}_b\Jk\c B_b\\
 &=&  \frac{4i(r^2+a^2)\cos\th}{|q|^4}B_b +\de_{cd}\frac{2(r^2+a^2)\cos\th}{|q|^4}\Big[\in_{cb}+i\de_{cb}\Big]B_d\\
&& - \frac{2(r^2+a^2)\cos\th}{|q|^4}\Big[\in_{bc}+i\de_{bc}\Big]B_c +r^{-1}\Ga_b\c B\\
 &=&  \frac{4i(r^2+a^2)\cos\th}{|q|^4}B_b +\frac{2(r^2+a^2)\cos\th}{|q|^4}\Big[\in_{db}+i\de_{db}\Big]B_d\\
&& - \frac{2(r^2+a^2)\cos\th}{|q|^4}\Big[\in_{bc}+i\de_{bc}\Big]B_c +r^{-1}\Ga_b\c B\\
 &=&  \frac{4i(r^2+a^2)\cos\th}{|q|^4}B_b -\frac{4(r^2+a^2)\cos\th}{|q|^4}\dual B_b +r^{-1}\Ga_b\c B\\
 \eeaa
 and hence, since $\dual B= -iB$, we obtain 
\beaa 
(\ov{\DD}\c\Jk) B_b+\de_{cd}\ov{\DD}_c\Jk_b B_d - \ov{\DD}_b\Jk\c B_b &=&  \frac{8i(r^2+a^2)\cos\th}{|q|^4}B_b  +r^{-1}\Ga_b\c B
\eeaa
as stated in \eqref{eq:LeStep4[B]_{ren}2}.

Next, to prove  \eqref{eq:LeStep4[B]_{ren}3}, we write
\bea
\lab{eq:LeStep4[B]_{ren}4}
\bsplit
& \Jk_b\ov{\DD}\c B+\de_{cd}\Jk_d\ov{\DD}_c B_b - \Jk\c\ov{\DD}_b B -\left(-\frac{1}{4}\ov{\Jk}\c(\DD\hot B)+\frac{3}{2}(\ov{\DD}\c B)\Jk\right)_b\\
&= -\frac{1}{2}(\ov{\DD}\c B)\Jk_b+\frac{1}{4}(\ov{\Jk}\c(\DD\hot B))_b +\de_{cd}\Jk_d\ov{\DD}_c B_b - \Jk\c\ov{\DD}_b B.
\end{split}
\eea
Next we  make use of  the following identity
\bea
 \lab{eq:LeStep4[B]_{ren}5} 
-\frac{1}{2}(\ov{\DD}\c B)\Jk_b+\frac{1}{4}(\ov{\Jk}\c(\DD\hot B))_b +\de_{cd}\Jk_d\ov{\DD}_c B_b - \Jk\c\ov{\DD}_b B &=& 2\Re(\Jk)^c\nab_cB_b.
\eea
Using   \eqref{eq:LeStep4[B]_{ren}5},  the  identity   \eqref{eq:LeStep4[B]_{ren}4}  becomes
\beaa
 \Jk_b\ov{\DD}\c B+\de_{cd}\Jk_d\ov{\DD}_c B_b - \Jk\c\ov{\DD}_b B -\left(-\frac{1}{4}\ov{\Jk}\c(\DD\hot B)+\frac{3}{2}(\ov{\DD}\c B)\Jk\right)_b= 2\Re(\Jk)^c\nab_cB_b
\eeaa
which we rewrite in the form
\beaa
 \Jk_b\ov{\DD}\c B+\de_{cd}\Jk_d\ov{\DD}_c B_b - \Jk\c\ov{\DD}_b B &=& \left(-\frac{1}{4}\ov{\Jk}\c(\DD\hot B)+\frac{3}{2}(\ov{\DD}\c B)\Jk+2\Re(\Jk)^c\nab_cB\right)_b
\eeaa
as stated in   \eqref{eq:LeStep4[B]_{ren}3}.

 It only  remains to check  the identity  \eqref{eq:LeStep4[B]_{ren}5}. Since $B=\b+i\dual\b$ and $\Jk_2=-i\Jk_1$, we have
\beaa
&&\left(-\frac{1}{2}(\ov{\DD}\c B)\Jk_b+\frac{1}{4}(\ov{\Jk}\c(\DD\hot B))_b +\de_{cd}\Jk_d\ov{\DD}_c B_b - \Jk\c\ov{\DD}_b B\right)_{b=1}\\
&=& -\Jk_1(\div\b+i\curl\b)+\frac{1}{4}\ov{\Jk_1}\Big((\nab_1+i\nab_2)(\b_1+i\b_2)-(\nab_2-i\nab_1)(\b_2-i\b_1)\Big)\\
&&+\frac{1}{4}\ov{\Jk_2}\Big((\nab_1+i\nab_2)(\b_2-i\b_1)+(\nab_2-i\nab_1)(\b_1+i\b_2)\Big)\\
&&+\Jk_1(\nab_1-i\nab_2)(\b_1+i\b_2)+\Jk_2(\nab_2+i\nab_1)(\b_1+i\b_2)\\
&&-\Jk_1(\nab_1-i\nab_2)(\b_1+i\b_2) -\Jk_2(\nab_1-i\nab_2)(\b_2-i\b_1)\\
&=& -\Jk_1(\div\b+i\curl\b)\\
&&+\frac{1}{2}\ov{\Jk_1}\Big(\nab_1\b_1-\nab_2\b_2+i(\nab_1\b_2+\nab_2\b_1)\Big)+\frac{i}{2}\ov{\Jk_1}\Big(\nab_1\b_2+\nab_2\b_1+i(\nab_2\b_2-\nab_1\b_1)\Big)\\
&&+\Jk_1\Big(\div\b+i\curl\b\Big) -i\Jk_1\Big(-\curl\b+i\div\b\Big) \\
&& -\Jk_1\Big(\div\b +i\curl\b\Big) +i\Jk_1\Big(\curl\b-i\div\b\Big)\\
&=& \ov{\Jk_1}\Big(\nab_1\b_1-\nab_2\b_2+i(\nab_1\b_2+\nab_2\b_1)\Big)+\Jk_1\Big(\div\b+i\curl\b\Big).
\eeaa
Since $\Im(\Jk)=\dual\Re(\Jk)$, we obtain 
\beaa
&&\left(-\frac{1}{2}(\ov{\DD}\c B)\Jk_b+\frac{1}{4}(\ov{\Jk}\c(\DD\hot B))_b +\de_{cd}\Jk_d\ov{\DD}_c B_b - \Jk\c\ov{\DD}_b B\right)_{b=1}\\
&=& \Big(\Re(\Jk)_1-i\Re(\Jk)_2\Big)\Big(\nab_1\b_1-\nab_2\b_2+i(\nab_1\b_2+\nab_2\b_1)\Big)\\
&&+\Big(\Re(\Jk)_1+i\Re(\Jk)_2\Big)\Big(\div\b+i\curl\b\Big)\\
&=& 2\Re(\Jk)_1\Big(\nab_1\b_1+i\nab_1\b_2\Big)+2\Re(\Jk)_2\Big(\nab_2\b_1+i\nab_2\b_2\Big)\\
&=& \Big(2\Re(\Jk)^c\nab_cB_b\Big)_{b=1}.
\eeaa
Since the complex tensors on the LHS and RHS verify $V_2=-iV_1$, the equality holds also for $b=2$ and hence we have obtained 
\beaa
-\frac{1}{2}(\ov{\DD}\c B)\Jk_b+\frac{1}{4}(\ov{\Jk}\c(\DD\hot B))_b +\de_{cd}\Jk_d\ov{\DD}_c B_b - \Jk\c\ov{\DD}_b B &=& 2\Re(\Jk)^c\nab_cB_b
\eeaa
as stated in \eqref{eq:LeStep4[B]_{ren}5}. This concludes the proof of Lemma \ref{Le:Step4[B]_{ren}}.
\end{proof}

{\bf Step 5.} Next, we prove the following corollary.
\begin{corollary}
\lab{corr:Step5[B]_{ren}}
We have
\beaa
&&\Big(\nab_4-a\Re(\Jk)^b\nab_b\Big)\left(r^4 [B]_{red} \right) \\
&=& \frac{r^4}{2}\ov{\DD}'\c\Big( A  -a(\Jk\hot B)\Big)  + O(1)\Xh+O(r)\dk^{\leq 1}B+O(r^2)\nab_3B+O(r)\dkb^{\leq 1}\Pc +O(1)\trXc\\
&&+O(r^{2})\nab_3A+O(r)\dk^{\leq 1}A+r^4\Ga_g\c(B,A) +r^4\Ga_b\c\nab_3 A+r^2\dk^{\leq 1}(\Ga_g\c\Ga_g).
\eeaa
\end{corollary}

\begin{proof}
 We combine  Corollary  \ref{Corr:Step3[B]_{ren}} with Lemma \ref{Le:Step4[B]_{ren}} to derive
\beaa
&&\Big(\nab_4-a\Re(\Jk)^b\nab_b\Big)[B]_{ren}  +\frac{4}{r}[B]_{ren}
= \frac{1}{2}\ov{\DD}'\c A -\frac{a}{2}\ov{\DD}\c(\Jk\hot B) + \frac{3a^2}{4r}(\ov{\Jk}\c B)\Jk \\
&& + O(r^{-4})\Xh+O(r^{-3})B+O(r^{-3})\dkb^{\leq 1}\Pc +O(r^{-4})\trXc\\
&&+O(r^{-2})\nab_3A+O(r^{-3})\dk^{\leq 1}A+\Ga_g\c(B,A)+\Ga_b\c\nab_3 A+r^{-2}\Ga_g\c\Ga_g.
\eeaa

Also, recalling \eqref{eq:Step2[B]_{ren}1}
\beaa
\DD_a' V_B &=& \DD_a V_B+\frac{1}{2}\underline{F}_a f^c\nab_c V_B+\frac{1}{2}\underline{F}_a\nab_4V_B+\left(\frac{1}{2}F_a+\frac{1}{8}|f|^2\underline{F}_a\right)\nab_3V_B\\
&&+O(r^{-3})V+r^{-1}\Ga_bV,
\eeaa
and applying  it to $V=\Jk\hot B$ we  obtain
\beaa
\ov{\DD}'\c(\Jk\hot B) &=& \ov{\DD}(\Jk\hot B)+\frac{1}{2}\ov{\underline{F}}^{c}f^c\nab_c(\Jk\hot B)_a+\frac{1}{2}\ov{\underline{F}}\c\nab_4(\Jk\hot B)\\
&&+\left(\frac{1}{2}\ov{F}+\frac{1}{8}|f|^2\ov{\underline{F}}\right)\c\nab_3(\Jk\hot B)+\Big(O(r^{-4})+r^{-2}\Ga_b\Big)B.
\eeaa
We deduce
\beaa
&&\Big(\nab_4-a\Re(\Jk)^b\nab_b\Big)[B]_{ren}  +\frac{4}{r}[B]_{ren}\\ 
&=& \frac{1}{2}\ov{\DD}'\c\Big( A  -a(\Jk\hot B)\Big) +\frac{a}{4}\ov{\underline{F}}^{c}f^c\nab_c(\Jk\hot B)_a+\frac{a}{4}\ov{\underline{F}}\c\nab_4(\Jk\hot B)\\
&& +\frac{a}{4}\left(\ov{F}+\frac{1}{4}|f|^2\ov{\underline{F}}\right)\c\nab_3(\Jk\hot B) \\
&& + O(r^{-4})\Xh+O(r^{-3})B+O(r^{-3})\dkb^{\leq 1}\Pc +O(r^{-4})\trXc\\
&&+O(r^{-2})\nab_3A+O(r^{-3})\dk^{\leq 1}A+\Ga_g\c(B,A)+\Ga_b\c\nab_3 A+r^{-2}\Ga_g\c\Ga_g.
\eeaa
Since $e_4(r)=1$ and $e_1(r)=e_2(r)=0$, we infer
\beaa
&&\Big(\nab_4-a\Re(\Jk)^b\nab_b\Big)\left(r^4 [B]_{ren} \right) \\ 
&=& \frac{r^4}{2}\ov{\DD}'\c\Big( A  -a(\Jk\hot B)\Big)  +\frac{ar^4}{4}\ov{\underline{F}}^{c}f^c\nab_c(\Jk\hot B)_a+\frac{ar^4}{4}\ov{\underline{F}}\c\nab_4(\Jk\hot B)\\
&& +\frac{ar^4}{4}\left(\ov{F}+\frac{1}{4}|f|^2\ov{\underline{F}}\right)\c\nab_3(\Jk\hot B) \\
&& + O(1)\Xh+O(r)B+O(r)\dkb^{\leq 1}\Pc +O(1)\trXc\\
&&+O(r^{2})\nab_3A+O(r)\dk^{\leq 1}A+r^4\Ga_g\c(B,A)+r^4\Ga_b\c\nab_3 A+r^2\Ga_g\c\Ga_g.
\eeaa
Also, we have
\beaa
f=\Big(-1+O(r^{-1})+r\Ga_b\Big)\nab u, \qquad \fb=\Big(-1+O(r^{-1})+r\Ga_b\Big)\nab u,
\eeaa
and hence
\bea
\lab{eq:asympt-ffb-a}
&& f=-a\Re(\Jk)+O(r^{-1})\Re(\Jk)+\Ga_b, \quad \fb=-a\Re(\Jk)+O(r^{-1})\Re(\Jk)+\Ga_b,  
\eea
which yields
\beaa
&&  \frac{ar^4}{4}\ov{\underline{F}}^{c}f^c\nab_c(\Jk\hot B)_a+\frac{ar^4}{4}\ov{\underline{F}}\c\nab_4(\Jk\hot B) +\frac{ar^4}{4}\left(\ov{F}+\frac{1}{4}|f|^2\ov{\underline{F}}\right)\c\nab_3(\Jk\hot B)\\
&=& O(r)\dk^{\leq 1}B+O(r^2)\nab_3B+r^2\Ga_b\c\dk^{\leq 1}B +r^3\Ga_b\nab_3B
\eeaa
and hence
\beaa
&&\Big(\nab_4-a\Re(\Jk)^b\nab_b\Big)\left(r^4 [B]_{ren} \right)\\
&=& \frac{r^4}{2}\ov{\DD}'\c\Big( A  -a(\Jk\hot B)\Big)  + O(1)\Xh+O(r)\dk^{\leq 1}B+O(r^2)\nab_3B+O(r)\dkb^{\leq 1}\Pc +O(1)\trXc\\
&&+O(r^{2})\nab_3A+O(r)\dk^{\leq 1}A+r^4\Ga_g\c(B,A) +r^3\Ga_b\nab_3B+r^4\Ga_b\c\nab_3 A+r^2\Ga_g\c\Ga_g.
\eeaa
Using the Bianchi identity for $\nab_3B$, i.e.
\beaa
\nab_3B -\DD\ov{\Pc} &=& \frac{2}{r}B+O(r^{-2}) B +O(r^{-2})\Pc+O(r^{-3})\Hc +O(r^{-4})\widecheck{\DD(\cos\th)}\\
&&+r^{-1}\Ga_b\c\Ga_g,
\eeaa
we infer $\nab_3B=r^{-2}\dkb^{\leq 1}\Ga_g$ and hence $r^3\Ga_b\nab_3B=r\Ga_b\dkb^{\leq 1}\Ga_g$ which is included in error terms of the type $r^2\dk^{\leq 1}(\Ga_g\c\Ga_g)$. We deduce
\beaa
&&\Big(\nab_4-a\Re(\Jk)^b\nab_b\Big)\left(r^4 [B]_{ren} \right)\\
&=& \frac{r^4}{2}\ov{\DD}'\c\Big( A  -a(\Jk\hot B)\Big) +    O(1)\Xh+O(r)\dk^{\leq 1}B+O(r^2)\nab_3B+O(r)\dkb^{\leq 1}\Pc +O(1)\trXc\\
&+&O(r^{2})\nab_3A+O(r)\dk^{\leq 1}A+r^4\Ga_g\c(B,A) +r^4\Ga_b\c\nab_3 A+r^2\dk^{\leq 1}(\Ga_g\c\Ga_g)
\eeaa
as desired. This concludes the proof of Corollary \ref{corr:Step5[B]_{ren}}.
\end{proof}

{\bf Step 6.} Next, we derive the following lemma.
\begin{lemma}
 \lab{Le:Step6[B]_{ren}}
For a complex horizontal 1-form $U$,  we have
\beaa
\, \left[\nab_4-a\Re(\Jk)^b\nab_b, r[\DD]_{ren} \right]U= O(r^{-3})\dk^{\leq 1}U+O(r^{-2})\nab_3U +\Ga_g\c\dk^{\leq 1}U,
\eeaa
where we recall that  $ [\DD]_{ren}=\left(\ov{\DD}\c -\frac{a}{2}\ov{\Jk}\c\nab_4 -\frac{a}{2}\ov{\Jk}\c\nab_3\right)$.
\end{lemma}

\begin{proof}
We compute, using  $e_1(r)=e_2(r)=0$, 
\beaa
&&\left[\nab_4-a\Re(\Jk)^b\nab_b, r[\DD]_{red}\right]\\
&=& [\nab_4, r\ov{\DD}\c] -ar[\Re(\Jk)^b\nab_b, \ov{\DD}\c] -\frac{a}{2}\nab_4(r\ov{\Jk})\c\nab_4 -\frac{a}{2}\nab_4(r\ov{\Jk})\c\nab_3 -\frac{a}{2}r\ov{\Jk}\c[\nab_4, \nab_3]\\
&&+\frac{a^2}{2}r\left[\Re(\Jk)^b\nab_b, \ov{\Jk}\c(\nab_4+\nab_3)\right].
\eeaa
Since $\nab_4(\ov{q}\ov{\Jk})=0$, and, for a complex 1-form $U$, 
\beaa
\,[\nab_4, \nab_3]U &=& 2\omb\nab_4U+2(\eta+\ze)\c\nab U+4\eta\hot(\ze\c U) -4\ze\hot(\eta\c U)+4i\dual\rho U,\\
\,[\nab_4, \ov{\DD}\c]U &=& -\frac{1}{\ov{q}}\ov{\DD}\c -\frac{1}{\ov{q}}\ov{Z}\c U+r^{-1}\Ga_g\c\dk^{\leq 1}U,
\eeaa 
we infer, for a complex 1-form $U$, 
\beaa
&&\left[\nab_4-a\Re(\Jk)^b\nab_b, r[\DD]_{red}\right]U\\
&=& \left(1-\frac{r}{\ov{q}}\right)\ov{\DD}\c U -\frac{r}{\ov{q}}\ov{Z}\c U -ar[\Re(\Jk)^b\nab_b, \ov{\DD}\c]U +\frac{a^2}{2}r\left[\Re(\Jk)^b\nab_b, \ov{\Jk}\c(\nab_4+\nab_3)\right]U\\
&& -\frac{a}{2}r\ov{\Jk}\c\Big(2\omb\nab_4U+2(\eta+\ze)\c\nab U+4\eta\hot(\ze\c U) -4\ze\hot(\eta\c U)+4i\dual\rho U\Big) +\Ga_g\c\dk^{\leq 1}U\\
&&-\frac{a}{2}\left(1-\frac{r}{\ov{q}}\right)\ov{\Jk}\c\nab_4 -\frac{a}{2}\left(1-\frac{r}{\ov{q}}\right)\ov{\Jk}\c\nab_3.
\eeaa
We infer
\beaa
&&\left[\nab_4-a\Re(\Jk)^b\nab_b, r[\DD]_{red}\right]U\\
&=& -\frac{ia\cos\th}{r}\ov{\DD}\c U -\frac{a}{r}\ov{\Jk}\c U -ar[\Re(\Jk)^b\nab_b,\ov{\DD}\c]U +O(r^{-3})\dk^{\leq 1}U+O(r^{-2})\nab_3U +\Ga_g\c\dk^{\leq 1}U.
\eeaa

Next, we compute
\beaa
[\Re(\Jk)^b\nab_b, \ov{\DD}\c]U &=& -\ov{\DD}_a\Re(\Jk)_b \nab_bU_a+\Re(\Jk)^b[\nab_b, \ov{\DD}_a]U_a.
\eeaa
Recall that we have
\beaa
\nab_a\Re(\Jk)_b=\frac{(r^2+a^2)\cos\th}{|q|^4}\in_{ab}+r^{-1}\Ga_b
\eeaa
and hence
\beaa
\ov{\DD}_a\Re(\Jk)_b &=& (\nab_a-i\in_{ac}\nab_c)\Re(\Jk)_b\\
&=& \frac{(r^2+a^2)\cos\th}{|q|^4}\Big[\in_{ab}-i\in_{ac}\in_{cb}\Big]
+r^{-1}\Ga_b\\
&=& \frac{(r^2+a^2)\cos\th}{|q|^4}\Big[\in_{ab}+i\de_{ab}\Big]+r^{-1}\Ga_b\\
&=& \frac{\cos\th}{r^2}\Big[\in_{ab}+i\de_{ab}\Big]+O(r^{-4})+r^{-1}\Ga_b.
\eeaa
We infer
\beaa
-\ov{\DD}_a\Re(\Jk)_b \nab_bU_a &=& -\frac{\cos\th}{r^2}\Big[\in_{ab}+i\de_{ab}\Big]\nab_bU_a+O(r^{-5})\dkb U+r^{-2}\Ga_b\dkb U\\
&=& -\frac{\cos\th}{r^2}\Big[-\curl(U)+i\div(U)\Big]+O(r^{-5})\dkb U+r^{-2}\Ga_b\dkb U\\
&=& -\frac{i\cos\th}{r^2}\ov{\DD}\c U+O(r^{-5})\dkb U+r^{-2}\Ga_b\dkb U.
\eeaa
Also, we have, see Proposition 2.58 in \cite{GKS1}, 
\beaa
[\nab_b, \nab_a]U_c &=& \frac{1}{2}\in_{ba}\left(\atrch\nab_3+\atrchb\nab_4\right)U_c -\frac{1}{2}E_{cdba}U_d -\in_{cd}\in_{ba}\rho U_d
\eeaa
where
\beaa
E_{abcd} &:=& \chi_{ac}\chib_{bd}+\chib_{ac}\chi_{bd} -\chi_{bc}\chib_{ad}  -\chib_{bc}\chi_{ad}.
\eeaa
Since
\beaa
\chi_{ab}=\frac{1}{r}\de_{ab}+O(r^{-2})+\Ga_g, \qquad \chib_{ab}=-\frac{1}{r}\de_{ab}+O(r^{-2})+\chibh_{ab},
\eeaa
we have\footnote{A priori, one also expects a contribution $r^{-1}\Ga_b$ coming from
\beaa
\de_{ac}\chibh_{bd}+\chibh_{ac}\de_{bd} -\de_{bc}\chibh_{ad}  -\chibh_{bc}\de_{ad},
\eeaa
but this tensor vanishes in fact identically.}
\beaa
E_{abcd} &=& -\frac{2}{r^2}\Big(\de_{ac}\de_{bd}-\de_{bc}\de_{ad}\Big) +O(r^{-3})+r^{-1}\Ga_g
\eeaa
and hence
\beaa
[\nab_b, \nab_a]U_c &=&  \frac{1}{r^2}\Big(\de_{cb}\de_{da}-\de_{db}\de_{ca}\Big) U_d \\
&&+O(r^{-3})\dk^{\leq 1}U+O(r^{-2})\nab_3U+r^{-1}\Ga_g\dk^{\leq 1} U+\Ga_g\nab_3U.
\eeaa
We infer
\beaa
[\nab_b, \ov{\DD}_a]U_a &=& \de_{ac}[\nab_b, \nab_a]U_c -i\in_{ca}[\nab_b, \nab_a]U_c\\
&=& \frac{1}{r^2}(\de_{ac}-i\in_{ca})\Big(\de_{cb}\de_{da}-\de_{db}\de_{ca}\Big) U_d \\
&&+O(r^{-3})\dk^{\leq 1}U+O(r^{-2})\nab_3U+r^{-1}\Ga_g\dk^{\leq 1} U+\Ga_g\nab_3U\\
&=& \frac{1}{r^2}\Big(-\de_{db}-i\in_{bd}\Big)U_d\\
&&+O(r^{-3})\dk^{\leq 1}U+O(r^{-2})\nab_3U+r^{-1}\Ga_g\dk^{\leq 1} U+\Ga_g\nab_3U\\
&=& -\frac{1}{r^2}\Big(U_b+i\dual U_b\Big) +O(r^{-3})\dk^{\leq 1}U+O(r^{-2})\nab_3U+r^{-1}\Ga_g\dk^{\leq 1} U+\Ga_g\nab_3U
\eeaa
and hence, using also $\Im(\Jk)=\dual\Re(\Jk)$, we obtain
\beaa
\Re(\Jk)^b[\nab_b, \ov{\DD}_a]U_a &=& -\Re(\Jk)^b\frac{1}{r^2}\Big(U_b+i\dual U_b\Big) \\
&&+O(r^{-4})\dk^{\leq 1}U+O(r^{-3})\nab_3U+r^{-2}\Ga_g\dk^{\leq 1} U+r^{-1}\Ga_g\nab_3U\\
&=& -\frac{1}{r^2}\Big(\Re(\Jk)\c U+i\Re(\Jk)\c\dual U\Big)\\
&&+O(r^{-4})\dk^{\leq 1}U+O(r^{-3})\nab_3U+r^{-2}\Ga_g\dk^{\leq 1} U+r^{-1}\Ga_g\nab_3U\\
&=& -\frac{1}{r^2}\Big(\Re(\Jk)\c U-i\dual\Re(\Jk)\c U\Big)\\
&&+O(r^{-4})\dk^{\leq 1}U+O(r^{-3})\nab_3U+r^{-2}\Ga_g\dk^{\leq 1} U+r^{-1}\Ga_g\nab_3U\\
&=& -\frac{1}{r^2}\ov{\Jk}\c U +O(r^{-4})\dk^{\leq 1}U+O(r^{-3})\nab_3U+r^{-2}\Ga_g\dk^{\leq 1} U+r^{-1}\Ga_g\nab_3U.
\eeaa
In view of the above, we deduce
\beaa
[\Re(\Jk)^b\nab_b, \ov{\DD}\c]U &=& -\ov{\DD}_a\Re(\Jk)_b \nab_bU_a+\Re(\Jk)^b[\nab_b, \ov{\DD}_a]U_a= -\frac{i\cos\th}{r^2}\ov{\DD}\c U -\frac{1}{r^2}\ov{\Jk}\c U\\
&&+O(r^{-4})\dk^{\leq 1}U+O(r^{-3})\nab_3U+r^{-2}\Ga_b\dk^{\leq 1} U+r^{-1}\Ga_g\nab_3U.
\eeaa
This yields
\beaa
&& -\frac{ia\cos\th}{r}\ov{\DD}\c U -\frac{a}{r}\ov{\Jk}\c U -ar[\Re(\Jk)^b\nab_b,\ov{\DD}\c]U\\
&=& -\frac{ia\cos\th}{r}\ov{\DD}\c U -\frac{a}{r}\ov{\Jk}\c U -ar\left(-\frac{i\cos\th}{r^2}\ov{\DD}\c U -\frac{1}{r^2}\ov{\Jk}\c U\right)\\
&&+O(r^{-3})\dk^{\leq 1}U+O(r^{-2})\nab_3U+r^{-1}\Ga_b\dk^{\leq 1} U+\Ga_g\nab_3U\\
&=& O(r^{-3})\dk^{\leq 1}U+O(r^{-2})\nab_3U+r^{-1}\Ga_b\dk^{\leq 1} U+\Ga_g\nab_3U.    
\eeaa
Coming back to 
\beaa
\left[\nab_4-a\Re(\Jk)^b\nab_b, r[\DD]_{red}\right]U
&=& -\frac{ia\cos\th}{r}\ov{\DD}\c U -\frac{a}{r}\ov{\Jk}\c U -ar[\Re(\Jk)^b\nab_b,\ov{\DD}\c]U \\
&&+O(r^{-3})\dk^{\leq 1}U+O(r^{-2})\nab_3U +\Ga_g\c\dk^{\leq 1}U,
\eeaa
we infer
\beaa
\left[\nab_4-a\Re(\Jk)^b\nab_b, r[\DD]_{red}\right]U= O(r^{-3})\dk^{\leq 1}U+O(r^{-2})\nab_3U +\Ga_g\c\dk^{\leq 1}U
\eeaa
as desired. This concludes the proof of Lemma \ref{Le:Step6[B]_{ren}}.
\end{proof}

{\bf Step 7.} We  are finally ready to prove Lemma  \ref{Le:mainpointwise-ell=1B}. We start with the identity of Corollary \ref{corr:Step5[B]_{ren}}
\beaa
&&\Big(\nab_4-a\Re(\Jk)^b\nab_b\Big)\left[r^4 [B]_{red} \right]\\
 &=& \frac{r^4}{2}\ov{\DD}'\c\Big( A  -a(\Jk\hot B)\Big)   + O(1)\Xh+O(r)\dk^{\leq 1}B+O(r^2)\nab_3B+O(r)\dkb^{\leq 1}\Pc +O(1)\trXc\\
&&+O(r^{2})\nab_3A+O(r)\dk^{\leq 1}A+r^4\Ga_g\c(B,A) +r^4\Ga_b\c\nab_3 A+r^2\dk^{\leq 1}(\Ga_g\c\Ga_g).
\eeaa
According to  Lemma  \ref{Le:Step6[B]_{ren}}  we have, with $U=  A  -a(\Jk\hot B)$,
\beaa
\, \left[\nab_4-a\Re(\Jk)^b\nab_b, r[\DD]_{ren} \right]U= O(r^{-3})\dk^{\leq 1}U+O(r^{-2})\nab_3U +\Ga_g\c\dk^{\leq 1}U.
\eeaa
We deduce
\bea
\lab{eq:Step7[B]_{ren}} 
\nn&&\Big(\nab_4-a\Re(\Jk)^b\nab_b\Big)       \Big(    r [\ov{\DD}\c]_{ren}  \big(r^4 [B]_{ren} \big)\Big)\\
\nn&=& r[\DD]_{ren}      \left[\frac{r^4}{2}\ov{\DD}'\c\Big( A  -a(\Jk\hot B)\Big)\right] \\
\nn&& + O(1)\dk^{\leq 1}\Xh+O(r)\dk^{\leq 2}B+O(r^2)\dk^{\leq 1}\nab_3B+O(r)\dk^{\leq 2}\Pc +O(1)\dk^{\leq 1}\trXc\\
\nn&&+O(r^{2})\dk^{\leq 1}\nab_3A+O(r)\dk^{\leq 2}A+r^4\dk^{\leq 1}\big(\Ga_g\c(B,A)\big) +r^3\dk^{\leq 1}\big(\Ga_b\nab_3B\big)\\
&&+r^4\dk^{\leq 1}\big(\Ga_b\c\nab_3 A\big)+r^2\dk^{\leq 1}\big(\Ga_g\c\Ga_g\big).
\eea
Next we note that
\bea
r\ov{\DD'}\c &=&  r[\DD]_{ren} +O(r^{-2})\dk^{\leq 1}+O(r^{-1})\nab_3+r^{-1}\Ga_b\dk^{\leq 1}+\Ga_b\nab_3.
\eea
This follows in view of the formula \eqref{eq:Step2[B]_{ren}1} 
\beaa
\DD_a' V_B &=& \DD_a V_B+\frac{1}{2}\underline{F}_a f^c\nab_c V_B+\frac{1}{2}\underline{F}_a\nab_4V_B+\left(\frac{1}{2}F_a+\frac{1}{8}|f|^2\underline{F}_a\right)\nab_3V_B\\
&&+\Big\{O(r^{-3})+ r^{-1} \Ga_b \Big\} V
\eeaa
from which we infer, using  also  \eqref{eq:asympt-ffb-a},
\beaa
r\ov{\DD'}\c &=& r\left(\ov{\DD}\c -\frac{a}{2}\ov{\Jk}\c\nab_4 -\frac{a}{2}\ov{\Jk}\c\nab_3\right)+O(r^{-2})\dk^{\leq 1}+O(r^{-1})\nab_3+r^{-1}\Ga_b\dk^{\leq 1}+\Ga_b\nab_3\\
&=& r[\DD]_{ren} +O(r^{-2})\dk^{\leq 1}+O(r^{-1})\nab_3+r^{-1}\Ga_b\dk^{\leq 1}+\Ga_b\nab_3
\eeaa
as stated above.

We deduce
\beaa
&& r[\DD]_{ren}\left[\frac{r^4}{2}\ov{\DD}'\c\Big( A  -a(\Jk\hot B)\Big)\right]\\
&=& r\ov{\DD'}\c\left[\frac{r^4}{2}\ov{\DD}'\c\Big( A  -a(\Jk\hot B)\Big)\right] + O(r)\dk^{\leq 2}A+O(r^2)\dk^{\leq 1}\nab_3A+r^2\Ga_b\dk^{\leq 1}A\\
&&+r^3\Ga_b\dk^{\leq 1}\nab_3A + O(1)\dk^{\leq 2}B+O(r)\dk^{\leq 1}\nab_3B+r\Ga_b\dk^{\leq 1}B+r^2\Ga_b\dk^{\leq 1}\nab_3B.
\eeaa
Inserting this on the right hand side of \eqref{eq:Step7[B]_{ren}}, since 
$e_1'(r)=e_2'(r)=0$,  we deduce
\beaa
&&\Big(\nab_4-a\Re(\Jk)^b\nab_b\Big)       \Big(    r [\ov{\DD}\c]_{ren}  \big(r^4 [B]_{ren} \big)\Big)\\
 &=& \frac{r^5}{2}\ov{\DD}'\c\ov{\DD}'\c\Big( A  -a(\Jk\hot B)\Big)    + O(1)\dk^{\leq 1}\Xh+O(r)\dk^{\leq 2}B+O(r^2)\dk^{\leq 1}\nab_3B+O(r)\dk^{\leq 2}\Pc \\
 &&+O(1)\dk^{\leq 1}\trXc +O(r^{2})\dk^{\leq 1}\nab_3A+O(r)\dk^{\leq 2}A+r^4\dk^{\leq 1}\big(\Ga_g\c(B,A)\big) +r^3\dk^{\leq 1}\big(\Ga_b\nab_3B\big)\\
&&+r^4\dk^{\leq 1}\big(\Ga_b\c\nab_3 A\big)+r^2\dk^{\leq 2}\big(\Ga_g\c\Ga_g\big).
\eeaa
Using, as before, that $\Ga_b\nab_3B$ in included in $r^{-1}\dk^{\leq 1}(\Ga_g\c\Ga_g)$, we deduce
\beaa
&&\Big(\nab_4-a\Re(\Jk)^b\nab_b\Big)       \Big(    r [\ov{\DD}\c]_{ren}  \big(r^4 [B]_{ren} \big)\Big)\\
&=& \frac{r^5}{2}\ov{\DD}'\c\ov{\DD}'\c\Big( A  -a(\Jk\hot B)\Big)   + O(1)\dk^{\leq 1}\Xh+O(r)\dk^{\leq 2}B+O(r^2)\dk^{\leq 1}\nab_3B\\
&&+O(r)\dk^{\leq 2}\Pc +O(1)\dk^{\leq 1}\trXc+O(r^{2})\dk^{\leq 1}\nab_3A+O(r)\dk^{\leq 2}A+r^4\dk^{\leq 1}\big(\Ga_g\c(B,A)\big) \\
&&+r^4\dk^{\leq 1}\big(\Ga_b\c\nab_3 A\big)+r^2\dk^{\leq 2}\big(\Ga_g\c\Ga_g\big)
\eeaa
as desired. This concludes the proof of the Lemma  \ref{Le:mainpointwise-ell=1B}.

%%%%%%%%%%%%%%%%%%%%%%%%%%%%%%%%%%%%%%

%%%%%%%%%%%%%%%%%%%%%%%%%%%%%%%%%%%%%%

%%%%%%%%%%%%%%%%%%%%%%%%%%%%%%%%%

\chapter{Proof of results in Chapter 9}

%%%%%%%%%%%%%%%%%%%%%%%%%%%%%%%%%

%%%%%%%%%%%%%%%%%%%%%%%%%%%%%%%%

\section{Proof of Proposition \ref{Prop:linearizedPTstructure1}}
\lab{appendix:ProofProp{Prop:linearizedPTstructure1}}

%%%%%%%%%%%%%%%%%%%%%%%%%%%%%%%%

Recall that we have
\beaa
\nab_4\tr X +\frac{1}{2}(\tr X)^2 &=& -\frac{1}{2}\Xh\c\ov{\Xh}=\Ga_g\c\Ga_g.
\eeaa
We infer
\beaa
\nab_4\left( \frac{2}{q}+\trXc\right) +\frac{1}{2}\left( \frac{2}{q}+\trXc\right)^2 &=& \Ga_g\c\Ga_g
\eeaa
and since $e_4(q)=1$, we deduce
\beaa
\nab_4\trXc +\frac{2}{q}\trXc &=& \Ga_g\c\Ga_g
\eeaa
as desired.

Next, recall 
\beaa
\nab_4Z +\frac{1}{2}\tr X Z &=&  \frac{1}{2}\tr X\Hb+\frac{1}{2}\widehat{X}\c(-\ov{Z}+\ov{\Hb}) -B.
\eeaa
We infer
\beaa
&&\nab_4\left(\frac{a\ov{q}}{|q|^2}\Jk +\Zc\right) +\frac{1}{2}\left( \frac{2}{q}+\trXc\right) \left(\frac{a\ov{q}}{|q|^2}\Jk +\Zc\right) \\
&=&  \frac{1}{2}\left( \frac{2}{q}+\trXc\right)\left(-\frac{a\ov{q}}{|q|^2}\Jk\right)+\frac{1}{2}\widehat{X}\c\left(-\ov{\frac{a\ov{q}}{|q|^2}\Jk}+\ov{-\frac{a\ov{q}}{|q|^2}\Jk}\right) -B+\Ga_g\c\Ga_g.
\eeaa
We have
\beaa
\nab_4\left(\frac{a\ov{q}}{|q|^2}\Jk\right) &=& \pr_r\left(\frac{a\ov{q}}{|q|^2}\right)\Jk+\frac{a\ov{q}}{|q|^2}\left(-\frac{1}{q}\Jk\right)=-\frac{2}{q}\frac{a\ov{q}}{|q|^2}\Jk
\eeaa
and hence
\beaa
\nab_4\Zc + \frac{1}{q}\Zc &=&   - \frac{a\ov{q}}{|q|^2}\trXc \Jk    -\frac{aq}{|q|^2}\ov{\Jk}\c\widehat{X}  -B+\Ga_g\c\Ga_g.
\eeaa

Next, recall  
\beaa
\nab_4H &=&  -\frac{1}{2}\ov{\tr X}(H-\Hb) -\frac{1}{2}\Xh\c(\ov{H}-\ov{\Hb}) -B.
\eeaa
We infer
\beaa
\nab_4\left(\frac{aq}{|q|^2}\Jk+\Hc\right) &=&  -\frac{1}{2}\left(\frac{2}{\ov{q}}+\ov{\trXc}\right)\left(\frac{aq}{|q|^2}\Jk+\Hc+\frac{a\ov{q}}{|q|^2}\Jk\right)\\
&& -\frac{1}{2}\Xh\c\left(\ov{\frac{aq}{|q|^2}\Jk}+\ov{\frac{a\ov{q}}{|q|^2}\Jk}\right) -B+\Ga_b\c\Ga_g.
\eeaa
Now, we have
\beaa
\nab_4\left(\frac{aq}{|q|^2}\Jk\right) &=& \pr_r\left(\frac{aq}{|q|^2}\right)\Jk+\frac{aq}{|q|^2}\left(-\frac{1}{q}\right)\Jk\\ 
&=& -\frac{a}{\ov{q}^2}\Jk -\frac{a}{|q|^2}\Jk 
\eeaa 
and hence
\beaa
\nab_4\Hc+ \frac{1}{\ov{q}}\Hc &=& -\frac{ar}{|q|^2}\ov{\trXc}\Jk  -\frac{ar}{|q|^2}\ov{\Jk}\c\Xh -B+\Ga_b\c\Ga_g.
\eeaa

Next, recall 
\beaa
\nab_4\omb  -(\eta-\etab)\c\ze +\eta\c\etab&=&   \rho
\eeaa
which we rewrite 
\beaa
\nab_4\omb  -\frac{1}{2}\Re\Big((H-\Hb)\c \ov{Z} -H\c\ov{\Hb}\Big) &=&   \Re(P).
\eeaa
We infer
\beaa
&& \frac{1}{2}\pr_r^2\left(\frac{\De}{|q|^2}\right)+\nab_4\ombc  -\frac{1}{2}\Re\left(\left(\frac{aq}{|q|^2}\Jk+\Hc+\frac{a\ov{q}}{|q|^2}\Jk\right)\c \ov{\frac{a\ov{q}}{|q|^2}\Jk+\Zc} +\left(\frac{aq}{|q|^2}\Jk+\Hc\right)\c\ov{\frac{a\ov{q}}{|q|^2}\Jk}\right)\\
 &=&   \Re\left(-\frac{2m}{q^3}+\Pc\right)
\eeaa
and hence
\beaa
 \nab_4\ombc  &=&  -\frac{1}{2}\pr_r^2\left(\frac{\De}{|q|^2}\right) +\frac{a^2}{2|q|^2}\Re\left(\Jk\c\ov{\Jk}\right) +\frac{a^2}{|q|^4}\Re\left(q^2\Jk\c\ov{\Jk}\right) + \Re\left(-\frac{2m}{q^3}\right)\\
 && + \Re\left(\Pc\right) +\frac{ar}{|q|^2}\Re\left(\Jk\c \ov{\Zc}\right)+\frac{2a}{|q|^2}\Re\left(q\ov{\Jk}\c \Hc\right)+\Ga_g\c\Ga_g.
\eeaa
Since $\Jk\c\ov{\Jk}=\frac{2(\sin\th)^2}{|q|^2}$, we infer
\beaa
 \nab_4\ombc  &=&   \Re\left(\Pc\right) +\frac{ar}{|q|^2}\Re\left(\Jk\c \ov{\Zc}\right)+\frac{2a}{|q|^2}\Re\left(q\ov{\Jk}\c \Hc\right)+\Ga_g\c\Ga_g.
\eeaa

Next, recall
\beaa
\nab_4\tr\Xb +\frac{1}{2}\tr X\tr\Xb  &=& \DD\c\ov{\Hb}+\Hb\c\ov{\Hb}+2\ov{P} -\frac{1}{2}\Xh\c\ov{\Xbh}.
\eeaa
We infer
\beaa
&&\nab_4\left(-\frac{2q\Delta}{|q|^4}+\trXbc\right)+\frac{1}{2}\left( \frac{2}{q}+\trXc\right)\left(-\frac{2q\Delta}{|q|^4}+\trXbc\right) \\
 &=& \DD\c\ov{\Hb}+\Hb\c\ov{\Hb}+2\left(-\frac{2m}{\ov{q}^3}+\ov{\Pc}\right) +\Ga_b\c\Ga_g.
\eeaa
Now, we have, using in particular $\Jk\c\ov{\Jk}=\frac{2(\sin\th)^2}{|q|^2}$, 
\beaa
\DD\c\ov{\Hb}+\Hb\c\ov{\Hb} &=& -\DD\c\left(\ov{\frac{a\ov{q}}{|q|^2}\Jk}\right)+\frac{a\ov{q}}{|q|^2}\Jk\c\ov{\frac{a\ov{q}}{|q|^2}\Jk}\\
&=& -\frac{aq}{|q|^2}\DD\c\ov{\Jk} +\frac{a}{\ov{q}^2}\DD(\ov{q})\c\ov{\Jk} +\frac{a^2}{|q|^2}\Jk\c\ov{\Jk}\\
&=& -\frac{aq}{|q|^2}\left(\frac{4i(r^2+a^2)\cos\th}{|q|^4}+\widecheck{\DD\c\ov{\Jk}}\right) +\frac{a}{\ov{q}^2}\DD(r)\c\ov{\Jk}\\
&& - \frac{ia^2}{\ov{q}^2}\Big(i\Jk+\widecheck{\DD(\cos\th)}\Big)\c\ov{\Jk} +\frac{a^2}{|q|^2}\Jk\c\ov{\Jk}\\
&=&  -\frac{4iaq(r^2+a^2)\cos\th}{|q|^6} +\frac{2(\sin\th)^2a^2(q^2+|q|^2)}{|q|^6} \\
&& -\frac{aq}{|q|^2}\widecheck{\DD\c\ov{\Jk}} +\frac{a}{\ov{q}^2}\DD(r)\c\ov{\Jk} - \frac{ia^2}{\ov{q}^2}\widecheck{\DD(\cos\th)}\c\ov{\Jk}. 
\eeaa
We deduce
\beaa
\nab_4\trXbc +\frac{1}{q}\trXbc  &=& 2\ov{\Pc} + \frac{q\Delta}{|q|^4}\trXc -\frac{aq}{|q|^2}\widecheck{\DD\c\ov{\Jk}} +\frac{a}{\ov{q}^2}\DD(r)\c\ov{\Jk} - \frac{ia^2}{\ov{q}^2}\widecheck{\DD(\cos\th)}\c\ov{\Jk} +\Ga_b\c\Ga_g.
\eeaa

Next, recall
\beaa
\nab_4\Xbh +\frac{1}{2}\tr X\, \widehat{\Xb}  &=& \DD\hot\Hb  +\Hb\hot\Hb -\frac{1}{2}\ov{\tr\Xb} \widehat{X}.
\eeaa
We have
\beaa
\DD\hot\Hb  +\Hb\hot\Hb &=& -\DD\hot\left(\frac{a\ov{q}}{|q|^2}\Jk\right)  +\frac{a\ov{q}}{|q|^2}\Jk\hot\frac{a\ov{q}}{|q|^2}\Jk\\
&=& -\frac{a\ov{q}}{|q|^2}\DD\hot\Jk +\frac{a}{q^2}\DD(q)\hot\Jk+\frac{a^2\ov{q}^2}{|q|^4}\Jk\hot\Jk\\
&=& -\frac{a\ov{q}}{|q|^2}\DD\hot\Jk +\frac{a}{q^2}\DD(r)\hot\Jk+\frac{ia^2}{q^2}\Big(i\Jk+\widecheck{\DD(\cos\th)}\Big)\hot\Jk+\frac{a^2\ov{q}^2}{|q|^4}\Jk\hot\Jk\\
&=& -\frac{a\ov{q}}{|q|^2}\DD\hot\Jk +\frac{a}{q^2}\DD(r)\hot\Jk+\frac{ia^2}{q^2}\widecheck{\DD(\cos\th)}\hot\Jk
\eeaa
and hence
\beaa
\nab_4\Xbh +\frac{1}{q}\widehat{\Xb}  &=&  -\frac{a\ov{q}}{|q|^2}\DD\hot\Jk +\frac{a}{q^2}\DD(r)\hot\Jk+\frac{ia^2}{q^2}\widecheck{\DD(\cos\th)}\hot\Jk +\frac{q\Delta}{|q|^4} \widehat{X}+\Ga_b\c\Ga_g.
\eeaa

Next, recall
\beaa
 \nab_4\Xib &=& \nab_3\Hb +\frac{1}{2}\ov{\tr\Xb}(\Hb-H) +\frac{1}{2}\Xbh\c(\ov{\Hb}-\ov{H}) -\Bb.
\eeaa
We infer
\beaa
 \nab_4\Xib &=& \nab_3\Hb +\frac{1}{2}\ov{-\frac{2q\Delta}{|q|^4}+\trXbc}\left(-\frac{2ar}{|q|^2}\Jk -\Hc\right) -\frac{ar}{|q|^2}\ov{\Jk}\c\Xbh -\Bb+\Ga_b\c\Ga_b.
\eeaa
Also, we have 
\beaa
\nab_3\Hb &=& \nab_3\left(-\frac{a\ov{q}}{|q|^2}\Jk\right)=-\frac{a\ov{q}}{|q|^2}\nab_3\Jk +\frac{a}{q^2}\nab_3(q)\Jk\\
&=& -\frac{a\ov{q}}{|q|^2}\left(\frac{\De q}{|q|^4}\Jk+\widecheck{\nab_3\Jk}\right)+\frac{a}{q^2}\Big(e_3(r)+iae_3(\cos\th)\Big)\Jk\\
&=&  -\frac{a\De(|q|^2+\ov{q}^2)}{|q|^6}\Jk  -\frac{a\ov{q}}{|q|^2}\widecheck{\nab_3\Jk}+\frac{a}{q^2}\left(\widecheck{e_3(r)}+iae_3(\cos\th)\right)\Jk.
\eeaa
We infer
\beaa
 \nab_4\Xib &=&  \frac{\ov{q}\Delta}{|q|^4}\Hc - \frac{ar}{|q|^2}\ov{\trXbc}\Jk  -\frac{ar}{|q|^2}\ov{\Jk}\c\Xbh -\Bb -\frac{a\ov{q}}{|q|^2}\widecheck{\nab_3\Jk}+\frac{a}{q^2}\left(\widecheck{e_3(r)}+iae_3(\cos\th)\right)\Jk +\Ga_b\c\Ga_b.
\eeaa
This concludes the proof of Proposition \ref{Prop:linearizedPTstructure1}.

%%%%%%%%%%%%%%%%%%%%%%%%%%%%%%%%%%

\section{Proof of Lemma \ref{lemma:linkPGandTframeinMext}}
\lab{sec:proofoflemma:linkPGandTframeinMext}

%%%%%%%%%%%%%%%%%%%%%%%%%%%%%%%%%%

In view of the initialization of the PG and PT structures of $\Mext$ on $\Si_*$, we have 
\beaa
\la=1, \qquad f=0, \qquad u'=u, \qquad r'=r, \qquad \th'=\th, \qquad \Jk'=\Jk\quad\textrm{on}\quad \Si_*,
\eeaa 
see in particular Remark \ref{rmk:comparisionchangeofframePGandPTonSigmastar}.

Next, let $F=f+i\dual f$. Since $\xi=\xi'=0$ and $\om=\om'=0$, we have in view of Corollary \ref{cor:transportequationine4forchangeofframecoeffinformFFbandlamba:PTcase} 
\beaa
\nab_{\la^{-1}e_4'}F+\frac{1}{2}\ov{\tr X} F+2\om F &=& -2\Xi -\chih\c F+E_1(f, \Ga),\\
\la^{-1}\nab_4'(\log\la) &=& 2\om+f\c(\ze-\etab)+E_2(f, \Ga),
\eeaa
where $E_1(f, \Ga)=O(f^2\Ga)$ and $E_2(f, \Ga)=O(f^2\Ga)$. Since $\la=1$ and $f=0$ on $\Si_*$, this propagates immediately to $\Mext$, and hence
\beaa
\la=1, \qquad f=0, \qquad e_4'=e_4\quad\textrm{on}\quad \Mext.
\eeaa 

Next, since $e_4'=e_4$ on $\Mext$, and since $e_4(r)=e_4'(r')=1$, $e_4(u)=e_4'(u')=0$, and $e_4(\th)=e_4'(\th')=0$, we infer $e_4(r'-r)=0$, $e_4(u'-u)=0$ and $e_4(\th'-\th)=0$. Hence, since $u'=u$, $r'=r$ and $\th'=\th$ on $\Si_*$, this propagates immediately to $\Mext$
\beaa
u'=u,\qquad r'=r, \qquad \th'=\th, \qquad q'=q\quad\textrm{on}\quad \Mext.
\eeaa 

Also, we compute $\nab_4(q\Jk')$. Since $e_4'=e_4$, $q'=q$, and $\nab_4'(q'\Jk')=0$, we have
\beaa
0 &=& \nab_4'(q'\Jk')_a = e_4'(q'\Jk'_a)-\g(\D_{e_4'}e_a', e_b')\Jk_b' = e_4(q\Jk_a')-\g(\D_{e_4'}e_a', e_b')\Jk_b'\\
&=& \nab_4(q\Jk')_a -\Big(\g(\D_{e_4'}e_a', e_b')-\g(\D_{e_4}e_a, e_b)\Big)\Jk_b'.
\eeaa
Since, using $\la=1$ and $f=0$, we have in view of the frame transformation 
\beaa
\g(\D_{e_4'}e_a', e_b') &=& \g\left(\D_{e_4}\left(e_a +\frac 1 2  \fb_a  e_4\right), e_b +\frac 1 2  \fb_b  e_4\right)\\
&=& \g(\D_{e_4}e_a, e_b) -\fb_b\xi_a+\fb_a\xi_b=\g(\D_{e_4}e_a, e_b)
\eeaa
where we used the fact that $\xi=0$, we infer $\nab_4(q\Jk')=0$ on $\Mext$. Hence, since $\nab_4(q\Jk)=0$, we deduce on $\Mext$
\beaa
\nab_4\big( q(\Jk-\Jk')\big)=0.
\eeaa
As $\Jk'=\Jk$ on $\Si_*$, this propagates immediately to $\Mext$, and we deduce 
\beaa
\Jk'=\Jk\quad\textrm{on}\quad \Mext.
\eeaa 

Next, we consider the transition coefficients $(f', \fb', \la')$  from the PT frame to the PG one, i.e. $(f', \fb', \la')$ corresponds to the inverse transformation of $(f, \fb, \la)$. Since $\la=1$ and $f=0$, we infer immediately from \eqref{relations:laffb-to-primes} that  
\beaa
\la'=1, \quad f'=0, \quad \fb'=-\fb.
\eeaa

Next, we derive the transport equation for $\fb$. We compute\footnote{We use here a more precise transformation formula for $\etab$ than the one derived in Proposition \ref{Proposition:transformationRicci}.} 
\beaa
2\etab_a' &=& \g(\D_{e_4'}e_3', e_a') = \g\left(\D_{e_4}\left(e_3 + \fb^b e_b  + \frac 1 4 |\fb|^2 e_4\right), e_a +\frac 1 2  \fb_a  e_4\right)\\
&=& \g\left(\D_{e_4}e_3, e_a +\frac 1 2  \fb_a  e_4\right)+e_4(\fb_a)+\fb^b\g\left(\D_{e_4}e_b, e_a +\frac 1 2  \fb_a  e_4\right)\\
&&+\frac 1 4 |\fb|^2\g\left(\D_{e_4}e_4, e_a +\frac 1 2  \fb_a  e_4\right)\\
&=& 2\etab_a -2\om\fb_a+\nab_4\fb_a -\fb_b\fb_a\xi_b+\frac{1}{2}|\fb|^2\xi_a
\eeaa
and since $\xi=0$ and $\om=0$, we infer
\beaa
\nab_4\fb &=& 2(\etab'-\etab).
\eeaa
Since  $(e_1, e_2, e_3, e_4)$ is a PG frame and $(e_1', e_2', e_3', e_4')$ is a PT frame, we have $\etab=-\ze$ and $\Hb'=-\frac{a\ov{q'}}{|q'|^2}\Jk'$. Since $q'=q$ and $\Jk'=\Jk$, we infer in view of the definition of $\Zc$ in Definition \ref{def:renormalizationofallnonsmallquantitiesinPGstructurebyKerrvalue},
\beaa
\Hb'-\Hb &=& -\frac{a\ov{q'}}{|q'|^2}\Jk'+Z=Z-\frac{a\ov{q}}{|q|^2}\Jk=\Zc
\eeaa
so that $\etab' -\etab=\widecheck{\ze}$
and hence 
\beaa
\nab_4\fb &=& 2\widecheck{\ze}.
\eeaa

Finally, we derive the identity for $\nab_4'\Fb'$ where $\Fb'=\fb'+i\dual\fb'$. To this end, we compute\footnote{We use here a more precise transformation formula for $\ze$ than the one derived in Proposition \ref{Proposition:transformationRicci}.} 
\beaa
2\ze_a &=& \g(\D_{e_a}e_4, e_3) = \g\left(\D_{e_a' +\frac 1 2  {\fb'}_a  e_4'}e_4', e_3' + {\fb'}^b e_b'  + \frac 1 4 |\fb'|^2 e_4'\right)\\
&=&  \g\left(\D_{e_a'}e_4', e_3' + {\fb'}^b e_b'\right)+\frac 1 2  \fb_a'\g\left(\D_{e_4'}e_4', e_3' + {\fb'}^b e_b'\right)\\
&=& 2\ze_a' +{\fb'}^b\chi_{ab}'+2\om'\fb_a' +(\fb'\c\xi')\fb_a'
\eeaa
and since $\xi'=0$ and $\om'=0$, we infer
\beaa
\ze &=& \ze' +\frac{1}{4}\trch'\fb'+\frac{1}{4}\atrch'\dual\fb'+\frac{1}{2}\fb'\c\chih'.
\eeaa
Since $q'=q$ and $\Jk'=\Jk$, and in view of the linearization of $\ze$ in Definition \ref{def:renormalizationofallnonsmallquantitiesinPGstructurebyKerrvalue} and the one for $\ze'$ in Definition \ref{def:renormalizationofallnonsmallquantitiesinPTstructurebyKerrvalue:chap9}, we have
\beaa
\ze'-\ze &=& \Re\left(\frac{a\ov{q'}}{|q'|^2}\Jk'-\frac{a\ov{q}}{|q|^2}\Jk\right)+\zec'-\zec=\zec'-\zec
\eeaa
and hence
\beaa
\zec &=& \zec' +\frac{1}{4}\trch'\fb'+\frac{1}{4}\atrch'\dual\fb'+\frac{1}{2}\fb'\c\chih'.
\eeaa
Together with the above equation for $\nab_4\fb$, this yields 
\beaa
\nab_4\fb &=& 2\widecheck{\ze}=2\zec' +\frac{1}{2}\trch'\fb'+\frac{1}{2}\atrch'\dual\fb'+\fb'\c\chih'.
\eeaa
Since we have obtained above
\beaa
e_4'=e_4, \qquad \fb'=-\fb, \qquad \g(\D_{e_4'}e_a', e_b') =\g(\D_{e_4}e_a, e_b),
\eeaa
we infer
\beaa
\nab_4'\fb_a' &=& e_4'(\fb_a')-\g(\D_{e_4'}e_a', e_b')\fb_b'\\
&=& -e_4(\fb_a) + \g(\D_{e_4}e_a, e_b)\fb_b =-\nab_4\fb_a
\eeaa
and hence
\beaa
\nab_4'\fb'+\frac{1}{2}\trch'\fb'+\frac{1}{2}\atrch'\dual\fb' &=& -2\zec' -\fb'\c\chih'.
\eeaa
With the notation $\Fb'=\fb'+i\dual\fb'$, we infer
\beaa
\nab_4'\Fb'+\frac{1}{2}\tr X'\Fb' &=& -2\Zc' -\Fb'\c\chih'
\eeaa
as desired. This concludes the proof of Lemma \ref{lemma:linkPGandTframeinMext}.

%%%%%%%%%%%%%%%%%%%%%%%%%%%%%%%%%%%%%%%%%%

\section{Proof of Proposition \ref{prop:propertiesoftauusefulfortheoremM8:chap9}}
\lab{sec:proofofprop:propertiesoftauusefulfortheoremM8:chap9}

%%%%%%%%%%%%%%%%%%%%%%%%%%%%%%%%%%%%%%%%%%

In this section, we prove Proposition \ref{prop:propertiesoftauusefulfortheoremM8:chap9} in Kerr, see Proposition \ref{prop:propertiesoftauusefulfortheoremM8:appendix} below. The general case 
follows by obvious modifications and is left to the reader.

%%%%%%%%%%%%%%%%%%%%%%%%%%%%%%%%%%%%%%%%%%%%%%%%%%

\subsection{Setting in Kerr}

%%%%%%%%%%%%%%%%%%%%%%%%%%%%%%%%%%%%%%%%%%%%%%%%%%

Let $r_0\gg m$ a fixed large constant. In the Boyer-Lindquist coordinates, we introduce 
\bea
u:=t-\int_{r_0}^r\frac{{r'}^2+a^2}{\Delta(r')}dr', \qquad \ub:=t+\int_{r_0}^r\frac{{r'}^2+a^2}{\Delta(r')}dr',
\eea
where the above normalization is such that we have
\bea
u=\ub\quad\textrm{on}\quad r=r_0.
\eea

The following lemma will be useful.
\begin{lemma}\lab{lemma:usefulexpansioninrofubminusu}
We have for $r\geq r_0$
\bea
\frac{1}{2}(\ub-u) &=& r-r_0+2m\log\left(\frac{r}{r_0}\right)-\frac{4m^2}{r}+\frac{4m^2}{r_0}+\int_{r_0}^rO\left(\frac{m^3}{{r'}^3}\right)dr'.
\eea
\end{lemma}

\begin{proof}
We have
\beaa
\frac{r^2+a^2}{\De} &=& \frac{r^2+a^2}{r^2-2mr+a^2}=\frac{1+\frac{a^2}{r^2}}{1-\frac{2m}{r}+\frac{a^2}{r^2}}\\
&=& \left(1+\frac{a^2}{r^2}\right)\left(1+\frac{2m}{r}-\frac{a^2}{r^2}+\frac{4m^2}{r^2}+O\left(\frac{m^3}{r^3}\right)\right)\\
&=& 1+\frac{2m}{r}+\frac{4m^2}{r^2}+O\left(\frac{m^3}{r^3}\right)
\eeaa
and hence 
\beaa
\frac{1}{2}(\ub-u) &=& \int_{r_0}^r\frac{{r'}^2+a^2}{\Delta(r')}dr' = \int_{r_0}^r\left(1+\frac{2m}{r'}+\frac{4m^2}{{r'}^2}+O\left(\frac{m^3}{{r'}^3}\right)\right)dr'\\
&=&  r-r_0+2m\log\left(\frac{r}{r_0}\right)-\frac{4m^2}{r}+\frac{4m^2}{r_0}+\int_{r_0}^rO\left(\frac{m^3}{{r'}^3}\right)dr'
\eeaa
as desired. 
\end{proof}

Also, let $\de_\HH>0$  a small enough constant, and $r_+$ is given by
\bea
r_+ &:=& m+\sqrt{m^2-a^2}. 
\eea
Moreover, let $r_*\gg r_0$. We define the spacelike hypersurfaces
\bea
\AA:=\Big\{r=r_+-\de_\HH\Big\}, \qquad \Si_*:=\Big\{u+r=u_*+r_*, \,\, 1\leq u\leq u_*\Big\}. 
\eea
Also, note that $\Mint'=\Mint$ in Kerr with 
\bea
\Mint:=\Big\{r_+-\de_\HH\leq r\leq r_0, \quad 1\leq\ub\leq u_*\Big\},
\eea
and that $\Mext$ is given by 
\bea
\Mext:=\Big\{r\geq r_0, \quad 1\leq u\leq u_*, \quad u+r\leq u_*+r_*\Big\}.
\eea

Next, we look for  the scalar function $\tau$ under the form 
\bea\lab{def:choiceoftau}
\tau &:=& \ub+f(r),
\eea
where $f$ will be carefully chosen such that:
\begin{enumerate}
\item The level sets of $\tau$ are spacelike, i.e. 
\bea\lab{eq:condf1}
\g(\D\tau, \D\tau)<0\quad \textrm{ on }r\geq r_+-\de_\HH.
\eea

\item $\tau=\ub$ at $(\ub=u_*, r=r_+-\de_\HH)$, i.e. 
\bea\lab{eq:condf2}
f(r_+-\de_\HH)=0.
\eea

\item $\tau=u$ at $(u=u_*, r=r_*)$, i.e. 
\bea\lab{eq:condf3}
f(r_*)=-2\int_{r_0}^{r_*}\frac{r^2+a^2}{\De(r)}dr.
\eea
\end{enumerate}

The choice of $\tau$ allows us to define ${}^{(top)}\Si$. 
\begin{definition}\lab{def:definitionofSigmatop}
Let $f(r)$ a function satisfying \eqref{eq:condf1} \eqref{eq:condf2} \eqref{eq:condf3},  and let $\tau$ given by \eqref{def:choiceoftau}. Then, we define the hypersurface ${}^{(top)}\Si$  by 
\bea
{}^{(top)}\Si &:=& \Big\{\tau=u_*, \,\, r_+-\de_\HH\leq r\leq r_*\Big\}.
\eea
\end{definition}

\begin{remark}
In view of the definition of $\tau$ and the fact that $f$ satisfies properties \eqref{eq:condf1} \eqref{eq:condf2} \eqref{eq:condf3}, the hypersurface ${}^{(top)}\Si$ given by Definition \ref{def:definitionofSigmatop} satisfies:
\begin{enumerate}
\item ${}^{(top)}\Si$ is spacelike,

\item ${}^{(top)}\Si\cap\AA=\{\ub=u_*\}\cap\AA$,

\item ${}^{(top)}\Si\cap\Si_*=\{u=u_*\}\cap\Si_*$.
\end{enumerate}
 \end{remark}

We then define $\Mtop$ as follows.
\begin{definition}\lab{def:definitionofMtop}
Let $f(r)$ a function satisfying \eqref{eq:condf1} \eqref{eq:condf2} \eqref{eq:condf3},  and let $\tau$ given by \eqref{def:choiceoftau}. Then, we define the spacetime region $\Mtop$  by 
\bea
\Mtop &:=& \{\tau\leq u_*\}\cap\{u\geq u_*\}\cap\{\ub\geq u_*\}.
\eea
\end{definition}

Finally, $\MM$ is given by
\bea
\MM &:=& \Mext\cup\Mint\cup\Mtop,
\eea
and we also introduce the region $\Mtop'$ given by
\bea
\Mtop' &:=& \{\tau\leq u_*\}\cap\{u\geq u_*'\}\cap\{\ub\geq u_*'\},
\eea
where $u_*-2\leq u_*'\leq u_*-1$.

We are now ready to state the analog of Proposition \ref{prop:propertiesoftauusefulfortheoremM8:chap9} in Kerr.
\begin{proposition}\lab{prop:propertiesoftauusefulfortheoremM8:appendix}
Let $\tau$ a scalar function given by \eqref{def:choiceoftau}. There exists a particular choice of $f(r)$ satisfying \eqref{eq:condf1} \eqref{eq:condf2} \eqref{eq:condf3}, and such that:
\begin{enumerate}
\item We have on  $\MM$
\bea
\g(\D\tau, \D\tau) &\leq& -\frac{m^2}{4|q|^2}<0,
\eea
so that the level sets of $\tau$ are spacelike and asymptotically null.

\item Denoting, on each level set of $\ub$ in $\,{}^{(top)}\MM'(r\geq r_0)$, by $r_+(\ub)$ the maximal value of $r$  and by $r_-(\ub)$ the minimal value of $r$, we have\footnote{Note that \eqref{eq:upperboundrpubminusrmubonMtop:App} depends on the choice of ${}^{(top)}\Si$ and hence on the choice of $\tau$.}
\bea\lab{eq:upperboundrpubminusrmubonMtop:App}
0\leq r_+(\ub)-r_-(\ub)\leq 2m+1.
\eea

\item In $\Mtop'(r\leq r_0)$, $\tau$ satisfies   
\bea
u_*-(m+2)\leq \tau \leq u_*.
\eea

\item In $\MM(r\leq r_0)$, $\tau$ satisfies 
\bea
e_4(\tau) = \frac{2(r^2+a^2) -\frac{m^2}{r^2}\De}{|q|^2}, \quad  e_3(\tau) = \frac{m^2}{r^2},\quad \nab(\tau) = a\Re(\Jk).
\eea
\end{enumerate}
\end{proposition}

In the next section, we construct the suitable function $f(r)$. This will allow us to prove Proposition \ref{prop:propertiesoftauusefulfortheoremM8:appendix} in section \ref{sec:proofofprop:propertiesoftauusefulfortheoremM8:appendix}.

%%%%%%%%%%%%%%%%%%%%%%%%%%%%%%%%%%%%%%%%%%%%%%%%%%

\subsection{Construction of a suitable function $f(r)$}

%%%%%%%%%%%%%%%%%%%%%%%%%%%%%%%%%%%%%%%%%%%%%%%%%%

In this section, we exhibit a function $f(r)$ satisfying in particular the properties \eqref{eq:condf1} \eqref{eq:condf2} \eqref{eq:condf3}.

We start with the following lemma.
\begin{lemma}
We have 
\bea
\g(\D\tau, \D\tau) &=& \frac{1}{|q|^2}\Big(\Delta(f'(r))^2+2(r^2+a^2)f'(r)+a^2(\sin\th)^2\Big).
\eea
\end{lemma}

\begin{proof}
We have, using  the fact that $\tau=\ub+f(r)$, $e_3(\ub)=0$, $e_3(r)=-1$, and $\nab(r)=0$, 
\beaa
\g(\D\tau, \D\tau) &=& -e_4(\tau)e_3(\tau)+|\nab(\tau)|^2 = \Big(e_4(\ub)+e_4(r)f'(r)\Big)f'(r)+|\nab(\ub)|^2.
\eeaa 
Since
\beaa
e_4(\ub)=\frac{2(r^2+a^2)}{|q|^2}, \qquad e_4(r)=\frac{\Delta}{|q|^2}, \qquad |\nab(\ub)|^2=\frac{a^2(\sin\th)^2}{|q|^2},
\eeaa
we infer
\beaa
\g(\D \tau, \D \tau) &=& \left(\frac{2(r^2+a^2)}{|q|^2}+\frac{\Delta}{|q|^2}f'(r)\right)f'(r)+\frac{a^2(\sin\th)^2}{|q|^2}
\eeaa 
as desired.
\end{proof}

Motivated by the above lemma, we consider, for $r>r_+$, the following second order polynomial 
\bea\lab{eq:secondorderpolynomialforfuturespacelikeboundarycondition}
P(X) &:=& \Delta X^2+2(r^2+a^2)X+a^2.
\eea

\begin{lemma}
For $r>r_+$, we have $P(X)<0$ if and only if
\beaa
X_-(r) <X<  X_+(r).
\eeaa
where
\beaa
X_-(r) &:=& -\frac{r^2+a^2}{\Delta}-\frac{\sqrt{r^4+a^2r^2+2a^2mr}}{\Delta},\\ 
X_+(r) &:=& -\frac{a^2}{r^2+a^2+\sqrt{r^4+a^2r^2+2a^2mr}}.
\eeaa
\end{lemma}

\begin{proof}
We compute 
\beaa
D_P(r) &=& 4(r^2+a^2)^2-4a^2\Delta\\
&=& 4\Big(r^4+2a^2r^2+a^4 -a^2(r^2-2mr+a^2)\Big)\\
&=& 4\Big(r^4+a^2r^2+2a^2mr\Big)
\eeaa
so that $D_P(r)>0$ for any $r>0$. In particular, for $r>r_+$, since $\Delta>0$ in this case, $P$ is a second order polynomial in $X$, and hence, $P$ has two distinct roots $X_{\pm}(r)$ given by 
\beaa
X_{\pm}(r) &=& -\frac{r^2+a^2}{\Delta}\pm\frac{\sqrt{D_P(r)}}{2\Delta}, \qquad X_-(r)<X_+(r)<0.
\eeaa
Note that we may rewrite $X_-(r)$ and $X_+(r)$ as
\beaa
X_-(r) = -\frac{r^2+a^2}{\Delta}-\frac{\sqrt{r^4+a^2r^2+2a^2mr}}{\Delta},\qquad X_+(r) = -\frac{a^2}{r^2+a^2+\sqrt{r^4+a^2r^2+2a^2mr}}.
\eeaa
Also, for $r>r_+$,  since $\Delta>0$, we have $P(X)<0$ if and only if $X_-(r)<X<X_+(r)$ as desired.
\end{proof}

\begin{lemma}
For $r$ large, we have the following expansions
\beaa
X_-(r) &=&  -2-\frac{4m}{r}-\frac{8m^2-\frac{1}{2}a^2}{r^2}+O\left(\frac{m^3}{r^3}\right),\\
X_+(r) &=& -\frac{a^2}{2r^2}+O\left(\frac{m^3}{r^3}\right).
\eeaa
\end{lemma}

\begin{proof}
For $r$ large, we have the expansions
\beaa
\sqrt{r^4+a^2r^2+2a^2mr} &=& r^2\sqrt{1+\frac{a^2}{r^2}+O\left(\frac{m^3}{r^3}\right)}\\
&=& r^2\left(1+\frac{a^2}{2r^2}+O\left(\frac{m^3}{r^3}\right)\right),
\eeaa
\beaa
&& -\frac{r^2+a^2}{\Delta}-\frac{\sqrt{r^4+a^2r^2+2a^2mr}}{\Delta}\\
 &=& -\frac{2+\frac{3a^2}{2r^2}+O\left(\frac{m^3}{r^3}\right)}{1-\frac{2m}{r}+\frac{a^2}{r^2}}\\
&=& -\left(2+\frac{3a^2}{2r^2}+O\left(\frac{m^3}{r^3}\right)\right)\left(1 +\frac{2m}{r}+\frac{4m^2-a^2}{r^2}+O\left(\frac{m^3}{r^3}\right)\right)\\
&=& -2-\frac{4m}{r}-\frac{8m^2-\frac{1}{2}a^2}{r^2}+O\left(\frac{m^3}{r^3}\right),
\eeaa
and
\beaa
-\frac{a^2}{r^2+a^2+\sqrt{r^4+a^2r^2+2a^2mr}} &=& -\frac{a^2}{2r^2}+O\left(\frac{m^3}{r^3}\right),
\eeaa
which concludes the proof of the lemma.
\end{proof}

\begin{lemma}\lab{lemma:defintionoff1andf2andtheirbasicpropertiesusefullater}
Let, for $r\geq r_+-\de_\HH$, 
\bea
\bsplit
f_1(r) &:= \frac{m^2}{r}-\frac{m^2}{r_+-\de_\HH},\\
f_2(r) &:=  -2(r-r_0)-4m\log\left(\frac{r}{r_0}\right)-\frac{m^2}{r}+c_{0,*},
\end{split}
\eea
with the constant $c_{0,*}$ chosen such that $f_2$ satisfies \eqref{eq:condf3}. Then: 
\begin{enumerate}
\item We have, for all $r\geq r_+-\de_\HH$,
\beaa
f_1'(r)>f_2'(r).
\eeaa

\item there exists a unique solution $r_1\geq r_+-\de_\HH$ of 
\beaa
f_1(r_1)=f_2(r_1).
\eeaa

\item We have $r_1\in(r_0, r_0+m)$.

\item For $r_+-\de_\HH<r<r_1$, we have $f_1(r)<f_2(r)$, and for $r>r_1$, we have $f_1(r)>f_2(r)$. 

\item The following holds:
\begin{itemize}
\item For all $r\geq r_+-\de_\HH$, we have
\beaa
P(f_1'(r))\leq -\frac{m^2}{4}.
\eeaa

\item For all $r\geq r_0$, we have
\beaa
P(f_2'(r)) &\leq& -8m^2.
\eeaa

\item For all $r\geq r_0$, we have
\beaa
P\big(\sigma f_1'(r)+(1-\sigma)f_2'(r)\big)\leq -\frac{m^2}{4}, \qquad 0\leq \sigma\leq 1.
\eeaa
\end{itemize}
\end{enumerate}
\end{lemma}

\begin{proof}
In view of the definition of $f_1(r)$ and $f_2(r)$, we have
\beaa
f_1'(r) &=&  -\frac{m^2}{r^2},\\
f_2'(r) &=& -2-\frac{4m}{r}+\frac{m^2}{r^2}.
\eeaa
In particular, note that, on $r>r_+-\de_\HH$, 
\beaa
f_1'(r)-f_2'(r) &=& -\frac{2m^2}{r^2}+\frac{4m}{r} +2 =\frac{2r^2+4mr-2m^2}{r^2}\\
&=& \frac{2}{r^2}(r+m-\sqrt{2}m)(r+m+\sqrt{2}m)\\
&\geq& \frac{2}{r^2}\left(r_+ -(\sqrt{2}-1)m\right)\left(r_++(1+\sqrt{2})m\right)\\
&\geq& \frac{2}{r^2}(2-\sqrt{2})(2+\sqrt{2})m^2>0.
\eeaa

Also, since $f_2$ satisfies \eqref{eq:condf3}, we have
\beaa
f_2(r_*)=-2\int_{r_0}^{r_*}\frac{r^2+a^2}{\De(r)}dr,
\eeaa
and hence, using in particular  Lemma \ref{lemma:usefulexpansioninrofubminusu}, 
\beaa
c_{0,*} &=&  \frac{9m^2}{r_*}-\frac{8m^2}{r_0}+\int_{r_0}^{r_*}O\left(\frac{m^3}{{r'}^3}\right)dr'\\
&=& -\frac{8m^2}{r_0}\left(1+O\left(\frac{r_0}{r_*}\right)+O\left(\frac{m}{r_0}\right)\right).
\eeaa
We infer
\beaa
 f_1(r_0) &=& -\frac{m^2}{r_+-\de_\HH}+\frac{m^2}{r_0}, \\ 
f_2(r_0) &=& O\left(\frac{m^2}{r_0}\right),\\
f_1(r_0+m) &=& -\frac{m^2}{r_+-\de_\HH}+\frac{m^2}{r_0+m}, \\ 
f_2(r_0+m) &=& -2m+O\left(\frac{m^2}{r_0}\right),
\eeaa
and hence
\beaa
f_1(r_0)<f_2(r_0), \qquad f_1(r_0+m)>f_2(r_0+m)
\eeaa
so that, there exists, by the mean value theorem, $r_1\in(r_0, r_0+m)$ such that 
\beaa
f_1(r_1)=f_2(r_1).
\eeaa
Note that, since $f_1'(r)>f_2'(r)$ for all $r>r_+-\de_\HH$, $r_1$ is the unique 0 of $f_1-f_2$ on $r>r_+-\de_\HH$, and we have $f_1(r)<f_2(r)$ for $r<r_1$, and $f_1(r)>f_2(r)$ for $r>r_1$. 

Next, we compare $f_1'(r)$ to the roots $X_\pm(r)$ of $P$. We have
\beaa
X_+(r) - f_1'(r) &=&  -\frac{a^2}{r^2+a^2+\sqrt{r^4+a^2r^2+2a^2mr}}+\frac{m^2}{r^2}\\
&=& \frac{m^2(r^2+a^2+\sqrt{r^4+a^2r^2+2a^2mr})-a^2r^2}{r^2(r^2+a^2+\sqrt{r^4+a^2r^2+2a^2mr})}\\
&=& \frac{(m^2-a^2)r^2+m^2a^2+m^2\sqrt{r^4+a^2r^2+2a^2mr}}{r^2(r^2+a^2+\sqrt{r^4+a^2r^2+2a^2mr})}\\
&\geq& \frac{(2m^2-a^2)r^2+m^2a^2}{r^2(r^2+a^2+\sqrt{r^4+a^2r^2+2a^2mr})}\\
&\geq& \frac{m^2r^2}{r^2(r^2+m^2+\sqrt{(r^2+m^2)^2})}\\
&\geq& \frac{m^2}{4r^2}
\eeaa
and
\beaa
f_1'(r) - X_-(r) &=& -\frac{m^2}{r^2}+\frac{r^2+a^2}{\Delta}+\frac{\sqrt{r^4+a^2r^2+2a^2mr}}{\Delta}\\
&=& \frac{r^2(r^2+a^2+\sqrt{r^4+a^2r^2+2a^2mr})-m^2\Delta}{r^2\Delta}\\
&=& \frac{r^4 -(m^2-a^2)r^2 +2m^3r-a^2m^2+r^2\sqrt{r^4+a^2r^2+2a^2mr}}{r^2\Delta}.
\eeaa
Note that we have, for $r>r_+$,
\beaa
(r^4 -(m^2-a^2)r^2 +2m^3r-a^2m^2)' &=& 4r^3 -2(m^2-a^2)r+2m^3\\
&=& 2r(2r^2-2m^2+2a^2)+2m^3\\
&\geq& 4ma^2+2m^3>0
\eeaa
which implies,  for $r>r_+$,
\beaa
f_1'(r) - X_-(r) &\geq& \frac{m^4 -(m^2-a^2)m^2 +2m^4-a^2m^2+r^2\sqrt{r^4+a^2r^2+2a^2mr}}{r^2\Delta}\\
&=& \frac{2m^4+r^2\sqrt{r^4+a^2r^2+2a^2mr}}{r^2\Delta}\\
&\geq& \frac{r^2}{\Delta}.
\eeaa
For $r>r_+$, we infer
\beaa
P(f_1'(r)) &=& \Delta(f_1'(r)-X_-(r))(f_1'(r)-X_+(r))\\
 &=& -\Delta(f_1'(r)-X_-(r))(X_+(r)-f_1'(r))\\
 &\leq& -\frac{m^2}{4}.
\eeaa
Also, we have for $r_+-\de_\HH\leq r\leq r_+$
\beaa
P(f_1'(r)) &=& -|\Delta|(f_1'(r))^2+2(r^2+a^2)f_1'(r)+a^2\\
&\leq& 2(r^2+a^2)\left( -\frac{m^2}{r^2}\right)+a^2\\
&\leq& -2m^2+a^2\\
&\leq& -m^2.
\eeaa
Thus, we have for $r\geq r_+-\de_\HH$
\beaa
P(f_1'(r))  &\leq& -\frac{m^2}{4}.
\eeaa

Next, we  derive an upper bound for $P(f_2'(r))$ when $r\geq r_0$. Recall that we have for $r$ large the following expansions
\beaa
X_-(r) &=&  -2-\frac{4m}{r}-\frac{8m^2-\frac{1}{2}a^2}{r^2}+O\left(\frac{m^3}{r^3}\right),\\
X_+(r) &=& -\frac{a^2}{2r^2}+O\left(\frac{m^3}{r^3}\right).
\eeaa
We infer
\beaa
X_+(r) - f_2'(r) &=& -\frac{a^2}{2r^2}+O\left(\frac{m^3}{r^3}\right) +2+\frac{4m}{r}-\frac{m^2}{r^2}\\
&=& 2+O\left(\frac{m}{r}\right)
\eeaa
and 
\beaa
f_2'(r) - X_-(r)   &=& -2-\frac{4m}{r}+\frac{m^2}{r^2} +2+\frac{4m}{r}+\frac{8m^2-\frac{1}{2}a^2}{r^2}+O\left(\frac{m^3}{r^3}\right)\\
&=& \frac{9m^2-\frac{1}{2}a^2}{r^2}+O\left(\frac{m^3}{r^3}\right).
\eeaa
Since
\beaa
\Delta &=& r^2\left(1+O\left(\frac{m}{r}\right)\right),
\eeaa
we infer, for $r\geq r_0$, 
\beaa
P(f_2'(r)) &=& \Delta(f_2'(r)-X_-(r))(f_2'(r)-X_+(r))\\
&=& -\left(1+O\left(\frac{m}{r}\right)\right)\left(2+O\left(\frac{m}{r}\right)\right)\left(9m^2-\frac{1}{2}a^2+O\left(\frac{m^3}{r}\right)\right)\\
& \leq & -\frac{17}{2}m^2\left(1+O\left(\frac{m}{r_0}\right)\right)\\
&\leq& -8m^2
\eeaa
for $r_0$ large enough compared to $m$. 

Finally, for $r>r_+$, we have $P''(X)=2\Delta>0$ and hence $P$ is convex in $X$ which implies, for any $0\leq \sigma\leq 1$, for $r\geq r_0$,
\beaa
P\big(\sigma f_1'(r)+(1-\sigma)f_2'(r)\big) &\leq& \sigma P(f_1'(r))+(1-\sigma)P(f_2'(r))\\
&\leq& -\frac{\sigma m^2}{4} -8m^2(1-\sigma)\\
&\leq&-\frac{m^2}{4}
\eeaa
as desired. This concludes the proof of Lemma \ref{lemma:defintionoff1andf2andtheirbasicpropertiesusefullater}.
\end{proof}

\begin{corollary}\lab{cor:choiceoffforthefurtureboundaryofthespacetime}
The exists a smooth function $f$ such that 
\begin{enumerate}
\item $f(r)=f_1(r)$ for $r\leq r_0$.

\item $f(r)=f_2(r)$ for $r\geq r_0+m$.

\item For all $r\geq r_+-\de_\HH$, $\tau:=\ub+f(r)$ verifies 
\beaa
\g(\D\tau, \D\tau) &\leq& -\frac{m^2+4a^2(\cos\th)^2}{4|q|^2}<0
\eeaa
so that the level sets of $\tau$ are spacelike in $r\geq r_+-\de_\HH$.
\end{enumerate}
\end{corollary}

\begin{proof}
Recall from Lemma \ref{lemma:defintionoff1andf2andtheirbasicpropertiesusefullater} that 
\begin{enumerate}
\item We have, for all $r\geq r_+-\de_\HH$,
\beaa
f_1'(r)>f_2'(r).
\eeaa

\item there exists a unique solution $r_1\geq r_+-\de_\HH$ of 
\beaa
f_1(r_1)=f_2(r_1).
\eeaa

\item We have $r_1\in(r_0, r_0+m)$.

\item For $r_+-\de_\HH<r<r_1$, we have $f_1(r)<f_2(r)$, and for $r>r_1$, we have $f_1(r)>f_2(r)$. 
\end{enumerate}
We deduce the existence of a smooth function $f$ on $r\geq r_+-\de_\HH$ such that 
\begin{enumerate}
\item $f(r)=f_1(r)$ for $r_+-\de_\HH\leq r\leq r_0$.

\item $f(r)=f_2(r)$ for $r\geq r_0+m$.

\item For all $r_0\leq r\leq r_0+m$, $f'(r)$ verifies
\beaa
f_2'(r)\leq f'(r)\leq f_1'(r).
\eeaa
\end{enumerate}

It remains to check the property of the level sets of $\tau=\ub+f(r)$. 
From the above properties of $f$, we have, for all $r\geq r_+-\de_\HH$
\beaa
P(f'(r)) &\leq& -\frac{m^2}{4}
\eeaa
and hence
\beaa
\Delta (f'(r))^2+2(r^2+a^2)f'(r)+a^2 &\leq&  -\frac{m^2}{4}.
\eeaa
We deduce, for $r\geq r_+-\de_\HH$, 
\beaa
\g(\D\tau, \D\tau) &=& \frac{1}{|q|^2}\Big(\Delta(f'(r))^2+2(r^2+a^2)f'(r)+a^2(\sin\th)^2\Big)\\
&=& \frac{1}{|q|^2}\Big(\Delta(f'(r))^2+2(r^2+a^2)f'(r)+a^2-a^2(\cos\th)^2\Big)\\
&\leq& -\frac{m^2+4a^2(\cos\th)^2}{4|q|^2}
\eeaa
as desired.
\end{proof}

\begin{corollary}\lab{cor:fstasfiesallthecorrectconditions}
The function $f$ in Corollary \ref{cor:choiceoffforthefurtureboundaryofthespacetime} satisfies \eqref{eq:condf1} \eqref{eq:condf2} \eqref{eq:condf3}. 
\end{corollary}

\begin{proof}
Since $\tau=\ub+f(r)$ verifies $\g(\D\tau, \D\tau)<0$ in view of Corollary \ref{cor:choiceoffforthefurtureboundaryofthespacetime}, $f$ satisfies \eqref{eq:condf1}. Also, since $f(r)=f_1(r)$ for $r_+-\de_\HH\leq r\leq r_0$, and since $f_1(r_+-\de_\HH)=0$, $f$ satisfies \eqref{eq:condf2}. Finally, since $f(r)=f_2(r)$ for $r\geq r_0+m$, and since 
\beaa
f_2(r_*)=-2\int_{r_0}^{r_*}\frac{r^2+a^2}{\De(r)}dr
\eeaa
by the choice of the constant $c_{0,*}$ appearing in the definition of $f_2$, $f$ satisfies \eqref{eq:condf3}. This concludes the proof of the corollary. 
\end{proof}

%%%%%%%%%%%%%%%%%%%%%%%%%%%%%%%%%%%%%%%%%%%%%%%%%%

\subsection{Proof of Proposition \ref{prop:propertiesoftauusefulfortheoremM8:appendix}}\lab{sec:proofofprop:propertiesoftauusefulfortheoremM8:appendix}

%%%%%%%%%%%%%%%%%%%%%%%%%%%%%%%%%%%%%%%%%%%%%%%%%%

Let $\tau$ be given by \eqref{def:choiceoftau}, i.e. $\tau=\ub+f(r)$. We choose $f(r)$ as in Corollary \ref{cor:choiceoffforthefurtureboundaryofthespacetime}. In particular, $f(r)$ satisfies \eqref{eq:condf1} \eqref{eq:condf2} \eqref{eq:condf3} in view of Corollary \ref{cor:fstasfiesallthecorrectconditions}. Also, we have, in view of Corollary \ref{cor:choiceoffforthefurtureboundaryofthespacetime}, for all $r\geq r_+-\de_\HH$,  
\beaa
\g(\D\tau, \D\tau) &\leq& -\frac{m^2}{4|q|^2}<0
\eeaa
so that the first property of Proposition \ref{prop:propertiesoftauusefulfortheoremM8:appendix} is satisfied. 

Next, we consider the second property of Proposition \ref{prop:propertiesoftauusefulfortheoremM8:appendix}, i.e. the upper bound for $r_+(\ub)-r_-(\ub)$. First, one easily sees that on each level set of $\ub$ in $\,{}^{(top)}\MM'(r\geq r_0)$, the maximal value $r_+(\ub)$ of $r$ corresponds to the value of $r$ on $\{u=u_*'\}$, and that the minimal value  $r_-(\ub)$ of $r$ corresponds to the value of $r$ on ${}^{(top)}\Si\cup\{r=r_0\}$. Since $r_+(\ub)$ is the value on $u=u_*'$, and since we have  for $r\geq r_0$
\beaa
\frac{1}{2}(\ub-u) &=& r-r_0+2m\log\left(\frac{r}{r_0}\right)-\frac{4m^2}{r}+\frac{4m^2}{r_0}+\int_{r_0}^rO\left(\frac{m^3}{{r'}^3}\right)dr',
\eeaa
we infer
\beaa
\frac{1}{2}(\ub-u_*') &=& r_+(\ub) -r_0+2m\log\left(\frac{r_+(\ub)}{r_0}\right)-\frac{4m^2}{r_+(\ub)}+\frac{4m^2}{r_0}+\int_{r_0}^{r_+(\ub)}O\left(\frac{m^3}{{r'}^3}\right)dr'.
\eeaa
We first focus on the case where $r_-(\ub)\geq r_0+m$. In that case, $r_-(\ub)$ is the value on $\ub+f(r)=u_*$ and hence
\beaa
\ub+f(r_-(\ub))=u_*. 
\eeaa
We infer
\beaa
&&\frac{1}{2}\Big(u_*-u_*'-f(r_-(\ub))\Big)\\ 
&=& r_+(\ub) -r_0+2m\log\left(\frac{r_+(\ub)}{r_0}\right)-\frac{4m^2}{r_+(\ub)}+\frac{4m^2}{r_0}+\int_{r_0}^{r_+(\ub)}O\left(\frac{m^3}{{r'}^3}\right)dr'.
\eeaa
For $r\geq r_0+m$, we have $f=f_2$ and hence
\beaa
&& \frac{1}{2}\big(u_*-u_*'\big)+r_-(\ub)-r_0 +2m\log\left(\frac{r_-(\ub)}{r_0}\right)+\frac{m^2}{2r_-(\ub)}-\frac{1}{2}c_{0,*}\\
 &=& r_+(\ub) -r_0+2m\log\left(\frac{r_+(\ub)}{r_0}\right)-\frac{4m^2}{r_+(\ub)}+\frac{4m^2}{r_0}+\int_{r_0}^{r_+(\ub)}O\left(\frac{m^3}{{r'}^3}\right)dr'.
\eeaa
This yields
\beaa
r_+(\ub) - r_-(\ub) +2m\log\left(\frac{r_+(\ub)}{r_-(\ub)}\right) &=&   \frac{1}{2}\big(u_*-u_*'\big)+O\left(\frac{m^2}{r_*}\right)+O\left(\frac{m^3}{r_0^2}\right)+  \frac{m^2}{2r_-(\ub)}+\frac{4m^2}{r_+(\ub)}.
\eeaa
and thus, for $r_-(\ub)\geq r_0+m$, and since $u_*'\geq u_*-2$, we infer
\beaa
0\leq r_+(\ub) - r_-(\ub) &\leq& 1+O\left(\frac{m^2}{r_0}\right)\leq 1+m.
\eeaa
In the other case, i.e. $r_0\leq r_-(\ub)\leq r_0+m$, denoting by $\ub_0$ the value of $\ub$ such that $r_-(\ub_0)=r_0+m$, we have $r_+(\ub)\leq r_+(\ub_0)$ and $r_-(\ub)\geq r_0$ and hence
\beaa
0\leq r_+(\ub) - r_-(\ub) &\leq& r_+(\ub_0)-r_0=r_+(\ub_0) -r_-(\ub_0)+m\leq 1+2m
\eeaa
in this case. Thus, in both cases, we have obtained on $\,{}^{(top)}\MM'(r\geq r_0)$
\beaa
0\leq r_+(\ub) - r_-(\ub) \leq 1+2m
\eeaa
so that the second property of Proposition \ref{prop:propertiesoftauusefulfortheoremM8:appendix} is satisfied. 

Next, we consider the third property of Proposition \ref{prop:propertiesoftauusefulfortheoremM8:appendix}, i.e. the lower bound for $\tau$ on $\Mtop'(r\leq r_0)$. Since $\tau=\ub+f(r)$, and since $\ub\geq u_*'\geq u_*-2$ on $\Mtop'$, we infer on $\Mtop'(r\leq r_0)$
\beaa
\tau &=& \ub+f(r)\geq u_*'+f(r)\geq u_*-2+\min_{r_+-\de_\HH\leq r\leq r_0}f(r).
\eeaa
In view of the definition of $f$, $f$ is strictly decreasing so that its minimum is achieved at $r=r_0$. Also, at $r=r_0$, $f=f_1$. We infer on $\Mtop'(r\leq r_0)$
\beaa
\tau &\geq& u_*-2+f_1(r_0)=u_*-2 -\frac{m^2}{r_+-\de_\HH}+\frac{m^2}{r_0}
\eeaa
and hence
\beaa
\tau &\geq& u_*-2 -m.
\eeaa
so that the third property of Proposition \ref{prop:propertiesoftauusefulfortheoremM8:appendix} is satisfied. 

Finally, we consider the last property of Proposition \ref{prop:propertiesoftauusefulfortheoremM8:appendix}. 
Since $\tau=\ub+f((r)$, since $f=f_1$ in for $r\leq r_0$, and since $f_1'(r)=-\frac{m^2}{r^2}$, we have in $\MM(r\leq r_0)$
\beaa
e_4(\tau) &=& e_4(\ub)+e_4(r)f'(r)=\frac{2(r^2+a^2)}{|q|^2}-\frac{m^2}{r^2}\frac{\De}{|q|^2},\\
e_3(\tau) &=& e_3(\ub)+e_3(r)f'(r)=\frac{m^2}{r^2},\\
\nab(\tau) &=& \nab(\ub)+f'(r)\nab(r)=a\Re(\Jk),
\eeaa
as desired. This concludes the proof of of Proposition \ref{prop:propertiesoftauusefulfortheoremM8:appendix}.

%%%%%%%%%%%%%%%%%%%%%%%%%%%%%%%%%%%%%%

%%%%%%%%%%%%%%%%%%%%%%%%%%%%%%%%%%%%%%


\begin{thebibliography}{99}

\bibitem{ABBMa2019} L.  Andersson, Thomas B\"ackdahl, P.  Blue and S.  Ma, \textit{Stability for linearized gravity
on the Kerr spacetime}, arXiv:1903.03859.

\bibitem{A-B} L. Andersson and P.  Blue, \textit{Hidden symmetries and decay for the wave equation on the Kerr spacetime}. Ann. of Math. (2) \textbf{182} (2015), 787-853.
 
\bibitem{AnArGa} Y. Angelopoulos, S. Aretakis and  D. Gajic, \textit{A vector field approach to almost-sharp decay for the wave equation on spherically symmetric, stationary spacetimes}, Ann. PDE \textbf{4} (2018), Art. 15, 120 pp.
 
 
\bibitem{Bar-Press} J. M. Bardeen and W. H. Press, \textit{Radiation fields in the Schwarzschild background}, J. Math. Phys. \textbf{14} (1973), 719.


\bibitem{B-S1}  P. Blue and A. Soffer, \textit{Semilinear wave equations on the Schwarzschild manifold. I. Local decay estimates}, Adv.
Differential Equations \textbf{8} (2003), 595--614.

\bibitem{B-S2}  P. Blue and A. Soffer, \textit{Errata for  ``Global existence and scattering for the nonlinear Schr\"odinger equation on
Schwarzschild manifolds'' , ``Semilinear wave equations on the Schwarzschild manifold I: Local Decay Estimates'',
and `` The wave equation on the Schwarzschild metric II: Local Decay for the spin 2 Regge Wheeler equation'' },
gr-qc/0608073, 6 pages.

           
\bibitem{B-St}  P. Blue and  J. Sterbenz, \textit{Uniform decay of local energy and the semi-linear wave
              equation on {S}chwarzschild space}, Comm. Math. Phys. \textbf{268} (2006), 481--504.    


\bibitem{Ca-Ni}  G. Caciotta   and F. Nicolo     \textit{The non linear perturbation of the Kerr spacetime in an external region}, arXiv:0908.4330.







\bibitem{Carter} B. Carter, \textit{Global structure of the Kerr family of gravitational fields}, Phys. Rev. \textbf{174}  (1968), 1559--1571. 

\bibitem{Chand2} S. Chandrasekhar, \textit{On the equations governing the perturbations of the Schwarzschild black hole}, P. Roy. Soc. Lond. A Mat. \textbf{343} (1975), 289--298.

 \bibitem{Chen}  P.-N. Chen, M.-T.  Wang, S.-T. Yau, \textit{Quasilocal angular momentum and center of mass in general relativity}, Adv. Theor. Math. Phys. 20 (2016), 671-682.  


           
\bibitem{Chand} S. Chandrasekhar, \textit{The mathematical theory of black holes}, 1983, Oxford Classic Texts in the Physical Sciences.


 \bibitem{Chr} D. Christodoulou, \textit{Global solutions of nonlinear hyperbolic equations for small initial data}, Comm. Pure  Appl. Math. \textbf{39} (1986), 267--282.
 
 
 
 \bibitem{Ch-SW} D.Christodoulou, \textit{The formation of shocks in 3-dimensional fluids,} EMS Monographs in Mathematics, European Mathematical Society (EMS), Z\"urich, 2007. 
 
 
  \bibitem{Chr-BH} D. Christodoulou, \textit{The formation of Black Holes  in General Relativity}, EMS Monographs in Mathematics, 2009.
             
\bibitem{Ch-Kl0} D. Christodoulou and S. Klainerman,  {\it Asymptotic properties of linear field theories in Minkowski space}, Comm. Pure Appl. Math.  \textbf{43} (1990), 137--199.  
           
\bibitem{Ch-Kl} D. Christodoulou and S. Klainerman, \textit{The global nonlinear stability of the Minkowski space}, Princeton University Press, 1993.


\bibitem{Daf}  M. Dafermos, \textit{The mathematical analysis of black holes in general relativity}, Proceedings of the International Congress of Mathematicians, Seoul 2014, Vol. III, 747--772, Kyung Moon Sa, Seoul, 2014.
              
\bibitem{D-H-R} M. Dafermos, G. Holzegel and I. Rodnianski, \textit{Linear stability of the Schwarzschild solution to gravitational perturbations}, Acta Math. \textbf{222} (2019), 1--214.


\bibitem{D-H-R-Kerr} M. Dafermos, G. Holzegel and I. Rodnianski, \textit{Boundedness and decay  for the Teukolsky equation on Kerr spacetimes I: The case $|a|\ll M$},  Ann. PDE (2019), 118 pp.

\bibitem{Daf-Luk} M. Dafermos and J. Luk, \textit{The interior of dynamical vacuum black holes I: The C0-stability of the Kerr Cauchy horizon}, arXiv:1710.01722.
              
              
\bibitem{DaRo1} M. Dafermos and I. Rodnianski, \textit{The red-shift effect and radiation decay on black hole spacetimes}, Comm. Pure Appl. Math. \textbf{62} (2009), 859--919.

\bibitem{Da-Ro3}  M. Dafermos and I. Rodnianski, \textit{A new physical-space approach to  decay for the wave equation with applications to black hole spacetimes}, XVIth International Congress on Mathematical Physics, World Sci. Publ., Hackensack, NJ, 2010, 421--432.  


\bibitem{DaRo2} M. Dafermos and I. Rodnianski, \textit{A proof of the uniform boundedness of solutions to the wave equation on slowly rotating Kerr backgrounds}, Invent. Math. \textbf{185} (2011), 467--559.


\bibitem{mDiRySR2014} M.  Dafermos, I. Rodnianski and Y. Shlapentokh-Rothman, \textit{Decay for solutions of the wave equation on Kerr exterior spacetimes   $iii$: The  full subextremal case $|a|<m$}, Ann. of Math. \textbf{183} (2016), 787--913.

\bibitem{DHRT} M.  Dafermos, G. Holzegel,   I. Rodnianski and M. Taylor,  \textit{The non-linear stability of the Schwarzschild family of black holes}, arXiv:2104.08222.

\bibitem{GKS1} E. Giorgi, S.  Klainerman and J.   Szeftel, \textit{A general formalism for the stability of Kerr}, arXiv:2002.02740.


\bibitem{GKS2}   E. Giorgi, S.  Klainerman and J.   Szeftel,  \textit{Morawetz estimates in perturbations of  Kerr}, in preparation. 
 
 
 \bibitem{HHV} D. H\"afner, P. Hintz and A. Vasy,   \textit{Linear stability of slowly rotating Kerr black holes},   arXiv:1906.00860. 
  
\bibitem{HVas} P. Hintz and A. Vasy, \textit{The global non-linear stability of the Kerr-de Sitter family of black holes}, Acta Math. \textbf{220} (2018), 1--206.
  
 \bibitem{Kerr} R.  P. Kerr, \textit{Gravitational field of a spinning mass as an example of algebraically special metrics},  Physical Review Letters \textbf{11} (1963), 237.

\bibitem{Kl-vectorfield}  S. Klainerman, \textit{Uniform decay estimates and the Lorentz invariance of the classical wave equations}, Comm. Pure Appl. Math. \textbf{38} (1985), 321--332.

 \bibitem{Kl-ICM} S. Klainerman, \textit{Long time behavior of solutions to nonlinear wave equations}, Proceedings of the International Congress of Mathematicians, Vol. 1, 2 (Warsaw, 1983), 1209--1215, PWN, Warsaw, 1984.

\bibitem{Kl-null}  S. Klainerman, \textit{The Null Condition and global existence to nonlinear wave equations}, Lect. in Appl. Math. \textbf{23} (1986),  293--326.

\bibitem{Kl-vect2}   S. Klainerman, \textit{Remarks on the global Sobolev inequalities}, Comm. Pure Appl. Math. \textbf{40} (1987), 111--117.
 
\bibitem{KlNi} S. Klainerman and F. Nicolo, The evolution problem in general relativity. Progress in Mathematical Physics \textbf{25}, Birkhauser, Boston, 2003, +385 pp.

\bibitem{KlNi2} S. Klainerman and F. Nicolo, Peeling properties of asymptotic solutions to the Einstein vacuum equations, Class. Quantum Grav. \textbf{20} (2003), 3215--3257.



\bibitem{KRS}  S. Klainerman,  I. Rodnianski and J. Szeftel, \textit{The bounded $L^2$  curvature conjecture}, Inventiones \textbf{202} (2015), 91--216.

\bibitem{KS} S. Klainerman and J.  Szeftel, \textit{Global Non-Linear Stability of Schwarzschild Spacetime under Polarized Perturbations}, Annals of Math Studies, 210. Princeton University Press, Princeton NJ, 2020, xviii+856 pp. 

\bibitem{KS-GCM1}  S. Klainerman and J. Szeftel, \textit{Construction of GCM spheres in perturbations of Kerr}, arXiv:1911.00697.

\bibitem{KS-GCM2}  S. Klainerman and J. Szeftel, \textit{Effective  results in uniformization and intrinsic GCM spheres in perturbations of Kerr}, arXiv:1912.12195.

\bibitem{KS:Kerr-B}  S. Klainerman and J. Szeftel, \textit{Stability of Kerr. Boundedness of highest derivatives of curvature}, in preparation.


\bibitem{Ma} S. Ma, \textit{Uniform energy bound and Morawetz estimate for extreme components of spin fields in the exterior of a slowly rotating Kerr black hole II: linearized gravity}, Comm. Math. Phys. \textbf{377} (2020), 2489--2551.



\bibitem{MaMeTaTo} J. Marzuola, J. Metcalfe, D. Tataru and M. Tohaneanu, \textit{Strichartz estimates on Schwarzschild black hole backgrounds}, Comm. Math. Phys. \textbf{293} (2010), 37--83.



\bibitem{NP} E.  Newman and R. Penrose, \textit{An approach to gravitational radiation by a method of spin coefficients},  J. Math. Phys. \textbf{3} (1962), 566--578.



\bibitem{P-T} W. Press and S. A. Teukolsky, \textit{Perturbations of a rotating black hole. II. Dynamical stability of the Kerr metric}, Astrophys. J. \textbf{185} (1973), 649--673.


\bibitem{Re-W} T. Regge and J. A. Wheeler, \textit{Stability of a Schwarzschild singularity}, Phys. Rev. (2), 108:1063--1069, 1957.



 \bibitem{Rizzi} A.  Rizzi,    \textit{Angular Momentum in General Relativity: A new Definition},          Phys. Rev. Letters, vol 81, no6, 1150-1153, 1998.



\bibitem{Vo} V. Schlue, \textit{Decay of the Weyl curvature in expanding black hole cosmologies}, arXiv:1610.04172.


\bibitem{Shen}  D. Shen,  \textit{Construction of GCM hypersurfaces in perturbations of Kerr}, in preparation.


\bibitem{Yacov} Y.  Shlapentokh-Rothman, \textit{Quantitative Mode Stability for the Wave Equation on the Kerr Spacetime},  Ann. Henri Poincar\'e. \textbf{16} (2015), 289--345.   



\bibitem{Y-R} Y.  Shlapentokh-Rothman and R. Teixeira da Costa, \textit{Boundedness and decay for the Teukolsky equation on Kerr in the full subextremal range $|a|< M$: frequency space analysis}, arXiv:2007.07211.   






\bibitem{Sz} L.  B. Szabados, \textit{Quasi-Local Energy-Momentum and Angular Momentum in General Relativity}, Living Rev. Relativity, 12, (2009), 4.



\bibitem{Ta-Toh} D. Tataru and M. Tohaneanu, \textit{A Local Energy Estimate on Kerr Black Hole Backgrounds},  Int. Math. Res. Not.,  2 (2011), 248--292.


\bibitem{Teuk} S. A. Teukolsky, \textit{Perturbations of a rotating black hole. I. Fundamental equations for gravitational, electromagnetic, and neutrino-field perturbations}, Astrophys. J. \textbf{185} (1973), 635--648.


\bibitem{Whit}  B. Whiting, \textit{Mode stability of the Kerr black hole},  J. Math. Phys. \textbf{30} (1989), 1301--1305.


\end{thebibliography}
\end{document}